\documentclass[10pt, any]{book}
\pdfoutput=1 

\usepackage[paperwidth=155mm,paperheight=235mm,textwidth=117mm,textheight=199mm,includeheadfoot,inner=19mm,outer=19mm]{geometry}

\usepackage{emptypage}

\usepackage{amscd,amsmath, amssymb, fancyhdr, color, enumitem,  relsize, mwe, epsfig,dutchcal,textcomp, graphicx}
\usepackage[all]{xy}
\usepackage{pdfpages, courier}

\input cyracc.def
\font\tencyr=wncyi10
\def\cyr{\tencyr\cyracc}

\usepackage{imakeidx}
\usepackage{makeidx}
\makeindex[name=persons,title=Name Index, intoc]
\makeindex[name=terms,title=Subject Index, intoc]

\usepackage{epigraph}
\makeatletter
\let\oldendpart\@endpart 
\newcommand\partepigraph[3][60pt]{
\renewcommand{\@endpart}{%
\vskip#1
\epigraph{#2}{#3}
\vfil}
}
\let\oldendpart\@endpart 
\newcommand\twopartepigraph[5][60pt]{
\renewcommand{\@endpart}{%
\vskip#1
\epigraph{#2}{#3}\epigraph{#4}{#5}
\vfil}
}
\newcommand\removeepigraph{%
  \let\@endpart\oldendpart} 
\makeatother

\usepackage{tocloft}
\usepackage[T1]{fontenc}
\usepackage{libertine}
\usepackage{tocbibind}

\usepackage[backref=page]{hyperref}
\hypersetup{
     colorlinks   = true,
     citecolor    = magenta,
}
\renewcommand*{\backref}[1]{}
\renewcommand*{\backrefalt}[4]{%
    \ifcase #1 (Not cited.)%
    \or        (Cited on page~#2.)%
    \else      (Cited on pages~#2.)%
    \fi}


\makeatletter

\@addtoreset{equation}{chapter}
\@addtoreset{Mycounter}{chapter}
\makeatother

\makeatletter
\def\x@arrow{\DOTSB\Relbar}
\def\xlongrightarrowfill@{\arrowfill@\relbar\relbar\longrightarrow}
\newcommand{\xlongrightarrow}[2][]{%
        \ext@arrow 0099\xlongrightarrowfill@{#1}{#2}}
\def\xlongleftarrowfill@{\arrowfill@\longleftarrow\relbar\relbar}
\newcommand{\xlongleftarrow}[2][]{%
	\ext@arrow 0099\xlongleftarrowfill@{#1}{#2}}
\newcommand*{\doublerightarrow}[2]{\mathrel{
  \settowidth{\@tempdima}{$\scriptstyle#1$}
  \settowidth{\@tempdimb}{$\scriptstyle#2$}
  \ifdim\@tempdimb>\@tempdima \@tempdima=\@tempdimb\fi
  \mathop{\vcenter{
    \offinterlineskip\ialign{\hbox to\dimexpr\@tempdima+1em{##}\cr
    \rightarrowfill\cr\noalign{\kern.5ex}
    \rightarrowfill\cr}}}\limits^{\!#1}_{\!#2}}}
\makeatother

\usepackage{diagrams}
\newarrow{Dots}{}...{>}
\newarrow{DDots}{<}---{>}

\newcommand{\version}{version 6.0,\ \ December 5, 2024}

\newcommand{\la}{\lambda}
\newcommand{\al}{\alpha}
\newcommand{\be}{\beta}
\newcommand{\ga}{\gamma}

\newcommand{\ka}{K\"ahler}
\newcommand{\f}{\varphi}
\newcommand{\e}{\varepsilon}

\newcommand{\ts}{\theta^\sharp}

\newcommand{\myhash}{\raisebox{\depth}{\Large\#}}

\newcommand{\Ll}{\operatorname{Lie}}
\newcommand{\Lie}{\operatorname{Lie}}

\newcommand{\Av}{\operatorname{Av}}

\newcommand{\Ae}{{\operatorname{Ae}}}
\newcommand{\w}{\frac{d}{dw}}

\newcommand{\fin}{{\operatorname{\sf fin}}}
\newcommand{\diam}{{\operatorname{\sf diam}}}
\newcommand{\sing}{{\operatorname{\sf sing}}}
\newcommand{\dev}{{\operatorname{\sf dev}}}
\newcommand{\alg}{{\operatorname{\sf alg}}}
\newcommand{\kod}{\operatorname{\sf kod}}
\newcommand{\tr}{\operatorname{\sf trace}}

\newcommand{\CC}{\mathbb{C}}
\newcommand{\HH}{\mathbb{H}}
\newcommand{\RR}{\mathbb{R}}
\newcommand{\ZZ}{\mathbb{Z}}
\newcommand{\NN}{\mathbb{N}}

\newcommand{\GL}{{\sf{GL}}}
\newcommand{\PGL}{{\sf{PGL}}}
\newcommand{\PSL}{{\sf{PSL}}}
\newcommand{\SU}{{\sf{SU}}}
\newcommand{\CU}{{\sf{CU}}}
\newcommand{\PU}{{\sf{PU}}}
\newcommand{\SO}{{\sf{SO}}}
\newcommand{\Sol}{{\sf{Sol}}}
\newcommand{\sol}{{\goth{sol}}}
\newcommand{\U}{{\sf{U}}}
\newcommand{\OO}{{\sf{O}}}
\renewcommand{\H}{{\sf{H}}}
\newcommand{\SL}{{\sf{SL}}}
\newcommand{\Sp}{{\sf{Sp}}}
\newcommand{\sconf}{{\sf{SConf}}}

\newcommand{\spp}{\operatorname{\sf span}}

\def\eqref#1{(\ref{#1})}

\newcommand{\goth}{\mathfrak}
\newcommand{\g}{{\mathfrak g}}
\newcommand{\ra}{{\:\longrightarrow\:}}
\newcommand{\N}{{\mathbb N}}
\newcommand{\Z}{{\mathbb Z}}
\newcommand{\C}{{\mathbb C}}
\newcommand{\R}{{\mathbb R}}
\newcommand{\Q}{{\mathbb Q}}
\renewcommand{\H}{{\mathbb H}}
\newcommand{\6}{\partial}
\def\1{\sqrt{-1}\,}

\newcommand{\gsl}{\mathfrak{sl}}
\newcommand{\restrict}[1]{{\left|_{{\phantom{|}\!\!}_{#1}}\right.}}

\newcommand{\cntrct}                
{\hspace{2pt}\raisebox{1pt}{\text{$\lrcorner$}}\hspace{2pt}}

\newcommand{\arrow}{{\:\longrightarrow\:}}

\newcommand{\calo}{{\cal O}}
\newcommand{\cac}{{\cal C}}
\newcommand{\caf}{{\cal F}}
\newcommand{\cad}{{\cal D}}
\newcommand{\cax}{{\cal X}}
\newcommand{\caw}{{\cal W}}

\newcommand{\cav}{{\cal V}}

\renewcommand{\bar}{\overline}
\renewcommand{\phi}{\varphi}
\renewcommand{\epsilon}{\varepsilon}
\renewcommand{\geq}{\geqslant}
\renewcommand{\leq}{\leqslant}

\newcommand{\ev}{{\operatorname{\sf even}}}
\newcommand{\odd}{{\operatorname{\sf odd}}}
\newcommand{\fl}{{\rm fl}}
\newcommand{\Char}{\operatorname{\sf char}}
\newcommand{\ri}{\sqrt{-1}}
\newcommand{\im}{\operatorname{\sf im}}
\renewcommand{\exp}{\operatorname{\sf exp}}
\renewcommand{\log}{\operatorname{\sf log}}
\renewcommand{\mod}{\operatorname{\sf mod}}
\newcommand{\End}{\operatorname{End}}
\newcommand{\Symp}{\operatorname{Symp}}
\newcommand{\Reff}{\operatorname{Ref}}
\newcommand{\Mat}{\operatorname{Mat}}

\newcommand{\Tot}{\operatorname{Tot}}
\newcommand{\Id}{\operatorname{Id}}
\newcommand{\id}{\operatorname{{\sf id}}}
\newcommand{\ind}{\operatorname{{\sf ind}}}
\renewcommand{\sup}{{\operatorname{{\sf sup}}}}
\newcommand{\Aff}{{\operatorname{{\sf Aff}}}}
\renewcommand{\inf}{{\operatorname{{\sf inf}}}}
\newcommand{\const}{\operatorname{{\sf const}}}
\newcommand{\Vol}{\operatorname{Vol}}
\renewcommand{\Vec}{\operatorname{Vec}}
\newcommand{\Hom}{\operatorname{Hom}}
\newcommand{\Hot}{\operatorname{Hot}}
\newcommand{\intHom}{\operatorname{{\cal H\!o\!\,m}}}
\newcommand{\Ext}{\operatorname{Ext}}
\newcommand{\Tor}{\operatorname{Tor}}

\newcommand{\Aut}{\operatorname{Aut}}
\newcommand{\Conf}{\operatorname{{\sf Conf}}}
\newcommand{\Map}{\operatorname{Map}}
\newcommand{\Inv}{\operatorname{Inv}}

\newcommand{\Diff}{\operatorname{{\sf Diff}}}

\newcommand{\Iso}{\operatorname{Iso}}

\newcommand{\Sym}{\operatorname{\sf Sym}}
\newcommand{\codim}{\operatorname{\sf codim}}
\renewcommand{\dim}{\operatorname{\sf dim}}
\newcommand{\supp}{\operatorname{\sf supp}}
\newcommand{\depth}{\operatorname{\sf depth}}

\newcommand{\coker}{\operatorname{\sf coker}}

\newcommand{\rk}{\operatorname{rk}}

\newcommand{\Hol}{\operatorname{Hol}}

\newcommand{\Tw}{\operatorname{Tw}}
\newcommand{\Tr}{\operatorname{Tr}}
\newcommand{\Spec}{\operatorname{{\sf Spec}}}
\newcommand{\Alt}{\operatorname{Alt}}

\newcommand{\Rep}{\operatorname{Rep}}
\newcommand{\Pos}{\operatorname{Pos}}
\newcommand{\Der}{\operatorname{Der}}
\newcommand{\Comp}{\operatorname{\sf Comp}}
\newcommand{\Teich}{\operatorname{\sf Teich}}
\newcommand{\Mass}{\operatorname{\sf Mass}}
\newcommand{\Ric}{\operatorname{Ric}}
\newcommand{\Pic}{\operatorname{Pic}}
\newcommand{\Scal}{\operatorname{\sf Scal}}
\newcommand{\Alb}{\operatorname{Alb}}
\newcommand{\Sh}{\operatorname{Sh}}
\newcommand{\Ob}{\operatorname{{\cal O b}}}
\newcommand{\Mor}{\operatorname{{\cal M or}}}

\renewcommand{\div}{\operatorname{\sf{div}}}
\newcommand{\Ad}{\operatorname{Ad}}
\newcommand{\ad}{\operatorname{ad}}
\newcommand{\St}{\operatorname{St}}
\newcommand{\Spin}{\operatorname{Spin}}
\newcommand{\Gr}{\operatorname{Gr}}

\newcommand{\aut}{\mathfrak{aut}}
\newcommand{\ham}{\mathfrak{ham}}

\newcommand{\orb}{{\operatorname{\sf orb}}}
\newcommand{\kah}{{\operatorname{\text{\sf kah}}}}
\newcommand{\sas}{{\operatorname{\text{\sf sas}}}}
\newcommand{\bas}{{\operatorname{\text{\sf bas}}}}
\newcommand{\hor}{{\operatorname{\text{\sf hor}}}}
\newcommand{\ver}{{\operatorname{\text{\sf ver}}}}
\newcommand{\mix}{{\operatorname{\text{\sf mix}}}}
\renewcommand{\Re}{\operatorname{Re}}

\renewcommand{\Im}{\operatorname{Im}}
\newcommand{\Arg}{\operatorname{Arg}}
\newcommand{\grad}{\operatorname{\sf grad}}

\newcommand{\dash}{\operatorname{{\textnormal{--}}}}

\newcommand{\vol}{\operatorname{vol}}

\newcounter{Mycounter}[chapter]
\newcounter{lemma}[chapter]
\renewcommand{\thelemma}{{Lemma \thechapter.\arabic{lemma}}}
\newcommand{\lemma}{%
     \setcounter{lemma}{\value{Mycounter}}
     \refstepcounter{lemma}
     \stepcounter{Mycounter}
\noindent
     {\bf \thelemma:\ }
}

\newcounter{claim}[chapter]
\renewcommand{\theclaim}{{Claim \thechapter.\arabic{claim}}}
\newcommand{\claim}{%
     \setcounter{claim}{\value{Mycounter}}
     \refstepcounter{claim}
     \stepcounter{Mycounter}
\noindent
     {\bf \theclaim:\ }
}

\newcounter{sublemma}[chapter]
\newcounter{corollary}[chapter]
\renewcommand{\thecorollary}{{Corollary \thechapter.\arabic{corollary}}}
\newcommand{\corollary}{%
     \setcounter{corollary}{\value{Mycounter}}
     \refstepcounter{corollary}
     \stepcounter{Mycounter}
\noindent
     {\bf \thecorollary:\ }
}

\newcounter{theorem}[chapter]
\renewcommand{\thetheorem}{{Theorem \thechapter.\arabic{theorem}}}
\newcommand{\theorem}{%
     \setcounter{theorem}{\value{Mycounter}}
     \refstepcounter{theorem}
     \stepcounter{Mycounter}
\noindent
     {\bf \thetheorem:\ }}

\newcounter{conjecture}[chapter]
\renewcommand{\theconjecture}{{Conjecture \thechapter.\arabic{conjecture}}}
\newcommand{\conjecture}{%
     \setcounter{conjecture}{\value{Mycounter}}
     \refstepcounter{conjecture}
     \stepcounter{Mycounter}
\noindent
     {\bf \theconjecture:\ }
}

\newcounter{proposition}[chapter]
\renewcommand{\theproposition} {{Proposition \thechapter.\arabic{proposition}}}
\newcommand{\proposition}{%
     \setcounter{proposition}{\value{Mycounter}}
     \refstepcounter{proposition}
     \stepcounter{Mycounter}
\noindent
     {\bf \theproposition:\ }
}

\newcounter{definition}[chapter]
\renewcommand{\thedefinition} {{Definition~\thechapter.\arabic{definition}}}
\newcommand{\definition}{%
     \setcounter{definition}{\value{Mycounter}}
     \refstepcounter{definition}
     \stepcounter{Mycounter}
\noindent
     {\bf \thedefinition:\ }}

\newcounter{example}[chapter]
\renewcommand{\theexample}{{Example \thechapter.\arabic{example}}}
\newcommand{\example}{%
     \setcounter{example}{\value{Mycounter}}
     \refstepcounter{example}
     \stepcounter{Mycounter}
\noindent
     {\bf \theexample:\ }}

\newcounter{remark}[chapter]
\renewcommand{\theremark}{{Remark \thechapter.\arabic{remark}}}
\newcommand{\remark}{%
     \setcounter{remark}{\value{Mycounter}}
     \refstepcounter{remark}
     \stepcounter{Mycounter}
\noindent
     {\bf \theremark:\ }}

\newcounter{problem}[chapter]
\newcounter{question}[chapter]
\renewcommand{\thequestion}{{Question \thechapter.\arabic{question}}}
\newcommand{\question}{%
     \setcounter{question}{\value{Mycounter}}
     \refstepcounter{question}
     \stepcounter{Mycounter}
\noindent
     {\noindent \bf \thequestion:\ }}

\def\blacksquare{\hbox{\vrule width 5pt height 5pt depth 0pt}}
\def\endproof{\blacksquare}
\def\shortdash{\mbox{\vrule width 4.5pt height 0.55ex depth -0.5ex}}

\newcommand{\proof}{{\bf Proof: \ }}
\newcommand{\pstep}{{\bf Proof. Step 1: \ }}

\begin{document}

\title{\huge \textbf{\textsc{ Principles of Locally Conformally K\"ahler  Geometry}}}

\author{\Large\bf Liviu Ornea and Misha Verbitsky}
\date{December 5, 2024}

\maketitle
%
%
%
%
%
%
%
%
%
%
%
%
%
%
%
\hypersetup{linkcolor=blue}

\vspace*{-2cm}

\tableofcontents

\chapter*{Introduction}
\addcontentsline{toc}{chapter}{Introduction}
\markboth{Introduction}{Introduction}

\epigraph{\it Writing long books is a laborious and impoverishing act of foolishness: expanding in five hundred pages an idea that could be perfectly explained in a few minutes. A better procedure is to pretend that those books already exist and to offer a summary, a commentary.}{\sc \scriptsize Jorge Luis Borges}

It was at the beginning of 1970s when Lieven
\index[persons]{Vanhecke, L.} Vanhecke gave a talk in Ia\c si (Romania) about conformal
K\"ahler manifolds. Izu \index[persons]{Vaisman, I.} Vaisman, a member of the Faculty of 
Mathematics  of the University of Ia\c si at
that time, asked what can be said about the local
situation. Vanhecke replied: ``That's for you to find
out''. And it is exactly what Vaisman did. In 1976, he
published the paper \cite{va_isr} in which he introduced
the  ``locally conformal (almost) K\"ahler manifolds''. It
was the birth of LCK geometry. In a long series of
papers, Vaisman clarified the notion by comparing it with
the K\"ahler manifolds, gave the first examples,  and
introduced the class of ``LCK manifolds with parallel Lee
form'' (which now bears his name). Other examples were\index[terms]{form!Lee!parallel}
given at the beginning of the 1980s by Franco \index[persons]{Tricerri, F.} Tricerri,
\cite{tric}. Starting in 1980, the notion was also
considered by several Japanese mathematicians, important
results being obtained by Toyoko \index[persons]{Kashiwada, T.} Kashiwada, Yoshinobu
\index[persons]{Kamishima, Y.} Kamishima, Kazumi \index[persons]{Tsukada, K.} Tsukada, etc.

In the first 20 years, this new kind of Hermitian 
geometry was studied mainly using differential geometry
methods. The focus was on Riemannian and conformal
properties. For example, a great deal of work was
dedicated to isometric submanifolds (totally geodesic,
totally umbilical, minimal, totally real, CR etc.). The
exception was the paper \cite{tric} in which the blow-up at
points was proved to preserve the LCK class.

Most of these findings were gathered in the monograph
\cite{do} written by Sorin Dragomir\index[persons]{Dragomir, S.} and the first author. The present book starts from where the previous
one ends.

In our book we combine the methods of algebraic
geometry, functional analysis and complex analysis
to study the LCK manifolds, which properly belong
to differential geometry. Similar to K\"ahler
geometry, the subject often transcends the boundary
of its domain, and we arrive at the point when
no differential geometry is involved. Still,
the background is always differential geometric:
we mention several works where the authors study
the $p$-adic versions of LCK manifolds 
(\cite{_Mustafin_, _Voskuil_,_Scholze:congress_}), but we never
pursue this direction.

We tried to be accessible to  beginner students;
this is one of the reasons why we cover several preliminary
topics, such as foliations, Frobenius theorem, and Ehresmann
connections.

We start with the definition of complex and K\"ahler manifolds, mainly
to fix the notation. For an introduction to complex
geometry, see \cite{demailly, griha,
	_moroianu:book_, vois, _Huybrechts_}. Note that the notation,
especially the signs, vary from one author to another.
Within this book, we try to be consistent, sometimes
without success. 

The preliminary requirements for this book vary from part
to part. The first part, expanded from several lecture
courses, is oriented toward advanced undergraduate and master
students. We assume a working knowledge of differential geometry
(Riemannian structures, connections, principal bundles, 
de Rham algebra), topology (de Rham cohomology, Poincar\'e duality, 
fundamental groups, local systems), Lie groups and
algebras, basic algebraic geometry, basic 
complex analysis and basic functional analysis. 
We use Hodge theory, citing, without proof, 
several key results, such as the Hodge 
decomposition. We also use, without proofs,
elements of the theory of Stein manifolds;
\index[persons]{Demailly, J.-P.} Demailly in \cite{demailly} gives all the 
tools that are necessary for our use.

The second part treats more advanced subjects
and the requirements are much higher. We try to 
give a basic introduction to several key notions,
such as derived functors and the \index[persons]{Grothendieck, A.} Grothendieck
spectral sequence, but a working knowledge of
algebraic geometry (in the scope of \index[persons]{Demailly, J.-P.} Demailly's
textbook \cite{demailly}) and homological algebra 
(\index[persons]{Grothendieck, A.}Grothendieck's Toh\^oku paper \cite{_Grothendieck:Tohoku_})
is necessary. A graduate student studying algebraic
geometry won't find it too specialized, and we tried
to lower the requirements by introducing each subject
within the text.

The third part is a survey of current research on LCK
geometry. We give new proofs of most results, using
the methods developed earlier in this book. The third
part is oriented towards the researchers who operate with the
knowledge of differential geometry required for the subject.

In the first two parts of this book, every chapter
is augmented by a sequence of exercises. Originally, 
the exercises were used in the lecture courses, but we 
expanded them to include the whole new series, giving an
introduction to some subjects (such as elliptic curves
and Galois theory). Near the end of the writing, we 
put a moratorium on adding new chapters; at that moment,
all new results we got went to the exercises. This is
why some of the exercises are difficult theorems in 
themselves. There is little consistency or system
in the exercises, but we hope that every reader would find
something interesting to herself or himself.

The last chapter (``Open questions'') is a logical
conclusion of the same approach: some of the questions
we mention are pretty much impossible to solve, but most
of them could be used by an early or advanced graduate
student as a part of her or his diploma work. We 
avoided using general ideas such as ``develop
a geometric flow'' in favour of concrete 
conjectures, and described the context 
whenever possible.

Every chapter of the first two parts 
starts with an introduction, which is more or less
independent of the body of the chapter. Its purpose
is to set forth the narrative of the main body, 
explain the history, and place the subject in
a wider mathematical context; often,
it tells the same story, with more flourish and less 
mathematical rigour. Most of the definitions given 
in the introductions are repeated in the respective
chapters, in more formal fashion.

With painful trepidation, 
we avoid several important subjects, most notably, 
stable bundles, stable coherent sheaves, and Yang--Mills
connections: the book is too big as it is. For an
introduction to Yang--Mills geometry, see 
\cite{_Lubke_Teleman_}. Also, most of the concepts 
relating hyperk\"ahler and hypercomplex geometry
to LCK geometry is relegated to ``open questions''.
Finally, we did little justice to Sasakian geometry,
which we view as an integral part of our subject. Almost
everything related to things Sasakian is beautifully
explained in the book \cite{bog},
by Charles P. Boyer\index[persons]{Boyer, C. P.} and Krzysztof \index[persons]{Galicki, K.} Galicki,
and we did not want to repeat their work.

\section*{Thanks}

In many of the tight places, we were helped by
F. \index[persons]{Ambro, F.} Ambro, E. \index[persons]{Amerik, E.} Amerik, D. \index[persons]{Angella, D.} Angella,
M. \index[persons]{Aprodu, M.} Aprodu, R. Bryant, M. Entov, P. \index[persons]{Gauduchon, P.} Gauduchon, N. \index[persons]{Istrati, N.} Istrati, D. \index[persons]{Kaledin, D.} Kaledin,
A. \index[persons]{Moroianu, A.} Moroianu, S. Nemirovski\index[persons]{Nemirovski, S.},  A. \index[persons]{Otiman, A.} Otiman, J. V. \index[persons]{Pereira, J. V.} Pereira, Yu. Prokhorov,\index[persons]{Prokhorov, Yu.}
K. Shramov\index[persons]{Shramov, K.} and M. \index[persons]{Toma, M.} Toma. 

Victor \index[persons]{Vuletescu, V.} Vuletescu was almost a bona fide
coauthor of this book, with his constant attention,
encouragement and many beautiful ideas.

M. Verbitsky is grateful to all his students, who
helped and inspired him during the lecture courses on
LCK geometry and related subjects:
A. Abasheva,\index[persons]{Abasheva, A.} R.~Deev,\index[persons]{Deev, R.} V. Gizatulin,\index[persons]{Gizatulin, V.} Iu. Gorginian,\index[persons]{Gorginian, Iu.} N. \index[persons]{Klemyatin, N.} Klemyatin,
D. Korshunov,\index[persons]{Korshunov, D.} N. \index[persons]{Kurnosov, N.} Kurnosov, S.~Makarova,\index[persons]{Makarova, S.} P. Osipov,\index[persons]{Osipov, P.} G. Papayanov,
V. \index[persons]{Rogov, V.} Rogov, A. Soldatenkov,\index[persons]{Soldatenkov, A.} L. Soukhanov,\index[persons]{Soukhanov, L.} A.~Viktorova\index[persons]{Viktorova, A.}. Their
ideas and questions enriched the narrative and
were often incorporated here.

Many thanks to the colleagues and students who read the
first draft of the book and corrected some of our errors:
C. P. Boyer, K. Broder, C. Ciulic\u a\index[persons]{Ciulic\u a, C.}, N. Istrati, T. Orban\index[persons]{Orban, T.}, A. Otiman, M. Stanciu\index[persons]{Stanciu, M.}, and V. Vuletescu.

\hfill

\hfill {Bucharest -- Rio de Janeiro, \ \ August 2022.}

\hfill

We thank the  referees and editors of Birkh\"auser for improving the content and the layout of the manuscript.
After the book was published, we received a great help from Ethan Addison, who sent us a big list of errors and corrections.


\markboth{intro}{Introduction}

{\setlength\epigraphwidth{0.5\linewidth}
\partepigraph{\em The Dao gives birth to unity,\\
Unity gives birth to duality,\\
Duality gives birth to trinity,\\
Trinity gives birth to a myriad of things.}{\sc\scriptsize Tao Te Ching, by  Laozi}
\part{Lectures in locally conformally K\"ahler geometry}
\removeepigraph
}

\chapter{K\"ahler manifolds}

{\setlength\epigraphwidth{0.9\linewidth}
\epigraph{
\it \qquad
The shift penetrates verses throughout (especially contemporary
verses); it is one of the most important parts of the verse. It modifies
the word, the stanza, the sounds. \\

\smallskip
\qquad One must establish a special "shift police force" for the
prompt capture of shifts, which would leave their authors
gaping in amazement. \\
\smallskip

\qquad The shift conveys movement and space.\\
\qquad The shift conveys multiplicity of meanings and images.\\
\qquad The shift is the style of our contemporary life.\\
\qquad The shift is a new discovery of America!..
}
{\sc\scriptsize Shiftology of Russian
Verse: An Offensive and
Educational Treatise, by 
A. Kruchenykh}}

\section{Complex manifolds}

\definition  Let $M$ be a smooth manifold. 
	An {\bf almost complex structure}\index[terms]{!almost complex structure} is a section
	$I\in \End(TM)$ that satisfies $I^2 = - \Id_{TM}$.
	
	The couple $(M,I)$ is called an {\bf almost complex manifold}.\index[terms]{manifold!almost complex}

\hfill

Slightly abusing the language, we denote 
the extension of $I$ to  $TM_\C$ by $I$ as well.
The eigenvalues of this operator are $\pm \1$.
Let $TM_\C=T^{0,1}M\oplus T^{1,0}M$ be the corresponding eigenvalue 
decomposition.

\hfill

\remark
In algebraic geometry, one should
always make clear that the vector bundle
is distinct from its space of sections.
One usually uses the notation $\Gamma(B)$
or $H^0(M, B)$ for the space of sections
of a vector bundle $B$. However, in 
differential geometry one could avoid
this distinction (see Section \ref{_Bundles_Intro_Section_}).
Throughout this book, we often
use the same letter for the bundle
over a smooth manifold and its space 
of sections.

\hfill

\definition 
	An almost complex structure is {\bf  integrable\index[terms]{almost complex structure!integrable}}
	if $[X,Y]\in T^{1,0}M$ for all $ X,Y \in T^{1,0}M$.
	In this case, $I$ is called {\bf  a complex structure\index[terms]{complex structure}}.
	
	A manifold with an integrable almost complex structure
	is called {\bf a complex manifold\index[terms]{manifold!complex}}.

\hfill

\theorem  {\rm\bf (Newlander--Nirenberg)\\ \index[persons]{Newlander, A.}\index[persons]{Nirenberg, L.}}\label{nn}
	An almost complex structure is integrable if and only if the manifold admits an atlas with charts taking values in $\C^n$, $n=\dim_\C M$, and with {\bf  holomorphic} changes of coordinates.\index[terms]{theorem!Newlander--Nirenberg}

\hfill

\remark  The  ${C}^{\infty} (M)$-linear map 
	$N:\;\Lambda^2(T^{1,0}M)\arrow T^{0,1}M$ defined by the commutator is  called {\bf                            
		the Nijenhuis tensor\index[terms]{tensor!Nijenhuis}} of $I$. One can represent $N$ as a section 
	of $\Lambda^{2,0}M \otimes T^{0,1}M$.

\section{Holomorphic vector fields}

\definition\label{holo_field} A real vector field 
$Z\in TM$ is called {\bf real holomorphic} if $Z(f)$ is
holomorphic for any holomorphic function $f$ defined
on some open subset of $M$.\index[terms]{vector field!real holomorphic}

\hfill

\remark Let  $X\in TM$ be a real vector field. Then
$$X=\frac 12 (X-\1 IX) + \frac 12 (X+\1 IX).$$
Clearly, $X-\1 IX\in T^{1,0}M$ (and is called the (1,0)
part of $X$), whereas $X+\1 IX\in T^{0,1}M$ (and is called
the (0,1) part of $X$). 

\hfill

\definition A (1,0)-vector field $X\in T^{1,0}M$ is called
            {\bf holomorphic} if its real part is
            real holomorphic in the sense of \ref{holo_field}.\index[terms]{vector field!real holomorphic}

\hfill

\remark \label{_holo_1,0_real_Remark_}
A (1,0)-vector field $X\in T^{1,0}M$ is
holomorphic if and only if $X(f)$ is
holomorphic for any local holomorphic function $f$ defined
on some open subset of $M$. Moreover, taking the (1,0)-part
gives a bijective correspondence between real holomorphic
vector fields and (1,0)-holomorphic vector fields.\index[terms]{vector field!holomorphic}

\hfill

The next proposition is clear.

\hfill

\proposition\label{realholofield} Let $X$ be a real vector field. The following are equivalent:
\begin{itemize}
	\item[(i)] $X$ is real holomorphic.\index[terms]{vector field!real holomorphic}
	\item[(ii)] $\Lie_XI=0$, where $\Lie_X$ denotes the Lie derivative along $X$. 
	\item[(iii)] The flow generated by $X$ consists of biholomorphic transformations.
\end{itemize}

\remark Let $(M,I)$ be a compact complex manifold. Then the group of biholomorphisms of $M$, $\Aut(M,I)$, is a Lie group whose Lie algebra is the space of real holomorphic vector fields. See also Exercise \ref{holofield}.

\hfill

\remark\label{_holo_field_commute_Remark_}
For any holomorphic vector field $X$, the fields
$X$ and $X^c:=I(X)$ commute.\index[terms]{vector field!holomorphic}
Indeed, 
\[ [X, X^c]=\Lie_X(X^c)= \Lie_X(IX)= \Lie_X(I) (X) + I(\Lie_X(X))= 0.
\]

\section{Hermitian manifolds}

\definition  A Riemannian metric $g$ on
	an almost complex manifold $(M,I)$ is called 
	{\bf Hermitian} if $g(IX, IY)= g(X,Y)$.
	In this case 
\[g(IX, Y)= g(I^2X, IY) = - g(X, IY)=-g(IY, X)
\]
	and hence $\omega(X,Y):= g(IX, Y)$ is skew-symmetric.\index[terms]{metric!Hermitian}

\hfill

\definition The differential 
	form $\omega\in \Lambda^2(M)$ is called
	{\bf  the Hermitian form} of $(M,I,g)$.\index[terms]{form!Hermitian}

\hfill

\remark $\omega$ is $\U(1)$-invariant, and hence  {\bf of Hodge type} (1,1).\index[terms]{Hodge!type}

\hfill

\definition Let $g, g'$ be two Riemannian metrics on the same manifold $M$. They are said to be {\bf conformal} if there exists $f\in C^\infty(M)$ such that $g'=e^{f}g$. The set $[g]=\{e^fg\,;\, f\in C^\infty(M)\}$ is called {\bf the conformal class} of $g$. If $(M,I)$ is a complex manifold and one metric in a conformal class is Hermitian with respect to $I$, then all metrics in the conformal class are so.

\section{K\"ahler manifolds}\index[terms]{manifold!K\"ahler}

\hfill

\definition A Hermitian manifold $(M,I,g,\omega)$ is called {\bf K\"ahler} if $d\omega=0$. $\omega$ is then called the K\"ahler form, and its cohomology class $[\omega]\in H^2(M,\R)$ is called the {\bf K\"ahler class}.\index[terms]{form!K\"ahler}\index[terms]{class!K\"ahler}

\subsection{Examples of K\"ahler manifolds}

\example\label{_standard_Kahler_structure_on_C^n_example_} $\C^n$ with the flat metric $g=\Re\left(\sum dz_i\otimes d\bar z_i\right)$, and with K\"ahler form $\omega=\1\sum dz_i\wedge d\bar z_i$.

\hfill

\example {\bf  The complex projective space.}  
Let $\C P^n$ be the complex projective
space, and $g$ a $\U(n+1)$-invariant Riemannian metric. It is called
{\bf the  Fubini--Study metric on $\C P^n$}.\index[terms]{metric!Fubini--Study} 
Note that the Fubini--Study metric is defined up to a constant factor.
In some textbooks this constant is fixed by requiring that 
the Hopf fibration $h:\; S^{2n+1} \arrow \C P^n$ takes a
tangent vector $v\in T_x S^{2n+1}$ orthogonal to a 
fibre of $h$ to a vector $dh(v)\in T_{h(x)} \C P^{n}$ of 
the same length.\footnote{This is equivalent to asking
that the sectional curvature $K(x,y)$ of the Fubini--Study metric
satisfies $1 \leq K(x,y) \leq 4$ for any pair of 
orthogonal unit tangent vectors; see \cite{besse}.}

The Fubini--Study metric can be obtained by taking an arbitrary Riemannian metric, then averaging it with $\U(n+1)$ using the Haar measure on $\U(n+1)$.\index[terms]{Haar measure}

\hfill

\remark  For any $x\in \C P^n$, the stabilizer for the $\U(n+1)$ action 
	$\St(x)$ is isomorphic to $\U(n)$. The  Fubini--Study metric on
	$T_x\C P^n= \C^n$ is $\U(n)$-invariant, and hence unique up to a constant. As  $d\omega\restrict x$ is a $\U(n)$-invariant 3-form
	on $\C^n$, it has to vanish, because $-\Id\in \U(n)$, and hence $d\omega=0$.

\hfill

\remark 
	With the same argument,  Hermitian symmetric spaces  are  K\"ahler.

\hfill 

\example The product of two \ka\ manifolds is \ka\ with respect to the product metric.

\hfill

\example The blow-up at points and along submanifolds of a K\"ahler manifold is again \ka .\index[terms]{blow-up}

\hfill

\example {\bf  Complex submanifolds.}

\hfill

\definition 
	{\bf An almost complex submanifold} $X \subset M$
	of an almost complex manifold $(M,I)$ is a smooth submanifold
	that satisfies $I(TX)= TX$.
	
	If $I$ is integrable, then $X$ is called {\bf  a complex submanifold} of $M$.\index[terms]{submanifold}

\hfill

\remark  
	Let $X \subset M$ be an almost complex submanifold of $(M,I)$,
	where $I$ is integrable. Then $(X, I\restrict{TX})$
	is a complex manifold.

\hfill

Since exterior differentiation commutes with pullback, the restriction of the \ka\ form is closed on each complex submanifold of a \ka\ manifold.

In particular, {\em every projective manifold \index[terms]{manifold!projective}
	(complex submanifold of \ $\C P^n$) is K\"ahler.}\index[terms]{submanifold}

\subsection{Menagerie of complex geometry}

Usually, in algebraic geometry\index[terms]{geometry!algebraic} one deals with
projective manifolds. There are two wider classes one 
should consider when studying projective ones.

\hfill

\example  A complex manifold that is birational
to a projective manifold is called {\bf Moishezon}.\index[terms]{manifold!Moishezon}\index[persons]{Moishezon, B. G.} 

{\em The transcendence degree\index[terms]{function!meromorphic} $a(M)$ of the field $k(M)$ of 
	global meromorphic functions on a compact complex manifold $M$
	satisfies $a(M)\leq \dim_\C M$}, as shown by Moishezon, 
\cite{moi}; equality here means 
that $M$ is Moishezon. 
The number $a(M)$ is also called {\bf the algebraic dimension} of $M$.


\hfill

\theorem   {({\bf Moishezon}, \cite{moi})} {Any K\"ahler Moishezon manifold is
		projective.}\index[terms]{theorem!Moishezon}\index[terms]{manifold!projective}

\hfill

\example {\bf  Small deformations of K\"ahler manifolds}
often result in non-projective K\"ahler ones (even for a 
torus and a K3 surface).

\hfill

The class that includes Moishezon and K\"ahler manifolds is called
{\bf    Fujiki class C}, \cite{fujiki}.\index[terms]{Fujiki class C}\index[persons]{Fujiki, A.} A manifold is {\bf of  Fujiki class C}
if it is bimeromorphic to a K\"ahler manifold. As shown in {\em loc. cit.},
\index[persons]{Fujiki, A.} Fujiki class C is closed under all natural operations that occur
in algebraic geometry\index[terms]{geometry!algebraic} (such as taking moduli spaces or images).

\hfill

\remark  The K\"ahler minimal model program\index[terms]{minimal model program} \cite{_Horing_Peternell_}, \cite{cdv}, would imply
	that any K\"ahler manifold admits a sequence of bimeromorphic
	fibrations with fibres that are either projective,\index[terms]{manifold!projective}\index[terms]{manifold!hyperk\"ahler} hyperk\"ahler
	or tori, and hence  the class of K\"ahler manifolds
	is probably very restricted.
	By contrast,  the class of complex manifolds 
	is  huge.

\hfill

The fundamental\index[terms]{fundamental group!of K\"ahler manifolds}
group of a K\"ahler manifold is very special. On the other hand:

\hfill

\theorem  {({\bf \index[persons]{Taubes, C. H.}Taubes} \cite{taubes}, {\bf \index[persons]{Panov, D.} Panov---\index[persons]{Petrunin, A.}Petrunin} \cite{pp})}
\label{taubes_fund_group}\index[persons]{Panov, D.}\index[persons]{Taubes, C. H.}\index[terms]{theorem!Taubes, Panov--Petrunin} 
For any finitely generated, finitely presented group $\Gamma$, 
	there exists a compact, complex 3-di\-men\-sio\-nal
	manifold $M$ with $\pi_1(M)=\Gamma$.

\hfill

\conjecture ({\bf \index[persons]{Yau, S.-T.} Yau} \cite{yau})
	Let $(M,I)$ be a compact almost complex manifold, $\dim_\C M\geq 3$.\index[terms]{conjecture!Yau}
	Then $I$ can be deformed to a complex structure.\index[terms]{almost complex structure}

\hfill

\remark The following important result of \index[persons]{Gromov, M.} Gromov (\cite[p. 103]{_Gromov_}) can be cited to
support this conjecture. Let $M$ be a non-compact almost complex
manifold, $\dim_\C M \geq 3$. Then $M$ admits a complex structure.

\hfill

\remark The (non-)existence of a complex structure     is 
	highly non-trivial even in the simplest cases, such as
	$S^6$. See  \cite{hkp}, \cite{cdv}. On the other hand, $S^6$ has a natural  almost complex structure constructed out of octonions,\index[terms]{octonions} see \cite{bal} for example, that admits a compatible  {\bf nearly-K\"ahler}\index[terms]{manifold!nearly-K\"ahler} metric, which means that $\nabla I$ is antisymmetric, see \cite{bfgk}, \cite{br}.
	
	\hfill
	
\remark It is known that non-K\"ahler complex manifolds are much more abundant
	than K\"ahler ones, except in complex dimension 2, where non-K\"ahler
	manifolds are few and much better understood than projective
	ones. However, it is very hard to
	come up with new examples of compact,  non-K\"ahler complex manifolds.

\section{Exercises}

\subsection{K\"ahler geometry and holomorphic vector fields}
\index[terms]{geometry!K\"ahler}
\index[terms]{vector field!holomorphic}
\begin{enumerate}[label=\textbf{\thechapter.\arabic*}.,ref=\thechapter.\arabic{enumi}]

\item\label{_abelian_compact_complex_Exercise_}
Let $G$ be a compact, complex, connected Lie group.
Prove that $G$ is abelian.

	\item\label{holofield}
Let $X$ be a holomorphic vector field on a complex manifold, that is, 
one that satisfies
\index[terms]{vector field!holomorphic} $\Lie_XI=0$.
Prove that $IX$ is also holomorphic.
	
	{\em Hint:}   $\Lie_XI=A(I)$, where $A=\nabla(X)$ acts by the formula
	$A(I)(v)=A(Iv)-IA(v)$. Therefore, {\em $X$ is holomorphic
		if and only if $\nabla(X)$ is complex linear.} Since
	$\nabla(I)=0$, one has $\nabla(IX)=I(\nabla(X))$,
	and hence $\nabla(X)$ is complex linear, implying that 
	$\nabla(IX)$  is complex linear. \index[terms]{vector field!holomorphic}

	\item Let $(M,I,g,\omega)$ be an almost complex
          Hermitian manifold, $\nabla$ a torsion-free
          connection. \index[terms]{connection!torsion-free}
Assume $\nabla\omega=0$ and
          $\nabla I=0$. Then $(M,I,g,\omega)$ is
          K\"ahler.

	{\em Hint:}  
	Use  $[X,Y]=\nabla_X Y - \nabla_Y X$, to show that $T^{1,0}$ is involutive, then use 
	$d\omega= \Alt(\nabla\omega)$, which holds for torsion-free connections.

\item 
Let $(M,I)$ be an almost complex manifold,
and $d^c:= I d I^{-1}$. Prove that $I$ is integrable
if and only if $dd^c=-d^c d$.

\item
Let $\omega$ be a non-degenerate 2-form on a Riemannian manifold,
and $\nabla$ its Levi--Civita connection. Assume that $\nabla(\omega)=0$.
Prove that $M$ admits a complex structure $I$ such that $\nabla(I)=0$.

\item
Let $(M,I)$ be an almost complex manifold, 
$\dim_\C M=n$, $U\subset M$ a dense, open subset,
and $\Omega \in \Lambda^{n,0}(U)$ a non-degenerate 
$(n,0)$-form. Assume that $d\Omega=0$.
Prove that the almost complex structure $I$ is integrable.


\item
A {\bf holomorphic differential}
on an almost complex manifold is a closed $(1,0)$-form.
Let $G$ be a finite group acting on $M=\C^n$ by
holomorphic maps. Prove that $M$ admits\index[terms]{differential!holomorphic}
a non-zero $G$-invariant holomorphic differential.

\item
Let $M=\C P^{2m}\times \C P^{2n}$, for some $m, n \in \Z^{>0}$.
Prove that $M$ does not admit a K\"ahler structure
with non-standard orientation.

{\em Hint:} Prove that all complex manifolds have a canonical
orientation. Prove that this orientation can be given by 
the top power of the K\"ahler form, if the manifold is K\"ahler.

\item Let $M:= \frac{\C^n \backslash 0}{\langle A \rangle}$
be the so-called linear Hopf manifold, with $A$ an invertible linear operator
with operator norm $\Vert A \Vert < 1$. Prove that $M$
admits no symplectic structures.

\item
Let $M$ be a complex manifold,
admitting a non-zero holomorphic vector field $\xi$ with zero set $Z$. Suppose that $Z$
is non-empty and zero-dimensional.\index[terms]{vector field!holomorphic}
Prove that the topological Euler\index[terms]{Euler--Poincar\'e characteristic}
characteristic of $M$ is non-negative.

\item 
Let $M$ be a compact K\"ahler manifold, and
 $\Alb(M):=\frac{H^0(\Omega^1(M))^*}{\Lambda}$,
where $H^0(\Omega^1(M))$ is the space of holomorphic differentials,
and $\Lambda$ the group generated by all 
integrals over integer homology classes.
The group $\Alb(M)$ is called {\bf the Albanese variety}
of $M$, denoted by $\Alb(M)$. Prove that:\index[terms]{variety!Albanese}
\begin{enumerate}
	\item  The lattice $\Lambda\subset H^0(\Omega^1(M))^*$
is discrete and cocompact, and hence $\Alb(M)$ is a compact complex torus. 
	\item Fix $m\in M$. Consider the map $\Psi:\; M \arrow \Alb(M)$,
taking $x \in M$ to the map $\theta\mapsto \int_\gamma \theta$,
where $\gamma$ is any path-connecting $m$ to $x$.
Prove that $\Psi$ is defined unambiguously and is holomorphic.
\end{enumerate}

\end{enumerate}

\subsection{The Lie algebra of holomorphic Hamiltonian Killing fields}\index[terms]{vector field!holomorphic}\index[terms]{vector field!Hamiltonian}\index[terms]{vector field!Killing}

The following exercises are not elementary;
we advise the less prepared reader to skip them.
We will sometimes refer to these statements later in this book.
The reader can find these results in \cite[Chapter 2.H]{besse}.

\begin{enumerate}[label=\textbf{\thechapter.\arabic*}.,ref=\thechapter.\arabic{enumi}]
\setcounter{enumi}{11}

\item \label{_ex_needed_in_chapter_21_}
Let $(M, \omega)$ be a compact K\"ahler manifold,
$G(M)$ the group of its holomorphic diffeomorphisms,
$\g\subset TM$ the Lie algebra of vector fields satisfying
$\nabla(X)=0$, where $\nabla$ is the Levi--Civita
connection, and $G$ its Lie group.
\begin{enumerate}
	\item Prove that $G\subset G(M)$.
	\item Prove that all its orbits are compact complex tori in $M$.
\item Prove that  $i_X \omega$ is the real
part of a holomorphic 1-form for any $X \in \g$.
	
\item Let $G_1\subset G(M)$ be the group of
          holomorphic diffeomorphisms
	acting trivially on $\Alb(M)$. 
Prove that for  $X\in \g$ there is a holomorphic
1-form $\eta$ such that $\langle \eta, X\rangle \neq 0$. Use this
to prove that $G_1 \cap G=\{\Id\}$.
\end{enumerate}


\item\label{_Killing_antisym_Exercise_}
Let $X$ be a vector field 
on a Riemannian manifold $(M,g)$. Prove that
$X$ is Killing\index[terms]{vector field!Killing} if and only 
if the operator $A:=\nabla(X)\in \End(TM)$ 
taking $Y$ to $\nabla_Y(X)$ is 
antisymmetric: $g(A(x),y)=- g(x, A(y))$
for all $x, y \in T_m M$.

\item\label{_antisym_I_lin_Exercise_} 
Consider the standard Euclidean form $g$
and the Hermitian form $\omega$ on a vector space
$V=\R^{2n}=\C^n$.
Let $A\in \End_\R(V)$ be a matrix
that satisfies $\omega(Ax, y)= - \omega(x, Ay)$ and
$g(Ax, y)=- g (x, Ay)$, and $\omega(IAx, y)= - \omega(x, IAy)$ 
for all $x, y\in V$.
Prove that $A=0$.

\item\label{_nabla_IX_Exercise_}
Let $X,Y\in TM$ be holomorphic vector fields
on a K\"ahler manifold, and $\nabla$ the Levi--Civita
connection.
\begin{enumerate}
\item Prove that $\nabla_{IX} Y = I\nabla_X Y$.
\item Let $A(\cdot):= \nabla_X \cdot$. Assume that 
$\nabla_X \omega= \nabla_{IX}\omega=0$.
Prove that $\omega(IA\cdot, \cdot)= - \omega(\cdot, IA\cdot)$ 
is equivalent to $\Lie_{IX}\omega=0$.
\end{enumerate}

{\em Hint:}
Prove that $\nabla_{IX} -\Lie_{IX}= IA$,
and show that $\nabla_{IX}\omega -\Lie_{IX}\omega=
\omega(IA\cdot, \cdot)+ \omega(\cdot, IA\cdot)$.

\item \label{_Killing_holomorphic_is_parallel_Exercise_}
Let $(M, I, \omega)$ be a K\"ahler manifold and let 
$\g\subset TM$ be the  subalgebra formed by vector fields
$X\in TM$ such that $\Lie_X I=0$, and 
$\Lie_X\omega=\Lie_{IX}\omega=0$.
\begin{enumerate}
\item Prove that for all $X\in \g$, one has $\nabla(X)=0$.
\item Let $Y$ be a holomorphic vector field that satisfies
$\nabla(Y)=0$. Prove that $Y\in \g$.
\item Prove that $\g$ is abelian, if $M$ is compact.
\end{enumerate}

{\em Hint:} Use Exercise \ref{_Killing_antisym_Exercise_}, 
Exercise \ref{_nabla_IX_Exercise_}
and Exercise \ref{_antisym_I_lin_Exercise_} for the first
part.

\item Let $\nabla$ be a torsion-free 
connection\index[terms]{connection!torsion-free} on
a complex manifold $(M,I)$, satisfying $\nabla(I)=0$.
Prove that a vector field
$V$ is holomorphic if and only if
$\nabla(V)\in \End(TM)$ commutes with 
$I\in \End(TM)$. 

{\em Hint:} Prove that $\nabla_{IZ}(V)= I(\nabla_{V}(Z))+ [V, IZ]$
for any vector fields $Z, V$ on $M$. Use
$\nabla_{V}(Z)+[V, Z]= \nabla_Z V$.

\item
Let $v$ be a vector field
on a K\"ahler manifold $(M,I, \omega)$, satisfying
$\Lie_v \omega=0$, and $\eta:=i_v(\omega)$
the corresponding 1-form, that is  closed by Cartan formula.\index[terms]{Cartan formula}
Let $\nabla$ be the Levi--Civita connection.
Prove that $\nabla \eta\in \Sym^2(T^*M)$
for any closed 1-form $\eta$.
Prove that $v$ is holomorphic if and only if the 1-form
$\eta:=i_v(\omega)$ satisfies
$\nabla \eta\in \Sym^{1,1}(T^*M)$, that is, 
$\nabla \eta$ is of type (1,1) with respect to the 
Hodge decomposition on $\Sym^*(T^*M)$.

{\em Hint:} Use the previous exercise.

\item
Let $\eta$ be a closed real 1-form on a compact K\"ahler manifold
that satisfies $\nabla \eta\in \Sym^{1,1}(T^*M)$.
Consider the decomposition $\eta= \alpha + \eta'$ 
where $\eta'$ is exact and $\alpha$ the real
part of a holomorphic 1-form. Prove that
$\nabla \eta'\in \Sym^{1,1}(T^*M)$ 
and $\nabla\alpha=0$.

{\em Hint:} Prove that the harmonic decomposition commutes
with the connection and use this to show that 
$\nabla \alpha\in \Sym^{1,1}(T^*M)$. Show that
$\nabla \alpha\in \Sym^{2,0}(T^*M)$
when $\alpha$ is  a closed holomorphic 1-form.

\definition
Let $(M, \omega)$ be a symplectic manifold.
We denote the contraction of a differential form $\eta$ with a vector 
field $v$ by $i_v(\eta)$.
A vector field $v$ on $M$ is called {\bf Hamiltonian}
if $\Lie_v \omega=0$ and the 1-form $i_v(\omega)$
is exact.\footnote{The form $i_v(\omega)$ is closed,
because $\Lie_v \omega= d(i_v(\omega))=0$.}\index[terms]{vector field!Hamiltonian}

\item\label{_killing=hamiltonian_plus_parallel_Exercise_}
Let $v$ be a holomorphic Killing vector
field\index[terms]{vector field!Killing}\index[terms]{vector field!holomorphic} on a compact K\"ahler manifold $M$.
Prove that there exists a unique Hamiltonian
holomorphic vector field $v'$\index[terms]{vector field!Hamiltonian}\index[terms]{vector field!holomorphic}
such that $\alpha:=i_{v+v'}(\omega)$
is the real part of a holomorphic 1-form.
Prove that $v$ is Hamiltonian
if and only if $\alpha=0$.

{\em Hint:} Use the previous exercise.

\item
Let $\goth s$ be the space of closed
1-forms on a compact K\"ahler manifold that satisfy
$\nabla \eta\in \Sym^{1,1}(T^*M)$, 
$\g \subset \goth s$ the space of 
parallel 1-forms in $\goth s$,
and $\goth h$ the space of all 
exact 1-forms in $\goth s$.
\begin{enumerate}
\item Show that any exact 1-form 
on a compact manifold $M$ vanishes somewhere in $M$.
\item Prove that $\g \cap \goth h=0$.
\item 
Use Exercise \ref{_killing=hamiltonian_plus_parallel_Exercise_}  
to show that $\g \oplus \goth h=\goth s$.
\end{enumerate}

\item
Denote by ${\goth h}$ the algebra of Hamiltonian
holomorphic vector fields (they are a posteriori Killing)
on a compact K\"ahler manifold $(M,I, \omega)$,
and by ${\goth s}$ the Lie algebra of 
holomorphic Killing vector fields.\index[terms]{vector field!Killing}\index[terms]{vector field!holomorphic}
Prove that $\goth h$ is normal in $\goth s$,
and the quotient ${\goth s}/{\goth h}$
is isomorphic to the Lie algebra $\g$
of holomorphic parallel vector fields\index[terms]{vector field!parallel}\index[terms]{vector field!holomorphic}
(Exercise \ref{_Killing_holomorphic_is_parallel_Exercise_}).
Prove that $\goth s$ can be obtained as
a semi-direct product of $\g$ and ${\goth h}$.

\item\label{_holo_Killing_fixed_point_Exercise_}
Prove that a holomorphic Killing vector field\index[terms]{vector field!Killing}\index[terms]{vector field!holomorphic} on a compact
K\"ahler manifold is Hamiltonian if and only if it has a fixed point.

{\em Hint:} Use the decomposition  $\g \oplus \goth h=\goth s$.

\end{enumerate}


\chapter{Connections in vector bundles and the Frobenius
  theorem}


{{\setlength\epigraphwidth{0.9\linewidth}
\epigraph{\it
\qquad 
    A man once said: Why such reluctance? If you only
    followed the parables you yourselves would become
    parables and with that rid of all your daily cares.
\smallskip

\qquad 
    Another said: I bet that is also a parable.
\smallskip

\qquad 
    The first said: You have won.
\smallskip

\qquad 
    The second said: But unfortunately only in parable.
\smallskip

\qquad 
    The first said: No, in reality: in parable you have
    lost.}{\sc\scriptsize Franz Kafka, ``On Parables'' 
(1922), translation by Willa and Edwin Muir}
}}

\section{Introduction}
\label{_Bundles_Intro_Section_}

This is the second chapter containing preliminaries.
We tell a few basic notions of differential geometry,
related to vector bundles, connections and foliations.

Vector bundles is one of the central notions of
differential geometry, though many textbooks on
differential geometry barely mention them.
The reason is clear: the subject is hard to adjust to 
and takes more effort to learn than most other concepts.
If one works locally, it is sometimes easier to
avoid mentioning bundles at all. 

A proper introduction to vector bundles would
take us a long while, so we just mention some basic
notions. 

A vector bundle over a smooth manifold can be 
defined in several ways: as a locally free sheaf of modules over
smooth functions, or locally with the charts, 
trivializations and gluing functions. Finally, one could apply 
 Serre--Swan theorem\index[terms]{theorem!Serre--Swan} and say that a vector bundle
on $M$ is a projective module over the ring $C^\infty M$
of smooth functions. This is the approach used
by some books on K-theory (\cite{_Karoubi_}).\index[terms]{K-theory}

The last definition motivates our convention,
that is  somewhat controversial: we use the
same letter to denote the bundle and the
space of its sections. Then the differential
operators on bundles, such as the connection,
are denoted by $\nabla:\; B \arrow B \otimes \Lambda^1 M$
and not $\nabla:\; \Gamma(B) \arrow \Gamma(B \otimes\Lambda^1 M)$.

This forces us to distinguish between ``the bundle''
and its ``total space''; we treat a bundle $B$ and
its total space $\Tot(B)$ as separate
objects, related by the procedures of taking ``total space
of a bundle $B$'' and ``the space of sections of the
projection $\pi:\; \Tot(B) \arrow M$''.

It is suggested that the reader takes a course of 
differential geometry focusing on bundles (\cite{_Taubes:Bundles_}) or that they 
read a separate textbook on bundles, such
as \cite{_Kobayashi_Bundles_}.

Connections on vector bundles (sometimes called
``covariant derivatives'') is one of the central
notions of this topic. We assume a basic knowledge
of this subject. However, we need a more general
notion, the {\em Ehresmann connection} (\ref{_Ehresmann_conne_Definition_}).
Let $\pi:\; M \arrow X$ be a smooth submersion.
The tangent bundle $TM$ contains the sub-bundle
$T_\pi M$ of fibrewise tangent vectors; an Ehresmann
connection is a choice of a direct sum decomposition 
$TM= T_\pi M\oplus T_\hor M$.

It is not hard to see that an Ehresmann
connection on a proper smooth fibration
defines diffeomorphisms between the fibres.
Indeed, given a smooth path in the base, the
tangent vector field to this path can be uniquely lifted
to a section of $T_\hor M$. This vector field integrates
to a flow of diffeomorphisms that identifies the fibres.

A connection on a vector bundle induces a special kind of 
Ehresmann connection on its total space, called {\em a
linear Ehresmann connection}\index[terms]{connection!Ehresmann}
(\ref{_linear_Ehresmann_Definition_}). \index[terms]{connection!Ehresmann!linear}
Conversely, a linear 
Ehresmann connection on the total space defines a
connection on a vector bundle. We explain this
correspondence in some detail, and give an 
interpretation of the curvature of a connection on a vector bundle
in terms of the corresponding Ehresmann connection.

This correspondence will be used in 
Chapter \ref{_holomorphic_bundles_Chapter_}
to compute $dd^c$ of the square length function on the
total space of an Hermitian holomorphic vector bundle.
In the present chapter, we use the curvature of the Ehresmann connections\index[terms]{connection!Ehresmann}\index[terms]{correspondence!Riemann--Hilbert}
to deduce the Riemann--Hilbert correspondence from the Frobenius
theorem.\index[terms]{theorem!Frobenius}

Frobenius theorem is a result of basic calculus on 
manifolds, characterizing smooth foliations in terms
of their tangent sub-bundle. We give an independent
elementary proof of this result, for completeness.

The Riemann--Hilbert correspondence is the functorial equivalence
between the category of local systems and the category
of flat vector bundles with connection. The term 
was used by M. \index[persons]{Kashiwara, M.} Kashiwara (\cite{_Kashiwara:RH_}) in the
context of the correspondence between perverse
sheaves and holonomic D-modules.

The category
of local systems is equivalent to the category of the
representations of the fundamental group.\index[terms]{fundamental group} The Riemann--Hilbert
correspondence  is often understood (\cite{_Simpson:RH_}) as a smooth map between the moduli space of flat
connections to the ``character variety''
(the moduli of representations of the fundamental 
group).\index[terms]{correspondence!Riemann--Hilbert}

Ehresmann connections\index[terms]{connection!Ehresmann} provide a very simple
way of producing the Riemann--Hilbert correspondence.
Indeed, the vanishing of the curvature of the
bundle is equivalent to the integrability of the
horizontal distribution of the corresponding
Ehresmann connection. Applying the Frobenius theorem,\index[terms]{theorem!Frobenius}
we obtain that the horizontal distribution is
a foliation, and hence any flat bundle is locally 
generated by its parallel sections.\index[terms]{section!parallel} The sheaf of parallel
sections is locally constant, that is, it is a 
local system. \index[terms]{sheaf!locally constant}

\section{Connections in vector bundles}

The notions we recall here are well-known and can be found in standard references such as \cite{griha}, \cite{_oss_}. For sheaves and their applications, use \cite{gode}.

\hfill

\noindent{\bf  Notation:}
Let $M$ be a smooth manifold, $TM$ its tangent bundle,
$\Lambda^iM$ the bundle of differential $i$-forms, and 
${C}^{\infty} (M)$ the algebra of smooth functions. 

\hfill

\definition \label{vecsheaf} A (smooth) {\bf vector bundle} on a smooth manifold
	is a locally free sheaf of ${C}^{\infty}(M)$-modules.\index[terms]{bundle!vector bundle}

\hfill

\remark  {{\bf The space of sections} of a bundle $B$ is 
denoted by $\Gamma(B)$, or, sometimes, by the same letter
$B$.  A section $b$ of a bundle $B$ is often 
		denoted as $b\in B$ instead 
of $b\in\Gamma(B)$}.\index[terms]{bundle!vector bundle!section of}
	We put 
	$\Lambda^iM$ for $\Gamma( \Lambda^iM)$.

\hfill

\definition 
	A {\bf connection} in a vector bundle $B$ is
	a map $\Gamma(B) \stackrel \nabla \arrow \Lambda^1 M \otimes \Gamma(B)$ which
	satisfies \[ \nabla(fb) = df \otimes b + f \nabla b\]
	for all $b\in \Gamma(B)$, $f\in {C}^{\infty} (M)$.
\index[terms]{connection!in a vector bundle}

\hfill

\remark  A connection $\nabla$ in $B$ induces
	a connection $B^* \stackrel {\nabla^*} \arrow \Lambda^1 M \otimes \Gamma(B^*)$
	in the dual bundle, by the formula
\begin{equation}\label{_dual_bundle_connection_Equation_}
	d(\langle b, \beta\rangle) = \langle \nabla b, \beta\rangle+
	\langle b, \nabla^*\beta\rangle
\end{equation}
	These connections are usually denoted {\em   by the same letter $\nabla$.}

\hfill

\definition
If $X\in TM$, we shall write $\nabla_Xb$ for $(\nabla b)(X)$. The operator $\nabla_X:\Gamma(B)\to\Gamma(B)$ is called {\bf the covariant derivative of $b$ along $X$}.

\hfill

\remark 
	For any tensor bundle 
	${\cal B}_1:=
	B^*\otimes B^* \otimes \cdots \otimes B^* \otimes B\otimes B \otimes\cdots \otimes B$
	{a connection on $B$ defines a connection on ${\cal B}_1$} that acts according to
	the Leibniz formula:
	\[
	\nabla(b_1 \otimes b_2) = \nabla(b_1) \otimes b_2 + b_1 \otimes \nabla(b_2).
	\]
	In particular, a connection in $B$ induces connections in the Grassmann algebra bundle of differential forms on $B$.

\section{Curvature of a connection}

\definition
Let $\nabla:\; B \arrow B \otimes \Lambda^1 M$ be a connection
on a vector bundle $B$. We extend $\nabla$ to an operator
\[
B \stackrel{\nabla}\arrow \Lambda^{1}(M)\otimes B
\stackrel{\nabla}\arrow \Lambda^{2}(M)\otimes B 
\stackrel{\nabla}\arrow \Lambda^{3}(M)\otimes B \stackrel{\nabla}\arrow \cdots
\]
using the Leibniz identity\index[terms]{Leibniz identity}
$\nabla(\eta \otimes b) = d\eta\otimes b + (-1)^{\tilde \eta} \eta \wedge \nabla b$.
Then the operator $\nabla^2:\; B \arrow B\otimes \Lambda^{2}(M)$
is called {\bf the curvature} of $\nabla$.

\hfill

\remark
{The algebra of differential forms
with coefficients in $\End B$ acts on
$\Lambda^* M \otimes B$} via
$\eta \otimes a (\eta' \otimes b) = \eta \wedge \eta'
\otimes a(b)$, where
$a\in \End(B)$, $\eta, \eta'\in \Lambda^* M$, and $b\in B$.

\hfill

\remark
$\nabla^2(fb) = d^2 f b + df \wedge \nabla b - df \wedge
\nabla b + f \nabla^2 b$, and hence  {
the curvature is a $C^\infty (M)$-linear operator.}
{We shall consider the curvature $B$ as 
a 2-form with values in $\End B$}.
Then $\Theta_B:= \nabla^2  \in \Lambda^2 M \otimes \End B$,
where  an $\End(B)$-valued form acts on $\Lambda^* M \otimes B$
as above.

\hfill

The following claim is well known. We need it to be able
to compute the curvature explicitly in terms of a covariant derivative.
The proof we give is very detailed, and the reader is advised to 
treat this claim as an exercise.

\hfill

\claim 
Let $X, Y\in TM$ be vector fields, 
$(B, \nabla)$ a bundle with connection, and $b\in B$ a section of $B$.
Consider the operator
\[ \Theta_B^* (X, Y, b):= \nabla_X\nabla_Yb-\nabla_Y\nabla_Xb-\nabla_{[X,Y]}b
\]
Then: 
\[ \Theta_B^* (X, Y, b)=\Theta_B(X,Y)(b),\] 
where $\Theta_B$ is the curvature of $B$.

\hfill

\pstep We prove that $\Theta_B^* (X, Y, b)$ is $C^\infty M$-linear in $X$ and $Y$. By antisymmetry, it is enough to prove the linearity in $X$. For $f\in C^\infty M$ and $Z\in TM$ we have 
$$\nabla_{fZ}b=(\nabla b)(fZ)=f(\nabla b)(Z)=f\nabla_Z b.$$ Then:
\begin{equation*}
\begin{split}
\nabla_{fX}(\nabla_Yb)&=f\nabla_X\nabla_Yb,\\
\nabla_Y(\nabla_{fY}b)&=\nabla_Y(f\nabla_Xb)=\nabla(f\nabla_Xb)(Y)=
(df\otimes\nabla_Xb+\nabla(\nabla_Xb))(Y)\\
&=\Lie_Y(f)\nabla_Xb+\nabla_Y\nabla_Xb,\\
\nabla_{[fX,Y]}b&=\nabla_{f[X,Y]-\Lie_Y(f)X}b=f\nabla_{[X,Y]}b-\Lie_Y(f)\nabla_Xb,
\end{split}
\end{equation*}
hence $\Theta_B^* (fX, Y, b)=f\Theta_B^* (X, Y, b)$.

\hfill

{\bf Step 2:}
Since $\Theta^*_B(X,Y,b)$ and $\Theta_B(X,Y)b$ are linear
in $X, Y$, it would suffice to prove this equality
for coordinate vector fields $X, Y$. In particular,
we may assume that $[X,Y]=0$.

Consider the operator
$i_X:\; \Lambda^iM \otimes B \arrow \Lambda^{i-1}M \otimes B$ 
of contraction with a vector field $X$.
Writing $\nabla = d+ A$, where $A\in \Lambda^1 M \otimes \End B$,
we obtain $\nabla_X = \Lie_X + A(X)$, which gives
$[\nabla_X, i_Y]= [\Lie_X, i_Y]= 0$ when $X, Y$ are coordinate
vector fields. Then:
\begin{equation*}
\begin{split}
\nabla_X\nabla_Yb&=\nabla_X(\Lie_Yb+A(Y)b)\\
&=\Lie_X(\Lie_Yb)+\Lie_X(A(Y)b)\\
&+A(Y)\Lie_Xb+A(Y)\Lie_Xb+A(X)\wedge A(Y)b,\\
\nabla_{[X,Y]}b&=\Lie_{[X,Y]}b+A([X,Y]).
\end{split}
\end{equation*}
Then: 
$$\Theta^*_B(X,Y,b)=\left(\Lie_X(A(Y))-\Lie_Y(A(X))+A(X)A(Y)-A(X)A(Y)\right)b$$
since $[\Lie_X,\Lie_Y]=\Lie_{[X,Y]}=0$ and $[\Lie_X,i_Y]=i_{[X,Y]}=0$ because $X,Y$ are coordinate vector fields.
On the other hand, $\nabla^2=dA+A\wedge A$, and hence
\begin{equation*}
\begin{split}
\Theta_B(X,Y)b&=(dA+A\wedge A)(X,Y)b\\
&=\left(\Lie_X(A(Y))-\Lie_Y(A(X))-A([X,Y])\right)b\\
&+(A(X) A(Y)-A(Y)A(X))b\\
&=\left(\Lie_X(A(Y))-\Lie_Y(A(X))+A(X)A(Y)-A(X)A(Y)\right)b	
\end{split}
\end{equation*}		
and the proof is complete.	
%
%
\endproof

\section{Ehresmann connections}

\subsection{Ehresmann connection on 
smooth fibrations}\index[terms]{connection!Ehresmann}

\definition\label{_Ehresmann_conne_Definition_}
Let $\pi:\;M \arrow Z$ be a smooth fibration, with
$T_\pi M$ {\bf the bundle of vertical tangent vectors}
(vectors tangent to the fibres of $\pi$).
An {\bf Ehresmann connection} on $\pi$
is a sub-bundle $T_\hor M\subset TM$ such that
$TM= T_\hor M\oplus T_\pi M$.
The {\bf parallel transport}\index[terms]{parallel transport} along the path $\gamma:\; [0, a]\arrow Z$
associated with the Ehresmann connection\index[terms]{connection!Ehresmann}
is a diffeomorphism $V_t:\;  \pi^{-1}(\gamma(0))\arrow  \pi^{-1}(\gamma(t))$
smoothly depending on $t\in [0, a]$ and 
satisfying $\frac {dV_t}{dt}\in T_\hor M$.

\hfill

The following claim is well known and easy to prove 
using basic results on solutions of ODEs.
  
\hfill

\claim
Let $\pi:\;M \arrow Z$ be a smooth fibration with compact fibres.
Then the parallel transport, associated with the Ehresmann
connection, always exists.\index[terms]{parallel transport}

\subsection{Linear Ehresmann connections on vector bundles}\index[terms]{connection!Ehresmann!linear}

\definition
Let $\pi_1:\;M_1 \arrow Z$, $\pi_2:\;M_2 \arrow Z$
be smooth fibrations. The {\bf fibred product} $M_1 \times_Z M_2$
is 
\[
M_1 \times_Z M_2:= \{(u, v)\in M_1 \times M_2 \ \ |\ \ \pi_1(u)=\pi_2(v)\}.
\]
It is also a smooth fibration over $Z$.

\hfill

\remark
Suppose that $\pi_1, \pi_2$ are smooth fibrations equipped with
the Ehresmann connections.\index[terms]{connection!Ehresmann}
Then $M_1 \times_Z M_2$
is equipped with an Ehresmann connection
defined as $T_\hor M_1 \times_Z M_2=T_\hor M_1 \times T_\hor M_2$.

\hfill

\definition\label{_linear_Ehresmann_Definition_}
Let $B$ be a vector bundle on $M$ and $\Tot B \stackrel \pi \arrow M$
its total space. An Ehresmann connection\index[terms]{connection!Ehresmann}
on $\pi$ is called {\bf linear}\index[terms]{connection!Ehresmann!linear} if it is preserved by 
the homothety map $\Tot B \arrow \Tot B$ 
mapping $v$ to $\lambda v$ and by the addition map
$\Tot (B\oplus B) \arrow \Tot B$,
that is, addition preserves horizontal vectors.

\hfill

\remark
Let $A, B$ be sets mapped to $C$.
Recall that {\bf the fibred product}  $A\times_C B$
is the set of all pairs $(a\in A ,b\in B)$ that are 
mapped to the same element of $C$. Then
$\Tot (B\oplus B) = \Tot B\times_M  \Tot B$.

\hfill

\definition
Let $X\in TM$ be a vector field on $M$.
Recall that the projection $\pi:\; \Tot B \arrow M$
defines an isomorphism $T_\hor \Tot B\restrict x\tilde\arrow T_{\pi(x)}M$,
where we denote with $T_\hor \Tot B\restrict x$  the fibre of
$T_\hor \Tot B$ in $x\in \Tot B$. Therefore, there
exists a unique vector field $X_\hor\in T_\hor \Tot B$
such that $d\pi(X_\hor)=X$, called {\bf the horizontal
lift} of $X$. 

\hfill

\claim\label{_Ehresmann_linear_via_flow_Claim_}
An Ehresmann connection\index[terms]{connection!Ehresmann}\index[terms]{connection!Ehresmann!linear} is linear
if and only if the flow $V_t$ of any 
horizontal lift $X_\hor\in T_\hor \Tot B$
is compatible with the structure of 
vector space on the fibres of $\pi$. 

\hfill

\proof
Indeed, since the Ehresmann connection is
invariant under the homothety map $\Tot B \arrow \Tot B$ 
mapping $v$ to $\lambda v$, any horizontal
lift is also invariant under this homothety.
Therefore, $V_t$ commutes with this action.
Similarly, the addition maps a 
horizontal lift of $X$ to $\Tot(B\oplus B)$
to the horizontal lift of $X$ to $\Tot B$,
hence commutes with the corresponding flow.
\endproof

\hfill

\remark\label{_dual_connec_Remark_}
The Ehresmann
connection on $\Tot B$ defines the {\bf dual
Ehresmann connection}\index[terms]{connection!Ehresmann!dual} on $\Tot B^*$
in such a way that the differential of the pairing map\index[terms]{pairing}
$\Tot B^*\times_M \Tot B \arrow \Tot C^\infty M= \C \times M$
maps horizontal vectors to horizontal vectors.
It is not hard to see that the dual Ehresmann
connection to an Ehresmann connection \index[terms]{connection!Ehresmann}\index[terms]{connection!Ehresmann!dual}
\begin{equation}\label{_Ehresmann_splitting_Equation_}
T\Tot B= T_\hor \Tot B\oplus T_\pi \Tot B
\end{equation} \index[terms]{connection!Ehresmann!linear}
is always well-defined, and is linear
when \eqref{_Ehresmann_splitting_Equation_} is linear.
 
\hfill

Let $B$ be a vector bundle, and 
$f$ a function on $\Tot B$ that is 
linear on all fibres of $\pi$. Such a function
is called {\bf fibrewise linear}.
If the Ehresmann connection\index[terms]{connection!Ehresmann!linear} on $\Tot B$ is linear,
the Lie derivative of a fibrewise linear function along
a horizontal lift of vector field on $M$ is again fibrewise linear.
Indeed, the corresponding diffeomorphism flow
preserves addition and multiplication by a number 
(\ref{_Ehresmann_linear_via_flow_Claim_}),
hence induces a linear map on fibres.

\hfill


\claim \label{_connection_via_fibrewise_linear_Claim_}
This procedure defines a connection on the bundle $B^*$
of fibrewise linear functions
on $\Tot B$.

\hfill

\proof
Recall that a vector bundle over a manifold $M$
is a locally free sheaf of $C^\infty(M)$-modules.
Let $\cal F$ be the sheaf of fibrewise linear functions on $\Tot(B)$,
and ${\cal F}_0:= \pi_* \cal F$ its pushforward to $M$.
This sheaf is isomorphic to $B^*$. Indeed,
any section of $B^*$ gives a fibrewise linear 
function on $\Tot B$ and vice versa. We identify
fibrewise linear functions on $\Tot B$ with the
sections of $B^*$.
Given a  section $\xi \in B^*$,
consider $\xi$ as a function on $\Tot B$. 
Let $X_\hor$ be the horizontal lift of a vector field $X\in TM$. 
As shown above, the Lie derivative $\Lie_{X_\hor}(\xi)$
is fibrewise linear. This is used to 
define the operator $\nabla_X \xi:= \Lie_{X_\hor}(\xi)$.
This map satisfies the Leibniz identity \index[terms]{Leibniz identity}
$\nabla_X (f\xi)=\Lie_X(f) \xi + f\nabla_X(\xi)$
by construction, 
and is additive because the Ehresmann 
connection is compatible with the addition.
Therefore, $\xi \mapsto \nabla_X \xi$ gives a connection on $B^*$.
\endproof

\hfill

\proposition\label{_Ehresmann_vs_vector_bundle_Proposition_}
Let $B$ be a vector bundle, $\pi:\; \Tot B \arrow M$
its total space, and $T\Tot B= T_\hor \Tot B\oplus T_\pi \Tot B$
a linear Ehresmann connection.\index[terms]{connection!Ehresmann!linear} Construct the dual connection on\index[terms]{connection!Ehresmann!dual}
$\Tot B^*$ as in \ref{_dual_connec_Remark_}. 
After we  identify $B$ with the fibrewise linear functions on $\Tot B^*$,
the Ehresmann connection on $\Tot B^*$ gives
 a connection on the vector bundle $B$ 
(\ref{_connection_via_fibrewise_linear_Claim_}).
Then this procedure defines a bijective correspondence
between linear Ehresmann connections on $\Tot B$
and connections on $B$.

\hfill

{\bf Proof:} Let $(B,\nabla)$ 
be a vector bundle over $M$ with a connection.
It is not hard to see that the derivations on $C^\infty(\R^n)$ are 
the same as $C^\infty(\R^n)$-valued 
derivations of the polynomial ring $\R[x_1, ..., x_n]$.
Therefore, a vector field on $\Tot(B)$ is the same
as a derivation on the algebra of fibrewise
polynomial functions on $\Tot(B)$. However,
a homogeneous fibrewise polynomial function
of degree $k$ is the same  as a section
of the symmetric tensor power $\Sym^k B^*$.
Given a vector field $X\in TM$,
we consider the operator 
$\nabla_X:\; \Sym^* B^*\arrow \Sym^* B^*$
as a derivation of the algebra of fibrewise
polynomial functions on $\Tot(B)$. This
defines a vector field $X_1\in T \Tot B$;
after we construct the Ehresmann connection,\index[terms]{connection!Ehresmann}
$X_1$ will become the horizontal lift of $X$.

The vector field $X_1 \in T \Tot B$ is mapped
to $X$ under $d\pi$, because it acts on
$\pi^* C^\infty M$ the same way as $X$.
Consider the sub-bundle $T_\hor \Tot B$
generated by all such vectors. By construction,
the map $d\pi$ induces an isomorphism of the
fibres of $T_\hor \Tot B$ and $TM$; therefore, 
$T_\hor \Tot B\subset T\Tot B$ defines an Ehresmann 
connection. 

This construction is clearly the inverse of the one in
\ref{_connection_via_fibrewise_linear_Claim_}; 
the proof is left to the reader.
 \endproof

\section{Frobenius form and Frobenius theorem}\index[terms]{form!Frobenius}

\definition
{\bf A distribution} on a manifold is 
a sub-bundle $B\subset TM$.\index[terms]{distribution}

\hfill

\remark
Let  $\Pi:\; TM \arrow TM/ B$ be the projection, and
$x, y \in B$ some vector fields. Then 
$[fx, y]= f[x,y] - \Lie_y (f) x$. This implies that
$\Pi([x,y])$  is $C^\infty(M)$-linear as a function of $x$ and $y$.

\hfill

\definition\label{_Frobenius_form_Definition_}
The map $[B,B]\arrow TM/B$, sending $x, y$ to $\Pi([x, y])$
is called the {\bf Frobenius bracket} (or {\bf Frobenius form}); \index[terms]{form!Frobenius}\index[terms]{bracket!Frobenius}
it is a skew-symmetric $C^\infty(M)$-linear form on $B$ with values in $TM/B$.

\hfill

\definition
A distribution is called {\bf  integrable},
 or {\bf holonomic}, or {\bf  involutive}, if
its Frobenius form vanishes.\index[terms]{distribution!integrable}\index[terms]{distribution!involutive}\index[terms]{distribution!holonomic}

\hfill

\theorem The {\bf (Frobenius theorem)}\index[terms]{theorem!Frobenius}
Let $B\subset TM$ be a sub-bundle. Then $B$ is involutive
if and only if each point $x\in M$ has a neighbourhood
$U\ni x$ and a smooth submersion $U\stackrel \pi \arrow V$ 
such that $B$ is its vertical tangent space: $B= T_\pi M$.

\hfill

We prove Frobenius theorem at the end of this section.

\hfill

\definition\label{_leaf_foliation_Definition_}
The fibres of $\pi$ are called {\bf leaves},
or {\bf integral submanifolds} of the distribution $B$.
Globally on $M$, {\bf a leaf of $B$} is a maximal
connected manifold $Z\hookrightarrow M$ that is  immersed to $M$ 
and tangent to $B$ at each point.
A distribution for which Frobenius theorem holds is called
{\bf integrable}. If $B$ is integrable, the set of 
its leaves is called {\bf a foliation}. The leaves\index[terms]{foliation}
are manifolds that are  immersed in $M$. \footnote{The leaves
are immersed, but not necessarily closed.
Quite often it occurs that some (or all) leaves 
of a foliation are dense in $M$.}

\hfill

\remark
Let $U\subset M$ be an open subset.
A leaf of the foliation ${\cal F}$
in $U\subset M$ is a connected component of 
$F\cap U$, where $F$ is a leaf of ${\cal F}$ 
in $M$. However, the intersection 
$F\cap U$ can have many connected
components. This means that the set of leaves
of ${\cal F}$ on $U$ is different from the
set of leaves of ${\cal F}$ on $M$.

\hfill

\remark
If $B$ is the tangent bundle to a foliation,
then $[B,B]\subset B$ is clear, because it is
true leaf-wise. To prove the Frobenius theorem we
need only to prove the existence of the foliation
tangent to $B$, for any involutive $B$.

\hfill

\remark
To prove the Frobenius theorem for $B\subset TM$, 
 it suffices to show that
each point is contained in an integral submanifold. In this case,
the smooth submersion $U\stackrel \pi \arrow V$  is the projection
to the leaf space of $B$.

\hfill

{\bf Proof of the Frobenius theorem. Step 1:}\index[terms]{theorem!Frobenius}
Suppose that $G$ is a Lie group
acting on a manifold $M$. Assume that the vector fields from 
the Lie algebra of $G$ generate a sub-bundle $B\subset TM$. 
Then $B$ is integrable,
that is, Frobenius theorem holds of $B\subset TM$.
Indeed, the  orbits of the $G$-action are tangent to
$B\subset TM$.

\hfill

{\bf Step 2:} 
Let $u, v$ be commuting vector fields on a manifold $M$,
and $e^{tu}$, $e^{tv}$ be corresponding diffeomorphism flows.
Then $e^{tu}$, $e^{tv}$ commute. This easily follows by
taking a coordinate system such that $u$ is the coordinate
vector field.

\hfill

{\bf Step 3:} 
The commutator of vector fields in $B$ belongs to $B$; 
however, this does not immediately produce any finite-dimensional
Lie algebra: it is not obvious that any subalgebra
generated by such vector fields is finite-dimensional. 
To produce a Lie group with orbits tangent to $B$, we need
to find a collection $\xi_1, ..., \xi_k \in B$ 
of vector fields generating $B$ and make sure that
the $\xi_1,..., \xi_k$ generate a finite-dimensional
Lie algebra.

The statement of the Frobenius theorem is
local, and hence  we may replace $M$ with a small neighbourhood of a given
point. We are going to show that $B$ locally has a basis of commuting
vector fields. By Step 2, these vector fields can be locally 
integrated to a commutative group action, and the Frobenius 
theorem follows from Step 1.

\hfill

{\bf Step 4:}  
Let $\sigma:\; M \arrow M_1$ be a smooth submersion, $d\sigma:\;
T_x M \arrow T_{\sigma(x)}M_1$ its differential, and $v\in
TM$ a vector field that satisfies
\begin{equation}\label{_projecting_vector_field_Equation_}
d\sigma(v)\restrict x= d\sigma(v)\restrict y  
\end{equation}
for any $x, y\in \sigma^{-1}(z)$ and any $z\in M_1$.
In this case, the vector field $d\sigma(v)$
is well-defined on $M_1$. Given two vector
fields $u$ and $v$ that satisfy \eqref{_projecting_vector_field_Equation_},
we can easily check that the commutator $[u, v]$
also satisfies \eqref{_projecting_vector_field_Equation_}, and, moreover,
$d\sigma([u, v]) = [d\sigma(u), d\sigma(v)]$.

\hfill

{\bf Step 5:} Now we can finish the proof of Frobenius theorem.
We need to produce, locally in $M$, a basis of commuting vector fields
$\xi_i \in B$. We start by producing (locally in $M$) an auxiliary
submersion $\sigma$, with the fibres that are  
complementary to $B$. To define such a submersion,
we put coordinates locally on $M$, identifying $M$ with
an open subset in $\R^n$, and take
a linear map $\sigma:\; M \arrow M_1=\R^{\dim B}$ such that
$d\sigma:\; B\restrict x \arrow T_{\sigma(x)}M_1$
is an isomorphism at some $x\in M$. 

Then $d\sigma:\; B\restrict x \tilde\arrow T_{\sigma(x)}M_1$ 
 is an isomorphism in a neighbourhood of $x$; replacing
$M$ with a smaller open set, we may assume that 
$d\sigma:\; B\restrict x \tilde\arrow T_{\sigma(x)}M_1$  is an
isomorphism everywhere on $M$. Let
$\zeta_1, ..., \zeta_k$ be the coordinate vector fields on $M_1$.

Since $d\sigma:\; B\restrict x \arrow T_{\sigma(x)}M_1$
is an isomorphism, there exist unique vector fields
$\xi_1, ..., \xi_k\in B\subset TM$ such that $d\sigma(\xi_i)=\zeta_i$.
By Step 4, $d\sigma([\xi_i, \xi_j])= [\zeta_i, \zeta_j]=0$. 
Since $B$ is involutive, 
the commutator $[\xi_i, \xi_j]$ is a section of $B$.
Now, the map $d\sigma:\; B\restrict x \arrow T_{\sigma(x)}M_1$
is an isomorphism, and therefore, the vanishing of
 $d\sigma([\xi_1, \xi_j])$
implies $[\xi_1, \xi_j]=0$. We have constructed
a basis of commuting vector fields in $B$ and
finished the proof of Frobenius theorem.
\endproof

\section{Basic forms}

\definition\label{_basic_form: definition_}
Let $B\subset TM$ be an involutive
sub-bundle, that is, a sub-bundle
that satisfies $[B, B]\subset B$.
A differential form $\alpha \in \Lambda^*(M)$
is called {\bf basic with respect to $B$}
if $\Lie_X\alpha= 0$ and $i_X \alpha=0$ 
for any vector field $X\in B$, where 
$i_X:\; \Lambda^k(M) \arrow \Lambda^{k-1}(M)$
is the contraction (interior product) with $X$. 
\index[terms]{form!basic}

\hfill

\remark\label{_basic_for_closed-Remark_}
By the Cartan formula, 
$\Lie_X\alpha= d (i_X \alpha)+ i_X (d \alpha)$.
Therefore, whenever $\alpha$ is closed, the condition
$\Lie_X\alpha=0$ is implied by $i_X \alpha=0$.
We obtain that a closed form is basic with respect to $B$
if and only if it vanishes on $B$.

\hfill

The following theorem is well known.

\hfill

\theorem\label{_basic_forms_Frobenius_lifted_Theorem_}
Let $B\subset TM$ be an involutive
sub-bundle. By the Frobenius theorem,\index[terms]{theorem!Frobenius}
locally $M$ admits a smooth submersion
$\pi:\; M \arrow U$ such that $B$ is the 
tangent bundle to the fibres of $\pi$.
Then a form $\alpha$ is basic with respect
to $B$ if and only if $\alpha= \pi^* \alpha_0$
for some form $\alpha_0$ on $U$.

\hfill

\proof 
Let $x_1, ..., x_k, x_{k+1},..., x_n$
be coordinates on $M$ such that $\pi$ is
the forgetful map, sending $(x_1, ..., x_k, x_{k+1},..., x_n)$
to $(x_{k+1},..., x_n)$. Since $\ker \alpha \supset \ker \pi$,
the form 
\[ \alpha= \sum_{\{i_1, ..., i_l\}\subset \{k+1,..., n\}} f_{i_1, ..., i_l} dx_{i_1}\wedge... \wedge
    dx_{i_l}
\]
is a sum of differential
monomials in $dx_{k+1}, ..., dx_{n}$.
However, for each $i = 1, ..., k$, $\Lie_{d/dx_i}\alpha=0$,
hence the coefficients $f_{i_1, ..., i_l}$ of $\alpha$ depend only 
on the coordinates $x_{k+1},..., x_n$.
Therefore, $\alpha$ is lifted from $U$.
\endproof

\hfill

\remark\label{_multilinear_forms_Remark_}
By \ref{_basic_forms_Frobenius_lifted_Theorem_},
a differential form $\eta$ on $M$ is basic with respect to
$B$ if and only if it is locally obtained as a pullback
of forms on the leaf space.
In the same way, we could define {\bf a basic multilinear form}
(not necessarily exterior) $\eta \in (T^*M)^{\otimes n}$
on $M$ as a form on $M$ that is  locally obtained as a pullback
of multilinear forms on the leaf space.
The same proof as above can be used to 
show that $\eta$ is basic if and only if
$\eta$ satisfies
$\Lie_X \eta= 0$, and $\eta$ contracted
with $X$ on any of its arguments vanishes for all $X \in B$.

\section{The curvature of an Ehresmann connection}\index[terms]{connection!Ehresmann}

\definition\label{_curvature_Ehresmann_Definition_}
Consider a smooth fibration  $\pi:\; E \arrow M$, and
let $T_\hor E \subset TE$ be the horizontal sub-bundle
defining an Ehresmann connection. For any
$X, Y\in TM$,  consider the horizontal lifts
$X_\hor, Y_\hor\in T_\hor E$. Evaluating the Frobenius form\index[terms]{form!Frobenius}
associated with  $T_\hor E\subset TE$ on $X_\hor, Y_\hor$,
we obtain a section $\Psi(X_\hor, Y_\hor)\in TE/T_\hor E$.
Since the bundle $TE/T_\hor E$ is identified with
the vertical (fibrewise) tangent bundle $T_\pi E$,
we can consider $X, Y \arrow \Psi(X_\hor, Y_\hor)$
to be  a map from $\Lambda^2 TM$ to $T_\pi E$.
This map is called {\bf the curvature form\index[terms]{form!curvature}
of the Ehresmann connection}.

\hfill

For the next step, we need the following definition.

\hfill

\definition
A vector field $v$ on $\R^n$ is called {\bf linear} if
it satisfies  $v\restrict{x+y}=v\restrict{x}+v\restrict y$
and $v\restrict{\lambda x}= \lambda v\restrict{x}$
for any $x, y\in \R^n$. 

\hfill

\claim \label{_linear_v_fields_Claim_}
Linear vector fields on $\R^n$ and $\End(\R^n)$ are
in bijective correspondence.
Given $E\in \End(\R^n)$, the corresponding 
vector field $v$ associates to $x\in \R^n$
the vector $E(x)\in T_x\R^n=\R^n$.

\hfill

\proof Left as an exercise to the reader.
\endproof

\hfill

 Applying this
construction to vertical tangent vectors
to a total space of a vector bundle, we arrive at the notion of
a linear tangent vector field, as follows.

\hfill

\definition
Let $\pi:\; \Tot B \arrow M$ be a total space of a vector
bundle. A vertical vector field  $v \in T_\pi \Tot(B)$ 
is called {\bf fibrewise linear} if it is linear
on all fibres, that is, satisfies
$v\restrict{x+y}=v\restrict{x}+v\restrict y$
and $v\restrict{\lambda x}= \lambda v\restrict{x}$
for any $x, y$ in the same fibre of $\pi$.

\hfill

\remark
From \ref{_linear_v_fields_Claim_}
we obtain that the fibrewise linear vector fields are in bijective
correspondence with sections of $\End(B)$.
For any $E\in \End(B)$, the corresponding
vector field $v$ associates to $x\in \Tot(B)$
the vector $E(v)\in T_x \Tot(B)$.

\hfill

\proposition\label{_Ehresmann_curvature_bundles_Proposition_}
 Let $(B, \nabla)$ be a vector bundle with connection over $M$, 
$\pi:\; \Tot B \arrow M$ its total space, and 
$T\Tot B= T_\hor \Tot B\oplus T_\pi \Tot B$
the linear Ehresmann connection\index[terms]{connection!Ehresmann} associated with $\nabla$
as in \ref{_Ehresmann_vs_vector_bundle_Proposition_}.
Given $X, Y\in TM$, consider the curvature vector 
$\Psi(X_\hor, Y_\hor)\in T_\pi \Tot(B)$ defined as in 
\ref{_curvature_Ehresmann_Definition_}.
Then $\Psi(X_\hor, Y_\hor)$ is a fibrewise linear 
vector field. Moreover, the corresponding
section of $\End(B)$ is $\Theta_B(X,Y)$,
where $\Theta_B\in \Lambda^2(M) \otimes \End(B)$ 
is the curvature of the connection $\nabla$ on the 
vector bundle $B$.

\hfill

\pstep
The curvature of the linear Ehresmann connection
is fibrewise linear because the addition
and homothety preserve the Ehresmann connection,\index[terms]{connection!Ehresmann}
hence preserve its curvature. 

\hfill

{\bf Step 2:} Now we need to identify the
endomorphism $\Theta_B(X,Y)$ with an action
of the linear vector field  $\Psi(X_\hor, Y_\hor)$.
It is slightly more comfortable to prove it
for $B^*$; the proof for $B$ will follow by duality.
Recall that sections of $B^*$ are in bijective
correspondence with the fibrewise
linear functions $\alpha\in C^\infty \Tot B$.

By \ref{_Ehresmann_vs_vector_bundle_Proposition_}, the
connection (covariant derivative) $\nabla_X$ acts
on sections of $B^*$, interpreted as
fibrewise linear functions, as $\Lie_{X_{\hor}}$,
hence the curvature acts as $[X_\hor, Y_\hor]$.

For any commuting
vector fields $X, Y\in TM$, the Lie derivative of a fibrewise
linear function $\alpha\in C^\infty \Tot B$ 
along $X_\hor, Y_\hor$  gives $\nabla_X \alpha$ and $\nabla_Y\alpha$
by \ref{_Ehresmann_vs_vector_bundle_Proposition_}.
Therefore, their commutator $[X_\hor, Y_\hor]$ 
acting on fibrewise linear functions gives an operator
$[\nabla_X, \nabla_Y]=\Theta_B(X,Y)$. This implies
that derivation of $\alpha$ along the
vector field $\Psi(X_\hor, Y_\hor)\in T_\pi \Tot B$
acts in the same way as $\Theta_B(X,Y)$.
\endproof

\section{The Riemann--Hilbert correspondence}
\label{_Riemann_Hilbert_Section_}\index[terms]{correspondence!Riemann--Hilbert}

\subsection{Flat bundles and parallel sections}\index[terms]{section!parallel}

\definition (\index[persons]{Cartan, E.}Cartan, 1923)
Let $(B,\nabla)$ be a vector bundle with connection over $M$.
For each loop $\gamma$ based in $x\in M$, let 
$V_{\gamma, \nabla}:\; B\restrict x \arrow B\restrict x$
be the corresponding parallel transport\index[terms]{parallel transport} along the connection.
The {\bf holonomy group} of $(B,\nabla)$
is a group generated by $V_{\gamma, \nabla}$,
for all loops $\gamma$. If one takes all contractible
(homotopically trivial) loops  instead, $V_{\gamma, \nabla}$ generates
{\bf  the local holonomy}, or {\bf 
the restricted holonomy} group.\index[terms]{holonomy}\index[terms]{holonomy!local}

\hfill

\definition
Let $B$ be a vector bundle, and $\Psi$ a section of its tensor power.
We say that {\bf connection $\nabla$ preserves $\Psi$}
if $\nabla(\Psi)=0$. In this case, we also say that the tensor $\Psi$
is {\bf  parallel} with respect to the connection.

\hfill

\remark
$\nabla(\Psi)=0$ is equivalent to $\Psi$ being a solution of
$\nabla(\Psi)=0$ on each path $\gamma$. This means that 
the parallel transport preserves $\Psi$.

\hfill

We obtained:

\hfill

\corollary
Let $B$ be a bundle with connection.
A section of a tensor power of $B$ 
is parallel if and only if it is holonomy invariant.
\endproof

\hfill

%
%
%
%
%
%
%
%

\definition
A bundle with connection is {\bf flat} if 
the connection has vanishing curvature.\index[terms]{bundle!vector bundle!flat}

\hfill

\remark The next theorem claims that
a flat bundle is locally generated by parallel
sections. Clearly, in this case the sheaf of
parallel sections\index[terms]{section!parallel} is locally constant.

\hfill

\remark\label{_local_system_repres_Remark_}
Recall that a locally constant \index[terms]{sheaf!locally constant}
sheaf is also called ``a local system''.\index[terms]{local system}
A local system on $M$ is uniquely determined by
a representation of the fundamental group \index[terms]{fundamental group}$\pi_1(M)$.
In particular, any local system on a 
simply connected space is trivial.

\hfill

\theorem
Let $(B, \nabla)$ be a vector bundle with connection. 
Then the following conditions are equivalent:
\begin{description}
\item[(i)] The connection $\nabla$ is flat.\index[terms]{connection!flat}
\item[(ii)] The restricted holonomy group of $\nabla$ 
is trivial. 
\item[(iii)] $B$ is locally generated by its parallel
sections.
\end{description}

\hfill

\pstep
We prove that (iii) is equivalent to (i). 
Let $B$ be a flat bundle on $M$, and $X, Y \in TM$ commuting vector fields.
Then $\nabla_X:\; B \arrow B$ commutes with $\nabla_Y$. Then 
the Ehresmann connection\index[terms]{connection!Ehresmann} bundle $E_\nabla$ is generated
by commuting vector fields $\tau_\nabla(X)$, $\tau_\nabla(Y)$, ..., 
hence it is involutive. 
By Frobenius theorem, every\index[terms]{theorem!Frobenius}
point $b\in \Tot(B)$ is contained in a leaf of 
the corresponding foliation, tangent to $E_\nabla$.
In other words, the horizontal distribution $T_\hor
\Tot(B)$ is a foliation.
By \ref{_Ehresmann_vs_vector_bundle_Proposition_}, 
a leaf of this foliation is a parallel section \index[terms]{section!parallel}of $B$.
Since every point of $\Tot(B)$ is contained
in a leaf of the horizontal foliation,
$B$ is locally generated by its parallel sections.

Conversely, if $B$ is locally generated by
parallel sections,\index[terms]{section!parallel} its curvature vanishes.

\hfill

{\bf Step 2:} 
(ii) implies (iii): if the restricted 
holonomy is trivial, $B$ is locally generated by
the parallel sections. Indeed, locally in a contractible 
open subset $U\subset M$ the parallel
translation of a vector $v\in B\restrict x$ along a path
$\gamma\subset U$ connecting $x$ to $y$ is independent on  the 
choice of $\gamma$. Therefore, such translation produces
a parallel section\index[terms]{section!parallel} passing through any $v\in B\restrict x$. 

\hfill

{\bf Step 3:} 
(iii) implies (ii). The sheaf of parallel sections  of $B$ is locally constant. Indeed, for any
$x\in M$, there exists a collection of parallel
sections of $B$ in a neighbourhood $U\ni x$.
Therefore, the sheaf of parallel sections is
trivial in $U$.  However, the monodromy of a
local system over a contractible loop is trivial.
By definition, this monodromy generates the restricted 
holonomy group.\index[terms]{monodromy}
\endproof

\hfill

\remark
This theorem is a special case of the famous
Ambrose--Singer theorem, claiming that the local\index[terms]{theorem!Ambrose--Singer}
holonomy group of $(B,\nabla)$ in $x\in M$ is a Lie group with the Lie algebra
generated by all ``curvature elements''
$\Theta_B(x, y)\in \End_m(B)$ translated from
$m\in M$ to $x$ by parallel transport.\index[terms]{parallel transport}

\subsection{Local systems}

The following theorem is often called ``the
Riemann--Hilbert correspondence''. \index[terms]{correspondence!Riemann--Hilbert} 

\hfill

\theorem\label{_Riemann--Hilbert_Theorem_}
The category of locally constant sheaves of vector spaces 
 is naturally equivalent to the category of vector bundles on $M$
equipped with a flat connection. \index[terms]{sheaf!locally constant}\index[terms]{connection!flat}

\hfill

\pstep
Consider the constant sheaf $\R_M$ on $M$. This is a sheaf of rings, and
any locally constant sheaf is a sheaf of $\R_M$-modules.

Let ${\mathbb V}$ be a locally constant sheaf, and $B:= {\mathbb
  V}\otimes_{\R_M} C^\infty M$.
Since ${\mathbb V}$ is locally constant, the sheaf $B$ is a
locally free sheaf of $C^\infty M$-modules,
that is, a vector bundle. Let $U\subset M$ be an
open set such that ${\mathbb V}\restrict U$ is constant.
If $v_1, ..., v_n$ is a basis in ${\mathbb V}(U)$, all sections of $B(U)$
have the form $\sum_{i=1}^n f_i v_i$, where $f_i \in C^\infty U$.
Define the connection $\nabla$ by 
$\nabla\left (\sum_{i=1}^n f_i v_i\right) = \sum df_i \otimes v_i$.
This connection is flat because $d^2=0$. 
It is independent on  the choice of $v_i$
because the basis $\{v_i\}$ is 
defined canonically up to a matrix with constant coefficients.
We have constructed a functor from locally constant sheaves
to flat vector bundles.

\hfill

{\bf Step 2:}
Now let  $(B, \nabla)$ be a flat bundle over $M$.
The functor to locally constant sheaves takes 
 $U\subset M$ and maps it to the space of parallel
sections of $B$ over $U$. This defines a sheaf
${\mathbb B}(U)$. For any simply connected $U$, and any $x\in M$,
the space ${\mathbb B}(U)$ is identified with a vector space
$B\restrict x$, and hence  ${\mathbb B}(U)$ is locally constant.
Clearly, $B= {\mathbb B}\otimes_{\R_M} C^\infty M$, and hence 
this construction gives an inverse functor to 
${\mathbb V}\mapsto {\mathbb V}\otimes_{\R_M} C^\infty M$.
\endproof

\hfill

\definition \label{dnabla}
Let $(B, \nabla)$ be a flat bundle.
	The {\bf $B$-valued de Rham differential}\index[terms]{differential!de Rham}
	on the complex $ \Lambda^iM\otimes B \arrow\Lambda^{i+1}M\otimes B$
	is $d_\nabla(\eta\otimes b):= 
	d\eta\otimes b + (-1)^{\tilde\eta-1}\eta\wedge \nabla b$. It is easy to check that $d_\nabla^2=0$.
	
\hfill

\claim
The cohomology of the complex
$(\Lambda^*M\otimes B, d_\nabla)$ is isomorphic to the
cohomology of the local system ${\cal B}:=\ker \nabla$.
\index[terms]{local system!cohomology of}

\hfill

\proof Indeed, the complex
\[
B \stackrel \nabla \arrow B \otimes \Lambda^1 M
\stackrel {d_\nabla} \arrow B \otimes \Lambda^2 M
 \stackrel {d_\nabla} \arrow ...
\]
gives a fine resolution of the sheaf of \index[terms]{resolution!fine}
constant sections of $\ker \nabla$.
\endproof

\section{Exercises}

\begin{enumerate}[label=\textbf{\thechapter.\arabic*}.,ref=\thechapter.\arabic{enumi}]


\item Build a 4-dimensional smooth manifold
not admitting rank 3 distributions.

\item Prove that a compact $n$-dimensional torus admits
a rank 1 foliation with all leaves dense.

\item
Find a rank 1 sub-bundle $B\subset TS^3$ 
such that the corresponding foliation has 
non-compact leaves.


\item
Let $B$ be a vector bundle. Prove that $\Lambda^2 B$
is oriented, or find a counterexample.

\item
Construct a rank 2 vector bundle
not admitting a non-degenerate
bilinear form of signature (1,1),
or prove that such bundle does not exist. 

\item
Let $M = \R$. Find two vector fields $X, Y$ on $M$
such that the successive commutators of $X,Y$ 
generate an infinite-dimensional Lie algebra.


\item
Construct a 4-manifold $M$ and a rank 2 distribution
$B\subset TM$ such that $[B, B]$ has rank 3 and
$[[B, B], B]$ has rank 4.

\item Let $\pi:\; M \arrow X$ be a continuous map of manifolds.
Prove that $\pi$ is proper (preimages of compacts are compact)
if and only if all its fibres are compact and have the
same number of connected components or find a counterexample.

\item  Let $\pi:\; M \arrow X$ be a smooth submersion of manifolds.
Prove that when $\pi$ is proper, 
the fibres of $\pi$ are diffeomorphic. 
Find an example when $\pi$ is not proper, the
fibres are not diffeomorphic.

\item Prove that the Ehresmann connection\index[terms]{connection!Ehresmann}
defined on $\Tot B$ as in \ref{_Ehresmann_vs_vector_bundle_Proposition_}
 is indeed linear.

 \item Prove Ehresmann's fibration theorem: Let $f:M\arrow
  N$ be a smooth, surjective and  proper submersion of manifolds. Then
  $f$ is a locally trivial fibration.\index[terms]{theorem!Ehresmann's fibration}

\item Consider the polynomial ring 
$\R[x_1, ..., x_n]$ as a subring of $C^\infty(\R^n)$. 
Prove that any derivation of $C^\infty(\R^n)$
is uniquely determined by its restriction to $\R[x_1, ..., x_n]$.
Prove that any $C^\infty(\R^n)$-valued derivation on $\R[x_1, ..., x_n]$
defines a derivation of $C^\infty(\R^n)$.

\end{enumerate}


\chapter{Locally conformally \ka\ manifolds}\label{lck}

{\setlength\epigraphwidth{0.7\linewidth}
\epigraph{\em The learned discovers the Dao, duly obliges; \\
The learning discovers the Dao, and questions its potency;\\
The unlearned discovers the Dao, and roars into laughter,\\
Without the laughs, it would not be the Dao.}{\sc\scriptsize Tao Te Ching, by  Laozi}
}


\section{Introduction}
\label{_LCK_chapter_Intro_}

In this chapter,
we introduce the protagonist of this
book, the LCK (locally conformally K\"ahler) manifolds.
The study of LCK manifolds originates with I. \index[persons]{Vaisman, I.} Vaisman (\cite{va_isr}),
although the notion appeared, at least implicitly, in
previous works by other authors. Vaisman used the
geometric definition (the existence of an atlas with local
K\"ahler structures conformally related on intersections)
and proved that it is equivalent to the tensorial one;
in subsequent papers he proved that this is equivalent
to the existence of a K\"ahler cover with the deck group
acting by holomorphic homotheties. 

The term
``locally conformally K\"ahler'' is motivated by
the following definition: a manifold
is locally conformally K\"ahler if
it is covered by charts that are  K\"ahler,
and the gluing functions are conformal and 
holomorphic. It is called ``globally conformally
K\"ahler'' (GCK) if there is a global choice
of K\"ahler structure compatible with the gluing.
Generally, all LCK manifolds are tacitly
assumed to be not GCK. As follows from \index[terms]{theorem!Vaisman} Vaisman's
\ref{vailcknotk}, any LCK structure\index[terms]{structure!LCK}
on a compact K\"ahler-type complex manifold is itself GCK. 

Note that it makes very little
sense to speak of ``locally conformally K\"ahler manifolds''
in complex dimension 1, and we always exclude this case.
This assumption is universal and has no exceptions
throughout this book.

This chapter is mainly about definitions.
We introduce four equivalent definitions of LCK
manifolds: the one in terms of charts and
gluing maps, the tensor one, the K\"ahler homothety 
covering and the one that uses the K\"ahler
form with coefficients in a local system.
We prove that all these diverse notions are equivalent.

These proofs are rather formal. In the earlier version
of this chapter we  chose an informal approach,
that is  sufficient for most readers. We expect that
the rigorous proofs are self-evident for many of the readers,
but we include them for the sake of completeness.
The reader is invited to skip the formalities, as long
as the matters are already clear.

In this chapter, we do not assume profound knowledge of 
K\"ahler, complex and symplectic geometry\index[terms]{geometry!symplectic}, but some understanding of calculus on manifolds is necessary.

\section{Locally conformally symplectic  manifolds}\index[terms]{manifold!LCS}

\definition
Let $L$ be an oriented real line bundle\footnote{A line
  bundle is a one-dimensional bundle; a real line bundle is
  a bundle with 1-dimensional real fibres.} on $M$,
equipped with a flat connection, and let $\omega\in \Lambda^2M\otimes L$ be 
an $L$-valued differential form. \index[terms]{connection!flat}
We say that $\omega$ is {\bf non-degenerate} if for some (and then, for any)
trivialization of the bundle $L$ on an open set $U$, the corresponding
2-form taking values in $L \restrict U \cong \C^\infty M\restrict U$ 
is non-degenerate.

\hfill

\definition \label{_weight_for_lcs_}
We say that $(M,\omega, L)$
is {\bf  locally conformally symplectic} (LCS) if $\omega$
is a closed, non-degenerate form with coefficients in the local\index[terms]{manifold!LCS} system\index[terms]{local system} defined by $L$, that is, 
$d_\nabla\omega=0$ (see \ref{dnabla}).

In this situation, $L$ is called {\bf  the weight bundle}
of $(M,\omega)$, or {\bf the bundle of conformal weights}.\index[terms]{bundle!weight}

\hfill

\claim   
An oriented real line bundle $L$ can be smoothly 
trivialized.

\hfill

Indeed, choose a Riemannian metric on $L$. Then the set 
of positive unit vectors determines a nowhere degenerate section of $L$.

\hfill

\claim  
Let $(M,\omega, L)$ be a LCS
manifold\index[terms]{manifold!LCS}, $\psi$  non-degenerate section
of $L$, and $-\theta\in \Lambda^1 M$ the corresponding connection
form, $\nabla(\psi)=-\theta\otimes \psi$. Then $d\omega_\psi=-\theta\otimes\omega_\psi$,
where $\theta$ is a closed 1-form,  
and $\omega_\psi\in \Lambda^2M$ is $\omega$ considered to be  a differential
form after the identification $L\simeq  C^\infty(M)$ provided by $\psi$.

\hfill

\proof After identifying $L$ with a trivial
line bundle, we obtain $0=d_\nabla(\omega)= d(\omega) -\omega \wedge \theta$.
\endproof

\hfill

Let $\psi$ be a non-degenerate section of $L$, and
$\alpha\in \Lambda^k(M,L)$ an $L$-valued $k$-form. 
The multiplication by $\psi$ induces an isomorphism
between $C^\infty M$ and $L$. Then, 
for any section $u$ of $L$, the quotient $\frac  u \psi$
is considered to be  a function. Now, $\frac \alpha\psi$
can be naturally considered to be  a section of $\Lambda^k(M)$.
If we apply this construction to the $L$-valued 1-form $\nabla\psi$,
we obtain the connection form $-\theta$ in $L$.
In this situation the connection can be expressed 
as $\nabla(f\psi) = df\otimes \psi - f \theta\otimes \psi$.

\hfill

\theorem
1. Let $(M, \omega, L)$ be an LCS manifold, and
$\psi$ a non-degenerate section of $L$.
Then the 2-form $\omega_\psi:= \frac\omega\psi\in \Lambda^2(M)$
satisfies the equation 
$d(\omega_\psi)=\omega_\psi\wedge\theta$, 
where $\theta:= -\frac {\nabla\psi}{\psi}\in \Lambda^1 (M)$.

2. Conversely, for any non-degenerate
2-form $\omega_0\in \Lambda^2(M)$ satisfying
$d\omega_0= \theta \wedge \omega_0$, and
a real line bundle with non-degenerate section $\psi$ and the connection
$\nabla(f\psi) = df\otimes \psi - f\theta\otimes \psi$,
the $L$-valued form $\omega_0\otimes \psi$ satisfies
$d_\nabla(\omega_0\otimes \psi)=0$. In these assumptions,
$\theta$ is closed if and only if $(L, \nabla)$ is flat.

\hfill

\proof
The first part is clear: 
\[ 
0=d_\nabla(\omega_\psi\otimes \psi)=d (\omega_\psi)
\otimes \psi + \omega_\psi \otimes d_\nabla\psi = 
(d \omega_\psi- \omega_\psi \wedge \theta) \otimes \psi,
\]
implying $d \omega_\psi= \omega_\psi \wedge \theta$.

For the converse, we write 
\[
d_\nabla(\omega_0\otimes \psi)= d\omega_0\otimes \psi+
\omega_0 \wedge \nabla\psi= (d\omega_0 -\omega_0\wedge
\theta)\otimes \psi=0.
\]
The curvature of $L$ is $-d\theta + \theta\wedge \theta$,
and it vanishes if and only if $\theta$ is closed.
\endproof

\hfill

\remark
We have just shown that
the following two definitions are equivalent:

1. An LCS manifold\index[terms]{manifold!LCS} is one
equipped with a non-degenerate 
2-form $\omega$ satisfying $d\omega=\theta\wedge\omega$,
where $\theta$ is a closed 1-form, 
called {\bf the Lee form}.\index[terms]{form!Lee}

2. An LCS manifold is one equipped with a non-degenerate,
\index[terms]{manifold!LCS}
closed 2-form $\omega$ taking values in a flat, oriented
line bundle.\index[terms]{bundle!line!flat}

\hfill

\remark
Locally, any closed form $\theta$ is exact: $\theta\restrict{U}= df$ on some open set $U$.
Then, on this open set, $e^{-f} \omega\restrict{U}$ is closed; indeed,
$d(e^{-f} \omega) = e^{-f} d\omega - df\wedge e^{-f}\omega=
e^{-f}(\theta \wedge \omega-\theta\wedge\omega)=0$ on $U$.
Any LCS form is locally conformally equivalent
to a symplectic form; this explains the term.

\hfill

\remark 
Let $\omega$ be a non-degenerate 2-form on a manifold. 
An easy calculation can be used to 
show that   $d\theta=0$ follows from $d\omega=\theta\wedge\omega$
	when $\dim_\R M>4$. Indeed, $0=d(d\omega)= \theta
        \wedge \theta\wedge\omega + d\theta \wedge
        \omega$, but $d\theta \wedge
        \omega=0$ implies $d\theta=0$ when  $\dim_\R M>4$
(Exercise \ref{_mult_by_symple_injective_Exercise_}).

\hfill

\remark As $d\theta=0$, the operator
$d_\theta:=d-\theta\wedge$ satisfies $d_\theta^2=0$, and
hence determines a ``twisted'' cohomology, also called
Morse--Novikov cohomology.
\index[terms]{cohomology!twisted}\index[terms]{cohomology!Morse--Novikov}

\section{Galois covers and the deck transform group}
\label{carac}

\definition\label{_Galois_cover_definition_}
Let $M$ be a manifold, and $\pi:\; \tilde M \arrow M$
a covering. {\bf An automorphism of the covering $\pi$}
is a diffeomorphism $\phi:\; \tilde M \arrow \tilde M$
commuting with the projection to $M$. The group of
automorphisms of the covering $\tilde M$ is denoted
$\Aut_M(\tilde M)$.
A connected covering $\pi:\; \tilde M \arrow M$ is
called {\bf a Galois cover} if its group of
automorphisms acts transitively on the fibres
of $\pi$. In this case, the  group
$\Aut_M(\tilde M)$ is called {\bf the deck transform group},
{\bf the monodromy group}, or {\bf the Galois group}
of the covering.\index[terms]{group!Galois}\index[terms]{group!monodromy}\index[terms]{group!deck transform}

\hfill

The main theorem of Galois theory\index[terms]{Galois theory!for coverings} of the coverings is
very similar to the results on Galois groups
of finite field extensions. In fact, there
is a treatment of both subjects using the
language of Galois categories (due to \index[persons]{Grothendieck, A.} Grothendieck),
where all results are deduced from the fundamental
principles at the same time (\cite{_Cadoret:Galois_}).

\hfill

\theorem \label{_Galois_covers_Theorem_}
(The main theorem of Galois theory of the coverings)\\
Let $M$ be a connected manifold.\index[terms]{Galois theory!for coverings}
Given a connected covering $\tilde M \arrow M$,
we obtain a subgroup $\pi_1(\tilde M) \subset \pi_1(M)$.
This defines a bijection between
the isomorphism classes of connected coverings
of $M$  and the subgroups $G\subset \pi_1(M)$.
Under this correspondence, the Galois
covers correspond to the normal subgroups\index[terms]{cover!Galois}
$G\subset \pi_1(M)$.

\proof \cite[Section 1.3]{_Hatcher_}. \endproof

\hfill

\remark
Under this correspondence, the universal cover  of
$M$ corresponds to the trivial subgroup
$\pi_1(\tilde M) =\{e\}$. It is obviously a Galois cover.

\hfill

\remark
Let $\tilde M \arrow M$ be a Galois cover. Then the
deck transform group of the covering $\tilde M$  
is identified with 
the quotient $\frac{\pi_1(M)}{\pi_1(\tilde M)}$.

\hfill

\definition\label{_minimal_cove_Definition_}
Let $L$ be a bundle over $M$ with a 
flat connection $\nabla$, and $x\in M$ a point. 
Its monodromy group (monodromy,\index[terms]{group!monodromy}
in this context, means the same as holonomy,
but applied to flat bundles)\index[terms]{bundle!vector bundle!flat}
defines a map $\chi:\; \pi_1(M,x) \arrow L\restrict x$.\index[terms]{holonomy!of a flat bundle}\index[terms]{connection!flat}
{\bf The minimal cover\index[terms]{cover!minimal} of $M$ associated with 
$(L, \nabla)$} is the connected covering $\tilde M$
that satisfies $\pi_1(\tilde M)= \ker \chi$.

\hfill

\remark
Let $\tilde L$ be the pullback of $L$ to the
corresponding minimal cover $\tilde M$. By construction,
the monodromy\index[terms]{monodromy} of $\tilde L$ is trivial, and hence  $\tilde L$
is naturally trivialized by its parallel sections.\index[terms]{section!parallel}
In other words, $\tilde M$ is the smallest cover\index[terms]{cover!minimal}
where the pullback of $L$ is trivialized by the
parallel sections.

\section{Locally conformally \ka\ manifolds}

In this section we give many equivalent definitions
of an LCK manifold, and prove their equivalence.

\subsection{LCK manifolds: the tensorial definition}

\definition\label{_LCK_via_form_Definition_}
Let $(M,I, \omega)$ be a Hermitian manifold, $\dim_\C M >1$.
The manifold $M$ is called {\bf  locally conformally K\"ahler}
(LCK) if $d\omega=\theta\wedge\omega$, where $\theta$ is a closed
1-form, called, as in the LCS case, {\bf  the Lee form}.\index[terms]{form!Lee}
\index[terms]{manifold!LCK}

\hfill

Clearly, every LCK manifold is locally conformally symplectic  as well.

%

\hfill

\remark
The LCK property is conformally invariant: if $g$ is LCK
with Lee form $\theta$, then $e^fg$ is LCK with the
 Lee form\index[terms]{form!Lee}
$\theta+df$. If $M$ is compact and LCK, one can choose to
work with the \index[persons]{Gauduchon, P.} 
Gauduchon metric (unique up to constant multiplier in
any given conformal class, see Section \ref{gauduchon_metric}). 
The form $\theta$ is harmonic if and only if
$\omega$ is Gauduchon (\ref{_Gauduchon_theta_coclosed_Claim_}). 
\index[terms]{metric!Gauduchon}

\hfill

\remark
Usually, one tacitly assumes that $\theta$ is not exact.
Indeed, if $\theta=d\phi$, then $d(e^{-\phi}\omega)=
e^{-\phi}d\omega - e^{-\phi}\theta\wedge\omega=0$,
and $e^{-\phi}\omega$ is K\"ahler. In this case
$(M,I, \omega)$ is called {\bf  globally conformally
K\"ahler}, sometimes abbreviated as GCK.

\hfill

\remark 
As pullback commutes with exterior differential, any
complex submanifold $\iota:N\hookrightarrow M$ of an LCK
manifold inherits an induced LCK structure\index[terms]{submanifold}\index[terms]{structure!LCK}
$(\iota^*\omega, \iota^*\theta)$.

\subsection{The weight bundle and the homothety character}
\index[terms]{form!automorphic}
\index[terms]{homothety character}
\label{_weight_automo_Subsection_}

\remark\label{_Kahler_on_cover_Remark_}
Let $(M, \omega, \theta)$ be an LCK manifold, and
$\pi:\; \tilde M\arrow M$ a co\-ve\-ring such that the pullback
$\pi^* \theta$ is exact, $\pi^* \theta= d \phi$.
Then the form $e^{-\phi} \pi^*\omega$ is closed,
\begin{align*}
d(e^{-\phi} \pi^*\omega) & = - e^{-\phi}d\phi \wedge \pi^*\omega+
 e^{-\phi} \pi^*\theta \wedge \pi^*\omega= 
\\ &= - e^{-\phi} \pi^*\theta \wedge \pi^*\omega+ 
e^{-\phi} \pi^*\theta \wedge \pi^*\omega=0.
\end{align*}
Therefore, $M$ admits a covering with a K\"ahler
form that is  conformally equivalent to $\pi^*\omega$.
The converse is also true, as shown below (Subsection 
\ref{_Kahler_cover_Subsection_})

\hfill

\definition\label{_homothety_character_Definition_}
Let $M$ be an LCK manifold and $(\tilde M, \tilde \omega)$
its K\"ahler cover.\index[terms]{cover!K\"ahler} Denote by $\Gamma$ the deck transform group
of $\tilde M$. For each $\gamma\in \Gamma$, we have
$\gamma^*(\tilde \omega)= \chi(\gamma) \tilde \omega$
where $\chi:\; \Gamma \arrow \R^{>0}$ is a character,
called the {\bf homothety character}.\index[terms]{homothety character}
Since $\Gamma$ is a quotient group of $\pi_1(M)$,
we consider $\chi$ to be a character on $\pi_1(M)$.
This is the character associated with the K\"ahler
structure on the universal cover of $M$.

\hfill

\claim 
Let $(M, \omega)$ be an LCK manifold, and
$\chi:\; \pi_1(M) \arrow \R^{>0}$ its
homothety character. Then $\chi$ is 
independent on  the choice of the
K\"ahler cover $(\tilde M, \tilde \omega)$.

\hfill

\proof
Let $(\tilde M_1, \tilde \omega_1)$ and
$(\tilde M_2, \tilde \omega_2)$ be two K\"ahler 
covers, and $\chi_1, \chi_2$ the corresponding
homothety characters\index[terms]{homothety character}. 
Passing to the universal covering, we may
assume that $\tilde M_1= \tilde M_2$ is the universal
covering $\tilde M$ of $M$. Indeed, for any $\gamma\in \pi_1(M)$,
the action of $\gamma$ scales the form $\tilde \omega_i$
by the same factor on $\tilde M_i$ as on $\tilde M$.
Since $\tilde \omega_1$ and $\tilde \omega_2$ 
are conformally equivalent to the pullback of $\omega$,
they  are conformally equivalent: $\tilde \omega_2= f \tilde \omega_1$.
This implies $0 = d(\tilde \omega_2)=df  \wedge \tilde \omega_1$,
giving $df=0$, see Exercise
\ref{_mult_by_symple_injective_Exercise_}. 
Since the forms $\tilde \omega_1$,
$\tilde \omega_2$ are proportional, the corresponding
scale factors are equal, which implies $\chi_1=\chi_2$.
\endproof

\hfill

\definition\label{_trivializing_on_cover_Definition_}
Let ${\Bbb L}$ be a local system on a manifold $M$,
and $\hat M$ a covering of $M$. We say that
{\bf the local system ${\Bbb L}$ is trivial
on $\hat M$} if its pullback to $\hat M$
is a trivial local system.

\hfill

\remark 
Let $(M, x)$ be a connected manifold with a marked point $x\in M$.
Recall that the category of local systems on $M$ 
is equivalent to the category of representations of $\pi_1(M)$
(\ref{_local_system_repres_Remark_}).
According to 
\ref{_Galois_covers_Theorem_}, 
 the coverings of $M$ are in bijective correspondence
with the subgroups $\Gamma \subset \pi_1(M)$.
The covering $M_\Gamma$, corresponding
to $\Gamma$, is obtained as $M_\Gamma:= \frac{\tilde M}{\Gamma}$,
with $\Gamma= \pi_1(M_\Gamma)$. Consider a local system ${\Bbb L}$
associated with a representation $\rho:\; \pi_1(M) \arrow \GL(n, \R)$.
Then ${\Bbb L}$ is trivial on the covering $M_\Gamma$ if and only if
$\Gamma \subset \ker \rho$. Therefore, the minimal cover\index[terms]{cover!minimal}
trivializing the local system ${\Bbb L}$
(\ref{_minimal_cove_Definition_}) is the
covering $M_\Gamma$, associated with $\Gamma= \ker \rho$.

\hfill

\definition\label{_weight_bun_min_Kahler_Definition_}
Let $M$ be an LCK manifold, and
$\chi:\; \pi_1(M) \arrow \R^{>0}$ its homothety
character. Consider the normal subgroup $G:= \ker \chi\subset \pi_1(M)$.
{\bf The minimal K\"ahler cover}\index[terms]{cover!minimal K\"ahler} is the Galois cover associated
with this subgroup. By construction, this is the minimal
cover that trivializes {\bf the weight bundle} $L$
(\ref{_weight_for_lcs_}). In this context, the weight\index[terms]{bundle!weight}
bundle is the flat bundle that is  associated with the\index[terms]{bundle!line!flat}
local system given by the character $\chi:\; \pi_1(M) \arrow \R^{>0}$
via the Riemann--Hilbert correspondence (\ref{_Riemann--Hilbert_Theorem_}).
\index[terms]{correspondence!Riemann--Hilbert}

\hfill

\definition\label{_weight_bundle_for_kahler_cover_Definition_}
Let $M$ be a complex manifold such that
its universal covering $\tilde M$ is equipped with a K\"ahler
form $\tilde \omega$. Suppose that the deck transform group acts on 
$(\tilde M, \tilde \omega)$ by K\"ahler homotheties.
Denote by $\chi$ the corresponding character
$\chi:\; \pi_1(M) \arrow \R^{>0}$.
Then the local system associated with this
character is called {\bf the weight bundle}.

\subsection{Automorphic forms related to the homothety character}

\definition\label{_automo_forms_Definition_}
Let $M$ be an LCK manifold and $(\tilde M, \tilde \omega)$
its K\"ahler cover. Denote by $\Gamma$ the deck transform group
of $\tilde M$. For each $\gamma\in \Gamma$, we have
$\gamma^*(\tilde \omega)= \chi(\gamma) \tilde \omega$,
where $\chi:\; \Gamma \arrow \R^{>0}$ is the homothety character.\index[terms]{homothety character}
The {\bf automorphic differential forms}
(or functions, that we understand as sections of
$\Lambda^0(\tilde M)$) on $\tilde M$ are differential forms
$\alpha \in \Lambda^*(\tilde M)$
that have the same automorphic properties
as $\tilde \omega$, namely, $\gamma^*\alpha=\chi(\gamma)\alpha$
for all $\gamma\in \Gamma$.

\hfill

\claim\label{_automo_sections_of_L_Claim_}
The space of automorphic $k$-forms on $\tilde M$
is naturally identified with the space of
sections of $\Lambda^k(M) \otimes L$.

\hfill

{\bf Proof:} Consider the pullback $\tilde L:= \pi^* L$,
where $\pi:\; \tilde M \arrow M$ is the covering map.
Since $\Gamma$ is the deck transform group of $\tilde M$,
we have $\Gamma=\pi_1(M)/\pi_1(\tilde M)$, and
the character $\chi$ vanishes on $\pi_1(\tilde M)\subset\pi_1(M)$.
Therefore, the bundle $\tilde L$ is flat and has trivial monodromy\index[terms]{monodromy} on $\tilde M$,
hence it is naturally trivialized by the flat connection.\footnote{For
the sake of clarity, we are deliberately inaccurate here. 
A flat line bundle with trivial monodromy is trivialized by 
a flat section. However, this section is unique only up to a constant
multiplier. This choice of a constant actually affects the
identification of the space of sections of $\Lambda^k(M) \otimes L$
with the space of automorphic forms. We ignore this ambiguity,
to avoid more confusion.}
This is used to interpret sections of $\pi^* L$ as functions
on $\tilde M$, and $\pi^*L$-valued differential forms
as differential forms on $\tilde M$. 

For any section $u$ of $L$, let $\tilde u$ be
its lift to $\pi^* L$. By construction, the bundle
$\pi^* L$ is $\Gamma$-equivariant, and $\Gamma$ 
acts on the pullbacks of sections of $L$ multiplying 
a section with the homothety character.\index[terms]{homothety character}
Therefore, $\tilde u$ is automorphic.

Conversely, let $\tilde f\in C^\infty (\tilde M)$ be an automorphic function,
and $\tilde u:= \pi^* u$, where $u$ is a nowhere degenerate section of $L$.
Then $\tilde f_1:= \tilde f\tilde u^{-1}$ is a $\Gamma$-invariant function on $\tilde M$,
that is,  $\tilde f_1 = \pi^* f_1$, where $f_1$ is 
a function on $M$. Therefore, $\tilde f= \tilde u \tilde f_1 = \pi^*(u f_1)$.

The same argument applies to $L$-valued differential forms as well.
\endproof

\hfill

\example
Let $(M, \omega)$ be an LCK manifold, and $\tilde \omega$
the K\"ahler form on its K\"ahler cover.
Clearly, the  K\"ahler form on $\tilde M$ is automorphic.

\hfill

\definition 
 Let $\chi$ be the homothety  character 
of $L$. An {\bf automorphic form of weight 
$\lambda$}\index[terms]{form!automorphic!of weight $\lambda$}
is an automorphic form that satisfies
$\gamma^*\tilde \eta=\chi(\gamma)^{\lambda}\eta$
for each $\gamma\in \pi_1(M)$.

\hfill

Hence, an automorphic form of weight $\lambda$ (on
$\tilde M$) can be viewed as a form on $M$ with values
in the tensor power $L^\la$. Note that $L^\la$ is endowed
with the connection form $\la\theta$ and is still
flat because $d\theta=0$.

\subsection{K\"ahler covers of LCK manifolds: the second
  definition}
\label{_Kahler_cover_Subsection_}

Let $M$ be a complex manifold.  Assume that
its universal covering $\tilde M$ is equipped with a K\"ahler
form $\tilde \omega$, and the deck transform group acts on 
$(\tilde M, \tilde \omega)$ by K\"ahler homotheties.
 Consider the local system $L$ on $M$ associated with
the homothety character\index[terms]{homothety character} 
$\chi:\; \pi_1(M)\arrow \R^{>0}$ (\ref{_homothety_character_Definition_})
as in \ref{_weight_bundle_for_kahler_cover_Definition_},
and let $\psi_0$ be a smooth, positive
section that trivializes $L$, considered to be  an
oriented line bundle (\ref{_Riemann--Hilbert_Theorem_}). 
The lift of $\psi_0$ to $\tilde M$ 
is a section of the local system $\pi^*L$, that is  trivial 
on $\tilde M$ by our assumption.
The bundle $\pi^*L$ is obtained from a constant sheaf
by tensor multiplication with $C^\infty \tilde M$;
let $u$ be the constant section of this constant sheaf.
We consider $u$ to be a  section of $\pi^* L$. 
Since all sections of $\pi^* L$ are
proportional, there exists a function $f$ on $M$
such that $\psi_0= e^f u$. For each
$\gamma\in \Aut_M(\tilde M)$, 
one has $1 = \frac{\gamma^*\psi_0}{\psi_0}=
\frac{\gamma^*(e^f u)}{e^f u}=
\frac{\gamma^*u}{u} \cdot \frac{\gamma^* (e^f)}{e^f}$.
Then $\chi(\gamma)=\frac{\gamma^*u}{u} = \frac{\gamma^* (e^{-f})}{e^{-f}}$.
By definition of the local system $L$, 
$\chi(\gamma)=\frac{\gamma^*u}{u}= \frac{\gamma^*(\tilde
  \omega)}{\tilde\omega}$. Therefore,
{\bf  $e^{-f}\tilde\omega$ is a
$\pi_1(M)$-invariant form on $\tilde M$}.
Denote by $\omega$ the corresponding form on $M$.
Then 
\begin{equation}\label{_LCK_from_sections_of_L_Equation_} 
d\omega= d(e^{-f}\tilde \omega)= -e^{-f} df\wedge
\tilde\omega=-df\wedge\omega.
\end{equation}
This shows that the form $\omega$ satisfies 
$d\omega=\theta\wedge\omega$, where $\theta=-df$.

\hfill

This argument is used to state the following theorem.

\hfill

\theorem \label{_LCK_in_terms_of_covers_Theorem_}
Let $M$ be a complex manifold such that
its universal covering $\tilde M$ is equipped with a K\"ahler
form $\tilde \omega$, and the deck transform group acts on 
$(\tilde M, \tilde \omega)$ by K\"ahler homotheties.
Consider a section $\psi_0$ trivializing the weight
bundle $L$ associated with $(\tilde M, \tilde \omega)$
as above, and let $\psi$ be the function
on $\tilde M$ such that $\gamma^* \psi = \chi(\gamma)\psi$.\footnote{In the notation introduced above, $\psi=e^f$.} 
Then the corresponding Hermitian
form $\psi^{-1}\tilde\omega$ on $M$ is LCK. Conversely, any
LCK form can be obtained in this way from a K\"ahler
form on a covering and a trivialization
of the corresponding weight bundle.

\hfill

\proof
The form $\psi^{-1}\tilde\omega$ is LCK by 
\eqref{_LCK_from_sections_of_L_Equation_}.
Conversely, starting from an LCK structure on
$M$, we obtain an automorphic K\"ahler 
form on $\tilde M$ by 
\ref{_Kahler_on_cover_Remark_}.
\endproof

\hfill

We have arrived at an alternative definition of an LCK
manifold; its equivalence to the first one is
\ref{_LCK_in_terms_of_covers_Theorem_}.

\hfill

\definition \label{_LCK_in_terms_of_covers_Definition_}
  An LCK manifold is a complex manifold such that
its universal covering $\tilde M$ is equipped with a K\"ahler
form $\tilde \omega$, and the deck transform group acts on 
$(\tilde M, \tilde \omega)$ by K\"ahler homotheties.

\hfill

\remark\label{_Weyl_conn_definition_Remark_}
Note that the Levi--Civita connection on a K\"ahler
manifold does not change if we multiply the K\"ahler
metric by a constant. Therefore, the deck transform maps,
acting by homotheties on $(\tilde M, \tilde\omega)$,
preserve the Levi--Civita connection, that can
 be obtained as the pullback of a connection
$D:\; TM \arrow TM \otimes \Lambda^1M$ on $M$.
The connection $D$ is called {\bf the  Weyl connection}\index[terms]{connection!Weyl}
on the LCK manifold (see Section \ref{_Weyl_conne_Appendix_}). Its holonomy group\index[terms]{group!holonomy} contains
multiplication by a constant, and hence  the Weyl connection  cannot
preserve any global Riemannian structure on $M$.
However, $D$ preserves the conformal structure.

\subsection{LCK manifolds via an $L$-valued
  K\"ahler form: the third definition}
\label{_L_valued_Kahler_Subsection_}

Let $(L,\nabla)$ be a real line bundle with a flat connection on a complex manifold 
$M$, and $\tilde M$ a covering of $M$ such that the pullback of $L$ is trivial on 
$\tilde M$ as a local system (\ref{_trivializing_on_cover_Definition_}).
Earlier, we  indentified $L$-valued differential 
forms with automorphic forms on $\tilde M$. 
\ref{_LCK_in_terms_of_covers_Theorem_} is used to
interpret LCK forms on $M$ as automorphic K\"ahler
forms on $\tilde M$.

Recall that $L$-valued differential forms are equipped
with a differential $d_\nabla:\; \Lambda^k(M, L) \arrow  \Lambda^{k+1}(M, L)$
(\ref{dnabla}). After identifying, as above, the sections of $\Lambda^k(M, L)$ with 
automorphic forms, this operator corresponds to the de Rham differential.

This leads to the third equivalent definition of LCK manifolds 
(\ref{_LCK_via_L_valued_Definition_}
below). We start by defining ``$L$-valued K\"ahler forms''.

\hfill

\definition
Let $(L, \nabla)$ be a flat, oriented 
real line bundle on a complex manifold $M$,
and $\omega\in \Lambda^{1,1}(M, L)$
a differential form. We say that
$\omega$ is {\bf an $L$-valued K\"ahler form}
if $d_\nabla \omega=0$ and, for any
$v\in TM$, the expression $\omega(v, Iv)$,
considered to be  a section of $L$, is non-negative,
and strictly positive in all points 
of $M$ where $v\neq 0$.

\hfill

Let $\phi, \psi$ be sections of the line bundle $L$, 
with $\psi$ non-degenerate.
The ratio $\frac \phi \psi$ of two sections of $L$
is a function. Given an $L$-valued differential
form $\eta \in \Lambda^k(M,L)$, the ratio
$\frac \eta \psi$ is considered to be  a differential
form, $\frac \eta\psi \in \Lambda^kM$.

\hfill

\theorem
Let $(L, \nabla)$ be a be a flat, oriented
real line bundle on a complex manifold $M$,
and $\psi\in L$ a smooth section that trivializes $L$. 
Then for any $L$-valued K\"ahler form
$\omega_L$, the form $\omega:=\frac {\omega_L}\psi \in \Lambda^{1,1}(M)$
is LCK. Conversely, for any LCK manifold
$(M, \theta, \omega)$, and a flat line bundle
$(L, \nabla)$ trivialized by a section $\psi$ 
such that the corresponding connection form is $-\theta$,
the form $\omega\otimes\psi\in \Lambda^{1,1}(M,L)$
is an $L$-valued K\"ahler form.

\hfill

\proof 
Let $\theta$ be a 1-form such that 
$d_\nabla \psi = - \psi \otimes \theta$.
We consider $\psi^{-1}$ as a section of the
dual bundle $L^*$; then $d_\nabla(\psi^{-1})= \psi^{-1}\otimes \theta$.
Given an $L$-valued K\"ahler form $\omega_L$, we 
obtain 
\[ 
 d \omega=d(\omega_L\otimes\psi^{-1})= \omega \wedge d_\nabla(\psi^{-1})=
\omega_L \wedge \theta \otimes \psi^{-1} = \omega\wedge \theta.
\]
Since $\nabla$ is flat, the 1-form $\theta$ is closed:
\[ d^2_\nabla(\psi)= d_\nabla(\psi\otimes \theta)= \psi \otimes d\theta.
\]
This implies that $(M, \omega, \theta)$ is an LCK structure.

Conversely, for any LCK structure $(M, \omega, \theta)$,
the form $\omega \otimes \psi$ is $d_\nabla$-closed:
\[
d_\nabla(\omega \otimes \psi)= d\omega\otimes \psi +
\omega\wedge d_\nabla(\psi) = \omega\wedge\theta\otimes \psi 
-\omega\wedge\theta\otimes \psi=0.\ \ \ \endproof
\]

\hfill

Now we can state an equivalent definition of LCK structures.

\hfill

\definition\label{_LCK_via_L_valued_Definition_}
An LCK manifold is a complex manifold equipped
with a flat line bundle $L$ and an $L$-valued K\"ahler form.

\subsection{Conformally equivalent K\"ahler forms}

Locally conformally K\"ahler manifolds are defined as Hermitian manifolds
$(M,\omega)$ that satisfy $d\omega =\omega\wedge \theta$, where
$d\theta=0$. However, the assumption $d\theta=0$
is automatic when $\dim_\C M>2$, as explained below.

Let $\omega$ be a non-degenerate 2-form on a 
smooth manifold $M$ of dimension $2n$.
Denote by $L_\omega$ the operator
of exterior multiplication by $\omega$, with 
$L_\omega(\alpha):= \omega\wedge\alpha$.
A simple linear-algebraic calculation
implies that $L_\omega:\; \Lambda^k(M)\arrow \Lambda^{k+2}(M)$
is injective for all $k\leq n-1$
(Exercise \ref{_mult_by_symple_injective_Exercise_}). 
This algebraic statement  is well known, because it is used in the proof of 
the Lefschetz theorem in Hodge theory.\index[terms]{theorem!Lefschetz}\index[terms]{Hodge!theory}

For us, it is important, because 
$d^2\omega=0$ implies
$d(\omega\wedge \theta)=  \omega\wedge \theta\wedge\theta 
+ \omega\wedge d\theta= \omega\wedge d\theta= 0$,
and in $\dim_\C M>2$ this gives $d\theta=0$.

This argument also implies that
any conformal map of connected symplectic manifolds
of dimension $>4$ is in fact a symplectic homothety.
Indeed, let $\psi^*(\omega) = f\omega$,
then 
\begin{equation}\label{_conformal_then_homothety_Equation_}
0=d(\psi^*(\omega))= d(f\omega) = df\wedge \omega,
\end{equation}
implies $df=0$, and hence  $f$ is a constant.
Therefore, any holomorphic conformal map between
K\"ahler manifolds in complex dimension $>2$ is
homothetic.

\subsection[LCK manifolds via charts and atlases: the fourth definition]{LCK manifolds via charts and atlases:\\ the fourth definition}

Instead of the tensorial definition, we can give a 
definition in terms of charts and atlases. 
Let $(M, \omega, \theta)$ be an LCK manifold.
Locally, the closed form $\theta$ is exact:
$\theta = d f$. Then $e^{-f}\omega$ is closed,
\begin{align*}
d(e^{-f} \omega) & = - e^{-f}df \wedge \omega+
 e^{-f} \theta \wedge \omega= 
\\ &= - e^{-f} \theta \wedge \omega+ 
e^{-f} \theta \wedge \omega=0.
\end{align*}
Hence, (up to a conformal factor) the manifold
$(M, \omega)$ is locally K\"ahler. This explains
the term ``Locally Conformally K\"ahler''. 

\hfill

\theorem
Let $M$ be a complex manifold
covered by charts $\{U_\alpha\}$ that are 
all K\"ahler, and the transition functions between
the charts are homotheties.\footnote{By 
\eqref{_conformal_then_homothety_Equation_}, for $\dim_\C M >2$, any
conformal equivalence between K\"ahler forms
is a homothety.} The transition functions
are homothety coefficients, and hence  they are
constant. They satisfy the cocycle condition;
let $L$ be the local system defined by these
transition functions.
Then the local K\"ahler forms
can be interpreted as a single K\"ahler form with
coefficients in $L$. Conversely, for any
$L$-valued K\"ahler form on $M$, an atlas
$\{U_\alpha\}$ and the set of trivializations
of $L$ over each $U_\alpha$,
the transition functions between the 
corresponding collection of K\"ahler
forms are homotheties.
\endproof

\hfill

This leads us to the fourth definition
of an LCK manifolds, equivalent 
to the rest.

\hfill

\definition
Let $M$ be a complex 
manifold covered by charts $\{U_\alpha\}$.
Suppose that each $U_\alpha$ is equipped with
a K\"ahler form $\omega_\alpha$, and the transition functions
between the charts are homotheties with respect to $\omega_\alpha$. Then
$(M, \{U_\alpha, \omega_\alpha\})$ is called 
an LCK manifold.
{\bf An LCK form} on $(M, \{U_\alpha, \omega_\alpha\})$
is a Hermitian form\index[terms]{form!Hermitian} that is  conformally equivalent
with each $\omega_\alpha$. 

%

\section{The LCK rank}\index[terms]{manifold!LCK}

Using the homothety 
character $\chi$, we can introduce another useful notion:

\hfill

\definition \label{lck_rank} 
Let $M$ be an LCK manifold, and $\chi$ its homothety
character. The rank of the image
$\im(\chi)\subset\R^{>0}$ is called the {\bf LCK rank} of
the LCK manifold $M$.\index[terms]{rank!LCK}

\hfill

An equivalent definition is provided in the following:

\hfill

\proposition (\cite{pv})
The LCK rank equals the rank of the monodromy group of the flat connection $\theta$ in the weight bundle $L\arrow M$.\index[terms]{bundle!weight}\index[terms]{group!monodromy}

\hfill

\proof The logarithm\index[terms]{logarithm!of a character} of the character is the morphism 
\[\log\chi:\pi_1(M)\arrow (\R,+).
\]
 As such, it can be identified with an element of $$\Hom(\pi_1(M)^{\text{ab}},\R)\simeq\Hom(H_1(M),\R)\simeq H^1(M,\R).$$ 
Under these isomorphisms, $[\theta]$ corresponds to $\log\chi$.

Now, the monodromy group of the connection $\theta$ is generated by all the integrals $\int_\gamma\theta$, where $\gamma$ are loops generating $H_1(M)$. But 
$\int_\gamma\theta=\log\chi(\gamma),$
and this completes the proof.
\endproof

\section{A first example}

\noindent {\bf The classical Hopf manifold} \index[terms]{manifold!Hopf!classical} 
is $H:=\C^n \backslash 0/ \Z$, where $\Z$ acts diagonally, generated by the holomorphic transformation $(z_i)\mapsto (\la\cdot z_i)$ 
for a fixed number $\lambda\in \C$, $|\lambda|>1$. 

Clearly,  the universal covering $\C^n \backslash 0$ has the flat \ka\ metric  $\frac  1 2 \sum(dz^i \otimes d \bar{z}^i  +d \bar{z}^i \otimes dz^i)$ with respect to which 
the deck group $\Z$ acts by homotheties.

The \ka\ form then reads $\frac{\1}{2}\6\bar\6|z|^2$, the LCK
form $\omega$ is $\frac{\1}{2|z|^2}\6\bar\6|z|^2$,  and the
Lee form \index[terms]{form!Lee}is $\theta =-d\log |z|^2$ (I. \index[persons]{Vaisman, I.} Vaisman, \cite{va_isr}).

It is easy to see that $H$ is diffeomorphic to $S^1 \times S^{2n-1}$,
and fibred over $\C P^{n-1}$ with fibre the torus $\C^*/\langle \lambda \rangle$.

\hfill

\remark 
For any complex submanifold $\iota:X\hookrightarrow \C P^{n-1}$,
its preimage in the Hopf manifold $H$ is a complex manifold.\index[terms]{submanifold}  
The induced Hopf bundle $M:=\iota^*H\rightarrow X$ is again an LCK manifold. In particular, if $X$ is a genus $g$ compact curve, then $b_1(M)=2g+1$, $b_2(M)=4g$, yielding compact LCK surfaces with arbitrary large $b_1$, \cite{va_gd}
(Chapter \ref{comp_surf}).

\section{Notes}\label{_first_notes_}

\begin{enumerate}

\item 
The Lee form\index[terms]{form!Lee} is named after Hwa-Chung Lee,\index[persons]{Lee, H.-C.} who introduced it in \cite{lee}, where he studied nondegenerate differential 2-forms. It appeared in the context of finding an integrating factor for such a 2-form $\omega$ (a scalar $f$ such that $f\omega$ is exact). The precise result proven by Lee is  (\cite[Theorem 6]{lee}): 

\smallskip

\begin{minipage}{10.8cm}
{\em A necessary and sufficient condition that a nonsingular exterior differential form of degree two in an even number ($>4$) of variables admit an integrating factor is that the form itself be a factor of its exterior derivative.} 
\end{minipage}

\smallskip

 That is precisely the equation $d\omega=\theta\wedge\omega$.

	\item LCS manifolds were introduced, with this
          name, in \cite{lib} (1954).\index[terms]{manifold!LCS}
          Their interest was renewed starting with
          \index[persons]{Vaisman, I.} Vaisman's paper
          \cite{va_lcs} (which he wrote some 10 years
          after he initiated the study of LCK manifolds,
          see below, Note no 5). See also the  papers by
          A. \index[persons]{Banyaga, A.} Banyaga
          (\cite{ban1}, \cite{ban2}),
          S. \index[persons]{Haller, S.} Haller
          (\cite{hal}), S. Haller \&
          T. \index[persons]{Rybicki, T.} Rybicki
          (\cite{hr1}, \cite{hr2}) etc. Recently, there
          was a growing interest in LCS geometry\index[terms]{geometry!LCS}, e.g.
          \cite{_Bertelson_Meigniez_}, \cite{bk1},
          \cite{bk2}, \cite{oti1}, \cite{ots}, \cite{cm}
          etc. \index[persons]{Westlake, W. J.} As for LCK
          manifolds, it seems that their first appearance
          can be implicitly  found in J. W. Westlake's
          paper  \cite{westlake} (1954).

\item There exist compact symplectic manifolds \index[terms]{manifold!symplectic} that cannot carry any K\"ahler metric (see \cite{_Thurston:Kodaira_, _Babenko_Taimanov_}). The obstructions are topological.
Similarly, one would like to have examples of compact LCS
manifolds carrying no LCK metric. The difficulty is that,
as we shall see further, there are no known topological
obstructions for the existence of LCK structures\index[terms]{structure!LCK}. In \cite{bm} an example of a compact LCS\index[terms]{manifold!LCS} nilmanifold\index[terms]{nilmanifold} with
no LCK metric is constructed in dimension $4$, using the
classification of compact complex surfaces. Blowing-up points on these two examples, \cite{cy}, one
can obtain new examples of compact LCS manifolds without
LCK structures, but still only in dimension 4.

On the other hand, in \cite[Theorem 1.11]{em}, as a byproduct of their
theory of overtwisted contact structures\index[terms]{contact!structure!overtwisted}, \index[persons]{Eliashberg, Y.} Eliashberg and
\index[persons]{Murphy, E.} Murphy prove the following result: {\em Let $(M,\eta)$ be a closed $2n$-dimensional almost symplectic manifold
	and suppose that $[\mu]\in H^1(M; \Z)$ is a non-zero cohomology class. Then there exists a
 Lee class \index[terms]{class!Lee} $[\theta]$ in the formal homotopy class of $\mu$ such that $\theta=c\cdot\mu$
	for some real $c\neq 0$.} 
See also \cite{_Bertelson_Meigniez_} for a sharper result.

\item Historically, LCK manifolds were introduced in
  \cite{va_isr}, with what we labelled as ``fourth
  definition'', as almost Hermitian manifolds $(M,J,g)$
  whose Hermitian metric is locally conformal to some
  locally defined almost K\"ahler metrics: there exists an
  open covering $\{U_\al\}$ and almost K\"ahler metrics
  $g_\al$ on each $U_\al$ such that
  $g\restrict{U_\al}=e^{f_\al}g_\al$. This definition was
  then shown to be equivalent to the tensorial
  one. Indeed, if setting $\theta_\al=df_\al$, we get
  $d\omega\restrict{U_\al}=\theta_\al\wedge
  \omega\restrict{U_\al}$ and on intersections $U_\al\cap
  U_\be$ one has $\theta_\al\wedge
  \omega\restrict{U_\al\cap U_\be}=\theta_\be\wedge
  \omega\restrict{U_\al\cap U_\be}$. This implies
  $\theta_\al=\theta_\be$ on $U_\al\cap U_\be$ since the
  wedge product with the almost Hermitian
  form\index[terms]{form!Hermitian} is injective if
  $\dim_\C M\geq 2$. The first mention of the
  characterization in  terms of K\"ahler covers
  (\ref{_LCK_in_terms_of_covers_Definition_}) appears in
  \cite[Remark, p. 236]{va_gd}.

\item The notion of  singular K\"ahler space\index[terms]{space!K\"ahler!singular} (i.\,e.   a K\"ahler structure on a singular complex analytic space) is not yet fixed. There exists various attempts. This makes the study of LCK geometry on singular spaces very difficult. The first results were obtained by G. Ionita, O. Preda  and M. Stanciu in\index[persons]{Ionita, G.}\index[persons]{Preda, O.} \index[persons]{Stanciu, M.}  \cite{preda} and \cite{preda_stanciu,preda_stanciu_V_}, working with the definition of a singular K\"ahler space proposed in \cite{_Varouchas_}.

\item Indefinite LCK metrics are introduced in \cite{dd}. The main focus of the paper is on indefinite Vaisman manifolds  (\ref{_Vaisman_Definition_}).

\end{enumerate}

\section{Exercises}

\begin{enumerate}[label=\textbf{\thechapter.\arabic*}.,ref=\thechapter.\arabic{enumi}]

	\item Let $(M,\omega, \theta)$ be a LCS manifold\index[terms]{manifold!LCS} and let $\xi$ be the $\omega$-dual of the Lee form\index[terms]{form!Lee} $\theta$ ($\xi$ is called {\bf the Reeb field}). Show that $\Lie_\xi\omega=0$.\index[terms]{vector field!Reeb}
	
	\item Show that for any differential forms $\al, \be$ on a LCS manifold, one has:
	$$d_\theta(\al\wedge\be)=d_\theta\al\wedge\be+(-1)^{\deg(\al)}\al\wedge d_\theta\be+\theta\wedge\al\wedge\be.$$
	In particular, for a smooth function $f$, one has:
	$$d_\theta(f\al)=df\wedge\al+fd_\theta\alpha.$$
	
	\item \label{dtheta_inj} Show that a 
LCS manifold $M$ with Lee form\index[terms]{form!Lee} $\theta$ is  GCS if and
only if $d_\theta:  C^\infty(M)\ra\Lambda^1M$ is not
injective.
	
\item \label{_mult_by_symple_injective_Exercise_}
Let $\omega$ be a non-degenerate 2-form
on a smooth manifold $M$ of dimension $2n$, and
$\alpha\in \Lambda^k M$, $k \leq n-1$.
Suppose $\alpha \wedge \omega=0$. Prove that
$\alpha=0$.

	\item Let $\omega$ be a non-degenerate 2-form
on a smooth manifold of dimension $2n>4$, and $\theta$
a 1-form such that $\theta\wedge\omega=d\omega$.
Prove that $d\theta=0$.

{\em Hint:} Use the previous exercise for $\alpha:=d\theta$.
	
	\item Let $\lambda >1$ be a real number.
	Define  a {\bf weight $\lambda$ homogeneous form} on $\R^n\backslash 0$ 
	as a differential form $\eta$ that satisfies  $\rho_t^*\eta=\lambda^t\eta$,
	where $\rho_t$ is the homothety  $z \mapsto tz$,
        $t>0$. Prove that a $d$-closed weight $\lambda$ homogeneous form
        $\eta\in \Lambda^*(\R^n\backslash 0)$, satisfies $\eta=d\mu$,
where $\mu$ is weight $\lambda$ homogeneous.

	\item Let $M =\R^n\backslash 0/(x\sim 2x)$ be a real Hopf manifold, \index[terms]{manifold!Hopf!real}
	$\theta$ a closed and non-exact 1-form, 
	$d_\theta=d-\theta\wedge\cdot$, and let $H^*_\theta(M)$ be the cohomology
	of the complex $(\Lambda^*M, d_\theta)$
        (``Morse--Novikov cohomology''). Using the previous\index[terms]{cohomology!Morse--Novikov}
        exercise, prove that $H^i_\theta(M)=0$ for all $i$.

\item Let $M=X\times S^1$, $\pi:\; M \arrow S^1$ the standard projection,
and $\theta:= \pi^* dt$. Prove that  $H^i_\theta(M)=0$ for all $i$.

	
		\item (Moser stability for LCS manifolds\index[terms]{manifold!LCS}, \cite{bk1}) Let $(M,\omega,\theta)$ be a compact 
	locally conformally symplectic manifold, \index[terms]{Moser stability}
	$\omega_t$ a continuous deformation of $\omega$
	satisfying $d\omega_t=\omega_t\wedge\theta$
	and $[\omega_t]=\const$, where
	$[\omega_t]\in H^2_\theta(M)$ is the Morse--Novikov class of $\omega$.\index[terms]{class!Morse--Novikov}
	Find a flow of diffeomorphisms mapping $\omega$ to $\omega_t$.
		

	\item Prove  the  relation between the Levi--Civita connections of two conformal metrics $g'=e^{2f}g$:
	$$\nabla'=\nabla - df\otimes\id-\id\otimes df+ g\cdot \grad f,$$
	where $\grad f$ is defined by $g(\grad f,
        X)=\langle df, X\rangle$. 

\item \label{_Riemannian_volume_via_Kahler_Exercise_}
Let $X, Y\subset M$ be homologous compact complex submanifolds
in a K\"ahler manifold. Prove that they have the same 
Riemannian volume.

\item\label{_preimaga_non_compact_family_Exercise_}
 Let ${\cal X}=\{X_t\}$ be a continuous 
family of compact, simply connected submanifolds in an LCK manifold 
$(M, \omega, \theta)$
parametrized by $t\in S^1$. Denote by $\gamma\in H_1(M)$
the homology class realized by the corresponding circle in ${\cal X}$.
Assume that $\langle \gamma, \theta\rangle\neq 0$.
Prove that all connected components
of the preimage of ${\cal X}$ on the K\"ahler
cover of $M$ are non-compact.

\item \label{_product_with_CP^1_Exercise_}
Let $M= X \times \C P^1$, where $X$ is a compact complex 
manifold. Prove that any LCK structure\index[terms]{structure!LCK} on $M$ is GCK (globally
conformally K\"ahler).

{\em Hint:} Use Exercise 
\ref{_preimaga_non_compact_family_Exercise_} to show that
the volume of different connected 
components of the preimage of $\C P^1$ in the K\"ahler cover of $M$
is different, and apply Exercise 
\ref{_Riemannian_volume_via_Kahler_Exercise_}.

\item (see also Section \ref{_Topological_Criterion_})
Let $M$ be a compact LCK manifold
equipped with a holomorphic projection $\pi$ to a complex
variety $X$. Assume that the general fibre of $\pi$ is 
positive-dimensional and simply connected. Prove
that $M$ is GCK. 

{\em Hint:} Use the same logic as in Exercise
\ref{_product_with_CP^1_Exercise_}.

\end{enumerate}

\chapter{Hodge theory on complex manifolds and Vaisman's theorem}

\epigraph{\it Cast thy bread upon the waters: for thou shalt find it after many days.}{\sc\scriptsize Ecclesiastes, \ \ King James' Bible}

\section{Introduction}

``Hodge theory'' refers to the way of doing complex
algebraic geometry\index[terms]{geometry!algebraic} started by W. V. D. Hodge and S. Lefschetz\index[persons]{Hodge, W. V. D} \index[persons]{Lefschetz, S}
in the 1920s and 1930s. More specifically, it refers to
the applications of the ``Hodge decomposition'' to de Rham
cohomology and the de Rham differential graded algebra.
It turns out that most of the
 basic principles of algebraic
geometry can be reduced to results about harmonic
forms on compact K\"ahler manifolds. 
In Chapter \ref{_Harmonic_forms_chapter_} we
use a modern approach based on Lie superalgebras
to deduce the foundational notions of Hodge theory.

Here we define the Hodge decomposition
on differential forms, which lies at the heart of 
Hodge theory, and use it to prove a theorem of Vaisman
that is  a cornerstone of LCK geometry.

This chapter is very basic and should be
accessible to readers with a background in geometry.
We assume basic knowledge of linear algebra
and calculus on manifolds, including the notions
introduced in \ref{lck}.

Let $(M, \omega, \theta)$ be a compact LCK manifold,
$d\omega=\omega\wedge \theta$, where $\theta$ is
a closed 1-form. \index[terms]{theorem!Vaisman} Vaisman's theorem states that
unless $\theta$ is exact,\footnote{In this case,
$\omega$ is conformally equivalent to a K\"ahler form.}
$M$ does not admit a K\"ahler structure.
 Vaisman's theorem shows that on
compact, complex  manifolds, K\"ahler and locally
conformally K\"ahler geometries are dichotomic.
 
Our proof of  Vaisman's theorem is essentially
the same as the original one in \cite{va_tr}.

Note that
several authors (e.g.    A. \index[persons]{Moroianu, A.} Moroianu and collaborators) include
K\"ahler metrics among the LCK ones. However, they mostly
work on a non-compact case, when this dichotomy is not
applicable.

\section{Preliminaries}

Recall that for a complex linear space $(V,I)$ the 
{\bf Hodge decomposition} of the complexified space \index[terms]{Hodge!decomposition}
$V\otimes_\R \C:= V^{1,0}\oplus V^{0,1}$ is defined in such a way that
$V^{1,0}$ is the $\1$-eigenspace of $I$, and $V^{0,1}$ its $-\1$-eigenspace.

\hfill

\remark \label{hodge}
	The Hodge decomposition $V_\C=V^{1,0}\oplus V^{0,1}$ induces
	$$\Lambda^*_\C(V) =
	\Lambda^*_\C V^{0,1}\otimes \Lambda^*_\C V^{1,0},$$ giving
	\[ 
	\Lambda^d V_\C= \bigoplus_{p+q=d} \Lambda^p V^{1,0} 
	\otimes \Lambda^q V^{0,1}.
	\]
This follows from an isomorphism that is  true
for Grassmann algebras of arbitrary vector spaces:
$\Lambda^*(V\oplus W)= \Lambda^*V\otimes \Lambda^*W.$
We denote $\Lambda^p V^{1,0} \otimes \Lambda^q V^{0,1}$ by $\Lambda^{p,q}V$.

\hfill

The resulting decomposition 
$$\Lambda^n V_\C= \bigoplus_{p+q=n}\Lambda^{p,q}V$$
is called {\bf the Hodge decomposition of the Grassmann algebra\index[terms]{Grassmann algebra}}.

\hfill

\remark\label{_Hodge_dual_Remark_}
The Hodge decomposition\index[terms]{Hodge!decomposition} can be
carried on fibrewise to the tangent bundle of an almost
complex manifold and, applied to $V^*$, to the algebra of
differential forms. However, the complex structure on
$V^*$ is by convention opposite to one on $V$:
$\langle I(X), \xi\rangle= - \langle X, I(\xi)\rangle$
for $X\in V, \xi\in V^*$. This is done so that
$\langle I(X), I(\xi)\rangle=  \langle X, \xi\rangle$.

\subsection{Hodge decomposition on complex manifolds}

The $(p,q)$-decomposition is defined on the differential forms
on a complex manifold, as in \ref{hodge}.

Extend the action of the operator $I$ to differential forms by
$$(I\eta)(X_1,\ldots,X_k)=(-1)^k\eta(IX_1,\ldots,IX_k),\quad
\text{for}\,\, \eta\in \Lambda^kM$$
as in \ref{_Hodge_dual_Remark_}.

\hfill

\definition 
	Let $(M,I)$ be a complex manifold. 
	A differential form $\eta\in \Lambda^*M$
	is {\bf of type (1,0)} if $I(\eta)=\1\eta$, and
	{\bf of type (0,1)} if $I(\eta)=-\1\eta$. The corresponding
	vector bundles are denoted by $\Lambda^{1,0}M$ and
	$\Lambda^{0,1}M$. 

\hfill

\remark  
	The Cauchy--Riemann equations\index[terms]{Cauchy--Riemann equations} can be written as $df\in \Lambda^{1,0}(M)$.
	That is, {\em a function $f\in {C}^{\infty}_\C(M)$ is holomorphic\index[terms]{function!holomorphic}
		if and only if $df \in \Lambda^{1,0}M$.}

\hfill

\remark  
	Let $(M,I)$ be a complex manifold, and $z_1, \ldots, z_n$ a 
	holomorphic coordinate system in $U\subset M$, with $z_i$ being 
	holomorphic functions on $U$. {\em Then $dz_1,\ldots, dz_n$ generate
		the bundle $\Lambda^{1,0}U$, and $d\bar z_1,\ldots, d\bar z_n$
		generate $\Lambda^{0,1}U$.}

\hfill

\definition 
	Let $(M,I)$ be a complex manifold, $\{U_i\}$ a holomorphic atlas,
	and  $z_1,\ldots, z_n$ holomorphic coordinate system
	on each  patch.
	{\em The bundle $\Lambda^{p,q}(M,I)$ 
		of $(p,q)$-forms on $(M,I)$} is generated
	locally on each coordinate patch by monomials 
	$dz_{i_1}\wedge dz_{i_2}\wedge \cdots\wedge dz_{i_p}
	\wedge d\bar z_{i_{p+1}}\wedge \cdots\wedge d\bar z_{i_{p+q}}$.
	{\bf The Hodge decomposition} is the decomposition
	of vector bundles: 
	\[ \Lambda^d_\C M=\bigoplus_{p+q=d} \Lambda^{p,q} M.
	\]

\remark 
	The  de Rham differential\index[terms]{differential!de Rham} on a complex manifold
	has only two Hodge components: 
	\[
	d\left(\Lambda^{p,q}M\right)\subset
	\Lambda^{p+1,q}M\oplus  \Lambda^{p,q+1}M.
	\]

\definition 
	Let $d= d^{0,1}+d^{1,0}$ be the Hodge decomposition
	of the de Rham differential on a complex manifold: 
	$$d^{0,1}:\; \Lambda^{p,q}M \arrow \Lambda^{p,q+1}M, \ \ \ 
	\text{and}\ \ \  d^{1,0}:\; \Lambda^{p,q}M \arrow \Lambda^{p+1,q}M.$$
	The operators $d^{0,1}$, $d^{1,0}$ are denoted $\bar\6$
	and $\6$ and are called {\bf the Dolbeault differentials}\index[terms]{differential!Dolbeault}.

\hfill

\remark 
	(i) Since $d^2=0$, we have 
 $\6^2=0$, $\bar\6^2=0$, $\6\bar\6+\bar\6\6=0$.
	
	(ii) 
From the Newlander--Nirenberg theorem it is possible to see that
  $\bar\6^2=0$ is equivalent to the integrability of the
	complex structure (Exercise \ref{_Nijenhuis_tensor_Exercise_}).

\hfill

\definition 
	The {\bf complex (twisted) differential}\index[terms]{differential!complex (twisted)} operator is defined as $d^c=IdI^{-1}$.

\hfill

\claim \label{dc}
	The complex differential is related 
to $\6$ and $\bar\6$ by the formulae:
	$$\6=\frac{d-\1 d^c}{2},\qquad \bar\6=\frac{d+\sqrt{-1} d^c}{2}.$$

\begin{lemma}
	On a complex manifold, one has 
	$d^c = [{\cal W}, d]$, where $\caw$ is the Weil operator\index[terms]{Weil operator} $\caw (\eta)=(p-q)\1$ for $\eta\in\Lambda^{p,q}M$.
\end{lemma}

\hfill

\proof Clearly, $[{\cal W}, d^{1,0}]= \1 d^{1,0}$ and
$[{\cal W}, d^{0,1}] = - \1 d^{0,1}$. Adding these equations, we 
obtain $d^c = [{\cal W}, d]$.
\endproof

\hfill


This brings the following corollary.

\hfill

\corollary  \label{_anticommu_Corollary_}
$\{d, d^c\} = \{d, [d, {\cal W}]\}=0$,
where $\{a,b\}:= ab +ba$ is the anticommutator of the operators $a, b$. 
\endproof

\hfill

\remark 
	Clearly, $d= \6 + \bar\6$, $d^c=-\1 (\6-\bar\6)$,
	$dd^c=-d^cd = 2\1\6\bar\6$.

\subsection{Holomorphic one-forms and first cohomology}

\begin{lemma}\label{_holo_exact_1-form_Lemma_}
	Let $\eta$ be an exact holomorphic 1-form
	on a compact manifold. Then $\eta=0$.
\end{lemma}

\proof 
By the maximum principle, a holomorphic function
cannot have local maxima. Therefore, on a compact
manifold, any holomorphic function is constant.
Clearly, $\eta=df$, where $f$ is a function satisfying
$\bar\6 f=0$, and hence  holomorphic. Then $f=\const$
by the maximum principle.\index[terms]{maximum principle}
\endproof

\hfill

\definition 
	A $(0,p)$-form $\eta$ is called {\bf antiholomorphic}
	if $\bar\eta$ is holomorphic.

\hfill

The following result is implied by  Hodge theory 
(see \cite{demailly, griha, vois}):\index[terms]{Hodge!theory}

\hfill

\theorem \label{hodge1}
	Let $(M,I)$ be a compact K\"ahler manifold,
	and let $[\eta]\in H^1(M,\C)$ be a cohomology class.
	Then $[\eta]$ can be represented by a form $\eta=\eta^{1,0}+\eta^{0,1}$,
	where $\eta^{1,0}$ is holomorphic and $\eta^{0,1}$
        antiholomorphic.

\subsection{Positive $(1,1)$-forms}

\definition 
	A {\bf positive} $(1,1)$-form $\eta$ is a real $(1,1)$-form
	on a complex manifold, that can be locally written as
	$\eta=\sum_i a_i \beta_i \wedge I(\beta_i)$,
	where $\beta_i$ are real 1-forms, and $a_i$ are
        non-negative functions.\index[terms]{form!positive $(1,1)$}

\hfill

\remark  Hermitian forms are clearly positive. 
	Moreover, the cone of positive forms is the closure
	of the cone of Hermitian forms.  One may think of
	positive forms as of positive semi-definite Hermitian forms.\index[terms]{form!Hermitian}

\hfill

\remark \label{_French-positive_Remark_}
A form $\eta$ is positive if and only if 
$\eta (x, Ix)\geq 0$ for any tangent vector $x\in T_m M$.
Usually, this would be called ``semi-positive'', but
in complex algebraic geometry\index[terms]{geometry!algebraic} the French convention
prevails. To avoid the usual confusion,
we would suggest calling these forms ``French positive''.
For the benefit of the reader, we do not follow this suggestion
anyway.

\hfill

We recall the following classical fact of linear algebra. 

\hfill

\theorem  
\label{_simultaneous_diag_Theorem_}
\index[terms]{normal form}
	Let $g$ be a Hermitian metric on the real vector space $V$, and $h$ a pseudo-Hermitian form.\index[terms]{form!pseudo-Hermitian}
	There exists a $g$-orthonormal  basis 
	\[ \{x_1, I(x_1), x_2, I(x_2), \ldots, x_n, I(x_n)\}\] in $V^*$ 
	such that $h=\sum a_i x_i \wedge I(x_i)$.

\hfill

One often refers to this theorem as to
``a theorem on simultaneous diagonalization of two Hermitian forms''.

\hfill

\remark
A differential $(1,1)$-form $\eta$ on a manifold is positive
if and only if it is semi-Hermitian.
Let $(M,I, g)$ be an Hermitian manifold,
and $\eta$ a positive form. By the theorem about
simultaneous diagonalization of Hermitian forms\index[terms]{form!Hermitian} (\ref{_simultaneous_diag_Theorem_}),
locally there exists an orthonormal basis 
$\xi_1, ..., \xi_n \in \Lambda^{1,0}M$ and some  real and non-negative  $a_i$ such that
$\eta=-\1\sum_i a_i \xi_i \wedge \bar\xi_i$.
Moreover, let $\zeta_i\in  \Lambda^{1,0}M$ be
an arbitrary collection of sections (not necessarily a basis). Then
any linear combination of type
$-\1\sum_i u_i \zeta_i \wedge \bar\zeta_i$
with non-negative real $u_i$ is (``French'') positive.

\hfill

\definition 
	Hermitian forms are called {\bf strictly positive}.\index[terms]{form!Hermitian!strictly positive}

\hfill

\claim 
	Let $(M,I)$ be a complex manifold, and let $\eta$ be a real $(1,1)$-form.
	Then, for each 2-dimensional real subspace $W\subset T_x M$
	such that $I(W)=W$, {the restriction of $\eta$ to $W$
		is proportional to its volume form with non-negative
		coefficient.} Conversely, {if $\eta\restrict W$
		is non-negative for all such $W$, then $\eta$ is positive.}

\hfill

\proof A (1,1)-form $\eta$ is Hermitian if and only if $\eta(x,I(x))>0$
for each $x\in M$;  and it is positive if and only if $\eta(x,I(x))\geq 0$.
\endproof

\hfill

\definition 
	Let $(M,I, \omega)$ be a Hermitian manifold, $\dim_\C M=n$.
	The {\bf mass} of a positive $(1,1)$-form $\eta$, denoted $\Mass(\eta)$,  is the volume form
	$\eta \wedge \omega^{n-1}$.\index[terms]{form!positive $(1,1)$!mass of}

\hfill

Now we can state.

\hfill

\claim  
	Let $(M,I, \omega)$ be a Hermitian manifold of complex dimension $n$, 
	let $\{x_1, I(x_1), x_2, I(x_2),\ldots, x_n, I(x_n)\}$
	be an orthonormal basis in $\Lambda^1(M,\R)$, 
	and let $\eta=\sum_i \alpha_i x_i \wedge I(x_i)$
	be a positive $(1,1)$-form.
	{Then $\eta \wedge \omega^{n-1}= \sum \alpha_i \omega^n$.}

\hfill

\corollary  
	The mass of a positive form is always a (not strictly) positive volume form.	{A positive form vanishes if and only if its mass vanishes.}
\endproof

\section{Vaisman's theorem}

\theorem {(\cite{va_tr})}\label{vailcknotk}
Let $(M,\omega, \theta)$ be a compact LCK manifold, not globally conformally K\"ahler. 
 Then $M$ does not admit a K\"ahler structure.
 \index[terms]{theorem!Vaisman}

\hfill

\pstep That $M$ is not globally conformally K\"ahler means that $\theta$ is not cohomologous with zero, that is $\theta$ is not $d$-exact. 

Let $d\omega=\omega\wedge\theta$,
$\theta'=\theta+d\phi$. Then
$$d(e^\phi\omega)= e^\phi\omega\wedge \theta+
e^\phi\omega\wedge d\phi =e^\phi\omega\wedge \theta'.$$
 This means that
 we can replace the triple $(M,\omega, \theta)$ by
$(M,e^\phi \omega, \theta')$ for any 1-form $\theta'$
cohomologous to $\theta$. 

\hfill

{\bf Step 2:} Assume that $M$ admits a K\"ahler structure.
Then, by Hodge theory, $\theta$ is cohomologous to the sum
of a closed holomorphic and a closed antiholomorphic
form.\footnote{A form $\alpha$ is called {\bf
		antiholomorphic} if its complex conjugate $\bar\alpha$
	is holomorphic.}
After a conformal transformation (which changes $\theta$ to $\theta+d\phi$) as in Step 1,  we may assume
that $\theta$ itself  is the sum of a holomorphic and an  antiholomorphic form.

\hfill

{\bf  Step 3:} Then $dd^c\theta= \1 d \bar\6\theta=0$,
giving 
\[ dd^c(\omega^{n-1})= 
(n-1) d(I(\theta)\wedge\omega^{n-1})= 
- (n-1)^2\omega^{n-1}\wedge \theta\wedge I(\theta).
\] 
Therefore, 
$0=\int_M dd^c(\omega^{n-1})=\int \Mass(\theta\wedge I(\theta))$.
A mass of a positive form is positive, and strictly positive
if this form is non-zero. Then $\Mass(\theta\wedge I(\theta))=0$
implies $\theta\wedge I(\theta)=0$, giving $\Vert\theta\Vert^2=0$ and hence 
the  initial metric is globally conformally K\"ahler.
\endproof

\hfill

The same argument proves the following useful corollary.

\hfill

\corollary\label{_theta_not_d^c_closed_Corollary_}
Let $(M, \theta, \omega)$ be a compact LCK manifold. Then the cohomology class
$[\theta]\in H^1(M, \R)$ cannot be represented by a form that is 
$d^c$-closed.

\hfill

\proof
Indeed, for each representative of $[\theta]$, this form can be realized
as the Lee form\index[terms]{form!Lee} for an LCK metric that is  conformally equivalent to $\omega$.\index[terms]{metric!LCK}
Therefore, it would suffice to show that $d^c\theta\neq 0$ for any compact LCK
manifold $(M, \theta, \omega)$. If $d^c\theta= 0$, we would have
 $dd^c(\omega^{n-1})= \omega^{n-1}\wedge \theta\wedge I(\theta)$,
giving, as above, $0=\int_M dd^c(\omega^{n-1})=\int \Mass(\theta\wedge I(\theta))$,
hence $\theta\wedge I(\theta)=0$, implying that $\theta=0$.
\endproof

\hfill

\remark 
The above result shows that, unlike symplectic structures, that can coexist with LCS\index[terms]{structure!symplectic} structures\index[terms]{structure!LCS} on the same compact manifold, once one adds a compatible complex structure, one arrives at a true dichotomy: LCK non-GCK structures cannot coexist with K\"ahler structures\index[terms]{structure!K\"ahler} on the same compact complex manifold.

\hfill

\remark
\index[terms]{manifold!Vaisman} Vaisman theorem does not hold for non-compact LCK manifolds:
there exists an LCK structure\index[terms]{structure!LCK} with non-exact Lee form\index[terms]{form!Lee} on
$\C^* \times \C^n$, as the following example shows.

\hfill

\example\label{_LCK_non_compact_Lee_non_exact_Example_}
Consider the complex manifold $M= \C^* \times \C^n$
with coordinates $z_0, ..., z_n$ and let $\tilde M= \C^{n+1}$ 
with coordinates $u_0, z_1, ..., z_n$ map to $M$
via $z_0 = e^{u_0}$. Then $M = \tilde M/\Z$,
with $\Z$ generated by the map
$\gamma(u_0, z_1, ..., z_n)= (u_0+2\1\pi, z_1, ..., z_n)$. 
Consider the K\"ahler form
$\tilde \omega=  dd^c\alpha$
on $\tilde M$, 
where $\alpha$ is a positive, strictly plurisubharmonic function
that satisfies $\alpha(x + 2\1\pi) = c\alpha(x)$ for all $x\in\C$,  with $c>1$.
For example, we could take $\alpha= |e^{-\1 u_0}|^2+ 
|e^{-\1 u_0}|^2\sum_{i=1}^n |z_i|^2$.
Since $\alpha$ is a sum of squares of absolute values of holomorphic functions,
it is plurisubharmonic; it is strictly plurisubharmonic, 
because 
\begin{multline*} dd^c\alpha= (-\1)\6(e^{-\1 u_0})\wedge \bar\6(e^{-\1 \bar u_0})+\\+
(-\1)\sum_{i=1}^n \6(e^{-\1 u_0}z_i)\wedge \bar\6(e^{-\1 \bar u_0}\bar z_i),
\end{multline*}
each of the terms 
\[ (-\1)\6(e^{-\1 u_0})\wedge \bar\6(e^{-\1 \bar u_0})\] and
\[ (-\1)\6(e^{-\1 u_0}z_i)\wedge \bar\6(e^{-\1 \bar u_0}\bar z_i)\]
is positive, and the 1-forms $\6(e^{-\1 u_0})$, 
$\6(e^{-\1 u_0}z_i)$ form a basis of the space $\Lambda^{1,0}(\C^{n+1})$.

Clearly, $\tilde \omega$ is a K\"ahler form
that satisfies $\gamma^*(\tilde\omega) = c\tilde \omega$,
where $\gamma$ is the generator of the 
$\Z$-action\index[terms]{action!$\Z$-} on $\tilde M$
defined above. The non-compact complex manifold $M$ is thus LCK. The corresponding weight bundle\index[terms]{bundle!weight}
has non-trivial monodromy,\index[terms]{monodromy} defined by the
character $\chi(\gamma) =c$,  because 
$\gamma^*(\tilde\omega) = c\tilde \omega$.
Therefore, the Lee form\index[terms]{form!Lee} is not exact. Obviously, $M$ also admits a (product) K\"ahler metric.\index[terms]{metric!K\"ahler}

\section{Exercises}

\begin{enumerate}[label=\textbf{\thechapter.\arabic*}.,ref=\thechapter.\arabic{enumi}]

\item
Let $\Omega$ be a holomorphic $(n-1)$-form on a compact
complex manifold $M$ with $\dim_\C M=n$. Prove that
$d\Omega=0$.

\item
Let $M$ be a compact complex manifold, $\dim_\C M=2$.
Prove that all holomorphic differential forms on $M$ are closed.

	\item 
	Prove \ref{dc}: $\6=\frac{d+\1 d^c}{2}$, and $\bar\6=\frac{d-\sqrt{-1} d^c}{2}.$

\item\label{_Nijenhuis_tensor_Exercise_}
Let $(M,I)$ be an almost complex manifold,
and $\bar\6=\frac{d-\sqrt{-1} d^c}{2}.$ Prove that $\bar\6^2=0$
if and only if the Nijenhuis tensor of $I$ vanishes. 

	\item  Prove that  on a compact K\"ahler manifold,
	any holomorphic form is closed.
	
	\item Prove \ref{hodge1} for complex curves,
without relying on Hodge theory: {\em Let $(M,I)$ be a compact K\"ahler manifold,
	$\dim_\C M=1$,	and  $[\eta]\in H^1(M,\C)$  a cohomology class.
		Then $[\eta]$ can be represented by a form $\eta=\eta^{1,0}+\eta^{0,1}$,
		where $\eta^{1,0}$ is holomorphic and $\eta^{0,1}$ antiholomorphic.}

\item Let $\eta$ be a positive (1,1)-form on a compact almost complex
manifold. Prove that there exists a collection of (1,0)-forms
$\xi_1,..., \xi_n$ and non-negative real functions $u_1, ..., u_n$
such that $\eta= -\1 \sum_{i=1}^n u_i \xi_i \wedge \bar\xi_i$, or find a
counterexample.

\item Let $M$ be an almost complex manifold, $\dim_\R M=2n$.
A differential form $\eta\in\Lambda_\R^{n-1,n-1}(M)$ is called {\bf positive}
if $\eta \wedge \omega$ is a non-negative volume form for
any positive form $\omega\in \Lambda_\R^{1,1}(M)$. Prove that
any positive $(n-1,n-1)$-form $\eta$ satisfies $\eta=\kappa^{n-1}$
for some positive (1,1)-form $\kappa$.

	\item Let $\omega$ be a non-degenerate (1,1)-form on 
a compact complex manifold,
satisfying $d\omega=\theta \wedge\omega$. 
\begin{enumerate}
\item Assume that $\dim_\C M >2$. Prove that $d\theta=0$.
\item Assume, moreover $d^c\theta=0$, and $\omega(\theta, I\theta)>0$. 
Prove that $\theta=0$.
\end{enumerate}
	
	\item Let $(M,\omega, \theta)$ be a compact LCK manifold, $\dim_\C M> 2$, satisfying
	$dd^c\omega=0$. Prove that $\theta=0$.
	
	\item \label{_Gauduchon_co-closed_Exercise_}
Prove that on a compact LCK manifold
          $(M,\omega,\theta)$, $\dim_\C M =n$, satisfying
          $dd^c\omega^{n-1}=0$\footnote{Hermitian metrics
that satisfy this condition are called {\bf \index[persons]{Gauduchon, P.} Gauduchon},
see Section \ref{gauduchon_metric}.}, \index[terms]{metric!Gauduchon}
the Lee form\index[terms]{form!Lee} $\theta$ is co-closed: $d^*\theta=0$. 
	
	\item Let $M$ be a complex manifold, $\dim_\C M=n$.
	A Hermitian metric on $M$ is called {\bf balanced}
	if $d\omega^{n-1}=0$ (see \cite{_Michelson_}, also \cite[Chapter 4]{_Popovici:book_}). Prove that a classical Hopf manifold $\C^n\setminus  0/x\sim \lambda x$
	does not admit  balanced metrics.\index[terms]{metric!balanced}

{\em Hint:} Find an $(n-1)$-dimensional Hopf submanifold in $M$
and use the Stokes' theorem.
	
	\item Let $\omega$ be a non-degenerate 2-form on a
	$2n$-dimensional smooth manifold, and $d(\omega^k)=0$ for
	some $k$ satisfying $0<k<n-1$. Prove that
        $d\omega=0$.

\item Let $X$ be a K\"ahler manifold.
Prove that $M= \C^* \times X$ admits
an LCK metric with non-exact Lee form.\index[terms]{form!Lee}

{\em Hint:} See \ref{_LCK_non_compact_Lee_non_exact_Example_}.

\item
Let $(M, I, \omega)$ be a Hermitian $n$-manifold,
and $L_{\omega^{n-1}}(\alpha):= \alpha\wedge \omega^{n-1}$.
Prove that $L_{\omega^{n-1}}:\; \Lambda^1(M)\arrow \Lambda^{2n-1}(M)$
is an isomorphism.

\item\label{_LCK_balanced_Kahler_Exercise_}
Prove that any LCK metric that is  balanced is also K\"ahler.

{\em Hint:} Use the previous exercise.

\item 
A Hermitian metric $\omega$ on a complex $n$-manifold $M$
is called {\bf SKT} (``strongly K\"ahler with torsion''), 
or {\bf pluriclosed}, if $dd^c \omega=0$.
Let $\omega$ be a metric on a complex $n$-manifold 
that is  both LCK and SKT. Prove that 
$dd^c(\omega^{n-1})=(n-1)^2 \theta \wedge \theta^c\wedge \omega^{n-1}$,\index[terms]{metric!pluriclosed (SKT)}
where $\theta$ is the Lee form.\index[terms]{form!Lee}

\item\label{_LCK_SKT_Kahler_Exercise_}
Let $(M, \omega, \theta)$ be a compact LCK manifold, $\dim_\C M > 2$.
Assume that $\omega$ is SKT. Prove that $\omega$ is K\"ahler.

{\em Hint:} Use the previous exercise.

\end{enumerate}

The following exercises assume some knowledge of Hodge theory,
chapter 0 of \cite{griha} is definitely sufficient;
for a more detailed account of the K\"ahler identities,
see \ref{_kah_susy_Theorem_}.

\begin{enumerate}[label=\textbf{\thechapter.\arabic*}.,ref=\thechapter.\arabic{enumi}]
\setcounter{enumi}{18}

\item
Let $(M, I, \omega)$ be a Hermitian $n$-manifold,
and $\langle L, H, \Lambda\rangle$ 
its Lefschetz $\goth{sl}(2)$-triple, \ref{_kah_susy_Theorem_}.
A form $\alpha\in \Lambda^p(M)$, $p \leq n$
is called {\bf primitive} if $\Lambda \alpha =0$.
\begin{enumerate}\index[terms]{form!primitive}
\item Prove that $\alpha \wedge \omega^{n-p+1}=0$
if and only if $\alpha$ is primitive.
\item Prove that $\alpha$ is primitive if and only 
if $I(\alpha)$ is primitive.

\end{enumerate}

\item
Let $(M, \omega)$ be a Hermitian $n$-manifold.
Prove that $\omega$ is balanced $\Leftrightarrow$  $d\omega$ is primitive 
$\Leftrightarrow$  $d^c\omega$ is primitive.

\item\label{_HR_3-form_Exercise_}
Let $\alpha$ be a primitive\index[terms]{form!primitive} 3-form of Hodge type $(2,1)$
on a Hermitian manifold. Prove that 
$\alpha \wedge \bar\alpha \wedge \omega^{n-2}= -\1 |\alpha|^2 C \omega^{n}$,
where $C$ is a positive constant depending only on the dimension of $M$.

{\em Hint:} This is a special case of ``Hodge--Riemann relations''.\index[terms]{relations!Hodge--Riemann}

\item
Let $(M, \omega, I)$ be an SKT $n$-manifold, and $\alpha:= \6 \omega$.
Prove that $\6\bar\6\omega^{n-1} =   \alpha \wedge \bar\alpha\wedge \omega^{n-3}$.

\item
Let $(M, \omega, I)$ be a balanced SKT manifold, and let
$\alpha:= \6 \omega$.
Prove that $\6\bar\6\omega^{n-1} = -\1 |\alpha|^2 C \omega^{n}$,
where $C$ is the same constant as in Exercise 
\ref{_HR_3-form_Exercise_}.

\item\label{_balanced_SKT_Kahler_Exercise_}
Let $(M, \omega, I)$ be a compact balanced SKT manifold.
Prove that $d\omega=0$.

{\em Hint:} Use the previous exercise.

\item
Let $D\subset M$ be a compact divisor in a balanced
manifold. Prove that $D$ cannot be homologous to zero.

\item\label{_balanced_exact_no_divisors_Exercise_}
Let $M$ be a compact balanced $n$-manifold
with $H^{2n-2}(M)=0$. Prove that $M$ has
no divisors.

{\em Hint:} Use the previous exercise.

\item\label{_Bismut_torsion_Exercise_}
Let $(M, I, \omega, g)$ be a
complex Hermitian manifold, and 
$\nabla^b$ a con\-nec\-tion on $TM$
such that $\nabla^b(I)=\nabla^b(\omega)=0$. Let 
$T^b\in \Hom(\Lambda^2TM, TM)$ be  the torsion form
of $\nabla^b$, and assume that it satisfies 
$g(T(X,Y), Z)= - g(T(X,Z), Y)$.
Prove that $\nabla^b$
can be expressed in terms of the 
Levi--Civita connection $\nabla$ of $g$ as
\[
g(\nabla^b_XY,Z)=g(\nabla_XY,Z)+\frac 12 d\omega(IX,IY,IZ). 
\]

\item\label{_LCK_p_pluriclo_Exercise_}
A Hermitian form $\omega$ on a complex $n$-manifold
is called {\bf $p$-pluriclosed}\index[terms]{form!$p$-pluriclosed} if $dd^c(\omega^p)=0$, 
where $0 < p < n-1$. 
\begin{enumerate} 
\item Let $\Lambda:= L^*$
be the Lefschetz operator, $A:= \Lambda^2(dd^c\omega)$ and 
$B:= \Lambda^3(d\omega\wedge d^c\omega)$. 
Prove that for any  $p$-pluriclosed\index[terms]{form!$p$-pluriclosed} $\omega$,
we have $p A= - p(p-1) B$.

{\em Hint:\ } Use
$dd^c(\omega^p)= p(p-1) d\omega\wedge d^c\omega \wedge \omega^{p-2}
+ p dd^c\omega \wedge \omega^{p-1}.$
\item 
Let $\omega$ be a form that is  $p$-pluriclosed\index[terms]{form!$p$-pluriclosed}
and $q$-pluriclosed, for $p\neq q$. Prove that $A=B=0$.
Prove that $\omega$ is closed if it is LCK.

{\em Hint:} Prove that 
$\Lambda^3(d\omega\wedge d^c\omega)= \const \cdot |\theta|^2$
whenever $d\omega=\theta\wedge\omega$.
\end{enumerate}

\end{enumerate}

\chapter{Holomorphic vector bundles} 
\label{_holomorphic_bundles_Chapter_}

\epigraph{\it {\bf Meander, vi.}  To proceed sinuously and aimlessly. The word is the ancient name of a river about one hundred and fifty miles south of Troy, which turned and twisted in the effort to get out of hearing when the Greeks and Trojans boasted of their prowess.}{\sc \scriptsize The Devil's Dictionary by Ambrose Bierce }

\section{Introduction}

In this chapter, we reproduce several calculations
related to holomorphic vector bundles and their
total spaces. Our interest is motivated by the
following construction: given a line bundle
$L$ on $M$, denote by $\Tot^\circ(L)$ the
space of non-zero vectors in its total space.
Let $\phi$ be a holomorphic map from $M$
to itself, such that $\phi^*(L)\cong L$.
Denote the corresponding automorphism
of $\Tot^\circ(L)$ by $\phi_L$. Assume that
$\phi_L$ multiplies the length of a vector
by a constant $c> 1$. It is easy to see
that the quotient $\Tot^\circ(L)/\langle \phi_L\rangle$
is a complex manifold. Later in this book, we 
are going to prove that whenever the curvature of $L$ is a K\"ahler
form (in this case $L$ is called  
``a positive line bundle''), the manifold
$\Tot^\circ(L)/\langle \phi_L\rangle$ is LCK.
This construction is quite important, because
it gives a special class of LCK manifolds,
so-called ``Vaisman manifolds'' (Chapter \ref{vaiman}).
Moreover, any complex manifold admitting
a Vaisman structure, can be obtained from
this construction when $M$ is an orbifold
(see Chapter \ref{orbif} for an introduction
to the orbifolds).\index[terms]{manifold!Vaisman}

On a holomorphic Hermitian vector bundle, one defines
{\em the  Chern connection}\index[terms]{connection!Chern}
(\ref{_Chern_conne_Definition_}),  that is  uniquely determined
by the holomorphic structure and the Hermitian form.\index[terms]{form!Hermitian}\index[terms]{connection!Chern}
We assume that the Chern connection is a familiar
notion for the reader (\cite{demailly}, \cite{griha}).

Write the corresponding Ehresmann connection\index[terms]{connection!Ehresmann}
explicitly, as follows. Let $l:\; \Tot B \arrow \R^{\geq  0}$
 be the square length function, $l(v)=|v|^2$.  
Define $T_\hor \Tot B$ as the orthogonal complement
to the fibrewise tangent space with respect to 
the form $\eta:= dd^cl$. It turns out that
this defines a linear Ehresmann connection,
and the corresponding connection in $B$
is precisely the  Chern connection
(\ref{_Ehresmann_Chern_Theorem_}).

This is used to prove the fundamental
formula that is  used everywhere in LCK geometry;
for the lack of a better name, we call it\index[terms]{Calabi formula (2.6)} 
``Calabi formula (2.6)'', from \cite{_Calabi:hk_};
see also \cite[(15.19)]{besse}.\footnote{Note that 
in \cite[(15.19)]{besse}, the sign of one term in 
(15.19) is wrong; this seems to be a common error.}

Let $B$ be a holomorphic Hermitian vector bundle,
and $l(v):=|v|^2$ the length function on $\Tot(B)$.
 Calabi formula (2.6) expresses the (1,1)-form
$dd^c l$ in terms of the curvature $\Theta_B$ of the Chern
connection on $B$. Essentially, it interprets
$\Theta_B$ as a (1,1)-form on $\Tot(B)$, and
writes $dd^c l$ as a sum of the curvature and the 
fibrewise Hermitian metric on $\Tot(B)$.

This is used to construct a K\"ahler
metric on a total space of a negative
holomorphic line bundle\footnote{A holomorphic line
bundle is negative if the curvature of its
 Chern connection\index[terms]{connection!Chern} is a strictly 
negative (1,1)-form; it is equivalent to $L^*$ being
positive, which implies ampleness by Kodaira theorem.}.
This K\"ahler metric has conical-type symmetry;
later on we use it to construct LCK structures\index[terms]{structure!LCK}
on the $\Z$-quotients of these total spaces.

\section{Holomorphic vector bundles}

\subsection{Holomorphic structure operator}

\definition A {\bf holomorphic
	vector bundle} on a complex manifold $M$
is a locally free sheaf of $\calo_M$-modules.\index[terms]{bundle!vector bundle!holomorphic}\index[terms]{bundle!vector bundle!smooth}

\hfill

\claim
Let $B$ be a holomorphic vector bundle.
Then the sheaf $B_{{C}^{\infty}}:=B \otimes_{\calo_M} {C}^{\infty} (M)$ is  locally free, and hence {$B_{{C}^{\infty}}$ is 
	a smooth vector bundle.}

\hfill

\definition
$B_{{C}^{\infty}}$ is called the {\bf smooth vector bundle underlying $B$}.

\hfill

\definition
Let $d= d^{0,1}+d^{1,0}$ be the Hodge decomposition
of the de Rham differential on a complex manifold, 
$$d^{0,1}:\; \Lambda^{p,q}M \arrow \Lambda^{p,q+1}M, \quad 
\text{and}\quad d^{1,0}:\; \Lambda^{p,q}M \arrow \Lambda^{p+1,q}M.$$
The operators $d^{0,1}$, $d^{1,0}$ are denoted $\bar\6$
and $\6$ and called {\bf the Dolbeault differentials}.\index[terms]{differential!Dolbeault}

\hfill

\remark\label{_6^2=0_Remark_}
(i) From $d^2=0$,
one obtains {$\bar\6^2=0$ and $\6^2=0$.}

(ii) {The operator $\bar\6$ is $\calo_M$-linear.}

\hfill

\definition
Let $B$ be a holomorphic vector bundle, and  
$\bar\6:\; B_{{C}^{\infty}}\arrow B_{{C}^{\infty}}\otimes \Lambda^{0,1}M$
the operator mapping $b \otimes f$ to $b\otimes \bar\6 f$,
where $b\in B$ is a holomorphic section, and $f$ a 
smooth function. This operator is called {\bf the
	holomorphic structure operator} on $B$. {It is 
	correctly defined, because $\bar\6$ is $\calo_M$-linear.}

\hfill

\remark  Clearly, the kernel of {$\bar\6$ coincides with the set 
		of holomorphic sections} of the bundle $B$.


\subsection{The $\bar\6$-operator on vector bundles}
\index[terms]{operator!$\bar\6$ on vector bundles}


\definition 
	{A $\bar\6$-operator} on a smooth bundle 
	is a map $V \stackrel {\bar\6}\arrow \Lambda^{0,1}M\otimes V$,
	satisfying $\bar\6(fb) = \bar\6(f)\otimes b + f\bar\6(b)$
	for all $f\in {C}^{\infty} (M), b\in V$.

\hfill

\remark  A  $\bar\6$-operator on $B$ can be extended
		to 
		\[ \bar\6:\; \Lambda^{0,i}M\otimes V \arrow \Lambda^{0,i+1}M\otimes V,\]
using $\bar\6 (\eta \otimes b) = \bar\6(\eta)\otimes b + 
	(-1)^{\tilde \eta}\eta\wedge\bar\6(b)$, 
	where $b\in V$,  $\eta \in \Lambda^{0,i}M$,
and $\tilde \eta=i$ denotes the degree of $\eta$.

\hfill

\remark  Generally, $\bar\6^2$ does not need to vanish;
however, if $\bar\6$ is a holomorphic structure 
operator,\index[terms]{operator!holomorphic structure}
	then $\bar\6^2=0$, by \ref{_6^2=0_Remark_}.

\hfill

\theorem \label{_Koszul--Malgrange_Theorem_}
{\bf (\index[persons]{Koszul, J.-L.}Koszul--Malgrange)}\\
 Let $V$ be a $C^\infty$-bundle over a complex
manifold $M$, and  $\bar\6:\; V \arrow \Lambda^{0,1}M\otimes V$
	 a $\bar\6$-operator, satisfying $\bar\6^2=0$. {Then
		$B:=\ker \bar\6\subset V$ is a holomorphic vector
		bundle of the same rank as $V$.}\index[terms]{theorem!Koszul--Malgrange}

\proof  See \cite[Chapter I, Proposition
  3.7.]{_Kobayashi_Bundles_},
\cite{_Koszul_Malgrange_}.
\endproof

\hfill

\remark  This statement is the vector bundle analogue
	of the New\-lan\-der-Nirenberg theorem (\ref{nn}).

\hfill

\definition 
	A $\bar\6$-operator $\bar\6:\; V \arrow \Lambda^{0,1}M\otimes V$
	on a $C^\infty$-bundle $V$ is called a {\bf 
		holomorphic structure operator} if 
$\bar\6^2=0$.\index[terms]{holomorphic structure}


\subsection{Connections and holomorphic structure operators}


\definition \label{_compatible_with_holo_Definition_}
	Let $B$ be a smooth complex vector bundle with 
a complex-linear connection $\nabla$
	and the holomorphic structure
operator $\bar\6:\; B \arrow \Lambda^{0,1}M\otimes B$. 
	Consider the Hodge decomposition of $\nabla$,\index[terms]{Hodge!decomposition}
	$\nabla= \nabla^{0,1} + \nabla^{1,0}$, defined as follows:
	\[
	\nabla^{0,1}:\; V \arrow \Lambda^{0,1}M\otimes V, \ \ \ 
	\nabla^{1,0}:\; V \arrow \Lambda^{1,0}M\otimes V.
	\]
	We say that $\nabla$ is {\bf  compatible 
		with the holomorphic structure} if $\nabla^{0,1}=\bar\6$.

\hfill

\definition 
	A {\bf Hermitian holomorphic vector bundle}
	is a smooth complex vector bundle equipped with a Hermitian
	metric and a holomorphic structure.\index[terms]{bundle!vector bundle!holomorphic Hermitian}

\hfill

\definition \label{_Chern_conne_Definition_}
	A {\bf \index[persons]{Chern, S.-S.} Chern} (or {\bf second canonical}, 
\cite{li2}, \cite{gau_tor})  connection on a
	holomorphic Hermitian vector bundle is a connection
	compatible with the holomorphic structure and 
preserving the metric.\index[terms]{connection!Chern}

\hfill

\theorem 
	On any holomorphic Hermitian vector bundle, {a 
		Chern connection\index[terms]{connection!Chern} exists, and is unique.} 

\proof 
\cite[Proposition 3.12]{vois}, \cite[Proposition 10.3]{demailly}.
\endproof

\hfill

\remark \label{chern_tor}
	Being unique, we call it  {\bf the Chern
          connection}.\index[terms]{connection!Chern} It is the analogue of the
        Levi--Civita connection from Riemannian geometry. 
Consider the  Chern connection in the holomorphic tangent
bundle  of a complex Hermitian manifold $(M,I, h)$. Unless 
the Hermitian form \index[terms]{form!Hermitian}$h$ is K\"ahler, the Chern connection\index[terms]{connection!Chern}
has torsion. However, if $h$ is K\"ahler, the
Levi--Civita connection coincides with the Chern
connection, as can be easily seen from the definitions.


\subsection{Curvature of holomorphic line bundles}


The term {``curvature of a holomorphic
	bundle''}  usually means the curvature of the
 Chern connection\index[terms]{connection!Chern} associated with some Hermitian metric.

Let $B$ be a holomorphic Hermitian line bundle, and $b$ 
a non-degenerate holomorphic section of $B$. Denote by $\eta$ a (1,0)-form
that satisfies $\nabla^{1,0} b=\eta\otimes b$.
Then: 
$$d|b|^2= \Re g(\nabla^{1,0} b, b) = \Re\eta|b|^2,$$
which gives 
$$\nabla^{1,0} b= \frac{\6 |b|^2}{|b|^2}b=
2\6\log|b| b.$$

Then: 
\begin{equation}\label{_Chern_conne_curvature_dd^c_log_Equation_}
\Theta_B(b)= 2\1 \bar\6\6\log|b| b,\,\, \text{that is,}\,\,\,
\Theta_B = -2\1 \6\bar\6\log|b|.
\end{equation}

\claim\label{_dd^c_log_independent_Claim_}
The 2-form $2\6\bar\6\log|b|$ is independent on  the choice of 
a holomorphic section $b$.

\hfill

\proof  Indeed, let $b_1=b f$, where $f$ is a non-vanishing holomorphic function.
Then: 
\begin{equation*}
	\begin{split} 
\6\bar\6\log(|b_1|^2)&= 2\6\bar\6\log(|b|)+ \6\bar\6\log (f\bar f)\\
&=2\6\bar\6\log(|b|)+ \6\bar\6\log f+ \6\bar\6\log \bar f.
	\end{split}
\end{equation*}

The last two terms vanish, because the logarithm of a holomorphic
function is also holomorphic, and the logarithm of an\index[terms]{logarithm} antiholomorphic
function is antiholomorphic. \endproof 

\subsection{K\"ahler potentials and plurisubharmonic functions}

\definition 
	A real-valued smooth function on a complex manifold
	is called {\bf plurisubharmonic (psh)} if the (1,1)-form $dd^c f$
	is positive, and {\bf strictly plu\-ri\-sub\-har\-mo\-nic} if $dd^c f$
	is an Hermitian form.\index[terms]{form!Hermitian}\index[terms]{function!plurisubharmonic!(strictly) plurisubharmonic}

\hfill

\remark 
	Since $dd^c f$ is always closed, {it is also K\"ahler when
		it is strictly positive.}

\hfill

\definition 
	Let $(M,I,\omega)$ be a K\"ahler manifold.
	A {\bf K\"ahler potential} is a smooth real-valued function $f$
	such that $dd^c f=\omega$.\index[terms]{potential!K\"ahler}

\hfill

\remark  {Locally, K\"ahler potentials exist for any
  K\"ahler form}.\index[terms]{form!K\"ahler}
	This is a non-trivial theorem that follows from the 
	Poincar\'e--Dolbeault--Grothendieck lemma (\cite{griha, vois}).\index[terms]{lemma!	Poincar\'e--Dolbeault--Grothendieck}

\hfill

\example  The function $z\mapsto |z|^2$ is a K\"ahler potential for the
	usual (flat) Hermitian metric on $\C^n$.

\subsection[Chern connection obtained from an Ehresmann connection]{Chern connection obtained from an Ehresmann\\ connection}\index[terms]{connection!Chern}\index[terms]{connection!Ehresmann}

Let $B$ be a holomorphic Hermitian vector bundle
on a complex manifold $M$, and $\pi:\; \Tot B \arrow M$ its total space.
Consider the function $l:\; \Tot B \arrow \R$ 
taking $v\in B\restrict x $ to $|v|^2$, and let
$dd^c  l$ be the corresponding (1,1)-form on $\Tot B$.
Clearly, $dd^c  l$ is a Hermitian form,\index[terms]{form!Hermitian} constant 
on the fibres of $\pi$.

\hfill

\remark \label{_Ehre_holo_defined_Remark_}
Let $W\subset T\Tot B$ be the orthogonal
complement to $T_\pi \Tot B$ with respect to the
form $dd^c l$. Since $dd^c l$ is non-degenerate on
the fibres of $\pi:\; \Tot B \arrow M$, we have
a direct sum decomposition
 $T_\pi \Tot B\oplus W = T\Tot B$.
Therefore, $T_\hor\Tot B:= W$ 
defines an Ehresmann connection\index[terms]{connection!Ehresmann} on $\Tot B$.
It is easy to see that this connection is linear.
The following theorem gives an alternative definition
of the  Chern connection.\index[terms]{connection!Chern}

\hfill

\theorem\label{_Ehresmann_Chern_Theorem_}
Let $B$ be a holomorphic Hermitian vector bundle
on a complex manifold $M$, $\pi:\; \Tot B \arrow M$ its total space,
and $\nabla$ the linear  Ehresmann connection defined above
in \ref{_Ehre_holo_defined_Remark_}.
Then $\nabla$ induces the  Chern connection on the vector 
bundle $B$.
 
\hfill

\pstep The  Chern connection\index[terms]{connection!Chern} is defined as the unique connection
${}^C\!\nabla$ that is  unitary and satisfies 
${}^C\!\nabla^{0,1}(\text{holomorphic sections})=0$, 
that is  clearly equivalent to $\bar\6={}^C\!\nabla^{0,1}$. 
We prove that $\nabla={}^C\!\nabla$
by checking that $\nabla$ satisfies these two conditions.

The second condition 
is easy to verify. Indeed, a holomorphic
section $b$ of $B$ gives a holomorphic function $\mathcal b$ on 
$\Tot(B^*)$. For any $(0,1)$-vector field $X$ on $M$,
its horizontal lift $X_\hor$ to $T \Tot(B^*)$ is also
a (0,1)-vector field, and hence  the Lie derivative
$\Lie_{X_\hor} \mathcal b$ vanishes. After dualizing, this implies
$\nabla^{0,1}$(holomorphic sections)$\ =0$,
where $\nabla$ is the connection on $B$ induced by
the Ehresmann connection.\index[terms]{connection!Ehresmann}

\hfill

{\bf Step 2:} It remains only to show that the
vector bundle connection induced by $\nabla$ is unitary.
As in \ref{_Ehresmann_vs_vector_bundle_Proposition_},
the covariant derivative operator $\nabla_X$ is
identified with the Lie derivative along the
horizontal lift $X_\hor$ of $X$.
To finish the proof it suffices to show that
the diffeomorphism flow $V_t$, associated with $X_\hor$, 
induces an isometry on the fibres of $B$.

Let $l(v)=|v|^2$ be the length function on $\Tot(B)$.
Since the metric on $B$ is determined by the corresponding
quadratic form $l$, to prove that $V_t$ induces an isometry, 
we need only to show that $\Lie_{X_\hor} l=0$. 
Taking the Hodge components, it will suffice to prove this assuming 
$X$ is a (0,1)-vector field.

Let $\vec r=\sum_i t_i \frac \6{\6t_i}\in T_\pi\Tot B$ be the radial vector field
(the one tangent to the homothety flow). 
Since the radial vector field $\vec r$ is holomorphic, it commutes with
$X_\hor$, and 
$dd^c l(X_\hor,\vec r) =  \Lie_{X_\hor}\Lie_{I(\vec r)}  l$
by the Cartan formula. Since $l$ is homogeneous of weight 2
with respect to the homothety, we have 
$\Lie_{I(\vec r)}  l=2l$. This gives
\begin{equation}\label{_dd^c_l_radial_Equation_}
dd^c l(X_\hor,\vec r) =  \Lie_{X_\hor}\Lie_{I(\vec r)}  l=2 \Lie_{X_\hor}l.
\end{equation}
However, by definition of horizontal tangent vectors,
they are $(dd^cl)$-orthogonal to the vertical tangent vectors,
giving $dd^c l(X_\hor,\vec r) = 0$. Therefore,
\eqref{_dd^c_l_radial_Equation_} brings $\Lie_{X_\hor} l=0$,
implying that the Ehresmann connection\index[terms]{connection!Ehresmann} $\nabla$ induces
an isometry.
\endproof

\subsection{Calabi formula (2.6).}\index[terms]{Calabi formula (2.6)}

\definition\label{_quadratic_from_Herm_Definition_}
Let $B$ be a complex vector bundle on $M$.
A skew-Hermitian endomorphism 
$a\in {\goth u}(B)$ can be interpreted
as a real-valued quadratic function 
$q_a$ on fibres of $B$ as follows:
for any $v \in \Tot(B)$, we define $q_a(v)$ as 
$\omega_B(v, a(v))$, where $\omega_B$ denotes the
symplectic form on fibres of $B$ associated 
with the Hermitian structure. 

\hfill

\remark \label{_curvature_on_total_space_Remark_}
Let $B$ be a holomorphic Hermitian vector bundle
on a complex manifold $M$, $\nabla$ its Chern connection\index[terms]{connection!Chern},
$\pi:\; \Tot B \arrow M$ its total space, 
and $\Theta_B\in \Lambda^{1,1}(M) \otimes {\goth u}(B)$
its curvature. 
Let $q:\; {\goth u}(B)\arrow C^\infty \Tot B$
be the map defined in 
\ref{_quadratic_from_Herm_Definition_}.
Then  $q(\Theta_B)$ is a section of 
$\pi^*(\Lambda^{1,1}(M))$.
Consider the decomposition
$T\Tot(B)= T_\hor \Tot(B) \oplus T_\pi \Tot(B)$
defined in \ref{_Ehre_holo_defined_Remark_}.
Identifying $T_\hor \Tot(B)$ and $\pi^*T M$,
we interpret $q(\Theta_B)$ as a 
section of $\Lambda^{1,1}(T_\hor^* \Tot(B))\subset\Lambda^{1,1}(\Tot B)$.

\hfill

\remark\label{_fibrewise_Hermitian_Definition-remark_}
Let $B$ be a holomorphic Hermitian vector bundle.
Then $\Tot B$ is equipped with a canonical
Ehresmann connection \index[terms]{connection!Ehresmann}(\ref{_Ehre_holo_defined_Remark_}).
Consider the corresponding decomposition:
\[
\Lambda^2(\Tot B)=\Lambda^2_\pi(\Tot B)\oplus 
\Lambda^1_\pi(\Tot B)\otimes\Lambda^1_\hor(\Tot B) 
\oplus \Lambda^2_\hor(\Tot B) 
\]
This allows us to realize $\Lambda^2_\pi(\Tot B)$ as a
sub-bundle in $\Lambda^2(\Tot B)$. Let $v\in \Tot B$.
Using the identification 
$T_\pi\Tot B\restrict v= B\restrict{\pi(v)}$,
we consider $T_\pi \Tot B$ as a Hermitian vector bundle.
The corresponding Hermitian form gives a section $\sigma$ of
$\Lambda^2_\pi(\Tot B)\subset \Lambda^2(\Tot B)$.
The fibrewise Hermitian form\index[terms]{form!Hermitian} on $\Tot B$ 
is the form $\sigma \in \Lambda^2(\Tot B)$ obtained in this way.

\hfill

Now we can prove a foundational formula,
known as ``Calabi formula (2.6)'' or 
``formula (15.19) from \cite{besse}.''\index[terms]{Calabi formula (2.6)}

\hfill

\theorem\label{_formula_15_19_besse_Theorem_}
In the assumptions of \ref{_curvature_on_total_space_Remark_},
let $l(v):=|v|^2$ be the length function,  
 $q(\Theta_B)\in \Lambda^{1,1}(\Tot B)$ the 2-form
obtained from the curvature as in 
\ref{_curvature_on_total_space_Remark_}, and
\[ \omega_\pi\in \Lambda^{1,1}_\pi (\Tot B)\subset \Lambda^{1,1}(\Tot B)
\] 
the fibrewise Hermitian form\index[terms]{form!Hermitian} on $\Tot B$,
obtained in \ref{_fibrewise_Hermitian_Definition-remark_}.
Then 
\begin{equation}\label{_formula_15_19_besse_}
dd^c l=-q(\Theta_B)+ \omega_\pi.
\end{equation}

\pstep
Let $(dd^c l)_\hor$, $(dd^c l)_\mix$ and 
$(dd^c l)_\ver$ be the three components of $dd^cl$
associated with the decomposition
\[
\Lambda^2(\Tot B)=\Lambda^2_\pi(\Tot B)\oplus 
\Lambda^1_\pi(\Tot B)\otimes\Lambda^1_\hor(\Tot B) 
\oplus \Lambda^2_\hor(\Tot B) 
\]
To prove \ref{_formula_15_19_besse_Theorem_},
we need to show that $(dd^c l)_\mix=0$,
$(dd^c l)_\ver=\sigma$ and $(dd^c l)_\hor=-q(\Theta_B)$.
The first equation 
is the easiest to show, because the
Ehresmann connection\index[terms]{connection!Ehresmann} on $\Tot B$ was defined
in such a way that $T_\hor \Tot B$ is the
orthogonal complement of $T_\pi \Tot B$ 
with respect to $dd^c l$. This gives $(dd^c l)_\mix=0$.

The equation $(dd^c l)_\ver=\sigma$
is also clear, because $l=|v|^2$, and hence  on fibrewise 
tangent vectors $dd^c l$ is equal to the standard
Hermitian form.\index[terms]{form!Hermitian}

\hfill

{\bf Step 2:}
It remains to show that  $(dd^c l)_\hor=-q(\Theta_B)$.
Let $X, Y\in T_\hor \Tot B$ be horizontal tangent vector fields.
By \ref{_Ehresmann_Chern_Theorem_}, Step 2, the Lie
derivative of $l$ along any horizontal vector field vanishes. Therefore,
the Cartan formula gives 
\begin{equation*}
\begin{split}
dd^c(l)(X,Y)&= -(d^cl)([X,Y])+ \Lie_X(d^cl(Y))-\Lie_Y(d^cl(X))\\
&= -(d^cl)([X,Y]) - \Lie_X(dl(I(Y)))+\Lie_Y(dl(I(X)))\\
&=-(d^cl)([X,Y])- \Lie_X\Lie_{IY}(l) + \Lie_Y\Lie_{IX}(l).
\end{split}
\end{equation*}
The last two terms vanish because $IX$ and $IY$ are
horizontal vector fields, giving $\Lie_{IY}(l) =\Lie_{IX}(l)=0$.
This gives $dd^c(l)(X,Y)=-(d^cl)([X,Y])$.
Applying \ref{_Ehresmann_Chern_Theorem_}, Step 2 again, we
obtain that the 1-form $d^cl$ vanishes on horizontal vector fields.
Therefore, $(d^cl)([X,Y])$ depends only on the vertical
part $[X,Y]_\ver$ of the commutator $[X,Y]$.
However, by \ref{_Ehresmann_curvature_bundles_Proposition_},
this vertical part\footnote{In the notation used in 
\ref{_Ehresmann_curvature_bundles_Proposition_}, the vertical
part of the commutator is $\Psi(X,Y)$.}  
 is expressed through the
usual curvature of $B$ as 
$[X,Y]_\ver\restrict b= \Theta_B(X,Y)(b)$,
where $b\in \Tot B$ is a point, and 
$\Theta_B(X,Y)(b)$ is understood as a
vector in $B\restrict{\pi(b)}= T_\pi \Tot B\restrict b$.

To finish the proof of $(dd^c l)_\hor=q(\Theta_B)$,
it remains to show that 
\begin{equation}\label{_curvature_tautology_Equation_}
(d^cl)(\Theta_B(X,Y)(b))=q(\Theta_B)(X,Y)\restrict b.
\end{equation}
Let $F$ be a fibre of $B$, and $\alpha$ the
linear vector field on $F$ associated with
the linear operator $\Theta_B(X,Y)\restrict F\in \End F$ as in 
\ref{_Ehresmann_curvature_bundles_Proposition_}.
Then \eqref{_curvature_tautology_Equation_}
is rewritten as 
\begin{equation}\label{_curvature_vector_space_tautology_Equation_}
d^cl(\alpha)(b)= \omega_B(\alpha, b),
\end{equation}
where $b\in F$ is a point, and $\omega_B$ the fibrewise
symplectic form on $\Tot B$ restricted to $F$.
Since $\omega_B= dd^c l$, and $b$ can be interpreted as
the radial vector field $\vec r$,
we have $\omega_B(\alpha, b) = -dd^cl(\alpha, \vec r)$.
On the other hand, $d^c l = \omega_B(\vec r, \cdot)$, and hence 
\eqref{_curvature_vector_space_tautology_Equation_}
becomes
$\omega_B(\vec r, \alpha)= - dd^cl(\alpha, \vec r)$,
that follows from $dd^c l=\omega_B$.
\endproof

\hfill

\corollary\label{_curvature_line_bundle_via_15.19_Corollary_}
Let $L$ be a holomorphic Hermitian line bundle
over $M$, $\Tot^\circ(L)$ the total
space of non-zero vectors in $L$,
and $l(v):=|v|^2$ the length function.
Then the form $dd^c\log l$ is basic
with respect to the fibration
$\Tot^\circ(L)\arrow M$, giving
$dd^c\log l= \pi^* \eta$.
Moreover, $\eta$ is equal to $-\Theta_L$,
where $\Theta_L$ is the curvature of the  Chern 
connection on $L$.

\hfill

\proof 
Clearly,
\[
dd^c \log l = \frac{dd^c l}{l} - \frac {dl}{l} \wedge \frac {I(dl)}{l}.
\]
The form $\frac {dl}{l} \wedge \frac {I(dl)}{l}$
is equal to $\frac{dd^c l}{l}$ in the fibres of $\pi$
(in the usual coordinates on $\C^*$, both of these forms
are equal to $\frac{dx\wedge dy}{x^2+ y^2}$).
Therefore, the form $dd^c \log l = \frac{dd^c l}{l}$
vanishes on the fibres of $\pi$. Since it is closed,
it is basic (\ref{_basic_for_closed-Remark_}).
It is equal to $-\Theta_L$ by 
\ref{_formula_15_19_besse_Theorem_}.
\endproof



\section{Positive line bundles}\label{posline}


Let $L$ be a holomorphic Hermitian line bundle,
and $\Theta_L\in \Lambda^{1,1}M$ the curvature of its
 Chern connection. \index[terms]{connection!Chern}

\hfill

\definition 
	A Hermitian holomorphic line bundle $L$ is called {\bf positive}
	if $\1 \Theta_L$ is a strictly positive form,
and {\bf negative} if it is dual to positive, that is,
if $-\1 \Theta_L$ is a strictly positive form.
\index[terms]{bundle!vector bundle!holomorphic Hermitian!positive}
\index[terms]{bundle!vector bundle!holomorphic Hermitian!negative}

\hfill

\remark  In this case, $M$  is forced to be K\"ahler: 
	indeed, the curvature form\index[terms]{form!curvature} $\Theta$ of 
a positive line bundle is identified
        with a $(1,1)$ form on $M$ and, being equal to
        $dd^c\log|h|$, it is closed. When $M$ is compact,
it is also projective by 
 Kodaira's embedding theorem
(\ref{_Kodaira_embedding_Theorem_}).\index[terms]{theorem!Kodaira's embedding}

\hfill

The next two statements will be used later on in this book.

\hfill

\claim 
Let $L$ be a holomorphic Hermitian line bundle on a
complex manifold $M$, 
and $l$ a function on its total space $\Tot(L)$ defined 
by $l(v):= |v|^2$. Then
{$dd^c\log l= -\1 \pi^* \Theta_L$,
		where $\Theta_L$ is the curvature of $L$.} 

\hfill

\proof From \ref{_formula_15_19_besse_Theorem_},
it follows that $dd^c l= -\1 |v|^2 \pi^* \Theta_L +\omega_\pi,$
where $\omega_\pi$ is the fibrewise Hermitian form \index[terms]{form!Hermitian}
(\ref{_fibrewise_Hermitian_Definition-remark_}).
Therefore, 
\[
dd^c\log l= \frac{-\1 |v|^2 \pi^* \Theta_L +\omega_\pi}{l}
- \frac{dl \wedge d^c l}{l^2}
\]
(Exercise \ref{_dd^c_log_Exercise_}). However,
$\omega_\pi = \frac{dl \wedge d^c l}{l}$
because $B$ is rank 1 bundle, and on a complex
line with holomorphic coordinate $z$ we have $l=z\bar z$ and
the standard Hermitian form $\omega_\C$ on $\C$ can be obtained as
\[
\frac{dl \wedge d^c l}{l}= 
- \frac{\6 l \wedge \bar \6 l}{-2\1 l}=
\frac{\bar z dz\wedge z d\bar z}{-2\1z\bar z}=
\frac{dz\wedge d\bar z}{-2\1}=\omega_\C. \quad \endproof
\]

\hfill

\remark\label{_signs_lost_Remark_}
It is easy to get lost in the signs.
However, there is a rule of thumb that leads to 
choosing signs correctly. A very ample line bundle $L$\index[terms]{bundle!line!ample}
over a compact manifold is generated by holomorphic
sections, and hence   a plurisubharmonic function on $\Tot(L)$ has 
to be constant. On the other hand, a holomorphic section
$\xi$ defines a holomorphic function $\sigma(\xi)$ on 
$\Tot(L^*)$ mapping  $v \in L^*$ to $v(\xi)$. Taking a sum
$\sum_i |\sigma(\xi_i)|^2$ for a set of
holomorphic sections $\xi_i$ generating $L$, 
we obtain a function on $\Tot(L^*)$ that is  strictly 
plurisubharmonic outside of  the zero section.
Moreover, by  Remmert's theorem (\ref{rem_red}), there
is a bimeromorphic, holomorphic map from
$\Tot(L^*)$ to a singular Stein variety,
blowing down the zero section to a point.
In this way one obtains the cone of $L^*$
(complete with the origin) as a complex variety.

\hfill

\corollary\label{_dd^c_log_l_Corollary_}
Let $L$ be a holomorphic Hermitian line bundle on a
complex manifold $M$, and $\Tot^\circ(L)\subset \Tot(L)$ the set of
non-zero vectors in $\Tot(L)$. Suppose that $L$ is negative.
Then the form $dd^c\log l$ is positive on
$\Tot^\circ(L)$. Moreover,  $dd^c\log l$  has one 
	zero eigenvalue along the fibres of $L$
and the remaining eigenvalues are 
(strictly) positive. 
\endproof

\hfill

Existence of positive line bundles is very restrictive:

\hfill

\theorem  ({\bf \index[persons]{Kodaira, K.} Kodaira}, \cite{griha})
\label{_Kodaira_embedding_Theorem_} 
	Let $M$ be a compact complex manifold. Then {$M$ is
		projective if and only if it admits a positive line bundle.}\index[terms]{theorem!Kodaira's embedding}\index[terms]{manifold!projective}
\endproof

\section{Exercises}

\begin{enumerate}[label=\textbf{\thechapter.\arabic*}.,ref=\thechapter.\arabic{enumi}]

\item \label{_dd^c_log_Exercise_}
Let $\psi$ be a $\R^{>0}$-valued function on a complex manifold.
Prove that 
$$dd^c (\log \psi)= \frac {dd^c \psi}{\psi} - \frac{d\psi \wedge d^c \psi}{\psi^2}.$$

\item
Let  $l(z_1, ..., z_n):=\sum |z_i|^2$.
Prove that $dd^c \log l$ is semi-positive.
Prove that for $n >1$, the form $dd^c \log l$ has at least one positive
eigenvalue.

\item
Consider the tautological fibration
$\C^{n+1}\backslash 0 \stackrel \pi \arrow \C P^n$.
We consider $\pi$ as a quotient map, $\C P^n=(\C^{n+1}\backslash 0)/\C^*$,
and take $r= \sum_{i=1}^{n+1} z_i \frac \6{\6z_i}$.
 Prove that any $\C^*$-invariant form $\eta$ such that 
$i_r\eta = i_{\bar r}\eta =0$ is the pullback of a
form on $\C P^n$.  Prove that $dd^c \log l$ is the pullback of a
form $\omega$ on $\C P^n$.

\item
Recall that a {\bf Fubini--Study form} on $\C P^n$ is an
$\U(n+1)$-invariant Hermitian form.\index[terms]{form!Fubini--Study}
Consider the form $\omega\in \Lambda^{1,1}(\C P^n)$ 
defined in the previous exercise.
Prove that this form is $\U(n+1)$-invariant. 
Prove that all $\U(n+1)$-invariant
2-forms on $\C P^n$ are proportional.

\item
Let $\C^n \subset \C P^n$ be an affine chart with affine
coordinates $z_1, ..., z_n$. Prove that the Fubini--Study form
on this chart is given by 
\[ \omega= 
\frac{\sum_{i=1}^n dz_i \wedge d\bar z_i}{1+ \sum_{i=1}^n|z_i|^2} -
\frac{\sum_{i=1}^n\bar z_i dz_i}{1+ \sum_{i=1}^n|z_i|^2}\wedge 
\frac{\sum_{i=1}^nz_i d\bar z_i}{1+ \sum_{i=1}^n|z_i|^2}.
\]

\item
Let $f_1, ..., f_n$ be holomorphic functions without common zeroes 
on   a complex manifold. Prove that $\log(\sum_i |f_i|)$ is
a plurisubharmonic function.

\item
Let $L$ be a Hermitian holomorphic bundle on 
a complex manifold $M$, and $\Tot L$ its total space.
Denote by $\Tot^\circ(L)\subset \Tot L$ the set of non-zero vectors in 
$\Tot L$. Clearly, $\Tot^\circ(L)$ is equipped with a free action of $\C^*$, and
$\Tot^\circ(L)/\C^* = M$.  Consider the function 
$l:\; \Tot^\circ(L) \arrow \R^{>0}$
mapping a vector $v\in L\restrict x$ to $|v|^2$.
Prove that $dd^c \log l$ is the pullback of a form on $M$.
Prove that the corresponding (1,1)-form on $M$ 
is equal to $-\1\Theta_L$, where $\Theta_L$
is the curvature of the Chern connection \index[terms]{connection!Chern}on $L$.

\item Let $(M,\omega,\theta)$ be an LCK manifold, and 
$\nabla^C$ the  Chern connection on its tangent bundle. 
Show that  $\nabla^C$ is related to the  Weyl connection (see \ref{_Weyl_conn_definition_Remark_}) $D$  by the formula:
$$\nabla^C = D +\frac 12 \left(\theta\otimes \id +I\theta\otimes I \right).$$
Find the relation between the curvature tensors of the Chern and the Weyl connections.
\end{enumerate}

\chapter{CR, Contact and Sasakian manifolds} \label{CR_and_such}\index[terms]{manifold!Sasaki--Einstein}

\epigraph{\it The problems are solved, not by giving new information, but by arranging what we have known since long.}{\sc\scriptsize 
	Ludwig Wittgenstein, Philosophical Investigations}

\section{Introduction}

This chapter is a survey and introduction
to geometric structures on odd-di\-men\-si\-o\-nal
manifolds, that are  intimately related to 
symplectic, complex and locally conformally
K\"ahler geometries. 

The contact structures on a manifold
go back to Sophus Lie\index[persons]{Lie, S.} who defined the
notion of contact transformation and studied
the structure of the group of local contact
transformations. We define a contact structure
on a manifold $M$ as a symplectic structure 
on its cone $C(M)=M \times \R^{>0}$, such that
the shift $(m, t) \mapsto (m, \lambda t)$
multiplies the symplectic form by $\lambda^2$
for any $\lambda \in \R^{>0}$.

The notion of contact
manifold was first properly studied
by \index[persons]{Boothby,  W. M.} Boothby and \index[persons]{Wang, H. C.} Wang in \cite{_Boothby_Wang_}, see also
\cite{_Gray:contact_} and \cite{_Reeb_}. At the time, the definition
of contact structures was more hands-on, in fact,
in most cases a contact manifold was defined in
terms of the contact form $\eta$, which was fixed,
or a bunch of local contact forms with gluing conditions
(\cite{_Gray:contact_}).
In this context, the notion of contact diffeomorphism
becomes very complicated.
Later, people realized that a contact structure
is uniquely determined by the distribution $\ker \eta$,
which simplifies the local computations.
The ``symplectization'' $C(M)$ acquired a definite
treatment in the appendix to \index[persons]{Arnol'd, V. I.} Arnol'd's ``Mathematical 
methods of classical mechanics'' (\cite{_Arnold_}). Arnol'd did not use
it as a definition, but he used symplectization to prove
the contact Darboux theorem.\index[terms]{symplectization}

Our treatment of contact geometry\index[terms]{geometry!contact} is influenced
by a more recent notion -- that of Sasakian manifold.
In modern terms, a Sasakian structure on $M$
is a K\"ahler structure on $C(M)= M \times \R^{>0}$
such that the shift $(m, t) \mapsto (m, \lambda t)$
is a holomorphic map that multiplies the K\"ahler form\index[terms]{form!K\"ahler} by $\lambda^2$
for any $\lambda \in \R^{>0}$.

In 1960-1961 \index[persons]{Sasaki, S.}  Sh. Sasaki
and Yo. \index[persons]{Hatakeyama, Y.} Hatakeyama (\cite{_Sasaki_}, \cite{_Sasaki_Hatakeyama_}) considered a class of almost contact structures compatible with a metric that they called  
$(\phi, \xi, \eta,g)$-structure. 
Sasaki and \index[persons]{Hatakeyama, Y.} Hatakeyama defined  $(\phi, \xi, \eta,g)$-structure on $M$
as a complicated collection of tensors on $M$. Without using this term, they defined
the symplectization $C(M)$ (they used the cylinder
instead of the cone, but this is essentially a matter
of convention). Then they checked the vanishing
of the Nijenhuis tensor on $C(M)$. In these two 
papers, they did not mention that $C(M)$ is K\"ahler. Soon after, \index[persons]{Sasaki, S.} Sasaki added the whole structure that we now know under his name (see Subsection \ref{_Sasaki:definition_}). 

The modern approach to Sasakian geometry,\index[terms]{geometry!Sasaki}
where it is defined in terms of the K\"ahler
structure on the cone $C(M)$, is due to 
C. P. Boyer \index[persons]{Boyer, C. P.}and K. \index[persons]{Galicki, K.} Galicki, who spent the 1990s
writing a series of papers on the 
subject, culminating in the book
\cite{bog}, which remains the main reference on the subject.

The Sasakian structures are related to the
contact structures in the same way as 
K\"ahler structures are related to the
symplectic ones. In fact, the links
of isolated singularities defined by
homogeneous polynomials are always
Sasakian. This allows one to study
the Sasakian manifolds in
terms of algebraic geometry\index[terms]{geometry!algebraic}. The algebraic analogue
of the Sasakian structure is called the Seifert bundle\index[terms]{bundle!Seifert}
(\cite{_Kollar:Seifert_}). However, the Seifert bundles
correspond to Sasakian manifolds that are 
quasi-regular (\ref{_regular_quasi_regular_sasakian_}), 
avoiding the more complicated
(and, arguably, more interesting) case of irregular
Sasakian manifolds.\index[terms]{manifold!Sasaki!regular}\index[terms]{manifold!Sasaki!quasi-regular}

The notion of CR-manifold\index[terms]{manifold!CR} is one of the basic subjects
of complex analysis. We give an account of its history and
motivation in Subsection \ref{_CR_history_Subsection_}.
Similar to the contact structures, a CR-structure\index[terms]{structure!CR} is
defined by fixing a sub-bundle $B\subset TM$; however,
we also fix a complex structure $I:\; B \arrow B$,
$I^2= - \Id_B$. A triple $(M, B,I)$ defines a
CR-structure if the sub-bundle $B^{1,0}\subset B\otimes_\R\C$ 
of $(1,0)$-vectors in $B$ is integrable, that is, it satisfies
$[B^{1,0}, B^{1,0}]\subset B^{1,0}$.

Historically, this structure was observed on
real hypersurfaces $S\subset M$ in complex manifolds. In this
case, $B$ is identified with the intersection $TS\cap I(TS)$.
The Frobenius form $\Lambda^2 B \arrow TS/B$ \index[terms]{form!Frobenius}
(\ref{_Frobenius_form_Definition_})
takes values in the real line bundle $TS/B$, and has 
Hodge type (1,1). There are many situations when 
the contact manifold $S$ is studied by identifying 
$S$ with the boundary of a Stein manifold (the so-called
``Stein filling'', \cite{_Eliashberg:filling_}).\index[terms]{Stein filling} 
In this case, the Frobenius form is
non-degenerate, and the CR-structure determines the contact structure.\index[terms]{structure!CR}\index[terms]{manifold!CR}
However, the notion of CR-structure is more rigid; indeed,
for any CR-manifold $(M,B,I)$, the bundle $B$ is equipped
with $TS/B$-valued Hermitian form\index[terms]{form!Hermitian}. In particular, there
is no hope of obtaining a normal form of CR-manifolds
or an analogue of Darboux theorem.\index[terms]{theorem!Darboux}

The notion of CR-manifold is by its nature intermediate
between the contact and Sasakian manifolds. Indeed, 
for any Sasakian manifold $S$, we realize $S$ as a 
real hypersurface in its cone $C(S)$, that is  
a complex (and even K\"ahler) manifold. The 
contact distribution on $S$ can be obtained by
taking $TS\cap I(TS)$, as in the case of CR 
hypersurfaces. In \ref{_Reeb_fie_from_CR_Theorem_}, we  characterize the
Sasakian manifolds as CR-manifolds equipped
with a certain one-parametric group of CR-automorphisms.

For the purposes of the present book, Sasakian
manifolds have utmost importance as a source of 
locally conformally K\"ahler manifolds. Indeed,
given a compact Sasakian manifold $S$,
its cone $C(S)= S \times \R^{>0}$ is equipped with a natural
action by holomorphic homotheties
$(s, t)\mapsto (s, \lambda t)$.
Taking a discrete subgroup $\Gamma=\Z \subset \R^{>0}$,
we obtain a compact quotient
$C(S)/\Gamma$, that is  by construction locally conformally
K\"ahler. It turns out that in this way we can construct a large
class of LCK manifolds (\ref{str_vai}).\index[terms]{manifold!LCK} 

The cone-type metrics are ubiquitous in Riemannian geometry.
Similar to \index[persons]{Berger, M.} Berger's classification of holonomies of Riemannian
manifolds, one can classify the manifolds with Riemannian cones \index[terms]{cone!Riemannian}
having prescribed holonomy.
This gives an interesting menagerie of
special manifolds: Sasakian, 3-Sasakian, nearly K\"ahler and
nearly parallel $G_2$ (Chapter \ref{other}).\index[terms]{manifold!nearly-K\"ahler}\index[terms]{manifold!Sasaki}\index[terms]{manifold!3-Sasaki}\index[terms]{manifold!nearly parallel $G_2$}
This correspondence was introduced
by C. B\"ar\index[persons]{B\"ar, C.} who used it to classify the Killing spinors on\index[terms]{spinor!Killing}
Riemannian manifolds; it is now known as {\bf B\"ar correspondence}
(\cite{_Bar:killing_,_BG:cones_}).\index[terms]{correspondence!B\"ar}

This chapter is based on the material introduced in the
previous chapters; however, the Sasakian and contact 
geometry are not the central subject of this book.
To fully comprehend these topics, the reader is encouraged to 
look elsewhere. 

Contact geometry\index[terms]{geometry!contact} is one of the central subjects of 
symplectic geometry;\index[terms]{geometry!symplectic} there is a comprehensive introduction
in \cite{_Cannas_da_Silva_} and \cite{_McDuff_Salamon_}.
 For an introduction in CR geometry\index[terms]{geometry!CR}, see 
\cite{dt} and \cite{_Jacobowitz_}, whereas for Sasakian geometry,\index[terms]{geometry!Sasaki}
the definitive reference is \cite{bog}.


\section{CR-manifolds}


The natural generalization of a complex structure on manifolds of  
odd dimensions is the  notion of CR-structure, which
was first defined on boundaries of domains in $\C^n$.

\hfill

\definition  \label{_CR_mfld_Definition_}
	Let $M$ be a smooth manifold,
$B\subset TM$ a sub-bundle in the tangent bundle,
and $I:\; B \arrow B$ an endomorphism satisfying
$I^2=-1$. Consider its $\1$-eigenspace $B^{1,0}(M)\subset
B\otimes \C \subset T_C M=TM\otimes \C$.
Suppose that $[B^{1,0}, B^{1,0}]\subset B^{1,0}$.
Then $(B,I)$ is called a {\bf CR-structure on $M$}.\index[terms]{structure!CR}

\hfill

\example  
\begin{itemize}
\item[1.] A complex manifold
is CR, with $B=TM$: by  \ref{nn}, the inclusion {$[T^{1,0}M, T^{1,0}M]\subset T^{1,0}M$	is equivalent to the integrability of the complex structure}.

\item[2.] Let $X$ be a complex manifold,
and $M\subset X$ a hypersurface. Then\index[terms]{hypersurface} 
$B:= TM \cap I(TM)$, and $\rk B=\dim_\C (TM \cap I(TM))=\dim_\C X-1=n-1$, {$M$ is a CR-manifold.}\index[terms]{manifold!CR} since $[T^{1,0}X,T^{1,0}X]\subset T^{1,0}X$.

\end{itemize}

%
%

The notions of holomorphic map and holomorphic vector field are naturally extended to the CR context:

\hfill

\definition \label{_CR_holomorphic_Definition_}
Let $(M,B,I)$ be a CR-manifold. A function $f$ on $M$ is
called {\bf CR-holomorphic} if for any vector field
$v\in B^{0,1}$, we have $\Lie_v f=0$. \index[terms]{function!CR-holomorphic}
A vector field
$v\in TM$ is called {\bf CR-holomorphic} if
the corresponding diffeomorphism flow preserves 
$B$ and $I$.\index[terms]{vector field!CR-holomorphic}

\section[Contact manifolds and pseudoconvex CR-ma\-ni\-folds]{Contact manifolds and pseudoconvex\\ CR-ma\-ni\-folds}

\subsection{Contact manifolds and symplectic cones}

The link between contact and symplectic geometry\index[terms]{geometry!symplectic}\index[terms]{geometry!contact} is made
through {\bf symplectic cones}.\index[terms]{cone!symplectic}

\hfill

\definition  Let $S$ be a manifold. Then $C(S):= S
  \times \R^{>0}$ is called {\bf the  cone over $S$}. The
  multiplicative group $\R^{>0}$ acts on $C(S)$ by
  $h_\la(x, t) \mapsto (x, \lambda t)$.

\hfill

We can now define the contact manifolds.

\hfill

\definition 
 A {\bf contact manifold} is \index[terms]{manifold!contact}
 a manifold $S$ such that its cone $C(S)$ is endowed with
 a symplectic form\index[terms]{form!symplectic} $\omega$ on
 which $\R^{>0}$ acts through
 $h_\la^*\omega=\la^2\omega$. 
A contact structure on $S$ is a symplectic form
on $C(S)$ that satisfies this relation.
A contactomorphism (morphism of contact manifolds)
is a diffeomorphism of $S$ that induces
a symplectomorphism of $C(S)$.

\hfill

\example \label{contact_sphere}
An odd-dimensional sphere $S^{2n-1}$ is contact.
Indeed, its cone is $\R^{2n}\backslash0$, with the
standard symplectic form $\omega:=\sum_i dp_i \wedge dq_i$
satisfying  $h_\la^*\omega=\la^2\omega$.
 
 \hfill
 
 One has the following characterizations.
 
\hfill

 \theorem  Let $S$ be a differentiable manifold. The
 following are equivalent:
\begin{description}
\item[(i)]  $S$ is contact.
\item[(ii)]   $S$ is odd-dimensional and there exists an oriented sub-bundle of codimension 1,
  $B\subset TS$,  
 	with non-degenerate Frobenius form $\Phi\in \Lambda^2 (B,TS/B)$. \index[terms]{form!Frobenius}
 	 This $B$ is called  {\bf the contact bundle}.
\item[(iii)] 
 $S$ is odd-dimensional and there exists an oriented sub-bundle of codimension 1,
 	$B\subset TS$, such that for any nowhere degenerate 1-form  $\eta\in
 	\Lambda^1 S$ that annihilates $B$, the form $\eta \wedge (d\eta)^k$ 
 	is a non-degenerate volume form (where $\dim S=2k+1$). Every such $\eta$ 
 	is called {\bf a contact form}.\index[terms]{form!contact}
\end{description} 
	
\remark The contact form $\eta$ is not unique. Moreover,
contactomorphisms of a manifold $S$ do not necessarily
preserve $\eta$.

\hfill

\pstep {\bf Equivalence of (ii) and (iii).} 
Since for a 1-form $\al$ one has 
$$d\al(X,Y)=X(\al(Y))-Y(\al(X))-\al([X,Y]),$$
for each $X, Y \in B$, $d\eta(X,Y)= \eta([X,Y])= \Phi(X,Y)$.
Therefore, the Frobenius form \index[terms]{form!Frobenius}$\Lambda^2 B \stackrel \Phi \arrow TS/B$
can be expressed as $\langle \Phi(X,Y), \eta\rangle= d\eta(X,Y)$.
On the other hand, the non-degeneracy of $\eta \wedge (d\eta)^k$ on $TS$
	is equivalent to the non-degeneracy of $d\eta=\Phi$ on $B=\ker \eta$.
Therefore, $\langle \Phi(X,Y), \eta\rangle= d\eta(X,Y)$
is of maximal rank if and only if $\eta \wedge (d\eta)^k$
is non-degenerate.

\hfill

{\bf Step 2:  (iii) implies (i).} Let $M\stackrel \pi\arrow S$ be the space of positive vectors in the oriented
1-dimensional bundle $L:=TS/B$, that is  trivialized by
the form $\eta$, $V\in T M$  a 
unit vertical vector field, and $t:\; M \arrow \R$
a map that associates $\eta(v)$ to a point $(s,v)\in M$,
$s\in S, v \in L\restrict s$. Let
$\tau:=t\pi^*\eta\in \Lambda^1 M$, and let $\omega:=d\tau$.
Consider the radial vector field $r=tV\in TM$. Clearly,\index[terms]{vector field!radial}
$\Lie_r \tau=2\tau$, giving $\Lie_r d\tau=2d\tau$.
{To prove that $M$ is a symplectic cone of $S$,
	it remains to show that $d\tau$ is non-degenerate.}

Indeed, since $\ker dt=\pi^{-1}(S)$, 
any vector field $x\in TS$ can be naturally lifted to a vector
field $\pi^*(x)\in \ker dt\subset TM$. Now let  $x,y\in B$ and $X=\pi^*(x)$, 
$Y:= \pi^*(y)$. Then one has $d\tau(X,Y)= -\tau([X,Y])=-\tau(\pi^*([x,y]))$.
Since the Frobenius form is non-degenerate,\index[terms]{form!Frobenius}
the form  $d\tau\restrict {\pi^*B}=-\tau(\pi^*([x,y]))$ 
is non-degenerate on $\pi^*B$. Also, the contraction
$i_V d\tau= \tau$, and $\ker \tau=  \pi^*B\oplus\langle V\rangle$; 
hence {$d\tau$ is non-degenerate on the symplectic
	orthogonal complement of $\pi^*B$, 
generated by $V$ and $\ker \tau \restrict {\pi^* B}$.}

\hfill

{\bf Step 3:  (i) implies (iii).}  Let $M=C(S)=S\times \R^{>0}$, and $t\in  C^\infty( M)$ the 
standard coordinate along $\R^{>0}$. Consider the vector field
$r:= t\frac d {dt}$, and the form $\eta:= i_r \omega$
obtained by contraction with $r$.
Since $i_r \eta=0$ and $d\eta=\Lie_r\omega=\omega$ by the 
Cartan formula and $h_\la^*\omega=\la^2\omega$, we have:
\[ \Lie_{r}(t^{-1}\eta)=i_r d(t^{-1}\eta)  + d(i_r t^{-1}\eta) 
= t^{-1}\eta-t^{-1} \eta=0,
\]
and hence the form $t^{-1}\eta$ is the pullback of a form on $S$, with
$t^{-1}\eta=\pi^* \eta_0$, under the projection
	$\pi:\; C(S)\arrow S$.  This gives a form $\eta_0$ on $S$.
Moreover, $(d\eta)^{k+1}$
is non-degenerate because $d\eta$ is symplectic.  Therefore,
	$i_r (d\eta)^{k+1}= (k+1)(d\eta)^{k}\wedge \eta$
	is non-degenerate on $S=S\times\{1\}$. This implies that $(d\eta_0)^k\wedge\eta_0$ is non-degenerate on $S$.
 \endproof
 
\hfill

\remark 
Let $S$ be a contact CR-manifold, with  the  bundle $B$
contact {\em and}\  CR, and a complex structure
operator $I$ on $B$. Since the Frobenius form vanishes \index[terms]{form!Frobenius}
when both arguments are from  $B^{0,1}$ and $B^{1,0}$, it is a pairing between $B^{0,1}$ and $B^{1,0}$. Indeed, $[B^{1,0}, B^{1,0}]\subset B^{1,0}$\index[terms]{pairing}
and $[B^{0,1}, B^{0,1}]\subset B^{0,1}$. This proves that the 
Frobenius form is of Hodge type (1,1), and hence  Hermitian.
\index[terms]{form!Frobenius}\index[terms]{form!Hermitian} 	
 	
\subsection{Levi form and pseudoconvexity}

\definition \label{_Levi_form_Definition_}
	The Frobenius form\index[terms]{form!Frobenius} of a CR contact 
manifold, considered to be  a Hermitian form, \index[terms]{form!Hermitian}
is called {\bf Levi form}.
\index[terms]{form!Levi}

\hfill

CR-manifolds of hypersurface type (meaning that
$\codim B=1$) with definite Levi form are closely related
to K\"ahler manifolds.\index[terms]{manifold!CR!of hypersurface type}

\hfill

\definition 
	Let $(M,B,I)$ be a CR-manifold,
	with $\codim B =1$. Then $M$ is called {\bf a strictly pseudoconvex
		CR-manifold} if its Levi form is positive definite. \index[terms]{manifold!CR!strictly pseudoconvex}

\hfill

Strictly pseudoconvex CR-manifolds usually appear as level sets
of strictly plu\-ri\-subharmonic functions. To prove this we need
the following preliminary calculation.

\hfill

\proposition \label{propos}
Let $M$ be a complex manifold,  $\phi\in  C^\infty( M)$ 
 a smooth function, and $s$  a regular value of $\phi$.
Consider $S:=\phi^{-1}(s)$ as a CR-manifold,
with $B=TS\cap I(TS)$,
and let $\Phi$ be its Levi form,\index[terms]{form!Levi} taking values
in 
\[ TS/B= \frac{\ker d\phi}{\ker d\phi \cap I(\ker d\phi)}.
\]
Then $d^c \phi:\; TS/B\arrow  C^\infty (S)$ trivializes $TS/B$.
Consider the tangent vectors $u, v \in B$.
{Then $-d^c \phi(\Phi(u,v))=dd^c\phi(u,v)$.}

\hfill

\proof
Indeed, extend $u, v$ to vector fields
$u, v\in B= \ker d\phi\cap I(\ker d\phi)$. 
Then 
$$-d^c\phi(\Phi(u,v))= -d^c\phi([u,v])=
dd^c\phi(u, v) - \Lie_v d^c\phi(u) + \Lie_u d^c\phi(v).$$
(the last equality follows from the Cartan formula).
However, $d^c\phi(u)=d^c\phi(v)=0$
because $v, u \in I(\ker d\phi)$.
This gives $-d^c\phi(\Phi(u,v))=dd^c\phi(u,v)$.
\endproof

\hfill

A general recipe for  constructing pseudoconvex CR-manifolds is as follows.

\hfill

\corollary\label{_psh_then_pseudoconvex_Corollary_}
Let $M$ be a complex manifold, $\phi\in  C^\infty (M)$ 
a strictly plurisubharmonic function, and $s$ a regular value of $\phi$.
{Then $S:=\phi^{-1}(s)$ is strictly pseudoconvex.}\index[terms]{function!plurisubharmonic!(strictly) plurisubharmonic}

\hfill

\proof By \ref{propos}, the Levi form\index[terms]{form!Levi}
of $S$ is expressed as $dd^c\phi(u, v)$, and hence  it is positive definite.
\endproof

\hfill

\definition 
	Given a contact manifold $S$ with a  contact form $\eta$,\index[terms]{form!contact} one associates a vector field, usually denoted $\xi$,  defined by the equations:
	$$i_\xi\eta=1,\qquad i_\xi d\eta=0.$$
The vector field $\xi$ is called {\bf the Reeb} or {\bf characteristic} field.\index[terms]{vector field!Reeb}

\subsection{Normal varieties}

Before we can proceed, we need to state a few definitions
of complex analytic geometry. We follow \cite{demailly}.

\hfill

\definition\label{_normality_Definition_}
Recall that a complex variety $X$ is called {\bf normal}
if any locally bounded meromorphic function on an open subset 
$U\subset X$ is holomorphic (\cite[Definition II.7.4]{demailly}). 
In algebraic geometry, a variety
is normal if all its local rings are integrally closed; these
two notions are equivalent for a complex variety obtained
from an algebraic one (\cite[Theorem II.7.3]{demailly},
\cite[Satz 4, p. 122]{_Kuhlmann_}).

\hfill

\remark
Each complex variety $A$ admits {\bf the normalization},
that is, a normal variety $B$ equipped with a finite
holomorphic map $n:\; B \arrow A$. The normalization
is unique up to an isomorphism. 
In the sequel, we will usually consider 
the normalization maps $n:\; B \arrow A$
that are  homeomorphisms. The concept may 
seem alien, but in fact it is very natural.
For example, consider the plane curve $C$
defined by the equation $x^2=y^3$. This
curve is clearly singular in zero.
It can be parametrized by a variable
$t$ with $t^2=y$ and $t^3=x$; this
parametrization defines a map $n:\; \C \arrow C$
that is  finite, and hence  $n$ is the
normalization of $C$. However, $n$ is clearly 
a homeomorphism.

\hfill

\claim\label{_quotient_normal_Claim_}
Let $G$ be a finite group holomorphically acting 
on a normal complex variety $Y$, and $X= Y/G$ the quotient 
variety. Then $X$ is normal.

\hfill

\proof Let $f$ be a locally bounded meromorphic function
on $X$. Its pullback $\tau^* f$ to $Y$ is clearly meromorphic and locally
bounded, and hence  holomorphic. Since $\tau^* f$ is $G$-invariant,
it is holomorphic on $Y$.
\endproof

\subsection{Stein completions and Rossi-Andreotti--Siu theorem}\index[terms]{completion!Stein}

We shall need the following deep result which characterizes the CR-pseudo-con\-vex manifolds in terms of their ring of CR
holomorphic functions and relates them to Stein manifolds
(see Section \ref{_Stein_Subsection_} below for the definition of and details on
Stein manifolds).\index[terms]{manifold!Stein}

\hfill

\theorem  ({\bf (\index[persons]{Rossi, H.}Rossi and \index[persons]{Andreotti, A.} Andreotti-\index[persons]{Siu, Y.-T.}Siu)}; \cite{andreotti_siu}, \cite{rossi})\label{asr} 
Let $S$ be a compact strictly pseudoconvex CR-manifold, \index[terms]{theorem!Rossi, Andreotti--Siu}
$\dim_\R S\geq 5$, and let $H^0(\calo_S)$ be the ring of
CR-holomorphic functions. {Then $S$ is the boundary of a
Stein variety $M$ with isolated singularities,} such that
$H^0(\calo_S) = H^0(\calo_M)$, where $H^0(\calo_M)$
denotes the ring of holomorphic functions on $M$, considered 
a compact complex manifold with boundary.

\hfill

\remark \label{_Forster_Stein_uniqueness_Remark_}
Let $A$ be a commutative Fr\'echet algebra over $\C$ and $\Spec(A)$\index[terms]{algebra!Fr\'echet}
the {\bf continuous spectrum} of $A$, defined as the\index[terms]{spectrum!continuous}
set of all continuous $\C$-linear homomorphisms $A\arrow \C$.
According to \cite{_Forster:1966_,_Forster:1967_}, $\Spec(H^0(\calo_X))=X$ for any Stein
variety $X$, where $H^0(\calo_X)$ is the algebra
of holomorphic functions equipped with the topology
of uniform convergence on compacts.
From \cite{_Forster:1967_} it
also follows that a Stein manifold is determined uniquely
from the Fr\'echet algebra $H^0(\calo_X)$.
This implies that the Stein variety obtained in
the  Rossi and  Andreotti--Siu theorem is unique.

\hfill

Forster's theorem, together with the Hartogs principle, implies that 
we can ``fill in'' a missing compact subset
of a Stein variety, and recover the Stein variety from the
complement $M\backslash K$; this procedure (usually
applied to the complement to a point) is called
``the Stein completion''.\index[terms]{completion!Stein} 

\hfill

\definition\label{_Stein_completion_Definition_}
Let $M$ be a Stein variety, $\dim M > 1$,
and $K\subset M$ a compact subset.
Following \cite{andreotti_siu},
we define {\bf the Stein completion} of $X:=M\backslash K$
as the continuous spectrum of the algebra 
$H^0(\calo_X)$, that is, $\Spec(H^0(\calo_X))=M$
equipped with the structure of a Stein variety.\index[terms]{variety!Stein}

\hfill

This notion is based on the following version of the Hartogs principle,
due to H. Rossi.\index[persons]{Rossi, H.}

\hfill

\theorem\label{_Hartogs_Stein_Theorem_}
Let $X$ be a normal Stein variety, $\dim_\C X >1$, 
and $K \subset X$ a compact subset. Then every holomorphic
function on $X \backslash K$ can be extended to $X$.

\proof \cite[Theorem 6.6]{_Rossi:fields_}. \endproof

\hfill

\remark
When $M$ is normal,\index[terms]{variety!normal} we can define the Stein completion \index[terms]{completion!Stein}as follows.
Let $M$ be a normal Stein variety, and $K\subset M$ a compact subset.
By \ref{_Hartogs_Stein_Theorem_}, the ring of functions on 
$M$ is identified with $H^0(\calo_{M \backslash K})$;
by \ref{_Forster_Stein_uniqueness_Remark_}, 
this ring with $C^0$ topology uniquely defines $M$.
Then $M$ is the Stein completion of
$M\backslash K$. By \ref{_Hartogs_Stein_Theorem_},
this definition is equivalent to \ref{_Stein_completion_Definition_}.

\hfill

\definition\label{_weak_Stein_Definition_}
{\bf A weak Stein completion}\index[terms]{completion!Stein!weak}
of $M \backslash K$ is any Stein variety
$M'$ containing a compact $K'$ such that 
$M'\backslash K'\cong M\backslash K$.

\hfill

\remark\label{_normalization_of_weak_Stein_comple_Remark_}
Clearly, if $X\backslash K$ is normal,\index[terms]{completion!Stein}
the normalization of a weak Stein completion\index[terms]{completion!Stein!weak}
of $X \backslash K$ is the Stein completion
of  $X \backslash K$. \index[terms]{normalization}

\hfill

\remark\label{_Stein_comple_unique_Forster_Remark_}
Given any normal complex variety $X$ and an open embedding
$X \hookrightarrow M$ to a Stein variety,
with $M \backslash X$ compact, we identify $M$
with the Stein completion\index[terms]{completion!Stein} of $X$. This can
be used as an alternative definition
of the Stein completion; by \ref{_Hartogs_Stein_Theorem_},
$M$ is uniquely determined by $X$.

\hfill

\corollary  
The Lie group $G:=\Aut(S)$ of CR-automorphisms
is identified with the group of complex automorphisms
of the corresponding Stein space $M$. {Its Lie algebra
 is the Lie
algebra of holomorphic vector fields on $M$.}

\section{Sasakian manifolds}\index[terms]{manifold!Sasaki}

Sasakian structures relate to contact ones as K\"ahler
structures relate to complex ones. Precisely:

\hfill

\definition
Let $S$ be a Riemannian manifold. 
The {\bf Riemannian cone}\index[terms]{cone!Riemannian} of $(S,g)$ is the
manifold $C(S):=\R^{>0}\times S$ with the metric
$dt^2+t^2g$.\index[terms]{cone!Riemannian}

\hfill

\example
Clearly, the flat space $\R^n \backslash 0$ can be
obtained as a Riemannian cone of the unit sphere
$S^{n-1}\subset \R^n$ with the
standard metric.  

\hfill

\definition \label{_Sasaki_definition_}  
A {\bf Sasakian structure} on a Riemannian manifold $(S, g)$ is
a complex structure on its Riemannian cone
\index[terms]{cone!Riemannian} $C(S):=\R^{>0}\times S$, 
 such that the cone metric on $C(S)$ is K\"ahler and the
homothety map $h_\lambda:\; C(S)\arrow C(S)$ sending
$(t,m)$ to $(\lambda t, m)$ is holomorphic. {\bf A Sasakian manifold}
is  a manifold equipped with a Sasakian structure. 
{\bf A morphism} of Sasakian manifolds is an isometric
immersion $S_1 \hookrightarrow S_2$ such that
the corresponding map of the cones $C(S_1)\arrow C(S_2)$
is holomorphic.
\index[terms]{manifold!Sasaki}

\hfill

\remark\label{_iR_holo_isometries_Remark_}
Let $r=t \frac d {dt}$ be the radial vector field on $C(S)$, acting
on $C(S)$ by holomorphic homotheties. It is also called
{\bf the Euler vector field}. \index[terms]{vector field!Euler} Since
$r$ is holomorphic, its complex conjugate $Ir$
is also holomorphic. Since $Ir$ is orthogonal to $r$,
it preserves the level sets of $t$, giving $\Lie_{Ir}(t)=0$. 
From the following claim, we obtain that $t^2$ is a K\"ahler
potential.  Therefore, $Ir$ 
preserves the metric, and acts on $C(S)$ by isometries.
Note that the vector field $Ir$ commutes with $r$, because it
is holomorphic (\ref{_holo_field_commute_Remark_}).

\hfill

\claim \label{_Sigma_defi_Claim_}
Let $S$ be a Sasakian manifold, $C(S)=S \times \R^{>0}$ its cone,
$t$ the coordinate along the second variable, and 
$r=t\frac{d}{dt}$ the Euler field. Then
$t^2$ is a K\"ahler potential on $C(S)$.\index[terms]{potential!K\"ahler} Moreover,
the form $dd^c \log t$
vanishes on $\langle r, I(r)\rangle$ and the rest of its
eigenvalues are positive.

\hfill

\proof
\[ 2\omega=\Lie_{r}\omega=  d(i_r \omega)= d(tIdt)=\frac
1 2 dd^c(t^2).
\]
Therefore, $t^2$ is a K\"ahler potential. Now, 
$dd^c \log t^2= \frac{\tilde \omega}{t^2}-\frac {dt\wedge
  Idt}{t^2}$, and hence  all eigenvalues of $dd^c \log t^2$
are positive except one, that vanishes.
\endproof

\hfill

\remark
A Sasakian manifold is clearly contact, because its cone
is symplectic and $h_\lambda$ acts by symplectic homotheties, that is,
this diffeomorphism maps the symplectic form $\omega$ to $\lambda^2 \omega$.
Indeed, $h_\lambda(I) = I$ and $h_\lambda(g) = \lambda^{2} g$.

\hfill

\remark \label{cylcon}
Sometimes, it is useful to view the above cone as a
cylinder $S\times\R$ ({ via}  the transformation
$u=\log t$) with the metric $e^{2u}g+ e^{2u}du^2$.

\hfill

Let $S$ be a Sasakian manifold, $\omega$ the K\"ahler form\index[terms]{form!K\"ahler}
on $C(S)$, and  $r=t\frac{d}{dt}$ the Euler field.\index[terms]{vector field!Euler}
Let $\eta\in \Lambda^1 S$ be the contact form,
$\eta= -i_r \omega\restrict S$.
Then $\Lie_{Ir}t= \langle dt, Ir\rangle=0$. Clearly,
$\eta(Ir)= -\omega (r, Ir)=t^2=1$ for $t=1$.
On the other hand, 
\begin{equation}\label{_Reeb_eta_Equation_}
i_{Ir} d\eta= \Lie_{Ir}(\eta)- i_{Ir} \eta = \Lie_{Ir}(\eta) = 
\Lie_{Ir}(i_r \omega)=0.
\end{equation}
The last term $\Lie_{Ir}(i_r \omega)$ 
vanishes because $Ir$ commutes with $r$, being
holomorphic, and preserves $\omega$ by \ref{_iR_holo_isometries_Remark_}.
This implies  

\hfill

\claim\label{_Reeb_via_cone_Claim_}
The  Reeb vector field of a Sasakian manifold is $\xi=Ir$,
where $r$ is the radial vector field.
\endproof 

\hfill

\theorem \label{_Reeb_Sasakian_Theorem_}
	The Reeb field acts on a Sasakian manifold
	by contact isometries, that is the flow of
 the Reeb field contains only isometries that preserve the
        contact form. Moreover, its action lifts
to a holomorphic action on its cone.\index[terms]{contact!sub-bundle}

\hfill

\proof Follows from \ref{_iR_holo_isometries_Remark_}
and \eqref{_Reeb_eta_Equation_}. \endproof

\hfill

\example 
The odd spheres are Sasakian.  Indeed, $C(S^{2n-1})=\C^n\setminus  0$. Note that $\C^n\setminus  0$ has many K\"ahler structures, and thus odd-dimensional spheres have many Sasakian structures. More generally, totally umbilical, oriented real hypersurfaces in K\"ahler manifolds are Sasakian, \cite{bl}. In particular, the product $S^2\times S^3$, that can be understood as the unit tangent bundle of $S^3$ is again Sasakian. Many examples are constructed and discussed in \cite{bog}.\index[terms]{hypersurface}\index[terms]{hypersurface!totally umbilical}

\hfill

There is a second way in which contact (resp. Sasakian) geometry is related to symplectic	(resp. K\"ahler) geometry. Namely, if the Reeb field is regular enough, one can look at its space of orbits that will  be a symplectic (resp. K\"ahler)  orbifold or even  manifold. This is called the {\bf \index[terms]{vector field!Reeb!regular} Boothby--Wang fibration}\index[terms]{Boothby--Wang fibration} (in circles), \cite{bl}.   The precise definitions are as follows.

\hfill

\definition\label{_regular_quasi_regular_sasakian_} 
A contact (resp. Sasakian ) manifold is called {\bf regular} if its \index[terms]{manifold!Sasaki!regular}\index[terms]{manifold!Sasaki!quasi-regular}
Reeb field generates a free action of $S^1$; it is called {\bf 
quasi-regular} if all orbits of the Reeb field are closed, and it is called 
{\bf irregular} otherwise.	

\hfill

\remark 
The Reeb field $\xi$ is  nowhere degenerate since, by definition, $\eta(\xi)=1$, where $\eta$ denotes the contact
1-form. Therefore, the Reeb action is always locally free.

\hfill

\corollary\label{regsas}
Let $S$ be a regular Sasakian manifold, and $\xi$ its Reeb
field. Then the space of $\xi$-orbits $X$ is K\"ahler.
Moreover, $X$ is equipped with a negative holomorphic
	Hermitian  line bundle\index[terms]{manifold!Sasaki!regular}
	$L$ such that $S$ is the space of 
unit vectors in $L$.\footnote{Therefore, 
whenever $X$ is compact, it is a complex projective manifold.}\index[terms]{manifold!projective}

\hfill

\proof   
Let $X=S/\xi$ be the space of orbits. This quotient is well
defined and smooth, because $\xi$ is regular.
Then $X=C(S)/\C^*$, where the $\C^*$-action is generated
by $ r=t\frac{d}{dt}$ and $ I(r)=\xi$, and hence it is holomorphic. {Therefore,
	$X$ is a complex manifold} (as  a quotient of a complex
manifold by a holomorphic action of a Lie group).
It is K\"ahler by \ref{_Sigma_defi_Claim_}, because
$dd^c(\log t)$ defines a K\"ahler form\index[terms]{form!K\"ahler} on the space
of orbits of $\C^*$.

 Consider now the  holomorphic
Hermitian line bundle
obtained from the $\C^*$-bundle $C(S)\arrow X$. Clearly,
$S$ is its space of unit vectors. Its curvature is
expressed by $dd^c\log |v|=-\omega_X$, and {hence this line 
	bundle is negative.} \index[terms]{bundle!vector bundle!holomorphic Hermitian!positive}
\endproof

\hfill

\remark We say ``$L$ is negative'' to keep with the
convention established by 
\ref{_curvature_line_bundle_via_15.19_Corollary_} 
and \ref{_signs_lost_Remark_}: the function
$dd^c |v|$ is plurisubharmonic on $\Tot(L)$ when
$L$ is a negative line bundle. However, this
distinction is irrelevant for the present purposes,
because as complex fibred spaces, the manifolds
$\Tot^\circ(L)$ and $\Tot^\circ(L^*)$ are isomorphic,
as shown in \ref{_C^*_bundle_isom_to_dual_Claim_} below.

\hfill

\claim\label{_C^*_bundle_isom_to_dual_Claim_}
Let $L$ be a holomorphic line bundle over $M$, and 
$\Tot^\circ(L)\mapsto M$ the total space of its
non-zero vectors, interpreted as a $\C^*$-bundle over
$M$. Then $\Tot^\circ(L)$ and $\Tot^\circ(L^*)$ 
are naturally isomorphic as principal holomorphic
bundles.\index[terms]{bundle!principal}

\hfill

\proof
Denote by $\underline L$ and $\underline L^*$ the
$\C^*$-bundles associated with $L$ and $L^*$. 
The corresponding total spaces are $\Tot^\circ(L)$ 
and $\Tot^\circ(L^*)$. We need
to show that $\underline L$ and $\underline L^*$ are isomorphic.
Let $\{U_i\}$ be a covering of $M$ and 
$\phi_{ij}:\; U_i\cap U_j\arrow \C^*$ 
the cocycles defining $\underline L$. 
Then $\underline L^*$ is defined by the cocycle
$\phi_{ij}^{-1}$. These cocycles are equivalent, because
the map $t\arrow t^{-1}$ defines an automorphism of $\C^*$
mapping $\phi_{ij}$ to $\phi_{ij}^{-1}$.
\endproof

\hfill

\remark 
The same statement and proof work for quasi-regular\index[terms]{manifold!Sasaki!quasi-regular}
Sasakian manifolds, only the base manifold will be a
K\"ahler orbifold \index[terms]{orbifold}(see Chapter
\ref{orbif}), and the  holomorphic line
bundles should also be considered in the orbifold
category. Hence, in the (quasi-)regular case, a contact
manifold $S$ sits between two symplectic
manifolds or orbifolds: the symplectic cone above
and the  Boothby--Wang quotient below; both
are K\"ahler when $S$ is Sasakian. \index[terms]{Boothby--Wang fibration}

\hfill

	\example 
		Let $X\subset \C P^n$ be a complex submanifold\index[terms]{submanifold}
		and $C(X) \subset \C^{n+1}\backslash 0$ the corresponding cone.
		The cone $C(X)$ is obviously K\"ahler, and hence  the
			intersection $C(X)\cap S^{2n-1}$,
                        called {\bf the link} of the singularity\index[terms]{link of singularity}
                        of $C(S)$, is Sa\-sa\-kian. This
		intersection  is an $S^1$-bundle over $X$. This construction
		gives many interesting contact manifolds,
		including \index[persons]{Milnor, J.} Milnor's exotic 7-spheres, that happen to be Sasakian, \cite{bog}.
				In other words, 
		{\em a link of a singularity of a cone
                  over a projective manifold is
			always Sasakian.}
		
		Moreover, it follows from \cite{ov_sas} that  every quasi-regular Sasakian manifold is obtained\index[terms]{manifold!Sasaki!quasi-regular}
		this way, for some  K\"ahler metric on $\C^{n+1}$

\hfill
	
\remark \label{ruki}
A 3-dimensional Sasakian manifold is a pseudoconvex 
hypersurface in a Hopf manifold or is 
quasi-regular
(\ref{_Sasakian_3-dim_qr_or_S^3_Corollary_}, \cite{geiges}, \cite{bel}).
In arbitrary dimension, every compact Sasakian manifold
can be deformed to a quasi-re\-gu\-lar one 
(\ref{_defo_qr_sas_}, \cite{_Rukimbira_deform_}, \cite{_Itoh_}).

\hfill

\example 
Let $M= \C^2\backslash 0=C(S^3)$ considered to be  a conical K\"ahler
manifold with the standard K\"ahler structure, 
and let $M_1=M/G$, where $G=\Z/4\Z$ is generated by 
$\tau(x,y)=(y,-x)$. Since $\tau^2=-1$, the action of $G$
on $M$ is free. {Therefore, $M_1$ is also a conical K\"ahler
	manifold.}

Then {\em $M_1$ is quasi-regular, but not regular.}\index[terms]{manifold!Sasaki!regular}\index[terms]{manifold!Sasaki!quasi-regular}

Indeed, the free orbits are those that satisfy
$(tx,ty)\neq (y,-x)$ for each $t\in \U(1)=\{\lambda\in
\C\ \ |\ \ |\lambda|=1\}$, $t\neq 1$. A non-free orbit
gives $tx=y, ty=-x =t^2x$, and hence  $t=\pm\1$ and 
$x=\pm y$.

\hfill

\remark 
	It was proven in \cite{bou} that all
        odd-dimensional tori bear contact structures (that cannot be regular \cite{bl}).\index[terms]{manifold!contact!regular} \index[terms]{contact!structure} However, odd-dimensional tori do not admit Sasaki structures, since the first Betti number of a compact Sasakian manifold should be even (see \cite{_Tachibana_} and Chapter \ref{_Harmonic_forms_chapter_}).

\hfill

\remark 
Given a contact structure on a manifold, there can be many compatible Sasakian structures. Concerning the relation among Sasakian and subjacent contact structures, see \cite{ov_geom_ded}.


\section{Notes}


\subsection{CR-structures and CR-manifolds}
\label{_CR_history_Subsection_}

The notion of CR-structure was motivated by a 
problem treated by \index[persons]{Poincar\'e, H.} Poincar\'e in 1907 (\cite{_Poincare_07_}). 
Clearly, all smooth codimension 1 subvarieties
in $\C$ are locally equivalent by a holomorphic
diffeomorphism. Poincar\'e noticed that this is
false in dimension 2, and proved that there
exist a pair of real hypersurfaces in $\C^2$ which
cannot be related by a biholomorphism even locally.
He asked whether it is possible to determine
a set of local invariants that can be used
to determine whether two real hypersurfaces 
are locally biholomorphic. \'Elie \index[persons]{Cartan, E.} Cartan
solved this problem in \cite{_Cartan_} for dimension
2, but it was much harder for arbitrary dimension.

This question is related to the famous Levi's problem
in complex analysis, named after Eugenio Elia Levi (1883--1917). 
Levi was interested in conditions for an open subset
of $\C^n$ to be a ``domain of holomorphy'' (\cite{_E_Levi_})
(in modern terms, to be Stein, see Subsection \ref{_Stein_Subsection_}).
He proved that the boundary of such a subset
cannot be arbitrary, and is characterized
by non-negativity of what is now called
the Levi form.\index[terms]{form!Levi} The ``Levi problem'' asks
whether the non-negativity of the Levi
form on the boundary is sufficient for
the domain to be Stein. It was solved
by \index[persons]{Oka, K.} Oka, \index[persons]{Bremermann, H. J.} Bremermann and \index[persons]{Norguet, F.} Norguet
(\cite{_Oka_, _Bremermann:Levi_, _Norguet_}, see also \cite{_Grauert_}).
The notion of strictly pseudoconvex CR-manifold
is in fact due to this circle of ideas and
precedes the appearance of CR-manifolds by
almost 50 years.

The problem of finding whether two hypersurfaces
are biholomorphic was treated by N. Tanaka (\cite{_Tanaka_}), who
used the method of prolongation going back to \'E. \index[persons]{Cartan, E.} Cartan.
In 1974, S. S.  Chern and J. K. \index[persons]{Moser, J. K.} Moser solved \index[persons]{Poincar\'e, H.} Poincar\'e's
problem in full generality (\cite{_Chern_Moser_}). 
Most of the modern treaties on CR-geometry
consider this to be a starting point for the research.

However, the modern notion of CR-manifold
(defined in terms of a complex structure on
a sub-bundle of arbitrary codimension, \ref{_CR_mfld_Definition_})
was treated much earlier in a thesis by
S. J. Greenfield\index[persons]{Greenfield, S. J.} (\cite{_Greenfield_}), who did his Ph. D. in 
Brandeis University under H. \index[persons]{Rossi, H.} Rossi.
Some of the people who treated the
subject (such as \cite{_Andreotti_Friedricks_})
refer to Greenfield directly.



\subsection[Sasakian manifolds: the tensorial definition by Sh. Sasaki]{Sasakian manifolds: the tensorial definition\\ by Sh. Sasaki}\label{_Sasaki:definition_}
Introduced by Shigeo \index[persons]{Sasaki, S.} Sasaki in the 1960s (see \cite{_Sasaki:notes_}), first called  ``normal metric contact manifolds'', and intensively studied by the Japanese school of geometry (see the references in the monograph \cite{bl}), what is known as Sasakian manifolds were defined locally, in terms of tensor fields. This structure was seen as a counterpart of K\"ahler geometry\index[terms]{geometry!K\"ahler} in odd dimension. Namely, a Sasakian structure on a manifold $S$ was given by the following tensor fields: a Riemannian metric $g$, a (1,1) tensor field $\f$, and a contact form $\eta$ with Reeb field $\xi$, 
that satisfied the relations:
\begin{enumerate}
	\item $g(\f X,\f Y)=g(X,Y)-\eta(X)\eta(Y)$, $\eta\circ\f=0$,
	\item $\eta(X)=g(X,\xi)$,
	\item $\f^2=-\id+\eta\otimes\xi$,
	\item $[\f,\f](X,Y):=\f^2[X,Y]-\f[\f X,Y]-\f[X, \f Y]-[\f X,\f Y]$ (the Nijenhuis tensor of $\f$) satisfies $[\f,\f]+2d\eta\otimes\xi=0$ (normality of $\f$),
	\item $\xi$ is Killing, and\index[terms]{vector field!Killing}
	\item $(\nabla^g_X\f)Y=g(X,Y)\xi-\eta(Y)X$
\end{enumerate}
for all $X,Y\in TS$.  The equivalence of this definition
with the one in terms of holonomy of the Riemannian cone\index[terms]{cone!Riemannian} (see also the Table ``Riemannian cones with
special holonomy" in Chapter \ref{other}) 
is implicit in the work of many authors who wrote on this
subject, but was explicitly proven and imposed in the
1990s by a long series of papers by C.P. Boyer\index[persons]{Boyer, C. P.},
K. \index[persons]{Galicki, K.} Galicki and their co-authors (see \cite{bog} and the
references therein). Ever since, Sasakian manifolds have proved
to be extremely important in string theories and also in
algebraic geometry\index[terms]{geometry!algebraic}. Various generalizations of these notions exist: {\bf K-contact manifolds} ($\xi$ is Killing,
but the almost complex structure on the cone is not necessarily
integrable), almost contact metric etc.\index[terms]{vector field!Killing}


\section{Exercises}


\begin{enumerate}[label=\textbf{\thechapter.\arabic*}.,ref=\thechapter.\arabic{enumi}]

\item
Let $S$ be a regular Sasakian manifold,\index[terms]{manifold!Sasaki!regular} and $X$ the space
of its Reeb orbits, considered to be  a smooth manifold.
Assume that $\pi_1(S)$ is finite. Prove that
$\pi_2(X)$ is infinite.

\item\label{_pullback_trivial_Exercise_}
Let $L$ be a complex Hermitian line bundle on 
a topological space $X$, and $S\stackrel \pi \arrow X$ the
space of unit vectors in $L$, that is  considered
as an $S^1$-fibration over $X$. Prove that
the pullback bundle $\pi^*L$ admits a trivialization.

\item
Let $S$ be a regular Sasakian manifold\index[terms]{manifold!Sasaki!regular}, and $X$ the space
of its Reeb orbits, considered to be  a smooth manifold. Assume 
that the universal cover of $S$ is contractible.
\begin{enumerate}
\item 
Prove that the  universal cover of $X$ is contractible.
\item
Construct the exact sequence
$0 \arrow \Z \arrow \pi_1(S) \arrow \pi_1(X) \arrow 0.$
\item
Prove that $\pi_1(S)$ contains a subgroup 
isomorphic to $\Z\oplus \Z$.\footnote{As follows from
\cite[Corollaire 7.3]{_CDP:geometrie_groupes_},
this also implies that $S$ does not admit a metric
of strictly negative sectional curvature.}
\item Let $\pi:\; S \arrow X$ be the
canonical projection. Using Exercise \ref{_pullback_trivial_Exercise_},
prove that there exists a complex line bundle
$L$ on $X$ that is  topologically non-trivial,
but $\pi^* L$ is topologically trivial.
\item
Using this observation, prove that the
exact sequence \\
$0 \arrow \Z \arrow \pi_1(S) \arrow \pi_1(X) \arrow 0$
does not split.
\end{enumerate}

\item 
Let $\phi$ be a real-valued smooth 
function on a complex manifold, $S$ its level set for a regular value,
that is  connected, and $(B,I)$ the standard CR-structure\index[terms]{structure!CR} on $S$,
with $B= TS \cap I(TS)$. Prove that
$dd^c\phi\restrict B$ is proportional to the
Levi form\index[terms]{form!Levi} of $S$ with positive or negative
coefficient (not necessarily constant).

\item Let $S\subset M$ be a compact real hypersurface in
a complex manifold. Assume that the CR-structure\index[terms]{structure!CR}
induced by the sub-bundle $(TS \cap I(TS))\subset TS$ is strictly pseudoconvex.
Prove that there exists a strictly plurisubharmonic
function defined in a neighbourhood of $S$ such 
that $S$ is its level set. 

{\em Hint:} Use the previous exercise.

\item
{(J. Gray, see  \cite{bl})}. 
Consider a smooth hypersurface $i:M\hookrightarrow\index[terms]{hypersurface}
\R^{2n+2}$, such that no tangent space of $M$ (viewed as
an affine hyperplane) contains the origin of
$\R^{2n+2}$. Prove that $M$ is contact with the contact form $\eta=i^*\al$, where \index[terms]{form!contact}
$$\al=x_1dx_2-x_2dx_1+\cdots+x_{2n+1}dx_{2n+2}-x_{2n+2}dx_{2n+1}.$$
{\em Hint:}
Let $\{X_1, \ldots, X_{2n+1}\}$ be a basis for $T_{x^0}M$. Then a  normal vector $W=(W_i)$ to $M$ in $x^0$ can be written as 
$$W_i=(*dx_i)(X_1,\ldots,X_{2n+1}),$$
where $*$ is the Hodge operator of the flat metric on $\R^{2n+2}$. Then 
$$\left(\al\wedge(d\al)^{n}\right)(X_1,\ldots,X_{2n+1})=\sum_i x^0_iW_i.$$
In particular, this is another proof that $S^{2n+1}$ is contact.

\item Consider {\bf the Penrose hypersurface} $S:=\{z\in \C^3 \ \ |\ \ |z_1|^2+ 
|z_2|^2 = 1+ |z_3|^2\}$. Prove that its Levi form\index[terms]{form!Levi} is non-degenerate\index[terms]{hypersurface!Penrose}
and find its signature.

\item A real-valued function $\phi$ on a complex manifold
is called {\bf pluriharmonic}
if $dd^c\phi=0$.
Let $\phi$ be a pluriharmonic function, $c$ a regular value, and $S:=
\phi^{-1}(c)$. Prove that the Levi form $\Phi$ of $S$  vanishes
(such $S$ is called {\bf Levi-flat}).\index[terms]{function!pluriharmonic}

\item 
The {\bf local $dd^c$-lemma}\index[terms]{lemma!$dd^c$} states that any closed
$(p,q)$-form, $p, q>0$, on an open ball in $\C^n$ is {\bf
  $dd^c$-exact}, that is, it belongs to the image of $dd^c$.
Deduce this statement from the
{\bf Poincar\'e--Dolbeault--Grothendieck lemma}, 
which says that any $\bar \6$-closed $(p,q)$-form, $q>0$,
is $\bar\6$-exact.

\item Prove that any pluriharmonic function on an 
open ball is a sum of holomorphic and antiholomorphic
functions (use the $dd^c$-lemma).

\item  Prove that any pluriharmonic function on  
$\C^*$ is a sum of holomorphic and antiholomorphic
functions, or find a counterexample.

\item Let $S\subset \C^n$ be a compact smooth real hypersurface,
such that the Levi form\index[terms]{form!Levi} $\Phi$ of $S$ is non-degenerate. Prove that 
$\Phi$ is sign-definite.

\item Let $S\subset \C^n$ be a smooth Levi-flat hypersurface.\index[terms]{hypersurface!Levi-flat}
Prove that $S$ is non-compact.

\item Let $M$ be a complex manifold, and $\phi:\; M \arrow \R$
a smooth function satisfying $d\phi \wedge d^c \phi \wedge  dd^c\phi=0$. 
Prove that for any regular value $c$ of $\phi$,
the preimage $\phi^{-1}(c)$ is Levi-flat.

\item Let $(S,g,\eta)$ be a Sasakian manifold with Reeb
  field $\xi$  and define a (1,1)-tensor field $\f$ by
  $\f=\nabla\xi\in TM \otimes \Lambda^1 M =\End(TM)$. Prove that:
	\begin{equation*}
	\begin{split}
	\f^2&=-\id +\eta\otimes\xi,\\
	\nabla\f&=g-\eta\otimes\eta.
	\end{split}
	\end{equation*} 
	 This $\f$ corresponds to the tangent to $S$ part of $I$ on the cone $C(S)$.
\item Recall that the isometry group \index[terms]{group!isometry}$\Iso(M,g)$	of a Riemannian manifold $(M,g)$ is a Lie group (this is Meyers-\index[persons]{Steenrod,  N. E.}Steenrod's theorem, \index[terms]{theorem!Myers--Steenrod} see \cite{_Kobayashi_Transformations_}). Prove that if $M$ is compact, then $\Iso(M,g)$ is compact. For a generalization to metric spaces,\index[terms]{space!metric} see Exercises \ref{_compact_isom_group_1_exercise_}, \ref{_compact_isom_group_2_exercise_}.

\item (\cite{gopp})\label{iso_con} Let $W$ be a compact
  Riemannian manifold, and $C(W)$ its Riemannian
  cone. Prove that the group of isometries $\Iso(C(W))\simeq
  \Iso(W)$. Precisely, every isometry of $C(W)$ has the\index[terms]{isometry}
  form $(w,t)\mapsto (\psi(w), t)$, where
  $\psi\in\Iso(W)$. More generally, a homothety of $C(W)$
  with scale factor $a\in \R$ is of the form 
$(w,t)\mapsto (\psi(w), a\cdot t)$, where $\psi\in\Iso(W)$.

{\em Hint:} Show first that the metric completion of $C(W)$ is obtained by adding a single point, called origin. Then show that  every isometry of $C(W)$ can be trivially extended to the origin, preserving the rays and the levels of the form $\{t\}\times W$. This proves the statement for isometries, the one for homotheties follows easily.

\end{enumerate}

\chapter{Vaisman manifolds}\label{vaiman}

{\setlength\epigraphwidth{0.5\linewidth}
\epigraph
{\it De la musique encore et toujours!\\
Que ton vers soit la chose envol\'ee\\
Qu'on sent qui fuit d'une \^ame en all\'ee\\
Vers d'autres cieux \`a d'autres amours.\\ \medskip

Que ton vers soit la bonne aventure\\
Eparse au vent crisp\'e du matin\\
Qui va fleurant la menthe et le thym...\\
Et tout le reste est litt\'erature.}
{\sc\scriptsize Paul Verlaine, ``Art po\'etique''}
}

\section{Introduction}\index[terms]{manifold!Vaisman}

\subsection{Many definitions of Vaisman manifolds}

Vaisman manifolds were an invention of Izu \index[persons]{Vaisman, I.} Vaisman
(\cite{va_rendiconti}). When he started working with the subject,
he just called them ``LCK manifolds\index[terms]{manifold!LCK} with parallel Lee\index[terms]{form!Lee!parallel}
form''. After a while, Vaisman coined the name ``generalized
Hopf manifolds'' for this class. The name was not very
appropriate, because, as it turned out, not all Hopf
manifolds are Vaisman. In \cite{do}, the authors
suggested the term ``Vaisman manifolds'', that is 
commonly used now.

The notion of \index[terms]{manifold!Vaisman} Vaisman manifold is still slightly
ambiguous, because ``Vaisman manifold'' can mean either
"a complex manifold admitting a Vaisman metric"
or ``an Hermitian complex manifold equipped with
a Vaisman metric''. The same ambiguity is present
in K\"ahler geometry.\index[terms]{geometry!K\"ahler} To distinguish between
these notions, people sometimes use the term
``manifold of Vaisman type'' (``manifolds of K\"ahler
type, in K\"ahler case''). 

Another ambiguity is present in LCK geometry,
but it is pretty much absent in the Vaisman case.
Speaking of LCK manifold\index[terms]{manifold!LCK}, one can fix the metric,
or its conformal class. It turns out that the
Vaisman metric is unique (up to a constant)
in its conformal class.

The first and most straightforward definition of Vaisman manifold\index[terms]{manifold!Vaisman}
was given by Vaisman\index[persons]{Vaisman, I.}: it is an LCK manifold $(M, I, g,
\theta)$ such that the Lee form\index[terms]{form!Lee} $\theta$ is parallel\index[terms]{form!Lee!parallel}
under the Levi--Civita connection. For the students
trained in algebraic geometry\index[terms]{geometry!algebraic}, several equivalent
versions are suggested.\index[terms]{manifold!LCK}

Let $(M, I, g, \theta)$ be a locally conformally K\"ahler manifold
equipped with a holomorphic vector field $v$. Suppose that
$v$ is Killing\index[terms]{vector field!Killing} (that is, its diffeomorphism flow acts by
isometries), and acts on the corresponding K\"ahler metric
by non-trivial homotheties. Then $(M, I, g, \theta)$ is
Vaisman (\ref{kami_or}). This result was proven in
\cite{kor} in compact situation; we give its simplified
proof in this chapter. Sometimes this is actually used
as a definition.

Another definition of complex manifolds of Vaisman type,
that is  valid only in the compact case, can be obtained as a consequence
of the structure theorem, proven in \ref{str_vai}.
Let $X$ be a projective orbifold (for the definition
of an orbifold, see Chapter \ref{orbif}), and $L$\index[terms]{bundle!line!ample}
an ample line bundle. Denote by $\Tot^\circ(L)$ the space
of non-zero vectors in the total space of $L$.
Suppose that $\Tot^\circ(L)$ is smooth (this depends on
the singularities of the orbifold $X$ and the
structure of the line bundle $L$ in the neighbourhood of
these singularities). 
Let $\phi_X$ be a holomorphic automorphism of $X$ acting
equivariantly on $L$; denote by $\phi$
its extension to $\Tot^\circ(L)$. Since $L$ is ample, it admits
a Hermitian metric $h$ with positive curvature. Suppose that $\phi$
acts on $\Tot^\circ(L)$ multiplying the length of each 
vector by a scalar $\lambda>1$, 
\[
|\phi(v)|_h = \lambda |v|_h.
\]
Then the quotient manifold
$\Tot^\circ(L)/\langle \phi\rangle$ is of Vaisman type,
and any compact complex manifold of Vaisman type
can be obtained in this way\index[terms]{manifold!Vaisman} (\ref{_Structure_of_quasi_regular_Vasman:Theorem_} and \ref{defovai}).

This suggests that the Vaisman geometry is
in fact "more algebraic" than the K\"ahler geometry;\index[terms]{geometry!K\"ahler}
indeed, any Vaisman manifold\index[terms]{manifold!Vaisman} is associated with
a projective orbifold, and can be interpreted
in terms of projective algebraic geometry\index[terms]{geometry!algebraic!projective}.

Recall that {\bf a linear Hopf manifold}
is a quotient of $\C^n \backslash 0$ by
an action of a cyclic group generated by
an invertible endomorphism with operator norm
$\Vert A\Vert < 1$.
In \ref{semihopf} we prove that 
a linear Hopf manifold is Vaisman
if and only if $A$ is diagonalizable.

Similarly to
projective manifolds, that are  complex manifolds
admitting a complex embedding to $\C P^n$,
compact Vaisman manifolds are manifolds
that admit a complex embedding to a
diagonal Hopf manifold (\ref{_diagonal_Hopf_Definition_}). This can be considered
as another definition of Vaisman manifolds.\index[terms]{manifold!Vaisman}

\subsection{Riemannian cones}\index[terms]{cone!Riemannian}

Recall that the holonomy group of a Riemannian 
manifold is the holonomy of its Levi--Civita
connection. The study of Riemannian geometry
from this point of view has a long and illustrious history
starting from de Rham and \index[persons]{Berger, M.} Berger. 
De Rham proved that any Riemannian manifold
splits locally as a Riemannian product of
manifolds with irreducible holonomy
(\ref{_de Rham_holonomy_Theorem_}),
and \index[persons]{Berger, M.} Berger gave a classification of 
Riemannian manifolds with irreducible
holonomies it terms of their holonomy
representations (Chapter \ref{other}). 

For Vaisman geometry, de Rham's observation
has a crucial importance: a Vaisman manifold\index[terms]{manifold!Vaisman}
is locally isometric to a product of a Riemannian
manifold and a real line. Indeed, on a Vaisman
manifold the Lee field\index[terms]{Lee field} is parallel, and hence \index[terms]{Lee field!parallel}
its holonomy group is reducible, and the
direction of this vector field corresponds to
the real line summand.

Another important property of the Vaisman
manifolds relates to the Riemannian 
structure on its K\"ahler cover.
It turns out that this K\"ahler
cover is locally isometric
to a Riemannian cone (\ref{vai_cone}).\index[terms]{cone!Riemannian}

Let $(S,g_S)$ be a Riemannian manifold.
Recall that the Riemannian cone 
over $S$ is the manifold $S\times \R^{>0}$
with the metric $t^2 g_S + (dt)^2$, where
$t$ is the parameter in $\R^{>0}$.
To prove that the coverings of the Vaisman manifolds
are locally conical, we use the following\index[terms]{manifold!Vaisman} 
local characterization of the cone geometry.

Let $(M,g)=(S\times \R^{>0}, t^2 g_S + (dt)^2)$
be a Riemannian cone,\index[terms]{cone!Riemannian} and $v=t \frac{d}{dt}$ the Euler vector field. \index[terms]{vector field!Euler}
It is easy to see that the corresponding
diffeomorphism flow acts on $(M,g)$
by homotheties, multiplying the metric
$g$ by a constant: $\Lie_v(g) = 2g$.
We give a characterization of Riemannian
cones in terms of the Euler vector field:
a Riemannian manifold is locally conical if
and only if it admits a vector
field $v$ such that $\Lie_v(g) = 2\lambda g$
and the dual form $\eta=v^\flat$ is closed
(\ref{_nabla(X)=Id_then_cone_Theorem_}).

This result makes sense even when $\lambda=0$,
in which case it characterizes manifolds
satisfying $\nabla(v)=0$, that is, the Riemannian
products with a line (``Riemannian cylinders'').

This observation makes it possible to pass
from the cone manifolds to cylinders;
this is used as a model for Vaisman geometry.

Let $(M,g)$ be a Riemannian manifold, and 
$v$ a tangent vector field such that $\Lie_v(g) = 2\lambda g$
and the dual form $\eta=v^\flat$ is closed.
Consider a function  $\phi$ 
such that $\eta= d\phi$. We choose $\phi$
in such a way that $\Lie_v(\phi) = 2\lambda \phi$,
and show that $\phi^{-1} g$ is cylindrical.
In Vaisman context, this function 
gives a K\"ahler potential for the K\"ahler
metric on the cone (\ref{_tilde_M_potential_Corollary_}).

\subsection{Basics of Vaisman geometry}
\label{_gauduchon_basics_Subsection_}

One of the most important tools in complex
geometry is the Gauduchon theorem. This is 
possibly the only way to approach a general non-K\"ahler
complex manifold, apart from the topological methods,
such as the index theorem. We prove the \index[persons]{Gauduchon, P.} Gauduchon theorem
in Section \ref{_appendix_gauduchon_}. 

Let $(M,I, g)$ be a compact complex Hermitian manifold,
$\dim_\C M=n$, and $\omega\in \Lambda^{1,1}(M)$ its Hermitian form.\index[terms]{form!Hermitian}
We say that $g$ is {\bf a Gauduchon metric} if it
satisfies {\bf the Gauduchon condition} $dd^c (\omega^{n-1})=0$.\index[terms]{condition!Gauduchon}
Gauduchon theorem claims that any Hermitian metric
is conformal to a Gauduchon metric, that is  unique
up to a constant multiplier.

Let $(M, I, g, \theta)$ be a Vaisman manifold,
$n=\dim_\C M$. Denote by $*$ the Hodge star operator,
e.g. \cite{demailly}.\index[terms]{manifold!Vaisman}
Since $\nabla(\theta)=0$, the form $\theta$
is closed, and co-closed, that is, satisfies
$d(*\theta)=0$. 
It is not hard to compute that
$*(\theta)= c I(\theta)\wedge \omega^{n-1}$,
where $c$ is a constant. However, $d^c(\omega) = I(\theta) \wedge \omega$, which gives
$d^c(\omega^{n-1})= (n-1) I(\theta)\wedge \omega^{n-1}$.
Then
\[
d(*\theta)= c d((I(\theta) \wedge \omega^{n-1}) = \frac{c}{n-1}d d^c (\omega^{n-1}).
\]
We obtain that the \index[persons]{Gauduchon, P.} Gauduchon condition $d d^c (\omega^{n-1})=0$
is equivalent to $d^*\theta=0$; the latter follows from the
Vaisman condition $\nabla \theta=0$ automatically. 

By Gauduchon's theorem, on each compact Hermitian manifold
there exists a Gauduchon metric in the
same conformal class. Moreover, the \index[persons]{Gauduchon, P.} Gauduchon metric is 
unique up to a constant. This fixes the conformal gauge
for Vaisman geometry: the conformal ambiguity, which
is always present in the LCK geometry, disappears when
one passes to Vaisman manifolds. When we say ``Vaisman
manifold'', we usually assume that the metric is fixed.

Examples of Vaisman manifolds abound.\index[terms]{manifold!Vaisman} 
Not all Hopf manifolds are Vaisman, however,
the Vaisman Hopf manifolds are dense in the
space of Hopf manifolds (\ref{def_lckpot2Vai}). In Chapter \ref{comp_surf}
we prove that any non-K\"ahler surface with $b_1 >1$
is Vaisman, too. 

In general, it's hard to say whether a given 
manifold with $b_1 >0$ admits an LCK structure.
In fact, we do not seem to have any example of
a compact manifold $M$ with 
$b_1(M) >0$ and $\dim_\R M > 4$ that does not admit a
complex structure of LCK type.
However, the situation for a Vaisman manifold is opposite:
the topology and geometry of Vaisman manifolds is very restrictive.
For example, the first Betti number of a compact Vaisman manifold is\index[terms]{manifold!Vaisman}
always odd (\ref{_odd_first_betti_}). Vaisman conjectured that $b_1(M)$ is odd\index[terms]{conjecture!Vaisman}
for all compact LCK manifolds\index[terms]{manifold!LCK}, but in 2004, a counterexample
was found by \index[persons]{Oeljeklaus, K.} Oeljeklaus and \index[persons]{Toma, M.} Toma (\ref{_Vaisman_conjecture_disproved_}).

Unlike more general LCK manifolds\index[terms]{manifold!LCK},
Vaisman manifolds have big automorphism groups;
indeed, the {\bf Lee field} $\theta^\sharp$ (the vector field dual to the Lee form)\index[terms]{Lee field}\index[terms]{form!Lee}
acts on a Vaisman manifold by holomorphic isometries.
The same is true for the {\bf anti-Lee field}\index[terms]{Lee field!anti-} $I(\theta^\sharp)$.\index[terms]{manifold!Vaisman}
These two vector fields generate a holomorphic, totally geodesic foliation on
any given Vaisman manifold, that is  called {\bf the canonical foliation}.\index[terms]{foliation!canonical}

The canonical foliation $\Sigma$ is a convenient object, because
it is determined by the complex structure, being 
independent on  the choice of the Vaisman metric
(\ref{_Subva_Vaisman_Theorem_}). Moreover, 
any Vaisman manifold admits a semi-positive,
closed (1,1)-form $\omega_0$ vanishing 
on $\Sigma$ and defining a K\"ahler structure
on the leaf space of $\Sigma$ (the leaf space
does not always exist globally, but locally
it is well-defined, and the K\"ahler structure
induced by $\omega_0$ is clearly global when
the leaf space is defined). 

This chapter is also introductory and should
be accessible to a reader with some background in geometry.
We use the content of the previous chapters and 
some basic notions of Riemannian geometry, such
as Killing fields\index[terms]{vector field!Killing} and Levi--Civita connection (see, for example, \cite{_Gallot_Hulin_Lafontaine_}).

\section{Holonomy and the de Rham split\-ting theorem}

The following exposition is standard in Riemannian
geometry. The reader is referred to \cite{besse} for 
the proofs and more details.

\hfill

\definition (\index[persons]{Cartan, E.}Cartan, 1923)
Let $(B,\nabla)$ be a vector bundle with connection over $M$.
For each loop $\gamma$ based in $x\in M$, let 
$V_{\gamma, \nabla}:\; B\restrict x \arrow B\restrict x$
be the corresponding parallel transport along the connection.
The {\bf  holonomy group} of\index[terms]{parallel transport} $(B,\nabla)$\index[terms]{group!holonomy}
is the group generated by $V_{\gamma, \nabla}$,
for all loops $\gamma$. If one takes all contractible
loops instead, $V_{\gamma, \nabla}$ generates
{\bf the local holonomy}, or {\bf 
the restricted holonomy} group.\index[terms]{group!holonomy!restricted}

\hfill

\remark A bundle is {\bf flat} (has vanishing curvature)
if and only if its restricted holonomy vanishes.\index[terms]{bundle!vector bundle!flat}

\hfill

\remark If $\nabla(\phi)=0$ for some tensor 
$\phi\in B^{\otimes i}\otimes (B^*)^{\otimes j}$,
the holonomy group preserves $\phi$.

\hfill

\definition {\bf The holonomy of a Riemannian manifold}
is the holonomy of its Levi--Civita connection.

\hfill

\example The holonomy of a Riemannian manifold lies in
$\OO(T_x M, g\restrict x)\simeq \OO(n)$.\index[terms]{group!holonomy}

\hfill

\example  The holonomy of a K\"ahler manifold lies in
$\U(T_x M, g\restrict x, I \restrict x)\simeq \U(n)$.

\hfill

\remark 
Let $x, y \in M$ be points, and $\Hol_x(M)\subset \GL(T_xM)$,
 $\Hol_y(M)\subset \GL(T_yM)$ the holonomy groups
associated with all loops in $x$ and in $y$.
Clearly, for any path $\gamma$ from $x$ to $y$,
the holonomy along $\gamma$ induces
an isomorphism $h_\gamma:\; \Hol_x(M)\arrow \Hol_y(M)$.
Abusing the language, this is sometimes
expressed as ``the holonomy group of a connected manifold
does not depend on the choice of a point $x\in M$.'' 
Note that the isomorphism $h_\gamma$ depends 
on the choice of $\gamma$, and different choices
give isomorphisms that are  equivalent up to a conjugation 
with an element of $\Hol_y(M)$.
Since we are interested in the holonomy group
up to conjugation, we will usually adopt this point of view.

\hfill

The following result is not very hard to prove.

\hfill

\claim
Let $M$ be a Riemannian manifold, 
and $\Hol_0(M)\stackrel \rho \arrow \End(T_xM)$
the reduced holonomy representation. Suppose that $\rho$ is reducible:
$T_xM = V_1\oplus V_2 \oplus ...\oplus V_k$.  Then $G=\Hol_0(M)$ also
splits: $G= G_1\times G_2 \times \cdots\times G_k$,
with each $G_i$ acting trivially on all $V_j$ with $j\neq i$.

\hfill

This leads to the following global result, that is  
easy to establish locally, but the global version is
rather hard (\cite{KoNFD2}).

\hfill

\theorem \label{_de Rham_holonomy_Theorem_}
(de Rham) 
A complete, simply connected  \index[terms]{theorem!de Rham splitting}
Riemannian manifold with non-irreducible holonomy 
splits as a Riemannian product onto factors
corresponding to irreducible components of the holonomy
representation.

\hfill

We will not need de Rham splitting  in whole generality, 
but we use the following corollary of the
local version of the de Rham theorem.

\hfill

\corollary\label{_parallel_field_de_Rham_Corollary_}
Let $X\in TM$ be a vector field satisfying $\nabla X=0$.
Then $M$ locally splits as a Riemannian manifold: $M=M_1 \times I$,
where $I\subset \R$ is an interval equipped with the standard metric.
\endproof

\section{Conical Riemannian metrics}

\proposition\label{_Killing_via_nabla_Proposition_}
Let $g\in \Sym^2T^*M$ be a Riemannian form on $TM$, 
$v\in TM$ a vector field, $\nabla$ the Levi--Civita connection, 
and $h:=\Lie_v(g)$.  Then 
$h(Y, Z)=g(\nabla_YX,Z)+ g(\nabla_Z X, Y)$.

\hfill

\proof
Since the Lie derivative is compatible with contraction,
\[ \Lie_X(g(Y,Z))=\Lie_X(g)(Y,Z)+ g([X,Y],Z)+ g([X,Z],Y).\]
Similarly, 
\[ \nabla_X(g(Y,Z))=\nabla_X(g)(Y,Z)+ g(\nabla_XY,Z)+ g(\nabla_XZ,Y).\]
However, $\nabla_X(g(Y,Z))=\Lie_X(g(Y,Z))$, giving
\begin{align*}
\Lie_X(g)(Y,Z)&=g(\nabla_XY,Z)+ g(\nabla_XZ,Y)-g([X,Y],Z)- g([X,Z],Y)=\\
& = g(\nabla_YX,Z)+ g(\nabla_Z X, Y)
\end{align*}
using $\nabla_XY-[X,Y]=\nabla_YX$, $\nabla_XZ-[X,Z]=\nabla_ZX$.
\endproof

\hfill

\corollary 
A vector field $X$ satisfies $\Lie_X(g)=0$
if and only if $\nabla X$ is an antisymmetric operator: $g(\nabla_YX, Z)=- g( Y, \nabla_ZX)$
for all vector fields $Y,Z$.
\endproof

\hfill

\definition
Such a vector field is called {\bf Killing}.
A vector field is Killing if and only if the corresponding
diffeomorphism flow acts by isometries.\index[terms]{vector field!Killing}

\hfill

\theorem \label{_A_X_defi_Theorem_}
Let $M$ be a Riemannian manifold, $\nabla$ the Levi--Civita connection,
and $X\in TM$ a vector field. Consider the operator 
$A_X(\psi):= \nabla_X(\psi)-\Lie_X(\psi)$ on tensors. Then
\begin{description}
\item[(i)]
$A_X$ is $C^\infty$-linear, satisfies the Leibniz rule,
and commutes with contraction.
 \item[(ii)]
On tensors $\psi \in TM^{\otimes k} \otimes T^*M^{\otimes l}$, 
the operator $A_X$ acts as 
\[ A_X(\psi)=\sum_{i=1}^{k+l} {\mathbb A}_i,
\] 
where
${\mathbb A}_i= A_X \in TM\otimes \Lambda^1 M=\End(TM)=\End(T^*M)$
is understood as an endomorphism of $TM$ or $\Lambda^1 M$,
acting on the $i$-th component of the tensor product.
If $A_X=\lambda\Id$, this gives $A_X(\psi)=(k+l) \lambda\psi$.
\end{description}

\hfill

{\bf  Proof. Step 1:} Linearity follows from 
$\nabla_X(f\psi)= f \nabla_X(\psi)+\Lie_X f(\psi)$
and $\Lie_X(f\psi)= f \Lie_X(\psi)+\Lie_X f(\psi)$,
and Leibniz and contraction identity from similar identities for $\Lie_X$
and $\nabla$. 

\hfill

{\bf Step 2:} On vector fields, $\nabla_X Y - [X,Y]=\nabla_Y X$
because $\nabla$ is torsion-free. \index[terms]{connection!torsion-free}

\hfill

{\bf Step 3:} On $\Lambda^1 M$, we have $A_X=-\nabla(X)$,
because $A_X$ is a derivation that preserves the
pairing between 1-forms and the vector fields:
\[
0=A_X (\langle v, \xi\rangle)= \langle A_X(v), \xi\rangle +
\langle v, A_X(\xi)\rangle.
\]
Then $A_X(\psi)=(k+l)\nabla(X)$ follows, because
$A_X$ satisfies the Leibniz rule with respect
to the tensorial multiplication.
\endproof

\hfill

We obtain the following characterization of homothetic vector fields\index[terms]{vector field!homothetic} on a Riemannian manifold.

\hfill

\proposition \label{homoth}
	Let $\eta$ be a  1-form on a Riemannian manifold $(M,g)$, 
	and $X=\eta^\sharp$ its $g$-dual vector field.  
Then the following are equivalent.
\begin{description}
\item[(i)]
 $\nabla(X)=\lambda \Id$,\ $\la\in\R$.
 \item[(ii)]$\nabla(\eta)=\lambda g$. 
\item[(iii)] $d\eta=0$ and $\Lie_X g=2\lambda g$.
\end{description}

\proof  Indeed,  
$\nabla(X)=\lambda \Id$ is clearly equivalent to
$\nabla\eta= \lambda g$ (one is obtained from another by applying $g^{-1}$,
that is  parallel).

Then  (ii) and (i) imply $d\eta=0$, as 
$d\eta=\Alt(\nabla\eta)=0$, since $g$ is symmetric.

From $\nabla(X)=\lambda \Id$ and \ref{_A_X_defi_Theorem_} we obtain 
$\nabla_X g - \Lie_X g= -2\lambda g$, giving $\Lie_X g=2\lambda g$.

Finally,  (iii) implies (i), 
since from $\Lie_X g=2\lambda g$ and $\nabla_X g=0$,
we obtain that $\nabla(X)(g)=2\lambda g$. 
This yields that the
	symmetric part of $\nabla(X)$, considered to be  a section
	of $\End(TM)$, is equal to $\lambda \Id$, whereas
the antisymmetric part of $\nabla(X)$ acts on $g$ trivially.
To obtain the antisymmetric part, it is more convenient
to replace $\nabla(X)$ by $\nabla(\eta)$.  Then 
	the antisymmetric part of $\nabla(\eta)$
	is equal to $\Alt(\nabla(\eta))=d\eta=0$.
 \endproof

\hfill

\theorem\label{_nabla(X)=Id_then_cone_Theorem_}
Let $(M,g)$ be a Riemannian manifold, and 
$X$ a nowhere vanishing vector field that satisfies 
$\nabla(X)=\lambda \Id$, $\lambda>0$. By \ref{homoth},
the dual 1-form $\eta= X^\flat$ is closed.
Then, locally, in a neighbourhood of each point
there exists a function $\phi\in C^\infty M$ 
such that $d\phi=\eta$ and $\Lie_X \phi=2\lambda\phi$.
Moreover, $X$ is the gradient vector field for $\phi$,
and $\phi^{-1} X$ is a parallel vector field on a Riemannian 
manifold $(M, \phi^{-1} g)$. Finally,\index[terms]{vector field!parallel}
$(M, g)$ is locally isometric to a Riemannian cone.\index[terms]{cone!Riemannian}

\hfill

\proof
It suffices to prove the theorem if we rescale $X$ in such a way that
$\lambda=1$. Since $\Lie_X X=[X,X]=0$, and $\Lie_X g= 2 g$,
we have $\Lie_X \eta=2\eta$. 
Let $\phi$ be a function such that $d\phi=\eta$. 
This gives $\Lie_X \phi=2\phi+\const$. Adding
a constant to $\phi$ if necessarily, we can assume that
$\Lie_X \phi=2\phi$. Since $X$ is dual to $d\phi$,
$X$ is the gradient vector field for $\phi$. Then $\Lie_X \phi\geq 0$.
Since $2\phi=\Lie_X \phi \geq 0$, a point $x$ where $\phi(x)=0$ 
is a local minimum of $\phi$, and $X= \grad(\phi)=0$ at this point.
Since $X$ is nowhere zero, this does not happen,
and we have $2\phi=\Lie_X \phi>0$.

Consider the metric $g_1 = \phi^{-1}g$, and let
$\eta_1= \phi \eta$ be the dual 1-form to $X$
with respect to $g_1$. Then $\Lie_X g_1 = 0$
and $d\eta_1= d\phi \wedge d\phi=0$.
This means that $(M, g_1, X)$ satisfies
the assumptions of \ref{homoth} (iii) with $\lambda=0$,
giving $\nabla_1(X)=0$ by \ref{homoth} (i).

From the de Rham holonomy decomposition theorem
(\ref{_parallel_field_de_Rham_Corollary_})
it follows that $(M, g_1)$ is a cylinder, that is, 
a product of $\R$ with a Riemannian manifold $(M', g')$.
Here $t=\log \phi$ is the parameter on $\R$, giving
$(M, g_1) = (M' \times \R, g' + (dt)^2)$. Clearly,
$d\phi = d e^{t} = e^{t} dt$. Let $\psi:= \sqrt \phi$.
Then $d\psi = 1/2 \psi dt$, and
$(d\psi)^2 = 1/4 \phi (dt)^2$. This identifies
$M$ with $M' \times \R^{>0}$, with parameter $\psi$ on $R^{>0}$,
giving 
\[ g=\psi^2 g_1=\psi^2 g' + 4(d\psi)^2,
\]
with $\psi= e^{\frac t 2}$ the parameter on $\R^{>0}$.
This metric is clearly conical.
\endproof


\section{Vaisman manifolds: local properties}\index[terms]{manifold!Vaisman}


Here is the motivating example for all that follows:
 
\hfill

\example  Recall (\ref{_Sasaki_definition_}) that we defined a Sasakian\index[terms]{manifold!Sasaki}
manifold as a manifold with its Riemannian cone\index[terms]{cone!Riemannian} equipped with a K\"ahler
structure. For any given $\lambda\in \R^{>1}$, {the quotient
$C(S)/h_\lambda$ of a conical K\"ahler manifold associated
    with a Sasakian manifold is locally
conformally K\"ahler.}

%

\hfill

\definition \label{_Vaisman_Definition_}
An LCK manifold $(M, g, \omega, \theta)$ is called
{\bf Vaisman} if $\nabla\theta=0$, where
$\nabla$ is the Levi--Civita connection associated
with $g$.\index[terms]{manifold!Vaisman}

\hfill

\remark I. Vaisman\index[persons]{Vaisman, I.} called the LCK manifolds with parallel
Lee form\index[terms]{form!Lee!parallel} ``generalized Hopf''. Unfortunately, not all Hopf
manifolds are Vaisman, and this name did not stick. The 
name ``Vaisman manifold'' was coined in \cite{do}, but 
both names still circulate.\index[terms]{manifold!Vaisman}

\hfill

\remark \label{gau_vai}
Being parallel, the Lee form\index[terms]{form!Lee} $\theta$
of a Vaisman metric is\index[terms]{metric!Vaisman}
co-closed, that is, satisfies $d^*\theta=0$ 
(see Subsection \ref{_gauduchon_basics_Subsection_})
and thus a Vaisman metric is a Gauduchon
metric, see Exercise \ref{_Gauduchon_co-closed_Exercise_}
and Chapter  \ref{_appendix_gauduchon_}. Therefore, it is unique in
its conformal class, up to a constant
multiplier.\index[terms]{metric!Gauduchon}

\hfill

\remark 
	As $\theta$ is parallel, its norm is constant on (connected) $M$, and hence, after rescaling the metric, we may assume $|\theta|=1$. In particular, the Lee field\index[terms]{Lee field} $\theta^\sharp$\footnote{The dual $\theta^\sharp$ is taken with respect to the LCK metric.\index[terms]{metric!LCK} Note that other authors prefer to denote the Lee field by $\xi$ or $T$.} has
        no zeros,thus, if $M$ is compact,
        $\chi(M)=0$, this being the first obstruction to
        the existence of Vaisman metrics on compact
        complex manifolds.

\hfill

\corollary \label{_Lee_field_Killing_Corollary_}
If an LCK manifold is Vaisman then its Lee field\index[terms]{Lee field} 
$\theta^\sharp$ is Killing.\index[terms]{vector field!Killing}

\hfill

\proof  If $M$ is Vaisman, then $\nabla\theta=0$ and, in
\ref{homoth} (ii), 
$\la=0$, and
hence, from (iii), $\Lie_{\theta^\sharp} g=0$.
\endproof

\hfill

\remark 
On compact K\"ahler manifolds, any Killing field\index[terms]{vector field!Killing}\index[terms]{vector field!holomorphic} is
holomorphic; this follows from Hodge theory, see
\cite{_moroianu:book_}. The property is also true on compact LCK
manifolds\index[terms]{manifold!LCK} that are  neither Hopf manifolds nor have a hyperk\"ahler universal covering (\cite{_Moro_Pilca_}), but its proof is too complicated to be presented here. 
On Vaisman manifolds\index[terms]{manifold!Vaisman}, the Lee field is
holomorphic (and the result is local). \index[terms]{vector field!holomorphic}	

\hfill

\proposition \label{_Lee_field_holo_Proposition_}
On a Vaisman manifold $(M,g,\omega,\theta)$, the Lee field\index[terms]{Lee field}
is holomorphic: $\Lie_X I=0$. Moreover, its K\"ahler cover
$(\tilde M,\tilde g,\tilde \omega)$ is locally isomorphic
to the K\"ahler cone over a Sasakian manifold.\index[terms]{manifold!Sasaki}

\hfill

\pstep  
Here we briefly repeat the 
argument from \ref{_nabla(X)=Id_then_cone_Theorem_}.
Since $\nabla\theta=0$, by
\ref{_parallel_field_de_Rham_Corollary_}, $M$ is locally a
product $S\times \R$, with $g$ being (locally) a product
metric $g_S+dt^2$ and $\frac{d}{dt}=\theta^\sharp$. Then
$\theta=dt$.  Consider a K\"ahler cover
$\tilde M$ on which $\theta$ is exact, $\theta=d\f$. 
After adding a constant, we may assume that $\f=t$.
Then the K\"ahler metric can be written as $\tilde
g=e^{-\f}g$ and the manifold
$(\tilde M, \tilde g)$ is a Riemannian cone\index[terms]{cone!Riemannian}
of $S=\phi^{-1}(c)$. The
vector field $X=\theta^\sharp$ acts  by homotheties with
respect to $\tilde g$, giving\footnote{Note that here, as everywhere else,
we dualize $\theta\mapsto\theta^\sharp$ using the metric $g$.}
 $\Lie_X\tilde g=\la\cdot \tilde g$ for some constant
$\la$.

\hfill

{\bf Step 2:}  We show that $\theta^\sharp$
is holomorphic.   Denote with $\tilde\nabla$ the Levi--Civita connection of
$\tilde g$. By \ref{homoth} and \ref{_A_X_defi_Theorem_} (ii),
$\tilde\nabla(X)I=0$.

On the other hand, $\tilde\nabla_X I=0$ because $\tilde g$
is K\"ahler and the complex structure on a  K\"ahler manifold 
is parallel. However,
$\Lie_YI=\tilde\nabla_YI-\tilde\nabla(Y)I$ for each vector
field $Y$, as
$\tilde\nabla$ is symmetric. This gives $\Lie_XI=0$.

\hfill

{\bf Step 3:}
Sasakian manifolds are defined as manifolds
with the K\"ahler structure on the Riemannian
cone such that the Euler vector fields is holomorphic.
Therefore, $\tilde M$ is by definition a cone over
a Sasakian manifold.
\endproof

\hfill

\definition \label{_conical_Kahler_Definition_}
Recall that {\bf a conical K\"ahler manifold}
is a K\"ahler manifold $\tilde M$ equipped with an action of 
$\R^{>0}$ acting by holomorphic K\"ahler homotheties.
By definition, this is equivalent to $\tilde M$ being
locally holomorphically isometric to 
a cone over a Sasakian manifold.

\hfill

\corollary \label{vai_cone}
Let $M$ be a Vaisman manifold\index[terms]{manifold!Vaisman} and let $\tilde M$ be a
K\"ahler cover on which $\theta$ is exact. Then
$\tilde M$ is a K\"ahler cone and $\theta^\sharp$ is a
radial vector field acting by holomorphic homotheties.\index[terms]{vector field!radial}
\endproof

\hfill

\corollary \label{prop_lee}
On a Vaisman manifold, $[\theta^\sharp, I\theta^\sharp]=0$
and the complex vector field $\theta^\sharp-\1
I\theta^\sharp$ is a holomorphic $(1,0)$-field.
Consequently, the foliation $\Sigma$ generated 
by ${\theta^\sharp}$ and $I\theta^\sharp$ is holomorphic. 

\hfill

\proof
The field $\theta^\sharp-\1	I\theta^\sharp$  is
holomorphic because its real part is holomorphic,
and the corresponding foliation is obtained by
exponentiating this vector field. \endproof

\hfill

\definition 
	The above foliation $\Sigma$ 
(also denoted $\caf$) is called the {\bf canonical} or {\bf vertical foliation}.
\index[terms]{foliation!canonical}

\hfill

\remark\label{_canon_foli_totally_geodesic_Remark_}
The real vector fields $\theta^\sharp$ and $I\theta^\sharp$ generating the
canonical foliation are Killing\index[terms]{foliation!canonical} (\ref{_Lee_field_Killing_Corollary_}, 
\ref{_Lee_field_holo_Proposition_} and
\ref{_Reeb_Sasakian_Theorem_}). Setting $A=X$ in the equation  
\[ 
g(\nabla_AX, B)= - g(A, \nabla_BX)
\]
defining the Killing fields\index[terms]{vector field!Killing}, we get 
\[
g(\nabla_XX, B)= - g(X, \nabla_BX)=-\frac 12 \Lie_B(g(X,X))=0,
\]
and hence $g(\nabla_XX, B)= 0$ for all $B$.
Then $\nabla_XX= 0$, thus  the 
trajectories of a Killing field of constant length are
geodesics. Therefore, the canonical foliation in
a Vaisman manifold\index[terms]{manifold!Vaisman} is totally geodesic.\index[terms]{foliation!totally geodesic}

\hfill

Finally, we show that the quotients of 
conical K\"ahler manifolds by holomorphic homotheties
are always Vaisman.

\hfill

\theorem  \label{halfstr}
	Let $(\tilde M, \tilde g, \tilde \omega) = S\times \R^{>0}$ 
	be a conical K\"ahler manifold,
	$\Z =\langle \gamma \rangle$ a group acting on
	$\tilde M = C(S)$ by K\"ahler homotheties, and $t:\; C(S)\arrow
	\R^{>0}$ the projection map. Then 
		the form $\omega:= t^{-2}\tilde\omega$ is an
		LCK form\index[terms]{form!LCK} on $M:= \tilde M/ \langle \gamma \rangle$,
		its Lee form\index[terms]{form!Lee} is $- t^{-1} dt$, and $\nabla\theta=0$.
		Here $\nabla$ is the Levi--Civita connection on 
		$(M, g)$, and $g=t^{-2} \tilde g$.
In particular, this implies that $(M, \omega, \theta)$ is Vaisman.

\hfill

\proof   Indeed, $t^{-2} \tilde g$ is the product
metric on $C(S)= S\times \R$ (see \ref{cylcon}), where $s=\log (t^2)$ is the
coordinate on $\R$ and $ds= t^{-1} dt$ the unit covector. 
To find $\theta$, notice that 
\[ d\omega=d(t^{-2}\tilde\omega)= -2t^{-3}dt\wedge \tilde
\omega= -2t^{-1}dt \wedge \omega= - 2 ds\wedge \omega.
\]
This gives $\theta = - 2 d(s) = - d\log(t^2)$.
Then $\nabla(ds)=0$, because it is the unit covector on 
the $\R$ component of $(C(S), g) =S \times \R$.
\endproof

\section{Vaisman metrics obtained from holomorphic automorphisms}

The existence of a Vaisman metric in a conformal class is
a very strong condition.  Note that a Vaisman metric is,\index[terms]{metric!Vaisman}
in particular, a Gauduchon metric (see \ref{gau_vai} and 
Subsection \ref{_gauduchon_basics_Subsection_}) and\index[terms]{metric!Gauduchon}
hence, if it exists,  it is unique (up to homothety) in a
given conformal class. 

The following theorem was stated and proven in 
\cite{kor} in the compact case, but most of the arguments
used here do not rely on compactness. Here we give a new
proof, which works for non-compact LCK manifolds too.\index[terms]{manifold!LCK}

\hfill

\theorem \label{kami_or}
Let $(M,\omega, \theta)$ be an LCK manifold\index[terms]{manifold!LCK} equipped with a 
holomorphic and conformal $\C$-action $\rho$ without fixed points,
which lifts to non-isometric homotheties on 
the K\"ahler cover $\tilde M$. {Then $(M,\omega, \theta)$
	is conformally equivalent to a Vaisman manifold.}\index[terms]{manifold!Vaisman}

\hfill

{\bf Proof. Step 1:}
Let $\vec r\in TM$ be a real vector field tangent to this $\C$-action,
and let $\tilde \omega$ be the K\"ahler form\index[terms]{form!K\"ahler} of $\tilde M$.
Then $\Lie_{\vec r} \tilde\omega=a\tilde\omega$, $a\neq 0$, and
$\Lie_{I\vec r} \tilde\omega=b\tilde\omega$, with $a, b \in \R$ . 
Replacing $\vec r$ with an appropriate linear combination
of $\vec r$ and $I\vec r$, we can assume that
	$\Lie_{I\vec r} \tilde\omega=0$ and 
	$\Lie_{\vec r} \tilde\omega=2 \tilde\omega$. 
To do this, let $u= I\vec r- b a^{-1} \vec r$.
Then 
\[ \Lie_u \tilde\omega= b\tilde \omega - b a^{-1}  a\tilde\omega= 0,
\]
and 
\[ 
\Lie_{Iu} \tilde\omega=  - \Lie_{\vec r} \tilde\omega - ba^{-1} \Lie_{I\vec r}
- a \tilde\omega  - b^2 a^{-1} \tilde\omega= K\tilde\omega,
\]
where $K=-a- b^2a^{-1}$.
If $K=0$, then $a^2+b^2=0$, that is  impossible because $a\neq 0$.
If we replace $\vec r$ with $2K^{-1} Iu$, we obtain,
precisely, that $\Lie_{I\vec r} \tilde\omega=0$ and 
	$\Lie_{\vec r} \tilde\omega=2 \tilde\omega$.

\hfill

{\bf Step 2:}
Let 
$\eta:= I(i_{\vec r}\tilde\omega )=i_{I\vec r}\tilde\omega$.
Then  $\Lie_{I\vec r} \tilde \omega=d\eta=0$.
Since $\Lie_{\vec r} \tilde \omega = 
d(i_{\vec r}\tilde \omega )= d^c\eta= 2 \tilde \omega$,
one has $d^c \eta = 2 \tilde \omega$.

Now \[
2\eta = \Lie_{\vec r} I(i_{\vec r}\tilde\omega )=
 \Lie_{\vec r}\eta=d(i_{\vec r} \eta)
=d(\tilde\omega(I\vec r, \vec r)).
\]
This gives $4\tilde \omega= dd^c \phi$, where 
$\phi=\tilde\omega(\vec r, I\vec r)=|\vec r|^2$. 

\hfill

{\bf Step 3:}
We obtained that the action of the holomorphic
vector field $\vec r$ multiplies the
K\"ahler potential $\phi$ by a constant, and the action by $I(\vec r)$ 
preserves $\phi$. It remains to show that 
$(\tilde M, \tilde \omega)$ is
locally isometric to a K\"ahler cone, or, equivalently,
that $(M, f\omega)$ is locally isometric to a cylinder
for some positive function $f\in C^\infty M$.

Notice that $\vec r$ multiplies $\phi$ by the same
constant factor as $\tilde\omega$. The function 
$\phi=\omega(\vec r, I\vec r)$ is
strictly positive everywhere because the $\C$-action $\rho$
has no fixed points. Then $g_0 :=
\phi^{-1}\tilde\omega(\cdot, I\cdot)$ is an Hermitian form
on $\tilde M$, that is  $\rho$-invariant.\index[terms]{form!Hermitian}

Therefore, 
$\rho$ acts on the Hermitian manifold $(\tilde M, g_0)$ by holomorphic
isometries. We are going to prove that $\nabla_0\eta=0$,
where $\nabla_0$ is the Levi--Civita connection 
associated with $g_0$.

Let $S_c= \phi^{-1}(c)$ be a level set of $\phi$. Using
Sard's theorem,\index[terms]{theorem!Sard} we can be sure that $c$ is a regular value
of $\phi$ and $S_c$ is smooth.
The vector field $\vec r$ is orthogonal to $S_c$ because
$g(\vec r, \cdot) = -\eta(\cdot)= -(d\phi)(\cdot)$. The action of $\vec r$
on $M$ maps $\phi$ to $\const\cdot \phi$. 
This gives a local decomposition $M= \R \times S_c$.
Since $\vec r$ is orthogonal to $S_c$, and
all $S_c$ are naturally isometric with respect
to $\omega_0$, to prove that $(M, g_0)$
is isometric to a product it remains to show
that $g_0(\vec r, \vec r)=\const$.
However, $g_0(\vec r, \vec r)=\phi^{-1}g(\vec r, \vec r)=
\phi^{-1}\omega(\vec r, I\vec r)=1$.
\endproof	

\hfill

The formula  $4\tilde \omega= dd^c \phi$, where
$\phi=\tilde\omega(\vec r, I\vec r)$ (Step 2),
has the following useful corollary.

\hfill

\corollary\label{_tilde_M_potential_Corollary_}
Let $(M, \omega,\theta)$ be a Vaisman manifold,\index[terms]{manifold!Vaisman}
and $(\tilde M, \tilde \omega)$ its K\"ahler cover
locally isometric to the cone $C(S)= S \times\R^{>0} $ over a Sasakian
manifold, where $t$ is the parameter on 
$\R^{>0}$. Then $4\tilde \omega= d d^c t^2$.

\hfill

\proof
Choose $\vec r = t \frac {d}{dt}$.
Then $\Lie_{\vec r}\tilde \omega=2\tilde \omega$
because $\tilde \omega$ is homogeneous of order 2.
Therefore, $\phi=\tilde\omega(\vec r, I\vec r)= t^2$,
which gives $4\tilde \omega= d d^c t^2$
as indicated in \ref{kami_or}, Step 2. \endproof

\hfill

\example
In Exercise \ref{_example_of_isometry_shells_Exercise_},
we construct an example of an isometric flow on an LCK Hopf manifold $H$\index[terms]{manifold!Hopf}
which lifts to an isometry of the K\"ahler metric on its universal cover.
This is an example of an isometric flow that does not satisfy the
assumptions of \ref{kami_or}. 
When $H$ is a classical Hopf manifold,
$H = \frac{\C^n \backslash 0}{\Lambda \Id}$, and the metric on $H$
is given by $\omega = \frac{dd^c l}{l}$, where $l(z)=|z|^2$,
any $u\in \U(n)$ will act on $H$ and its universal cover by isometries.


\section{The canonical foliation on compact Vaisman manifolds}
\label{canf}\index[terms]{manifold!Vaisman}\index[terms]{foliation!canonical}


On compact Vaisman manifolds, the canonical foliation
$\Sigma$ has  specific properties. In the next
theorem we collect (and rephrase) results from \cite{tsu},
\cite{tsu2}, \cite{_Verbitsky:Vanishing_LCHK_}.\index[terms]{foliation!canonical}

\hfill

\theorem \label{_Subva_Vaisman_Theorem_} 
Let $M$ be a compact Vaisman manifold,\index[terms]{manifold!Vaisman} and 
$\Sigma\subset TM$ its canonical foliation.\index[terms]{foliation!canonical} 
We rescale the metric on $M$ in such a way that $|\theta|=1$.
\begin{description}
\item[(i)] $\Sigma$ is
	independent on  the choice of the Vaisman metric.
\item[(ii)] There exists a semi-positive, exact
	(1,1)-form $\omega_0$ with $\Sigma=\ker \omega_0$.
\item[(iii)] For any complex compact subvariety\index[terms]{subvariety}
	$Z\subset M$, $Z$ is tangent to $\Sigma$. 
\item[(iv)] For any compact complex subvariety
	$Z\subset M$, the set of smooth points of $Z$ is
  Vaisman\footnote{One could define the notion of {\em
      Vaisman variety} when the variety is possibly
    singular. In that case, all subvarieties of a Vaisman 
variety are Vaisman.}.
\end{description}


\proof Let $\tilde M = C(S)$ be the
conical K\"ahler manifold that covers $M$
(\ref{vai_cone}), and 
$\phi$  the K\"ahler potential\index[terms]{potential!K\"ahler} on $C(S)$, written
explicitly as $t^2$ (\ref{_Sigma_defi_Claim_})
 Then $\tilde \omega_0:=dd^c\log\phi$ is a pseudo-Hermitian form which\index[terms]{form!pseudo-Hermitian}
	vanishes on $\Sigma$ and is positive on
        $TM/\Sigma$ (\ref{_Sigma_defi_Claim_}).
The deck transform group acts on $\log\phi$ 
by adding a constant. Therefore, $d\log\phi$
is deck transform invariant, and the form
$\tilde \omega_0=d(I d\log\phi)$
is the pullback of a form $\omega_0\in \Lambda^{1,1}M$. This proves (ii).

To prove (i), note that the zero foliation (annihilator) of
$\omega_0$ is independent on  the choice of the Vaisman
metric. Indeed, if there are two Vaisman metrics
with fundamental forms $\omega_0$ and $\omega_0'$ vanishing on different 
1-dimensional complex foliations, the sum 
$\omega_0+\omega'_0$ would be positive definite.
However, $\int_M( \omega_0+\omega_0')^{\dim_\C M}$
vanishes, because $\omega_0$ and $\omega_0'$ are exact.
{Since $\Sigma=\ker\omega_0$, $\Sigma$ 
	is then independent on  the Vaisman structure.}

Now, for any compact subvariety
$Z\subset M$, the integral $\int_Z\omega_0^{\dim_\C Z}$
vanishes, because $\omega_0$ is exact. {Therefore,
	$\omega_0\restrict {TZ}$ has one zero eigenvalue at each
	point of $Z$.} This means precisely that $\Sigma\subset TZ$
at this point, and hence (iii).

As for (iv), since the Lee field\index[terms]{Lee field} is tangent
to $Z$, the covering $\tilde Z\subset C(S)$
is preserved by the homotheties. Therefore, it is also a
conical K\"ahler manifold. Then its image $Z\subset M$
	is Vaisman (\ref{kami_or}). \endproof

\hfill

\corollary\label{subvai}
A closed complex submanifold of a compact Vaisman manifold\index[terms]{manifold!Vaisman} has an induced Vaisman structure. Moreover, it contains the leaves of the canonical foliation intersecting it.\index[terms]{submanifold}
\endproof\index[terms]{foliation!canonical}

\hfill

\remark\label{_Sigma_transv_K_} 
In the terminology of \cite{to} (see also \cite{bog}),
from  \ref{_Subva_Vaisman_Theorem_} (ii) it follows that $\Sigma$ is a {\bf transversally
  K\"ahler} foliation.\index[terms]{foliation!transversally K\"ahler}

%
%
%
%

\section{Exercises}

\begin{enumerate}[label=\textbf{\thechapter.\arabic*}.,ref=\thechapter.\arabic{enumi}]

\item Let $\theta$ be a 1-form on a Hermitian manifold $M$,
with $\dim_\C M=n$. Prove that
$*(\theta)= c I(\theta)\wedge \omega^{n-1}$,
where $c$ is a constant. Find this constant.

\item 
Let $S$ be a Riemannian manifold, $C(S)$ its cone.
Assume that $C(S)$ is a flat Riemannian manifold.
Prove that $S$ is locally isometric to a 
sphere of constant curvature.
Is the converse true?

\item
Let $L$ be a positive line bundle on
a complex projective manifold $X$.
Denote by $M$ the quotient of the
space $\Tot^\circ(L)$  of all non-zero vectors in $\Tot(L)$
by a homothety $v \mapsto \lambda v$, $\lambda \in \C$,
$|\lambda| > 1$. Prove that $M$ is Vaisman.

\item Prove that all compact complex curves in a Vaisman manifold\index[terms]{manifold!Vaisman}
are elliptic curves.

\item Prove that a blow-up of a Vaisman manifold 
cannot be Vaisman.

\item\label{_Hopf_diffe_S^1_x_S^2n-1_Exercise_}
Denote by $\psi^N$ the $N$-th iteration of $\psi$.
A {\bf holomorphic contraction} \index[terms]{holomorphic contraction}
of $\C^n$ is an invertible map $\psi:\; \C^n \arrow \C^n$
such that for each compact subset $K\subset \C^n$
and open subset $U\subset  \C^n$ containing 0,
there exists $N$ such that $\psi^N(K)\subset U$.
A {\bf Hopf manifold}\footnote{This is the {\em most\index[terms]{manifold!Hopf}\index[terms]{manifold!Hopf!linear}
    general} definition of a Hopf manifold, that is 
due to \index[persons]{Kodaira, K.} Kodaira, \cite[p. 694]{_Kodaira_Structure_II_}. In this book, we shall sometimes use
a less general notion, such as {\bf a linear Hopf
  manifold}.} is a quotient of $\C^n \backslash 0$
by a holomorphic contraction. Prove that a Hopf
manifold is diffeomorphic to $S^{2n-1}\times S^1$.

\item \label{_Hopf_surface_explicit_Exercise_}
Let $A\in \End(\C^2)$
be an Hermitian matrix with two eigenvalues 
$\alpha,\beta$ such that $|\alpha| = |\beta| > 1$.
Prove that the function $\phi(z):=|z|^2$ on $\C^2$ 
satisfies $A^*(\phi)= |\alpha|^2 \phi$.
Prove that the corresponding Hopf manifold
$\frac{\C^2\backslash 0}{\langle A \rangle}$
is Vaisman. 

{\em Hint:} Prove that $dd^c \phi$ defines a conical
metric on $\tilde M :=\C^2\backslash 0$, and the linear 
vector field $v\mapsto \log A(v)$ acts on $\tilde M$ by
holomorphic homotheties.

\item
In the assumptions of Exercise
\ref{_Hopf_surface_explicit_Exercise_},
prove that $\ker (dd^c \log \phi)$
defines the canonical foliation \index[terms]{foliation!canonical}$\Sigma$ on $M$.
Show that for any elliptic curve $C\subset M$,
its preimage $\tilde C \subset  \C^2\backslash 0$
is $\C^*$-invariant, where $\C^*$ acts on
$ \C^2\backslash 0$ in a standard way,
$x \mapsto t x$.

\item
In the assumptions of Exercise
\ref{_Hopf_surface_explicit_Exercise_},
let $C\subset M$ be an elliptic curve,
and $\tilde C \subset  \C^2\backslash 0$
its preimage. Prove that  $\tilde C=\C^2\backslash 0\cap W$,
where $W$ is an eigenspace of $A$. Prove that
$M$ contains precisely 2 or infinitely many elliptic curves.

{\em Hint:}
Prove that $\tilde C$ is $A$-invariant.

\item
Let $g$ be a Vaisman metric on a manifold 
$M$ with $\pi_1(M)=\Z$, and $\theta^\sharp$ the corresponding
Lee field. Prove that the anti-Lee field\index[terms]{Lee field} $I(\theta^\sharp)$\index[terms]{Lee field!anti-}
has precompact\footnote{Precompact means ``having a
  compact closure''.} orbits on the universal 
covering of $M$. Prove that a real linear 
combination $a\theta^\sharp + b I(\theta^\sharp)$
has precompact orbits if and only if $a=0$.

\item
Let $g_1, g_2$ be Vaisman metrics on a Hopf manifold 
$M= \frac{\C^n\backslash 0}{\Z}$,
and $\theta_1^\sharp, \theta_2^\sharp$  the corresponding
Lee fields\index[terms]{Lee field}. Prove that $\theta_1^\sharp, \theta_2^\sharp$
are proportional with a real coefficient.

{\em Hint:} Use the previous exercise.

\item (\cite{va_gd}) Let $\f(z)=|z|^2$ on $\C^n$ and
  $A(z):=\lambda z$, $\la\in\C$, $|\la|\neq
  0,1$. Show that the function $\f^{t+1}$ is
  plurisubharmonic for $t>-1$ and the 2-form
  $\omega_t:=\frac{dd^c\f^{1+t}}{\f^{1+t}}$ induces an LCK
  form on the Hopf manifold $\frac{\C^n\backslash
    0}{\langle A \rangle}$, that is  Vaisman only for
  $t=0$. Show that the Lee field \index[terms]{Lee field}of each LCK form
  $\omega_t$ does not depend on $t$.\index[terms]{form!LCK}



\item Show that on a Vaisman manifold\index[terms]{manifold!Vaisman} $M$, the distribution
  $\Sigma^\perp$ is not integrable, that is,
  $[\Sigma^\perp, \Sigma^\perp]\not\subset \Sigma^\perp$.
  Prove that $\rk [\Sigma^\perp, \Sigma^\perp]= \dim_\R M -1$.

\item Prove that a compact Vaisman
manifold is not a product of two positive-dimensional
complex manifolds.

{\em Hint:} Indeed, was it a product of two
complex manifolds, each one should contain the leaves of
the canonical foliation, contradiction.\index[terms]{foliation!canonical}

\item Let $M$ be a complex manifold with universal
covering biholomorphic to $\C^n$. Prove that any holomorphic 
map from a Hopf manifold to $M$ is trivial.

{\em Hint:} Use Hartogs' extension theorem.\index[terms]{theorem!Hartogs}

\item
A real $(p,p)$-form $\eta$ on a complex vector space is called
{\bf weakly positive} if for any collection
of Hermitian forms $\omega_1, ..., \omega_{n-p}$,
the volume form $\eta \wedge \omega_1 \wedge ... \wedge \omega_{n-p}$
is non-negative. These forms clearly form a cone; a form that belongs
to its interior is called {\bf strictly weakly positive}.
A differential form on a complex manifold is 
strictly weakly positive if it is strictly weakly positive pointwise.
Let $\eta$ be a strictly weakly positive 
$(p,p)$-form on a compact complex manifold $M$, and $\rho$ 
a non-zero $(n-p,n-p)$-form that can be obtained as a
product of positive (1,1)-forms. Prove that
$\int_M \eta \wedge \rho>0$.

\item
Let $\eta$ be a strictly weakly positive 
$(p,p)$-form on a compact Vaisman $n$-manifold $(M, \omega, \theta)$,
with $p < n-1$. Prove that 
\[
\int_M dd^c \eta \wedge \theta \wedge \theta^c\wedge \omega_0^{n-p-2}
=-\int_M \eta \wedge \omega_0^{n-p}<0.
\]

{\em Hint:} Use the formula $d^c\theta=\omega_0$ and
apply the previous exercise.

\item\label{_Vaisman_no_dd^c_closed_pos_p_p_Exercise_}
Let $(M, \omega, \theta)$ be a compact Vaisman $n$-manifold.
Prove that $M$ does not admit $dd^c$-closed
strictly weakly positive $(p,p)$-forms, $1 < p < n-1$.

{\em Hint:} Use the previous exercise.

\end{enumerate}


\chapter{The structure of compact Vaisman manifolds}\label{str_th}\index[terms]{manifold!Vaisman}


{\setlength\epigraphwidth{0.6\linewidth}
	\epigraph{\it Tels que les excr\'ements chauds d'un vieux colombier,\\
	Mille R\^eves en moi font de douces br\^ulures :\\
	Puis par instants mon coeur triste est comme un aubier\\
	Qu'ensanglante l'or jeune et sombre des coulures.}{\sc\scriptsize Arthur Rimbaud, \ \ Oraison du soir}
}
\section{Introduction}

The structure theorem for Vaisman manifolds
claims that any compact Vaisman manifold, considered
as a complex manifold, is a quotient
of a Riemannian cone\index[terms]{cone!Riemannian} over a compact Sasakian manifold.

This is a global version of \ref{_Lee_field_holo_Proposition_},
which claims that the K\"ahler cover of 
any Vaisman manifold (compact or not) is 
locally isomorphic to the cone over a Sasakian manifold.

The original proof of this result (\cite{ov_str})
had an error that we fixed in \cite{ov_jgp_16}.

The actual statement of the structure theorem
and its proof use the notion of ``LCK rank'' of\index[terms]{rank!LCK}
an LCK manifold. Recall that an LCK manifold\index[terms]{manifold!LCK}
is a manifold $M$ equipped with a positive, closed
(1,1)-form taking values in an oriented rank 1 local system.
Denote the monodromy group\index[terms]{group!monodromy} of this local
system by $\Gamma$; by definition, $\Gamma$
is a subgroup of the group $\R^{>0}$
of automorphisms of the fiber.
Therefore, $\Gamma$ is a torsion-free\index[terms]{group!torsion-free}
abelian group. Its rank is called
{\bf the LCK rank} of $M$.
We gave the same notion 
in different wording in \ref{lck_rank}.

Let $C(S)$ be the cone over
a Sasakian manifold, considered
as a K\"ahler manifold, and let $\Lambda\simeq \Z$
act on $C(S)$ properly discontinuously 
by holomorphic homotheties. 
Such an action is always induced by
an isometry action on $S$ and a shift along
the $\R^{>0}$ direction in $C(S) = S\times \R^{>0}$ 
(see the exercises to this chapter).
By \ref{halfstr}, the quotient 
$C(S)/\Lambda$ is Vaisman.

It is not very hard to prove that 
a compact Vaisman manifold \index[terms]{manifold!Vaisman}of LCK
rank 1 is always holomorphically
isometric to a manifold obtained 
using this construction (\ref{str_vai}).

In this chapter we deal with 
compact Vaisman manifolds of LCK rank\index[terms]{manifold!LCK} $>1$
proving that the set of Vaisman metrics
of rank 1 on a given complex manifold
of Vaisman type is dense in the set
of all Vaisman metrics
(\ref{_Vaisman_defo_transve_Proposition_}).
This result uses a harmonic decomposition\index[terms]{decomposition!harmonic}
(\ref{_harmo_deco_1-form_Proposition_}),
that is  a special case of the harmonic decomposition theorem of
\index[persons]{Kashiwada, T.} Kashiwada and \index[persons]{Vaisman, I.} Vaisman, proved in Chapter \ref{_Harmonic_forms_chapter_}.

\section[The Vaisman metric expressed through
the Lee form]{The Vaisman metric expressed in terms of
  the Lee form}\index[terms]{form!Lee}

From now on, we denote $I(\theta)$ by $\theta^c$.
Further on, we shall use the following formula.

\hfill

\proposition
Let $(M, \omega, \theta)$ be a Vaisman manifold.\index[terms]{manifold!Vaisman}
Then for an appropriate constant $u\in \R^{>0}$, we have
$u \omega= d^c \theta + \theta \wedge \theta^c$.

\hfill

\proof
Let $(\tilde M, \tilde \omega)$ be a K\"ahler cover of $M$.
Then $\tilde M= C(S)$ is locally a cone of a Sasakian
manifold, with $\theta=-t^{-2}dt^2$ and
$4 \tilde\omega= dd^c t^2$
(\ref{_Lee_field_holo_Proposition_}, \ref{_tilde_M_potential_Corollary_}).
By construction, we have $\omega = \lambda t^{-2}\tilde \omega$,
where $\lambda$ is a constant.
Then  
\[ d^c\theta = - d^c (t^{-2} dt^2)= t^{-2} d d^c (t^2) - t^{-4} dt^2\wedge
d^c t^2=  4 t^{-2} \tilde \omega - \theta \wedge \theta^c.\]
This gives $4 \lambda^{-1}\omega = d^c\theta +\theta \wedge \theta^c$.
\endproof

\hfill

\remark \label{_scaling_Vaisman_remark_}
If we replace $\omega$ by $u \omega$, 
the metric remains Vaisman and the Lee form\index[terms]{form!Lee} remains the same
(however, its length changes).
Then 
\begin{equation}\label{_omega_via_theta_Chapter_8_Equation_}
\omega= d^c \theta +\theta \wedge \theta^c.
\end{equation}
This is equivalent to $|\theta| =1$, as the following
result implies. 

\hfill

\proposition\label{_preferred_gauge_Vaisman_theta=1_Proposition_}
Let $(M, \omega, \theta)$ be a Vaisman manifold,\index[terms]{manifold!Vaisman}
with the metric rescaled in such a way that
$\omega= d^c \theta + \theta \wedge \theta^c.$
Then $|\theta|=1$. Conversely, $|\theta|=1$
implies $\omega= d^c \theta + \theta \wedge \theta^c.$

\hfill

\pstep
Let $(\tilde M, \tilde \omega)$ be a K\"ahler cover of $M$,
and $\phi=4^{-1}t^2$ the K\"ahler potential. By 
\ref{_tilde_M_potential_Corollary_}, $dd^c t^2 =4\tilde \omega$, and hence 
$dd^c \phi =\tilde \omega$.
By \ref{_Subva_Vaisman_Theorem_}, the (1,1)-form
$\omega_0:=dd^c\log\phi$ is semi-positive, vanishes
on the leaves of the canonical foliation $\Sigma$,
and is positive in the transversal directions.\index[terms]{foliation!canonical} 
Since $\omega= \const \phi^{-1}\tilde\omega$,
one has 
\[ d \omega= \const d(\phi^{-1}) \wedge \tilde\omega=
-\const \frac{d\phi}{\phi^2} \wedge \tilde\omega=
- \frac{d\phi}{\phi}\wedge \omega.
\]
This gives 
\begin{equation}\label{_theta_d_log_Equation_}
\theta=-d\log\phi.
\end{equation}
Now, 
\begin{equation}\label{_omega_0_Equation_}
\omega_0 = dd^c\log\phi= - d^c d\log\phi = d^c\theta.
\end{equation}
because $\theta=-d\log\phi$.

\hfill

{\bf Step 2:} Suppose that $\omega=  d^c \theta + \theta \wedge \theta^c.$
Then $\omega= \theta \wedge \theta^c$ on $\Sigma$.
Since $\Sigma$ is generated by the Lee\index[terms]{Lee field} and the anti-Lee fields\index[terms]{Lee field!anti-}
$\theta^\sharp$ and $I \theta^\sharp$, this gives
\[ |\theta|^2=\omega(\theta^\sharp, I\theta^\sharp) = 
\theta \wedge I\theta (\theta^\sharp,  I\theta^\sharp) =
|\theta|^4,
\]
implying $|\theta|=1$. The converse implication is proven
by the same argument.
\endproof


\section[Decomposition for harmonic 1-forms on Vaisman manifolds]{Decomposition for harmonic 1-forms on\\ Vaisman manifolds}\index[terms]{manifold!Vaisman}


We prove a theorem that is  a particular case of
\index[persons]{Kashiwada, T.} Kashiwada and \index[persons]{Vaisman, I.} Vaisman's harmonic decomposition theorem (see Chapter 
\ref{_Harmonic_forms_chapter_}, for the general statement; 
for a different proof in degree 1, see \cite{_Madani_Moroianu_Pilca:toric_}).\index[terms]{decomposition!harmonic}

\hfill

We will need the following lemma, that is  almost trivial.

\hfill

\lemma \label{_invariance_for_cohomology_Vaisman_Lemma_}
Let $M$ be a compact Vaisman manifold,\index[terms]{manifold!Vaisman} and $G$
the closure of the Lie group generated by the
action of the Lee and the anti-Lee fields\index[terms]{Lee field} on $M$.\index[terms]{Lee field!anti-}
Then $G$ is a compact torus, and any harmonic form
on $M$ is $G$-invariant.

\hfill

\proof Since the Lee and anti-Lee fields\index[terms]{Lee field} are acting on $M$\index[terms]{Lee field!anti-}
by isometries (\ref{_canon_foli_totally_geodesic_Remark_}), 
the group $G$ is a closed subgroup of the
isometry group of $M$, and the latter is compact.
Then $G$ is compact; since it is the closure of an
abelian group, it is abelian.

Now let  $\eta$ be a harmonic form. Then for 
any $g\in G$, the form $g^* \eta$
is also harmonic, because $g$ acts on $M$ by isometries.
However, a path-connected topological group 
acts trivially on the cohomology of any topological space,
hence $g^*\eta$ is cohomologous to $\eta$.
Then $\eta = g^*\eta$ because both forms are
harmonic representatives of the same cohomology class.
\endproof

\hfill

By \cite[Corollary 2.9]{_Verbitsky:G2_forma_}
a product of a harmonic form and a parallel
form on any Riemannian manifold is harmonic,
see also
\cite{_Chern:Chern_,_Lichnerowicz:colloque_,_Weil:Chern_,_Sampson:Chern_}.\footnote{We are very grateful to Robert Bryant for 
these references and for his insightful comments  on the Mathoverflow topic
 {\scriptsize\url{https://mathoverflow.net/questions/396600/reference-to-a-theorem-about-a-product-of-harmonic-and-parallel-forms/396647}}
which shed more light on the history of this result.}
We give a simple proof of this result for 
parallel 1-forms on a compact Riemannian manifold.

\hfill

\lemma\label{_product_harmonic_parallel_Lemma_}
Let $\eta$ be a harmonic form on 
a compact Riemannian manifold, and
$\beta$ a parallel 1-form.
Then $\eta \wedge \beta$ is also harmonic.

\hfill

\proof
Let $L_\beta:\; \Lambda^i(M) \arrow \Lambda^{i+1}(M)$
be the map $\alpha\stackrel {L_\beta} \mapsto \beta \wedge \alpha$,
and $\{\cdot, \cdot\}$ denote the anticommutator.
Then 
\[ d^*(\eta \wedge \beta)= - L_\beta(d^*\eta) + \{L_\beta, d^*\}(\eta)
\]
hence $d^*(\eta \wedge \beta)= 0$ would follow if we prove
$\{L_\beta, d^*\}(\eta)=0$. The vector field
$\beta^\sharp$ is parallel, and hence \index[terms]{vector field!parallel} Killing.\index[terms]{vector field!Killing}
Since harmonic forms are 
invariant with respect to a connected group acting by
isometries (\ref{_invariance_for_cohomology_Vaisman_Lemma_}),
the equality $\{L_\beta, d^*\}(\eta)=0$ would follow if we show that 
\begin{equation}\label{_d^*_L_eta_anticommu_Equation_}
\{L_\beta, d^*\}=-\Lie_{\beta^\sharp}.
\end{equation}
Let $i_{\beta^\sharp}: \Lambda^i(M) \arrow \Lambda^{i-1}(M)$ 
be the contraction operator, $\alpha \mapsto i_{\beta^\sharp} \alpha$.
Then $\{i_{\beta^\sharp}, d\}= \Lie_{\beta^\sharp}$
by the Cartan formula. However, $(i_{\beta^\sharp})^* = L_\beta$, 
which gives 
\[
(\Lie_{\beta^\sharp})^* = (\{i_{\beta^\sharp}, d\})^*= \{L_\beta, d^*\}.
\]
The Lie derivative along a Killing field\index[terms]{vector field!Killing} is anti-selfdual, 
because $[\Lie_{\beta^\sharp}, *]=0$,
hence $(\Lie_{\beta^\sharp})^*= -\Lie_{\beta^\sharp}$, which
implies \eqref{_d^*_L_eta_anticommu_Equation_}.
\endproof

\hfill

The following theorem is a version of
the harmonic decomposition theorem proven by
\index[persons]{Kashiwada, T.} Kashiwada and \index[persons]{Vaisman, I.} Vaisman in \cite{kashiwada_kodai,va_gd}.
However,  the Hodge decomposition
on the transversally harmonic forms is due to
\index[persons]{Tsukada, K.} Tsukada, \cite{tsuk}.

\hfill

\proposition\label{_harmo_deco_1-form_Proposition_}
Let $(M, \omega, \theta)$ be a compact Vaisman manifold and 
$\alpha$ a harmonic 1-form on $M$. Then\index[terms]{manifold!Vaisman} 
$\alpha = a\theta + \rho$, where $a\in \C$ is a constant,
and $\rho$ a 1-form that vanishes on 
the canonical foliation\index[terms]{foliation!canonical} $\Sigma$ and 
satisfies $d\rho=d^c\rho=0$.\footnote{In other words,
$\rho$ is a basic form with respect to $\Sigma$, see 
\ref{_basic_form: definition_}.}

\hfill

{\bf Proof. Step 1:} 
 From \ref{_product_harmonic_parallel_Lemma_}, we obtain that
the contraction operator 
$i_{\theta^\sharp}$ maps harmonic forms to harmonic forms.
Then $i_{\theta^\sharp}(\alpha)=\const$. Replacing $\alpha$ 
by $\rho=\alpha - \const \theta$, we can assume that
$\alpha$ is orthogonal to $\theta$. To finish the proof of
\ref{_harmo_deco_1-form_Proposition_} it remains to 
show the following. Let $\alpha$ be a harmonic 1-form such that 
the scalar product $(\alpha, \theta)$ vanishes pointwise.
Then  $d^c\alpha=0$ and $\alpha$ vanishes on the canonical foliation. \index[terms]{foliation!canonical}

\hfill

{\bf Step 2:} As usual, we write $\theta^c:= I(\theta)$.
As in \ref{_scaling_Vaisman_remark_}, we assume that $|\theta|=1$.
 We prove now that $(\alpha, \theta^c)=0$
for any harmonic 1-form $\alpha$. By Step 1, we may
assume, without restricting the generality, 
that $(\alpha, \theta)=0$. Let 
$\alpha= f \theta^c + \alpha_0$, where
$(\alpha_0, \theta^c)=0$, and 
$f= (\alpha, \theta^c)\in C^\infty M$.

Let $G$ be the group defined in
\ref{_invariance_for_cohomology_Vaisman_Lemma_}.
By this lemma, $\alpha$ is $G$-invariant;
the decomposition $\alpha= f \theta^c + \alpha_0$
is also $G$-invariant.
Then $\alpha_0$ is a $G$-invariant 1-form
that satisfies $(\alpha_0, \theta^c) = (\alpha_0, \theta)=0$,
hence vanishes on $\Sigma$. 
This implies that $\alpha_0$ is basic\index[terms]{form!basic}, and hence  $d\alpha_0$  also
vanishes on $\Sigma$. From 
$d(f \theta^c) = df \wedge \theta^c + f d\theta^c$ we obtain
\[
0 = d(\alpha) = d (f \theta^c + \alpha_0) = 
df \wedge \theta^c + f d\theta^c + d\alpha_0.
\]
Since $d \theta^c = \omega_0$ by \eqref {_omega_0_Equation_},
the last two terms vanish on $\Sigma$; this gives
$df \wedge \theta^c=0$, because this term is orthogonal
to forms vanishing on $\Sigma$. Then $df$ is proportional
to $\theta^c$, $df= f_1 \theta^c$ which implies that 
$d \theta ^c= d(f_1^{-1} df))$  is proportional  to
$\theta^c \wedge d(f_1^{-1})$ outside of  the
zero set of $f$. However, this is impossible
because $d\theta^c=  \omega_0$
is not divisible by $\theta^c$ anywhere on $M$.
We proved that $f=0$ and $\alpha$ is orthogonal to $\theta^c$.

\hfill

{\bf Step 3:}
It remains to show that $d^c \alpha=0$ for
any harmonic form $\alpha\in \ker \Delta$ that vanishes on 
the canonical foliation $\Sigma$.\index[terms]{foliation!canonical}

Since $\alpha$ is also $G$-invariant
(\ref{_invariance_for_cohomology_Vaisman_Lemma_}),
it is a basic form with respect to
$\Sigma$. Denote by $\Lambda^*_b(M)$ the algebra of 
basic forms on $M$.
By \ref{_basic_forms_Frobenius_lifted_Theorem_}, 
$\alpha$ is locally lifted from a form $\underline\alpha$
on the leaf space $U$ of $\Sigma$, giving (locally in $M$)
$\pi^*\Lambda^*(U)= \Lambda^*_b(M)$. The manifold $U$
is K\"ahler because $\omega_0$ defines a 
K\"ahler structure on the leaf space of $\Sigma$. 
Let $\Delta_b$ be the Laplacian on $\Lambda^*(U)$
lifted to $\Lambda^*_b(\pi^{-1}(U))$. 
Since the metric $\omega_0$ is constant
along the leaves of $\Sigma$, the gluing maps
between two local instances of the leaf space
are holomorphic isometries. Therefore, the operators
$\Delta_b$ agree on open patches
in $M$ where the leaf space exists.
 This allows us to extend it to the
map $\Delta_b:\; \Lambda^*_b(M)\arrow \Lambda^*_b(M)$.
We call this map {\bf the basic 
Laplacian}  (see also Chapter \ref{_Harmonic_forms_chapter_}).
We are going to show that
$\ker \Delta\cap \Lambda^*_b(M) \subset \ker \Delta_b$
and $\ker \Delta_b\subset \ker d^c$.
This would imply that $d^c \alpha=0$.

\hfill

{\bf Step 4:} We prove 
$\ker \Delta\cap \Lambda^*_b(M) \subset \ker \Delta_b$

Let $*_b:\; \Lambda^*_b(M)\arrow \Lambda^*_b(M)$
be the Hodge star operator on the leaf space $U$ lifted to $\Lambda^*_b(M)$.
An orthonormal basis $\zeta_1, ..., \zeta_{2n-2}$ in 
$(\Sigma^\bot)^*= \pi^*(\Lambda^1 U)$ 
gives 
 an orthonormal basis $\zeta_1, ..., \zeta_{2n-2}, \theta, \theta^c$
in $\Lambda^1 M$. Therefore, for any basic form $\beta$, we have
\begin{equation}\label{_*_0_explicit_Equation_}
*\beta = (*_b\beta)\wedge \theta\wedge\theta^c, \ \ 
\text{and} \ \  *_b\beta=*(\beta\wedge \theta\wedge\theta^c)
\end{equation}
Denote by $d_b^*:\; \Lambda^*_b(M)\arrow \Lambda^*_b(M)$
be the operator $d^*$  on $\Lambda^*U$ lifted to $\Lambda^*_b(M)=\Lambda^*U$,
and let $\Pi_b$ be the orthogonal projection to the space of
forms vanishing on $\Sigma$.
By \eqref{_omega_0_Equation_},
$d(\theta\wedge\theta^c)= - \theta \wedge \omega_0$,
where $\omega_0 = \omega- \theta\wedge\theta^c$ by 
\eqref{_omega_via_theta_Chapter_8_Equation_}.
Then 
\[ 
(-1)^{\tilde \beta} d^* \beta = *d* \beta=  d^*_b(\beta) - *(\theta \wedge \omega_0\wedge
\beta)
\]
Therefore, $\Pi_b(d^* \beta)= d^*_b(\beta)$
for any basic form $\beta$. This is why
$(d^*\beta, d^*\beta)=0$ implies that $(d^*_b\beta, d^*_b\beta)=0$,
and $\ker \Delta \subset \ker d \cap \ker d^* \subset 
\ker d \cap \ker d_b^*\subset \ker \Delta_b$.

\hfill

{\bf Step 5:} We prove that any basic form $\alpha$
that satisfies $\alpha\in \ker \Delta_b$ also satisfies
$d^c\alpha=0$. Since $\Delta_b$ is the pullback of the Laplacian
on a K\"ahler manifold $U$ (the leaf space), 
it is $I$-invariant. Therefore, $I(\alpha)\in \ker \Delta_b$.
It remains to show that any basic form $\beta \in \ker \Delta_b$
is closed. By definition
$\Delta_b= d^*_b d+ d d^*_b$.
From \eqref{_*_0_explicit_Equation_} we obtain
\begin{equation*} 
	\begin{split}
d_b^*(x)&=(-1)^{\tilde x} *_b d *_b (x)= (-1)^{\tilde x} *_b (d * (x \wedge \theta \wedge \theta^c))\\
 &=(-1)^{\tilde x} *\bigg(d( * (x \wedge \theta \wedge \theta^c)) \wedge \theta \wedge \theta^c\bigg)
 	\end{split}
\end{equation*} 
By definition of $*$, we have $(a, b) = \int_M *a\wedge b$
for any differential forms of the same order.
Then for any basic forms $x, y \in \Lambda^*_b(M)$, we have
\begin{equation*} 
	\begin{split}
(d^*_b(x), y) &= 
(-1)^{\tilde x} \int_M d\big( * (x \wedge \theta \wedge \theta^c)\big) \wedge \theta \wedge \theta^c \wedge y\\
&= \int d\big( * (x \wedge \theta \wedge \theta^c)\big) \wedge ( y\wedge \theta \wedge \theta^c)\\
&=(-1)^{\tilde x} \int *(x \wedge \theta \wedge \theta^c) \wedge  d( y\wedge \theta \wedge \theta^c) \\
&=
\big(x \wedge \theta \wedge \theta^c, d(y \wedge \theta \wedge \theta^c)\big).
	\end{split}
\end{equation*}
By \eqref{_omega_0_Equation_}
again, we have $d(\theta \wedge \theta^c) = -\theta \wedge \omega_0$.
Then 
\[
d(y \wedge \theta \wedge \theta^c)=
dy \wedge \theta \wedge \theta^c - (-1)^{\tilde y}y \wedge \theta \wedge \omega_0.
\]
Since $x$ and $y$ are basic, and $\theta$ and $\theta^c$
are orthogonal to each other and to the basic forms,
the last term is orthogonal to $x \wedge \theta \wedge \theta^c$,
which gives 
\[ 
(d^*_b(x), y) = \big(x \wedge \theta \wedge \theta^c, d(y \wedge \theta \wedge \theta^c)\big) = (x \wedge \theta \wedge \theta^c, dy \wedge \theta \wedge \theta^c) = (x, dy).
\]
We proved that $d$ is adjoint to $d_b^*$ on basic forms. Therefore,
\[
(\Delta_b x, x) = (dd^*_b x, x) + (d^*_bd x, x)  =
(dx, dx) + (d^*_b x, d^*_b x),
\]
 and
whenever $\Delta_b x=0$ for a basic form $x$, we also have
$dx = d^*_b x=0$.
We proved \ref{_harmo_deco_1-form_Proposition_}.
\endproof


\section{Rank 1 Vaisman structures}


Recall  that the {\bf LCK rank}\index[terms]{rank!LCK} (\ref{lck_rank}) of an LCK
manifold\index[terms]{manifold!LCK} is the rank of the subgroup $\im(\chi)\subset
\R^{>0}$ or, equivalently, the rank of the monodromy
group\index[terms]{group!monodromy}\index[terms]{rank!LCK}
of the weight bundle, or, equivalently, \index[terms]{bundle!weight}
the rank of the group generated by the periods of
the Lee form $\theta$.\index[terms]{form!Lee} 

We now prove that every compact Vaisman manifold is
biholomorphic to a Vaisman manifold with LCK rank 1.\index[terms]{rank!LCK} The
precise statement is:\index[terms]{manifold!Vaisman}

\hfill

\theorem \label{_rk_1_approximated_Theorem_}
{\bf
  (\cite{ov_jgp_16})}\label{_Vaisman_defo_transve_Proposition_}
Let $(M,\theta,\omega)$ be a compact Vaisman manifold. 
Then $\omega$ can be approximated by a sequence
of Hermitian forms\index[terms]{form!Hermitian} that are  conformally equivalent
to Vaisman metrics of LCK rank 1.\index[terms]{rank!LCK}

\hfill

\pstep
As above, let $G$ be the closure of the group generated by the
Lee and the anti-Lee actions on $(M,\theta,\omega)$.
This group acts on $(M,\theta,\omega)$ by holomorphic
isometries.

Consider a 
cohomology class $[\alpha]\in H^1(M, \R)$ such that the 
class $[\theta]+ [\alpha]$ is rational, and let
$\alpha$ be its harmonic representative. Since
$G$ acts on $M$ by isometries, $\alpha$ is $G$-invariant.
Decomposing $\alpha$ as $a\theta + \beta$,
we may assume that $d^c\beta=0$ (\ref{_harmo_deco_1-form_Proposition_}).

Then the one-form $\theta':= \theta+\alpha$ has rational cohomology class.
Consider the (1,1)-form $\omega':= d^c\theta'+\theta' \wedge I\theta' $ 
obtained as a deformation of 
$\omega= d^c\theta+ \theta \wedge I\theta$
(\ref{_scaling_Vaisman_remark_}).  Assume
that $\alpha$ is chosen sufficiently small in such a way
that $\omega'$ is positive definite. In the next two steps, we prove that
$\omega'$ is conformally equivalent to a Vaisman metric,
using \ref{kami_or}.

\hfill

{\bf Step 2:}
We start by proving that $\omega'$ is an LCK
form,\index[terms]{form!LCK} with the Lee form\index[terms]{form!Lee} $\theta'$. 
Let $d_{\theta'}:\; \Lambda^i (M) \arrow \Lambda^{i+1}(M)$
map $\lambda$ to $d\lambda - \theta' \wedge \lambda$.
Then the LCK condition can be written down as 
$d_{\theta'}(\omega')=0$. Now, the form $d^c\theta'$
is proportional to $d^c \theta$ by \ref{_harmo_deco_1-form_Proposition_}, and hence  it
satisfies $d^c \theta' \in \Lambda^{1,1}(M)$.
This gives $Id^c\theta'= d^c\theta'$, and
\[ d_{\theta'} (I\theta') = - I dI^{-1} \theta'   
 - \theta' \wedge I(\theta')  = - d^c \theta'
 - \theta' \wedge I(\theta')=
-\omega',
\] and hence $\omega'$ is $d_{\theta'}$-closed.

\hfill

{\bf Step 3:}
Consider the holomorphic flow $F$ generated by the Lee field\index[terms]{Lee field}
$\theta^\sharp$. It fixes $\omega$, $\theta$ and $\alpha$,
and its lift $\tilde F$ to the universal cover  $\tilde M$
acts by non-trivial homotheties with respect to the
K\"ahler metric $\tilde \omega'$, because $\tilde \omega'$ 
approximates $\tilde \omega$.  Then \ref{kami_or}
implies that $\omega'$ is conformally equivalent to a Vaisman metric. 
\endproof 

\hfill

The same argument also proves the following corollary.

\hfill

\corollary\label{_Lee_classes_for_Vaisman_open_Corollary_}
Let ${\cal L}_v\subset H^1(M,\R)$ be the set
of all Lee classes\index[terms]{class!Lee} $[\theta]\in H^1(M, \R)$ for
all Vaisman structures on a given compact
complex manifold $M$. Then ${\cal L}_v$ is open in $H^1(M, \R)$.
\endproof 

\hfill

\remark
In fact, the set of Lee classes on $M$ is a half-space
(\ref{_Lee_cone_on_Vaisman_Theorem_}),
that is  not very hard to see from the same argument.

\hfill

\remark\label{_First_appearance_of_LCK_potential_equation_}
The proof of \ref{_Vaisman_defo_transve_Proposition_}
is motivated by the notion of ``LCK manifold with potential''
introduced in Chapter \ref{lckpotchapter}.\index[terms]{manifold!LCK!with potential} Consider 
the operator $d_\theta$ defined as
$d_\theta(\alpha)= d\alpha-\theta\wedge\alpha$, and
let $d_\theta^c:= I d_\theta I^{-1}$. Then
the equation $\omega= d^c\theta+ \theta \wedge I\theta$
(\ref{_scaling_Vaisman_remark_}) becomes
$\omega=d_\theta d_\theta^c(\phi)$, where $\phi=1$.
We express this by saying that $\phi$ {\bf is an 
LCK potential for $\omega$.} A small\index[terms]{potential!LCK}
deformation of $\theta$ and $\phi$ leads
to a new form $\omega'=d_{\theta'} d_{\theta'}^c(\phi')$
that is  also positive, because positivity is
an open condition, and LCK, because
$d_{\theta'}\omega'=d_{\theta'}^2 d_{\theta'}^c(\phi')=0$.

\section{The structure theorem}

\remark
From \ref{_rk_1_approximated_Theorem_} it follows
that any compact Vaisman manifold is biholomorphic to a Vaisman\index[terms]{manifold!Vaisman}
manifold of rank 1. Therefore, the following theorem
implies that any  compact Vaisman manifold
is biholomorphic to a $\Z$-quotient of a Riemannian
cone over a Sasakian manifold.

\hfill

\theorem \label{str_vai} 
({\bf Structure Theorem}, \cite{ov_str}) \\
Every compact Vaisman manifold of LCK rank 1 is\index[terms]{rank!LCK}
biholomorphically isometric to a Vaisman manifold
obtained as $C(S)/\Z$ (\ref{halfstr}), 
where $S$ is Sa\-sa\-kian, \index[terms]{manifold!Sasaki}
$\Z= \bigg\langle (x, t) \mapsto (\phi(x), q t)\bigg\rangle$, $q>1$,
 $\phi$ is a Sasakian automorphism of $S$, and $C(S)$ is the Sasakian cone considered to be  a complex manifold.
Moreover, the triple $(S, \phi, q)$ is unique.

\hfill

 {\bf Proof:} 
As the Lee field\index[terms]{Lee field!parallel} $\theta^\sharp$ is parallel, $M=S\times\R$ locally
(\ref{_parallel_field_de_Rham_Corollary_}). 

The submersion over $S^1$ is obtained  as follows. By \ref{vai_cone},
$\tilde M=C(S)$ is a K\"ahler cone over a Sasakian manifold 
whenever $\theta$ is exact.

Since $M$ has LCK rank 1\index[terms]{rank!LCK}, the deck group $\Gamma$ of its
minimal K\"ahler cover is cyclic, $\Gamma\cong \Z$.
The manifold $M$ is obtained from $\tilde M$
as the quotient $M=\tilde M/\Gamma$, where $\Gamma$
acts on $\tilde M=C(S)= \R^{>0}\times S$ by holomorphic homotheties. Let $\gamma$
be the generator of $\Gamma$.
Since $\gamma$ acts on the first component of $C(S)=
\R^{>0}\times S$ as multiplication by a number $c\neq 1$,
the submersion $M \arrow S^1$ is given by 
$\sigma:\; \tilde M/\Gamma \arrow \R^{>0}/\langle c \rangle=S^1$.

The fibres of $\sigma$ are Sasakian by construction.
The composition of the monodromy map $\gamma:\; C(S) \arrow C(S)$\index[terms]{map!monodromy}
and the homothety $(t, m) \arrow (c^{-1}t, m)$ acts on 
$C(S)=\R^{>0}\times S$ by K\"ahler isometries; 
hence this action descends to a Sasakian isometry on $S$.\index[terms]{isometry!Sasakian}
The corresponding mapping torus \index[terms]{mapping torus}is identified with the
quotient $M=\tilde M /\Gamma$ by construction.
\endproof


\section{Exercises}
\label{_Exercises_for_Vaisman_chapter_Subsection_}

\begin{enumerate}[label=\textbf{\thechapter.\arabic*}.,ref=\thechapter.\arabic{enumi}]

\item Let $\nabla$ be the Levi--Civita connection on a Riemannian
	manifold $(M,g)$
	and  $X$ a vector field satisfying $\nabla(X)=\Id_{TM}$.
	\begin{enumerate}
	\item Prove that $(M,g)$ is locally isometric to a Riemannian cone.\index[terms]{cone!Riemannian}
	\item Let $R\in \Lambda^2 M \otimes \End(TM)$ be the curvature tensor of the Levi--Civita connection.\index[terms]{tensor!curvature}
	Prove that $\Lie_X(R)=0$ and $\nabla_X (R)=-2R$.
	\item Suppose that $(M,g)$ is Einstein: $\Ric(M)=c g$.\index[terms]{manifold!Einstein}
	Prove that $c=0$.
	\end{enumerate}

\item 
Let $X$ be a vector field on a Riemannian manifold $M$ satisfying
$\nabla(X)=\Id_{TM}$, and $V_t$ its diffeomorphism flow. Suppose
that $V_t$ is defined for all $t\in \R$.
Denote by $d$ the Riemannian distance.
\begin{enumerate}
\item
Prove that $d(V_t(a), V_t(b)) = e^{2t} d(V_t(a), V_t(b))$ for
any $a, b \in M$.

\item Let $\bar M$ be the metric completion of $M$.
Prove that there exists a point $p\in \bar M$
such that $\lim\limits_{t \to -\infty} V_t(a)=p$
for all $a\in M$. We call this limit point
{\bf the origin}.

\item Prove that $X$ is the gradient vector
field, $X= \grad(\phi)$, where $\phi$ is the square
of the distance to the origin.
\end{enumerate}

\item
Let $X$ be a compact Riemannian manifold, $C(X)$ its cone,
and $M:=\overline {C(X)}$ its metric completion. 
Prove that $M$ is non-singular (that is, isometric
to a Riemannian manifold) if and only if $C(X)$ is flat 
as a Riemannian manifold.

{\em Hint:} A homothety multiplies the norm of the sectional 
curvature by a constant. Express this constant in terms of
the homothety.

\item
We say that a vector field $v$ is {\bf complete}
if the corresponding diffeomorphism flow $V_t$
is defined for all $t\in \R$.
Let $S$ be a complete Riemannian manifold, 
	$M= C(S)$ its Riemannian cone\index[terms]{cone!Riemannian} and 
$X_1$, $X_2$ complete vector fields satisfying $\nabla(X_i)=\Id_{TM}$.
Prove that $X_1=X_2$.

{\em Hint:} Use the previous exercise and 
the gradient expression for $X_i$. Treat separately the case when $M$ is flat.

\item 
Let $S$ be a complete Riemannian manifold, 
	$M= C(S)$ its Riemannian cone,\index[terms]{cone!Riemannian}
and $\phi:\; M \arrow M$ an isometry. Prove\index[terms]{isometry!of a Riemannian cone}
	that $\phi$ is induced by an isometry of $S$.

{\em Hint:} Use the previous exercise.

\item
Let $M$ be a compact complex manifold,
and $\theta$ a closed 1-form such that
$\omega:=d(I\theta)+ \theta \wedge I\theta$ is an Hermitian form.\index[terms]{form!Hermitian}
Prove that $\omega$ is LCK.\index[terms]{form!LCK}

\item 
Let $M$ be a compact complex manifold,
and $\theta$ a closed 1-form such that
the dual vector field $\theta^\sharp$ is holomorphic and Killing.\index[terms]{vector field!Killing}\index[terms]{vector field!holomorphic}
Assume that $\omega:=d(I\theta)+ \theta \wedge I\theta$ is 
an Hermitian form.
Prove that $M$ admits a Vaisman metric.

\item
Let $(M, \omega, \theta)$ be a compact Vaisman manifold.
Prove that $b_1(M)$ is odd.\index[terms]{manifold!Vaisman}

\item
Let $(M, \omega, \theta)$ be a compact Vaisman manifold,
and $k$ an integer between 1 and $b_1(M)$.
Prove that $\omega$ can be obtained as a limit 
of a sequence of LCK forms\index[terms]{form!LCK} $\omega_i$ of LCK rank $k$,\index[terms]{rank!LCK}
such that each $\omega_i$ is conformal to Vaisman.

{\em Hint:} Use the same argument as in the proof
of \ref{_Vaisman_defo_transve_Proposition_}.

\item\label{_omega_hermitian_on_Tot_circ_L_}
Let $L$ be a positive line bundle
on a projective manifold, and $\Tot^\circ(L)$ the space of 
non-zero vectors in the total space $\Tot(L)$.
Let $\theta = - d\log (l^2)$, where $l$ is the length
function. Prove that 
$\omega:=d(I\theta)+ \theta \wedge I\theta$ is 
an Hermitian form on $\Tot^\circ(L)$.\index[terms]{form!Hermitian}

{\em Hint:} Use \ref{_dd^c_log_l_Corollary_}.

\item
In the assumptions of Exercise \ref{_omega_hermitian_on_Tot_circ_L_},
let $q:\; \Tot^\circ(L)\arrow \Tot^\circ(L)$
be a holomorphic homothety, generating a 
properly discontinuous $\Z$-action.\index[terms]{action!$\Z$-}
Prove that $\Tot^\circ(L)/\Z$ is Vaisman.

\item Construct a compact Vaisman manifold
with $b_1 >1$.\index[terms]{manifold!Vaisman}

{\em Hint:} Use the previous exercise.

\item 
Let $(M, \omega)$ be a complex Hermitian manifold,
and $\theta$ a closed 1-form such that
$\omega=d(I\theta)+ \theta \wedge I\theta$.
\begin{enumerate}
\item Let $\omega_0 =d(I\theta)$.
Prove that 
$*(\omega_0) = \frac1 {(n-2)!} \omega_0^{n-2} \wedge\theta \wedge I\theta$.
\item Prove that $*(I\theta)= \frac1 {(n-2)!} \omega_0^{n-2}\wedge\theta$.
\item Prove that $d * d(I\theta)= -\frac 1 {(n-2)!} \omega_0^{n-1} \wedge\theta$.
\item Prove that 
$d^*d(I\theta)= -(n-1) \theta$.
\end{enumerate}

\item\label{_holo_curva_Heisenberg_Exercise_}
Let $(T,I, \omega)$ be a compact flat K\"ahler torus, 
and $L$ a holomorphic Hermitian line bundle 
such that the curvature of its  Chern connection\index[terms]{connection!Chern}
is equal to $\omega$. Denote by $S$
the total space of the corresponding $S^1$-bundle
of unit vectors, fibred over $T$ with the fibre $S^1$.
\begin{enumerate}
\item Prove that $S$ is Sasakian. Prove that the group of 
Sasakian automorphisms of $S$ is 1-dimensional.
\item 
The manifold $S$ is the total space of a principal $S^1$ bundle
$L$, equipped with the Chern connection that induces an
Ehresmann connection\index[terms]{connection!Ehresmann}
$TS = T_\hor S \oplus T_\ver S$, where
$T_\ver S$ is the bundle of vectors tangent to the 
fibres of $\pi: S \arrow T$.
Consider the metric $g$ on $S$ 
compatible with the decomposition $TS = T_\hor S \oplus T_\ver S$,
and such that $g$ on $T_\hor S =\pi^* TT$ is
induced from the metric on $T$ and $g$ on the trivial
bundle $T_\ver S$ is constant.
Prove that $g$ is the Sasakian metric on $S$.
\item
Let $R\in \Sym^2(\Lambda^2 TS)$ be the Riemannian curvature tensor
on $S$. Prove that
$R(X,Y,Z,T)= - 3 \omega_0(X,Y)\omega_0(Z,T)$, where
$\omega_0\in \Lambda^2 S$ is the pullback of the symplectic form
on the torus.
\item
Let $M$ be a Vaisman manifold associated with $S$,
locally isometric to $S\times S^1$. Such a Vaisman manifold
is called {\bf a Vaisman manifold of Heisenberg type}.\index[terms]{manifold!Vaisman}
Define {\bf the holomorphic sectional curvature}
of $M$ as a map $X\arrow R_M(X,IX,X,IX)$,
where $R_M$ is the Riemannian curvature of $M$, and $X\in TM$.
Prove that the holomorphic sectional curvature
of $M$ is negative when $X$ is orthogonal to the Lee
field, and vanishes on the Lee field\index[terms]{Lee field} direction.
\end{enumerate}

\item
Recall that {\bf a nilmanifold}\index[terms]{nilmanifold}
is a quotient $G/\Gamma$, where $G$
is a simply connected nilpotent Lie group,
and $\Gamma$ a discrete, cocompact lattice.\index[terms]{lattice!cocompact}
We usually consider the left quotients.
Note that a nilmanifold is not left homogeneous,
but it is right homogeneous.
For more details on nilmanifolds,\index[terms]{nilmanifold} see
Chapter \ref{_nil_and_solv_} in the third part of this book.
Let $\goth g$ be a nilpotent Lie algebra
generated by $x_i, y_i, z, t$, with 
the only non-vanishing commutators 
$[x_i, y_i]=z$, $G$ the corresponding
Lie group, $\Gamma$ the lattice \index[terms]{lattice}
generated by the exponents of the
generators. Define a left-invariant almost
complex structure on $G$ by taking the left translations
of the complex structure operator $I\in \End(\g)$
with $I(x_i) = y_i, I(y_i)=-x_i, I(z)=t, I(t)=-z$.
Prove that this almost complex structure is 
integrable, and that the nilmanifold\index[terms]{nilmanifold} $(G/\Gamma, I)$
is biholomorphic to the Vaisman manifold of the Heisenberg type
defined in Exercise \ref{_holo_curva_Heisenberg_Exercise_}.\index[terms]{manifold!Vaisman}

\item\label{_nilma_Heisenberg_forms_Exercise_}
Let $M:=G/\Gamma$ be a nilmanifold\index[terms]{nilmanifold} obtained as a left quotient
of a simply connected nilpotent Lie group.
A tensor field $\Phi$ on $G/\Gamma$ is called {\bf locally
invariant} if its pullback to $G$ is left invariant.
\begin{enumerate}
\item Let $\g$ be a Lie algebra. Define the differential
on $\Lambda^*(\g^*)$ by taking a 1-form $\eta\in \Lambda^1(\g^*)$
to $d\eta(x,y):= -\eta([x, y])$ and extending it to 
$\Lambda^*(\g^*)$ by the Leibniz rule. 
Prove that $d^2=0$, and show that, conversely,
$d^2=0$ implies the Jacobi identity in $(\g, [\cdot, \cdot])$.
The complex  $(\Lambda^*(\g^*), d)$ is called {\bf the \index[persons]{Chevalley, C.} Chevalley-Eilenberg
complex}. We consider it as a differential graded algebra 
(\ref{_DGA_Definition_}).
\item
Prove that the differential graded algebra of locally invariant differential
forms on $M=G/\Gamma$ is isomorphic to the \index[persons]{Chevalley, C.} Chevalley-Eilenberg
complex of $\g$. 
\item
Let $M:=G/\Gamma$ be the  Vaisman nilmanifold\index[terms]{nilmanifold} of Heisenberg type
obtained in the previous exercise, $x_i, y_i, z, t$ the basis
in $\g$, and $\xi_i, \upsilon_i, \zeta, \tau$ the dual basis
in $\g^*$, extended to a basis in the space of 
locally invariant forms on $M$. Prove that
$d(\xi_i)=d(\upsilon_i)=d(\tau)=0$, and
$d(\zeta) =  \sum_i \xi_i\wedge\upsilon_i$.
\item
Suppose that $\dim_\C M=2$; in this case, 
the basis in $\g^*$ is denoted 
$\xi, \upsilon, \zeta,\tau$, and the only non-trivial differential is
$d(\tau)= \xi \wedge\upsilon$. Prove that
$\tau^{1,0}\wedge \xi^{0,1} + \xi^{1,0}\wedge \tau^{0,1}$
is a non-degenerate, closed, pseudo-Hermitian form of 
signature (2,2). \index[terms]{form!pseudo-Hermitian}
\item Prove that any closed
locally invariant 2-form on $M$ is degenerate, if 
$\dim_\C M >2$.
\item
 K. \index[persons]{Nomizu, K.} Nomizu (\cite{nomizu}) has shown that the
natural embedding from the \index[persons]{Chevalley, C.} Chevalley-Eilenberg differential
graded algebra $(\Lambda^*(\g^*), d)$ to the de Rham
algebra of $M$ induces an isomorphism in cohomology
for any nilmanifold \index[terms]{nilmanifold}$M$. Now let  $M$ be
a Vaisman nilmanifold of Heisenberg type. Using
Nomizu theorem,\index[terms]{theorem!Nomizu} prove that any closed 2-form $\eta$ on $M$
satisfies $\int_M\eta^{\dim_\C M}=0$, if $\dim_\C M >2$.
\end{enumerate}

\end{enumerate}


\chapter{Orbifolds}\label{orbif}

\epigraph{\it Consider the following subtraction problem: 342 minus 173. Now, remember how we used to do that: Three from two is nine, carry the one, and if you're under 35 or went to a private school, you say seven from three is six, but if you're over 35 and went to a public school, you say eight from four is six ...and carry the one, so we have 169. But in the new approach, as you know, the important thing is to understand what you're doing, rather than to get the right answer.}{\sc\scriptsize Tom Lehrer,  New Math}

\section{Introduction}

Orbifolds are one of the most straightforward
generalizations of a manifold. In fact, most statements
about de Rham algebra, harmonic forms, vector bundles
etc. generalize to the orbifold category immediately.
Roughly speaking, an orbifold is a variety with 
quotient singularities. However, for some purposes
this ``rough'' notion is not sufficient, and we 
need to explain the notion in more depth.

The difference between orbifold and the ``naive'' notion
of ``variety with quotient singularities'' is that
the orbifold remembers which group action was used.
One can define an ``orbispace'' as a space obtained
locally as a quotient $U_i/G_i$ of charts $U_i$ by finite group action,
with each point $x\in U_i$ equipped with an ``automorphism group''
$\Mor(x,x):= \St_{G_i}(x)$. The (complex) orbifold is
an orbispace with all $U_i$ open balls in $\R^n$ ($\C^n$
for the complex case) and all gluing maps (holomorphic)
smooth.

In algebraic geometry\index[terms]{geometry!algebraic} the orbifolds are known as
``Deligne--Mumford  stacks'';\index[terms]{Deligne--Mumford,  stack} indeed, if we drop the
requirement that all $G_i$ are finite,  the notion of an orbifold
becomes the notion of a ``smooth stack''.

Almost every result and construction that is  valid
for complex manifolds can be trivially extended to orbifolds.
This includes the notion of complex and Riemannian structure,
that we shall use further on for K\"ahler orbifolds.
It takes a little more care to define a holomorphic
vector bundle and a principal $G$-bundle over an
orbifold; we give the relevant definitions for
completeness.\index[terms]{bundle!principal}

The aim of this chapter is simple, we define an
orbifold, and prove that any quotient of a manifold 
by a compact torus action is an orbifold, as long
as this torus has orbits of the same dimension.  
In the smooth category, the same is actually true for any 
compact Lie group, but we work in the holomorphic category, and the only
compact connected complex Lie group is the torus.
We give a proof that is  adapted for a torus.
However, if one generalizes \ref{_limits_torus_Proposition_} 
to arbitrary groups,  the general proof will be just as easy.

It is interesting that these torus quotients
are locally trivial bundles 
in the orbifold category (\ref{torus_quotient} (ii)).

A reader who is acquainted with the stacks%
\footnote{\url{https://stacks.math.columbia.edu/}} won't find
anything new in this chapter, but readers who
are trained in differential geometry might find 
the exposition too compressed. We recommend
\cite{bog}, where the same basic results and definitions
are given in more detail. 

The reason why we need the orbifold is basically the same
as why \cite{bog} need them. Any Sasakian manifold is equipped
with an action of a one-parametric group generated by the
Reeb field; this action preserves the Sasakian structure
(both the CR-structure \index[terms]{structure!CR}and the metric). The Sasakian
manifold is {\bf quasi-regular}\index[terms]{manifold!Sasaki!quasi-regular} if all orbits of
the Reeb action are compact. By the result proven
in this chapter, the space of Reeb orbits is
actually an orbifold. As shown in \ref{regsas},
this orbifold is actually K\"ahler (the proof
in \ref{regsas} does not invoke the orbifolds,
but the orbifold proof is just the same).
Moreover, this orbifold is also projective.
This allows one to reduce most questions
about quasi-regular Sasakian manifolds to
projective algebraic geometry.\index[terms]{geometry!algebraic!projective}

This observation is especially useful, because
any Sasakian manifold can be deformed to a
quasi-regular one, and, moreover, the quasi-regular
Sasakian manifolds are dense in the appropriate
moduli space (\ref{_defo_qr_sas_}).

For the Vaisman manifolds, the situation is even
better. Recall that the Vaisman manifold $M$ is\index[terms]{manifold!Vaisman}
called {\bf quasi-regular} \index[terms]{manifold!Vaisman!quasi-regular}when its canonical
foliation has compact leaves. 
The leaves of the canonical foliation \index[terms]{foliation!canonical}can be identified 
with the orbits of the action of the Lie group
$G$ generated by the Lee and anti-Lee fields.\index[terms]{Lee field}\index[terms]{Lee field!anti-}
This group acts on $M$ by holomorphic isometries. 
A compact Vaisman manifold is quasi-regular
if and only if $G$ is a compact torus.
In this case, the quotient $M/G$
is a complex orbifold. In 
\ref{_Structure_of_quasi_regular_Vasman:Theorem_} we prove that $M/G$ is actually
projective, and, moreover, $M$ can be recovered
from the complex projective data.
Again, the compact quasi-regular Vaisman manifolds\index[terms]{manifold!Vaisman!quasi-regular}
are dense in the set of all compact Vaisman manifolds,
hence we can reduce most topological and complex-analytic
questions about Vaisman manifolds to the projective
algebraic geometry.\index[terms]{geometry!algebraic!projective}

Orbifolds were introduced in \cite{satake} as
``V-manifolds''\index[terms]{manifold!Satake (``V-'')}
and have appeared with this name in many papers since then.
In 1970-ies, W. \index[persons]{Thurston, W.} Thurston used the V-manifolds in his
geometrization program\index[terms]{geometrization program} for
3-manifolds. It was Thurston
who renamed the V-manifolds to orbifolds.
Orbifolds are generalizations of manifolds and naturally appear as leaf
spaces of quasi-regular foliations, such as the
canonical foliation\index[terms]{foliation!canonical} of a quasi-regular Vaisman manifold. We 
provide here a brief introduction to the subject, sending
the reader to \cite{bog} and the references therein for a
more detailed one. \index[terms]{manifold!Vaisman}\index[terms]{manifold!Vaisman!quasi-regular}

\section{Groupoids and orbispaces}\index[terms]{groupoid}\index[terms]{orbispace}

\definition 
A {\bf groupoid} is a category with all morphisms invertible.

\hfill

\definition 
An action of a group on a manifold is called  {\bf rigid} if
the set of points with trivial stabilizer is dense.\index[terms]{action!rigid}

\hfill

\definition 
An {\bf orbispace} is a topological space $M$, 
equipped with:
\begin{enumerate}
	\item  A structure of a groupoid\index[terms]{groupoid} (the points of $M$
are objects of the groupoid category);
\item  An open cover $\{U_i\}$;
 \item Continuous maps $\phi_i:\; V_i\arrow U_i$, where each
$V_i$ is equipped with a rigid action of a finite group $G_i$,
satisfying the following properties: \index[terms]{group!finite}
\begin{enumerate}
\item $\phi_i:\; V_i \arrow V_i/G_i=U_i$ is the quotient map.
\item For each $x\in M$ and $U_i\ni x$, 
the group $\Mor(x,x)$ is equal to the stabilizer $\St_{G_i}(x)$ of $x$ in $G_i$.
\end{enumerate}
\end{enumerate}
The atlas $\{(V_i, G_i)\}$  defined by these data is called
{\bf an orbifold cover.}

\hfill

All notions that are  known from the geometry of
manifolds can be defined in the orbifold setting as well.

To define 
the {\bf gluing maps}, associated with an orbifold,
we take an intersection $U_{ij}:=U_i \cap U_j$ realized as
a subset of $U_i$ and  $U_{ji}=U_j \cap U_i$ realized as as a
subset in $U_j$. Let $V_{ij}\subset V_i$ and
$V_{ji} \subset V_j$ be their preimages.
By construction, the stabilizer $\St_{G_i}(V_{ij})$ of $V_{ij}$ in $G_i$ is
identified with the stabilizer $\St_{G_j}(V_{ji})$ of
$V_{ji}$ in a way that is  compatible with the chain of
identifications
\[
U_{ij}= \frac{V_{ij}}{\St_{G_i}(V_{ij})} =
\frac{V_{ji}}{\St_{G_j}(V_{ji})} = U_{ji}.
\]
The {\bf gluing map} is a homeomorphism
$\psi_{ij}:\; V_{ij}\arrow V_{ji}$ constructed above. By
construction, it satisfies {\bf the
cocycle condition} $\psi_{ij} \psi_{jk} =\psi_{ik}$,
similar to that for the manifolds.

We identify the orbispaces defined by an orbifold cover and
its $G_i$-invariant refinement, similarly to how it is
done for ``charts and atlases'' definition of the manifolds.

\hfill

\definition 
{\bf A morphism} $\Phi:\; X \arrow Y$ of orbispaces is 
given by the following data:
\begin{enumerate}
\item A continuous map $\Phi$ of the underlying topological spaces,
together with the functors $\Phi_x:\; \Mor(x,x)\arrow \Mor(\Phi(x),\Phi(x))$
of the corresponding grou\-po\-ids.

\item For any orbifold cover $\{(W_i, H_i)\}$ of $Y$, there
  is an orbifold cover $\{(V_j, G_j)\}$
of $X$, such that for any chart $(W_i, H_i)$ on $Y$
there exists a chart $(V_i, G_i)$ on $X$ such that
$\Phi(V_i/G_i) \subset W_i/H_i$.

\item For any charts $(V_i, G_i)$, $(W_i, H_i)$ such that
$\Phi_\top(V_i/G_i) \subset W_i/H_i$
there is a group homomorphism $\Phi_{(V_i,W_i)}:\; G_i \arrow H_i$
and a continuous map $\tilde \Phi_{(V_i,W_i)}:\; V_i
\arrow W_i$ making the following diagram commutative
for any $g\in G_i$
\[\begin{CD}
V_i @>{\tilde \Phi_{(V_i,W_i)}}>> W_i \\
@V{g}VV @VV{\Phi_{(V_i,W_i)}(g)}V\\
V_i @>>{\tilde \Phi_{(V_i,W_i)}}> W_i
\end{CD}
\]
Moreover, the liftings $\tilde \Phi_{(V_i,W_i)}:\; V_i
\arrow W_i$ are compatible with the relevant transition
maps between the charts in $X$ and $Y$.
\end{enumerate}
The atlases $\{(W_i, H_i)\}$, $\{(V_j, G_j)\}$
satisfying these conditions are called {\bf compatible}.

\hfill

\remark 
An orbispace is a topological space, locally obtained as
a quotient, {with the quotient structure remembered
	via the groupoid structure.}\index[terms]{groupoid}

\hfill

\definition 
A point $x\in M$ is called  {\bf orbipoint} if \index[terms]{orbipoint}
 $\Mor(x,x)$ is non-trivial. The {\bf order} of an orbipoint
is $|\Mor(x,x)|$.  The group $\Mor(x,x)$
is called {\bf the monodromy group}\index[terms]{group!monodromy} of the orbipoint.

\section{Real orbifolds}

\definition 
 An {\bf orbifold} is an orbispace\index[terms]{orbifold}
$(M, \{\phi_i:\; V_i \arrow V_i/G_i=U_i\})$,
where all $V_i$ are diffeomorphic to open balls in $\R^n$.
The {\bf local coordinates} on an orbifold are
coordinates on $V_i$ defined by the embeddings
to $\R^n$. A {\bf morphism of orbifolds} is a morphism of the
corresponding orbispaces defined by smooth maps
in local coordinates.

\hfill

\example 
Let $M=\C P^1/((x,y)\sim(x, -y))$. {This quotient
	is homeomorphic to $\C P^1$.} {However, it is a different orbifold}
if we consider the cover induced from $\C P^1/G$, with $G=\{\pm 1\}$
and the groupoid\index[terms]{groupoid} structure where $\Mor(x,x)=\St_G(x)$.

\hfill


Orbifolds where introduced to treat ``manifolds with
singularities''. However, on a smooth orbifold one can
introduce all the machinery of differential geometry:

\hfill

\definition 
 A {\bf smooth orbifold} is an orbifold $M$ equipped with
a sheaf of functions $ C^\infty (M)$ in such a way that for each 
$U_i=V_i/G_i$, the corresponding ring of sections
$ C^\infty (U_i)$ is identified with the  
ring of $G_i$-invariant smooth functions on $V_i$.

\hfill

\definition 
 A {\bf Riemannian metric} on a smooth orbifold is a
$G_i$-invariant metric on each $V_i$, compatible with the
gluing maps.

\section{Complex orbifolds}

All notions in complex geometry can be generalized to orbifolds.

\hfill

\definition 
 A {\bf complex orbifold} is an orbifold $M$ equipped with
a sheaf of functions $\calo_M$ in such a way that each $V_i$
is an open ball in $\C^n$, and for each 
$U_i=V_i/G_i$, the corresponding ring of sections
$\calo_{U_i}$ is identified with the
ring of $G_i$-invariant holomorphic functions on $V_i$.

\hfill

\definition 
The {\bf underlying complex variety} of a 
complex orbifold is a complex variety with the topological
space $M$ and the structure sheaf $\calo_M$.

\hfill

\example \label{weight_CP}
Let $\C^*$ act on $\C^n$ as \[ h_t(x_1, ..., x_n)=
(t^{a_1}x_1, t^{a_2}x_2, ..., t^{a_n}x_n).\]
The quotient $\C P^{n-1}(a_1, ..., a_n):=(\C^n\backslash 0)/\C^*$ is called a
{\bf weighted projective space}. It is an orbifold, as
follows from \ref{torus_quotient}.\index[terms]{weighted projective space}

\hfill

 \definition 
 A {\bf  projective orbifold} is a complex
 orbifold with the underlying complex
 variety projective.\index[terms]{orbifold!projective}

\hfill

 \definition 
 A {\bf holomorphic vector bundle}, also known as
{\bf V-bundle} on a complex 
orbifold\index[terms]{bundle!vector bundle!holomorphic}
 is a $G_i$-equivariant vector bundle on each $V_i$, equipped
 with the $G_i$-invariant gluing maps satisfying the cocycle 
condition.\footnote{We give a definition of {\bf principal
orbifold bundle} in Section \ref{_orbibundles_Section_}.}
 
A holomorphic Hermitian line bundle $L$ 
on an orbifold is called {\bf positive} if the curvature
of $L$ is positive definite on each  $V_i$.

\hfill

 We shall need the following 
extension of  Kodaira embedding theorem to orbifolds:
 
\hfill

 \theorem   {\bf (W. L. \index[persons]{Baily, W. L.} Baily)} \label{baily}
 Let $M$ be a compact complex orbifold equipped with 
 a positive holomorphic Hermitian vector bundle $L$. Then
 $M$, considered to be  a complex variety, admits an
embedding to $\C P^n$.\index[terms]{theorem!Kodaira's embedding!for
 orbifolds (Baily)}\index[terms]{orbifold!projective}\index[terms]{theorem!Baily}

\proof
It is proven in the same way as the 
usual  Kodaira embedding theorem; see
\cite{baily}.
\endproof


\section{Quotients by tori}

 
In the following, we shall need the notion of
stratification associated with a group action. A good
reference is \cite{dui}. 

\hfill

\definition\label{_Hausdorff_distance_Definition_}
Let $X, Y$ be subsets of a metric space $M$.\index[terms]{space!metric} Recall
        that the
        {\bf Hausdorff distance} is defined as
\begin{equation}\label{_Hausdorff_distance_Equation_}
\delta(X, Y):= \max\left(\sup_{x\in X}\:d(x, Y),\; \sup_{y\in
   Y}\:d(y, X)\right).
\end{equation}
The Hausdorff distance is clearly metric-dependent.
However, it is not hard to see that the
topology defined by this distance
on the set of compact subsets of $M$
is independent of the choice of the metric.
Further on, we shall speak of {\bf Hausdorff
convergence} of a sequence of compact
subsets in $M$ referring to this topology.

\hfill

The next result essentially says that if a subgroup $G_1$
is close enough to a closed subgroup $G$ of a torus, in the
\index[persons]{Gromov, M.} Gromov--Hausdorff distance, then $G_1$ is a subgroup of
$G$. The precise statement is as follows.\index[terms]{Gromov--Hausdorff distance}

\hfill

 \begin{lemma}\label{_GH_limit_subgroups_Lemma_}
Fix a flat Riemannian
 	metric on a compact torus $T^n$. Then for any compact
 		subgroup $G\subset T^n$ there exists a positive number $\epsilon(G)$
 		such that $\delta (G_1, G)> \epsilon(G)$
                unless $G_1\subset G$.
 \end{lemma}
 
\hfill

 \proof 
Let $\epsilon(G): = \frac 2 3 R$, where
 $R$ is the metric diameter of the smallest circle in the decomposition
 $T^n/G= (S^1)^k$, where  $T^n/G$ is considered with the  
 flat metric induced from $T^n$. This number is the
 smallest possible diameter of a non-trivial subgroup
 $F\subset T^n/G$. Since the quotient map $T^n \arrow T^n/G$ is
1-Lipschitz, this gives  $\epsilon(G) \leq
 \delta(0, G'/G'\cap G)$ for any  subgroup $G'\subset T^n$
 with 
 $G'\not \subset  G$. Therefore, 
 $\delta(G_1, G\cap G_1)\geq \delta(0, G_1/(G_1\cap G))\geq
 \epsilon(G)$ unless $G_1 \subset G$.
 \endproof

\hfill

The Hausdorff distance defines a topology
on the set of closed subsets of a metric space.\index[terms]{space!metric} 
It is well known that 
the set of compact subsets of a compact metric space
with Hausdorff topology is compact (\cite{_Burago_Burago_Ivanov_}).
Clearly, the limit of closed subgroups of a torus
is again a closed subgroup, that is, a subtorus.
\ref{_GH_limit_subgroups_Lemma_} implies the following.

\hfill

 \proposition \label{_limits_torus_Proposition_}
 For any converging 
sequence of compact subgroups of a torus, $T_i \subset T^n$, the
 limit $T_\infty:=\lim_i T_i$ contains all $T_i$, except a finite number.	
 
\hfill

 \proof
Clearly, $T_\infty$ is a closed subgroup in $T$, that is, a subtorus.
Consider the number $\epsilon(T_\infty)$ defined in 
\ref{_GH_limit_subgroups_Lemma_}.
Then $\delta(T_{m}, T_\infty) < \epsilon(T_\infty)$,
 for all $m$, except finitely many. 
From $\delta(T_m, T_\infty)< \epsilon(T_\infty)$ and
\ref{_GH_limit_subgroups_Lemma_} we obtain that
$T_m \subset T_\infty$ for all
$T_m$ except finitely many (in \cite{_Mozes_Shah_}
a similar result is proven in the setting of measure theory). \endproof
 
\hfill

\corollary 
	Let $T^n$ be a compact torus acting on a connected
        topological space.
	Consider the function $x \stackrel \Psi \mapsto \St(x)$.
	Then there exists a stratification of $M$ by closed
	strata $M_i$ such that {$\Psi$ is constant on the complement
	$M_i \backslash \bigcup_{M_j \subsetneq M_i} M_j$, to the
smaller strata
and $\Psi(M_i)\supset \Psi(M_j)$
	whenever $M_j \supset M_i$.}\index[terms]{stratification}

\hfill

{\proof} Consider the set ${\goth A}$ of all compact
	subgroups of $T^n$, and let 
	$M_\alpha:= \{x\in M \ \ |\ \ \Psi(x)\supset\alpha\}$,
	where $\alpha \in {\goth A}$. The map $\Psi$ is
        semi-continuous\index[terms]{map!semi-continuous} in the following sense: given
a sequence $\{x_i\}$ converging to $x$, 
the \index[persons]{Gromov, M.} Gromov--Hausdorff limit of $\St(x_i)$ contains
$\St(x)$. This follows from the closedness of the
set 
\[ Z=\big\{(g, x)\in T^n \times M\ \ |\ \ g(x)=x\big\}.
\]
By semi-continuity of
        $\Psi$, the set $M_\alpha$ is closed
	for each $\alpha$. The relation $\Psi(M_i)\supset \Psi(M_j)$
	for smaller strata follows from the semi-continuity
        of $\Psi$  and \ref{_limits_torus_Proposition_}. 
	\endproof

\hfill

\definition
Let $\Phi:\; X \arrow Y$ be an orbifold morphism,
$\{(V_i, G_i)\}$ and $\{(W_i, H_i)\}$ the compatible atlases
on $X$ and $Y$ and $\tilde\Phi_{(V_i,W_i)}:\; V_i \arrow W_i$ the
lifts of $\Phi$ to the orbifold covers. We say that
$\Phi$ is a {\bf locally trivial orbifold fibration}
if all maps $\tilde\Phi_{(V_i,W_i)}:\; V_i \arrow W_i$
are locally trivial fibrations.

\hfill

We now state the main result of this section (we shall use
it further, in order to understand the structure of
compact quasi-regular Vaisman\index[terms]{manifold!Vaisman!quasi-regular} manifolds):\index[terms]{manifold!Vaisman}

\hfill

\theorem \label{torus_quotient}
Let  $G= T^n$ be a torus acting on a manifold $M$
by biholomorphic maps, and let $M_i$ be the corresponding
stratification. Assume that all orbits of $G$ 
have the same dimension and the action of $G$ is 
effective.
Then:
\begin{description}
\item[(i)] The quotient $M/G$ is equipped with an orbifold structure,
and for each $x\in M/G$, the corresponding group of morphisms $\Mor(x,x)$
is equal to $\St(x)$.
\item[(ii)] There is an orbifold structure
on $M$, with $\Mor(x,x)=H_x$, such that
the quotient map $\pi:\; M\arrow M/G$ is an orbifold
morphism. Moreover, $\pi$ is a locally trivial
orbifold fibration.
\end{description}

{\bf Proof of (i):}
We first show that all orbits of $G$ are smooth. Indeed,
each orbit $O$ is diffeomorphic to $T^n/\St(x)$ for any $x\in O$.

{\bf A local section} (also known as {\bf a slice}) of the action 
of $G$ at $x\in M$ is a smooth submanifold $S\ni x$ defined locally
in some neighbourhood of $x$, transversal to the orbit
$G\cdot x$ and having complementary dimension to this orbit.
 Clearly, a local section exists at each $x\in M$.
 
Observe that a local section at $x$ can always be chosen
$\St(x)$-invariant.  To see this, choose a $G$-invariant
Riemannian metric, let $W\subset T_x M$ be the orthogonal complement
of the tangent space to $G\cdot x$, and let $S$ be the union of 
sufficiently short segments of all geodesics passing
through $x$ and tangent to $W$.

Finally, let $H_x:= \St(x)$. Let $S$ be an $\St(x)$-invariant section
(slice) of the $G$-action in $x$.
Take a tubular neighbourhood $U$ of an orbit $G\cdot x$ given by
\[ U:=\bigcup\limits_{g\in G/H_x} gS.\] 
For $S$ sufficiently small, 
the quotient map $U/G\arrow S/H_x$ is well-defined. 
Indeed, the stabilizer is semi-continuous, and hence 
for an appropriate choice of $S$, the stabilizer of any point 
$y\in S$ is contained in $H_x$. 

\hfill

{\bf Proof of (ii):} 
Let $S$ be the local slice defined above, and
$U\subset M$ the corresponding tubular neighbourhood. 
Denote by $\tilde U=S\times G$ the trivial $G$-bundle
over $S$. The natural evaluation map $\tilde U \arrow U$
mapping $(s, g)$ to $g(s)\in U$ is finite;
indeed, the preimage of a point $u\in U$ is identified
with its stabilizer in $G$. We are going to define
the orbifold structure on $M$ using $\tilde U$,
for various $U$, as a chart. The corresponding
orbipoints are points $m$ with non-trivial stabilizer
$\St_G(m)$. By construction, the  map
$\tilde U \arrow U$ is a quotient by the $\St_G(m)$-action.
This takes care of the definition of an orbifold;
the map $\pi:\; M \arrow M/G$ is compatible with this
orbifold structure, because it maps a point $m$
with $\Mor(m,m)= \St_G(m)$ to a point
$\pi(m)$ with $\Mor(\pi(m),\pi(m))= \St_G(m)$,
consistent with taking the local quotient
on the appropriate charts.

In orbifold charts, the map $\pi:\; M \arrow M/G$
is written as the quotient $S\times G \arrow S$,
where $S$ is an appropriate slice, and hence  it is
locally trivial.
\endproof

\section{Principal orbifold bundles}
\label{_orbibundles_Section_}

\definition
Let $X$ be an orbifold, and
$\pi:\; E \arrow X$ a locally trivial orbifold fibration.
Assume that $E$ is equipped
with a fibrewise transitive action of a group $G$.
Suppose that for any local
model $U= \tilde U/\Gamma$, there exists a 
$\Gamma$-equivariant 
principal bundle $E_{\tilde U}\arrow \tilde U$  with fibre $G$
and an isomorphism $E\cong E_{\tilde U}/\Gamma$.\index[terms]{bundle!principal}
Suppose, moreover, that for any two local models
$\tilde U_1 \arrow U_1$ and $\tilde U_2 \arrow U_2$
there exist transition maps of the corresponding principal 
fibrations
\[ 
\tilde \phi_{12}:\;  E_{\tilde U_1}\restrict {U_1 \times U_2}
\tilde \arrow E_{\tilde U_2}\restrict {U_1 \times U_2} 
\]
that satisfy the cocycle condition.
Then $E$ is called {\bf an orbifold
principal bundle} over $X$.

\hfill

\remark
Sometimes an orbifold principal bundle
is called {\bf a V-bundle}.

\hfill

\example
Consider an action of a finite group $\Gamma$
on a complex manifold $X$, and let $L$ be a
$\Gamma$-equivariant holomorphic line bundle on $X$.
Suppose that the stabilizer $\St_\Gamma(x)$ of any point
$x\in X$ in $\Gamma$ acts freely on
the fibre of $L$ in $x$. This defines
a structure of principal orbifold $\C^*$-bundle 
over $X/\Gamma$. Moreover, any
principal orbifold $\C^*$-bundle on 
a complex orbifold is locally
modeled on this example.

\hfill

\claim\label{_quotient_principal_Claim_}
Let $S$ be a manifold equipped with an action of a torus $T$
with finite stabilizers, and $X:=S/T$ the
corresponding orbifold quotient. Then
$\pi:\; S \arrow X$ is a principal orbifold
$T$-bundle.

\hfill

{\bf Proof:}  By \ref{torus_quotient},
the manifold $S$ is equipped with an orbifold
structure in such a way that $\pi$ is a
locally trivial orbifold fibration. 
Moreover, there is a free and transitive
$T$-action on fibres of $\pi$ in orbifold charts.
This puts a structure of a $G$-principal
orbifold bundle on $\pi:\; S \arrow X$.
\endproof

\section{Exercises}
 
\begin{enumerate}[label=\textbf{\thechapter.\arabic*}.,ref=\thechapter.\arabic{enumi}]

\item Prove that a 1-dimensional
complex orbifold is uniquely defined by the \index[terms]{orbifold!complex}
following data: a smooth complex curve $M$, some orbipoints
$x_i$ marked on $M$, and order $p_i\in \Z^{>1}$ of the
monodromy group\index[terms]{group!monodromy} $\Z/p_i \Z$ at each $x_i$.

\item Let $S$ be a quasi-regular Sasakian\index[terms]{manifold!Sasaki!quasi-regular} manifold,\index[terms]{manifold!Sasaki} and $X:= S/\xi$
the corresponding orbifold. Prove that monodromy group of each
orbipoint is a cyclic group. 

\item Prove that the weighted projective space defined in \ref{weight_CP} is a complex orbifold.\index[terms]{weighted projective space}

\item Let $\tilde M:=\ \C^2 \backslash 0$, and $\C^*$ act on
$\tilde M$ as $h_t(x,y)=(tx, t^2 y)$. 
Prove that $\tilde M/\C^*$ is $\C P^1$ with one orbipoint of order 2.

 \item Define  coverings in the orbifold category, and prove the 
 existence of a universal covering.
 
 \item Let $M$ be a complex curve of genus 1 with at least
 one orbipoint. Prove that the universal cover of $M$
 (in the orbifold category) is the Poincar\'e disc $\Delta$.

\item 	\label{_orbifold_funda_Exercise_}
Define the {\bf orbifold fundamental group}\index[terms]{fundamental group!orbifold} $\pi_1^\orb(M)$
	as the group of au\-to\-mor\-phisms of the universal covering compatible
	with the projection to $M$. 
	In this context, we denote by 
$\pi_1^{top} (M)$ the {\em topological fundamental group}.\index[terms]{fundamental group!topological}
	Construct a monomorphism $\pi_1^{top}(M)\arrow \pi_1^\orb(M)$.
Find an orbifold $M$ such that $\pi_1^{top}(M)$ is trivial, and
$\pi_1^\orb(M)$ is non-trivial.

\item Let $M$ be a 1-dimensional complex orbifold with
a complete Hermitian metric of constant negative 
curvature, and $\pi_1^\orb(M)=0$. Prove that
$M$ is equivalent to the Poincar\'e disk $\Delta$ and
has no orbipoints.

\item
Let $M$ be $\C P^1$ with two orbipoints of order $p$ and $q$.
\begin{enumerate}
\item Find $\pi_1^\orb(M)$.
 	
\item  Find a quasi-regular Sasakian\index[terms]{manifold!Sasaki!quasi-regular}
manifold $S$ with $S/\xi=X$.
\end{enumerate}

\item
Let $G$ be a finite group acting on a complex curve $X_1$,
$X= X_1/G$ be its quotient equipped with an 
orbifold structure. Consider a $G$-equivariant
holomorphic line bundle $L$ on $X_1$.
Let $\Tot^\circ(L)$ be the space of non-zero vectors in $\Tot(L)$.
Construct a $G$-equivariant line bundle $L$ with free action of
$G$ on $\Tot^\circ(L)$.
Prove that the associated map $\Tot^\circ(L)\arrow X$
is a locally trivial orbifold fibration. 

 \end{enumerate}

 
\chapter{Quasi-regular  foliations}
\label{_qr_foliation_Chapter_}\index[terms]{foliation!quasi-regular}


{\setlength\epigraphwidth{0.6\linewidth}
\epigraph{ 
{\cyr Byli ochi ostree tochimo{\u i} kosy -\\
Po zegzice v zenice i po kaple rosy, -\\ \medskip
 I edva nauchilis{\cprime} oni vo ves{\cprime} rost\\
Razlichat{\cprime} odinokoe mnozhestvo zvezd.} \\ \medskip
\font\tenit = cmssi17 at 10pt
\tenit
(Eyes once keener than a sharpened scythe -\\
In the pupil a cuckoo, a drop of dew -\\ \medskip

Now barely able to pick out, in full magnitude,\\
The lonely multitude of stars.)\\  \medskip

8--9 February 1937}{\sc\scriptsize Osip Mandelstam, translated by James Greene} 
}
\section{Introduction}

The theory of foliations is an important field of mathematics,
and we do not pretend to even start introducing it.
There is much literature on foliations, for example
\cite{_Molino_,_Lawson:Foliations_,_Camacho_Neto_,to, _Candel_Conlon_}.
In this chapter we give self-contained proofs
of some basic results, mostly targeted to the Reeb foliation
on Sasakian manifolds and the canonical\index[terms]{foliation!Reeb} foliation\index[terms]{foliation!canonical} on Vaisman
manifolds.

The divide between dynamics and foliations is quite thin.
A 1-dimensional oriented foliation is the same as a 
vector field up to a conformal multiplier. This allows
one to speak of the ``dynamics'' of a foliation. The dynamics
of a Reeb field is a huge subfield in contact geometry.\index[terms]{geometry!contact}
However, the dynamics of the Reeb field on a Sasakian
manifold is more predictable (Section \ref{_Closed_Reeb_orbits:Section_}).
The main reason is that a Sasakian manifold is
``transversally Riemannian''\index[terms]{foliation!transversally Riemannian}, and this puts 
severe restrictions on its ``holonomy'', defined below.

We describe the main objective of the present chapter.

We are interested in the {\bf regular foliations.}\index[terms]{foliation!regular}
Usually, a foliation ${\cal F}$ on $M$ is called {\bf regular}
if there exists a smooth submersion $\pi:\; M \arrow B$ such that
all leaves of ${\cal F}$ are fibres of $\pi$ (\cite{palais}).
It is not hard to see that this is equivalent to 
the existence of a base of open sets of $M$ intersecting
each leaf of ${\cal F}$ at most once. 


The regular contact manifolds are 
defined as thoses where the Reeb foliation is tangent
to the orbits of a free circle action \index[terms]{action!$S^1$-}(\ref{_regular_quasi_regular_sasakian_} ). To
obtain such an action on a manifold $M$ admitting
a smooth submersion $\pi:\; M \arrow B$ with circle fibres, we need to put
a metric on $M$ with all fibres of $\pi$ of the same
length, and arrange that the flow goes around the circles
with constant speed.

This argument shows that the dynamics of regular foliations is 
very simple. To work with dynamics of more complicated
foliations, one uses the notion of holonomy\index[terms]{foliation!holonomy of}
due to C. \index[persons]{Ehresmann, C.} Ehresmann (\cite{_Ehresmann_Shih:holonomy_}).
Given a leaf $F$ of a foliation ${\cal F}$
and a sufficiently small open set $U$ in it,
the foliated structure identifies the 
germs of the leaf space of ${\cal F}$ 
at each $x\in U$. Globally on $F$
this defines a homomorphism from
$\pi_1(F)$ to the group of diffeomorphisms of the germs of the
leaf space (\ref{_holonomy_folia_Definition_}).
This homomorphism is called {\bf the holonomy
of the foliation}. 

The holonomy of a foliation is related to the monodromy
of a flat connection (that is  actually\index[terms]{monodromy!of a connection} \index[terms]{connection!flat}
identified with its full holonomy group). Indeed, this is how
\index[persons]{Ehresmann, C.} Ehresmann originally discovered this notion
in 1940s. Let $(B, \nabla)$ be a flat vector
bundle over $M$. A connection $\nabla$ on $B$ defines
the linear Ehresmann connection \index[terms]{connection!Ehresmann}
$T\Tot(B) = T_\pi (\Tot B) \oplus T_\hor(\Tot B)$
(\ref{_Ehresmann_vs_vector_bundle_Proposition_}).
By \ref{_Ehresmann_curvature_bundles_Proposition_},
$\nabla$ is flat if and only if the distribution
$T_\hor(\Tot B)$ is integrable. We call the
corresponding foliation (defined whenever $\nabla$ 
is flat) {\bf the horizontal foliation}.

Let $F\subset \Tot(B)$
be the zero section of $B$, considered to be  a leaf
of the horizontal foliation. The holonomy
of the horizontal foliation in the leaf $F$ is 
equal to the full holonomy (that is, the monodromy) 
of $\nabla$.

We apply this notion to Vaisman and Sasakian
manifolds, that are  equipped with {\em transversally
Riemannian} foliations. A transversally Riemannian foliation
(\ref{_tra_Riemannian_foli_Definition_}) is a foliation
${\cal F}$ such that the bundle $TM/T{\cal F}$ 
is equipped with a metric $g$, defining a Riemannian metric
on each local leaf space (in other words, the metric
$g$ on $TM/T{\cal F}$ defines a basic semi-positive form on $TM$,
see \ref{_multilinear_forms_Remark_}). It is easy to see
that the holonomy of a transversally Riemannian foliation
preserves the Riemannian metric on the local leaf space.
This is the running assumption that we use.

Under this assumption, we prove a simplified version of
the Reeb stability theorem, showing that any transversally
Riemannian foliation that has trivial holonomy and all
leaves compact is actually regular 
(\ref{_regular_holonomy_finite_Proposition_}).
This is used further on to study the regular Vaisman and
Sasakian manifolds.\index[terms]{manifold!Vaisman!regular}

The main notion in this chapter, however, is not regularity,
but {\em qua\-si-re\-gu\-la\-ri\-ty}. We define {\bf quasi-regular foliations}
as those with all leaves compact. This is different from\index[terms]{foliation!quasi-regular}
the definition used in \cite{bog}, who defined 
quasi-regular foliations as ones where each leaf
meets a sufficiently small open subset a finite
number of times, globally bounded from above.
However, in Sasakian geometry\index[terms]{geometry!Sasaki} compactness
of all leaves is sometimes used as
a definition of quasi-regularity, conflicting with
\cite{bog} (see, for example, 
\cite{_Belgun_Moroianu_Semmelmann_} or \cite{_Sparks:SE_}).

When, in addition,
the leaves of the quasi-regular foliation are tangent to a torus action,
the quotient space is an orbifold (\ref{torus_quotient}).

To keep up with this notion, in this book we actually redefine the
notion of a regular foliation. For us, a foliation
is {\bf regular} if it is quasi-regular, and the projection
to its leaf space is a smooth submersion. When the
ambient manifold $M$ is compact, all the standard
notions of regularity (quasi-regularity) agree;\index[terms]{foliation!quasi-regular}
when it is not, then our definition of quasi-regularity 
is neither weaker nor stronger than the one in \cite{bog}.

The main purpose of this chapter is to prove that
all quasi-regular Sasakian manifolds are isometric to circle
bundles of unit vectors in positive holomorphic line
bundles over projective orbifolds. The local proof
is straightforward and obtained from the transversal
K\"ahler structure inherent in each Sasakian manifold.
However, to make it global we need to prove that
the Reeb action on a quasi-regular Sasakian manifold
factorizes through a circle\index[terms]{manifold!Sasaki!quasi-regular}
(\ref{_Reeb_action_factorizes_Corollary_}). 
This theorem is usually
deduced in a roundabout way from the work of \index[persons]{Wadsley, A. W.} Wadsley (\cite{_Wadsley_});
we give a new direct proof adapted for Sasakian geometry.\index[terms]{geometry!Sasaki}

 
\section{Quasi-regular foliations and holonomy}\index[terms]{foliation!quasi-regular}



\hfill

\definition
Recall that {\bf a leaf} of a foliation
is a maximal connected immersed submanifold
that is  everywhere tangent to the corresponding sub-bundle
(\ref{_leaf_foliation_Definition_}).
Let ${\cal F}$ be a foliation on $M$ 
and $U, U'\subset M$ open subsets. Denote by
 $U_{\cal F}$ and $U_{\cal F}'$ the sets of
leaves of ${\cal F}$ on $U$, $U'$. We say that a connected leaf
$F\in U_{\cal F}$  is {\bf an extension}
of $F'\in U'_{\cal F}$ if the union
$F\cup F'$ is a connected leaf in $U\cup U'$.

\hfill

\definition \label{_diffeo_of_germ_Definition_}
Let $x\in M$ be a point in a manifold, and
$u_1, u_2:\; U_x \arrow V_x$
diffeomorphisms mapping a neighbourhood of $x$ to
another neighbourhood of $x$. Write $u_1 \sim u_2$
if $u_1=u_2$ in some neighbourhood of $x\in M$
where both $u_1, u_2$ are defined. An
equivalence class of maps $u:\; U_x \arrow V_x$
is called {\bf a germ of diffeomorphisms of $M$ in $x$},
or {\bf a diffeomorphism of the germ}. 

\hfill

\remark
It is not very hard to show that one could compose germs
of diffeomorphisms, and that this composition defines
the group structure on the set of such germs.

\hfill

\definition\label{_holonomy_folia_Definition_}\index[terms]{foliation!holonomy of} 
{\bf (\index[persons]{Ehresmann, C.}Ehresmann)} \\ 
Let ${\cal F}$ be a smooth foliation on a 
manifold $M$, and $F$ a leaf that is  closed in $M$ (we call
such leaves ``closed leaves'').
Consider a loop $\gamma$ in $F$ passing through 
a point $x\in F$, and let $U$ be a neighbourhood 
of $x$ such that the leaf space $U_{\cal F}$ of ${\cal F}$ in
$U$ is well-defined. Choosing $U$ sufficiently
small, we may assume that $\gamma$ is covered by
a  finite collection $\{U_i\}$   of open sets 
intersecting sequentially, such that the
leaf space of ${\cal F}$ exists for all $U_i$,
and any leaf $F_1\subset U$ 
admits an extension to $U_1$, this
extension admits an extension to $U_2$ and so forth. Given a 
leaf $F_1 \in U_{\cal F}$, extend this leaf to all $U_i$
sequentially in the same order as $\gamma$ passes through $U_i$.
This gives an extension of $F_1$ in a tubular
neighbourhood $U_\gamma$ of $\gamma$. 
An extension of the leaf $F_1$ could, in theory, diverge 
from $U$ far enough so that the result of
successive extensions does not intersect $U$.
However, starting from $F_1$ sufficiently close to $F$
and passing along $\gamma$, we 
obtain a leaf $F'_1\in U_{\cal F}$ 
which intersects $U$ twice, before
the extension and after it.

This correspondence gives a diffeomorphism from
a neighbourhood of $F\in U_{\cal F}$ to $U_{\cal F}$ 
called {\bf holonomy  of ${\cal F}$ along $\gamma$.} Usually, the holonomy
is interpreted as a diffeomorphism of the germ of
the leaf space of ${\cal F}$ in $x$.

\hfill

\remark The same way we can define holonomy
for a path on a closed leaf $F$ connecting 
$x\in F$ to $y\in F$. The holonomy 
is a germ of diffeomorphisms from 
the germ at
$x\in M$ of the space of leaves of ${\cal F}$ to 
the germ at $y\in M$ of the space of leaves of ${\cal F}$.

\hfill

\claim\label{_holonomy_foliation_topological_Claim_}
The holonomy of the foliation along $\gamma$
in the closed leaf $F$
depends only on the homotopy class of 
the loop $\gamma$ in $F$.
This defines a homomorphism from $\pi_1(F, x)$
to the group of diffeomorphisms of the germ of
the leaf space in $x\in F$.

\hfill

\proof We leave the proof as an exercise.
One could also use the definition
of holonomy in terms of locally constant sheaves,
given below. In this interpretation, \index[terms]{sheaf!locally constant}
the holonomy of a foliation is the monodromy
of a local system on $F$, making\index[terms]{monodromy} 
\ref{_holonomy_foliation_topological_Claim_} trivial.
\endproof

\hfill

Alternatively, the holonomy of a foliation can be defined
in terms of locally constant sheaves. Let $F\subset M$
be a closed leaf of a foliation ${\cal F}$. We define a sheaf 
${\cal  G}$ over $F$ as follows. 

%

Let $U\subset F$ be an
open set in $F$ obtained as $U = U_M \cap F$, where 
$U_M\subset M$ is an open set. Denote by $\hat U$ the germ of
the closed subset $U=U_M \cap F\subset U_M$. Since any neighbourhood of $U$ in 
$U_M$ is foliated, the germ $\hat U$ is foliated.
We say that two functions defined in a neighbourhood of 
$U\subset U_M$ {\bf have the same germ in $U$}
if they are equal on a smaller neighbourhood of $U$;
this is the same as to say that these functions
agree on $\hat U$. We denote the ring of such 
functions by $C^\infty(\hat U)$.
Let ${\cal G}(U)$ be the ring of all 
$f\in C^\infty(\hat U)$  that are  constant
on the leaves of ${\cal F}$. Since any connected
leaf of ${\cal F}$ is uniquely determined
by its intersection with a smaller open set
$V_M \subset U_M$, the sheaf ${\cal G}$ is
locally constant, with the fibres identified with
the germs of functions constant on the leaves 
of ${\cal  F}$.

It is possible to prove that any automorphism
of this ring induced from a diffeomorphism of
open sets defines a germ of diffeomorphisms
of the leaf space of ${\cal F}$ in a neighbourhood
of $x\in F$ (see Exercise \ref{_germs_diffeo_Exercise_}).
This defines a map $\pi_1(M)\arrow \Diff ({\cal L}_x)$,
where ${\cal L}_x$ is the germ of the leaf space.
This is an alternative (more conceptual)
way of defining the holonomy of a foliation.

\hfill


\definition\label{_tra_Riemannian_foli_Definition_}
A {\bf transversally Riemannian foliation}\index[terms]{foliation!transversally Riemannian}
is a foliation ${\cal F}$ such that its
leaf space is locally equipped with a Riemannian
metric that is  compatible with the transition maps.
In other words, there is a Riemannian 
structure on the bundle $TM/T{\cal F}$, where
$T{\cal F}\subset TM$ is the bundle of vector
fields tangent to the leaves of ${\cal F}$,
that is  basic with respect to $T{\cal F}$.

\hfill

\example\label{_Sasakian_transversal_Example_}
Let $g$ be the Riemannian metric on
a Sasakian manifold, and $\eta$ the contact form.
The form $\eta$ is dual to the Reeb field, 
that has constant length 1 by \ref{_Sigma_defi_Claim_}.
Therefore, the form $g_0 := g - \eta\otimes \eta$ vanishes in the
Reeb direction and is invariant under
the Reeb flow. We obtain that $g_0$ is 
basic with respect to the Reeb foliation,
and defines the transversally Riemannian
structure (\ref{_multilinear_forms_Remark_}).

\hfill

\example\label{_Vaisman_transversal_Example_}
Let $M$ be a Vaisman manifold, and $\Sigma$
its canonical foliation.\index[terms]{manifold!Vaisman}
The form $\omega_0$ is transversally\index[terms]{foliation!canonical}
K\"ahler on the leaf space of $\Sigma$
(\ref{_Subva_Vaisman_Theorem_}),
hence $\Sigma$ is also tranversally 
Riemannian.

\hfill

\remark\label{_Riema_holono_folia_Remark_}
Clearly, the holonomy of a transversally  Riemannian
foliation preserves the transversal Riemannian metric. Therefore, it maps
a sufficiently small $\epsilon$-ball around a point 
$x$ in the leaf space of ${\cal F}$ in $U\ni x$
to itself. We will consider the holonomy of a 
transversally Riemannian foliation as 
a diffeomorphism of an
$\epsilon$-ball. For arbitrary foliations
this will not work, because the holonomy 
may contract or expand an open neighbourhood
in arbitrary (and unpleasant) ways, and 
it may act only on the  germs
of the leaf spaces.

\hfill

It is easy to see that any isometry of a connected Riemannian
manifold $M$ that fixes a point $x\in M$ and acts
trivially on $T_x M$ is trivial. This gives an
injective map $\St_x(\Iso(M)) \arrow \OO(T_x M)$,
where $\Iso(M)$ is the isometry group and
$\St_x(\Iso(M))$ the stabilizer of $x\in M$. 
This motivates the following proposition.

\hfill

\proposition
Let ${\cal F}$ be a transversally Riemannian
foliation, $F$ a closed leaf and $H$ the
holonomy group of ${\cal F}$ acting on the germ of the
space $M/{\cal F}$ of leaves of ${\cal F}$.
Taking the differential of the holonomy diffeomorphism  in 
$T_x (M/{\cal F})$, we obtain a group homomorphism
$r:\; H \arrow O\left (TM/T{\cal F}\restrict x\right)$
to the appropriate orthogonal group. 
Then, $r$ is injective.

\hfill

\proof
Let $U\ni x$ be a sufficiently small open set, and $X:= U/{\cal F}$
the leaf space of ${\cal F}$ in $x$. Since ${\cal F}$
is transversally Riemannian, $X$ is equipped with a
natural Riemannian structure.
By definition, $r$ maps a germ of diffeomorphisms
of $X$ to $\OO(T_x X)$. The kernel of such 
a map is an isometry of $X$ fixing $x$ and acting trivially
on $T_x X$; this map is trivial because it acts
trivially on all geodesics in $X$ passing through $x$.
\endproof 

\hfill

\definition\label{_qr_and_reg_foliation_Definition_}
A foliation ${\cal F}$ on $M$ is called {\bf quasi-regular}
if all its leaves are compact, and {\bf regular}\index[terms]{foliation!regular}\index[terms]{foliation!quasi-regular}
if it is quasi-regular, and the projection
$M \arrow M/{\cal F}$ to the leaf space is a locally
trivial fibration. In this case, the leaf
space is {\em a posteriori} a manifold as well.

\hfill

\proposition\label{_regular_holonomy_finite_Proposition_}
Let ${\cal F}$ be a transversally Riemannian, quasi-regular foliation.\index[terms]{foliation!transversally Riemannian}
Then it is regular if and only if its holonomy is trivial.

\hfill

{\bf Proof:}
Let $F$ be a closed leaf of ${\cal F}$, and
$U_\epsilon$ an $\epsilon$-neighbourhood of $[F]$
in the leaf space $X:= M/{\cal F}$. 
Given $x\in F$, consider a submanifold $S\subset U_\epsilon$
transversal and complementary dimension to $F$
(for example, we could take a 
union of geodesics passing through $x$ in a certain 
direction transversal to $F$). Choosing $\epsilon$ sufficiently
small, we can always assume that $S\cap U_\epsilon$
is a closed submanifold. 

Locally around $S$ we can identify $S$ with the leaf space
of ${\cal F}$, because $S$ intersects each leaf transversally.
Therefore, we can consider the holonomy of the foliation
 a diffeomorphism of $S$ (\ref{_Riema_holono_folia_Remark_}).

Consider the projection $\pi:\; S \arrow U_\epsilon/{\cal F}$.
The preimage of a point $s\in U_\epsilon/{\cal F}$ is the intersection
of the corresponding leaf $F_s$ and $S$, that is  finite because 
$F_s$ is compact and $S$ is closed. 

For each point $s \in S \cap F_s$, we can reach
$r$ by going along a path $\gamma_s$ in $F_s$. 
Projecting this path to $F$ and connecting its
ends within a small neighbourhood of $s$, 
we obtain a loop $\gamma$ in $F$.
By definition of holonomy, $r$
is the image of $s$ under the 
holonomy action associated with $\gamma$.
Therefore, $S \cap F_s$ is the orbit of
$s$ under the holonomy action.

Clearly, $S$ intersects
each leaf in the $\epsilon$-neighbourhood $U_\epsilon$ of $F$ 
once if and only if the projection
$U_\epsilon \arrow U_\epsilon/{\cal F}$ is locally trivial, for
sufficiently small $\epsilon$. This is equivalent
to the triviality of holonomy, because the
image of $S$ in $U_\epsilon/{\cal F}$
is the quotient by the holonomy action.
\endproof

\hfill

\remark
The famous Reeb stability theorem (\cite{_Reeb_stability_}) claims that\index[terms]{foliation!quasi-regular}
any codimension one foliation admitting a compact leaf with finite holonomy
is quasi-regular. We do not use this result, because
we always work with transversally Riemannian foliations,
which makes things much easier.\index[terms]{theorem!Reeb stability}

 
\section{Circle bundles over Riemannian orbifolds}


Recall that {\bf a Riemannian submersion}
is a smooth map $\pi:\; M\arrow M_1$
of Riemannian manifolds
such that for any $m\in M$ the
metric on the orthogonal complement
to the fibres of $\pi$ is equal to the
pullback $\pi^* g_1$ of the Riemannian \index[terms]{submersion!Riemannian}
metric on $M_1$.

\hfill

\definition
Let $(M,g)$ be a Riemannian manifold
equipped with a 1-dimen\-sional foliation ${\cal F}$.
We say that the Riemannian metric $g$
is {\bf compatible with a transversally
Riemannian structure} if the projection
$M \arrow M/{\cal F}$ to the leaf space (locally defined)
is a Riemannian submersion, for each open
subset of $M$ where the leaf space is defined.

\hfill

\example
The standard metric on a Sasakian manifold
is compatible with the transversally Riemannian
structure on the Reeb  foliations,
that is  clear from \ref{_Sasakian_transversal_Example_}.
The same can be established for
the canonical foliation on a \index[terms]{foliation!canonical}
Vaisman manifold.\index[terms]{manifold!Vaisman}

\hfill

\lemma\label{_Killing_transversal_Lemma_}
Let $(M,g)$ be a Riemannian manifold,
and $V_t$ a flow of isometries, $t\in \R$, without
fixed points.
Denote by ${\cal F}$ the 1-dimensional foliation whose leaves
are orbits of $V_t$. Then $g$ is compatible with a transversally
Riemannian structure on $(M, {\cal F})$.

\hfill

\proof
Let $T{\cal F}^\bot$ be the orthogonal complement
to the tangent bundle $T{\cal F}$ to the leaves of ${\cal F}$.
Since $V_t$ acts by isometries, the metric induced on $T{\cal F}^\bot$ 
is $V_t$-invariant. On the other hand, $V_t$ acts transitively on the
leaves of ${\cal F}$, and hence  the metric on $T{\cal F}^\bot$ 
can be obtained as a pullback from the leaf space (locally defined).
\endproof

\hfill

%

Let $(M, g, R)$ be a connected Riemannian manifold
equipped with a  vector field $R$,
and  ${\cal F}$  the foliation defined
by the orbits of the corresponding diffeomorphism
flow $e^{tR}$. Generally speaking, we have no control
over the dynamics of the action of $R$.
Even if all orbits are circles, it could
possibly happen that the rate of rotation
is different, and the map $t\mapsto e^{tR}$ does not factorize
through the circle (see the notes in Section \ref{_foliations_Notes_} for many
counterexamples). 
It turns out that when the vector
field  $R$ is Killing\index[terms]{vector field!Killing}
the rate of rotation is constant
(\ref{_Killing_flow_compact_orbts_Theorem_}).
We start the proof of 
\ref{_Killing_flow_compact_orbts_Theorem_}
by showing that the leaf space of a
transversally Riemannian rank 1 foliation
is an orbifold (this is a theorem due
to P. \index[persons]{Molino, P.} Molino in greater generality).

\hfill

\proposition\label{_S^1_qr_orbifold_Proposition_}
Let ${\cal F}$ be a quasi-regular rank 1 foliation
on a smooth manifold $M$. Assume that ${\cal F}$
is transversally Riemannian and\index[terms]{foliation!quasi-regular}\index[terms]{foliation!transversally Riemannian}
that the tangent bundle $T{\cal F}$ is oriented.
Then for each leaf $F$ of ${\cal F}$ there
is a foliated neighbourhood $W_F\subset M$ of $F$ such that
$T{\cal F}$ is tangent to a smooth circle
action on $W_F$. Moreover, the
quotient map $\pi:\; M \arrow M/{\cal F}$ is 
a locally trivial fibration in an open,
dense set of $M/{\cal F}$, and the leaf space
$M/{\cal F}$ is an orbifold.\index[terms]{orbifold} 

\hfill

\pstep
When $T{\cal F}$ is tangent to a smooth circle
action, $M/{\cal F}$ is an orbifold by
\ref{torus_quotient}. Since the notion of 
an orbifold is local, the local existence of
the circle action\index[terms]{action!$S^1$-} guarantees that 
$M/{\cal F}$ is globally an orbifold.

A leaf $F\in M/{\cal F}$ 
is called {\bf submersive} if 
$\pi:\; M \arrow M/{\cal F}$ is a locally trivial fibration
in a neighbourhood of $F$.
Clearly, a leaf $F$ is submersive
if and only if its holonomy is trivial.

\hfill

{\bf Step 2:}
We prove that the set of submersive leaves is open and 
dense in $M/{\cal F}$.

Let $F$ be a leaf of ${\cal F}$, and $U$ an open 
neighbourhood of $F$ in $M$.
Consider 
a local section  $S\subset M$ of ${\cal F}$, that is,
a closed submanifold of $U$ of complementary dimension
and transversal to $F$. Let $U_S$ be neighbourhood of $S$ in $U$.
Choosing $U$ and $U_S$ sufficiently small, we could 
identify $S$ with the leaf space
of ${\cal F}$ on $U_S$.

Under this identification, for any closed leaf
$F'$ in $U$, the set $F' \cap S$ can be identified
with the holonomy orbit of a point in $F'\cap S= U_0/{\cal F}$.

The number of intersections with $S$ is semi-continuous
as a function on $M/{\cal F}$, in the Hausdorff topology on
$M/{\cal F}$ considered to be  a set of closed subsets in $M$
(\ref{_Hausdorff_distance_Definition_}). Indeed, 
if a  leaf $F_\infty$ is a Hausdorff limit of a sequence
of leaves $F_1, F_2, ...$,
then $F_\infty\cap S$ is the Hausdorff limit
of $F_i\cap S$, and hence  the number of points
$\myhash(F_\infty\cap S)$ in
$F_\infty\cap S$ is bounded by 
$\overline\lim_i \myhash(F_i\cap S)$.

Submersive leaves
are the leaves where this intersection number is 
a local maximum. Indeed, in a neighbourhood $U_1$
of a submersive leaf $F'$, the set $S$ intersects it only once,
because $\pi$ is locally trivial and $S\cap U_1$ is
a section of $\pi$. Let $s_1, ..., s_k$ be the intersection
points of $F'$ and $S$, and $U_1, ..., U_k$
sufficiently small neighbourhoods of these points.
Any submersive leaf close enough to $F'$ intersects 
each of $U_i$ exactly once; hence  the number
$\myhash(F'\cap S)$ is locally constant
(and, by semi-continuity, maximal) on
the set of submersive leaves.

Conversely, if $\myhash(F'\cap S)$ 
is a local maximum on a leaf $F'$,
each leaf in a neighbourhood of $F'$
intersects $S\cap U_1$ exactly once,
defining a locally trivial projection
from the neighbourhood of $F'$ in the leaf
space to $S\cap U_1$.

By semi-continuity the 
set where $\myhash(F'\cap S)$
is a local maximum is
open and non-empty in $U$. Therefore, the
set of submersive leaves is dense in $M/{\cal F}$.

\hfill

{\bf Step 3:} 
Let $F$ be a leaf of ${\cal F}$,  $x\in F$ a point and $U\ni x$
a sufficiently small neighbourhood. 
Then the leaf space $U/{\cal F}$ is
a smooth Riemannian manifold.
The holonomy action induces 
an isometry on $U/{\cal F}$.
Since all leaves of ${\cal F}$ on $M$
are compact, the holonomy acts
on $U/{\cal F}$ with finite orbits.
Therefore, its action on the tangent space
$T_{F}(U/{\cal F})$ cannot have infinite
order (indeed, if the holonomy action has infinite order on
$T_{F}(U/{\cal F})$, some point 
$p\in U/{\cal F}$ close to $F$ 
would have an infinite orbit). Let $o$ be the order of
the action of the holonomy of ${\cal F}$ on $T_{F}(U/{\cal F})$.
An isometry of a Riemannian manifold 
preserving a point and acting trivially
on its tangent space is trivial. 
Denote by $\gamma$ the generator
of the holonomy  action 
$\Z=\pi_1(F)\arrow \Iso(U/{\cal F})$. Since $\gamma^o$
acts trivially on $T_{F}(U/{\cal F})$
and preserves the origin, it acts trivially
on $U/{\cal F}$, and the isometry $\gamma$ has finite order.

Define {\bf the index} $\ind F$ of a leaf $F$ as
the order of holonomy in $F$. Clearly, 
$\ind(F)=\myhash(F'\cap S)$, 
where $F'$ is a submersive leaf in 
a neighbourhood of $F$ and $S$ a
submanifold transversal to $F$ (Step 2).

\hfill

{\bf Step 4:} 
Choose a Riemannian metric $g$ on $M$, and let $|F|$ denote
the length of a leaf $F$. Then, for any
leaf $F$ obtained  as a limit of a sequence
of submersive leaves $F_1, F_2, ...$, 
we have $|F|= \ind(F)^{-1}\lim_i |F_i|$.
To see this, consider a sufficiently
small neighbourhood $U$ of $F$.
Then $F_i$ converges to $F$ together
with the unit tangent vector field.
Locally in a neighbourhood $V$ of a point $x\in F$,
$F_i\cap V$ is a union of connected segments 
$D_1(i), ..., D_k(i)$,  and by construction 
each segment is obtained from $D_1(i)$
by the holonomy action. Therefore, there
are $\ind(F)$ such segments. As $i$ goes to infinity,
each segment 
$D_l(i) \subset F_i\cap V$ converges to $F\cap V$,
and there are precisely $\ind(F)$ segments,
giving  $\ind(F) |F\cap V| = \lim_i |F_i \cap V|$.

\hfill

{\bf Step 5:} 
Let $l:\;M/{\cal F} \arrow \R^{>0}$
map a leaf $F$ to $|F| \cdot\ind(F)$.
As shown in Step 4, this function is
continuous on $M/{\cal F}$ with 
respect to the Hausdorff topology
induced by the Hausdorff metric
on the leaf space. 

Let $f:= \pi^* l$, where 
   $\pi:\; M \arrow M/{\cal F}$ is the projection map.
Rescale the Riemannian metric $g$ in $M$ 
replacing it with $f^{-1}g$. Then each 
leaf of ${\cal F}$ has length $\ind(F)^{-1}$.
Since $M/{\cal F}$ is locally compact,
and $\ind(F)$ is semi-continuous (Step 2), 
$\ind(F):\; M/{\cal F}\arrow \Z^{>0}$
is locally bounded. Indeed,
by semi-continuity, the set
$X_A:=\{H\in M/{\cal F}\ \ |\ \ \ind(H) \geq A\}$
is closed for any $A\in \Z$; hence  the
sequence $M/{\cal F}=X_1\supset X_2\supset X_3 \supset...$
terminates by Cantor's intersection theorem.

Let $W\subset M/{\cal F}$ be a subset with
compact closure, and $k:= \max_{F\in W} \ind(F)$.
Then the length of each $F\in W$ with respect
to the metric $k! g$ is integer, and 
the unit vector field tangent to ${\cal F}$
generates a circle action\index[terms]{action!$S^1$-} on $\pi^{-1}(W)$.
\endproof

\hfill

The following theorem will be applied to Sasakian
manifolds, but the statement is very general.

\hfill

\theorem\label{_Killing_flow_compact_orbts_Theorem_}
\index[terms]{vector field!Killing}
Let $(M, g, R)$ be a connected Riemannian manifold
equipped with a Killing vector field $R$,
and ${\cal F}$ the foliation defined
by the orbits of the corresponding diffeomorphism
flow $e^{tR}$. Assume that all orbits of $e^{tR}$ are compact. 
Then $e^{tR}$ defines a circle action\index[terms]{action!$S^1$-}, that is,
the action $t \mapsto e^{tR}$ factorizes through $S^1$.

\hfill

{\bf Step 1:}
Let $F$ be an orbit of $R$; since $F$
is compact, it is  diffeomorphic
to a circle. Since $e^{tR}$ acts on $F$ 
as rotation with a constant speed, 
its action on $F$ is factorized through a 
$S^1$ acting on a circle $F$ by isometries.
Let $\sigma(F)$ be the smallest $t\in \R^{>0}$ such that
$e^{t R}$ acts trivially on $F$ (``the first return 
function''). Clearly, $\sigma(F)$ is the generator
of the kernel of the corresponding surjective map $\R \arrow S^1$.
By construction, $F \mapsto \sigma(F)$ defines
a function on the set of all orbits of $R$.
We need to show that for a dense subset
of orbits, $\sigma(F)$ is constant,  $\sigma = s$.
Then $e^{s R}$ acts as identity on $M$, because
it is the identity on a dense subset.

\hfill

{\bf Step 2:}
Let $M_0\subset M$ be the union of all
submersive leaves of ${\cal F}$
(see the proof of 
\ref{_S^1_qr_orbifold_Proposition_} for the definition).
By \ref{_S^1_qr_orbifold_Proposition_},
$M_0$ is open and dense in $M$ 
and the projection $M_0 \arrow M_0/{\cal F}$
is a locally trivial fibration.

The function $y \mapsto \sigma(F_y)$ is continuous on 
$M_0/{\cal F}.$ Let $F$ be an orbit, and $U_\epsilon$ 
an $\epsilon$-neighbourhood of $F$ in $M_0$. Since 
$R$ is Killing,\index[terms]{vector field!Killing} the  set $U_\epsilon$
is foliated by the circles. Unless
$\sigma$ is constant,  we may choose 
an orbit $F_y\subset U$ such that the quotient
$\frac{\sigma(F)}{\sigma(F_y)}$ is irrational.

We rescale $R$ by  a constant
multiplier in such a way that $\sigma(F)=1$.
Since $\frac{\sigma(F)}{\sigma(F_y)}$ is irrational, 
the set of all points $e^{n R}(y)$, with $n \in \Z$,
is dense in the circle $F_y$ (a rotation with irrational
period has dense orbit in a circle). 

Then there exists a time $n\in \Z\subset \R$ such that
$e^{n R}(y)$ belongs to any given open subset 
in $F_y \cap U_\epsilon$, and  $e^{n R}(x)=x$.
Since the map $e^{n R}$ 
acts by isometries and $\diam(S)< 2\epsilon$,
we have $d(e^{n R}(y), x)< 2\epsilon$ for
all $x, y\in S$, and all $n\in \Z$.
Applying the triangle inequality to the\index[terms]{inequality!triangle}
triangle $x, y, e^{n R}(y)$, we obtain
that $d(y, e^{n R}(y)) < d(y,x) + d(x, e^{n R}(y)) < 4\epsilon$.
Since $e^{n R}(y)$ can be chosen in any
given small open set of $F_y$, this implies
that  $\diam(F_y) < 4 \epsilon$.

Since $F_y$ can be chosen in a given open, dense
sets of orbits, this argument actually implies
that diameter of any orbit of $R$ is bounded by
$4\epsilon$, where $\epsilon$ is chosen arbitrarily
small; this is clearly impossible.
\endproof

\hfill

\corollary\label{_Reeb_action_factorizes_Corollary_}
Let $S$ be a connected quasi-regular Sasakian manifold
(i.\,e.   with all Reeb orbits compact).\index[terms]{manifold!Sasaki!quasi-regular}
Then the Reeb action $e^{tR}$ factorizes through the
circle.\index[terms]{vector field!Reeb}

\hfill

\proof
The Reeb flow acts on $S$ by isometries, and hence
\ref{_Killing_flow_compact_orbts_Theorem_}
applies. \endproof

\hfill

\remark
In Chapter \ref{CR_and_such}, we defined 
regular and quasi-regular Sasakian manifold
in terms of the Reeb action (\ref{_regular_quasi_regular_sasakian_}).
From \ref{_Reeb_action_factorizes_Corollary_}
it follows that this definition is consistent
with \ref{_qr_and_reg_foliation_Definition_}. 

\hfill

Further on, we will need a version of this theorem
for 2-dimensional action.

\hfill

\theorem\label{_Killing_flow_Vaisman_orbits_Theorem_}
Let $(M, g)$ be a connected Riemannian manifold
en\-dowed with a pair of commuting Killing vector\index[terms]{vector field!Killing} fields
$R_1, R_2$, and ${\cal F}$ the foliation defined
by the orbits of the corresponding diffeomorphism
flow $e^{t_1R_1+t_2R_2}$. Assume that all orbits of 
this 2-dimensional flow are compact, and $R_1$, $R_2$
are linearly independent everywhere. 
Then $e^{t_1R_1+t_2R_2}$ defines a torus action, that is,
the action $(t_1, t_2) \mapsto e^{t_1R_1+t_2R_2}$ 
factorizes through $S^1\times S^1$.

\hfill

\proof The proof is almost identical to the proof of 
\ref{_Killing_flow_compact_orbts_Theorem_}.
For more details, see Exercise 
\ref{_Killing_flow_Vaisman_orbits_Exercise_}.
\endproof

 
\section{Quasi-regular Sasakian manifolds}\index[terms]{manifold!Sasaki!quasi-regular}


Recall that a Sasakian manifold is called 
{\bf quasi-regular} if its Reeb foliation has compact
leaves and {\bf regular} if it is regular and the corresponding 
quotient map is a locally trivial fibration.

\hfill

\remark
Let $S$ be a (quasi-)regular Sasakian manifold.\index[terms]{manifold!Sasaki!regular}
By \ref{torus_quotient} and \ref{_Reeb_action_factorizes_Corollary_}, 
the Reeb orbit space
of $S$ is naturally equipped with an orbifold
structure; it is smooth\footnote{That is,
its orbifold structure is trivial. Saying that an
orbifold is ``smooth'' is ambiguous, because
the quotient space can be smooth but equipped
 with a non-trivial orbifold structure;
this is what happens for some Sasakian
structures on $S^3$, for example; see 
Exercises \ref{_qr_Sasakian_S^3_Exercise_}
and \ref{_orbifold_CP^1_Exercise_}.} 
when $S$ is regular.

\hfill

\theorem\label{_quasireg_Sasakian_orbibundles_Theorem_}
Let $S$ be a (quasi-)regular Sasakian manifold,\index[terms]{manifold!Sasaki!quasi-regular}\index[terms]{manifold!Sasaki!regular}
$X$ the leaf space of its Reeb fibration,
and $\pi:\; S \arrow X$ the quotient map.
Then: 
\begin{description}
\item[(i)] The space $X$ is equipped with a natural structure
of a K\"ahler orbifold. Moreover, $S$ 
is identified with  the space of unit vectors in
an orbifold holomorphic line bundle $L$ over $X$.

\item[(ii)] There exists a Hermitian
metric on $L$ such that its curvature is 
 the K\"ahler metric on $X$. 

\item[(iii)] 
Moreover, the Sasakian manifold
$S$ is isometric to the space $S^1(L)$ of 
vectors of length 1 in the total 
space $\Tot(L)$. Here we consider 
the total space $\Tot(L)$
to be a metric space\index[terms]{space!metric} with the K\"ahler metric $dd^c l$
defined in \ref{_dd^c_log_l_Corollary_}.\footnote{Note that, 
in order to apply  \ref{_dd^c_log_l_Corollary_}, we need $L^*$ to be a 
positive holomorphic line bundle; this is provided by
\ref{_quasireg_Sasakian_orbibundles_Theorem_} (ii).}
In this case, the cone $C(S)$ is holomorphically
isometric to the space $\Tot^\circ(L)\subset \Tot(L)$
of non-zero vectors in $L$.

\item[(iv)] Conversely, if $L$ is a negative holomorphic line bundle
over an orbifold and $S^1(L)$ is smooth, then $S^1(L)$ is Sasakian.
\end{description}

{\bf Proof of \ref{_quasireg_Sasakian_orbibundles_Theorem_} (i):} 
Let $C(S)$ be the Riemannian cone\index[terms]{cone!Riemannian} of $S$,
endowed with the natural K\"ahler
structure. The Reeb action on $S$
is extended to the action of $\C$ on
$C(S)$ generated by the Reeb field
$R$ and $I(R)$. Since $R$ has compact
orbits, this action is actually 
factorized through $\C^*$ (\ref{_Reeb_action_factorizes_Corollary_}). Indeed,
the total space of the corresponding $\U(1)$-bundle
is $S$. 

Let $X$ be the space of Reeb orbits in $S$.
This space is identified with the space
of $\C^*$-orbits in the cone $C(S)$,
hence it is a complex manifold. 
The transversal K\"ahler structure on $X$
is given by $dd^c\log t$, where $t$ is 
the parameter on the generator of $C(S)$
(\ref{_Sigma_defi_Claim_}).

By \ref{_quotient_principal_Claim_},
the quotient map $S \arrow X$
is a principal orbifold $S^1$-fibration.\index[terms]{bundle!principal}
Multiplying this $S^1$-fibration with $\R^{>0}$, we
obtain a principal orbifold $\C^*$-bundle; it is
clear that $C(S)$ is its total space. Denote the
corresponding orbifold line bundle over $X$ by $L$.
It is holomorphic because the $\C^*$-quotient map
$C(S) \arrow X$ is holomorphic.  This 
identifies $C(S)$ with the total space of an 
orbifold holomorphic line bundle $L$ on $X$.

\hfill

{\bf Proof of \ref{_quasireg_Sasakian_orbibundles_Theorem_} (ii):} 

The standard metric $t^2 g + (dt)^2$ on the cone $C(S)$ is 
homogeneous of order 2 in $t$ and invariant with respect
to $I(R)= I \frac d {dt}$. 
Therefore, it is restricted to a flat Hermitian metric
on each orbit of the corresponding $\C$-action.
This implies that $t^2 g + (dt)^2$ is 
proportional to an 
Hermitian metric on the $\C^*$-bundle $C(S)\arrow X$.
and its length function $l(v)= |v|^2$ on $\Tot(L)$ is equal to
$t^2$.

The transversally K\"ahler metric on $C(S)$
is expressed as $dd^c \log (t)$ by \ref{_Sigma_defi_Claim_}.
However, $t^2$ is the length function on $\Tot^\circ (L)$.
By \ref{_curvature_line_bundle_via_15.19_Corollary_},
the curvature $\Theta_L$ of $L$ can be written
(up to a multiplier of $\1$) as $dd^c \log (t^2)$. Then
$-\1\Theta_L$ is equal to the 
K\"ahler metric on $X$.

\hfill

{\bf Proof of \ref{_quasireg_Sasakian_orbibundles_Theorem_} (iii):} 
The map $S \arrow C(S)$ sending $s$ to $(s, 1)$
is by construction an isometry. However, the
metric on $C(S)$ is written as $dd^c(l)= dd^c(t^2)$
by \ref{_dd^c_log_l_Corollary_} and \ref{_Sigma_defi_Claim_}.

\hfill

Finally, \ref{_quasireg_Sasakian_orbibundles_Theorem_} (iv)
is directly implied by \ref{_curvature_line_bundle_via_15.19_Corollary_}.
Indeed, the total space $\Tot^\circ(L)$ of $L$ is K\"ahler.
The homothety action on $\Tot^\circ(L)$ is generated by the
radial vector field, that is  dual to $d l$,
where $l(v)=|v|^2$.  By \ref{homoth}, $\Tot^\circ(L)$
is a Riemannian cone\index[terms]{cone!Riemannian} over $S^1(L)$.
\endproof

\section{Notes}
\label{_foliations_Notes_}

A foliation on a Riemannian manifold
is called {\bf Riemannian foliation}\index[terms]{foliation!Riemannian}
if the restriction of the Riemannian
metric to the normal bundle $T{\cal F}^\bot$
is transversally Riemannian (\cite[Chapter 5]{_Tondeur:Geometry_}).
The study of Riemannian foliations was initiated by
B. \index[persons]{Reinhart, B.} Reinhart (\cite{_Reinhart:1959_}). We refer
to \cite{to,_Tondeur:Geometry_,_Molino_,_Reinhart:book_}
for the background and references. It is clear
that any Riemannian foliation is equipped with
a transversally Riemannian metric. Conversely,
a transversally Riemannian foliation ${\cal F}$ on $M$ can be equipped
with a Riemannian foliation structure by taking
an arbitrary direct sum decomposition
$TM= T{\cal F} \oplus N{\cal F}$, and putting
a metric on $TM$ in such a way that this
decomposition is orthogonal. 

For us, the main utility of the 
Riemannian (and  transversally Riemannian)
foliations is a theorem due to \cite{_Molino_},\index[terms]{theorem!Molino}
that states that a leaf space of a Riemannian
foliation with compact fibres is always an
orbifold. We prove this result for foliations
of rank 1 (\ref{_S^1_qr_orbifold_Proposition_}).

\index[persons]{Molino, P.} Molino's theorem is related to the following
question. Suppose that ${\cal F}$ is a rank 1 foliation
with compact fibres. Is it true that 
there exists a circle action\index[terms]{action!$S^1$-} tangent to ${\cal F}$?
Of course, existence of the circle action
implies that the leaf space is an orbifold
(\ref{torus_quotient}). 

When the foliation is Riemannian, a positive 
answer was obtained by A. W. \index[persons]{Wadsley, A. W.} Wadsley (\cite{_Wadsley_}).\index[terms]{theorem!Wadsley}

Without the Riemannian assumption, this
question was the subject of much research
(\cite{_Epstein_}, \cite{_Wadsley_}, 
\cite{_Sullivan:counterexamples_}).
D. \index[persons]{Sullivan, D.} Sullivan, who produced a counterexample
in dimension 5 (\cite{_Sullivan:counterexamples_})
calls it ``the periodic flow conjecture''.
\index[persons]{Sullivan, D.} Sullivan found a rank 1 foliation with compact \index[terms]{conjecture!Sullivan (periodic flow)}
leaves on a product of $\R P^3$ and a disk
and (after doubling this example to hide
the boundary) on $S^5$.

Later, Epstein\index[persons]{Epstein, D. B. A.} and Vogt\index[persons]{Vogt, E.} constructed a counterexample
on a compact manifold of dimension 4
(\cite{_Epstein_Vogt_}).

In this chapter we have re-proved \index[persons]{Wadsley, A. W.} Wadsley's theorem and 
used it to obtain the orbifold structure on the leaf space,
with an application to Sasakian geometry.\index[terms]{geometry!Sasaki}

We are grateful to Jorge Vit\'orio \index[persons]{Pereira, J. V.} Pereira
and participants of Mathoverflow for directing
us to the relevant reference.\footnote{\tiny \url{https://mathoverflow.net/questions/401470/smooth-rank-1-foliations-with-closed-leaves/401568}.}

\section{Exercises}

\begin{enumerate}[label=\textbf{\thechapter.\arabic*}.,ref=\thechapter.\arabic{enumi}]	

\item
Let $(B, \nabla)$ be a flat vector
bundle over $M$. Consider the corresponding Ehresmann connection \index[terms]{connection!Ehresmann}
$T\Tot(B) = T_\pi (\Tot B) \oplus T_\hor(\Tot B)$
(\ref{_Ehresmann_vs_vector_bundle_Proposition_}).
By \ref{_Ehresmann_curvature_bundles_Proposition_},
$\nabla$ is flat if and only if the distribution
$T_\hor(\Tot B)$ is integrable. 
Let $F\subset \Tot(B)$
be the zero section of $B$, considered to be  a leaf
of the horizontal foliation. Show that the holonomy
of the horizontal foliation in the leaf $F$ 
acts by linear automorphisms on each fibre of $B$. Show 
that the holonomy of the horizontal foliation is
equal to the monodromy of $\nabla$.\index[terms]{monodromy!of a connection}

\item (\cite{_Thurston:slithering_}, Example 2.8)\\
Let $\phi:\; M \arrow M$ be an automorphism of the
smooth manifold $M$, and $M\times [0,1]/(0,m)\sim(1, \phi(m))$
the corresponding mapping torus (some authors call  
this construction  {\bf the suspension} of $\phi$, see for example   
\cite{_Molino_,_Richardson:suspension_}).\index[terms]{foliation!suspension}
The {\bf suspension foliation} is a rank 1 foliation
with leaves tangent to $\{m\}\times [0,1]$.
Let $m_0\in M$ be a fixed point of $\phi$,
and $F=[0,1]/(0,m_0)\sim(1, m_0)$ the
corresponding closed leaf of the suspension foliation ${\cal  F}$.
Prove that the holonomy of ${\cal F}$ in $F$ 
maps the generator of $\pi_1(F)=\Z$ to $\phi$.

\item Construct a foliation with no transversally  Riemannian structure.

{\em Hint:} Use the previous exercise.

\item Construct a rank 1 foliation that admits
a transversally Riemannian structure, and
has a compact leaf, but is not quasi-regular.\index[terms]{foliation!quasi-regular}

\item
Let ${\cal F}$ be a regular foliation
on $M$ admitting an action\index[terms]{foliation!regular}
of $T^n$ that is  free and transitive on the leaves.
Denote by $\pi:\; M \arrow X$ the projection to the leaf space.
Consider the local system $R^1\pi_*(\C_M)$
of fibrewise first cohomology 
(``Gauss--Manin local system'').
Prove that it is trivial.

\item Construct a regular foliation
with a fibre $T^n$ and not admitting an action of $T^n$
that is  free and transitive on the leaves.

{\em Hint:} Use the previous exercise.

\item Let ${\cal F}$ be a transversally Riemannian foliation
with at least one closed fibre $F$ diffeomorphic to a sphere $S^n$, $n>1$.
Prove that ${\cal F}$ is regular in a neighbourhood 
of $F$.

\item (a version of Reeb stability theorem)
Let ${\cal F}$ be a foliation
with a compact fibre $F$ with the fundamental group \index[terms]{fundamental group}
$\pi_1(F)$ finite. Prove that ${\cal F}$ is quasi-regular 
in a neighbourhood of $F$.\index[terms]{theorem!Reeb stability}

\item 
Let $M$ be a compact manifold, and $\goth m$ a maximal ideal
in the ring $C^\infty M$ of all smooth real functions. 
Prove that $\goth m$ is the ideal
of all functions vanishing in a point.

\item
Let $M$ be a manifold, not necessarily compact, and $\goth m$ a maximal ideal
in the ring $C^\infty M$. Assume that the natural map
$\R \arrow C^\infty M/\goth m$ is an isomorphism.
Prove that $\goth m$ is the ideal
of all functions vanishing in a point.
Further on, we denote by $\Spec_\R(C^\infty M)$
the set of maximal ideals of this type; clearly,
$\Spec_\R(C^\infty M)$ is identified with $M$.

\item
Let $M$ be a manifold and $A$ an $\R$-linear 
automorphism of $C^\infty M$.
Then the corresponding map 
$A:\; \Spec_\R(C^\infty M)\arrow \Spec_\R(C^\infty M)$
is a homeomorphism with respect to the topology
induced from $M=\Spec_\R(C^\infty M)$.

\item
Let $M$ be a smooth manifold, and $\Aut_\R(C^\infty M)$
the group of all $\R$-linear automorphisms of the ring
$C^\infty M$. Prove that the natural embedding
of  $\Diff(M)$ in $\Aut_\R(C^\infty M)$
is an isomorphism.

{\em Hint:} Use the previous exercise.

\item\label{_germs_diffeo_Exercise_}
Let $x\in M$ be a point in a manifold,
and $C^\infty_x(M)$ the algebra of germs of
$C^\infty M$ in $x$. Denote by $\widetilde \Aut_\R(C^\infty_x (M))$ 
the group of automorphisms of $C^\infty_x (M)$ induced by
diffeomorphisms $U\arrow V$ of open neighbourhoods $U, V$ of $x\in M$.
Prove that $\widetilde \Aut_\R(C^\infty_x (M))$ is isomorphic
to the group of diffeomorphisms of the
germ of $M$ in $x$ defined in 
\ref{_diffeo_of_germ_Definition_}.

\item Prove that the 3-dimensional torus $T^3$ 
does not admit a regular Sasakian structure.\index[terms]{manifold!Sasaki!regular}

\item Let $M$ be a regular 3-dimensional Sasakian manifold
with infinite fundamental group.\index[terms]{fundamental group} Prove that
$\pi_1(M)$ is not abelian.

\item Construct a regular Sasakian structure on $\R P^3$.

\item Let $M$ be a simply connected regular Sasakian
manifold, $\dim S=3$. Prove that $M$ is diffeomorphic to $S^3$.

\item Prove that any regular Sasakian 3-manifold with  finite fundamental group 
is a lens space.\index[terms]{lens space}

\item
Construct a positive orbifold line bundle $L$
on the weighted projective line (\ref{weight_CP}). Prove that the
Sasakian manifold obtained
as the space of length one vectors in $L$
is diffeomorphic to $S^3$.

\item \label{_qr_Sasakian_S^3_Exercise_}
Construct a quasi-regular but non-regular 
Sasakian structure on $S^3$.\index[terms]{manifold!Sasaki!regular}

{\em Hint:} Use the previous exercise.

\item\label{_orbifold_CP^1_Exercise_}
Apply \ref{torus_quotient} to the Reeb action on a quasi-regular Sasakian\index[terms]{manifold!Sasaki!quasi-regular}
structure on $S^3$ defined in \ref{_qr_Sasakian_S^3_Exercise_}.
Prove that the quotient orbifold
is an orbifold $\C P^1$ with two orbifold points with
cyclic isotropy groups.

\item\label{_Killing_flow_Vaisman_orbits_Exercise_}
Let $M$ be a Riemannian manifold equipped with
a pair of commuting vector fields $R_1$, $R_2$.
 Assume that $R_1$, $R_2$
are linearly independent everywhere  and generate a  
2-dimensional flow with compact orbits.
Denote by ${\cal F}$ the corresponding
2-dimensional foliation, and let $\pi:\; M \arrow M/{\cal F}$
be the projection to its leaf space. The following
sequence of exercises  proves
\ref{_Killing_flow_Vaisman_orbits_Theorem_},
using the proof of 
\ref{_Killing_flow_compact_orbts_Theorem_}
as a framework.
\begin{enumerate}
\item Call an orbit ${\goth O}$ {\bf submersive}
if the projection $M \arrow M/{\cal F}$
is locally trivial in a neighbourhood of
${\goth O}\in M/{\cal F}$.
Prove that the set of submersive orbits
is dense in $M/{\cal F}$.

\item Define {\bf the index} of an orbit
${\goth O}$ as the order of the
holonomy of the foliation ${\cal F}$
in  ${\goth O}$. Prove that 
the index is finite and semi-continuous
in $M/{\cal F}$.

\item Replace $M$ by an open subset
that is  locally trivially $T^2$-foliated with all orbits submersive.
For any $u\in M/{\cal F}$, let $\Lambda_u \subset \R^2$
be a rank 2 sublattice such that the
$\R^2$-action on $\pi^{-1}(u)$ factorizes through
$\R^2/\Lambda_t$. Since the family
$\pi:\; M \arrow M/{\cal F}$ is locally trivial, the lattice\index[terms]{lattice}
$\Lambda_t$ depends continuously on $t$. Fix an orbit
${\goth O} \subset M$.
Suppose that for a general $t$ in any neighbourhood
of ${\goth O}$, one has $\Lambda_t \cap \Lambda_0=0$.
Prove that $\diam({\goth O}_t) \leq 4 d({\goth O}_t,
{\goth O})$.

\item Suppose that for the general $t$ in any neighbourhood
of $0$ one has $\rk \Lambda_t \cap \Lambda_0=1$.
Let $R$ be the generator of the rank 1 group
 $\Lambda_t \cap \Lambda_0$ (for general $t$)
considered  a Killing vector field \index[terms]{vector field!Killing}on $M$.
Prove that an orbit of $R$ is a circle.
Prove that ${\goth O}_t$ lies in a 
$4\epsilon$-neighbourhood of an orbit
of $R$, where $\epsilon= d({\goth O}_t,
{\goth O})$. 

\item Deduce \ref{_Killing_flow_Vaisman_orbits_Theorem_}
from these arguments.
\end{enumerate}

\end{enumerate}

 
\chapter{Regular and quasi-regular  Vaisman
 manifolds}\index[terms]{manifold!Vaisman}\index[terms]{manifold!Vaisman!quasi-regular}\index[terms]{manifold!Vaisman!regular}

\epigraph{\it
The heavens and earth have great beauty, but no words. The
four seasons have clearly defined patterns, but no
debates. All living things have basic principles, but no
explanations. A wise person begins with the beauty of the
heavens and earth and arrives at the principles of all
living things. This is what is called the Root and Origin:
to be able to observe from the viewpoint of the heavens!}
{\sc\scriptsize Chuang Tz\v u, Chapter 22: Knowledge Wanders North}


\section{Introduction}
\label{_qr_chapter_Intro_Section_}

{\bf Quasi-regular} Vaisman manifolds are Vaisman manifolds
equipped with an elliptic fibration (possibly with multiple fibres).\index[terms]{manifold!Vaisman}\index[terms]{manifold!Vaisman!quasi-regular}
They are called {\bf regular} when this fibration is smooth,
that is, when it has no multiple fibres.\footnote{For another,
equivalent version of these definitions see \ref{_qr_Definition_}.}\index[terms]{manifold!Vaisman!regular}

Every Vaisman manifold $M$ is equipped with a foliation,
called ``the canonical foliation'' $\Sigma$ (Section
\ref{_gauduchon_basics_Subsection_}).\index[terms]{foliation!canonical}
This foliation is transversally K\"ahler;
moreover, the transversal K\"ahler form\index[terms]{form!K\"ahler!transversal}
is globally defined on $M$ and exact
(\ref{_Subva_Vaisman_Theorem_}).
This implies that any complex 
subvariety of $M$ is tangent 
to $\Sigma$. In particular,
whenever $M$ admits an elliptic fibration, all
leaves of $\Sigma$ are elliptic curves.

Since $\Sigma$ is transversally K\"ahler, the
base of this elliptic fibration is a K\"ahler
orbifold. In fact, it is projective
(\ref{_Structure_of_quasi_regular_Vasman:Theorem_}).

This is related to an explicit construction
of Vaisman manifolds, defined in \cite[Section 2.3]{ov_imm_vai}. \index[terms]{manifold!Vaisman} 
Let $L$ be an anti-ample bundle over a projective orbifold $X$
and $\Tot^\circ(L)$ the space of non-zero vectors in 
its total space. Then $\Tot^\circ(L)$ is the cone of
$M$, in the sense of algebraic geometry. \index[terms]{geometry!algebraic}

Assume that $\Tot^\circ(L)$ is smooth.
This assumption imposes some restrictions on the
singularities of $X$, but some orbifold singularities
are still allowed; for example, the cone of a weighted
projective space is $\C^n \backslash 0$, that is  smooth.

When $L$ is equipped with a metric with negative curvature,
the space $\Tot^\circ(L)$ inherits the conical K\"ahler 
metric $dd^c l$, where $l(v)=|v|^2$ (\ref{_formula_15_19_besse_Theorem_}).
Taking a quotient of $\Tot^\circ(L)$ by a $\Z$-action\index[terms]{action!$\Z$-}
generated by a homothety $\phi$, we obtain an LCK manifold $M$.\index[terms]{manifold!LCK}
It turns out that this manifold is Vaisman. The simplest
example is obtained when $\phi=h_\lambda$ is the homothety
mapping $v$ to $\lambda v$, where $\lambda$ 
is a complex number,  $|\lambda|>1$.
In this case, $M$ is clearly elliptically fibred over $X$, 
hence quasi-regular. \index[terms]{manifold!Vaisman!quasi-regular}


Suppose  the homothety $\phi:\; \Tot^\circ(L)\arrow \Tot^\circ(L)$
preserves the fibres of the projection $\Tot^\circ(L) \arrow X$,
and acts with finite order on $X$. It turns out that
then the manifold $M=\Tot^\circ(L)/\langle \phi\rangle$
is a quasi-regular Vaisman manifold, and, moreover,
any quasi-regular Vaisman manifold is obtained in this way
(\ref{_Structure_of_quasi_regular_Vasman:Theorem_}). \index[terms]{manifold!Vaisman}\index[terms]{manifold!Vaisman!quasi-regular}

This theorem has many applications owing to another result,
also proven in \cite{ov_imm_vai}: any compact
Vaisman manifold $M$ admits a deformation that is  quasi-regular,
and, moreover, $M$ can be approximated by quasi-regular Vaisman
manifolds (\ref{defovai}). This allows one to consider only
quasi-regular (that is, elliptic) Vaisman manifolds when
studying the topology.

Regular Vaisman manifolds were first discussed in
\cite{va_rendiconti, va_gd}. Vaisman\index[persons]{Vaisman, I.}  showed that the
leaves of the foliation $\ker\theta$ carry Sasakian
structures then proved that a compact Vaisman manifolds
with regular Lee flow\index[terms]{Lee flow} fibres in circles over a compact
Sasakian manifold. He also showed that the fibration is a
Riemannian submersion. Then he showed that if the\index[terms]{submersion!Riemannian}
canonical foliation is regular, the manifold fibres in
complex one-dimensional tori over a projective
manifold. Finally, he identified this projective manifold
as the base of the  Boothby--Wang fibration\index[terms]{Boothby--Wang fibration} associated with 
the Sasakian manifold above, thus obtaining a commutative
triangle of principal fibre bundles and Riemannian
submersions (\ref{triangle}).\index[terms]{bundle!principal}

 
\section[Quasi-regular Vaisman manifolds as cone quotients]{Quasi-regular Vaisman manifolds\\ as cone quotients}\index[terms]{manifold!Vaisman}\index[terms]{manifold!Vaisman!quasi-regular}


\definition\label{_qr_Definition_}
A Vaisman manifold $M$ is called {\bf regular} if
the leaves of the canonical foliation $\Sigma$ 
are orbits of\index[terms]{foliation!canonical}
the group $T^2=(S^1)^2$ freely acting on $M$,  {\bf 
	quasi-regular} if these leaves are compact, and 
{\bf irregular} otherwise.
\index[terms]{manifold!Vaisman!regular}\index[terms]{manifold!Vaisman!quasi-regular}

\hfill

\remark\label{_qr_Vaisman_defi_Remark_}
In Chapter \ref{_qr_foliation_Chapter_},
we defined the notion of regular and quasi-regular 
foliation in a general situation
(\ref{_qr_and_reg_foliation_Definition_}). A Vaisman manifold
is regular (quasi-regular) if its canonical foliation
is regular (quasi-regular); this follows, for example,
from\index[terms]{foliation!canonical} \ref{_Killing_flow_Vaisman_orbits_Theorem_}.\index[terms]{manifold!Vaisman}

\hfill

Let $M$ be a Vaisman
manifold of LCK rank 1\index[terms]{rank!LCK}. By the Structure theorem\index[terms]{rank!LCK}
(\ref{str_vai}), $M$ is obtained 
as a quotient of a cone $C(S)$ over a Sasakian manifold
by an automorphism $(s, t) \mapsto (\phi(s), \lambda t)$
where $\lambda >1$ is a real constant, and $\phi:\; S \arrow S$
is a Sasakian automorphism.

\hfill

\theorem\label{_qreg_Vaisman_via_Sasakian_Theorem_}
Let $M$ be a Vaisman manifold\index[terms]{manifold!Vaisman}
obtained as a $\Z$-quotient of a cone $C(S)= S\times \R^{>0}$ over a Sasakian
manifold  $\pi:\; C(S)\arrow M= C(S)/\Z$.
Then:
\begin{description}
\item[(i)]
The leaves of the canonical foliation are obtained\index[terms]{foliation!canonical}
as $\pi({\goth O}\times \R^{>0})$, where $\goth O$ is a Reeb orbit on $S$.
\item[(ii)]
Moreover, the following are equivalent
\begin{description}
\item[(a)]
$M$ is quasi-regular 
\item[(b)]
 $S$ is a quasi-regular Sasakian manifold, and 
the action of $\phi$ on the space of Reeb orbits in 
$S$ has finite order.\index[terms]{manifold!Sasaki!quasi-regular}
\end{description}
\end{description}

{\bf Step 1. \ref{_qreg_Vaisman_via_Sasakian_Theorem_} (i):}
By \ref{halfstr}, the Lee form \index[terms]{form!Lee}lifted to the cone 
is $- \frac{dt}{t}$. Therefore, the Lee field\index[terms]{Lee field} is proportional
to $\frac d {dt}$, and the canonical foliation\index[terms]{foliation!canonical} is generated
by $\frac d {dt}$ and $I\left(\frac d {dt}\right)$, which
is proportional to the Reeb field by \ref{_Reeb_via_cone_Claim_}.

\hfill

{\bf Step 2: quasi-regularity of $M$ implies
  quasi-regularity of $S$ and the fi\-niteness condition on $\phi$.}

We identify $S\times \{1\}\subset C(S)$ with its image in
$M$, and denote it by $S \subset M$.
The leaves of $\Sigma$ intersect with $S$ in 
orbits of the Reeb field $I\theta^\sharp$ (\ref{_Reeb_via_cone_Claim_})
because $TS\cap \Sigma = \langle I\theta^\sharp\rangle$.
If the intersection of a leaf of $\Sigma$ and $S$ 
is compact, the Reeb orbits are compact.
Therefore, the quasi-regularity of $M$ implies 
quasi-regularity of $S$. By \ref{_Reeb_action_factorizes_Corollary_},
 the Reeb action factorizes through the circle.

The holomorphic diffeomorphism flow 
$V_\mu(s, t) = (s, \mu t)$
along the Lee field \index[terms]{Lee field}acts on $M$ by
automorphisms that preserve each leaf $\Sigma_x$ of the canonical
foliation. Since $M= C(S)/ (\phi \times \lambda)$,
the map $V_\lambda$ sends each Reeb orbit 
${\goth O}\subset \Sigma_x$ to $\phi({\goth O})$.
The same argument actually gives 
\begin{equation}\label{_intersection_leaf_Sigma_Equation_}
\Sigma_x \cap S= \bigcup_{i\in \Z} \phi^i({\goth O})
\end{equation}
where ${\goth O}$ is a single Reeb orbit that is  obtained
as an intersection of $\Sigma_x$ and $S \subset M$.
If $\Sigma_x$ is compact, this set is also compact,
hence contains finitely many orbits.

\hfill

{\bf Step 3: Proof of 
\ref{_qreg_Vaisman_via_Sasakian_Theorem_} (ii):}
Let $S/R$ be the space of Reeb orbits in $S$.
It remains to show that $M$ is quasi-regular
if $S$ is quasi-regular and the induced action 
$\phi\restrict{S/R}$ is finite.


Since $S$ is quasi-regular
(\ref{_Reeb_action_factorizes_Corollary_}), the $\R^2$-action generated by the Reeb and
the Lee vector fields factorizes through
$S^1 \times \R$.
As  $M= C(S)/ (\phi \times \lambda)$,
the subgroup $S^1 \times \lambda\Z$ fixes the image 
of $S\times \{t\}$ in $M$. 
Therefore, the leaf $\Sigma_x$  is compact if and only if
 $\Sigma_x$ intersects
an image of $S\times \{t\}$ in a union of finitely many circles.
However, this intersection is equal to 
$\bigcup_{i\in \Z} \phi^{i}({\goth O})$ by 
\eqref{_intersection_leaf_Sigma_Equation_}.

We obtain that 
the leaf $\Sigma_x$  is compact if and only if the set 
$\bigcup_{i\in \Z} \phi^{i}({\goth O})$ is finite.
This is equivalent to $\phi$ acting on the set
of Reeb orbits with finite order.
\endproof

 
\section{Regular Vaisman manifolds}\index[terms]{manifold!Vaisman}\index[terms]{manifold!Vaisman!regular}


The following claim gives an example of a regular
Vaisman manifold. Later, we will prove that all
regular Vaisman manifolds are obtained in this way.

\hfill

\claim
Consider the total space $\Tot^\circ(L)$ 
of all non-zero vectors in a
positive line bundle $L$ over a projective manifold\index[terms]{manifold!projective}
$X$. Consider the action of $\Z$ on $\Tot^\circ(L)$ mapping $v$ to
$\lambda v$, with $|\lambda|>1$. Then 
the quotient $\Tot^\circ(L)/\Z$ is equipped with
a natural Vaisman structure.

\hfill

\proof
Let $\tilde \omega:= dd^c l^2$ be the K\"ahler form\index[terms]{form!K\"ahler} on 
$\Tot^\circ(L)$, defined in \ref{_dd^c_log_l_Corollary_},
$V_t(x) = t x$ the homothety flow on $\Tot^\circ(L)$, 
and $X= \frac{d V_t}{dt}$ the corresponding vector field.
Clearly, $\Lie_X \tilde g=2\tilde g$, where $g$ is 
the K\"ahler metric on $\tilde M$. Moreover, the dual 1-form 
$X^\flat=d l$ is closed. Let 
$\nabla_{\tilde g}$ denote the Levi--Civita connection of $\tilde g$.
By \ref{homoth}, $\nabla_{\tilde g}(X)= \Id$,
and $X$ is the Euler field of the cone metric.\index[terms]{vector field!Euler}
Replacing $\tilde g$ by $g= \frac{\tilde g}{l^2}$,
we obtain the structure described in  \ref{homoth} (iii):
a vector field $X$ that satisfies 
$d(X^\flat)=0$,\footnote{We use $X^\flat$ to denote the 1-form
dual to $X$ with respect to the metric $g$.}
and $\Lie_X g= \lambda g$, with $\lambda=0$.
Then $\nabla_g X = \lambda \Id =0$ by \ref{homoth} (ii),
that is, $X$  is parallel\index[terms]{vector field!parallel} with respect to the
Levi--Civita connection associated with the
metric $g= \frac{\tilde g}{l^2}$. 
\endproof

\hfill

\remark Recall that the canonical
foliation on a compact Vaisman manifold 
is determined by its complex structure
(\ref{_Subva_Vaisman_Theorem_}).
In the next theorem, we give
a complex geometric characterization
of regular Vaisman manifolds.\index[terms]{manifold!Vaisman!regular}
The same theorem can be stated
and proven for quasi-regular Vaisman manifolds if\index[terms]{manifold!Vaisman!quasi-regular}
we replace ``manifolds'' with ``orbifolds''
(Section \ref{_quasireg_over_orbifold_Section_}).

\hfill

\theorem \label{_regular_total_space_Theorem_}
Let $M$ be a compact  Vaisman manifold.
Then $M$ is regular if and only if
it is biholomorphic to \index[terms]{manifold!Vaisman!regular}
a smooth elliptic fibration over 
a projective manifold $X$, \index[terms]{elliptic fibration}\index[terms]{manifold!projective}
obtained as a quotient of the total space $\Tot^\circ(L)$ 
of non-zero vectors in a
negative line bundle $L$ over $X$ by the action of $\Z$ mapping $v$ to
$\lambda v$, with $|\lambda|>1$.

\hfill

 {\bf Proof. Step 1: total spaces of elliptic fibrations
are regular.}
Let $M$ be a Vaisman manifold \index[terms]{manifold!Vaisman}obtained as a quotient
of $\Tot^\circ(L)$ by $\Z = \langle \lambda  \rangle$.
By \ref{_Subva_Vaisman_Theorem_}, all complex curves on $M$ are tangent to the
canonical foliation. Therefore, the leaves of the
canonical foliation are fibres of the elliptic fibration
$\pi:\; M \arrow X$, and the canonical foliation is
clearly regular.\index[terms]{foliation!canonical}

\hfill

{\bf Step 2: regularity implies that $M$ is a quotient
of a cone over a projective manifold.} 

Deforming the Vaisman metric as in
\ref{_Vaisman_defo_transve_Proposition_}, 
we can assume $M$ has LCK rank 1. \index[terms]{rank!LCK}
By \ref{halfstr}, 
$\Sigma$ is generated by  $\frac d {dt}$ and 
$I(\frac d {dt})$.

The Reeb foliation on a Sasakian manifold $S\subset M$ 
is obtained by intersection with the leaves of the canonical
foliation $\Sigma$. Therefore, the Reeb foliation
is quasi-regular when the corresponding Vaisman manifold\index[terms]{manifold!Vaisman}\index[terms]{manifold!Vaisman!quasi-regular}
is quasi-regular. Moreover, the holonomy of the Reeb
foliation is a subgroup of the holonomy of $\Sigma$.
By \ref{_regular_holonomy_finite_Proposition_},
the holonomy of a quasi-regular foliation is trivial if 
and only if the foliation is regular. Therefore, $S$ is 
a regular Sasakian manifold when $M$ is regular.

Let $X_1:= S/R$ be the space of Reeb orbits on $S$.
By construction, the leaf space of the
canonical foliation\index[terms]{foliation!canonical} is identified with
the quotient of $X_1/\phi$. Since the leaf
space $X_1/\phi$  is Hausdorff and $X_1$ is compact, $\phi$ acts on $X_1$
with finite order $d$.

A leaf $F$ of $\Sigma$ corresponds to an orbit of
$\C^*$-action on $C(S)$. It is mapped to 
$M$ by identifying the points $(s, t)$
and $(\phi(s), \lambda t)$. Therefore,
the intersection of $F$ with $S$
is the union 
\[ 
{\goth O} \cup \phi({\goth O}) \cup \phi^2({\goth O}) \cup ...
\cup \phi^{d-1}({\goth O}),
\]
where $\goth O$ is a Reeb orbit in $S$.
We need to show that $\phi$ acts freely on $X_1$.
Consider a point $x\in X_1$ where $\phi$ acts with 
order $k< d$, and let $F$ be the corresponding leaf of $\Sigma$.
The map $\phi^k:\; X_1 \arrow X_1$
takes a Reeb orbit ${\goth O} \subset S$
and sends it to another Reeb orbit that belongs
to the same leaf of $\Sigma$. Therefore,
$\phi^k$ generates the holonomy of $\Sigma$ 
near $F$. Since $\Sigma$ is regular, the holonomy of $\Sigma$ 
is trivial, $\phi^k=\Id$, and $k=d$.\index[terms]{foliation!regular}

We have shown that $S_1:=S/\langle \phi\rangle$ is a 
regular Sasakian manifold. Then\index[terms]{manifold!Sasaki!regular} \index[terms]{manifold!Sasaki} 
\ref{_quasireg_Sasakian_orbibundles_Theorem_} (iii)
implies that $S_1$ it is the space of unit 
vectors in a negative line bundle $L$ over 
$X:=S_1/\langle I\theta^\sharp\rangle$.

By construction, $X$ is the space of
leaves of $\Sigma$, and hence  $C(S)$ is a principal holomorphic 
$\C^*$-fibration over $X$.\index[terms]{bundle!principal} 
\endproof

\hfill

\corollary
Let $M$ be a Vaisman manifold obtained as
$M= C(S)/\lambda \times \phi$, where $S$ is regular,
and $\phi$ is a Sasakian automorphism of $S$ of 
constant finite order. Then $M$ is a regular
Vaisman manifold.\index[terms]{manifold!Vaisman!regular}
\endproof

\hfill

\example  The {\bf classical Hopf manifold} 
$\C^n\backslash 0 /(x\sim \lambda x)$
is  a regular Vaisman manifold, with $X=\C P^{n-1}$.
\index[terms]{manifold!Vaisman!regular}

\section[Quasi-regular Vaisman manifolds are or\-bi\-fold elliptic fibrations]{Quasi-regular Vaisman manifolds are\\ or\-bi\-fold elliptic fibrations}\index[terms]{manifold!Vaisman!quasi-regular}
\label{_quasireg_over_orbifold_Section_}\index[terms]{elliptic fibration}

If we work in the category of orbifolds, 
instead of manifolds, and speak  of quasi-regular
foliations, instead of the regular ones, the
analogue of \ref{_regular_total_space_Theorem_} is straightforward.

We defined quasi-regular Vaisman manifolds as ones
where the action induced by the Lee and anti-Lee fields\index[terms]{Lee field}\index[terms]{Lee field!anti-} (tangent to the\index[terms]{manifold!Vaisman!quasi-regular}
leaves of the canonical foliation) is factorized through
a compact torus. However,
by \ref{_Killing_flow_Vaisman_orbits_Theorem_},
a Vaisman manifold is quasi-regular if and only if
all leaves of the canonical foliation are compact
(\ref{_qr_Vaisman_defi_Remark_}).\index[terms]{foliation!canonical}

\hfill

\theorem \label{_Structure_of_quasi_regular_Vasman:Theorem_}
	Let $M$ be a Vaisman manifold,
$\Sigma$ its canonical foliation, and $M/\Sigma$\index[terms]{foliation!canonical}
the quotient space. Then the following statements are equivalent.
\begin{description} 
\item[(i)] $M$ is quasi-regular.
\item[(ii)] 
There exists a negative holomorphic orbifold
line bundle $L$ over $X$, such that
$M$ is biholomorphic to a $\Z$-quotient
of the space $\Tot^\circ(L)$ of non-zero vectors in $L$.
The leaves of the canonical foliation are compact,
and  their preimages in $\Tot^\circ(L)$ coincide
with the fibres of $L$.
\index[terms]{orbifold!projective}
\end{description}
Moreover, for any negative line bundle $L$
over a projective orbifold $X$, the quotient of
$\Tot^\circ(L)$ by the action of $\lambda \in \C$,
$|\lambda| >0$ is Vaisman.

\smallskip

\proof We prove the implication (i) $\Rightarrow$ (ii).
The leaf space $X=M/\Sigma$ of the canonical foliation\index[terms]{foliation!canonical}
is an orbifold\footnote{This also follows from
  \cite[Proposition 3.7]{_Molino_} because $\Sigma$ is a
  Riemannian foliation, but we prefer a direct proof.} as
proven in \ref{torus_quotient}.\index[terms]{foliation!Riemannian}
Since it is the  quotient of a complex space by a complex group
action, $X$ is a complex orbifold. 
Replacing the Vaisman metric with a
Vaisman metric of LCK rank 1 as in \index[terms]{rank!LCK}
\ref{_Vaisman_defo_transve_Proposition_}, we may assume
that $M$ admits a $\Z$-covering\index[terms]{cover!K\"ahler $\Z$-} $\tilde M=C(S)$
obtained as a cone over a compact Sasakian manifold
$S$. Since the Reeb orbits of $S$ are compact, 
$S$ is quasi-regular\index[terms]{manifold!Sasaki!quasi-regular} (\ref{_qreg_Vaisman_via_Sasakian_Theorem_}).
By \ref{_qreg_Vaisman_via_Sasakian_Theorem_} (iii) and
\ref{_quasireg_Sasakian_orbibundles_Theorem_} (ii)--(iii),
$\tilde M$ is the  total space
of a negative $\C^*$-bundle $L$ (in the orbifold sense).
  Now, \index[persons]{Baily, W. L.} Baily's \ref{baily} implies that $X$ is projective.\index[terms]{theorem!Baily}

The opposite implication is clear, because
all compact curves in a Vaisman manifold are
tangent to the canonical foliation (\ref{_Subva_Vaisman_Theorem_}),\index[terms]{manifold!Vaisman}
hence any elliptic fibration on a Vaisman manifold
coincides with the canonical foliation.\index[terms]{foliation!canonical} Therefore,
the leaves of the canonical foliation are compact, and
the manifold $M$ is quasi-regular
(\ref{_qr_Vaisman_defi_Remark_}).

Finally, \ref{_formula_15_19_besse_Theorem_}
implies that the function $l(v):= |v|^2$ is 
plurisubharmonic on $\Tot^\circ(L)$. The $\Z$-action\index[terms]{action!$\Z$-}
generated by $\lambda \in \C$,
$|\lambda| >0$  maps $l$ to $ |\lambda|^2 l$, and hence 
$(\Tot^\circ(L), dd^cl)/\langle \lambda\rangle$ is
locally conformally K\"ahler.
By \ref{kami_or}, it is Vaisman.
\endproof

\smallskip

\corollary \label{triangle}
Let $M$ be a compact quasi-regular Vaisman manifold. \index[terms]{manifold!Vaisman!quasi-regular}
Then there exists the\index[terms]{manifold!Vaisman}
following diagram of principal (orbi-)bundles. It is defined in the manifold
category when $M$ is regular, and in the orbifold\index[terms]{bundle!principal}
category for the quasi-regular case.
In this diagram, $S$ is a Sasakian orbifold, $X$ is a projective orbifold,
and the arrows correspond to the projections to 
the leaf spaces of transversally Riemannian 
foliations.
%
%
\begin{equation}\label{_t2_diagram_Vaisman_Equation_}
	\begin{minipage}{0.85\linewidth}
	\xymatrix{
		&M\ar[ld]_{S^1}\ar[rd]^{T^2}&\\
	S\ar[rr]_{S^1}&&X
	}	
\end{minipage}
\end{equation}
\index[terms]{orbifold}


\proof
The orbifold elliptic fibration $M \arrow X$ is 
constructed in \ref{_Structure_of_quasi_regular_Vasman:Theorem_}.
To construct an $S^1$-fibration over a Sasakian orbifold, we use
\ref{_qreg_Vaisman_via_Sasakian_Theorem_}. Here it is shown that
a $\Z$-covering \index[terms]{cover!K\"ahler $\Z$-}of $M$ is biholomorphic to a cone 
$C(S_1)=S_1\times \R^{>0}$ over 
a quasi-regular Sasakian manifold $S_1$\index[terms]{manifold!Sasaki!quasi-regular}, with the $\Z$-action\index[terms]{action!$\Z$-}
given by $(s, t)\arrow (\phi(s), \lambda t)$, 
where $\lambda>0$ is a constant, and $\phi$ is
a Sasakian automorphism of finite order. This
gives an $S^1$-fibration $M \arrow S_1/\langle \phi\rangle$.
We consider the quotient $S:=S_1/\langle \phi\rangle$
as a Sasakian orbifold. The leaf space of its Reeb
foliation is the projective orbifold $X$, and the
$T^2$-fibration $M \arrow X$ is factorized through $S$.
\endproof


\section{Density of  quasi-regular Vaisman manifolds}\index[terms]{manifold!Vaisman!quasi-regular}


We prove that every compact Vaisman manifold can be
deformed to a quasi-regular one, and hence it is enough to
understand the topology in the quasi-regular case. This is
in contrast with the fact that the LCK class is not stable
at small deformations, see \cite{bel} and Chapter
\ref{comp_surf}.

Recall that a conical K\"ahler manifold\index[terms]{manifold!K\"ahler!conical}
is a K\"ahler manifold locally isomorphic to a cone
over a Sasakian manifold (\ref{_conical_Kahler_Definition_}).
We start with a  closer look at the local structure of
conical K\"ahler manifolds.
 
\hfill

 \theorem \label{conloc}
 Let $(M,  \omega)$ be a conical K\"ahler manifold,
 and $X$ a vector field acting on $M$ by holomorphic, non-isometric 
 homotheties, such that $IX$ also acts by homotheties,
 and $e^{t X}$ is defined for any real $t$. Then:
 
 (i) $dd^c \phi=  2\omega$, where $\phi=|X|^2$.
 
 (ii) Let $S_X:= \phi^{-1}(1)$. Then $S_X$ is Sasakian,\index[terms]{manifold!Sasaki}
 and	$M$ is isometric to the Riemannian cone \index[terms]{cone!Riemannian}$C(S_X)$.
 
 (iii) $S_X$ is quasi-regular if and only if the action of $X$
 	integrates to a holomorphic $\C^*$-action.

\hfill

\proof The same argument was used in the proof of \ref{kami_or}.
Since both $X$ and $IX$ act by homotheties, there exists 
 a character $\chi:\; \langle X, IX\rangle\arrow \R$
 such that $\Lie_Z\omega=\chi(Z)\omega$. Replacing $X$ by
 	some linear combination of $X, IX$ if necessary, we 
 	may assume that $IX$ acts by isometries.
 Rescaling, we may assume that $\Lie_X g=2g$.
 
 Now define  $X^\flat:= g(X, \cdot)$ {(``the metric dual 
 	1-form'').} By Cartan formula, \index[terms]{Cartan formula} $dX^\flat= \Lie_{IX}\omega=0$ and
 $2 X^\flat = \Lie_X(X^\flat)=d\langle X, X^\flat\rangle=d|X|^2$.
 
 Since  $\Lie_X \omega= 2\omega$, we have  
 $2\omega = d(i_X \omega)= d(IX^\flat)=  dId|X|^2$, thus  
proving (i).

 Let $M\arrow S_X$ map $m$ to the intersection
$\left(e^{tX}\cdot m\right)\cap S_X$ 
of  $e^{tX}\cdot m$ (the orbit of $X$ through $m$) with  $S_X$. 
	By \ref{halfstr}, this gives a decomposition $M=S_X \times \R^{>0}$,
	compatible with the conical metric on $S_X \times \R^{>0}=C(S_X)$.

Finally, let $C$ be the group generated by 
$e^{tX}, e^{tIX}$. Clearly, $C=\R^{>0}\times \{e^{tIX}\}$.
The Reeb orbits on $S_X$  
are orbits of $e^{tIX}$, {hence they are compact if and
	only if $\{e^{tIX}\}$ is compact, equivalently,
	if and only if $C=\C^*$ (\ref{_Reeb_action_factorizes_Corollary_}).}
\endproof

\hfill

%
%
%

The key result in the proof of the density of
quasi-regular Vaisman manifolds among compact Vaisman ones
is as follows. \index[terms]{manifold!Vaisman!quasi-regular}

\hfill

\theorem \label{_approximate_on_cone_by_C^*_Theorem_}
Let $C(S)$ be a conical K\"ahler manifold, with $S$ compact, $h_t$
the corresponding homothety action, and $X$ its vector field.
Then there exists a vector field $X'$ arbitrarily close to $X$
acting on $C(S)$ by holomorphic homotheties, with $IX'$
also acting by homotheties, such that the action of $X'$ integrates
	to a $\C^*$-action on $C(S)$.

\hfill

\proof Fix some $\lambda\in \R^{>1}$, and let $M:= C(S)/\langle h_\lambda\rangle$
be the corresponding compact Vaisman manifold, on which the flow $h_t$ acts\index[terms]{manifold!Vaisman}
isometrically. 

Consider the Lie group $G\subset \Iso(M)$ 
obtained as the closure of $\{h_t\}$. Since
the isometry group is compact\footnote{\cite{_Myers_Steenrod_,_Kobayashi_Transformations_}; 
see also the exercises
  to this chapter.}, $G$ is compact;
it is commutative because it contains a dense
commutative subgroup. Therefore, $G$ is a compact torus.
Using the same argument as in \ref{conloc}, we choose
$X$ in such a way that $I(X)$ acts by isometries on 
$C(S)$, that is, it preserves the subsets
$S \times \{t\}\subset C(S)$. Denote by
$G_0\subset G$ the group of all isometries
of $M$ acting by isometries on the cone.
The manifold  $M$ is obtained as a quotient
of the cone $C(S)=S\times \R^{>0}$
by an action of $\Z$ generated
by $(s,t) \mapsto (s, \lambda t)$.
Therefore, $M$ is naturally identified
as a Riemannian manifold with $S \times S^1$. 
The metric on $M$ and $C(S)$ are conformally
equivalent, with the conformal factor of $t^2$.
Since $G_0$ is a subgroup of
$G$ acting by isometries on the cone $C(S)$ 
and on the product $S \times S^1$, it preserves
the conformal factor $t^2$, defined on the
$\Z$-covering\index[terms]{cover!K\"ahler $\Z$-} of  $S \times S^1$. Therefore,
$G_0$ maps a point $(s, t)$ to $(s', t)$
(with the same $t$).

We obtained that $G_0$ is a subgroup
of $\Iso(S\times S^1)$ acting
trivially on the $S^1$-factor.
In particular, $G_0$ is a closed subtorus of $G$.

Each vector field
$X'\in \Lie(G)$ acts on $M$ by holomorphic isometries,
hence it acts on $C(S)$ by homotheties; non-isometrically when
$X'$ is sufficiently close to $X$. 

Since $G_0$ is a torus, its Lie algebra
$\Lie(G_0)$ is equipped with a rational lattice.\index[terms]{lattice!rational}
All closed subgroups of $G$ are Lie groups by 
\index[persons]{Cartan, E.} the Cartan theorem (see for example  \cite{_Knapp_}).\index[terms]{theorem!Cartan}
Closed, connected subgroups of $G_0$
are in bijective correspondence with 
rational Lie subalgebras in $\Lie(G_0)$.

Choosing $Y\in \Lie (G_0)$ rational and sufficiently close
to $I(X)$, we obtain a flow of isometries of $C(S)$ which
acts on $S$ and factorizes  through $S^1$.
Then the Lie algebra generated by $Y, I(Y)$
integrates to a holomorphic $\C^*$-action
that approximates the action by $I(X), X$.
 \index[terms]{isometry}
\endproof

\hfill

\remark 
By \ref{conloc} this gives a new cone structure on $C(S)$.

\hfill

The proof of \ref{_approximate_on_cone_by_C^*_Theorem_}
implies a statement that is  used quite often.

\hfill

\corollary\label{_Vaisman_Lee_action_contains_monodromy_Corollary_}
Let $(M, \theta, \omega)$ be a Vaisman manifold of LCK
rank 1, $\theta^\sharp$ its Lee field\index[terms]{Lee field}, and \index[terms]{manifold!Vaisman}
$\tilde M$ the K\"ahler $\Z$-cover\index[terms]{cover!K\"ahler $\Z$-} of $M$.
Denote by $\gamma$ the generator of the
$\Z$-action\index[terms]{action!$\Z$-} on $\tilde M$.
Denote by $G\subset \Iso(M)$ the closure
of the one-parametric subgroup generated by
$e^{t\theta^\sharp}$, and let $\tilde G$
be the connected group obtained as the
lift of the $G$-action to  $\tilde M$.
Then $\tilde G$ contains $\gamma^k$, for some
$k\in \Z^{>0}$. 

\hfill

\proof
The surjective map $\tau:\; \tilde G \arrow G$ is 
an infinite cover, because $\tilde G$ acts
on $\tilde M$ conformally. Clearly, any
element in $\ker \tau$ belongs to the monodromy\index[terms]{group!monodromy}
group. Therefore, $\gamma^k \in \ker \tau$.
\endproof

\hfill

\remark
A stronger statement was made in 
\cite[Theorem 2.1]{ov_mz}, where we claimed that 
$\tilde G$ contains the deck transform group of
the covering $\tilde M \arrow M$.
This would imply, in particular,
the existence of the logarithm of the
generator of the deck transform group 
if it is cyclic, that is  false by
\ref{_only_for_gamma^k_log_exists_Example_}. 
Therefore, \cite[Theorem 2.1]{ov_mz}
is false, and \ref{_Vaisman_Lee_action_contains_monodromy_Corollary_}
is the strongest possible statement.

\hfill

\corollary 
Let $C(S)$ be a cone over a compact Sasakian
manifold $S$. {Then $C(S)$ is holomorphically 
	isometric to the total space of non-zero sections 
	of a negative line bundle over a 
	projective orbifold.}\index[terms]{orbifold!projective}

\proof Follows from \ref{_approximate_on_cone_by_C^*_Theorem_}
and \ref{_quasireg_Sasakian_orbibundles_Theorem_} (iii). 
\endproof

\hfill

\theorem {(\cite{ov_imm_vai})}\label{defovai}
Any compact Vaisman manifold $(M,I)$ admits a complex deformation\index[terms]{manifold!Vaisman}
$(M,I')$, that is  Vaisman and quasi-regular.\index[terms]{manifold!Vaisman!quasi-regular} Moreover, $I'$ can be chosen
arbitrarily close to $I$. 

\hfill

\proof
We can always approximate the Vaisman metric on $M$ by 
another Vaisman metric of LCK rank 1 (\ref{_rk_1_approximated_Theorem_}).
Therefore, it will suffice to prove \ref{defovai}
when $M$ has LCK rank 1.\index[terms]{rank!LCK}

In this case, $M = C(S)/\langle \gamma \rangle$,
where $\gamma$ is the generator of the $\Z$-action\index[terms]{action!$\Z$-}, acting
on $C(S)$ by holomorphic homotheties.

The idea of the proof is straightforward. By
\ref{_approximate_on_cone_by_C^*_Theorem_}, 
we can choose a $\C^*$-action by homotheties on $C(S)$.
We approximate $\gamma$ by a homothety $\gamma'$ which
preserves the orbits of a $\C^*$-action on $C(S)$ and acts 
with finite order on the corresponding leaf space.
Then the quotient $C(S)/\langle \gamma' \rangle$
is a quasi-regular Vaisman manifold\index[terms]{manifold!Vaisman!quasi-regular} 
(\ref{_Structure_of_quasi_regular_Vasman:Theorem_})
which  can be chosen arbitrarily close to $M$.

\hfill

{\bf Step 1:} This is an argument we have already
used (and will use more) in this book 
(\ref{_Vaisman_Lee_action_contains_monodromy_Corollary_}).

Let $G$ be the closure of the group generated
by the Lee and anti-Lee field\index[terms]{Lee field} in $\Iso(M)$,\index[terms]{Lee field!anti-}
and $\tilde G$ the group of conformal 
transforms of $C(S)$ generated by the
vector fields lifted from the vector fields
in $\Lie(G)\subset H^0(TM)$. 
Since $M = C(S)/\Z$, the group 
$\tilde G$ is mapped to $G$ with a
kernel $\Z$. Let $G_0\subset \tilde G$
be the subgroup acting on $C(S)$
by isometries. Clearly, it is
projected to its image in $G$ 
bijectively. 

Let $\chi:\ \tilde G  \arrow \R$
map $\alpha \in \tilde G$ to
its conformal factor such that
$\alpha^*(\tilde \omega) = \chi(\alpha) \tilde \omega$,
where $\tilde \omega$ is the K\"ahler form\index[terms]{form!K\"ahler} on $C(S)$.
Then $G_0= \ker\chi$. Therefore, $G_0$
is a codimension 1 subtorus in $G$.
Since $G_0$ preserves the Vaisman
metric on $M$ and the K\"ahler metric
on the torus, it can be obtained as 
the closure of the group of conformal 
automorphisms of $C(S)$ generated by the lifts 
of the Lee and the anti-Lee fields.\index[terms]{Lee field}\index[terms]{Lee field!anti-}
Therefore, $G_0$ commutes with the action
of $\gamma$ on $C(S)$.

\hfill

{\bf Step 2:} The Structure Theorem
(\ref{str_vai}) implies that 
 $\gamma$ acts as $(s, t) \mapsto (\phi(s), \lambda t)$,
where $\lambda>1$ and $\phi$ is a Sasakian automorphism
of $S$. The group $\Aut(S)$ of Sasakian automorphisms 
is a closed subgroup of the group of isometries of $S$,
hence it is compact.

On the conical K\"ahler manifold $C(S)$, we 
replace the Euler vector field $X$ by $X'$ in such\index[terms]{vector field!Euler}
a way that the $\C$-action generated by 
$X'$ and $I(X')$ factorizes through $\C^*$
(\ref{_approximate_on_cone_by_C^*_Theorem_}).
Moreover, $X'$ is chosen within $\Lie(G_0)$.
As follows from Step 1, $\gamma$ preserves $X'$ and the
$\C^*$-action $\rho$ generated by $X'$ and $I(X')$.

Since $\Aut(S)$ is compact, we can
always approximate $\phi$ by
$\phi'$ that satisfies $(\phi')^n=\Id$.
Let $\gamma'(s, t) := (\phi'(s), \lambda t)$.
This map commutes with $\rho$
and $(\phi')^n$ acts trivially on $S$, and hence 
$(\gamma')^n$ fixes each orbit of $\rho$.
Then $C(S)/\langle \gamma'\rangle$ 
is holomorphically fibred over the leaf space of
$\rho$ with the fibre $\C^*/\Z= T^2$.
Any elliptic curve on a Vaisman 
manifold is a leaf of the canonical
foliation (\ref{_Subva_Vaisman_Theorem_}), and hence 
the leaves of the canonical foliations are compact.
Now, \ref{_qr_Vaisman_defi_Remark_} \index[terms]{foliation!canonical}
implies that $C(S)/\langle \gamma'\rangle$
is quasi-regular.
\endproof

\hfill

\remark
The proof of \ref{defovai}
 gives in fact a stronger result than stated.
We choose a cone $C(S)$ such that
$M=C(S)/\langle \gamma\rangle$,
then approximate $\gamma$ by $\gamma'$
in the group of holomorphic homotheties
of $C(S)$ such that $M'=C(S)/\langle \gamma'\rangle$
is quasi-regular. Of course, $M'$ is a complex
deformation of $M$, but not all deformations
of $M$ can be realized this way.

\hfill

We also obtain the density of quasi-regular Sasakian
manifolds among compact Sasakian manifolds\index[terms]{manifold!Sasaki!quasi-regular}
(previously proven differently in \cite{ru}, see \ref{ruki}):

\hfill

\theorem \label{_defo_qr_sas_}
Any compact
Sasakian manifold $S$ admits a deformation\index[terms]{manifold!Sasaki}
$S'$ that is  quasi-regular. Moreover, the Sasakian 
structure on $S'$ can be chosen
arbitrarily close to $S$. 
\endproof

\hfill

For a more explicit form of this statement, see 
 \ref{qr_reeb}.

\hfill

\remark
Now the diagram from \ref{triangle}
can be constructed for any Vaisman manifold,
not necessarily quasi-regular.

\section{Immersion theorem for Vaisman manifolds}\label{immersion_vaisman}\index[terms]{manifold!Vaisman}

One can use \ref{defovai} to prove an analogue of the Kodaira
embedding theorem. However, in this case, we only obtain a
holomorphic {\bf immersion}. For embedding, another
argument is used (Chapter \ref{_Embedding_for_LCK_with_potential:Chapter_}).
Instead of the projective
space, the model space will now be the Hopf
manifold.\index[terms]{manifold!Hopf}
The first step is as follows.

\hfill

\theorem \label{_Immersion_for_cones_Theorem_}
Let $C(S)$ be the conical K\"ahler manifold of a compact Sasakian manifold $S$.
There exists a holomorphic immersion $C(S)\arrow C(S^{2n-1})$ that is  
	equivariant under homothety, where $C(S^{2n-1})=\C^n \backslash 0$
is the standard (flat) cone.

\hfill

\proof
By \ref{_approximate_on_cone_by_C^*_Theorem_},
$C(S)$ admits a biholomorphic isometry to
$C(S')$, where $S'$ is quasi-regular Sasakian manifold.
Therefore, we can assume that $S$ is quasi-regular from the start.\index[terms]{manifold!Sasaki!quasi-regular}

When $S$ is quasi-regular, $C(S)$ is  the 
space of non-zero vectors in the total space of a 
negative line bundle $L$ over a projective orbifold \index[terms]{orbifold!projective}
$X$ (\ref{_quasireg_Sasakian_orbibundles_Theorem_} (iii)).
However, we are free to replace $L$ with $L^*$ 
(\ref{_C^*_bundle_isom_to_dual_Claim_}),
hence we may assume that $L$ is positive.
By \index[persons]{Baily, W. L.} Baily's \ref{baily}, $L^N$ is very ample,\index[terms]{theorem!Baily}\index[terms]{bundle!line!ample}
and there exists an embedding $X \stackrel j \hookrightarrow \C P^{n-1}$
such that $L^N=j^*(\calo(1))$. Here we consider $X$ as a 
complex variety; the embedding $j$ does not respect or
care about the orbifold structure. However,
we may always assume that 
$H^0(X, L^N)=H^0(\C P^{n-1},\calo(1))$.

Consider the holomorphic map
$\psi:\; C(S) \arrow \Tot(L^N)$ mapping $v$ to $v^N$.
It is an $N$-sheeted covering. 
Now, define $\Psi:\; C(S)\arrow C(S^{2n-1})=\Tot^\circ(\calo(1))$
as $\Psi(v):= j(\psi(v))$.  Since $\psi$ is \'etale
	and $j$ is an embedding, $\Psi$ is an immersion.

We obtained an embedding of $C(S)$ to the space 
$\Tot^\circ(L^N)$ of non-zero vectors in $L^N$,
embedded to $\Tot^\circ(\calo(1))$. The latter space
can be identified with $H^0(X, L^N)\backslash 0$,
because $H^0(X, L^N)=H^0(\C P^{n-1}, \calo(1))$.
Hence, the Sasakian sphere $S^{2n-1}$ is a sphere in the
vector space $H^0(X, L^N)$, and $\Psi$ acts
as 
\begin{equation}\label{_Sasakian_cone_embedding_Equation_}
\Psi:\; C(S)\arrow H^0(X, L^N)
\end{equation}
\endproof

From the same argument as used in the proof of
this theorem, we also obtain the following
useful corollary.

\hfill

\corollary\label{_cone_Sasakian_total_space_C^*_Corollary_}
Let $S$ be a compact Sasakian manifold. Then 
there exists a projective orbifold $X$ and a positive 
holomorphic Hermitian line\index[terms]{orbifold!projective}
bundle $L$ over $X$ such that the
conical K\"ahler manifold $C(S)$ admits a holomorphic 
isometry with $\Tot^\circ(L)$.\footnote{As elsewhere, $\Tot^\circ(L)$
denotes the space of non-zero vectors in the total space 
$\Tot(L)$.} \endproof

\hfill

\definition \label{_linear_Hopf_Definition_}
Recall that {\bf a linear Hopf manifold} is a quotient of
$\C^n \backslash 0$ by a linear automorphism with all
eigenvalues $|\alpha_i|<1$ (see also the Exercises at 
the end of this chapter).\index[terms]{manifold!Hopf!linear}

\hfill

\theorem \label{_immersions_Theorem_}
 { (\cite{ov_imm_vai})}
Let $M$ be a compact quasi-regular Vaisman manifold. {Then
	$M$ admits a holomorphic immersion into a linear Hopf manifold.}\index[terms]{manifold!Vaisman!quasi-regular}

\hfill

\proof 
Let $C(S)$ be the conical K\"ahler 
cover of $M$, and $\gamma:\ \Z \arrow \Aut(C(S))$
be the homothety action. Consider  the immersion 
$\Psi:\ C(S)\arrow C(S^{2n-1})$ constructed in
\ref{_Immersion_for_cones_Theorem_}. 

By \ref{_cone_Sasakian_total_space_C^*_Corollary_},
$C(S)$ admits a holomorphic
isometry to $\Tot^\circ(L)$, where $L$ is a 
negative line bundle over a K\"ahler orbifold $X$.

We consider $\Tot^\circ(L)$ as a principal $\C^*$-bundle.%
\footnote{To be more precise, we should call it a principal
$\C^*$-orbi-bundle; however, later we invoke \index[persons]{Baily, W. L.} Baily's\index[terms]{bundle!principal} theorem\index[terms]{theorem!Baily}
which ignores the orbifold structure on $X$ and $L$.}
By \ref{defovai}, Step 2, the isomorphism $C(S)=\Tot^\circ(L)$
can be chosen in such a way that $\gamma$ commutes
with this $\C^*$-action.

Since the action of $\gamma$ on $\Tot^\circ(L)$
commutes with the $\C^*$-action, the corresponding
line bundle $L$ over $X$ is $\gamma$-equivariant.
Therefore,
$\gamma$ actually induces a linear automorphism $\Gamma$ of the vector space
$\C^n=H^0(L^N)$. Since $\gamma$ acts as a contraction on
$\Tot(L)$, it uniformly decreases the norm on $H^0(L^N)$.
Therefore,  the eigenvalues of $\Gamma$ are all
$|\alpha_i|<1$. 
This gives a commutative square with the vertical
lines generated by metric contractions.

\begin{equation*}
\begin{CD}
C(S) @>{\Psi}>> C(S^{2n-1}) \\
@VV{/\gamma}V  @VV{/\Gamma}V              \\
M@>>> C(S^{2n-1})/\langle\Gamma\rangle
\end{CD}
\end{equation*}
By \eqref{_Sasakian_cone_embedding_Equation_}, $C(S^{2n-1})$ is identified with 
$H^0(L^N)\backslash 0$, and hence  this diagram
is equivalent to
\begin{equation*}
\begin{CD}
C(S) @>{\Psi}>> H^0(L^N)\backslash 0 \\
@VV{/\gamma}V  @VV{/\Gamma}V              \\
M@>>> \bigg(H^0(L^N)\backslash 0\bigg)\bigg/\langle\Gamma\rangle
\end{CD}
\end{equation*}
with the bottom arrow being the holomorphic immersion.
This proves \ref{_immersions_Theorem_}. 
\endproof

\hfill

\remark 
In fact, as we shall see in Chapter \ref{lckpotchapter},
for each compact Vaisman manifold\index[terms]{manifold!Vaisman}
{there exists an embedding
	into a linear Hopf manifold.} This will also solve
the embedding problem for compact Sasakian manifolds.

\section{Notes}

\begin{enumerate}
\item The compactness of the leaves of the canonical foliation\index[terms]{foliation!canonical} played a central role in this chapter. When $\Sigma$ is not quasi-regular, a natural question is to determine how many compact leaves it can have. The first result in this direction belongs to K. \index[persons]{Tsukada, K.} Tsukada who  proved  in \cite{tsu} that on compact
         Vaisman manifolds, $\Sigma$ must have at least
         one compact leaf. This can also be seen 
by the following chain of arguments. Any \index[terms]{manifold!Vaisman}
compact Vaisman manifold is a submanifold in
a Vaisman Hopf manifold, \ref{embedding}. Any submanifold
in a Hopf manifold contains at least two elliptic curves,
 \cite{kato2}. In Chapter \ref{_Elliptic_curves:Chapter_}, we show that any
Vaisman manifold of complex dimension $n$ contains
at least $n$ elliptic curves.
 \index[persons]{Tsukada, K.} Tsukada  used the existence of a compact leaf of $\Sigma$ to show that the product of two compact Vaisman manifolds cannot bear any LCK metric.\index[terms]{metric!LCK}

\item On diagonal Hopf surfaces $H_{\al,\be}$, the\index[terms]{surface!Hopf!diagonal}  quasi-regularity of the canonical foliation\index[terms]{foliation!canonical}    is equivalent to the Hopf surface being elliptic and this, in turn, is equivalent to the existence of $n,m\in\NN^*$ such that $\al^m=\be^n$ (Kodaira; see also  \cite{parton}). 

\item For compact regular Vaisman manifolds, the commutative diagram \eqref{_t2_diagram_Vaisman_Equation_} was exhibited in \index[terms]{manifold!Vaisman}\index[terms]{manifold!Vaisman!regular}
\cite{va_gd}; 
see also \cite{do}. In this case, 
the \index[persons]{Leray, J.} Leray spectral sequence associated with  the
fibration 
can be used to obtain 
relations between the Betti
numbers of $M$ and its K\"ahler base,
\cite{va_gd}. \index[terms]{Betti numbers!rational}
\end{enumerate}

\section{Exercises}

 \begin{enumerate}[label=\textbf{\thechapter.\arabic*}.,ref=\thechapter.\arabic{enumi}]	

\item \label{_compact_isom_group_1_exercise_}
Let $M$ be a metric space,\index[terms]{space!metric} and ${\goth C}$
the set of compact subsets of $M$. Define {\bf the
  Hausdorff distance} on ${\goth C}$ as in 
\eqref{_Hausdorff_distance_Equation_}.
\begin{enumerate}
\item
 Prove that
${\goth C}$ is a complete metric space.
\item Prove that ${\goth C}$ is compact if $M$
is compact.
\end{enumerate}

\item \label{_compact_isom_group_2_exercise_} Let $M$ be a compact metric space,
and $\Iso(M)$ the group of its isometries.
Define the metric on $\Iso(M)$ by
setting $d(f_1, f_2) := d_H(\Gamma_{f_1}, \Gamma_{f_2})$,
where $d_H$ denotes the Hausdorff distance in the
set of compact subsets of $M\times M$, and 
$\Gamma_{f_i}\subset M\times M$
denotes the graph of $f_i$. Prove that
$\Iso(M)$ is compact.

{\em Hint:} Use the previous exercise.

\item Prove that any isometry of a Riemannian
manifold, considered a metric space,\index[terms]{space!metric} preserves the
geodesics.

\item Prove that any isometry of a Riemannian
manifold is smooth and preserves the Riemannian metric
tensor $g\in \Sym^2(T^*M)$.

{\em Hint:} Use the previous exercise.

\item Let $M$ be a {\bf classical Hopf surface},\index[terms]{surface!Hopf!classical}
that is, the quotient $\C^2 \backslash 0$ by $x \mapsto \lambda x$,
where $\lambda$ is a complex number, $|\lambda| >1$.

\begin{enumerate}
\item Construct a holomorphic action of a finite group $G$ on $M$
such that the quotient is a smooth complex surface.
Prove that $M/G$ is a quasi-regular Vaisman manifold,
but not a regular one. \index[terms]{manifold!Vaisman!quasi-regular}

\item Let $X$ be the leaf space of the canonical
foliation on $M/G$. Prove that $X$ is $\C P^1$
with a non-trivial orbifold structure.

\item
Let $G$ be a finite group acting freely
on a Hopf surface. The quotient is called
 {\bf a secondary Hopf surface}. \index[terms]{surface!Hopf!secondary}
Construct a free action of a finite cyclic
group on a classical Hopf surface\index[terms]{surface!Hopf!classical}
and find the orbifold fundamental group\index[terms]{fundamental group}
of the quotient (see Exercise
\ref{_orbifold_funda_Exercise_}
for the definition of the orbifold fundamental group).

\item Recall that a {\bf dihedral group} is the group of
  symmetries of a regular plane polygon. It is generated
  by rotations and reflections.  Let $\epsilon$ be a root
of unity. Prove that 
$(x, y) \arrow (\epsilon x, y)$, $(x, y) \arrow (y, x)$
generates an action of a dihedral group
on a classical Hopf surface, and
determine the orbifold points on the corresponding
leaf space $X$. 

\item Construct an action of the alternating\index[terms]{surface!Hopf!classical}
group $A_4$ on a classical Hopf surface $M$, and
determine the orbifold points on the leaf space $X$.

\end{enumerate}

\item
Construct  a quasi-regular (but not regular) Sasakian
structure on $S^3$.\index[terms]{manifold!Sasaki!quasi-regular}

{\em Hint:} Use the previous exercise.

\item \label{_reg_Sasakian_pi_1_topology_Exercise_}
Let $S$ be a regular Sasakian manifold, and $X=S/\xi$
the corresponding K\"ahler manifold. 
\begin{enumerate}
\item
Construct an
exact sequence of homology groups:
\[
H_2(X)\stackrel {\delta}\arrow  \Z \arrow H_1(S) \arrow H_1(X) \arrow 0.
\]
Prove that there exists an integer K\"ahler class
$[\omega] \in H^2(X)$ such that\index[terms]{class!K\"ahler}
$\delta(x)= \langle[\omega], x\rangle$,
where $\langle\cdot, \cdot\rangle$
denotes the pairing between $H_2$ and $H^2$.\index[terms]{pairing}

{\em Hint:} Interpret $\delta$ as the first  Chern
class of the ample line bundle associated with the Reeb\index[terms]{bundle!line!ample}
foliation on $S$.
\item
Construct an
exact sequence of homotopy groups:
\[
\pi_2(X)\stackrel {\delta'}\arrow  \Z \arrow \pi_1(S) \arrow \pi_1(X) \arrow 0.
\]
Let $[\omega] \in H^2(X)$ be $c_1(L)$, where 
$L$ is the line bundle on $X$ associated with $S$ as in 
\ref{_quasireg_Sasakian_orbibundles_Theorem_}.
Prove that $\delta'$ is non-zero if and only if
$\langle A, [\omega]\rangle\neq 0$, where $A$ denotes  
the image of the Hurewicz map $\pi_2(X)\arrow H_2(X)$.
\end{enumerate}

 \item Let $S$ be a regular Sasakian manifold, and $X=S/\xi$. Suppose that\index[terms]{manifold!Sasaki!regular}
	$\pi_1(X)=0$. Prove that $\pi_1(S)$
	is a finite cyclic group.

{\em Hint:} Use 
Exercise \ref{_reg_Sasakian_pi_1_topology_Exercise_}.

\item Let $G$ be a finite cyclic group acting on a projective
manifold $X$. Find a Sasakian manifold $S$ such that\index[terms]{manifold!projective}
$S/\xi=X/G$.

\item Let $S$ be a regular Sasakian manifold such that $S/\xi=T^2$
(compact complex torus of dimension 1). Prove that
$\pi_1(S)$ is non-abelian.

{\em Hint:} Use 
Exercise \ref{_reg_Sasakian_pi_1_topology_Exercise_}.

\item More generally, let $S$ be a regular Sasakian manifold,
and $X:=S/\xi$. Assume that $\pi_1(X)=\Z^2$, $\pi_2(X)=0$.
Prove that $\pi_1(S)$ is non-abelian.

\item Let $S$ be a Sasakian manifold, $H$ the group of
Sasakian automorphisms, and $H_\C$ the corresponding
complex Lie group acting on $C(S)$ by holomorphic
automorphisms. Prove that $\dim_\R H=\dim_\C H_\C$.

{\em Hint:} Use Exercise \ref{_Killing_holomorphic_is_parallel_Exercise_}.

 	\item Define a {\bf contraction} of a manifold $M$ to a point
 	$x\in M$ as a diffeomorphism $\phi$ such that for any 
 	open subset $U\ni x$ with compact\index[terms]{contraction}
 	closure there exists $N>0$ such that for all
 	$n>N$, the map $\phi^n$ maps $U$ to a compact subset $K\subset U$ (for the local expression of contractions, see \cite{_Sternberg_contraction_}).
 	Let $\phi:\ M \arrow M$ be a contraction to $x$.
 	Prove that $\lim_n \phi^n(m)=x$, for each $m\in M$.
 	
 	\item Let $\phi:\ \C^n \arrow \C^n$ be a polynomial map
 	preserving $0$, and with $D\phi\restrict 0$ invertible
 	with all eigenvalues $|\alpha_i|<1$.
 	\begin{enumerate}
 		\item Prove that $\phi$ defines a contraction on an open
 		ball $B_\epsilon(0)$. 
 		\item Consider the equivalence
 		relation on $B_\epsilon(0)\backslash 0$ generated by $x\sim\phi(x)$.
 		Prove that $(B_\epsilon(0)\backslash 0)/\sim$ is a complex manifold
 		diffeomorphic to $S^1\times S^{2n-1}$.
 	\end{enumerate}
 	 The manifold $(B_\epsilon(0)\backslash 0)/\sim$ defined above is called
 {\bf a Hopf manifold}. When $\phi$ is linear, it is called
 {\bf a linear Hopf manifold}. In complex dimension 2, a Hopf manifold 
 is also called {\bf a Hopf surface}, or {\bf a primary Hopf surface}.\index[terms]{surface!Hopf!primary}\index[terms]{surface!Hopf!secondary} 
 A {\bf secondary Hopf surface} is the quotient of a Hopf surface
 by a finite group freely acting
 on\index[terms]{group!finite} it.\index[terms]{surface!Hopf!primary, 
   secondary}\footnote{This definition of a Hopf
 manifold, due to \index[persons]{Kodaira, K.} Kodaira,  \cite[p. 694]{_Kodaira_Structure_II_},
is the most general one; in this book, we almost never invoke
it, using the linear Hopf manifolds instead.}
 
 	\item Find secondary Hopf surfaces with $\pi_1(H)=\Z
 	\oplus \Z/n\Z$, for each $n\in \Z^{>0}$
 	
 	\item Find a secondary Hopf surface with $\pi_1(H)=\Z
 	\oplus S_3$, where $S_3$ is a symmetric group.\index[terms]{surface!Hopf!secondary}
 	
  
 	

\item
Let $M$ be a Vaisman manifold, $T\Sigma\subset TM$
the tangent bundle of the canonical foliation.\index[terms]{foliation!canonical}
Prove that the bundle $T\Sigma$ is\index[terms]{manifold!Vaisman}
trivial.

\item
Let $M$ be a quasi-regular Vaisman manifold,\index[terms]{manifold!Vaisman!quasi-regular}
and $\sigma:\; M \arrow X$ the projection to the
leaf space of the canonical foliation, considered
as a K\"ahler orbifold. Prove that
$K_M \cong \sigma^*K_X$.

{\em Hint:} Use the previous exercise.

 \end{enumerate}


 \chapter{LCK manifolds with potential}\index[terms]{manifold!LCK!with potential}\label{lckpotchapter}

{\setlength\epigraphwidth{0.5\linewidth}
\epigraph{\it Come slowly -- Eden!
	
	Lips unused to Thee --
	
	Bashful -- sip thy Jessamines --
	
	As the fainting Bee --

\smallskip
	
	Reaching late his flower,
	
	Round her chamber hums --
	
	Counts his nectars --
	
	Enters -- and is lost in Balms.}{\sc \scriptsize Emily Dickinson, \ \ Come slowly -- Eden!}}


 \section{Introduction}\label{_LCK_pot_Intro_Section_}

Now we pass to the notion that will  be  central to 
the rest of this book: that of an LCK manifold with
potential. The class of LCK manifolds with potential
is intermediate between the Vaisman manifolds and
the general LCK manifolds. \index[terms]{manifold!Vaisman}

We define an LCK manifold with potential\index[terms]{manifold!LCK!with potential} by
requesting that its K\"ahler cover
has a positive global K\"ahler potential, transforming
by the same automorphic law as the K\"ahler metric
(\ref{_LCK_potential_Definition_}). 

There are two flavours of this notion,
which correspond to LCK rank 1 and LCK rank $>1$\index[terms]{rank!LCK}.
Lets us call an LCK potential\index[terms]{potential!LCK} on the K\"ahler cover
$\tilde M$ {\bf proper} if it is proper as a\index[terms]{potential!LCK!proper}
map $\tilde M \arrow \R^{>0}$. By \ref{proper_equiv_Z},
\ref{_proper_potential_LCK_rank_1_Remark_},\index[terms]{rank!LCK}
an LCK potential on $(M, \omega, \theta)$
is proper if and only if $M$ has LCK rank 1,
or, equivalently, when the cohomology class of
the Lee form\index[terms]{form!Lee} $\theta$ 
is proportional to rational.

By \ref{defor_improper_to_proper}, 
LCK metrics with potential\index[terms]{metric!LCK!with potential}
on a given compact complex manifold can be
approximated by LCK metrics with proper potential.\index[terms]{potential!LCK!proper}
This allows us to assume the properness of potential
in most applications.

\index[persons]{Kodaira, K.} Kodaira and \index[persons]{Spencer, D. C.} Spencer have proven that a small
deformation of a compact K\"ahler manifold is
again K\"ahler. The proof, which takes most
of the paper \cite{_Kod-Spen-AnnMath-1960_},
can be greatly simplified if we use the 
semi-continuity of the relevant terms
of the Fr\"olicher spectral sequence
(\cite[Theorem 2.3]{_Bell_Narasimhan_}; see also Section \ref{_KS_stability_Section_}).

Unfortunately, LCK
metrics\index[terms]{metric!LCK} do not enjoy this property. Indeed, \index[persons]{Belgun, F. A.} Belgun has 
shown that a particular type of Inoue surfaces cannot
bear any LCK metric, while being a small complex
deformation of another type of Inoue surfaces\index[terms]{surface!Inoue} admitting LCK
metrics (see Chapter \ref{comp_surf}). It is then natural
to look for some subclass of LCK structures\index[terms]{structure!LCK} that is 
stable under small deformations. 

This is how we arrived at
defining the LCK manifolds with potential\index[terms]{manifold!LCK!with potential} in
\cite{ov_lckpot} (the preprint was posted on arXiv in
2004): the class of LCK manifolds with potential
is stable under small deformations.

All Vaisman manifolds are LCK with potential\index[terms]{manifold!LCK!with potential}, but the
inclusion is strict, even for Hopf\index[terms]{manifold!Vaisman} surfaces.\index[terms]{surface!Hopf}

We define ``a linear Hopf manifold'' \index[terms]{manifold!Hopf!linear}as 
a quotient of $(\C^n\backslash 0)$ by a 
linear operator $A$ with all eigenvalues $|a_i|>1$.
A linear Hopf manifold is called ``a diagonal Hopf
manifold'' (\ref{_diagonal_Hopf_Definition_})
if $A$ can be diagonalized (such operators
are also called ``semisimple''). 

The diagonal Hopf manifolds 
form a dense subset\index[terms]{manifold!Hopf!diagonal} 
in the deformation space of linear Hopf manifolds,
because the linear operators with distinct eigenvalues
are dense in the space of all linear operators,
and the linear operators with distinct eigenvalues
are diagonalizable. Moreover, for each non-semisimple
operator $A$, there exists a sequence $\{g_i\}$ of
linear operators such that $\lim_i g_i A g_i^{-1}= B$,
where $B$ is diagonalizable 
(Exercise \ref{_closure_of_non_ss_orbit_contains_ss_Exercise_}).

By \ref{semihopf},  a Hopf manifold
is Vaisman if and only if it is diagonal.
Taking $A$, $B$ and $\{g_i\}$ as above,
we see that the Vaisman Hopf manifold
$(\C^n\backslash 0)/\langle B \rangle$ contains in each open neighbourhood
of its  Kuranishi deformation space a
Hopf manifold isomorphic to \index[terms]{space!Kuranishi deformation}
$(\C^n\backslash 0)/\langle A \rangle$.
Since the property of admitting an LCK
potential is open, this implies that
all linear Hopf manifolds are LCK
with potential.

The notion of LCK manifold with potential\index[terms]{manifold!LCK!with potential}
is, in itself, conformally invariant.
By contrast, the standard definition
of Vaisman manifold, $\nabla \theta=0$,
is {\em not} conformally invariant.\index[terms]{manifold!Vaisman} 
If we rescale the Vaisman metric by
a constant, we arrive at the equation
\eqref{_omega_via_theta_Chapter_8_Equation_}:
\begin{equation}\label{_omega_via_theta_Equation_}
\omega= d^c \theta + \theta \wedge I(\theta),
\end{equation}
where $\theta$ is the Lee form,\index[terms]{form!Lee} and
$\omega$ the Vaisman metric.

It is not hard to see that
\eqref{_omega_via_theta_Equation_}
implies that $\omega$ has a potential
(\ref{dcthetaex}).
Moreover, every LCK metric with a
potential admits a conformal
representative such that
\eqref{_omega_via_theta_Equation_} is satisfied
(\ref{dcthetaex}).
We call the metrics that satisfy
\eqref{_omega_via_theta_Equation_}
``metrics with preferred conformal gauge''
(\ref{_preferred_gauge_Definition_}).

In the literature, \eqref{_omega_via_theta_Equation_}
is often taken as a definition of LCK manifolds
with potential\index[terms]{manifold!LCK!with potential} (e.g.    \cite{_Madani_Moroianu_Pilca:holo_}).

The Vaisman condition fixes the metric
in its conformal class, up to a constant,
and the condition \eqref{_omega_via_theta_Equation_}
fixes the constant. However, the preferred gauge
condition in itself does not fix the metric
in its conformal class (\ref{_non-unique-omega_via_theta_Example_}).

There are other ways to characterize LCK manifolds with
potential. By \ref{embedding} and \ref{_linear_LCK_pot_Corollary_}, 
a compact complex manifold
$M$ admits an LCK metric with potential \index[terms]{metric!LCK!with potential}if and
only if $M$ admits a complex embedding to a 
linear Hopf manifold. By \ref{_S^1_potential_Theorem_} and
\ref{logar}, an LCK manifold\index[terms]{manifold!LCK} admits an
LCK metric with potential if and only
if the $\Z$-action\index[terms]{action!$\Z$-} on its K\"ahler
cover is an exponent of a vector
field that is  lifted from $M$.
These theorems could be considered
alternative definitions as well.

\hfill

Further in
this chapter, we determine the structure of the K\"ahler
covers admitting proper potential \index[terms]{potential!LCK!proper}and show that it is very
close to being a cone. In fact, its metric completion is a
Stein variety with only one singular point, precisely the
one added to complete it, see \ref{potcon}. 

The proof we give here for \ref{potcon} is different
from the original one in \cite{ov_lckpot}. 
The original proof was based on the  Poincar\'e--Dulac\index[terms]{theorem!Poincar\'e--Dulac}
theorem, but it assumed that the deck transform action
on the K\"ahler cover is an exponent of a vector field
action. 

The present proof is based on techniques from complex
geometry. In Section
\ref{_Stein_and_Montel_preliminaries_} we present the
minimal necessary background on Stein varieties,\index[terms]{variety!Stein} 
on normal families of holomorphic sections,
and the Montel theorem.

\hfill

This chapter serves as an introduction to two important
notions of complex geometry, that of Stein varieties, and
the normal families. 

A complex space, or complex analytic space,\index[terms]{space!complex}
is a complex-analytic analogue of a scheme.
One defines a complex space as a ringed\index[terms]{space!ringed}
space locally isomorphic to the zero
set $Z_J$ of a coherent\index[terms]{sheaf!coherent}
ideal sheaf $J\subset \calo_B$, where\index[terms]{sheaf!ideal}
$B\subset \C^n$ is an open ball.
The ring of functions on $Z_J$
is $\calo_B/J$. 

A complex space is called
{\bf a complex variety} if its structure
sheaf has no nilpotents. Conversely,\index[terms]{variety!complex}
one can understand a complex space
as a complex variety which allows for
nilpotents in its structure sheaf.
We are not going to use Stein spaces
(only varieties), but almost all
literature that deals with this
subject is written about Stein spaces,\index[terms]{space!Stein}\index[terms]{variety!Stein} 
and not Stein varieties (\cite{_Grauert_Remmert:Stein_}).

Similar to affine varieties in algebraic
geometry, Stein spaces are the building blocks
of complex varieties and complex spaces.

One can characterize Stein spaces as
complex spaces that admit a closed
holomorphic embedding to $\C^n$, or
as spaces such that higher cohomology
of any coherent sheaf vanish; these two
conditions are equivalent by 
\index[persons]{Cartan, H.} the Cartan theorems A and B \index[terms]{theorem!the Cartan A, B}
(\cite{_Gunning_Rossi_}).

The third condition, due to \index[persons]{Oka, K.} Oka \cite{_Oka_},
is the existence of a strictly
plurisubharmonic exhaustion function. 

Further on, we introduce the normal
families, used to state the Montel theorem.\index[terms]{theorem!Montel}
A family $\{f_i\}\subset H^0(\calo_M)$ 
of holomorphic functions on a complex manifold
$M$ is called {\bf a normal family}\index[terms]{normal family} if 
it is uniformly bounded  on 
each compact $K\subset M$.
Montel theorem states that any
normal family is precompact, that is,
its closure in the space $H^0(\calo_M)$ of holomorphic
functions equipped with $C^0$-topology\index[terms]{topology!$C^0$} is compact.

This theorem has the same flavour as the  Arzel\`a--Ascoli\index[terms]{theorem!Arzel\`a--Ascoli}
theorem for the Lipschitz functions on a metric space.\index[terms]{function!Lipschitz}
Recall that a function $f$ on a metric space\index[terms]{space!metric}
is called Lipschitz, or 1-Lipschitz, if $|f(x)-f(y)|\leq d(x,y)$.
Any family of Lipschitz functions bounded on compacts
is precompact in $C^0$-topology;\index[terms]{topology!$C^0$} indeed, it is
precompact in Tychonoff topology (topology of
pointwise convergence), because the space of
bounded functions is compact in Tychonoff topology.\index[terms]{topology!Tychonoff}
However, for Lipschitz functions, pointwise
convergence implies $C^0$-convergence.

Every complex manifold is equipped with
a natural metric structure, called
\index[persons]{Kobayashi, S.} the Kobayashi pseudo-metric (see the exercises
to this chapter). Any bounded holomorphic
function can be interpreted as a Lipschitz
function with respect to the Kobayashi\index[terms]{pseudometric!Kobayashi}
pseudometric. Then, Montel theorem
becomes the special case of the 
Arzel\`a-Ascoli theorem for the 
Lipschitz functions.

Recall that a map $f:\; M \arrow M$
is called {\bf a contraction} of $M$
to an open subset $U\subset M$ if for each compact
subset $K\subset M$, a sufficiently
big iteration of $f$ gives $f^n(K) \subset U$.

Using Montel theorem, we prove that for any holomorphic
function $u$ on $M$, and any contraction $f:\; M \arrow M$
of $M$ to $U$, with $U$ having compact closure, 
the limit $\lim_n (f^n)^* u$ is constant.
This is the idea of the proof of \ref{potcon},
where this argument is applied to the
Stein completion\index[terms]{completion!Stein} of the K\"ahler $\Z$-cover\index[terms]{cover!K\"ahler $\Z$-} 
$\tilde M$ of an LCK manifold with potential.\index[terms]{manifold!LCK!with potential}

This is used to show that the Stein completion\index[terms]{completion!Stein}
(\ref{_Stein_completion_Definition_})
of $\tilde M$ coincides with its metric completion.

This is remarkable, because
the Stein completion itself is independent
from the metric; of course, the metric completion
is independent on  the complex structure.


 \section{Deformations of LCK structures}\index[terms]{structure!LCK}


 The following theorem is well-known.
 
\hfill

 \theorem\label{_Kodaira_Spencer:Theorem_}  ({\bf  Kodaira--Spencer stability theorem}, 
\cite{_Kod-Spen-AnnMath-1960_})\\
A small deformation of a compact K\"ahler manifold is K\"ahler.
 \index[terms]{deformation}\index[terms]{theorem!Kodaira--Spencer stability}

\hfill

We are not going to use this theorem; we recall it as a motivation for
two similar results in LCK geometry (\ref{defvai}, 
\ref{_LCK_pot_stable_Theorem_}). 
However, the LCK analogue of  Kodaira--Spencer stability 
is more elementary and straightforward.\index[terms]{Kodaira--Spencer stability}

For the sake of completeness,
we give the proof of  Kodaira--Spencer 
stability in the Appendix to this chapter (Section 
\ref{_KS_stability_Section_}).

\hfill

\definition 
	Let ${\goth A}$ be a property of compact complex manifolds,
	and ${\cal X}\stackrel \pi \arrow B$ a smooth, proper map, that is,
	a smooth family of compact manifolds. We say that
	${\goth A}$ is {\bf stable under small deformations}
	if the set of all $z\in B$ such that $X_z:=\pi^{-1}(z)$
	has ${\goth A}$ is open in $B$.

\hfill

This is why \ref{_Kodaira_Spencer:Theorem_}
is called ``the  Kodaira--Spencer stability theorem''.

\hfill

Besides K\"ahler, another stable 
property is being  {\bf Hermitian symplectic}, i.\,e.      
 admitting a symplectic form $\omega$ such that its\index[terms]{manifold!Hermitian symplectic}
$(1,1)$-part is Hermitian (Exercise 
\ref{_Hermitian_symplectic_Exercise_}).\index[terms]{form!symplectic}
However, other important properties are not stable:
a small deformation of a projective manifold,
such as a K3 surface or a compact torus, is\index[terms]{surface!K3}
not necessarily projective.

Unlike the K\"ahler class,\index[terms]{class!K\"ahler} the LCK class is not stable
under small deformations, see \cite{bel} or see Chapter
\ref{comp_surf}. Moreover, the Vaisman property is not
stable (\cite{go}, \ref{def_lckpot2Vai}). 
However, all small deformations of
Vaisman manifolds are LCK (\ref{defvai}).\index[terms]{manifold!Vaisman}

\medskip

Let $(\tilde M, \tilde \omega)$ be a K\"ahler cover
of an LCK manifold\index[terms]{manifold!LCK}, $\Gamma$ its deck transform group,
and $\chi:\; \Gamma \arrow \R^{>0}$
the homothety character\index[terms]{homothety character} that maps $\gamma\in \Gamma$
to the coefficient $c_\gamma$ such that 
$c_\gamma\tilde\omega=\gamma^*\tilde\omega$.
Recall that {\bf an automorphic function}\index[terms]{function!automorphic}
on $\tilde M$ is a function $f$ that satisfies
$c_\gamma f=\gamma^*f$, and the space of automorphic
functions is naturally identified with the 
space of sections of the weight bundle $L$ \index[terms]{bundle!weight}
(\ref{_automo_sections_of_L_Claim_}).

\medskip

\theorem \label{defvai}
Let ${\cal X} \stackrel\pi\arrow B$ be a smooth,
proper, holomorphic map, and $z\in B$ a point.
Assume that the fibre $X_z:= \pi^{-1}(z)$  admits a Vaisman metric.
{Then there exists a neighbourhood
	$W\ni z$ such that for each $y\in W$,
	the fibre $X_y:= \pi^{-1}(y)$ is LCK.}\index[terms]{manifold!LCK}

\smallskip

\proof 
We start by replacing the Vaisman metric
by a Vaisman metric of LCK rank 1\index[terms]{rank!LCK} (\ref{_rk_1_approximated_Theorem_}). 
Let $\tilde X_z=C(S_z)$ be the conical K\"ahler\index[terms]{rank!LCK}
cover of $X_z$. By \index[persons]{Ehresmann, C.} Ehresmann fibration theorem (see
\cite{ehr}), \index[terms]{theorem!Ehresmann's fibration}
$\pi$ is a locally trivial fibration. Replacing
$B$ by a sufficiently small open neighbourhood of
$z$, we may assume that $\pi$ is trivial
	as a smooth fibration: ${\cal X}=X_z\times B$.
Consider a cover $\tilde {\cal X}\arrow {\cal X}$ 
with $\tilde {\cal X}=\tilde X_z\times B$,
and let $\tilde X_y$ denote the fibres
of the projection $\tilde {\cal X}\stackrel {\tilde \pi} \arrow B$.

Recall from the proof of \ref{regsas} that the K\"ahler
form on $C(S_z)$ has global potential $\phi$ that is 
automorphic in the sense of \ref{_automo_forms_Definition_}. Now extend 
$\phi$ to $\tilde {\cal X} = \tilde X_z\times B$
using the projection $\tilde {\cal X}\arrow \tilde X_z$.
Restricting $\phi$ to $\tilde X_y\subset \tilde {\cal X}$,
we obtain an automorphic function $\phi_y$ on any $\tilde X_y$.

The form $dd^c \phi_y$
is closed, automorphic and of type (1,1). Therefore, $\tilde X_y$
is LCK whenever the pseudo-Hermitian form\index[terms]{form!pseudo-Hermitian} $dd^c \phi_y$
is positive definite. However, the complex structure
on $X_y$ depends smoothly  on $y\in B$, and hence 
the function $y \mapsto dd^c \phi_y$
is continuous, and its eigenvalues 
depend continuously on $y\in B$. 
Since $\phi_y$ is automorphic, it 
would suffice to check the positivity on
the fundamental domain of the $\Z$-action\index[terms]{action!$\Z$-} on $\tilde X_y$.
The closure of the fundamental domain of the $\Z$-action can be identified
with $[a, b] \times S_z \subset C(S_z) = \R^{>0}\times S_z$  
since $\tilde X_z= C(S_z)$ and the deck transform 
group acts by shifting the $\R^{>0}$-coordinate
(\ref{str_vai}). Since $S_z$ is compact, the
fundamental domain of $\Z$-action\index[terms]{action!$\Z$-} has compact closure.
Therefore, for $y$ sufficiently
	close to $z$, these eigenvalues remain positive,
and $\phi_y$ gives an automorphic K\"ahler potential.
\endproof
 
\smallskip

This result motivates the definition of LCK manifolds with
potential. This is the class of manifolds for which the
argument of \ref{defvai} can be applied.


\section{LCK manifolds with potential}\index[terms]{manifold!LCK!with potential}
\label{_Definition_of_LCK_potential:Section_}


\subsection{LCK manifolds with potential, proper and improper}\index[terms]{potential!LCK!proper}

\definition\label{_LCK_potential_Definition_} 
Let $M$ be an LCK manifold, and 
$(\tilde M, \tilde \omega)$ a K\"ahler cover.
We say that $M$ has  {\bf LCK potential}
if $\tilde M$ admits a {\bf global, automorphic} K\"ahler potential
$\phi:\; \tilde M \arrow \R^{>0}$, $dd^c \phi=\tilde
\omega$.  
\index[terms]{potential!LCK}\index[terms]{potential!K\"ahler}\index[terms]{potential!automorphic}
If the automorphic potential $\f$ is a proper map (the
preimages of compact subsets are compact),
then we say that $M$ admits {\bf a proper LCK potential},
otherwise it is {\bf an improper LCK potential}.\index[terms]{potential!LCK}
\index[terms]{potential!LCK!proper}

\hfill

\remark\label{_Oka_Stein_Remark_}
Oka's theorem claims that any manifold $M$
admitting a positive, proper, strictly
plurisubharmonic map $f:\; M \arrow \R$
is Stein (\ref{_Oka_Stein_Theorem_}).
We want to warn the reader that
the automorphic potential $\phi$
does not satisfy this condition:
a proper map $f:\; M \arrow \R^{>0}$
can be non-proper, if considered
as a map $f:\; M \arrow \R$, because
the natural embedding $\R^{>0}\arrow \R$
is non-proper. In other words, 
a proper function $f:\; M \arrow \R^{>0}$
can be non-proper if considered to be  a 
function  $f:\; M \arrow \R$. For example,
let $M=\C^n\backslash 0$, and let
$f(z)=|z|^2$ be the standard K\"ahler
potential on $M$. As a map
$f:\; M \arrow \R^{>0}$, this map
is proper, and hence  it gives an
automorphic potential for
the deck transform $z\mapsto 2z$.
However, it is not proper
as a map $f:\; M \arrow \R$,
and, indeed, $M=\C^n\backslash 0$
is not Stein.

\hfill

\example\label{_Vaisman_has_potential_Example_}
By \ref{_tilde_M_potential_Corollary_}, all Vaisman manifolds admit an
LCK potential:\index[terms]{potential!LCK} $4\tilde \omega = dd^ct^2$.\index[terms]{manifold!Vaisman} 
Indeed, by \ref{str_vai}  a Vaisman manifold is covered
by a cone $C(S)= S \times \R^{>0}$ of a Sasakian manifold
$S$, and the K\"ahler form\index[terms]{form!K\"ahler} $\tilde \omega$ on this cone
is given by $\tilde\omega = dd^c(t^2)$ (up to a constant
multiplier). 
By \ref{str_vai} again, the deck transform group acts on $C(S)$
by mapping $(s, t)$ to $(\phi(s), \lambda t)$, and hence 
the function $t^2$ is automorphic.

\hfill

\remark \label{_potential_norm_Lee_Remark_}
The formula
 $\tilde\omega = dd^c(t^2)$ from
\ref{_Vaisman_has_potential_Example_}
can be written
in a more invariant way as 
$\tilde \omega= \const \cdot dd^c(\tilde \omega(\theta^\sharp, I\theta^\sharp))$, 
where $\theta^\sharp$ denotes the Lee field.\index[terms]{Lee field}
Indeed, for a Vaisman manifold,\index[terms]{manifold!Vaisman} 
$\omega(\theta^\sharp, I\theta^\sharp)=|\theta|^2 =\const$ because
$\theta$ is parallel. However, on the cone 
$C(S)= S \times \R^{>0}$, we have $\tilde \omega= t^2 \omega$,
hence $\tilde \omega(\theta^\sharp, I\theta^\sharp))=\const \cdot t^2$.

\hfill

\remark
Existence of a K\"ahler potential 
for $\tilde \omega$ on $\tilde M$
does not always imply the existence of
LCK potential,\index[terms]{potential!LCK} because the automorphicity
condition is not guaranteed. For an example,
take an Inoue surface or \index[terms]{surface!Inoue}
an LCK OT-manifold, that is  a generalization
of an Inoue surface to higher dimension (see the proof of \ref{_Existence_of_LCK_metric_on_OT_}). The K\"ahler form \index[terms]{form!K\"ahler}
on a K\"ahler cover of an Inoue surface has
a potential (Chapter \ref{inoue_lck}; in fact, this K\"ahler cover is\index[terms]{manifold!Oeljeklaus--Toma (OT)}
Stein, so existence of the potential is guaranteed),
but this potential is never automorphic.

\hfill

\remark Another situation when an LCK
potential never exists is a blow-up of 
an LCK manifold in a point. As shown in
Sections \ref{_Blow_up_at_points_section_} and  \ref{_Blow_up_on_submanifolds_section_} a blow-up of an LCK manifold
is LCK. However, its K\"ahler cover contains
infinitely many copies of the exceptional
divisor, that is  isomorphic to $\C P^n$.
On each of these $\C P^n$, the K\"ahler
form cannot have a potential, because it is
not exact.

\subsection{LCK manifolds with potential and preferred
   gauge}\index[terms]{manifold!LCK!with potential and preferred gauge}

\definition\label{_preferred_gauge_Definition_}
Let $(M, \omega, \theta)$ be an LCK manifold,\index[terms]{manifold!LCK}
and $\tilde M\stackrel \pi \arrow M$ its K\"ahler cover
admitting an LCK potential\index[terms]{potential!LCK} $\phi:\; \tilde M\arrow \R^{>0}$.
In that case, we say that $(M, \omega, \theta)$ is {\bf an LCK manifold
with potential}\index[terms]{manifold!LCK!with potential}. If, in addition, 
$\pi^*\omega= \frac{dd^c
  \phi}{\phi}$, we say that $(M, \omega, \theta)$ is {\bf an LCK manifold
with potential and preferred conformal gauge},
or just {\bf LCK with potential and preferred gauge}.\index[terms]{manifold!LCK!with potential}
\index[terms]{manifold!LCK!with potential and preferred gauge}

\hfill

Locally conformally K\"ahler manifolds  with potential with the preferred gauge
can be characterized intrinsically in terms of a
useful equation \eqref{dctheta}.

\hfill

\proposition  { (\cite{ov_jgp_09}
} \label{dcthetaex}
Let $(M,I,\omega,\theta)$ be an LCK manifold, and 
$\tilde M$ its K\"ahler cover. Then $(M,\omega, \theta)$ admits
an LCK potential\index[terms]{potential!LCK} with preferred conformal gauge if and
only if \index[terms]{manifold!LCK!with potential and preferred gauge}
\begin{equation}\label{dctheta}
\omega = d^c\theta+\theta\wedge I\theta.
\end{equation}	

\proof Let $\pi:\tilde M\arrow M$ be the K\"ahler cover
with global automorphic potential $\phi=e^{-\nu}$. Then
the pullback of the Lee form \index[terms]{form!Lee}is $\pi^*\theta=d\nu$.
Indeed, $\omega= \phi^{-1} dd^c\phi$, which gives
$d\omega= d(\phi^{-1}) dd^c \phi$, and hence 
\[ \pi^*\theta= \phi d(\phi^{-1})= - \phi
\frac{d\phi}{\phi^2}=- d \log \phi=d \nu.
\]
This gives
\begin{equation}\label{_potential_pref_gauge_computation_Equation_}
\begin{aligned}
\pi^*\omega&=\phi^{-1}\tilde\omega=e^\nu dd^c e^{-\nu}\\& =e^\nu d(-e^{-\nu}d^c\nu)=e^\nu(e^{-\nu}d\nu\wedge d^c\nu-e^{-\nu}dd^c\nu)\\
&=d\nu\wedge d^c\nu-dd^c\nu=d\nu\wedge d^c\nu+d^cd\nu\\
&=\pi^*\theta\wedge\pi^*(I\theta)+d^c\pi^*\theta\\
&=\pi^*(d^c\theta+\theta\wedge I\theta),
\end{aligned}
\end{equation}
since $\pi^*$ commutes  with both $d$ and $I$.
As $\pi^*$ is injective, the result is proven.  

Conversely, let $(M, \omega, \theta)$
satisfy \eqref{dctheta}. Choose a cover 
$\tilde M \stackrel \pi\arrow M$ such that 
$\pi^*\theta= d\nu$, and let $\phi:= e^{-\nu}$.
Then $dd^c\phi= \phi \pi^*(d^c\theta+\theta\wedge
I\theta)$, as follows from 
\eqref{_potential_pref_gauge_computation_Equation_}, and hence  
$\phi$ is an LCK potential\index[terms]{potential!LCK} that satisfies $\omega= \frac{dd^c\phi}{\phi}$.
\endproof

\hfill

\remark
In the literature, LCK manifolds with potential\index[terms]{manifold!LCK!with potential} are usually
tacitly equipped with the metric with preferred gauge.
However, it is not unique in its conformal class
(\ref{_non-unique-omega_via_theta_Example_}).

\hfill

\remark 
Then equation \eqref{dctheta} is precisely equation
\eqref{_omega_via_theta_Chapter_8_Equation_}, valid on Vaisman
manifolds after fixing the conformal gauge in such a way that
 $|\theta|=1$.
This implies that the Vaisman manifolds\index[terms]{manifold!Vaisman}
that satisfy \eqref{_omega_via_theta_Equation_} are
LCK with potential and preferred conformal gauge\index[terms]{manifold!LCK!with potential and preferred gauge}.
However, this property is not invariant under
a homothety: if we multiply a Vaisman metric
by a number, the Lee form\index[terms]{form!Lee} will not change,
and the equation \eqref{dctheta} becomes invalid.

\hfill

\remark
The Vaisman metric  satisfying  
the equation 
$\omega= d^c \theta + \theta \wedge I(\theta)$
is unique in its conformal class. 
Indeed, this equation implies that
$|\theta|=1$ 
(\ref{_preferred_gauge_Vaisman_theta=1_Proposition_}).
The Vaisman condition implies Gauduchon, 
which gives uniqueness up to a scalar multiplier,
and $|\theta|=1$ fixes the scalar multiplier.
However, the representative that satisfies \eqref{dctheta}
is not unique without the Vaisman condition, as the following
example implies. 

\hfill

\example\label{_non-unique-omega_via_theta_Example_}
Let $(M, \omega)$ be a compact LCK manifold with proper potential,\index[terms]{potential!LCK!proper}
and $\phi$ an LCK potential\index[terms]{potential!LCK} on the K\"ahler cover
$(\tilde M, \tilde \omega)$. Consider an
automorphic holomorphic function $f$ on $\tilde M$.
Such a function is known to exist. For example, on a classical
Hopf manifold $(\C^n \backslash 0)/\langle \lambda\Id\rangle$,
any quadratic polynomial is automorphic.
Multiplying $f$ by a sufficiently small scalar
coefficient $a>0$, we can always assume that $|f| < \phi$.
Indeed, it suffices to check that $|f| < \phi$ on a compact
closure of the fundamental domain $D$, and hence 
it suffices to take $a:= \frac{\inf_{x\in D} \phi(x)}{\sup_{x\in D}
	|f(x)|}$. This gives $\phi + \Re f>0$. Since $dd^c \Re f=0$, the
function
$\phi_1:= \phi + \Re f$ is also an LCK potential\index[terms]{potential!LCK}
for the same K\"ahler metric $\tilde \omega=dd^c \phi=dd^c \phi_1$,
and $\frac {dd^c \phi_1}{\phi_1}$ is an
LCK metric with potential\index[terms]{metric!LCK!with potential} on $M$ satisfying the preferred gauge
condition \eqref{dctheta}, as follows from \ref{dcthetaex}. This metric is
conformally equivalent to $\omega= \frac {dd^c \phi}{\phi}$,
but the conformal multiplier is clearly non-trivial.

%

\hfill

\remark\label{_preferred_gauge_elliptic_Remark_}
Let $(M, \omega,\theta)$ be an LCK manifold with
potential and preferred conformal gauge,
and $(M, e^f \omega, \theta + df)$ be an 
LCK structure\index[terms]{structure!LCK} in the same conformal class.
Then $(M, e^f \omega, \theta + df)$ has
preferred conformal gauge if and only if
\[
e^f \omega= d^c\theta + dd^c f + (\theta + df)\wedge I(\theta + df).
\]
Clearly, this condition is elliptic in $f$.
Indeed, 
\[ e^f \omega^n = 
\omega^{n-1} \wedge (d^c\theta + dd^c f + (\theta + df)\wedge I(\theta + df))
\]
is equivalent to the equation
\[
\frac{\omega^{n-1} \wedge (dd^cf)}{\omega^n}=
\frac{\omega^{n-1} \wedge (e^f \omega^{n-1} - 
(\theta + df)\wedge I(\theta + df))}{\omega^n}.
\]
The left-hand side of this equation
is the Laplacian $\Delta(f)$, and the right-hand side
is a first order non-linear differential operator.

\hfill

\remark
As follows from \ref{_preferred_gauge_elliptic_Remark_},
the preferred conformal gauge condition is elliptic.
In other words, there is at most a finite-dimensional
space of conformal multipliers of a given form $\omega$
that have preferred conformal gauge.

\subsection[The monodromy of LCK manifolds with proper potential]{The monodromy of LCK manifolds with\\ proper potential}\index[terms]{monodromy}\index[terms]{potential!LCK!proper}

Let $M$ be an LCK manifold\index[terms]{manifold!LCK} and $\tilde M\stackrel
\pi\arrow M$ its smallest
K\"ahler cover. Recall that {\bf the monodromy group}\index[terms]{group!monodromy}
of $M$ is the deck transform group of $\tilde M\stackrel
\pi\arrow M$. Since the LCK form can be considered
as a closed (1,1)-form with values in the weight bundle\index[terms]{bundle!weight}
$L$, the monodromy group coincides with the monodromy
group of $L$. This gives a natural embedding
$\Gamma\subset \R^{>0}$  of the monodromy group.
We interpret the properness of the LCK potential \index[terms]{potential!LCK!proper}
in terms of the monodromy as follows.

\hfill

\proposition
\label{proper_equiv_Z}
Let $M$ be a compact LCK manifold, $\Gamma\subset \R^{>0}$ 
its monodromy group, and\index[terms]{group!monodromy} 
$(\tilde M, \tilde \omega)$ its smallest K\"ahler cover, with
$\tilde M/\Gamma=M$. Assume that $\tilde \omega$ admits a 
$\Gamma$-automorphic, positive  K\"ahler potential $\phi$. 
The map  $\phi$ is proper if and only if
	$\Gamma\simeq\Z$.

\hfill

\proof Suppose $\Gamma\simeq\Z$ and 
let $\gamma$ be a generator of $\Z$, such that
$\gamma^*\phi=\lambda\phi$, and $\pi:\; \tilde M \arrow
\tilde M/\Gamma =M$ the quotient map.
{Then $\phi^{-1}([1, \lambda[)$ is a fundamental domain\index[terms]{cover!fundamental domain}
	of the $\Gamma$-action.} Therefore, 
$\pi:\;\phi^{-1}([1,\sqrt  \lambda])\arrow M$ is bijective
onto its image, that is  compact, and hence  
$\phi^{-1}([1,\sqrt  \lambda])$ is also compact.
This implies that 
{the  preimage of any closed interval is compact}.

Conversely, suppose that $\phi$ is proper and, by contradiction,  
assume that $\Gamma\not\simeq \Z$. Then {
	$\Gamma$ is a dense subgroup of $\R^{>0}$.}
Fix $x\in \tilde M$ and a nonempty interval
$]a,b[\subset \R^{>0}$, and let 
$${\goth H}:=\left\{\gamma\in \Gamma\,;\, \phi(\gamma(x))\in ]a, b[\right\}.$$ 
Since $\Gamma$ is dense, this ${\goth H}$ is infinite.
However, $\phi({\goth H}\cdot x)\subset [a, b]$. But
{${\goth H}\cdot x\subset \Gamma\cdot x$ that is  discrete, as $\Gamma$ is discrete. And hence the infinite, discrete set ${\goth H}\cdot x$
	is contained in the compact set $\phi^{-1}([a,b])$.}
This contradiction ends the proof. \endproof

\hfill

\remark \label{_proper_potential_LCK_rank_1_Remark_}
Let $\tilde M\stackrel \pi \arrow M$ be the smallest
K\"ahler cover of $M$, that is, the smallest cover
where $\theta$ becomes exact. Then  its monodromy\index[terms]{group!monodromy}
is equal to $\Z^k$, where $k$ is the rank of the smallest rational
subspace $V\subset H^1(M, \Q)$ such that $V\otimes_\Q \R$
contains $[\theta]$.
In particular,  the  class $[\theta]\in H^1(M, \R)$ is
proportional to a rational one precisely when $\tilde
M\stackrel \pi \arrow M$ is a $\Z$-covering\index[terms]{cover!K\"ahler $\Z$-}, i.\,e.     \ the
potential is proper. In other words, $(M, \omega, \theta)$
has proper potential \index[terms]{potential!LCK!proper}if and only if its LCK rank is 1.\index[terms]{rank!LCK}

\hfill

As we already saw, Vaisman manifolds are examples of  LCK manifolds\index[terms]{manifold!Vaisman}
with potential.\index[terms]{manifold!LCK!with potential}
The proof of \ref{defvai} can be repeated to obtain the following.

\hfill

\theorem \label{_LCK_pot_stable_Theorem_}
{ (\cite{ov_lckpot})}
	 The property of being LCK with potential
	is stable under small deformations. \index[terms]{deformation}
\blacksquare

\hfill

\remark \label{lck_pot_subm}
	Any complex submanifold $Z\subset M$\index[terms]{submanifold}
of an LCK manifold with (proper) potential is\index[terms]{potential!LCK!proper} 
an LCK manifold with (proper) potential.

\hfill

\remark\label{_proper_potential_S^1_Remark_}
Let $(M, \omega, \theta)$ be an LCK manifold with
an automorphic potential $\phi:\; \tilde M \arrow \R^{>0}$,
and $\Gamma$ the deck group of the minimal K\"ahler cover\index[terms]{cover!minimal K\"ahler} $\tilde M \arrow M$.
Denote by $\Gamma_0\subset \R^{>0}$ the monodromy of the
action of $\Gamma$ on $\phi$, that is, the image of the
homothety character\index[terms]{homothety character} $\chi:\; \Gamma \arrow \R^{>0}$. 
The potential $\phi$ is proper if and only if $\Gamma_0$
is a closed subgroup of $\R^{>0}$, or, equivalently,
when $\Gamma_0 = \Z$. Then $\phi$ defines a map
$\underline \phi:\; M \arrow S^1= \frac{\R^{>0}}{\Gamma_0}$.
The fibres of $\underline \phi$ are all strictly
pseudoconvex submanifolds in $M$, just like those of $\phi$ (in fact, these fibres are the 
same, up to an obvious identification).



\subsection{$d_\theta d^c_\theta$-potential}
\label{_potential_weight_Subsection_}
\index[terms]{bundle!weight}

In this section, we shall consider real flat line bundles.
Notice that all oriented real line bundles are trivial
in the smooth category; however, the choice of trivialization
is not canonical.\index[terms]{bundle!line!flat}

\hfill

Let $\theta$ be a closed form, and 
$L_\theta:\; \Lambda^*(M) \arrow \Lambda^{*+1}(M)$
be the operator of multiplication by $\theta$.
We define an operator $d_\theta$, also called 
{\bf the Morse--Novikov differential}, {\bf twisted differential},\index[terms]{differential!twisted} \index[terms]{differential!Morse--Novikov} \index[terms]{differential!Lichnerowicz}
or {\bf \index[persons]{Lichnerowicz, A.} Lichnerowicz differential}, as $d_\theta= d - L_\theta$.
The cohomology of this differential is equal to the
cohomology of the local system defined by\index[terms]{connection!flat}
the flat connection $\nabla_\theta:=\nabla_0 - L_\theta$
on the trivial line bundle $L$, where $\nabla_0$ is
the trivial connection. We can consider $d_\theta$
as the de Rham operator with coefficients in
the flat bundle $(L, \nabla_0 - L_\theta)$
(Section \ref{_Riemann_Hilbert_Section_}).

For the applications to the LCK geometry,\index[terms]{geometry!LCK}
$\theta$ is the Lee form\index[terms]{form!Lee}, and $(L, \nabla_\theta)$
the weight bundle. In this case, we can
also define the operator $d_\theta^c := I d_\theta I^{-1}$.
The operators $d_\theta, d_\theta^c$ behave in many respects
similar to the usual $d, d^c$; they can be understood
as $d, d^c$ taking values in the weight bundle.\index[terms]{bundle!weight}

\hfill

\definition
Let $\rho:\; \pi_1(M) \arrow \R^{>0}$ be a character.
Differential forms
that satisfy $\gamma^*\eta= \rho(\gamma)\eta$
are called  {\bf $\rho$-automorphic}. 
When $\rho$ is the homothety character,\index[terms]{homothety character}
$\Lambda_\rho^* M$ are automorphic forms\index[terms]{form!automorphic} in the sense
of \ref{_automo_forms_Definition_}.

\hfill

\definition\label{_d_theta_d^c_theta-potential_Definition_}
Let $\phi_0:\; M \arrow \R^{>0}$ 
be a positive smooth function  satisfying 
\begin{equation}\label{_LCK_pot_MN_Definition_}
d_\theta d^c_\theta(\phi_0)=\omega. 
\end{equation}
Then $\phi_0$ is called {\bf a $d_\theta
  d^c_\theta$-potential}
for the LCK metric.\index[terms]{metric!LCK}

\hfill

\proposition \label{_LCK_pot_via_d_theta_d^c_theta_Proposition_}
Let $(M,\omega,\theta)$ be a compact LCK manifold.
Then $M$ has an LCK potential \index[terms]{potential!LCK}if and only if 
$M$ has $d_\theta d^c_\theta$-potential.

\hfill

\proof
Consider
the smallest covering $\pi:\; \tilde M \arrow M$ such that
$\pi^*\theta$ is exact, and take a function $\nu$ on $\tilde M$ 
satisfying $d\nu=\pi^*\theta$. Since $\pi^*\theta$ is invariant
under the deck transform group $\Gamma$,
for each $\gamma\in \Gamma$ one has $\gamma^*\nu=\nu +
c_\gamma$, where $c_\gamma$ is a constant. Consider the
multiplicative character $\chi:\; \Gamma\arrow \R^{>0}$
given by $\chi(\gamma)= e^{-c_\gamma}$. 
The K\"ahler form\index[terms]{form!K\"ahler} $\tilde \omega$ 
on the K\"ahler cover of $M$ can be written
as $\tilde \omega=e^{-\nu}\pi^*\omega$. Indeed,
\begin{equation*}
	\begin{split} 
d(e^{-\nu}\pi^*\omega)&= - e^{-\nu} d\nu \wedge \pi^*\omega
+ e^{-\nu}\pi^*d\omega\\
& = - e^{-\nu} d\nu \wedge \pi^*\omega
+ e^{-\nu}\pi^*\theta\wedge \pi^*\omega =0
	\end{split}
\end{equation*}
because $d\nu=\pi^*\theta$.
Then $\chi(\gamma)$ is the homothety character\index[terms]{homothety character}
of $\tilde\omega$, that is,
$\gamma^*\tilde \omega= \chi(\gamma)\tilde \omega$.

Let $\Lambda_\chi^* M$ denote the space 
of $\chi$-automorphic forms on $\tilde M$
(same as the automorphic forms\index[terms]{form!automorphic} in the sense
of \ref{_automo_forms_Definition_}).

The map $\Lambda^* M \stackrel \Psi \arrow \Lambda_\chi^* M$
mapping $\eta$ to $e^{-\nu}\pi^*\eta$ makes the following
diagram commutative:
\[
\begin{CD}
\Lambda^* M @>\Psi>> \Lambda_\chi^*M\\
@V {d_\theta} VV @V {d} VV\\
 \Lambda^*M @>\Psi>> \Lambda_\chi^*M
\end{CD}
\]
Then $\Psi$ maps a $d_\theta d^c_\theta$-potential
 to a K\"ahler potential $e^{-\nu}$ on $\tilde M$,
and {\em vice versa}, making these notions equivalent.
\endproof

\hfill

\remark\label{_potential_section_of_wb_Remark_}
The equation $\omega= d_\theta d^c_\theta (\phi_0)$
can be interpreted in a more straightforward way
if we view $\phi_0$ as a section of the
weight bundle $(L, \nabla)$. After passing to the K\"ahler
cover, the monodromy of $L$ becomes trivial,\index[terms]{monodromy}
and it becomes trivialized by parallel sections.\index[terms]{section!parallel}
Under this trivialization, the lifts of 
$d_\theta$ and $d_\theta^c$ correspond to the
usual $d, d^c$, the sections of $L$ become
automorphic functions, and the equation
for $d_\theta d^c_\theta$-potential
$\omega= d_\theta d^c_\theta (\phi_0)$
becomes the equation 
$\tilde \omega= dd^c (\psi)$
for the K\"ahler potential.

\hfill

\remark
The equation for the preferred 
conformal gauge (\ref{_preferred_gauge_Definition_}) can be written 
as \[ \omega= d^c\theta + \theta \wedge \theta^c =d_\theta
d^c_\theta(1):
\] 
the constant function
is a $d_\theta d^c_\theta$-potential for $\omega$.
If we interpret this equation as in
\ref{_potential_section_of_wb_Remark_}, we can 
consider ``1'' as a non-degenerate section of $L$,
denoting it as ${\goth f}$.
This choice fixes the choice of $\theta$, giving
$\theta = -\frac{\nabla({\goth f})}{\goth f}$.
Using the same trivialization, the $L$-valued 2-form 
$d_\theta d^c_\theta(1)$ becomes a 2-form on $M$,
that can be written as 
$\omega=\frac{d_\theta d^c_\theta{\goth f}}{\goth f}$.
However, for most practical purposes the equation
$\omega= d_\theta d^c_\theta(1)= d^c\theta+ \theta\wedge \theta^c$
is just as convenient.

\hfill

\remark\label{_preferred_gauge_via_d_theta_Remark_}
Let $\omega=d_\theta d^c_\theta(\phi_0)$ be an LCK form with
$d_\theta d^c_\theta$-potential $\phi_0$. \index[terms]{form!LCK}
When $\phi_0=1$, this gives an equation for a
preferred gauge; for non-constant $\phi_0$, the form $\omega$
is usually  not a form of preferred conformal gauge.
Indeed, $\omega$ is of preferred conformal gauge
if and only if 
$\omega=d^c\theta+ \theta\wedge \theta^c$,
and $\omega=d_\theta d^c_\theta(\phi_0)$
is equivalent to 
\[ \omega= dd^c\phi_0+ d\phi_0 \wedge \theta^c + 
\theta\wedge d^c\phi_0 + \phi_0\cdot (d^c\theta+ \theta\wedge \theta^c),
\]
and the equation
\[
d^c\theta+ \theta\wedge \theta^c =
dd^c\phi_0+ d\phi_0 \wedge \theta^c + 
\theta\wedge d^c\phi_0 + \phi_0\cdot(d^c\theta+ \theta\wedge \theta^c)
\]
is elliptic in $\phi_0$, and hence  it has only finite-dimensional space of
solutions. This implies that $\omega=d_\theta d^c_\theta(\phi_0)$
is of preferred gauge for at most a finite-dimensional space of
$d_\theta d^c_\theta$-potentials 
$\phi_0$. \footnote{Compare this with \ref{_preferred_gauge_elliptic_Remark_},
where a non-linear elliptic equation was written for
preferred gauge metrics that are  conformally 
equivalent to a given one.}

\hfill

\example 
Consider 
the classical Hopf manifold $\frac{\C^n\setminus  0}{\Z}$,
with $\Z$ acting by homotheties. The
LCK form\index[terms]{form!LCK} and Lee form\index[terms]{form!Lee} can be written on $\C^n\setminus  0$ as 
  $\omega=|z|^{-2}\sum dz_i\wedge d\bar z_i$,
and respectively $\theta=-d\log|z|^2$. In this case, the $d_\theta
d^c_\theta$-potential  $\f_0\in\C(M)$ is the constant function 
$f_0=1$, whereas the potential function on $\C^n\setminus  0$ is $|z|^2$.

\subsection{Deforming an LCK potential to a proper potential}\index[terms]{potential!LCK}\index[terms]{potential!LCK!proper}

For any Vaisman manifold of LCK rank\index[terms]{rank!LCK} $>1$,\index[terms]{rank!LCK}\index[terms]{manifold!Vaisman}
its Vaisman metric can be approximated by a
Vaisman metric of LCK rank 1 (\ref{_rk_1_approximated_Theorem_}).
The proof of this result involved harmonic decomposition\index[terms]{decomposition!harmonic}
for 1-forms on Vaisman manifolds. A similar result is true
for LCK manifolds with potential,\index[terms]{manifold!LCK!with potential} but its proof is much
more elementary.

\hfill

\theorem {(\cite{{ov_jgp_09}, ov_jgp_16})} 
\label{defor_improper_to_proper}
Let $(M, \omega, \theta, \phi)$ be a compact LCK manifold
with improper LCK potential. Then $(\omega, \theta, \phi)$
can be approximated in the ${C}^\infty$-topology by an LCK
structure with proper LCK potential.\index[terms]{potential!LCK!proper}
\index[terms]{potential!improper}\index[terms]{topology!$C^\infty$}

\hfill

\proof 
Using \ref{_LCK_pot_via_d_theta_d^c_theta_Proposition_},
we write the $d_\theta d^c_\theta$-potential 
for $\omega$: 
\[ \omega=d_\theta d^c_\theta (\phi_0).
\]
 Replace $\theta$ with a form $\theta'$ with rational
cohomology class $[\theta']$ in a sufficiently small
${C}^\infty$-neighbourhood of $\theta$, and 
let $\omega':=d_{\theta'} d_{\theta'}^c(\phi_0)$.
Then $\omega'$ approximates $\omega$
in ${C}^\infty$-topology, and hence  for $\theta'$
sufficiently close to $\theta$, the form $\omega'$
is positive. It is $d_{\theta'}$-closed, because
$d_{\theta'}^2=0$, and hence  $0=d_{\theta'}\omega'=d\omega'-\theta'\wedge\omega'$.  
This implies that $(\omega', \theta')$
is an LCK structure.\index[terms]{structure!LCK} The LCK rank of\index[terms]{rank!LCK}
an LCK manifold\index[terms]{manifold!LCK} is the dimension of the smallest
rational subspace $W\subset H^1(M, \Q)$ such that
$W\otimes_\Q \R$ contains the cohomology class of the
Lee form.\index[terms]{form!Lee} Since $[\theta']$
is rational,  $W=\langle[\theta']\rangle$, and 
$(M,\omega', \theta')$ has LCK rank 1. \endproof

\hfill

Any LCK manifold with potential\index[terms]{manifold!LCK!with potential} that satisfies 
$b_1(M)>1$ admits a metric with improper potential.
This is the easiest way of producing LCK metrics
of arbitrary LCK rank.\index[terms]{rank!LCK}

\hfill

\proposition \label{_impro_pote_exists_Proposition_}
Let $(M, \omega, \theta)$ be a compact 
LCK manifold with $d_\theta d^c_\theta$-potential
$\phi_0$, and suppose 
$b_1(M)>1$. Then $M$ admits an
LCK metric $(M, \omega', \theta')$ with improper
potential, the same $d_\theta d^c_\theta$-potential
$\phi_0$, and arbitrary LCK rank between 2 and $b_1(M)$.
Moreover, $(\omega', \theta')$ can be chosen in an arbitrary
${C}^\infty$-neighbourhood of $(\omega, \theta)$.

\hfill

\proof 
Choose a closed $\theta'$ in a sufficiently
small neighbourhood of $\theta$, and let $W_\theta$ be the 
smallest rational subspace of $H^1(M, \R)$ such that
$W_\theta\otimes_\Q \R$ contains $[\theta]$. Since the
choice of the cohomology class $[\theta']$ is arbitrary
in a neighbourhood of $[\theta]$, the dimension of
$W_\theta$ can be chosen in an arbitrary way. Choosing 
$\theta'$ sufficiently close to $\theta$, we can
assume that the (1,1)-form $\omega':=d_{\theta'}
d^c_{\theta'}(\phi_0)$ is positive definite. 
Then $(M, \omega', \theta')$ is an LCK manifold
with improper potential and arbitrary LCK rank.\index[terms]{rank!LCK}
\endproof


\section{Stein manifolds and normal families}
\label{_Stein_and_Montel_preliminaries_}


In this section, 
we give the definitions and properties that 
we shall use, see for example \cite{demailly,_Gunning_Rossi_}.

\subsection{Stein manifolds}
\index[terms]{variety!Stein}\index[persons]{Stein, K.}
\label{_Stein_Subsection_}

\hfill

\definition 
	A complex variety\index[terms]{variety!complex} with isolated singularities $M$ is called {\bf holomorphically convex}\index[terms]{variety!complex!holomorphically convex}
if for any infinite discrete subset $S\subset M$,
	there exists a holomorphic function $f\in \calo_M$ which
	is unbounded on $S$.

\hfill

\remark
This definition is equivalent to the one given in 
most textbooks, such as \cite{demailly}.
Indeed, it is often used as the standard definition, see
for instance \cite{_Kapovich:Kleinian_,_Bogomolov_Oliveira_}.

\hfill

Clearly, every compact complex variety is holomorphically convex.
Holomorphically convex varieties can be characterized in terms
of their boundary, as follows.

\hfill

\theorem (\index[persons]{Grauert, H.}Grauert's solution of the Levi problem) 
\label{_Grauert_Levi_Theorem_}
\\
Let $M$ be a complex manifold with smooth boundary $S$.
Assume that the Levi form\index[terms]{form!Levi} on $S$ is strictly pseudoconvex.
Then $M$ is holomorphically convex.

\proof \cite[\S 2, Theorem 1]{_Grauert_}. \endproof

\hfill

\definition 
	A complex variety $M$ is called {\bf Stein}
	if it is holomorphically convex, and 
	for each distinct $x, y\in M$ 
there is a holomorphic function $f$ separating them,
 $f(x)\neq f(y)$.\index[terms]{variety!Stein}

\hfill

\remark 
	Equivalently, a complex variety is Stein
		if it admits a closed holomorphic embedding into $\C^n$
(\cite{_Gunning_Rossi_}).

\hfill

Recall that a map between topological spaces is {\bf proper} if the preimage of any compact set is compact. The following is a characterization of Stein manifolds:

\hfill

\theorem  \label{_Oka_Stein_Theorem_}
({\bf K. \index[persons]{Oka, K.} Oka}, \cite{_Oka_})
	A complex manifold $M$ is Stein \index[terms]{theorem!Oka}\index[terms]{manifold!Stein}
	if and only $M$ admits a K\"ahler metric with 
	a {K\"ahler potential that is  positive
		and proper}.\index[terms]{potential!K\"ahler}

\hfill

Recall that a continuous function $f:\; M \arrow \R$ is called
{\bf an exhaustion function} if $f^{-1}(]-\infty, a])$\index[terms]{function!exhaustion}
is compact for all $a\in \R$. \index[persons]{Oka, K.} Oka's theorem is often
stated in terms of exhaustion functions: a manifold
is Stein if and only if it admits a strictly
plurisubharmonic exhaustion function.\index[terms]{function!plurisubharmonic}
A version of this theorem is valid on complex varieties,
but the definition of strictly plurisubharmonic functions
is complicated, see \cite[Theorem 6.1]{_Fornaess_Narasimhan_}.

\hfill

We repeat the caution that we have already
taken in \ref{_Oka_Stein_Remark_}. Though
$f:\; M \arrow \R$ is proper and positive,
the condition of being proper as
a map $f:\; M \arrow \R^{>0}$
is more relaxed, and does not need
to imply that $M$ is Stein;
take, for example, the K\"ahler
potential $f(z) = |z|^2$, that is  proper as a map
$f:\; M \arrow \R^{>0}$ on $M=\C^n \backslash 0$.
It is not hard to see that $M=\C^n \backslash 0$
is not Stein unless $n=1$.



%
%
%
%
%
%

\hfill

In terms of cohomology, the Stein property is characterized by the following.

\hfill

\theorem  {\bf (H. \index[persons]{Cartan, H.} Cartan, 1951)}\index[persons]{Cartan, H.} \index[terms]{theorem!Cartan}
	{A complex variety $M$ is Stein}
	if and only if for any coherent sheaf $F$ on $M$,\index[terms]{sheaf!coherent}
	{its cohomology $H^i(F)$ vanishes for all $i>0$.}

\hfill

We shall also need the following result, that is  well known
for smooth Stein manifolds (\cite{_Behnke_Stein_}). For Stein
spaces, it is proven in \cite{_Stein:58_}; see also
\cite{_Andreotti_Narasimhan_}.
Recall that a subset of a topological space is
{\bf precompact} if its closure is compact.

\hfill

\proposition \label{_union_Stein_Proposition_}
Let $M= \bigcup_i M_i$ be the union
of an increasing sequence of Stein varieties, with
each $M_i$ precompact and open in $M_{i+1}$. Then $M$ is Stein.
\endproof

\hfill

\remark
In \cite{_Stein:58_,_Andreotti_Narasimhan_}, this result is stated
in the assumption that each $M_i$ is {\bf relatively Runge} in $M_{i+1}$,
that is, each holomorphic function defined on the closure \index[terms]{variety!relatively Runge}
of $M_i$ can be approximated by restrictions of 
holomorphic functions on $M_{i+1}$. \index[persons]{Oka, K.} The Oka--\index[persons]{Weil, A.}Weil theorem\index[terms]{theorem!Oka--Weil}
(\cite[Theorem 18]{_Fornaess_Forstneric_Wold_}) 
shows that precompact open subsets
$M_i \subset M_{i+1}$ are always relatively Runge.

\hfill

\corollary\label{_union_hol_convex_Corollary_}
Let $X= \bigcup_i X_i$ be the union
of an increasing sequence of holomorphically convex manifolds, with
each $X_i$ precompact and open in $X_{i+1}$. Then $X$ is
holomorphically convex.

\hfill

\proof \index[persons]{Remmert, R.}The  Remmert reduction theorem (\ref{rem_red}) implies
that each $X_i$ admits a proper holomorphic map
$X_i \arrow M_i$, with $M_i$ a Stein variety.
Then the natural projection
$\bigcup_i X_i \arrow \bigcup_i M_i$
is proper and holomorphic.
By \ref{_union_Stein_Proposition_},
$\bigcup_i M_i$ is Stein, and hence 
$\bigcup_i X_i$ is holomorphically convex.
\endproof

\subsection{Normal families of functions} 

\definition \label{normal_family}
Let $M$ be a complex manifold, 
and ${\cal F}\subset H^0(\calo_M)$ a 
family of holomorphic functions. We call ${\cal F}$ 
{\bf a normal family} if for each compact
$K\subset M$ there exists $C_K>0$ such that
for each $f\in {\cal F}$, $\sup_K |f| \leq C_K$.\index[terms]{normal family}

\hfill

\begin{lemma}\label{_f'_uniformly_bounded_Lemma_}
Let $M$ be a complex Hermitian manifold,
${\cal F}\subset H^0(\calo_M)$ a normal family,\index[terms]{normal family}
and $K\subset M$ a compact subset. {Then there exists
	a number $A_K>0$ such that $\sup_K |f'| \leq A_K$.}
\end{lemma}

\hfill

\proof
Assume that the contrary holds. Then there exists $x\in K$,
$v\in T_xM$, and a sequence $f_i \in {\cal F}$
such that $\lim_i |D_v f_i|=\infty$. Pick a disk
$\Delta\stackrel j \hookrightarrow M$ with compact
closure in $M$, tangent to $v$ in $x$, 
such that $j(0)=x$. Let $w=tv$ be the unit tangent vector.
Then $\sup_\Delta |f_i| < C_\Delta$ by the normal family\index[terms]{normal family}
condition. {By the Schwarz lemma (see \cite{jost}), this implies $|D_w f_i|<C_\Delta$.}\index[terms]{lemma!Schwarz}
However, $t^{-1}\lim_i |D_w f_i| =\lim_i |D_v f_i|=\infty$, yielding a 
contradiction. \endproof

\subsection{The $ C^0$ - topology on spaces of functions}
	
\definition \label{_compact_open_topology_Definition_}
\index[terms]{topology!compact-open, $C^1$}
Let $C(M)$ be the space of functions on a topological space.
	The {\bf topology of uniform convergence on compacts} (also known
	as {\bf compact-open topology}, usually denoted as $ C^0$) 
	is the topology on $C(M)$ whose  base of open sets
	is given by 
	\[ U(X, C):= \{f\in C(M)\ \ |\ \ \sup_K |f|< C\},
	\] for all
	compacts $K\subset M$ and $C>0$. 
	
	A sequence $\{f_i\}$ of functions converges
	to $f$ if it converges to $f$ uniformly on all compacts.

\hfill
	
\remark \label{convc0}
	When $M$ is locally compact, 
	any sequence of continuous functions converging in
	$C^0$ topology \index[terms]{topology!$C^0$} {converges to a continuous function.}

\hfill

\remark  Similarly, one defines 
	{the  $ C^0$-topology on the space of sections of a bundle.}

\subsection{The $ C^1$ - topology on spaces of sections}\index[terms]{topology!$C^1$}

\definition 
Let $B$ be a vector bundle on a smooth manifold
$M$, and $\nabla:\; B \arrow B\otimes \Lambda^1 M$ a connection. 
The {\bf $ C^1$-topology} on the space of sections of 
$B$ (denoted, as usual, by the same letter $B$)
is the one whose sub-base of open sets is given by 
$ C^0$-open sets on $B$ and $\nabla^{-1}(W)$,
where $W$ is an open set in $ C^0$-topology in
$B\otimes \Lambda^1 M$.

\hfill

\remark 
A sequence $\{f_i\}$ {\bf converges in the $ C^1$-topology} if it converges 
uniformly on all compacts, and the first derivatives $\{f_i'\}$ also converge
uniformly on all compacts. {This can be seen as an equivalent definition
	of the  $C^1$-topology.}\index[terms]{topology!$C^1$}

\subsection{Montel theorem for normal families}

\theorem  (Montel, \cite[Lemma 1.4]{_Wu:Montel_}).\label{montel} 
Let $M$ be a complex manifold and ${\cal F}\subset H^0(\calo_M)$ 
a normal family of functions. Denote by $\bar{\cal F}$\index[terms]{normal family}
its closure in the $C^0$-topology.\index[terms]{topology!$C^0$} {Then $\bar{\cal F}$
	is compact and contained in $H^0(\calo_M)$.}\index[terms]{theorem!Montel}

\hfill

\proof 
Let $\{f_i\}$ be a sequence of functions in ${\cal F}$.
By Tychonoff's theorem\index[terms]{theorem!Tychonoff}, for each compact $K$,
there is a subsequence of $\{f_i\}$
which converges pointwise on a dense countable
subset $Z\subset K$. Taking a diagonal
subsequence, we find  a subsequence $\{f_{p_i}\}\subset \{f_i\}$
which converges pointwise on a dense 
countable subset $Z\subset M$. By \ref{_f'_uniformly_bounded_Lemma_}, 
$\{|f'_i|\}$ is uniformly bounded on compacts.
Therefore, the limit $f:= \lim_i f_i$ is  Lipschitz \index[terms]{function!Lipschitz}
on all compact subsets of $M$. Then it is 
	continuous, because a pointwise limit of Lipschitz 
functions is again Lipschitz.

Since $\{|f'_i|\}$ is uniformly bounded on compacts, 
after passing to a subsequence,  
we can assume that $f'_i$ also converges pointwise in $Z$,
and $f:= \lim_i f_i$ is smooth. A limit
of complex linear operators is complex linear, and hence  
$Df$ is complex linear, and $f$ is holomorphic.
Therefore, $\bar{\cal F}\cap H^0(\calo_M)$ 
	is compact.  \endproof


\section[The Stein completion of the K\"ahler cover ]{The Stein completion of the K\"ahler cover for LCK
  manifolds with potential}\index[terms]{completion!Stein}\index[terms]{manifold!LCK!with potential}


Compact LCK manifolds with proper potential do not have
conical K\"ahler covers, as the Vaisman manifolds do,
because the holomorphic homothety is missing.
Still, the covering is very close to a cone.\index[terms]{manifold!Vaisman} 

\hfill

\theorem { (\cite{ov_lckpot})}\label{potcon}
Let $M$ be an LCK manifold with proper potential,
and $\tilde M$ its K\"ahler $\Z$-covering\index[terms]{cover!K\"ahler $\Z$-}.
If  $\dim_\C M\geq 3$, then the metric completion 
$\tilde M_c$\index[terms]{completion!metric}
admits a structure of a complex variety, 
compatible with the complex structure on
$\tilde M \subset\tilde M_c$, and 
the complement $\tilde M_c\setminus  \tilde M$
is just one point. Moreover, $\tilde M_c$ is 
Stein.\index[terms]{variety!Stein}

\hfill

The rest of this section is devoted to the proof of \ref{potcon}.

\hfill

\claim \label{claim1}
Let $\tilde M$ be a K\"ahler $\Z$-cover \index[terms]{cover!K\"ahler $\Z$-}
of an LCK manifold with potential,\index[terms]{manifold!LCK!with potential} and
$\tilde M_c$ the metric completion of $\tilde M$. 
The complement $\tilde M_c\backslash \tilde M$ is just one
	point, called {\bf the origin}.\footnote{This claim
can be stated and proven, in a more general setting, for 
Riemannian manifolds; see Exercise \ref{_metric_completion_Exercise_}.}

\hfill

\proof 
Consider a Cauchy sequence $\{z_i\}$ in $\tilde M$ not converging to a point in $\tilde M$.
Let $\phi:\; \tilde M \arrow \R^{>0}$ be the LCK potential,\index[terms]{potential!LCK}
and $M_n:=\phi^{-1}([\lambda^n,\lambda^{n+1}[)$ its fundamental domain.
The closure of $M_n$ is compact, because $\phi$ is proper
(\ref{proper_equiv_Z}).
Clearly, the distance between two fundamental domains
\[
M_n=\gamma^{n}\phi^{-1}([1,\lambda[)=\phi^{-1}([\lambda^n,\lambda^{n+1}[)
\] and 
\[ M_{n+k+2}=\gamma^{n+k+2}\phi^{-1}([1,\lambda[) \]
satisfies
\begin{equation}\label{star}
\sum_{i=0}^{[k/2]} \lambda^{n+2i} v\leq
 d(M_n,M_{n+k+2})\leq \sum_{i=0}^k \lambda^{n+i} v, 
\end{equation}
 where $v$ is the distance between $M_0$ and $M_2$. 
This implies that for any Cauchy sequence $\{z_i \in M_{n_i}\}$,
the sequence $\{n_i\}$ converges to $-\infty$, or 
all $z_i$ except finitely many 
belong to the set $M_p \cup M_{p-1}$ for some $p$. 
The second case is irrelevant, because each $M_p$ is compact,
and in this case $\{z_i\}$ converges to a point in $M_p$ or $M_{p-1}$.
By \eqref{star}, any sequence $\{w_i \in M_{n_i}\}$ such that
$\lim_i n_i = - \infty$ is a Cauchy sequence,
and all such sequences are equivalent. Therefore, they
converge to the same point in the metric
	completion.
\endproof

\hfill

We now return to the proof of \ref{potcon}.

\hfill

{\bf  Step 1:} 
The boundary of $\phi^{-1}([a, \infty[)$
is strictly pseudoconvex, because $\phi$ is
strictly plurisubharmonic (\ref{_psh_then_pseudoconvex_Corollary_}).
Applying the  Rossi and  Andreotti--Siu 
theorem \index[terms]{theorem!Rossi, Andreotti--Siu}
(\ref{asr}) to 
$\phi^{-1}(]0, a])$, we obtain a Stein variety
$\hat M_c^a$ containing $\phi^{-1}(]0, a])$.
Since $\hat M_c^a$ contains $\phi^{-1}(]0, a_1])$
for any $a_1<a$, and the  Rossi and Andreotti-\index[persons]{Siu, Y.-T.}Siu variety
is unique, one has a natural embedding
$\hat M_c^{a_1}\hookrightarrow \hat M_c^{a}$, for any $a_1<a$.
Let $\hat M_c:=\cup_{a\in \R^{>0}}\hat M_c^a$ be the union of
all $\hat M_c^a$.

By \ref{_union_Stein_Proposition_} this is
a Stein variety with isolated singularities.

To prove \ref{potcon}, it remains 
to identify $\hat M_c$ 
	with the metric completion of $\tilde M$.
By \ref{claim1}, {this is equivalent to the complement 
	$\hat M_c\backslash \tilde M$ being a singleton.}

\hfill

{\bf  Step 2:} 
By \ref{asr} again, any holomorphic function
 $\phi$ can be extended from $\tilde M$ to $\hat M_c$
uniquely. Indeed, $\phi$ is CR-holomorphic on the boundary
of $\tilde M_a$, and any CR-holomorphic function on this
boundary can be extended uniquely to $\tilde M_a$.
Therefore, the monodromy group\index[terms]{group!monodromy} $\Gamma=\Z$ maps
holomorphic functions on $\hat M_c$ to holomorphic
functions. This implies that
 $\Gamma$ acts on $\hat M_c$ by holomorphic
	automorphisms.

\hfill

{\bf  Step 3:} Denote by $\gamma$ the generator of $\Gamma$ which 
decreases the metric by a factor 
$\lambda<1$, and let $\hat M^a_c\subset \hat M_c$ 
be  as above. Since
$\gamma(\hat M^a_c)=\hat M^{\lambda a}_c$, for any
holomorphic function $f$ on  $\hat M_c$,
one has 
\[ 
\sup_{z\in\hat M^a_c} |f(\gamma^i(z))|= 
\sup_{z\in\hat M^{\lambda^i a}_c} |f(z)|\leq
\sup_{z\in\hat M^{a}_c} |f(z)|.
\]
{Therefore, $\{(\gamma^n)^*f\}$ is a normal family.}\index[terms]{normal family}

\hfill

 {\bf  Step 4:}
Let $f_{\sf lim}(z)$ be a limit point of the sequence $\{f(\gamma^n(z))\}$.
Since the sequence 
$t_i:=\sup_{z\in\hat M^{\lambda^i a}_c} |f(z)|$
is non-increasing, it converges, and $\sup_{z\in\hat M^a_c} f_{\sf lim} = \lim t_i$.
Similarly, $\sup_{z\in\hat M^{\lambda a}_c} f_{\sf lim} = \lim t_i$.
By the strong maximum principle,
\index[terms]{maximum principle} \ref{hopf_theorem}, {a non-constant holomorphic function 
	on a complex manifold with boundary cannot
	have local maxima (even non-strict) outside of  the boundary.}
Since $\hat M^{\lambda a}_c$ does not intersect the boundary of
$\hat M^a_c$, the function $f_{\sf lim}$ must be constant.

\hfill

{\bf  Step 5:} Consider now the complement
$Z:=\hat M_c\backslash \tilde M$, and suppose it has two distinct points
$x$ and $y$. Let $f$ be a holomorphic function
that satisfies $f(x)\neq f(y)$. Replacing $f$ by an exponent of
$\mu f$, $\mu\in \C$, if necessary, we may assume that $|f(x)|< |f(y)|$.
Since $\gamma$ fixes $Z$, that is  compact, {for any limit $f_{\sf lim}$ 
	of the sequence $\{f(\gamma^n(z))\}$, the supremum\index[terms]{supremum}
	$f_{\sf lim}$ on $Z$ is not equal to the infimum of $f_{\sf lim}$ on $Z$.}
This is impossible, and hence  $f=\const$ on $Z$, and $Z$ is one point.

This finishes the proof of \ref{potcon}.
\endproof

\hfill

\remark
Steps 3--4 above
give a special case of the following
general observation. Let 
$f:\; Z \arrow Z$ be a holomorphic 
self-map of a complex variety
such that for some precompact open subset $U\subset Z$
and each compact $K\subset Z$
there exists an iteration $f^N$
satisfying $f^N(K) \subset U$.
In this case, we say that {\bf $f$ 
contracts $Z$ to $U$.}\index[terms]{contraction}
Then for any holomorphic function
$u\in H^0(\calo_Z)$, the limit
$\lim_n (f^n)^* u$ is constant.

Clearly, the generator of the $\Z$-action\index[terms]{action!$\Z$-}
on $\hat M_c$ contracts $\hat M_c$
to $\hat M^{\lambda a}_c$. 
Using Montel theorem as above, one can prove
$\lim_n (f^n)^* u=\const$  for all contractions.

\hfill

\corollary\label{_gamma_on_T^*_contraction_Corollary_}
Let $M$ be an LCK manifold with potential,\index[terms]{manifold!LCK!with potential}
$\tilde M_c$ the Stein completion of its $\Z$-cover\index[terms]{cover!K\"ahler $\Z$-}\index[terms]{completion!Stein}
constructed in \ref{potcon}, and $\gamma$ 
the generator of the $\Z$-action\index[terms]{action!$\Z$-} acting
on $\tilde M_c$ by holomorphic contractions.
Let $\gamma^*\in\Aut(A)$ 
be the action of $\gamma$ on the ring of holomorphic functions
$A:=H^0(\calo_{\tilde M_c})$.
Then $\gamma^*$ acts on the 
Zariski cotangent space $T_c^* \tilde M_c$
with eigenvalues $|\alpha_i|< 1$.

\hfill

\proof
Let $\goth m\subset A$ be the maximal ideal of the origin 
$c\in \tilde M_c$. By definition, 
$T_c^* \tilde M_c = \frac{\goth m}{{\goth m}^2}$.
Let $u\in \frac{\goth m}{{\goth m}^2}$ be an eigenvector
of $\gamma^*$, with eigenvalue $\alpha$ and $f\in \goth m$ its representative.
Since $f(c)=0$, we have $\lim_N(\gamma^*)^N f=0$
(\ref{potcon}, Step 4). Then 
$\lim_N(\gamma^*)^Nu= \lim_N \alpha^N u=0$,
hence $|\alpha|<1$.
\endproof

\hfill

\remark 
Let $\phi:\; Z \arrow Z$ be a continuous map 
fixing a point $x\in Z$. We say that $\phi$ is
{\bf a contraction centred in $x$} if
for any open subset $U\ni x$ and any
compact $K \subset Z$, for a sufficiently
big iteration $\phi^n$ of $\phi$, we would
have $\phi^n(K) \subset U$. Clearly, any
metric homothety with scale factor $< 1$
is a contraction. Therefore, we can always
choose the generator of the $\Z$-action
in \ref{potcon} to be a holomorphic contraction
of the Stein variety $\tilde M_c$.

\section{Appendix 1: another construction of the Stein completion}

The proof of \ref{potcon} relies on the Rossi\index[terms]{theorem!Rossi, Andreotti--Siu} and Andreotti--Siu theorem,
that is  invalid in complex dimension 2. However, when applied to 
a K\"ahler $\Z$-cover of a Vaisman manifold, it can be proven directly
by means of projective algebraic geometry.
We will use this result in the study of Vaisman\index[terms]{manifold!Vaisman} and LCK structures on 
complex surfaces. 

\hfill

We start by recalling the Remmert--Stein theorem.\index[terms]{theorem!Remmert--Stein}

\hfill

\theorem\label{_Remmert_Stein_Theorem_}
(Remmert--Stein theorem)\\ 
Let  $B$ be a complex analytic variety, and $C \subset B$ and 
 $A \subset B \backslash C$
 complex analytic subvarieties. Assume that
all irreducible components $A_i$ of $A$ satisfy
$\dim A_i > \dim C$. Then the closure of $A$ is
complex analytic in $B$.

\proof
\cite[\S II.8.2]{demailly}.
\endproof

\hfill

\remark \label{_Stein_Remmert_completion_Remark_}
We will use \ref{_Remmert_Stein_Theorem_}
in one special case, that is  also the case used
in the proof of Chow's theorem\index[terms]{theorem!Chow}
 (\cite[Theorem II.8.10]{demailly}). Suppose that
$A \subset \C^n \backslash 0$ is a complex
analytic subset without 0-dimensional connected components; then
its closure $\bar A$ in $\C^n$ is complex analytic. A posteriori,
$\bar A$ is isomorphic to a weak Stein completion of
$A$, \ref{_weak_Stein_Definition_}.\index[terms]{completion!Stein!weak}

\hfill

\proposition (\cite{ov_ac}\label{_Stein_completion_cone_Proposition_}
Let $P$ be a projective orbifold, and $L$ 
an ample line bundle on $P$. Assume that 
the total space $\Tot^\circ(L)$\index[terms]{orbifold!projective}
of all non-zero vectors in $L$
is smooth. Then the Stein completion\index[terms]{completion!Stein}
of $\Tot^\circ(L)$ is obtained by adding a
point (called ``the origin'' or ``the apex'' elsewhere).\index[terms]{apex}

\hfill

\proof
 The space
$\Tot^\circ(L)$ can be interpreted
algebraically as follows. 
Consider the standard $\C^*$-action
on $\Tot^\circ(L)$ understood
as a principal $\C^*$-bundle.
A finite-dimensional $\C^*$-representation
is a direct sum of 1-di\-men\-si\-o\-nal irreducible representations,
with $\C^*$ acting by $\rho(t) = t^w$; the number $w$ is called
{\bf the weight} of the representation. Therefore, the
category of $\C^*$-representations is equivalent
to the category of graded vector spaces. A 
ring with $\C^*$-action is the same as a graded ring.
In particular, the ring of fibrewise polynomial
holomorphic functions on $\Tot^\circ(L)$ is graded,
with the functions of weight $w$ being polynomials of
degree $w$ on each $\C^*$-orbit. We identify
the space $(\calo_{\Tot^\circ(L)})_w=(\calo_{\Tot^\circ(L^*)})_w$ of functions
of degree $w$ with the space of sections of $L^{\otimes w}$.
The projective variety $P\subset {\Bbb P} (H^0(L^N)^*)$ 
is identified with the graded spectrum\index[terms]{spectrum!graded}
of the graded ring 
$\bigoplus_w H^0(X, L^{\otimes wN})$, whenever $L^N$ is very ample.
Then, the manifold $\Tot^\circ(L)$ is the preimage of $P$ under the
natural projection 
$H^0(L^N)^*\backslash 0 \arrow {\Bbb P} (H^0(L^N)^*)$.
The Stein completion of $\Tot^\circ(L)$\index[terms]{completion!Stein}
is its closure in $H^0(L^N)^*$, obtained by 
adding the zero (called ``the origin'' or 
``the apex'' elsewhere).\index[terms]{apex} It is
complex analytic by Remmert--Stein theorem\index[terms]{theorem!Remmert--Stein}
(\ref{_Stein_Remmert_completion_Remark_}), and normal by
\ref{proj_normality_Proposition_} below. \index[terms]{variety!normal}
\endproof

\section{Appendix 2: the proof of the Kodaira--Spencer stability theorem}
\label{_KS_stability_Section_}

In this book, we do not often require the reader to 
follow complicated complex geometry arguments, but now
we cannot avoid using one; this is why we present 
the proof of  Kodaira--Spencer stability in this Appendix.

\hfill

First, we state the $dd^c$-lemma for compact K\"ahler manifolds.
For its proof, the reader is referred to any of the textbook
on complex geometry (\cite{griha,demailly});
for some of its many applications, see
\cite{_DGMS:Formality_}.

\hfill

\theorem\label{_dd^c_lemma_Theorem_}
Let $\eta$ be a form on a compact K\"ahler manifold,
satisfying one of the following conditions.
\begin{enumerate}
\item $\eta$ is an exact $(p,q)$-form. 
\item $\eta$ is $d$-exact, $d^c$-closed. 
\item $\eta$ is $\bar\6$-exact, $\6$-closed. 
\end{enumerate}
Then $\eta\in \im dd^c$.
\endproof

\hfill

\theorem {\bf (Kodaira-\index[persons]{Spencer, D. C.}Spencer stability theorem)}\\
Let ${\cal X} \stackrel\pi\arrow B$ be a smooth
family of compact complex manifolds.
Assume that the fibre $X_z:= \pi^{-1}(z)$ is K\"ahler
(that is, it admits a K\"ahler structure).
Then there exists a neighbourhood
$W\ni z$ such that for each $y\in W$,
the fibre $X_y:= \pi^{-1}(y)$ is also K\"ahler.\index[terms]{theorem!Kodaira--Spencer stability}

\hfill

The idea of the proof is easy. Let $X_t$ be a family of complex
manifolds, with K\"ahler central fibre. A small deformation of
a K\"ahler form\index[terms]{form!K\"ahler} is also positive 
definite, but we need it to remain closed. To deform the K\"ahler
form to the neighbouring fibres, we first deform its cohomology class, 
in such a way that it can be represented by a closed (1,1)-form.
Then we apply the $dd^c$-lemma (\ref{_dd^c_lemma_Theorem_}), 
to deform the representative
of this class along with the class itself. 
To use $dd^c$-lemma\index[terms]{lemma!$dd^c$}
and to be able to deform the closed (1,1)-form, we need
to show that the Fr\"olicher (a.k.a. Hodge--de Rham)\index[terms]{spectral sequence!Hodge--de Rham-Fr\"olicher}
spectral sequence degenerates for small deformations
of the K\"ahler central fibre. This is implied by
semi-continuity of the cohomology (\cite[Theorem 2.3]{_Bell_Narasimhan_}).

\hfill

{\bf Proof. Step 1:} 
Consider the relative Fr\"olicher 
spectral sequence \index[terms]{spectral sequence!relative Fr\"olicher}
\begin{equation} \label{_relative_spectral_seq_Inequality_}
R^i\pi_*(\Omega^j_B{\cal X})\Rightarrow R^{i+j}\pi_*(\C_{\cal X})
\end{equation}
Here $R^\bullet\pi_*(\C_{\cal X})$ is the derived pushforward 
of a constant sheaf (that is, a graded local system over $B$ with the fibres
of grading $k$ in $y\in B$ identified with $k$-th cohomology of $X_y$; see also Subsection \ref{_Derived_functors_Subsection_}),
and the $E_2$ term $R^i\pi_*(\Omega^j_B{\cal X})$ is 
a coherent sheaf obtained as a derived direct image\index[terms]{sheaf!coherent}
of the fibrewise de Rham algebra 
$\Omega^j_B{\cal X}=\Omega^j({\cal X}/B)$.

It is a relative (over $B$) version of the usual
Fr\"olicher (Hodge--de Rham) spectral sequence 
\[H^i(\Omega^j M)\Rightarrow H^{i+j}(M,\C).\]
This spectral sequence gives the inequality
\begin{equation} \label{_spectral_seq_Inequality_}
\sum_{i+j=k}\dim H^i(\Omega^jX_y)\geq  \dim H^{i+j} (M,
\C) =\sum_{i+j=k}\dim
H^i(\Omega^jX_z).
\end{equation}
for each $y\in B$.

\hfill

{\bf Step 2:}
Since $X_z$ is K\"ahler, the Fr\"olicher spectral
sequence for $X_z$ degenerates in $E_2$, giving 
$\sum_{i+j=k}\dim H^i(\Omega^jX_z)= b_k(X_z)$.
By upper semi-continuity (\cite[Theorem 2.3]{_Bell_Narasimhan_}), 
\[ \sum_{i+j=k}\dim H^i(\Omega^jX_y)\leq \sum_{i+j=k}\dim H^i(\Omega^jX_z)\]
in a sufficiently small neighbourhood $U$ of $z$.
Comparing this with \eqref{_spectral_seq_Inequality_}, 
we find that the rank of
$H^i(\Omega^jX_y)$ is constant in $U$, and hence 
the inequality  is equality in $U$,
and the Fr\"olicher spectral sequence  for $X_y$
degenerates. The same argument proves that
$X_y$ satisfies the $dd^c$-lemma.

\hfill
 
{\bf Step 3:} We extend the K\"ahler form $g_z$ in the K\"ahler
fibre $X_z$ of ${\cal X}$ smoothly to a family $g_y$ of Hermitian forms,
$y\in B$. Let $I_y$ denote the
associated complex structure in the family ${\cal X}$, and let
$\bar\6_y:\; \Lambda^{p,q}(M, I_y) \arrow \Lambda^{p,q+1}(M, I_y)$ be
the corresponding Dolbeault operator.
 Consider the Bott--Chern cohomology 
$H^{p,q}_{BC}(M, I_y):= \frac {\ker \bar\6_y\cap \ker \6_y}{\im \6_y\bar\6_y}$
of $(M, I_y)$ (it is introduced in Section \ref{_Bott_Chern_Section_}
later in this book). We identify $H^{p,q}_{BC}(M, I_y)$
with the kernel of the corresponding Laplacian operator $\boxdot$
(\cite{sch}; see also Section \ref{_MN_BC_Intro_Section_}). 
As follows from Step 2, the dimension of $H^{p,q}_{BC}(M, I_y)$
is constant in a sufficiently small neighbourhood $U$ of $z\in B$.
To simplify notations, we replace $B$ with $U$, and assume
that $H^{p,q}_{BC}(M, I_y)$ has constant rank for all $y\in B$.
Then, the space $\ker \boxdot$ smoothly depends on $y\in U$.
Let ${\cal H}$ be the $\ker \boxdot$ considered to be  a 
vector bundle over $U$. We identify its fibre in $y\in B$ with
the Bott--Chern cohomology of $(M, I_y)$.

Let $\Lambda^{1,1}_{cl}({\cal X}/B)$ be the sheaf
of  closed fibrewise (1,1)-forms on ${\cal X}$,
and \[ \pi_*\Lambda^{1,1}_{cl}({\cal X}/B)\arrow {\cal H}\]
the natural projection sending a form to its Bott--Chern 
cohomology class. The kernel ${\cal L}$ of this projection
is the sheaf of fibrewise exact (1,1)-forms.  
The fibre of ${\cal L}$ in $y\in B$ is identified with the 
$\im \6_y\bar\6_y= \frac{C^\infty X_y}{\const}$
(the kernel of $\6\bar\6$ on a compact complex manifold is the space of constant
functions, as follows from the maximum principle, \ref{hopf_theorem}).\index[terms]{maximum principle}
We obtain that ${\cal L}$ is also locally trivial.
This implies that the sheaf $\pi_*\Lambda^{1,1}_{cl}({\cal X}/B)= {\cal L}\oplus {\cal H}$
is a locally trivial infinite--dimensional vector bundle on $B$.

\hfill

{\bf Step 4:} Let $\tilde \omega$ be a smooth section of 
$\pi_*\Lambda^{1,1}_{cl}({\cal X}/B)$
satisfying $\tilde\omega\restrict z=[\omega_z]$.
Then $\tilde \omega$ is a family of
closed $(1,1)$-forms $\omega_y\in \Lambda^{1,1}_{cl}(X_y)$,
depending smoothly on $y\in B$.
Since all eigenvalues of $\omega_z$ 
are positive, the same is true for
$\omega_y$ for $y$ sufficiently close to $z$.
However, a closed, positive $(1,1)$-form is K\"ahler.
\endproof

\section{Notes}

\begin{enumerate}
	
\item Examples of LCK metrics with potential\index[terms]{metric!LCK!with potential} appeared in the
  literature long before this notion was formally
  introduced and the actual name was coined. Indeed,
  potentials appear implicitly in Vaisman's papers (e.g.   
  \cite{va_gd}) in the construction of LCK metrics on
  Vaisman manifolds of the type $(\C^n\backslash
  0)/\langle A\rangle$ where $A$ is the transformation
  $(z_i)\mapsto \lambda (z_i)$, $\lambda\in\C$,
  $|\lambda|\neq 0, 1$. Also, for the construction of
  Vaisman metrics on general diagonal Hopf
  surfaces, in  \cite{go}  the existence of an
  automorphic (in our sense) plurisubharmonic function
  $\f$ on $\C^2\backslash 0$ is proven, defined by an implicit
  equation with unique solution. Then the (pull back of
  the) LCK metric was written as $\frac
  14\f^{-1}dd^c{\f}$. On the other hand, in
  \cite{_Verbitsky:Vanishing_LCHK_}, the squared length of
  the Lee field,\index[terms]{Lee field}
  $\tilde\omega(\theta^\sharp,\theta^\sharp)$,  
is identified with a potential for $\tilde \omega$. 
The automorphicity of this potential
is not stated, but it is implicit in the presentation.

\item 
In \cite{_Brunella:Kato_}, 
M. \index[persons]{Brunella, M.} Brunella uses a similar construction to
obtain an LCK metric on \index[persons]{Kato, Ma.} Kato surfaces\index[terms]{surface!Kato}
(\ref{_Brunella_Theorem_}). To construct
this metric, he expresses a K\"ahler metric
on a blown-up annulus in $\C^2$ in terms
of a potential, that is  singular in the
blown-up points, and glues this 
singular K\"ahler potential to
its rescaled image under the $\Z$-action.\index[terms]{action!$\Z$-}
This way, he obtains an automorphic LCK
potential on the Kato surface with
exceptional divisors removed; after
gluing the exceptional divisors back,
this gives an LCK metric on the surface.
This construction was generalized to
Kato manifolds\index[terms]{manifold!Kato} of arbitrary dimension
in \cite{iop}.

\item
Montel theorem was interpreted in 
a more abstract way by \index[persons]{Grothendieck, A.} Grothendieck
(\cite{_Grothendieck:Montel_,_Golovin_}).
V. D. \index[persons]{Golovin, V. D.} Golovin defined {\bf the Montel--Fr\'echet sheaves}\index[terms]{sheaf!Montel--Fr\'echet}
as those sheaves of metrizable (that is, Fr\'echet) topological 
vector spaces such that the restriction map\index[terms]{topological vector space!metrizable}
is a compact operator (\ref{_compact_operator_Definition_}).\index[terms]{operator!compact}
In this language, Montel theorem\index[terms]{theorem!Montel}  
can be stated as ``the sheaf $\calo_M$ of
holomorphic functions with $C^0$-topology\index[terms]{topology!$C^0$}
is Montel'' (\ref{idcomp}). 

\index[persons]{Grothendieck, A.} Grothendieck proved that a 
 sheaf ${\cal F}$ of Fr\'echet--Montel\index[terms]{topological vector space!Fr\'echet-Montel} spaces on a compact space $M$
has finite-dimensional cohomology, assuming that
$M$ has a base of open sets $\{U_i\}$ such that
$H^{>0}(U_i, {\cal F})=0$.
He applied this result to reprove the  Cartan--Serre theorem \index[terms]{theorem!Cartan--Serre}
stating that coherent sheaves on compact complex manifolds have
finite-dimensional cohomology;\index[terms]{sheaf!coherent}
see Exercise \ref{bounded}, which 
is based on this idea.

We use a similar approach in Chapter 
\ref{_Embedding_for_LCK_with_potential:Chapter_} 
when we prove an embedding theorem 
for LCK manifolds with potential.\index[terms]{manifold!LCK!with potential}

\item Levi's problem on holomorphic convexity of manifolds
with strictly pseudoconvex boundary was solved by \index[persons]{Oka, K.} Oka (who called
it ``the inverse Hartogs problem'')
for domains in $\C^2$ in 1942 (\cite{_Oka:pseudoconvex_1942_})
and then in 1953 independently by \index[persons]{Oka, K.} Oka, \index[persons]{Bremermann, H. J.} Bremermann and \index[persons]{Norguet, F.} Norguet
for $\C^n$ (\cite{_Oka_,_Norguet_,_Bremermann:Levi_}).
For an excellent survey of the history of Levi problem, see 
\cite{_Siu:Levi_}. More recently, J. \index[persons]{Noguchi, J.} Noguchi published a survey of 
a series of unpublished papers of K. Oka, hand-written in Japanese
in 1943 \cite{_Noguchi:Oka_congress_,_Noguchi:Oka_arxiv_}. 
The originals, that were sent to T. Takagi, were lost,
but the complete set of first-draft manuscripts was found
posthumously. The first of these papers contains another solution
of Levi problem for domains in $\C^n$. 

\end{enumerate}

\section{Exercises}

\begin{enumerate}[label=\textbf{\thechapter.\arabic*}.,ref=\thechapter.\arabic{enumi}]

\item\label{_conjugate_of_Jordan_cells_Exercise_}
Let $\Mat(n, \C)$ denote the algebra of square matrices.
Let $A\in \Mat(n, \C)$ be a Jordan cell,\index[terms]{Jordan cell}
{\[\small 
A =\begin{pmatrix} \alpha & 1 & 0 & \ldots & 0\\
0 & \alpha & 1 & \ldots & 0\\
\vdots &\vdots &\vdots & \cdots & \vdots \\
0&0&0 & \ldots &1\\
0&0&0 & \ldots &\alpha
\end{pmatrix}.
\]}
Prove that this matrix is conjugate to
{\[\small
A_\lambda =\begin{pmatrix} \alpha & \lambda  & 0 & \ldots & 0\\
0 & \alpha & \lambda  & \ldots & 0\\
\vdots &\vdots &\vdots & \cdots & \vdots \\
0&0&0 & \ldots &\lambda \\
0&0&0 & \ldots &\alpha
\end{pmatrix}
\]}
for each $\lambda\in \C\backslash 0$.

\item\label{_closure_of_non_ss_orbit_contains_ss_Exercise_}
Let $A\in \Mat(n, \C)$
be a matrix (not necessarily diagonalizable).
Prove that there exists a sequence $\{g_i\} \in \GL(n, \C)$
such that $\lim_i g_i A g_i^{-1}=B$,
where $B$ is a diagonalizable matrix with the same eigenvalues.

{\em Hint:} Use the previous exercise.

\item \label{_LCK_pot_conf_invariant_Exercise_}
Let $(M, \omega)$ be an
LCK manifold with potential,\index[terms]{manifold!LCK!with potential} and
$(M, f\omega)$ a conformally equivalent manifold.
Prove that $(M, f\omega)$  is also LCK with potential.
Prove that $(M, f\omega)$ has proper potential if and
only if $(M, \omega)$ has.\index[terms]{potential!LCK!proper}

\item\label{_metric_completion_Exercise_}
Let $(M,g)$ be a compact Riemannian manifold, 
$(\tilde M, \tilde g)$ its $\Z$-cover\index[terms]{cover!K\"ahler $\Z$-},
and $\rho:\; \Z\times \tilde M \arrow \tilde M$ the 
corresponding $\Z$-action.\index[terms]{action!$\Z$-}

\begin{enumerate}
\item Consider the cohomology class
$\theta_\Z \in H^1(M, \Z)$ associated with
the homomorphism $\pi_1(M) \arrow \Z$ obtained from
the $\Z$-cover\index[terms]{cover!K\"ahler $\Z$-}, and let $\theta$ a closed 1-form
representing $\theta_\Z$ in $H^1(M, \R)$.
Prove that the pullback $\tilde \theta$ of $\theta$
to $\tilde M$ is exact, $\tilde \theta = df$.
Prove that it satisfies 
$\rho(n)^*f = f- n c$ for some non-zero 
number $c\in \R$.

\item Let $\lambda \in \R^{>1}$.
Show that there exist
smooth functions $\phi: \tilde M \arrow \R^{>0}$ 
satisfying $\rho(n)^*\phi = \lambda^n \phi$
for all $n\in \Z$. 

\item
Denote by $\tilde M_c$ the metric completion of \index[terms]{completion!metric}
the Riemannian manifold \\ $(\tilde M, \phi\tilde g)$.
Prove that the complement $\tilde M_c \backslash \tilde M$ is always
one point.
\end{enumerate}

\item
Let $M$ be a compact LCK manifold of LCK rank 1, and\index[terms]{rank!LCK}
$\tilde M$ its K\"ahler $\Z$-cover\index[terms]{cover!K\"ahler $\Z$-}.
Denote by $\tilde M_c$ the metric completion of $\tilde M$.
Prove that the complement $\tilde M_c \backslash \tilde M$ is always\index[terms]{completion!metric}
one point.

{\em Hint:} Use the previous exercise.

	
\item \label{_Hermitian_symplectic_Exercise_}
Call a complex manifold {\bf Hermitian symplectic}\index[terms]{manifold!Hermitian symplectic}
if it admits a symplectic form $\omega$ such that its
$(1,1)$-part is Hermitian.
Prove that a small deformation of an Hermitian
  symplectic manifold\index[terms]{manifold!symplectic!Hermitian} is again Hermitian symplectic.

\item  Prove that $C^1$-topology\index[terms]{topology!$C^1$} on the space of sections of a bundle $B$
is independent on  the choice of
	a connection in $B$.
	
\item Prove that the topological vector space  $C^1(M)$ of 
	1-differentiable functions on a manifold 
	is complete in $C^1$-topology.\index[terms]{topology!$C^1$}

\item {\bf A pseudometric} is a function\index[terms]{pseudometric}
$d:\; M \times M \arrow \R^{\geq 0}$ that satisfies
the triangle inequality, $d(x,x)=0$ and $d(x, y) = d(y, x)$\index[terms]{inequality!triangle}
for all $x,y\in M$. This is the same as the metric,
but without the positivity of $d(x,y)$ for $x\neq y$.
Let $\{d_\alpha\}$ be a family of pseudometrics
on a set. Prove that $d_\sup(x,y):= \sup_\alpha d_\alpha(x,y)$
is a pseudometric.

\item \label{_Kobayashi_pseudometric_Exercise_}
Let $M$ be a complex manifold.
Define {\bf the Kobayashi pseudometric} on $M$\index[terms]{pseudometric!Kobayashi}
as $\sup_\alpha d_\alpha(x,y)$, where
$\{d_\alpha\}$ is the set of all metrics
that satisfy \[ d_\alpha(f(x), f(y)) \leq d_p(x, y)\]
for all holomorphic maps $f:\; \Delta \arrow M$ 
from the Poincar\'e disk $\Delta$, where $d_p$ 
is the Poincar\'e distance on $\Delta$.
\begin{enumerate}
\item Prove that any holomorphic
map is 1-Lipschitz with respect to the\index[terms]{map!1-Lipschitz}
\index[persons]{Kobayashi, S.} Kobayashi pseudometric.\index[terms]{distance!Poincar\'e}

\item
Using Schwarz lemma, prove that the\index[terms]{lemma!Schwarz}
Kobayashi pseudometric on $\Delta$ is equal to $d_p$.

\item Prove that the Kobayashi pseudometric is continuous
as a function $M \times M \arrow \R$.\index[terms]{pseudometric!Kobayashi}
\end{enumerate}

\item Let $M$ be a compact complex manifold.
Prove that the \index[persons]{Kobayashi, S.} Kobayashi pseudometric is bounded from
above by another metric induced by a Riemannian
structure.

{\em Hint:} By compactness, cover $M$ with a finite number of open
balls, then estimate the Kobayashi pseudometric in a compact
subset of each ball in terms of an appropriate Riemannian metric.

\item Let $M$ be a complex manifold, and
$\pi:\; \tilde M \arrow M$ a covering map. Prove that
$\pi$ is locally an isometry with respect to the
Kobayashi pseudometric.

\item Let $X$ be a curve of genus $g>1$. Prove that
its Kobayashi pseudometric is induced by the
Poincar\'e metric on its universal cover.\index[terms]{metric!Poincar\'e}

\item Let $Y$ be a compact complex manifold with 
a Kobayashi pseudometric equal to a metric induced
by some Riemannian structure. 
\begin{enumerate}
\item Let $X$ be a compact complex manifold. Prove
that there exists a constant
$C>0$ such that $\sup_{x\in X}|dF_x| < C$ for any holomorphic
map $F:\; X \arrow Y$.

{\em Hint:} Choose on $X$ a Riemannian metric
that bounds the Kobayashi pseudometric.

\item Using Montel theorem, \index[terms]{theorem!Montel}
prove that the space of holomorphic maps
$F:\; X \arrow Y$ is compact with respect to $C^0$-topology.\index[terms]{topology!$C^0$}
\end{enumerate}

\item Let $X,Y$ be complex manifolds, and let 
$M(X,Y)$ be the space of all holomorphic maps from $X$ to $Y$,
with the $C^0$-topology.\index[terms]{topology!$C^0$}
\begin{enumerate}
\item Prove that $M(X, \C P^1)$ is not compact for any
compact complex curve $X$.
\item Let $X$ be a complex manifold, and $Y$ a compact complex
curve of genus $>1$. Prove that $M(X,Y)$ is compact.

\end{enumerate}

\item\label{bounded} 
Recall that  
a subset $K$ of a  topological vector space $V$ is called {\bf bounded} if for
any open neighbourhood $U\ni 0$, there exists
$a\in \R$ such that $aU\supset K$. 

A topological vector space $V$ is called {\bf Montel}
if any closed bounded subset $K\subset V$\index[terms]{topological vector space!Montel}
is compact.

Now let $M$ be a complex manifold, $B$ a vector bundle,
and $V=H^0(B)$ the space of holomorphic sections, 
with the topology of uniform convergence on compacts.
\begin{enumerate}
\item Prove that bounded subsets of normed spaces are subsets which
are contained in a ball of a sufficiently large radius.
\item Prove that $V$ is a Montel space.
\item Prove that any normed Montel space is finite-dimensional.
\item Show that $H^0(B)$ is always finite-dimensional
if $M$ is compact.
\end{enumerate}

\item Let $(M, \omega, \theta)$ be a Vaisman manifold,\index[terms]{manifold!Vaisman}
$(\tilde M, \tilde \omega)$ its K\"ahler cover, and $\phi\in C^{\infty}M$ its
potential,
$\phi= \tilde \omega(\theta^\sharp, I\theta^\sharp)$ 
(\ref{_potential_norm_Lee_Remark_}).
\begin{enumerate}
\item
Prove that $\phi^\lambda$ is an LCK potential \index[terms]{potential!LCK}for all $\lambda >0$.
\item Prove that the potential $\phi^\lambda$ defines a Vaisman
structure on $M$.
\item Prove that the corresponding LCK form $\omega_\lambda$ on\index[terms]{form!LCK}
$M$ can be written as 
$\omega_\lambda= \lambda d^c\theta+ \lambda^2 \theta\wedge \theta^c$.
\item Prove that the Lee form\index[terms]{form!Lee} of $\omega_\lambda$ is equal to
$\lambda\theta$.
\end{enumerate}

\item\label{_uniqueness_of_potential_Exercise_}
Let $(M, \theta, \omega)$ be a Vaisman
manifold, $\tilde M$ its K\"ahler $\Z$-cover\index[terms]{cover!K\"ahler $\Z$-},
and $\phi\in C^\infty \tilde M$ its LCK potential.\index[terms]{potential!LCK}
Prove that for all $\alpha>0$ the manifold
$(M, \alpha \theta, 
\frac{d_{\alpha\theta} d^c_{\alpha\theta}\phi^\alpha}{\phi^\alpha})$
is Vaisman. Assume that $\tilde M$ is simply
connected. Prove that for $\alpha$ sufficiently small,
the LCK potential for 
$(M, \alpha \theta,
\frac{d_{\alpha\theta} d^c_{\alpha\theta}\phi^\alpha}{\phi^\alpha})$  is unique 
up to a constant multiplier.

\item
Find an LCK manifold with proper potential
such that the LCK potential\index[terms]{potential!LCK!proper} with the preferred gauge
is unique.

{\em Hint:} Use the previous exercise.

\item
Let $(M, \theta, \omega)$ be an LCK manifold,\index[terms]{manifold!LCK}
and $\phi\in C^\infty M$.
Prove that $\omega_1 :=\omega + d_\theta d^c_\theta \phi$
is an LCK metric\index[terms]{metric!LCK} with the same Lee form,\index[terms]{form!Lee} if $\phi$
sufficiently small in $C^2$-norm. Prove that
$\omega_1$ is not conformally equivalent to $\omega$
for general $\phi$.

\end{enumerate}

\chapter{Embedding LCK manifolds with potential in Hopf
  manifolds}\index[terms]{manifold!LCK!with potential}
\label{_Embedding_for_LCK_with_potential:Chapter_}
{\setlength\epigraphwidth{.8\textwidth}
\epigraph{\it Tout mouvement cr\'eateur implique un rien de prostitution. Cela s'applique \`a Dieu, comme \`a quiconque est pourvu d'un talent quelconque. On ne devrait pas s'ext\'erioriser, si on veut rester pur. Rentrer en soi, en
toute rencontre, --  tel appara\^it le devoir de l'homme <<int\'erieur>>. L'autre,  l'<<ext\'erieur>>, ne compte gu\`ere : il fait partie de l'humanit\'e.}{\sc \scriptsize Emil Cioran,\ \ Cahiers (1957-1965)}}
 
\section{Introduction}

An important feature of LCK manifolds with potential
(\cite{ov_lckpot})\index[terms]{manifold!LCK!with potential}
is the existence of a ``logarithm'' of the monodromy action.
The logarithm was introduced in the original paper,\index[terms]{action!monodromy}
where the LCK manifolds with potential were defined,
but the concept turned out to be more delicate than
we expected.

Let $M$ be an LCK manifold with proper potential.\index[terms]{potential!LCK!proper}
From \ref{potcon} it follows that, in $\dim\geq 3$, the metric completion
of  its K\"ahler $\Z$-cover\index[terms]{cover!K\"ahler $\Z$-} $\tilde M$ is  a Stein variety\index[terms]{variety!Stein}
with a singular point in the origin.

In this chapter we prove that this Stein variety
is  in fact an algebraic cone over a projective
variety. Indeed, let $\gamma$ be the generator
of the monodromy $\Z$-action\index[terms]{action!$\Z$-} on $\tilde M$,
and let $\calo_{\tilde M}^\gamma$
be the ring of {\bf $\gamma$-finite functions},
that is, holomorphic functions\index[terms]{vector!$F$-finite}
$f\in \calo_{\tilde M}$ such that
the space $\langle f, \gamma^* f, (\gamma^2)^* f, ...\rangle$
is finite-dimensional. We prove that the ring
$\calo_{\tilde M}^\gamma$ is finitely generated, and
$\calo_{\tilde M}$ is its completion with respect
to the $C^0$-topology\index[terms]{topology!$C^0$} (\ref{_cone_cover_for_LCK_pot_Theorem_}).
This result will be used in Chapter \ref{_embe_theorem_Vaisman_Chapter_}
to prove that,  in $\dim\geq 3$, the complex structure on a compact
LCK manifold with potential \index[terms]{manifold!LCK!with potential} can be approximated
by complex structures of Vaisman type (\ref{def_lckpot2Vai}).

Consider the action of $\gamma$ on a finite-dimensional 
space generating the ring $\calo_{\tilde M}^\gamma$.
Since the exponential $\exp:\; {\goth{gl}}(n, \C) \arrow \GL(n, \C)$
is surjective, one would expect that $\gamma$ is
the exponential of a derivation 
$\delta:\; \calo_{\tilde M}^\gamma\arrow \calo_{\tilde M}^\gamma$,
that is  represented by an algebraic vector field
$\vec r \in T \tilde M$.

This argument is heuristically attractive, but it is false.
It breaks down because $\vec r$ must preserve the relations
between the generators of $\calo_{\tilde M}^\gamma$,
hence we need $\gamma$ to be the exponential of
an element of the Lie algebra of the group
$G$ of automorphisms of $\calo_{\tilde M}^\gamma$.
However, this group is not necessarily connected,
hence the logarithm does not always exist
(see \ref{_only_for_gamma^k_log_exists_Example_}
and Exercise \ref{_secondary_Kodaira_Exercise_}).

However, the set of connected components
of an algebraic group is finite, and hence 
for some power of $\gamma$ the logarithm
exists. We give an analytic proof of this
fact, independent on  the algebraization
(\ref{logar}). This proof is based on
the Riesz--Schauder theorem applied to the
Banach ring of bounded holomorphic functions.

In the Vaisman case, the existence
of the logarithm can be deduced from
the holomorphic isometry action generated by the
Lee field (\ref{defovai}, Step 1).\index[terms]{Lee field} 
In that case, the logarithm vector field
generates a conformal $\C$-action. 
However, an LCK manifold\index[terms]{manifold!LCK}
admitting such an action is {\em a posteriori}
Vaisman (\ref{kami_or}), and hence  for non-Vaisman
LCK manifold, the logarithm of monodromy \index[terms]{monodromy}
cannot generate a conformal $\C$-action\index[terms]{action!$\C$-}:
it is a conformal $\R$-action\index[terms]{action!$\R$-} in the best case.

The famous  Kodaira embedding theorem 
\index[terms]{theorem!Kodaira's embedding}
states that any compact K\"ahler manifold with a
rational K\"ahler class\index[terms]{class!K\"ahler} can be holomorphically
embedded to $\C P^n$. This theorem presents a geometer
with two major difficulties. First, it is
hard to find a rational K\"ahler class
on a given manifold $M$, unless $H^{2,0}(M)=0$
(in the latter case, rational classes are
dense in the K\"ahler cone\index[terms]{cone!K\"ahler}). Second,
we have no control over the dimension
of the model $\C P^n$.

In locally conformally K\"ahler geometry,\index[terms]{geometry!K\"ahler}
a similar theorem can be proven. Instead of
$\C P^n$, for the model target variety one takes
the linear Hopf manifold $\C^n \backslash 0 /\langle A\rangle$,\index[terms]{manifold!Hopf!linear}
where $A$ is a linear endomorphism with all eigenvalues $\al_i$ 
satisfying $|\alpha_i|>1$. Later on in this book, we prove that
all linear Hopf manifolds admit an LCK potential\index[terms]{potential!LCK}
(\ref{_linear_LCK_pot_Corollary_}).\footnote{This 
result was proven for diagonal Hopf surfaces in 
\cite{go}; this argument is generalized for \index[terms]{surface!Hopf!diagonal}
greater dimension in \cite[Theorem 3.1]{kor}.}
This implies that all complex submanifolds
in Hopf manifolds admit LCK potential.
The embedding theorem we prove here states,
conversely, that all LCK manifolds with
potential admit  holomorphic embeddings 
into Hopf manifolds.

To embed a compact K\"ahler manifold in $\C P^n$,
one needs to have a rational K\"ahler class;\index[terms]{class!K\"ahler}
in LCK geometry, this obstruction does not exist.
Indeed, a small deformation of an LCK manifold
with potential is again LCK with potential, and
can be embedded into a Hopf manifold.\index[terms]{manifold!LCK!with potential}

In this sense, one could say that LCK manifolds
with potential are ``more algebraic'' than
K\"ahler manifolds. Indeed, they can be described
algebraically as $\Z$-quotients of algebraic cones\index[terms]{cone!algebraic}
over projective varieties 
(\ref{_cone_cover_for_LCK_pot_Theorem_}).

However, we have even less control over
the target manifold than in the projective case.
Indeed, for a projective manifold, we lose
only control over the dimension of $\C P^n$,
and here we cannot control which Hopf manifold
will play the role of the target space.

The embedding result belongs to complex geometry.
The metric structures are auxiliary; indeed,
if an LCK structure\index[terms]{structure!LCK} on $M$
can be obtained by restricting an LCK structure 
from the Hopf manifold, it has LCK rank 1,\index[terms]{rank!LCK}
and the corresponding LCK potential\index[terms]{potential!LCK} is proper.
If an LCK structure does not have proper
potential, this LCK structure cannot be obtained 
as a restriction from the ambient Hopf manifold.\index[terms]{potential!LCK!proper}

The proof of our embedding theorem is, in a sense,
much simpler than the proof of  the Kodaira embedding theorem.
Indeed, the major problem we face (and solve)
is to construct the Stein completion of the\index[terms]{completion!Stein}
$\Z$-cover\index[terms]{cover!K\"ahler $\Z$-} of an LCK manifold with potential
(\ref{potcon}). \index[terms]{manifold!LCK!with potential}

After that, the main obstacle one faces is to 
find a vector field that integrates to the $\Z$-action\index[terms]{action!$\Z$-}
(``the logarithm'' of the $\Z$-action, \ref{logar}).
In \cite{ov_lckpot}, the logarithm was used to
construct the embedding. \index[terms]{logarithm!of a $\Z$-action}

The proof we give is different from the original in
\cite{ov_lckpot}. Instead of using the logarithm,
we apply the Riesz--Schauder theorem to the Banach algebra of 
bounded holomorphic functions on a Stein variety.
This argument requires several notions of
functional analysis that we recall in the next section
(see, for example \cite{_Bourbaki:TVS_}, 
\cite[Section 5.2]{friedman}).

\section{Preliminaries on functional analysis}

\subsection{The Banach space of holomorphic functions}\index[terms]{space!Banach}

We briefly recall some notions of functional analysis, that are 
standard.

\hfill

\definition
A {\bf Banach space}\index[terms]{space!Banach} is a topological vector space
that is  complete with respect to a norm. 
A {\bf Fr\'echet space}\index[terms]{space!Fr\'echet} is a locally convex topological vector space
that is  complete with respect to a translation-invariant metric.

\hfill

\remark
Let $C^0(X)$ be the space of continuous functions on a compact
space, with the norm $\|f\|_\sup:= \sup_{x\in X} |f(x)|$. It is not
hard to see that this space is Banach. Similarly,
if $X$ is a compact smooth manifold equipped with a connection $\nabla$,
the space $C^k(M)$ of $k$ times differentiable functions 
with the norm $\| f\|_{C^k}:=\sum_{i=0}^k \|\nabla^i(f)\|_\sup$
is also Banach.

\hfill

\remark
 {\bf A seminorm} on a vector space is a function
$\nu:\; V \arrow \R^{\geq 0}$ that satisfies the triangle inequality
$\nu(x+y) \leq \nu(x) + \nu(y)$ and $\nu(\lambda v) = |\lambda| \nu(v)$
for any $\lambda$ in the ground field (usually $\R$ or $\C$).
A seminorm is the same as a norm that is  allowed to take value
0 on non-zero vectors. Let $\{\|\cdot \|_i\}$ be a countable family of 
semi-norms on a vector space $V$.
Consider the topology with sub-base given by open balls in $\|\cdot \|_i$.
Suppose that for each $x$ there exists $i$ such that $\|x \|_i\neq 0$.
Then this topology is Hausdorff. This is the same topology as given by the
translation-invariant metric $d(x-y)= \sum_{i=1}^\infty \frac 1 {2^i} \max(1, \|x-y\|_i)$.
If this metric is complete, we say that {\bf the Fr\'echet topology\index[terms]{topology!Fr\'echet}
on $V$ is induced by the family of seminorms $\{\|\cdot \|_i\}$.}

\hfill

\example
The space $C^\infty M$ of smooth functions on a compact manifold
equipped with a family of norms $\|\cdot \|_{C^k}$
is Fr\'echet because each norm $\|\cdot \|_{C^k}$
is complete; however, it is not Banach.

\hfill

\theorem\label{_Banach_bounded_holo_Theorem_}
Let $M$ be a complex manifold, and let $H^0_b(\calo_M)$  
be the space of all bounded holomorphic functions, equipped
with  the sup-norm $|f|_{\sup} := \sup_M |f|$.
 Then $H^0_b(\calo_M)$ is a Banach space.\index[terms]{space!Banach}

\hfill

\proof 
Let $\{f_i\}\in H^0_b(\calo_M)$ be a Cauchy sequence in
the $\sup$-norm. {Then $\{f_i\}$ converges to a continuous
	function $f$} in the $\sup$-topology.\index[terms]{topology!$\sup$}

Since $\{f_i\}$ is a normal family\index[terms]{normal family} (see \ref{normal_family}), it
has a subsequence which converges in $ C^0$-topology
to $\tilde f\in H^0(\calo_M)$, by Montel \ref{montel}.\index[terms]{theorem!Montel} However, the {$ C^0$-topology
	is weaker than the $\sup$-topology,\index[terms]{topology!$\sup$} hence $\tilde f=f$.}
Therefore, $f$ is holomorphic. \endproof

\hfill

\corollary\label{_holo_Frechet_Corollary_}
Let $M$ be a complex manifold, not necessarily compact,
and $K_0\subset K_1 \subset ...$ a family of compact submanifolds with smooth 
boundary such that $M = \bigcup K_i$. 
Consider the family $\|\cdot \|_{K_i}$ of norms on $H^0(M, \calo_M)$,
with $\|f \|_{K_i}:= \sup_{x\in K_i} |f(x)|$. This family of norms defines a Fr\'echet topology\index[terms]{topology!Fr\'echet} on 
$H^0(M, \calo_M)$.\footnote{This topology on $H^0(M, \calo_M)$
is called {\bf the compact-open topology}, or {\bf the $C^0$-topology.}}\index[terms]{topology!$C^0$}

\hfill

\proof Let $\{ f_i\in  H^0(M, \calo_M)\}$ be
a sequence of holomorphic functions. Assume that $\{f_i\}$
is a Cauchy sequence in each norm $\|\cdot \|_{K_i}$.
To prove that $H^0(M, \calo_M)$ is Fr\'echet, it would suffice
to show that the sequence $\{f_i\}$ converges in the $C^0$-topology \index[terms]{topology!$C^0$}
to a holomorphic function. However, it converges to a holomorphic
function on each submanifold $K_m$, by Montel theorem
(\ref{_Banach_bounded_holo_Theorem_}), and 
$M = \bigcup K_m$. \endproof

\subsection{Compact operators}

Recall that 
a  subset $X$ of a topological space $Y$ is called {\bf 
	precompact},
or {\bf relatively compact in $Y$}, 
if its closure is compact. \index[terms]{subspace!precompact}

\hfill

\definition \label{_bounded_set_Definition_}
A subset
$K\subset V$ of a topological vector space
is called {\bf bounded} (Exercise \ref{bounded})  if for any 
open set $U\ni 0$, there exists a number $\lambda_U\in
\R^{>0}$ such that $\lambda_U K \subset U$.\index[terms]{subspace!bounded}

\hfill

\definition \label{_compact_operator_Definition_}
Let $V, W$ be topological vector spaces, A continuous operator $\phi:\; V
\arrow W$ is called {\bf  compact} if the image of
any bounded set is precompact.
\index[terms]{operator!precompact, compact, bounded}

\hfill

Montel \ref{montel} can be restated as follows:\index[terms]{theorem!Montel}

\hfill

\claim \label{idcomp}
Let $V=H^0(\calo_M)$ be the space of holomorphic functions on a complex
manifold $M$ with $C^0$-topology.\index[terms]{topology!$C^0$} Then any 
	bounded subset of $V$ is precompact. In this case, the
identity map is a compact operator.\index[terms]{operator!compact}

\hfill

\remark
It is not hard to deduce from Montel theorem \index[terms]{theorem!Montel}
that the space of bounded holomorphic functions on $M$
is Banach (that is, complete as a metric space) \index[terms]{space!metric}
with respect to the $\sup$-norm;
we leave it as an exercise to the reader.

\hfill

\remark  \label{_Riesz_theorem_Remark_}
According to the Riesz theorem (\cite[Chapter 1, Theorem 4]{_Diestel_}, 
Exercise \ref{_Riesz_Exercise_}), a closed ball in a normed
	vector space $V$ is never compact, unless $V$ is\index[terms]{theorem!Riesz}
finite-dimensional. This means that $(H^0(\calo_M),  C^0)$
does not admit a norm. A topological vector space in which any bounded subset is
precompact is called a {\bf Montel space}.\index[terms]{topological vector space!Montel}

\subsection{Holomorphic contractions}

\definition 
A {\bf contraction} of a topological space $X$ to a point
$x\in X$ is a continuous map $\phi:\; X \arrow X$ such that for any 
compact subset $K\subset X$ and any open set $U\ni x$, 
there exists $N>0$ such that for all
$n>N$, the map $\phi^n$ maps $K$ to $U$.\index[terms]{contraction!holomorphic}

\hfill

\theorem \label{_contra_compact_Theorem_}
Let $X$ be a complex variety, and 
$\gamma:\; X \arrow X$ be a holomorphic contraction
to $x\in X$ such that $\gamma(X)$ is precompact. 
Consider the Banach space\index[terms]{space!Banach} $V=H^0_b(\calo_X)$ of bounded holomorphic
functions with the sup-norm, and let $V_x\subset V$
be the space of all $v\in V$ vanishing in $x$. Then $\gamma^*:\; V \arrow V$
	is compact, and the eigenvalues of its restriction to 
	$V_x$ are strictly less than 1 in absolute value.\footnote{Since
$\gamma^*$ maps constants to constants identically, one cannot expect
that $\|\gamma^*\|< 1$ on $V$. However, if we add a condition
that excludes constants, such as $v(x)=0$, we immediately
obtain $\|\gamma^*\|< 1$.}

\hfill

\pstep For any $f\in H^0(\calo_X)$ we have
\[|\gamma^* f|_{\sup}= \sup_{m\in \overline{\gamma(X)}}
|f(m)|.
\]
This implies that $\gamma^*(f)$ is bounded.
Therefore, {for any sequence $\{f_i\in H^0(\calo_X)\}$ converging in the 
	$ C^0$-topology, the sequence $\{\gamma^* f_i\}$ converges
	in the $\sup$-to\-po\-lo\-gy.}\index[terms]{topology!$\sup$}
The set $B_C:=\{v\in V\ \ |
\ \  |v|_{\sup} \leq C\}$ is precompact in the 
$ C^0$-topology, because it is a normal family.\index[terms]{normal family}
Then $\gamma^* B_C$ is precompact in the 
$\sup$-topology.\footnote{The $C^0$-convergence for holomorphic
functions is strictly weaker that the $\sup$-convergence; see
Exercise \ref{_C_0_but_not_sup_Exercise_} for an example.} 
This proves that the operator $\gamma^*:\; V\arrow V$ is compact.

It remains to show that its operator norm is $<1$ on $V_x$.

\hfill

{\bf Step 2:}
Since $\sup_X |\gamma^* f|= \sup_{\gamma(X)} |f| \leq \sup_X |f|$,
one has $\|\gamma^*\|\leq 1$. If this inequality is not
strict, for some sequence $\{f_i\}$ of holomorphic
functions $f_i\in V_x$ with $\sup_X |f_i|\leq 1$ (that is, $f_1\in B_1$)
one has $\lim_i \sup_{m\in \gamma(X)} |f_i(m)|=1$.
Since $B_1$ is a normal family,\index[terms]{normal family} $f_i$ has a subsequence
converging in $ C^0$-topology to $f$. Then $\{\gamma(f_i)\}$
converges to $\gamma(f)$ in $\sup$-topology,\index[terms]{topology!$\sup$} giving 
$$\lim_i \sup_{m\in \gamma(X)} |f_i(m)|= \sup_{m\in  \gamma(X)} |f(m)|=1.$$
This supremum is realized somewhere on the closure 
$\overline{\gamma(X)}\subset X$, because $\overline{\gamma(X)}$ is compact.
Since $f(x)=0$, $f$ is non-constant.
By the maximum modulus principle,\index[terms]{maximum principle} 
a non-constant holomorphic function has no local maxima; this means
	that $|f(m)| >1$ somewhere on $X$. Then $f$ cannot
be the $C^0$-limit of $\{f_i\}\subset B_1$. \endproof

\subsection{The Riesz--Schauder theorem}

The following result  is a Banach analogue of 
the usual spectral theorem for compact operators\index[terms]{operator!compact} on
Hilbert spaces. It will be the central piece in our
argument.

\hfill

\theorem ({\bf Riesz--Schauder}, \cite[Section 5.2]{friedman})\\
\label{_Riesz_Schauder_main_Theorem_}
Let $F:\; V \arrow V$ be a compact operator\index[terms]{operator!compact} on a Banach space\index[terms]{space!Banach}.
Then for each non-zero $\mu \in \C$, there exists a sufficiently
large number $N\in \N$ such that for each $n>N$ {one has 
	$V= \ker(F-\mu\Id)^n \oplus \overline{\im (F-\mu\Id)^n}$,
	where $\overline{\im (F-\mu\Id)^n}$ is the closure of the image.}
Moreover, $\ker(F-\mu\Id)^n$ is finite-dimensional and independent on $n$.\index[terms]{theorem!Riesz--Schauder}\index[terms]{operator!compact}
\endproof

\hfill

\remark\label{_root_space_RS_Remark_}
Define  {\bf the root space of an operator \index[terms]{root space}
$F\in \End(V)$, associated with an eigenvalue $\mu$,}
as $\bigcup_{n\in \Z} \ker(F-\mu\Id)^n$. 
In the finite-dimensional case, the root spaces
coincide with the  Jordan cells of the corresponding 
matrix. Then 
\ref{_Riesz_Schauder_main_Theorem_} can be reformulated
by saying that any compact operator\index[terms]{operator!compact} $F\in \End(V)$ admits a 
Jordan cell decomposition, with $V$ identified with\index[terms]{Jordan cell}
a completed direct sum of the root spaces, that are 
all finite-dimensional; moreover, the eigenvalues
of $F$ converge to zero.

\hfill

We need the following corollary of the Riesz--Schauder theorem, that is 
obtained using the same arguments as in the finite-dimensional case.

\hfill

\theorem\label{rs}
Let  $F:\; V \arrow V$ be a compact operator\index[terms]{operator!compact} on a Banach space.\index[terms]{space!Banach}
We say that $v$ is {\bf a root vector}\index[terms]{vector!root} for $F$ if $v$ lies in 
a root space of $F$, for some eigenvalue $\mu\in \C$. 
Then the space generated by the root vectors is dense in $V$.
\endproof

\hfill

\definition\label{_F_finite_Definition_}
Let $F\in \End(V)$ be an endomorphism of a vector space.
A vector $v\in V$ is called {\bf $F$-finite}\index[terms]{vector!$F$-finite}
if the space generated by $v, F(v), F(F(v)),  ...$
is finite-dimensional. 

\hfill

\remark Clearly, a vector is 
$F$-finite if and only if it is obtained as a combination
of root vectors. Then  \ref{rs} 
implies the following.

\hfill

\corollary\label{_finite_RS_Corollary_}
 Let  $F:\; V \arrow V$ be a compact operator\index[terms]{operator!compact} on a Banach space,\index[terms]{space!Banach}
and $V_0\subset V$ the space of all $F$-finite vectors.
Then $V_0$ is dense in $V$.
\endproof


\section{The embedding theorem}


We recall the definition of what will be our model space:

\hfill

\definition 
Let $A\in \End(\C^n)$ be an invertible linear endomorphism
with all eigenvalues $|\alpha_i| <1$. The quotient
$H:= (\C^n \backslash 0)/\langle A\rangle$ is called
{\bf a linear Hopf manifold}.\index[terms]{manifold!Hopf!linear}

\hfill

We now state the central result of this chapter.

\hfill

\theorem \label{embedding}
Let $M$ be a compact LCK manifold with potential\index[terms]{manifold!LCK!with potential}, $\dim_\C M\geq 3$.
{Then $M$ admits a holomorphic embedding into a linear 
	Hopf manifold.}

\hfill

\ref{embedding} is implied by \ref{_gamma_finite_dense_Theorem_} below.
Let $F\in \End(V)$ be an endomorphism of a vector space.
Recall that a vector $v\in V$ is called {\bf $F$-finite}\index[terms]{vector!$F$-finite}
(\ref{_F_finite_Definition_})
if it belongs to a finite-dimensional $F$-invariant subspace.

\hfill

\theorem \label{_gamma_finite_dense_Theorem_}
	 Let $M$ be an LCK manifold with proper 
potential,\index[terms]{potential!LCK!proper} $\dim_\C M\geq 3$, and $\tilde M$ its K\"ahler $\Z$-covering.\index[terms]{cover!K\"ahler $\Z$-}
Consider the metric completion $\tilde M_c$ with its complex
structure and a contraction $\gamma:\; \tilde M_c \arrow \tilde M_c$ 
generating the
$\Z$-action.\index[terms]{action!$\Z$-} Let $H^0(\calo_{\tilde M_c})_\fin$ be the space of
functions that are  $\gamma^*$-finite.\index[terms]{vector!$F$-finite} {Then
	$H^0(\calo_{\tilde M_c})_\fin$ is dense in the $\sup$-topology\index[terms]{topology!$\sup$}
	on each compact subset of $\tilde
        M_c$.}\footnote{This is the same as the $C^0$-topology\index[terms]{topology!$C^0$}
  on $H^0(\calo_{\tilde M_c})$.}

\hfill

\remark
Further on, we prove \ref{embedding} in assumption that
$M$ has a proper LCK potential. This assumption is harmless,
because we can always replace an improper LCK potential\index[terms]{potential!LCK!proper}
with a proper one (\ref{defor_improper_to_proper})
without affecting the complex geometry
of the manifold.

\subsection{Density implies the embedding theorem}
\label{_density_implies_Subsection_} 

We deduce the embedding 
\ref{embedding} from 
\ref{_gamma_finite_dense_Theorem_} about the density of 
the space of $\gamma^*$-finite functions.
\index[terms]{function!$\gamma^*$-finite}
\hfill

{\bf Step 1:} 
Let $W\subset H^0(\calo_{\tilde M_c})_\fin$ be an
$m$-dimensional $\gamma^*$-invariant subspace $W$ with basis
$\{w_1, \ldots, w_m\}$. Then the following diagram is
commutative:
\begin{equation*}
\begin{CD}
\tilde M@>\Psi>> \C^m \\
@V{\gamma}VV  @VV{\gamma^*} V \\
\tilde M@>\Psi >>  \C^m
\end{CD}
\end{equation*}


where $\Psi(x)=(w_1(x), w_2(x), \ldots, w_m(x))$.

Suppose that the map $\Psi$ associated with a given 
$W\subset  H^0(\calo_{\tilde M_c})_\fin$ is injective
on a precompact fundamental domain of $\gamma$. Then the
	quotient map gives an embedding 
	$\Psi:\; \tilde M/\Z \arrow (\C^m\backslash 0)/\gamma^*$.
All eigenvalues of $\gamma^*$ are $<1$ because its operator
norm is $<1$, by \ref{_contra_compact_Theorem_}, and hence  the quotient
$ (\C^m\backslash 0)/\gamma^*$ is a Hopf manifold.

\hfill

{\bf Step 2:} To find an appropriate 
$W\subset  H^0(\calo_{\tilde M_c})_\fin$, choose a holomorphic embedding
$\Psi_1:\; \tilde M_c\hookrightarrow \C^n$, which exists because
$\tilde M_c$ is Stein. Let $\tilde w_1, \ldots, \tilde w_n$
be the coordinate functions on $\tilde M_c$ induced by $\Psi_1$. 
\ref{_gamma_finite_dense_Theorem_} allows one to approximate $\tilde w_i$ by 
$w_i \in H^0(\calo_{\tilde M_c})_\fin$ in $ C^0$-topology.
{Choosing $w_i$ sufficiently close to $\tilde w_i$
	in a precompact fundamental domain of the $\Z$-action\index[terms]{action!$\Z$-}, we obtain
	that $x\mapsto (w_1(x), w_2(x), ..., w_n(x))$ is injective
	in a precompact fundamental domain of the $\Z$ - action.} 

Finally,
take $W\subset H^0(\calo_{\tilde M_c})_\fin$ generated by 
the action of $\gamma^*$ on $w_1, \ldots, w_n$, and apply Step 1.
\endproof

\subsection{Density of $\gamma^*$-finite functions on
the minimal K\"ahler cover} \index[terms]{function!$\gamma^*$-finite}

Now we prove  \ref{_gamma_finite_dense_Theorem_} about the density of 
the space of $\gamma^*$-finite functions.\index[terms]{function!$\gamma^*$-finite}
Let $\tilde M$ be the minimal K\"ahler cover of $M$,
$\phi$ its K\"ahler potential, and $\tilde M_c$ the metric
completion, equipped with the natural complex structure as
in \ref{potcon}. Since $\phi$ is automorphic,
one has $\gamma^*(\phi)=\lambda \phi$, where
$\gamma$ is the contraction generating the deck transform
group of $\tilde M$, and $\lambda\in ]0,1[$.

Let $\tilde M_c^a:=\phi^{-1}([0, a])\subset \tilde M_c$,
and let $M_i:= M_c^{\lambda^{-i}}$,
where $i\in \Z$. Clearly, $M_i \subset M_{i+1}
\subset M_{i+2}\subset  ...$, with $\tilde M_c =\bigcup_i M_i$, and
$M_i\backslash M_{i-1}$ is
a fundamental domain of the action of $\Z=\langle
\gamma\rangle$, with $\gamma^k$ 
inducing an isomorphism $(\gamma^k)^*:\; M_i\arrow M_{i-k}$.

The operator $\gamma^*:\; H^0_b(\calo_{M_i}) \arrow
H^0_b(\calo_{M_{i}})$ is compact by \ref{_contra_compact_Theorem_}. Applying
the Riesz--Schauder theorem, we obtain a dense subspace\index[terms]{theorem!Riesz--Schauder}
\[ {\cal W}_i:=\bigoplus_{\mu\in \C} \ker
(\gamma^*-\mu\Id)^N\subset H^0_b(\calo_{M_i}). 
\]

Clearly, for any $f\in {\cal W}_i$, there is a
finite-dimensional space $W\ni f$ such that
$\gamma^*(W)= W$.  Let $A:\; W \arrow W$
be the linear map induced by $\gamma^*$. Consider the
function $A^{-k}(f)\in W$. Since
$(\gamma^k)^*A^{-k}(f)=f$, the function
$(\gamma^{-k})^*(f) \in H^0_b(\calo_{M_k})$
is equal to $f$ on $M_0\subset M_k$.
This gives a holomorphic extension of
$f$ to $M_k$ and to $\tilde M_c = \bigcup_i M_i$.

We have proven that ${\cal W}_i= H^0(\calo_{\widetilde
  M_c})_\fin\restrict {M_i}$.

Since ${\cal W}_i$ is $\sup$-dense in
$H^0_b(\calo_{M_k})$, it follows that the space $H^0(\calo_{\widetilde
  M_c})_\fin$ is $ C^0$-dense in $H^0(\calo_{\widetilde
  M_c})$. \endproof

\hfill

\remark The proof of the embedding theorem in $\dim_\C
M=2$ fails because the Rossi and Andreotti--Siu theorem fails\index[terms]{theorem!Rossi, Andreotti--Siu}
in this dimension. However, this theorem seems to be true
in  $\dim_\C M=2$. In Section \ref{_embedding_surfaces_pot_}, we deduce the embedding theorem
for $\dim_\C M=2$ from the GSS conjecture on complex
surfaces, widely accepted to be true.\index[terms]{conjecture!GSS}

\section{Notes} 
\begin{enumerate}
\item  The embedding theorems in LCK geometry\index[terms]{geometry!LCK} have been
extended to LCS geometry\index[terms]{geometry!LCS} in
\cite{torres} and
\cite[Theorem 1]{mmp}. They prove that a compact  LCS\index[terms]{manifold!LCS}
manifold  with $d_\theta$-exact LCS form and integral Lee\index[terms]{form!LCS}
class can be embedded into a Hopf manifold, in such a way
that the LCS form of the Hopf manifold is compatible with
that on M. This, moreover, implies a CR embedding result in
contact geometry,\index[terms]{geometry!contact} where the model space is a contact
sphere, see also Chapter \ref{sasemb}.

\item For Vaisman manifolds,\index[terms]{manifold!Vaisman} one can also give an isometric and holomorphic embedding into a Hopf manifold, see \ref{_Isometric_embedding_for Vaisman_Theorem_}. 
\end{enumerate}

\section{Exercises}

\begin{enumerate}[label=\textbf{\thechapter.\arabic*}.,ref=\thechapter.\arabic{enumi}]

\item
\begin{enumerate}
\item
Let $H:= \C^n\backslash 0 /\langle A\rangle$, $n\geq 2$, be a Hopf
manifold, where $A=\lambda \Id$, with $|\lambda|< 1$. 
Consider an automorphism $\phi:\; H \arrow H$.
Denote by $\tilde \phi$ the lift of $\phi$ to
$\C^n\backslash 0$. Prove that $\tilde \phi$ 
can be holomorphically extended to 0.

\item
Let $\lambda$ be a complex number, $0<|\lambda|< 1$,
and $\tilde \phi:\; \C^n \arrow \C^n$ be a holomorphic
map satisfying $\tilde\phi(\lambda z) = \lambda \tilde
\phi(z)$. Prove that $\tilde \phi$ is linear.

\item 
Prove that any holomorphic automorphism of 
$H=\C^n\backslash 0 /\langle \lambda \Id\rangle$,  $n\geq 2$
is induced by a linear automorphism of $\C^n$.
\end{enumerate}

\item
Let $H=\C^n\backslash 0 /\langle A\rangle$,
$n\geq 2$ be a linear Hopf manifold, with $A$
a linear map having all eigenvalues $0< |\alpha_i| <1$.
Prove that any holomorphic automorphism of 
$H$  is induced by a linear automorphism of $\C^n$
commuting with $A$.

{\em Hint:} Use the same logic as 
in the previous exercise.

\item Let $M$ be a complex manifold,
and $\phi:\; M \arrow M$ a holomorphic
self-map taking $M$ to a precompact subset.
Prove that $\phi$ has a fixed point which 
is unique.

\item\label{_attr_basin_Exercise_}
Let $\phi:\; M \arrow M$,  be a 
holomorphic diffeomorphism on a complex
manifold $M$. A fixed point $p$ of $\phi$ is called
{\bf attraction point} if $\|D\phi\|<1 $
at $p$. An {\bf attraction basin}
for $\phi$ is the set of all $x\in M$
such that $\lim_{n\to \infty} \phi^n(x)=p$.
Prove that the attraction basin $B$ of $\phi$ is biholomorphic
to $\C^n$.

{\em Hint:} Show that $B$ is a union of
a sequence of open balls.

\item
Consider the map\footnote{This is 
a special case of a {\bf H\'enon map}\index[terms]{map!H\'enon}, studied in complex dynamics.} 
$\Psi:\; \C^2 \arrow C^2$ given by 
$\Psi(x,y) = (\frac 1 2 y + x^2, \frac 1 2 x)$.
\begin{enumerate}
\item
Prove $\Psi$ is a diffeomorphism and 
$(0,0)$ is an attraction point for $\Psi$.

\item
Prove that $\Psi^2(x,y)=(\frac 1 4 x + yx^2 + \frac 1 4 y^2 + x^4,  
\frac 1 4 y + \frac 1 2 x^2)$.
Use this to show that there exists an open ball $B$
centred in $(100,100)\in \C^2$ such that 
all elements $(a, b) \in \Psi^{2n}(B)$ satisfy
$\Re(a) > 100$, for all $n \in \Z^{>0}$.
\end{enumerate}

\item Prove that there exists 
an open holomorphic embedding 
$\sigma:\; \C^n \arrow \C^n$, $n\geq 2$
such that $\C^n \backslash \im \sigma$ has a non-empty
interior.

{\em Hint:} Use the previous exercise.

\item
Let $\C \hookrightarrow \C$ be an open holomorphic embedding.
Prove that it is surjective.

\item
Let $X, Y$ be complex manifolds,
and $\Map(X,Y)$ the set of all holomorphic
maps. Denote by $W(K,U)\subset \Map(X,Y)$ the set of all
maps putting a given compact set $K\subset X$
in an open subset $U\subset Y$.
We equip $\Map(X,Y)$ with {\bf compact-open
  topology}, with the sub-base given by all
$W(K,U)$. 
\begin{enumerate}
\item Consider the following sequence
$\{f_n\}\in \Map(\C,\C)$: 
\[ f_n(z)=\epsilon\sum_{i=0}^n
  (\epsilon z)^i.
\]
 Prove that each open neighbourhood of 
$0\in \Map(\C,\C)$ contains all $\{f_n\}$, for $\epsilon$
sufficiently small.

\item Prove that $\{f_n\}$ is discrete in $\Map(\C, \C)$,
that is, it has no converging subsequences.

\item Prove that $\Map(\C,\C)$ is not locally compact.

\item Prove that $\Map(X,Y)$ is locally compact
for any compact complex manifolds $X,Y$.
\end{enumerate}

\end{enumerate}


\chapter{Logarithms and algebraic cones}

\epigraph{\it I must confess, contemplation of imbecility is my vice, my sin... Just like this: contemplation...Giulio Cesare Vanini, who was burned as an erethic, used to recognize God's greatness while contemplating a clod; others while contemplating the sky. I recognize it in the imbecile. Nothing is deeper, more abyssal, more dizzying, more untouchable... Only, one shouldn't contemplate it for too long.}{\sc\scriptsize Leonardo Sciascia, \ \ Todo modo}


\section{Introduction}

Locally conformally K\"ahler manifolds
with potential\index[terms]{manifold!LCK!with potential} are in a sense very ``algebraic''.
The difference between complex projective and K\"ahler geometry\index[terms]{geometry!K\"ahler}\index[terms]{geometry!complex projective}
is understood through the Kodaira embedding theorem:
a K\"ahler manifold with rational K\"ahler class\index[terms]{class!K\"ahler}
admits a complex embedding into a projective space.
In this setup, one can say that the geometry of
complex projective manifolds is part of the K\"ahler geometry,
and this part is certainly very algebraic.

The analogue of the Kodaira embedding theorem is already known:
an LCK manifolds with potential\index[terms]{manifold!LCK!with potential} is precisely a complex manifold
that admits a holomorphic embedding to a Hopf manifold.
This notion is already very close to algebraic geometry.\index[terms]{geometry!algebraic}
For example, there exists a $p$-adic version of the theory
of Hopf manifolds and their complex subvarieties,
within the framework of the rigid analytic spaces
(\cite{_Scholze:congress_}). The $p$-adic Hopf manifolds are known 
and well-studied (\cite{_Mustafin_,_Voskuil_}).

The passage from LCK geometry to complex algebraic geometry
is based on the notion of the ``algebraic cone''.
 In Chapter \ref{linear_hopf_shells},
we prove that an algebraic cone is an affine cone\index[terms]{cone!affine} over a 
projective orbifold.
Indeed, one can define algebraic cones as complex varieties
which underlie affine algebraic cones that are  non-singular
outside of  the origin (\ref{_cone_is_total_space_Proposition_}). In this chapter, we
define algebraic cones in a less explicit fashion.

Let $\tilde M_c$ be a Stein variety admitting
a holomorphic $\C^*$-action $\rho$, free outside of
a fixed point $c\in \tilde M_c$. Assume that $\tilde M_c$
is smooth outside of $c$. Suppose that
for all $\lambda\in \C^*$, $|\lambda|<1$, 
the action of $\rho(\lambda)$ contracts
$M$ to the unique fixed point $c$ of 
$\C^*$-action in $\tilde M_c$.
Then $(\tilde M_c, \rho)$ is called {\bf a closed algebraic cone}.
The complement $\tilde M= \tilde M_c \backslash \{c\}$ is called
{\bf an open algebraic cone}. In \ref{_closed_cone_normal_from_open_Corollary_},
we show that the open algebraic cone determines the
closed algebraic cone, up to normalization.\index[terms]{normalization}

Let $M$ be an LCK manifold with proper potential.\index[terms]{potential!LCK!proper}
From \ref{potcon} it follows that, in $\dim\geq 3$, the metric completion
of  its K\"ahler $\Z$-cover \index[terms]{cover!K\"ahler $\Z$-}$\tilde M$ is  a Stein variety
with a singular point in the origin.

In this chapter
we prove that this Stein variety
is in fact an algebraic cone.
Moreover, the algebraic structure
on $\tilde M_c$ is uniquely determined
by its complex analytic structure
(\ref{_same_algebra_structure_on_cone_Theorem_}, 
\ref{_cone_cover_for_LCK_pot_Theorem_}).

Let $\gamma$ be the generator
of the monodromy\index[terms]{monodromy} $\Z$-action\index[terms]{action!$\Z$-} on $\tilde M$,
and let $\calo_{\tilde M}^\gamma$
be the ring of {\bf $\gamma$-finite functions},\index[terms]{function!$\gamma^*$-finite}
that is, holomorphic functions
$f\in \calo_{\tilde M}$ such that
the space $\langle f, \gamma^* f, (\gamma^2)^* f, ...\rangle$
is finite-dimensional. We prove that the ring
$\calo_{\tilde M}^\gamma$ is finitely generated, and
that $\calo_{\tilde M}$ is its completion with respect
to the $C^0$-topology\index[terms]{topology!$C^0$} (\ref{_cone_cover_for_LCK_pot_Theorem_}).
Moreover, $\calo_{\tilde M}^\gamma$ admits an embedding
to the polynomial ring $\C[V]$, taking $\gamma$ to a linear
endomorphism $A\in \End(V)$ and the corresponding map
of complex manifolds $\tilde M \arrow \C^n \backslash 0$ induces
a holomorphic embedding from $M =  \tilde M/\Z$ to
the Hopf manifold $\C^n \backslash 0$.

Since the image of $\tilde M_c$ is defined by a graded ideal in $\C[V]$,
it admits a $\C^*$-action contracting $\tilde M_c$ to the origin.

This result is used in Chapter \ref{_embe_theorem_Vaisman_Chapter_}
to prove that,  in $\dim\geq 3$, the complex structure on a compact
LCK manifold with potential\index[terms]{manifold!LCK!with potential}  can be approximated
by complex structures of Vaisman type (\ref{def_lckpot2Vai}).

The main utility of the algebraic cones is
exposed in Chapter \ref{linear_hopf_shells}, where
we use the cones to prove that every Hopf manifold
(and, indeed, every $\Z$-quotient of an open 
algebraic cone, with $\Z$ acting by contractions
on the corresponding closed cone) admits an
LCK metric with potential. This is done 
by taking ``a pseudoconvex shell'' \index[terms]{pseudoconvex shell}(a pseudoconvex
hypersurface of the cone with certain properties
which mimics the properties of a Sasakian manifold
inside its cone), and solving a certain differential equation
that uses the logarithm of the monodromy action.\index[terms]{action!monodromy}

The logarithm was introduced 
in the original paper \cite{ov_lckpot},
where the LCK manifolds with potential\index[terms]{manifold!LCK!with potential} were defined,
but the concept turned out to be more delicate than
we expected.

For our purposes, {\bf a logarithm} of the 
monodromy action is a vector field $\vec r$ on
the K\"ahler $\Z$-cover\index[terms]{cover!K\"ahler $\Z$-} of an LCK manifold\index[terms]{manifold!LCK}
such that the automorphism $e^{\vec r}$ generates
the monodromy action.

Consider the action of $\gamma$ on a finite-dimensional 
space generating the ring $\calo_{\tilde M}^\gamma$.
Since the exponential $\exp:\; {\goth{gl}}(n, \C) \arrow \GL(n, \C)$
is surjective, one would expect that $\gamma$ is
the exponential of a derivation 
$\delta:\; \calo_{\tilde M}^\gamma\arrow \calo_{\tilde M}^\gamma$,
that is  represented by an algebraic vector field
$\vec r \in T \tilde M$.

This argument is heuristically attractive, but it is false.
It breaks down because $\vec r$ must preserve the relations
between the generators of $\calo_{\tilde M}^\gamma$,
hence we need $\gamma$ to be the exponential of
an element of the Lie algebra of the group
$G$ of automorphisms of $\calo_{\tilde M}^\gamma$.
However, this group is not necessarily connected,
hence the logarithm does not always exist
(see \ref{_only_for_gamma^k_log_exists_Example_}
and Exercise \ref{_secondary_Kodaira_Exercise_}).

However, the set of connected components
of an algebraic group is finite, and hence  the logarithm
exists for some power of $\gamma$. 
We give an analytic proof of this
fact, independent on  the algebraic structure 
on the cone $\tilde M_c$ (\ref{logar}). This proof is based on
the Riesz--Schauder theorem applied to the
Banach ring of bounded holomorphic functions.

In the Vaisman case, the existence
of the logarithm can be deduced from
the holomorphic isometry action generated by the
Lee field \index[terms]{Lee field}(\ref{defovai}, Step 1). 
In that case, the logarithm vector field
generates a conformal $\C$-action. 
However, an LCK manifold\index[terms]{manifold!LCK}
admitting such an action is {\em a posteriori}
Vaisman (\ref{kami_or}), and hence  for non-Vaisman
LCK manifolds, the logarithm of monodromy \index[terms]{action!monodromy}
cannot generate a conformal $\C$-action\index[terms]{action!$\C$-}:
it is a conformal $\R$-action\index[terms]{action!$\R$-} in the best case.

We finish this introduction by a brief mention of 
formal schemes and formal geometry; a reader
who does not know schemes can skip it. 
Let $J\subset A$ be an ideal of a ring,
$S$ the set of all prime ideals in $A/J$, equipped with
the Zariski topology. Each prime ideal in $A/J$
gives a prime ideal in the $J$-adic completion
$\hat A_J$. This allows one to define a sheaf
of rings $F$ over $S$, with $\hat A_J= F(S)$,
and appropriate localizations as $F(U)$ for
Zariski open subsets $U\subset S$.
The ringed space $(S, F)$ is called
{\bf the formal spectrum of $A$, completed at $J$}\index[terms]{spectrum!formal}

A {\bf formal scheme} is a ringed topological space
locally isomorphic to the formal spectrum.

{\bf An algebraization} of a formal scheme $S$
is a scheme $X$ such that $S$ is homeomorphic
to a closed subset $S_0\subset X$, and its
structure sheaf is identified with the $S_0$-adic
completion of $\calo_X$. 

Not every formal scheme admits an algebraization;
however, an algebraization can be obtained from a $\C^*$-action
on a formal scheme, acting with weights $|\alpha_i|>1$
in the formal direction. Then the algebraization
is obtained by taking $\C^*$-finite functions
in the structure sheaf. 

Applying this observation to the adic completion
of the closed algebraic cone in the origin, we obtain
that this ring is a completion of a finitely--generated
ring; this argument motivated the original proof of
\ref{_cone_cover_for_LCK_pot_Theorem_}.
This explains why we used the $\C^*$-action in the
definition of an algebraic cone, instead of
the $\C$-action.

This approach could be used to show that
an open algebraic cone is always a total space
of an ample $\C^*$-bundle, instead of using
the LCK metrics, as we did in\index[terms]{metric!LCK} \ref{_cone_is_total_space_Proposition_}.

\section{The logarithm of an automorphism}\label{logaritm}

\subsection{Logarithms of an automorphism of a Banach 
ring}\index[terms]{logarithm!of an automorphism}

We start with the construction of the logarithm of an
automorphism of a finitely generated Banach ring. It is a
general construction needed in several
applications.

\hfill

\definition 
A {\bf Banach ring} is a Banach space\index[terms]{space!Banach}
equipped with a commutative, associative, continuous bilinear
product map  satisfying $\|xy\|\leq \|x\|\|y\|$.
A Banach ring is {\bf finitely generated}
if it is the topological closure of a finitely generated ring
and {\bf finitely presented} when it is a completion
of a normed, finitely presented ring.\index[terms]{Banach ring}

\hfill

\newcommand{\I}{\goth{I}}
\example 
The  ring of bounded holomorphic functions with $\sup$-norm on a complex
variety is a Banach ring.

\hfill

\proposition \label{log}
Let $R$ be a finitely generated, finitely 
presented commutative Banach ring, and let 
$V\subset R$ denote a finite-dimensional subspace containing the 
unit, which generates $R$ multiplicatively. We write 
$R=\overline{\C[V]/\I}$,
where $\I$ is an ideal. Let
 $R_N=\bigoplus_{i=0}^N\Sym^i(V)$ 
be the subspace of $\C[V]$ generated by 
powers of $V$ of degree $\leq N$.
Since $R$ is finitely presented,
there exists a number $N$ such that $\I\cap R_N$ generates $\I$. 
Consider  an automorphism $A$ of $\C[V]$ preserving
$\I$ such that on $R_N$ one has $\|A-\Id\|< 1$, 
where $\|\cdot\|$ is the operator norm.
For each $x\in R_N$, define {\bf the logarithm} by:
\begin{equation}\label{_log_power_series_Equation_}
\log(A)(x) := \sum_{i=1}^\infty (-1)^i \frac{(\Id-A)^i}{i}(x)
\end{equation}
(the series converges, because $(\|A-\Id\|< 1$ on $R_N$). 
{Then $\log A$ can be extended to a derivation on $R$,
	that satisfies $e^{\log A}=\Id$.}

\hfill

\proof For each $x,y, xy \in R_N$, one has the Leibniz identity\index[terms]{Leibniz identity}
$\log(A)(xy)=\log(A)(x) y + x \log(A)(y)$ by standard formal
identities with logarithms (\cite[\S 6.1]{_Bourbaki:Lie_}).  
Since all relations
are generated by elements of $R_N\cap \I$, and $\log(A)$
preserves $\I$, the operator $\log(A)$ can be extended
to $R=\overline{\C[V]/\I}$ using  Leibniz identity.
\endproof

\subsection{Logarithms of the homothety of the cone}

\definition
Let $M$ be an LCK manifold\index[terms]{manifold!LCK}, $\tilde M$ its K\"ahler
cover, and $\gamma$ an element of $\pi_1(M)$. 
We say that a vector field $\vec r$ on $\tilde M$
is {\bf a logarithm of $\gamma^k$} if $e^{\vec r} = \gamma^k$.

\hfill 

\remark From 
\ref{_S^1_potential_Theorem_} it follows that the existence of a logarithm
of $\gamma\in \pi_1(M)$, when $\gamma$ acts
on $\tilde M$ by non-trivial homotheties,
 implies the existence of an LCK potential.\index[terms]{potential!LCK}

\hfill

It turns out that for an LCK manifold with proper potential\index[terms]{potential!LCK!proper}
and $\gamma$ the generator of the monodromy group\index[terms]{group!monodromy}, the logarithm of $\gamma^k$
always exists, for $k$ sufficiently large;
however, the logarithm of $\gamma$ does not
necessarily exist, even for a Vaisman manifold \index[terms]{manifold!Vaisman}
(\ref{_only_for_gamma^k_log_exists_Example_}). 

\hfill

\theorem {(\cite{ov_pams})}\label{logar} 
Let $M$ be an LCK manifold with proper potential, $\tilde M$
its K\"ahler  $\Z$-covering\index[terms]{cover!K\"ahler $\Z$-}, and $M=\widetilde
M/\langle\gamma\rangle$.
Then there exists $k\in \Z^{>0}$ and a holomorphic
	vector field $\vec r$ on $\tilde M$ such that
	$\gamma^k=e^{\vec r}$.

\hfill

\proof
Let $a\in \R^{>0}$. By
\ref{_Stein_comple_unique_Forster_Remark_}, 
the automorphisms of $\tilde M$ are the same as automorphisms of its
Stein completion\index[terms]{completion!Stein} $\tilde M_c$ preserving the origin.
We apply \ref{log} to a power of the automorphism
$\gamma\in \Aut(\tilde M)=\Aut(\tilde M_c)$ acting
on the Banach ring of holomorphic functions
on $\phi^{-1}([0, a])\subset \tilde M_c$, where
$\phi$ is the LCK potential.\index[terms]{potential!LCK}

\hfill

{\bf Step 1:} We start by defining the logarithm of $\gamma^k$
on the ring $\calo_{\tilde M_c}^\gamma$  
of $\gamma$-finite functions, continuous \index[terms]{function!$\gamma^*$-finite}
in the appropriate Banach norm. The choice of $k$ depends on the
eigenvalues of $\gamma$ on the set of multiplicative
generators of $\calo_{\tilde M_c}^\gamma$.

 As follows from \ref{_gamma_finite_dense_Theorem_}, 
the ring $\calo_{\tilde M_c}^\gamma$
is finitely generated and dense in
$H^0(\calo_{\tilde M_c})$. To define the action of the logarithm $\vec r$
on $\tilde M_c$, it would suffice to define a continuous
derivation $\log (\gamma^k)$ on the topological ring
$H^0(\calo^\gamma_{\tilde M_c})$ in such a way that
$e^{\log \gamma^k}=\gamma^k$. Identifying derivations
of $H^0(\calo_{\tilde M_c})$ and vector fields,
we obtain a vector field $\vec r$ on $\tilde M_c$
satisfying $\gamma^k=e^{\vec r}$.

We define $\log (\gamma^k)$ as the sum of a power series,
$\log(\gamma^k) =
\sum_{i=0}^\infty(-1)^i\frac{(1-\gamma^k)^i}{i}$
(to simplify the notation, we write $1$ for
the identity operator in the algebra of endomorphisms). 
To extend this series to $\calo_{\tilde M_c}$, we need
this series to converge absolutely in the operator norm
$\|\cdot \|$, that is  equivalent to
$\|1-\gamma^k\| < 1$. However, we only know
that $\|\gamma\|< 1$, because $\gamma$ is a holomorphic 
contraction (\ref{_gamma_on_T^*_contraction_Corollary_}). 
The open ball $\|\gamma\|< 1$ centred in 0
does not lie in the open ball $\|1-\gamma\| < 1$ centred
in $1=\Id\in \End(V)$, and hence  convergence of this
series is an additional condition on $\gamma$.

We prove that for any operator $\gamma\in \End(V)$
with $\|\gamma\|< 1$, there exists a power $k\in \Z^{>0}$ such that
$\|1-\gamma^k\| < 1$. This explains why we cannot
take the logarithm of $\gamma$, but only of $\gamma^k$.


Observe now that for any finite set of complex numbers $\{a_1, \ldots, a_n\}$ 
that satisfy $0<|a_i|<1$, there exists an integer $k>0$ such
	that $|a_i^k-1|< 1$. Indeed, let $a_i = b_i u_i$, where
	$|u_i|=1$, $b_i \in \R$. We consider $z:=(u_1, ..., u_n)$
as a point in the torus $T:= (S^1)^n$. Then the closure of the sequence
$\{z, z^2, z^3, ...\}$ in $T$ is a compact torus in $T$
(see, for example,    \cite[section 5]{_Wilkinson_}). Therefore, there 
is a subsequence of $\{z, z^2, z^3, ...\}$ converging to 1.
For any given $\epsilon$ one can 
	find $k\in \Z^{>0}$ such that $\sum_i |\arg(u_i^k)| <\epsilon$.
	Take $\epsilon := \frac \pi 3$.
Any operator $A$ with eigenvalues $\{a_1, \ldots, a_n\}$
with $0<|a_i| < 1$ and $|\Arg(a_i)| < \frac \pi 3$
satisfies $\|A-\Id\| < 1$.\footnote{Here $\Arg$ denotes {\bf the argument} of
a complex number $\lambda$, that is, a number  $u \in ]-\pi, \pi]$ such that
$\frac{\lambda}{|\lambda|}= e^{\1u}$.}

\centerline{\includegraphics[width=0.6\linewidth]{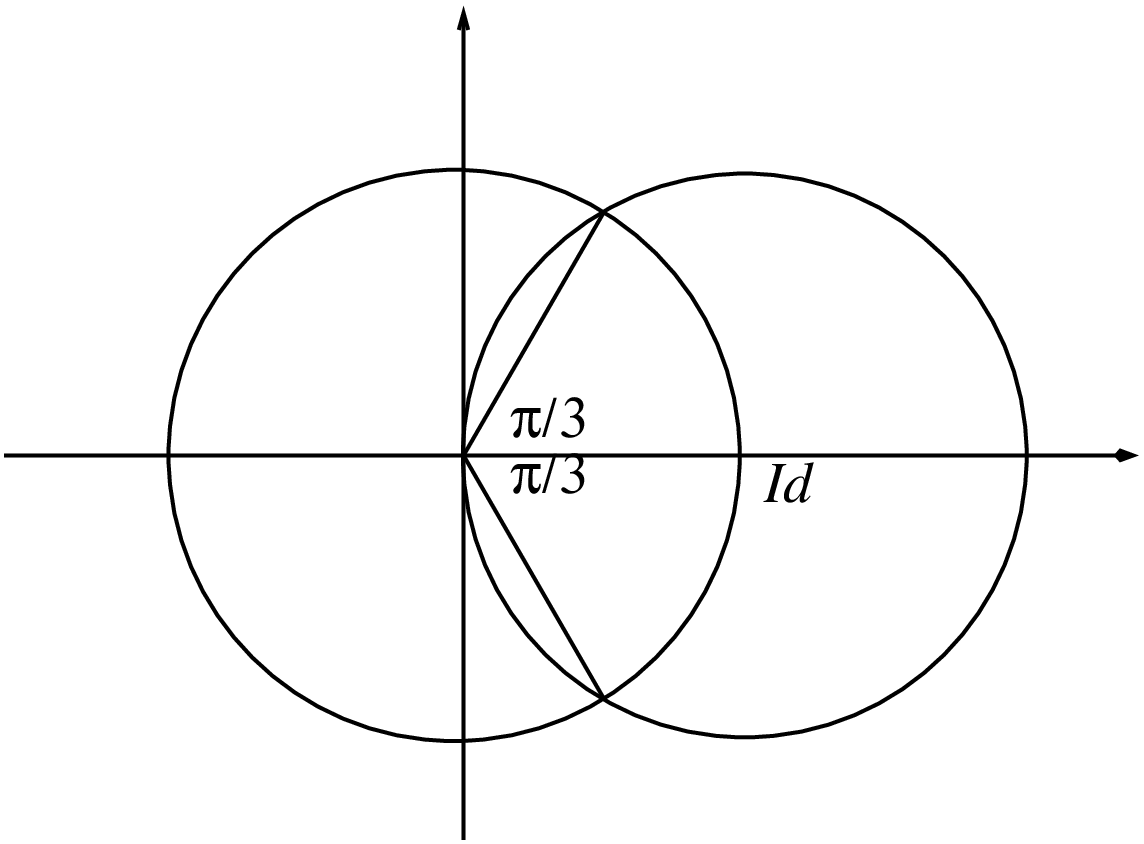}}
	
This allows one to choose $k\in \Z^{>0}$ such that on 
the space $V$ generating $\calo_{\tilde M_c}^\gamma$,
one has $|\gamma^k-\Id| <1$. Let $\log \gamma^k$ be the 
logarithm defined in \ref{log}. Then $e^{\log \gamma^k}=\gamma^k$,
hence we can take $\vec r :=\log \gamma^k$.  

\hfill

{\bf Step 2:}
In \ref{embedding} we constructed an embedding
of a K\"ahler $\Z$-cover\index[terms]{cover!K\"ahler $\Z$-} of an LCK manifold with potential\index[terms]{manifold!LCK!with potential} to $\C^n$.
Consider the closed embedding $u:\; \tilde M_c\arrow V^*$ taking $m \in \tilde M_c$ and
$v \in V$ to $v(m)$, where $V\subset \calo_{\tilde M_c}^\gamma$
is an appropriate space of generators. Clearly, $u$ is $\gamma$-equivariant.
Since the ideal of $\tilde M_c\subset V^*=\C^n$
is generated by $R_N \cap \I$, it is preserved
by $\log \gamma^k$. Since $\gamma$ is a linear 
endomorphism of $V^*$, the vector field
$\vec r =\log \gamma^k$ is linear on 
$V^*$. Then its restriction to 
$\tilde M_c\subset V^*$ is a holomorphic
vector field. For any function on 
$\tilde M_c$ obtained as a restriction of a
holomorphic function $f$ on $V^*$,
the function $\log \gamma^k(f)$
is well-defined on $\tilde M_c$,
and the map 
$f\restrict{\tilde M_c} \mapsto \log \gamma^k(f)\restrict{\tilde M_c}$
is continuous by construction.
\endproof

\hfill

\remark 
When $M$ is a Vaisman manifold, \index[terms]{manifold!Vaisman}
\ref{logar}  follows directly from 
\ref{_Vaisman_Lee_action_contains_monodromy_Corollary_}.

\hfill

\example\label{_only_for_gamma^k_log_exists_Example_}
The statement of \ref{logar}  is optimal: it is
possible to find a Vaisman manifold
$M= \tilde M /\langle \gamma\rangle$ such that
the logarithm of $\gamma$ does not exist.\index[terms]{logarithm!of $\gamma$}
For example, take an ample line\index[terms]{bundle!line!ample}
bundle $L$ over a curve $C$ of genus $g >1$,
and let $\sigma$ be a non-trivial automorphism of $C$ of order $k$.
We can assume that $\sigma$ is not homotopic
to a trivial map (\cite{_Eberlein_}).
Replacing $L$ with $\bigotimes_{i=0}^{k-1} \sigma^i(L)$,
we may assume that $L$ is $\sigma$-equivariant,
and the equivariant action is isometric with
respect to the Hermitian metric on $L$.
Then $\sigma$ extends to an isometry
of the total space $\Tot^\circ(L)$
of non-zero vectors in $L$. 
The space $\Tot^\circ(L)$
is a K\"ahler cone over a Sasakian
manifold (\ref{_quasireg_Sasakian_orbibundles_Theorem_} (iv)).
Let $\gamma$ map $v\in \Tot^\circ(L)$
to $\lambda \sigma(v)$, where $\lambda$ is
a complex number with $|\lambda|>1$.
Since $\sigma$ acts on $\Tot^\circ(L)$
by an isometry, the quotient
$M:= \Tot^\circ(L)/\langle \gamma\rangle$ is
Vaisman. If the logarithm $\vec r$ of $\gamma$ exists, it defines
a homotopy of $\gamma$ to the identity map;
indeed, $e^{0\cdot\vec r}=\Id$ and $e^{1\cdot\vec r}=\gamma$. Let
$S:= \{ v\in  \Tot^\circ(L)\ \|\ \ |v|=1\}$
be the Seifert manifold\index[terms]{manifold!Seifert} associated with $C$ and $L$;
it is clearly Sasakian.
The corresponding \index[terms]{exact sequence!Serre} Serre exact sequence gives
\[
0 \arrow \Z \arrow \pi_1(S) \arrow \pi_1(C) \arrow 0.
\]
A diffeomorphism of $C$ is homotopy equivalent to $\Id$
if and only if it acts trivially on $\pi_1(C)$, because
$C$ is 
$K(\pi_1(C), 1)$.\footnote{An \index[persons]{Eilenberg, S.} Eilenberg-\index[persons]{MacLane, S.}MacLane space\index[terms]{space!Eilenberg-MacLane}
$K(\Gamma, 1)$ is a CW-space $X$ with $\pi_1(X)=\Gamma$
and $\pi_i(X)=0$ for $i>1$. \index[persons]{Eilenberg, S.} Eilenberg-MacLane spaces are unique
up to a homotopy, and the group of homotopy auto-equivalences
of $K(\Gamma, 1)$ is $\Aut(\Gamma)$ (see, for example, 
\cite{_Fuks_Fomenko_}). A 
CW-space is $K(\Gamma, 1)$ if and only if its
universal cover is contractible.}
Therefore, for any automorphism of $C$
that is  not homotopy equivalent to 0,
the corresponding automorphism of $S$ is also not
homotopy equivalent to 0.
However, $\Tot^\circ(L)$ is homotopy equivalent
to $S$, and $\gamma$ acts on $S$ as $\sigma$,
that is  not homotopy equivalent to a trivial map.
Therefore, $\gamma$ does not have a logarithm.

This construction makes sense for a K\"ahler cone
over an arbitrary Sasakian manifold $S$. Consider an automorphism
$\phi$ of a compact Sasakian manifold $S$ not homotopy equivalent to zero,
and let $M:= \frac{C(S)}{(s, t) \sim (\phi(s), q t)}$
be the corresponding Vaisman manifold (\ref{str_vai}).\index[terms]{manifold!Vaisman}
Since the map $(s, t) \arrow (\phi(s), q t)$ is not
homotopy equivalent to identity, it has no logarithm.
However, for any automorphism $\phi$ of a Sasakian manifold,
$\phi^k$ is homotopy equivalent to identity for some 
$k \in \Z^{>0}$. Indeed, the group $\Aut(S)$ is compact, 
as a closed subgroup of an isometry group (\cite{_Myers_Steenrod_,_Kobayashi_Transformations_}), and hence 
it has finitely many connected components. Therefore,
a certain power $\gamma^k$ of $\gamma$ is homotopic
to identity, and this obstruction to the existence
of the logarithm of $\gamma^k$ disappears.


\section{Algebraic structures on Stein completions}\index[terms]{completion!Stein}


\subsection{Ideals of the embedding to a Hopf manifold}

Let $M\subset H$ be a subvariety of a linear Hopf manifold.
In this section we prove that $M$ can be defined by a set 
of polynomial equations.

Let $\gamma\in \Aut(\C^n)$ be a holomorphic diffeomorphism.
Recall that a function $f\in \calo_{\C^n}$ is 
{\bf $\gamma$-finite} if the space $\langle f, \gamma^*f, 
{\gamma^2}^* f, ...\rangle$ is finite-dimen\-sio\-nal 
(see \ref{_F_finite_Definition_}).\index[terms]{function!$\gamma^*$-finite}

We start from the following lemma.

\hfill

\lemma\label{_gamma_finite_on_C^n_Lemma_}
Let $\gamma$ be an invertible linear contraction of $\C^n$.
A holomorphic function on $\C^n$
is $\gamma$-finite if and only if it is polynomial.

\hfill

\proof
Clearly, a polynomial function is $\gamma$-finite.
The operator $\gamma^*$
acts on homogeneous polynomials of degree $d$
with eigenvalues $\alpha_{i_1}\alpha_{i_2}... \alpha_{i_d}$,
where $\alpha_{i_j}$ are the eigenvalues of $\gamma$ on $\C^n$.
Since $\gamma$ is a contraction, all $\alpha_{i_j}$
satisfy $|\alpha_{i_j}| <1$. Therefore,
any sequence $\{\alpha_{i_1}\alpha_{i_2}... \alpha_{i_d}\}$
converges to 0 as $d$ goes to infinity. This implies that
the Taylor decomposition of a $\gamma$-finite\index[terms]{function!$\gamma^*$-finite} function $f$ can   have 
only finitely many components, otherwise the eigenspace
decomposition of $f$ with respect to the
action of $\gamma$ is infinite. 
\endproof

\hfill

We give a new and self-contained proof of the 
following basic result, stated in \cite{ov_pams}.
The proof in \cite{ov_pams} is a 4-line sketch referring 
further to \cite{ov_imm_vai}, where this statement is
actually proven, but not stated anywhere as is.

\hfill

\theorem \label{_cone_cover_for_LCK_pot_algebraic_Theorem_}
{\cite{ov_lee}} for the present proof.)\\
Let $M \hookrightarrow H$ be a complex subvariety of a linear Hopf manifold,
and $\tilde M_c \arrow \C^N$ the corresponding map of Stein completions\index[terms]{completion!Stein},
with $\tilde M_c$ obtained as the closure of $\tilde M \subset \C^N$
by adding the zero. Then $\tilde M_c$ is an algebraic subvariety, that is,
a set of common zeroes of a system of polynomial equations.

\hfill

\pstep  
Let $\tilde M_c$ be the Stein\index[terms]{completion!Stein}
completion of $\tilde M$ equipped with an $A$-equivariant embedding to
$\C^N$, where $A$ acts as a linear operator with all
eigenvalues $|\alpha_i|< 1$ (\ref{_gamma_on_T^*_contraction_Corollary_}). 
The $A$-finite functions\index[terms]{function!$A$-finite} are polynomial and
vice versa by \ref{_gamma_finite_on_C^n_Lemma_}.

\hfill

We want to produce an explicit fundamental domain $U_0$
for the action of 
$\Z\cong \langle A\rangle =
\{ ..., A^{-n}, A^{-n+1}, ... A^{-1}, \Id_{\C^n}, A, A^2, ... \}$
on $\C^N$, in such a way that $U_0 = V \backslash A(V)$,
where $V$ is Stein. Let $B\subset \C^n$ be 
the unit ball. When the operator norm $\|A\|$ of $A$ 
is less than 1, one has $A(B) \Subset B$, and 
$B \backslash A(B)$ is the fundamental domain that we can use.
This would hold, for example, when $A$ is diagonalizable.
On the other hand, the operator norm of a contraction
can be greater than 1. Consider for example 
the matrix
$A={\small \begin{pmatrix} \frac 1 2 & 1000 \\ 0 & \frac 1 2\end{pmatrix}}$;
its norm is at least 1000. Therefore, one should take more
care when choosing the fundamental domain. Recall that any
matrix over $\C$ admits a Jordan decomposition, and
every Jordan cell 
{\begin{equation}\label{_Jordan_cell_Equation_}
		{\scriptsize
		\begin{pmatrix} \alpha & 1 & 0 & \ldots & 0\\
			0 & \alpha & 1 & \ldots & 0\\
			\vdots &\vdots &\vdots & \cdots & \vdots \\
			0&0&0 & \ldots &1\\
			0&0&0 & \ldots &\alpha
		\end{pmatrix}}
	\end{equation}
}
is conjugate to 
{\begin{equation}\label{_Jordan_cell_epsilon_Equation_}
		{\scriptsize
		\begin{pmatrix} \alpha & \epsilon & 0 & \ldots & 0\\
			0 & \alpha & \epsilon & \ldots & 0\\
			\vdots &\vdots &\vdots & \cdots & \vdots \\
			0&0&0 & \ldots &\epsilon\\
			0&0&0 & \ldots &\alpha
		\end{pmatrix}}
\end{equation}}(see {e.g.} \ref{_semisimple_operator_approx_Proposition_} below).
Writing $A$ in a Jordan basis and replacing each cell 
\eqref{_Jordan_cell_Equation_}
with \eqref{_Jordan_cell_epsilon_Equation_}, for $\epsilon$
sufficiently small, we obtain a contraction with an operator
norm $<1$ conjugate to $A$; then $A(B)\Subset B$
is a fundamental domain for the action of $\langle A\rangle$.

\hfill

{\bf Step 3:}
Let $U_0$ be a fundamental domain for $A$ acting
on $\C^N$. As indicated in Step 2,
every linear contraction is conjugate to an operator
$A$ with operator norm $<1$. 
The fundamental domain $U_0$ for $A$ with operator norm
$<1$ can be obtained by taking an open ball
$B\subset \C^n$ and removing $A(B)$ from $B$. Denote by
$U_n$ a copy of this domain obtained as $U_n:= A^{-n}(U_0)$,
and let $V_n:= \{0\}\cup \bigcup_{i>-\infty}^n U_i$.
Since  $V_n= A^{-n}(B)$, it is a Stein domain in $\C^N$.

Let $A:\; H^0_b(\calo_{V_n}) \arrow H^0_b(\calo_{V_n})$
be the operator on the ring of bounded holomorphic functions
induced by the action of $A$ on $\C^N$.
Clearly, this map is compatible with the map
$A:\; H^0(\calo_{\tilde M_c}) \arrow H^0(\calo_{\tilde M_c})$
constructed above; this is what allows us to denote
them by the same letter. By \ref{_contra_compact_Theorem_}, the operator
$A:\; H^0_b(\calo_{V_n}) \arrow H^0_b(\calo_{V_n})$
is compact.

\hfill

{\bf Step 4:} 
Let $\goth{I}(V_n)$ be the ideal of $\tilde M_c\cap V_n$
in $H^0_b(\calo_{V_n})$.
By \ref{_finite_RS_Corollary_},  a compact endomorphism of 
a Banach space\index[terms]{space!Banach} admits a Jordan decomposition.
Then the $A$-finite vectors\index[terms]{function!$A$-finite} are finite linear combinations
of the vectors from the Jordan cells. This implies that
the set of $A$-finite 
functions in $\goth{I}(V_n)$
is dense in $\goth{I}(V_n)$ with respect to the $\sup$-topology.\index[terms]{topology!$\sup$}
On the other hand, all $A$-finite functions can be holomorphically
extended to $\C^N$ by automorphicity. 

The base of $C^0$ (that is, compact-open) topology on $H^0(\calo_{\C^N})$ 
is formed by translations of open sets consisting of all
functions $f\in H^0(\calo_{\C^N})$ 
that satisfy $|f|< C< \infty$ on a given compact set, for some
positive $C\in \R$. Therefore,
it is the weakest topology such that
its restriction to $H^0_b(\calo_{V_n})$
with $\sup$-topology\index[terms]{topology!$\sup$} is continuous for all $n$.
This implies that any set of functions $S\subset H^0(\calo_{\C^N})$
that is  bounded on compacts and dense in $H^0_b(\calo_{V_n})$,
for all $n$, is dense in $H^0(\calo_{\C^N})$.

Since the space of $A$-finite functions
is dense in $\goth{I}(V_n)$,
the space $\goth{I}^A$  of $A$-finite functions\index[terms]{function!$A$-finite} in 
$\goth{I}$ is dense in $\goth{I}$ with respect to the
$C^0$-topology.\index[terms]{topology!$C^0$} In particular,
the set of common zeros of $\goth{I}^A$
coincides with $\tilde M_c\subset \C^N$.

\hfill

{\bf Step 5:}
The $A$-finite functions are polynomials, as shown in Step 1.
By Hil\-bert's basis theorem, \index[terms]{theorem!Hilbert's basis}
any ideal in the ring of polynomials is finitely generated.
Therefore, the ideal $\goth{I}^A$ is finitely generated
over polynomials. Let $f_1, ..., f_n$ be the
set of its generators. By Step 2, the set of
common zeros of $\goth{I}^A$ is 
$\tilde M_c\subset \C^N$; therefore,
$\tilde M_c\subset \C^N$ is given by polynomial
equations $f_1=0, f_2=0, ...,  f_n=0$.
\endproof

\subsection[Algebraic structures on Stein completions: the existence]{Algebraic structures on Stein completions:\\ the existence}\index[terms]{completion!Stein}

Let $X\subset M \subset \C P^n$ be
projective subvarieties of $\C P^n$.\index[terms]{variety!quasi-projective}
Then $Z= M\backslash X$ is called {\bf a quasi-projective variety}.
The {\bf Zariski topology} on $Z$ is the topology such that the
closed subsets of $Z$ are closed algebraic subvarieties.\index[terms]{function!regular}
The {\bf regular functions on a Zariski open subset $W\subset \C P^n$} 
are rational functions on $\C P^n$ that have no
poles on $W$. The {\bf sheaf of regular functions} on $Z$
is a subsheaf of the sheaf of holomorphic functions
in the Zariski topology\index[terms]{topology!Zariski} on $Z$, constructed as follows.
Let  $U\subset Z$ be an open subset, 
obtained as the intersection of $Z$ and a
Zariski open subset $W\subset \C P^n$.
A holomorphic function on $U$ is
{\bf regular} if it is the restriction
of a regular function on $W$.

{\bf An algebraic structure} on a \index[terms]{structure!algebraic}
complex analytic variety $Z$ is a subsheaf of the
sheaf of holomorphic functions that can be
realized as a sheaf of regular functions 
for some biholomorphism between $Z$ and
a  quasi-projective variety $Z= M\backslash X$.\index[terms]{variety!quasi-projective}

Note that the algebraic structure is not
unique; indeed, even the complex manifold
$(\C^*)^2$ has more than one algebraic structure
(\cite{_Simpson:rank_one_}). In \cite{_Jelonek:exotic_},
Z. \index[persons]{Jelonek, Z.} Jelonek has constructed an uncountable set 
of pairwise non-isomorphic algebraic structures on
certain Stein manifolds.\index[terms]{manifold!Stein}

It turns out that 
any K\"ahler $\Z$-cover\index[terms]{cover!K\"ahler $\Z$-} of an LCK manifold with potential\index[terms]{manifold!LCK!with potential}
is equipped with a distinguished affine structure, that is 
uniquely determined by the complex structure  
(\ref{_algebraic_str_on_cone_existence_Theorem_},
\ref{_same_algebra_structure_on_cone_Theorem_}).

\hfill

The main result of this section is the following theorem.

\hfill

\theorem (\cite{ov_ac}) \label{_algebraic_str_on_cone_existence_Theorem_}
Let $M$ be an LCK manifold with proper potential,
and  $\tilde M_c$ the Stein variety
obtained as a Stein completion\index[terms]{completion!Stein} of its K\"ahler $\Z$-cover\index[terms]{cover!K\"ahler $\Z$-} $\tilde M$.
Denote by $\gamma\in \Aut(\tilde M_c)$ 
the generator of the $\Z$-action\index[terms]{action!$\Z$-} on $\tilde M$ extended to 
$\tilde M_c$. This map is a bijective holomorphic contraction
centred in the apex\index[terms]{apex} $c$ of $\tilde M_c$. 
Then $\tilde M_c$ admits an algebraic structure such that
the regular functions are precisely the $\gamma$-finite functions.\index[terms]{function!$\gamma^*$-finite}

\hfill

\pstep
\ref{_gamma_finite_on_C^n_Lemma_}
implies \ref{_algebraic_str_on_cone_existence_Theorem_}
when $\gamma$ is a linear contraction.

\hfill

{\bf Step 2:} By \ref{embedding}, $M$ admits a holomorphic
embedding to a linear Hopf manifold. 
Consider the corresponding embedding 
$\tilde M_c \hookrightarrow \C^n$.
By \ref{_cone_cover_for_LCK_pot_algebraic_Theorem_}, 
$\tilde M_c$ is an algebraic subvariety of $\C^n$ 
given by a finite set of polynomial equations.

To finish the proof of \ref{_algebraic_str_on_cone_existence_Theorem_},
it remains to show that a holomorphic function on $\tilde M_c$ is 
polynomial if and only if it is $\gamma$-finite.\index[terms]{function!$\gamma^*$-finite}
It is clear that any holomorphic function on $\tilde M_c$
obtained from a polynomial function on $\C^n$
is $\gamma$-finite. It remains to show the converse.
Consider a $\gamma$-finite function on $\tilde M_c$. By 
\ref{_gamma_finite_on_C^n_Lemma_},
to prove that it is polynomial, it would suffice
to prove that it can be extended to 
a $\gamma$-finite function on $\C^n$.

Consider the exact sequence of sheaves of holomorphic
functions on $\C^n$
\[ 
0 \arrow I_{\tilde M_c} \arrow \calo_{\C^n}\arrow \calo_{\tilde M_c}\arrow 0,
\]
where $I_{\tilde M_c}$ is the ideal sheaf of $\tilde M_c$.
The corresponding long exact sequence gives
\[
0 \arrow H^0(\C^n,I_{\tilde M_c}) \arrow H^0(\C^n,\calo_{\C^n})\arrow 
H^0(\C^n,\calo_{\tilde M_c})\arrow H^1(\C^n,I_{\tilde M_c}).
\]
Since $\C^n$ is Stein, $H^1(\C^n,I_{\tilde M_c})=0$,
hence the function $f$ is the restriction of
a holomorphic function $\tilde f \in H^0(\C^n,\calo_{\C^n})$.
It remains to show that $\tilde f$ can be chosen
$\gamma$-finite.\index[terms]{function!$\gamma^*$-finite}

\hfill

{\bf Step 3:} This assertion is implied by the following lemma,
which finishes the proof.

\hfill

\lemma
Consider the commutative square of continuous 
operators of Banach spaces\index[terms]{space!Banach}
\[ \begin{CD}
W_1 @> K_1 >> W_1\\
@V R VV @VV R V \\
W_2 @> K_2 >> W_2
\end{CD}
\]
with $R$ a surjective map, and both $K_i$, $i=1,2$, are compact.
Restricting $R$ to the space $W_i^{K_i}$ of $K_i$-finite vectors, we
obtain a map $R_f:\; W_1^{K_1}\arrow  W_2^{K_2}$.
Then $R_f$ is surjective.

\hfill

\proof
Let $U_2\subset W_2$ be a finite-dimensional space
of $K_2$-finite vectors, and $U_1\subset W_1$ its preimage.
Since $R$ is surjective, the natural map
$U_1 \to U_2$ is surjective. Since $U_1^{K_1}$ is dense,
and $U_2$ is finite-dimensional, the restriction
of $R$ to $U_1^{K_1}$ is also surjective.
\endproof

\subsection[Algebraic structures on Stein completions: the  uniqueness]{Algebraic structures on Stein completions:\\ the  uniqueness}\index[terms]{completion!Stein}

In this subsection, we prove that the 
affine algebraic structure constructed on $\tilde M_c$ in 
\ref{_algebraic_str_on_cone_existence_Theorem_}
is uniquely defined by the complex geometry of the
Stein variety $\tilde M_c$.

We start with the following lemma, showing that
a function is $\gamma$-finite\index[terms]{function!$\gamma^*$-finite} (that is, polynomial)
if and only if it admits a certain growth condition.
This is similar to a well-known result of complex
analysis that states that an entire holomorphic
function on $\C^n$ is polynomial if and only if
it has polynomial growth. \index[terms]{function!of polynomial growth}

\hfill

\lemma (\cite{ov_ac})\label{_poly_growth=gamma_finite_Lemma_}
Let $M$ be an LCK manifold with proper potential,
and  $\tilde M_c$ the Stein variety,
obtained as a Stein completion \index[terms]{completion!Stein}of its K\"ahler 
$\Z$-cover \index[terms]{cover!K\"ahler $\Z$-}$\tilde M$. Denote by $\gamma\in \Aut(\tilde M_c)$ the
holomorphic contraction generating the $\Z$-action\index[terms]{action!$\Z$-} on
$\tilde M_c$. Choose a compact set $K\subset \tilde M_c$
containing an open neighbourhood of the apex $c$. Denote by
${\cal B}\subset H^0(\calo_{\tilde M_c})$ 
the ring of functions on $\tilde M_c$
\begin{equation}\label{_polynomial_growth_Equation_}
{\cal B}:= \{ f \in H^0(\calo_{\tilde M_c})\ |\ 
\exists\; C>0 \text{\ such that\ } \forall i \ \ \sup_{K}|(\gamma^*)^{-i}|f| < C^i\}.
\end{equation}
%
Then ${\cal B}$ coincides with the space of $\gamma$-finite \index[terms]{function!$\gamma^*$-finite}
functions.\footnote{We call the function satisfying
\eqref{_polynomial_growth_Equation_} {\bf the function of polynomial growth}.
This terminology is justified because for $\gamma$ a linear contraction 
of $\C^n$, \eqref{_polynomial_growth_Equation_} is equivalent to having 
polynomial growth.}

\hfill

\pstep
Let $f$ be a $\gamma$-finite function\index[terms]{function!$\gamma^*$-finite}, and $W$ the
space generated by $f, (\gamma^*)f, (\gamma^*)^2 f, ...$.
Let $\|\cdot \|_K$ be the norm on $W$ defined by
$\|f\|:= \sup_{K} |f|$, and let
$C:= \sup_{\|f\|=1}\|(\gamma^*)^{-1}f\|$ be
the operator norm of the operator
$(\gamma^*)^{-1}\in \End(W)$ in this norm.
Then $\sup_{K}|(\gamma^*)^{-i}f|\leq C^i \sup_{K} |f|$,
and $f$ has polynomial growth, in the sense of 
\eqref{_polynomial_growth_Equation_}.

\hfill

{\bf Step 2:} 
Suppose that $f$ has polynomial growth, in the sense 
of \eqref{_polynomial_growth_Equation_},
and let $W$ be the space of all functions
generated by $f, (\gamma^*)f, (\gamma^*)^2 f, ...$.
Then all elements of $W$ have the same growth as $f$,
with the same constant $C$, and hence  the closure 
$\bar W$ of $W$ in the norm $\|f\|:= \sup_{K} |f|$ consists
of functions with polynomial growth. 

The condition \eqref{_polynomial_growth_Equation_}
holds for all $f \in \bar W$ if and only if 
$(\gamma^*)^{-1}$ has finite norm on
$\bar W$. To finish the proof of \ref{_poly_growth=gamma_finite_Lemma_},
it remains to show that the norm of $(\gamma^*)^{-1}$
is infinite on $\bar W$ if $\bar W$ is infinite-dimensional.

The operator $\gamma^*$ on $\bar W$ is compact;
by the Riesz--Schauder theorem, it has the Jordan cell
decomposition with eigenvalues converging to 0.\index[terms]{theorem!Riesz--Schauder}
The norm of a linear operator $A$ with eigenvalues
$\alpha_i$ satisfies $\|A\| \geq \sup |\alpha_i|$.
Therefore, a compact operator\index[terms]{operator!compact} cannot be 
invertible on an infinite--dimensional Banach space:\index[terms]{space!Banach}
the inverse operator has infinite norm.
\endproof

\hfill

\theorem (\cite{ov_ac})\label{_same_algebra_structure_on_cone_Theorem_}
Let $M$ be an LCK manifold with proper potential,\index[terms]{potential!LCK!proper}
and  $\tilde M_c$ the Stein variety,
obtained as a Stein completion\index[terms]{completion!Stein} of its K\"ahler $\Z$-cover\index[terms]{cover!K\"ahler $\Z$-} $\tilde M$.
Consider the affine algebraic structure defined in 
\ref{_algebraic_str_on_cone_existence_Theorem_}.
Then this affine algebraic structure is uniquely determined
by the complex structure on $\tilde M_c$.

\hfill

\pstep
Let $\gamma_1, \gamma_2 \in \Aut(\tilde M_c)$
be holomorphic contractions centred in the
apex $c$\index[terms]{apex}, associated with two LCK manifolds
with proper potential such that their  K\"ahler 
$\Z$-covers\index[terms]{cover!K\"ahler $\Z$-} have the same Stein
completion\index[terms]{completion!Stein} $\tilde M_c$.
We prove \ref{_same_algebra_structure_on_cone_Theorem_},
using \ref{_poly_growth=gamma_finite_Lemma_}.
To show that the algebraic structures induced by
$\gamma_1, \gamma_2$ are the same, it would suffice
to show that the spaces of $\gamma_1$- and $\gamma_2$-finite
functions on $\tilde M_c$ coincide. \index[terms]{function!$\gamma^*$-finite}
By \ref{_poly_growth=gamma_finite_Lemma_}, 
this would follow
if the growth estimates 
\eqref{_polynomial_growth_Equation_}
for $\gamma_1$ and $\gamma_2$ are equivalent.

We have reduced \ref{_same_algebra_structure_on_cone_Theorem_},
to the following statement. Let $K\subset \tilde M_c$
be a compact subset that contains an open neighbourhood of 
the apex \index[terms]{apex}$c\in \tilde M_c$. Then 
a function $f \in H^0(\calo_{\tilde M_c})$ that satisfies
$\sup_K(\gamma_1)^{-n}|f| < C_1^n$, also
satisfies  $\sup_K(\gamma_2)^{-n}|f| < C_2^n$
for some other constant $C_2 >0$ and all $n>0$.

\hfill

{\bf Step 2:} Suppose that there
exists an integer $d>0$ such that
$\gamma_1^d(K) \subset \gamma_2(K)$.
Then the polynomial growth estimate \eqref{_polynomial_growth_Equation_}
for $\gamma_2$ follows from the growth estimate for
$\gamma_1^d$, that is  equivalent to the 
growth estimate for $\gamma_1$.
Therefore, \ref{_same_algebra_structure_on_cone_Theorem_}
would follow if we prove that the integer $d$
always exists.

Recall that a continuous map $f:\; M \arrow M$ fixing $x\in M$
is called {\bf a contraction centred in $x\in M$}
if for each compact subset $K\subset M$, 
and each open neighbourhood of $x$, a sufficiently
big iteration of $f$ gives $f^N(K) \subset U$.
The maps $\gamma_1, \gamma_2$ are contractions centred in 
the apex $x$\index[terms]{apex}, and hence  for some $N>0$, the set 
$\gamma_1^N(K)$ lies in the interior
of $\gamma_2(K)$, which is  a neighbourhood of $x$.
\ref{_same_algebra_structure_on_cone_Theorem_} is proven.
\endproof


\section{Algebraic cones}


\subsection{Algebraic cones defined in terms of $\C^*$-action}

As we already mentioned, the K\"ahler cover of an LCK
manifold with potential\index[terms]{manifold!LCK!with potential} is not a cone in the Riemannian sense
(it lacks the homothety), but is very close to the idea of
cone, because it can be completed with just one (singular)
point. Here we give the precise notion of cone that we
shall need:

\hfill

\definition\label{_alge_cone_Definition_}
A {\bf closed algebraic cone}  is an
affine variety $\cac$ admitting a $\C^*$-action $\rho$
with a unique fixed point $x_0$, called {\bf the origin} or {\bf the apex}\index[terms]{apex},
and satisfying the following: \index[terms]{cone!algebraic}

(i) $\cac$ is smooth outside of  $x_0$, and $\C^*$ acts freely 
on $\cac \backslash x_0$.

(ii) $\rho$ acts on the Zariski tangent
space $T_{x_0}\cac$ with all eigenvalues
$|\alpha_i|<1$.\index[terms]{Zariski tangent space}

An {\bf open algebraic cone} is a complex manifold
that is  biholomorphic to a closed algebraic
cone without the origin.

\hfill
 
\theorem (\cite{ov_ac})\label{_cone_cover_for_LCK_pot_Theorem_}
Let $M$ be an LCK manifold with proper potential,\index[terms]{potential!LCK!proper}
and $\tilde M$ be its K\"ahler $\Z$-cover.\index[terms]{cover!K\"ahler $\Z$-}
Then $\tilde M$ is an open algebraic cone.

\hfill

\proof
Let $\gamma$ be the generator of the $\Z$-action\index[terms]{action!$\Z$-} on $\tilde M_c$.
Using \ref{embedding}, we can 
embed $M$ into a linear Hopf manifold.
This produces a holomorphic map
$\tilde M \arrow \C^n\backslash 0$ taking
$\gamma$ to a linear contraction $A\in \End(\C^n)$.
Let $\tilde M_c$ be the closure of $\tilde M$
in $\C^n$. By \ref{_cone_cover_for_LCK_pot_algebraic_Theorem_},
$\tilde M_c$ is equipped with a natural algebraic structure.
To prove \ref{_cone_cover_for_LCK_pot_Theorem_}
it remains only to show that 
this variety admits a holomorphic $\C^*$-action, free
outside of the origin $x$, and acting
with eigenvalues $|\alpha_i|< 1$ on 
the Zariski tangent space $T_x \tilde M_c$.

Let ${\cal G}_A$ be the Zariski closure of
$\langle A \rangle$ in $\GL(\C^N)$. This is a commutative
algebraic group, acting on the variety $\tilde M_c\subset \C^N$.
To prove that $\tilde M_c$ is an algebraic cone, we need to find
a $\C^*$-action contracting $\tilde M_c$ to the origin.
Let $A=SU$ be the Jordan--Chevalley decomposition for
$A$, with $S, U\in {\cal G}_A$. Since ${\cal G}_A$ preserves $\tilde M_c$, 
the endomorphisms $S$ and $U\in \End(\C^n)$ also act on $\tilde M_c$.
Since the eigenvalues of $S$ are the same
as the eigenvalues of $A$, the map $S$ is a linear contraction.
Let $G_S\subset {\cal G}_A$, $G_S= e^{\C \log S}$ be the one-parametric subgroup
containing $S$. We prove that $G_S$ can be approximated by
subgroups of ${\cal G}_A$ isomorphic to $\C^*$; then these subgroups
also contain a contraction, and we are done.

Consider the map taking any $A_1\in {\cal G}_A$ to its unipotent component
$U_1$. Since ${\cal G}_A$ is commutative, this map is a group
homomorphism. Therefore, its kernel $G_s$ (that is,
the set of all semisimple elements in ${\cal G}_A$) is an
algebraic subgroup of ${\cal G}_A$. A connected commutative
algebraic subgroup of $\GL(\C^N)$ consisting of semisimple elements
is always isomorphic to $(\C^*)^k$ 
(\cite[Proposition 1.5]{_Borel_Tits:Groupes_Reductifs_}).\footnote{%
See also Exercise \ref{_semisimple_commu_Exercise_}.}
The one-parametric subgroups $\C^*\subset (\C^*)^k$ 
are dense in $(\C^*)^k$ because one-parametric
complex subgroups $\C^*\subset (\C^*)^k$ can be obtained as 
complexification of subgroups $S^1\subset \U(1)^k \subset (\C^*)^k$,
and those are dense in $\U(1)^k$. Therefore,
the contraction $S\in G_s=(\C^*)^k$ can be approximated
by an element of $\C^*$ acting on $\tilde M_c$.
\endproof

\subsection{Algebraic cones and Hopf manifolds}

An open algebraic cone is always a $\Z$-cover\index[terms]{cover!K\"ahler $\Z$-}
of a Vaisman manifold; we prove this result in the next
chapter. As a preparation, we embed a $\Z$-quotient of an open
algebraic cone to a Hopf manifold.\index[terms]{manifold!Vaisman}

\hfill

\lemma\label{_properly_discont_contra_Lemma_}
Let $\gamma$ act on a complex variety 
$\tilde M_c$ by contractions, contracting $\tilde M_c$ to a point $c$.
Then the corresponding $\Z$-action\index[terms]{action!$\Z$-} on $\tilde M:= \tilde M_c\backslash \{c\}$ 
is properly discontinuous, and hence  $\tilde M/\Z$ is Hausdorff; it is 
a manifold when $\tilde M$ is a manifold.

\hfill

\proof
By definition, the $\Z$-action is properly discontinuous if every
point has a neighbourhood $U$ such that the set 
$\{g\in \Z \ \ |\  \ g(U)\cap  U \neq \emptyset\}$
is finite. Let $x\in \tilde M$ and 
$K$ be the compact closure of an open neighbourhood
of $U\subset \tilde M$ containing $x$. Since $\tilde M_c$ is Hausdorff,
there exists a neighbourhood $W\ni c$ whose closure does not intersect $K$. By definition of contractions,
there exists $N>0$ such that 
$\gamma^n(K)\subset W$ for all $n\geq N$.
This implies that $\gamma^n(K)\cap K = \emptyset$, and also
 $K \cap \gamma^{-n}(K)=\emptyset$.
We have shown that $\gamma^n(K)\cap K= \emptyset$
for all $n \notin [-N, N]$.
\endproof

\hfill

\theorem (\cite{ov_ac}) \label{_alge_cone_to_Hopf_Theorem_}
Let $(\tilde M, \rho:\; \C^* \arrow \Aut(\tilde M))$ be an open algebraic
cone, and $\lambda \in \C^*$, $|\lambda|<1$ a complex number.
Denote by $\langle \gamma \rangle$ the $\Z$-action\index[terms]{action!$\Z$-} on $\tilde M$
generated by $\gamma=\rho(\lambda)$.  By \ref{_properly_discont_contra_Lemma_},
$\tilde M/\langle \gamma \rangle$ is a complex manifold.
Then there exists a holomorphic embedding $j:\; \tilde M/\langle \gamma \rangle\hookrightarrow H$
to a Hopf manifold.  

\hfill

\proof
Let $R$ be the ring of $\gamma$-finite holomorphic functions on\index[terms]{function!$\gamma^*$-finite}
$\tilde M_c$. By \ref{_cone_cover_for_LCK_pot_algebraic_Theorem_},
Step 3, this ring is finitely generated. Let $V\subset R$
be a $\gamma$-invariant space containing all generators 
of the maximal ideal. 
Since $R$ is dense in $\calo_{\tilde M_c}$, the
functions in $V$ separate the points in $\tilde M$. 
By Montel theorem,\index[terms]{theorem!Montel} $R$ is dense in $C^1$-topolog\index[terms]{topology!$C^1$}y whenever it is
dense in $C^0$-topology;\index[terms]{topology!$C^0$} therefore, 
the derivatives of the functions $f\in R$ generate $T^*_x\tilde M$
for all $x\in \tilde M$. This implies that
the tautological map $\tilde M \arrow V^*$
is a holomorphic embedding.
This map is by construction $\gamma$-equivariant,
hence it induces a holomorphic embedding
$\tilde M/\langle \gamma \rangle\hookrightarrow H=V^*/\langle \gamma\rangle$
(see also Subsection \ref{_density_implies_Subsection_}).
\endproof



%
%
%


\section{Exercises}


\begin{enumerate}[label=\textbf{\thechapter.\arabic*}.,ref=\thechapter.\arabic{enumi}]

\item
Let $A=-\Id\in \GL(2n+1, \R)$
be the diagonal $(-1)$-matrix. Prove that $A\neq e^B$
for any other $B\in \GL(n, \R)$.

\item
Let $A_1, ..., A_n\in \End(\C^n)$ be a collection
of commuting endomorphisms of $\C^n$. 
Prove that there exists a basis $x_1, ..., x_n \in \C^n$
such that all $A_i$ are upper triangular.

\item
Let $A_1, ..., A_n\in \End(\C^n)$ be a collection
of commuting diagonalizable endomorphisms of $\C^n$. 
Prove that there exists a basis $x_1, ..., x_n \in \C^n$
such that all $A_i$ are diagonal.

\item
Let $G\subset \GL(\C^n)$ be a commutative Lie subgroup.
Prove that there exists a basis $x_1, ..., x_n \in \C^n$
such that all $g\in G$ are upper triangular.

\item\label{_semisimple_commu_Exercise_}
Let $G\subset \GL(\C^n)$ be a commutative Lie subgroup.
Assume that all $g\in G$ are diagonalizable.
Prove that there exists a basis $x_1, ..., x_n \in \C^n$
such that all $g\in G$ are diagonal in this basis.

\item
Let $G$ be a discrete group acting on a locally compact 
Hausdorff topological space $M$.
Its action is {\bf properly discontinuous} if every
point has a neighbourhood $U$ such that the set 
$\{g\in G \ \ |\  \ g\cdot  U\cap   U \neq \emptyset\}$
is finite.

\begin{enumerate}
\item Prove that $M/G$ is Hausdorff when
the action of $G$ is properly discontinuous.
\item Prove that any orbit of a properly
discontinuous action is closed and discrete.

\end{enumerate}

\item {\bf The codiscrete topology} is a topology 
on $M$ with any open set equal to $M$ or $\emptyset$.
Find an example of $\Z$-action\index[terms]{action!$\Z$-} on $S^1$
such that $S^1/\Z$ is codiscrete.

\item
Let $M$ be a metric space,\index[terms]{space!metric} and $G$ a discrete group
acting on $M$ by isometries. Assume that all orbits
of $G$ are closed and discrete, and that the stabilizers
of all points in $M$ are finite. Prove that
the action of $G$ is properly discontinuous.

\item
Let $G$ act on $M$ freely and properly discontinuously.
Prove that the map $M \arrow M/G$ is a covering.

\item
Let $C$ be a compact 2-dimensional manifold, and
$\sigma:\; C\arrow C$ a diffeomorphism, that is 
not homotopy equivalent to identity. 
Prove that $\sigma$ acts non-trivially on
$\pi_1(C)$.

\item Let $X$ be a projective manifold, $L$ an ample line bundle,
and $\phi$ an automorphism of $X$
that is  non-trivial on some homotopy group
$\pi_i(X)$, $i\in \Z^{>0}$. Assume that $L$ is 
$\phi$-equivariant and admits a $\phi$-invariant
Hermitian metric with positive curvature. Prove that 
the total space $S$ of the corresponding $S^1$-bundle 
is Sasakian. Prove that $\phi$ acts on $S$ as a Sasakian
automorphism. Prove that the action of $\phi$ on $S$
is not homotopy equivalent to identity.

\item
Let $T^n$ be a compact torus, considered to be  a Lie group,
and $\Gamma_a\subset T^n$ a cyclic subgroup generated by $a\in T^n$.
Prove that $\Gamma_a$ is dense in $T^n$ for $a$ outside
of a measure zero subset of $T^n$.

\item Let $R$ be a Banach ring, and $A:\; R \arrow R$
a continuous derivation. Prove that the operator 
$e^A:= \sum_{k=0}^\infty \frac {A^k}{k!}$ is continuous.
Prove that it is an automorphism of $R$.

\item Let $H^0_b(\calo_X)$
be the space of bounded holomorphic
functions on a complex manifold, equipped with the $\sup$-norm.
Prove that it is a Banach ring (that is, complete
with respect to the metric induced by the norm).

\item\label{_C_0_but_not_sup_Exercise_}
Consider the sequence of bounded holomorphic
functions on an open disk, $f_i:= z^i$. Prove that
$f_i$ converges in $C^0$-topology,\index[terms]{topology!$C^0$} and find its limit.
Prove that it does not converge in the $\sup$-topology.\index[terms]{topology!$\sup$}
	
\item Prove \ref{idcomp}: 
let $V=H^0(\calo_M)$ be the space of holomorphic functions on a complex
manifold $M$ with $ C^0$-topology. Prove that any 
	bounded subset of $V$ is precompact. Prove that 
the identity map $\Id_V$ is a compact operator.\index[terms]{operator!compact}

\item
Prove that the ring of holomorphic functions on a disk $\Delta\subset \C$
is non-Noetherian.

\item
Prove that the ring  $\calo_1$ of germs of holomorphic functions 
in one variable is Noetherian.\index[terms]{ring!Noetherian}

\item
Let $A:\; M \arrow M$ be a holomorphic
contraction of a complex variety,
fixing $x\in M$ and  acting on 
the Zariski tangent space $T_x M$\index[terms]{Zariski tangent space}
with all eigenvalues $|\alpha_i|<1$.
Suppose that $A$ is invertible 
as an automorphism of the ring of germs
$\calo_{M, x}$. Let $\calo^A\subset\calo_{M, x}$
be the ring of all $A$-finite elements.\index[terms]{function!$A$-finite}
\begin{enumerate}

\item
Prove that all eigenvalues of $A$ acting on
$\calo^A$ are of form $\alpha_1^{k_1} \alpha_2^{k_2}...\alpha_n^{k_n}$,
where $\alpha_1, \alpha_2,...,\alpha_n$ are the eigenvalues of $A$ on
$T_x M$. 

\item Prove that $\calo^A$ is Noetherian.

{\em Hint:} Use the embedding  \ref{_alge_cone_to_Hopf_Theorem_}.

\item Let $R$ be a Noetherian 
local ring, with the maximal ideal $\goth m$,\index[terms]{ring!local}\index[terms]{ring!Noetherian}
and $R_0\subset R$ its subring. Suppose that 
$\frac{R_0\cap \goth m}{(R_0\cap \goth m)^2} = \frac{\goth m}{\goth m^2}$
and $\frac{R}{\goth m}= \frac{R_0}{R_0\cap \goth m}$.
Prove that $R=R_0$.

{\em Hint:} Use Nakayama lemma.\index[terms]{lemma!Nakayama}

\item For a linear operator $F\in \End(V)$,
a {\bf root space} (``Jordan cell'')\index[terms]{root space}
with eigenvalue $\alpha$\index[terms]{Jordan cell}
is defined as $\bigcup_k \ker (F-\alpha\Id_V)^k\subset V$.
Let $W_k$ is the direct sum of all 
root spaces of $A$
with eigenvalues $\alpha_1, \alpha_2,...,\alpha_n$
in $\calo^A/{\goth m}^k$, where $\goth m$ is the
maximal ideal of $\calo^A$.
Prove that the natural maps
$W_1\leftarrow W_2 \leftarrow W_3 \leftarrow ...$
are all surjective, and for $k$ sufficiently big
they are all isomorphisms.

\item
Let $W\subset \calo^A$ be the direct sum of all
root spaces of $A\in \End(\calo^A)$
with eigenvalues $\alpha_1, \alpha_2,...,\alpha_n$.
Prove that $W$ is finite-dimensional and generates $\calo^A$.

\item
Let ${\goth m}_x\subset \calo_{M, x}$ be 
the maximal ideal of $x$. 
We equip $\calo_{M, x}$ with {\bf ${\goth m}_x$-adic \index[terms]{topology!${\goth m}$-adic}
topology,} with the base ${\goth m}_x^k$, for all
$k\in \Z^{>0}$. Prove that $\calo^A$ is  
dense in $\calo_{M, x}$ in ${\goth m}_x$-adic 
topology.
\end{enumerate}

\item\label{_grading_via_automo_Exercise_}
Let $R$ be a ring over $\C$
equipped with an automorphism $A$,
and $V$ an $A$-invariant finite-dimensional space
of generators of $R$. Assume that all ei\-gen\-va\-lues
$\alpha_i$ of $A$ satisfy $|\alpha_i| >1$.
Prove that $R$ admits a multiplicative grading
$R = \bigoplus_{i\in \Z^{\geq i}} R^i$,
with all $R^i$ finite-dimensional.

\item\label{_affine_cone_Exercise_}
Let $X$ be projective variety. The {\bf affine cone}\index[terms]{cone!affine} of  $X$
is the spectrum of the ring\index[terms]{spectrum}
$\bigoplus_{i\in \Z^{\geq i}} H^0(X, \calo(i))$.
Let $R=\bigoplus_{i\in \Z^{\geq i}} R^i$
be a finitely generated graded ring, with all $R^i$ finite dimensional.
Prove that there exists a projective variety such that
$R$ is the ring of regular functions on its affine cone.\index[terms]{cone!affine}

\item\label{_cone_total_space_Exercise_}
Let $\tilde M_c$ be a closed algebraic cone.
Prove that $\tilde M_c$ is biholomorphic to
an affine cone\index[terms]{cone!affine} of a projective variety.\footnote{In Chapter \ref{linear_hopf_shells},
we prove this result using Vaisman geometry.}

{\em Hint:} Apply Exercise \ref{_grading_via_automo_Exercise_}
to the ring of $\rho$-finite functions\index[terms]{function!$\rho$-finite} on $\tilde M_c$
and use Exercise \ref {_affine_cone_Exercise_}.


\end{enumerate}


\chapter{Pseudoconvex shells and LCK metrics on Hopf manifolds}\index[terms]{metric!LCK}
\label{linear_hopf_shells}

{\setlength\epigraphwidth{0.6\linewidth}
\epigraph{\it
This book is most famous for its sandpaper cover. An
auto-destruction feature that enabled it to damage not
only the book it might be standing next to in the
bookshelf, but also the person who would be reading it. An
anti-book to destroy all other books.}{\sc\scriptsize
``BOOKS OF WARFARE:
THE COLLABORATION BETWEEN
GUY DEBORD \& ASGER JORN
FROM 1957-1959'', by Christian Nolle}
}


\section{Introduction}

\subsection{LCK metrics on Hopf manifolds}\index[terms]{metric!LCK}

The question of finding LCK metrics on Hopf manifolds
goes back to I. Vaisman, who originally called what is 
now known as \index[terms]{manifold!Vaisman} Vaisman manifold ``generalized Hopf manifolds''.\index[terms]{manifold!Vaisman}
Later, (\cite{bel})
F. \index[persons]{Belgun, F. A.} Belgun realized that not all Hopf manifolds are
``generalized Hopf'', which was the main reason why
we call them Vaisman manifolds now.

In fact, the term
``generalized Hopf manifold" \index[terms]{manifold!generalized Hopf} was already used by \index[persons]{Brieskorn, E.} Brieskorn and
van de Ven (see \cite{_Brieskorn_Ven_}) for some products 
of homotopy spheres that do
not bear an LCK structure.\index[terms]{structure!LCK}

The LCK metrics\index[terms]{metric!LCK} on all Hopf surfaces were
constructed in \cite{go}; this construction was generalized
for higher dimension in \cite{kor}. \index[terms]{surface!Hopf}
However, even for the Vaisman Hopf surfaces, 
this construction was very difficult.

An explicit LCK metric on a diagonal Hopf manifold
was given in \cite{_Verbitsky:Vanishing_LCHK_};
however, as observed by Matei \index[persons]{Toma, M.} Toma and Ryushi \index[persons]{Goto, R.} Goto 
later, this metric is singular.
An easy way of building LCK structures on 
Hopf manifolds was suggested in \cite{ov_pams};
in this chapter we present an alternative 
approach to the same construction.

We can define an {\bf open algebraic cone}\index[terms]{cone!algebraic}
as a smooth manifold obtained as the  space
of non-zero vectors in an ample $\C^*$-bundle over a projective
variety (these varieties are {\em a posteriori} 
orbifolds, by
\ref{_cone_Sasakian_total_space_C^*_Corollary_}).
By \ref{_cone_is_total_space_Proposition_}, this
definition is equivalent to \ref{_alge_cone_Definition_} given before.

As follows from \ref{_cone_cover_for_LCK_pot_Theorem_},
a K\"ahler $\Z$-cover\index[terms]{cover!K\"ahler $\Z$-} over an LCK manifold with
proper potential\index[terms]{potential!LCK!proper} is an open algebraic cone.
A {\bf closed algebraic cone} is an affine variety
obtained from an open algebraic cone by 
adding the origin.

Now let  $\cac$ be an open algebraic cone
equipped with a holomorphic action $\rho_t$ of $\R^{>0}$
that is  extended to a holomorphic action
on the closed cone, and contracts $\cac$ to the
origin. A {\bf pseudoconvex shell}\index[terms]{pseudoconvex shell}
is a compact strictly pseudoconvex 
hypersurface $S\subset \cac$ which
intersects each orbit of the $\R^{>0}$-action
exactly once. For any pseudoconvex shell
on an algebraic cone, we obtain a
K\"ahler potential $\phi$ that is  automorphic
in the following sense: $\rho_t^*(\phi)=t^\alpha \phi$,
where $\alpha>0$ is a fixed real number, and  $t\in \R^{>0}$.
Then $\phi$ defines an LCK potential \index[terms]{manifold!LCK!with potential}on the compact manifold 
$\cac/\Z$, where the $\Z$-action\index[terms]{action!$\Z$-} is
generated by $\rho_{t_0}$, for a given $t_0 >1$.
This allows one to define LCK metrics with
potential on any quotient $\cac/\Z$ by contraction
(\ref{_LCK_on_quotients_Corollary_}), in 
particular, on all linear Hopf manifolds.
For the Hopf manifold, one could take the
usual sphere as the pseudoconvex shell\index[terms]{pseudoconvex shell}, and
the usual linear logarithm for the $\vec r$-action.

Let $M$ be an LCK manifold with proper potential,
and $\tilde M$ its $\Z$-covering.\index[terms]{cover!K\"ahler $\Z$-}
The logarithm of $\gamma:\; \tilde M \arrow \tilde M$
can be interpreted as a holomorphic $S^1$-action on 
$M$ which lifts to an $\R$-action\index[terms]{action!$\R$-} on $\tilde M$.
It is proven in \ref{_S^1_potential_Theorem_} 
that the existence of such an action on an LCK manifold implies the existence
of LCK potential.\index[terms]{potential!LCK}

Not all LCK metrics with potential\index[terms]{metric!LCK!with potential} are obtained
from pseudoconvex shells\index[terms]{pseudoconvex shell}; however, we can fully
describe the ones that are  obtained in this way.
In \ref{_LCK_pot_via_d_theta_d^c_theta_Proposition_},
we characterized LCK metrics with potential in terms
of $d_\theta d^c_\theta$-potential 
(\ref{_d_theta_d^c_theta-potential_Definition_}),
proving that an LCK metric $\omega$ on an
LCK manifold $(M, \omega,  \theta)$ admits a potential
if and only if $\omega=d_\theta d^c_\theta(\phi_0)$,
for a positive function $\phi_0\in C^\infty(M)$.

The LCK metrics with potential\index[terms]{metric!LCK!with potential} obtained from pseudoconvex shells
can be characterized in terms of $d_\theta d^c_\theta$-potential
as follows.

Let $(M, \theta, \omega)$ be an LCK 
manifold with proper potential admitting a holomorphic 
$S^1$-action obtained as the logarithm of the monodromy.
Then $\omega$ is obtained from a\index[terms]{action!monodromy} 
pseudoconvex shell\index[terms]{pseudoconvex shell} if and only if 
$\omega=d_\theta d^c_\theta(\phi_0)$, where $\phi_0$
is an $S^1$-invariant $d_\theta d^c_\theta$-potential. 
Moreover, every LCK manifold with
such an $S^1$-action admits a 
metric with $S^1$-invariant $d_\theta d^c_\theta$-potential
(\ref{_averaging_LCK_potential_proposition_}).

\subsection{Affine cones of projective varieties}\index[terms]{cone!affine}

We shall use the structure theorem for Vaisman manifolds (\ref{str_vai}) to\index[terms]{manifold!Vaisman}
give an equivalent definition of an algebraic cone.
Let $X$ be a projective variety, and
$ R:= \bigoplus_{i\in \Z^{\geq 0}} H^0(X, \calo(i))$ 
the corresponding graded ring. Since $R$ is finitely
generated, its spectrum is an affine variety,\index[terms]{variety!affine}
called {\bf the affine cone} of $X$.\index[terms]{cone!affine}
We prove that a variety $Z$ with at most\index[terms]{spectrum}
one isolated singularity is a closed algebraic cone\index[terms]{cone!algebraic}
if and only if it is biholomorphic to an 
affine cone over a projective orbifold.

This theorem is purely algebraic, and can be proven
algebraically, but the proof using LCK geometry is more 
geometric and intuitive.

The notion of an open algebraic cone\index[terms]{cone!algebraic} is 
defined in terms of the closed algebraic cone:
it is a smooth manifold $\tilde M$ admitting a free $\C^*$-action
and embedded to a closed algebraic cone $\tilde M_c$,
such that $\tilde M_c \backslash \{c\}= \tilde M$,
and $\C^*$ contracts $\tilde M_c$ to $c$, called
{\bf an origin}, or {\bf an apex}\index[terms]{apex}. The procedure called 
``the Stein completion''\index[terms]{completion!Stein}
(\ref{_Stein_completion_Definition_})
allows  a closed
algebraic cone \index[terms]{cone!algebraic} to be reconstructed from an open one. However,
there can be many closed cones associated
with a given open cone $\tilde M$.

Recall that a complex variety $X$ is called \index[terms]{variety!normal}
{\bf normal} if any bimeromorphic finite \index[terms]{map!bimeromorphic}
map $X' \arrow X$ is an isomorphism.\index[terms]{map!finite}
\footnote{We actually use an equivalent
definition (Subsection
\ref{_normality_of_cone_Subsection_}): a variety is normal if all
locally bounded meromorphic functions
are holomorphic.\index[terms]{function!meromorphic}}

We use the Stein completion\index[terms]{completion!Stein} to prove that
any open algebraic cone $\tilde M$ can be associated
with a unique closed algebraic cone\index[terms]{cone!algebraic} $\tilde M_c$, that is 
normal. Other closed algebraic cones associated
to $\tilde M$ are homeomorphic to $\tilde M_c$,
because they are all obtained by adding a point
in the origin. We prove that the natural homeomorphism map 
$\sigma:\; \tilde M_c\arrow \tilde M_c'$ to another closed 
algebraic cone\index[terms]{cone!algebraic} associated with $\tilde M$
is holomorphic. Moreover, $\sigma$ is the
normalization map.\index[terms]{normalization} 

This is related to the notion of ``projective normality''
that is much studied in classical algebraic geometry\index[terms]{geometry!algebraic}.\index[terms]{variety!projectively normal}
A projective variety $X \subset \C P^n$ is called 
{\bf projectively normal} if the ring
$R= \bigoplus_{i\in \Z^{\geq 0}} H^0(X, \calo(i))$ of the
regular functions in its affine cone is integrally closed;\index[terms]{function!regular}
this is equivalent to the affine cone being 
normal. Projective normality is a very subtle notion,
because it greatly depends on the choice of
projective embedding (\cite[ Exercise I.3.18]{_Hartshorne:AG_}).
Before the advent of the Hodge theory,
projective normality was used to prove things such as the Riemann--Roch formula.
O. Zariski and his student H. T. \index[persons]{Muhly, H. T.} Muhly, who were the first to study
projective normality in a systematic manner,
called it ``arithmetic normality'' 
(\cite{_Muhly_,_Zariski:complete_});
however,  ``projective normality'' stuck 
because it was used in the textbook by Hodge and
\index[persons]{Pedoe, D.} Pedoe, \cite{_Hodge_Pedoe_}.


\section{Pseudoconvex shells}\label{shells}\label{pseudo_shell}


\subsection{Pseudoconvex shells in algebraic cones}\index[terms]{cone!algebraic}

The notion we need is the following.

\hfill

\definition 
Let $\tilde M$ be an open algebraic cone,\index[terms]{cone!algebraic} $\tilde M_c$
the corresponding closed cone, 
and $\vec r\in T\tilde M_c$ a holomorphic vector field
such that for all $t>0$ the diffeomorphism $e^{t\vec r}$
is a holomorphic contraction of $\tilde M_c$ to the origin.
In this situation, the 1-parameter 
family $e^{t\vec r}$ is called {\bf the flow of contractions} (or {\bf contraction flow}).\index[terms]{contraction!flow of}
A strictly pseudoconvex hypersurface $S\subset \tilde M$ is called \index[terms]{hypersurface!strictly pseudoconvex}
a {\bf  pseudoconvex shell} if $S$ intersects each orbit
of $e^{t\vec r}$, $t\in \R$, exactly once.\index[terms]{pseudoconvex shell}

\hfill

Pseudoconvex shells are used to obtain 
automorphic plurisubharmonic functions, as follows.\index[terms]{function!plurisubharmonic}

\hfill

\theorem  { (\cite{ov_pams})}\label{shell_char}
	 Let $\tilde M$
be an algebraic cone,\index[terms]{cone!algebraic} $e^{t\vec r}$ a contraction flow, and
$S\subset \tilde M$ a pseudoconvex shell.\index[terms]{pseudoconvex shell} Then for each $\lambda\in \R$
there exists a unique function $\phi_\lambda$ such that 
$\Lie_{\vec r}\phi = \lambda \phi$ and $\phi_\lambda\restrict S=1$.
Moreover,
{such $\phi_\lambda$  is strictly 
	plurisubharmonic when $\lambda \gg 0$.}

\hfill

\proof For each $\lambda$,
{$\phi_\lambda$ is uniquely determined on each orbit of
	$e^{t\vec r}$,} $t\in \R$, because $\phi_\lambda$
restricted to this orbit is $e^{\lambda t}$.

 Let $B:= e^{\R \vec r}\cdot (TS \cap
I(TS))\subset T\tilde M$ be the sub-bundle obtained from
$TS \cap I(TS)$ by translations along the flow $e^{t\vec r}$.
{Then $dd^c \phi\restrict B$ is the Levi form\index[terms]{form!Levi} of $S$,}
hence it is positive definite (\ref{propos}).

 Replacing $\phi$ by $\phi^{a}$
amounts to replacing $\lambda$ by $a\lambda$.
Then
\begin{equation}\label{_power_of_potential_Equation_}
dd^c \phi^{a} = \phi^{a-2} (a \cdot \phi dd^c \phi + a(a-1)
d\phi\wedge d^c\phi).
\end{equation}
To finish the proof, it would suffice to show that 
$dd^c \phi^{a}\restrict S >0$ for $a$ sufficiently big.
However, $S$ is compact, and hence  this is implied by the
following elementary linear algebra lemma applied to $V=TM$, $W=B$, 
$h_1 = \phi dd^c \phi$, $h_2= d\phi\wedge d^c\phi$.

\hfill

\lemma\label{sum_herm_forms} (\cite[Lemma 2.17]{ov_pams},
\cite[Lemma 2]{_Tomberg_})\\
Let $h_1, h_2$ be pseudo-Hermitian forms\index[terms]{form!pseudo-Hermitian} on a complex 
vector space $V$, and let $W\subset V$ be a subspace of
codimension 1. Assume that $h_1 \restrict W$ is
strictly positive, $h_2 \restrict W=0$, and
$h_2 \restrict {V/W}$ is also strictly positive.
{Then there exists a number $T_0\in \R$ that depends
	continuously on $h_1, h_2$ such that
	$h_T:=h_1 + T h_2$ is positive definite for all $T>T_0$.}

\hfill

\proof We think of $h_1$, $h_2$ as of 
real-valued bilinear symmetric forms. Let $y\in V$ be a
vector that satisfies $h_2(y,y)=1$. Then any $x\in V$ 
can be written as $x=ay+z$, $z\in W$. This gives
\begin{equation}\label{bu}
h_T(x,x)= Ta^2+ a^2 h_1(y,y)+ h_1(z,z) +2 a h_1(z, y) 
\end{equation}
Consider \eqref{bu} as a polynomial in the variable $a$. Then $h_T(x,x)$ is positive
definite for all $a$ if and only if 
\begin{equation}\label{bubu} 
(h_1(z,y))^2 - (T+ h_1(y,y))\cdot h_1(z,z) <0.
\end{equation}
Since $h_1 \restrict W$ is non-degenerate, there exists $y' \in W$
satisfying $h_1(z, y')=h_1(z,y)$. Take $T$ such that
$T> h_1(y',y')- h_1(y,y)$.
Then \eqref{bubu} becomes
\[ (h_1(z,y'))^2 - h_1(y',y') h_1(z,z) <0
\]
that is  true by Cauchy--Buniakovsky--Schwarz inequality\index[terms]{inequality!Cauchy--Buniakovsky--Schwarz} for the scalar product defined by $h_1$ in $W$.
Then \eqref{bubu} and \eqref{bu} are satisfied for this choice of $T$.
\endproof

\hfill

\remark\label{_power_of_pote_Remark_}
The equation \eqref{_power_of_potential_Equation_}
can be used to show that $\phi^a$ is an LCK
potential for all $a\in \R^{>1}$, if $\phi$ is
an LCK potential.\index[terms]{potential!LCK} The corresponding Lee form\index[terms]{form!Lee} is
$\theta_a:= - d \log(\phi^a)= -a d\log\phi=a \theta$.
In other words, 
if the Lee form of $\phi$ is equal to $\theta$,
the Lee form of $\phi^a$ is $a\theta$.

\subsection{All linear Hopf manifolds are LCK with potential}\index[terms]{manifold!LCK!with potential}

\corollary \label{_linear_LCK_pot_Corollary_}
{ (\cite{ov_pams})} 
{All linear Hopf manifolds are LCK with potential.}

\hfill

\proof Let $M=(\C^n \backslash 0)/\langle A\rangle$,
where $A\in \GL(n, \C)$ is a linear map with all 
eigenvalues $|\alpha_i|<1$. Since the exponential
map is surjective on $\GL(n, \C)$, we can always
choose a logarithm matrix $\log A$, with\index[terms]{logarithm!of a matrix}
$e^{\log A}=A$. Since the eigenvalues $\beta_i=\log \alpha_i$
of $\vec r:=\log A$ all satisfy $\Re \beta_i <0$, the family
$e^{t\vec r}$ is a contraction flow.\index[terms]{contraction!flow of}

If $A$ is diagonal, and the Hermitian metric
on $\C^n$ is diagonalized in the same coordinates,
the contraction flow restricted on the unit sphere 
$S\subset \C^n$, $S=\{(z_1, ..., z_n)\ \ |\ \ \sum_i|z_i|^2=1\}$
is directed towards the inside
of $S$, and hence  $S$ is a pseudoconvex shell.\index[terms]{pseudoconvex shell}

Otherwise, it is possible that
$A(S)$ intersects $S$, and hence  some trajectories
of the contraction flow intersect $S$ several times.
To avoid this, we do a coordinate change, as follows.
First, use the Jordan normal form theorem
to choose the basis in $\C^n$ in such a way
that $A$ is a sum of several Jordan cells.\index[terms]{Jordan cell} 
Using the standard linear algebraic argument
(\ref{_semisimple_operator_approx_Proposition_})
we change the basis once more,
replacing all Jordan cells 
{\[\scriptsize 
A =\begin{pmatrix} \alpha_i & 1 & 0 & \ldots & 0\\
0 & \alpha_i & 1 & \ldots & 0\\
\vdots &\vdots &\vdots & \cdots & \vdots \\
0&0&0 & \ldots &1\\
0&0&0 & \ldots &\alpha_i
\end{pmatrix}
\]}
with 
{\[\scriptsize
A_\lambda =\begin{pmatrix} \alpha_i & \epsilon  & 0 & \ldots & 0\\
0 & \alpha_i & \epsilon  & \ldots & 0\\
\vdots &\vdots &\vdots & \cdots & \vdots \\
0&0&0 & \ldots &\epsilon \\
0&0&0 & \ldots &\alpha_i
\end{pmatrix}.
\]}
For $\epsilon$ sufficiently small,
the flow $e^{t\vec r}$ restricted to $S$ is directed 
towards the center of the sphere, in the same way
as it happens for the diagonal matrix $A$. 

Then $S$ is a pseudoconvex shell\index[terms]{pseudoconvex shell}, and for $\lambda$ sufficiently big 
a plurisubharmonic function $\phi_\lambda$ gives an LCK potential.\index[terms]{potential!LCK}\index[terms]{function!plurisubharmonic}
\endproof

\hfill

\remark 
	For Hopf surfaces, the above result was proven in \cite{go}, 
using an {\sl ad hoc} deformation argument.\index[terms]{surface!Hopf}

\subsection{Pseudoconvex shells in algebraic cones}\index[terms]{cone!algebraic}

Using the embedding to $\C^n$ constructed in 
\ref{_alge_cone_to_Hopf_Theorem_}, we will build a pseudoconvex shell in 
any algebraic cone.\index[terms]{pseudoconvex shell}

\hfill

\proposition\label{_shell_in_alge_cone_Proposition_}
Let $(\tilde M_c, \rho:\; \C^* \arrow \Aut(\tilde M))$ 
be a closed algebraic cone,\index[terms]{cone!algebraic} and $\vec r\in \Lie(\C^*)$ a vector field
tangent to the $\C^*$-action and contracting $\tilde M_c$ to the
origin. Then there exists a pseudoconvex shell\index[terms]{pseudoconvex shell} $S \subset \tilde M$
in the corresponding open algebraic cone $\tilde M$, that is,
a strictly pseudoconvex hypersurface\index[terms]{hypersurface!strictly pseudoconvex}
meeting each orbit of $e^{\R \vec r}$ exactly once.

\hfill

\proof
Consider the $\Z$-action\index[terms]{action!$\Z$-} generated by $\gamma \in \rho(\C^*)$, 
$\gamma = e^{\lambda \vec r}$
that acts on $\tilde M_c$ by contractions. By \ref{_alge_cone_to_Hopf_Theorem_},
the quotient $\tilde M/\langle \gamma\rangle$ admits
an embedding to a Hopf manifold $(\C^n\backslash 0)/\langle \gamma\rangle$.
This gives a $\gamma$-invariant map $\tilde M \arrow \C^n\backslash 0$.
By construction, this map takes $\vec r$ to $\log \gamma\in T(\C^n\backslash 0)$.

Choose a pseudoconvex shell $S\subset \C^n\backslash 0$
which meets each orbit of $e^{\R \vec r}$ once. Then 
$S$ is a pseudoconvex shell\index[terms]{pseudoconvex shell} in $\tilde M$.
\endproof

\subsection{LCK manifolds admitting an $S^1$-action}
\label{_S^1_on_LCK_Subsection_}\index[terms]{manifold!LCK}

The argument proving \ref{_linear_LCK_pot_Corollary_}
can in fact be generalized to obtain LCK structures\index[terms]{structure!LCK!with potential}
with potential on $M$ whenever $M$ admits
an $S^1$-action that is  lifted to a free
$\R$-action \index[terms]{action!$\R$-}on a $\Z$-cover\index[terms]{cover!K\"ahler $\Z$-} $\tilde M$, 
and $\tilde M$ contains a pseudoconvex shell\index[terms]{pseudoconvex shell}. 

\hfill

\definition
Let $M$ be a manifold and $\tilde M$ its $\Z$-cover,\index[terms]{cover!K\"ahler $\Z$-}
with the monodromy\index[terms]{monodromy} generated by $\gamma\in \Diff(\tilde M)$.
We call an $\R$-action $e^{t\vec r}$ on $\tilde M$ 
{\bf a logarithm of $\gamma$} if \index[terms]{logarithm!of $\gamma$}
$\gamma= e^{t_0\vec r}$ for some $t_0\in \R$.
This is equivalent to having an $S^1$-action on $M$,
that is  lifted to a free $\R$-action\index[terms]{action!$\R$-} on $\tilde M$.

\hfill

\theorem\label{_LCK_with_pot_from_shell_Theorem_}
Let $M$ be a complex manifold, 
and $\tilde M$ its $\Z$-cover\index[terms]{cover!K\"ahler $\Z$-}, with the $\Z$-action\index[terms]{action!$\Z$-}
generated by $\gamma$. Suppose that $\gamma$
admits a logarithm $e^{t\vec r}$, and 
$\tilde M$ contains a pseudoconvex shell,\index[terms]{pseudoconvex shell}
that is, a compact pseudoconvex hypersurface
$S\subset \tilde M$ which intersects
each orbit of the flow $e^{t\vec r}$
precisely once. As above, we consider 
$e^{t\vec r}$ to be a holomorphic $S^1$-action on $M$.
Then $M$ admits an $e^{t\vec r}$-invariant LCK structure\index[terms]{structure!LCK!with potential}
with potential constructed using the same
procedure as in \ref{shell_char}. Moreover, every
$e^{t\vec r}$-invariant LCK structure
with potential and preferred gauge\index[terms]{structure!LCK!with potential and preferred gauge} 
is obtained in this way.

\hfill

\proof
Let $\phi_\lambda$ be a function such that 
$\Lie_{\vec r}\phi_\lambda = \lambda \phi_\lambda$ 
and $\phi_\lambda\restrict S=1$.
This function exists and is unique (\ref{shell_char}).
Moreover, $\phi:=\phi_\lambda^\alpha$ is an LCK
potential, for $\alpha\in \R$ sufficiently big
(see the proof of \ref{shell_char}).
The corresponding LCK metric\index[terms]{metric!LCK}
with preferred gauge $\frac{dd^c \phi}{\phi}$
is $e^{t\vec r}$-invariant, because 
$dd^c\phi$ and $\phi$ have the same automorphy
multipliers. 

Conversely, let $(\omega, \theta)$
be an $e^{t\vec r}$-invariant LCK structure\index[terms]{structure!LCK!with potential and preferred gauge}
with potential $\phi\in C^\infty M$ and preferred gauge on $M$.
Then $e^{t\vec r}$ fixes the Lee form \index[terms]{form!Lee}
$\theta =-d\log\phi$, and hence  it multiplies 
the potential $\phi$ by a constant.
Consider the hypersurface $S=\{ x\in \tilde M \ \ |\ \ \phi(x)=1\}$;
it is strictly pseudoconvex because $\phi$ is 
strictly plurisubharmonic. Clearly, the function
$\phi$ is recovered from $S$ and $e^{t\vec r}$-action 
as above.
\endproof

\hfill

\definition\label{_metric_from_shell_Definition_}
An LCK structure with potential \index[terms]{structure!LCK!with potential}obtained as in
\ref{_LCK_with_pot_from_shell_Theorem_} is called
{\bf an LCK structure with potential obtained from a
pseudoconvex shell}.\index[terms]{pseudoconvex shell}

\subsection{Existence of $S^1$-action on an LCK manifold with potential}\index[terms]{manifold!LCK!with potential}

As seen from \ref{_only_for_gamma^k_log_exists_Example_}, 
an LCK manifold with potential does not always admit an $S^1$-action
associated with the logarithm of the monodromy\index[terms]{monodromy} as in 
\ref{_LCK_with_pot_from_shell_Theorem_}. In this subsection,
we  prove that it always has a 
finite cover that admits an $S^1$-action.

We often need to average the 
Lee form\index[terms]{form!Lee} and the LCK metric\index[terms]{metric!LCK} with respect to 
a holomorphic $S^1$-action. When the $S^1$-action
exists on a finite cover, this is still possible 
by averaging on the cover and making sure the
result is monodromy--invariant.

\hfill
 
\proposition\label{_averaging_LCK_potential_proposition_}
Let $\phi_0$ be a $d_\theta d^c_\theta$-potential
on an LCK manifold\index[terms]{manifold!LCK} $(M, \omega, \theta)$, and
$\tilde M$ its K\"ahler $\Z$-cover\index[terms]{cover!K\"ahler $\Z$-} with 
the deck transform action generated by $\gamma$.
Let $\vec r$ be the logarithm of $\gamma^k$
constructed above. Since $\vec r$ is 
$\gamma^k$-invariant, we can consider it
as a vector field on the quotient 
$M_k:=\tilde M / \langle \gamma^k\rangle$.
Denote by $\sigma:\; M_k \arrow M$ the covering map.
Then there exists an $e^{t\vec r}$-invariant 
Lee form\index[terms]{form!Lee} $\theta'\in \Lambda^1(M_k)$ in the same cohomology
class as $\sigma^*\theta$ and an $e^{t\vec r}$-invariant
$d_{\theta'} d^c_{\theta'}$-potential 
$\phi_0'$ on $M_k$.

\hfill

\proof
Since the action of $e^{t\vec r}$  on $M_k$ is
factorized through a circle $S^1$, we have an $e^{t\vec
  r}$-invariant closed 1-form $\theta'$ in the
same cohomology class as $\sigma^*\theta$.
Let $\omega_k\in \Lambda^{1,1}(M_k)$
be an LCK form \index[terms]{form!LCK}in the same conformal class
as $\sigma^* \omega$ that satisfies
$d_{\theta'}\omega_k=0$, and
$\phi_k$ its $d_{\theta'} d^c_{\theta'}$-potential.
Averaging $\phi_k$ with the $S^1$-action, we obtain
an $S^1$-invariant $d_{\theta'} d^c_{\theta'}$-potential
$\phi_0'$ such that $d_{\theta'} d^c_{\theta'}\phi_0'=\omega'$,
where $\omega'$ is obtained by averaging $\omega_k$
with $S^1$.
\endproof

\hfill

\remark The $d_{\theta'} d^c_{\theta'}$-potential
$\phi_0'$ on $\tilde M / \langle \gamma^k\rangle$
constructed in \ref{_averaging_LCK_potential_proposition_}
is obtained from the pseudoconvex shell\index[terms]{pseudoconvex shell} $\phi_0'^{-1}(1)$
and the vector field $\vec r$ as in \ref{_metric_from_shell_Definition_}. 

\hfill

\theorem\label{_LCK_potential_from_shells_Theorem_} 
 { (\cite{ov_pams})} 
Let $M$ be a compact LCK manifold $M$ with 
proper potential, $\tilde M$ its K\"ahler\index[terms]{manifold!LCK!with potential} $\Z$-cover\index[terms]{cover!K\"ahler $\Z$-},
$\gamma$ the generator of $\Z$-action,\index[terms]{action!$\Z$-}
and $\vec r$ the logarithm of some power $\gamma^k$,
which exists by \ref{logar}. Then 
$\tilde M$ admits a $\gamma$-invariant
LCK potential\index[terms]{potential!LCK} $\phi$, that is  also
$e^{t\vec r}$-automorphic: 
$e^{t\vec r}(\phi)= e^{a t}\phi$, 
for some $a> 1$. In particular,
$\phi$ is obtained from 
a pseudoconvex shell\index[terms]{pseudoconvex shell} as in 
\ref{_metric_from_shell_Definition_}. 

\hfill

\proof 
Let $M_k:= \tilde M/\langle \gamma^k\rangle\stackrel p\arrow M$ 
be the $k:1$ cover of $M$, and $\phi\in C^\infty M$ the LCK potential
obtained from an $S^1$-invariant $d_{\theta'} d^c_{\theta'}$-potential
$\phi_0'$  (\ref{_averaging_LCK_potential_proposition_}).
By construction, the LCK structure\index[terms]{structure!LCK} $(M_k, \theta',\omega')$ 
has the same Lee cohomology class as the original LCK
structure $(M_k, \sigma^*\theta,\sigma^*\omega)$.

Then $\theta'=-d\log\phi$ can be understood as the Lee form\index[terms]{form!Lee} of $M_k$.
The action of $\gamma$ on $\tilde M$  generates
an action of $\Z/k\Z$ on $M_k=\tilde M/\langle\gamma^k\rangle$.


Since the projection $\sigma:\; M_k \arrow M$ is compatible
with the LCK structure,\index[terms]{structure!LCK} the cohomology class of the
Lee form\index[terms]{form!Lee} $\theta'$ is lifted from the projection $M_k \arrow M$.
Therefore, the class $[\theta']\in H^1(M_k)$ is
$\Z/k\Z$-invariant, or, equivalently, it is 
$\gamma$-invariant.

Averaging the conformal multiplier with $\Z/k \Z$,
we can replace the LCK structure $(M, \omega', \theta')$
with a $\Z/k\Z$-invariant LCK structure
$(M, \omega'', \theta'')$ in the same conformal class.
The corresponding metric
is by construction also $e^{t\vec  r}$-invariant
(it is obtained by averaging $e^{t\vec  r}$-invariant 
metrics over $\Z/k \Z$).

This new LCK structure\index[terms]{structure!LCK} also admits an LCK potential\index[terms]{potential!LCK}
$\phi_1:\; \tilde M \arrow \R^{>0}$
(Exercise \ref{_LCK_pot_conf_invariant_Exercise_}).
Then the $e^{t\vec r}$-invariant LCK metric with preferred gauge
$\omega_1:= \frac{dd^c\phi_1}{\phi_1}$ satisfies
$\omega_1=d^c\theta''+ \theta''\wedge (\theta'')^c$ (\ref{dcthetaex}),
and hence it is $\Z/k\Z$ invariant on $M_k$, and
can be lifted from $M$. However, any
$e^{t\vec r}$-invariant LCK metric on $M_k$ is obtained from 
the pseudoconvex shell\index[terms]{pseudoconvex shell} construction (\ref{_LCK_with_pot_from_shell_Theorem_});
thus the automorphic K\"ahler potential on $\tilde M$ is
determined by a pseudoconvex shell as in \ref{shell_char}.

We have shown that $M$ admits an LCK metric
with potential\index[terms]{metric!LCK!with potential} $\phi_1$ (and preferred conformal gauge)
obtained from a pseudoconvex shell on its $\Z$-cover\index[terms]{cover!K\"ahler $\Z$-}
$\tilde M$.
\endproof

\hfill

\remark 
In fact, \ref{_LCK_potential_from_shells_Theorem_} 
fully characterizes the class of complex manifolds
admitting a metric with LCK potential.\index[terms]{potential!LCK}
Indeed, by \ref{_S^1_potential_Theorem_}, any LCK manifold admitting
a holomorphic $S^1$-action admits an LCK potential,
if this action does not lift to an $S^1$-action
on its K\"ahler cover. In other words, the existence
of a logarithm of the monodromy action\index[terms]{action!monodromy} on the 
K\"ahler $\Z$-cover $\tilde M$
{\em is equivalent} to the existence of an
LCK potential.

\hfill

Another application of the existence of the logarithm
gives a criterion for an LCK manifold with potential\index[terms]{manifold!LCK!with potential} to be
Vaisman.

\hfill

\theorem \label{_Vaisman_via_logarithm_Theorem_}
 Let $(M,\omega)$ be an LCK-manifold with potential,
 $\tilde M$ its algebraic cone\index[terms]{cone!algebraic}, $\tilde M/\langle
\gamma\rangle =M$, and  $\phi$  the K\"ahler potential on $\tilde M$.
Then {$\omega$ is conformally equivalent to a Vaisman
metric if and only if there exists the logarithm 
$\vec r$ of $\gamma$ such that $\Lie_{I\vec r}\phi =0$.}\index[terms]{logarithm!of $\gamma$}

\hfill

\proof  If $M$ is Vaisman, then $\tilde M=C(S)$,
where $S$ is Sasakian and $\vec r := t\frac{d}{dt}$ is the \index[terms]{manifold!Sasaki}
logarithm.  Moreover, $I\vec r$ is the Reeb field on $S$,
acting on $C(S)$ by holomorphic isometries.\index[terms]{vector field!Reeb}

Conversely, if $\tilde M$ admits a logarithm $\vec r$ 
with such properties, then the corresponding holomorphic
flow acts on $\tilde M$ by homotheties, and $M$ is
Vaisman by \ref{kami_or}.  
\endproof 

\subsection{Quotients of algebraic cones are LCK}\index[terms]{cone!algebraic}

As one of the applications of the formalism developed in
this chapter, we construct an LCK structure\index[terms]{structure!LCK!with potential} with potential
on a quotient of an algebraic cone.\index[terms]{cone!algebraic}

\hfill

\corollary\label{_LCK_on_quotients_Corollary_}
Let $\tilde M$ be an open algebraic cone,
and $\tilde M_c$ the corresponding closed cone,
obtained by adding the origin. Consider a
holomorphic contraction $\gamma:\; \tilde M_c \arrow \tilde M_c$
of $\tilde M_c$ to the origin.\footnote{Unless $\tilde M_c=\C^n$,
the origin is the only singular point of $\tilde M_c$, and hence 
any holomorphic automorphisms of $\tilde M_c$ preserves the origin.}
Then the quotient $\tilde M/\langle \gamma\rangle$
admits an LCK metric with potential.\index[terms]{metric!LCK!with potential}

\hfill

\proof
By \ref{logar}, the map $\gamma^k$ admits a a logarithm
for some $k\in \Z^{>0}$. Using \ref{_shell_in_alge_cone_Proposition_}, we construct a
pseudoconvex shell\index[terms]{pseudoconvex shell} in $\tilde M_c$. Then 
 \ref{_LCK_potential_from_shells_Theorem_} 
implies that $\tilde M/\langle \gamma\rangle$ admits
an LCK potential.\index[terms]{potential!LCK}
\endproof

\subsection[Holomorphic isometries of LCK manifolds with potential]{Holomorphic isometries of LCK manifolds\\ with potential}\index[terms]{manifold!LCK!with potential}

Using the metric obtained from
\ref{_LCK_potential_from_shells_Theorem_}, we 
construct Killing \index[terms]{vector field!Killing}\index[terms]{vector field!holomorphic}holomorphic vector fields on any LCK manifold with 
potential. Recall that {\bf a homothetic vector field} on a Riemannian
manifold $M$ is a vector field whose flow acts on $M$ by homotheties.

\hfill

\theorem\label{_holomorphic_action_on_LCK_Theorem_}
Let $M$ be an LCK manifold with proper potential. 
Then $M$ admits another LCK metric with potential\index[terms]{manifold!LCK!with potential} 
and a holomorphic vector field $\vec {\underline r}$,
which is  lifted to a non-isometric homothetic vector field 
on its K\"ahler $\Z$-covering\index[terms]{cover!K\"ahler $\Z$-} \index[terms]{vector field!holomorphic}
$(\tilde M, \tilde \omega)$.

\hfill

\proof By \ref{_averaging_LCK_potential_proposition_},
there is a finite $k:1$ covering $M_1$ of $M$ 
such that the $\Z$-action\index[terms]{action!$\Z$-} on the $\Z$-covering\index[terms]{cover!K\"ahler $\Z$-}
$\tilde M$ of $M_1$ admits a logarithm $\vec r$.
As we explain in Subsection
\ref{_S^1_on_LCK_Subsection_}, this logarithm 
induces a holomorphic $S^1$-action on $\tilde M$.
Moreover, $M_1$ admits an $S^1$-invariant
LCK metric $g_1$ with potential,\index[terms]{metric!LCK!with potential} obtained
as the pullback of an LCK metric $g$ with potential
on $M$ (\ref{_LCK_potential_from_shells_Theorem_}).

Let $\chi:\; \pi_1(M) \arrow \Z$ be the homothety character,\index[terms]{homothety character} that is, the character
associated with the weight bundle of the LCK structure\index[terms]{structure!LCK} (\ref{_homothety_character_Definition_}).
Using the Galois theory\index[terms]{Galois theory!for coverings} for coverings (Section \ref{carac}), we obtain
that the coverings of any manifold $X$ are in bijective
correspondence with subgroups of the fundamental group.\index[terms]{fundamental group}
In particular, $M_1$ is associated with the
kernel of the composition of $\chi:\; \pi_1(M) \arrow \Z$ 
and the projection $\Z\arrow \Z/k\Z$. Since this subgroup
is normal, the group $\Z/k\Z$ acts on $M_1$ freely
by holomorphic isometries, giving an isomorphism
$M= \frac{M_1}{\Z/k\Z}$.

Averaging $\vec r$ with the $\Z/k\Z$-action, we obtain
a $\Z/k\Z$-invariant holomorphic vector field $\vec r_1$
on $M_1$. Then $\vec r_1$ is the pullback of the
vector field  $\vec {\underline r}$ on $M$.
Since $\vec r_1$ is the average of $\vec r$
and its $\Z/k\Z$-translates $\vec r'$ that satisfy
$\Lie_{\vec r'}\tilde \omega = \tilde \omega$,
the same is true about $\vec {\underline r}$.
Moreover, $\vec r$ and its $\Z/k\Z$-translates
are Killing with respect to the LCK metric $g_1$
on $M_1$, and hence  $\vec {\underline r}$ is Killing
on $(M, g_1)$.\index[terms]{vector field!Killing} 
\endproof

\hfill

\corollary\label{_S^1_action_exists_Corollary_}
Let $M$ be a compact complex manifold admitting an LCK metric with potential.\index[terms]{metric!LCK!with potential}
Then $M$ admits a holomorphic $S^1$-action $\rho$ preserving an LCK structure\index[terms]{structure!LCK!with potential} $(M, \theta, \omega)$ with potential.
Moreover, the infinitesimal generator vector field for this action 
acts on the K\"ahler cover $\tilde M$ of $(M, \theta, \omega)$ by non-isometric homotheties.

\hfill

\proof
Let $\vec {\underline r}$ be the holomorphic Killing vector field\index[terms]{vector field!Killing}\index[terms]{vector field!holomorphic}
constructed in \ref{_holomorphic_action_on_LCK_Theorem_}, and $G$ be the closure
of the Lie group generated by $e^{t\vec{\underline r}}$. Since $\vec{\underline r}$ is Killing, 
$G$ is compact; clearly, it is commutative and acts on $M$ by holomorphic
isometries. Therefore, $G$ is a compact torus, and the one-parametric
subgroup $e^{t\vec{\underline r}}\subset G$ can be approximated by compact
one-parametric subgroups $S^1 \subset G$. This gives a holomorphic
isometric circle action\index[terms]{action!$S^1$-}  on $M$. Lifted to $\tilde M$,\index[terms]{action!$\R$-}
this circle action produces an $\R$-action by holomorphic homotheties,
that is  non-isometric if we chose an one-parametric subgroup 
sufficiently close to $e^{t\vec{\underline r}}$.
\endproof

\hfill

\remark
The Killing vector field\index[terms]{vector field!Killing} $\vec {\underline r}$
obtained in \ref{_holomorphic_action_on_LCK_Theorem_}
is not necessarily a logarithm of the $\Z$-action\index[terms]{action!$\Z$-} on $\tilde M$.
To illustrate this, consider the following example.
Let $\rho$ be a $\Z/k\Z$-action on a Riemannian
surface $S$ of genus $\geq 2$, generated by $\sigma \in \Aut(S)$ 
and $L$ a $\Z/k\Z$-equivariant ample bundle on $S$. Consider 
the Vaisman manifold $\Tot^\circ(L)/\langle \gamma\rangle$
where $\gamma$ is generated by $(x, l) \mapsto (\sigma(x), a \sigma(l))$\index[terms]{manifold!Vaisman}
and $a\in \C$, $|a|>1$, $x\in S, l\in L\restrict x$,
and $(x, l) \mapsto (\sigma(x), \sigma(l))$ denotes the
equivariant action of the group $\Z/k\Z= \langle \sigma\rangle$
on $\Tot^\circ(L)$. As we have already explained (\ref{_only_for_gamma^k_log_exists_Example_}),
the $\Z$-action\index[terms]{action!$\Z$-} on $\Tot^\circ(L)$ does not admit a logarithm.
However, its $\Z/k\Z$-cover is $M_1=\Tot^\circ(L)/\langle \gamma^k\rangle$,
that is  a regular Vaisman manifold,\index[terms]{manifold!Vaisman!regular} because
$\gamma^k(x,l) = (x, a^k l)$. The logarithm $\vec r$ of
the $\Z$-action\index[terms]{action!$\Z$-} on $\Tot^\circ(L)$ is a constant
linear vector field tangent to the fibres of $L$.
On $M_1$, it is the constant vector field tangent
to the leaves of the elliptic fibration.
In particular, it is $\Z/k\Z$-invariant.
The quotient $M= \frac{M_1}{\Z/k\Z}$
is also an elliptic surface, and $\vec {\underline r}$
is also the constant vector field tangent
to the leaves of the elliptic fibration;
however, its exponent preserves the
fibres of the projection $\Tot^\circ(L)\arrow S$,
and the $\Z$-action $(x, l) \arrow (\sigma(x), a \sigma(l))$ 
exchanges some of these fibres. Therefore, $\vec {\underline r}$
is no longer a logarithm of the $\Z$-action.\index[terms]{action!$\Z$-}

\hfill

\remark
Let $M$ be an LCK manifold\index[terms]{manifold!LCK} equipped with a holomorphic 
conformal vector field $X$, such that its lift to the K\"ahler cover acts
by non-isometric homotheties. In \ref{_holomorphic_action_on_LCK_Theorem_},
such a holomorphic vector field is constructed on every LCK manifold with potential,\index[terms]{manifold!LCK!with potential}
for an appropriate choice of LCK metric with potential.\index[terms]{metric!LCK!with potential}
By \ref{kami_or}, $M$ is Vaisman, if 
the 2-dimensional Lie group generated by the vector fields $X, I(X)$
acts on $M$ conformally. When $M$ is non-Vaisman, this
last condition is not satisfied: the holomorphic vector field $I(X)$
cannot be conformal, unless $M$ is Vaisman. If $M$ is compact,
the action of $X, I(X)$ can be integrated to a holomorphic action of 
$\C$, considered to be an abelian complex Lie group; however, this action
is, generally speaking, not conformal.


\section{Algebraic cones as total spaces of $\C^*$-bundles}\index[terms]{cone!algebraic}


\subsection{Algebraic cones: an alternative definition}

The following proposition
can be used as an alternative, more explicit, definition of an open algebraic cone.

Recall that {\bf a closed algebraic cone}\index[terms]{cone!algebraic} is an affine variety
equipped with a $\C^*$-action acting by contractions, free outside
an isolated singularity (the origin) fixed by the $\C^*$-action. 
The {\bf open algebraic cone} is the same space without the origin.\index[terms]{singularity!isolated}

\hfill

\proposition (\cite{ov_ac}) \label{_cone_is_total_space_Proposition_}
An open algebraic cone\index[terms]{cone!algebraic} $\cac_0$ is biholomorphic to the space
$\Tot^\circ(L)$ of non-zero vectors in
an ample line bundle $L$ over a projective orbifold.
Conversely, the total space $\Tot^\circ(L)$ of non-zero
vectors in an ample line bundle $L$ over a projective
orbifold, equipped with the standard $\C^*$-action,
is always an open algebraic cone, if it is smooth.\index[terms]{bundle!line!ample}

\hfill

\proof
Fix the contraction map associated with the
$\C^*$-action $\rho$ on $\cac_0$. 
By \ref{_shell_in_alge_cone_Proposition_}, $\cac_0$
admits a pseudoconvex shell;\index[terms]{pseudoconvex shell} let $\phi$ be the
plurisubharmonic function associated with this shell.
Averaging $\phi$ with the $S^1$-action associated with 
$\rho\restrict{S^1\subset \C^*}$,
we obtain a plurisubharmonic function  $\phi_1$ such that
$\rho(t)^*\phi_1 = \const\cdot \phi_1$. Let $\Z\subset \C^*$ 
be generated by a non-isometric homothety.\index[terms]{function!plurisubharmonic}
By \ref{kami_or}, the corresponding $\Z$-quotient of 
the K\"ahler manifold $(\cac_0, dd^c\phi_1)$ is Vaisman. 
From \ref{str_vai}, we obtain that $\cac_0$ is the
Riemannian cone\index[terms]{cone!Riemannian} of a Sasakian manifold, equipped
with the intrinsic K\"ahler structure. From
\ref{_cone_Sasakian_total_space_C^*_Corollary_} 
we obtain that $\cac_0$ is the Riemannian cone of 
another Sasakian manifold, that is  quasi-regular.\index[terms]{manifold!Sasaki!quasi-regular}
From \ref{_quasireg_Sasakian_orbibundles_Theorem_} 
it follows that $\cac_0$ is the space of 
non-zero vectors in an ample line bundle over 
a K\"ahler orbifold.

Conversely, for any ample line bundle $L$,
the space $\Tot^\circ(L)$ is equipped
with the $\C^*$-action that acts fibrewise and contracts
the Stein completion (\ref{_Stein_completion_Definition_})\index[terms]{completion!Stein}
to the origin, \index[terms]{completion!Stein}
hence it is an algebraic cone.\index[terms]{cone!algebraic}
Indeed, by
\ref{_Structure_of_quasi_regular_Vasman:Theorem_}, 
$\Tot^\circ(L)$ 
is a $\Z$-cover \index[terms]{cover!K\"ahler $\Z$-}of a Vaisman manifold, 
and its Stein completion\index[terms]{manifold!Vaisman}
is given by adding the origin to all
orbits of the $\C^*$-action.
\endproof

\hfill

\remark The proof of
\ref{_cone_is_total_space_Proposition_}
uses the structure theorem for quasi-regular \index[terms]{manifold!Vaisman!quasi-regular}
Vaisman manifolds. However, there exists a purely
algebraic proof of this statement
(Exercise \ref{_cone_total_space_Exercise_}).

%

\subsection{Closed algebraic cones and normal varieties}\index[terms]{cone!algebraic}
\label{_normality_of_cone_Subsection_}

Recall that a complex variety $X$ is called {\bf normal}\index[terms]{variety!complex!normal}\index[terms]{variety!normal}
if any locally bounded meromorphic function on an open subset \index[terms]{function!meromorphic}
$U\subset X$ is holomorphic. In algebraic geometry\index[terms]{geometry!algebraic}, a variety
is normal if all its local rings are integrally closed; these\index[terms]{ring!local}
two notions are equivalent for a complex variety obtained\index[terms]{ring!local!integrally closed}
from an algebraic one (\cite[Theorem II.7.3]{demailly},
\cite[Satz 4, p. 122]{_Kuhlmann_}).\footnote{We are grateful to Francesco Polizzi\index[persons]{Polizzi, F.}
for this reference; please see an
excellent Mathoverflow thread 
{\tiny\url{ https://mathoverflow.net/questions/303406/algebraic-vs-analytic-normality/}}
for more details about the complex analytic and complex algebraic normality.}

We are interested in the normality of a closed algebraic cone.\index[terms]{cone!algebraic}

\hfill

Let $M$ be an LCK manifold with potential,
embedded to a Hopf manifold $H$, and $\tilde M\subset \C^n\backslash 0$
the corresponding $\Z$-coverings. By the Remmert--Stein theorem, 
the closure of $\tilde M$ is complex analytic in $\C^n$;
however, it is not necessarily normal (\ref{_cone_not_normal_Example_}).
This may lead to serious contradictions
if not properly understood. We explain the algebro-geometric ramifications
of this notion; readers who are not interested in algebraic geometry\index[terms]{geometry!algebraic}
are invited to skip the rest of this section.

Now let  $X$ be a projective variety and $L$ an ample line bundle\index[terms]{variety!projective}
on $X$. The {\bf homogeneous coordinate ring} of $X$ is\index[terms]{ring!homogeneous coordinate}\index[terms]{Euler characteristic!holomorphic} 
the ring $\bigoplus_i H^0(X, L^{\otimes i})$. This ring is finitely
generated whenever $L$ is ample (\cite[Theorem 2.3.15]{_Lazarsfeld:1_}). 
Indeed, for an ample line bundle\index[terms]{bundle!line!ample}
$L$ and $i$ sufficiently big, the cohomology $H^{>0}(X,  L^{\otimes i})$ vanishes,
hence the number $\dim H^0(X,  L^{\otimes i})$ is equal to the holomorphic Euler 
characteristic $\chi(L^{\otimes i})$. The latter is polynomial in $i$
by the Riemann--Roch--Hirzebruch theorem, which implies that finitely
many generators is enough to generate this ring.\index[terms]{theorem!Riemann--Roch--Hirzebruch}

A projective variety $(X, L)$ for which the ring $\bigoplus_i H^0(X, L^{\otimes i})$ is integrally
closed is called {\bf projectively normal}. By definition, the affine variety associated with this ring
is {\bf the affine cone of $X$},\index[terms]{cone!affine} that is  a closed algebraic cone in our \index[terms]{cone!algebraic}
parlance. \index[terms]{variety!projectively normal}

The closed algebraic cone, obtained from an open algebraic
cone by taking the Stein completion,\index[terms]{completion!Stein} is always normal.
However, the projective normality is a very tricky condition,
depending on the choice of the bundle $L$.

This was the subject of Exercise I.3.18 from \index[persons]{Hartshorne, R.} Hartshorne's ``Algebraic Geometry''
(\cite{_Hartshorne:AG_}),
where Hartshorne considers the same smooth complex curve with two different
projective embeddings; the affine cone of the first is normal, and the 
affine cone of the second is not. \index[terms]{cone!affine} 

In Exercise II.5.14, \cite{_Hartshorne:AG_},
Hartshorne gives a criterion for projective normality.

\hfill

\proposition\label{proj_normality_Proposition_}
A projective variety $X\subset \C P^n$ is projectively normal if and
only if $X$ is normal, and the restriction map 
\[ H^0(\C P^n, \calo(i)) \arrow H^0\left(X, \calo(i)\restrict X\right)\]
is surjective for all $i \geq 0$. 
\endproof 

\hfill

\remark
This proposition can be applied to algebraic cones,
because the quotient singularities are normal 
(\ref{_quotient_normal_Claim_}),
and the orbifolds have only quotient singularities.

\hfill

\example\label{_cone_not_normal_Example_}
\footnote{We are grateful to Yu. Prokhorov\index[persons]{Prokhorov, Yu.} for this example.}
Let $X\subset \C P^n$ be a projective manifold,
that is  not contained in a projective subspace
$\C P^k \subset \C P^n$ of smaller dimension.
Consider a linear projection $p:\; \C P^n \arrow \C P^{n-1}$
centred in a point $z\notin X$. Then the restriction $p\restrict X$ is holomorphic.
Assume that $p:\; X \arrow \C P^{n-1}$ is also injective
(this can be always achieved if $2\dim X < n-1$ for 
an appropriately general choice of $z$; indeed, $p$ is not injective
if and only if the center $z$ does not belong to the secant 
variety\footnote{A secant variety of $X\subset \C P^n$ is the Zariski closure of the union of all\index[terms]{variety!secant}\index[terms]{Zariski closure}
lines $\C P^n \subset \C P^n$ intersecting $X$ at at least two points.} of $X$,
that has dimension $\leq 2\dim X +1$). Consider the natural map of affine cones\index[terms]{cone!affine}
$u:\; C(X) \arrow C(p(X))$. Outside  the origin $c$, this map is biholomorphic,
hence it is birational and finite.\index[terms]{map!birational}\index[terms]{map!finite} By definition, a variety $Z$ is normal
if any birational finite map $Z_1 \arrow Z$ is an isomorphism.
Were $C(p(X))$ normal, this would imply\index[terms]{Zariski tangent space}
that $u$ is an isomorphism. However, the Zariski tangent space
$T_c C(X)$ is $n+1$-dimensional, because $C(X)\subset C(\C P^n)$ generates 
the vector space $C(\C P^n)= \C^{n+1}$ (otherwise $X$ would have been contained
in a smaller dimension projective subspace). On the other hand,
$\dim T_c C(p(X)) \leq n$, and hence  $u:\; C(X) \arrow C(p(X))$
is not an isomorphism.

\hfill

This example makes sense, if one considers the following statement.

\hfill

\proposition (\cite{ov_ac})\label{_closed_cone_normal_from_open_Corollary_}
Let $\tilde M$ be an open algebraic cone,\index[terms]{cone!algebraic} and
${\goth S}$ the set of all closed algebraic cones
$\tilde M_c$ obtained by adding the origin to $\tilde M$.
Then there exists only one closed algebraic cone $Z\in {\goth S}$ which
is normal, and for any other $Z' \in {\goth S}$, 
its normalization is $Z$.\index[terms]{normalization}

\hfill

\proof
By definition, the Stein completion\index[terms]{completion!Stein} of
$\tilde M$ (\ref{_Stein_completion_Definition_}), denoted by $Z$,
is normal; by \ref{_Forster_Stein_uniqueness_Remark_}, \index[terms]{completion!Stein}
it is unique. By construction, $Z$ is equipped with
a finite, bijective, birational map to any\index[terms]{algebraic cone}
closed algebraic cone $Z'$ associated with $\tilde M$.
Indeed, the holomorphic functions on $Z'$ are
holomorphic on $\tilde M$, and hence  they can be
extended to its Stein completion.
This implies that $Z$ is the normalization of all other 
$Z' \in {\goth S}$.
\endproof

\section{Exercises}

\definition\label{_equivariant_structure_Definition_}
\footnote{See also \ref{_equiv_object_Definition_}.}
Let $B$ be a vector bundle over a manifold
$M$ equipped with an action of a group $\Gamma$.
Recall that the bundle $B$ is called {\bf $\Gamma$-equivariant},
if the group action can be extended from $M$ to
the total space of $B$, preserving the
structure of vector spaces on the fibres.\index[terms]{equivariant structure}
This extension is called {\bf a $\Gamma$-equivariant
structure on $B$}. The equivariant structure is defined by
a family of isomorphisms $R_g:\; g^*B \arrow B$, for all
$g\in G$, which 
are associative in the following sense:
\[ g_2^*(R_{g_1}) \circ R_{g_2}= R_{g_1 g_2},\] where
$R_{g_1 g_2}$ is considered an isomorphism
${(g_1g_2)}^*B \arrow B$, and $g_2^*(R_{g_1}) \circ R_{g_2}$
as a composition of the isomorphism
$g_2^*(R_{g_1}):\; g_2^*g_1^* (B)\arrow g_2^*(B)$
and $R_{g_2}:\; g_2^*(B) \arrow B$.

\begin{enumerate}[label=\textbf{\thechapter.\arabic*}.,ref=\thechapter.\arabic{enumi}]

\item
Let $A\in \End (\C^n)$ be a linear contraction,
that is, an endomorphism with all eigenvalues satisfying
$|\alpha_i| < 1$, and $S\subset \C^n$ a unit sphere.
Prove that $S \cap e^{tA}(S)=\emptyset$ for $t$ sufficiently big.
Find examples where $S \cap e^{tA}(S)\neq \emptyset$
for $t>0$.

\item
Let $U_1, U_2\subset \R^n$ be bounded neighbourhoods of $0\in \R^n$,
such that the boundaries $\6 U_i$ are connected.
Assume that $\6 U_1\cap \6 U_2=\emptyset$. 
\begin{enumerate}
\item  
For any $x\in \R^n \backslash \6 U_1$,
consider a path $\gamma$ transversal to $\6 U_1$ and
connecting $0$ to $x$. 
Prove that the parity of the number of intersection points
between $\gamma$ and $\6 U_1$ is independent on  the choice of $\gamma$.

\item Prove that the hypersurface $\6 U_1$ is oriented.

\item
Let $\ind_{U_1}(0, x)$ be the number
of intersection points $\sharp(\gamma\cap \6 U_1)$
coun\-ted with the sign depending on the orientation. Prove that
 $\ind_{U_1}(0, x)$ is independent on  the choice of $\gamma$.

\item
Let $\gamma$ be a smooth path-connecting
$0$ to $x$, and $R_1, ..., R_n$ the connected segments
of $\gamma$ outside of $U_1$. Each $R_i$
intersects $\6 U_1$ twice, in (say) $p_i$ and $q_i$.
Prove that one can connect $p_i$ to $q_i$ by a segment that belongs to $\6 U_1$,
and this can be used to replace the path $\gamma$
with another path which intersects $\6 U_1$ at most once.

\item
Prove that the function $\ind_{U_1}(0, x)$ is constant on $\6 U_2$.
Prove that if $\ind_{U_1}(0, x)=0$ on $\6 U_2$,
then $U_2 \subset U_1$.

\item Prove that
either $U_1 \subset U_2$ or $U_2 \subset U_1$.

\item Let $X$ be a complex variety,
$z\in X$ a point, and $U_1, U_2\subset X$  connected  neighbourhoods of $z\in X$
such that the boundaries $\6 U_i$ are smooth and connected. 
Assume that $\6 U_1\cap \6 U_2=\emptyset$ and that the closures $\bar U_i$
are compact. Prove that either $U_1 \subset U_2$ or $U_2 \subset U_1$.
\end{enumerate}

\item
Let $\phi\in C^\infty M$ be a positive smooth function
 without critical points on a
compact complex Hermitian manifold with boundary.
Assume that all its level sets are compact
pseudoconvex hypersurfaces. Assume, moreover,
that $|\nabla (d\phi)|$ is globally bounded 
from above, and $|d\phi|$ is globally bounded
from below. Prove that $\phi^a$ is plurisubharmonic
for $a\in \R^{>0}$ sufficiently big.

\item 
Let $q_1, q_2$ be quadratic forms on a vector space
$W$. Suppose that $q_1$ is positive definite on
a subspace $W\subset V$, and $q_2$ vanishes on $W$
and is positive definite on $V/W$. Prove that
a linear combination $q_1 + T q_2$ is positive
definite for $T\in \R$ sufficiently big.

\item 
Let $X$ be a projective manifold, $L$ a positive line bundle
\index[terms]{bundle!line!positive}
over $X$, and $S\subset \Tot^\circ(L)$ the corresponding Sasakian
manifold, considered to be the space of unit vectors in $\Tot^\circ(L)$.
\begin{enumerate}
\item Let $\tau:\; S^1 \arrow \Aut(X)$ be a holomorphic
 isometric action of $S^1$. Assume that $L$ is $\tau$-equivariant.
Prove that the action of $\tau$ is extended to an action
of $S^1$ on $S$ by Sasakian automorphisms.

\item
Let $\gamma$ be a non-isometric homothety of 
$\Tot^\circ(L)$ obtained from the usual
$\C^*$-action, and $\rho$ the corresponding action of
$\C^*$ on the Vaisman manifold $\Tot^\circ(L)/\langle \gamma\rangle$.\index[terms]{manifold!Vaisman}
Since $L$ is $\tau$-equivariant, its action is naturally
extended to an action on $\Tot^\circ(L)$ by holomorphic isometries.
Consider the map $\rho_1(t):= \tau(t) \rho(t)$ on 
$\Tot^\circ(L)/\langle \gamma\rangle$.
Prove that $\rho_1$ is an action of $S^1$ by 
holomorphic isometries on the Vaisman manifold
$\Tot^\circ(L)/\langle \gamma\rangle$.
Prove that its lifting to $\Tot^\circ(L)$
defines a logarithm of $\gamma$, different from the standard 
(fibrewise) logarithm.

\end{enumerate}

\item\label{_Albanese_Picard_Exercise_}
\begin{enumerate}
\item
Let $X$ be a compact K\"ahler manifold.
We define {\bf the Picard group} $\Pic(X)$ of $X$ as the
group $H^1(X, \calo_X^*)$ of all holomorphic line bundles, up to isomorphism.
Using the exponential exact sequence,
\[
\arrow H^1(X, \calo_X) \arrow H^1(X, \calo_X^*)\arrow H^2(X, \Z)\arrow  
\]
 we put the structure
of a complex Lie group on $\Pic(X)$.
Prove that the connected component $\Pic_0(X)$ of the Picard group
is dual to the Albanese variety of $X$.\index[terms]{variety!Albanese}
\item
Let $X$ be a projective variety and $\sigma\in \Aut(X)$ its
automorphism. Assume $\sigma$ acts trivially on
the Albanese variety\index[terms]{variety!Albanese} of $X$ and on $H^2(X, \Z)$. Prove that any line bundle
$L$ on $X$ satisfies $\sigma^*(L)\cong L$.
\item
Let $\tau$ be a holomorphic  $S^1$-action on $X$, which has a fixed point.
Prove that any line bundle on $X$ is $\tau$-equivariant.
\item
Let $\tau$ be a holomorphic $S^1$-action on a projective manifold
$X$ that has a fixed point. Prove that there exists
a K\"ahler form\index[terms]{form!K\"ahler} $\omega$ on $X$ such that $\tau$ acts on
$(X, \omega)$ by isometries.
\end{enumerate}

\item
Let $X$ be a projective manifold admitting an action $\tau$ of $S^1$,
that fixes a point $x\in X$. Prove that any ample line bundle
on $X$ admits a $\tau$-equivariant structure (\ref{_equivariant_structure_Definition_}) 
and an $S^1$-equivariant 
Hermitian metric with positive curvature.
Consider a regular Vaisman manifold\index[terms]{manifold!Vaisman!regular} $M=\Tot^\circ(L)/\Z$
associated with the standard homothety action $\gamma$ on $\Tot^\circ(L)$.\index[terms]{manifold!Vaisman}
Prove that $\gamma$ has two different logarithms,
that is, there are two non-equal one-pa\-ra\-me\-tric subgroups
$\R \subset \Aut(\Tot^\circ(L))$ which contain $\gamma$.

{\em Hint:} Use the previous exercise.

\item\label{_Kodaira_surface_Exercise_}
Let 
$M:=\Tot^\circ(L)/\langle u \rangle$,
where $L$ is an ample line bundle over an elliptic curve $E$,
and $u \in \C$, $|u|>1$. Then $M$ is called {\bf a Kodaira surface}.\index[terms]{surface!Kodaira}
\begin{enumerate}
\item Prove that the canonical foliation\index[terms]{foliation!canonical} $\Sigma$ on $M$
is the foliation tangent to the fibres of $L$.

\item Prove that any holomorphic automorphism 
of $M$ is compatible with the projection to $E$.

\item Let $G$ be the image of $\Aut(M)$ in the group of
  automorphisms of $E$. Prove that the bundle $L$
is $G$-equivariant.

\item
Prove that $G$ preserves a point $[L]$ in
$\Pic_d(E)$, where $d=\deg L$ and
$\Pic_d(E)$ is the space of line bundles
of degree $d$, considered  a connected
component of the Picard variety.

\item
Prove that the subgroup of $\Aut(E)$
fixing a point $[L]\in \Pic_d(E)$, $d>0$,
is finite.
\end{enumerate}

{\em Hint:} Use Exercise \ref{_Albanese_Picard_Exercise_}.

\item
Let $M$ be a  Kodaira surface. Prove that the
group $\Aut(M)$ of holomorphic automorphisms of $M$ is compact.\index[terms]{surface!Kodaira}

{\em Hint:} Use the previous exercise.

\item\label{_secondary_Kodaira_Exercise_}
Let $E$ be an elliptic curve.
Consider a finite cyclic group $G= \Z/k\Z$ acting on 
$E$, $k=2, 3, 4, 6$, fixing a point,
and let $\sigma$ be its generator.
\begin{enumerate}
\item Prove that for any ample line bundle $L$ on $E$,\index[terms]{bundle!line!ample}
the bundle $\bigotimes_{g\in G} g^*(L)$ is $G$-equivariant and ample.

\item 
Let $L$ be a $G$-equivariant ample bundle on $E$, and let 
$s:\; \sigma^* L \arrow L$ be the line bundle isomorphism associated with\index[terms]{equivariant structure}
the equivariant structure (\ref{_equivariant_structure_Definition_}). Given a number $\lambda \in \C$, $|\lambda|>1$,
consider the automorphism $\gamma\in \Aut(\Tot^\circ(L))$ taking
$(z, v)$ to $(\sigma^{-1}(z), \lambda s(v))$, where $z\in E$ and $v\in L\restrict z$.
Prove that $L$ admits a $G$-invariant Hermitian metric
with positive curvature. Let $l(v):= |v|^2$ be the length
function on $\Tot^\circ(L)$. Prove that $l$ is strictly
plurisubharmonic, and that the Hermitian metric $\frac{dd^c l}{l}$
is $\gamma$-invariant. Prove that the corresponding
metric on $M:=\Tot^\circ(L)/\langle\gamma\rangle$ 
is Vaisman. 

\item Prove that a $G$-cover of $M$ is 
 $M_1:=\Tot^\circ(L_1)/\langle u \rangle$
where $L_1= L^{\otimes k}$ is a positive line bundle, and
$u \in \C$, $|u|>1$ is the standard homothety.
Prove that $M_1$ is a regular Vaisman\index[terms]{manifold!Vaisman!regular} 
manifold, equipped with a free $G$-action,
and that $M_1/G = M$.\footnote{A Vaisman manifold \index[terms]{manifold!Vaisman}
$M_1:=\Tot^\circ(L_1)/\langle u \rangle$,
where $L$ is an ample line bundle over an elliptic curve
and $u \in \C$, $|u|>1$ 
is called {\bf a primary Kodaira surface},
or just {\bf Kodaira surface}. For any finite group $G$
freely acting on $M_1$, the quotient\index[terms]{surface!Kodaira!primary}
$M_1/G$ is {\bf a secondary Kodaira surface}\index[terms]{surface!Kodaira!secondary}
if the leaf space of the canonical foliation on $M_1/G$
is $\C P^1$.\index[terms]{surface!Kodaira!secondary}}

\item
Prove that 
the $\Z$-action\index[terms]{action!$\Z$-} on $\Tot^\circ(L)$
is not homotopy equivalent to zero,
and the corresponding Vaisman manifold
does not admit a logarithm.
\end{enumerate}

\item
Prove that the complex surface obtained 
in Exercise \ref{_secondary_Kodaira_Exercise_}
is a $T^2$-fibration over $\C P^1$ and show
that it has the same de Rham cohomology as the Hopf surface.
Prove that it is not homeomorphic to a Hopf surface.\index[terms]{surface!Hopf}
\index[terms]{cohomology!de Rham}

{\em Hint:} Compute its fundamental group.\index[terms]{fundamental group}

\item 
Construct two LCK structures with potential\index[terms]{structure!LCK!with potential}
$(\theta, \omega)$ and $(\theta', \omega')$ on a 
Kodaira surface in such a way that the cohomology
classes $[\theta], [\theta']\in H^1(M,\R)$ are\index[terms]{surface!Kodaira}
not proportional. 

\item
Construct an LCK structure with improper potential
on a  Kodaira surface.\index[terms]{surface!Kodaira}

\item
Let $(M, \theta, \omega)$ be an LCK manifold with
potential and preferred gauge, 
$\pi:\; \tilde M\arrow M$ a K\"ahler $\Z$-cover\index[terms]{cover!K\"ahler $\Z$-} of $M$,
and $\vec r\in T\tilde M$ a logarithm of the $\Z$-action\index[terms]{action!$\Z$-}.
Assume that $\theta$ is $\vec r$-invariant. Prove that
the LCK metric on $M$ is obtained from\index[terms]{logarithm!of a $\Z$-action}
a pseudoconvex shell\index[terms]{pseudoconvex shell} as in \ref{_metric_from_shell_Definition_}.

\item \label{_LCK_proper_logarithm_inva_confo_Exercise_}
Let $(M, \theta, \omega)$ be an LCK manifold with proper
potential, and $\rho$ an $S^1$-action on $M$ obtained
from the logarithm of the monodromy \index[terms]{action!monodromy}as in 
\ref{_LCK_with_pot_from_shell_Theorem_}.
\begin{enumerate}
\item Prove that $\omega$ is conformally equivalent to 
a metric obtained from a pseudoconvex shell
if and only if there exists a $\rho$-invariant
Hermitian metric in the same conformal class
as $\omega$.

\item Consider an LCK potential\index[terms]{potential!LCK} on the $\Z$-cover\index[terms]{cover!K\"ahler $\Z$-} $\tilde M$
that is  not $e^{t\vec r}$-auto\-mor\-phic, and let $(M, \omega, \theta)$ be
the corresponding LCK structure with preferred gauge.
Assume that any auto\-mor\-phic LCK-potential $\phi$ with preferred gauge
and the Lee form \index[terms]{form!Lee}cohomologous to $\theta$
is uniquely determined by the LCK structure\index[terms]{structure!LCK}
(see Exercise \ref{_uniqueness_of_potential_Exercise_} for an example)
Prove that $\theta$ and $\omega$ are not
$S^1$-invariant.
\end{enumerate}

\item\label{_exist_unique_preferred_gauge_by_Lee_class_Exercise_}
Using Exercise
\ref{_uniqueness_of_potential_Exercise_}, prove that there exists an
LCK manifold\index[terms]{manifold!LCK} with the Lee cohomology class $[\theta]$
such that for any LCK form $\omega$ with potential and\index[terms]{form!LCK!with potential}
the same Lee cohomology class, the potential $\phi$ with preferred
gauge is uniquely determined by $\omega$.

\item
Find an LCK manifold $(M, \omega, \theta)$ with
proper potential, such that its K\"ahler $\Z$-cover\index[terms]{cover!K\"ahler $\Z$-} $\tilde M$
admits the logarithm of the $\Z$-action,\index[terms]{action!$\Z$-} but
$\omega$ is not conformally equivalent to 
a metric obtained from a pseudoconvex shell.\index[terms]{pseudoconvex shell}

{\em Hint:} Use Exercises \ref{_LCK_proper_logarithm_inva_confo_Exercise_}
and   \ref{_exist_unique_preferred_gauge_by_Lee_class_Exercise_}.

\item
Let $\tilde M$ be a K\"ahler $\Z$-cover\index[terms]{cover!K\"ahler $\Z$-} of an LCK manifold 
$(M, \theta, I)$ with potential,\index[terms]{manifold!LCK!with potential}
 $\rho_t$ the logarithm flow, and $S= \psi^{-1}(1)$ the
corresponding pseudoconvex shell. Consider a holomorphic map
$\Phi:\; \tilde M \arrow \tilde M$ commuting with $\rho_t$
and preserving $S$. 
\begin{enumerate}
\item Prove that $\Phi$ commutes with the
$\Z$-action\index[terms]{action!$\Z$-} on $\tilde M$.
\item Prove that $\Phi$ is induced by a holomorphic
isometry of $(M, \theta, I)$.
\end{enumerate}

\item\label{_example_of_isometry_shells_Exercise_}
Let $H$ be a linear 
Hopf manifold, $H=\frac{\C^n \backslash 0}{\langle A\rangle}$
and $V_t:\; \C^n \arrow \C^n$ a flow of isometries of the standard
metric on $\C^n$ commuting with $A$.
Assume that the sphere in $\C^n$ is a shell for the
logarithmic action $e^{t\log A}$. \index[terms]{action!logarithmic}
\begin{enumerate}
\item Prove that $V_t$ acts on the 
corresponding LCK manifold\index[terms]{manifold!LCK} by holomorphic isometries.
\item Find an example of such a flow $V_t$ when $A$ is
non-diagonalizable, and the  sphere in $\C^n$ is a shell for the
 action $e^{t\log A}$. 
\end{enumerate}

\end{enumerate}


\chapter{Embedding theorem for Vaisman manifolds}\index[terms]{manifold!Vaisman}
\label{_embe_theorem_Vaisman_Chapter_}

\epigraph{\it Gliding o'er all, through all,\\
	Through Nature, Time, and Space,\\
	As a ship on the waters advancing,\\
	The voyage of the soul -- not life alone,\\
	Death, many deaths I'll sing.}{\sc\scriptsize Walt Whitman, \ \ Gliding O'er All}

\section{Introduction}

In this chapter, we work only with LCK manifolds\index[terms]{manifold!LCK}
of complex dimension $>2$. We take $\dim_\C M>2$ as
a running assumption.

Previously, we have established that the compact
LCK manifolds with potential \index[terms]{manifold!LCK!with potential}are precisely the
manifolds that admit a complex embedding to a Hopf
manifold (\ref{embedding}, \ref{_linear_LCK_pot_Corollary_}).

The Vaisman manifolds are an important subclass of
LCK manifolds with potential. It turns out that 
 the Vaisman manifolds are precisely the LCK
manifolds with potential that can be embedded
to a ``diagonal Hopf manifold'', that is,
a linear Hopf manifold $\frac{\C^n\backslash 0}{\langle A \rangle}$,
where $A\in \GL(\C, n)$ is a diagonalizable linear
operator.

We start by proving this result for a Hopf manifold,
showing that any diagonal Hopf manifold is Vaisman
(\ref{semihopf}). Conversely, 
\ref{vaisman_embed} implies that a non-diagonal
Hopf manifold is not Vaisman.

In this chapter we focus on proving that Vaisman 
manifolds are embeddable to a diagonal Hopf manifold.

In \ref{_cone_cover_for_LCK_pot_Theorem_}, 
we associated an algebraic cone\index[terms]{cone!algebraic}
to any compact LCK manifold with potential:\index[terms]{manifold!LCK!with potential}
this manifold can be obtained as a $\Z$-quotient
of this cone. In fact, the $\Z$-action\index[terms]{action!$\Z$-} is sufficient
to recover the algebraic structure. Let $\tilde M$
be a $\Z$-cover\index[terms]{cover!K\"ahler $\Z$-} of an LCK manifold\index[terms]{manifold!LCK}, and $\gamma$
the generator of $\Z$-action. Then the algebraic
functions on $\tilde M$ are precisely holomorphic
functions, which are  {\bf $\gamma$-finite},\index[terms]{function!$\gamma^*$-finite} that
is, are contained in a finite-dimensional
$\gamma$-invariant subspace (\ref{_cone_cover_for_LCK_pot_Theorem_}).

We call a compact LCK manifold with potential\index[terms]{manifold!LCK!with potential}
{\bf semisimple} if the action of $\gamma$\index[terms]{manifold!LCK!with potential}
on the space of $\gamma$-finite functions on
$\tilde M$ is diagonalizable. Clearly,
a linear Hopf manifold is semisimple if
and only if it is diagonal.

In \ref{vaisman_embed}, 
we prove that $M$ is semisimple
if and only if it is Vaisman. Since the embedding
to a Hopf manifold is constructed using a
sufficiently big finite-dimensional space
of $\gamma$-finite functions,\index[terms]{function!$\gamma^*$-finite} this also shows
that a manifold is Vaisman if and only if it
can be embedded to a diagonal Hopf manifold 
(\ref{_non-Vaisman_Corollary_}). 

The term ``semisimple'' comes from the 
Jordan--Chevalley theorem, that is  a part\index[terms]{decomposition!Jordan--Chevalley}
of algebraic group theory. Recall that
{\bf an algebraic group} is a subgroup of $\GL(n)$
given by polynomial equations. 
An element $g$ of an algebraic group is
called {\bf semisimple} if it is
diagonalizable over the algebraic 
closure of the base field, and
{\bf unipotent} if $g-\Id$ is nilpotent.

The Jordan--Chevalley theorem\index[terms]{theorem!Jordan--Chevalley decomposition}
claims that an element $g$ of an 
algebraic group can be decomposed
as $g=su$, where $s$ is semisimple,
$u$ unipotent, and $su=us$.
Moreover, this decomposition is
unique, and compatible with the
morphisms of algebraic groups.

Consider the adjoint action
of $\GL(n)$ on itself, and
let $\Ad_{\GL(n)}(g)$ be the adjoint orbit
of $g$, with $g=su$ its Jordan--Chevalley
decomposition. It is not hard to see that
$s$ belongs to the closure of $\Ad_{\GL(n)}(g)$ 
(\ref{_semisimple_operator_approx_Proposition_}).

This observation is used to show that
any LCK manifold with potential \index[terms]{manifold!LCK!with potential}can
be deformed to a Vaisman manifold\index[terms]{manifold!Vaisman}
(\ref{def_lckpot2Vai}).

\section[Embedding Vaisman manifolds to Hopf
manifolds]{Embedding Vaisman manifolds to\\ Hopf
  manifolds}\label{embedding_vaisman_hopf}

In this section, we prove the following equivalence: a compact LCK
manifold with proper potential is Vaisman if and only if
it can be embedded in a diagonal Hopf manifold.
\index[terms]{manifold!LCK!with potential!LCK!proper}

\subsection{Semisimple Hopf manifolds are Vaisman}

Recall that the linear operator $A\in\End(V)$ is called {\bf semisimple}
	if it is diagonalizable over the algebraic 
closure \index[terms]{operator!semisimple}
        of $k$. 

\hfill

\definition\label{_diagonal_Hopf_Definition_}
A {\bf diagonal} (or {\bf semisimple}) 
{\bf Hopf manifold} is a quotient of
$\C^n \backslash 0$ by a semisimple (that is, 
diagonalizable) linear automorphism with all
eigenvalues $|\alpha_i|<1$. \index[terms]{manifold!Hopf!diagonal}

\hfill

\definition 
	Let $M$ be an LCK manifold with potential,\index[terms]{manifold!LCK!with potential} 
	and $j:\; M \arrow H$ be a holomorphic embedding to a Hopf
	manifold $H := \frac{\C^n \backslash 0}{\langle A\rangle}$.
	Then $M$ is called {\bf semisimple} if $A$ is
	semisimple.

\hfill

\theorem \label{semihopf}
A diagonal linear Hopf manifold admits a Vaisman 
metric. \index[terms]{manifold!Hopf!linear!semisimple}

\hfill

\proof Let $V:=\C^n$ and 
$H:=\frac{V\backslash0}{\langle A\rangle}$ be a semisimple Hopf manifold,
let	$\{e_i\}$ be an eigenvalue basis in $V$.
Representing the eigenvalues as a product of 
an absolute value and the argument, we write
$A(e_i)=\alpha_i u_i e_i$,
	with $\alpha_i \in ]0, 1[$ and $u_i \in \U(1)$.
	Consider the unit sphere $S\subset V $, and 
	let $\rho(t)(e_i):= \alpha_i^t e_i$. 
Clearly, $S$ is a pseudoconvex shell\index[terms]{pseudoconvex shell} for the
action of $\rho(t)$.
By \ref{_LCK_potential_from_shells_Theorem_}, there exists
	a $\rho$-automorphic K\"ahler potential $\phi_\rho$
	on $V\setminus 0$ with $\phi_\rho\restrict S=1$. 
Since $A\circ \rho(-1)$ preserves $S$
	and $A$ commutes with $\rho(t)$,
	the function $\phi_\rho$ is $A$-automorphic.
	
By \ref{_Vaisman_via_logarithm_Theorem_}, to prove that
$\phi_\rho$ defines a Vaisman structure, it is
enough to check that $\phi_\rho$ is invariant
with respect to the action of $I(\vec r)$,
where $\vec r$ is the vector field tangent to
$\rho(t)$.  

This is clear, because 
$I(\vec r)$ preserves the sphere $S$ and
commutes with $\rho(t)$, and $\phi_\rho$ is determined
by $S$ and $\rho(t)$.
To see that the exponent of $I(\vec r)$ preserves $S$,
observe that
$e^{t I\vec r}(e_i)= e^{\1 t\log \alpha_i} e_i$,
	and this operator is unitary.
	\endproof

\hfill

\remark 
For primary diagonal Hopf surfaces, this result was proven
directly, by exhibiting an automorphic potential on
$\C^2\backslash 0$, in \cite{go}. See also \cite{bel},
\cite{kor} and \cite{ov_pams}.	\index[terms]{surface!Hopf!diagonal!}

\hfill

\remark \ref{semihopf} has a stronger form (\ref{_Hopf_Vaisman_diagonal_Corollary_})
stating that a linear Hopf manifold is Vaisman type if and
only if it is diagonal.


\subsection[Algebraic groups and Jordan--Chevalley de\-com\-po\-si\-tion]{Algebraic groups and Jordan--Chevalley\\ de\-com\-po\-si\-tion}\label{alg_gr_jc}


We recall the following definitions:

\hfill

\definition 
	\begin{enumerate}
	\item The {\bf category of affine varieties} is\index[terms]{variety!affine}
          the category of finitely generated rings without
          nilpotent elements\index[terms]{nilpotent element}  with arrows inverted.
		\item An {\bf algebraic group} is a group object in the
category of affine varieties.\index[terms]{group!algebraic}\index[terms]{group!pro-algebraic}

\item A {\bf pro-algebraic group} is an inverse limit
of algebraic groups. 
\end{enumerate}

Further on, all algebraic groups are
considered over $\C$. For more reference
on algebraic groups, please see \cite{hum}.

\hfill

\definition 
An element of an algebraic group $G$ is called
{\bf semisimple} if its image is semisimple
for some exact\index[terms]{semisimple element}
algebraic representation of $G$, and is called 
{\bf unipotent} if its image is unipotent\index[terms]{unipotent element}
(that is, exponential of a nilpotent) \index[terms]{nilpotent element}
for some exact algebraic representation of $G$.

\hfill

\remark 
For any algebraic representation of an algebraic group
$G$, the image of any semisimple element is a 
semisimple operator, and the image of any unipotent
element is a unipotent operator (\cite[\S 15.3]{hum};
for an alternative argument, see the exercises to 
this chapter).

\hfill

\theorem  ({\bf Jordan--Chevalley  decomposition}, 
\cite[\S 15.3]{hum}) \label{jcdec}\\
Let $G$ be an algebraic group, and $A\in G$.
{Then there exists a unique decomposition $A= S U$\index[terms]{Jordan--Chevalley decomposition}
	of $A$ in a product of commuting elements $S$ and $U$,
	where $U$ is unipotent and $S$ semisimple.}\index[terms]{theorem!Jordan--Chevalley decomposition}

\hfill

\remark 
Since this decomposition is unique, it is functorial.
{Therefore, it is also true for all pro-algebraic groups.}

\subsection{The algebraic cone of an LCK
  manifold with potential}\index[terms]{manifold!LCK!with potential} 
\remark 
	Let $M$ be an LCK manifold with proper potential, $\tilde M$
	its K\"ahler $\Z$-covering\index[terms]{cover!K\"ahler $\Z$-}, $M= \tilde M/\langle
	A\rangle$, $\tilde M_c$ its completion, and $j:\; M \hookrightarrow H$ a holomorphic embedding to a Hopf
	manifold $H = \C^n \backslash 0/\langle A\rangle$ as in \ref{embedding}. 
	
	Let $R$ be the  $\goth
	m$-adic completion of $\calo_{\tilde M_c}$ in
	the maximal ideal of the origin $x_0 \in \tilde M_c$,
	$R= \lim\limits_{\longleftarrow}\ \calo_{\tilde M_c}/ {\goth m}^k$,
	and $G:= \Aut(R)$. Then  
	\[ 
	G= \lim\limits_{\longleftarrow}\ \Aut(\calo_{\tilde M_c}/ {\goth m}^k),
	\]
	hence $G$ is a pro-algebraic group. \index[terms]{completion!adic}

\hfill

Using this notion, we  interpret the semisimplicity of an
LCK manifold with proper potential as follows.

\hfill

\proposition \label{_LCK_pro_algebraic_linear_Proposition_}
Let  $j: M \hookrightarrow H$ be a holomorphic embedding of an LCK manifold\index[terms]{manifold!LCK}
$M= \tilde M/\langle A\rangle$
to a Hopf manifold $H = V\setminus 0/\langle A\rangle$,
such that $V$ is an $A$-invariant subspace in $\calo_{\tilde M_c}$.
{Then the action of $A$ on $V$ is semisimple if and only if
	$A$ is semisimple as an element of the pro-algebraic group \index[terms]{group!pro-algebraic}
	$G= \lim\limits_{\longleftarrow}\ \Aut(\calo_{\tilde M_c}/ {\goth m}^k)$.}

\hfill

\proof If $A$ is semisimple as an element of $G$, its action
on $V$, considered to be an $A$-invariant subspace in
$R\supset \calo_{\tilde M_c}$, is also semisimple.

Conversely, if $A$ is semisimple on $V$, then $\calo_{\tilde M_c}$
is a subring in $R$, that is  a quotient ring of
$\C[[V]]$. The latter is the ${\goth m}$-adic completion of the polynomial
ring $\C[V]$, where $A$ is  clearly semisimple.
\endproof

\hfill

\definition
Let $M$ be an LCK manifold with proper potential\index[terms]{manifold!LCK!with potential}
and $\tilde M_c$ the corresponding algebraic cone.
The manifold $M$ is called {\bf semisimple}
if the $\Z$-action \index[terms]{action!$\Z$-}on $\tilde M_c$ is semisimple,
that is, its image in the pro-algebraic
group of the automorphisms of the completion of the cone
is semisimple.

\hfill

We are going to show that semisimplicity is equivalent 
with the existence of a Vaisman metric.

\hfill

\theorem { (\cite{ov_lckpot})}\label{vaisman_embed} 
Let $M$ be a compact LCK manifold with proper potential. Then $M$
	is semisimple if and only if  it is Vaisman.

\hfill

\proof 
Let $M$ be a Vaisman manifold of LCK rank 1, \index[terms]{rank!LCK}
$\tilde M$ its $\Z$-cover\index[terms]{cover!K\"ahler $\Z$-},\index[terms]{manifold!Vaisman}
and $\tilde M_c$ its completion, equipped
with the K\"ahler potential 
$\phi:\; \tilde M_c \arrow \R^{> 0}$. 
Consider a compact subset $\tilde M_c^a:= \phi^{-1}([0,a])$.
We introduce an $L^2$-structure\index[terms]{structure!$L^2$} on  the ring 
$H^0(\calo_{\tilde M^a_c})$ of holomorphic functions by 
$|f|^2 = \int_{\tilde M^a_c}|f|^2 \tilde\omega^n$.

Let $\theta^\sharp$ be the Lee field\index[terms]{Lee field} of $M$ lifted to $\tilde M$.
Since the field
$I\theta^\sharp$ is Killing,\index[terms]{vector field!Killing} it acts on $H^0(\calo_{\tilde M^a_c})$
	by isometries, and hence  its action on each finite-dimensional
subspace of $H^0(\calo_{\tilde M^a_c})$ is semisimple.

Since $\theta^\sharp$ is holomorphic, its action 
on $H^0(\calo_{\tilde M^a_c})$ is given by
\[ \Lie_{\theta^\sharp} f = -\1 \Lie_{I\theta^\sharp} f;\]
hence, it is also semisimple.

Let $G\subset \Iso(M)$ be the closure of the complex
one-parametric group of isometries generated by 
the flow of the Lee field\index[terms]{Lee field} $\theta^\sharp$ and $I\theta^\sharp$ 
on $M$, and $\tilde G$ its lift to $\tilde M$.
Then the action of all
vector fields in $\Lie(\tilde G)$ on $H^0(\calo_{\tilde M^a_c})$
is semisimple. Indeed,  any 
$v\in \Lie(\tilde G)$ acts on $\tilde M$ by homotheties;
composing $e^{v}$ with $e^{t\theta^\sharp}$ for
appropriate $t$, we obtain an element
$g= e^v e^{t\theta^\sharp}$ acting on $\tilde M$
by isometries, and hence  it is semisimple on 
$H^0(\calo_{\tilde M^a_c})$. Then 
$e^v= g e^{-t\theta^\sharp}$ is semisimple as a product
of two commuting semisimple operators.

Let $\gamma$ be the generator
of the monodromy action\index[terms]{action!monodromy} on $\tilde M$,
with $M = \tilde M/\langle \gamma\rangle$.
To prove that $M$ is semisimple, 
we need to show that the action of $\gamma$
on $\calo_{\tilde M}$ is semisimple.
By \ref{_Vaisman_Lee_action_contains_monodromy_Corollary_},
the group $\tilde G$ contains a power of $\gamma$,
and its action is semisimple, and hence  $\gamma^k$
is semisimple; however, a power of a non-semisimple
operator cannot be semisimple. This implies that
$\gamma$ is also semisimple.

To prove that semisimple implies Vaisman, recall that 
 all compact complex subvarieties of Vaisman manifolds are
 again Vaisman. Then the statement follows from
 \ref{semihopf}.\index[terms]{manifold!Vaisman}
 \endproof

\hfill

\corollary\label{_non-Vaisman_Corollary_}
Let $M\subset H$ be an LCK manifold with
potential embedded to a linear Hopf manifold $H$.
Then $M$ is Vaisman if and only if $H$ is diagonal.

\hfill

\proof
If $H$ is diagonal, $M$ is Vaisman by \ref{semihopf}.
Conversely, suppose that the Hopf manifold 
$H=\frac{\C^n\backslash 0}{\langle  A\rangle}$  
is {\em not} diagonal, and let $\gamma= A\restrict {\tilde M}$
be the generator of the $\Z$-action \index[terms]{action!$\Z$-}on the corresponding
K\"ahler $\Z$-covering\index[terms]{cover!K\"ahler $\Z$-} $\tilde M\subset \C^n \backslash 0$ of $M$.

The space $V^*$
of linear functions on $\C^n=V$ gives
a finitely--dimensional $\gamma$-invariant
subspace in the ring $\calo^\gamma_{\tilde M}$
of $\gamma$-finite functions \index[terms]{function!$\gamma^*$-finite}on $\tilde M$.
If $H$ is not diagonal, the action of
$\gamma$ on $V^*$ is not semisimple, and 
hence $M$ is not semisimple.
\endproof

\hfill

\corollary \label{_Hopf_Vaisman_diagonal_Corollary_}
A linear Hopf manifold admits a Vaisman structure if and only if it is diagonal.
\endproof

\hfill

\remark
For complex surfaces this was proven by
F. \index[persons]{Belgun, F. A.} Belgun (\cite{bel}). Note that none of the
arguments given in this chapter can be applied
to $\dim =2$, because these arguments rely on 
the structure of an algebraic cone on $\tilde M$,
and this, in turn, relies on the  Rossi and  Andreotti--Siu
theorem, that is  known to be false in $\dim =2$.

\hfill

As an application, we give an explicit description
of the canonical foliation\index[terms]{foliation!canonical} on the diagonal Hopf manifold.\index[terms]{manifold!Hopf!diagonal}
Since we are not going to use this result, in lieu of the proof
we refer to Exercise \ref{_linear_Hopf_canonical_explicitly_Exercise_}
below, where we suggest a way of proving it, using a step-by-step argument.

\hfill

\proposition\label{_canonical_explicit_Hopf_Proposition_}
Let $M$ be a linear Hopf manifold, 
$M=\frac{\C^n \backslash 0}{\langle A\rangle}$,
where $A(e_i) =\alpha_i e_i$, with
$\alpha_i \in \C$ and $e_i$ the basis in $\C^n$.
Choose a Vaisman metric on $M$, and let $\Sigma\subset TM$
be the canonical foliation. 
Let $\log A_\R$ be the linear map given by 
$\log A_\R(e_i) =\log |\alpha_i| e_i$. Then the leaves of 
$\Sigma$ are precisely the orbits of $e^{\C \log A_\R}$.
\endproof

\section[Deforming an LCK manifold with proper potential to Vaisman manifolds]{Deforming an LCK manifold with proper\\ potential to Vaisman manifolds}\index[terms]{manifold!Vaisman}\index[terms]{manifold!LCK!with potential}

As a first application of \ref{vaisman_embed}, we prove
that, in dimension greater than 3, the complex structure
of a compact LCK manifold with potential\index[terms]{manifold!LCK!with potential} can be deformed
to one admitting a Vaisman metric. Therefore, any
LCK manifold with potential is diffeomorphic to a Vaisman
manifold. We start with the following algebraic
notions.

\hfill

\definition
Let $U \subset X$ be a subset of an affine algebraic variety.
{\bf The Zariski closure} $\bar U$ of $U$ is the smallest closed algebraic
subvariety of $X$ containing $U$. In other words,
$\bar U$ is the set of all common zeros for all\index[terms]{Zariski closure}
regular functions $f\in \calo_X$ vanishing  in $U$.

\hfill

\remark 
Let $G$ be an algebraic group.
Then the Zariski closure of a subgroup $\Gamma\subset G$
is an algebraic subgroup of $G$ (see Exercise 
\ref{_Zariski_closure_of_subgroups_Exercise_}).

\hfill

\theorem  { (\cite{ov_imrn_10})} \label{def_lckpot2Vai}
  Let $(M,\omega,\theta)$, 
$\dim_\C M \geq 3$, be a compact LCK manifold
with proper potential.\index[terms]{manifold!LCK!with potential} Then there exists a complex
analytic deformation of $M$ which
admits a Vaisman metric. \index[terms]{deformation}

\hfill

\proof The idea  is to embed $M$ into a Hopf manifold
defined by a
linear operator $A$, then, using the Jordan--Chevalley decomposition of $A$
to show
that its semisimple part preserves the subvariety \index[terms]{decomposition!Jordan--Chevalley}
$\tilde M \subset\C^n\setminus 0$, thus
yielding an embedding of a small  deformation of $M$ into a  diagonal Hopf
manifold,
the Vaisman metric of that can be pulled back on $M$. We now provide the
details.\index[terms]{metric!Vaisman}

\hfill

{\bf Step 1.} Let $V= \C^n$, $A\in \End(V)$ be an invertible
linear operator
with all eigenvalues $|\alpha_i|<1$, and $H := \frac{V\setminus  0}{\langle
A\rangle}$  the corresponding Hopf manifold. 
One may see that the complex submanifolds
of $H$ are identified with complex subvarieties
$Z$ of $V$, that are  smooth outside of  $0$ and are fixed by $A$.
Indeed, by  Remmert--Stein theorem 
(\cite[Chapter 2, \S 8.2]{demailly}\index[terms]{theorem!Remmert--Stein}), 
for every complex subvariety $X\subset H$,
the closure of $\pi^{-1}(X)$ is complex analytic in $V=\C^n$,
where $\pi:\; V\setminus 0 \arrow H$ is the natural projection.

\hfill

{\bf Step 2.}
We are going to prove that any such $Z$ is fixed by the 
group $G_A$ defined as the Zariski closure of 
$\langle A\rangle$ in $\GL(V)$. Denote by ${\goth m}\subset \calo_V$ 
the maximal ideal of 0. Let $I_Z$ be \index[terms]{Zariski closure}
the ideal of $Z$, and let 
 $\hat I_Z$ be the corresponding ideal in the ${\goth m}$-adic
completion 
$\hat \calo_{V}\supset \calo_V$ of the structure ring $\calo_{V}$ in $0$. 
To prove that $I_Z$ is fixed by $G_A$,
it is enough to show that $\hat I_Z$ is fixed by $G_A$.
Indeed, $I_Z= \hat I_Z \cap \calo_V$, and
$\calo_V$ is fixed by $G_A$ by construction.

Recall that, by definition (see for example \cite[Chapter 10]{_Atiyah_MacDonald_}),   
$\hat I_Z$ is the inverse
limit\index[terms]{inverse limit}
of the projective system $\displaystyle\frac{I_Z}{I_Z \cap {\goth m}^k}$ :
\[ \hat I_Z= \lim\limits_{\longleftarrow} \frac{I_Z}{I_Z \cap {\goth m}^k}.
\]
To prove that $\hat I_Z$ is fixed by $G_A$
it only remains to show that $\displaystyle\frac{I_Z}{I_Z \cap {\goth m}^k}$
is fixed by $G_A$. However,
$\displaystyle\frac{I_Z}{I_Z \cap {\goth m}^k}$
is a subspace in the  vector space
$\calo_{V}/{\goth m}^k$, finite-di\-men\-si\-o\-nal by
\cite[Corollary 6.11 and Exercise 8.3]{_Atiyah_MacDonald_}.
A finite-dimensional subspace that is  fixed by $A$ is
automatically fixed by the Zariski closure  $G_A$.\index[terms]{Zariski closure}

\hfill

{\bf Step 3.} Now, for any linear operator there
exists a unique decomposition $A:= S U$ in a product
of commuting operators, with $S$ semisimple (diagonal),
and $U$ unipotent, i.\,e.     its spectrum contains only\index[terms]{spectrum}
the number $1$. Then $S$ acts on $V$ by linear contractions.
By \ref{jcdec}, for any finite-dimensional representation of $\GL(n)$,
any vector subspace that is  fixed by $A$ is also fixed
by $S$.
By the argument in Step 2, this proves that  $S$
fixes the ideal $\hat I_Z$, and the subvariety
$Z\subset V$.

\hfill

{\bf Step 4.} 
The eigenvalues of $S$ are the same as the eigenvalues of $SU=A$;
hence, the manifold $H_S:= (V\setminus  \{0\})/\langle S\rangle$ 
is diagonal Hopf. By construction, $H_S$ contains
$M_0 := (Z\setminus  0)/\langle S\rangle$, that is 
Vaisman because it is a complex submanifold of a
Vaisman manifold $H_S$ (\ref{semihopf}).\index[terms]{manifold!Vaisman} 
Let $U^t$ denote $e^{t\log U}$; this makes
sense because $U$ is unipotent.
Since $Z$ is $S$-invariant and $U$-invariant (Step 3),
it is $SU^t$-invariant. Then
$M_t := (Z\setminus  0)/\langle SU^t\rangle$
gives a continuous family of LCK manifolds
with potential\index[terms]{manifold!LCK!with potential} connecting $M_0$ (that is  Vaisman)
to $M$. \endproof


\hfill

\remark \label{_Vaisman_defo_same_cone_Remark_}
The Vaisman deformation of an LCK manifold $M$ with proper potential\index[terms]{manifold!LCK!with potential} 
 has the same cone covering $\tilde M$,\index[terms]{cover!cone} but with a different $\Z$-action.\index[terms]{action!$\Z$-}

\section{Exercises}

\begin{enumerate}[label=\textbf{\thechapter.\arabic*}.,ref=\thechapter.\arabic{enumi}]


%




\item Let $A\in\GL(n,\C)$ be an endomorphism
 of finite order. Prove that $A$ is semisimple.

\item Let $G$ be an algebraic group over $\C$ \index[terms]{group!algebraic}
that does not contain unipotent elements. Prove that
$G$ is commutative, or find a counterexample.

\item\label{_Zariski_closure_of_subgroups_Exercise_}
Let $G$ be an algebraic group, and $\Gamma\subset G$
a subgroup, that is, a subset that is  closed under the group
operations. 
\begin{enumerate}
\item Prove that ``being closed under the group operations''
is an algebraic condition, that is, a condition given by
polynomial equations.
\item
Prove that the Zariski closure of $\Gamma$
in $G$ is an algebraic subgroup.
\end{enumerate}

\item
Let $G$ be an algebraic group. Given $u\in G$,
consider the element $u_0:=\Id-u$ in the group algebra
$\C[G]$. 
\begin{enumerate}
\item
Prove that $u$ is unipotent if and only if
$u_0^k=0$, for $k$ sufficiently big.
\item Let $G\arrow G_1$ be a group homomorphism.
Prove that the image of an unipotent element is unipotent.
\end{enumerate}

\item\label{_closure_of_jordan_cell_Exercise_}
Let {\[\scriptsize 
p =\begin{pmatrix} \alpha & 1 & 0 & \ldots & 0\\
0 & \alpha & 1 & \ldots & 0\\
\vdots &\vdots &\vdots & \cdots & \vdots \\
0&0&0 & \ldots &1\\
0&0&0 & \ldots &\alpha
\end{pmatrix}
\]}
be a Jordan cell in $\GL(V)$, with $\alpha$ of infinite order,\index[terms]{Jordan cell}
and $\C[p]\subset \Mat(V)$ the group algebra generated by $p$.
\begin{enumerate}
\item Prove that the closure of $\C[p]$ (in the usual topology)
contains the diagonal matrix $\alpha \Id$.
\item Prove that the Zariski closure of $\langle p\rangle$
is the group of all matrices
 \[\scriptsize  \begin{pmatrix} \beta & \gamma & 0 & \ldots & 0\\
0 & \beta & \gamma & \ldots & 0\\
\vdots &\vdots &\vdots & \cdots & \vdots \\
0&0&0 & \ldots &\gamma\\
0&0&0 & \ldots &\beta
\end{pmatrix}
\]
where $\beta\in \C^*$ and $\gamma\in \C$.
\end{enumerate}

\item
Let $G$ be an algebraic group over $\C$, $s\in G$ and
let $A_s$ be the Zariski closure of $\langle s\rangle$ in $G$.\index[terms]{Zariski closure}
Prove that $s$ is semisimple if and only if 
$A_s$ does not contain unipotent elements.

{\em Hint:} Use Exercise 
\ref{_closure_of_jordan_cell_Exercise_}.

\item
Let $p\subset \GL(n, \C)$.
Prove that there exists {\bf the Jordan decomposition}\index[terms]{Jordan decomposition}
  $p=su$,
where $s$ is a semisimple operator and $u$
a unipotent operator commuting with $s$.
Prove that it is unique.

\item\label{_Jordan_Chevalley_Exercise_}
Let $\rho:\; G\arrow G'$ be a homomorphism
of algebraic groups, $p \in G$ a Jordan cell, as in 
Exercise \ref{_closure_of_jordan_cell_Exercise_}, \index[terms]{Jordan cell}
and $\rho(\C[p])$ the image of the group 
algebra generated by $p$. 
\begin{enumerate}
\item Prove that $\rho$ sends the unipotent
part of $p$ to the unipotent part of $\rho(p)$.
\item Prove that $\rho$ maps
semisimple elements to semisimple elements.
\item Prove that $\rho$ preserves
the Jordan decomposition $p=su$, that
is, $\rho(s)$ is semisimple and $\rho(u)$ is
unipotent.
\end{enumerate}

\item Let $V_0 \leftarrow V_1 \leftarrow V_2 \leftarrow ...$
be a sequence of surjective maps of vector spaces, 
and $A_i \in \End(V_i)$ a sequence of endomorphisms 
making the following diagram commutative
\[ \begin{CD}
V_0 @<<< V_1 @<<< V_2 @<<< ...\\
@VV{A_0}V @VV{A_1}V @VV{A_2}V\\
V_0 @<<< V_1 @<<< V_2 @<<< ...
\end{CD}
\]
Denote by $A$ the corresponding endomorphism
of the inverse limit $\lim\limits_\leftarrow V_i$.\index[terms]{inverse limit}
Prove that $A$ is semisimple if and only if
all $A_i$ are semisimple, and $A$ is unipotent
if and only if all  $A_i$ are unipotent.

\item \label{_Hopf_cone_Exercise_}
Let $V=\C^n$, and $A\in \End(V)$ be an automorphism with 
all eigenvalues $|\alpha_1|<1$, and ${\cal C}\subset
\C^n$  an $A$-invariant algebraic subvariety of $V$, not
necessarily smooth.
Prove that there exists $a\in \End V$ such that $A= e^a$
and ${\cal C}$ is $e^{ta}$-invariant for all $t\in \R$.

\item  In the assumptions of Exercise \ref{_Hopf_cone_Exercise_},
prove that the group of automorphisms of ${\cal C}$
commuting with the action of $A$ is an algebraic group.\index[terms]{group!algebraic}

\item  In the assumptions of Exercise
  \ref{_Hopf_cone_Exercise_}, prove
that the set of such $a$ is either finite or has
cardinality of continuum.

\item
Let $\Gamma\subset (\C^*)^n$ be an infinite
cyclic subgroup. Prove that the quotient
$(\C^*)^n/\Gamma$ is not isomorphic to an algebraic group.

{\em Hint:} Use induction in $n$.

\item
Let $\Z^{m+k}\subset \C^m\subset \GL(m, \C)$ be a discrete abelian
subgroup. Prove that $\C^m/\Z^{m+k}$ is not isomorphic
to an algebraic group if $k>0$.

{\em Hint:} Use the previous exercise.

\item
Let $G\subset (\C^*)^n\subset \GL(n, \C)$ be an
algebraic subgroup, and $\goth g$
its Lie algebra. Prove that the group
$\Lambda_{\goth g}:=\goth g \cap \Z^n$ 
satisfies $\rk \Lambda_{\goth g}=\dim_\C \goth g$.

{\em Hint:}
Identifying $(\C^*)^n$ with $\C^n/\Z^n$,
we obtain that $(\C^*)^n/G = \C^{n-k}/\Lambda$,
where $k= \dim_\C \goth g$, and $\Lambda$ a 
discrete subgroup of rank $n-\rk \Lambda_{\goth g}$.
Use the previous exercise to show that 
this is possible only if $\rk \Lambda=n-k$.

\item
Let  $(\C^*)^n\subset \GL(n, \C)$ be the diagonal
subgroup. Prove that any algebraic subgroup of $(\C^*)^n$ is isomorphic to
$(\C^*)^k$.

{\em Hint:} Use the previous exercise.

\item 
Let $G$ be the algebraic closure 
of a cyclic subgroup $\langle g\rangle$ in the diagonal
subgroup $(\C^*)^n\subset \GL(n, \C)$.
Prove that $G$ is $(\C^*)^n$ for a 
general $g\in (\C^*)^n$.

{\em Hint:} Prove that the number of algebraic
subgroups of $(\C^*)^n$ is countable, and
show that a general $g\in (\C^*)^n$ does not
belong to a proper algebraic subgroup.

\item Let $A\in \GL(n, \C)$ be a
diagonalizable element, and $G$ the connected component of the
Zariski closure of the cyclic subgroup $\langle A \rangle
\subset \GL(n, \C)$.
Prove that $G= (\C^*)^k$, where $k\in \Z^{\geq 0}$.
Prove that any $k\in \{0, 1, ..., n\}$ can be realized this way.
 
\item \label{_algebraic_closure_Hopf_Exercise_}
Let $M$ be a linear Hopf manifold,
$M= \frac{\C^n\backslash 0}{\langle \gamma\rangle}$,
where $\gamma= \lambda U$, $\lambda \in \C\backslash 0$, $|\lambda|< 1$
and $U\in \U(\C^n)$ a unitary matrix.
\begin{enumerate}
\item Prove that  the LCK\index[terms]{metric!LCK}
metric $\omega= \frac{dd^c |z|^2}{|z|^2}$ on $M$ is Vaisman.

\item Clearly, the radial vector field $\tilde{ r}$ on $\C^n$
is $\gamma$-invariant. Prove that $\tilde{ r}$
is the pullback of a vector field $\vec r$ on $M$; we call
$\vec r$ {\bf the radial vector field}.

\item Prove that the canonical foliation\index[terms]{foliation!canonical} of the LCK
metric $\omega= \frac{dd^c |z|^2}{|z|^2}$ on $M$ 
is tangent to the radial vector field on $M$.\index[terms]{metric!LCK}

\item
Let $G_\C$ be the smallest complex Lie subgroup 
in $\Aut(M)$ that contains $e^{t\vec r}$, for all $t\in \R$.
Prove that $G_\C$ contains a compact
subgroup $G$ such that $\Lie G_\C$
is the smallest complex Lie algebra
containing $\Lie G$.

\item 
Choose the coordinates on $\C^n$ such that
the unitary endomorphism $U$ is diagonal.
Using the structure theorem for reductive
Lie groups, prove that $G_\C= (\C^*)^k$.
Prove that $G_\C$ is the smallest Lie subgroup of $(\C^*)^n$
acting diagonally on $\C^n$ and containing $U$.

\item 
Let $\tilde G_\C$ be the lifting of $G_\C$ to
$\Aut(\C^n)$. Prove that $\tilde G_\C$
is the algebraic closure of $\langle \gamma\rangle$.

\end{enumerate}

\end{enumerate}


\chapter[Non-linear Hopf manifolds]
{Holomorphic contractions on Stein varieties and non-linear Hopf manifolds}
\label{_non-linear_Hopf_Chapter_}
{\setlength\epigraphwidth{0.9\linewidth}
\epigraph{\it\footnotesize Conclusions such as these open a wide field
  for speculation and exciting conjecture. They should be
  regarded, perhaps, in connexion with some of the most
  faintly-detailed incidents of the narrative; although in
  no visible manner is this chain of connexion
  complete. Tekeli-li! was the cry of the affrighted
  natives of Tsalal upon discovering the carcass of the
  white animal picked up at sea. This also was the
  shuddering exclamation of the captive Tsalalian upon
  encountering the white materials in possession of
  Mr. Pym. This also was the shriek of the swift-flying,
  white, and gigantic birds that is sued from the vapoury
  white curtain of the South. Nothing white was to be
  found at Tsalal, and nothing otherwise in the subsequent
  voyage to the region beyond. It is not impossible that
  "Tsalal," the appellation of the island of the chasms,
  may be found, upon minute philological scrutiny, to
  betray either some alliance with the chasms themselves,
  or some reference to the Ethiopian characters so
  mysteriously written in their windings.}{\sc\scriptsize
The Narrative of Arthur Gordon Pym by Edgar Allan Poe}}


\section{Introduction}

So far in this book we have discussed only the 
linear Hopf manifolds, defined as quotients\index[terms]{manifold!Hopf!linear}
of $\C^n \backslash 0$ by a cyclic group $\langle A \rangle$ 
generated by a linear contraction $A$. If we drop the linearity
assumption, we arrive at a general notion of ``Hopf manifold''
(\ref{_Hopf_gene_Definition_}).\index[terms]{manifold!Hopf!non-linear}

This definition is motivated by \index[persons]{Kodaira, K.} Kodaira's
research on complex surfaces. In \cite[Theorem 34]{_Kodaira_Structure_II_}
Kodaira proved that a complex surface with $b_1=1$ and $b_2=0$
that contains a curve and has no non-constant meromorphic
functions is a Hopf surface,\index[terms]{surface!Hopf} that is, a Hopf manifold of 
complex dimension 2. Using the Poincar\'e classification
of contractions of $\C^n$, and the further works by 
\index[persons]{Dulac, H.} Dulac, Latt\`es and \index[persons]{Sternberg, S.} Sternberg 
\cite{_Lattes_,_Dulac_,_Sternberg_contraction_},
\index[persons]{Kodaira, K.} Kodaira wrote a normal form 
for a Hopf surface, expressing a given
contraction as a polynomial map
in appropriate complex coordinates.

Following Poincar\'e, Dulac and others, 
\index[persons]{Kodaira, K.} Kodaira made a distinction between 
{\bf resonant} and {\bf non-resonant}\index[terms]{manifold!Hopf!(non-)resonant}
Hopf manifolds. Let $H= \frac{\C^n \backslash 0}{\langle \gamma \rangle}$
be a Hopf manifold. Using the Hartogs theorem, we extend\index[terms]{theorem!Hartogs}
$\gamma$ to a holomorphic automorphism of $\C^n$; let
$\alpha_1, ..., \alpha_n$ be the eigenvalues of 
the differential of $\gamma$ in 0.

A Hopf manifold is
called {\bf non-resonant} if the numbers
$\alpha_i$ do not satisfy a relation
of form $\alpha_i = \prod_{j=1}^d {\alpha_{l_j}}$
for any $d>1$. \index[persons]{Kodaira, K.} Kodaira showed that
all non-resonant Hopf surfaces are linear,\index[terms]{surface!Hopf!non-resonant}
and classified the resonant ones using the
normal form provided by  \index[persons]{Sternberg, S.} Sternberg's 
version of the Latt\`es--Poincar\'e--Dulac
theorem.

We re-prove a part of this classification theorem in 
Section \ref{_flat_Hopf_Section_}, showing that 
any non-resonant Hopf manifold is linear.
The full strength of the Sternberg (and  Poincar\'e--Dulac)
classification can be found in \cite{_Arnold:ODE+_}
and \cite{_Sternberg_contraction_}; its application
to the Hopf manifold is straightforward.

In this chapter, we prove that a Hopf manifold
(considered in this generality) is LCK with potential.\index[terms]{manifold!LCK!with potential}
For linear Hopf manifolds, this is proven in 
\ref{_linear_LCK_pot_Corollary_}. We prove that
all non-linear Hopf manifolds are  LCK with potential
by embedding them into linear Hopf manifolds. \index[terms]{manifold!linear Hopf}

In fact, the result we prove is even stronger.
Let $V$ be a Stein variety with an isolated\index[terms]{variety!Stein}
singularity $x$, and $\gamma:\; V \arrow V$ an
invertible holomorphic contraction centred in $x\in V$.
Then $M:=\frac{V\backslash x}{\langle \gamma\rangle}$
is a compact, smooth complex variety. Using the same argument
as used for the LCK manifolds with potential\index[terms]{contraction!holomorphic}
(\ref{embedding}), we prove that $M$ can be
embedded to a linear Hopf manifold. Then 
\ref{_linear_LCK_pot_Corollary_} implies that
$M$ is LCK with potential \index[terms]{manifold!LCK!with potential}whenever it is smooth.

In the last section, we discuss holomorphic embeddings 
$\Psi:\; H \hookrightarrow H_1$ from Hopf manifolds to Hopf manifolds. 
A special case is {\bf a linear embedding}, that is, a map
\[ \frac{\C^n \backslash 0}{\langle A \rangle}\hookrightarrow 
\frac{\C^m \backslash 0}{\langle A_1 \rangle}
\]
induced by a linear map $u:\; \C^n \hookrightarrow \C^m$
which makes the diagram
\begin{equation*}
\begin{CD}
\C^n @> u >> \C^m \\
@V{A}VV  @VV{A_1} V \\
\C^n @> u >> \C^m.
\end{CD}
\end{equation*}
commutative.

We say that
$\Psi:\; M \hookrightarrow H_1$ is {\bf a minimal embedding} 
to a linear Hopf manifold $H_1$ if it cannot be factorized through 
a chain of embeddings
$\Psi:\; M \hookrightarrow H_2 \hookrightarrow H_1$,
where $H_2 \hookrightarrow H_1$ is a proper linear embedding.

Let $\Psi:\; M \hookrightarrow H_1$ be a
holomorphic embedding of a complex manifold to a linear Hopf manifold,
$H_1=\frac{\C^n \backslash 0}{\langle A \rangle}$
and $\tilde M_c \subset \C^n$ the closure of the
pullback of $M$. 
We characterize the minimal embeddings using the
Zariski tangent space $T_c \tilde M_c$,\index[terms]{Zariski tangent space}
showing that $\Psi$ is minimal if and
only if $T_0 \tilde M_c= T_0 \C^n=\C^n$.

This result, surprisingly, does not imply
that any minimal embedding of Hopf manifolds
is an isomorphism. The reason is,
the closure $\tilde M_c$ does not need to be
normal. When $M$ is a Hopf manifold, this
closure is homeomorphic, but not necessarily
isomorphic to $\C^m$, and its Zariski tangent
space can be arbitrarily big.
We produce an example of a minimal
embedding of a classical
Hopf surface to a 7-dimensional 
linear Hopf manifold in \index[terms]{surface!Hopf!classical}
Exercise \ref{_minimal_Hopf_in_C^7_Exercise_}.

\section[Hopf manifolds and  holomorphic contractions]{Hopf manifolds\\ and  holomorphic contractions}

Further on in this section, we use the following 
elementary lemma.

\hfill

\lemma\label{_embeddings_are_open_Lemma_}
Let $K \subset M$ be a precompact open
subset of a complex manifold, $X$ another complex manifold, and
$\Map(M, X)$ the space of holomorphic maps from $M$ to $X$, equipped
with a $C^0$-topology.
\index[terms]{topology!$C^0$} Denote by $E\subset \Map(M, X)$ 
the space of holomorphic maps that are  injective immersions on $K$.
Then $E$ is open in $\Map(M, X)$.

\hfill

\pstep
Put some Hermitian metrics on $X$ and $M$, and 
define the map $u:\;\Map(M, X)  \arrow \R^{\geq 0}$ by 
\[ \Map(M, X)\ni f\mapsto u(f)=\inf_{x, y\in K, x\neq y} 
\frac{d(f(x), f(y))}{d(x,y)},
\] where $d$ denotes the 
Riemannian distance on $M$ and $X$.
We are going to show that $u(f)>0$ for all $f\in E$.
Around the diagonal $\Delta \subset K^2$,
the quantity $\frac{d(f(x), f(y))}{d(x,y)}$ is bounded from below by
the number $\inf_{v\in T_x M\backslash 0} |Df(v)|$, where $Df$
denotes the differential of $f$. Indeed,
\begin{equation}\label{_distance_and_deri_Equation_}
\lim_{y\rightarrow x} \frac{d(f(x), f(y))}{d(x,y)}=
   \lim_{y\rightarrow x}
   \left(\frac{\delta(x',y')}{d(x,y)}\right)=
   \lim_{y\rightarrow x}
   \left \|Df\restrict x\left(\frac{x'-y'}{d(x,y)}\right)\right\|,
\end{equation}
where $x', y'$ are images of $x$ and $y$ in a
chart defined in a neighbourhood of $x\in M$,
and $\delta$ the Riemannian distance in this chart.
Taking a sufficiently small precompact neighbourhood $U$
of the diagonal $\Delta\subset K^2$, we can assume that
\[ 
\frac{d(f(x), f(y))}{d(x,y)}\geq \inf_{v\in T_x M\backslash 0} |Df(v)| >0
\]
for all $(x, y) \in U\backslash \Delta$ 
(the number $|Df(v)|$ is positive, because the differential
$Df$ is injective, $f$ being an immersion).
Since $f$ is injective, we also have 
$\frac{d(f(x), f(y))}{d(x,y)}> \epsilon >0$
for $(x, y) \in K^2 \backslash U$. This implies that $u(f)>0$ for all $f\in E$. 

\hfill

{\bf Step 2:}
We have shown that $u:\;  E \arrow \R^{\geq 0}$ is bounded from below and 
strictly positive everywhere on $E$; moreover, if $u(f)>0$, 
the map $f$ is an injective immersion.

In a sufficiently small neighbourhood of $f\in E$, $u$ is positive, and hence 
this neighbourhood belongs to $E$.
\endproof

\subsection{Holomorphic contractions on Stein varieties}

Recall that {\bf a holomorphic contraction}
of a variety $V$ with center in $x\in V$
is a map $\gamma:\; V \arrow V$, $\gamma(x)=x$, such that
for each compact subset $K \subset V$ and
any open neighbourhood $U \ni x$, there
exists $N >0$ such that $\gamma^n(K) \subset U$
for any $n>N$. In this case, 
$\gamma$ is called {\bf the contraction centred in $x$}.

\hfill

\ref{_Stein_by_contract_to_linear_Hopf_Theorem_} below
is a generalization, in a sense, of \ref{embedding},
and their proofs are very similar. 
We decided to include both theorems
because the argument proving \ref{_Stein_by_contract_to_linear_Hopf_Theorem_}
is more abstract and harder to digest. 
Moreover, the embedding for LCK manifolds
with potential\index[terms]{manifold!LCK!with potential} has a cohomological interpretation
(\ref{_embedding_to_Hopf_via_BC_class_Corollary_})
which the seemingly more general 
\ref{_Stein_by_contract_to_linear_Hopf_Theorem_}
lacks. Therefore, it
makes sense to discuss these two theorems separately.

\hfill

\theorem\label{_Stein_by_contract_to_linear_Hopf_Theorem_} (\cite{ov_non_linear}) 
Let $V$ be a Stein variety, and $\gamma:\; V \arrow V$ 
an invertible holomorphic contraction centred in $x$.
Assume that the group $\langle \gamma\rangle$ acts on $V\backslash x$ properly
discontinuously. Then the quotient variety 
$\frac{V\backslash x}{\langle \gamma\rangle}$ admits
a holomorphic embedding to a linear Hopf manifold.

\hfill

\proof 
The proof is similar to the proof of \ref{_contra_compact_Theorem_}.
We start from an embedding $\Phi:\; V \hookrightarrow \C^n$, and
consider $\Phi$ as a collection of holomorphic functions on $V$.
We approximate each of these functions by a $\gamma$-finite\index[terms]{function!$\gamma^*$-finite} function.
This allows us to obtain the map $\Phi:\; V \hookrightarrow \C^n$
as a limit of $\Phi_i:\; V \hookrightarrow \C^n$,
with all coordinate components of $\Phi_i$ $\gamma$-finite.
Then we add to the coordinate component functions $f_{i,j}$, $j=1, ..., n$ of
$\Phi_i$ the functions $\gamma^* f_{i,j}, (\gamma^2)^* f_{i,j}, ...$
obtaining a finite-dimensional $\gamma^*$-invariant space of functions.
This gives a sequence of $\gamma^*$-equivariant maps 
$\Psi_i:\; V \arrow \C^{N_i}$, with $\Phi_i$ obtained from
$\Psi_i$ by a coordinate projection $\C^{N_i}\arrow \C^n$.

From \ref{_embeddings_are_open_Lemma_}, it follows that the space
of injective holomorphic maps is open in the uniform topology.
Let ${\cal U}$ be an open neighbourhood of $\Phi$
consisting of injective immersions. Since $\Phi= \lim_i \Phi_i$,
for $i$ sufficiently big,  $\Phi_i$ is injective. 
We have constructed a $\gamma$-finite 
injective immersion $V \arrow \C^{N_i}$. After removing the origins
and taking the $\gamma^*$-quotients, we obtain 
an holomorphic embedding 
\[ 
\frac{V\backslash x}{\langle \gamma\rangle } \arrow \frac{\C^{N_i}\backslash 0} {\langle \gamma^*\rangle }
\]
This was the scheme of the proof; we follow it with a detailed argument.

\hfill

\pstep
This step is more or less identical to 
\ref{_contra_compact_Theorem_}, Step 1.
Let $V_0 \ni x$ be a precompact Stein neighbourhood of $x\in V$.
By the definition of a contraction, for $k$ sufficiently big,
we have $\gamma^k(V_0) \Subset V_0$.\footnote{Here, as elsewhere,
the notation $A \Subset B$ means that $A$ is {\em relatively compact} in $B$,
that is, $A$ is a subset of $B$ and its closure is compact.}
Let $H^0_b(V_0, \calo_V)$ be the space of bounded holomorphic
functions on $V_0$, equipped with the $\sup$-metric,
and  $B(C)\subset H^0_b(V_0, \calo_V)$ the set of all
functions $f$ with $|f|\leq C$. This is a normal family\index[terms]{normal family}
(\ref{normal_family}), and hence  it is precompact in $C^0$-topology\index[terms]{topology!$C^0$}
(that is, in the topology of uniform convergence on compacts).
However, it is not compact in the $\sup$-topology,\index[terms]{topology!$\sup$} because
the closed ball in an infinite-dimensional Banach space\index[terms]{space!Banach} is never compact 
(this statement is called ``the Riesz theorem'', see
\ref{_Riesz_theorem_Remark_} or Exercise \ref{_Riesz_Exercise_}). 
\index[terms]{theorem!Riesz}

Consider $(\gamma^k)^*:\; H^0_b(V_0, \calo_V)\arrow H^0_b(V_0, \calo_V)$
as a linear map of Banach spaces.\index[terms]{space!Banach}
The map $(\gamma^k)^*$ takes any sequence $\{f_i\}$ converging in
$C^0$-topology\index[terms]{topology!$C^0$} to a sequence converging in the $\sup$-topology.\index[terms]{topology!$\sup$}
Indeed, $\{f_i\}$ uniformly converges on the closure
$\overline{\gamma^k(V_0)}$ that is  compact because $\gamma^k(V_0) \Subset V_0$.
Since $(\gamma^k)^*(f_i)(z)= f_i(\gamma^k(z))$, the
function $f_i$ takes the same values on $\gamma^k(V_0)$ as
$(\gamma^k)^*(f_i)$ takes on $V_0$.
Therefore, the sequence $\{(\gamma^k)^*(f_i)\}$ uniformly converges
on $V_0$.

Since an open ball in $H^0_b(V_0, \calo_V)$ is precompact
in the $C^0$-topology,\index[terms]{topology!$C^0$} the operator
$(\gamma^k)^*:\; H^0_b(V_0, \calo_V)\arrow H^0_b(V_0, \calo_V)$
takes an open ball to a precompact set in the Banach space
$H^0_b(V_0, \calo_V)$. Therefore, the operator 
$(\gamma^k)^*$ is compact as a map of Banach spaces.\index[terms]{space!Banach}

\hfill

{\bf Proof. Step 2:} Recall that a function $f \in H^0_b(V_0, \calo_V)$
is called {\bf $\gamma^k$-finite} if\index[terms]{function!$\gamma^k$-finite} 
the sequence 
$\{f, (\gamma^k)^*(f), (\gamma^{2k})^*(f), 
 (\gamma^{3k})^*(f),...\}$ belongs to a finite-dimensional
subspace of $H^0_b(V_0, \calo_V)$.

Let $R$ be a compact operator\index[terms]{operator!compact} without kernel
on a Banach space $B$\index[terms]{space!Banach}. Recalling the statement of \index[terms]{theorem!Riesz--Schauder}
the Riesz--Schauder theorem (\ref{_Riesz_Schauder_main_Theorem_}, 
\ref{_root_space_RS_Remark_}), we see that $B$ admits a\index[terms]{Schauder basis} Schauder basis%
\footnote{A {\bf Schauder basis} in a Banach space $W$\index[terms]{space!Banach}
is a set $\{x_i\}$ of vectors such that the closure of the
space generated by $x_i$ is $W$, and the closure of the space generated 
by all of $\{x_i\}$ except one of them does not contain the last one.}
such that the matrix for $R$ in this basis\index[terms]{function!$\gamma^*$-finite}\index[terms]{function!$\gamma^*$-finite}
is a sum of finite-dimensional Jordan blocks. This implies
that the set of $\gamma^k$-finite vectors in $H^0_b(V_0, \calo_V)$
is dense (see also \ref{_finite_RS_Corollary_}).

\hfill

{\bf Proof. Step 3:} 
Every $\gamma^k$-finite function\index[terms]{function!$\gamma^k$-finite} $f$ can be extended
to a holomorphic function on $V$. Indeed,
let $W_f$ be the finite-dimensional subspace
of $H^0_b(V_0, \calo_V)$ generated by 
$\{f, (\gamma^k)^*(f), (\gamma^{2k})^*(f), 
 (\gamma^{3k})^*(f),...\}$. Since 
$(\gamma^k)^*$ takes $W_f$ to itself, and $W_f$ is finite-dimensional,
the map $(\gamma^k)^*$ is invertible on $W_f$.
However, $((\gamma^k)^*)^{-1}$ takes functions
defined on $V_0$ to functions defined on $\gamma^{-k}(V_0)$,
hence any $w\in W_f$ can be holomorphically extended to
$\gamma^{-k}(V_0)$. Repeating this argument, we extend
$w\in W_f$ from $V_0$ to the sets
$\gamma^{-mk}(V_0)\supset \gamma^{(-m+1)k}(V_0)\supset ... \supset 
\gamma^{-2k}(V_0)\supset \gamma^{-k}(V_0)\supset V_0$.
Since $\gamma$ is a contraction, we have
$\bigcup_m \gamma^{-mk}(V_0)=V$, and hence  all elements of $W_f$  
are holomorphically extended to $V$; this includes $f\in W_f$.

\hfill

{\bf Proof. Step 4:} 
It is clear that every  $\gamma$-finite \index[terms]{function!$\gamma^*$-finite}function $f_1\in H^0(V, \calo_V)$ 
is $\gamma^k$-finite.\index[terms]{function!$\gamma^k$-finite} Now we prove that, conversely,
every $\gamma^k$-finite function $f\in H^0(V, \calo_V)$ 
is in fact $\gamma$-finite. Indeed, let 
$W_f\subset H^0(V, \calo_V)$ be 
the subspace generated by 
$\{f, (\gamma^k)^*(f), (\gamma^{2k})^*(f), 
 (\gamma^{3k})^*(f),...\}$, and $W'_f\supset W_f$ the subspace
generated by $\{f, (\gamma)^*(f), (\gamma^{2})^*(f), 
 (\gamma^{3})^*(f),...\}$. 

Until the rest of this step, we
use $\sum_{i=0}^\infty z_i$ as a shorthand 
for the subspace generated by the
vectors $\{z_i\}$.
Since $\Z^{\geq 0}= \bigcup_{i=0}^{k-1} \big[i+(k\Z^{\geq 0})\big],$
the space
$W'_f= \sum_{i=0}^{k-1} (\gamma^i)^* \sum_{l\geq 0}^\infty  (\gamma^{lk})^*(f)$ 
is equal to the subspace generated by the finite-di\-men\-si\-o\-nal
spaces $W_f, \gamma^*(W_f), ..., (\gamma^{k-1})^*(W_f)$,
hence it is also finite-dimensional.

\hfill

{\bf Proof. Step 5:} 
Since the space $V$ is Stein, it admits
a closed holomorphic embedding to $\C^d$.
In other words, there exists a 
$d$-dimensional subspace $A\subset H^0(V, \calo_V)$ 
such that the corresponding map $P:\; V \arrow \C^d$
is a closed embedding.
Since the space of $(\gamma^k)$-finite\index[terms]{function!$\gamma^k$-finite} functions is dense in 
$H^0_b(V_0, \calo_V)$, the space of $\gamma$-finite \index[terms]{function!$\gamma^*$-finite}functions
is dense in $H^0(V, \calo_V)$ with $C^0$-topology.\index[terms]{topology!$C^0$} Let $\Omega$ be a
precompact fundamental domain for the action of $\langle \gamma\rangle $ on $V$.
Since $H^0(V, \calo_V)^\gamma$ is $C^0$-dense in $H^0(V, \calo_V)$,
there exists a finite-dimensional
space $A_1= \C^d$ of $\gamma$-finite functions such that the
corresponding map $P_1:\; V \arrow \C^d$ approximates
the embedding $P:\; V \arrow \C^d$ uniformly on $\Omega$.
By \ref{_embeddings_are_open_Lemma_}, for an appropriate choice of
$A_1$, the map $P_1$ restricted to the closure $\bar\Omega$ of $\Omega$ is injective.
Let $W\subset H^0(V, \calo_V)$ be the smallest
$\gamma^*$-invariant subspace generated by $A_1$;
since all vectors in $A_1$ are $\gamma$-finite,
the space $W$ is finite-dimensional.
This makes the following diagram commutative
\begin{equation*}
\begin{CD}
V@>\Psi>> W \\
@V{\gamma}VV  @VV{\gamma^*} V \\
V @>\Psi >>  W
\end{CD}
\end{equation*}
where $\Psi:\; V \arrow W$ is a closed embedding
of $V$ to $W=\C^N$. Removing the zero and taking
the quotient over $\gamma$, respectively $\gamma^*$, this gives a
holomorphic map
$\frac{V\backslash x}{\langle \gamma\rangle}
\hookrightarrow \frac{W\backslash 0}{\langle \gamma^*\rangle}$,
where $\frac{W\backslash 0}{\langle \gamma^*\rangle}$ is by construction
a linear Hopf manifold. It is an injective immersion, 
because $P_1\restrict {\bar \Omega}$ is injective.
\endproof

\hfill

\corollary
Let $V$ be a Stein variety, and $\gamma:\; V \arrow V$ 
be an invertible holomorphic contraction centred in $x$.
Assume that the group $\langle\gamma\rangle$ acts on $V\backslash x$ properly
discontinuously, and $V\backslash x$ is smooth. Then the quotient manifold
$\frac{V\backslash x}{\langle \gamma\rangle}$ admits
an LCK metric with potential.\index[terms]{metric!LCK!with potential}

\hfill

\proof 
By \ref{_Stein_by_contract_to_linear_Hopf_Theorem_},
$\frac{V\backslash x}{\langle \gamma\rangle}$ admits
a holomorphic embedding to a linear Hopf manifold $H$.
By \ref{_linear_LCK_pot_Corollary_},
$H$ admits an LCK metric with potential.
Any closed complex submanifold of an LCK manifold with potential\index[terms]{manifold!LCK!with potential}
is again LCK with potential, and hence  
$\frac{V\backslash x}{\langle \gamma\rangle}$
is LCK with potential.
\endproof

\subsection{Non-linear Hopf manifolds are LCK}

So far in this book, we used only the linear Hopf manifolds.
However, a more general notion is worth studying.

\hfill

\definition\label{_Hopf_gene_Definition_}
Let $\gamma:\; \C^n \arrow \C^n$ be an invertible holomorphic
contraction centred in 0. Then the quotient
$\frac{\C^n \backslash 0}{\langle \gamma\rangle}$
is called {\bf a Hopf manifold}.\index[terms]{manifold!Hopf}

\hfill

\remark 
As follows from Exercise \ref{_Hopf_diffe_S^1_x_S^2n-1_Exercise_},
a Hopf manifold is diffeomorphic to $S^1 \times S^{2n-1}$.
In particular, it is never K\"ahler.

\hfill

The following corollary is directly implied
by \ref{_Stein_by_contract_to_linear_Hopf_Theorem_}.

\hfill

\corollary\label{_gene_Hopf_is_LCK_Corollary_} (\cite{ov_non_linear}) 
Let $H$ be a Hopf manifold. Then $H$ admits an LCK
metric with potential.\index[terms]{manifold!Hopf}\index[terms]{metric!LCK!with potential}
\endproof

\hfill

\remark
Most of the results about the LCK manifolds with potential\index[terms]{manifold!LCK!with potential}
rely upon the  Andreotti--Siu and  Rossi theorem used to construct
the Stein completion\index[terms]{completion!Stein} of the K\"ahler $\Z$-cover\index[terms]{cover!K\"ahler $\Z$-}. This is why we
restrict most of the statements to manifolds of dimension $\geq 3$.
However, a Hopf manifold $H$ is covered by $\C^n \backslash 0$,
and this manifold already has the Stein completion $\C^n$.
Therefore, in the results about the Hopf manifolds\index[terms]{theorem!Rossi, Andreotti--Siu}\index[terms]{completion!Stein}
we do not need the restriction $\dim_\C H \geq 3$.


\section{Minimal Hopf embeddings}


\definition
Let $M$ be an LCK manifold with potential\index[terms]{manifold!LCK!with potential},
admitting a holomorphic embedding $\Psi:\; M \hookrightarrow H$ 
to a Hopf manifold. We say that $\Psi$ is {\bf minimal}
if $\Psi$ does not admit a factorization
$M \hookrightarrow H_1 \stackrel \psi \hookrightarrow H$,
with $H_1$ a Hopf manifold of dimension less than $H$,
and $\psi$ is induced by a linear embedding of the
universal covering spaces.

\hfill

The following claim is trivial.

\hfill

\claim
Let  $\Psi:\; M \hookrightarrow H$ be a subvariety
of a Hopf manifold, and $\tilde \Psi:\;\tilde M \arrow \C^n \backslash 0$
the corresponding map of the covering spaces. Then $\Psi$ is minimal
if and only if the linear closure of $\im \tilde \Psi$ is $\C^n$.
\endproof

\hfill

Let $M$ be an LCK manifold with potential\index[terms]{manifold!LCK!with potential}, $W$ a vector space over $\C$,
$H= \frac{W\backslash 0}{\langle A\rangle}$ a linear Hopf
manifold, and $\Psi:\; M \hookrightarrow H$ a holomorphic embedding.
Passing to a $\Z$-covering\index[terms]{cover!K\"ahler $\Z$-}, we obtain the following diagram
\begin{equation*}
\begin{CD}
\tilde M_c @>{\tilde \Psi}>> W \\
@V{\gamma}VV  @VV{A} V \\
\tilde M_c @>{\tilde \Psi} >>  W
\end{CD}
\end{equation*}
where $\tilde M_c$ is 
a weak Stein completion\index[terms]{completion!Stein!weak} of 
the $\Z$-covering\index[terms]{cover!K\"ahler $\Z$-} of $M$,
and $\gamma$ the generator of the deck group of this covering.
Generally speaking, the weak Stein completion is defined only up
to normalization (\ref{_normalization_of_weak_Stein_comple_Remark_}).
However, the subvariety in $W$ is obtained as the closure
of $\tilde M \subset W \backslash 0$, and hence  $\tilde M_c$ is
defined uniquely.\index[terms]{variety!normal}\index[terms]{normalization}

\hfill

\theorem
In these assumptions,
let $T_c \tilde M_c$ be the Zariski tangent space of $\tilde M_c$.
Consider the map $\tilde \Psi:\; \tilde M_c \hookrightarrow \C^n$
obtained from the pullback of  $\Psi$, and let 
$D_c\tilde \Psi:\; T_c \tilde M_c \arrow T_0\C^n$
be its differential in the origin $c\in \tilde M_c$.
Then $\Psi$ is minimal if and only if
the differential $D_c\tilde \Psi:\; T_c \tilde M_c \arrow T_0\C^n$ is an isomorphism.

\hfill

\proof
By \ref{_algebraic_action_on_LCK_with_pot_Corollary_}, 
$\tilde\Psi(\tilde M_c)$ is 
invariant with respect to the Zariski closure of $\langle A\rangle$;
by \ref{_Zariski_closure_of_Jordan_cell_Claim_}, this
Zariski closure contains the group $\{\C^*\cdot \Id\}$.\index[terms]{Zariski closure} 
Since $\tilde\Psi(\tilde M_c )$ is $\C^*$-invariant,
it coincides with its Zariski tangent cone in 0, that is 
defined as the set of all 1-jets of the curves $S\subset W$ 
tangent to $\tilde\Psi(\tilde M_c)$ in $c$. \index[terms]{Zariski tangent cone}

It is not hard to see that the Zariski tangent space
of a variety is the linear closure of the Zariski tangent cone
(\cite[Exercise 20.1]{_Harris:first_course_}).
Therefore, the Zariski tangent space $T_0 \tilde\Psi(\tilde M_c)$ of
$\tilde \Psi(\tilde M_c)$ coincides with the minimal linear subspace 
containing $\tilde\Psi(\tilde M_c)$, that is  equal to $W$ if and only if
$\Psi$ is minimal.
\endproof

\hfill

\remark
Let $\Psi:\; H \arrow H'$ be a minimal 
embedding of Hopf manifolds. 
Then $\Psi$ is not necessarily an isomorphism
(Exercise \ref{_minimal_Hopf_in_C^7_Exercise_}).


\section{Poincar\'e--Dulac normal forms}


 Kodaira's classification of Hopf surfaces\index[terms]{surface!Hopf}
(\cite{_Kodaira_Structure_II_}) is based on \index[persons]{Sternberg, S.} Sternberg's version of Poincar\'e--Dulac
normal form theorem, that we explain in more detail in Chapter \ref{_Mall_bundles_Chapter_}.
Sternberg proves that any holomorphic contraction $A:\; \C^n \arrow \C^n$ can be expressed in a normal form
by a polynomial matrix, written in the upper diagonal form.
\begin{align*}
A_1(z_1, ..., z_n) &= \alpha_1 z_1 + P_1(z_2, ..., z_n)\\
A_2(z_1, ..., z_n)& = \alpha_2 z_2 + P_2(z_3, ..., z_n)\\
... \ \ \ \ \ \ \ \ \ \ & \ \ \ \ \ \ \ \ \ \ ...\\
A_n(z_1, ..., z_n) & = \alpha_n z_n.
\end{align*}
Let $A=US$ be the corresponding Jordan--Chevalley decomposition
(\ref{_LCK_pro_algebraic_linear_Proposition_}).\index[terms]{Jordan--Chevalley decomposition}
If $A$ is non-linear, that is, $P_i \neq 0$, 
the action of $A$ on polynomials is upper
triangular in an appropriate monomial
basis, with non-trivial terms above
the diagonal. It is not hard
to see that a matrix that is  upper
triangular and has non-trivial terms
above the diagonal cannot be made
diagonal in any other basis. Therefore,
the action of $A$ on the polynomials 
is not semisimple, and $U\neq \Id$.

As shown in \ref{_non-Vaisman_Corollary_}, 
a Vaisman manifold admits an embedding to a 
diagonal Hopf manifold $H$, that is, 
to a manifold $H=\frac{\C^n \backslash 0}{\langle \gamma\rangle}$,
where $\gamma$ is diagonal.
Recall  that $\gamma$-finite functions on $\C^n$ are
always polynomial (\ref{_gamma_finite_on_C^n_Lemma_}). Therefore, the monodromy
action on $\gamma$-finite functions on the $\Z$-cover of 
a diagonal Hopf manifold (and hence, on the Vaisman manifolds)
is semisimple.

If $U$ is non-trivial, the action of the
monodromy on $A$-finite functions on $\C^n$ cannot
be semisimple, and hence  $M$ cannot be Vaisman.

This brings the following corollary.

\hfill

\corollary (\cite{ov_non_linear}) 
Let $H$ be a non-linear Hopf manifold. Then $H$ does not
admit a Vaisman metric.
\endproof

\section{Exercises}

\begin{enumerate}[label=\textbf{\thechapter.\arabic*}.,ref=\thechapter.\arabic{enumi}]

\item
Let $U \subset V$ be a closed subspace of a normed vector
space.
\begin{enumerate}
\item Let $z\in V \backslash U$
and $y\in U$ be a point that satisfies
\[ d(z,y) \leq (1-\epsilon) d(z, U).\]
Let $x:= z-y$.
Prove that $|x| \geq (1-\epsilon)d(x, U).$

\item
Prove {\bf the Riesz lemma:} for any $\epsilon >0$
there exists $z\in V$ such that the distance\index[terms]{lemma!Riesz}
$d(z, U)$ satisfies $d(z, U) \geq (1-\epsilon) |z|$.
\end{enumerate}

\item
Let $B$ be a unit ball in an infinite-dimensional
normed vector space. Prove that $B$ contains an infinite sequence
$x_1, x_2, ..., x_n, ...$ such that $d(x_i, x_j) >1/2$
for all $i\neq j$.

{\em Hint:} Use the previous exercise when $U$
is finite-dimensional and apply induction on $\dim U$.

\item\label{_Riesz_Exercise_}
Let $V$ be a normed vector space such that the
closed unit ball in $V$ is compact. Prove that $V$
is finite-dimensional.

{\em Hint:} Use the previous exercise.

\item
Construct an infinite-dimensional
topological vector space $V$ admitting
a translation-invariant metric $d$
such that the closed unit ball in $V$ is compact.

{\em Hint:} Find a Fr\'echet space with
all bounded subsets precompact. Use Montel
theorem.\index[terms]{theorem!Montel}\index[terms]{topological vector space!Fr\'echet}

\item
Let $A\in \GL(n, \C)$ be a linear contraction acting on
$\C^n$, and $V= \bigoplus_d V_d$ the space of polynomial functions on $\C^n$,
with $V_d$ denoting the space of homogeneous polynomials of degree $d$.
Recall that the root space of $A$ associated with an eigenvalue $\alpha$ 
is $\bigcup_{i=1}^\infty \ker(A-\alpha\Id)^i$.
\begin{enumerate}
\item
Let $U\subset H^0(\C^n \calo_{\C^n})$ be a finite-dimensional
$A$-invariant subspace. Prove that all functions $f\in U$ are polynomial.

\item
Prove that all the root spaces $V_\alpha$ of the action of $A$ on $V$ are finite-dimensional.\index[terms]{root space}

\item Prove that any root space 
is contained in a sum $\oplus_{i=0}^m V_{i}$, for some $m=1, 2, 3, ...$

\item 
Prove that 
$V = \bigoplus_{\alpha_i} W_{\alpha_i}$, where $\alpha_i$ are the eigenvalues and 
$W_{\alpha_i}$ are  all root spaces for the action of $A$ on $W$.

\item
Prove that the projection
$V \arrow W_{\alpha_i}$ along the other root spaces is
continuous in $C^0$-topology \index[terms]{topology!$C^0$}on $V$.
\end{enumerate}

\item\label{_Dense_finite_Hopf_Exercise_}
Let $A\in \GL(n, \C)$ be a linear contraction acting on
$\C^n$, $V= \bigoplus_d V_d$ the space of polynomial functions on $\C^n$,
with $V_d$ denoting the space of homogeneous polynomials of degree $d$,
and $W\subset V$ an $A$-invariant subspace that is  dense in $C^0$-topology.\index[terms]{topology!$C^0$}

\begin{enumerate}
\item Prove that $W=\bigoplus_d (W\cap V_d)$.

\item Prove that $W=V$.
\end{enumerate}

{\em Hint:} Use the previous exercise.

\item
Let $\Psi:\; H_1 \arrow H_2$ be a holomorphic embedding of linear Hopf manifolds.
\begin{enumerate}
\item Prove that $\Psi$ induces a holomorphic embedding 
$\tilde \Psi:\;  \C^m \backslash 0 \arrow \C^n \backslash 0$
of the universal covering spaces.

\item Prove that $\tilde \Psi$ can be extended to 
a holomorphic map $\tilde \Psi_c:\;  \C^m \arrow \C^n$.

\item
Assume that $H_1$ and $H_2$ are linear Hopf manifolds.
Prove that the map $\tilde \Psi_c$ is polynomial.

\item 
Assume that $H_2$ is diagonal. 
Using Exercise \ref{_Dense_finite_Hopf_Exercise_}, prove 
that $H_1$ is also diagonal.

\item 
Let $H_1 = (\C^2\backslash 0)/\Z$ be a classical
Hopf manifold, with the $\Z$ action generated
by $\lambda\Id$, and $\psi:\; \C^2 \arrow \C^3$
the {\bf Veronese map} taking $(x, y)$ to $(x^2, xy, y^2)$
Prove that $\psi$ can be extended to a holomorphic map
$\Psi:\; H_1 \arrow H_2$ where $H_2$ is a classical
Hopf manifold, $H_2= (\C^3\backslash 0)/\langle \lambda^2 \Id\rangle$.
Prove that $\Psi$ is a 2 to 1 covering map of its image.\index[terms]{map!Veronese}

\item
Let $\Psi:\; H_1 \arrow H_2$ be a holomorphic embedding of classical
Hopf manifolds, and 
$\tilde \Psi_c:\;  \C^m \backslash 0 \arrow \C^n \backslash 0$
its lift to the universal covers.
Prove that $\tilde \Psi$ is a homogeneous polynomial map.
Prove that $\Psi$ is not an embedding, unless it is linear.

\end{enumerate}

\item
Consider the ring
$\C[x, y]_d$ of polynomials of degree $\geq d$ plus constants.
Prove that $\C[x, y]_d$ is finitely generated over $\C$.

\item
\begin{enumerate}
\item Prove that $\C[x, y]_2= \C[t_1, t_2, ..., t_7]/I$
where $t_1=x^2, t_2=xy, t_3=y^2, t_4=x^3, t_5 = x^2 y, t_6= xy^2, t_7=y^3$
and $I$ is the ideal generated by all relations between $t_i$ 
in $\C[x,y]$. 
\item Prove that $I$ is finitely generated.
\end{enumerate}

\item
Consider the algebraic variety $V= \Spec(\C[x, y]_2)$.
Prove that $\C^2$ is the normalization of $V$, and the
normalization map $\C^2 \arrow V$ is bijective and
defines an isomorphism outside of zero.\index[terms]{map!normalization}

\item
Consider a matrix $A$ on the space $\langle t_1, ..., t_7\rangle$
with $A(t_i) = 4 t_i$ for $i=1, 2, 3$ and $A(t_i) = 8 t_i$ for
$i=4,5,6,7$. Prove that the action of $A$ on the generators
$t_1, ..., t_7$ of $\C[x, y]_2$ can be extended to an 
automorphism $\alpha$ of this algebra.

\item
Prove that the
normalization map $\nu:\; \C^2 \arrow V$ induces a commutative diagram
\begin{equation*}
\begin{CD}
\C^2 @>{\nu}>> V \\
@V{A_0}VV  @VV{\alpha} V \\
\C^2 @>{\nu} >>  V
\end{CD}
\end{equation*}
where $A_0(z)=2z$.
\begin{enumerate}
\item Prove that the corresponding holomorphic map
$\frac{\C^2\backslash 0}{\langle A_0\rangle} \arrow 
\frac{V\backslash \nu(0)}{\langle \alpha \rangle}$
is an isomorphism.
\item Prove that $V$ is a closed algebraic
cone, and that $\C^2 \backslash 0$ is the
corresponding open algebraic cone.
\end{enumerate}

\item\label{_minimal_Hopf_in_C^7_Exercise_}
Let $\Psi:\; V \arrow \C^7$ take $z\in V$ to $t_1(z) , ..., t_7(z)$,
where $t_i\in \C[x,y]_2$ are considered to be functions on $V=\Spec(\C[x,y]_2)$.
\begin{enumerate}
\item 
Prove that the following diagram is commutative
\begin{equation*}
\begin{CD}
V @>{\Psi}>> \C^7 \\
@V{\alpha}VV  @VV{A} V \\
V @>{\Psi} >>  \C^7
\end{CD}
\end{equation*}
\item 
Prove that $T_{\nu(0)} ^*V= \langle dt_1, ..., dt_7\rangle$,
where $T_{\nu(0)}^* V= \frac{\goth m}{\goth m^2}$ 
is the Zariski cotangent space to the origin $\nu(0)\in V$, and $\goth m$ 
the maximal ideal of the origin.
\item 
Prove that $\Psi$ induces a minimal embedding
of linear Hopf manifolds
$\frac{\C^2\backslash 0}{\langle A_0\rangle} \arrow 
\frac{\C^7\backslash 0}{\langle A \rangle}$.
\end{enumerate}

\item
Let $A= {\scriptsize \begin{pmatrix} \lambda & 1 \\ 0 & \lambda\end{pmatrix}}$
be an invertible linear contraction on $\C^2$. 
We say that a function $f$ on $\C^*$ is {\bf $A^*$-automorphic}
if $A^*(f)= \const \cdot f$.
\begin{enumerate} 
\item Prove that
any $A^*$-automorphic polynomial on $\C^2$ belongs to $\C[z_2]$,
where $z_2$ is an $A^*$-automorphic coordinate function.
\item Prove that the Hopf surface $\frac{\C^2 \backslash 0}{\langle A \rangle}$\index[terms]{surface!Hopf}
contains only one complex curve.
\end{enumerate}

\item
Let 
\[ A= \small \begin{pmatrix} 
\lambda & 1 & 0 \\ 0 & \lambda & 1 \\
0 & 0 & \lambda 
\end{pmatrix}
\]
be an invertible linear contraction on $\C^3$
that has a unique Jordan cell.
\begin{enumerate} 
\item Prove that the space of $A^*$-automorphic
quadratic polynomials on $\C^3$ is 2-dimensional.
\item Prove that the Hopf manifold 
$\frac{\C^2 \backslash 0}{\langle A \rangle}$
contains a non-constant family of divisors.
\end{enumerate}

\end{enumerate}


\chapter{Morse--Novikov and Bott--Chern cohomology of LCK manifolds} \index[terms]{manifold!LCK}
\label{_MN_and_BC_Chapter_}

\epigraph{\it\footnotesize
-- ``Oh, they're there all right,'' Orr had assured him about
the flies in Appleby's eyes after Yossarian's fist fight
in the officers' club, ``although he probably doesn't even
know it. That's why he can't see things as they really
are.''\\ \medskip

-- ``How come he doesn't know it?'' inquired Yossarian.\\ \medskip

-- ``Because he's got flies in his eyes,'' Orr explained with
exaggerated patience. ``How can he see he's got flies in
his eyes if he's got flies in his eyes?''}{\sc\scriptsize Catch-22, by Joseph Heller}



\section{Introduction}
\label{_MN_BC_Intro_Section_}

The geometry of a K\"ahler manifold
is strongly influenced by an important cohomological
invariant, the K\"ahler class,\index[terms]{class!K\"ahler} that is, the de Rham 
cohomology class of its K\"ahler form.\index[terms]{form!K\"ahler}

Its analogue in LCK geometry is called
{\bf the twisted Bott--Chern class}. Many results\index[terms]{class!Bott--Chern!twisted}
of K\"ahler geometry \index[terms]{geometry!K\"ahler} that interpret geometric
properties of a manifold in terms of the 
K\"ahler class\index[terms]{class!K\"ahler} can be translated to LCK geometry this way.

The first such result is an LCK analogue of the
 Kodaira theorem, that states that a K\"ahler
manifold is projective if and only if it admits
a rational K\"ahler class. The model manifold
for K\"ahler geometry is $\C P^n$, and
this means that a manifold
with rational K\"ahler class\index[terms]{class!K\"ahler} can be embedded
to $\C P^n$. The model manifold for LCK
geometry is the Hopf manifold.\index[terms]{manifold!Hopf} However,
an LCK analogue of  the Kodaira theorem\index[terms]{theorem!Kodaira embedding} is
even easier to state: an LCK manifold is embeddable
to Hopf if and only if its twisted Bott--Chern class vanishes 
(\ref{_embedding_to_Hopf_via_BC_class_Corollary_}).\index[terms]{class!Bott--Chern!twisted} 

Let $(M, \theta, \omega)$ be an LCK manifold,
and $d_\theta(\alpha):= d\alpha - \theta\wedge\alpha$
be the Morse--Novikov differential (Section
\ref{_potential_weight_Subsection_}). 
The LCK equation $d\omega=\theta \wedge \omega$
can be reformulated as $d_\theta\omega=0$.
The cohomology of the complex $(\Lambda^*(M), d_\theta)$ is 
identified with the cohomology of $M$
with coefficients in the weight local system $L$;\index[terms]{bundle!weight}
the class of $\omega$ in $H^2(M, L)$ is
called {\bf the Morse--Novikov class} of $(M, \theta, \omega)$.
\index[terms]{class!Morse--Novikov}
It is an important cohomological invariant of an LCK
manifold\index[terms]{manifold!LCK}; however, our attempts to deduce geometric
information from the Morse--Novikov class failed, so far.
We conjecture that the vanishing of the Morse--Novikov class implies
the existence of an LCK potential,\index[terms]{potential!LCK} but this conjecture
seems to be very hard.

However, the converse statement is easy to prove;
indeed, the existence of LCK potential immediately
implies that $\omega\in \im (d_\theta)$.
Moreover, the cohomology with coefficients in the weight
local system vanishes for all Vaisman manifolds,
because the Vaisman manifolds are diffeomorphic\index[terms]{manifold!Vaisman}
to mapping tori over a circle, and the weight local
system is the pullback of a non-trivial local system
on the circle.
By the deformation theorem (\ref{def_lckpot2Vai}),
the Morse--Novikov cohomology vanishes for the
LCK manifolds with potential \index[terms]{manifold!LCK!with potential}as well
(\ref{mn_van}).

However, the Morse--Novikov class\index[terms]{class!Morse--Novikov} has an important refinement,
called the twisted Bott--Chern class.\index[terms]{class!Bott--Chern!twisted} Let us recall the definition
of the {\bf Bott--Chern cohomology} of a complex manifold.

A compact complex manifold has two
important cohomology theories, that are  (almost)
equal for  K\"ahler manifolds: the Dolbeault
cohomology and the de Rham cohomology. For 
complex manifolds, they are hard to relate,
because a Dolbeault cohomology class cannot\index[terms]{cohomology!Dolbeault}\index[terms]{cohomology!de Rham}
always be realized by a closed form, and the
de Rham cohomology class cannot always be 
realized by a Dolbeault closed form.
The Fr\"olicher spectral sequence
with $E_2^{p,q}= H^{p,q}(M)$ converges to the
de Rham cohomology, but there are no
direct maps between these groups.\index[terms]{cohomology!Bott--Chern}

Bott--Chern cohomology, defined 
as $\frac{\ker \6 \cap \ker \bar\6}{\im \6\bar\6}$,
was introduced precisely to rectify this problem.
Indeed, all forms in $\ker \6 \cap \ker \bar\6$
are closed, and Dolbeault closed. Moreover,
the group $\im \6\bar\6$ is contained in 
$\im\6\cap \im \bar \6$, and hence  one has
natural homomorphisms from the Bott--Chern
cohomology to both Dolbeault and de Rham 
cohomology groups. By $dd^c$-lemma
(\ref{_dd^c_lemma_Theorem_}), on
a compact K\"ahler manifold, the
Bott--Chern cohomology is equal to 
de Rham and Dolbeault cohomologies.
In the more general context, the 
Bott--Chern cohomology is related to
de Rham and Dolbeault cohomologies through 
several foundational exact sequences
discovered by \index[persons]{Angella, D.} Angella and \index[persons]{Tomassini, A.} Tomassini
(\cite{angella,_Angella_Tomassini_}).

The notion of Bott--Chern cohomology is
less straightforward than that of Dolbeault
and de Rham cohomology.\index[terms]{cohomology!de Rham} First of all,
the Bott--Chern cohomology\index[terms]{cohomology!Bott--Chern} has no multiplicative
structure; in fact, a product of two $dd^c$-exact
forms is no longer $dd^c$-exact. Second, there
is no second-order Laplacian operator with the
kernel that can be identified with the Bott--Chern
cohomology. A fourth-order ``Laplacian'' for
the Bott--Chern cohomology was defined in 
\cite{sch}; however, to prove that the Bott--Chern
cohomology is finite-dimensional, it is easier
to use the notion of an elliptic complex.

In this chapter, we recall the notions
of elliptic differential operators and
elliptic complex, and apply
these notions to the Bott--Chern cohomology.

This formalism is applied to the notion of ``twisted Bott--Chern
cohomology'', defined as the Bott--Chern cohomology with
values in a flat line bundle $L$, usually identified with the
weight local system:\index[terms]{cohomology!Bott--Chern}
\[
H^{p, q}_{BC}(M,L)
:=\frac{\ker \6_\theta \cap \ker \bar\6_\theta}{\im \6_\theta\bar\6_\theta}.
\]
Using a trivialization of $L$, we write its connection
as $\nabla:= \nabla_0 + \theta$, where $\nabla_0$
is the standard connection on a trivial bundle.
Then $\6_\theta(\alpha) := \6(\alpha)-\theta^{1,0}\wedge \alpha$
and $\bar\6_\theta(\alpha) := \bar\6(\alpha)-\theta^{0,1}\wedge \alpha$.

Since the symbol of $\6_\theta\bar\6_\theta$ is the same as the symbol of
$\6\bar\6$ and the symbols of $\6_\theta$, $\bar\6_\theta$ are
the same as the symbols of $\6, \bar\6$, the group $H^{p, q}_{BC}(M,L)$,
identified with the cohomology of an elliptic complex, is
finite-dimensional.

The twisted Bott--Chern class\index[terms]{class!Bott--Chern!twisted} $[\omega]_{BC}$ of an LCK form\index[terms]{form!LCK}
vanishes when $\omega= \6_\theta\bar\6_\theta\phi$,
hence it is zero for LCK manifolds with potential\index[terms]{manifold!LCK!with potential}.
The converse assertion looks trivial (and, indeed,
in \cite{ov_jgp_09} we stated it without a proof).
However, the cohomological condition $[\omega]_{BC}=0$ 
does not imply that $\phi$ is positive, and without the positivity
it is hard to establish the existence of an LCK potential.\index[terms]{potential!LCK}

In Chapter \ref{_posi_pote_Chapter_}, we prove that 
$[\omega]_{BC}=0$ implies the existence of an LCK
potential for an LCK metric\index[terms]{metric!LCK}, possibly in a different
conformal class; this result is highly non-trivial.

An important special case of this theorem was established in
a paper of N. \index[persons]{Istrati, N.} Istrati \cite{nico2}. One of the applications of the 
twisted Bott--Chern class\index[terms]{class!Bott--Chern!twisted} vanishing was given in
\cite{ov_mz}, where we proved that any
LCK manifold admitting a holomorphic conformal
flow, acting by non-isometric homotheties on
the K\"ahler cover, admits an LCK metric
$\omega$ with potential\index[terms]{metric!LCK!with potential}. \index[persons]{Istrati, N.} Istrati proved that the
equation $\omega= \6_\theta\bar\6_\theta\phi$
has a positive solution $\phi$ in this special case.

\section{Preliminaries on differential operators}

\subsection{Differential operators}

We recall the basic notions about differential operators
on  vector bundles (see, for example, \cite{_LM_}), to be used when
considering different types of cohomology on complex
manifolds. We omit most of the proofs, referring to
\cite{_LM_}, \cite{griha}, \cite{demailly},
\cite{_Berline_Getzler_Vergne_} and other
books where the subject is treated.

\hfill

\definition ({\bf Grothendieck}, \cite{_co_}) 
\index[persons]{Grothendieck, A.}
\label{_diff_ope_Gr_Definition_}\\
Let $R$ be a commutative ring over a field $k$, and let 
and $A, B$ be $R$-modules.
A {\bf differential operator of order 0}\index[terms]{operator!differential}
from $A$ to $B$ is a $R$-linear map $\phi\in \Hom_R(A,B)$.
Differential operators of order $i>0$ are defined inductively:
$\alpha \in \Diff^i(A,B)$ if for any $r\in \R$,
the commutator $[\al,L_r]=\alpha L_r-L_r\alpha$ belongs to $\Diff^{i-1}(A,B)$,
where $L_r(x)=rx$.

\hfill

Now let  $F$, $G$ be vector bundles on a same smooth
manifold $M$. As usual, we denote the space of sections
of a smooth vector bundle $F$ by the same letter $F$.
This space is considered to be a $C^\infty M$-module.
Notice that, the by Serre--Swan theorem, the category of\index[terms]{theorem!an}
vector bundles over $M$ is equivalent to the
category of projective $C^\infty M$-modules
(\cite{_Karoubi_}); this explains our convention.

\hfill

\definition
A {\bf differential operator} 
on vector bundles $F$, $G$ is a differential
operator between the corresponding spaces of sections
in the sense of \ref{_diff_ope_Gr_Definition_}. The set of
all such differential operators is denoted $\Diff^i(F,G)$.
Any element of  
$\Diff^i(M):= \Diff^i({C}^{\infty}( M),{C}^{\infty}( M))$
is called a {\bf differential operator on
  $M$}.\index[terms]{operator!differential}

\hfill

\remark
{This definition is equivalent to the usual one:}
locally (in coordinates) any differential operator
is expressed as a composition of derivations and
multiplications by $f\in {C}^{\infty} (M)$.

\hfill

To define {\em the symbol}
 of a differential operator, we consider the filtration $$\Diff^0(M)\subset
\Diff^1(M)\subset \Diff^2(M)\subset \cdots$$ 
For a vector bundle $V$, let  $\Sym^i(V)$ denote the $i$-th symmetric power (i.\,e.     the symmetric part
of the $i$-th tensor power) of $V$. We recall:

\hfill

\theorem
The  associated graded ring to this filtration is isomorphic to the ring 
$\bigoplus_i \Sym^i(TM)$, identified with the
ring of fibrewise polynomial functions on $T^*M$.

\hfill

\corollary
For vector bundles $F,G$ on $M$, we have:
\[
\Diff^i(F,G)/\Diff^{i-1}(F,G)=\Sym^i(TM)\otimes\Hom(F,G).
\]

Now we can give:

\hfill

\definition
Let $F, G$ be vector bundles on $M$, 
and $D\in \Diff^i(F,G)$  a differential operator.
Consider its class in $\Diff^i(F,G)/\Diff^{i-1}(F,G)$ 
as a $\Hom(F,G)$-valued function on $T^*(M)$ (polynomial
of order $i$ on each cotangent space).
This function is called the {\bf symbol} \index[terms]{operator!symbol of}
of $D$.

\hfill

We can now introduce an important class of differential operators.

\hfill

\definition
Let $F$, $G$ be vector bundles of the same rank 1  $M$.
A differential
operator $D:\; F \arrow G$ is called {\bf elliptic}\index[terms]{operator!elliptic}
if its symbol $\sigma(D)\in \Hom(F,G)\otimes \Sym^i(TM)$
is invertible at each non-zero $\xi \in T^*M$.

\hfill

\example
Consider a second-order differential operator
$D:\; C^\infty M \to C^\infty M$. Locally 
in coordinates $z_1,..., z_n$, the operator $D$ can be written
as
\[ D(f) = t + u f + \sum_{i=1}^n v_i \frac{\6 f}{\6 z_i} 
+ \sum_{i,j= 1}^n w_{ij} \frac{\6^2 f}{\6 z_i\6 z_j}
\]
Clearly, the symbol of $D$ is 
$ \sum_{i,j= 1}^n w_{ij} \frac{d}{dz_i}  \frac{d}{dz_j}$
considered to be  a function on $T^* M$.
It is elliptic if and only if the matrix
$w_{ij}$ is positive definite or negative
definite everywhere on $M$. For example,
the Laplacian $\Delta(f) = \sum_{i=1}^n \frac{\6^2 f}{\6 z_i^2}$
is always elliptic.

\hfill

Before giving the essential properties of elliptic
operators, we need to recall the definition of 
$L^2$-topology\index[terms]{topology!$L^2$} on the sections of a vector bundle.

\hfill

\definition
Let $F$ be a vector bundle on a compact manifold,\index[terms]{topology!$L^2$}
and $\nabla$ a connection on $F$. Fix a connection on
$\Lambda^1 M$. These two connections define a connection
on the tensor product $F\otimes (\Lambda^1 M)^{\otimes
  k}$, that we denote by the same letter $\nabla$.
Consider the operator
$\nabla^k:\; F \arrow F\otimes  (\Lambda^1 M)^{\otimes k}$
obtained by iterating $\nabla$.
The {\bf $L^2_p$-topology} on the space of sections
of $F$ is the topology defined by the quadratic form
$$\|f\|^2=\sum_{i=0}^p\, \int_M |\nabla^if|^2 \Vol_M,$$
for some  scalar product on $F$ and a Riemannian 
form on $M$.

Let $L^2_p(F)$ denote the completion of the
space of sections of $F$ with the above topology.

\hfill

\remark
One can prove that the above topology is independent of the choice of the connection and scalar product used in the definition.

\hfill

\definition
A continuous operator $\psi:\; A \arrow B$ on topological
vector spaces is called {\bf Fredholm}\index[terms]{operator!Fredholm}
if 
\begin{enumerate}
	\item $\ker D$ is finite-dimensional;
	\item $\im D$ is closed and has finite codimension.
\end{enumerate}

The next result is foundational for the entire Hodge theory (see, for example  \cite{_Berline_Getzler_Vergne_}, \cite{demailly}). 

\hfill

\theorem
Let $D:\; F \arrow G$ be an elliptic operator of order
$d$. Clearly, $D$ defines a continuous map 
$L^2_{p}(F) \arrow L^2_{p-d}(G)$. {Then this map is
	Fredholm.}

\subsection{Elliptic complexes}

\definition
Let $F$, $G$, $H$ be vector bundles,
and $F\stackrel D \arrow G\stackrel D \arrow H$ a complex
of differential operators (that is, $D^2=0$).
It is called {\bf an elliptic complex} if its symbols 
$F\stackrel {\sigma(D)} \arrow G\stackrel  {\sigma(D)} \arrow H$
give an exact sequence at each non-zero $\xi \in T^*M$.\index[terms]{complex!elliptic}

\hfill

\definition
Let $A$, $B$, $C$ be topological vector spaces
and $A\stackrel D \arrow B\stackrel D \arrow C$ a complex
of continuous maps. It is called {\bf a Fredholm complex} \index[terms]{complex!Fredholm}
if $\im D$ is closed, 
and $\frac{\ker D}{\im D}$ is finite-dimensional.

\hfill

With all these preparations we can state the following.

\hfill

\theorem
Let $F\stackrel {D_1} \arrow G\stackrel {D_2} \arrow H$ be an 
elliptic complex of differential operators, with
$D_1$ of order $d_1$ and $D_2$ of order $d_2$.
Then the complex
	$L^2_{p}(F) \stackrel{D_1}\arrow 
	L^2_{p-d_1}(G)\stackrel{D_2}\arrow L^2_{p-d_1-d_2}(H)$
	is Fredholm.

\hfill

\corollary\label{ell_com}
The cohomology modules of any elliptic complex on a
compact manifold are finite-dimensional.

\section{Bott--Chern cohomology}
\label{_Bott_Chern_Section_}

\definition
Let $M$ be a complex manifold, and $H^{p,q}_{BC}(M)$
the space of closed $(p,q)$-forms modulo $dd^c(\Lambda^{p-1,q-1} M)$.
Then $H^{p,q}_{BC}(M)$ is called {\bf the Bott--Chern
  cohomology} of $M$. \index[terms]{cohomology!Bott--Chern}

\hfill

\remark
A $(p, q)$-form $\eta$ is closed if and only if
$\6\eta=\bar\6\eta=0$. 
Using $2\1 \6\bar\6= dd^c$, we  could define the
Bott--Chern cohomology $H^*_{BC}(M)$ as 
$H^*_{BC}(M):=\frac{\ker \6 \cap \ker \bar\6}{\im \6\bar\6}$.

\hfill

\remark
There are natural (and functorial) maps from the
Bott--Chern cohomology to the Dolbeault cohomology\index[terms]{cohomology!Dolbeault}\index[terms]{cohomology!Bott--Chern}
$H^*(\Lambda^{*,*} M, \bar\6)$ and to the de Rham
cohomology, but no morphisms between the de Rham and
the Dolbeault cohomology (see also \cite{angella}).

\hfill

\theorem\label{_BC_f_dim_Theorem_}
Let $M$ be a compact complex manifold. Then
$H^{p,q}_{BC}(M)$ is finite-dimensional.

\hfill

\proof Since
$H^*_{BC}(M)=\frac{\ker \6 \cap \ker \bar\6}{dd^c}$,
it would suffice to show that the complex
\[
\Lambda^{p-1,q-1} M\stackrel{dd^c} \arrow
\Lambda^{p,q} M \stackrel{\6+\bar\6} \arrow \Lambda^{p+1,q}M
\oplus \Lambda^{p,q+1}M
\]
is elliptic. At $\xi\in T^*M=\Lambda^{1,0}M$, the symbol
$\sigma(dd^c)$ of $dd^c$ is 
equal to the operator of multiplication of a form by $\xi\wedge \bar\xi$,
 the symbol $\sigma(\6)$ of $\6$ is multiplication by $\xi$ and
the symbol $\sigma(\bar\6)$ of $\bar \6$ is multiplication by $\bar\xi$.
Therefore, $\ker \sigma(\6) = \im \sigma(\6)$,
$\ker \sigma(\bar\6) = \im \sigma(\bar\6)$ (this
proves finite-dimensionality of Dolbeault cohomology),\index[terms]{cohomology!Dolbeault}
and $\ker \sigma(\bar\6)\cap \ker \sigma(\6)=\im\sigma(\6\bar\6)$.
\endproof

\section{Morse--Novikov cohomology}

\subsection{Morse--Novikov class of an LCK manifold}\index[terms]{manifold!LCK}\index[terms]{class!Morse--Novikov}

Let $(M,I,\omega,\theta)$ be an LCK manifold and let $L$
be its weight bundle. \index[terms]{bundle!weight}As usual, let $\pi:\tilde M\arrow M$ be the universal cover of $M$.  Recall the
definitions and results in Section \ref{_Riemann_Hilbert_Section_}
and Subsections \ref{_weight_automo_Subsection_}
and \ref{_L_valued_Kahler_Subsection_}.  

\hfill

\definition 
 The cohomology of the complex $(\Lambda^*M, d_\theta:=d-\theta)$
is called {\bf Morse--Novikov cohomology}, or
{\bf \index[persons]{Lichnerowicz, A.} Lichnerowicz cohomology}, or {\bf twisted cohomology}; the corresponding
complex is the {\bf Morse--Novikov complex}.\index[terms]{complex!Morse--Novikov}
It computes the cohomology of the local system $L$
associated with the weight bundle.  \index[terms]{cohomology!Morse--Novikov}
\index[terms]{cohomology!twisted}\index[terms]{cohomology!Lichnerowicz}

\hfill

\remark Recall that the local system $L$ associated with the weight bundle of an LCK manifold\index[terms]{manifold!LCK} that is  not GCK is never trivial, and hence it cannot have global sections. We conclude that for any compact LCK manifold that is  not GCK, $H^0(L)=0$.

\hfill

\definition
Let $(M,I,\omega,\theta)$ be an LCK manifold, and $(\tilde
M, \tilde \omega)$ its
K\"ahler cover, with the deck group $\Gamma$.
Given $\gamma\in \Gamma$, we denote by $\chi(\gamma)\in \R^{>0}$
the scale factor given by $\gamma^*\tilde
\omega=\chi(\gamma)\tilde \omega$. Clearly, this
defines a multiplicative character $\chi:\;
\Gamma\arrow\R^{>0}$. Let $\lambda\in \R^{>0}$.
A function $f\in  C^\infty(M)$ is called {\bf
  $\lambda$-automorphic} if $\gamma^*(f)=
\chi(\gamma)^\lambda f$ for all $\gamma\in \Gamma$.
In a similar way, one defines $\lambda$-automorphic\index[terms]{form!automorphic!of weight $\lambda$}
differential forms.\index[terms]{function!automorphic!of weight $\lambda$}

\hfill

\proposition \label{_automorphic_forms_sections_Proposition_}
 Let $(L,\nabla)$ be the weight bundle on $M$, 
$\phi$ its trivialization, and $\theta$ the connection form in $L$.
Consider  a non-zero parallel section\index[terms]{section!parallel} $\Phi$ of $\pi^* L$.
Then:

(i) {$\frac {\pi^*\phi}\Phi$ is an automorphic function of
weight $\lambda$.}

(ii) For each $L^\lambda$-valued
differential form $\eta$, the differential form
$\frac {\pi^*\eta}{\pi^*\phi}\in \Lambda^*\tilde M$ is automorphic of weight 
$\lambda$, giving an {equivalence $\Xi$ between the space of 
sections of $\Lambda^*M\otimes L$ and
$\Lambda^*M_\lambda$.} \index[terms]{form!automorphic!of weight $\lambda$}

(iii) Under this equivalence, the {de Rham differential on
$\Lambda^*M_\lambda$ corresponds to $d_\nabla$.}

\hfill

\proof (i) is clear, because the monodromy \index[terms]{action!monodromy} 
acts on $\pi^*\phi$ trivially and acts on $\Phi$ with weight $\lambda$.

(ii) is clear by the same reason: any section of 
$\Lambda^*M\otimes L$ produces an automorphic form,
and any automorphic form\index[terms]{form!automorphic} $\rho$ gives a section
$\rho \Phi^{-1}$ of $\pi^* L^\lambda$, which is  fixed
by the monodromy, and hence it is obtained as a pullback.

To prove (iii), take any $\rho \in \Lambda^*M\otimes L$. 
Then 
\[ d(\Xi(\rho))= d(\pi^*\rho \Phi^{-1})=
\pi^* d_\nabla \rho \Phi^{-1}-  \rho\wedge\nabla \Phi\cdot \Phi^{-1}=
\Xi(d_\nabla \rho)- \Xi(\rho \wedge \theta), 
\]
which completes the proof.
\endproof

\hfill

\corollary
 The Morse--Novikov complex on $M$ 
is identified with the de Rham 
complex of automorphic forms\index[terms]{form!automorphic} on $\tilde M$.
\endproof

\hfill

\remark 
 The above \ref{_automorphic_forms_sections_Proposition_} remains true in the more general context of a manifold endowed with a local system. Indeed, in the proof we did not make use of the complex structure.

\hfill

For an LCK manifold\index[terms]{manifold!LCK}, $d_\theta\omega=0$, and hence we can give:

\hfill

\definition 
 The cohomology class $[\omega]_{MN}$ of $\omega$ in the Morse--Novikov
cohomology is called {\bf Morse--Novikov class} of $M$.\index[terms]{class!Morse--Novikov}

\hfill

\proposition 
The Morse--Novikov class $[\omega]_{MN}$ vanishes for LCK manifolds with
potential and hence for Vaisman manifolds.\index[terms]{manifold!Vaisman}

\hfill

\proof Since an LCK form $\omega$ with potential 
is $d_\theta d^c_\theta$-exact \index[terms]{form!LCK!with potential}
(\ref{_LCK_pot_via_d_theta_d^c_theta_Proposition_}), $[\omega]_{MN}$ vanishes.
\endproof

\hfill

\remark  The class $[\omega]_{MN}$ 
is known to be non-zero for some other LCK manifolds. \index[terms]{manifold!LCK}
Indeed, a blow-up of an LCK manifold at a point is again LCK\index[terms]{blow-up}
(see Section \ref{_Blow_up_at_points_section_}), 
and the LCK form on the blow-up is 
actually globally conformally K\"ahler when restricted
to the exceptional divisor.
Therefore, it cannot be exact.\index[terms]{divisor!exceptional}
In  \cite{ban2} it was shown that the standard
LCK structure\index[terms]{structure!LCK} on an Inoue surface (see Chapter \ref{inoue_lck}) has non-exact
Morse--Novikov class. \index[terms]{surface!Inoue}
In fact, all known examples of compact LCK manifolds
with vanishing Morse--Novikov class\index[terms]{class!Morse--Novikov} admit an LCK
metric with potential.\index[terms]{metric!LCK!with potential}

\hfill

On compact LCK manifolds with potential\index[terms]{manifold!LCK!with potential}, we have a very sharp result:

\hfill

\theorem { (\cite{ov_jgp_09})} \label{mn_van}
The Morse--Novikov cohomology of a compact 
LCK manifold $M$ with potential, $\dim_\C M >2$, vanishes.	

\hfill

\proof
By \ref{_MN_vanishes_Lee_class_Corollary_}, 
it is enough to take $M$ to be Vaisman. 
Indeed, by \ref{_MN_vanishes_Lee_class_Corollary_}, 
the set of Lee classes on an LCK\index[terms]{class!Lee}
manifold with potential\index[terms]{manifold!LCK!with potential} is the same as the set
of Lee classes on its Vaisman deformation.

Then \ref{mn_van} would follow from the Structure theorem
\ref{str_vai}.

Indeed, by \ref{str_vai}, a 
Vaisman manifold is diffeomorphic to a mapping torus\index[terms]{mapping torus}\index[terms]{manifold!Vaisman}
of a Sasakian automorphism $q:\; S\arrow S$,
for a compact Sasakian manifold $S$. \index[terms]{manifold!Sasaki}However, the group of
automorphisms of a Sasakian manifold is contained
in its isometry group, and hence  it is compact 
(\cite{_Myers_Steenrod_,_Kobayashi_Transformations_}). Therefore,\index[terms]{group!isometry}
for any automorphism $q:\; S \arrow S$ of 
a Sasakian manifold, $q^d$ is isotopic to the identity
for some integer $d>0$. We obtain that 
a finite covering of a Vaisman manifold is 
diffeomorphic to $S^1 \times S$.
Since the harmonic forms on $M$ are contained
in the harmonic forms in its finite covering,
it suffices to prove vanishing of Morse--Novikov
cohomology for a finite cover, that is  a Vaisman manifold 
diffeomorphic to $S^1 \times S$.

The Morse--Novikov cohomology of the $d_\theta$-complex is
isomorphic to the cohomology of the local system
$L$.  By construction, $L=p^*L'$, with $L$ being a local system on
$S^1$ and $p$ the submersion $M \arrow S^1$ 
(\ref{_proper_potential_S^1_Remark_}). Write the
cohomology of $L$ as the derived direct image $R^iP_*(L)$,
where $P$ is the projection to a point (see also Subsection \ref{_Derived_functors_Subsection_}). Decomposing
$P= \pi \circ p$ onto a projection $\pi$ 
to $S^1$ and projection $p$ from $S^1$ to a point, we obtain
$R^*P_*(L)=R^*p_*\pi_*(L')$, that is 
equal to the hypercohomology  ${\Bbb H}^*(S^1, R^*\pi_*(L'))$ of the
complex of sheaves $R^*\pi_*(L')$. 
The \index[persons]{Leray, J.} Leray spectral\index[terms]{spectral sequence!Leray}
sequence has $E_2$-term  $H^p(H^q(S)\otimes_\C L')$ and 
converges to  ${\Bbb H}^*(S^1,R^*\pi_*(L'))=H^*(M, L')$. 

However each sheaf $H^q(S)\otimes L'$ has
vanishing cohomology, being a locally trivial sheaf on
$S^1$ with non-trivial constant monodromy.\index[terms]{monodromy} Therefore, this
Leray spectral sequence \index[terms]{spectral sequence!Leray}
vanishes in $E_2$ and converges to $0$.

For another proof of this result, using a theorem from
\cite{llmp}, see \ref{_MN_vanishes_Lee_class_Corollary_}.
\endproof

\subsection{Twisted Dolbeault cohomology}\index[terms]{cohomology!Dolbeault!twisted}

\definition 
 Let $M$ be a complex manifold,
and $(L,\nabla)$ a flat, oriented, real line bundle on $M$.
Identifying sections of $L$ with automorphic\index[terms]{bundle!line!flat}
forms of weight 1 on $\tilde M$ as above, we consider {the
Hodge decomposition $d_\theta=\6_\theta+\bar\6_\theta$,}\index[terms]{form!automorphic!of weight $\lambda$}
where $d_\theta$ is the de Rham differential on 
automorphic forms,\index[terms]{form!automorphic} and $\6_\theta$, $\bar\6_\theta$ are its Hodge components.\index[terms]{cohomology!Dolbeault!twisted}

\hfill

\remark 
Consider the map $\Xi$ identifying $\lambda$-automorphic
differential forms and the forms with values in
$L^\lambda$. Clearly, this equivalence
is compatible with the Hodge decomposition. This proves the following.

\hfill

\claim\label{_del_del_bar_Claim_}
 The following identities hold on an LCK manifold:\index[terms]{manifold!LCK}
 $$\6_\theta = \6 - \theta^{1,0},  
\qquad  \bar\6_\theta = \bar\6 - \theta^{0,1}.$$

This leads to the following commutation relations (see also \cite{va_comp}):

\begin{equation}\label{comm_form}
 \begin{split}
  \{\6_\theta, \bar\6_\theta\}&= \{\6_\theta, \6_\theta\}
=\{\bar\6_\theta, \bar\6_\theta\}=0,\\ 
 \{d_\theta, d_\theta^c\}&=0, \quad -2\1\6_\theta\bar\6_\theta= d_\theta d_\theta^c,
 \end{split}
\end{equation}
where
$d_\theta^c= I d_\theta I^{-1}=
-\1(\6_\theta-\bar\6_\theta)=d^c - I(\theta)$,
and $\{\cdot, \cdot\}$ denotes the anticommutator
(\ref{_anticommu_Corollary_}).

\section{Twisted Bott--Chern cohomology}\label{twisted_BC}

On LCK manifolds,\index[terms]{manifold!LCK} we need to consider Bott--Chern cohomology with values in the weight bundle (and hence associated with the $d_\theta$). The precise definition is as follows.\index[terms]{bundle!weight}

\hfill

\definition \label{_twi_BC_class_Definition_}
 Let $M$ be a complex manifold, and $L$ a flat vector
bundle. Consider the corresponding differential
$d_\nabla=d_\theta$, and let $\6_\theta$, $\bar\6_\theta$
be its Hodge components.
{\bf The weighted (twisted) Bott--Chern cohomology groups} 
are defined as 
\[ H^{p,q}_{BC}(M,L):= \frac{\ker
   d_\theta\restrict{\Lambda^{p,q}M\otimes L}}{\im
  \6_\theta\bar\6_\theta}.\index[terms]{cohomology!Bott--Chern!twisted}
\]

\hfill

The following theorem is completely analogous to \ref{_BC_f_dim_Theorem_};
its proof is the same.\index[terms]{bundle!vector bundle!flat}

\hfill

\theorem \label{finite_dim_BC}
 Let $M$ be a compact complex manifold, and $L$ a flat vector
bundle.  Then the spaces  $H^{p,q}_{BC}(M,L)$ are finite-dimensional.
\endproof

\hfill

\definition 
 Let $(M, \omega,\theta)$ be an LCK manifold,\index[terms]{manifold!LCK} and $L$ its
weight bundle. The cohomology class $[\omega]_{BC}$ of $\omega$ in 
$H^{1,1}_{BC}(M,L)$ is called {\bf the Bott--Chern
 class of $M$}.\index[terms]{class!Bott--Chern}

\hfill

Note that this is the closest analogue of the K\"ahler class,\index[terms]{class!K\"ahler}
and the following theorem (together with the Hopf
embedding result) are an LCK analogue of the Kodaira embedding
theorem.\index[terms]{theorem!Kodaira embedding} We state it here, for convenience,
but its proof refers to  \ref{_main_positive_intro_Theorem_} about the positivity
of LCK potential\index[terms]{potential!LCK} proven in the next chapter.

\hfill

\theorem \label{vanishing_BC_class_implies_potential_Theorem_}
 Let $(M, \omega,\theta)$ be an LCK manifold.
Suppose that $[\omega]_{BC}=0$.
{Then $M$ admits an LCK structure $(\omega', \theta')$ with proper  potential.}\index[terms]{structure!LCK!with potential}

\hfill

\proof 
Automorphic K\"ahler potential on the K\"ahler cover
$(\tilde M, \tilde \omega)$ is the
same as a positive function $\phi_0\in C^\infty M$ such that 
$d_\theta d^c_\theta\phi_0=\omega$ 
(\ref{_LCK_pot_via_d_theta_d^c_theta_Proposition_}).
By \ref{_main_positive_intro_Theorem_},
an LCK manifold
$(M, \omega, \theta)$ with $\omega\in \im d_\theta d^c_\theta$
admits an LCK metric $\omega'$ with
positive automorphic potential $\phi_0$.

 The properness of the potential
 is equivalent to $(M,\theta, \omega')$ having LCK rank 1.\index[terms]{rank!LCK}
However, we can deform the LCK structure\index[terms]{structure!LCK} $(\theta, \omega')$
retaining the complex structure on $M$ 
in such a way that the LCK rank of the deformation is 1
(\ref{defor_improper_to_proper}). \endproof 

\hfill

\corollary\label{_embedding_to_Hopf_via_BC_class_Corollary_}
Let $(M, \theta, \omega)$ be a compact
LCK manifold,\index[terms]{manifold!LCK} $\dim_\C M \geq 3$,
and $[\omega]_{BC}\in H^{1,1}_{BC}(M,L)$ 
its twisted Bott--Chern class.\index[terms]{class!Bott--Chern!twisted} Suppose that $[\omega]_{BC}=0$.
Then $M$ admits a holomorphic embedding to a Hopf manifold.
\endproof

\section[Bott--Chern classes and Morse--Novikov cohomology]{Bott--Chern classes and Morse--Novikov\\ cohomology}\index[terms]{class!Bott--Chern}

Let $L$ be a holomorphic line bundle over the complex manifold $M$. The holomorphic cohomology of $L$ can be realized as cohomology of the complex
\begin{equation}
C^\infty(L)\stackrel{\bar\6}{\arrow}\Lambda^{0,1}(M,L)\stackrel{\bar\6}{\arrow}\Lambda^{0,2}(M,L)\stackrel{\bar\6}{\arrow}...
\end{equation}
Let $L$ be equipped with a flat connection $\nabla$. Then $\6: \Lambda^{0,1}(M,L)\arrow\Lambda^{1,1}(M,L)$ induces a natural map
\begin{equation}\label{_del_morphism_Equation_}
H^1(L)\stackrel{\6}{\arrow} H^{1,1}_{BC}(M,L)
\end{equation}
Let $L'$ be the bundle $L$ with a holomorphic structure defined by the complex conjugate of the  $\nabla^{1,0}$-part of the connection. The complex
\begin{equation}\label{_cohomology_of_L'_Equation_}
	C^\infty(L)\stackrel{\nabla^{1,0}}{\arrow}\Lambda^{0,1}(M,L)\stackrel{\nabla^{1,0}}{\arrow}\Lambda^{0,2}(M,L)\stackrel{\nabla^{1,0}}{\arrow}...
\end{equation}
computes the holomorphic cohomology of $L'$. 

If the bundle $L$ is real, then $L\equiv L'$ and the cohomology of  \eqref{_cohomology_of_L'_Equation_} is naturally equivalent to $\bar{H^*(L)}$. The map $\bar\6:\Lambda^{0,1}(M,L)\arrow \Lambda^{1,1}(M,L)$ induces a homomorphism
\begin{equation}\label{_del_bar_morphism_Equation_}
\bar{H^1(L)}\arrow H^{1,1}_{BC}(M,L).
\end{equation}
that is  entirely similar to \eqref{_del_morphism_Equation_}. 

The following result allows one to compute the Bott--Chern classes\index[terms]{class!Bott--Chern} in terms of holomorphic cohomology and
Morse--Novikov cohomology.

\hfill

\theorem \label{_BC_exact_sequence_via_Dolbeault_Theorem_}
{ (\cite{ov_jgp_09})}
Let $M$ be a complex manifold, $L_\R$ a trivial real line bundle with flat connection $d-\theta$, where $\theta$ is a real closed\index[terms]{connection!flat}
1-form. Denote by $L$ its complexification viewed as a holomorphic bundle. Then there exists the following exact sequence:
\begin{equation}\label{_exact_sequence_BC_MN_Equation_}
	H^1(L)\oplus\bar{H^1(L)}\stackrel{\6+\bar\6}{\arrow} H^{1,1}_{BC}(M,L)\stackrel{\nu}{\arrow} H^2(M,L)
\end{equation}
where $H^2(M,L)$ is the Morse--Novikov cohomology, $\nu$ is the tautological map and the first arrow is obtained as the direct sum of \eqref{_del_morphism_Equation_} and \eqref{_del_bar_morphism_Equation_}.

\hfill

\proof
Let $\eta$ be a (1,1)-form such that $[\eta]=0\in H^2(M,L)$. Then $$\eta=d_\theta\alpha=\bar\6_\theta\alpha^{1,0}+\6_\theta\alpha^{0,1},$$
where $\6_\theta=\6-\theta^{1,0}$ and $\bar\6_\theta=\bar\6-\theta^{0,1}$ (see \ref{_del_del_bar_Claim_}). However, these operators are precisely those that are used to define the first arrow in \eqref{_exact_sequence_BC_MN_Equation_}. Moreover, since $\eta$ has Hodge type (1,1), we have
\[
\bar\6_\theta\alpha^{1,0}=0, \qquad \6_\theta\alpha^{0,1}=0,
\]
and hence $\alpha^{1,0}$ and $\alpha^{0,1}$ represent cohomology classes in $H^{1,0}(L)$,and $\bar{H^{1,0}(L)}$. Now, the Bott--Chern class \index[terms]{class!Bott--Chern}of $\eta$ is obtained as $\6[\alpha^{0,1}]+\bar\6[\alpha^{1,0}]$, and thus the sequence \eqref{_exact_sequence_BC_MN_Equation_} is exact.
\endproof

\hfill

Together with \ref{vanishing_BC_class_implies_potential_Theorem_}, the above result immediately implies the following.

\hfill

\proposition
Let $(M,\omega,\eta)$ be a compact LCK manifold\index[terms]{manifold!LCK}, and $L$ its weight bundle seen as a holomorphic line bundle. Assume $H^1(M,L)=0$ and $H^2_\theta(M)=0$. Then $H^{1,1}_{BC}(M,L)=0$ and $M$ admits an LCK structure\index[terms]{structure!LCK!with potential} with proper potential. \endproof


\section{Exercises}

\begin{enumerate}[label=\textbf{\thechapter.\arabic*}.,ref=\thechapter.\arabic{enumi}]

\item
Let $R$ be a ring, and 
$D^i \in \Diff^i(R)$, $D^j \in \Diff^j(R)$
differential operators on $R$ defined inductively as 
in \ref{_diff_ope_Gr_Definition_}. Prove that the composition
$D^i D^j$ lies in $\Diff^{i+j}(R)$. 

{\em Hint:}
Use induction and the identity
$[ v, D^i D^j] = [v, D^i] D^j + D^i [v, D^j]$

\item 
Let  $D^i \in \Diff^i(R)$, $D^j \in \Diff^j(R)$
be differential operators. Prove that the commutator
$[D^i, D^j]$ lies in $\Diff^{i+j-1}(R)$.

{\em Hint:} Use induction and Jacobi identity
\[
[ v, [D^i, D^j]] = [[v, D^i], D^j] + [D^i, [v, D^j]].
\]

\definition
Let $R$ be a $k$-algebra, and
$D:\; R \arrow A$ a $k$-linear map from
$R$ to an $R$-module. It is called
{\bf a $k$-derivation}, or just {\bf derivation}
 if it satisfies the
Leibniz rule: $D(xy) = y D(x) + x D(y)$.

\item
\begin{enumerate}
\item Prove that $D(k)=0$ for any $k$-derivation
on a $k$-algebra (we assume $\Char k=0$).
\item
Let $R$ be a finite extension of a field $k$ of characteristic
0. Prove that the space $\Der_k(R,R)$ of derivations
vanishes.
\end{enumerate}

\item 
Let $x_1, ..., x_n$ be coordinates on $\R^n$.
Prove that any derivation on $C^\infty \R^n$ is
written as coordinates as $D(f) = \sum_{i=1}^n f_i \frac \6
{\6 x_i}$, where $f_i\in C^\infty M$.

{\em Hint:} Use the Hadamard lemma\index[terms]{lemma!Hadamard}: for each point
$x\in \R^n$, and any function $f\in C^\infty \R^n$,
$f$ is a sum of an affine function and a function
$g$ with $g(x)=0$ and $dg\restrict x=0$.

\item 
Let $D\in \Diff^1(R)$ be a differential operator of
first order on a ring $R$. Prove that $D- D(1)$ is a derivation of $R$.
Prove that $\Diff^1(R)/\Diff^0(R)$ is isomorphic to the
space of derivations of $R$.

\item Find an example of a compact LCK (non-K\"ahler)
manifold admitting symplectic structure.\index[terms]{manifold!LCK}

\item  Let $M$ be a compact complex manifold,
$H^2(M)=H^1(M)=0$. Prove that $\dim H^{1,1}_{BC}(M)= \dim H^1(\calo_M)$,
where $H^{p,q}_{BC}(M)$ denotes the Bott--Chern cohomology.

\item\label{124} 
Assume that the standard map $\bigoplus_{p,q} H^{p,q}_{BC}(M)\arrow H^*(M)$
to the de Rham cohomology \index[terms]{cohomology!de Rham}is injective. Prove that it is surjective.

\item  Assume that the standard map 
$\bigoplus_{p,q} H^{p,q}_{BC}(M)\arrow \bigoplus_{p,q}H^{p,q}_{\bar \6}(M)$
to the Dolbeault cohomology\index[terms]{cohomology!Dolbeault} is injective. Prove that it is surjective.

\item
Let $B$ be a vector bundle on a manifold $M$.
Consider a first order operator $D:\; B \arrow \Lambda^1M \otimes B$,
and let $\sigma(D)$ be its symbol.
Denote by $\pi:\; T^*M \arrow M$ the projection.
We consider $\sigma(D)$ as an
$\pi^*(\End(B) \otimes \Lambda^1 M)$-valued function on $T^*M$,
linear on fibres of $T^*M$, that is, as a
section of $\End(B) \otimes \Lambda^1 M\otimes TM= \End(B) \otimes \End(TM)$.
Prove that $D$ is a connection if and only if 
$\sigma(D) = \Id_E \otimes \Id_{TM}$.

\item  Let $M$ be a compact Riemannian 4-manifold, and 
$\Lambda^2M=\Lambda^+M\oplus\Lambda^-M$ the decomposition
on eigenspaces of the Hodge $*$-operator. 
Consider the map $d^+:\; \Lambda^1M \arrow \Lambda^+M$
obtained by projecting the de Rham differential to $\Lambda^+M$.
\begin{enumerate}
\item Prove that the complex
\[
C^\infty M \stackrel d\arrow\Lambda^1 M \stackrel {d^+}\arrow { \Lambda^+M}
\]
is elliptic.
\item
Prove that $\ker d^+/\im d$ is finite-dimensional.
\end{enumerate}

\item
Let $D_1:\ F \arrow G$, $D_2:\; G \arrow H$
be differential operators, with $F, G, H$
vector bundles of the same rank. 
\begin{enumerate}
\item
Prove that the
composition $D_1 \circ D_2$ is elliptic if $D_2, D_2$ are elliptic.
\item Prove that $D_1$ and $D_2$ are elliptic if $D_1\circ D_2$ is elliptic.
\end{enumerate}

\item
Let $\nabla:\; B \arrow B\otimes \Lambda^1 M$
be a connection on a Hermitian
vector bundle on a Riemannian manifold, and 
$\nabla^*:\; B\otimes \Lambda^1 M\arrow B$ its Hermitian
adjoint.
Prove that $\nabla^*\nabla$ is elliptic, and its symbol
is equal to $-g\otimes \Id_B$, where $g\in \Sym^2 TM$
is the metric tensor. 

\item
Consider the operator
$d+d^*:\; \Lambda^* M \arrow
\Lambda^* M$. Prove that it is elliptic.

\item Prove that the equations
\eqref{comm_form} follow from the usual \index[terms]{Kodaira identities} Kodaira identities
\[
  \{\6, \bar\6\}= \{\6, \6\}
=\{\bar\6, \bar\6\}=0, \quad
 \{d, d^c\}=0, \quad -2\1\6\bar\6= d d^c,
\]
twisted by the operator of multiplication
by $e^\psi$, where $\theta=d\psi$.

\item
Let $(E, \bar\6_0)$
be a trivial holomorphic line bundle, $\rho^{0,1}$ a 
$\bar\6$-closed (0,1)-form,
and $\bar\6:= \bar\6_0 + \rho^{0,1}$
the corresponding holomorphic structure operator.
Let $\alpha\in \R$ be a real number.
{\bf The tensor power} $E^\alpha$ is the trivial line
bundle with the holomorphic structure operator
$\bar\6_0 + \alpha\rho^{0,1}$.
Prove that $E^\alpha = E^{\otimes \alpha}$
when $\alpha$ is integer.

\item
Let $(M, \theta, \omega)$ be an LCK manifold\index[terms]{manifold!LCK} with the Lee
form non-vanishing everywhere on $M$. Prove that 
$H^*(\Lambda^*(M), d_\theta)=0$.

\item
Let $(M, \theta, \omega)$ be an
LCK manifold with vanishing 
Morse--Novikov cohomology class
$[\omega]_{MN}\in H^2(M, L)$.
Denote by $(\tilde M, \tilde \omega)$ the K\"ahler
cover of $M$. 
\begin{enumerate}
\item
Prove that $\tilde \omega=d\eta$, where $\eta$ 
is an automorphic 1-form.

\item Prove that there exists a 
holomorphic Hermitian line bundle $E$
over $\tilde M$ such that $\tilde \omega$ is the
curvature of its Chern connection.\index[terms]{connection!Chern}

{\em Hint:} Take a trivial line bundle $E$ 
with $\bar \6_0$ its operator of holomorphic
structure, and define the holomorphic structure
operator $\bar\6$  on $E$
as $\bar\6 = \bar\6_0 + \eta^{0,1}$.

\item
Prove that the deck transform maps $E$ to a
line bundle that is  isomorphic to the tensor power 
$E^\alpha$, for some $\alpha \in \R^{>0}$.
\end{enumerate}

\item Let $M$ be a complex manifold,
and $\alpha$ a closed 1-form. Prove that $d^c(\alpha)$
is a (1,1)-form.

\item
Let $M$ be a complex manifold.
Consider the map $H^1(M, \C) \stackrel\delta \arrow H^{1,1}_{BC}(M)$
mapping $\alpha$ to $d^c\alpha$, where $H^{1,1}_{BC}(M)$ denotes the Bott--Chern cohomology. 
\begin{enumerate}
\item Prove that $\ker \delta=H^{1,0}(M)\oplus
\overline{H^{1,0}(M)}$, where $H^{1,0}(M)$
denotes the space of closed holomorphic 1-forms.
\item Prove that the following sequence is exact
\begin{multline*}
0 \arrow H^{1,0}(M)\oplus
\overline{H^{1,0}(M)} \arrow H^1(M, \C)
\stackrel \delta \arrow \\ H^{1,1}_{BC}(M) \stackrel \tau \arrow H^2(M, \C),
\end{multline*}
where $\tau$ is the tautological map taking a
class represented by a form $\alpha$ to a class represented
by the same form $\alpha$.
\end{enumerate}

\item
Let $M$ be a compact complex manifold, $\dim_\C M=n$, equipped
with a Gauduchon metric $\omega$, and $L$ a line bundle on $M$.
Define {\bf the degree} of $L$ by $\deg_\omega (L):= \int_M \Theta_L \wedge \omega^{n-1}$,
where $\Theta_L$ is the curvature of the  Chern connection\index[terms]{connection!Chern} associated
with a Hermitian metric $h$ on $L$. Prove that $\deg_\omega (L)$
is independent on  $h$.\index[terms]{bundle!vector bundle!degree of}

\item Let $(M, \theta, \omega)$
be a compact Vaisman manifold, and $L$ its
weight bundle. Prove that $\omega$ is Gauduchon.
Prove $\deg_\omega (L)>0$, and $L$ is non-trivial
as a holomorphic line bundle.\index[terms]{manifold!Vaisman}

\item \label{_weight_bun_non-tri_Exercise_}
Let $(M, \theta, \omega)$ be
a compact LCK manifold\index[terms]{manifold!LCK}, and $L$ its weight bundle.
\begin{enumerate}

\item
Prove that there exists a Hermitian structure on $L$
such that the curvature of its Chern connection\index[terms]{connection!Chern}
is equal to $d^c\theta$.

\item
Suppose that $L$ is trivial as a
holomorphic line bundle. Prove that $d^c\theta_1=0$
for some LCK structure\index[terms]{structure!LCK} $(\theta_1, \omega_1)$
 conformally equivalent to $(\theta, \omega)$

\item
Prove that any compact LCK manifold that satisfies
$d^c\theta=0$ is actually  K\"ahler.

\item Prove that $L$ is non-trivial as a
 holomorphic line bundle.
\end{enumerate}

{\em Hint:} Prove that $\int_M dd^c(\omega^{n-1})>0$
whenever $d^c\theta=0$.

\end{enumerate}


%
%


\chapter{Existence of positive potentials}
\label{_posi_pote_Chapter_}


\epigraph{\it Have you not done tormenting me with your accursed time! It's abominable! When! When! One day, is that not enough for you, one day he went dumb, one day I went blind, one day we'll go deaf, one day we were born, one day we shall die, the same day, the same second, is that not enough for you? They give birth astride of a grave, the light gleams an instant, then it's night once more..}{\sc\scriptsize Samuel Beckett, \ \ Waiting for Godot}


\section{Introduction}


The existence of an LCK potential\index[terms]{potential!LCK} is essentially
a cohomological property. Indeed, an LCK form $\omega$ that admits
an LCK potential has vanishing class 
$[\omega]_{BC}$ in the
twisted Bott--Chern cohomology\index[terms]{cohomology!Bott--Chern!twisted} group $H^{1,1}_{BC}(M,L)$
(\ref{_LCK_pot_via_d_theta_d^c_theta_Proposition_}). 
In \cite{ov_jgp_09}, we stated this cohomological
equivalence as obvious. However, 
this observation does not take
into account the positivity of an LCK potential,\index[terms]{potential!LCK}
which makes the equivalence more difficult.

In this chapter, we resolve this problem.

Recall that an LCK metric\index[terms]{metric!LCK} on an LCK manifold\index[terms]{manifold!LCK}
$(M,\omega, \theta)$ can be understood
as a K\"ahler metric taking values in a 
flat line bundle $L$, called {\em the weight bundle}.
We denote the de Rham operator on $L$-valued\index[terms]{bundle!weight}
forms by $d_\theta$; it acts as 
$d_\theta(\alpha)= d\alpha -\theta\wedge \alpha$,
where $\theta $ is the Lee form.\index[terms]{form!Lee}
Then the equation for LCK potential\index[terms]{potential!LCK}
can be written as $\omega = d_\theta d^c_\theta \phi_0$.
A function that satisfies this equation is called
{\bf a $d_\theta d^c_\theta$-potential}.
When one does not fix the left-hand side,
one arrives at the definition of 
{\bf $d_\theta d^c_\theta$-plurisubharmonic
functions}, that is, real-valued functions
$\psi$ with $d_\theta d^c_\theta \psi$ a 
positive (1,1)-form. When this form is
strictly positive, we say that $\psi$ is
{\bf strictly $d_\theta d^c_\theta$-plurisubharmonic};
this is equivalent to being a $d_\theta d^c_\theta$-potential
for some LCK form.\index[terms]{form!LCK}

Notice that the term ``LCK potential''
is reserved for positive automorphic functions on 
the K\"ahler cover $\tilde M$; this is equivalent
to having a positive $d_\theta d^c_\theta$-potential.
We use the term ``$d_\theta d^c_\theta$-potential''
for $\phi_0$ that satisfies $\omega = d_\theta d^c_\theta \phi_0$
to emphasize that the function is not necessarily positive.

Originally, we expected that all 
$d_\theta d^c_\theta$-potentials
be positive. Indeed, this seems to be true in
many examples, such as in Exercise 
\ref{_LCK_pote_unique_Exercise_}.
In 2016, Victor \index[persons]{Vuletescu, V.} Vuletescu produced
an example of a $d_\theta d^c_\theta$-potential
on a Hopf manifold $H$ that is  {\em not} positive
(Section \ref{_counterexample_positive_Section_}).
However, in his example, $H$ admits another
 $d_\theta d^c_\theta$-potential that is  positive.

For a long time we thought that this situation
is universal, and that any LCK metric\index[terms]{metric!LCK} $\omega$ 
with $[\omega]_{BC}=0$ admits a positive
$d_\theta d^c_\theta$-potential.
However, even if this were true,
we could not prove it. Instead,
we prove a weaker theorem, which
is completely sufficient for all our purposes,
showing that any LCK manifold\index[terms]{manifold!LCK}
$(M,\omega, \theta)$ with $\omega$
$d_\theta d^c_\theta$-exact
admits another LCK metric
$\omega'$ such that
$\omega' = d_\theta d^c_\theta \phi_0$,
for a positive $\phi_0 \in C^\infty M$.

The argument is based on the following observation,
that we learned from Matei \index[persons]{Toma, M.} Toma. Suppose that $\phi_0$ is a negative
$d_\theta d^c_\theta$-potential on an LCK manifold.
Then $dd^c(\log (-\phi_0))$ is K\"ahler. 

Indeed, 
the function $t \mapsto \log(-t)$ is convex\index[terms]{function!convex} and monotonous
on the half-line $]-\infty, 0[$, and a composition
of a plurisubharmonic function and a monotonous convex
function is again plurisubharmonic. Applying this 
argument to a negative automorphic K\"ahler potential $\phi$
on the K\"ahler cover $\tilde M$, we obtain a
function $u:=\log (-\phi)$ that satisfies $\gamma^* u= u +\const$
for any $\gamma$ in the deck group $\Gamma$ of $\tilde M$, and hence 
$dd^c u$ is a $\Gamma$-invariant K\"ahler form\index[terms]{form!K\"ahler} on $\tilde M$.

The form $d^c u$ is also $\Gamma$-invariant.
This implies, in particular, that $M$ admits an exact K\"ahler
form, and hence  $M$ cannot be compact (unless it has a boundary). 

This also implies that $M$ is Stein whenever $M$ \index[terms]{manifold!Stein!with boundary} 
is a compact manifold with smooth boundary \index[terms]{manifold!LCK!with boundary}
$\6 M$ and $\phi_0=0$ on $\6 M$. Indeed, in this case the function
$u$ is constant on the boundary of $\tilde M$.
Since $u$ is strictly plurisubharmonic, the
boundary of $\tilde M$ is strictly pseudoconvex
(\ref{_psh_then_pseudoconvex_Corollary_}).
However, strict pseudoconvexity
is a local property, and hence  $\6 M$ is also 
strictly pseudoconvex. 
By \index[persons]{Grauert, H.} Grauert's solution of the Levi problem\index[terms]{Levi problem}
(\ref{_Grauert_Levi_Theorem_}), this implies that
$M$ is holomorphically convex. By the Remmert
reduction theorem (\ref{rem_red}),\index[terms]{theorem!Remmert reduction}
a holomorphically convex manifold which\index[terms]{variety!complex!holomorphically convex}
has no positive-dimensional compact
complex subvarieties is Stein. However,
$M$ has no compact complex subvarieties,
because it admits an exact K\"ahler form,\index[terms]{form!K\"ahler}
hence it is Stein.

We spend some time discussing Stein LCK\index[terms]{manifold!LCK!Stein}
manifolds, and prove that all such manifolds
admit a positive $d_\theta d^c_\theta$-potential
(\ref{_D_has_pos_d_theta_potential_Proposition_}).

Return now to a compact LCK manifold
$(M, \theta, \omega)$ without boundary
equipped with a $d_\theta d^c_\theta$-potential $\phi_0$.
If $\phi_0$ is non-negative everywhere, the
function $\phi_0+\epsilon$ is an LCK potential\index[terms]{potential!LCK}
for $\epsilon$ sufficiently small, and we are done.
Otherwise, $M$ is a union
of $D:=\phi_0^{-1}(]-\infty, \epsilon])$ and $M \backslash 
D$.\footnote{The set $\phi_0^{-1}(]-\infty, \epsilon])$
is compact, because $\phi_0$ is a continuous
function on a compact set $M$, and hence  it is bounded from below.}
We use the gluing formalism, developed in 
K\"ahler geometry\index[terms]{geometry!K\"ahler} (mostly due to J.-P. \index[persons]{Demailly, J.-P.} Demailly)
to glue the positive $d_\theta d^c_\theta$-potential on $D$ 
(that is  Stein) to the positive $d_\theta d^c_\theta$-potential 
$\phi_0$ on $M \backslash D$ to obtain a 
positive $d_\theta d^c_\theta$-potential on $M$.

The gluing of K\"ahler metrics is a very useful construction
based on the notion of {\em regularized maximum} of K\"ahler potentials.
Let $\max_\epsilon:\; \R \times \R \arrow \R$ be \index[terms]{regularized maximum}
a convex, smooth function, monotonous in both arguments,
and satisfying $\max_\epsilon(u, v) = \max (u, v)$ whenever 
$|u-v|>\epsilon$. It is not hard to see that
$\max_\epsilon(\phi, \psi)$ is strictly plurisubharmonic
whenever $\phi$ and $\psi$ are strictly plurisubharmonic
functions on a complex manifold.
To write the function $\max_\epsilon$ explicitly,
we note that $\max(u, v) = \max(u-v,0) + v$.

To produce $\max_\epsilon$, we replace
the absolute value function by its approximation
$R$ that is  smooth, convex, and satisfies
$R(u) = |u|$ whenever $|u| > \epsilon$, 
and write $\max_\epsilon(u, v) = R(u-v) + v$.

\centerline{\includegraphics[width=0.3\linewidth]{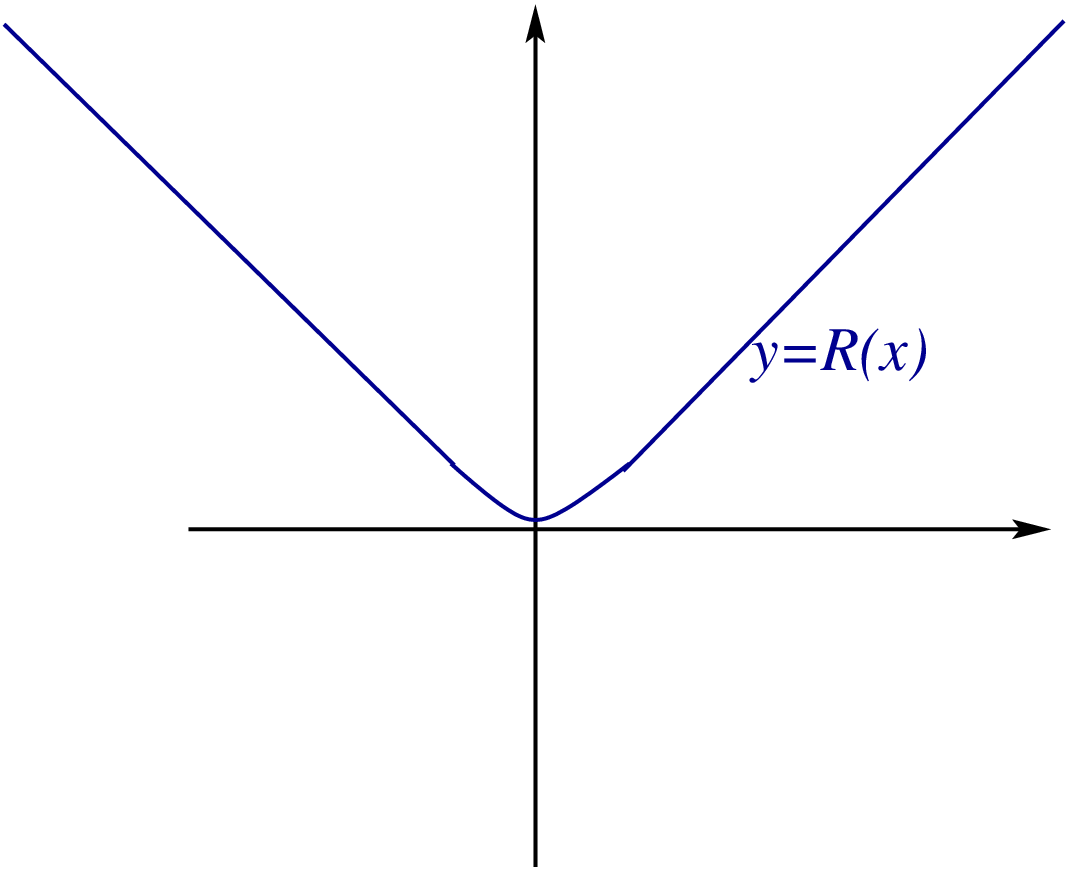}}

The regularized maximum is used when two K\"ahler
metrics $\omega_1, \omega_2$ on $M$
have potentials $\phi_1, \phi_2$ in a neighbourhood of a hypersurface
$S\subset M$, with $\phi_1\restrict S = \phi_2 \restrict S$
and $\Lie_X \phi_1 < \Lie_X \phi_2$ for a vector field
$X \in TM\restrict S$ transversal to $S$.
In this case, for $\epsilon$ sufficiently small,
the regularized maximum $\max_\epsilon(\phi_1, \phi_2)$
is equal to $\phi_1$ on one side of $S$ outside
of its neighbourhood, and equal to $\phi_2$
on the other side, and hence  $dd^c(\max_\epsilon(\phi_1, \phi_2))$
can be glued to $\omega_1$ on one side and to $\omega_2$ on another
side. This gives a K\"ahler form\index[terms]{form!K\"ahler} that is  equal to $\omega_1$
on the one side of $S$ and to $\omega_2$ on another side.

In the LCK situation, we notice that
a regularized maximum of two strictly $d_\theta d^c_\theta$-plurisubharmonic 
functions is again strictly $d_\theta d^c_\theta$-plurisubharmonic.
Now, suppose that $M$ is a compact LCK manifold\index[terms]{manifold!LCK}, $D\subset M$ 
a Stein subset with smooth boundary equipped with 
a positive $d_\theta d^c_\theta$-potential $\psi$, and
$\phi_0$ is a $d_\theta d^c_\theta$-potential which
satisfies $\phi_0 > 3\epsilon$ outside of $D$.
Then $\max_\epsilon(\delta \psi, \phi_0)$
is positive,  $d_\theta d^c_\theta$-plurisubharmonic 
and equal to $\phi_0$ outside of a sufficiently 
small neighbourhood of $D$, whenever $\delta$
is a positive number that satisfies 
$\delta < \epsilon (\sup_{x\in D}\psi(x))^{-1}$
(\ref{_main_positive_Theorem_}).

This proof is different from the one
published in \cite{ov_poz_pot}, which was much more
complicated and relied on classical (but very difficult)
results from complex analysis, due to H. J. \index[persons]{Bremermann, H. J.} Bremermann.


\section{A counterexample to the positivity of the
  potential}\label{_counterexample_positive_Section_}


Automorphic potentials are not necessarily positive,
as seen in the following.

\hfill

\example \label{_non_posi_potential_Victors_Example_}
{(courtesy of V. \index[persons]{Vuletescu, V.} Vuletescu):} 
Take a Hopf manifold $\C^n \backslash 0/\Z$
where $\Z$ acts by multiplication with $\lambda >1$
and let $L$ be the weight bundle seen as a local system\index[terms]{local system} with the same\index[terms]{bundle!weight}
monodromy\index[terms]{monodromy}. Then the usual flat K\"ahler form\index[terms]{form!K\"ahler}
on $\C^n$ can be considered to be  a closed
Hermitian form\index[terms]{form!Hermitian} with values in the bundle $L^2$.
Its LCK potential \index[terms]{potential!LCK}is a function $l(z):= |z|^2$.
Any quadratic polynomial on $\C^n$ gives a holomorphic
section of $L^2\otimes_\R \C$; let $v$ be its real part.
Then $d_\theta d^c_\theta (v)=0$, because $v$ is
the real part of a holomorphic section of $L^2$,
and $d_\theta d^c_\theta (l + Av)= d_\theta d^c_\theta(l)$
is the LCK form\index[terms]{form!LCK} on $M$, for any real constant
$A$. However,
for $A$ large enough, the $d_\theta d^c_\theta$-potential $l+Av$
is non-positive.

\hfill

Note that the obvious solution, which one would use
in the K\"ahler case to make a $d_\theta d^c_\theta$-potential positive (adding a constant) does not
work, because the operator $d_\theta d^c_\theta$ 
does not vanish on constants. We shall need 
to find a positive function $h$ such that  
$d_\theta d^c_\theta (h)$ is non-negative; adding
$C d_\theta d^c_\theta (h)$ to $\omega$, $C\gg 0$, gives 
us an LCK form\index[terms]{form!LCK} with positive potential in
the same cohomology class. The precise result is.

\hfill

\theorem\label{_main_positive_intro_Theorem_}{( \cite{ov_poz_pot})} 
Let $M$ be a compact LCK manifold\index[terms]{manifold!LCK}  with a K\"ahler cover 
admitting an automorphic K\"ahler potential (not
necessarily positive).
Then $M$ also admits an LCK metric with a
positive automorphic potential.\index[terms]{potential!positive}

\hfill

\proof See \ref{_main_positive_Theorem_}. \endproof

\hfill

We start by showing that any $d_\theta d^c_\theta$-potential 
for an LCK metric is positive somewhere.


 \section{A $d_\theta d^c_\theta$-potential on a compact LCK manifold
is positive somewhere}\label{_nega_}


\theorem\label{_strictly_negative_pots_Theorem_}
Let $(M, \theta, \omega)$ be an LCK manifold\index[terms]{manifold!LCK}
that is  not K\"ahler,
and $\omega= d_\theta d^c_\theta(\phi_0)$ for some
smooth  $d_\theta d^c_\theta$-plurisubharmonic function
$\phi_0\in C^\infty (M)$.
Then $\phi_0 >0$ at some point of $M$.

\hfill

\proof By absurd, suppose $\phi_0\leq 0$
everywhere on $M$. Let $\pi:\; \tilde M\arrow M$ be the minimal 
K\"ahler cover of $M$, and $\rho$ a function
on $\tilde M$ such that $d\rho=\theta$.
Then the $d_\theta d^c_\theta$-plurisubharmonicity of
$\phi_0$ is equivalent to plurisubharmonicity of
$\phi:=e^{-\rho}\pi^*\phi_0$ on $\tilde M$.

Strict $d_\theta d^c_\theta$-plurisubharmonicity is stable under $C^2$-small deformations
of $\phi_0$; therefore, the function $\phi_0-\epsilon$ is also strictly
$d_\theta d^c_\theta$-plurisubharmonic. Therefore, we may 
assume that $\phi <0$ everywhere.

Define
$$\psi:=-\log(-\phi).$$ 
Since $x\mapsto -\log(-x)$ is strictly
monotonous and convex, the function
$\psi$ is strictly plurisubharmonic.  Moreover,
for every $\gamma\in\Gamma$, we have
\[
\gamma^*\psi=-\log(-(\phi\circ\gamma))=-
\log(\chi(\gamma))-\log(-\phi)=\const+\psi.
\]
Therefore, the K\"ahler form\index[terms]{form!K\"ahler}
$dd^c\psi$ is $\Gamma$-invariant on $\tilde M$, and $M$ is
K\"ahler.
\endproof

\hfill

\remark We are grateful to Matei \index[persons]{Toma, M.} Toma for this beautiful argument.

\section{Stein manifolds with negative $d_\theta d^c_\theta$-potential}

Let $(M, \theta, \omega)$ be a compact LCK manifold admitting
a $d_\theta d^c_\theta$-potential $\phi_0$.  By 
\ref{_strictly_negative_pots_Theorem_}, $\phi_0$
is strictly positive somewhere in $M$.
Denote by $D$ the set of points where $\phi_0$ is
negative (adjusted to make its boundary smooth). 
It turns out that the manifold $D$
is Stein and admits a strictly positive 
$d_\theta d^c_\theta$-potential.

We start by stating a classical result about
holomorphically convex manifolds.\index[terms]{variety!complex!holomorphically convex}

\subsection{Remmert theorem and 1-jets on Stein manifolds}

\theorem ({\bf  Remmert \index[persons]{Remmert, R.}reduction}, \cite{rem})\label{rem_red} 
Let $X$ be a holomorphically convex
manifold. Then there exist a Stein variety $Y$ and a proper, surjective, holomorphic map\index[terms]{theorem!Remmert reduction}\index[terms]{manifold!holomorphically convex}
$f : X \ra Y$ such that
\begin{enumerate}
	\item[(i)] $f_*\calo_X = \calo_Y$.
	
	Moreover, the fact that $Y$ is Stein and (i) imply:
	\item[(ii)] $f$ has connected fibres.
	\item[(iii)]  The map $f^* : \calo_Y (Y)\ra \calo_X (X)$ is an isomorphism.
	\item[(iv)] The pair $(f, Y)$ is unique up to biholomorphism, i.\,e.   for any other pair $(f_0, Y_0)$
	with $Y_0$ Stein and property (i), there exists a biholomorphism $g : Y \ra Y_0$ such that
	$f_0 = g \circ f$. 
\end{enumerate}

We will also need the following result about 
Stein manifolds.\index[terms]{manifold!Stein} Recall that {\bf the $k$-jet
of a section of the bundle $B$ in $x$} is 
an element of the space $H^0(B/ {\goth m}_x^{k+1})$,
where ${\goth m}_x$ is the maximal ideal of $x$.

\hfill

\lemma\label{_Stein_globally_genera_1-jet_Lemma_}\index[terms]{function!jet of}
Let $D$ be a compact Stein manifold with smooth boundary,
and $(B, \nabla)$ a complex vector bundle equipped
with a flat connection. We consider $B$\index[terms]{connection!flat}
as a holomorphic vector bundle, using the
holomorphic structure operator $\bar\6=\nabla^{0,1}$.
Then there exists a collection of holomorphic sections
$f_1, ..., f_k$ of $B$ such that $\nabla^{1,0} f_1, ..., \nabla^{1,0} f_k$
generate the bundle $\Lambda^{1,0}(D) \otimes B$.

\hfill

\proof
For each $x\in D$, and any holomorphic section $f$ of $B$,
the form $\nabla^{1,0} f$ is identified with the 1-jet
of $f$ in $x$. On the other hand, 1-jets can be understood
as elements of $H^0(B/ B\otimes {\goth m}_x^2)$, where ${\goth m}_x$
is the maximal ideal of $x$ in $\calo_D$.

Since $D$ is Stein, the natural map
$H^0(B) \arrow H^0(B/ B\otimes {\goth m}_x^2)$
is surjective. Indeed, the exact sequence of coherent sheaves\index[terms]{sheaf!coherent}
\[
0 \arrow B\otimes {\goth m}_x^2 \arrow B \arrow B/ B\otimes {\goth m}_x^2
\arrow 0
\]
gives a long exact sequence
\[
... \arrow H^0(B) \arrow H^0(B/ B\otimes {\goth m}_x^2)
\arrow H^1(B\otimes {\goth m}_x^2)\arrow ...
\]
and the group $H^1(B\otimes {\goth m}_x^2)$ 
vanishes because $D$ is Stein.
Therefore, there exists a holomorphic
section of $B$ with any prescribed 1-jet.

This implies that for each $x\in D$ there
exists a collection of sections $f_1, ..., f_l$ of $B$
such that the restrictions $\nabla^{1,0} f_1, ..., \nabla^{1,0} f_l$
to $x\in D$ generate the space $\Lambda^{1,0}_x(D) \otimes B$.
Therefore, the bundle $\Lambda^{1,0}(D) \otimes B$
is generated by  $\nabla^{1,0} f_1, ..., \nabla^{1,0} f_l$
in a neighbourhood $U_x$ of $x$. 

Consider a cover of $D$ by such $U_x$. Since $D$ is
a compact manifold with boundary, it has a finite subcover. The 1-forms
$\nabla^{1,0} f_1, ..., \nabla^{1,0} f_k$
corresponding to this
collection of holomorphic sections 
of $B$ generate $\Lambda^{1,0}(D) \otimes B$.
\endproof

\subsection{Negative sets for $d_\theta d^c_\theta$-potentials are Stein}

\lemma\label{_level_set_for_d_theta_potential_Lemma_}
Let $\phi_0$ be a $d_\theta d^c_\theta$-potential
on an LCK manifold\index[terms]{manifold!LCK} $(M, \theta, \omega)$.
Assume that 0 is a regular value of $\phi_0$.
Then the CR-manifold $\phi_0^{-1}(0)\subset M$
is strictly pseudoconvex.\index[terms]{manifold!CR!strictly pseudoconvex}

\hfill

\proof
The strictly pseudoconvexity of $\phi_0^{-1}(0)\subset M$
is a local property, and hence  we can check it on a K\"ahler
cover $\tilde M \stackrel \pi \arrow M$. The pullback of
$\phi_0^{-1}(0)$ to $\tilde M$ is the level set of a K\"ahler
potential, and hence  it is strictly pseudoconvex 
(\ref{_psh_then_pseudoconvex_Corollary_}).
\endproof

\hfill

\proposition\label{_negative_set_for_d_theta_potential_Proposition_}
Let $M$ be a compact LCK manifold,\index[terms]{manifold!LCK}
$\phi_0$ a $d_\theta d^c_\theta$-potential,
$c$ a regular value of $\phi_0$, and $D:= \phi_0^{-1}(]-\infty, c])$.
Then $D$ is Stein for $c$ sufficiently
close to 0.

\hfill

\proof
A small deformation of a $d_\theta d^c_\theta$-potential
is again a $d_\theta d^c_\theta$-potential, and hence 
$\phi_0-c$ is a $d_\theta d^c_\theta$-potential for
$c$ sufficiently small. 
Indeed,
\[
d_\theta d^c_\theta(u) = dd^c u + \theta \wedge d^c u +
du\wedge \theta^c + u \theta \wedge \theta^c;
\]
positivity of this form is preserved if $u$
is $C^2$-close to $\phi_0$. Then, by
\ref{_level_set_for_d_theta_potential_Lemma_},
the level set $\phi_0^{-1}(c)= (\phi_0-c)^{-1}(0)$ is
strictly pseudoconvex.

Note that a compact manifold $D$ with smooth,
strictly pseudoconvex boundary is holomorphically convex
by the solution of Levi problem (\cite{_Grauert_}).
Then the  Remmert reduction (\ref{rem_red}) implies that $D$
admits a proper, surjective  and holomorphic \index[terms]{variety!Stein}
map $\pi:\; D \arrow D_0$ with connected fibres to a Stein
variety $D_0$ with isolated singularities.

Denote by $\tilde D$ the smallest cover of
$D$ where the Lee form\index[terms]{form!Lee} $\theta$
is exact. Let $\psi\in C^\infty \tilde D$ satisfy $d\psi=\theta$, and let
$\Phi:= e^{-\psi} (\phi_0-c)$ be
the negative K\"ahler potential on 
$\tilde D$ obtained from $\phi_0-c$.
Then the function $\Psi:=-\log(-\Phi)$ is strictly plurisubharmonic
on $\tilde D$ (see \ref{_strictly_negative_pots_Theorem_}),
hence the 1-form $d^c\Psi$ is the pullback of a
1-form $\Theta$ on  $D$.
On the other hand, $d\Theta$ is a K\"ahler form.\index[terms]{form!K\"ahler}
This implies that $D$ has no compact subvarieties (without boundary),
and the map $\pi:\; D \arrow D_0$ is bijective.
Then $D$ is Stein.
\endproof

\subsection{Stein LCK manifolds
admit a positive $d_\theta d^c_\theta$-potential}\index[terms]{manifold!LCK!Stein}

\proposition\label{_D_has_pos_d_theta_potential_Proposition_}
Let $(D, \theta, \omega)$ be a compact LCK manifold\index[terms]{manifold!LCK!with boundary} with smooth boundary.
Assume that $D$ is Stein.
Then $D$ admits a positive $d_\theta d^c_\theta$-potential.

\hfill

\proof 
Let $L$ be the weight bundle on $D$.\index[terms]{bundle!weight}
The bundle $L$ can be considered to be  a trivial
real line bundle $L_0$ with the connection $\nabla:=\nabla_0+\theta$,
where $\nabla_0$ is the trivial connection.
Then $\bar\6_\theta $ is the holomorphic structure 
operator for sections of $L\otimes_\R \C$.
Here $L\otimes_\R \C$ is understood as
a holomorphic line bundle on $D$, with the holomorphic
structure induced from the flat connection $\nabla$.
However, its sections are identified with the sections
of the trivial bundle $L_0\otimes_\R \C$, and hence 
with complex-valued functions on $M$.

Every holomorphic section $l$ of $L$ 
satisfies $\bar\6_\theta l =0$, which gives
\[d_\theta d^c_\theta(|l|^2)= - 2\1 \6_\theta \bar\6_\theta(|l|^2)=
- 2\1 \6_\theta(l) \wedge \bar\6_\theta(\bar l).
\]
This implies that $|l|^2$ is
$d_\theta d^c_\theta$-plurisubharmonic.

Since $D$ is Stein,
the space of sections of a holomorphic bundle is 
globally generated; this assures the existence of
sufficiently many holomorphic sections of $L$.
Then, for a sufficiently big collection $f_1, ..., f_k$ of sections
of $L$, the sum $\sum |f_i|^2$ is positive 
and strictly $d_\theta d^c_\theta$-plurisubharmonic
everywhere on $D$. To assure this, we need
$f_i$ to have no common zeros, and
$\6_\theta f_i$ to generate $\Lambda^1 D\otimes L$ at every
point of $D$, that is  possible
by \ref{_Stein_globally_genera_1-jet_Lemma_}.

We have proven that $D$ admits a positive
LCK potential.\index[terms]{potential!LCK} 
\endproof

\hfill

To find a positive $d_\theta d^c_\theta$-potential
on a manifold $(M, \theta, \omega)$ admitting
an $d_\theta d^c_\theta$-potential $\phi_0$,
we take a positive $d_\theta d^c_\theta$-potential
on the Stein manifold $\phi_0^{-1}(]-\infty, c])$
and glue it to $\phi_0$, as indicated below.
 

\section{Gluing the LCK forms}\label{_gluing_}\index[terms]{form!LCK}


\subsection[Regularized maximum of plurisubharmonic functions]{Regularized maximum of plurisubharmonic\\ functions}

\label{_reg_max_subsection_}
In \cite{_Demailly_1982_}, the notion
of a {\em regularized maximum} of two functions was
defined as follows.  \index[terms]{regularized maximum}\index[persons]{Demailly, J.-P.}

\hfill

\definition (\cite{_Demailly_1982_}) Choose $\epsilon >0$, and let 
$\max_\epsilon:\; \R^2\arrow \R$ be a smooth, convex\index[terms]{function!convex}
function, monotonous in both variables,
that satisfies $\max_\epsilon(x, y) = \max(x,y)$
whenever $|x-y|>\epsilon$. Then $\max_\epsilon$ is called
{\bf a regularized maximum}. 

\hfill

\remark The regularized maximum exists for any $\epsilon >0$.
Take any convex, monotonous function $\rho_\epsilon:\; \R \arrow
\R^{\geq 0}$ such that $\rho_\epsilon(x)=0$ for $x < -\frac \epsilon
2$ and $\rho_\epsilon(x)=x$ for $x > -\frac \epsilon 2$.
Then $\max_\epsilon(x, y) := \rho_\epsilon(y-x)+x$
is a regularized maximum.

\hfill

A simple computation brings the following.

\hfill

\claim (\cite{_Demailly_1982_}) \label{reg_max} 
A regularized maximum\index[terms]{theorem!Demailly}
of two plurisubharmonic functions is
again plurisubharmonic.
\endproof

\hfill

Taking regularized maxima of K\"ahler potentials
to glue the K\"ahler forms\index[terms]{form!K\"ahler} is a well-known procedure; it was 
much used by J.-P. \index[persons]{Demailly, J.-P.} Demailly (see for example   
\cite{dp}), and, in the LCK context,
in our paper \cite{ov_sas} and in \cite{_Brunella:Kato_} 
(see also Chapter \ref{sasemb}). Here
we apply the same argument to $d_\theta d^c_\theta$-potentials.

\hfill

\definition
Let $L$ be a line bundle over 
a complex manifold $M$ with a flat connection, and
$\psi$ a section of $L$. Assume that $L$ is trivialized,
with $\theta$ its connection form; we may consider
$\theta$ as a closed 1-form on $M$. Then the de Rham differential
on $L$-valued forms can be written as $d_\nabla=d-\theta=d_\theta$,
and the twisted de Rham differential as $Id_\nabla I^{-1}= d_\theta^c$. 
We say that a function $\psi$ is {\bf
  $d_\theta d_\theta^c$-plurisubharmonic}
if $d_\theta d_\theta^c(\psi)$ is a semi-positive Hermitian
form on $M$. This is equivalent to plurisubharmonicity
for $L$-valued functions. When $d_\theta d_\theta^c(\psi)$
is Hermitian, we say that $\psi$ is {\bf 
strictly $d_\theta d_\theta^c$-plurisubharmonic}.

\hfill

\remark
In Section
\ref{_potential_weight_Subsection_} we defined {\em $d_\theta
  d^c_\theta$-potential}; a smooth function is a
$d_\theta  d^c_\theta$-potential
if and only if it is strictly $d_\theta d_\theta^c$-plurisubharmonic.

\hfill

Later on, we apply this construction to the weight bundle
on an LCK manifold.\index[terms]{manifold!LCK}
We extend \index[persons]{Demailly, J.-P.} Demailly's method of gluing plurisubharmonic\index[terms]{bundle!weight}
functions to this context.

\hfill

\claim\label{_d_theta_reg_max_Claim_}
Let $\theta$ be a closed form on a complex manifold $M$,
and $\phi, \psi$ two $d_\theta d_\theta^c$-plurisubharmonic
functions. Then $\max_\epsilon(\phi, \psi)$ is also
$d_\theta d_\theta^c$-plurisubharmonic; it is
strictly $d_\theta d_\theta^c$-plurisubharmonic if
$\phi$ and $\psi$ are strictly $d_\theta d_\theta^c$-plurisubharmonic.

\hfill

\proof Since this result is local,  
we may always assume that $\theta= d\rho$ for some
function  $\rho\in C^\infty M$.
Then $d_\theta(\eta) = e^{-\rho}d(e^\rho\eta)$
and 
\[ 
d_\theta d_\theta^c(\max_\epsilon(\phi, \psi))=e^{-\rho} dd^c (\max_\epsilon(e^\rho \phi,e^\rho \psi)).
\]
Since $e^\rho \phi$, $e^\rho \psi$ are plurisubharmonic, the form
$dd^c (\max_\epsilon(e^\rho \phi,e^\rho  \psi))$ is positive, by \ref{_d_theta_reg_max_Claim_}.
\endproof

\subsection{Gluing of LCK potentials}\index[terms]{potential!LCK}

\proposition \label{_gluing_LCK_Proposition_}
{\bf (Gluing of LCK metrics)}  \\
Let $(M, \omega)$ be an LCK manifold,\index[terms]{metric!LCK}
with $\omega = d_\theta d^c_\theta \phi_0$.
Consider a number $\epsilon >0$ and let 
$\mu$ be a regular value of $\phi_0$ satisfying
$0>\mu > -\epsilon$. Let $D:= \phi_0^{-1}(]-\infty, \mu]$.
By \ref{_negative_set_for_d_theta_potential_Proposition_} 
and \ref{_D_has_pos_d_theta_potential_Proposition_}, 
$D$ is Stein and admits a positive LCK potential\index[terms]{potential!LCK} $\psi$. 
Then, for a sufficiently small
$\epsilon>0$ and  $\delta>0$,
the function 
\begin{equation*}
	\xi=
	\begin{cases}
		\max_\epsilon(\psi\delta, \phi_0+3\epsilon) \quad \text{on $D$},\\
		\phi_0+3\epsilon \quad \text{outside of $D$}
    \end{cases}
\end{equation*}		
is a positive strictly $d_\theta d^c_\theta$-plurisubharmonic
function.

\hfill

{\bf Proof:} 
Since a small deformation of $d_\theta d^c_\theta$-potential
is a $d_\theta d^c_\theta$-potential, the function $\phi_0+3\epsilon$
is strictly $d_\theta d^c_\theta$-plurisubharmonic for $\epsilon$
sufficiently small. 
Choose $\delta:= (\max_{D} \psi)^{-1} \epsilon$. Then
$ \psi\delta < \epsilon$ on $D$, and hence  near
the boundary of $D$, one has 
$\phi_0+3\epsilon\geq \mu + 2\epsilon > \epsilon > \psi\delta $,
and $\xi= \phi_0+3\epsilon$ near the boundary of $D$.
Therefore, $\xi$ is a positive 
$d_\theta d^c_\theta$-potential on $M$.
\endproof

\hfill

We can now finish the proof of the main result of this chapter.

\hfill

\theorem\label{_main_positive_Theorem_}
Let $(M, \omega, \theta)$ be a compact LCK manifold\index[terms]{manifold!LCK}
admitting a $d_\theta d^c_\theta$-potential $\phi_0$.
Then $M$ also admits an LCK metric with LCK 
potential.\footnote{By our convention, {\em an LCK
potential} is a positive K\"ahler potential on a
K\"ahler cover of $M$. The existence of
LCK potential\index[terms]{potential!LCK} is equivalent to 
the existence of $d_\theta d^c_\theta$-potential 
that is  positive.}

\hfill

\proof
By \ref{_strictly_negative_pots_Theorem_}, 
$\phi_0$ is positive somewhere on $M$.
Replacing $\phi_0$ by $\phi_0+\epsilon$,
for $\epsilon$ sufficiently small, we
can always assume that 0 is a regular value
of $\phi_0$. By \ref{_negative_set_for_d_theta_potential_Proposition_},
the manifold $D:=\phi_0^{-1}(]-\infty, 0])$ is
a Stein manifold with smooth boundary.
By \ref{_D_has_pos_d_theta_potential_Proposition_},
the manifold $D$ admits a positive 
strictly $d_\theta d^c_\theta$-plurisubharmonic function
$\psi$. Using \ref{_gluing_LCK_Proposition_},
we glue $\phi_0$ to $\psi$, obtaining
a positive, strictly $d_\theta d^c_\theta$-plurisubharmonic 
function $\xi$ on $M$.
\endproof

\section{Exercises}

\begin{enumerate}[label=\textbf{\thechapter.\arabic*}.,ref=\thechapter.\arabic{enumi}]

\item Let $\phi:\; \R \arrow \R$ be a monotonous, convex function,
and $\psi$ a plurisubharmonic function on a complex manifold $M$.
Prove that $x\mapsto \phi(\psi(x))$ is also plurisubharmonic.\index[terms]{function!convex}

\item Prove that any convex function on $\C^n$ is
plurisubharmonic.

\item
Let $P:\; \R^n \arrow\R$ be a convex smooth function
that is  strictly mo\-no\-to\-nously increasing in all arguments.
Prove that $P(u_1, ..., u_n)$ is strictly plurisubharmonic
for any collection $u_1, ..., u_n$ of strictly plurisubharmonic
functions on a complex manifold.

\item\label{_Reg_max_Exercise_}
Let $R:\; \R \arrow \R$ be a convex function
that is  strictly monotonously increasing on $[0, \infty[$
and satisfies $R(u) = |u|$ whenever $|u|>\epsilon$.
Prove that $\max_\epsilon(x, y):= \frac{R(x-y) + R(x+y)}{2}$
is strictly monotonous in both arguments and convex
as a function from $\R^2$ to $\R$.

\item Let $B\subset \C^n$ be an open ball.
Construct two
plurisubharmonic functions $f\neq g$ on $\C^n$ that are 
equal in $B$.

\item
Let $D$ be a Stein manifold.
Prove that there exists a holomorphic function
$f$ in the interior of $D$ that satisfies
$\lim_{i} |f(x_i)|=\infty$ for any
sequence $\{x_i \in D\}$ that has no 
limit points in $D$.

\item\label{_k_jets_diagonal_Exercise_}
Let $M$ be a complex manifold, and $J_\Delta$
the ideal of the diagonal in $M\times M$.\index[terms]{function!jet of}
Identifying $M$ with the diagonal, we can
consider $\calo_{M\times M}/ J_\Delta^k$
as a sheaf on $M$. Prove that the fibre of
$\calo_{M\times M}/ J_\Delta^{k+1}$ at $x\in M$
is naturally identified with the space of $k$-jets
of functions on $M$ in $x$.

\item Let $D$ be a compact Stein manifold with
smooth boundary.
Prove that there exists a collection
$\{f_i\}$ of holomorphic functions on $D$
such that the $k$-jets of $f_i$ generate
the bundle of $k$-jets of 
functions on $D$.\footnote{To define the
bundle of $k$-jets, use Exercise \ref{_k_jets_diagonal_Exercise_}.}

\item
A real-valued function $u$ is called {\bf pluriharmonic}\index[terms]{function!pluriharmonic} if $dd^cu=0$.
Show that $u$ is pluriharmonic if and only if $u$
is locally equal to the real part of a holomorphic function.

\item
Let $f$ be a pluriharmonic function on a K\"ahler
manifold. Prove that $f$ is harmonic.

\item
Let $f$ be a function on an open ball $B$ in $\C^n$.
Prove that $f$ is pluriharmonic if and only if 
for each complex line $l\subset \C^n$, the restriction
$f\restrict{l\cap B}$ is harmonic.

\item
Let $u_n$ be a sequence of pluriharmonic functions that are 
locally uniformly bounded. 
Show that there is a subsequence converging 
locally uniformly to a pluriharmonic function.

\item Let $f$ be a pluriharmonic function on $X=\C^2 \backslash 0$.
Prove that $f$ is the real part of a holomorphic function.
Find a counterexample to this statement for 
$X=\C\backslash 0$.

\item
Let $\lambda, \mu \in \C$, $\lambda\neq 0$.
Call a function $\rho$ on $\C^n$ 
{\bf $(\lambda, \mu)$-ho\-mo\-ge\-neous} if
$\rho(\lambda z)= \mu \rho(z)$.
Prove that $\C^n\backslash 0$ does not admit
non-zero $(\lambda, \mu)$-ho\-mo\-ge\-neous pluriharmonic 
functions unless $\mu=\lambda^n$ for some
$n\in \Z^{\geq 0}$.

\item\label{_LCK_pote_unique_Exercise_}
Let $H= \C^n\backslash 0/{\langle A\rangle}$
be a classical Hopf manifold, with $A = \alpha \Id$,
with $\alpha>1$ irrational, and let 
$\omega=\frac{dd^c (|z|^{2\alpha})}{|z|^{2\alpha}}$
be the LCK metric on $H$.\index[terms]{metric!LCK}
Prove that the $d_\theta d^c_\theta$-potential
for $(H, \omega)$ is unique and positive.

{\em Hint:} Use the previous exercise.

\item Let $f_1, ..., f_n$ be holomorphic functions on a
  complex manifold. Prove that $\log\left(\sum_i |f_i|^2\right)$ is
  plurisubharmonic.

{\em Hint:} Compute $dd^c\log \left(\sum |z_i|^2\right)$
on $\C^n$.

\item Let $f_1, ..., f_n$ be holomorphic functions on a
  manifold without common zeroes. 
Prove that $\log\left(\sum_i |f_i|^p\right)$ is
  plurisubharmonic, for any $p>0$.

\item (\cite{_moroianu:book_}) 
Let $(M,\theta, \omega)$ be a compact LCK manifold with
potential. Assume that $\omega$ is a Gauduchon metric,
that is, satisfies $d^* \theta=0$,
\index[terms]{metric!Gauduchon} and $\phi$ its $d_\theta d^c_\theta$-potential,
 $\omega=d_\theta d_{\theta}^c\phi$.
\begin{enumerate}
\item Show that $d_\theta^c\phi=-\phi I\theta$.

\item Rewrite explicitly the condition $d_\theta d^c_\theta(\phi)=\omega$ as:
\[
  \omega=dd^c\phi-\phi dI\theta+i\theta\wedge d\phi -
  \theta\wedge Id\phi  +\phi\theta\wedge I\theta.
\]

\item Take the scalar product with $\omega$ to obtain:
\begin{equation}\label{_omega_multiplied_by_dd^c_phi_Equation_}
n=g(dd^c\phi,\omega)-g(d(\phi I\theta),\omega)-g(\theta,
  d\phi)+\phi |\theta|^2.
\end{equation}

\item Prove that $d^*\omega=(1-n)I\theta$.

\item Taking into account that 
$d^*\omega=(1-n)I\theta$, integrate
  \eqref{_omega_multiplied_by_dd^c_phi_Equation_} 
on $M$ and prove that
 $$\vol(M)=\int_M\f |\theta|^2.$$
  Use this to obtain another proof of  
\ref{_strictly_negative_pots_Theorem_}. 
\end{enumerate}

\item {\bf (Gluing of K\"ahler metrics)} \\
Let $(M, \omega)$ be a K\"ahler manifold, and 
$D\subset M$ a submanifold of the same dimension with smooth
compact boundary such that $\omega = dd^c\phi$ in
a smooth neighbourhood of $D$, with $\phi$ a plurisubharmonic
function. Let $\psi$ be another plurisubharmonic function
with $\psi=\phi$ on $\6 D$. Consider a vector field 
$X\in TM\restrict{\6 D}$ that is  normal and
outward-pointing everywhere in $\6 D$. 
Assume that $\Lie_X \psi < \Lie_X \phi$ everywhere
on $\6 D$. Let $D_-$ be an open subset of $D$
that does not intersect a neighbourhood $U$ of $\6 D$,
and $D_+$ an open subset of $M \backslash D$ which
does not intersect $U$. 

\begin{enumerate}
\item
Prove that $\max_\epsilon(\psi, \phi)$
is equal to $\psi$ in
$D_-$ outside of the neighbourhood of 
$\6 D$ and equal to $\phi$ in $D_+$ 
outside of the neighbourhood of 
$\6 D$.

\item  Prove that there exists a K\"ahler form\index[terms]{form!K\"ahler}
$\omega_1$ that is  equal to $\omega$ on $D_+$
and to $d_\theta d^c_\theta \psi$ on $D_-$.
\end{enumerate}

\item 
Let $(\tilde M, \tilde \omega)$ be a K\"ahler 
manifold, equipped with a properly discontinuous
free action of a group $\Gamma$ by holomorphic
homotheties with compact quotient $M$
(a posteriori, this quotient is an LCK manifold\index[terms]{manifold!LCK}).
Take an automorphic K\"ahler potential
$\phi$ on $\tilde M$, and let
$u$ be another automorphic function on $\tilde M$.
Then $u^{-1} \phi$ is $\Gamma$-invariant, and obtained
as the pullback of $\phi_0\in C^\infty M$.
Prove that $u^{-1} \phi$ is a $d_\theta d^c_\theta$-potential
for an appropriate LCK structure\index[terms]{structure!LCK} $(\theta, \omega)$ on $M$.

\item
Let $(M, \theta, \omega)$ be an LCK manifold with 
potential. Prove that for any $\phi_0\in C^\infty M$,
there exists an LCK structure $(M, \theta', \omega')$
such that $\omega'= d_{\theta'} d^c_{\theta'} \phi_0$.

\item
Find an LCK manifold  $(M, \theta, \omega)$ 
with a $d_\theta d^c_\theta$-potential
$\phi_0$ such that the level set $\phi_0^{-1}(c)$ 
of $\phi_0$ is not strictly pseudoconvex for
some $c>0$.\footnote{By contrast, $\phi_0^{-1}(0)$ 
is {\em always} strictly pseudoconvex
(\ref{_level_set_for_d_theta_potential_Lemma_}).}

{\em Hint:} use the previous exercise.

\end{enumerate}


\chapter{Holomorphic $S^1$ actions on LCK manifolds}\index[terms]{manifold!LCK}\index[terms]{action!$S^1$-}

{\setlength\epigraphwidth{0.6\linewidth}
\epigraph{\it I am a man of cultivation; I have studied various remarkable books, but I cannot fathom the direction of my preferences; do I want to live or do I want to shoot myself, so to speak? But in order to be ready for all contingencies, I always carry a revolver in my pocket.}{\sc\scriptsize Anton Checkov, \ \ The Cherry Orchard}}

\section{Introduction}

In \ref{kami_or}, the Vaisman manifolds were\index[terms]{manifold!Vaisman}
characterized in terms of their automorphism groups.
It was shown that an LCK manifold is Vaisman if and only
if it admits an action of $\C$ by conformal holomorphic
automorphisms acting by non-isometric homotheties
on its K\"ahler cover. 

It turns out that compact LCK manifolds with potential\index[terms]{manifold!LCK!with potential}
admit a similar characterization in terms of 
automorphisms. Let $M$ be a compact LCK manifold  admitting
a holomorphic circle action\index[terms]{action!$S^1$-}. Assume that the Lee form\index[terms]{form!Lee}
is not exact on some orbits of this action.
Then $M$ admits the structure of LCK with potential.\index[terms]{structure!LCK!with potential}

In fact, $M$ admits an $S^1$-invariant LCK structure,
and the corresponding metric is automatically 
LCK with potential.

To prove this result, we start from an arbitrary
LCK metric and an $S^1$-action. First, we average
the Lee form with respect to the $S^1$-action, and
obtain an LCK structure\index[terms]{structure!LCK} in the same conformal class
and an $S^1$-invariant Lee form.\index[terms]{form!Lee}

Then we notice that a sum of two LCK metrics\index[terms]{metric!LCK} with the
same Lee form is again an LCK metric with the same
Lee form. This allows one to average an LCK metric
with $S^1$-action when the Lee form\index[terms]{form!Lee} is $S^1$-invariant.
In the end we obtain an $S^1$-invariant 
LCK structure.\index[terms]{structure!LCK}

It turns out that an $S^1$-invariant LCK 
structure has vanishing twisted Bott--Chern class\index[terms]{class!Bott--Chern!twisted}
whenever the $S^1$-action lifts to an action by non-isometric
homotheties on its K\"ahler cover. This implies the
existence of a $d_\theta d^c_\theta$-potential for the
metric (see \ref{_d_theta_d^c_theta-potential_Definition_}
for $d_\theta d^c_\theta$-potential). 

However, this potential is not {\em a priori}
positive; we need positivity for the notion of
LCK potential.\index[terms]{potential!LCK} Positivity would follow if we apply
\ref{_main_positive_Theorem_}. Alternatively, one could
prove directly that the $d_\theta d^c_\theta$-potential
obtained from the circle action\index[terms]{action!$S^1$-} is indeed 
positive, as \index[persons]{Istrati, N.} Istrati did in
\cite{nico2}.

\section{$S^1$-actions on compact LCK manifolds}\index[terms]{action!$S^1$-}\index[terms]{manifold!LCK}

The aim of this chapter is to prove the following:

\hfill

\theorem  \label{_S^1_potential_Theorem_}
{ (\cite{ov_imrn_12})} \label{s1action}
 Let $M$ be a compact complex manifold, equipped
with a holomorphic $S^1$-action and an LCK metric\index[terms]{metric!LCK}
(not necessarily $S^1$-invariant). Suppose that
the weight bundle $L$, restricted to a general\index[terms]{bundle!weight}
orbit of this $S^1$-action, is non-trivial
as a 1-dimensional local system. Then $M$ admits
an LCK metric with LCK potential.\index[terms]{potential!LCK}

\hfill

\remark 
The converse statement is also true, by \ref{_S^1_action_exists_Corollary_}.

\hfill

The rest of this chapter 
is devoted to the proof of \ref{s1action}, that is  separated in several lemmas.

\subsection{The averaging procedure}\label{_averaging_subsection_}

We first describe a general procedure, valid  on any compact LCK manifold.\index[terms]{manifold!LCK}

\hfill

\begin{lemma} { (\cite{ov_jgp_09})} 
 Let $M$ be a compact complex manifold, equipped
with a holomorphic $S^1$-action and an LCK metric
(not necessarily $S^1$-invariant). Then there exists
an LCK metric in the same conformal class with an
$S^1$-invariant Lee form.\index[terms]{form!Lee}
\end{lemma}

\hfill

\proof As $S^1$ is compact, averaging the Lee form\index[terms]{form!Lee} $\theta$ on $S^1$, we obtain a
closed $1$-form $\theta'$ that is  $S^1$-invariant and stays in the same
cohomology class as $\theta$: $\theta'=\theta+df$. Then
$\omega'=e^{-f}\omega$ is an LCK form\index[terms]{form!LCK} with Lee form $\theta'$ and conformal
to $\omega$.
\endproof

\hfill

\begin{lemma}\label{_LCK_averaging_Lemma_}
 { (\cite{ov_jgp_09})}
 Let $(M,\omega, \theta)$ be a compact complex manifold, equipped
with a holomorphic $S^1$-action and an LCK metric.\index[terms]{metric!LCK} {Then $M$ admits
an $S^1$-invariant LCK metric.}\index[terms]{action!$S^1$-}
\end{lemma}

\hfill

\proof  Using the previous lemma, we 
choose a metric in the same conformal class with $S^1$-invariant
Lee form.\index[terms]{form!Lee} 
Therefore, we may assume that $\theta$ is $S^1$-invariant.

Now, for each $t \in S^1$, let
$\omega_t:=\rho(t)^*\omega$, where $\rho$ denotes the action. Then 
$d(\omega_t) =\omega_t\wedge \theta$.
Averaging $\omega_t$ with respect to $t$,
{we obtain a positive, $S^1$-invariant form
$\omega_{av}$ satisfying $d(\omega_{av}) =\omega_{av}\wedge \theta$.}
\endproof

\hfill

As a result of the above lemmas, we can suppose from the
beginning that $S^1$ acts by holomorphic isometries with
respect to the fixed LCK metric\index[terms]{metric!LCK}. Such an action lifts to
the universal covering $\tilde M\arrow M$ in an action by
holomorphic conformalities with respect to the K\"ahler
metric of $\tilde M$. But, in real dimension greater than $2$,
conformalities of a symplectic form\index[terms]{form!symplectic}
are homotheties. We then have the following.

\hfill

\claim \label{claim}
 Let $(M,\omega, \theta)$ be a compact complex manifold, equipped
with a holomorphic $S^1$-action by isometries.\index[terms]{action!$S^1$-}
Then the covering $\R$ of $S^1$ acts on the K\"ahler cover $\tilde M \arrow M$
 by holomorphic homotheties.
 \endproof

\subsection{Holomorphic homotheties on a K\"ahler manifold}

The main technical step of the proof is the following local result.

\hfill

\proposition \label{prop}
 Let $A$ be a vector field acting on a
K\"ahler manifold $(\tilde M, \tilde \omega)$ by
holomorphic homotheties: $\Lie_A\tilde \omega
=\lambda\tilde\omega$, with $\lambda\neq 0$.
{Then
\[
dd^c|A|^2= \lambda^2\tilde \omega  +\Lie_{A^c}^2\tilde \omega,
\]
where $A^c = I(A)$.}

\hfill

\proof 
Let $\eta:= i_A\tilde\omega= I(A)^\flat$
be the dual form to $A^c$, where $i_v(\alpha)$ denotes the
contraction of a form $\alpha$ with a vector field $v$.
Replacing  $A$ by $\lambda^{-1}A$, we may
assume that $\lambda=1$. {By the Cartan formula,\index[terms]{Cartan formula}
\[
\tilde\omega=\Lie_A\tilde \omega  =d(i_A(\tilde\omega))=d\eta.
\]}
Since $A$ and $A^c$ are holomorphic,
$\Lie_{A^c}$ commutes with $I$. This gives
\[ \Lie_{A^c}\tilde\omega=
\Lie_{A^c} I\tilde\omega=  I \Lie_{A^c}\tilde\omega= I d (i_{A^c}\tilde\omega)=
- IdI^{-1}(i_A\tilde\omega)=-d^c \eta.
\]
 Since $\Lie_A$ commutes with $I$,
one has 
\begin{equation}\label{_Cartan_twisted_Equation_}
\{d^c, i_{A^c}\}= I\Lie_A I^{-1}=\Lie_A.
\end{equation}
The Cartan formula and $\Lie_{A^c}\tilde\omega=-d^c\eta$ give:
\begin{equation}\label{*} 
\Lie_{A^c}^2 \tilde\omega= -\Lie_{A^c} d^c\eta= 
-i_{A^c}dd^c \eta- d i_{A^c} d^c \eta 
\end{equation}
The first summand vanishes because  $dd^c \eta= -d^c d\eta= d^c \tilde\omega$. The second summand gives
\begin{equation}\label{**}
d i_{A^c} d^c \eta= dd^c \langle I(A), I(A)^\flat\rangle  -
d \{d^c, i_{A^c}\} \eta 
\end{equation}
Finally, \eqref{_Cartan_twisted_Equation_} gives
$$d \{d^c, i_{A^c}\} \eta= d\Lie_A\eta=\Lie_Ad\eta= \tilde\omega.$$
Therefore, \eqref{*} and \eqref{**} give
$\Lie_{A^c}^2 \tilde\omega=dd^c|A|^2- \tilde\omega$, for $\lambda=1$.
\endproof

\subsection{Vanishing of the twisted Bott--Chern class on
  manifolds endowed with an $S^1$-action}\index[terms]{class!Bott--Chern!twisted}

To finish the proof of \ref{s1action}, we may already assume that the LCK
metric\index[terms]{metric!LCK} on $M$ is $S^1$-invariant. Let $\tilde \omega$ be
the K\"ahler metric on the covering $\tilde M$. Denote 
the action of $S^1$ lifted to $\tilde M$ by $\tilde\rho$. From
\ref{claim}, $\tilde \rho$ acts on $(\tilde M, \tilde
\omega)$ by homotheties.

The restriction of the flat connection in the weight bundle\index[terms]{connection!flat}
$L$ to a loop has trivial monodromy\index[terms]{monodromy} whenever this\index[terms]{bundle!weight}
loop lifts to a homeomorphic loop in $\tilde M$.
Since $L$ is non-trivial on orbits of $\rho$,
the lift $\tilde \rho$ is an $\R$-action\index[terms]{action!$\R$-} by non-trivial homotheties. 
Rescaling, {we may assume that the vector field $A$
tangent to $\tilde\rho$ satisfies $\Lie_A \tilde\omega=\tilde\omega$.}

Now \ref{prop} gives:
\begin{equation}\label{***} 
\tilde \omega=  dd^c|A|^2 -\Lie_{A^c}^2\tilde \omega,
\end{equation}
where $A^c=I(A)$. Let $\mu_t:= \rho^c(t)^*[\omega]_{BC}$
be the Bott--Chern class\index[terms]{class!Bott--Chern} of $e^{tA^c}(\omega)$.
By \eqref{***}, $\mu_t$ satisfies the differential
equation $\mu_t''= -\mu_t$, and hence 
{$\mu_t = a\sin(t)+ b \cos(t)$,
for some $a, b \in H^{1,1}_{BC}(M,L)$.}

It then follows that $\int_0^{2\pi} e^{t A^c}[\tilde \omega]dt=0.$
Consider the K\"ahler form \index[terms]{form!K\"ahler}
$$\tilde \omega_0:=
\int_0^{2\pi} e^{t A^c}(\tilde \omega) dt$$
on $\tilde M$. This form is an average
of automorphic forms\index[terms]{form!automorphic} of the same character of automorphy,
because $e^{tA^c}$ commutes with $e^{t'A}$.
The Bott--Chern class of $\omega_0$ vanishes, because
$$\int_0^{2\pi}\sin(t) dt= \int_0^{2\pi}\cos(t) dt=0.$$ {Therefore, 
$\tilde \omega_0$ admits an automorphic potential.}

By \ref{_main_positive_Theorem_}, $M$ also admits an LCK metric with {\em positive} automorphic potential  which, by  \ref{defor_improper_to_proper} can be deformed to a proper one. The proof is complete. \index[terms]{metric!LCK!with potential}
\endproof

\hfill

\remark 
An alternative proof that directly produces a positive potential was recently presented in \cite{nico2}.

\section{Exercises} 

  Let $(M, \omega, \theta)$ be an LCS manifold and $L$ its weight bundle (see \ref{_weight_for_lcs_}). Denote with \index[terms]{manifold!LCS}\index[terms]{bundle!weight}
  $\sconf(M)$ the group of {\bf conformal symplectomorphisms} of $(M,\omega)$, i.\,e.      diffeomorphisms of $M$ that map $\omega$ into $e^f\omega$. Let $\mathfrak{sconf}(M)$ be its Lie algebra. It consists of vector fields $v$ whose flows $e^{tv}$ are formed by conformal symplectomorphisms.

\begin{enumerate}[label=\textbf{\thechapter.\arabic*}.,ref=\thechapter.\arabic{enumi}]
 
 \item Let $\mu\in\sconf(M)$. Prove that $\mu$ acts as a
   symplectic homothety on the symplectic covering
   $(\tilde M, \tilde\omega)$,
   $\mu^*\tilde\omega=\Theta(\mu)\tilde\omega$, where
   $\Theta:\sconf(M)\rightarrow\R^*$ is a group homomorphism.

 \item Prove that the weight bundle $L$ of an LCS manifold\index[terms]{manifold!LCS} is
 $\sconf(M)$-e\-qui\-va\-riant. Construct an $\sconf(M)$-e\-qui\-va\-riant symplectic structure\index[terms]{bundle!weight}
 on the space of non-zero vectors in $L\otimes\C$.
 
\item Suppose that an LCS manifold $M$ admits a vector
  field  $v\in \mathfrak{sconf}(M)$, with $\chi(v)\neq 0$,
  such that $e^{tv}$ induces a circle
  action\index[terms]{action!$S^1$-}. Prove that the
quotient $M/\langle e^{tv}\rangle$ is a contact orbifold.\index[terms]{orbifold!contact}

\item Let $(M,\omega,\theta)$ be an LCS manifold and $v\in
  \mathfrak{sconf}(M)$ inducing a circle action\index[terms]{action!$S^1$-}
  $\rho(t):=e^{tv}$. Prove that $\chi(v)=\int_S\theta$,
  where $S=S^1$ is an orbit of $\rho.$

\item Let $(M, \omega, \theta)$ be an LCK 
manifold\index[terms]{manifold!LCK}, $\pi:\tilde M\to M$ a K\"ahler cover on which 
$\pi^*\theta=d\phi$. Let $\rho(t)$ be a flow of
holomorphic isometries on $M$, $A$ its infinitesimal
generator. Denote by $\tilde \rho(t)$ the lift of
$\rho(t)$ to $\tilde M$, $\tilde A$ its infinitesimal
generator. 
\begin{enumerate}
	\item Show that $\tilde\rho(t)$ acts by isometries if and only if $\Lie_{\tilde A}\phi=0$.
	\item Show that if $A$ is not orthogonal to the
          Lee field\index[terms]{Lee field} $\theta^\sharp$, then $\tilde\rho(t)$
          does not act by isometries.
\end{enumerate}

\item
Let $(M, \theta, \omega)$ be a compact LCK manifold\index[terms]{manifold!LCK}
equipped with a circle actio\index[terms]{action!$S^1$-}n. Assume that
$\theta$ is not exact on some of the circle orbits.
Prove that $M$ admits
the structure of LCK with potential.

\item
Let $(M, \theta, \omega)$ be a compact homogeneous LCK manifold,\index[terms]{manifold!LCK!homogeneous}
that is, an LCK manifold with a transitive action of the
group of holomorphic isometries.
\begin{enumerate}
\item Prove that $M=G/G_1$, where $G$ is a compact Lie group
acting on $M$ by holomorphic isometries.
\item Given that $H=S^1\subset G$, prove that $H$ acts 
isometrically on the K\"ahler cover of $M$ if and only if
$\theta\restrict{H\cdot x}=0$ for all orbits 
$H\cdot x$ of $H$.
\item Prove that there exists an $S^1$-orbit such that
$\theta\restrict{S^1\cdot x}\neq 0$. Deduce that $M$ admits
the structure of LCK with potential.\index[terms]{structure!LCK!with potential}
\end{enumerate}

\item
Let $(M,\theta, \omega)$ be an LCK manifold\index[terms]{manifold!LCK}, and
$\tilde M$ its K\"ahler $\Z$-cover\index[terms]{cover!K\"ahler $\Z$-}. Assume that the
group of holomorphic automorphisms of $\tilde M$
is connected. Prove that $M$ admits the structure 
of LCK with potential.\index[terms]{group!holomorphic automorphisms}

\item
Let $(M,\theta, \omega)$ be an LCK manifold, and
$\tilde M$ its K\"ahler $\Z$-cover.\index[terms]{cover!K\"ahler $\Z$-} Assume that the
group of holomorphic conformal automorphisms of $\tilde M$
is connected. Prove that $M$ is Vaisman.\index[terms]{group!holomorphic conformal automorphisms}

\item
Produce an example of a compact LCK manifold\index[terms]{manifold!LCK}
such that its universal cover $\tilde M$ has a disconnected
 group of holomorphic automorphisms.

\item
Let $G$ be the connected component of the
group of holomorphic automorphisms
of a compact LCK manifold $M$, and $\tilde G$ its lift to
the minimal K\"ahler cover $\tilde M$. 
\begin{enumerate}
\item Prove that the kernel $K$ of the projection
$\tilde G \arrow G$ acts on $\tilde M$ by non-isometric
homotheties.
\item Prove that  $M$ admits the structure 
of LCK with potential if and only if $K$ is non-trivial.\index[terms]{structure!LCK!with potential}
\end{enumerate}

\item\label{_non-LCK+pot_isometry_on_cover_Exercise_}
Let $M$ be a compact LCK manifold that does not admit
an LCK metric with potential,\index[terms]{metric!LCK!with potential} and $X$ a holomorphic
Killing vector field on $M$. \index[terms]{vector field!holomorphic}\index[terms]{vector field!Killing}
\begin{enumerate}
\item Prove that $X$ is tangent to a compact group $G$ acting
on $M$ by holomorphic isometries.
\item Prove that the lift $\tilde X$ of $X$ acts on
the K\"ahler cover $\tilde M$ of $M$ by isometries.
\end{enumerate}
{\em Hint:} Suppose that $X$ acts on $\tilde M$
by non-isometric homotheties. Find a circle 
subgroup $S^1\subset G$ that acts on $\tilde M$
by non-isometric homotheties and apply 
\ref{_S^1_potential_Theorem_}.

\end{enumerate}


\chapter{Sasakian submanifolds in algebraic cones}\label{sasemb}

\epigraph{\it I would give my life for a man who is looking for the truth. But I would gladly kill a man who thinks that he has found the truth.}{\sc\scriptsize Luis Bu\~ nuel}
\section{Introduction}

Vaisman manifolds are closely\index[terms]{manifold!Vaisman}
related to Sasakian manifolds: a mapping torus of
a Sasakian manifold admits Vaisman structure (\ref{halfstr}), and every
Vaisman manifold contains a Sasakian pseudoconvex shell\index[terms]{pseudoconvex shell}
that determines its complex geometry 
(\ref{_Vaisman_via_logarithm_Theorem_}).

Traditionally, one would prove theorems in Vaisman
geometry using results about Sasakian manifolds. 
However, we can reverse this chain of arguments
to prove things about Sasakian manifolds using
Vaisman geometry.\index[terms]{geometry!Vaisman}

We start this chapter from a characterization
of Sasakian manifolds in terms of their CR-structures\index[terms]{structure!CR}.
Recall that a CR-manifold is a manifold
$M$ equipped with a sub-bundle $B\subset TM$
and a complex structure operator $I\in \End B, I^2=-\Id_B$,
such that $[B^{1,0}, B^{1,0}]\subset B^{1,0}$.
In many aspects, the CR geometry\index[terms]{geometry!CR} offers a good
substitute of complex geometry. This is
especially true about the strictly pseudoconvex CR-manifolds,
that is, the CR-manifolds for which $\codim B=1$ and the
Levi form\index[terms]{form!Levi} $[B, B] \arrow TM/B$ is positive definite.

Define CR-holomorphic functions as being 
functions $f$ such that $df$ vanishes on $B^{0,1}\subset TM\otimes \C$.\index[terms]{function!CR-holomorphic}
It is not hard to see that CR-holomorphic functions form
a ring; for a strictly pseudoconvex CR-manifold $S$ of dimension $>3$,
the ring of CR-holomorphic functions is identified with a ring
of holomorphic functions on a Stein variety $M$ (smooth or
with isolated singularities). Then $S$ is identified with the
boundary of $M$ (\ref{asr}).

Define a CR-holomorphic vector field\index[terms]{vector field!CR-holomorphic}
on a CR-manifold as a vector field with its diffeomorphism
flow preserving the CR-structure.\index[terms]{structure!CR} For strictly
pseudoconvex CR-manifolds, this is the same as a derivation
of the ring of CR-holomorphic functions; in particular,
a CR-holomorphic vector field acts on the associated
Stein variety.

It is easy to see that a hypersurface $S\subset M$ in a complex
manifold is equipped with a natural CR-structure,
$B= TS\cap I(TS)$. Sasakian manifolds are defined
in terms of the K\"ahler structure on their Riemannian
cones, and hence  they are realized as hypersurfaces in
complex manifolds. Therefore, every Sasakian manifold
carries a CR-structure; this CR-structure\index[terms]{structure!CR} is pseudoconvex,
because $S$ is the level set of a plurisubharmonic function
(\ref{_psh_then_pseudoconvex_Corollary_}).

We characterize CR-manifolds admitting Sasakian 
structure as follows: a stric\-tly pseudoconvex CR-manifold
$(S, B,I)$ is Sasakian if it admits a CR-holomorphic vector field
$\xi$ transversal to $B$ everywhere. Then $\xi$ or $-\xi$,
depending on orientation, is the Reeb field of this Sasakian
manifold (\ref{_Reeb_fie_from_CR_Theorem_}). 
The vector field determines the metric
uniquely, and the set of Sasakian structures
on a given CR-manifold is finite-dimensional.

It is known that the group of CR-automorphisms of a compact
strictly pseudoconvex CR-manifold $S$ is
compact unless $S$ is the standard sphere,
\cite{_Burns_Snider_,_Lee:CR_,_Schoen:CR_}. The Lie
algebra $\g$ of this group is the one which 
contains the Reeb fields, and the Reeb fields
of possible Sasakian structures form a 
convex cone in $\g$. 

This characterization of Sasakian manifolds
can be further used to apply results from 
complex geometry to Sasakian geometry.\index[terms]{geometry!Sasaki}
Imitating an argument of \index[persons]{Demailly, J.-P.} Demailly and
Demailly-\index[persons]{Paun, M.}Paun, we prove a K\"ahler embedding theorem,
showing that any K\"ahler metric on a 
submanifold $X\subset M$ can be extended
to a K\"ahler metric on $M$, if the possible 
cohomological obstructions vanish 
(\ref{_embed_Kahler_Theorem_}). 
This is used to prove the existence
of CR-holomorphic embeddings in
Sasakian geometry.\index[terms]{geometry!Sasaki} We prove that
every Sasakian manifold admits a
Sasakian embedding to a sphere 
equipped with an appropriate Sasakian
structure.\index[terms]{embedding!CR-holomorphic}

The dynamics of the Reeb  field action on a Sasakian
manifold can be understood in terms of the group
of CR-automorphisms, that is  either ${\rm PU}(1,n)$
(Exercise \ref{_ball_auto_Exercise_}) or compact.
Since the set of compact 1-parametric subgroups
of a compact Lie group is dense, for a
Sasakian manifold that is  not a sphere, the
Reeb fields acting with compact orbits are dense 
in the Reeb fields for all Sasakian structures.

For a sphere, the group of CR-automorphisms
${\rm PU}(1,n)$ is not compact; however, every Sasakian
structure is metrizable, and hence  any Reeb field
lies in a compact subgroup of ${\rm PU}(1,n)$,
and can also be approximated by a Reeb field
with compact orbits. This result implies, 
in particular, that every Sasakian manifold
is a limit of quasi-regular Sasakian manifolds
(\cite{_Rukimbira_deform_}, \cite{_Itoh_}).\index[terms]{manifold!Sasaki!quasi-regular}

The topology of Sasakian manifolds is significantly
different from that of general contact manifolds;
indeed, every Sasakian manifold is diffeomorphic
to an $S^1$-bundle over a K\"ahler orbifold,
which immediately follows from the quasi-regularity
(\ref{_quasireg_Sasakian_orbibundles_Theorem_}).
In real dimension 3, the classification of Sasakian
manifolds can be obtained from the classification
of Vaisman surfaces 
(\cite{bel,_Belgun:Sasakian1_,_Belgun:Sasakian2_}), or directly
(\cite{geiges}).
\index[persons]{Geiges, H.}\index[terms]{theorem!Geiges}
In dimension 5 or more, there is no classification,
but the topology of Sasakian manifolds can be deduced
from the topology of projective orbifolds using the
circle bundle.

\section{Sasakian structures on CR-manifolds}\index[terms]{manifold!Sasaki}

\remark Let $(S, B, I)$ be a 
strictly pseudoconvex CR-manifold. \index[terms]{manifold!CR!}
An orientation on $TS/B$ is called {\bf a co-orientation}
of a CR-manifold. The Levi
form $h$ (\ref{_Levi_form_Definition_}) 
defines the co-orientation as follows.
For every non-vanishing vector field $X\in B$, 
the Levi form\index[terms]{geometry!CR} applied to $X, I(X)$
gives a non-degenerate section
$h(X, IX)\in TS/B$, defining an
orientation on $TS/B$. Since
$h$ is positive definite, this orientation
is independent on  the choice of $X$.

\hfill

\definition \label{_pos_Reeb_}
Let $(S, B, I)$ be a 
strictly pseudoconvex CR-manifold.
A vector field $v\in TS$ is called {\bf positive}
if it is transversal to $B$ everywhere, and
its projection on $TS/B$ is positive.

\hfill

\example  \label{_Reeb_positive_holo_Example_}
The Reeb field $\xi$ of a Sasakian manifold $S$ is always positive 
(or negative, depending on the choice of orientation). 
Indeed, consider $S$ as a hypersurface in the
corresponding conical K\"ahler manifold $(M,I)$. Then
the radial vector field
$I\xi\in TM\restrict S$ is always normal to $S$, and hence 
$\xi\notin B=TS\cap I(TS)$. The Reeb field is also 
CR-holomorphic (\ref{_CR_holomorphic_Definition_}),
because the Reeb field acts on the cone of the Sasakian
manifold by holomorphic isometries (\ref{_Reeb_Sasakian_Theorem_}).

\hfill

The following theorem is stated for $\dim_\R S \geq 5$.
For $\dim_\R S =3$ it is also true, as 
follows from the classification
of compact 3-dimensional Sasakian  manifolds 
due to H. \index[persons]{Geiges, H.} Geiges (\cite{geiges}).
\index[persons]{Geiges, H.}\index[terms]{theorem!Geiges}

\hfill

\theorem  \label{_Reeb_fie_from_CR_Theorem_}
{ (\cite{ov_geom_ded,_BGS:cone_,_BGS:canonical_})}
 Let $S$ be a strictly pseudoconvex compact CR-manifold,
$\dim_\R S \geq 5$. {Then $S$ admits a Sasakian structure
if and only if $S$ admits a positive CR-holomorphic vector field.}\index[terms]{vector field!CR-holomorphic}
This vector field becomes the Reeb field of the Sasakian manifold.

\hfill

\pstep By \ref{_Reeb_positive_holo_Example_}, the Reeb
field on a Sasakian manifold is positive and CR-holomorphic.
Conversely, 
let $S$ be a strictly pseudoconvex compact CR-manifold,
$\dim_\R S \geq 5$ and let $v\in TS$ 
a positive CR-holomorphic vector field. 

By the Rossi and  \index[persons]{Andreotti, A.} Andreotti-\index[persons]{Siu, Y.-T.}Siu Theorem (\ref{asr}),
$S=\6 M$, where $M$\index[terms]{theorem!Rossi, Andreotti--Siu}
is a Stein variety with isolated singularities.
The variety $M$ can be realized as 
$M=\Spec(H^0(\calo_S))$, where $\Spec$ denotes the
continuous spectrum (that is, the set of continuous $\C$-linear \index[terms]{spectrum!continuous}
homomorphisms to $\C$)
of the commutative Banach algebra $H^0(\calo_S)$ of
CR-holomorphic functions on $S$. Since
$v$ is CR-holomorphic, the corresponding 
diffeomorphism flow acts on $M$ by 
holomorphic automorphisms.\index[terms]{function!CR-holomorphic}

Since $v$ is positive, $Iv$ is transversal
to $\6 M$. Replacing $v$ by $-v$ if needed, we can always assume that
$Iv$ points towards the interior of $M$. Then for 
small $\epsilon>0$ the holomorphic self-map 
$A_\epsilon:=e^{\epsilon I v}$
maps $M$ to a subset $A_\epsilon(M)\subset M$ 
inside the interior of $M$.

\hfill 

{\bf Step 2:}
 Consider now the ring ${\cal H}= H^0(\calo_M)_b$ of 
bounded holomorphic functions on the interior of $M$, with the sup-metric.
Then ${\cal H}$ is a Banach ring. Since $A_\epsilon(M)$
has compact closure in the interior of $M$,
the space $A_\epsilon^* {\cal H}$ is a normal
family, and $A_\epsilon^*$ is a compact operator\index[terms]{operator!compact}
(\ref{_contra_compact_Theorem_}).

By the strong maximum principle, \ref{hopf_theorem}, \index[terms]{maximum principle} for any
non-constant $f\in {\cal H}$, one has 
\begin{equation}\label{_contraction_maximum_Equation_}
\sup_{A_\epsilon(M)} |f|<\sup_M |f|.
\end{equation}
Indeed, if this inequality is not strict,
$|f|$ reaches its maximum on some point
in the closure of $A_\epsilon(M)$. This
is impossible, because $f$ is holomorphic
on an open subset containing $A_\epsilon$,
and a holomorphic function cannot have
maxima in the interior of its domain of definition.

\hfill

{\bf Step 3:}
Since the operator $A_\epsilon^*$ is compact,
the sequence $(A_\epsilon^k)^*(f)$ 
 has a limit point $f_{\lim}$.  This function  satisfies
$\sup_{A_\epsilon(M)} |f_{\lim}|=\sup_M |f_{\lim}|$;
applying \eqref{_contraction_maximum_Equation_},
we obtain that it  is constant. 

This implies that  the limit point $z_{\lim}$ 
of the sequence $\{A^k_\epsilon z\}_{k\geq 0}$ is unique and 
independent of $z\in M$. Indeed, $f_{\lim}(z)= f(z_{\lim})$, but $f_{\lim}=\const$.
This means that {$A_\epsilon$ is a holomorphic contraction
contracting $M$ to the origin point $x_0\in M$.}\index[terms]{contraction!holomorphic}

\hfill

{\bf Step 4:}
Let $C$ denote the union of $\Z^{>0}$ copies  $M_i$ of $M$, glued together according to the equivalence relation:
$$x\sim y \ \text{for}\, x\in M_i, \ y\in M_{i-1}\ \text{and}\ x=A_\epsilon(y).$$
\begin{figure}[htp]\centering{
\includegraphics[scale=0.43]{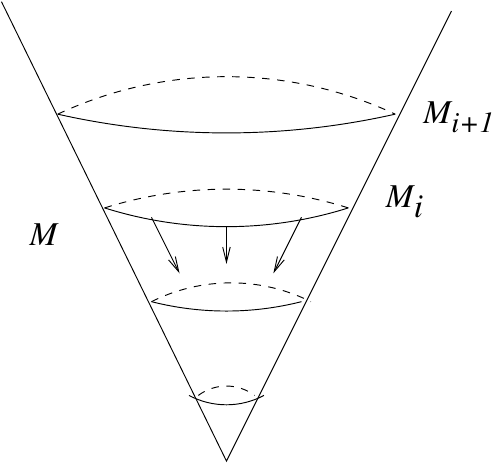}}
\end{figure}

Then $A_\epsilon$ is invertible on $C$ because $A_\epsilon$ moves $y\in M_i$ to $A_\epsilon(y)\in M_i$, and hence $y\in M_{i-1}$. This shows that $A_\epsilon^{-1}$ moves $s\in M_i$ to $x\in M_{i+1}$.

\hfill

{\bf Step 5:}
Since the vector field $v\restrict {S_\epsilon}=A_\epsilon(v)$
is nowhere vanishing for each $\epsilon$, the vector field 
$\vec r:= Iv$ is transversal to $S_\epsilon:= A_\epsilon(S)$ and  
points to the origin. Therefore, each point of $S$ 
is contained in the image of the map $\rho:\; \R \arrow M$, where $\rho$ 
is the solution  of the  equation 
$\frac{d\rho(t)}{dt}=\vec r$.

\hfill

{\bf Step 6:}
Let $\phi_\lambda$ be a $\rho$-automorphic
K\"ahler potential associated with the pseudoconvex shell\index[terms]{pseudoconvex shell} $S$ and 
the flow of diffeomorphisms $\rho$  (\ref{shell_char}).
The Lie algebra $\langle v, I v\rangle$ acts on $(M, dd^c \phi_\lambda)$
by holomorphic homotheties, and hence  $C$ is a conical K\"ahler manifold
(\ref{kami_or} and \ref{_Lee_field_holo_Proposition_}). 
Therefore, $S$ is Sasakian.
\endproof


\section{Isometric embeddings of K\"ahler and Vaisman ma\-nifolds}


The following general result will be used further on 
to construct isometric embeddings of Vaisman manifolds.
We are indebted to J.-P. \index[persons]{Demailly, J.-P.} Demailly who suggested the scheme\index[terms]{manifold!Vaisman}
of the proof.

\hfill

\theorem { (\cite{ov_sas})} \label{_embed_Kahler_Theorem_} 
 Let $(M, \omega)$ be a compact K\"ahler manifold, and $Z\subset M$ a  closed
complex submanifold. Denote by $[\omega]\in H^2(M)$ the K\"ahler class\index[terms]{class!K\"ahler} of $M$. Consider a
K\"ahler form\index[terms]{form!K\"ahler} $\omega_0$ on $Z$ such that $[\omega_0]=[\omega]\restrict{Z}$. 
Then there exists a K\"ahler form $\omega'$ on $M$, with $[\omega']=[\omega]$, and such that $\omega'\restrict{Z}=\omega_0$.

\hfill

\proof Let 
$$\omega_0=\omega\restrict{Z}+ dd^c(u),$$
for some function $u$ on $Z$. We extend $u$ to
$M$, add the square of the distance function $\delta(z)=d(z, Z)^2$ and smoothen the
result outside of the neighbourhood  $U$ of $Z$ where $\delta(z)$ is smooth.
This gives a smooth function $u$ on $M$
such that 
\begin{equation}\label{_bound_on_u_Equation_}
(1-\epsilon)\omega+dd^c(u)>0
\end{equation} 
on some neighbourhood
$U\supset Z$, for $\e<\frac 12$.

By \cite[Lemma 2.1]{dp}, there exists a function $\f$ with the following properties:
\begin{itemize}
\item[(i)] $\f$ is smooth outside of $Z$;
\item[(ii)] $dd^c(\f)>-C\omega$ for some $C<1$;
\item[(iii)] $\f$ has logarithmic poles on $Z$, that is,
locally in $Z$ one has
$$\phi= \const\cdot\log\left(\sum |g_i|^2\right),$$ where $g_i$ is
a collection of holomorphic functions vanishing on $Z$.
\end{itemize}
 Then $\f$ is very negative near $Z$. The regularized maximum (\ref{reg_max})\index[terms]{regularized maximum}
$$\psi:=\max{}_\delta\left(\frac{\e}{\delta}\f+A,u\right)$$
equals $u$ in some open neighbourhood of $Z$ for any constant $A\in\R$, and 
$$\omega':=\omega+dd^c(\psi)$$
satisfies $\omega'\restrict{Z}=\omega_0$. 

It remains to show that $\omega'$ is positive. First,
\eqref{_bound_on_u_Equation_} implies
$$dd^c(\psi)\restrict{U}\geq -(1-\e)\omega.$$
This immediately implies that $\omega'\restrict U$ is positive.
Now, with $A$ big enough, we may assume that $\psi=\f+A$
outside of $U$, giving
$$\omega'\restrict{M\setminus  U}=dd^c(\psi)+\omega=
dd^c(\phi)+\omega \geq (1-C) \omega.$$
This finishes the proof. \endproof

\hfill

We use this result to construct an isometric embedding
of a Vaisman manifold to a Hopf manifold.\index[terms]{manifold!Vaisman}

\hfill

\theorem\label{_Isometric_embedding_for Vaisman_Theorem_} 
Let $\Psi: M\hookrightarrow H$ be the
holomorphic embedding of a quasi-regular Vaisman manifold
to a quasi-regular Hopf manifold. Let $\tilde M$ be \index[terms]{manifold!Vaisman!quasi-regular}
a K\"ahler $\Z$-cover\index[terms]{cover!K\"ahler $\Z$-} of $M$, and $\phi$ an\index[terms]{embedding!holomorphic}
automorphic, positive K\"ahler potential on $\tilde M$.
We lift $\Psi$ to the holomorphic map of the
respective $\Z$-covers, $\tilde \Psi:\; \tilde M
\hookrightarrow \C^N\setminus 0$. 
 Then there exists a K\"ahler potential 
on $\C^N\setminus 0$ which extends the
potential $\f$ on $\tilde M$.

\hfill

\proof By the structure theorem for 
quasi-regular Vaisman manifolds \index[terms]{manifold!Vaisman}
(\ref{_Structure_of_quasi_regular_Vasman:Theorem_}),
$\tilde M$ can be obtained as the total space
$\Tot^\circ(L)$ of non-zero vectors in an ample line bundle 
over a projective orbifold $X$. 
 

Since all compact complex curves in $H$ are
leaves of the canonical foliation\index[terms]{foliation!canonical} (\ref{_Subva_Vaisman_Theorem_}),
the map $\Psi: M\hookrightarrow H$
maps the leaves of the canonical
foliation to the leaves of the canonical
foliation. On the other hand, under the
identification between $\tilde M$
and $\Tot^\circ(L)$, the leaves of the
canonical foliation correspond to the
fibres of $L$; a similar identification
is valid for $H$ as well.

Since $\Psi: M\hookrightarrow H$ is compatible with the
canonical foliation, it induces the holomorphic map 
 $f:X\hookrightarrow \C P^{N-1}(i_1, i_2, ..., i_N)$ of the
corresponding leaf spaces, where
$\C P^{N-1}(i_1, i_2, ..., i_N)$ denotes the 
weighted projective space.

Since $\Pic(\C P^{N-1}(i_1, i_2, ..., i_N))=\Z$ is generated\index[terms]{group!Picard}
by the tautological bundle \index[terms]{bundle!vector bundle!tautological} $\calo(1)$, the
weight bundle $L$ on $M$ satisfies $L=f^*\calo(d)$.
\index[terms]{weighted projective space}
The Vaisman metrics on $M$ and $H$ are obtained from
the Hermitian metric on the corresponding line 
bundles
(\ref{_Structure_of_quasi_regular_Vasman:Theorem_}).
Therefore, $\Psi$ is an isometric embedding if the
natural map $L\arrow f^*(\calo(d))$ of line bundles
is an isometry.

Using \ref{_embed_Kahler_Theorem_}, we may
choose the K\"ahler metric on $X$ in such a way that 
$f:X\hookrightarrow \C P^{N-1}(i_1, i_2, ..., i_N)$
is an isometry. Moreover, we may also assume
that the associated Hermitian form\index[terms]{form!Hermitian} is the curvature of the
line bundle $L$, and the Hermitian form on 
$\C P^{N-1}(i_1, i_2, ..., i_N)$ is the curvature of
$\calo(d)$. Since the curvature of the  Chern connection\index[terms]{connection!Chern}
on a line bundle determines its Hermitian structure
up to a constant multiplier, this implies, after adjusting
a constant, that the natural map
$L\arrow f^*(\calo(d))$ is an isometry.
Applying
\ref{_Structure_of_quasi_regular_Vasman:Theorem_},
we obtain that $\Psi: M\hookrightarrow H$
is also an isometry.
\endproof


\section{Embedding Sasakian manifolds in spheres}\index[terms]{manifold!Sasaki}


\subsection{Kodaira-like embedding for Sasakian manifolds}

 In this section, we construct a CR embedding  for  Sasakian manifolds of
 dimension greater than 5.  This result parallels the
 Kodaira embedding for K\"ahler manifolds and is optimal
 as shown in Subsection \ref{opti}
 below. 

\hfill

\theorem \label{_Sasa_embedding_Theorem_}
{ (\cite{ov_sas})}\label{emb_sas} Any compact Sasakian manifold of dimension 
$2n-1\geq 5$  admits a CR-embedding into a
Sasakian manifold $\mathbb{S}$ that is  diffeomorphic to a sphere, and such that the corresponding conical 
K\"ahler manifold $\mathbb{S} \times \R^{>0}$ is biholomorphic to $\C^N\setminus 0$. 

\hfill

\proof Let $S$ be a compact Sasakian manifold. Then  $M =
S \times S^1$  is  Vaisman, and hence it holomorphically
embeds in a diagonal Hopf manifold $H$\index[terms]{embedding!holomorphic}
(\ref{vaisman_embed}). Let $\Psi:M\hookrightarrow H$ be
this embedding and denote with $\widetilde \Psi$ the
natural lift of $\Psi$ to the universal covers,
$\widetilde\Psi:\tilde M\hookrightarrow \C^N\setminus 0$.

The manifold $\tilde M$ is a cone $S\times \R^{>0}$
endowed with a global automorphic potential  such that 
$S$ is its 1-level set. By 
\ref{_Isometric_embedding_for Vaisman_Theorem_}, we can extend the
automorphic K\"ahler potential on $S\times \R^{>0}$
to an automorphic K\"ahler potential on $\C^N\setminus 0$.
To prove \ref{_Sasa_embedding_Theorem_}
it remains to show that the corresponding 
level sets are always diffeomorphic to the
spheres.

\hfill

\claim \label{_level_sets_spheres_Claim_}
Let $\f$ be a positive K\"ahler potential on $V = \C^N\setminus 0$, satisfying
$X_A(\f) = 2\f$ for a vector field $X_A$ on $V$, $(X_A)\restrict v=A(v)$, 
 where $A$ is a linear operator
on $V$. Then the level sets $\f^{-1}(c)$ are  diffeomorphic to a sphere, for all  $c\in \R^{>0}$.

\hfill

\proof Since $X_A(\f) = 2\f$, $e^A$ acts on $V$ as a homothety, with respect to the
K\"ahler metric defined by
$dd^c\f$. Therefore, $e^{-A}$ is a contraction.

Let $\cax_A$ be the 1-dimensional foliation   generated on
$V$ by $X_A$, and  
\[ r(t, x):= e^{tA}\cdot x, \text{\ for \ } x\in V,
t\in\R.
\] Also, let $S_x:=\im(t\mapsto r(t,x))\subset V$.
Clearly, the sets $S_x$ are leaves of $\cax_A$. The leaf
space $V/\cax_A$ is
naturally isomorphic to $S^{2N-1}$,  because $e^{-A}$ is a
contraction.

On each leaf $S_x$, the function $\f$ becomes an
exponential after a coordinate change. Therefore,
$S_x$ intersects with $\f^{-1}(c)$ precisely once. This
defines the diffeomorphism between $\f^{-1}(c)$ and
$V/\cax_A\simeq S^{2N-1}$. We finished the proof of
\ref{_level_sets_spheres_Claim_} and \ref{_Sasa_embedding_Theorem_}.
\endproof

\subsection{Optimality of the embedding result}\label{opti}

Consider a CR-holomorphic immersion of CR-manifolds $X \hookrightarrow Y$, with $Y$ Sa\-sa\-kian,
and $X$ preserved by the Reeb field of $Y$. This Reeb action is positive in the sense of \ref{_pos_Reeb_} and hence defines a Sasakian structure on $X$.  This is analogous to the fact that the 
restriction of a K\"ahler metric to a complex subvariety
is again a K\"ahler metric. We claim that such a map is
compatible with a whole set of Sasakian structures
on $X$ considered to be  a CR-manifold. Moreover,\index[terms]{manifold!Sasaki!quasi-regular}
some of these Sasakian structures are quasi-regular.

\hfill

\remark\label{qr_reeb} Let $S_1, S_2$ be compact Sasakian manifolds and suppose there exists
a CR-embedding $S_1\hookrightarrow S_2$ which commutes with the flows of the respective Reeb
fields $\xi_1, \xi_2$. The group $G$ of Sasakian automorphisms of $S_2$ that preserve $S_1$ is, by
assumption, non-empty, as it contains the flow of $\xi_2$. It thus   contains a whole torus $T$ in which
one can deform the initial Reeb flow to a quasi-regular one (it is the same method employed to deform a Vaisman structure to a quasi-regular one, \ref{defovai}). Hence, there exists a quasi-regular Sasakian structure on $S_2$ whose Reeb
flow preserves  $S_1$.

\hfill

Now restrict this quasi-regular vector field
$\xi_2$ to $S_1$. It remains CR-ho\-lo\-mor\-phic and positive,
and hence it is still a Reeb field and defines a new
Sasakian structure on $S_1$, that is  by construction
quasi-regular. Then the above CR-embedding passes to the
K\"ahler orbifold quotients $M_i:=S_i/\xi_i$ and produces
a K\"ahler (i.\,e.     holomorphic and {\em isometric}) embedding\index[terms]{embedding!holomorphic}\index[terms]{embedding!isometric}
$M_1\hookrightarrow M_2$.

Suppose that in the embedding theorem we could have a
precise model space. Denote it by $S_2$. Then the
quotient K\"ahler orbifold
$M_2:=S_2/\xi$ is determined
by a choice of a quasi-regular Reeb field $\xi$. Such a
choice  corresponds 
to a Lie group homomorphism
$S^1\hookrightarrow T$, where $T=\overline{\langle
  e^{t\xi_2}\rangle}$ is the closure of the group
generated by the Reeb field $\xi_2$.  
 But there only exists a countable number of
Lie group homomorphisms of the Lie group $S^1$ into a compact torus, and
$S_1$  is an $S^1$-bundle over an arbitrary projective orbifold. One cannot\index[terms]{orbifold!projective}
expect to have an isometric embedding of an arbitrary K\"ahler orbifold into a model
K\"ahler orbifold, or a countable number of such orbifolds: obviously, the model space
must have continuous (and even infinite-dimensional)
moduli. This contradiction shows that \ref{emb_sas} is the
best that one can aim at\footnote{Unlike the Harvard people, we do
  end sentences by prepositions.}. 
 
\section{Notes}

\begin{enumerate}
	
	\item A more general isometric embedding result for LCK manifolds with potential\index[terms]{manifold!LCK!with potential} was obtained by \index[persons]{Angella, D.} Angella and Zedda on the model of K\"ahler geometry\index[terms]{geometry!K\"ahler}, where isometric immersions in $\C
	P^n$ are studied using \index[persons]{Calabi, E.} Calabi's
	diastasis\index[terms]{diastasis} function. This method was
	extended to LCK geometry in \cite{az}, allowing the
	ambient space to be infinite-dimensional.  They arrived at the
	following partial result: {\em Let $(M,\omega)$ be an LCK
		manifold with proper potential.\index[terms]{manifold!LCK!with potential} Assume that the metric
		completion of the K\"ahler cover is
                smooth, that is, biholomorphic to $\C^n$
                (Exercise \ref{_attr_basin_Exercise_}). Then
		$(M,\omega)$ admits a holomorphic and isometric
		immersion into a Hopf manifold if the
		following conditions are fulfilled:
		(1) a K\"ahler cover of $(M,\omega)$ is induced by an immersion into $l^2(\C)$;
		(2) the diastasis centred at the origin of the metric completion of  the cover is an automorphic
		potential.}\index[terms]{diastasis}
	
	\item  \ref{_Sasa_embedding_Theorem_} was reproven
	with different arguments by \index[persons]{van Coevering, C.} van Coevering, entirely within the
	framework of Sasakian geometry\index[terms]{geometry!Sasaki}, in \cite[Theorem
	3.1]{coev}, also identifying the corresponding weighted
	Sasakian structure on the ambient sphere.\index[terms]{weighted sphere}
	
	\item A still more general extension was obtained in
	\cite[Theorem 1 (i)]{torres} by Mart\'{\i}nez Torres, for
	compact contact manifolds (with globally defined contact
	form), that are  shown to admit approximate holomorphic
	embeddings in contact spheres. 
	
	\item  In \cite[Theorem 5.2]{my},
	\index[persons]{Marinescu, G.} Marinescu and Yeganefar prove a CR embedding of compact
	Sasakian manifolds to $\C^N$
	viewed as CR-manifold. Their methods are of
	pseudoconvex geometry.
	
\end{enumerate}

\section{Exercises}

\begin{enumerate}[label=\textbf{\thechapter.\arabic*}.,ref=\thechapter.\arabic{enumi}]

\item
Let $S$ be a compact strictly pseudoconvex CR-manifold, and $V$
the space of all CR-holomorphic vector fields.
Prove that $V$ is finite-dimensional.\index[terms]{vector field!CR-holomorphic}

\item
Let $(S, B, I)$ be a compact CR-manifold, and ${\goth S}$
the set of all Sasakian structures inducing the CR
structure on $S$. Prove that ${\goth S}$
is finite-dimensional.

{\em Hint:} Use the previous exercise.

\item
Let $U$ be a smooth complex manifold
equipped with a holomorphic contraction.\index[terms]{contraction!holomorphic}
Prove that $U$ is biholomorphic to $\C^n$.

\item Let $M$ be a Vaisman manifold,\index[terms]{manifold!Vaisman}
$\tilde M$ its K\"ahler $\Z$-cover,\index[terms]{cover!K\"ahler $\Z$-} and $\phi$
the automorphic potential. Prove that
all level sets of $\phi$ are diffeomorphic.

\item
Let $M$ be an LCK manifold with potential,\index[terms]{manifold!LCK!with potential}
$\tilde M$ its K\"ahler $\Z$-cover\index[terms]{cover!K\"ahler $\Z$-}, and $\phi$
the automorphic potential. Suppose that
$M$ admits an automorphic pluriharmonic function $\psi$.
Prove that there exists a linear combination $\phi + \lambda \psi$
that has a negative minimum somewhere on $\tilde M$.

\item
Let $M$ be an LCK manifold with potential,
$\tilde M$ its K\"ahler $\Z$-cover, and $\phi$
the automorphic potential. Adding a pluriharmonic
automorphic function to $\phi$ as in \index[terms]{function!pluriharmonic}
\ref{_non_posi_potential_Victors_Example_}, 
find an example of $M$ such that $\phi$ can 
be modified to acquire singularities, in
particular, some Morse-type singularities.

\item
Find an example of LCK manifold with
potential such that not all smooth level
sets of the potential are diffeomorphic.

{\em Hint:} Use the previous exercise.

\item
Let $M, M_1$ be Vaisman manifolds\index[terms]{manifold!Vaisman}
$\tilde M, \tilde M_1$ their $\Z$-covers\index[terms]{cover!K\"ahler $\Z$-},
and $\phi, \phi_1$ the corresponding automorphic potentials on
$\tilde M, \tilde M_1$. Denote by $S, S_1$ the 
level sets of $\phi, \phi_1$, considered to be  Sasakian manifolds.
Assume that there exists a CR embedding 
$S \hookrightarrow S_1$ mapping the Reeb field of $S$
to the Reeb field of $S^1$. Prove that it
can be extended to a holomorphic immersion
$M \arrow M_1$.

\item
Let $j:\; S \hookrightarrow S_1$\index[terms]{embedding!CR-holomorphic}
be a CR embedding of Sasakian manifolds
mapping the Reeb field of $S$
to the Reeb field of $S_1$. 
\begin{enumerate}
\item
Prove that $j$
induces a holomorphic isometric embedding 
$C(S)\hookrightarrow C(S_1)$ of
the corresponding conical K\"ahler manifolds.

\item Prove that $j$ is an isometric immersion.
\end{enumerate}

\item\label{_ball_auto_Exercise_}
Let $B\subset \C^n$ be an open ball.
\begin{enumerate}
\item Consider $B$ as a subset of $\C P^n\supset \C^n$.
Prove that there exists a pseudo-Hermitian metric of 
signature $(1,n)$ on $\C^{n+1}$ such that $B = {\mathbb   P}\Pos$,
that is, the projectivization of the positive cone $\Pos
\subset \C^{n+1}$.

\item 
Denote by $\PU(1,n)$ the quotient of the group 
  $\U(1,n)$ of pseudo-Her\-mi\-tian isometries by its center.
Prove that the group $\PU(1,n)$ acts on the ball
$B\subset \C^n$ by holomorphic automorphisms.
Prove that this action is transitive.

\item
Let $A\in \Aut(B)$ be a holomorphic automorphism of
a ball preserving the origin. Using Schwarz
lemma, prove that $A$ is linear.\index[terms]{lemma!Schwarz}

\item Prove that $\PU(1,n)= \Aut(B)$, that is,
$\PU(1,n)$ is the group of holomorphic automorphisms of a
  ball.
\end{enumerate}

\item
Consider the sphere $S^{2n-1}\subset \C^n$ with the standard CR-structure.\index[terms]{structure!CR}
Prove that the group of CR-automorphisms of $S^{2n-1}$
does not act on $S^{2n-1}$ by isometries.

{\em Hint:} Use the previous exercise.

\item
Let $S$ be a CR-manifold, and 
$g, g_1$ two Sasakian metrics compatible
with the CR-structure. Prove that 
$g$ is not necessarily proportional to $g_1$.

{\em Hint:} Use the previous exercise.

\item Prove that every real vector bundle $B$
over a manifold $M$ is a direct summand of a 
trivial vector bundle.

{\em Hint:} When $B$ is a direct summand of a tangent
bundle, this follows from Whitney's embedding theorem.
For general $B$, replace $M$ by the sphere bundle 
$M_1 \subset B \oplus C^\infty M$, and show that 
$B$ is a direct summand of $TM_1 \restrict M$.

\item Prove that every quaternionic vector bundle $B$
over a manifold $M$ is a direct summand of a 
trivial quaternionic vector bundle.

{\em Hint:} Let $B_0$ be $B$ considered to be  a real
vector bundle. Show that $B$ is a direct summand
of $B_0 \otimes_\R {\Bbb H}$, and apply the
previous exercise to $B_0$.

\item Let $B$ be a rank $p$ quaternionic vector bundle
over a manifold $M$. Prove that there exists
a quaternionic Grassmannian space $\Gr(p, n, {\Bbb H})$
and a map $M \arrow \Gr(p, n, {\Bbb H})$ such that $B$
is isomorphic to the pullback of the tautological bundle.

{\em Hint:} Use the previous exercise.

\item
Let $P\arrow M$ be a principal $\SU(2)$-bundle.\index[terms]{bundle!principal}
Prove that there exists a map $M\arrow {\Bbb H} P^\infty$
such that $P = M\times_{{\Bbb H} P^\infty} S^\infty$.

{\em Hint:} Use the previous exercise.

\item \label{_S^3_oer_S^2_Exercise_}
Prove that every principal $\SU(2)$-bundle
on $S^2$ is trivial.

{\em Hint:} Use the previous exercise.

\item
Let $L=\calo(1)$ be the standard ample bundle
on $\C P^1$, and $S$ the total space of the
corresponding $S^1$-bundle. Prove that 
$S$ is equipped with a canonical free action
of $\SU(2)$.

\item
Let $L=\calo(1,1)$ be the standard ample line bundle on 
$Q:= \C P^1 \times \C P^1$, and
$Z$ the total space of the corresponding
$S^1$-bundle. 
\begin{enumerate}
\item Prove that the projection $Z\arrow \C P^1$
is a principal $\SU(2)$-fibration.

\item
Prove that $Z$ is diffeomorphic to $S^2\times S^3$.
\end{enumerate}

{\em Hint:} Use the previous exercise and
Exercise \ref{_S^3_oer_S^2_Exercise_}.

\item
Construct a Sasakian embedding $S^5 \hookrightarrow S^7$,
with the standard Sasakian structures. Construct a
Sasakian structure on $S^2\times S^3$ and find
a Sasakian embedding $S^2\times S^3 \hookrightarrow S^7$.

\end{enumerate}


\chapter{Oeljeklaus--Toma manifolds}\label{OT_manifolds}

{\setlength\epigraphwidth{0.9\linewidth}
\epigraph{\it Appliquer le m\^eme traitement \`a un po\`ete et \`a un penseur me semble une faute de go\^ut. Il est des domaines auxquels les philosophes
	ne devraient pas toucher. 
	
	D\'esarticuler un po\`eme comme or d\'esarticule un syst\`eme est un d\'elit, voire un sacril\`ege.
	 
	Chose curieuse : les po\`etes exultent quand ils ne comprennent pas ce qu’on d\'ebite sur eux. Le jargon les flatte, et leur donne
	l’illusion d'un avancement. Cette faiblesse les rabaisse au niveau de leurs glossateurs. }{\sc\scriptsize Emil Cioran, \ \ Aveux et anath\`emes}}
\section{Introduction}

\subsection{Many species of Inoue surfaces}

The  Oeljeklaus-\index[terms]{manifold!Oeljeklaus--Toma (OT)}Toma manifolds, also known as OT-manifolds,
are multi-dimen\-si\-onal generalizations of the \index[persons]{Inoue, Ma.} Inoue surfaces of type $S_M$,
that are  a special class of surfaces constructed by Masahisa Inoue\index[terms]{surface!Inoue}
in his first papers (\cite{_Inoue:announce_,inoue}). To distinguish
these Inoue surfaces from the complex surfaces discovered by \index[persons]{Inoue, Ma.} Inoue later,
his earliest surfaces are sometimes called {\bf the Inoue--Bombieri} surfaces.\index[terms]{surface!Inoue--Bombieri}
They are named so because \index[persons]{Bombieri, E.} Bombieri mentioned a construction
similar to Inoue in his private communication to Inoue, \cite{inoue}.
Since the early 1970s, four new kinds of surfaces
were constructed, also known as {\em \index[terms]{surface!Inoue} Inoue surfaces}. The surfaces which
were later called ``parabolic Inoue surfaces''
were introduced  in \cite{_Inoue:Proc_}. The
``hyperbolic Inoue surfaces'' and ``half Inoue surfaces'' are
collectively known  as
``the Inoue--\index[persons]{Hirzebruch, F.}Hirzebruch surfaces''. They are described
in \cite{_Inoue:Surfaces_II_}, with the notion
traced back to \index[persons]{Hirzebruch, F.} Hirzebruch's paper \cite{_Hirzebruch:Hilbert_};
the terminology (half, parabolic and hyperbolic Inoue) is due to
I. \index[persons]{Nakamura, I.} Nakamura, \cite{_Nakamura:curves_}.\index[terms]{surface!Inoue!half}

The three kinds of Inoue surfaces introduced by 
\index[persons]{Nakamura, I.} Nakamura \cite{_Nakamura:curves_}
are only tangentially
related to the \index[persons]{Inoue, Ma.} Inoue surfaces defined in Inoue's early work.\index[terms]{surface!Inoue}
In \cite{_Nakamura:curves_}, these three classes are defined as follows.

Recall that the class VII${}_0$ surfaces are minimal surfaces with  Kodaira\index[terms]{dimension!Kodaira}
dimension $-\infty$ and $b_1=1$. 
The parabolic Inoue surfaces are class VII${}_0$ surfaces
with an elliptic curve and at least one cycle of rational curves.
The hyperbolic Inoue surfaces are class VII${}_0$ surfaces
with two cycles of rational curves. The half Inoue surfaces 
are class VII${}_0$ surfaces obtained as $\Z/2$-quotients
of the hyperbolic\index[terms]{surface!Inoue} Inoue surfaces by the $\Z/2$-action
that exchanges the two cycles of rational curves.

These surfaces were constructed by Inoue by taking a quotient
of a germ of $\C^2$ by a rational automorphism, gluing in the cusp and 
resolving the singularities; this construction was later generalized
by G. \index[persons]{Dloussky, G.} Dloussky \cite{_Dloussky:Kato_},
who proved that all such constructions lead to \index[persons]{Kato, Ma.} Kato surfaces
(Chapter \ref{comp_surf}).

Another construction of parabolic Inoue surfaces was given by 
I. \index[persons]{Enoki, I.} Enoki (\cite{_Enoki_,_Enoki:compactifiable_}),
and used by M. \index[persons]{Brunella, M.} Brunella to construct the LCK metrics on these
surfaces (\cite{_Brunella:Enoki_}). The Enoki surfaces are obtained
by taking a flat affine bundle with infinite monodromy \index[terms]{monodromy}over an
elliptic curve, and compactifying it with a cycle of rational
curves. When this flat affine bundle has a holomorphic
section, we obtain a surface with an elliptic curve
and a cycle of rational curves, that is, the parabolic
\index[terms]{surface!Inoue} Inoue surface. Later, \index[persons]{Brunella, M.} Brunella generalized his construction
of LCK metrics\index[terms]{metric!LCK} to all Kato surfaces\index[terms]{surface!Enoki} (\cite{_Brunella:Kato_}).\index[terms]{surface!Kato}

The fifth kind of \index[terms]{surface!Inoue} Inoue surfaces, not related to the
other four in any way, is the general type projective
surface constructed by Inoue in \cite{_Inoue:general_type_}.
This construction was generalized by \index[persons]{Bauer, I.} Bauer and \index[persons]{Catanese, F.} Catanese
to arbitrary dimension (\cite{_BC:Inoue_manifolds_}).
Bauer and \index[persons]{Catanese, F.} Catanese have proved in \cite{_BC:Burniat_I_}
that the surfaces constructed\index[terms]{surface!Burniat}
by \index[persons]{Inoue, Ma.} Inoue are precisely the \index[persons]{Burniat, P.} Burniat surfaces,
defined by Burniat in \cite{_Burniat_} 
(see also \cite{_Peters:Burniat_}).

\subsection{Class VII${}_0$ surfaces with $b_2=0$}
 
The earliest of Inoue surfaces (the ``Inoue--Bombieri  surfaces'')\index[terms]{surface!Inoue--Bombieri}
were constructed explicitly by \index[persons]{Inoue, Ma.} Inoue in \cite{inoue}; very soon,
F. \index[persons]{Bogomolov, F. A.} Bogomolov showed that any class VII${}_0$ surface with $b_2=0$
is of this type. \index[persons]{Inoue, Ma.} Inoue obtained his surfaces as quotients of
$\C \times {\Bbb H}$ by a cocompact group of affine transformations
(here ${\Bbb H}$ denotes the upper half-plane, considered to be 
an open complex submanifold in $\C$). \index[persons]{Bogomolov, F. A.} Bogomolov has shown that
any class VII${}_0$ surface with $b_2=0$ without curves is equipped with
a flat torsion-free connection, \index[terms]{connection!torsion-free}\index[terms]{connection!flat}and classified surfaces
admitting a flat torsion-free connection using the Galois
theory\index[terms]{Galois theory} (\cite{_Bogomolov:VII_76_}). The last step of Bogomolov's argument
was written up in a rather cursory manner.
Later on, \index[persons]{Teleman, A. } Teleman and Li-\index[persons]{Yau, S.-T.}Yau--Zhang gave another proof of this
part of \index[persons]{Bogomolov, F. A.} Bogomolov's argument \cite{_Li_Yau_Zang:VII_,_Teleman:bogomolov_},
using the  Donaldson--Uhlenbeck--Yau theorem about the existence
of the instanton metrics on stable bundles.\index[terms]{theorem!Donaldson--Uhlenbeck--Yau} 
Also, \index[persons]{Bogomolov, F. A.} Bogomolov has rewritten his paper (\cite{_Bogomolov:VII_82_})
to make it more readable.

Inoue has constructed the \index[persons]{Inoue, Ma.} Inoue surface\index[terms]{surface!Inoue} explicitly in terms
of linear algebra; since the \index[persons]{Oeljeklaus, K.} Oeljeklaus-\index[persons]{Toma, M.}Toma paper appeared, his construction
has always been interpreted in terms of number fields.

We start by giving a number-theory-free  account of the Inoue construction.
Let $T^3$ be a 3-dimensional torus, and $A$ an infinite order
automorphism of $T^3$ preserving the flat connection and the orientation.
Flat oriented automorphisms of the torus are in one-to-one
correspondence with elements $A\in \SL(3, \Z)$.
A characteristic polynomial $P_A(t)$ of $\SL(3, \Z)$ 
might have 3 real roots or one real root; \index[persons]{Inoue, Ma.} Inoue 
considers the case when $P_A(t)$ has only one real root,
and takes a mapping torus $M$ of $T^3$ over $S^1$
associated with $A$. Then he puts a complex structure
on $M$, using an explicit expression for the
complex structure.

The modern version of this construction uses the
number theory. Consider a degree 3 number field $K= \Q[\zeta]$
obtained from $\Q$ by adding the real root $\zeta$ of $P_A(t)$,
and let $\calo_K\subset K$ be the integral closure of $\Z$ in $K$.
It is not hard to see that $\calo_K$ is a lattice\index[terms]{lattice} in $K$,
that is, the additive group of $\calo_K$ is isomorphic to $\Z^3$.
Let $\tau:\; K \arrow \R$ denote the field embedding
mapping $\zeta$ to the corresponding real root, and $\sigma:\; K \arrow \C$
one of the two complex embeddings. 
By the Dirichlet unit theorem, the group $\calo_K^*$ of units in
$\calo_K$ is $\Z\times \{\pm 1\}$, with the negative
part $-\Z$ mapped to $\R^{< 0}$ under $\tau$, and
$\Z$ mapped to $\R^{>0}$. Let $U\subset \calo_K^*$ 
be the group $\calo_K^*\cap \tau^{-1}(\R^{>0})$, isomorphic to $\Z$.

Take $\zeta$ as a generator of $U$, and
consider the semidirect product $\calo_K\rtimes U$
associated with the multiplicative action of $U$ on
$\calo_K$. The action of $A:=\tau(\zeta)$
on $\calo_K$ belongs to $\SL(\calo_K)= \SL(3,\Z)$.
The corresponding characteristic polynomial $P_A(t)$
is equal to the minimal polynomial of $\zeta$,
and it is a cubic polynomial with integer
coefficients and one real root. This construction
gives an element of $\SL(3,\Z)$ from a 
cubic number field. The inverse construction
is also straightforward; any $A \in \SL(3,\Z)$
defines an algebra $\Q[A] \subset \Mat(3, \Q)$,
that is  identified with the number field
$\Q[\zeta]$. Since $A$ is invertible in $\Z[A]$,
any $A$ comes from an element of $U$, and hence 
we have a correspondence between 
$A\in \SL(3, \Z)$ and the units in the
ring of cubic algebraic integers (\ref{_Units_SL(n,Z)_Remark_}).

Clearly, $\pi_1(M) = \calo_K\rtimes U$,
where $M$ is the Inoue surface constructed above as a mapping torus.
Indeed, by \index[persons]{Serre, J.-P.} Serre, one has the following exact sequence of fundamental groups\index[terms]{fundamental group}
\begin{equation}\label{_mapping_torus_Inoue_Equation_}
0 \arrow \Z^3= \pi_1(T^3)\arrow \pi_1(M) \arrow \Z= \pi_1(S^1)\arrow 0.
\end{equation}
This exact sequence has a section, because any surjective
map to a free group has a section. Therefore, $\pi_1(M) = \pi_1(T^3)\rtimes \pi_1(S^1)$

This group is isomorphic to $\calo_K\rtimes U$,
where $K=\Q[A]$ is the number field associated with $A$.
Indeed, the mapping torus $M$, considered to be  a fibration
over $S^1$, is equipped with the flat Ehresmann connection,\index[terms]{connection!Ehresmann}
and its monodromy action defines a homomorphism $\Z\arrow \Diff(T^3)$ \index[terms]{action!monodromy}
taking $\Z$ to $\langle A \rangle$. Therefore, the action
of $\pi_1(S^1)$ on $\pi_1(T^3)$ provided by
the exact sequence \eqref{_mapping_torus_Inoue_Equation_} 
is generated by $A$. However, the action of
$A$ on $\calo_K= \Z^3$ is given by $A=\tau(\zeta)$,
as shown above. This implies $\pi_1(M) = \calo_K\rtimes U$.

Consider a complex manifold $\C \times {\Bbb H}$
as a product $\R^3 \times \R^{>0}$.
Taking the real and one of two complex embeddings
of the cubic field $K$, we can identify $K$
with a lattice\index[terms]{lattice} in $\C \times \R$. The group
$U$ acts on $\C \times \R = K \otimes_\Q \R$
(\ref{_tensor_product_with_R_Corollary_})
by automorphisms preserving $\calo_K \subset \C \times \R$.
Its action $u \arrow \sigma(u)$ on $\C$ is complex linear, and the action
on $\R$ is by real homotheties $u \arrow \tau(u)$. Extend the action of
$U$ on $\C \times \R$ to $\C \times {\Bbb H}$
by taking $u$ to $\sigma(u)\times \tau(u)$,
where $\tau(u) \in \R^{>0}$ acts on ${\Bbb H}$
by real homotheties; this action is clearly
affine and complex linear. We defined the
action of the multiplicative group 
$U$ on $\C \times {\Bbb H}$
that is  equivariant with respect
to the action of the additive group $\calo_K$.
Taken together, these two actions define
an action of $\calo_K\rtimes U$ on
$\C \times {\Bbb H}$, which is  complex affine.

The projection $\C \times {\Bbb H}\arrow \R^{>0}$
taking $(z, h)$ to $\Im (h)$ induces a map
$\pi:\; \frac{\C \times {\Bbb H}}{\calo_K\rtimes U}\arrow S$ to the
circle $S^1 = \R^{>0}/\langle \tau(\zeta)\rangle$.
The fibres of this projection are isomorphic to
$\frac{\C \times \R}{\calo_K}= T^3$. This allows one
to identify the projection $\pi$ with the mapping 
torus $M\arrow S^1$ constructed above. We obtain
$\frac{\C \times {\Bbb H}}{\calo_K\rtimes U}=M$.

This complex surface is called {\bf an \index[terms]{surface!Inoue} Inoue surface of class $S_M$}, 
sometimes also {\bf Inoue surface of class $S^0$.}
Two other classes of Inoue surfaces, called  $S^+$ and $S^-,$
are obtained from quadratic number fields using
a similar construction. These surfaces are also
cocompact quotients of $\C \times {\Bbb H}$
by a discrete, cocompact group of affine transforms.

The  Oeljeklaus-\index[terms]{manifold!Oeljeklaus--Toma (OT)}Toma manifolds are generalizations
of \index[terms]{surface!Inoue} Inoue surfaces of class $S_M$; the construction
of the LCK metric\index[terms]{metric!LCK} on some of Oeljeklaus--Toma manifolds,
given in \ref{_Existence_of_LCK_metric_on_OT_}, 
also gives an LCK metric on Inoue surfaces of class $S_M$.

The LCK structure\index[terms]{structure!LCK} on \index[terms]{surface!Inoue} Inoue surfaces was discovered
by F. \index[persons]{Tricerri, F.} Tricerri (\cite{tric}). He produced LCK metrics on\index[terms]{surface!Inoue--Bombieri}
all Inoue--Bombieri surfaces except a subclass of $S^-$.
Then F. \index[persons]{Belgun, F. A.} Belgun has proven that the rest of Inoue--Bombieri surfaces
do not admit LCK structures (\cite{bel}).

\subsection{Oeljeklaus--Toma manifolds and LCK geometry}
\label{_OT_LCK_Intro_Subsection_}

The Oeljeklaus--Toma manifolds are constructed from the
same number-theoretic data as the \index[terms]{surface!Inoue} Inoue surfaces of class $S_M$.
Let $K$ be a number field, and $K \otimes_\Q \R$ be its
tensor product with $\R$. The tensor product $K \otimes_\Q \R$
is isomorphic to a direct sum of $t$ copies of $\C$
and $s$ copies of $\R$ (\ref{_tensor_product_with_R_Corollary_};
see also the exercises in Subsection \ref{_idempotents_Exercise_}).
The projections $K \otimes_\Q \R\arrow \R$ and
$K \otimes_\Q \R\arrow \C$ define {\bf the real embeddings}
and {\bf the complex embeddings} of $K$.
We say that $s$ is {\bf the number of real embeddings of $K$},
and $2t$ is {\bf the number of complex embeddings};
we write $2t$ because each complex embedding has
a complex conjugate, and hence  the complex embeddings come in pairs.

Fix a number field $K$ with $s> 0, t >0$, and let
$\calo_K$ be the corresponding field of integers, that is,
the integral closure of $\Z$ in $K$. By Dirichlet's units 
theorem, the group $\calo_K^*$ of units is isomorphic
to $\{\pm 1\} \oplus \Z^{s+t-1}$ (\ref{_Dirichlet_units_Theorem_}).
A unit is {\bf positive} if it is mapped to $\R^{>0}$ 
under all real embeddings.

The group of positive  units $\calo_K^{+,*}$ acts on $\calo_K$ multiplicatively.
The corresponding semi-direct product
$\calo_K \rtimes \calo_K^{+,*}$ is similar to the group
$\Z^3 \rtimes \Z$ as used in Inoue's work.

To construct the  Oeljeklaus-\index[terms]{manifold!Oeljeklaus--Toma (OT)}Toma manifold,
 we need to make a choice to replace
$\calo_K^{+,*}\cong \Z^{s+t-1}$ by a free abelian
subgroup $U$ of rank $s$.

Consider the additive action of $\calo_K$ on 
$K \otimes_\Q \R= \R^s\times \C^t$.
Embedding each copy of $\R$ into
the  Poincar\'e half-plane ${\Bbb H}$,
we obtain an additive action of
$\calo_K= \Z^{s+2t}$ on ${\Bbb H}^s\times  \C^t$.
The quotient
${\Bbb H}^s\times \C^t/\calo_K^+$ can be naturally identified 
with a torus $\R^s\times \C^t/\calo_K^+$
multiplied by $(\R^{>0})^s$.

Let $\sigma_1,..., \sigma_s:\; K \arrow \R$
be all the real embeddings.
Consider the action of the group $\calo_K^{+,*}$
of positive units on $(\R^{>0})^s$ taking
$\zeta \in \calo_K^{+,*}$ to a map
\[ 
(t_1, ..., t_s)\arrow (\sigma_1(\zeta) t_1, ..., \sigma_s(\zeta) t_s).
\] 
The group $U\subset \calo_K^{+,*}$ is chosen in such a way
that the restriction of this action to  $U\subset \calo_K^{*}$
is cocompact and properly discontinuous on $(\R^{>0})^s$.
The existence of such $U$ is non-trivial; it is proven by
 \index[persons]{Oeljeklaus, K.} Oeljeklaus and \index[persons]{Toma, M.} Toma (\cite{ot}).

The action of $\calo_K \rtimes U$ 
on ${\Bbb H}^s\times  \C^t$ is defined as follows:
the additive part $\calo_K$ acts by adding 
the elements of $\C$ and $\R$ associated with the
corresponding complex and real embeddings, and the multiplicative
part $U$ acts multiplicatively on the $\C$-components
and by multiplication by $\sigma_i(\zeta)\in \R^{>0}$
on the $i$-th ${\Bbb H}$-component.
Since $\R^s\times \C^t/\calo_K= T^{s+2t} \times (\R^{>0})^s$
and $U$ acts properly discontinuously and cocompactly
on the $(\R^{>0})^s$-component, the whole action of
$\calo_K \rtimes U$ on ${\Bbb H}^s\times \C^t$ is cocompact,
free and properly discontinuous. The quotient
$\frac{{\Bbb H}^s\times \C^t}{\calo_K \rtimes U}$
is called {\bf an  Oeljeklaus-\index[terms]{manifold!Oeljeklaus--Toma (OT)}Toma manifold};
it is  a compact, complex manifold of algebraic
dimension 0 (\cite{ot}).

It is not hard to see that an OT manifold does
not admit a K\"ahler metric (\ref{_OT_non-Kahler_Theorem_}).
When $t=1$, the OT-manifolds carry a natural
LCK metric\index[terms]{metric!LCK}, that is  defined as follows.\index[terms]{manifold!Oeljeklaus--Toma (OT)}

Let $\omega_P$ be the Poincar\'e metric on ${\Bbb H}$,
and $\Psi:\; {\Bbb H}^s\arrow \R^{>0}$
take $z_1, ..., z_s$ to $\prod_{i=1}^s \Im z_i$.

Consider the metric $\omega_C$ induced by the archimedean 
norm on the $\C$-compo\-nent, and let
$\omega_1$ be the product metric associated with $\omega_C$ and $\omega_P$
on each of the components of ${\Bbb H}^s\times \C$. 
We claim that the metric $\Psi \cdot \omega_1$ is $\calo_K \rtimes U$-invariant;
since $\omega_1$ is K\"ahler, this implies that the OT-manifold is LCK.\index[terms]{manifold!Oeljeklaus--Toma (OT)}

The parallel transport\index[terms]{parallel transport}
associated with the $\calo_K$-action
preserves both the Poin\-car\'e metric on ${\Bbb H}$
and the metric $\omega_C$ on the $\C$-component.
It remains only to show that the action of $U\subset \calo_K^{+,*}$
on ${\Bbb H}^s\times \C$ leaves $\Psi\cdot  \omega_1$ invariant.

 The action of $\zeta \in U$ multiplies
$\omega_C$ by $|\tau(\zeta)|^2$, where $\tau:\; K \arrow \C$ 
is the complex embedding of $K$, and leaves $\omega_P= dd^c \log(\Im z)$
invariant. However, it multiplies $\Psi$ by $\prod_{i=1}^s |\sigma_i(\zeta)|$,
hence the whole metric $\Psi\cdot \omega_1$  is multiplied by the number
$|\tau(\zeta)|\prod_{i=1}^s |\sigma_i(\zeta)|$. This product
is equal to the absolute value of the norm of $\zeta$,
that is  equal to 1 because $\zeta$ is a unit (\ref{_units_norm_+-1_Remark_}).
Therefore, $\Psi\cdot  \omega$ is $U\subset \calo_K^{+,*}$-invariant.

In \ref{_Existence_of_LCK_metric_on_OT_} we give
an alternative construction of the same LCK metric
in terms of a K\"ahler potential \index[terms]{metric!LCK}
on the universal cover $\tilde M ={\Bbb H}^s\times \C$.

Before the work of \index[persons]{Oeljeklaus, K.} Oeljeklaus and \index[persons]{Toma, M.} Toma, all known LCK manifolds
were either \index[terms]{surface!Inoue} Inoue surfaces, LCK manifolds with potential,\index[terms]{manifold!LCK!with potential} or blow-ups of those.
Based on this evidence, \index[persons]{Vaisman, I.} Vaisman conjectured that all compact LCK manifolds\index[terms]{manifold!LCK}
that are  non-K\"ahler have $b_1$ odd. The first counterexample
to this conjecture is constructed by \index[persons]{Oeljeklaus, K.} Oeljeklaus and \index[persons]{Toma, M.} Toma using
the OT-manifolds; they constructed an OT-manifold with all
odd Betti numbers even.\index[terms]{manifold!Oeljeklaus--Toma (OT)}

The LCK geometry of OT-manifold is in a certain sense
very similar to Vaisman manifolds.\index[terms]{manifold!Vaisman} Just like the Vaisman
manifolds, any LCK OT-manifold $M$\index[terms]{manifold!Oeljeklaus--Toma (OT)} is equipped with
a canonical holomorphic foliation, preserved by all holomorphic\index[terms]{foliation!holomorphic}
automorphisms and embeddings of OT-manifolds.
This foliation is transversally K\"ahler, and
the transversal K\"ahler metric is exact and
semi-positive on $M$.

However, the difference with Vaisman geometry is
pronounced. On a Vaisman manifold $M$, $\dim_\C M=n$,
the canonical foliation always has at least
$n$ compact leaves\index[terms]{foliation!canonical} (\ref{_Vaisman_bound_from_Lfsch_Corollary_}). 
By contrast, the canonical foliation on an OT manifold has no
compact leaves, and in fact an OT-manifold
cannot contain a compact complex curve
(\ref{_nocurves_in_OT_Theorem_}).

\subsection{Subvarieties in the OT-manifolds}\index[terms]{manifold!Oeljeklaus--Toma (OT)}

In this chapter, we explain the construction of OT-manifolds,
and prove several basic results about complex subvarieties
in OT-manifolds, using the exact transversally K\"ahler form\index[terms]{form!K\"ahler!transversal}
$\omega_0$ defined in Subsection \ref{_LCK_OT_Subsection_}.
Since $\omega_0$ is exact, the integral $\int_Z \omega_0^{\dim_\C Z}$
always vanishes, for any compact complex subvariety $Z\subset M$.
This implies that $Z$ is nowhere transversal to the
canonical foliation\index[terms]{foliation!canonical} $\ker \omega_0$.  The canonical
foliation is tangent to the $\C^t$-component
in the decomposition $\tilde M = \C^t\times {\Bbb H}^s$,
where $\tilde M$ is the universal cover of the OT-manifold.
Therefore, its rank is $t$.

When $\dim_\C Z=1$,
this implies that $Z$ is contained in a leaf of
a canonical foliation. Using the result of 
S. \index[persons]{Verbitskaya, S. M.} Verbitskaya, we prove that any leaf \index[terms]{foliation!canonical}
of the canonical foliation is biholomorphic to 
$\C^t$, and hence  it cannot contain complex curves.

When $t=1$ (this is precisely the case when
the OT-manifold is LCK), the condition
$\int_Z \omega_0^{\dim_\C Z}=0$\index[terms]{manifold!Oeljeklaus--Toma (OT)}
implies that $Z$ contains a leaf $L$
of the canonical foliation whenever
$Z$ intersects $L$. Using the ``strong
approximation theorem'' (a powerful
result of algebraic number theory),
we prove that the closure $\bar L$ of $L$ in $M$
is the image of $\C \times \prod_{i=1}^s l_i$
under the map $\C\times {\Bbb H}^s\arrow M$,
where $l_i \subset {\Bbb H}$ is the  line
$\Im z=\const$. It is clear that
$\bar L$ is not contained in any complex subvariety;
indeed, the minimal complex vector space
containing $T_x \bar L$ is $T_x M$.
This implies that an LCK \index[terms]{manifold!Oeljeklaus--Toma (OT)} Oeljeklaus--Toma manifold
has no complex subvarieties of positive dimension.

\section{Number theory: local and global fields}

This chapter is devoted to a generalization in higher
dimensions of the Inoue surfaces $S_M$. The construction
and the proof of the main result (the non-existence of
subvarieties) make use of several facts from number theory
and Galois theory,\index[terms]{Galois theory} see \cite{_Cassels_Frolich_},
\cite{bs}, \cite{st}. We present
them briefly in this section.

For a systematic view of number theory we suggest using  \cite{_Cassels_Frolich_,_Neukirch_}; here we recall some
results and definitions that are  relevant for the
geometry of \index[terms]{manifold!Oeljeklaus--Toma (OT)} Oeljeklaus--Toma manifolds.

\subsection{Normed fields}
\definition 
 A {\bf  valuation} on a field $k$ is a function $|\cdot |:\; k
\arrow \R^{\geq 0}$, with the following properties:
\begin{enumerate}\index[terms]{valuation}
 \item Non-degeneracy: $|x|=0$ $\Leftrightarrow$ $x=0$.
 \item Multiplicativity: $|xy|= |x||y|$.
 \item There exists
$c>0$ such that $|\cdot |^c$ satisfies the triangle inequality.\index[terms]{inequality!triangle}
\end{enumerate}

\example 
\begin{enumerate}
\item Standard examples are the usual  valuation on $\Q$, $\R$, $\C$.
 \item Let $p$ be a prime number, and let $m,n\in \Z$ be  coprime with $p$.
Define the {\bf $p$-adic valuation} on $\Q$ by 
$|\frac m n p^k|:= p^{-k}$.\index[terms]{valuation!$p$-adic}
\end{enumerate}

 \remark 
 The   $p$-adic valuation satisfies an additional 
``non-archimedean axiom'': $|x+y| \leq \max(|x|, |y|)$.
Such valuations are called {\bf non-archimedean}.\index[terms]{valuation!(non-)archimedean}

One can see that any power of a non-archimedean valuation is again
non-archimedean, and satisfies the non-archimedean triangle inequality.\index[terms]{inequality!non-archimedean triangle}
 
\hfill

\definition 
 Let $|\cdot |$ be a valuation on a field $K$.
Consider the topology on $K$ with open sets generated by
\[ B_\epsilon (x) := \{y \in k \ \ |\ \ |x-y| < 
   \epsilon\}.
\]
Two valuations are called {\bf equivalent}
if they induce the same topology.

\hfill

The equivalence of valuations can be characterized as follows:

\hfill

\theorem  
 Two valuations $|\cdot |_1, |\cdot |_2$
are equivalent if and only if
$|\cdot |_1= |\cdot |_2^c$
for some $c>0$.

\proof \cite[Chapter II, Section 4]{_Cassels_Frolich_}. \endproof

\hfill

\theorem  {\bf (Ostrowski)}\index[terms]{theorem!Ostrowski}\\
 Every valuation on $\Q$
is equivalent to the usual  (``archimedean'') 
one or to a $p$-adic one.

\proof \cite[Chapter II, Section 1]{_Cassels_Frolich_}. \endproof

\hfill

\definition 
 The completion of a field $k$ under a valuation
 $|\cdot |$ is the completion of $k$ in the metric $|\cdot |^c$,
where $c>0$ is a constant such that $|\cdot |^c$ 
satisfies the triangle inequality. It is clearly a field.
The completion of $\Q$ under the $p$-adic
valuation is called the {\bf  field of $p$-adic
numbers}, denoted $\Q_p$.\index[terms]{field!of $p$-adic numbers}

%

\subsection{Local fields}

Let
$\mathbb{F}_p(t)$ denote the field of rational functions
over a finite field $\mathbb{F}_p$.

\hfill

\definition 
 A {\bf finite extension} $[K\!\!:\!\!k]$ of fields
is a field $K \supset k$ that is  finite-dimensional
as a vector space over $k$. \index[terms]{extension!finite} 

A {\bf number field} \index[terms]{field!number}
is a finite extension of  $\Q$. 

A {\bf functional field}\index[terms]{field!functional}
is a finite extension of the field  .

A {\bf global field} is a number or a functional field. \index[terms]{field!global}

A {\bf local field} is the completion of a global field
under a non-trivial valuation.\index[terms]{field!local}\index[terms]{valuation}

\hfill

The next result shows how to produce local fields.

\hfill

\theorem 
 Let  $\bar k$ be a field that is  complete and locally
compact under some valuation. {Then $\bar k$
is a local field}.

\proof \cite[Chapter II, Section 6]{_Cassels_Frolich_}. \endproof

\hfill

Now let   $[K\!\!:\!k]$ be a finite extension\index[terms]{extension!finite}, and 
$x\in K$. Consider the multiplication by $x$ as a
$k$-linear endomorphism of $K$. As the extension is finite, the determinant of this endomorphism is well-defined.

\hfill

\definition 
 The determinant of the above endomorphism is called a {\bf norm} 
associated with the extension and is denoted $N_{K/k}(x)$.

\hfill

\remark 
 A norm defines a homomorphism $K^* \arrow k^*$ of multiplicative groups.
For Galois extensions, the norm  $N_{K/k}(x)$ {is the 
product of all elements conjugate to $x$.}\index[terms]{extension!Galois}

\hfill

Norms can be used to extend valuations.

\hfill

\theorem 
 Let $ [K\!:\!  k]$ be a finite extension of 
local fields, of degree $n$. {Then a  valuation 
on $ k$ is uniquely extended to $ K$.} Moreover,
{this extension is expressed as
$|x|:= \left|N_{K/k}(x)\right|^{\frac 1 n}$.}

\hfill

\proof The only thing to check is the triangle inequality.\index[terms]{inequality!triangle}
For non-archi\-me\-de\-an norms, it is an easy exercise in
divisibility. For archimedean norms, we embed $ K$ to 
$\C$ or $\R$ and use the triangle inequality 
for the usual (archimedean) norm on $\R$ or $\C$.
\endproof

\subsection{Valuations and extensions of global fields}

Recall that if  $A,B$ are extensions of a same field $k$, $\Char k=0$,
where $A\!\!:\!\!k$ is finite, and we consider $A\otimes_k B$ 
as a $k$-algebra, then $A\otimes_k B$ 
is a direct sum of fields containing $A$ and $B$
(Exercise  \ref{_tens_product_direct_sum_Exercise_}). 
Now we can state the following.

\hfill

\theorem 
 Let $k$ be a number field,
$|\cdot|$ a valuation, $K\!\!:\!\!k$ a finite extension, and  
$\bar k$ its completion. Consider the decomposition of $K
\otimes_k \bar k$ into the direct sum of fields 
$K \otimes_k \bar k:= \bigoplus_i \bar K_i$.
{Then each extension of the  valuation
$|\cdot|$ from $k$ to $K$ is induced from some  
$\bar K_i$, and all such extensions are non-equivalent.}

\hfill

When $k=\Q$, and $|\cdot|$ is the usual
(archimedean) valuation, we obtain that
all $K_i$ are extensions of $\R$, 
that is, isomorphic to $\R$ or $\C$. This gives the following.

\hfill

\corollary\label{_tensor_product_with_R_Corollary_} 
 {For each number field $K$ of degree $n$
over $\Q$, there exists only a finite number
of different homomorphisms $K\hookrightarrow \C$},
all of them injective. 
Denote by $s$ the number of embeddings
whose image lies in $\R \subset \C$
(such an embedding is called {\bf ``real''}),
and $2t$ the number of embeddings whose image
does not lie in $\R$ ({\bf ``complex embeddings''}). 
Then $s+2t =n$.

\subsection{Dirichlet's unit theorem}

\definition 
 Let $[K\!:\!\Q]$ be a number field of degree $n$.
 The {\bf ring of integers} $\calo_K \subset K$ 
is the integral closure of $\Z$ in $K$. It can be
defined as the set of all roots $\alpha$ of monic polynomials
$P(t)= t^n + a_{n-1} t^{n-1} +a_{n-2}t^{n-2} + \cdots + a_0$,
$a_i \in \Z$, which lie in $K$.

\hfill

\claim 
The additive group $\calo_K^+$ of a number field $K$ of
degree $n$
is a finitely generated abelian group of rank $n$.

\hfill

\definition 
 A {\bf unit} of the ring $\calo_K$ is an invertible element.\index[terms]{unit}

\hfill

\remark \label{_units_norm_+-1_Remark_}
Note that $u \in \calo_K$ is a unit if and only if 
the norm $N_{K/\Q}(x)\in \Z$ is invertible, that is,
$N_{K/\Q}(x)=\pm 1$.

\hfill

The construction of \index[terms]{manifold!Oeljeklaus--Toma (OT)} Oeljeklaus--Toma manifolds makes essential use of the following: 

\hfill

\theorem   {\bf (Dirichlet's unit theorem)}
\label{_Dirichlet_units_Theorem_}
\\ \index[terms]{theorem!Dirichlet's unit}
 Let $K$ be a number field that has $s$ real
embeddings and $2t$ complex ones. Then {the group of
units  $\calo_K^*$ is isomorphic to $G \times \Z^{t+s-1}$,}
where $G$ is the finite group\index[terms]{group!finite} of roots of unity contained 
in $K$. Moreover, if $s>0$,  one has $G= \{\pm 1\}$.

\hfill

\remark 
 For a quadratic field, the group of units is 
the group of solutions of Pell's equation. For instance, if $K=\Q(\sqrt{3})$ then any solution $(a, b)\in \ZZ^2$ of Pell's equation\index[terms]{Pell equation}
$$x^2-3y^2=1$$
will define a unit $a+b\sqrt{3}\in \calo_K^*$. By
contrast,  in $K=\Q(\sqrt{-3})$ the only units are $\pm1$
and $ \pm \varepsilon$, where $\epsilon$ is a non-real
root of unity of order $3$. In the general case, there are
no easy explicit descriptions of the group of units.


\section{Oeljeklaus--Toma manifolds}


We now follow the original presentation in \cite{ot}, see also \cite{ovu}.

Let $K$ be a number field that has $2t$
complex embeddings denoted $\tau_i, \bar \tau_i$ and 
$s$ real ones denoted $\sigma_i$,  $s>0$, $t>0$.  
We denote the group $\calo_K^*\cap \bigcap_i
\sigma^{-1}_i(\R^{>0})$ by  $\calo_K^{+,*}$.
Clearly, $\calo_K^{+,*}$ is a finite index
subgroup of the group of units of $\calo_K$. The elements
of  $\calo_K^{+,*}$ are called {\bf totally positive
  units} of the number field $K$.\index[terms]{unit!totally positive}

 Denote by $\H=\{y\in \C\,;\, \Im y>0\}$ the upper half-plane. 
For any $\zeta\in \calo_K$  denote by $T_\zeta$ the
automorphism of $\C^t\times\H^s$ given by
$$T_\zeta(x_1,...,x_t, y_1,...,
y_{s})=\bigg(x_1+\tau_1(\zeta),...,
x_t+\tau_t(\zeta),y_1+\sigma_1(\zeta),...,
y_{s}+\sigma_{s}(\zeta)\bigg).
$$

Similarly, for any totally positive unit $\xi$, let $R_\xi$ be the automorphism of $\C^t\times \H^s$ defined by
$$R_\xi(x_1,\ldots,x_t, y_1,\ldots, y_{s})=\bigg (\tau_1(\xi)x_1,..., \tau_t(\xi)x_t,\sigma_1(\xi) y_1,  ..., \sigma_t(\xi) y_s\bigg).$$
Note that the total positivity of $\xi$ is needed for $R_\xi$ to act on $\C^t\times\H^s.$

For any subgroup $U\subset \calo_K^{*, +}$, the above maps define a free action of the
semidirect product $U\ltimes \calo_K$ on $\C^t\times\H^s$ (see below for an explanation of semidirect
products).

It is proven in \cite{ot} that one can always find
subgroups $U\subset \calo_K^{*, +}$ such that the above
action is also discrete and cocompact; such subgroups are
called {\bf admissible subgroups}. Note that if $U$ is an
admissible subgroup then necessarily one has
$\rk_\ZZ(U)+\rk_\ZZ(\calo_K)=2(s+t),$
hence $\rk_\ZZ(U)=s.$ This explains why the condition
$t>0$ is needed: otherwise we would have  $\rk
\calo_K^*=s-1,$ strictly less than $s,$ and hence
admissible subgroups could not exist.

Let then $U$ be an admissible subgroup and let 
$\Gamma:=U\ltimes \calo_K$.

\hfill

\definition  { (\cite{ot})}
{\bf An  Oeljeklaus--Toma manifold} (or OT-manifold for short) is a quotient
$\C^t \times {\H}^s/\Gamma$. In \cite{ot} this manifold
was  denoted $X(K,U)$.\index[terms]{manifold!Oeljeklaus--Toma (OT)}

\hfill

\remark \label{inoue}
For $s=t=1$, one recovers a well-known construction of the
 Inoue surface $S_M$ (see Chapter \ref{inoue_lck} or \cite{inoue}). \index[terms]{surface!Inoue}

\subsection{The solvmanifold structure}\label{sol_OT}
Topologically, OT-manifolds are\index[terms]{manifold!Oeljeklaus--Toma (OT)} solvmanifolds.\index[terms]{solvmanifold}
Recall that {\bf a solvmanifold} is a quotient of a connected,
solvable Lie group by a cocompact discrete subgroup.
An easy way of obtaining a solvable group is to consider
{\bf the semidirect product} $A \rtimes B$ of abelian Lie groups.
Let $A, B$ be abelian Lie groups, and $\phi:\; B \arrow
\Aut(A)$ a group homomorphism. Then $A \rtimes B$ 
is the set of pairs $(a, b)\in A\times B$, with
the product defined as $(a, b)(a_1, b_1)= (a \phi(b) a_1, bb_1)$.

Now consider $\Gamma=\calo_K\rtimes U$, and let
$R_*:= U\otimes_\Z \R$ and $R_+:= \calo_K\otimes_\Z \R$
be the corresponding abelian Lie groups. The group $U$ acts on
$R_+$ as a family of commuting endomorphisms, and the
eigenvalues of these maps are given by $\sigma_i(u)$,
hence they are positive. Since $U$ is commutative, we can
simultaneously diagonalize all $u\in U$, and 
take logarithms and then exponents of the corresponding
operators acting on $R_+$. This extends the $U$-action on
$R_+$ to an $R_*=U\otimes_\Z \R$-action. Then
$\Gamma$ is contained as a cocompact, discrete lattice\index[terms]{lattice!cocompact}
in $R_+\rtimes R_*$, and the corresponding quotient
is naturally identified with the corresponding
OT-manifold.\index[terms]{manifold!Oeljeklaus--Toma (OT)}

The complex structure on the universal covering
$R_+\rtimes R_*$ is by construction left-invariant,
but it is not right-invariant. Therefore, the
right action of $R_+\rtimes R_*$ on the
quotient $R_+\rtimes R_*/\Gamma=X(K,U)$
is not holomorphic, and the OT-manifold is
{\em not} a complex homogeneous manifold.\index[terms]{manifold!homogeneous}

\subsection{The LCK metric}\index[terms]{metric!LCK}
\label{_LCK_OT_Subsection_}

\theorem \label{_OT_non-Kahler_Theorem_}
{( \cite{ot})} \index[terms]{manifold!Oeljeklaus--Toma (OT)}
The  OT-manifolds cannot bear any K\"ahler metric.

\hfill

\proof We give a direct proof; for a topological
argument, see \cite{ot}. Denote by $\zeta_i$ the
coordinates on the $\H^s$ part of the universal covering of
an OT-manifold $X(K,U)$. Consider the 2-form
$\omega_0:=dd^c\sum\log\im(\zeta_i)$ and the 1-form
$\eta:=d^c\sum\log\Im(\zeta_i)$.
The form $\omega_0$ is equal to the Poincar\'e 
metric on $\H^s$, and therefore, it is
positive definite along the $\H^s$-directions. 
It is equal to zero in the $\C^t$-directions.

The homotheties
in $U$ will multiply $\Im(\zeta_i)$ by constants, adding
a constant additive term to $\sum\log\Im(\zeta_i).$
This term is annihilated by $d^c$. The
translations in $U$ do not affect the imaginary parts of
$\zeta_i$. This proves that $\eta$ is invariant at the
action of $U$ and hence $\eta$ and $\omega_0$ can be
considered to be  forms on $X(K,U)$. We have constructed
a (semi-)positive, exact $(1,1)$-form $\omega_0=d\eta$ on 
the OT-manifold $X(K,U)$. If $X(K,U)$ admits a K\"ahler metric $\omega$, then\index[terms]{manifold!Oeljeklaus--Toma (OT)}
$\omega^{n-1}\wedge\omega_0$ is strictly positive and
exact and hence  $0<\int_{X(K,U)}
\omega^{n-1}\wedge\omega_0=0$, contradiction.
 \endproof

\hfill

  However, some OT-manifolds can have LCK metrics\index[terms]{metric!LCK}.

\hfill

\theorem { (\cite{ot})}\label{_Existence_of_LCK_metric_on_OT_} \\
Let $K$ be a number field with $2t=2$ complex embeddings. Then, for any admissible group of totally positive units $U$, the  manifold $X(K, U)$ has an LCK metric\index[terms]{metric!LCK}.

\hfill

\proof We construct an automorphic 
K\"ahler metric on $\C \times {\H}^s$ of the form  $dd^c\phi$, where 
$\phi(x, \zeta_1, ..., \zeta_s)= |x|^2 + \prod_{i=1}^s \Im(\zeta_i)^{-1}$ (compare with the construction in the proof of \ref{_OT_non-Kahler_Theorem_}).
The function $\phi$ is clearly 
plurisubharmonic (its $dd^c$ is the 
Poincar\'e metric on each $\H$, and Euclidean on $\C$),
and $dd^c\phi$ is $\calo_K$-invariant. Any $\xi\in \calo_K^{+,*}$
multiplies $|x|^2$ by $A:=\tau(\xi)^2$ and $\prod_{i=1}^s \Im(\zeta_i)^{-1}$
by $B$, where $B:=\prod_{i=1}^s \sigma_i(\xi)^{-1}$. However, $AB=N(\xi)=1$,
because the norm $N(\xi)$ is integer and invertible, and
$\xi$ is totally positive (\ref{_units_norm_+-1_Remark_}). 
This proves that the K\"ahler metric $dd^c\phi$ is automorphic.

An alternative construction of the same LCK structure\index[terms]{structure!LCK} is given in
Subsection \ref{_OT_LCK_Intro_Subsection_}.
\endproof

\hfill

\remark 
For $s=t=1$, the LCK metric\index[terms]{metric!LCK} coincides with the warped product 
metric in \cite{tric} on the \index[terms]{surface!Inoue} Inoue surfaces of type  $S_M$.

\hfill

\remark
The potential $\f$ {\bf is not} acted on by homotheties. Moreover, it was proven in \cite{Kas}, using the solvmanifold structure,\index[terms]{solvmanifold} that OT-manifolds cannot be Vaisman. It can be shown that, more generally, OT-manifolds cannot have LCK potential\index[terms]{potential!LCK}, and, still more generally, cannot have $d_\theta$-exact LCK forms, see \cite{oti2}.\index[terms]{manifold!Oeljeklaus--Toma (OT)}

\hfill

\remark The above constructed LCK metric (and any other
possible LCK metric)\index[terms]{metric!LCK} on $X(K,U)$ is not GCK because
$X(K,U)$ is not of K\"ahler type (\ref{_OT_non-Kahler_Theorem_}).

\hfill

\remark \label{_Vaisman_conjecture_disproved_}\index[terms]{conjecture!Vaisman}
 It was conjectured in \cite{va_tr} that any compact LCK
 manifold\index[terms]{manifold!LCK} $X$ should have at least one odd rational Betti number of\index[terms]{Betti numbers!rational}
 odd degree. The following counterexample in complex
 dimension $3$ was given in  \cite{ot}. Let $K$ have
 $s=2, t=1$ and take any admissible subgroup $U$. Then
 $X=X(K, U)$ is LCK, by  \ref{_Existence_of_LCK_metric_on_OT_}, and $\dim_\C
 X=s+t=3$. The Betti numbers can be computed (see
 \ref{topinv} for details) as $b_1(X)=s=2-b_2(X)$.  To
 compute $b_3(X)$, observe that $X$ carries a global,
 nowhere vanishing 1-form (see Exercise \ref{form}), and
 hence its Euler--Poincar\'e characteristic vanishes, so
 $b_3(X)$ has to be  even too.\index[terms]{Euler--Poincar\'e characteristic}

\section[Non-existence of complex subvarieties in OT-manifolds]
{Non-existence of complex subvarieties\\ in OT-manifolds with LCK structure}\index[terms]{structure!LCK}

We now prove that LCK OT-manifolds are {\bf simple},
meaning that they have no subvarieties. They are the only known\index[terms]{manifold!Oeljeklaus--Toma (OT)}
example of LCK manifolds\index[terms]{manifold!LCK} with this very strong property.
Our proof will make use of  the Strong Approximation
Theorem, a high-tech version of the {\em Chinese remainder
theorem.}\index[terms]{theorem!Chinese remainder}

\hfill

\theorem  { (\cite{ov_mrl_11})} \label{no_subv}  %
 \ 
 Let $K$ be a number field that has $s$ real embeddings
and $2t$ complex ones, $t=1$, $s>0$.
{Then the corresponding  Oeljeklaus--Toma manifold
has no non-trivial complex subvarieties}.\index[terms]{manifold!Oeljeklaus--Toma (OT)}

\hfill

\proof 
Consider an OT-manifold $X=\C\times \H^s/\Gamma$. 
Let $x$ be the coordinate on $\C$ and 
$\zeta_1,..., \zeta_s$ the coordinates on $\H^s$. 
The following construction already appeared in the proof
of \ref{_OT_non-Kahler_Theorem_}.
Consider
the function 
$$\phi(x, \zeta_1, \ldots, \zeta_s):= \prod_i \Im(\zeta_i)$$
on $\C \times {\H}^s$.
Since the action of $\Gamma$ multiplies $\Im(z_i)$
by a number, the form $d\log \phi$ is
$\Gamma$-invariant.
Let $\theta$ denote the corresponding 1-form on $X$. Then the Lee form\index[terms]{form!Lee} of the LCK metric\index[terms]{form!LCK}\index[terms]{metric!LCK} of $X$ is $\theta=-d\log\f$ (see \ref{_Existence_of_LCK_metric_on_OT_}).

The 2-form $\omega_0:=d(I\theta)= dd^c \log \phi$ has Hodge type (1,1) and is 
positive definite on the leaves of the foliation
$\{x\}\times {\H}^s\subset \C \times {\H}^s$ 
\[
\omega_0 = \1\6\bar\6\log \phi = 
\1 \sum_i |\Im z_i|^{-2} dz_i\wedge d\bar z_i .
\]
Observe that  $\omega_0 \geq 0$ (the reader should notice
the close resemblance with the similar form on LCK
manifolds with potential).\index[terms]{manifold!LCK!with potential}

\hfill

\remark\label{_canoni_fol_Remark_}
The form $\omega_0:= \1\6\bar\6\log \prod_i \Im(\zeta_i)$
is well-defined on any OT-manifold \index[terms]{manifold!Oeljeklaus--Toma (OT)}
$M= \frac{\C^t\times\H^s}{\calo_K\rtimes U}$; it is 
exact and semi-positive of rank $s$. Its kernel is
the tangent space to the $\C^t$-foliation in $TM$;
it is a holomorphic foliation of constant rank.
We call $\ker \omega_0$ {\bf the canonical foliation}
of an OT-manifold.\index[terms]{foliation!holomorphic}\index[terms]{foliation!canonical}

\hfill

Let $\Sigma \subset TX$
be the null-foliation of $\omega_0$ (the foliation
generated by its null eigenspace). {
It is a holomorphic, involutive, smooth 
complex 1-dimensional foliation,}
with the leaves obtained from
$\C\times \{(\zeta_1, \ldots, \zeta_s)\}\subset \C \times {\H}^s$.

Now suppose there exists
a compact complex subvariety $Z\subset X$. 
Let $k=\dim_\C Z$. Then the integral $\int_Z \omega_0^k=0$, because
$\omega_0$ is exact. Therefore, $Z$ is at each point
tangent to a leaf of $\Sigma$.
{Since $\Sigma$ is 1-dimensional, this
means that $Z$ contains all leaves of $\Sigma$ intersecting $Z$.}

{It remains to show that any variety
that contains a leaf of $\Sigma$ coincides with $X$.}
To this end, let $\Sigma_0$ be a leaf of $\Sigma$.
Its preimage in $\C \times {\H}^s$
contains the set
\[
 \widetilde \Sigma_0(\zeta_1,\ldots, \zeta_s) :=\bigcup_{x\in \C,\  \xi\in \calo_K^+}
 \bigg (x, (\zeta_1+\sigma_1(\xi), \ldots, \zeta_s + \sigma_s(\xi))\bigg)
\]
where $(\zeta_1, \ldots, \zeta_s)\in {\H}^s$ is some fixed point
and $\sigma_i:\; K \arrow \R$ are the real embeddings.

Let 
\[
X_{\alpha_1, \ldots, \alpha_s}:=
\{(x, \zeta_1, \ldots, \zeta_s)\ \ 
|\ \ \Im \zeta_i=\alpha_i, \ i = 1, \ldots, s\}
\]
where $\alpha_i = \Im z_i>0$.

Since the smallest complex subspace containing
$T_pX_{\alpha_1, \ldots, \alpha_s}$ is $T_pX$, for all $p\in X_{\alpha_1, \ldots, \alpha_s}$, the theorem is implied by the following.

\hfill

\claim 
 The closure of $\tilde\Sigma_0(\zeta_1,\ldots, \zeta_s)$ 
contains the set $X_{\alpha_1, \ldots, \alpha_s}$, where
$\alpha_i = \Im \zeta_i$.

\hfill

This claim is immediately implied by the
following statement, applied to the set  
$\sigma_1, ..., \sigma_m$ of all real embeddings.

\hfill

\theorem \label{dens} 
 Let $[K\!:\!\Q]$ be a number field,
$\{\tau_i, \bar\tau_i\, ;\, i=1,\ldots,t\}$ its
 complex embeddings  and
  $\{\sigma_j \, ;\, j=1,\ldots, s\}$ the real embeddings, $t>0$,   $m=2s+t$. Consider
a set $\rho_1, \ldots, \rho_m$ of embeddings of $K$ 
to $\C$ or $\R$, with each of $\tau_i$
and $\sigma_i$ appearing  once, except for one. Consider the map
$R:\; K \arrow \R^a \times \C^b$,
 $R(\xi):=(\rho_1(\xi), \ldots, \rho_m(\xi))$.
{Then the image of $\calo_K$ is dense in $\R^a \times \C^b$.}

\hfill

The proof of this result is based on the {\bf Strong Approximation Theorem}. 
To be able to state it, we recall the following.

\hfill

\definition 
Let $K$ be a number field.
Consider the product $\prod_\nu K_\nu$ 
of all completions of $K$, where $\nu$ runs through the set
of equivalence classes of norms  on $K$, both archimedean and non-archimedean.
The {\bf adelic group} ${\cal A}_K\subset \prod_\nu K_\nu$
is the set of all sequences\index[terms]{group!adelic}
$(x_{\nu_1}, ..., x_{\nu_n}, ...)\in \prod_\nu K_\nu$
where $|x_{\nu_i}|\leq 1$ for all $i$
except a finite number, equipped with the product topology.

\hfill

\remark 
Tikhonov's theorem {implies that ${\cal A}_K$
is locally compact.} \index[terms]{theorem!Tikhonov}

\hfill

Now we can state:

\hfill

\theorem  ({\bf Strong Approximation Theorem}, \cite{st})  \\
Consider the natural embedding\index[terms]{theorem!strong approximation}
$K\hookrightarrow {\cal A}_K$. {Then its
image is a discrete, cocompact subgroup.}
Moreover, the projection of
${\cal A}_K \stackrel{P_{\nu_0}}\arrow \prod_{\nu\neq \nu_0} K_\nu$
to the product of all completions except for one
{maps $K$ to a dense subset of
$R_{\nu_0}({\cal A}_K)$.}

\hfill

We use the Strong Approximation Theorem for the projection
to the product of all completions except for one archimedean
completion, as follows.

\hfill

{\bf Proof of \ref{dens}:}  Consider $K$  as a subring of
${\cal A}_K$ and let $\calo_{A_K}$ be the ring of all {integer
  adeles}, that is, all $n$-tuples  
$(x_{\nu_1}, \ldots, x_{\nu_n},...)\in {\cal A}_K$, such that
$|x_{\nu_i}|\leq 1$ for each non-archimedean valuation.
By construction of the topology on ${\cal A}_K$, the subring
$\calo_{A_K}\subset {\cal A}_K$ is open in ${\cal A}_K$.
Clearly, 
 \begin{equation}\label{u}
 \calo_K= K \cap \calo_{A_K}.
\end{equation}
Now let  $P:\; {\cal A}_K \arrow {\cal A}_1$
be the projection of ${\cal A}_K$ to the product of all 
completions except for one archimedean. As projections are open mappings, {
and  $\calo_{A_K}$ is open in ${\cal A}_K$,
its projection to ${\cal A}_1$ is open in ${\cal A}_1$}.

By the Strong Approximation Theorem, the image $P(K)$ is dense in ${\cal A}_1$. 
Therefore, the image $P(K)\cap P(\calo_{A_K})$ 
is dense in the image $P(\calo_{A_K})$ of an open subset $A_K\subset {\cal A}$. 
{From \eqref{u}, we infer that $P(K)\cap P(\calo_{A_K})$ coincides with
$P(\calo_K)$.}

We proved that
$P(\calo_K)$ is dense in
${\cal A}_1\cap P(\calo_{A_K})$. {
Therefore, its projection to the product
of all archimedean completions except for one is also dense.}
\endproof

\hfill

The following corollary is immediate.

\hfill

\corollary \index[terms]{function!meromorphic!}
An LCK Oeljeklaus--Toma manifold does not admit
non-cons\-tant meromorphic functions.\footnote{This is known for all OT-manifolds,\index[terms]{manifold!Oeljeklaus--Toma (OT)}
but the proof is more complicated (\cite{_Battisti_Oe_}).}
\endproof




\section{Non-existence of curves on OT-manifolds}


Here we reproduce a result of S. \index[persons]{Verbitskaya, S. M.} Verbitskaya
\cite{sima}, who has proven that an OT-manifold has no
compact curves. This is a generalization of a similar
result for \index[terms]{surface!Inoue} Inoue surfaces, proven by Inoue \cite{inoue}.\index[terms]{manifold!Oeljeklaus--Toma (OT)}

\hfill

\theorem \label{_nocurves_in_OT_Theorem_}
Let $X$ be an OT-manifold. Then $X$ has no complex
subvarieties of dimension 1.

\hfill

\proof Let  $X=\frac{\C^t \times \mathbb{H}^s}{U\ltimes \calo_K}$ be an OT-manifold, and
$\tilde X$ its $U$-covering,
with $\tilde X = \C^t \times \mathbb{H}^s/\calo_K$.
Let $\zeta_i$, $i=1, ..., s$, be the complex coordinates on each factor 
${\H}$,
and  $\phi:= \log \left(\prod_{i=1}^s \Im \zeta_i\right)$.
As in the proof of \ref{no_subv}, the function $\phi$ is plurisubharmonic, 
 but not strictly
plurisubharmonic. Also, $\theta:=d^c\phi$  is $\Gamma$-invariant,
hence well-defined on $X$, and $\omega_0:=d\theta$ is a semi-positive
Hermitian form\index[terms]{form!Hermitian} (\ref{_canoni_fol_Remark_}). 
Since this form is exact, it has to vanish
on any 1-dimensional subvariety $Z\subset X$. This means
that $Z$ is everywhere tangent to the sub-bundle
$\Sigma:=\ker \omega_0$ which consists of vectors tangent
to $\C^t$. Therefore,
$Z$ is contained in a single leaf of the foliation
$\Sigma\subset TX$.

To finish the proof of \ref{_nocurves_in_OT_Theorem_},
it remains to show that the leaves of the foliation
$\Sigma\subset TM$ are biholomorphic to $\C^t$, and hence  cannot contain a compact
complex subvariety. However, each leaf $L$ is affine and
complete, and hence  it is obtained as a quotient of a
flat affine subspace $\tilde L= \C^t\subset \C^t \times \mathbb{H}^s$
by the stabilizer of $\tilde L$  in $\Gamma$. 

Let $\gamma\in \Gamma$ be an element fixing a leaf
$\tilde L$, with $\gamma=(u, a)\in U\ltimes \calo_K$.
Then for each real embedding $\sigma_i:\; K \arrow \R$, 
we have $\sigma_i(u)\zeta_i + \sigma_i(a)=\zeta_i$, $i= 1, ..., s$,
where $\zeta_i$ are coordinates on the $i$-th component $\mathbb{H}$.

These equations imply that $\zeta_i =
\frac{\sigma_i(a)}{1-\sigma_i(u)}$,
when $u\neq 1$. Hence 
 $\zeta_i\in\R$, but ${\H}$ has no real
elements, contradiction. We obtained that
$u=1$,  which gives  $\zeta_i +  \sigma_i(a)=\zeta_i$, implying
$a=0$ as well. We have shown that the stabilizer of
$\tilde L$ in $\Gamma$ is trivial.
\endproof

\hfill

The argument used in the proof of
\ref{_nocurves_in_OT_Theorem_}
also proves the following useful result.

\hfill

\claim\label{_cano_foli_OT_leaves_Claim_}\index[terms]{manifold!Oeljeklaus--Toma (OT)}
Let $M$ be an OT-manifold, $\Sigma$ the canonical
foliation on $M$ (\ref{_canoni_fol_Remark_}),
and $S$ a leaf of $\Sigma$ (this leaf is immersed
in $M$; the manifold structure on $S$ is obtained
from this immersion). Then $S$ is biholomorphic to $\C^t$.
\endproof


\section{Exercises} 


\subsection{Idempotents in tensor products}
\label{_idempotents_Exercise_}

\definition
Let $R$ be a commutative algebra over a field $k$.
$R$ is {\bf finitely--generated Artinian ring over $k$}
if $R$ is finite-dimensional as a vector space over $k$.
For the present purposes, we will tacitly abbreviate 
``finitely--generated Artinian ring'' as ``Artinian ring''.
An Artinian ring is called {\bf semi-simple}\index[terms]{ring!Artinian}
if it does not contain non-zero nilpotents.\index[terms]{ring!Artinian!semi-simple}

\begin{enumerate}[label=\textbf{\thechapter.\arabic*}.,ref=\thechapter.\arabic{enumi}]

\item
Let $R$ be an Artinian ring with unity and without zero divisors.
Prove that $R$ is a field.

{\em Hint:}
Prove that any injective endomorphism of a finite-dimensional
space is invertible. Use this to find $x^{-1}$ for any given $x\in R$.

\item
Prove that any prime ideal in an Artinian ring is 
maximal.

{\em Hint:}
Use the previous exercise.

\item
Let $v\in R$ be an element of a finite-dimensional
algebra $R$ over $k$. Consider the subspace 
$k[v]\subset R$ generated by $1, v, v^2, v^3, \dots$.
Suppose that $\dim k[v]=n$. Prove that
$P(v)=0$ for a polynomial
$P= t^{n} + a_{n-1} t^{n-1} + \dots$ with coefficients in $k$.
Prove that this polynomial is unique.\footnote{
This polynomial is called {\bf the minimal polynomial}
of $v\in \R$.}\index[terms]{minimal polynomial}

\item
Let $v\in R$ be an element of an Artinian ring with unity over $k$, and
$P(t)$ its minimal polynomial. Consider the subalgebra
$k[v]\subset R$ generated by $v$ and $k$. Prove that
$R[v]$ is isomorphic to the ring  $k[t]/(P)$ of residues modulo $P(t)$.

\end{enumerate}

\noindent
\definition
Suppose that $v\in R$ satisfies $v^2=v$. Then 
$v$ is called {\bf an idempotent}.  Idempotents  $e_1, e_2 \in R$
{\bf are orthogonal} if $e_1e_2=0$. An idempotent $e\in R$
is called {\bf indecomposable}
if there are no non-zero orthogonal 
idempotents $e_2, e_3$ such that $e = e_2 + e_3$.

\begin{enumerate}[label=\textbf{\thechapter.\arabic*}.,ref=\thechapter.\arabic{enumi}]
\setcounter{enumi}{4}

\item
Let $e\in R$ be an idempotent in a ring.
Prove that $1-e$ is also an idempotent.
Prove that a product of idempotents is an idempotent.

\item
Let $e\in R$ be an idempotent in a ring.
Consider the space  $eR\subset R$
(image of the multiplication by $e$.
Prove that $eR$ is a subalgebra in $R$, $e$ is 
unity in $eR$, and  $R=eR \oplus (1-e)R$.

\item
Let  $e_2,e_3 \in R$ be orthogonal idempotents.
Prove that  $e_1:= e_2 +e_3$ is also an idempotent satisfying
$e_2,e_3 \in e_1R$ and $e_1R=e_2R \oplus e_3R$.

\item
Let $\Char k \neq 2$, and  $e_1, e_2, e_3$ 
idempotents in an algebra $R$ over $k$.  Suppose that $e_1 = e_2 + e_3$. 
Prove that $e_2$, $e_3$ are orthogonal.

\item
Let $R$ be a semi-simple Artinian algebra with unity, and $e\in R$ 
a non-decomposable idempotent. Prove that $eR$ is a field.

\item
Let $R$ be a semi-simple Artinian ring over a field $k$,  $\Char k \neq 2$.
Prove that $1$ can be decomposed to a sum of indecomposable orthogonal
idempotents, $1=\sum_{i=1}^r e_i$. Prove that such a decomposition is unique.

{\em Hint:}
To prove existence, take an idempotent $e\in R$,
decompose $R$ to a direct sum of two subrings, 
$R=eR \oplus (1-e)R$, and use induction over $\dim_k R$.
For uniqueness, take two different orthogonal decompositions,
$1=\sum_{i=1}^r e_i$, and $1=\sum_{j=1}^s f_j$,
and prove that $e_i= \sum_{j=1}^s e_i f_j$ 
is an orthogonal decomposition.

\item
Let $R$ be a semi-simple Artinian ring with unity over a field 
$k$, $\Char k \neq 2$. Prove that $R$ is isomorphic to a direct sum of fields.
Prove that this decomposition is unique.

{\em Hint:}
Use the previous exercise.
\end{enumerate}

\noindent
\definition
Let $R$ be an algebra over a field $k$. A bilinear
symmetric form  $g$ on $R$ is called {\bf invariant} if
$g(x, yz) = g(xy, z)$ for all  $x$, $y$, $z\in R$.

\begin{enumerate}[label=\textbf{\thechapter.\arabic*}.,ref=\thechapter.\arabic{enumi}]
\setcounter{enumi}{11}

\item
Let $R$ be an Artinian ring equipped with a bilinear invariant
form $g$, and  ${\mathfrak m}$ an ideal in $R$.
Prove that its orthogonal complement  ${\mathfrak m}^\bot$ is also an ideal.

\end{enumerate}

\noindent
\definition
Let $R$ be an Artinian ring over $k$.
Consider the bilinear form $a, b \arrow \Tr(ab)$,
where $\Tr(ab)$ is the trace of the endomorphism
$L_{ab}\in \End_k R$, $x \stackrel {L_{ab}}\arrow abx$.
This form is called {\bf the trace form}, denoted
 $\Tr_k(ab)$.\index[terms]{form!trace}

\begin{enumerate}[label=\textbf{\thechapter.\arabic*}.,ref=\thechapter.\arabic{enumi}]
\setcounter{enumi}{12}

\item
Let $A$ be a linear operator on an $n$-dimensional
vector space of characteristic 0, such that
$\Tr A= \Tr A^2 = ... = \Tr A^n=0$.
Prove that $A$ is nilpotent.

\item
Let $[K:k]$ be a finite field extension in characteristic 0.
Prove that the trace form is always non-degenerate.

{\em Hint:}
Prove that $\Tr_k(x, x^{-1})= \dim_k K$ for any $x\in K$.
\end{enumerate}

\noindent
\definition
A finite field extension $[K:k]$ with non-degenerate
trace form is called {\bf separable}.\index[terms]{extension!separable}

\begin{enumerate}[label=\textbf{\thechapter.\arabic*}.,ref=\thechapter.\arabic{enumi}]
\setcounter{enumi}{14}

\item
Find an example of non-separable finite field extension
in characteristic $p$.

\item \label{_Tr_semisimple_Exercise_}
Let $R$ be an Artinian ring over $k$
with non-degenerate trace form. Prove that $R$ is semi-simple.
Prove that for $\Char k =0$, the trace form is non-degenerate on any
semi-simple Artinian ring.

\item
Let $A$, $B$ be rings over a field $k$.
\begin{enumerate}
\item Prove that there exists a multiplicative operation
$(A\otimes_k B)\times (A\otimes_k B) \arrow A\otimes_k B$,
mapping $a\otimes b, a'\otimes b'$ to $aa'\otimes bb'$.
\item Prove that this operation defines
the ring structure on $A\otimes_k B$.
\end{enumerate}

\item
Let $R, R'$ be Artinian rings over $k$, and
$g, g'$ the trace forms on $R, R'$. Consider the tensor product
$R \otimes_k R'$, and the bilinear symmetric form $g\otimes g'$ on
$R \otimes R'$, acting as $g\otimes g'(a\otimes a', b\otimes b'):= g(a, a')g'(b, b')$.
Prove that  $g\otimes g'$ is equal to the form $a, b \arrow \Tr(ab)$.

\item
Prove that the tensor product of semi-simple Artinian rings is semi-simple
if $\Char k=0$.

{\em Hint:}
Use the previous exercise.

\item\label{_tens_product_direct_sum_Exercise_}
Let  $[K_1:k]$, $[K_2:k]$ be finite extensions, $\Char k=0$.
Prove that the $k$-algebra $K_1\otimes_k K_2$ is semi-simple.
Prove that $K_1\otimes_k K_2$ is a direct sum of finite
extensions of $k$.

\end{enumerate}

\subsection{OT-manifolds}\index[terms]{manifold!Oeljeklaus--Toma (OT)}

\begin{enumerate}[label=\textbf{\thechapter.\arabic*}.,ref=\thechapter.\arabic{enumi}]
\setcounter{enumi}{20}

\item
Let $[K:\Q]$ be a finite extension, 
and $\calo_K$ the integral closure of $\Z$ in $K$ (it is called
{\bf the ring of integers} in $K$).
Consider the set ${\goth U}$ of all non-archimedean valuations on 
$K$. 
\begin{enumerate}
\item
Prove that $\calo_K= \bigcap_{\nu\in {\goth U}} \{ x\in K \ \ |\ \ \nu(x) \leq 1\}$.
\item Consider $\calo_K$ as an abelian group.
Prove that $\calo_K$ is torsion-free\index[terms]{ring!torsion-free} and finitely generated.
Prove that $\calo_K \cong \Z^n$, where $n$ is the degree of the extension
$[K:k]$.
\end{enumerate}

\item
Let $u \in \calo_K$ be an element in the ring of integers.
Consider the multiplicative action of $u$ on $\calo_K=\Z^n$.
\begin{enumerate}
\item Prove that the norm of $u$ is equal to 
the determinant of the action of $u$ on $\Z^n$. 

\item Prove that $u$ is a unit of $\calo_K$ if and only if
the determinant of the action of $u$ on $\Z^n$ is 1.
\end{enumerate}

\item 
Let $A\in \SL(n, \Z)$.  Prove that the 
characteristic polynomial $\det(t \Id -A)$ is irreducible
if and only if $A$ does not preserve proper rational
subspaces $V \subset \Q^n$.

\item
\begin{enumerate}
\item
Let $A\in \SL(n, \Z)$ be an operator such that
$P(t):=\det(t \Id -A)$ is irreducible. Consider the
ring $\Z[A] \subset \End(\Z^n)$.  Prove that
$\Z[A]$ has no zero divisors and its fraction field
is $\frac{Q[t]}{(P(t))}$.

\item
Consider the natural map $\Z[A]\arrow K$, where
$K:=\frac{Q[t]}{(P(t))}$ is the corresponding finite
extension of $\Q$. Prove that the image
of $A$ belongs to $\calo_K$. Prove that
it is invertible in $\calo_K$.
\end{enumerate}

\remark\label{_Units_SL(n,Z)_Remark_}
The previous exercises explain that units in number fields
correspond to irreducible matrices $A\in \SL(n, \Z)$,
and, conversely, irreducible matrices produce units
in number fields.

 \item Let $\zeta_k$, $k=1,\ldots, s$ be the coordinates on $\H^s$.  Show that the 1-forms
 $$\omega_k=\frac{1}{\Im \zeta_k}d\Im\zeta_k,\quad k=1,\ldots, s,$$
 defined on $\C^t\times\H^s$ are $U\ltimes \calo_K-$invariant, and hence  descend to 1-forms on $X(K, U)$.\label{form} 

 
\item
Let $M$ be an OT-manifold, $\Sigma=\ker \omega_0$ the canonical\index[terms]{manifold!Oeljeklaus--Toma (OT)}
foliation (\ref{_canoni_fol_Remark_}), and $X\subset M$ a compact K\"ahler
submanifold. Prove that $T_xX\subset \Sigma$ for each point
$x\in X$.

\item\label{_OT_LCK_explicit_Exercise_}
Consider an LCK OT-manifold $M= \frac{\tilde M}{\Gamma}$, obtained as a \index[terms]{manifold!LCK}
quotient of $\tilde M={\Bbb H}^n \times \C$,
and $\omega_\C$ be the constant K\"ahler form\index[terms]{form!K\"ahler} on
the $\C$-component in ${\Bbb H}^n \times \C$.
Prove that the LCK-form\index[terms]{form!LCK} on $M$ can be written as
$\omega= \omega_0 + \frac{\omega_\C}{\psi}$ 
where $\omega_0= dd^c\sum_{i=1}^n\log \Im\zeta_i$,
where $\zeta_i$ are complex coordinates on 
the $\H$-components, and $\psi= \prod_{i=1}^n\Im\zeta_i$.

\item
Prove that a K\"ahler manifold does not admit a positive,
exact (1,1)-form.

\item 
Prove that an OT-manifold cannot admit compact K\"ahler
submanifolds of positive dimension.\index[terms]{manifold!Oeljeklaus--Toma (OT)}

{\em Hint:} Use the previous exercise and \ref{_cano_foli_OT_leaves_Claim_}.

\item
Prove that a Hopf manifold does not admit holomorphic
embeddings in an OT-manifold.\index[terms]{manifold!Hopf}

{\em Hint:} Show that any Hopf manifold contains an
elliptic curve, and that an OT-manifold does not.\index[terms]{manifold!Hopf}

\item
Let $M= \frac{\C\times\H^s}{\calo_K\rtimes U}$
be an OT manifold equipped with the standard LCK structure\index[terms]{structure!LCK}
$(M, \theta, \omega)$
and $\chi:\; \pi_1(M) \arrow \R^{>0}$
the homothety character\index[terms]{homothety character} of the LCK structure
(\ref{_homothety_character_Definition_}).
Prove that $\ker \chi =\calo_K\rtimes U$, where $U_0\subset U$
is the group of all units $\zeta\in U$ such that 
$|\tau(\zeta)|=1$, where $\tau$ is the complex embedding
of $K$ to $\C$.

\item Let $X\subset M$ be a complex surface in an OT-manifold,
and $\omega_0$ the canonical (1,1)-form on $M$ defined in 
\ref{_canoni_fol_Remark_}.\index[terms]{surface!complex}
 Prove that $\omega_0\restrict X$ is
a semi-positive form of constant rank 1 on $X$. 

{\em Hint:} Use \ref{_cano_foli_OT_leaves_Claim_}.

\item
\begin{enumerate}
\item Prove that any OT-manifold is equipped with a flat, torsion-free\index[terms]{connection!torsion-free}\index[terms]{connection!flat}
connection $\nabla:\; TM \arrow TM \otimes \Lambda^1 M$ 
such that $\nabla^{0,1}=\bar\6$ (the holomorphic structure operator on $TM$).\index[terms]{operator!holomorphic structure}

\item Prove that the local system associated with $(TM, \nabla)$
is a direct sum of 1-dimensional local systems $L_i$.

\item Prove that all $L_i$ are non-trivial and non-isomorphic.
\end{enumerate}

\item
Let $(M, \omega, \theta)$ be an OT-manifold equipped with 
an LCK form\index[terms]{form!LCK} $\omega= \Psi^{-1} \omega_1$ as in Subsection \ref{_OT_LCK_Intro_Subsection_}.\index[terms]{manifold!Oeljeklaus--Toma (OT)}
Consider the action of the Lie group $R_+\rtimes R_*$ on 
the universal cover $\tilde M$ (Subsection \ref{sol_OT}).
Prove that the pullback of $\omega$ to $\tilde M$ 
is ($R_+\rtimes R_*$)-invariant.

\item\label{_anti_Lee_Killing_OT_Exercise_}
Let $(M, \omega, \theta)$ be an OT-manifold equipped with 
an LCK form $\omega= \Psi^{-1} \omega_1$ as in Subsection \ref{_OT_LCK_Intro_Subsection_}.
\begin{enumerate}
\item Prove that the Lee form\index[terms]{form!Lee} of $M$ is equal to $-d\log\Psi$.
\item Prove that $|\theta|_{\omega_1}=\const$, where $|\cdot|_{\omega_1}$
denotes the metric associated with $\omega_1$.
\item Prove that $dd^c(\omega^{n-1})=0$, that is, $\omega$ is Gauduchon.
\item Prove that the Lee field $\theta^\sharp$, dual to $\theta$ with respect
to the metric associated with $\omega$, is Killing.\index[terms]{vector field!Killing} Prove that this
Lee field is not parallel.\index[terms]{vector field!parallel}\index[terms]{Lee field}
\item Prove that the anti-Lee field\index[terms]{Lee field!anti-} $I(\theta^\sharp)$
is not Killing.\index[terms]{vector field!Killing} Prove that it preserves the Riemannian volume on $M$.
\end{enumerate}

{\em Hint:} Use Exercise \ref{_OT_LCK_explicit_Exercise_}.

\end{enumerate}

\chapter{Appendices}

\epigraph{\it La chair est triste, h\'elas! et j'ai lu tous les livres.}{\sc\scriptsize St\'ephane Mallarm\'e, \ \ Brise marine}

\hfill

In the next three appendices  we give a couple of important 
notions of conformal geometry that were used, intermittently, in the body of the lectures. 


\section{Appendix A. Gauduchon metrics}\label{_appendix_gauduchon_}


\label{gauduchon_metric}\index[persons]{Gauduchon, P.}

On a Hermitian manifold $(M,g,I)$,  with fundamental form
$\omega$, the  {\bf Chern connection}\index[terms]{connection!Chern} is the unique
connection $\nabla^C$ in the tangent bundle satisfying the
following three properties:
\begin{itemize}
\item $\nabla^CI=0$,
\item $\nabla^Cg=0$,
\item its torsion $T^C$ is of type (0,2), i.\,e.   $T^C(IX,Y)=IT^C(X,Y)$.
\end{itemize}	
If $g$ is a  K\"ahler metric, the Chern connection\index[terms]{connection!Chern}
coincides with its  Levi--Civita connection, but in general
the two connections are different.

\hfill

In \ref{_Chern_conne_Definition_}, the Chern connection \index[terms]{connection!Chern}was defined on an\index[terms]{bundle!vector bundle!holomorphic Hermitian}
arbitrary holomorphic Hermitian vector bundle $B$. This
definition coincides with the one above when $B=TM$.

\hfill

The {\bf torsion 1-form} $\eta$ of the Chern connection\index[terms]{connection!Chern} is defined as
$$\eta(X)=\Tr(Y\mapsto T(X,Y)), \quad X\in T_xM, \quad x\in M.$$
It is shown in \cite{gau_st} that the torsion 1-form satisfies the relation
\begin{equation}\label{torsion_form1}
\eta=Id^*\omega,
\end{equation}
where $d^*$ is the codifferential (i.\,e.  , up to a constant
multiplier, $d^*=*d*$, where $*$ is the Hodge star
operator defined by the metric $g$). 
Equation \eqref{torsion_form1} is equivalent to:
\begin{equation}\label{torsion_form2}
d\omega^{n-1}=\eta\wedge\omega^{n-1}.
\end{equation}
Here $n=\dim_\C M$. 

Note that the map
$\Lambda^{1}M\xlongrightarrow{\wedge\omega^{n-1}}\Lambda^{2n-1}M$
is an isomorphism, and hence the equation $d\omega^{n-1}=\theta\wedge\omega^{n-1}$
has a unique solution $\theta$. From
\eqref{torsion_form1}, this equation also determines the
torsion 1-form.
Clearly, for an LCK manifold\index[terms]{manifold!LCK}, the Lee form\index[terms]{form!Lee} is proportional
with the constant factor to the torsion 1-form.

\hfill

\definition \label{_Gauduchon_definition_}
Let $(M,I)$ be a complex manifold, $\dim_\C M=n$. A Hermitian metric is called {\bf Gauduchon} if its fundamental form $\omega$  satisfies \[ dd^c\omega^{n-1}=0.
\]
We call the $dd^c$-closed form $\omega^{n-1}$ 
{\bf the  Gauduchon form of $(M, \omega)$}.\index[terms]{form!Gauduchon}
It is not hard to see that this form determines the Hermitian structure uniquely.

\hfill

Gauduchon metrics\footnote{Originally, Gauduchon called
  these metrics  {\bf standard}.} are essential in
defining stability of vector bundles over non-K\"ahler
manifolds (\cite{_Lubke_Teleman_}). They also appear very naturally in locally
conformally K\"ahler geometry\index[terms]{geometry!K\"ahler}. \index[persons]{Gauduchon, P.} Gauduchon's theorem 
(\ref{standard_metric}  below)
states
that such metrics exist in every conformal class of
Hermitian metrics on a compact manifold. 

\hfill

\claim\label{_Gauduchon_theta_coclosed_Claim_}
A Hermitian metric is Gauduchon if and only if  its Lee form\index[terms]{form!Lee} is coclosed: $d^*\theta=0$.

\hfill

\proof Indeed, if $g$ is Gauduchon, using \eqref{torsion_form1} we have (up to constant multipliers):
\begin{equation}\label{prima}
0=dd^c\omega^{n-1}=dIdI\omega^{n-1}=dId\omega^{n-1}=dI\big(I\theta\wedge\omega^{n-1}\big)=d\big((I\theta)\wedge\omega^{n-1}\big).
\end{equation}
On the other hand, by \eqref{torsion_form2} we obtain:
\begin{equation}\label{_*theta_Equation_}
\begin{split}
\theta=&I*d*\omega=*Id*\omega; \ \text{hence,}\\
*\theta&=Id*\omega=Id\omega^{n-1}=I(\theta\wedge\omega^{n-1})\\
&=(I\theta)\wedge\omega^{n-1}
\end{split}
\end{equation}
By equation \eqref{prima} this implies $d*\theta=0$, and finally, since $*$ is an isomorphism, $*d*\theta=0$, i.\,e.   $\theta$ is coclosed. \endproof

\hfill

The main result of this section is as follows.

\hfill

\theorem { (\cite{gau_st})}\label{standard_metric} 
\index[terms]{metric!Gauduchon}\index[terms]{theorem!Gauduchon}
In each conformal class of Hermitian metrics on a compact
complex manifold of complex dimension $n\geq 2$ there
exists a  Gauduchon
metric.\index[persons]{Gauduchon, P.}
Moreover, this Gauduchon metric is unique, up to 
a constant multiplier.

\hfill

For the proof, we shall need the following consequence of the maximum principle.

\hfill

\theorem ({\bf
  Strong Maximum Principle; E. Hopf, 1927}; \cite[Theorem 8.19]{gt})\\\label{hopf_theorem}
\index[terms]{theorem!E. Hopf (maximum principle)}\index[persons]{Hopf,
  E.} Let $M$ be a compact manifold with boundary, and $L$ a second-order  elliptic operator acting on functions and such that
$L(1)=0$. Then all local maxima of the 
solutions of the equation $L(f)\geq 0$ 
occur on the boundary of $M$.
In particular, the only solutions of the inequality
$L(f)\geq 0$ on $M$ without boundary (or $L(f)\leq 0$) are the
constant functions.

\hfill

{\bf Proof of Gauduchon's theorem. Step 1:} 
Let $\omega\in \Lambda^{1,1}M$ be an Hermitian
form on $M$. Then it is enough to find a function $ \phi> 0 $
such that $ \6 \bar \6 (\phi \omega ^ {n-1}) = 0 $. Consider the differential
operator
\[
L(\phi) = \frac{\6\bar\6 (\phi\omega^{n-1})}{\omega^{n}}
\]
{Then $ L $ is an elliptic operator,
	with the same symbol as the Laplace operator}.

{\bf Step 2:} On the space $C^\infty(M)$ we consider
the $L^2$-metric\index[terms]{metric!$L^2$} given by formula
$(x,y) = \int_M xy \omega^n.$
Let $ \alpha \in C^\infty(M)$. From the Stokes formula we obtain
\[
\int_M L(\phi)\alpha\omega^{n} =
\int_M  \6\bar\6 (\phi\omega^{n-1}) \alpha=
\int_M \phi\omega^{n-1}\wedge \6\bar\6\alpha,
\]
hence the adjoint of $L$ is $ L^*\alpha =\frac{\omega^{n-1}\wedge
		\6\bar\6\alpha}{\omega^{n}}$.

{\bf Step 3:} Since $ L^* $ vanishes on constants,
its kernel is one-dimensional, by the Strong Maximum Principle
(\ref{hopf_theorem}).\index[terms]{maximum principle}
By the index formula for elliptic
operators of the second order, $ \ind L^* = 0 $, and hence 
its cokernel is one-dimensional. {We have proven that $ \dim \ker L = 1 $.}
 
{\bf Step 4:} {For the existence of the 
	Gauduchon metric it is thus enough to make sure that any
	 non-zero function $ u \in \ker L $ is everywhere positive, or
	is everywhere negative.} Then $ \6 \bar \6 (u \omega^{n-1}) = 0 $, and hence, $ \phi^{1 / (n-1)} \omega $ is the Gauduchon metric.

{\bf Step 5:} Let $f \in \im L^*$.
Suppose that $f\geq 0$ or $f\leq 0$ everywhere.
Then $\pm f=L(g) \geq 0$ for some $g\in C^\infty(M)$.
Applying the Strong Maximum Principle again, we obtain that
$g=\const$ and $f=0$.

{\bf Step 6:} If $\ker L$ contains  a function $ u $ which takes both positive and negative values, then 
we can construct the function $ \alpha \in C^\infty(M) $,
that is  everywhere positive, and satisfies
$\int_M u \alpha \omega^n=0$. Then $\alpha\in (\ker L)^\bot=\im L^*$,
that is  impossible owing to the assertion of Step 5.
Therefore, each $ u \in \ker L $ is either non-positive,
or non-negative.

{\bf Step 7:} Recall the Harnack inequality 
(\cite[Corollary 8.21]{gt}). \index[terms]{inequality!Harnack's}
Let $ D $ be an elliptic operator on a domain
$\Omega'\supset \Omega$, with the closure of $\Omega$ in
$\Omega'$ compact. Then
there exists a constant $ C $ such that for any
$u\in C^\infty(\Omega')$, $u\geq 0$, $D(u)=0$,
we have $ \sup_\Omega u \leq C \inf_\Omega u $.
This immediatelyimplies that $ u \neq 0 $
everywhere on $\Omega$. \endproof

\hfill

\remark (i) From the above argument it is clear that
on a compact complex  manifold, the Gauduchon metric is
unique (up to a constant multiplier) in each conformal class of
Hermitian metrics. See also the generalization of \ref{standard_metric} to conformal geometry (without compatible complex structure) in \cite{gau_tor}.

(ii) Let $(M,I,[g])$ be a compact LCK manifold\index[terms]{manifold!LCK}. Then two
Gauduchon metrics in $[g]$ give rise to the same Lee form,\index[terms]{form!Lee}
which is  coclosed with respect to any of these Gauduchon
metrics.\index[terms]{metric!Gauduchon}

\hfill

The following is a generalization of the notion of a Gauduchon metric.

\hfill

\definition (\cite{fww13})\label{_k_gauduchon_defi_}
A Hermitian metric $g$ is called {\em $k$-Gauduchon} if its fundamental form $\omega$ satisfies $dd^c\omega^{k} \wedge \omega^{n-k-1}=0$. \index[terms]{metric!$k$-Gauduchon}

\hfill

\proposition (\cite{_Ornea_Otiman_Stanciu_})
	\label{_k_gauduchon_equiv_prop_}
	Let $(M, I)$ be a compact complex manifold of complex dimension $n>2$ and $\omega$ a Hermitian metric. Then:
\begin{enumerate}
		\item[(i)] $\omega$ is $k$-Gauduchon for any $1\leq k \leq n-1$ if and only if
		\begin{equation*}
			(n-k-1)dId\omega \wedge \omega^{n-2}=-(k-1)d* \theta,
		\end{equation*}
		\item[(ii)] If $\omega$ is $k_1$-Gauduchon and $k_2$-Gauduchon for some $1 \leq k_1 \neq k_2 \leq n-1$, then $\omega$ is $k$-Gauduchon, for any $1\leq k \leq n-1$ (in particular, it is Gauduchon).
\end{enumerate}
\proof
The metric $\omega$ is $k$-Gauduchon if and only if:
\begin{equation}\label{eq:aceeasiForma1}
	\begin{split}
				0 &= dd^c \omega^k \wedge \omega^{(n-k-1)} = dId (\omega^k) \wedge \omega^{n-k-1} \\ 
				&= kdI(d\omega \wedge \omega^{k-1}) \wedge \omega^{n-k-1} = kd(Id\omega \wedge \omega^{k-1}) \wedge \omega^{n-k-1} \\
				&= kdId\omega \wedge \omega^{k-1} \wedge \omega^{n-k-1} - k(k-1)Id\omega \wedge d\omega \wedge \omega^{k-2} \wedge  \omega^{n-k-1} \\
				&= kdId\omega \wedge \omega^{n-2} - k(k-1)Id\omega \wedge d\omega \wedge \omega^{n-3}; 
	\end{split}
\end{equation}
		thus,
\begin{equation}\label{eq:k-Gaud}
			dId\omega \wedge \omega^{n-2} = \frac{k-1}{n-2}Id\omega \wedge d(\omega^{n-2}), 
\end{equation}
		which gives
\begin{equation*}
			dId\omega \wedge \omega^{n-2} = \frac{k-1}{n-2}\left( dId\omega \wedge \omega^{n-2}-d(Id\omega \wedge \omega^{n-2})\right),
\end{equation*}
and furthermore,
\begin{equation*}
			(n-k-1) dId\omega \wedge \omega^{n-1}=-(k-1)dId\omega^{n-1}.
\end{equation*}
Now, since $d\omega^{n-1}=\theta \wedge \omega^{n-1}$ and $* \theta = \frac{I\theta \wedge \omega^{n-1}}{(n-1)!}$, we conclude this is equivalent to
\begin{equation*}
			(n-k-1)dId\omega \wedge \omega^{n-2}=-(k-1)d* \theta.
\end{equation*}	
In particular, if $\omega$ is $k$-Gauduchon with $k \neq n-1$, then:	
\begin{equation*}\label{eq:1G-integrala}
			\int_M dd^c\omega \wedge \omega^{n-2} = 0,
\end{equation*}
since $k \neq n-1$. This proves (i).
		
If $\omega$ is $k_1$ and $k_2$-Gauduchon, from \eqref{eq:k-Gaud} we get
		\[
		dId\omega \wedge \omega^{n-2} = \frac{k_1-1}{n-2}Id\omega \wedge d(\omega^{n-2}) = \frac{k_2-1}{n-2}Id\omega \wedge d(\omega^{n-2}).
		\]
		Note that \eqref{eq:k-Gaud} also holds for $k=n-1$. Now, since $k_1 \neq k_2$,
\[
		dId\omega \wedge \omega^{n-2} = Id\omega \wedge d(\omega^{n-2}) = 0,
\]
and, for any $1 \le k < n-1$, using \eqref{eq:aceeasiForma1},
\[
		0 = kdId\omega \wedge \omega^{n-2} - k(k-1)Id\omega \wedge d\omega \wedge \omega^{n-3} = dd^c \omega^k \wedge \omega^{n-k-1}
\]
i.\,e.    $\omega$ is $k$-Gauduchon. Moreover, $dId\omega^{n-1}=dId\omega \wedge \omega^{n-2}-(n-2)Id\omega \wedge \omega^{n-2}=0$, and hence  $\omega$ is also Gauduchon, and hence  (ii).
\endproof


\section{Appendix B. An explicit formula of the Weyl connection }\label{_Weyl_conne_Appendix_}


The  Weyl connection\index[terms]{connection!Weyl} on an LCK manifold\index[terms]{manifold!LCK} was defined 
in \ref{_Weyl_conn_definition_Remark_} in terms 
of the K\"ahler cover. Here we give an independent explicit
construction of Weyl connections in the spirit of
conformal geometry. 

\hfill

\definition
	Given a conformal class $c$ on a manifold $M$, a connection $\nabla$ in $TM$ is called {\bf compatible} with $c$ if $\nabla$ preserves the class $c$.	Explicitly, this means that for any $g\in c$ there exists a 1-form $\eta_g$ such that 
	\begin{equation}\label{compcon}
	\nabla g=\eta_g\otimes g.
	\end{equation}
A torsion-free, compatible connection on $(M,c)$ is called
a {\bf Weyl connection}. \index[terms]{connection!torsion-free}\index[terms]{connection!Weyl}

\hfill

\theorem
Let $c$ be a conformal class of a Riemannian metric
on a manifold $M$, and $\theta$ a 1-form. Then 
there exists a unique Weyl connection $\nabla$
that satisfies  $\nabla g=\theta\otimes g$
for some $g$ in the conformal class $c$.
Moreover, such $g$ is unique up to a constant
multiplier.

\proof \cite{gau_weyl}. \endproof

\hfill

The form $\theta$
is called {\bf the Lee form}\index[terms]{form!Lee} of the pair $(\nabla, g)$.

\hfill

If $g$ is an LCK metric\index[terms]{metric!LCK} on $(M,I)$, with the Lee form\index[terms]{form!Lee}
$\theta$, then any conformally equivalent Hermitian
metric $e^fg$ is again LCK with the Lee
form $\theta+df$.  

For an LCK manifold\index[terms]{manifold!LCK} with the Lee form $\theta$, we denote
by $D$ the Weyl connection on $TM$ associated with $\theta$.
Then, for an open set $U$ on which $\theta=df$ and
$g'=e^{-f}g$ is K\"ahler, we have 
$$Dg'=-dfe^{-f}g+e^{-f}Dg=-\theta\otimes g'+\theta\otimes g'=0,$$
and hence $D\restrict U$ is the Levi--Civita connection of
the K\"ahler metric $g'$. Clearly, this means that the
pullback of the Weyl connection $D$ coincides with the
Levi--Civita connection of the K\"ahler metric.

In particular, $D$ preserves the complex structure: $DI=0$.

Let $\nabla$ be the Levi--Civita connection of the LCK
metric $g$\index[terms]{metric!LCK}. Being the Levi--Civita connections of two
conformal metrics, $D$ and $\nabla$ are related by the
equation (\cite[Theorem 1.159 (a)]{besse}):
\begin{equation}\label{concon}
D=\nabla - \frac 12\big\{\theta\otimes\id+\id\otimes \theta- g\cdot \theta^\sharp\big\}.
\end{equation}

Conversely, if $\theta$ is a closed one-form on the Hermitian manifold $(M,I,g)$ and if $D$ is the connection given by \eqref{concon}, then $Dg=\theta\otimes g$, and then 
$$D\omega=Dg\circ I+g\circ DI=\theta\otimes g\circ I+ g\circ DI,$$
and if $DI=0$, then $D\omega=\theta\otimes\omega$. This immediately implies $d\omega=\theta\wedge\omega$.

We have proven the following intrinsic
characterization of LCK 
manifolds.\index[terms]{manifold!LCK}

\hfill

\proposition 
The Hermitian manifold $(M,I,g)$ is LCK if and only if $M$ admits a closed one-form such that the Weyl connection associated with it preserves the complex structure.\index[terms]{connection!Weyl}

\hfill

The following lemma contains useful computational formulae.

\hfill

\lemma\label{formulas}
Let $\nabla$ be the Levi--Civita connection\index[terms]{connection!Levi--Civita} on an LCK
manifold. Then the covariant derivative 
of the complex structure is given by the formula:
\begin{equation}\label{DJ}
(\nabla_XI)Y=\frac 12 \left(\theta(IY)X-\theta(Y)IX+g(X,Y)I\theta^\sharp+
\omega(X,Y)\theta^\sharp\right).
\end{equation}
The exterior derivative of $\theta^c$ is given by:
\begin{equation}\label{dthetac}
d\theta^c(X,Y)=\left(-|\theta|^2\omega+\theta\wedge\theta^c\right)(X,Y)+\big((\nabla_Y\theta)(IX)-(\nabla_X\theta)(IY)\big).
\end{equation}

\proof  Equation \eqref{DJ} is obtained from the relation \eqref{concon} between the Levi--Civita connection $\nabla$ of $g$ and the Weyl connection associated with $\theta$:
and using $DI=0$.

For \eqref{dthetac}, one writes (recall that $\theta^c(X)=-\theta(IX)$):
\begin{equation*}
\begin{split}
(\nabla_X\theta^c)(Y)&=\nabla_X(\theta^c(Y))-\theta^c(\nabla_XY)\\
&=-\nabla_X(\theta(IY))+\theta(I\nabla_XY)\\
&=-(\nabla_X\theta)(IY)-\theta((\nabla_XI)(Y)),
\end{split}
\end{equation*}
then uses \eqref{DJ} and the fact that $d\theta^c$ is the
anti-symmetrization of $\nabla\theta^c$.\index[terms]{anti-symmetrization}
\endproof

\hfill

\remark 
From \eqref{DJ} also follows:
\begin{equation}\label{domega}
\nabla_X\omega=-\frac 12\{X^\flat\wedge \theta^c+(IX)^\flat\wedge\theta\}.
\end{equation}

\remark  Equation \eqref{DJ} shows that  LCK
  manifolds\index[terms]{manifold!LCK} belong to the $W_4$ class in \index[persons]{Gray, A.} Gray--Hervella's
  classification of almost Hermitian manifolds, \cite{gh}.

{\setlength\epigraphwidth{0.5\linewidth}
\partepigraph{\it Every number is infinite; there is no difference.}
	{\sc\scriptsize Liber AL, I. 4.}
\part{Advanced LCK geometry}
\removeepigraph
}

\chapter{Non-K\"ahler elliptic surfaces}
\label{_elliptic_Chapter_}

\epigraph{\it Nomenclature, the other foundation of botany, should provide the names as soon as the classification is made... If the names are unknown knowledge of the things also perishes... For a single genus, a single name.}
	{\sc\scriptsize Carl Linnaeus}

\section{Introduction}


The main focus of this chapter is  elliptic non-K\"ahler
surfaces.  We give a new proof that they are principal 
elliptic\index[terms]{bundle!principal}\index[terms]{bundle!principal!elliptic}
bundles in the orbifold category (\ref{_elli_bundle_Theorem_}) -- 
a result originally proven by \index[persons]{Br\^inz\u anescu, V.} 
 Br\^inz\u anescu in \cite{_Brinzanescu:manuscripta_} 
(see also \cite{_Brinzanescu:bundles_}) --, and locally 
conformally K\"ahler (\ref{_elli_then_Vaisman_Theorem_}) -- a result 
due to \index[persons]{Belgun, F. A.} Belgun, \cite{bel}. We prove that all these surfaces are Vaisman, 
giving a new proof\index[terms]{orbifold} of Belgun's classification of Vaisman surfaces
(\cite{bel}).

This chapter, together with Chapter \ref{comp_surf},
is aimed at giving a comprehensive rendition of the classification theory
for non-K\"ahler complex surfaces. It turns out that the
non-K\"ahler complex surfaces are actually easier to classify than the
K\"ahler ones. A version of these two chapters appeared as a paper
\cite{_ovv:surf_}, jointly with Victor \index[persons]{Vuletescu, V.} Vuletescu; we are grateful to Victor
for allowing us to use this content.

The split between Chapters \ref{_elliptic_Chapter_} and
\ref{comp_surf} is somewhat artificial; in fact, there are
a few places in  Chapter \ref{_elliptic_Chapter_}  where we refer to results
proven in  Chapter \ref{comp_surf}. In Chapter \ref{_elliptic_Chapter_}
we explain the classification of non-K\"ahler elliptic surfaces, and prove that
all non-K\"ahler elliptic surfaces are Vaisman, and in 
Chapter \ref{comp_surf} we use the classification of elliptic
surfaces to classify all non-K\"ahler complex surfaces.

The main result of this chapter is 
\ref{_elli_then_isotri_Theorem_}, together with \ref{_isotrivial_Proposition_},
which claims that any minimal non-K\"ahler surface that admits
an elliptic fibration is in fact a principal elliptic fibration.
This gives a sharp contrast with the elliptic fibrations on K\"ahler-type 
surfaces, which come in a great many flavours and usually have 
non-trivial monodromy\index[terms]{monodromy} and picturesque special fibres.

In lieu of ``non-K\"ahler'' we use an equivalent
assumption, that is  justified in Chapter \ref{comp_surf}.
Given an elliptic fibration $\pi:\; M \arrow S$, we assume that
the pullback of the volume form $\Vol_S$ to $M$ is exact.
This is equivalent to $M$ being non-K\"ahler, as follows from
\ref{_elli_bundle_Theorem_}.

To explain the structure of elliptic fibration, we introduce two key 
notions of complex geometry: the Gauss--Manin connection\index[terms]{connection!Gauss--Manin},
together with the variations of Hodge structure\index[terms]{structure!Hodge!variation of}, and
Bishop--Gromov's theorem\index[terms]{theorem!Bishop--Gromov's compactness} on the  compactness of the  Barlet space\index[terms]{space!Barlet} of 
complex curves in a complex manifold. 

We refer specifically to \index[persons]{Gromov, M.} Gromov, even if Bishop proved this
result for complex varieties much earlier, because Gromov's
result gives an explicit description of degenerations for
the limits of complex curves. We apply it to describe
all degenerations which an elliptic curve might have
(\ref{_limits_of_elliptic_Corollary_}). This result
is used to show that the only degenerate fibre that might 
occur in a non-K\"ahler surface is a multiple fiber.

In \ref{_elli_bundle_Theorem_}
we prove that the monodromy\index[terms]{monodromy!of a connection} of the Gauss--Manin\index[terms]{connection!Gauss--Manin}
connection of the elliptic fibration on a non-K\"ahler manifold is 
trivial.

The multiple fibres are the subject of Section 
\ref{_multiple_fibres_Subsection_}, where 
we give a local description of a multiple
fibre of a minimal elliptic surface with 
trivial monodromy\index[terms]{monodromy!of a connection} of the Gauss--Manin connection\index[terms]{connection!Gauss--Manin}
(\ref{_Albanese_around_special_multiplicities_Proposition_}).
In a neighbourhood of a special fiber, such a
surface can always be obtained as
a finite quotient of a trivial fibration
over a disk by an action of a cyclic group. 

After we prove that all
non-K\"ahler elliptic surfaces are principal
fibrations, we use topological results about the
principal fibrations to prove that a $\Z$-cover
of any non-K\"ahler elliptic surface $M$ is the total
space of a principal $\C^*$-bundle over an orbifold.
A classical argument due to \index[persons]{Blanchard, A.} Blanchard
(\ref{_Blanchard_Theorem_}) implies that
this $\C^*$-bundle is associated with an
ample line bundle. Then \ref{_formula_15_19_besse_Theorem_}
and \ref{kami_or} implies that $M$ is Vaisman.

In Chapter \ref{comp_surf}, we build
upon this result to show that (unless
the Global Spherical Shell Conjecture is false)
all non-K\"ahler complex surfaces, except one of 
three species of Inoue surfaces,\index[terms]{conjecture!GSS}
are locally conformally K\"ahler.

\section{Gauss--Manin local systems and variations of Hodge structure}\index[terms]{structure!Hodge!variation of}

\subsection{The Gauss--Manin connection}

For the sequel, we recall some basic facts about the\index[terms]{connection!Gauss--Manin}
Gauss--Manin connection. We may refer to \cite{_Griffiths:transcendental_} 
or \cite{vois}. Let $\pi:\; M \arrow B$ be
a smooth, proper map of smooth manifolds. By \index[persons]{Ehresmann, C.} Ehresmann's\index[terms]{theorem!Ehresmann's fibration}
theorem, the fibres of $\pi$ are diffeomorphic. Then
$\pi$ is a locally trivial fibration, and hence  for any fixed $k$ the
$k$-th cohomology of its fibres form a local system, called
{\bf Gauss--Manin local system}.\index[terms]{local system!Gauss--Manin}
By the Riemann--Hilbert correspondence, the category of local systems
is equivalent to the category of vector bundles equipped with
\index[terms]{correspondence!Riemann--Hilbert}
a flat connection (\ref{_Riemann--Hilbert_Theorem_}).

The bundle associated with the  Gauss--Manin local system
is called {\bf the Gauss--Manin bundle}, and the connection -- 
{\bf the Gauss--Manin connection}. 
It can be constructed as follows.

Let $T_\pi M\subset TM$ be the bundle of fibrewise tangent vectors.
By definition, an Ehresmann connection $e$ on $M$ \index[terms]{connection!Ehresmann}
is a decomposition $TM= T_\pi M \oplus T_\hor M$,
that is, a choice of such a splitting.
Identifying $T_\hor M$ and the pullback $\pi^*TB$,
we may consider the pullback of a vector field
$X\in TB$ as a vector field $X_e\in T_\hor M$.

A section of a vector bundle associated
with the fibrewise cohomology of $M$
is given by a fibrewise closed differential form 
$\eta\in\Lambda^k M$. The Lie derivative 
$\Lie_{X_e}\eta$
is closed on fibres of $\pi$ for any vector field
on $B$ lifted to a horizontal vector field $X_e$ on $M$.
Indeed, the corresponding diffeomorphisms
map fibres to fibres and fibrewise closed
forms to fibrewise closed forms.

Since different choices $e, e'$ of the
Ehresmann connection result in the vector fields
$X_e, X_{e'}$ that satisfy $Y:=X_e-X_{e'}\in T_\pi M$,
and the form $\Lie_{Y}\eta$ is fibrewise exact,
the cohomology class of the restriction of
$\Lie_{X_e}\eta$ is independent on the
choice of the Ehresmann connection.\index[terms]{connection!Ehresmann}

Now let  $[\eta]$ be the collection of the
cohomology classes of $\eta$ on all fibres of
$\pi$, considered to be  a section of the Gauss--Manin
bundle. Define\index[terms]{bundle!Gauss--Manin} 
\begin{equation}\label{_GM_conne_Equation_}
\nabla_X[\eta]:= [\Lie_{X_e}\eta],
\end{equation}
where $[\Lie_{X_e}\eta]$ is the collection
of the cohomology classes of $\Lie_{X_e}\eta$
on all fibres of $\pi$. This formula
defines the Gauss--Manin connection $\nabla$.

\hfill

We needed this observation to prove the following
lemma.

\hfill

\lemma\label{_Gauss_Manin_d_1_vanishes_Lemma_}
Let $\pi:\; M \arrow B$ be a smooth fibration, and
$\eta$ a $p$-form that is  closed on the fibres of $\pi$.
Using the pullback map, we consider $\pi^*(\Lambda^*B)$ 
as a subspace in $\Lambda^*M$.
Assume that $d\eta$ belongs to 
$\pi^*(\Lambda^2B)\bigwedge \Lambda^{p-1}M\subset \Lambda^{p+1}M$.
Then the section of the Gauss--Manin bundle, corresponding
to $[\eta]$, is parallel.

\hfill

\proof
We use the formula \eqref{_GM_conne_Equation_}: 
$\nabla_X[\eta]= [\Lie_{X_e}\eta]$. The Cartan formula gives 
$\Lie_{X_e}\eta= i_{X_e} d\eta + d(i_{X_e}\eta)$, the second term
on the right-hand side is exact, and the first vanishes on fibres
because $d\eta\in \pi^*(\Lambda^2B)\bigwedge \Lambda^{p-1}M$.
\endproof

\subsection{Variations of Hodge structures}\index[terms]{structure!Hodge!variation of}

\definition
Let $V_\R$ be a real vector space.
{\bf A  Hodge structure\index[terms]{structure!Hodge} of weight $w$} 
on a vector space $V_\C=V_\R \otimes_\R \C$ 
is a decomposition $V_\C =\bigoplus_{p+q=w} V^{p,q}$, satisfying 
$\overline{V^{p,q}}= V^{q,p}$. It is called {\bf an integer 
Hodge structure} if one fixes an integer  lattice\index[terms]{lattice!integer} $V_\Z$ or $V_\Q$
in $V_\R$. A Hodge structure is \index[terms]{structure!Hodge}
equipped with a $\U(1)$-action, with $u\in \U(1)$
acting as $u^{p-q}$ on $V^{p,q}$. {\bf  A morphism}
of  Hodge structures is a linear 
map that is  $\U(1)$-invariant and preserves 
the Hodge structure.

\hfill

\claim\label{_tori_Hodge_structures_Claim_}
The category of complex tori with a group 
structure (that is, the zero point fixed) is equivalent to the category of
integer Hodge structures of weight 1.

\hfill

\proof Indeed,
a Hodge decomposition $V^{1,0} \oplus V^{0,1}$ on a space $V= V_\R \otimes_\R \C$
is the same as a complex structure on $V_\R$. The translation
invariant complex structure on a torus $T$ is the same as a complex structure 
on the tangent space $T_x T$, for some $x\in T$.
This implies that
a translation invariant complex structure on a torus $T = V_\R / V_\Z$ is the
same as the Hodge decomposition $V^{1,0} \oplus V^{0,1}=V_\R \otimes_\R \C$.

To pass from the torus $T$ to
the Hodge structure,\index[terms]{structure!Hodge} we consider the Hodge decomposition on 
$H^1(T, \C) = H^{1, 0}(T) \oplus H^{0,1}(T)$ and $H^1(T, \Z)\subset H^1(T, \C)$, 
defining a Hodge structure on $H^1(T, \C)=H^1(T, \R)\otimes_\R \C$. For the converse,
we start with a Hodge structure $V= V^{1,0} \oplus V^{0,1}$,
$V_\Z\subset V$, pass to the dual Hodge structure on the dual space
$W=V^*$, and take $T:= W^{1,0}/\pi(W_\Z)$, where
$\pi:\; W\arrow W^{1,0}$ denotes the projection along $W^{0,1}$.
It is easy to see that $T$ is compact and $H_1(T,\C)= W$, with the
Hodge decomposition and the integer lattice \index[terms]{lattice!integer}the same as
we started from. 

Indeed, $T= H_1(T, \R)/H_1(T,\Z)$;
this can be seen from the integration of parallel 1-forms over
the paths connecting a given point of $T$ to a fixed reference
point in $T$. After this isomorphism is fixed, the complex structure on $T$ 
is obtained from the projection $H_{1,0}(T) \tilde \arrow H_1(T,\R)$,
where $H_{1,0}(T)\subset H_1(T,\C)$ is the dual space to 
$H^{1,0}(T)\subset H^1(T,\C) =  H_1(T,\C)^*$. \endproof

\hfill

\definition
Let $M$ be a complex manifold.
A {\bf variation of Hodge structures}\index[terms]{structure!Hodge!variation of}
on $M$ is a complex vector bundle $(B, \nabla)$ with a flat connection
equipped with a parallel anti-complex involution \index[terms]{condition!Griffiths transversality}
and Hodge structures at each point, $B= \bigoplus_{p+q=w} B^{p,q}$
that satisfy the  \index[persons]{Griffiths, P.} Griffiths transversality condition:
\begin{equation}\label{_Griffiths_tranversality_Equation_}
\begin{aligned} 
\nabla^{1,0}(B^{p,q}) & \subset \big[B^{p,q} \oplus
  B^{p-1,q+1}\big]\otimes \Lambda^{1,0}(M), \\ \nabla^{0,1}(B^{p,q}) & \subset \big[B^{p,q} \oplus
  B^{p+1,q-1}\big]\otimes \Lambda^{1,0}(M).
\end{aligned}
\end{equation}

\remark 
\index[persons]{Griffiths, P.} Griffiths transversality conditions can
be deduced using \ref{_Gauss_Manin_via_Hattori_Remark_} below.

\hfill

\remark
Notice that this section completely omits the notion
of ``polarization'' and the ``polarized Hodge structures'',
that we did not even define. Most of the results about
the variations of Hodge structures\index[terms]{structure!Hodge} require the polarization
to work.

\hfill

\example\label{_VHS_Example_}
Let $\pi:\; M \arrow X$ be a proper holomorphic 
submersion with K\"ahler fibres. Consider the bundle
$V := R^k\pi_*(\C_M)$ with the fibre in $x$ the $k$-th cohomology
of $\pi^{-1}(x)$, the Hodge decomposition coming
from the complex structure on $\pi^{-1}(x)$,
and the Gauss--Manin connection. This defines
a variation of Hodge structures.\index[terms]{connection!Gauss--Manin}
The Griffiths transversality condition for this local system
immediately follows from the explicit description of the Gauss--Manin 
connection using the differential in the Hattori spectral sequence 
(\ref{_Gauss_Manin_via_Hattori_Remark_}).\index[terms]{spectral sequence!Hattori}

\hfill

We need the variations of Hodge structures\index[terms]{structure!Hodge!variation of}
for only one observation, which allows one
to interpret the variations of Hodge structures
of weight 1 as toric fibrations, and vice versa.

\hfill

\claim\label{_VHS_weight_1_torus_Claim_}
The integer variations of
Hodge structures of weight 1 are in a functorial bijective correspondence
with holomorphic toric fibrations with a holomorphic section. More precisely,
the variations of Hodge structures of weight 1 correspond
to fibrations in complex tori $X \arrow Y$ with a holomorphic
section, giving the fibrewise choice of zero in the torus.

\hfill

\proof
Let $V$ be an integer variation of Hodge structures\index[terms]{structure!Hodge!variation of} of weight one over $X$,
and $V^{1,0}\subset V_\Z\otimes_\Z \C$ the corresponding bundle of $(1,0)$-vectors.
By the  Griffiths transversality condition,\index[terms]{condition!Griffiths transversality}
$\nabla^{0,1}(V^{1,0})\subset \nabla^{0,1}(V^{1,0})$.
Since $\nabla^2=0$, one also has $(\nabla^{0,1})^2(V^{1,0})=0$,
hence $V^{1,0}$ is a holomorphic vector bundle by \index[persons]{Koszul, J.-L.} Koszul-\index[persons]{Malgrange, B.}Malgrange theorem
(\ref{_Koszul--Malgrange_Theorem_}). The quotient of\index[terms]{theorem!Koszul--Malgrange}
the fibre of $V^{1,0}$ by the image of $V_\Z$ is a complex torus.
Therefore, the fibration $\Tot((V^{1,0})^*/V_\Z^*)\arrow X$
is a smooth, holomorphic fibration with the fibre a complex torus
and a section provided by the zero section of $(V^{1,0})^*$.

Conversely, each torus fibration $\pi: M \arrow X$ 
with a section is equipped with a group structure, because its fibres
are identified with their Albanese varieties;\index[terms]{variety!Albanese} this construction
also provides the fibres of $\pi$ with a flat structure.
The corresponding Gauss--Manin local system of first cohomology
gives a variation of Hodge structures (\ref{_VHS_Example_});
the corresponding torus fibration is by construction identified\index[terms]{map!Albanese!relative}\index[terms]{local system!Gauss--Manin}
with the relative Albanese fibration of $\pi$ (see 
\ref{_Alb_pi_Definition_} below). 
\endproof

\section{Gromov's compactness theorem}

In the sequel, we will use the standard results on the degeneration
of complex curves. In projective geometry, this theorem is classical;
however, we apply it to non-K\"ahler surfaces. For our present purposes,
it is easier to invoke \index[persons]{Gromov, M.} Gromov's compactness theorem 
(\cite{_Gromov:pseudoholomorphic_}), which works
for any almost complex manifold. 

Let $M$ be a compact metric
space. Recall that on the set of closed subsets of $M$
one can define the Hausdorff distance (\ref{_Hausdorff_distance_Definition_}).
The topology it induces on the set of closed subsets is called
{\bf the \index[persons]{Gromov, M.} Gromov--Hausdorff topology}; it is independent\index[terms]{topology!Gromov--Hausdorff}
from the choice of the ambient metric inducing the same topology on $M$. 
The set of closed subsets of a compact metric space\index[terms]{space!metric}
$M$ is compact in Gromov--Hausdorff topology
(\cite{_Gromov:Metric-book_}).

\hfill

\definition
{\bf A pseudoholomorphic curve} on an almost complex\index[terms]{curve!pseudoholomorphic}
manifold $(M,I)$ is a closed compact subset $\Sigma \subset (M,I)$
that is  a smooth submanifold of dimension 2 outside of a finite subset,
and satisfies $I(T_x\Sigma) = T_x \Sigma$ for any smooth point $x\in \Sigma$.

\hfill

\index[persons]{Gromov, M.} Gromov's compactness theorem asserts that on an almost complex
symplectic manifold, the
set of pseudoholomorphic curves  is compact in the 
Gromov--Hausdorff topology.

\hfill

As follows from \cite[Theorem 4.2.1]{_McD_Salamon:J-holomo_}, 
$\Sigma$ admits a smooth, complex analytic pa\-ra\-me\-tri\-zation,
and can be interpreted as a ``cusp-curve'' in the sense of 
\cite[Chapter VIII]{_Audin_Lafontaine_}. 

\hfill

\definition
Let $(M,J)$ be an almost complex manifold.
A {\bf cusp-curve} is a collection $\{\Sigma_i\}$ of Riemann surfaces
with certain nodes $x_j\in \Sigma_i$ identified, and a $J$-holomorphic
map $\bigcup \Sigma_i \stackrel \phi \arrow (M,J)$ compatible with this identification
(that is, the nodes that are  identified are mapped to the same
points in $M$).

\hfill

A pseudoholomorphic curve $\Sigma$ is complex
analytic outside of the singular set $\Sigma_\sing$, and the complex
structure on $\Sigma$ can be extended to $\Sigma_\sing$ using the pseudoholomorphic
version of the Riemann removable
singularities theorem, \cite[Theorem 4.2.1]{_McD_Salamon:J-holomo_}.

\index[persons]{Gromov, M.} Gromov's compactness theorem is a theorem about the limits
of pseudoholomorphic curves. It turns out that the limit
is also a pseudoholomorphic curve, as long as the volume of the
curves is bounded. We refer to \cite[Chapter VIII, Theorem 2.2.1]{_Audin_Lafontaine_},
where this theorem is proven in much greater generality.

\hfill

\theorem\label{_Gromov_compactness_Theorem_}
\index[terms]{theorem!Gromov's compactness}
Let $(M, I, \omega)$ be a compact almost complex Hermitian manifold,
$V>0$ a constant, and $P_V$ the set of pseudoholomorphic curves $S\subset (M,I)$
with the Riemannian volume $\int_S \omega$ bounded by $V$. Then $P_V$ is compact
with respect to the \index[persons]{Gromov, M.} Gromov--Hausdorff topology. Moreover, the limit
of a sequence of diffeomorphic cusp-curves $\Sigma_i \subset M$
is a cusp-curve $C\stackrel \phi \arrow (M,I)$ that can be realized as follows.
Take a finite collection  of disjoint simple loops avoiding the nodes on a
cusp-curve $\Sigma$ diffeomorphic to $\Sigma_i$. Then $C$ is homeomorphic to
the cusp-curve obtained from $\Sigma$ by contracting each of these
loops to a point.
\endproof

\hfill

\corollary 
Let $S_i$ be irreducible components of the limit of a sequence of 
smooth pseudoholomorphic curves
 of genus $g$. Then $\sum_i g(S_i) < g$.

\proof \cite[Chapter VIII, Theorem 2.2.1]{_Audin_Lafontaine_}. \endproof

\hfill

A typical situation that can occur is the following.
Take a contractible simple loop $S$ \index[terms]{curve!rational}
in a smooth curve $\Sigma$, and several disjoint loops $S_i$ in 
the interior of $S$. Contracting $S$, we obtain a rational curve
$C_1=\C P^1$ attached to $S$ in the point of $\Sigma$ obtained
by contracting $S$. Contracting $S_i$, we attach a sequence
of rational curves to $C_1$. The incidence graph\index[terms]{incidence graph}
of the sequence of rational curves obtained by contracting
a nested family of disjoint simple loops in a plane\index[terms]{tree}
is a tree, that can be seen from the following picture.

\centerline{\includegraphics[width=0.75\linewidth]{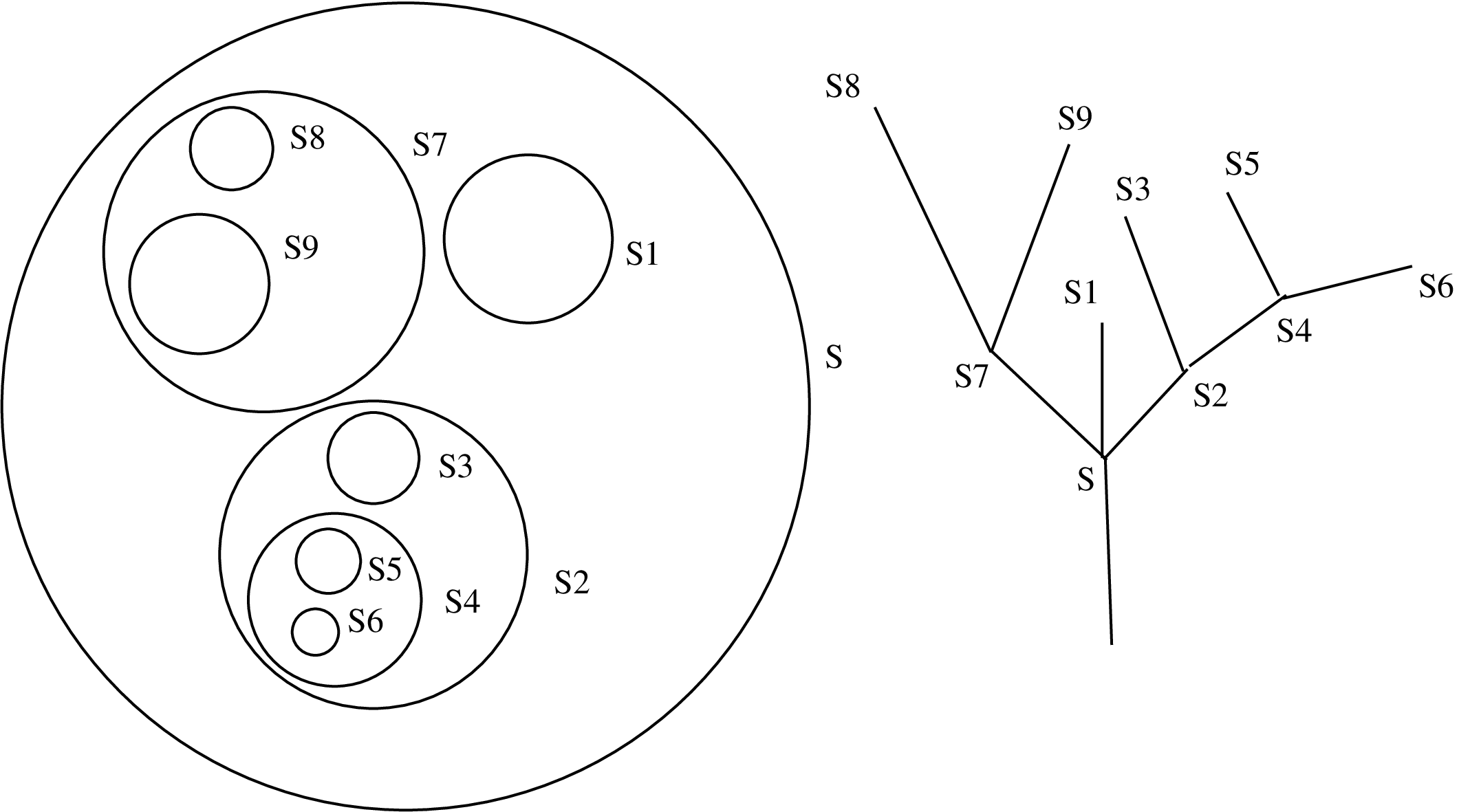}}

This geometric picture is called ``the tree bubbling phenomenon'';
for more details and ramifications of tree bubbling, see
\cite{_Parker_Wolfson:bubble_}.

\hfill

Let $C$ be a cusp-curve obtained as a limit of 
cusp-curves diffeomorphic to a smooth curve $\Sigma$, and $C_i$ its
smooth components, mapped to $(M,I)$. It might happen
that one of $C_i$ is contracted to a point. 

This explains the situation when 
the irreducible components of the pseudoholomorphic
curve obtained as a limit admit triple and higher intersections.
An example is found in \cite[Example V.3.5]{_Hummel:compactness_};
this triple intersection can be interpreted as an
intersection of three curves with a fourth one,
that is  contracted to a point.

\hfill

\centerline{\includegraphics[width=0.3\linewidth]{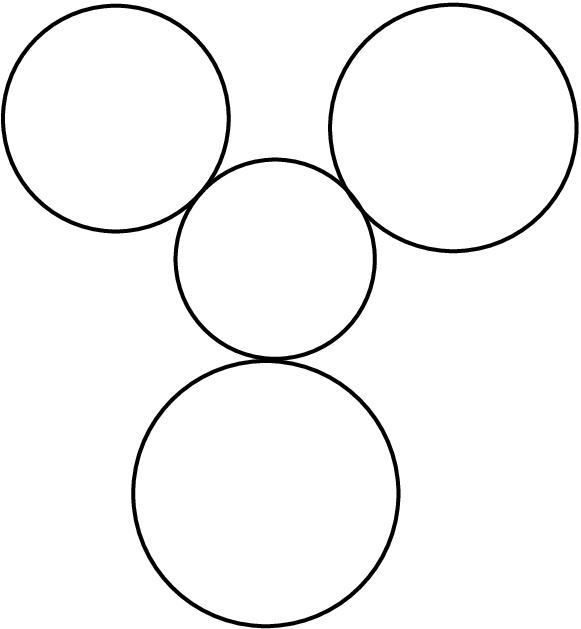}}

\hfill

\remark\label{_attach_tree_Remark_}
In the sequel, we describe attaching a tree-type sequence
of rational curves  to a curve $\Sigma$
by saying that we ``attach a tree of rational curves to $\Sigma$'',
or ``attach a bubble tree\index[terms]{tree!bubble} to $\Sigma$'', and the tree of rational
curves ``a bubble tree''.\index[terms]{tree!bubble}

\hfill

\centerline{\includegraphics[width=0.8\linewidth]{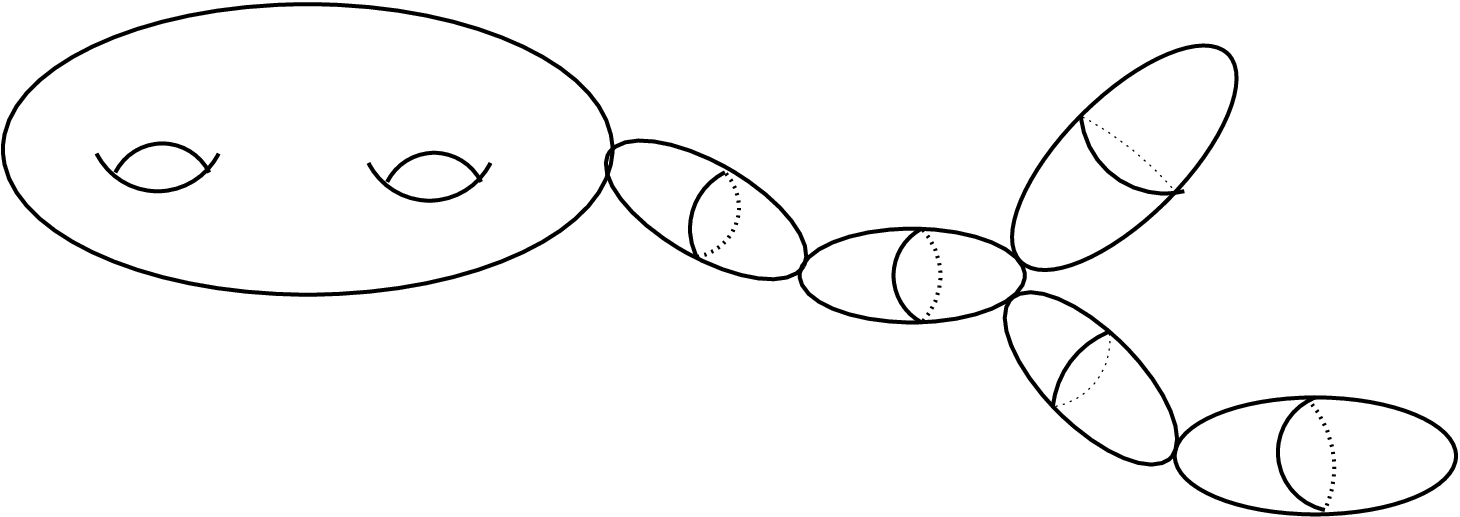}}

\hfill

Since we are going to study the elliptic surfaces,\index[terms]{surface!elliptic} 
we are most interested in the Gromov--Hausdorff limit of elliptic curves.\index[terms]{curve!elliptic}

\hfill

\corollary\label{_limits_of_elliptic_Corollary_}
Let $\Sigma$ be a Gromov--Hausdorff limit of a sequence
of smooth pseudoholomorphic curves of genus 1. Then $\Sigma$
is 
\begin{description}
\item[(i)] A curve of genus 1 with several bubble trees\index[terms]{tree!bubble} attached.

\centerline{\includegraphics[width=0.9\linewidth]{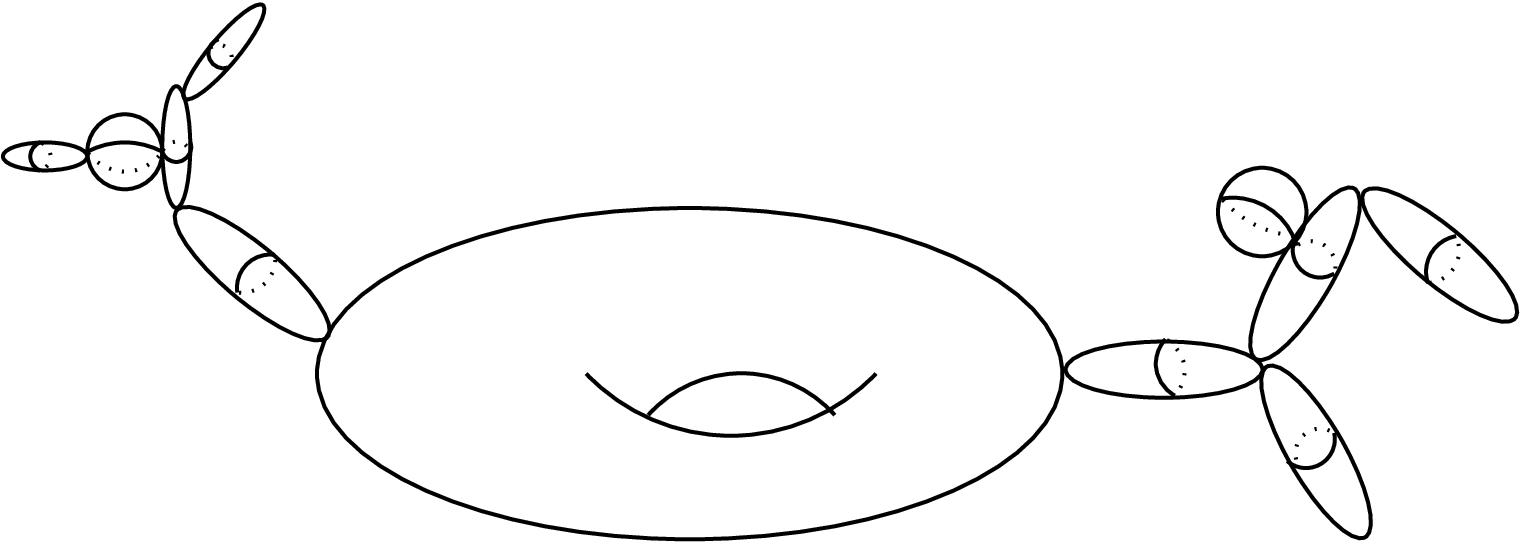}}
\item[(ii)] A  curve of genus 0 with a single
  self-intersecting node (``nodal rational curve'')\index[terms]{curve!rational!nodal}
with several bubble trees attached.

\centerline{\includegraphics[width=0.6\linewidth]{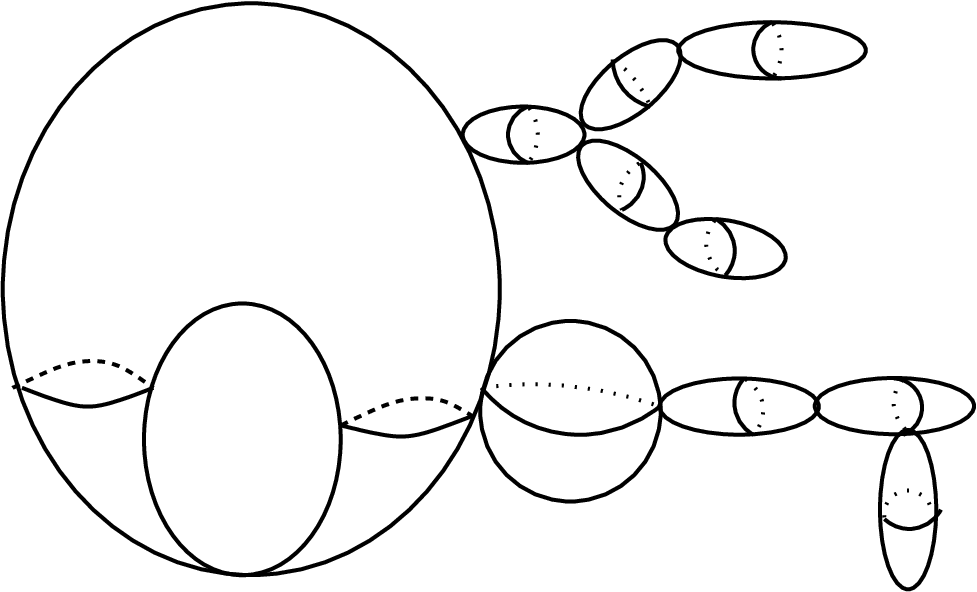}}
\item[(iii)] A sequence of curves of genus 0 arranged in a circle\footnote{This arrangement
of curves of genus zero is often called ``a wheel''.}  with  
several bubble trees attached.\index[terms]{tree!bubble}

\centerline{\includegraphics[width=0.6\linewidth]{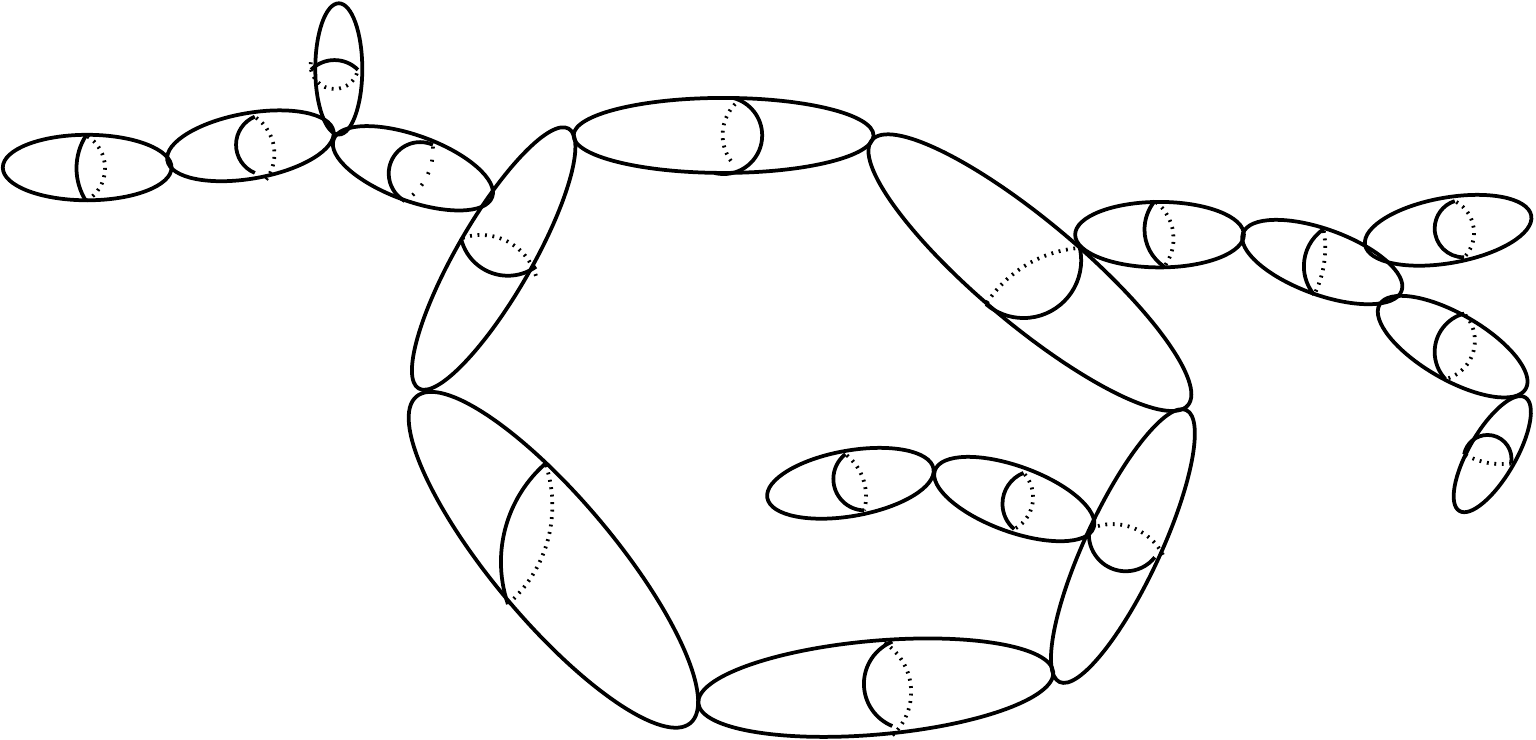}}
\end{description}

\proof
By \ref{_Gromov_compactness_Theorem_},
the curve $\Sigma$ is obtained from an elliptic curve $E$
by contracting a collection of disjoint simple loops $S_1, ..., S_n$.
If all $S_i$ are homotopic to zero in $E$, the contraction
gives an elliptic curve with several trees
of rational curves attached (\ref{_attach_tree_Remark_}).
If not all $S_i$ are homotopic to zero in $E$,
let us renumber $S_i$ in such a way that
$S_1, ..., S_d$ are not homotopic to zero, and
$S_{d+1}, ..., S_n$ are. Since $S_1, ..., S_d$ are
disjoint, they are homotopic to each other; indeed,
any two non-homotopic simple curves on a 2-torus intersect.
Replacing parallels by meridians if necessary,
we may assume that $S_1, ..., S_d$ is a collection
of meridians.

\centerline{\includegraphics[width=0.2\linewidth]{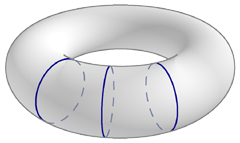}}

After contracting these loops, we obtain ``a circle of 
sausages'', for $d>1$ (a wheel of rational curves, arranged in a circle)

\centerline{\includegraphics[width=0.4\linewidth]{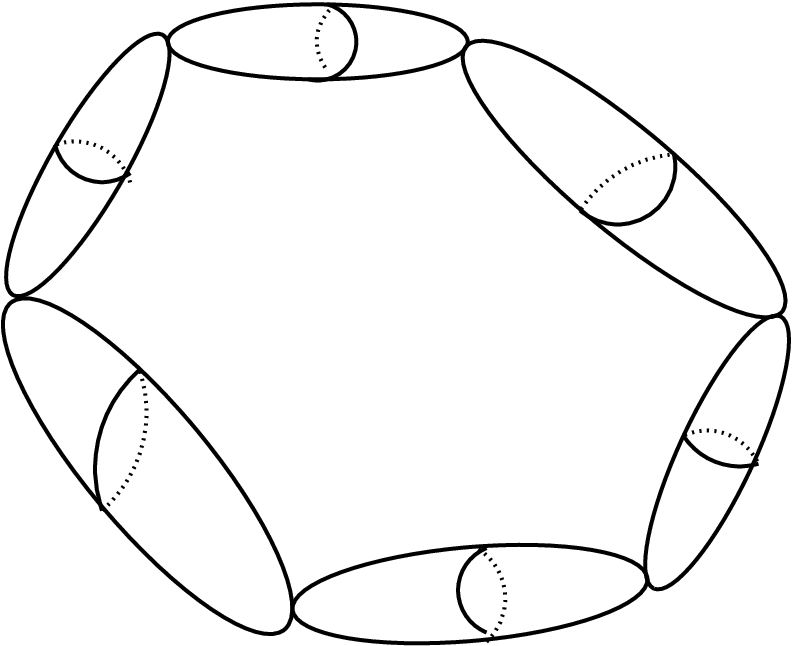}}

\noindent or a single rational curve with a self-intersection, if $d=1$.

\centerline{\includegraphics[width=0.3\linewidth]{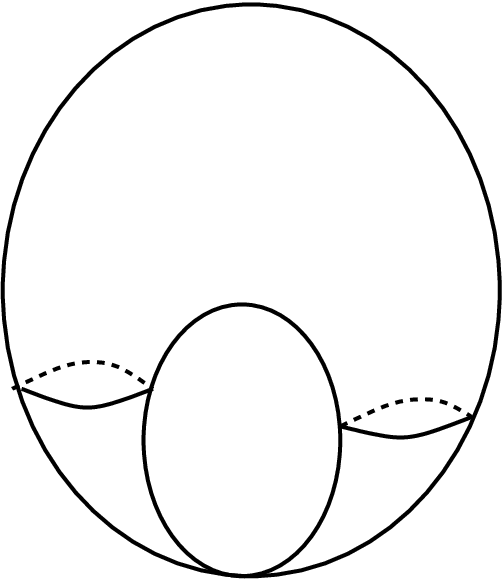}}

After attaching to these types of curves the bubble trees\index[terms]{tree!bubble}
associated with the rest of simple disjoint loops,
we arrive at the statement of
\ref{_limits_of_elliptic_Corollary_}.
\endproof


\section{Barlet spaces}\index[terms]{space!Barlet}
\label{_Barlet_Section_}

\index[persons]{Barlet, D.} Barlet spaces are spaces of cycles, that is, 
closed complex analytic subvarieties
of a given dimension in a given complex manifold
with multiplicities (positive integers)
assigned to their irreducible components.
They are similar but distinct from the \index[persons]{Douady, A.} Douady spaces, that are \index[terms]{space!Douady}
spaces of closed complex analytic subspaces (possibly with nilpotents
in the structure sheaf). 

For more details on Barlet spaces\index[terms]{space!Barlet} and their properties, please see the book
\cite{_Barlet_Magnusson:book_}.

Let $M$ be a metric space.\index[terms]{space!metric} Recall from \ref{_GH_limit_subgroups_Lemma_} that 
the (Gromov--)Hausdorff metric on \index[terms]{Gromov--Hausdorff distance}
the set ${\cal C}$ of closed subsets of $M$
is defined as follows: $d(X,Y)$ is the infimum
of all $\epsilon$ such that $X$ belongs to an
$\epsilon$-neighbourhood of $Y$, and
$Y$ belongs to an $\epsilon$-neighbourhood of $X$.
When $M$ is compact, the corresponding topology
on ${\cal C}$ is independent on  the choice
of metric on $M$ as long as the topology of $M$ 
remains the same. It is called 
{\bf the (Gromov-)Hausdorff topology}.

It is important for the sequel that the
\index[persons]{Barlet, D.} Barlet space of cycles is complex analytic, with
the topology that is  compatible with the (Gromov-)Hausdorff
topology on the set of all closed subvarieties.

For each $k \geq 0$, the Barlet space ${\goth B}_k(M)$ of $k$-cycles
on a manifold $M$ comes equipped with a universal family of cycles,
and their supports form a closed com\-plex-analytic subvariety  ${\goth B}^m_k(M) \subset
M \times {\goth B}(M)$, the ``marked Barlet space''
of pairs ``a complex analytic cycle and
a point in its support''. The
forgetful map ${\goth B}^m_k(M)\arrow {\goth B}_k(M)$ has all
fibres of the same dimension $k$, and the forgetful map $\Psi:{\goth B}^m_k(K)
\to M$ is complex-analytic. In particular, for any compact irreducible
component $Z$ in a \index[persons]{Barlet, D.} Barlet space, and its 
marked counterpart $Z^m$, the image $\Psi(Z^m)$ is a complex
analytic subvariety by Remmert's proper mapping\index[terms]{theorem!Remmert's proper mapping}
theorem (\cite[Theorem V.C.5]{_Gunning_Rossi_}, \cite[\S II.8.2]{demailly}). Geometrically,
this means that for any compact complex-analytic family $Z$ of
cycles in $M$, the union of supports of all these cycles
is a complex subvariety in $M$.

Further on, we shall use the following result, which
was proven for curves in \cite[Corollary 2.19]{_Verbitsky:twistor_},
and for divisors in \cite{_Barlet:divisors_}.\index[terms]{divisor}

\hfill

\theorem\label{_compact_Barlet_Theorem_}
Let $(M,I)$ be a compact complex manifold, and
$\omega$ a Hermitian form.\index[terms]{form!Hermitian} Assume that
$dd^c(\omega^k)=0$ for some integer $k>0$. Then any connected component
$Z$ of the  Barlet space ${\mathfrak B}_k(M)$ 
of $k$-cycles in $M$ is compact.

\smallskip

\proof
By Bishop's theorem (\cite{_Bishop:conditions_}),
a (Gromov--)Hausdorff limit of a family of compact complex
subvarieties is complex if its Hermitian volume
stays bounded. Since the space of closed
subvarieties in $M$ with (Gromov--)Hausdorff topology is compact,
this implies that the set of closed 
compact $k$-dimensional subvarieties
with volume bounded by a constant $C\in \R$
is compact.

Therefore, to prove compactness, it would suffice to show
that the volume $\Vol(S)$ is constant as a function of 
$[S]\in Z$.

Let $X$ denote the marked family associated with $Z$,
and $\pi_M:\; X \arrow M$, $\pi_X:\; X \arrow Z$
the forgetful maps. Then the volume function
$\Vol:\; Z\arrow \R^{>0}$ can be expressed as
$\Vol = (\pi_X)_*\pi_M^* \omega^k$, where
$(\pi_X)_*$ denotes the pushforward of a
differential form (generally speaking, the
pushforward is not a form, but it is well
defined as a current).\index[terms]{current}

Let $k$ be the dimension of the cycles parametrized by $Z$.
Since pullback and pushforward of differential forms
commute with $d$, $d^c$, this gives 
$$dd^c
\Vol=(\pi_X)_*\pi_M^* (dd^c\omega^k),$$ see for example 
\cite[(8.12)]{_NHYM_}, \cite[Theorem 2.10]{_Verbitsky:S^6_}
or \cite[Proposition 1.9]{_Ivashkovich:Annals_}.
Therefore, $\Vol$ satisfies $dd^c(\Vol)=0$
whenever $dd^c(\omega^k)=0$.

Functions that satisfy $dd^c f=0$ are called\index[terms]{function!pluriharmonic}
{\bf pluriharmonic}. Using the local $dd^c$-lemma,\index[terms]{lemma!$dd^c$} it
is easy to see that any pluriharmonic function
is locally the sum of a holomorphic and an antiholomorphic
function.

By Bishop's theorem, the set\index[terms]{theorem!Bishop}
$\Vol^{-1}(]-\infty, C])$ is compact for all $C\in \R$,
hence  $-\Vol$ has a maximum somewhere in $X$.
However, a pluriharmonic function that admits a 
maximum  is necessarily
constant by E. Hopf's strong maximum principle (\ref{hopf_theorem}).\index[terms]{maximum principle}
Therefore, $\Vol$ is constant on each
connected component of the Barlet space.
Now, each of these components is compact by Bishop's theorem.
\endproof

\hfill

Applying this result to the
space of curves on a complex surface
and using the Gauduchon form\index[terms]{form!Gauduchon} $\omega$
(Subsection \ref{_BC_degree_Subsection_}),
we obtain the following useful corollary.

\hfill

\corollary\label{_Barlet_on_surf_Corollary_}
Let $M$ be a compact complex surface,
${\goth B}_1(M)$ the  Barlet space\index[terms]{space!Barlet} of 1-cycles on $M$
and $Z$ its connected component. Then $Z$ is compact.
\endproof


\section{Elliptic fibrations with multiple fibres}
\label{_multiple_fibres_Subsection_}

\subsection{Multiple fibres of elliptic fibrations and the
  relative Al\-ba\-ne\-se map}
\label{_relative_Albanese_Subsection_}

Further on, we need the following general observation, which
is much used in algebraic geometry.\index[terms]{geometry!algebraic}

\hfill

\claim\label{_general_fiber_retracted_to_special_Claim_}
Let $\pi:\; M \arrow S$ be a proper, flat morphism
of complex varieties, with $M$ smooth, and $s\in S$ 
a point. Let  $C_s = \pi^{-1}(s)$ be the special fibre of $\pi$
Then there exists a small neighbourhood $U \ni s$ 
such that $M_U = \pi^{-1}(U)$ admits a retraction onto $C_s$.

\proof
\cite{_Persson:degene_}, \cite{_Clemens:degene_},
see also  \cite{_Morrison:Clemens_Schmid_}.
\endproof

\hfill

\definition\label{_Alb_pi_Definition_}
Let $U\ni s$ be a disk in $S$ such
that $M_U = \pi^{-1}(U)$ is homotopy equivalent to $C_s$.
Assume that $\pi:\; M_U \arrow U$ has a holomorphic section
$\tau:\; U \arrow M_U$. Let $R$ be the space
of 1-forms on $M_U$ that are  holomorphic on the fibres of $\pi$,
and  $C_x= \pi^{-1}(x)$  a fibre of $\pi$ in $x$.
The {\bf relative Albanese map}\index[terms]{map!Albanese!relative}
$\Alb_\pi:\; C_x\arrow \frac{H^{1,0}(C_x)}{H_1(C_x, \Z)}$
takes a fibrewise holomorphic 1-form $\eta$
to $\int_\gamma \eta$, where $\gamma$ is a path 
in $C_x$ connecting $y\in C_x$ to $\tau(x)\in C_x$.

\hfill

\remark\label{_rel_alb_mono_trivial_Remark_}
Let $\pi:\; M_U \arrow U$ 
be an elliptic fibration over a disk which 
admits a section $\tau$. Assume that the 
monodromy\index[terms]{monodromy} of the Gauss--Manin local system \index[terms]{local system!Gauss--Manin}
is trivial, and there exists a fibrewise holomorphic
1-form $\eta^{1,0}\in \Lambda^1(M)$ giving a parallel section of the \index[terms]{section!parallel}
Gauss--Manin local system\index[terms]{local system!Gauss--Manin} on $U\backslash s$. 
Then the relative Albanese map\index[terms]{map!Albanese!relative} restricted to the space
$A:=\langle \eta^{1,0}\rangle\subset H^{1,0}(C_x)$
defines a map $C_x \arrow A^*/\Lambda_x$,
where $\Lambda_x$ denotes the subgroup generated by
$a\mapsto \int_{t\in H^1(C_x, \Z)} a$.
The lattice\index[terms]{lattice} $\Lambda_x\subset  A^*$ is the same
for all $x\in U$, because the Gauss--Manin\index[terms]{local system!Gauss--Manin}
local system $V$ associated with the first
cohomology of the fibres has no monodromy\index[terms]{monodromy}, and the Hodge
decomposition on $V$ is constant.
Therefore, $\Alb_\pi:\;  M_U \arrow\frac{H^{1,0}(C_x)}{H_1(C_x, \Z)} $
takes a fibre of $\pi$ to the elliptic curve 
$A^*/\Lambda=\frac{H^{1,0}(C_x)}{H_1(C_x, \Z)}$ that is  independent on  $x$. 
Let $\widetilde \Alb_\pi:\; M_u \arrow U \times (A^*/\Lambda)$ 
take $x\in M$ to the pair $(\pi(x), \xi)$, where $\xi= \Alb_\pi(x)$
is an element of $A^*/\Lambda=\frac{H^{1,0}(C_x)}{H_1(C_x, \Z)}$ associated with
$x$ as above.

\hfill

\proposition\label{_Albanese_around_special_fiber_Proposition_}
Let $\pi:\; M \arrow S$ be an elliptic fibration
and $s\in S$ a singular value of $\pi$. 
Consider  a disk $U\ni s$ in $S$ such that $M_U = \pi^{-1}(U)$ 
contains no other critical values of $\pi$. 
Assume that the 
monodromy \index[terms]{monodromy}of the Gauss--Manin local system around $s$\index[terms]{local system!Gauss--Manin}
is trivial, and there exists a fibrewise holomorphic
1-form $\eta\in \Lambda^{1,0}(M)$ giving a parallel section of the \index[terms]{section!parallel}
Gauss--Manin local system on $U\backslash s$. 
\begin{description}
\item[(i)] Then the 
special fibre of $C_s$ is an elliptic curve with 
several bubble trees attached. \index[terms]{tree!bubble}
\item[(ii)] If, moreover,
$\pi:\; M_U \arrow U$ admits a holomorphic
section $\tau$, then the map 
$M_U \xlongrightarrow{\widetilde \Alb_\pi} U \times (A^*/\Lambda)$
(\ref{_rel_alb_mono_trivial_Remark_})
blows down all bubble trees.\index[terms]{tree!bubble}
\end{description}

\smallskip

\proof
We start by showing that $\pi_1(M_U) = \Z^2$.

Clearly, every
special fibre is a limit of general fibres. By
\ref{_limits_of_elliptic_Corollary_},
$C_s$ is an elliptic curve, a nodal rational\index[terms]{curve!rational!nodal}
curve, or a wheel of rational curves, with
several bubble trees attached.\index[terms]{tree!bubble}

Suppose that $C_s$ is a nodal rational\index[terms]{curve!rational!nodal}
curve or a chain of rational curves with bubble trees attached. Then the
normalization $\tilde C_s$ of $C_s$ is a union \index[terms]{normalization}
of several copies of $\C P^1$. The pullback of
the fibrewise holomorphic 1-form $\theta^{1,0}$ has to vanish identically 
on each of these $\C P^1$, because the sheaf of
holomorphic vector fields on $\C P^1$
is globally generated. This implies that
$\theta^{1,0}\restrict{C_s}=0$, and hence  $\theta^{1,0}$
is exact on $C_s$. 

By \ref{_general_fiber_retracted_to_special_Claim_},
$C_s$ is a deformation retract of $\pi^{-1}(U)$. 
If the form $\theta^{1,0}$ is exact on $C_s$, it is
exact on $\pi^{-1}(U)$, that contains a copy
of a general fibre of $\pi$. Then $\theta^{1,0}$ is
exact on a general fibre of $\pi$, which
is impossible. Therefore,
$C_s$ contains an elliptic curve $E$.

The relative Albanese map\index[terms]{map!Albanese!relative} blows down all bubble trees\index[terms]{tree!bubble}
because the fibrewise holomorphic 1-form $\theta^{1,0}$
vanishes on all rational curves.
\endproof 

\subsection{Structure of a neighbourhood of a multiple fiber}
\label{_mult_fibres_Subsection_}

\definition
Let $M \arrow S$ be an elliptic surface,
that is, a proper holomorphic map with general fibres 1-dimensional
complex tori. Let $C_s$ be a special fiber, and
$x\in C_s$ a point of its smooth locus.
Denote by $C$ the irreducible component
of $C_s$ containing $x$. Let $V\subset M$
be a holomorphic disk transversal to $C$ in $x$. We say that
$C_s$ {\bf has multiplicity $d$ in $C$} if $V$ intersects 
each general fibre that is  sufficiently close to $C$ in 
precisely $d$ points.

\hfill

\example
Let $E_1\times E_2$ be a product of two elliptic
curves, $\epsilon \in E_2$ a $2$-torsion point
and $\tau$ an involution taking $(x, y) \in E_1\times E_2$
to $(-x, y+\epsilon)$. This action is clearly free.
The quotient $M:= \frac{E_1\times E_2}{\langle \tau\rangle}$ 
is called {\bf the bielliptic surface}.\index[terms]{curve!elliptic}
By Exercise \ref{_E/+-_biholo_to_CP^1_Exercise_}, \index[terms]{surface!bielliptic}
$S:= E_1/\{\pm 1\}$ is biholomorpic to $\C P^1$.
The natural projection $\pi:\; M \arrow S$
is an elliptic fibration with multiple fibres of multiplicity 2
in the ramification points of $E_1\arrow  E_1/\{\pm 1\}=S$.
Indeed, let $D$ be a holomorphic disk transversal
to a special fibre $\frac{(\{0\}\times E_2)}{\langle \tau\rangle}$
of $\pi$, and $D'$ its preimage
in $E_1\times E_2$. Since $E_1\times E_2\arrow M$
is a 2-sheeted cover, $D'$ is the union of two copies
of a holomorphic disk, intersecting
$\{0\}\times E_2$ in $(0, t)$ and $(0, t+\epsilon)$;
each of these disks intersects the neighbouring
fibres of the projection $E_1 \times E_2 \arrow E_1$
once,  hence $D$ intersects each neighbouring 
fibre twice.

For an example of a Hopf surface with\index[terms]{surface!Hopf}
multiple fibres, see \ref{kodaira.rem}.

\hfill

Let $\pi:\; M \arrow S$ be an elliptic fibration,
and $C_s$ a singular fibre in $s\in S$. Choose a neighbourhood
$U$ of $s\in S$, and let $M_U:= \pi^{-1}(U)$. We assume that
$\pi\restrict{M_U}$ is smooth outside of $C_s$.

\hfill

\definition\label{_local_multisection_Definition_}
A {\bf local multisection of\index[terms]{multisection}
$\pi$ ramified in $s$} is a smooth holomorphic disk
$U'\subset M_U$ intersecting $C_s$ transversally
in a smooth point of $C_s$, such that the
restriction $\pi\restrict {U'} :\; U' \arrow U$
is proper. Let $d$ be the multiplicity of 
$C_s$ in the irreducible component containing
the point $U'\cap C_s$. 

\hfill 

The map 
$\pi\restrict {U'} :\; U' \arrow U$
is a finite, non-ramified covering outside of $s$,
because $\pi:\; M_U \arrow U$ is smooth
outside of $C_s$. Then 
$\pi\restrict {U'} :\; U' \arrow U$
is a $d$-sheeted finite covering 
ramified in $s$ and non-ramified
outside of $s$. A multisection can be
constructed as follows. A smooth
holomorphic disk meeting $C_s$ transversally
in its smooth point can be always constructed;
replacing $U$ with a smaller neighbourhood, we
can make sure that 
the restriction $\pi\restrict {U'} :\; U' \arrow U$
is proper.

\hfill

Suppose that $E\subset C_s$
has multiplicity $d$. Choose a local multisection
$\tau:\; U'\arrow M_U$ of $\pi:\; M_U \arrow U$ ramified in $s$ as in
\ref{_local_multisection_Definition_}.
We are going to define a ramified covering $\pi':\; M_U'\arrow U'$,
where $\sigma:\; U'\arrow U$ is the $d$-sheeted covering
of the disk $U$ associated with the multisection $\tau$,
equipped with an action of $\Z/d\Z$, in such a way that
$\frac{M_U'}{\Z/d\Z}=M_U$.

Let $M_U':= U' \times_U M_U$.
This is an elliptic fibration over $U'$ 
equipped with a section $\tau':\; U'\arrow M_U'$
associated with $\tau$ in a natural way.
The fibres of $\pi'$ are identified with the corresponding
fibres of $\pi$, because $M_U':= U' \times_U M_U$.

%

\hfill

The following proposition is a reformulation of definitions.

\hfill

\proposition\label{_Albanese_around_special_multiplicities_Proposition_}
Let $\pi:\; M_U \arrow U$ be an elliptic
fibration over a disk $U$ with a single singular fiber
$C_s$ and trivial monodromy\index[terms]{monodromy!of a connection} of the Gauss--Manin connection.\index[terms]{connection!Gauss--Manin}
Let $E$ its smooth component of multiplicity $d$,
$\tau:\; U' \arrow M_U$ the corresponding 
multisection ramified in $s$,
and $\pi':\; M_U' \arrow U'$ the ramified cover
of $M_U$ constructed from this multisection as above, 
with $M_U= \frac{M_U'}{\Z/d\Z}$.\index[terms]{multisection}
Consider the relative Albanese map \index[terms]{map!Albanese!relative}
(\ref{_Albanese_around_special_fiber_Proposition_})
$M_U' \stackrel{\widetilde\Alb_{\pi'}}\arrow U' \times A^*/\Lambda$.
Then, after taking the $\Z/d\Z$-quotient, this will give the
holomorphic map
$M_U \stackrel{\widetilde\Alb_{\pi}}\arrow U \times \frac{A^*/\Lambda}{\Z/d\Z}$,
also called {\bf the relative Albanese map}.
\endproof

\hfill

\remark Unlike the relative Albanese map
defined in \ref{_Albanese_around_special_fiber_Proposition_},
the map in \ref{_Albanese_around_special_multiplicities_Proposition_}
is not injective on the fibres; indeed, on each fibre of 
$\pi$ this map is a non-ramified $d$-sheeted covering.

\hfill

\remark
There is another interesting way of constructing the
relative Albanese map defined in\index[terms]{map!Albanese!relative} \ref{_Albanese_around_special_multiplicities_Proposition_}.
Let $\tau:\; U' \arrow M_U$ be a multisection,
$C_x$ a general fiber, and $\{x_1, ..., x_d\}\subset  C_x$ 
the intersection of $C_x$ and $\tau(U')$.
For each $y\in C_x$, let $\gamma_1, ..., \gamma_d$ be paths
connecting $y$ to $x_1, ..., x_d$, and let $\Alb_\pi$ take
$y$ to $\sum_{i=1}^d \int_{\gamma_i}\theta^{1,0}$.
This sum is well-defined up to addition of the periods
of $\theta^{1,0}$, that is, the lattice\index[terms]{lattice} $\int_{H^1(C_x, \Z)}\theta^{1,0}\subset A^*$.
The dependence of each $\int_{\gamma_i}\theta^{1,0}$
from the point $y$ is not holomorphic, because
$\tau$ is a multisection, and we cannot separate the
branches; however, the sum $\sum_{i=1}^d \int_{\gamma_i}\theta^{1,0}$
depends holomorphically on $y$. It is easy to see that
this operation gives the same map as 
\ref{_Albanese_around_special_multiplicities_Proposition_}.

\section{Non-K\"ahler elliptic surfaces}


\subsection{Structure of elliptic fibrations on  non-K\"ahler surfaces}

For a version of the following theorem (with a different proof),
see \cite[Theorem 3.17]{_Brinzanescu:bundles_} and \cite{_Brinzanescu:manuscripta_}.

\hfill

\theorem \label{_elli_then_isotri_Theorem_} 
Let $M$ be a non-K\"ahler compact complex surface
admitting an elliptic fibration $\pi:\; M \arrow S$. 
Assume that the pullback $\pi^*\Vol_S$ of the volume form
is equal to the differential of a $d^c$-closed
form.\footnote{This assumption is always satisfied
  when $M$ is non-K\"ahler, by
  \ref{_elli_bundle_Theorem_}.}
Then the general fibres of  $\pi$ are isomorphic elliptic curves.
If, moreover, $M$ is minimal, then all fibres of
$\pi$ are elliptic curves.\index[terms]{curve!elliptic}

\hfill

\pstep
Let  $\Theta:= \pi^*\Vol S$ be the pullback of the volume form,
and $\theta$ a closed 1-form on $M$ that satisfies
$d^c\theta = \Theta$. 
Denote by $S_0\subset S$ the set of regular values of $\pi$,
and let $C\subset \pi^{-1}(S_0)$ be a general fibre of $\pi$.
We are going to show that $\theta \restrict C$ is non-exact.

Assume, by absurd, that $\theta \restrict C$ is exact.

Let $\chi:\; \pi_1(M) \arrow \R^*$ be the multiplicative
character associated with the 1-form $\theta$; this character
can be understood as the monodromy\index[terms]{monodromy} of the flat line bundle $L$
with the connection form $\theta$. We need to show that 
the restriction $\chi\restrict{\pi_1(C)}$ is non-trivial.
Let $M_0:=  \pi^{-1}(S_0)$. 

Clearly, we may consider $\chi$ as a character on $\pi_1(M_0)$.
Since $S_0$ is the smooth locus of $\pi$, the natural map
$\pi:\;M_0 \arrow S_0$ is locally trivial by Ehresmann theorem.
Consider the \index[terms]{exact sequence!Serre} Serre exact sequence
\[
0 \arrow \pi_1(C) \arrow \pi_1(M_0) \arrow \pi_1(S_0) \arrow \pi_2(C)=0.
\]
Since $\theta \restrict C$ is exact, the character $\chi$ vanishes on
$\pi_1(C)$. This implies that the local system $L$ is trivial
on a general fibre $C$ of $\pi$, and $L\restrict {M_0}$
is the pullback of a local system on $S_0$.


Since $\pi_1(C_s)$ is a quotient of $\pi_1(C)$ and  
the local system $L$ is trivial on the general
fibre $C$, it is trivial on special fibres as well.
This implies that $L$ is the pullback of a local system on $S$,
and $\theta$ is cohomologous to a 1-form $\pi^*(\theta_0)$.
This is impossible, because 
\[ 
\deg_\omega(d^c(\theta)):=\int_M d^c(\theta) \wedge \omega=\int_M \Theta \wedge \omega  >0,
\]
where $\omega$ is the Gauduchon form.\index[terms]{form!Gauduchon}
However, the 2-form $d^c\theta_0$ is exact on $S$ and hence it is 
Bott--Chern cohomologous to 0, which brings 
$\int_M d^c\pi^*(\theta_0) \wedge \omega=\int_M d^c(\theta) \wedge \omega=0$
(as shown in Subsection \ref{_BC_degree_Subsection_} below, 
the degree $\deg_\omega$ vanishes for any form that is  Bott--Chern
cohomologous to zero).

\hfill

{\bf Step 2:}
Return to the map $\pi:\; M \arrow S$ constructed above.
The smooth fibres of $\pi$ are all elliptic curves by Step 1.
At this step we are going to prove that all its smooth fibres are isomorphic
as complex curves.

Let again $S_0\subset S$ be the smooth locus of $\pi$,
and let $H$ be the Gauss-\index[persons]{Manin, Yu. I.}Manin bundle
associated with the first cohomology of the fibres of $\pi$.
By Step 1, $\theta$ is not exact on the general fibres of $\pi$.
Since $\theta$ is closed, 
the section of $H$ represented by $\theta$ is
constant with respect to the Gauss--Manin connection.

The form $d^c\theta= \Theta$ belongs to $\pi^*\Lambda^2 B$.
Then $I(\theta)$ is closed on the fibres of $\pi$. By \ref{_Gauss_Manin_d_1_vanishes_Lemma_},
$I(\theta)$, and hence the Hodge components of $\theta$,
define a parallel section\index[terms]{section!parallel} of the Gauss--Manin 
bundle.\footnote{The same conclusion is implied 
	by \index[persons]{Schmid, W.} Schmid's fixed part theorem
	(\cite[Theorem 7.22]{_Schmid:singularities_}):
	for any polarizable variation of Hodge structures\index[terms]{structure!Hodge!variation of} over
	a quasi-projective base,
	the (1,0) and (0,1)-parts of a parallel section $\theta$ are also
	constant with respect to the Gauss--Manin connection.}

We obtain a basis $\theta^{1,0}$ and $\theta^{0,1}$
for $H$. Therefore, the variation of Hodge structures \index[terms]{structure!Hodge!variation of}
induced by periods of the elliptic curves on
$H$ is trivial. This implies that the
corresponding elliptic fibration is {\bf isotrivial},
that is, all its smooth fibres are isomorphic as complex curves.%
\footnote{Here we implicitly use the equivalence between the
torus fibrations admitting a holomorphic section and the variations
of Hodge structures, given in \ref{_VHS_weight_1_torus_Claim_}.}

Note that isotriviality implies the local triviality
by a theorem of \index[persons]{Grauert, H.} Grauert and Fischer (\cite{fg}), and hence 
the fibration $\pi$ is locally trivial over $S_0$.
We will give an independent proof of this
observation in \ref{_isotrivial_Proposition_}.

\hfill

{\bf Step 3:} 
By \ref{_Albanese_around_special_fiber_Proposition_}, the
special fibre of $\pi$ contains
an elliptic curve,  possibly with some bubble trees
attached. We need to show that the bubble trees
are empty if $M$ is minimal.\index[terms]{tree!bubble}

Let $d$ be the multiplicity of the elliptic curve component
of the special fiber.
\ref{_Albanese_around_special_multiplicities_Proposition_}
implies that the relative Albanese map\index[terms]{map!Albanese!relative}
$M_U \arrow U\times \frac{A^*/\Lambda}{\Z/d\Z}$ is a holomorphic
map to a smooth manifold which blows down all bubble
trees, and hence  $M_U$ cannot be minimal.
\endproof

\subsection{Isotrivial elliptic fibrations}

Let $\pi:\; M \arrow S$ be an elliptic fibration, 
and $m$ a point on a fibre $E$ of $\pi$. Since $E$ is an
elliptic curve, we represent $E$ as a quotient
$E= \C/\Lambda$, with $\Lambda =\Z^2$. However,
there is no natural choice of the origin point in $E$.
Let $\Aut_0(E)$ denote the group of translations of $E$;
this is the same elliptic curve, with a fixed choice\index[terms]{curve!elliptic}
of 0. The curves $\Aut_0(\pi^{-1}(x))$ form 
a smooth holomorphic family over the smooth locus 
$S_0 \subset S$ of $\pi$, called {\em the relative Albanese 
fibration} of $\pi:\; M \arrow S$\index[terms]{map!Albanese!relative} (\ref{_Alb_pi_Definition_}). It is fibrewise
isomorphic to $\pi$, but, unlike $\pi$,
the relative Albanese fibration is equipped
with a natural holomorphic section.

\hfill

%
%
%

Note that the assumptions of the following
proposition are contained in the conclusion of 
\ref{_elli_then_isotri_Theorem_}, and hence  
\ref{_isotrivial_Proposition_} can be applied in the same generality.

\hfill

\proposition\label{_isotrivial_Proposition_}
\index[terms]{elliptic fibration!isotrivial}
Let $M$ be a  non-K\"ahler compact complex 
surface, admitting an elliptic fibration $\pi:\; M \arrow S$. 
Assume that the pullback $\pi^*\Vol_S$ of the volume form
is equal to the differential of a $d^c$-closed
form, and all special fibres of $\pi$ are 
elliptic curves. Denote by $E$ a general fibre of
$M$; we use the same letter $E$ for the complex Lie group
associated with $E$ as an elliptic curve.
Then there exists a natural holomorphic action
of the group $E$ on $M$, transitive on each fibre and free on 
any non-multiple fibre. 

\hfill

\proof
The variation of Hodge structures\index[terms]{structure!Hodge!variation of} associated with $\pi$ is
trivial, as shown in \ref{_elli_then_isotri_Theorem_},
Step 1; it is trivialized by the form $\theta^{1,0}$,
that is  closed on general fibres of $\pi$. 
Let $M_0 \arrow S_0$ be the 
smooth locus of $\pi$, and
$\Alb_\pi(M_0)$ the relative Albanese fibration,
that is, the bundle of fibrewise translations of $M_0$.
Since the variation of Hodge structures associated
with $\pi$ is trivial, the family $\Alb_\pi(M_0)$ 
is also trivial. 


Since $\Alb_\pi(M_0)= S_0\times E$, where $E= \Aut_0(\pi^{-1}(x))$ is an elliptic curve, 
the family $\pi:\; M_0 \arrow S_0$ is a principal
elliptic bundle, and $E$ acts on $M_0$ freely
and transitively on the fibres. To prove 
\ref{_isotrivial_Proposition_} it remains
to show that this action can be holomorphically extended
to the special fibres; since the special fibres are limits
of the general ones, this action will be transitive.

Let $M_U:= \pi^{-1}(U)$ be a sufficiently small
neighbourhood of the multiple fibre $C_s$ of $\pi$.
By \ref{_Albanese_around_special_multiplicities_Proposition_},
after passing to the quotient $E\arrow \frac E{\Z/d\Z}$,
the relative Albanese map\index[terms]{map!Albanese!relative} $M_0\arrow E$ can be extended
to the map $M_U\arrow \frac E{\Z/d\Z}$, inducing an isomorphism of
the singular fibre and $\frac E{\Z/d\Z}$. 

To finish the proof, we need to show  that the
vector fields associated with the action of $E$ on $M_0$
can also be extended to the central fibre $C_s$.
After taking the finite quotient, the fibration becomes
trivial, and the vector fields that we want to extend become smooth.
This implies that these vector fields were bounded with respect to some
smooth Hermitian metric on the tangent bundle. Now, the Riemann
extension theorem implies that these vector fields
admit a holomorphic extension from $M_0$ to $M$.
\endproof

\hfill

\definition
Let $M$ be a complex manifold equipped with
an action of a compact complex torus $T$, 
with the orbits of the same dimension. Assume that
the quotient map $\pi:\; M \arrow S$ is well-defined. Then 
$\pi:\; M \arrow S$ is called {\bf a principal
	toric fibration}, and {\bf a principal
	elliptic fibration} if $\dim_\C T =1$.

\hfill

\remark
The ``isotrivial elliptic fibration''
is a fibration with general fibres isomorphic
to the same elliptic curve. This is a more general notion.
In \cite{_Brinzanescu:manuscripta_,_Brinzanescu:bundles_}, 
an isotrivial elliptic fibration is 
called ``a quasi-bundle''.\index[terms]{quasi-bundle}

\hfill

\remark\label{_orbifold_Remark_}
Let $\pi:\; M \arrow S$ 
be a principal elliptic fibration over a curve $S$, with
fibre $C$. We associate with $\pi$ the following orbifold \index[terms]{orbifold}
structure on $S$. For any multiple fibre $R=\pi^{-1}(s)$ of $\pi$, consider
the group $\Gamma_R\subset C$, obtained as the kernel of the
natural action of $C$ on $R$. Taking the quotient
$M/\Gamma_R$, we obtain another fibration, that is 
locally trivial in a neighbourhood of $s$. A smooth section
of this fibration in a neighbourhood of $s$ gives 
a $\Gamma_R$-invariant multisection\index[terms]{multisection} $\tilde U \arrow U$ 
of $\pi:\; \pi^{-1}(U) \arrow U$. Then $U$ is obtained
as a finite quotient of $\tilde U$, and the lift
of $\pi$ to $\tilde U$ is locally trivial.

\hfill


\remark
We obtained that any principal 
elliptic bundle $\pi:\; M \arrow S$  over a curve 
defines an orbifold structure on this curve,
and $\pi$ is locally trivial in this orbifold structure.\index[terms]{orbifold}
The same argument would work for any base $S$
regardless of dimension.

\hfill

\remark
To simplify the notation and the arguments,
we shall deal with the smooth orbifold fibrations in the same
way as for the smooth ones, and use the standard
terminology instead of adding ``orbifold smooth''
everywhere. Most of the standard results
and constructions in the smooth category are extended 
to the orbifold category in the usual way. For the only exception of
importance to us, see below \ref{_c_1_isotri_Remark_}.\index[terms]{orbifold}

\hfill

The topology and geometry of 
principal toric fibrations was explored in some\index[terms]{bundle!principal!toric}
depth in \cite{_Hofer:remarks_}.
For our present purposes, we shall need the following
theorem. Notice that the principal toric
fibration can be considered to be  a locally trivial principal
toric fibration if we work in the orbifold category.
Following \cite{_Hofer:remarks_}, we define the
Chern classes associated with principal toric
fibrations as follows.\index[terms]{class!Chern}

\hfill

\theorem\label{_c_1_toric_Theorem_}
Let $T$ be an $n$-dimensional 
torus and $TS$ be the sheaf of $T$-valued smooth functions on
$S$. Let $\tilde T=\R^n$ the universal covering of $T$, and
denote by $\tilde TS$ the sheaf of smooth $\tilde T$-valued
functions on $S$. Denote by $\Z^n(S)$ the constant sheaf
with fibre $\Z^n$.
Then the following exact sequence of sheaves
\[
0 \arrow \Z^n(S) \arrow \tilde TS\arrow  TS\arrow 0
\]
gives an exact sequence in cohomology
\[
0 = H^1(S,\tilde TS) \arrow H^1(S,TS) \arrow H^2(S,\Z^n(S))
\arrow H^1(S,\tilde TS) =0,
\]
and a bijective correspondence between
the set $H^1(TS)$ of smooth principal toric fibrations 
and the cohomology group $H^2(\Z^n(S))= H^2(M, \Z)^n$.

\hfill

\proof
Clear (see \cite{_Hofer:remarks_} for details). \endproof

\hfill

\definition\label{_c_1_toric_Definition_}
Let  $\pi:\; M \arrow S$ be a principal toric fibration,
with fibre $T=\R^n/\Z^n$. The $n$ cohomology classes
associated with $\pi$ as in \ref{_c_1_toric_Theorem_}
are called {\bf the  Chern classes of the principal toric
  fibration}.\index[terms]{class!Chern}

\hfill

\remark \label{_c_1_isotri_Remark_}
In the case when $\pi:\; M \arrow S$ 
is a principal elliptic fibration over a 
compact curve, one has $H^2(S, \Z)=\Z$. Therefore,
a toric fibration is uniquely determined by the
vector $c_1(\pi) \in \Z^2= H^2(S, \Z)^2$. Recall that a vector
$v$ in a lattice\index[terms]{lattice} $\Z^n$ is called {\bf primitive}\index[terms]{vector!primitive}
if $v$ is not divisible by an integer $n>1$.
For any primitive
vector  $v$ in a lattice  $\Lambda=\Z^2$, 
there exists $w\in \Lambda$ such that 
$\Lambda = \langle w, v\rangle$. This implies that
the Chern classes\index[terms]{class!Chern} of a
principal elliptic fibration over a curve
are determined (up to an automorphism of $\Lambda$) 
by the biggest $n \in \Z$ such that $c_1(\pi)$ is
divisible by $n$. Using this observation, we shall understand
$c_1(\pi)$ as a non-negative integer; it is equal to zero when
$v=0$ and equal to the largest integer divisor of $v$ when $v\neq 0$. Note that if $S$ is an orbifold, then the first Chern class is\index[terms]{divisor}
still perfectly well-defined and classifies principal elliptic\index[terms]{orbifold}
fibrations, but it is no longer necessarily true
that $H^2(S,\Z) = \Z$: it might also have some torsion.

\hfill

Now we shall prove that a principal elliptic fibration 
$\pi:\; M \arrow S$ over a curve with $c_1(\pi)\neq 0$ 
is always obtained as a quotient of the total space
of an ample $\C^*$-bundle. We need the following preliminary
topological statement. \index[terms]{bundle!line!ample}

\hfill

\claim\label{_c_1_via_covering_Equation_}
Let $\pi:\; M \arrow S$ 
be a principal elliptic fibration over a compact curve $S$,
and $C$ its general fibre. Assume that $c_1(\pi)$ is not torsion.
Then the rank of the natural map $H_1(C) \arrow
H_1(M)$ is 1. Moreover, there exists
a normal subgroup $G\subset \pi_1(M)$ with $\pi_1(M)/G=\Z$ such that
the corresponding $\Z$-covering is infinite on $C$.

\smallskip

\proof
Consider the exact sequence
\begin{equation}\label{_Leray_exact_Equation_}
H_2(S,\Q)\stackrel \delta \arrow  H_1(C,\Q) \stackrel \psi \arrow H_1(M,\Q) \arrow
H_1(S,\Q)\arrow 0
\end{equation}
obtained from the \index[persons]{Leray, J.} Leray spectral sequence of the fibration.\index[terms]{spectral sequence!Leray}
It is easy to see that $\delta$ is dual to the  Chern
class of $\pi$. This implies that the map $\psi:\; H_1(C) \arrow
H_1(M)$ has rank 1 when $c_1(\pi)$ is not torsion and 2 when
it is. Therefore, there exists an element
$v\in \pi_1(C)$ such that its image in $\pi_1(M)$ 
has infinite order.  Take a homomorphism
$H_1(M) \arrow \Z$, which is  non-zero on $\psi(H_1(C))$.
Since $H_1(M)= \frac{\pi_1(M)}{[\pi_1(M),\pi_1(M)]}$,
the map $v$ defines a group homomorphism
$\pi_1(M) \arrow \Z$, which is  non-trivial on
$\pi_1(C)= H_1(C)$.
\endproof

\hfill

\proposition\label{_elli_via_line_bundle_Proposition_}
Let $\pi:\; M \arrow S$ 
be a principal elliptic fibration over a compact curve $S$,
and $C$ its general fibre. Consider a $\Z$-covering $\phi:\;\tilde M
\arrow M$ such that $\phi^{-1}(C)\arrow C$ is an infinite covering.
Then $\tilde M$ is the total space of a principal $\C^*$-bundle\index[terms]{bundle!principal}
associated with a line bundle $L$ on $S$. Moreover, the surface $M$ is obtained
as the quotient of the total space $\Tot^\circ(L)$ of non-zero
vectors in $L$ by a holomorphic automorphism  $q:\; S \arrow S$
acting on $L$ equivariantly and linearly on all fibres.
Finally, the automorphism $q:\; S \arrow S$
has finite order, and $|c_1(L)|$ is equal to $c_1(\pi)$.

\hfill

\proof
Since $\phi^{-1}(C)\arrow C$ is an infinite cyclic covering, 
where $C$ is an elliptic curve, the 
space $\tilde M$ is a $\C^*$-bundle over $S$. Denote the
associated vector bundle by $L$.  Then $\tilde M=\Tot(L^*)$,
and $M$ is obtained from $\Tot^\circ(L)$ as a quotient
by a cyclic group of automorphisms generated by
an automorphism $\tilde q:\; \tilde M \arrow \tilde M$
commuting with the $\C^*$-action. Then $\tilde q$ defines
a holomorphic automorphism of $S$.
However, all automorphisms of $\C^*$ commuting with the 
$\C^*$-action are given by  multiplication by a number. 
Therefore, $\tilde q$ gives a holomorphic section 
$\lambda\in \Hom(L, q^*(L))$. The fibres of
$\pi:\; M \arrow S$ are compact, and hence  $q$ has to be of finite order.

The Chern class\index[terms]{class!Chern} of 
$L$ is obtained from the $\U(1)$-bundle associated with
$L$, and this bundle is homotopy equivalent to $\tilde M$,
hence $c_1(L)$ is equal to $c_1$ of the circle bundle
associated with the $\C^*$-bundle $\tilde M \arrow S$.
However, $c_1$ of $\tilde M \arrow S$ is by construction
equal with $\pm c_1(\pi)$ (the $\pm$ sign is due to the
ambiguity of the definition of $c_1$ of the toric 
fibration, see \ref{_c_1_isotri_Remark_}).
\endproof

\section[The Blanchard theorem]{The Blanchard theorem and the LCK\\ structure \index[terms]{structure!LCK}on 
non-K\"ahler elliptic surfaces}\label{_LCK_on_elli_Subsection_}\index[terms]{theorem!Blanchard}

The main result of the present section
is the next theorem. In the regular case, it was proven in 
\cite{va_gd} (Theorem 3.5 and the subsequent Remark); 
see also \cite{_Vuli:arxiv_}.

\hfill

\theorem\label{_elli_then_Vaisman_Theorem_}
(\cite[Theorem 1]{bel})
Let $M$ be a compact, non-K\"ahler, minimal surface
admitting an elliptic fibration. Then $M$ is 
Vaisman.

\hfill

We give the proof of \ref{_elli_then_Vaisman_Theorem_} at
the end of this section. We shall use \index[persons]{Blanchard, A.} Blanchard theorem
(see \cite{_Blanchard_} or \cite{_Rogov:Iwasawa_}), which 
is applied to the present case as follows.

\hfill

\theorem\label{_Blanchard_Theorem_}
Let $\pi:\; M \arrow S$ be a principal 
elliptic fibration. Then $M$ is K\"ahler if and only if
$c_1(\pi)$ is torsion.

\hfill

\pstep
Suppose that $c_1(\pi)$ is torsion, and consider the exact sequence
\[
0 \arrow H^1(S,\Q) \arrow H^1(M,\Q) \arrow H^1(C,\Q) \stackrel
{\delta^*} \arrow H^2(S,\Q)
\]
dual to \eqref{_Leray_exact_Equation_}.
As shown in the proof of \ref{_c_1_via_covering_Equation_},
one has $\delta^*=0$ if and only if $c_1(\pi)$ is torsion.
However, if $\delta^*\neq 0$, the pullback
$\pi^*(\omega_S)$ of the volume form is exact,
that is  impossible when $M$ is K\"ahler by 
\ref{_exact_pos_non-Kahler_Remark_}.\footnote{This argument 
can be generalized to all dimensions, see Exercise 
\ref{_Blanchard_generalization_Exercise_}.}

\hfill

{\bf Step 2:}
Assume now that $c_1(\pi)$ is torsion, or equivalently, that
$\delta^*=0$. Using \ref{_elli_via_line_bundle_Proposition_},
we obtain $M$ as a quotient of the total bundle
$\Tot^\circ(L)$ by an automorphism acting on the base $S$ as
$\phi:\; S \arrow S$ and on
the fibres as a constant endomorphism
$\lambda:\; L \arrow \phi^* L$, $|\lambda|\neq 0$. 
Since $c_1(L)=\pm c_1(\pi)$ is torsion,
this bundle admits a flat Hermitian connection.
Denote by $\rho$ its monodromy action.\index[terms]{action!monodromy}
Since $\rho$ acts on the fibres of $L$ isometrically,
we can choose a $\rho$-invariant K\"ahler metric $\omega_F$
on  the fibres of $\Tot^\circ(L)/\langle \lambda\rangle$, that are  elliptic curves.
Extending $\omega_F$ to $M$ using the flat Ehresmann
connection on $\pi:\; M \arrow S$, we obtain a \index[terms]{connection!Ehresmann}
closed, positive (1,1)-form $\hat \omega_F$ on $M$ that is 
strictly positive on the fibres of $\pi$;
the sum $\omega:=\hat\omega_F+ \pi^*(\omega_S)$ is
a positive, closed, strictly positive almost
everywhere (1,1)-form. Then $\omega^2>0$.
However, the intersection form\index[terms]{form!intersection} on $H^{1,1}(M)$ 
is negative definite when $M$ is non-K\"ahler
by \ref{_neg_def_h^11_Proposition_} below.
This implies that $M$ is K\"ahler when
$c_1(L)$ is torsion. 
\endproof

\hfill

Let us prove \ref{_elli_then_Vaisman_Theorem_}.
If $M$ is non-K\"ahler, it is obtained as a quotient
of $\Tot^\circ(L)$  by $\lambda:\; L \arrow \phi^* L$,
where $c_1(L)\neq 0$ by \ref{_Blanchard_Theorem_}, Step 2. 
Replacing $L$ by its dual
bundle if necessary, we may assume that $L$ is ample.\index[terms]{bundle!line!ample}
Choose a  Hermitian structure on $L$ such that
its curvature is a positive (1,1)-form on $S$,
and let $\psi\in C^\infty(\Tot^\circ(L))$ be the  
function $\psi(v) = |v|^2$. 

We want the isomorphism $\lambda:\; L \arrow \phi^* L$
to have constant length. Since $\phi$ is of finite order,
we can always represent the curvature $\Theta$ of $L$ by
a $\phi$-invariant form. Using the $dd^c$-lemma\index[terms]{lemma!$dd^c$} as usual,
we choose a metric on $L$ with curvature $\Theta$.
Then the curvature of $\Hom(L, \phi^*L)$ is $\phi^*(\Theta)-\Theta=0$.
 Therefore, $\Hom(L, \phi^*L)$ is a flat unitary bundle; any section
$\lambda\in \Hom(L, \phi^*L)$  gives a holomorphic
section of this flat unitary bundle, and hence  it
is of constant length.

Then $dd^c(\psi)$
is a K\"ahler form\index[terms]{form!K\"ahler} on $\Tot^\circ(L)$
(\ref{_formula_15_19_besse_Theorem_}).
Therefore, $\lambda$ acts on $\Tot^\circ(L)$ by non-isometric holomorphic
homotheties, and hence  the quotient $M = \Tot^\circ(L)/\langle \lambda\rangle$
is LCK. Also, this quotient is Vaisman, because $\C^*$ acts on $\Tot(L^*)= \tilde M$ by holomorphic homotheties (\ref{kami_or}). 

Recall that \ref{_formula_15_19_besse_Theorem_} was stated for
manifolds, and note that the base of the fibration $\Tot^\circ(L)\arrow S$
can be an orbifold, with $L$ a principal $\C^*$-bundle in the orbifold category.
This is not an issue, for the following reason.
Locally in a neighbourhood $U$ of a special point of $S$, the
fibre bundle $\Tot^\circ(L)/\langle \lambda\rangle\arrow S$ 
is obtained as a finite quotient of the smooth fibration
$M_U' \arrow U'$
over a ramified cover $U'$ of $U$
(Subsection \ref{_mult_fibres_Subsection_}).
In the notation of Subsection \ref{_mult_fibres_Subsection_}, 
we express this as $M_U= \frac{M_U'}{\Z/d\Z}$. On $M_U'$, the
form $dd^c(\psi)$ is K\"ahler by \ref{_formula_15_19_besse_Theorem_}.
Since $M_U$ is the quotient of $M_U'$ by the group $\Z/d\Z$ acting 
freely and isometrically, the metric $dd^c(\psi)$ is K\"ahler on
$M_U$ as well.
\endproof

\section{Exercises}

\begin{enumerate}[label=\textbf{\thechapter.\arabic*}.,ref=\thechapter.\arabic{enumi}]

\subsection{Group structure on a curve of genus 1}

\item
Let $E$ be a projective manifold. An
integer divisor in $E$ is a formal linear
combination of irreducible codimension 1 subvarieties 
in $E$ with integer coefficients.
Two divisors $D_1, D_2$ are called {\bf rationally
equivalent} if there exists a rational function on
$E$ with $D_1-D_2$ its zero divisor. Prove that
the group of divisors up to rational equivalence
is $\Pic(E)= H^1(E, \calo^*_E)$.

\item
Let $D= \sum_i n_i x_i$ be a divisor on a compact smooth complex curve $E$ (that is,
a Riemann surface). Irreducible divisors on a curve are points, and hence 
$x_i$ can be understood as points on $E$. {\bf The degree} of the divisor $D$
is the number $\sum_i n_i$. Prove that 
$\deg(D)= \int_E c_1(L)= \frac 1 {2\pi \1} \int_E \Theta_E$, where
$L$ is the line bundle associated with $D$, and $\Theta_E$ the curvature
of a  Chern connection\index[terms]{connection!Chern} on $E$.\footnote{The degree of a bundle $L$ 
on a complex curve $E$ is defined as $\int_E c_1(L)$.}

\item
Let $L$ be a degree 3 holomorphic line bundle on a smooth 
complex curve $E$ of genus 1, and $x_1, x_2\in E$ distinct points.
Denote by $k_{x_1}, k_{x_2}$ the skyscraper sheaves, $k_{x_i}:= \calo_E/{\goth m}_{x_i}$
(quotient by the maximal ideal),  and let
$$q:\; L \arrow k_{x_1}\otimes_{\calo_E} L \oplus k_{x_2}\otimes_{\calo_E}L$$
be the natural quotient map. 
Prove that the induced map $$H^0(E, L) \arrow 
H^0(E, k_{x_1}\otimes_{\calo_E} L \oplus k_{x_2}\otimes_{\calo_E}L)$$
is surjective.

{\em Hint:} Write the long exact sequence
\begin{multline}
0 \arrow H^0(E, \ker q) \arrow H^0(E, L) \arrow\\ 
\arrow H^0(E, k_{x_1}\otimes_{\calo_E} L \oplus k_{x_2}\otimes_{\calo_E}L)
\arrow H^1(E, \ker q)\arrow ...
\end{multline}
Prove that $\ker q$ is a line bundle of 
degree 1, and use the   Kodaira vanishing theorem.\index[terms]{theorem!Kodaira vanishing}

\item Let $L$ be a degree 3 bundle on a  smooth 
complex curve $E$ of genus 1. Prove that $L$ is very ample, 
and $H^0(C, L)$ is 3-dimensional. Prove that $E$
is biholomorphic to a smooth cubic in $\C P^2$.

{\em Hint:} Use the previous exercise.

\item\label{_lines_intersecting_cubic_Equation_}
Let $E\subset \C P^2$ be a smooth
curve given by a cubic equation $P(x,y,z)=0$, 

\begin{enumerate}
\item Let $l_1, l_2\subset \C P^2$ be two lines.
Prove that $(E\cap l_1)- (E\cap l_2)$ is
rationally equivalent to 0, where the
intersection divisor $E\cap l_i$ is
the sum of intersection points counted with their
multiplicities.

\item
Let $e\in E$ be any point, and $l\subset \C P^2$ the tangent
line to $E$ in $e$. Prove that $E\cap l=3e$

\item
Fix $e\in E$, and let $l$ be a line which intersects $E$ in $e, x, y$.
Prove that $2e$ is rationally equivalent to $x+y$.

\item
Fix $e\in E$. For any $x\in E$ distinct from $e$,
let $l_x$ be the line in $\C P^2$ connecting $x$ to $e$.
Then $l_x$ intersects $E$ in three points, $e, x$ and $y$.
Let $\sigma(x):=y$. Given $z, t\in E$, take the line
passing through $z$ and $t$ and let $\mu(z,t)$ denote the
third point of the intersection of this line and $E$.
Prove that the operation $x, y \mapsto \sigma(\mu(x,y))$
takes $(x, y)$ to the only point $z\in E$ such that
$z+e$ is rationally equivalent to $x+y$.
Prove that this operation defines a structure of
a commutative group on $E$, with $e$ its unit,
and $\sigma(x)$ the opposite.

\item Denote by $\Pic_i(E)$ the group of degree $i$ line bundles on $E$,
and let $L_e:\; \Pic_0(E)\arrow \Pic_1(E)$ take a line bundle
associated with a divisor $D$ to a line bundle associated with 
$D+e$. Prove that $\Pic_1(E)=E$,
and this map defines a group isomorphism $\Pic_0(E)\arrow E$,
with respect to the group structure on $E$ defined above.

\end{enumerate}

\item
Let $E$ be a curve of genus 1, and $e\in E$
a fixed point.
\begin{enumerate}
\item Prove that the divisor $3x$ is equivalent to $3e$
for any point $x\in E$.

\item Prove that for any
$x, y\in E$ there exists a unique $z\in E$
such that $x+y$ is rationally equivalent
to $z+e$. 

\item
Prove that the operation
$x, y \arrow z$ defines a group structure
on $E$. 

\item
Construct a natural isomorphism
of this group and $\Pic_0(E)$
(the connected component of $\Pic(E)$.
\end{enumerate}

{\em Hint:} Represent $E$ as a smooth cubic
in $\C P^2$ and apply Exercise \ref{_lines_intersecting_cubic_Equation_}.

\item\label{_E/+-_biholo_to_CP^1_Exercise_}
Let $E\subset \C P^2$ be a smooth curve of genus 1,
and $e\in E$ a point. Consider the group structure on $E$
defined as above, with $e$ its unit. 
\begin{enumerate}

\item Construct a bijective
correspondence between the lines in $\C P^2$ passing through $e$
and unordered pairs $\{x, -x\}\subset E$.

\item Prove that the quotient $E/\{\pm 1\}$ is biholomorphic
to $\C P^1$.
\end{enumerate}

\subsection{Elliptic fibrations}

\item
Let $M$ be a non-K\"ahler complex surface,
smoothly fibred over an elliptic curve
(such a surface is called {\bf a \index[terms]{surface!Kodaira} Kodaira surface}).
Compute $b_1(M)$.

\item 
Let  $\pi:\; M \arrow X$ be a holomorphic
map of complex manifolds.
\begin{enumerate}
\item Let $\nu\in \Lambda^{\dim_\R X} X$ be a positive top degree form on $X$.
Assume that $\pi^* \nu$ is exact. Prove that $M$ is non-K\"ahler.
\item Let $\eta\in \Lambda^k X$ be a non-exact form on $X$.
Assume that $\pi^* \eta$ is exact. Prove that $M$ is non-K\"ahler.
\end{enumerate}

\item\label{_Blanchard_generalization_Exercise_}
Let $\pi:\; M \arrow X$ be a principal
elliptic fibration with $E$ the general fibre. 
Consider the first  Chern class\index[terms]{class!Chern} $c_1(\pi)$ as a
map $H_2(M) \arrow H^1(E)$. 
\begin{enumerate}
\item Suppose that
$c_1(\pi)$ is not torsion. Prove that there exists
a non-exact 2-form $\eta$ on $X$ such that $\pi^*\eta$ is exact.

\item Suppose that
$c_1(\pi)$ is not torsion. Prove that $M$ is non-K\"ahler.
\end{enumerate}

{\em Hint:} Use the previous exercise.

\item
Let $M$ be a compact LCK manifold, 
$\pi:\; M \arrow B$ a complex fibration,
and $p:\; \tilde M \arrow M$ a K\"ahler cover,
equipped with an automorphic K\"ahler form\index[terms]{form!K\"ahler}
$\tilde \omega$. Let $S\subset B$ be a complex
submanifold, and $\{R_i\}$ the set of 
connected components of $p^{-1}(\pi^{-1}(S))$.
Prove that either all $R_i$ are non-compact,
or the set of the volumes
$\{ \int_{R_i} \tilde \omega^{\dim_\C R_i}\}$
is infinite.

\item
Let $\pi:\; M \arrow X$ be a principal
elliptic fibration with $E$ the general fibre. 
Consider the first  Chern class\index[terms]{class!Chern} $c_1(\pi)$ as a
map $H_2(M) \arrow H^1(E)$. Assume that this
map has rank 2 over $\Q$. 
\begin{enumerate}
\item Let $C$ be a general fibre of $\pi$.
Prove that the image of $\pi_1(C)$ in $\pi_1(M)$ is finite.
\item Let $p:\; \tilde M \arrow M$ be any covering of $M$.
Prove that $p^{-1}(C)$ is an infinite union of elliptic curves.
\item
Let $p:\; \tilde M \arrow M$ be a K\"ahler cover\index[terms]{cover!K\"ahler} of $M$,
$\tilde \omega$ its K\"ahler form,\index[terms]{form!K\"ahler}
and $\{C_i\}$ the set of connected components of
$p^{-1}(C)$. Prove that $\int_{C_i}\tilde \omega$ is
the same for all $C_i$.

\item 
Prove that $M$ does not admit an LCK structure\index[terms]{structure!LCK}.\footnote{This is a theorem of V. \index[persons]{Vuletescu, V.} Vuletescu, see \cite{_Vuli:arxiv_}.}
\end{enumerate}

{\em Hint:} Use the previous exercise.

\item
Let $S$ be a K\"ahler manifold, and
$\pi:\; M\arrow S$ a principal toric fibration given by a holomorphic map.
Assume that the Chern class\index[terms]{class!Chern} of $\pi$,
considered to be  a map $H_2(S, \Z) \arrow H^1(T, \Z)$,
where $T$ is the fibre of $\pi$, are torsion.
\begin{enumerate}
\item Prove that the fibration $\pi:\; M\arrow S$
admits a flat, holomorphic Ehresmann connection.\index[terms]{connection!Ehresmann}

\item Prove that $M$ is K\"ahler.
\end{enumerate}

\item
A complex Hermitian manifold $(M, I, \omega)$, $\dim_\C M=n$, is called
{\bf balanced} if $d(\omega^{n-1})=0$.
Let $M$ be a balanced manifold equipped with
a proper holomorphic submersion map $M \arrow S$.
Prove that $S$ is also balanced.

\item
Let $(M, I, \omega)$ be a complex Hermitian manifold,
and $X$ the connected component of the \index[persons]{Barlet, D.} Barlet space (Section \ref{_Barlet_Section_}) of
curves on $M$. Consider the function $\Vol:\; X \arrow \R^{>0}$
taking a point $[C]\in X$ associated with a curve $C$ to
$\int_C \omega$. Using the \index[persons]{Gromov, M.} Gromov compactness theorem, prove
that $\Vol$ reaches its minimum on $X$.

\item
Let $(M, I, \omega)$ be a complex Hermitian manifold,
and $X$ the connected component of the Barlet space of
curves on $M$. Assume that $dd^c \omega=0$.
\begin{enumerate}
\item Consider the function $\Vol:\; X \arrow \R^{>0}$
taking a point $[C]\in X$ associated with a curve $C$ to
$\int_C \omega$. Prove that $dd^c\Vol =0$.
\item Prove that the function $[C] \mapsto -\Vol(C)$ does not have maximum anywhere on
  $X$, unless $\Vol$ is constant.
\item From the \index[persons]{Gromov, M.} Gromov compactness theorem, deduce that 
$\Vol$ is constant and $X$ is compact. \footnote{This is \cite[Corollary 2.19]{_Verbitsky:twistor_}.}
\end{enumerate}

{\em Hint:} Use the previous exercise.

\end{enumerate}


\chapter{Kodaira classification for 
non-K\"ahler complex surfaces}\label{comp_surf}


{\setlength\epigraphwidth{0.8\linewidth}
\epigraph{\em
 \qquad These ambiguities, redundancies and deficiencies remind us of those that doctor Franz Kuhn attributes to a certain Chinese encyclopaedia entitled 'Celestial Empire of benevolent Knowledge'. In its remote pages it is written that the animals are divided into: (a) belonging to the emperor, (b) embalmed, (c) tame, (d) sucking pigs, (e) sirens, (f) fabulous, (g) stray dogs, (h) included in the present classification, (i) frenzied, (j) innumerable, (k) drawn with a very fine camelhair brush, (l) et cetera, (m) having just broken the water pitcher, (n) that from a long way off look like flies.
 }{\sc\scriptsize J. L. Borges, The Analytical Language of John Wilkins }
}
\section{Introduction}
\label{_Intro_Section_}

\subsection{An overview of this chapter}

This chapter, together with Chapter
\ref{_elliptic_Chapter_},
 appeared as a distillation of a lecture course
on complex surfaces given in 2008 and 2012 in Moscow Independent
University. The main reference on complex surfaces is the
great book by Barth, \index[persons]{Hulek, K.} Hulek, Peters and Van de Ven 
\cite{_Barth_Peters_Van_de_Ven_}. This book offers
a powerful narrative, but (as  often happens with
great books) some plots are nested within more
plots, and it sometimes becomes hard to disjoin a 
particular strain from the polyphonic discourse.

Our aim is to give a classification of LCK structures\index[terms]{structure!LCK}
on (a posteriori, non-K\"ahler) surfaces. We thought it was easier
(and more enlightening) to prove the non-K\"ahler part
of the \index[terms]{Kodaira--Enriques classification} Kodaira--Enriques classification directly along
with the classification of LCK structures.



We prove that all non-K\"ahler non-elliptic surfaces are 
of class VII (\ref{_elli_or_class_VII_Theorem_}). There is
enough good literature (\cite{_DOT:Kato_surfaces_},
\cite{_Dloussky:Kato_}, \cite{_Nakamura:curves_}, \cite{_Nakamura:towards_},
 \cite{_Teleman:instantons_}) on class VII surfaces for us 
to give less attention to this case. We gave more references
on the early work on class VII surfaces, in particular on \index[persons]{Inoue, Ma.} Inoue's work,
in the introduction to Chapter \ref{OT_manifolds}.

In this chapter we state the known classification 
results up to the GSS conjecture
and give a new proof of \index[persons]{Brunella, M.} Brunella's theorem on the  existence of LCK metrics
on \index[persons]{Kato, Ma.} Kato surfaces. Together with the GSS conjecture (still not fully proven)
this would imply that all non-K\"ahler complex surfaces are LCK,\index[terms]{conjecture!GSS}\index[terms]{surface!Inoue}
except for some of the Inoue surfaces\index[terms]{surface!Inoue}
(\cite{bel}). 

Results on LCK manifolds with potential\index[terms]{manifold!LCK!with potential} are usually
stated in the assumption that  dimension is $\geq 3$.
The reason is the failure of the \index[persons]{Andreotti, A.} Rossi  and \index[persons]{Siu, Y.-T.} Andreotti--Siu
theorem in dimension 2. In the closing pages
of the chapter, we explain how this restriction
can be avoided, if the GSS conjecture is assumed.
We prove that any minimal LCK complex surface, except
Inoue, admits an LCK metric with potential and
an embedding to a linear Hopf manifold.

\subsection{The Buchdahl--Lamari theorem}

In \cite{_Buchdahl:surfaces_,_Lamari_},\index[terms]{theorem!Buchdahl--Lamari}
N. \index[persons]{Buchdahl, N. P.} Buchdahl and A. \index[persons]{Lamari, A.} Lamari have proven a result previously
known only from the \index[terms]{Kodaira--Enriques classification} Kodaira--Enriques classification of complex surfaces.

\hfill

\theorem\label{_B-L-intro_Theorem_}
Let $M$ be a compact complex surface. Then the first Betti number 
$b_1(M)$ is odd if and only if $M$ is non-K\"ahler.\index[terms]{Betti numbers!first}
\endproof

\hfill

Its direct proof, however, simplifies this 
classification considerably. In this paper we attempt
to recover most of the  Kodaira--Enriques classification\index[terms]{Kodaira--Enriques classification}
for non-K\"ahler surfaces using the \index[terms]{theorem!Buchdahl--Lamari} Buchdahl--Lamari theorem
and an intermediate result 
(stated subsequently as
\ref{_Lamari-current-intro_Theorem_}), which was used\index[terms]{current}
by \index[persons]{Lamari, A.} Lamari to prove \ref{_B-L-intro_Theorem_}. Our proof is different
from the classical one, found in \cite{_Barth_Peters_Van_de_Ven_}, in a few aspects:
we do not rely on birational arguments and the classification
of the elliptic surfaces due to \index[persons]{Kodaira, K.} Kodaira. Also, we 
aim to classify the locally conformally K\"ahler 
structures on complex surfaces.

\hfill

The proof used by \index[persons]{Lamari, A.} Lamari is based on 
the Harvey--Lawson duality theorem (\cite{hl}), which
gives a criterion of K\"ahlerianity in terms of positive currents\index[terms]{current!positive}
and on \index[persons]{Demailly, J.-P.} Demailly's regularization of positive currents
(\cite{_Demailly:regularization_,_Demailly:regularization_flow_}).
\index[terms]{theorem!Demailly regularization}\index[terms]{theorem!Harvey--Lawson duality}

For an introduction to currents and their applications
in algebraic geometry\index[terms]{geometry!algebraic}, see \cite{demailly,_Demailly:analytic_}.
Recall that ``currents'' on $M$ are functionals on the space of 
differential forms on $M$ with compact support that are  continuous in $C^\infty$-topology.\index[terms]{topology!$C^\infty$}
A differential form $\alpha$ defines a current $\tau \arrow \int_M \tau\wedge \alpha$.
This allows us to consider the space of differential forms as a subspace of 
the space of currents.

\hfill

\remark\label{_Current_cohomo_Remark_}
The usual operators and constructions of K\"ahler geometry, \index[terms]{geometry!K\"ahler}
e.g. $d, d^*, \6, \bar\6$, Laplacian, the Hodge decomposition)
extend from differential forms to currents \index[terms]{current}in a natural way.
The corresponding cohomology (de Rham, Dolbeault, Bott--Chern) for
currents is equal to those of differential forms
(\cite{demailly}).

\hfill

A (1,1)-form on a complex manifold $M$, $\dim_\C M=n$ is {\bf positive}
if  defined by a pseudo-Hermitian form \index[terms]{form!pseudo-Hermitian}with non-negative
eigenvalues, and {\bf strictly positive} when it is Hermitian. 
An $(n-1, n-1)$-form is {\bf positive} if it is
a product of $n-1$ positive forms. A (1,1)-current is called\index[terms]{current!positive}
{\bf positive} if it is non-negative on any positive
$(n-1, n-1)$-form.\footnote{For historical reasons, ``positivity''
	for differential forms is understood in the French sense:
	0 is ``positive''. We idly suggest the term ``French-positive'', 
	to avoid confusion.} The notion of positivity for forms is 
compatible with that of currents (\cite{demailly}).

\hfill

\theorem\label{_Lamari-current-intro_Theorem_}
Let $M$ be a compact complex non-K\"ahler surface. Then
there exists a positive, exact\index[terms]{current} (1,1)-current $\Theta$ on $M$.

\proof \cite[Theorem 6.1]{_Lamari_}. \endproof

\hfill

\remark\label{_exact_pos_non-Kahler_Remark_}
The existence of an exact, positive (1,1)-current $\Theta$
immediately implies that the surface $M$ is non-K\"ahler.
Indeed, suppose that $M$ admits a K\"ahler form\index[terms]{form!K\"ahler} $\omega$.
Then $\int_M \omega\wedge \Theta >0$, that is  impossible
because $\Theta$ is exact.

\hfill

In \cite{hl}, Harvey and Lawson proved
a weaker version of \ref{_Lamari-current-intro_Theorem_}
for any $n$-dimensional compact complex manifold $M$.
They have shown that $M$ is non-K\"ahler if 
and only if $M$ admits an exact $(2n-2)$-current $\Theta$ with
its $(n-1,n-1)$-part positive. \index[persons]{Lamari, A.} Lamari\index[terms]{current}
proves that for $n=2$, the current $\Theta$ can be chosen 
of type (1,1).

\subsection{Locally conformally K\"ahler surfaces}
\label{_LCK_Subsection_}

\definition
Let $M$ be a compact complex surface with $b_1(M)=1$.
It is called {\bf a class VII surface} if its \index[terms]{dimension!Kodaira} Kodaira
dimension is $\kappa(M)=-\infty$, and {\bf a class VII${}_0$ surface} if it is also minimal.\index[terms]{surface!class VII}

\hfill

From the Kodaira--Enriques classification
(\cite{_Barth_Peters_Van_de_Ven_}) it follows that
all minimal non-K\"ahler surfaces not of class VII
are elliptic. 
We prove this result (independently from the
rest of the Kodaira--Enriques classification) in \ref{_elli_bundle_Theorem_}.
We prove that all elliptic surfaces are Vaisman when they are minimal
(\ref{_elli_then_Vaisman_Theorem_}).

\hfill

\theorem\label{_LCK_surface_Theorem_}
Let $M$ be a compact, non-K\"ahler surface, 
that is  not of class VII. Then $M$ admits an
LCK structure\index[terms]{structure!LCK}, and a Vaisman one if $M$ is\index[terms]{structure!Vaisman}
minimal.

\hfill

\proof By \ref{_elli_or_class_VII_Theorem_}, $M$ is elliptic,
when it is minimal, and by
\ref{_elli_then_Vaisman_Theorem_}, it is Vaisman. 
By \ref{_blow_up_is_LCK_Theorem_}  
(see also \cite{tric} and \cite{vu1}), the blow-up 
of an LCK manifold in points is again LCK. 
In particular, a blow-up of an LCK
surface remains LCK. Therefore, a surface is LCK\index[terms]{blow-up}
if its minimal model is LCK, for example, elliptic. 
\endproof

\hfill

For class VII surfaces, a complete classification
is not known, but it would follow from the so-called 
``Global Spherical Shell (GSS) conjecture'',\index[terms]{conjecture!GSS}
which claims that any class VII surface $M$ with $b_2>0$
contains an open complex subvariety $U\subset M$
biholomorphic to a neighbourhood of the standard sphere 
$S^3 \subset \C^2$, and $M \backslash U$ is connected.
Surfaces that satisfy this conjecture are called
{\bf Kato surfaces}. \index[terms]{surface!Kato}

\index[persons]{Kato, Ma.} Kato constructed the GSS surfaces explicitly, gluing
a boundary of a small ball excised from a blown-up complex
ball to the boundary of this ball. We give Kato's construction
in Section \ref{_Brunella_proof_}. From this construction
it is clear that any Kato surface $M$ has precisely $b_2(M)$
distinct rational curves. In \cite{_DOT:Kato_surfaces_},
it was shown that the converse is also true: if a class
VII surface $M$ has $b_2(M)$ distinct rational curves, it is Kato.

Another condition, equivalent to the GSS condition, was 
found in \cite{_Dloussky:anticanonical_}. In this paper,
G. \index[persons]{Dloussky, G.} Dloussky proves that a surface $M$ is Kato if and only if it 
admits a line bundle $L$ with $c_1(L)= -b_2(M) c_1(K_M)$ and 
$H^0(M,L)\neq 0$.

\index[persons]{Brunella, M.} Brunella has shown that all \index[persons]{Kato, Ma.} Kato surfaces 
are LCK (\cite{_Brunella:Kato_} and Section \ref{_Brunella_proof_}).
\index[persons]{Bogomolov, F. A.} Bogomolov's theorem on class VII
surfaces with $b_2(M)=0$ 
(see
\cite{_Bogomolov:VII_76_,_Bogomolov:VII_82_,
_Li_Yau_Zang:VII_,_Teleman:bogomolov_}) 
implies that they are
either Inoue surfaces or Hopf surfaces.\index[terms]{surface!Hopf}\index[terms]{surface!Inoue}
The Hopf surfaces are LCK (\cite{go,bel,ov_pams}), and
among the three classes of Inoue surfaces, two\index[terms]{surface!Inoue}
are LCK, and the third contains a subclass that does not admit
an LCK structure (\cite{tric, bel}).\index[terms]{structure!LCK} 

The modern proof of \index[persons]{Bogomolov, F. A.} Bogomolov's classification theorem
(due to A. \index[persons]{Teleman, A.} Teleman and Li-Yau-Zhang) is based on gauge
theory. Using gauge-theoretic methods, A. \index[persons]{Teleman, A.} Teleman
was able to prove the GSS conjecture for minimal
class VII surfaces with $b_2=1$
(\cite{_Teleman:b2=1_}). Extending this
approach to $b_2>1$, \index[persons]{Teleman, A.} Teleman was also able
to prove that any class VII surface with
$b_2(M)=2$ contains a cycle of rational curves,
hence can be smoothly deformed to a blown-up Hopf surface
(\cite{_Teleman:instantons_}). \index[terms]{blow-up}\index[terms]{surface!Hopf}

Once the GSS conjecture is proven, this finishes the
classification of LCK surfaces. If it is true, all non-K\"ahler
surfaces are LCK, except that particular case  of Inoue surfaces.

\hfill

\corollary\label{_Vaisman_qr_or_Hopf_}
All Vaisman complex surfaces are 
quasi-regular or Hopf.

\hfill

{\bf Proof:} By \ref{_LCK_surface_Theorem_},
a Vaisman surface is elliptic or class VII.
By \ref{_Subva_Vaisman_Theorem_}, non-K\"ahler elliptic surfaces
are quasi-regular; indeed, all positive-di\-men\-si\-o\-nal
subvarieties of a Vaisman manifold are tangent
to the leaves of the canonical foliation $\Sigma$, and hence 
the elliptic fibres of a Vaisman elliptic
surface coincide with the leaves of $\Sigma$.
By \ref{_Vaisman_is_Hopf_or_elli_Proposition_}
below, a Vaisman surface is Hopf or elliptic (or both).
\endproof

\hfill

\corollary\label{_Sasakian_3-dim_qr_or_S^3_Corollary_}
A 3-dimensional Sasakian manifold $S$
is quasi-regular or realized as a strictly
pseudoconvex hypersurface in a Hopf surface.

\hfill

\proof
Let $h:\; C(S) \arrow C(S)$ be
a non-trivial homothety of the cone, $(t, s) \mapsto (\lambda t, s)$,
and $M:= C(S)/\langle h \rangle$ the corresponding Vaisman manifold.
From \ref{_Vaisman_qr_or_Hopf_}
it follows that $M$ is a Hopf surface or
$M$ is quasi-regular (or both). In the first case,
$S$ is a hypersurface in $C(S)/\langle h \rangle= S\times S^1$, 
and in the second case, $S$ is quasi-regular
by \ref{_Vaisman_qr_or_Hopf_}.
\endproof


\section{Cohomology of non-K\"ahler surfaces}
\label{_Cohomology_Section_}

\subsection{Bott--Chern cohomology of a surface}
\label{_BC_degree_Subsection_}

In this section, we repeat several  
arguments about the Bott--Chern co\-ho\-mo\-lo\-gy previously given 
in \cite{_Teleman:cone_} and \cite{_Angella_Tomassini_V_}.
Recall that any compact complex manifold admits
a {\bf Gauduchon metric} in any conformal class of Hermitian
metrics (\cite{gau_tor} and Appendix \ref{gauduchon_metric}). A Gauduchon metric\index[terms]{metric!Gauduchon}
on an $n$-dimensional complex manifold is, by definition, a Hermitian metric
with Hermitian form\index[terms]{form!Hermitian} $\omega$ satisfying $dd^c(\omega^{n-1})=0$,
where $d^c=I dI^{-1}$ is the twisted differential.

\hfill

\definition
The {\bf Bott--Chern cohomology group} of a complex manifold
is \[ H^{p,q}_{BC}(M):=\displaystyle\frac{\ker d\restrict{\Lambda^{p,q}M}}{\im dd^c}.\]
\index[terms]{cohomology!Bott--Chern}

\remark
By the $dd^c$-lemma,\index[terms]{lemma!$dd^c$} on a compact K\"ahler manifold the Bott--Chern cohomology groups 
are equal to de Rham and (hence) the Dolbeault
cohomology groups. It is also easy to see that the
complex
\[
\Lambda^{p-1,q-1}M\stackrel{dd^c} \arrow \Lambda^{p,q}M \stackrel d\arrow \Lambda^{p+q+1}M
\]
is elliptic (\ref{_BC_f_dim_Theorem_}; 
see also \cite[Proposition 5]{_Kod-Spen-AnnMath-1960_}), and hence  the Bott--Chern cohomology is  finite-dimensional
on any compact manifold.

\hfill

\theorem\label{_BC_Degree_Theorem_}
Let $M$ be a compact non-K\"ahler surface. Then the kernel of the
natural map $P:\;H^{1,1}_{BC}(M) \arrow H^2(M)$ is 1-dimensional.

\hfill

\pstep
Let $\omega$ be a Gauduchon metric on $M$.
Consider the differential operator $D:\; f\mapsto dd^cf \wedge \omega$ mapping
functions to 4-forms. Clearly, $D$ is elliptic and its index is the
same as  the index of the Laplacian: $\ind D = \ind \Delta=0$, and hence  $\dim \ker D=\dim \coker D$. 
The Hopf maximum principle (\ref{hopf_theorem}) \index[terms]{maximum principle}implies that $\ker D$ only
contains constants, and hence  $\coker D$ is 1-dimensional.
However,
$\int_M D(\alpha) = \int_M dd^cf \wedge \omega = \int_M fdd^c \omega=0$.
This implies that a 4-form $\kappa$ belongs to $\im D$ if and only if $\int_M\kappa=0$.

\hfill

{\bf Step 2:}
Let $\alpha$ be a closed $(1,1)$-form. Define {\bf the degree} $\deg_\omega \alpha:= \int_M \omega \wedge\alpha$.
Since $\int_M dd^cf \wedge \omega=0$, this defines a map $\deg_\omega:\; H^{1,1}_{BC}(M,\R)\arrow \R$.
Given a closed $(1,1)$-form $\alpha$ of degree 0, 
the form $\alpha':=\alpha- dd^c(D^{-1}(\alpha\wedge \omega))$ 
satisfies $\alpha'\wedge \omega=0$; in other words, it is an $\omega$-primitive
(1,1)-form. For $\omega$-primitive forms,\index[terms]{form!primitive} one has $\alpha'\wedge \alpha' = -|\alpha'|^2\omega\wedge\omega$,
giving
\begin{equation} \label{_primitive_square_Equation_}
\int_M \alpha'\wedge \alpha'= -\|\alpha'\|^2_\omega,
\end{equation}
that is  impossible when $\alpha'$ is a non-zero class in $\ker P$,
because $\alpha'$ is exact. 
Therefore, any class of zero degree in $\ker P\subset H^{1,1}_{BC}(M,\R)$ vanishes. This implies that
any two classes in $\ker P$ are proportional.
\endproof

\subsection{First cohomology of non-K\"ahler surfaces}

Recall that the Dolbeault cohomology  $H^{p,0}(M)$ of a compact complex manifold $M$ \index[terms]{cohomology!Dolbeault}
coincides with the space of holomorphic $p$-forms on $M$.
On a compact K\"ahler manifold, holomorphic
$p$-forms are closed because they are harmonic.
On non-K\"ahler manifolds, this is generally false.
However, this is true on compact complex surfaces.

\hfill

\lemma\label{_holo_closed_Lemma_}
All holomorphic forms on a compact complex surface are closed.

\hfill

\proof
For 2-forms, this statement is trivial.
Indeed, the de Rham differential of a (2,0)-form $\beta$ on a
complex surface $M$ is equal to $\bar\6\beta=0$.

Let $\alpha\in \Lambda^{1,0}M$ be a holomorphic 1-form. 
Then $\bar\6\alpha=0$, because it is holomorphic, and
by the same reason $d\alpha$ is a holomorphic, exact (2,0)-form. 
Then $d\alpha\wedge d\bar\alpha$ is a positive (2,2)-form,
giving $0=\int_M d\alpha\wedge d\bar\alpha =\|d\alpha\|^2$.
Then $d\alpha=0$, and $\alpha$ is closed.
\endproof

\hfill

\theorem\label{_H^1_odd_Theorem_}
Let $M$ be a non-K\"ahler manifold and
$\Theta$ a non-zero exact positive (1,1)-current (\ref{_Lamari-current-intro_Theorem_}). \index[terms]{current}
Then  $\Theta$ is cohomologous in $H^{1,1}_{BC}(M)$ to
$d^c\theta$, where $\theta$ is a closed
1-form.
Moreover,
\begin{description}
	\item[(i)] 
Consider the natural embedding
$H^{1,0}(M)\oplus \overline{H^{1,0}(M)} \arrow H^1(M, \C)$
(\ref{_holo_closed_Lemma_}). Then 
	\begin{equation}\label{_Hodge_on_H^1_surface_Equation_}
	H^1(M) = H^{1,0}(M)\oplus \overline{H^{1,0}(M)} \oplus
	\langle[\theta]\rangle,
	\end{equation}
	that is, $H^1(M)$ is generated by the cohomology classes
	of holomorphic and antiholomorphic forms and $\theta$.

	\item[(ii)] Since all holomorphic forms are closed
          (\ref{_holo_closed_Lemma_}),
	the antiholomorphic forms are $\bar\6$-closed and have Dolbeault classes.
	This gives a natural map from
        $\overline{H^{1,0}(M)}$ to $H^{0,1}(M)$.
	Then  $H^{0,1}(M)$ is generated by $\overline{H^{1,0}(M)}$ and
	the Dolbeault class $[\theta^{0,1}]$, giving
	$H^{0,1}(M)= \overline{H^{1,0}(M)}\oplus \langle [\theta^{0,1}] \rangle$.
\end{description}

{\bf \pstep Existence of $\theta$:}
Let $\Theta_0$ be a smooth, closed, real $(1,1)$-form representing
the same Bott--Chern class \index[terms]{class!Bott--Chern}as the current $\Theta$. This\index[terms]{cohomology!Bott--Chern}
form exists because the Bott--Chern cohomology of $(p,q)$-currents\index[terms]{current!of integration}
is equal to the Bott--Chern cohomology of smooth $(p,q)$-forms.
We are going to show that $\Theta_0= d^c \theta$, where $\theta$ is 
a closed 1-form.

Let $H^{1,1}_{BC}(M) \stackrel \tau \arrow H^2(M, \C)$
be the tautological map.
The following exact sequence is apparent 
if we apply \ref{_BC_exact_sequence_via_Dolbeault_Theorem_}
in the case where the 1-dimensional local system $L$ is trivial:
\begin{equation}\label{_Dolbeault_BC_exact_Equation_}
H^{0,1}(M) \oplus \overline{H^{0,1}(M)} \xlongrightarrow{\6+\bar\6} 
H^{1,1}_{BC}(M)\stackrel \tau \arrow H^2(M, \C).
\end{equation}
By \ref{_BC_Degree_Theorem_}, the kernel of $\tau$ is
1-dimensional. Therefore, $\6(H^{0,1}(M))=\bar\6(\overline{H^{0,1}(M)})=\ker \tau$.
Choose
$[\alpha^{0,1}]\in H^{0,1}(M)$ and $[\beta^{1,0}]\in \overline{H^{0,1}(M)}$
in such a way that 
$[\alpha^{0,1}]= [\overline{\beta^{1,0}}]$
and $[\Theta_0] = \6[\alpha^{0,1}] = \bar\6[\beta^{1,0}]$,
where $[\Theta_0]$ is the Bott--Chern class\index[terms]{class!Bott--Chern} of $\Theta_0$.

An appropriate choice of the forms $\alpha^{0,1}$ and $\beta^{1,0}$
in the corresponding cohomology classes would give
\begin{equation}\label{_Theta_0_differential_Equation_}
\Theta_0= \6\alpha^{0,1}= \bar\6\beta^{1,0}, \text{\ and \ } 
\overline{\beta^{1,0}}=\alpha^{0,1}.
\end{equation}
Indeed, $\Theta_0 = \6\alpha^{0,1}+ \6\bar\6 f$,
because $\Theta_0$ is Bott--Chern-cohomologous to $\6\alpha^{0,1}$;
replacing $\alpha^{0,1}$ by $\alpha^{0,1}+\bar\6 f$,
that belongs to the same Dolbeault class, we obtain
a representative that satisfies $\Theta_0 = \6\alpha^{0,1}$.
Since $\Theta_0$ is real, this also implies
$\Theta_0 = \bar\6\alpha^{1,0}$, where
$\alpha^{1,0}= \overline{\alpha^{0,1}}$.
This gives 
\begin{equation}\label{_alpha_diffe_conju_Equation_}
d(\alpha^{1,0}-\alpha^{0,1})=\Theta_0 - \Theta_0 =0.
\end{equation}

Let $\theta:= -\1(\alpha^{1,0}-\alpha^{0,1})$.
Then \eqref{_alpha_diffe_conju_Equation_}
gives $d(\theta)=0$, and 
\eqref{_Theta_0_differential_Equation_} implies
$d(I\theta)= \Theta_0$.

\hfill

{\bf Step 2. Proof of \ref{_H^1_odd_Theorem_} (i):}
A non-zero linear combination of holomorphic and antiholomorphic 
forms is closed, but never exact (\ref{_holo_closed_Lemma_}).
It follow that  the natural map 
$H^{1,0}(M)\oplus \overline{H^{1,0}(M)} \stackrel \kappa \arrow H^1(M,\C)$ 
is injective. To prove \ref{_H^1_odd_Theorem_} (i), it remains to show that
its image has codimension 1, and prove that $H^1(M,\C)$ 
is generated by the class of  $\theta$ and the image of $\kappa$. 

Clearly, the kernel of the natural map $\6:\; H^{0,1}(M)\arrow H^{1,1}_{BC}(M)$
coincides with the subgroup $\overline{H^{1,0}(M)}\subset H^{0,1}(M)$. Indeed,
the cohomology class of a (0,1)-form $\alpha\in \ker \bar\6$ 
belongs to the kernel of this map if $\6\alpha =
\6\bar\6f$, and this implies that $\alpha -\bar\6 f$ is 
antiholomorphic.

The image of $\6:\; H^{0,1}(M)\arrow H^{1,1}_{BC}(M)$ is contained in 
the one-dimensional kernel of the map $P:\;H^{1,1}_{BC}(M) \arrow H^2(M)$ 
(\ref{_BC_Degree_Theorem_} and \eqref{_Dolbeault_BC_exact_Equation_}). 
Therefore, the image of $\kappa$ has codimension 
at most 1. The class $\6 \alpha^{0,1}\in H^{1,1}_{BC}(M)$ 
has non-zero degree in the sense of \ref{_BC_Degree_Theorem_},  Step 1,
because $\6 \alpha^{0,1}$ is Bott--Chern cohomologous
to $\Theta$, and $\Theta$ is a positive current.\index[terms]{current!positive} Therefore
$\6 \alpha^{0,1}$ generates $\ker P=  \langle [\Theta]\rangle$. 
We proved \eqref{_Hodge_on_H^1_surface_Equation_}.

\hfill

{\bf Step 3. Proof of  (ii):}
Consider the map $\6:\; H^{0,1}(M)\arrow H^{1,1}_{BC}(M)$.
By Step 2, $\6([\theta^{0,1}])$ is non-zero in
$H^{1,1}_{BC}(M)$ (in the notation of Step 2, we have $\theta^{0,1}=\alpha^{0,1}$)
Since $\6:\; H^{0,1}(M)\arrow H^{1,1}_{BC}(M)$
vanishes on antiholomorphic forms, 
the Dolbeault cohomology  class of $\theta^{0,1}$
does not belong to $\overline{H^{1,0}(M)}$,
giving an injective map 
\begin{equation}\label{_holo_to_Dolbeault_Equation_}
\overline{H^{1,0}(M)}\oplus \langle [\theta^{0,1}]
\rangle\arrow H^{0,1}(M).
\end{equation}
Now, the kernel of $\6:\; H^{0,1}(M)\arrow H^{1,1}_{BC}(M)$
is generated by antiholomorphic forms.
This gives an exact sequence
\begin{equation}\label{_exact_Dolbeault_Equation_}
0\arrow \overline{H^{1,0}(M)}\arrow H^{0,1}(M)\arrow
H^{1,1}_{BC}(M). 
\end{equation}
The image of $\6:\; H^{0,1}(M)\arrow H^{1,1}_{BC}(M)$
is at most 1-dimensional by \ref{_BC_Degree_Theorem_},
hence it is generated by $\6([\theta^{0,1}])$.
Using \eqref{_exact_Dolbeault_Equation_}, 
we obtain that the injective map 
\eqref{_holo_to_Dolbeault_Equation_}
is actually surjective. We proved 
\ref{_H^1_odd_Theorem_} (ii).
\endproof

\hfill

\remark 
As a corollary, we deduce the \index[persons]{Buchdahl, N. P.} Buchdahl--Lamari theorem
(\ref{_B-L-intro_Theorem_}). Indeed, from
\ref{_H^1_odd_Theorem_} (ii) it follows that $b_1(M)$ is odd for
$M$ non-K\"ahler. This way, one obtains
\ref{_B-L-intro_Theorem_} from the existence
of an exact, positive (1,1)-current\index[terms]{current!positive} on each complex surface
with $b_1$ odd (\ref{_Lamari-current-intro_Theorem_}).

\subsection{Second cohomology of non-K\"ahler surfaces}

\claim\label{_H^2,0_Claim_} 
Let  $M$ be a compact complex surface. Then: 
$$H^{2,0}(M) = \overline{H^{0,2}(M)}= H^0(K_M),$$
where $K_M$ denotes the canonical bundle.\index[terms]{bundle!vector bundle!canonical}

\hfill

\pstep
By Serre  duality,\index[terms]{duality!Serre} $H^{0,2}M = H^0(K_M)^*$, that has the same
dimension as $H^0(K_M)=H^{2,0}(M)$.

\hfill

{\bf Step 2:}
The natural map $R:\;H^{2,0}(M) \arrow \overline{H^{0,2}(M)}$ 
is injective, because its kernel is formed by $\6$-exact holomorphic
forms $\alpha= \6\beta$, but for such $\alpha$ one has 
$0=\int_Md\beta\wedge \bar \alpha = \int_M\alpha\wedge \bar\alpha= \|\alpha\|^2$,
that is  impossible unless $\alpha=0$. As shown in Step 1, these spaces
have the same dimension, and hence  $R$ is an isomorphism.
\endproof

\hfill

\corollary \label{_E_2_degenerates_Corollary_}
(\cite[Theorem IV.2.8]{_Barth_Peters_Van_de_Ven_})\\
The Hodge--de Rham--Fr\"olicher spectral sequence of a 
compact complex\index[terms]{spectral sequence!Hodge--de Rham--Fr\"olicher}
surface converges in $E_1$. In other words, $d_2$ of this
spectral sequence vanishes.

\hfill

\proof This argument is standard. The Hodge--de Rham-Fr\"olicher 
degenerates if $\sum_{p,q} \dim H^{p,q}(M)=  \sum_i \dim H^i(M)$.

We have proven the degeneration 
for $H^1(M)$ in \ref{_H^1_odd_Theorem_}.
 Serre  and Poincar\'e duality give the degeneration of this spectral sequence
in the page $E_1$ for the $H^3(M)$
term.\index[terms]{duality!Poincar\'e}\index[terms]{duality!Serre} 
The degeneration of this spectral sequence
for the second cohomology is implied by the
following general observation.

\hfill

\lemma\label{_E_2_all_lines_degenerate_Lemma_}
Let $(E_k^{p,q}, d_k)$ be a spectral sequence
such that $E_1^{p,q}= E_\infty^{p,q}$ for
$p+q \neq s$ (in our case, $s=2$).
Then $d_2=0$.

\hfill

\proof
By definition, the $k$-th page $\bigoplus_{p,q} E_{k+1}^{p,q}$ is equal to the cohomology
of $d_k$ on  $\bigoplus_{p,q} E_k^{p,q}$.
Since $ E_1^{p,q}= E_\infty^{p,q}$, for each $k \geq 2$,
the image of $d_k$ and the kernel of $d_k$ on  $E_1^{p,q}$
is equal to zero unless $p+q=s$. However, the $d_k$
originating in the space $E_{k-1}^{p,q}$ with $p+q=s$
has image in $E_k^{p-k+1,q+k}=E_1^{p-k+1,q+k}=E_\infty^{p-k+1,q+k}$, and it 
has to vanish because
$E_1^{p-k+1,q+k}=E_\infty^{p-k+1,q+k}$.
Similarly, the differential $d_k$ which arrives 
at  $E_k^{p,q}$ with $p+q=s$ has to originate in
$E_{k-1}^{p+k-1,q-k}$, and $d_k\restrict{ E_{k-1}^{p+k-1,q-k}}=0$
for the same reason.

This proves \ref{_E_2_all_lines_degenerate_Lemma_} and
\ref{_E_2_degenerates_Corollary_}.
\endproof

\hfill

\remark
We have shown that 
\[
H^2(M, \C)  = H^{2,0}(M) \oplus H^{1,1}(M) \oplus H^{0,2}(M)\ \text{
and}\ \ H^{2,0}(M)= \overline {H^{0,2}(M)},
\]
 same as for K\"ahler surfaces
(\ref{_H^2,0_Claim_}). 

\hfill

\lemma\label{_BC_to_H^11_Lemma_}
Let $M$ be a compact complex surface. Then the natural
map $$H^{1,1}_{BC}(M) \arrow H^{1,1}(M)$$ is surjective.

\hfill

\proof
Let $\alpha\in \Lambda^{1,1}M$ be a $\bar\6$-closed form representing
a Dolbeault cohomology class $[\alpha]\in H^{1,1}(M)$.
Clearly, the class $[\alpha]\in H^{1,1}(M)$ can be represented by a closed (1,1)-form
if and only if the Bott--Chern class\index[terms]{class!Bott--Chern} $[\beta]\in
H^{2,1}_{BC}(M)$ of $\beta:=\6 \alpha$ 
vanishes. Indeed, if $\beta=\6\bar\6\eta$, then
$\alpha - \bar\6\eta$ is $\6$-closed.

Since the Hodge--de Rham--Fr\"olicher 
spectral sequence degenerates in
$E_1$, and $\beta$ is 
$\6$-exact and $\bar\6$-closed,
one has $\beta= \bar\6\gamma$, for some $\gamma\in \Lambda^{2,0}M$.
However, the cokernel of the map
$\6:\; \Lambda^{1,0}M\arrow \Lambda^{2,0}M$
is generated by holomorphic 2-forms, and they are all
closed (\ref{_H^2,0_Claim_} ),
hence $\bar\6(\Lambda^{2,0}M) = \bar\6\6(\Lambda^{1,0}M)$,
and the Bott--Chern class\index[terms]{class!Bott--Chern} of $\beta\in \bar\6(\Lambda^{2,0}M)$ vanishes.
\endproof

\hfill

\proposition\label{_neg_def_h^11_Proposition_}
{(\cite[Theorem IV.2.14]{_Barth_Peters_Van_de_Ven_})} \\
Let $M$ be a compact non-K\"ahler surface. Then the
intersection form on $H^{1,1}(M)$ is negative definite.\index[terms]{form!intersection}

\hfill

\proof 
Fix a Gauduchon metric $\omega$ on $M$.\index[terms]{metric!Gauduchon}
Let $[\alpha]\in H^{1,1}(M)$ be a cohomology class.
By \ref{_BC_to_H^11_Lemma_}, we can represent $[\alpha]$
by a closed (1,1)-form $\alpha$. Consider the degree
functional $\deg_\omega:\; H^{1,1}_{BC}(M,\R) \arrow \R$
defined in Subsection \ref{_BC_degree_Subsection_}.
Since $\deg_\omega(\Theta)>0$ for an exact (1,1)-current $\Theta$,\index[terms]{current}
any cohomology class $[\alpha]\in H^{1,1}(M)\subset H^2(M)$
can be represented by a closed (1,1)-form $\alpha$ with
$\deg_\omega\alpha=0$. Acting as in the proof of
\ref{_BC_Degree_Theorem_}, we find 
$f\in C^\infty(M)$ such that $\alpha-dd^c f$ is primitive.\index[terms]{form!primitive}
Replacing $\alpha$ by $\alpha-dd^c f$, we
obtain $\int_M \alpha\wedge \alpha= -\|\alpha\|^2_\omega<0$
as in \eqref{_primitive_square_Equation_}.
\endproof

\hfill

\remark
In \cite[Theorem 2.37]{_Brinzanescu:bundles_} (see also \cite{_Brinzanescu_Flondor:1_, _Brinzanescu_Flondor:3_}) it was shown that
the space $H^{1,1}(M) \cap H^2(M, \Q)$ (the rational Neron--Severi group) 
of a complex surface with algebraic
dimension 0 has negative definite intersection form. This important result also follows from \ref{_neg_def_h^11_Proposition_}.\index[terms]{form!intersection}




\subsection{Vanishing in multiplication of holomorphic 1-forms}

\ref{_neg_def_h^11_Proposition_} has an interesting corollary.

\hfill

\proposition \label{_prod_hol_1_forms_Proposition_}
Let $M$ be a  non-K\"ahler, compact, complex surface. Then for any 
holomorphic  1-forms $\alpha, \beta$, the product
$\alpha\wedge \beta$ vanishes, and the product $\alpha\wedge \bar\beta$ is exact.

\hfill

\pstep
Let $\alpha, \beta$ be holomorphic 1-forms.
Then $\alpha\wedge \beta$ is a holomorphic 2-form, with
$\int_M \alpha\wedge \beta\wedge \bar \alpha\wedge \bar \beta>0$
unless $\alpha\wedge \beta=0$. 

\hfill

{\bf Step 2:}
The intersection form on $H^{1,1}(M)$ is negative definite by
\ref{_neg_def_h^11_Proposition_}. \index[terms]{form!intersection}
Therefore, $\int_M \eta\wedge \bar\eta=0$ 
implies $\eta=0$ for any $\eta\in H^{1,1}(M)$.  
Take $\eta=\alpha\wedge \bar \alpha$
and $\rho=\beta\wedge \bar \beta$. Then clearly 
$\int_M \eta\wedge \bar\eta=\int_M \rho\wedge \bar\rho=0$,
which implies 
$\int_M \alpha\wedge \beta\wedge \bar \alpha\wedge \bar \beta=0$,
hence $\alpha\wedge \beta=0$ (Step 1).

\hfill

{\bf Step 3:}
To see that $\eta:=\alpha\wedge \bar\beta$ is exact,
we use \ref{_neg_def_h^11_Proposition_}
again. Unless $\mu:=\alpha\wedge \bar \beta$ is exact, one would have
$\int_M \mu\wedge \bar\mu<0$,
but $\int_M \mu\wedge \bar\mu=0$ 
as we have already shown in Step 2. Therefore, $\mu$ is exact.
\endproof

\subsection{Structure of multiplication in
de Rham cohomology of non-K\"ahler surfaces without curves}\index[terms]{cohomology!de Rham}

Further on, we shall need the following lemma.

\hfill

\lemma\label{_A_defi_Lemma_}
Let $M$ be a compact non-K\"ahler surface,
and $\Theta$ an exact real (1,1)-form. Then there exists
a closed form $\tilde \theta \in \Lambda^1M$  such that 
$d^c(\tilde \theta)=\Theta$. 

\hfill

\proof
Without restricting the generality, we may assume that $\Theta$ 
is an exact real (1,1)-form representing a non-zero class in
$H^{1,1}_{BC}(M)$. Then, by \ref{_H^1_odd_Theorem_} (ii), 
the cohomology class of $\Theta$ belongs to the
image of the natural map $\6:\; H^{0,1}(M) \arrow H^{1,1}_{BC}(M)$,
and can be expressed as $\6(\theta^{0,1})$, where
$\theta$ is a closed form. Therefore,
$\Theta- d^c(\theta)= d^cdf$  for some function 
$f\in C^\infty M$. Take $\tilde \theta:= \theta+ d f$.
Then $\tilde\theta$ is a closed 1-form such that $d^c\tilde\theta=\Theta$.
\endproof

\hfill

In this subsection we prove the following structure theorem
about multiplication in $H^1(M)$. It is a posteriori true
in the general situation, but we need it only for surfaces without
complex curves.

\hfill

In the following theorem, we no longer make a distinction between
the closed forms and the de Rham cohomology\index[terms]{cohomology!de Rham} classes they represent; this abuse
of notation makes the text less unreadable.

\hfill

\theorem \label{_multiplica_H^1_Theorem_}
Let $M$ be a compact non-K\"ahler surface
without curves, and $\theta \in \Lambda^1M$
a closed 1-form such that $d^c(\theta)$ is non-zero 
in $H^{1,1}_{BC}(M, \R)$ (\ref{_H^1_odd_Theorem_} (i)).  Denote by $W \subset \Lambda^1M$
the subspace generated by holomorphic and antiholomorphic forms.
Using \ref{_H^1_odd_Theorem_}, we identify $W$ with its image in 
$H^1(M)= W \oplus \langle \theta \rangle$.
Then: 
\begin{description}
	\item[(i)] The multiplication $W \wedge W \arrow H^2(M)$ vanishes.
	The multiplication $W \wedge \theta \arrow H^2(M)$ is injective.
	\item[(ii)] Let $\theta^c:= I(\theta)$.
	Then for any non-zero $x \in W$, the form $x\wedge \theta^c$ is 
	closed and represents a non-zero element of $H^2(M)$.
	\item[(iii)] The Poincar\'e pairing on the spaces \index[terms]{pairing!Poincar\'e}
	$W \wedge \theta$ and $W \wedge \theta^c\subset H^2(M)$ vanishes. The 
	Poincar\'e pairing between these two subspaces is non-degenerate.
\end{description}

\pstep Multiplication $W \wedge W \arrow H^2(M)$ vanishes
by \ref{_prod_hol_1_forms_Proposition_}.
Let us prove that for any non-zero $x\in W$,
the form $x\wedge \theta^c$ is 
closed and represents a non-zero element of $H^2(M)$.

Without restricting the generality, we may assume that
$W\neq 0$. For any $x\in W$, the closed (1,1)-form $x\wedge\bar x$ 
is cohomologous to zero
(\ref{_prod_hol_1_forms_Proposition_}). 

For any non-zero holomorphic form $x$, the form
$\Theta:=x\wedge \bar x$ is positive, non-zero and exact.
We fix such $\Theta$, and fix a closed form 
$\theta$ such that $d^c\theta=\Theta$ as in 
\ref{_A_defi_Lemma_}. Choose any $y\in W$.
Then $d(y\wedge \theta^c) = y \wedge x\wedge \bar x=0$
(\ref{_prod_hol_1_forms_Proposition_}) hence
$y\wedge \theta^c$ is closed. 

\hfill

{\bf Step 2:}
The injectivity of $W \xlongrightarrow{x\wedge \theta} H^2(M)$
and $W \xlongrightarrow{x\wedge \theta^c} H^2(M)$ would follow if
we prove that the formula 
\[
x, y \arrow \int_M x\wedge y \wedge \theta\wedge \theta^c,
\ \ x, y\in W
\]
defines a non-degenerate pairing on $W$,
and hence (i) and (ii) follow from (iii).

The Poincar\'e pairing on the images $W\wedge \theta$ and
$W\wedge\theta^c$ in $H^2(M)$ 
vanishes because $\theta\wedge\theta = \theta^c \wedge\theta^c=0$.
To prove that the Poincar\'e pairing between $W\wedge \theta$\index[terms]{pairing!Poincar\'e}
and $W\wedge \theta^c$ is non-degenerate, take a holomorphic 
1-form $x\in W$. To finish the proof of \ref{_multiplica_H^1_Theorem_} (iii) it 
would suffice to show that the integral 
$\int_M \1 x\wedge \bar x \wedge \theta \wedge \theta^c$
is positive.  The form $\1 x\wedge \bar x \wedge \theta
\wedge \theta^c$ is positive, because it is
the product of a $(2,0)$-form and its complex conjugate.
It is  non-zero if $x$ is not
proportional to the (1,0)-part of $\theta$, denoted 
$\theta^{1,0}$. 

\hfill

{\bf Step 3:}
It remains to show that $\theta^{1,0}$
is not proportional to a holomorphic form; this is
where we use that $M$ has no holomorphic curves.
In this case, the zero set of any 1-form $x\in W$ is 
a finite set.

Suppose that $x\in W$ is holomorphic
and proportional to $\theta^{1,0}$.
Then there exists a smooth function $\alpha$
defined outside of the zero set $S$ of $x$
such that $\theta^{1,0}= \alpha x$.
Without restricting generality, we may assume that
$d\theta^{1,0}= x \wedge \bar x$ (Step 1). 
Since $\theta^{1,0}=\alpha x$,
this implies  $d(\alpha x)= d\alpha \wedge x= - \bar x \wedge x$, which gives
$\bar \6\alpha = -\bar x$. Therefore, $dd^c \alpha=0$
outside of $S$. Then $\alpha$ is locally
the real part of a holomorphic function, defined
on $M \backslash S$, where $S\subset M$ is a finite set.
Using the Hartogs extension theorem,\index[terms]{theorem!Hartogs} we obtain that
$\alpha$ is smooth and defined globally on $M$. Then $\alpha=0$
by the maximum principle (\ref{hopf_theorem}), because $dd^c$ is elliptic. 
\endproof

\hfill

\corollary\label{_W_rational_Corollary_}
Let $M$ be a compact complex surface
without curves, and $W \subset H^1(M, \C)$
the subspace generated by holomorphic and antiholomorphic forms.
Then $W$ is a rational subspace, that is, there
exists a subspace $W_\Q \subset H^1(M, \Q)$
such that $W = W_\Q \otimes_\Q \C$.

\hfill

\proof
For $x\in H^1(M)$, denote by $L_x:\; H^1(M)\arrow H^2(M)$
the map $y \mapsto x\wedge y$.
 Without restricting
generality, we may assume $W\neq 0$; then $\dim W\geq
2.$ By \ref{_multiplica_H^1_Theorem_}, $\rk L_\theta= \dim
W\geq 2$
and $\rk L_x =1$ for $x\in W$. Therefore,  
$W$ is the space of all $x\in  H^1(M)$ such that the
$\rk L_x=1$. Since the multiplication in
$H^*(M)$ is defined over $\Q$, the subspace $W\subset H^1(M)$
is rational.

This can be explained as follows. Consider the $\C$-valued de Rham 
cohomology $H^*(M,\C)$ as a tensor product $H^*(M,\Q)\otimes_\Q \C$,
and let $G:= \Aut_\Q \C$. We extend the action of $G$ from $\C$
to $H^*(M,\Q)\otimes_\Q \C$ by acting on the second factor.
A space $W\subset H^*(M, \C)$
is rational if and only if it is $G$-invariant; this
follows because the set $\C^G$ of $G$-invariant elements
of $\C$ is $\Q$, that is  not hard to prove.

Now, the subspace $W\subset H^1(M,\C)$ is determined 
as
\[
W = \{ x\in H^1(M,\C)\ \ |\ \ \rk L_x=1\};
\]
since $G$ preserves the multiplication in $H^*(M,\C)$,
the space $W$ is $G$-invariant.
\endproof

\hfill

\remark
\ref{_W_rational_Corollary_}
is true for all surfaces. We leave its proof
as an exercise to the reader.


\section{Elliptic fibrations on non-K\"ahler surfaces}


For a version of the following theorem (with a different proof),
see \cite[Theorem 3.17]{_Brinzanescu:bundles_} and \cite{_Brinzanescu:manuscripta_}.

\hfill

\theorem \label{_elli_bundle_Theorem_} 
Let $M$ be a non-K\"ahler compact complex surface
admitting a non-constant positive dimensional 
family $D_t$ of complex curves. Then $M$ 
admits a holomorphic, surjective map $\pi:\; M \arrow S$ 
to a curve, and the pullback $\pi^*\Vol_S$
is equal to the differential of a $d^c$-closed form. Moreover, 
all fibres of $\pi$ are elliptic curves, homologous to 0,
and the general fibres of $\pi$ are isomorphic.\index[terms]{curve!elliptic}

\hfill

\pstep Let $D=\bigcup D_i$ be an irreducible decomposition
of the general member of the family $D_t$. If all
$D_i$ are rigid (have no deformations in $M$), the curve $D$ 
is also rigid. Therefore, one of the irreducible
components of $D$ belongs to a family of positive dimension.
Replacing $D$ with this component, we can assume that
the general curve in our family is irreducible.

Let $C_1, C_2$ be two general curves in the same family $F$
of irreducible curves. Since $C_1$ is homologous to $C_2$,
and the intersection form\index[terms]{form!intersection} on $H^{1,1}(M)$ is negative definite,
one has $(C_1 \cdot C_2)=0$ and the curves $C_i$ are
homologous to 0. This implies that through each point of $M$
passes no more than one curve from $F$, and $F$ is 1-dimensional.

Let $S$ be the base of the family. 
This is a reduced irreducible complex space of dimension
$1$, and the curves in the family with multiplicity $1$ form an
analytic family of $1$-cycles on $M$. By the universal property of
the \index[persons]{Barlet, D.} Barlet spaces, the family is induced from the universal one via
a map $S \to {\goth B}_1(M)$, where ${\goth B}_1(M)$
is the Barlet space of complex curves on $M$. Replace $S$ with the irreducible\index[terms]{space!Barlet}
component of ${\goth B}_1(M)$ that contains the image of the
map. Then, as shown in \ref{_Barlet_on_surf_Corollary_}, each
connected component of ${\mathfrak B}_1(M)$ is compact. Therefore
$S$ is compact.  Let $S^m\subset M \times S$ be the corresponding incidence family,
that is, the set of pairs $(x, C)$ with $[C]\in S$ is a curve and $x\in C$.
Then the variety $S^m \subset M \times S$ 
is irreducible, and its image $\Psi(S^m) \subset M$ under the
projection $\Psi:M \times S \to M$ is an irreducible closed
subvariety. Since the family is not constant, $\Psi(S^m)$ cannot be
a curve, so $\Psi(S^m)=M$.

All complex curves $C\subset M$ with the fundamental class $[C]\neq
0$ have negative self-intersection by
\ref{_neg_def_h^11_Proposition_}. This implies that the general
curves in the family $S$ are homologous to 0.
For any family of varieties, the 
pullback of the volume form on the base is always 
cohomologous to the fundamental class of the fiber.
This implies that $\pi^*\Vol_S$ is exact.

Let $U \subset S$
be the dense open subset parametrizing cycles that are irreducible
with multiplicity one, and let $U^m \subset S^m$ be it preimage in
$S^m$. Since all curves in the family parametrized by $U$
do not intersect, the proper map $\Psi:S^m \to M$ is injective on
$U^m$. First of all, this means that $S$ is of dimension exactly
$1$, and $S^m$ is of dimension $2$. Secondly, $\Psi$ is one-to-one
over the complement $M \setminus \Psi(S^m \setminus U^m)$, and since
$S^m \setminus U^m$ is a finite disjoint union of ``bad'' fibres of
the projection $S^m \to S$, $\Psi(S^m \setminus U^m) \subset M$ is
at most one-dimensional. Last, note that although these ``bad''
fibres are different as cycles, they might have common irreducible
components, so they might intersect. However, if a curve $C \subset
S^m$ is contracted by $\Psi$, then it must lie in $S^m \setminus U^m$, 
so each of its connected components lies in a bad fibre. On the other hand,
$\Psi$ is injective on each of the fibres including
bad ones, so $C$ cannot exist.

We have shown that $\Psi:S^m \to M$ is a proper bijective map that is
an isomorphism over the complement to a subvariety of positive
codimension. A bijective bimeromorphic map to a normal variety is biholomorphic.\index[terms]{variety!normal}\index[terms]{map!bimeromorphic}
Since $M$ is smooth, thus normal, $\Psi$ is an
isomorphism, so we obtain a holomorphic map $\pi:M \cong S^m \to S$
to an irreducible compact curve. Taking its Stein factorization, we
may assume that the fibres of $\pi$ are connected (\cite{_Stein:factorization_,_Hartshorne:AG_}).\index[terms]{Stein factorization}

This implies that $S$ is normal. Indeed, a complex variety $S$
is normal if any locally bounded meromorphic function $f$ on $S$ is holomorphic\index[terms]{function!meromorphic}
(Subsection \ref{_normality_of_cone_Subsection_}). Since $\pi^*(f)$ is locally 
bounded and meromorphic on $M$, and $M$ is normal, this function
is holomorphic on $M$. Since $\pi^*(f)$ is constant on the fibres
of $\pi$ and they are connected, it is a lift of a holomorphic
function on $S$. Therefore, $f$ is also holomorphic.

Since $\dim S=1$, normality of $S$ implies the smoothness.

\hfill

{\bf Step 2:}
Now let  $C$ be a general fibre of $\pi$.
Using the standard isomorphism
$\pi^*(\Lambda^1S)= N^*_\pi M$ between the pullback of the
cotangent sheaf and the conormal bundle to the fibres of
$\pi$ we obtain an exact sequence
\begin{equation}\label{_adjunction_Equation_}
0 \arrow \pi^*(\Lambda^1S)\arrow \Lambda^1 M \arrow
\Lambda^1_\pi M \arrow 0,
\end{equation}
where $\Lambda^1_\pi M$ is the bundle of 
holomorphic differentials on fibres of $\pi$. 
This sequence makes sense on the smooth locus
of $\pi$, but we need it only in a neighbourhood
of the generic fibre $C$.

The exact sequence \eqref{_adjunction_Equation_}
is a special case of the {\bf conormal exact sequence}, that is  valid\index[terms]{conormal exact sequence}
for any smooth curve $C\subset M$
\[
0 \arrow N^* C \arrow \Lambda^1 M \arrow \Lambda^1 C \arrow 0
\]
It implies the {\bf adjunction formula}\index[terms]{adjunction formula}
\begin{equation}\label{_adjunction_formula_Equation_}
K_M \restrict C = N^* C \otimes_{\calo_C} K_C
\end{equation}
Since the normal bundle $NC$ is trivial, 
this implies $K_M\restrict C= K_C$, 
that is, 
the canonical bundle of $M$ restricted
to $C$ gives the canonical bundle to $C$.
Since the self-intersection of $C$ is zero,
$C$ is homologous to 0 (\ref{_neg_def_h^11_Proposition_}).
This gives $\int_C c_1(K_M)=0$,
hence the degree of $K_C$ is equal 0.
This implies that $C$ is an elliptic curve.

\hfill

{\bf Step 3:}
The form $\Theta:=\pi^* \Vol_S$ is exact, as shown in Step 1.
By \ref{_A_defi_Lemma_}, $\Theta$ is the differential of a
$d^c$-closed 1-form.
This means that we have arrived at the assumptions of 
\ref{_elli_then_isotri_Theorem_}, which implies
that $\pi$ is a principal elliptic fibration,
with all fibres elliptic curves, and general
fibres isomorphic.
\endproof

\hfill

\corollary\label{_fibra_Corollary_}
Let $M$ be a non-K\"ahler, compact, complex 
surface. Assume that $b_1(M)>3$. Then $M$ admits an elliptic fibration.

\hfill

\proof
\ref{_H^1_odd_Theorem_} implies that $\dim H^{1,0}(M)= \frac{b_1(M)-1}2$, and hence
we have $\dim H^0(\Omega^1M)\geq 2$.
By \ref{_prod_hol_1_forms_Proposition_},
all globally defined holomorphic 1-forms on $M$ are pointwise proportional.
Consider the rank 1 sheaf $L\subset \Omega^1M$ generated by
globally defined holomorphic 1-forms. If $\dim H^{1,0}(M)=\dim H^0(L) > 1$, the zero
divisors of the sections of this sheaf form a continuous family
of curves on $M$, and \ref{_elli_bundle_Theorem_}
can be applied.
\endproof

\hfill

\remark\label{kodaira.rem}
Note that the result we have proven does not imply that the 
elliptic fibration $\pi$ 
is a submersion. Indeed, $\pi$ may have multiple fibres. For
an example, take a classical Hopf surface $H:=\frac{\C^2 \backslash 0}{(z\sim 2z)}$\index[terms]{surface!Hopf!classical}
and let $\tau$ takes $(x, y)$ to $(\sqrt 2 x, -\sqrt 2 y)$. 
This surface is equipped with an elliptic fibration
$\pi:\; H \arrow \C P^1$ sending $(x, y)\in H$ to $[x:y]\in \C P^1$.
Clearly, $\tau$ is an
involution. This action is free on $H$, but the corresponding action
on $\C P^1$ has two fixed points, $[0:1]$ and $[1:0]$,
and $\tau$ defines an involution on the corresponding fibres.
These two fibres become multiple fibres of the
elliptic fibration $\pi_\tau:\; \frac H{\langle \tau\rangle} \arrow 
\frac {\C P^1}{\langle \tau\rangle}$. To see that $\pi_\tau$ is
not submersive in these two fibres, we use Ehresmann theorem,
claiming that any proper submersion is locally trivial.\index[terms]{theorem!Ehresmann's fibration}
To see that $\pi_\tau$ is not locally trivial, we 
compute the Riemannian volumes of its fibres.
By construction, the special fibres of $\pi_\tau$ 
are obtained  quotients of the fibres of $\pi$ by an involution,
and the general fibres are identified with the corresponding fibres of $\pi$.
Therefore, the Riemannian volume of the special fibres of $\pi_\tau$ 
is two times less than the volume of general fibres,
hence $\pi_\tau$ cannot be smooth in these fibres.


\section{Class VII surfaces}


\subsection{The Riemann--Roch formula for embedded curves}
\label{_RR_curves_Subsection_}

Here we give some basic results about singular curves
embedded to a complex surface. We follow \cite[II.11]{_Barth_Peters_Van_de_Ven_}.

\hfill

Let $\chi(F) = \dim H^0(F) - \dim H^1(F)$ denote {\bf the Euler characteristic}
of a sheaf $F$ on a curve $C$, and $\deg_C(L):= \int_C c_1(L)$. We shall use
 {\bf the Riemann--Roch formula} \index[terms]{Euler characteristic!of a sheaf}
\begin{equation}\label{_RR_Equation_}
\chi(L) = \deg_C(L) + \chi(\calo_C)
\end{equation}
(Exercise \ref{_RR_curve_Exercise_}).\index[terms]{Riemann--Roch formula}

\hfill

We need the following version of the adjunction formula \eqref{_adjunction_formula_Equation_}, valid for all
curves $C\subset M$, even singular. Define {\bf the arithmetical genus}\index[terms]{arithmetical genus}\index[terms]{adjunction formula}
of a curve $C$ by $g(C):= 1- \chi(\calo_C)$.
When $C$ is smooth, this number is clearly equal to its genus.

\hfill

\claim
Let $C\subset M$ be a complex curve on a complex surface. Then 
\begin{equation} \label{_adjunction_sing_Equation_}
g(C)=1 + \frac 1 2 \deg_C(K_M\otimes_{\calo_M} \calo(C)) = 1+ \frac {(C\cdot K_M) + (C\cdot C)} 2,
\end{equation}
where $\calo(C)$ denotes the line bundle 
associated with the divisor $C\subset M$.\footnote{The last equation of
\eqref{_adjunction_sing_Equation_} is vacuous, it is the same quantity 
written in different notation.}

\hfill

\proof
To prove \eqref{_adjunction_sing_Equation_}, we start by defining
{\bf the dualizing sheaf of $C$} as $\omega_C:= K_M\restrict C \otimes_{\calo_M} \calo(C)$.\index[terms]{sheaf!dualizing}
Serre's duality (\cite[II.6.1]{_Barth_Peters_Van_de_Ven_})
implies the existence of a surjective map
$H^1(\omega_C) \stackrel \delta \arrow \C$ such for any locally free sheaf $F$
the pairing $H^1(F) \times H^0(\omega_C \otimes_{\calo_C} F^*) \arrow H^1(\omega_C)  \stackrel \delta \arrow \C$
is non-degenerate. By definition, $\deg_C\omega_C=\deg_C(K_M\otimes_{\calo_M} \calo(C))$.
However, the Riemann--Roch formula \eqref{_RR_Equation_}
gives 
\begin{equation}\label{_degree_omega_C_Equation_}
\deg_C \omega_C = \chi(\omega_C) - \chi(\calo_C) = - 2 \chi(\calo_C),
\end{equation}
(the last equation follows from \index[terms]{duality!Serre} Serre duality).
Then 
\begin{equation}\label{_genus_via_degree_Equation_}
g(C)=1- \chi(\calo_C)= 1 + 1/2 \deg_C \omega_C.
\end{equation}
Again, $\deg_C\omega_C=\deg_C(K_M\otimes_{\calo_M} \calo(C))$
together with \eqref{_genus_via_degree_Equation_}
 immediately implies \eqref{_adjunction_sing_Equation_}.
\endproof

\hfill

This result can be used to obtain a geometric description of the genus
$g(C)= 1- \chi(\calo_C)$ when the curve $C$ is singular. 
Recall that the incidence graph of a curve $C$ is the graph with 
vertices enumerated by irreducible components of $C$, and
edges the intersection points, with multiple edges of multiplicity $r=\sum_{i=1}^k r_i$ \index[terms]{incidence graph}
between components of $C$ which intersect in $k$ points with
multiplicities $r_1, ..., r_k$.

\hfill

\proposition\label{_g(C)=0_tree_Proposition_}
Let $C$ be a connected, compact complex curve (possibly singular).
Then $g(C) \geq 0$. Moreover, $g(C)=0$ if
and only if $C$ is a tree of smooth rational curves,
that is, a union of smooth rational curves joined
in points of multiplicity one, with the 
incidence graph a tree.\index[terms]{tree}

\hfill

\proof 
Since $C$ is connected and compact, we have $H^0(C)=\C$,\index[terms]{normalization}
hence $\chi(C) \leq 1$. This implies that  $g(C)= 1- \chi(\calo_C) \geq 0$.
Let $\tilde C\stackrel \nu \arrow C$ be the normalization of $C$, that is,
the curve obtained by taking apart all irreducible components
of $C$ and resolving the singularities. Since $\nu$ is finite, 
the pushforward functor is exact on coherent sheaves (\cite{_Hartshorne:AG_}). 
The tautological map $\tau:\; \calo_C \arrow \nu_* \nu^* \calo_C$ is injective
because $\nu$ is surjective; outside of the singularities, $\tau$ is an isomorphism.
Consider the exact sequence
\begin{equation}\label{_normalization_exact_sequence_Equation_}
0 \arrow \calo_C \arrow \nu_* \nu^* \calo_C \arrow R \arrow 0
\end{equation}
where $R:= \coker \tau$. The sheaf $R$ is a direct sum of skyscraper
sheaves located at the singular points of $C$; its dimension at each
of these points is equal to the multiplicity of the singularity.
Since $\nu_*$ is exact, and $\nu^* \calo_C = \calo_{\tilde C}$,
we have $H^i(\nu_* \nu^* \calo_C) = H^i(\nu^* \calo_C)= H^i(\calo_{\tilde C})$.
Then \eqref{_normalization_exact_sequence_Equation_} gives:
$$\chi(\nu_* \nu^* \calo_C)= \chi(\calo_{\tilde C})= \chi(\calo_C) + \chi(R).$$
Since the support of $R$ is 0-dimensional, we have $H^1(R)=0$, and
$\chi(R)= \dim H^0(R)=e + \mu$, where $e$ 
is the number of edges of the incidence graph $\Gamma$ of $C$, and
$\mu$ is the sum of multiplicities of the singular points
of all irreducible components of $C$.

A connected graph $\Gamma$ with $n$ vertices is a tree if and only if
it has $n-1$ edges.

Let $n$ be the number of irreducible components $C_1, ..., C_n$ of $C$.
Since $C$ is connected, the number $e$ of edges of its incidence graph $\Gamma$ is
at least $n-1$, and it is equal to $n-1$ if and only if $\Gamma$ is a tree. 
This gives
\[ \chi(\calo_C)= \chi(\calo_{\tilde C})- \chi(R)\leq  
\chi(\calo_{\tilde C})- e = \sum_{i=1}^n \chi(C_i) - e \leq n-(n-1) =1.
\]
The last inequality is not strict if and only if $e=n-1$
(that is, the graph $\Gamma$ is a tree) and $\sum_{i=1}^n \chi(C_i)=n$,
that is, all $C_i$ are smooth rational curves.
If the inequality is strict, we have $\chi(\calo_C)<1$,
which implies that $g(C)= 1- \chi(\calo_C) >0$.
\endproof

\hfill

\remark
Gromov's theorem on the compactness of the moduli 
of pseudo-holo\-morphic curves (\ref{_Gromov_compactness_Theorem_}) implies\index[terms]{theorem!Gromov's compactness}
that a limit of complex curves is obtained
by contracting disjoint simple loops to points. When these
simple curves are homotopy equivalent to zero,
the limit acquires bubble trees\index[terms]{tree!bubble} formed by rational
curves. Parts of this result can be recovered
using the same argument as used in the proof
of \ref{_g(C)=0_tree_Proposition_}.\index[terms]{tree!bubble}
In particular, the appearance of trees of
rational curves in both 
\ref{_g(C)=0_tree_Proposition_}
and \ref{_Gromov_compactness_Theorem_} 
is not a random coincidence: \ref{_g(C)=0_tree_Proposition_}
can be seen as an algebraic version of the bubble
tree part of \ref{_Gromov_compactness_Theorem_}.

\subsection{(-1)-curves}

\definition
A {\bf (-1)-curve} on a complex surface
is a smooth rational curve of self-intersection -1.

\hfill

By the famous result of \index[persons]{Grauert, H.} Grauert (\cite{_Grauert:Modifikazionen_}), any 
curve of negative self-intersection on a surface can be
blown down (in analytic category). This result is a generalization
of a classical criterion of Castelnuovo.\index[terms]{Castelnuovo criterion}

\hfill

\theorem\label{_Catselnuovo_Theorem_} (Castelnuovo criterion)\\
Let $C \subset M$ be a (-1)-curve on a complex surface.
Then there exists a holomorphic, bimeromorphic map
$M \arrow M_1$ bijective outside of $C$ and contracting
$C$ to a point, with $M_1$ smooth.

\hfill

\proof \cite[V.5.7]{_Hartshorne:AG_}.
In \index[persons]{Hartshorne, R.} Hartshorne this result is proven for projective surfaces.
However, the smoothness part of the argument works in a complex
analytic setup just as well, and the contractibility of $C$
follows directly from \cite{_Grauert:Modifikazionen_}.
 \endproof

\hfill

Further on, we will use the following proposition, 
taken from \cite[III.2.2]{_Barth_Peters_Van_de_Ven_}.

\hfill

\proposition\label{_-1_via_intersections_Proposition_}
Let $C$ be a complex curve on a complex surface $M$.
Assume that $(C\cdot C)<0$ and $(K_M \cdot C)<0$.
Then $C$ is smooth, and satisfies $(C \cdot C)=-1$
(in other words, $C$ is a (-1)-curve).

\hfill

\proof
If $(C\cdot C)<0$ and $(K_M \cdot C)<0$, then 
$\deg_C \omega_C = \deg_C K_M + (C\cdot C) < 0$
(Subsection \ref{_RR_curves_Subsection_}).
Then $\chi(\calo_C) = 1 - \dim H^0(\omega_C) =1$,
hence $g(C)=0$, and $C$ is a smooth genus 0
curve by \ref{_g(C)=0_tree_Proposition_}.
\endproof

\subsection{Non-K\"ahler surfaces are either class VII or elliptic}

In this subsection we prove the following theorem.

\hfill

\theorem\label{_elli_or_class_VII_Theorem_}
Let $M$ be a non-K\"ahler compact minimal complex surface. 
Then $M$ is principal elliptic or belongs to class VII
(or both).\index[terms]{surface!elliptic}\index[terms]{surface!class VII} 

\hfill

\pstep 
We observe that any minimal non-K\"ahler surface admitting
a positive-dimensional family of complex curves is principal
elliptic (\ref{_elli_bundle_Theorem_}). To prove  \ref{_elli_or_class_VII_Theorem_}, we need
to show that any non-K\"ahler minimal surface is of class VII
or admits a positive-dimensional family of complex curves.

As $b_1(M)>0$, and $H_1(M, \Z)= \pi_1(M)/[\pi_1(M),  \pi_1(M)]$,
there exists a subgroup of any given finite  index $r$ in $\pi_1(M)$.
Denote the corresponding covering by $\sigma:\; M_1\arrow M$.
By construction, $M_1$ is a compact complex surface.
Since $\chi(\calo_{M})= \frac{c_1^2+c_2}{12}$ is expressed through the
curvature, one has $|\chi(\calo_{M_1})|= d \cdot |\chi(\calo_{M})|\geq d$
for any $d$-sheeted covering $\sigma:\; M_1\arrow M$.
Unless $\chi(\calo_M)=0$, we can find a covering $M_1$ with $\chi(\calo_{M_1}) < -3$
or $\chi(\calo_{M_1}) 3$. In the first case, $\dim H^{1,0}(M)>1$ and $M_1$ is elliptic by
\ref{_fibra_Corollary_}. 
In the second case, the canonical bundle\index[terms]{bundle!vector bundle!canonical}
$K_{M_1}$ satisfies $\dim H^0(K_{M_1})\geq 2$, and hence  it
has sections which give a continuous family of divisors.
Therefore, for $\chi(\calo_{M})\neq 0$, the finite covering $M_1$ has a 
continuous family of divisors. Then
$M$ also has such a continuous family, and it is elliptic by
\ref{_elli_bundle_Theorem_}.

\hfill

{\bf Step 2:} Assume that $b_1(M)=3$, but $M$ is not elliptic. 
Then, by Step 1, $\chi(\calo_{M})=0$, and hence  $h^{0,1}(M)=2$ and $h^{0,2}(M)=1$,
and the same is true for all non-ramified coverings of $M$.
We are going to prove that $M$ contains no curves.

First, we prove that $c_1(M)^2=(K_M\cdot K_M)=0$. Otherwise, $c_1(M)^2 < 0$,
so $K_M$ is non-trivial, and since $h^{0,2}(M)=1$, we have $K_M =
\calo(D)$ for some effective divisor $D=\sum_ia_iD_i$. If for some\index[terms]{divisor!effective}
$i$ we have $(D_i\cdot D_i) < 0$ and $(K_M \cdot D_i) < 0$, then $D$ is a
$(-1)$-curve (\ref{_-1_via_intersections_Proposition_}) 
and this is not possible since $M$ is minimal. If
$(D_i\cdot D_i)=0$, for all $i$, then $D_i$ is homologous to $0$ by
\ref{_neg_def_h^11_Proposition_}. In this case, $(K_M \cdot
D_i)=0$, and then $(K_M)^2 = \sum_ia_i(K_M\cdot D_i)=0$.

Using the Riemann--Roch formula for complex surfaces, from 
 $c_1(M)^2=0$ we derive \index[terms]{Riemann--Roch formula}
$$
0=\chi(\calo_{M})= \frac{c_1(M)^2+c_2(M)}{12}= \frac{c_2(M)}{12} =
\frac{e(M)}{12},
$$
where $e(M)$ is the topological Euler characteristic.
Then $b_2(M)= b_1(M) + b_3(M) - b_0(M) - b_4(M) = 6-2=4$. Moreover,
the same holds for any unramified cover $M_1$ of $M$.

Now assume that $M$ contains an irreducible curve $C$. Then by the
adjunction formula \eqref{_adjunction_sing_Equation_}, its arithmetical genus is $g(C) =
\frac{(C\cdot C)+(K_M \cdot C)}{2}+1$, and since $(C\cdot C) \leq 0$ and\index[terms]{adjunction formula}
$(K_M\cdot K_M)=0$, we have $(K_M \cdot C)=0$, and hence  $g(C) \leq 1$. Then
$C$ is smooth rational, or rational with a single node, or
rational with a single cusp, or smooth elliptic, and in any case,
$b_1(C) \leq 2$. Since $b_1(M)=3 > b_1(C)$, we have an unramified
covering $M_1 \to M$ of any degree $d$ that splits over $C$, and thus
contains non-intersecting curves $C_1,\dots,C_d$ with the same
self-intersection $l=(C_i\cdot C_i)$, $i=1,\dots,d$.

Then if $l < 0$, the classes of the curves $C_i$ span a
$d$-dimensional subspace in $H^2(M_1,\Q)$, so that $d \leq b_1(M_1)
= 4$, and this is a contradiction since $d$ was arbitrary.

Otherwise, $l=0$, and all the curves $C_i$ are homologous to $0$. Let
$D$ be their union. Since the curves $C_i$ are disjoint, one has
$\calo_D= \bigoplus_i \calo_{C_i}$. This gives the short exact sequence
\begin{equation}\label{_exact_O_curves_Equation_}
0 \arrow \calo(-D) \arrow \calo_{M_1} \arrow \bigoplus_i \calo_{C_i}\arrow 0.
\end{equation}
The corresponding long exact sequence
\[ \C = H^0(\calo_{M_1}) \arrow \C^d = \bigoplus_i H^0(\calo_{C_i})
\arrow H^1(\calo(-D)).\]
shows that $\dim H^1(M_1,L) \geq d-1$, where $L:=\calo(-D)$. However,
$\chi(L)=\chi(\calo_{M_1})=0$ because $c_1(L)=[D]=0$, and for any line
bundle $L'$ on $M_1$, we have $\dim H^0(M_1,L') \leq 1$ -- otherwise
$M_1$;  hence, $M$ also carries a non-trivial family of divisors and is elliptic by
\ref{_elli_bundle_Theorem_}. Thus $\dim H^0(M_1,L)$ and $\dim
H^2(M_1,L) = \dim H^0(M_1,K_{M_1} \otimes L^*)$ are at most $1$, and
$d \leq 3$. This is again a contradiction, and hence  $M$ contains
no complex curves.

\hfill

{\bf Step 3:}
We can prove now that all surfaces 
with $b_1(M)=3$ are elliptic.
Arguing by absurd, we may assume $M$ is non-elliptic,
but in this case, $M$ has no complex curves (Step 2).

Consider the space $W\subset H^1(M)$ generated by
holomorphic and antiholomorphic forms.
By \ref{_H^1_odd_Theorem_} (i), $H^{1,0}(M) =\C$, $H^{0,1}(M) =\C^2$,
and $W$ is 2-dimensional. 
Let $\alpha\in \Lambda^{1,0}M$ be a holomorphic 1-form
generating the space $H^{1,0}M$. Fix $x\in M$
and define the Albanese map by taking\index[terms]{map!Albanese}
$y\in M$ to $\int_\gamma \alpha$, where
$\gamma$ is a path-connecting $x$ to $y$.
This map depends on the choice of the path $\gamma$,
hence it gives a map $\Alb:\; M \arrow \C /\Lambda$,
where $\Lambda\subset \C$ is the period lattice\index[terms]{lattice!period} of $\alpha$,
that is, the set  $\{\int_v\alpha \ \ |\ \ v \in H_1(M,
\Z)\}$.  By 
\ref{_W_rational_Corollary_}, $W$ is rational whenever $M$
has no complex curves.
Therefore, the integral $\int_v\alpha$
vanishes on one of the generators in $H_1(M,\Z)$,
and  $\Lambda$ is a rank 2 lattice in $\C$.
It is easy to see that it is discrete. 
Then the Albanese map\index[terms]{map!Albanese} gives a fibration
$M \arrow \C/\Lambda$, and $M$ contains
a continuous family of curves, and hence  it is
elliptic.

\hfill

{\bf Step 4:} It remains to show that $M$ is 
of class VII or elliptic if $\chi(\calo_{M})=0$ (when $\chi(\calo_{M})\neq 0$,
$M$ is elliptic by Step 1).
When $b_1(M)=3$, the surface $M$ is elliptic (Step 3), and the same
holds if $b_1(M) \geq 5$ by \ref{_fibra_Corollary_}. 
Therefore, we may assume that $b_1(M)=h^{0,1}(M)=1$. 
The same is true for all coverings $M_1$ of $M$, ramified or
unramified, because otherwise $M_1$ would be elliptic, and
$M$ would contain a positive-dimensional family of complex curves. 
Since $\chi(\calo_M)=0$, we have $h^{0,2}(M)=0$.
To prove that $M$ is class VII one has
to check that $H^0(K_M^n)=0$ for any $n>0$.

If $K^n_M$ admits a non-zero section, then $K_{M_1}$
admits a non-zero section for some
finite ramified covering $\sigma:\; M_1\arrow M$. Indeed, given a 
section $\alpha$ of $K^n_M$, with zero divisor $C$,
locally around a point $x\in M\backslash C$, one can take
its $n$-th degree root $\sqrt[n]{\alpha}$ and obtain a section
of $K_M$. The sheaf of such sections is finite and locally constant
on $M \backslash C$, \index[terms]{sheaf!locally constant}
hence it becomes trivial after passing to a finite covering 
$\widetilde {M\backslash C} \arrow M\backslash C$. 
In a neighbourhood $U$ of $x\in C$, we interpret
$\sqrt[n]{\alpha}$ as a multivalued function with values in $K_M$,
that can be trivialized in $U$. The graph of $\sqrt[n]{\alpha}\subset U\times C$
can be glued to $\widetilde {M\backslash C}$, giving a 
ramified covering $\sigma:\; M_1\arrow M$ and a section of $K_{M_1}$.
Then $h^{0,2}(M_1)\neq 0$, and hence  $M_1$ is elliptic (Step 3).
Therefore, $M=\sigma(M_1)$ contains a continuous family
of holomorphic curves; it is elliptic by \ref{_elli_bundle_Theorem_}.
\endproof


\section{Brunella's theorem: all Kato surfaces are LCK}
\label{_Brunella_proof_}\index[terms]{theorem!Brunella}


The Kato surfaces are also called ``Global spherical shell surfaces''.
By definition, a \index[persons]{Kato, Ma.} Kato surface is a surface $M$ which possesses a 
{\bf global spherical shell}, that is, an open subset $U\subset M$
that is  biholomorphic to a neighbourhood of $S^3$ in $\C^2$
and such that $M \backslash U$ is connected.\index[terms]{surface!Kato}
In this section we prove \index[persons]{Brunella, M.} Brunella's theorem,
showing that all Kato surfaces are locally
conformally K\"ahler. 

For the original definition and early works on Kato surfaces see
\cite{_Kato:announce_,_Kato:Kinokunya_,_Kato:sugaku_,_Dloussky:Kato_},
and for \index[persons]{Brunella, M.} Brunella's original proof,
\cite{_Brunella:Enoki_} and \cite{_Brunella:Kato_}.

It is possible to see that all Kato surfaces are 
deformations of a blow-up of a Hopf surface. To see
this, we consider the following explicit \index[terms]{surface!Hopf}
construction of Kato surfaces.\index[terms]{blow-up}

Let $M$ be a Kato surface, and $S\subset U \subset M$ 
the corresponding 3-sphere. Consider the map $\chi:\; \pi_1(M) \arrow \Z$
mapping a path $\gamma$ to the intersection index $\gamma \cap S$.
Clearly, $\chi$ is a group homomorphism. Denote the
corresponding $\Z$-covering by $\tilde M$, and the preimages
of $S$ in $\tilde M$ by $S_i$, $i\in \Z$
(these preimages can be enumerated, because
the deck transform group $\Z$ acts on the
set of preimages of $S$ freely and transitively).
Denote by $M_i$ the subsets of $\tilde M$
situated between $S_i$ and $S_{i-1}$.
Clearly, each $M_i$ is a fundamental domain of
the deck transform action.

Each $M_i$ has two boundary components, with 
$S_i$ pseudo-convex and $S_{i-1}$ pseudo-concave.
Gluing a ball $B$ to $S_{i-1}$, we obtain a manifold
$\hat M_i= B \coprod_{S_{i-1}} M_i$ that is  compact and has a strictly
pseudoconvex boundary. By the solution of the Levi problem
(\index[persons]{Oka, K.}Oka, \index[persons]{Grauert, H.} Grauert, \ref{_Grauert_Levi_Theorem_}), 
$\hat M_i$ is holomorphically convex.\index[terms]{Levi problem}
Using  Remmert's reduction theorem\index[terms]{theorem!Remmert reduction} \cite{rem},\index[terms]{manifold!holomorphically convex} 
we obtain a proper map with connected fibres
$p:\; \hat M_i \arrow X$, where $X$ is Stein.
Since the neighbourhood of the boundary of
$X$ is biholomorphic to a neighbourhood of $S^3\subset B$,
Hartogs theorem \index[terms]{theorem!Hartogs}implies that $X$ is in fact biholomorphic
to $B$. Then $p:\; \hat M_i \arrow X$ is a bimeromorphic,
holomorphic map to a ball, and is obtained by a sequence
of  blow-ups.\index[terms]{blow-up}

Following \index[persons]{Dloussky, G.} Dloussky \cite{_Dloussky:Kato_}, we define the 
{\bf Kato data}\index[terms]{Kato data} on a closed ball $B\subset \C^2$
as a bimeromorphic, holomorphic map $\hat B \arrow B$,
together with an open subset $B_0\subset \hat B$
and a biholomorphism from $B_0$ to an open ball.
Then the complement $\hat B \backslash B_0$
has two smooth boundary components, that are  both isomorphic
to $S^3\subset B$ and can be glued together to obtain
a compact complex surface. 

\begin{figure}[ht]
	\begin{center}\includegraphics[width=0.3\linewidth]{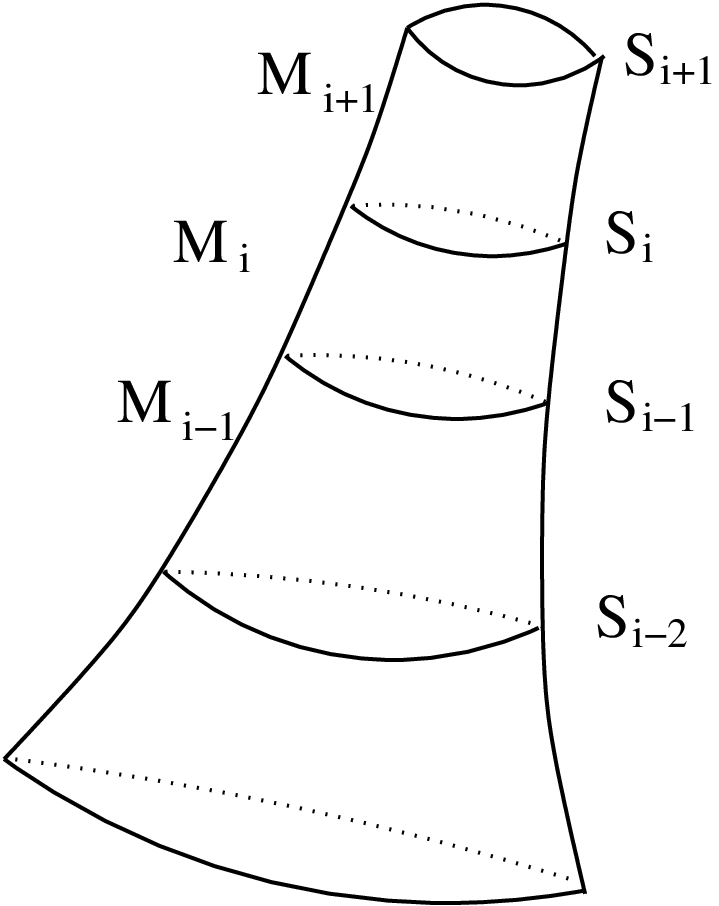}\\[.1in]
		{\em Action of $\Z$ on $\tilde M$, spherical shells and
			fundamental domains}
	\end{center}
\end{figure}

\vspace{-.1in}

We have just proven the following result, originally
due to Ma. \index[persons]{Kato, Ma.} Kato. (\cite{_Kato:announce_}). 

\hfill

\theorem
Let $M$ be a Kato surface, and $S\subset M$ its global
spherical shell. Then $M$ can be obtained
from the Kato data by gluing two boundary components
of $\hat B \backslash B_0$ as above.
\endproof

\hfill

\remark\label{_center_Kato_Remark_}
Using the same arguments, it is possible to show
that the \index[persons]{Kato, Ma.} Kato surface $M$ is not minimal unless
$\hat B$ is obtained by blowing up a point $x_0\in B$,
then blowing up points which lie on exceptional divisors.\index[terms]{divisor!exceptional}
Since a blow-up of an LCK manifold is LCK, we can
always assume that $M$ is minimal and that $\hat B$ is 
obtained from $B$ by successive blow-ups in $0\in B$.\index[terms]{blow-up}

\hfill

\remark\label{_domain_C^2_Remark_}
Suppose that the Kato data satisfy the
assumptions in \ref{_center_Kato_Remark_}.
Then, in the construction of the \index[persons]{Kato, Ma.} Kato surface $M$ from
the Kato data, we can always replace the open ball $B$ of radius
1 by a ball $B(r)\subset B$ of radius $r$ and
the ball $B_0$ by its image $\Psi(B(r))$.
Clearly, the resulting Kato surface is
biholomorphic to $M$. Similarly, the
ball $B$ can be replaced by any holomorphically
convex domain $U\subset \C^2$ with smooth boundary containing the origin $0$.

\hfill

Now we can prove \index[persons]{Brunella, M.} Brunella's theorem.

\hfill

\theorem\label{_Brunella_Theorem_}
\index[terms]{theorem!Brunella}
Let $M$ be a \index[persons]{Kato, Ma.} Kato surface. Then $M$ is locally conformally
K\"ahler.

\hfill

\pstep Let $\pi:\; \hat B \arrow B$ and $B_0\subset B$ be the
Kato data. Denote by $\Psi:\; B\arrow B_0$ the corresponding
biholomorphic equivalence. Choose a K\"ahler metric $\hat \omega$ on
$\hat B$.

\index[persons]{Brunella, M.} Brunella's theorem is proved by finding a metric
$\hat\omega$ on $\hat B$ with the following automorphic
condition. Consider the space $\tilde M$ obtained by
gluing $\Z$ copies of $\hat B \backslash B_0 = M_i$ as
above. A K\"ahler metric $\tilde \omega$ on $\tilde M$ is called
{\bf $\Z$-automorphic} if the deck transform group
mapping $M_i$ to $M_j$ acts on $(\tilde M, \tilde\omega)$
by homotheties. To obtain such a form we need to find
a K\"ahler form\index[terms]{form!K\"ahler} $\hat \omega$ on $\hat B$ 
such that $\hat\omega\restrict{B_0}$ is equal
to $\Psi^*\hat\omega$ in a neighbourhood of
the boundary of $B_0$. If this is true, $\Psi$ acts
by homotheties in a neighbourhood of
$S$. Then the restriction of $\hat\omega$
to $\hat B \backslash B_0 = M_i$
can be extended to a $\Z$-automorphic K\"ahler form on
$\tilde M = \bigcup_{i\in \Z} M_i$.

This is the strategy we follow, except that we replace
$B$ with  another strictly holomorphically convex domain as\index[terms]{domain!holomorphically convex}
in \ref{_domain_C^2_Remark_}.

\hfill

{\bf Step 2:}
Using the local $dd^c$-lemma,\index[terms]{lemma!$dd^c$} we can find a smooth
function $\phi$ on $B_0$ such that
$dd^c\phi=\hat\omega\restrict{B_0}$.
Adding an appropriate affine function, we may assume that
$\phi$ reaches its minimum in a point $x\in B_0$.
Adding an appropriate constant, we can assume that
$U_0:= \phi^{-1}(-\infty, 0)$ has compact closure
with smooth boundary $B_0$.
Let $U$ be the closure of $\Psi^{-1}(U_0)$,
and $\hat U$ the preimage of $U$ under the
bimeromorphic contraction $\hat B \arrow B$.
We obtain $M$ by gluing two boundary
components of $\hat U\backslash U_0$
as in \ref{_domain_C^2_Remark_}.

\hfill

{\bf Step 3:}
The map $\pi:\; \hat U \arrow U$ is proper
map of manifolds of the same dimension.
Therefore, the pushforward of a positive
$(p,p)$-form is a positive $(p,p)$-current.\index[terms]{current}

Let $\pi_*\hat \omega$ be the pushforward of
$\hat \omega$, considered to be  a current on $U$. 
Using the $dd^c$-lemma for
currents, we obtain $\pi_*\hat\omega= dd^c f$, where
$f$ is a plurisubharmonic function on $U$ that is 
smooth outside of the singularities of $\pi$.

Denote by $R$ the boundary of $U\subset \C^2$.
Then $f$ is smooth and strictly plurisubharmonic
in a neighbourhood of $R$.\index[terms]{function!plurisubharmonic}
Another plurisubharmonic function in a neighbourhood of
$R$ is obtained by taking $f_1:=(\Psi^{-1})^*\phi$,
where $\phi$ is the K\"ahler potential on
$U_0$ constructed in Step 2.

Rescaling $f$ if necessary and adding
a constant, we may assume that $-\epsilon <f\restrict R<0$
and $|df|\restrict R \ll \epsilon$. Let $A$ be a sufficiently
big positive number, and $0 < \delta \ll \epsilon$.
Then the regularized maximum\index[terms]{regularized maximum} (see Subsection \ref{_reg_max_subsection_}) $\max_{\delta}$ of
$f$ and $A f_1$
is equal to $A f_1$ in a very small neighbourhood of $R$
(because $f$ is negative on $R$ and $f_1=0$ on $R$),
and equal to $f$ in a neighbourhood $V$ of
$R_\epsilon:=Af_1^{-1}(-2\epsilon)$ because $|df| \ll A |df_1|$
and as $Af_1$ goes to $-2\epsilon$, $f$ does not go
much below $-\epsilon$.

Replacing $\hat \omega$ by
$dd^c \max_{\delta}(f, f_1)$ on the annulus between $R$
and $R_{2\epsilon}$, we obtain a K\"ahler form\index[terms]{form!K\"ahler} $\hat \omega_1$.
Since $\max_{\delta}(f, f_1) = f_1$ in a neighbourhood of $R$,
the map $\Psi:\; (U, \hat\omega_1) \arrow (U_0, \hat \omega_1)$ acts
by homothety mapping a neighbourhood of $R$ with the
metric $\hat \omega_1$ isometrically to a neighbourhood of $\Psi(R)$
with the metric $A\hat \omega_1$.
We have constructed an LCK metric on any Kato surface.
\endproof

\hfill

\remark Prior to publishing his general result
in \cite{_Brunella:Kato_},  M. Brunella constructed examples
of LCK metrics on Enoki surfaces\index[terms]{surface!Enoki}
(\cite{_Brunella:Enoki_}), which
are special cases of class VII surfaces admitting 
a global spherical shell. Even before this,  LCK metrics on some
Kato surfaces (hyperbolic Inoue and parabolic  Inoue)\index[terms]{surface!Kato} also\index[terms]{surface!Inoue!parabolic}\index[terms]{surface!Inoue!hyperbolic}
appeared as anti-self-dual bihermitian metrics\index[terms]{metric!bihermitian}, obtained
by twistor methods in the work of A. \index[persons]{Fujiki, A.} Fujiki and
M. \index[persons]{Pontecorvo, M.} Pontecorvo, see  \cite{_Fujiki_Max:JDG_,
	_Fujiki_Max:SIGMA_}, and \cite{_Fujiki_Max:arxiv016_}, 
where the possibilities for the Lee class (cohomology class of the\index[terms]{class!Lee}
Lee form\index[terms]{form!Lee} $\theta\in \Lambda^1M$) of these
structures is discussed (see also
\cite{_Apostolov_Dloussky_}). As M. \index[persons]{Pontecorvo, M.} Pontecorvo pointed out 
\cite{_Pontecorvo:letter_},
it is implicitly shown there that these Lee classes are different from the
ones constructed by \index[persons]{Brunella, M.} Brunella.

\hfill

\remark\label{_Kato_mflds_Remark_}
Compact complex manifolds with global spherical
shells exist in any dimension (if the complex dimension is
$n$, then the shell has dimension $2n-1$). They are called
{\bf Kato manifolds},\index[terms]{manifold!Kato}
\cite{_Kato:Kinokunya_}. Their construction can be
described in terms of Kato data, as for the surfaces,
see \cite{_Dloussky:Kato_}. Brunella's construction
explained in \ref{_Brunella_Theorem_} can be carried out
in any dimension proving that {\em all Kato manifolds
  arising as a sequence of blow-ups\index[terms]{blow-up} along
  smooth centers admit non-K\"ahler, LCK metrics}; these
LCK metrics are not exact (in particular, they cannot have
potential), unless the manifold is a primary Hopf. This
was first proved in \cite{iop} where examples of
Kato manifolds of any complex dimension $n\geq 4$ carrying
no LCK metric are also constructed (elaborating on 
\index[persons]{Hironaka, H.} Hironaka's example of modification; see also
\cite{shi}). Moreover, in the same paper a special class
of Kato manifolds, associated with so-called Kato matrices
is studied. This particular class of Kato manifolds
generalizes the  Inoue-Hirzebruch
surfaces. The\index[terms]{surface!Inoue-Hirzebruch} following
properties of a  Kato manifold $M_A$ associated with a
certain matrix $A\in\GL(n,\Z)$ are proven in \cite{iop}:
$\kod(M_A)= -\infty$; $M_A$ does not have holomorphic
forms; $M_A$ has holomorphic vector fields; the algebraic
dimension\index[terms]{dimension!algebraic} $a(M_A)$ is greater than 
or equal to the multiplicity of the eigenvalue 1 in the
matrix $A$; $M_A$ does not admit LCK metrics with
potential, in particular they do not admit Vaisman
metrics. 


\section[The embedding theorem in complex dimension $2$]{The embedding theorem in complex\\ dimension $2$}
\label{_embedding_surfaces_pot_}


The proof of the  embedding of LCK manifolds with
potential, \ref{embedding}, relies on an application of the 
 Rossi and  Andreotti--Siu theorem that only works in complex\index[terms]{theorem!Rossi, Andreotti--Siu}
dimension at least $3$. In the following, we give an
argument that works on compact complex surfaces.

\hfill

\theorem (\cite{ov_indam})\label{_sphe_implies_dim2_Theorem_}
Assume that the spherical shell conjecture is true.
Then \ref{potcon} and \ref{embedding} are true in\index[terms]{conjecture!GSS}
dimension 2: for any LCK surface $M$ with (proper) LCK\index[terms]{potential!LCK}
potential, the 1-point completion of its K\"ahler
$\Z$-covering is Stein and $M$ admits a holomorphic
embedding into a linear Hopf manifold. Moreover,
any Vaisman surface also admits such an embedding.

\hfill

\proof The only part of the proof missing for
dimension 2 is the Rossi and  Andre\-o\-tti-Siu Theorem (\ref{asr}), 
\index[terms]{theorem!Rossi, Andreotti--Siu}
which claims that any strictly pseudoconvex CR-manifold
can be realized as the boundary of a Stein variety with isolated
singularities. Were this result true,
we would have used it to prove \ref{_pot_then_ras_Conjecture_} below.

\ref{_pot_then_ras_Conjecture_} is 
stated here as a conjecture, because we do not know
how to prove it for class VII${}_0$ non-Kato surfaces.

For Vaisman surfaces the existence 
of the weak Stein completions
follows directly from \ref{_Stein_completion_cone_Proposition_}.
Indeed, in  \ref{_elli_then_Vaisman_Theorem_}, we prove that
any Vaisman surface is a $\Z$-quotient of an
algebraic cone\index[terms]{cone!algebraic} $\Tot^\circ(L)$
over a projective orbifold, where $L$ is an ample
bundle. The corresponding closed algebraic
cone is the Stein completion of $\tilde M$.
\endproof

\hfill

\begin{conjecture}\label{_pot_then_ras_Conjecture_}
Let $M$ be a compact LCK complex surface with proper
potential, and $\tilde M$ its K\"ahler $\Z$-covering.
Then the metric completion of $\tilde M$ is a
Stein variety.\index[terms]{variety!Stein}\index[terms]{completion!metric}
\end{conjecture}

\hfill

\remark 
For $\dim M\geq 3$, this is \ref{potcon}.

\hfill

\ref{_pot_then_ras_Conjecture_}
follows from the spherical shell conjecture and the
classification of surfaces. We give its proof in this situation as follows;
this also proves \ref{_sphe_implies_dim2_Theorem_}.

\hfill

{\bf Step 1:} Notice that an LCK surface $M$ with an LCK potential\index[terms]{potential!LCK} (be it proper or not) 
cannot contain rational curves. Indeed, if $M$ contains
rational curves, by homotopy lifting $\tilde M$ would also
contain rational curves, but $\tilde M$ is embedded
to a Stein variety. This implies that $M$ cannot
be a Kato surface, and that $M$ is minimal.

{\bf Step 2:} If the spherical shell conjecture is true,\index[terms]{conjecture!GSS}
class VII surfaces that are  not Kato have $b_2=0$. From 
\index[persons]{Bogomolov, F. A.} Bogomolov's theorem (\cite{_Bogomolov:VII_76_,_Bogomolov:VII_82_,
_Li_Yau_Zang:VII_,_Teleman:bogomolov_}), it follows that 
 any class VII surface with $b_2=0$
is biholomorphic to a Hopf surface or to an Inoue surface.\index[terms]{surface!Inoue}\index[terms]{surface!Hopf}

{\bf Step 3:} Inoue surfaces do not have 
LCK potential\index[terms]{potential!LCK} (\cite{oti2}). The Hopf surfaces are quotients
of $\C^2\setminus  0$ by an action of $\Z$, and hence  their 1-point completions
are affine.

{\bf Step 4:} The only non-K\"ahler minimal surfaces that are  not of class
VII are non-K\"ahler elliptic surfaces.
These surfaces are obtained as follows (\ref{_elli_then_Vaisman_Theorem_}).
Let $X$ be a 1-dimensional compact complex orbifold,\index[terms]{orbifold}
and $L$ an ample line bundle on $X$. Consider\index[terms]{bundle!line!ample}
the space $\tilde M$ of all non-zero vectors in the
total space of $L^*$, and let $\Z$ act on $\tilde M$
as $v \mapsto \alpha v$, where $\alpha\in \C$
is a fixed complex number, $|\alpha|>1$. 
Any non-K\"ahler elliptic surface is isomorphic
to $\tilde M/\Z$ for appropriate $\alpha$, $X$ and $L$.
However, the sections of $L^{\otimes n}$ 
define holomorphic functions on $\tilde M\subset
\Tot(L^*)$, identifying $\tilde M$ 
and the corresponding cone over $X$. The completion of this
cone is $\tilde M_c$, and it is affine
(\cite{_EGA2_}, \S 8).
This finishes the proof of \ref{_sphe_implies_dim2_Theorem_}.
\endproof

\hfill

The classification of LCK surfaces with potential\index[terms]{surface!LCK!with potential} would follow
immediately from the classification of Vaisman surfaces, which
is due to F. \index[persons]{Belgun, F. A.} Belgun, if we could prove that any LCK surface\index[terms]{surface!Vaisman}
with potential is a deformation of a Vaisman manifold.
This statement is true in dimension $\geq 3$, but
its proof in dimension 2 is unknown. We finish this
section with the proof of Belgun's theorem.

\hfill

\proposition\label{_Vaisman_is_Hopf_or_elli_Proposition_}
(\cite{bel})\\
Let $M$ be a complex surface. Then $M$ is Vaisman if
and only if $M$ is a Vaisman Hopf surface or an elliptic surface.\index[terms]{surface!Hopf}

\hfill

\proof
By \ref{_elli_or_class_VII_Theorem_}, a non-K\"ahler surface is either principal elliptic
or in class VII or both. Principal elliptic surfaces are Vaisman by \ref{_elli_then_Vaisman_Theorem_}.

If $M$ is Vaisman, it has a deformation which
is principal elliptic (\ref{defovai}). 
A principal elliptic surface $M$ cannot be of class VII,
unless it is a Hopf surface. Since class VII is deformation
invariant, and the class of Hopf surfaces is also\index[terms]{surface!Hopf}
deformation invariant, $M$ is either Hopf or elliptic.
\endproof


\section{Inoue surfaces}\index[terms]{surface!Inoue}


The Inoue surfaces constitute an important class of surfaces
constructed by \index[persons]{Inoue, Ma.} Inoue in 1974 (\cite{inoue}). 
They are classified into three classes: S${}^0$, S${}^+$
and  S${}^-$. We describe the construction of Inoue
surfaces in Chapter \ref{inoue_lck}, where we present 
the results of \index[persons]{Belgun, F. A.} Belgun and \index[persons]{Tricerri, F.} Tricerri on existence and
non-existence of LCK metrics on these surfaces.
The Inoue surfaces of class S${}^0$ can be also
considered to be  \index[terms]{manifold!Oeljeklaus--Toma (OT)} Oeljeklaus--Toma manifolds of lowest
dimension; see Chapter \ref{OT_manifolds}.

\section{Exercises}

The Riemann--Roch formula for the curve is
$\chi(F) = \deg(F) + \chi(\calo_C)\rk F$.
Here we deduce this formula, together with
$\deg(B) = \int_C c_1(B)$ for any vector bundle $B$ on $C$.
However, both the degree and $c_1$ are defined in such a way that 
the Riemann--Roch formula becomes a part of their definition.

\begin{enumerate}[label=\textbf{\thechapter.\arabic*}.,ref=\thechapter.\arabic{enumi}]

\subsection{Riemann--Roch formula for a curve}

\item (``Invariant factors theorem'').\index[terms]{Riemann--Roch formula}
Let $R$ be a principal ideal ring. Prove that any
finitely--generated $R$-module is a direct sum\index[terms]{theorem!invariant factors}
of cyclic $R$-modules. Use this result to deduce
the Jordan normal form theorem, and to classify
the finitely--generated abelian groups.

\item
Let $C$ be a complex curve,
and $x\in C$ a smooth point. Prove that
any coherent sheaf on $C$ supported in $x$
is isomorphic to $\bigoplus_{i=1}^k \calo_C/{\goth m}^{d_i}$,
where ${\goth m}$ is the maximal ideal of $x$, and
$d_1, ..., d_k$ is a collection of positive integers.

{\em Hint:} Use the previous exercise.

\item
Let $C$ be a complex curve, possibly singular,
and $V$ an abelian group, freely generated by
isomorphism classes of coherent sheaves on $C$.
{\bf The Gro\-then\-dieck K-group} $K_0(C)$
is the quotient of $V$ by its subgroup
generated by relations $[F_1]+ [F_3]=[F_2]$
for all exact sequences of coherent sheaves
$0 \arrow F_1 \arrow F_2 \arrow F_3 \arrow 0.$
\index[terms]{group!Grothendieck K}

\begin{enumerate}
\item Let $L$ be a line bundle, and
$0 \arrow \calo_C \arrow L \arrow R \arrow 0$
be an exact sequence associated with a section
$l\in H^0(C, L)$. Prove that $[L]-[\calo_C]= \sum a_i [x_i]$,
 where $a_i \in \Z^{>0}$, $[x_i]$ are classes of skyscraper sheaves
$\calo_C/{\goth m}_{x_i}$, and ${\goth m}_{x_i}$
is the maximal ideal of a point $x_i$.
\item Prove that $K_0(C)$ is generated
by $\calo_C$ and the classes of skyscraper sheaves
$\calo_C/{\goth m}_x$.
\end{enumerate}

\item
Let $C$ be a complex curve, possibly singular,
and $F$ a coherent sheaf on $C$. 
We define {\bf the Euler characteristic of $F$}
as $\chi(F):= \dim H^0(C,F) - \dim H^1(C, F)$.
Prove that $\chi$ defines a group homomorphism
$K_0(C) \arrow \Z$.\index[terms]{Euler characteristic!of a sheaf}

\item
We consider line bundles on a complex curve $C$, not necessarily smooth.
\begin{enumerate}
\item
Let $L$ be a line bundle, admitting a holomorphic section,
and 
\begin{equation}\label{_section_exact_sequence_Equation_} 
0 \arrow \calo_C \arrow L \arrow R \arrow 0
\end{equation}
the corresponding exact sequence.
Define the {\bf degree} $\deg L$ 
as $\chi(L) - \chi(\calo_C)$.
Prove that $\deg(L) = \dim H^0(C, R)$,
when $R$ is defined in the exact sequence 
\eqref{_section_exact_sequence_Equation_}.
\index[terms]{bundle!line!degree of}

\item
Prove that the degree is multiplicative,
$\deg(L_1 \otimes L_2) = \deg(L_1) + \deg(L_2)$,
whenever $L_1$ and $L_2$ are very ample line bundles.
\item Prove
that the function $\deg(L):= \chi(L) -\chi(\calo_C)$
satisfies $\deg(L_1 \otimes L_2) = \deg(L_1) + \deg(L_2)$
for all line bundles on $C$.
\end{enumerate}

\item \label{_RR_curve_Exercise_}
Let $C$ be a complex curve, possibly singular,
and $F$ a coherent sheaf on $C$. 
\begin{enumerate}  
\item Prove that the restriction of $F$ to an open set $U\subset C$
is isomorphic to a vector bundle. Prove that the rank of this vector
bundle is independent on  $U$ when $C$ is irreducible. This number
is called {\bf the rank} of $F$.

\item
We define {\bf the first Chern class}\index[terms]{class!Chern}
$c_1(F)\in H^2(C,\Z)$ as 
\[ c_1(F):= \left (\chi(F)-\sum_i \rk F\restrict {C_i}\chi(\calo_{C_i})\right)\cdot [x],
\]
where $[x]$ is the fundamental class of a point and $C_i$ irreducible components of $C$.
Prove that the first Chern class $c_1$  defines a group homomorphism
$c_1:\; K_0(C) \arrow H^2(C, \Z)$.

\item Prove that $c_1$ satisfies {\bf the Whitney formula}: for any two vector bundles\index[terms]{Whitney formula}
$B_1, B_2$ on a curve, $c_1(B_1\oplus B_2) = c_1(B_1) + c_1(B_2)$.

\item
Prove that $c_1(L)= \deg_C L\cdot [x]$ for any line bundle $L$ on $C$.

\item
Let $B_1, B_2$ be vector bundles on $C$.
Prove that $c_1(B_1 \otimes B_2) = \rk B_1 \cdot c_1(B_2) +
\rk B_2 \cdot c_1(B_1)$.
\end{enumerate}

\subsection{Complex surfaces}

\item
Let $C$ be a smooth complex curve on a non-K\"ahler surface $M$.
Assume that $C\subset M$ admits a  positive-dimensional
family of deformations. Prove that the $C$ is a genus 1 curve.

\item 
Choose a Gauduchon metric $\omega$ on a complex surface $M$, and let
$B$ be a holomorphic vector bundle on $M$. Fix a Hermitian metric $h$ on $B$, and
let $\Theta_B$ be the curvature of its Chern connection.\index[terms]{connection!Chern}
 Define {\bf the degree} $\deg_\omega B:= \frac1 {2\pi\1}\int_M \omega \wedge \Tr_B\Theta_B$,
where $\Tr_B:\; \Lambda^2(M) \otimes \End B\arrow \Lambda^2(M)$
is the trace map. 
\begin{enumerate}
\item Prove that $\deg_\omega B$ is independent
from the choice of $h$.
\item Prove that $\deg_\omega B$ is integer when 
$\omega$ is K\"ahler and its cohomology class is integer.
\item Let $M$ be a Hopf surface. Prove that
for any $\lambda \in \R$ there exists a line bundle $L$
such that $\deg_\omega L=\lambda$.\index[terms]{surface!Hopf}
\end{enumerate}

\item 
Choose a Gauduchon metric $\omega$ on a complex surface $M$, and let
$L$ be a holomorphic line bundle on $M$. Assume that $H^0(M,L)\neq 0$.\index[terms]{metric!Gauduchon}
\begin{enumerate}
\item
Prove that $\deg_\omega(L)>0$.
\item 
Suppose that $L$ admits a holomorphic
section with zero divisor $Z$. Assume that
$Z$ has multiplicity 1. Prove that 
$\deg_\omega(L) = \int_Z \omega$.
\end{enumerate}

\item
Let $M$ be a linear Hopf surface,\index[terms]{surface!Hopf!linear}
that is, a quotient $\frac{\C^2 \backslash 0}{\langle A\rangle}$
with $A$ a linear contraction. Denote by $G$ the connected component
of the group $\Aut(M)$ of holomorphic automorphisms.
\begin{enumerate}
\item Prove that $G$ contains a subgroup isomorphic to
$\C^* \times \C^*$ if $A$ is diagonal.
\item Prove that
$G$ contains a subgroup isomorphic to
$\C^* \times \C$ if $A$ is a non-diagonal Jordan block.
\end{enumerate}

\item
\begin{enumerate}
\item Let $f:\; \C^m\backslash 0 \arrow \R$ take
$z$ to $\log|z|$. Prove that $f$ is plurisubharmonic.
\item Let $\phi:\; \C^n \arrow \C^m$ be a holomorphic
map, and $g(z):= \log |\phi(z)|$. Prove that $g$ is plurisubharmonic.
\item
Let $f_1, ..., f_m$ be a collection of holomorphic functions on $\C^n$.
Prove that $\log \left(\sum_{i=1}^m |f_i|\right)$ is plurisubharmonic.
\item  Let $\beta_1, ..., \beta_n \in \R^{\geq 1}$ 
be real numbers, and $h$ a smooth function on
$\C^n \backslash 0$, defined as 
$h(z_1, ..., z_n):=\log \left(\sum_{i=1}^n |z_i|^{\beta_i}\right)$.
Prove that $h$ is plurisubharmonic.
\end{enumerate}

\item \label{_linear_Hopf_canonical_explicitly_Exercise_}
Let $M$ be a linear Hopf manifold, 
$M=\frac{\C^n \backslash 0}{\langle A\rangle}$,
where $A(e_i) =\alpha_i e_i$, with
$\alpha_i \in \C$ and $e_i$ the basis in $\C^n$.
Recall that a function $f\in C^\infty(\C^n\backslash 0)$
is called {\bf automorphic} if $A^* f= \lambda f$ for some
$\lambda\in \R^{\geq 0}$.
\begin{enumerate}
\item Let $\phi(z_1, ..., z_n)= \sum_{i=1}^n |z_i|^{\beta_i}$,
with $|\alpha_i|^{\beta_i} = |\alpha_1|^{\beta_1}$ for $i=2, ..., n$,
and $\beta_i \in \R^{>1}$. Prove that such $\beta_i$ exist.
Prove that $\phi$ is an automorphic plurisubharmonic function
on $\C^n\backslash 0$. 
\item Prove that $dd^c\log \phi$ is $A^*$-invariant, 
positive (1,1)-form\footnote{French-positive, see 
\ref{_French-positive_Remark_}.}  on $\C^n\backslash 0$.
\item Choose a Vaisman metric on $M$, and let $\Sigma\subset TM$
be the canonical foliation.\index[terms]{foliation!canonical} Prove that $\ker dd^c\log \phi\supset \Sigma$.
\item Let $\log A_\R$ be the linear map given by 
$\log A_\R(e_i) =\log |\alpha_i| e_i$. Prove that the leaves of 
$\Sigma$ are the orbits of $e^{\C \log A_\R}$.\footnote{This argument
is referenced in \ref{_canonical_explicit_Hopf_Proposition_}.}
\end{enumerate}

\item
Let $M$ be a linear Hopf surface,\index[terms]{surface!Hopf!linear}
$M=\frac{\C^2 \backslash 0}{\langle A\rangle}$. Suppose that $A$ is diagonalizable.
\begin{enumerate}
\item Prove that $M$ is Vaisman.
Let $S$ be a leaf of the canonical foliation
on $M$, and $\tilde S\subset \C^2 \backslash 0$ its 
preimage. Prove that $\tilde S$ is an orbit of the group
$e^{\C \log A_\R}$, where the action of $e^{\C \log A_\R}$
is defined in Exercise \ref{_linear_Hopf_canonical_explicitly_Exercise_}.

\item Let $E\subset M$ be an elliptic curve,
and $\tilde E\subset \C^2 \backslash 0$ its pullback.
Prove that $\tilde E$ is an orbit of $e^{\C \log A_\R}$.

\item Prove that $M$ is elliptic if the eigenvalues
$\alpha_1, \alpha_2$ of $A$ satisfy $\alpha_1^n = \alpha_2^m$
for some integers $m, n \neq 0$. 

\item 
Consider the action of the group $\Z=\langle A\rangle$ on the set
$S$ of  orbits of the $e^{\C \log A_\R}$-action.
Prove that any elliptic curve on $M$ is obtained
as a quotient of an orbit $s\in S$ of the $e^{\C \log A_\R}$-action
by its stabilizer $\St_{\langle A\rangle}(s)$.

\item Prove that all orbits of $e^{\C \log A_\R}$ on 
$\C^2 \backslash 0$ are biholomorphic to $\C$ or $\C^*$.
Prove that an elliptic curve $E\subset M$  
is a quotient of an orbit $\tilde E$ of $e^{\C \log A_\R}$
by $\St_{\langle A\rangle}(\tilde E)$. Prove that
$\tilde E\cong \C^*$. Prove that its closure 
in $\C^2$ is an eigenspace of  $A^d$, for some $d \in \Z^{>0}$.

\item Prove that $M$ is non-elliptic  if
the eigenvalues $\alpha_1, \alpha_2$ of $A$ 
do not satisfy $\alpha_1^n = \alpha_2^m$
for some integers $m, n \neq 0$.  

\end{enumerate}

\item 
Let $M$ be a linear Hopf surface,\index[terms]{surface!Hopf!linear}
$\frac{\C^2 \backslash 0}{\langle A\rangle}$. 
Suppose that $A$ is non-diagonalizable.
\begin{enumerate} 
\item Prove that $M$ is not Vaisman.
\item Prove that $M$ is not elliptic. Prove
that $M$ contains an elliptic curve.
\end{enumerate}

\item
Let $M$ be a linear Hopf surface,
$\frac{\C^2 \backslash 0}{\langle A\rangle}$.
\begin{enumerate}
\item Suppose that $\Aut(M)$ acts on $M$ transitively.
Prove that $M$ is elliptic. \index[terms]{surface!Hopf!elliptic}
\item Prove that an elliptic Hopf surface $M=\frac{\C^2 \backslash 0}{\langle A\rangle}$ has special fibres 
if and only if $A\neq\const\cdot \Id$.
\item Prove that the group $\Aut(M)$ acts on $M$ transitively
if and only if $A=\const\cdot \Id$.
\end{enumerate}

\item
Choose a Gauduchon metric $\omega$ on a complex surface $M$, and let
$K_M$ be its canonical bundle. Assume that $M$ admits two non-proportional 
holomorphic vector fields. Prove that $\deg_\omega(K_M)<0$.

\item
Choose a Gauduchon metric $\omega$ on a complex surface $M$, and let
$K_M$ be its canonical bundle. Prove that
\begin{enumerate}
\item $\deg_\omega(K_M)<0$ if $M$ is a Hopf surface.
\item $\deg_\omega(K_M)=0$  if $M$ \index[terms]{surface!Hopf} \index[terms]{surface!Kodaira} is a Kodaira surface.
\item $\deg_\omega(K_M)>0$ if $M$ is a non-K\"ahler elliptic
surface with $b_1 > 3$ or its non-ramified finite quotient.
\end{enumerate}

\item Let $C$ be a smooth complex curve on a complex surface $M$
with $b_1(M)=1$ and $b_2(M)=0$. Prove that $C$ is a genus
1 curve.

\item
Let $\alpha\in \Lambda^{1,0}(M)$ be a  $\6$-closed (1,0)-form on a complex surface.
Prove that $\alpha= \alpha_0 + \6 f$, where $\alpha_0$ is holomorphic,
and $f\in C^\infty M$ is a function.

\item
Let $M$ be a complex manifold, 
and $L$ a holomorphic line bundle on $M$ with $c_1(L)=0$. 
\begin{enumerate}
\item Prove that $L$ admits a smooth section $\beta$ without zeros.
\item Let $h$ be the Hermitian metric on $L$ such that $|\beta|=\const$,
and $\nabla_\beta$ the corresponding  Chern connection\index[terms]{connection!Chern}. Prove that
$\nabla_\beta(\beta) = \bar\6(\beta) +  \beta\otimes \alpha$,
where $\alpha\in \Lambda^{1,0}(M)$ is a $\6$-closed (1,0)-form on $M$.

\item Let $f$ be a function on $M$, and $\beta_1:=e^f \beta$.
Prove that $\nabla_{\beta_1}(\beta_1)= \bar\6(\beta_1)+ \beta_1\otimes (\alpha + \6 f)$.

\item Let $\nabla$ be a Chern connection on $M$ such that
$\nabla_\beta(\beta) = \bar\6(\beta) +  \beta\otimes \alpha$,
where $\alpha$ is a holomorphic 1-form. Prove that $\nabla$ is flat.

\item Suppose that $M$ is a complex surface.
Prove that $L$ admits a flat  Chern connection\index[terms]{connection!Chern} compatible with the holomorphic structure.
\end{enumerate}

{\em Hint:} Use the previous exercise.

\item
 Suppose that $M$ is a complex surface, 
$\Pic_\fl$ the group of isomorphism classes of complex 
line bundles with flat connection on $M$, and\index[terms]{connection!flat}
$\Pic_0$ the group of isomorphism classes of holomorphic line
bundles with $c_1=0$. A flat connection $\nabla$ on $L$ defines a 
holomorphic structure $\bar\6:=\nabla^{0,1}$.
Consider the natural forgetful map
$\tau:\; \Pic_\fl\arrow \Pic_0$ taking $\nabla$ to $\bar\6=\nabla^{0,1}$.

\begin{enumerate}
\item Prove that $\Pic_\fl= H^1(M, \C_M^*)$,
where $\C_M^*$ is a constant sheaf that satisfies
$\C_M^*(U)=\C^*$ for each connected open set $U\subset M$.

\item Let $\Z_M$ and $\C_M$ be the standard 
constant sheaves on $M$.
Consider the exponential exact sequences
\[ \begin{CD}
0 @>>> \Z_M  @>e>> \C_M  @>\exp>> \C_M^* @>>> 0\\
&& @VVV @VVV @VVV &\\
0 @>>> \Z_M  @>e>> \calo_M   @>\exp >> \calo_M^* @>>> 0
\end{CD}
\]
where $e(1)= 2\pi\1$.
Prove that the following diagram is commutative and has exact rows:
\[
\begin{CD}
@>>> H^1(M, \Z) @>e>> H^1(M, \C) @>>> \Pic_\fl @>>> 0\\
&& @VVV @VVV @VVV &\\
@>>> H^1(M, \Z) @>e>> H^1(M, \calo_M) @>>> \Pic_0 @>>> 0
\end{CD}
\]

\item Prove that the kernel of $\tau:\;  \Pic_\fl\arrow \Pic_0$ 
is equal to the kernel of the natural
homomorphism $H^1(M)\arrow H^{0,1}(M)$.
\item Prove that $\tau$ is an isomorphism when $b_1(M)=1$,\footnote{This result
was originally proven by \index[persons]{Kodaira, K.} Kodaira, \cite[page 57]{_Kodaira_Structure_III_}.}
and not an isomorphism when $b_1(M) >1$.
\end{enumerate}

\item
Let $L$ be a flat line bundle on a complex surface $M$.
Prove that the dimension of the kernel of 
the natural map $H^{1,1}_{BC}(M,L) \arrow H^2(M, {\Bbb L})$
is $H^0(M, {\Bbb L})$, where ${\Bbb L}$ is the
local system associated with $L$.

{\em Hint:} Extend the argument used to prove
\ref{_BC_Degree_Theorem_} to (1,1)-forms with coefficients in a line bundle.

\hspace{-3em} The following sequence of exercises is based on \index[persons]{Bogomolov, F. A.} Bogomolov's
work on the classification of class VII surfaces 
\cite{_Bogomolov:VII_76_,_Bogomolov:VII_82_,_Bogomolov_Buonerba_Kurnosov_}.

\item
Let $M$ be a complex surface with $b_1=1$,
and $L$ a line bundle. 
\begin{enumerate}
\item Prove that $\chi(\calo_M)=0$.
\item Prove that $\chi(L)=0$ for any line bundle with $c_1(L)=0$.
\item Assume that $b_2=0$. Prove that $H^0(M,L)=0$
if and only if $H^0(L^*\otimes K_M)=0$.
\end{enumerate}

\item Let $M$ be a complex surface with $b_1=1$ and
$b_2=0$. 
\begin{enumerate}
\item Prove that $c_2(M)$ is torsion.
\item Let $F$ be a tensor bundle on $M$,
$F= TM^{\otimes n} \otimes T^*M^{\otimes m}$.
Prove that $H^0(M,F)=0$
if and only if $H^0(F^*\otimes K_M)=0$.
\end{enumerate}

{\em Hint:} Prove that the  Chern classes\index[terms]{class!Chern} of $F$ are zero,
and find the Euler characteristic of $F$. Then use the Serre duality\index[terms]{duality!Serre}
to find $\chi(F)=\dim H^0(M,F)-\dim H^0(F^*\otimes K_M)$.


\item Let $M$ be a complex surface. Assume that the 
holomorphic bundle $TM$ has no proper coherent subsheaves.
Let $L$ be a line bundle on $M$.

\begin{enumerate}
\item Prove that $H^0(\Omega^1 M \otimes L)=0$
for any line bundle $L$.

\item 
Prove that the determinant $\det l\in L^{\otimes 2}$ 
of any non-zero section $l \in H^0(L\otimes \End(TM))$ 
is a non-zero holomorphic section of $L^{\otimes 2}$ .

\item Assume that $M$ is a surface of class VII.
Prove that \[ H^0(K_M\otimes \End(TM))=0.\]
\end{enumerate}

\item\label{_currents_basic_Exercise_}
Let $\Theta$ be a positive, closed (1,1)-current\index[terms]{current} on a 
Vaisman $n$-manifold $M$, and $\omega_0$ the transversally
K\"ahler form\index[terms]{form!K\"ahler} (\ref{_Subva_Vaisman_Theorem_}).
\begin{enumerate}
\item Assume that $\omega_0^{n-1}\wedge \Theta=0$.
Prove that $\Theta$ is basic with respect to the
canonical foliation $\Sigma$.\index[terms]{foliation!canonical}
\item Let $\Theta$ be a current that is  basic with respect to the
canonical foliation. Prove that $\Theta$ is invariant with respect
to the action of the Lie group generated by the Lie field.
\item Assume that $M$ is a compact Vaisman manifold.\index[terms]{manifold!Vaisman}
Prove that $\int_M \omega_0^{n-1}\wedge \Theta=0$.
Deduce that $\omega_0^{n-1}\wedge \Theta=0$.
\end{enumerate}

\end{enumerate}

\chapter{Cohomology of holomorphic bundles on Hopf
  manifolds}
\label{_cohomo_on_Hopf_Theorem_}
\epigraph{\it A definition worth its name - mathematicians learned this from Alexander Grothendieck - is not a concise wording of what everyone knows, but a pointer
	toward unknown. Many unexpected fruits have been harvested from what grew
	from the seeds of ideas in his definitions.}
{\sc\scriptsize Misha Gromov, Great Circle of Mysteries}

\section{Introduction}

In this chapter, we relate the results of D. \index[persons]{Mall, D.} Mall 
who computed the cohomology of the tangent bundle
on Hopf manifolds (\cite{_Mall:Contractions_}; see also 
\cite{_Libgober_,_Gan_Zhou_}). Eventually this will
allow us to express the \index[persons]{Kuranishi, M.} Kuranishi space of 
Hopf manifolds in explicit terms (Chapter \ref{_Kuranishi_Chapter_}).

Throughout this book,
we mostly used {\bf the linear Hopf manifolds}, that are 
the quotients of $\C^n\backslash 0$ by a linear contraction.
Since this chapter, we relax this condition by allowing
any invertible holomorphic contraction.\index[terms]{manifold!Hopf}
Note that the ``non-linear'' Hopf manifolds 
are LCK with potential \index[terms]{manifold!LCK!with potential}(\ref{_gene_Hopf_is_LCK_Corollary_})
and admit a holomorphic embedding to a linear Hopf manifold
(\ref{_Stein_by_contract_to_linear_Hopf_Theorem_}). 

Most of this chapter is taken by an introduction
to homological algebra and to $G$-equivariant sheaves;
in both cases, we follow \index[persons]{Grothendieck, A.} Grothendieck's brilliant
exposition in the T\^ohoku paper (\cite{_Grothendieck:Tohoku_}). 
However, the most difficult result here is not Mall's theorem
or Grothendieck spectral sequence\index[terms]{spectral sequence!Grothendieck}, but the calculation
of the cohomology groups $H^*(\C^n\backslash 0, \calo_{\C^n})$.
These cohomology groups were computed in \cite{_Grauert_Remmert:Stein_}
by an explicit calculation, using induction and
the \v Cech resolution; we give a new proof, based 
on Serre duality.

The cohomology of an annulus
$C_{a,b}:=\{z\in \C^n \ \ |\ \ a< |z| <b\}$ can be easily computed using  Serre duality
and the exact sequence
\begin{multline}\label{_exact_annulus_intro_Equation_}
... \arrow H^i_c(B_a, \calo_{\C^n}) \arrow \\
\arrow H^i_c(B_b, \calo_{\C^n}) \arrow 
H^i(C_{a,b}, \calo_{\C^n}) \arrow  H^{i+1}_c(B_a, \calo_{\C^n}) \arrow...
\end{multline}
where $B_a, B_b$ are open balls of radii $a, b$, and
$H^i_c(B_a, \calo_{\C^n})$ is the cohomology with compact support.
Since $H^i_c(B_a, \calo_{\C^n})$ is dual to $H^{n-i}(B_a, \calo_{\C^n})$,
and $B_a$ is Stein, this group vanishes for all $i\neq n$.
Then, by \eqref{_exact_annulus_intro_Equation_}, 
the group $H^i(C_{a,b}, \calo_{\C^n})$ vanishes
for all $i \neq 0, n-1$. 

The passage from the annulus to $\C^n\backslash 0$ is much less
trivial, and involves the \index[persons]{Grothendieck, A.} Grothendieck spectral sequence
associated with the map $\C^n \backslash 0 \xlongrightarrow{x\mapsto |x|} \R^{>0}$.
We show that the pushforward of $\calo_{\C^n}$ under this map
is a special kind of sheaf on $\R^{>0}$, called {\bf a directed sheaf},
and compute the cohomology of directed sheaves.\index[terms]{sheaf!directed}

The reader can skip the homological algebra and proceed directly to
Section \ref{_cohomo_C^n_without_0_Section_}, where the 
proof of \index[persons]{Mall, D.} Mall's theorem starts. However, the Grothendieck
spectral sequence is referred to in several other chapters in this book.

Another subject we touch is\index[terms]{cohomology!group}
the group cohomology. We give two definitions
of group cohomology $H^*(G,V)$ with coefficients in a $G$-module $V$.
First, we define $H^*(G,V)$ as the cohomology of the corresponding
local system on the \index[persons]{Eilenberg, S.} Eilenberg-\index[persons]{MacLane, S.}MacLane space $K(G, 1)$.\index[terms]{space!Eilenberg-MacLane}
Then we show that this definition is equivalent to the 
definition obtained if we consider $V$ as a module
over the group ring $\R[G]$.\index[terms]{group ring} In this setup, the 
cohomology groups are defined as $H^i(G, V):= \Ext^i_{\R[G]}(\R, V)$.
We use the equivalence of these two definitions to compute
the cohomology of $\Z$ with coefficients in a $\Z$-module.
This computation is used to finish the proof of Mall's theorem;
however, the general result about the group cohomology is also very useful.

\section[Derived functors and the Grothendieck spectral sequence]{Derived functors and the Grothendieck\\ spectral sequence}
\label{_Derived_functors_Subsection_}
\index[terms]{derived functor}
\index[terms]{derived category}\index[terms]{spectral sequence!Grothendieck}

For the benefit of the reader, we briefly introduce the derived
categories, which serve as a uniform way of obtaining almost all
cohomology functors on abelian categories. We introduce
the \index[persons]{Grothendieck, A.} Grothendieck spectral sequence, that is  used to describe
the derived functors of the composition of two left or right
exact functors $\Psi, \Phi$ in terms of derived  functors of
$\Psi$ and $\Phi$. This spectral sequence is used further on
in two independent aspects: we invoke \index[persons]{Grothendieck, A.} Grothendieck's spectral
sequence to describe the $G$-equivariant cohomology in this
chapter, and to compute cohomology of a sheaf in terms
of higher derived images (Section \ref{_Hopf_cohomo_Section_}).

We assume the functional knowledge of abelian categories;\index[terms]{category!abelian}
like with  set theory, most people are using the abelian
categories without thinking of its axiomatic definition,
and this attitude is sufficient for the present purposes.
The best (at least, the shortest) introduction to abelian
categories is found in Grothendieck's T\^ohoku paper,
\cite{_Grothendieck:Tohoku_}, 
who introduced them in this
paper, proved the basic results and gave the necessary 
counterexamples.

We also assume the definitions of projective and injective 
objects.

The notion of derived categories is due to 
A. \index[persons]{Grothendieck, A.} Grothendieck and J.-L. \index[persons]{Verdier, J.-L.} Verdier, who introduced them
in his Ph. D. thesis (\cite{_SGA41/2_,_Verdier:thesis_});
for a modern introduction to the subject, we recommend
\cite{_Gelfand_Manin_}. 

We start with defining  localization of
  categories. It is a notion very similar to the\index[terms]{localization}
localization in algebras, that is, formally taking
the fractions with denominators in a multiplicatively
closed set. Similar to the localization
in non-commutative algebras, to localize a category
you need {\bf the Ore condition}. In algebras,\index[terms]{condition!Ore}
one says that the multiplicatively closed
subset $S$ in an algebra $A$ {\bf satisfies the
Ore condition} if for each $a\in A$ and each $s\in S$,
there exists $a_1\in A$, and $s_1\in S$ 
such that $as_1=sa_1$. Then one can permute the
formal fractions, using the relations such as 
$s^{-1} a= a_1s_1^{-1}$, allowing one to multiply the
fractions. Each pair $(A, S\subset A)$ of an 
algebra and the multiplicatively closed subset
satisfying the Ore condition gives rise
to {\bf the right localization} $A[S^{-1}]$,
that is  an algebra generated by $A$ and with
the right multiplication by $s\in S$ inverted
(\cite{_Ranicki:Localization_}).

Localization of categories is a similar concept,\index[terms]{category!localization of}
where one can invert a class $S$ of morphisms obtaining
a category with the same objects and a greater set
of morphisms, in such a way that all $S$ become isomorphisms.

We follow 
\cite[Section 4.27, tag/04VB]{_Stacks_Project_}.
Let ${\cal C}$ be a category, and $S$ a set of morphisms
of ${\cal C}$, closed under composition and containing all
isomorphisms. We say that
$S$ is a {\bf right multiplicative system} if
the following two conditions are satisfied:
the right Ore condition RMS2, and the right 
cancellability property, RMS3. Here the RMS2 means that
every solid diagram \index[terms]{cancellability property}
\begin{equation}\label{_solid_dia_Equation_} 
\begin{diagram}[size=2em]
X & \rTo^{g} & Y \\
\dTo^t & & \dDots^s\\
	Z & \rDots_f & W
\end{diagram}
\end{equation}
(that is, the part of \eqref{_solid_dia_Equation_} 
that has $t, g, X, Y,$ and $Z$)
where $t \in S$, can be completed to a commutative diagram
by the dotted arrows, with $s\in S$. The RMS3 condition 
declares that for any 
pair of morphisms $a, b\in \Mor(X,Y)$ and any
$s\in \Mor(Y, Z)$ that belongs to $S$,
satisfying $a\circ s=b\circ s$, there exists
a morphism $s' \in \Mor(W, X)$ such that
$s' \circ a= s' \circ b$:
\[
W \stackrel {s'} \arrow X \doublerightarrow{a}{b}Y
\stackrel {s} \arrow Z.
\]
When these conditions are satisfied,
the category $\cac$ admits a functor
to the {\bf right localized category}
$\cac[S^{-1}]$, with the same set of objects,
and all morphisms from $S$ mapped to 
isomorphisms. The localized category
can be defined as the universal category
among all categories that satisfy this
property.\index[terms]{category!localized}

To define the derived categories, we start
with the category of complexes in an abelian category.
Given a morphism $\Psi:\; (A^*, d_A) \arrow (B^*, d_B)$ of complexes,
one defines {\bf the cone} $C(\Psi)$, that is  the\index[terms]{cone!of a morphism}
complex $A^{*-1} \oplus B^{*}$ with the differential
$-d_A+d_B+\Psi$; it is easy to see that 
$(-d_A+d_B+\Psi)^2=0$ if and only if $\Psi$ commutes
with the differentials, that is, defines a morphism
of complexes.  The sequence of complex 
morphisms $... \arrow A^* \arrow B^* \arrow C(\Psi)\arrow A^{*-1}\arrow ...$
is called {\bf a distinguished triangle};
the composition of any two consecutive maps in
a distinguished triangle is zero.

There exists a long exact sequence
of cohomology
\begin{equation} \label{_cone_long_exact_Equation_}
... \arrow H^i(A^*) \arrow H^i(B^*) \arrow H^i(C(\Psi))
\arrow H^{i+1}(A^*) \arrow ...
\end{equation}
For any complex $K$, denote by $\Hom(K^*, A^*)$ the
space of morphisms of complexes. Then, in addition to
\eqref{_cone_long_exact_Equation_}, there is a 
sequence
\begin{multline} \label{_cone_hom_long_exact_Equation_}
... \arrow \Hom(K^*, A^*) \arrow \Hom(K^*, B^*)\arrow\\
 \arrow \Hom(K^*, C(\Psi))
\arrow \Hom(K^*, A^*[1]) \arrow ...
\end{multline}
where $C^*[a]$ denotes the complex $C^*$ shifted by $a\in \Z$.
We leave the proof of exactness of these sequences
as an easy exercise.

A morphism $\Psi$ is called {\bf quasi-isomorphism}
if it induces an isomorphism on cohomology; \index[terms]{quasi-isomorphism}
from \eqref{_cone_long_exact_Equation_} it follows
that $\Psi$ is a quasi-isomorphism if and only if
$C(\Psi)$ is {\bf acyclic}, that is, has zero cohomology.\index[terms]{sheaf!acyclic}

The set of quasi-isomorphisms in the category of complexes
is clearly multiplicatively closed.
The right (and left) Ore condition RMS2 can be obtained
if we take $W= C(g\oplus t)$ in the diagram \eqref{_solid_dia_Equation_}.
Then $W=X^{*-1}\oplus Y^{*} \oplus Z^{*}$,
and the natural map $Y \arrow W$ is a quasi-isomorphism
because $X^*\oplus Z^{*-1}= C(t)$ has zero cohomology,
being a cone of a quasi-isomorphism.

To prove the right cancellation condition RMS3,
it suffices to check that for any morphism 
of complexes $f:\; X \arrow Y$ and a quasi-isomorphism
$s:\; Y \arrow Z$ such that $s\circ f =0$,
there exists a quasi-isomorphism $s':\; W \arrow X$
such that $f\circ s'=0$. To obtain such a quasi-isomorphism,
we consider the natural morphism 
$p:\; C(s)[-1]\arrow Y$, and the  
exact sequence \eqref{_cone_hom_long_exact_Equation_}:
\begin{multline}\label{_long_exact_for_cancellation_Equation_}
... \arrow \Hom(X, C(s)[-1]) \stackrel p \arrow \Hom(X, Y) \stackrel s \arrow 
\Hom(X, Z) \arrow \\ \arrow \Hom(X, C(s)) \arrow ...
\end{multline}
Consider $f \in \Hom(X,Y)$ as an element in the second
term of this exact sequence. Since $s\circ f=0$, 
the image of $f$ in $\Hom(X, Z)$ vanishes. From the exactness
of \eqref{_long_exact_for_cancellation_Equation_}, 
we obtain that $f= u\circ p$, where $u\in \Hom(X,C(s)[-1])$.
Since $s$ is a quasi-isomorphism, the complex
$C(s)$ is acyclic, and the natural map 
$s':\; C(u)\arrow X$ is a quasi-isomorphism.
However, $f\circ s'=0$ because $f\circ s'=p\circ u\circ s'$,
and the composition $u\circ s'$ of two consecutive arrows
in a distinguished triangle vanishes.
\begin{equation}
	\begin{minipage}{0.85\linewidth}
\xymatrix{
	C(u)\ar[rd]_{s'}&&&\\
	&X\ar[r]^f \ar[dr]^u &Y\ar[r]^s &Z\\
	&&C(s)[-1]\ar[u]_p&
}
\end{minipage}	
\end{equation}
This implies that quasi-isomorphisms define a right
(also, a left) multiplicative system.

Let $\psi, \phi:\; (A^*, d_A) \arrow (B^*, d_B)$
be morphisms of complexes. They are called {\bf homotopic}
if there exists an operator $h:\; A^*\arrow B^{*-1}$
such that $ h d_A + d_Bh= \psi-\phi$. In this case,
$h$ is called {\bf a homotopy operator}.
Denote by $\Hot(A^*, B^*)$ the group of morphisms
from $A^*$ to $B^*$ modulo homotopy.\index[terms]{category!homotopy}
The {\bf homotopy category} $\Hot({\cal A})$ over an abelian category
${\cal A}$ is the category with objects complexes
in ${\cal A}$, and morphisms 
$\Mor_{\Hot({\cal A})}(A^*, B^*):=\Hot(A^*, B^*)$.

Let $\cac$ be an abelian category, and $\Hot(\cac)$ its
homotopy category.
The {\bf derived category} $D(\cac)$ is obtained from $\Hot(\cac)$
category by inverting (localizing) the quasi-isomorphisms.\index[terms]{quasi-isomorphism}
This makes sense, because the quasi-isomorphisms
satisfy the right multiplicative system axioms
RMS2 and RMS3, as we have shown. 
Typically, one would consider a smaller
category $D^+(\cac)$ or $D^-(\cac)$
that is  called {\bf bounded from the left}\index[terms]{category!bounded from the right/left}
and {\bf bounded from the right}; these are
derived categories of complexes $C^*$ that have
vanishing cohomology $H^i(C^*)=0$ for $i\ll 0$
(bounded from the left) and $H^i(C^*)=0$ for $i\gg 0$
(bounded from the right) 
\cite[Section 13.11, tag/05RR]{_Stacks_Project_}.

To explain the meaning of this category, recall that
any complex $(P^*, d_P)$ of projective objects of an
abelian category, bounded from the right
and with the vanishing cohomology, admits a homotopy map
$h:\; P^*\arrow P^{*-1}$ satisfying $hd_P+ d_Ph= \Id_{P^*}$.%
\footnote{A complex $...\arrow C^i \arrow C^{i+1} \arrow ...$
is called {\bf bounded from the right} if all $C^i$ vanish for
$i \gg 0$, and {\bf bounded from the left} if all $C^i$ vanish for
$i \ll 0$.}
This implies that $(P^*, d_P)$ splits onto a direct sum
of length two complexes 
$\arrow Q_i  \tilde \arrow Q_{i+1} \arrow 0$,
with $Q_i \cong Q_{i-1}$.
The same statement is true for a complex of injective 
objects, bounded from the left and with the vanishing
cohomology. Then the identity morphism from $(P^*, d_P)$ 
to itself is homotopic to zero; in the corresponding
homotopy category, $(P^*, d_P)$ corresponds to the 
trivial object.

If we apply this argument to the cone $C(s)$ of a quasi-isomorphism,
we obtain that the identity map in $C(s)$ is homotopic to zero,
which implies that $s$ is an isomorphism in the homotopy category.

From this argument, we see that the derived category 
of appropriately bounded complexes is equivalent to the
homotopy category of complexes of projective (or injective)\index[terms]{sheaf!projective/injective}
objects, assuming that any object has a projective (or injective)
resolution. However, the derived categories are more flexible, and
can be used even if such resolutions do not exist.\index[terms]{resolution!injective/projective}

Let $F:\; \cac \arrow \cac_1$ be a functor of abelian
categories, that is  left exact (i.e. it maps an exact sequence to
a sequence that is  exact in the first two terms)
or right exact  (maps an exact sequence to
a sequence that is  exact in the last two terms).
The corresponding {\bf derived functors} (the {\bf left derived
functor} for the right exact functor $F$ and vice versa)
are functors on the derived categories 
$R^* F:\; D(\cac) \arrow D(\cac_1)$. The derived functors are 
defined whenever we can replace a given complex $K$ by
a quasi-isomorphic complex $I(K)$ of injective objects (for the left exact $F$)
or the quasi-isomorphic complex $P(K)$
of projective objects for the right exact $F$.
Then  $R^*F(K) = F(I(K))$ for $F$ left exact, and
$L^*F(K)= F(P(K))$ for $F$ right exact.\index[terms]{derived functor}

Since the derived category (for appropriately bounded complexes)
is isomorphic to the homotopy category of projective or injective complexes,
this construction actually defines a functor on the corresponding
derived category.\index[terms]{derived category}

The traditional way of defining derived functors is the
following. Let $X$ be an object of an abelian category,
and $C^*$ the complex of projective or injective modules,
quasi-isomorphic to the complex $0 \arrow X \arrow 0$.
This complex is obtained by taking a projective or\index[terms]{resolution!injective/projective}
injective resolution of $X$. Given (say) a left exact
 functor $F$, the right derived functor $R^iF(X)$
is defined as the $i$-the cohomology of the complex 
$F(C^*)$. This definition is due to \index[persons]{Cartan, H.} Cartan and Eilenberg 
(\cite{_Cartan_Eilenberg:Derived_}), who invented the term
and the concept. In his T\^ohoku paper \cite{_Grothendieck:Tohoku_}, 
\index[persons]{Grothendieck, A.} Grothendieck  created derived categories, making the definition 
of derived functors truly functorial.\index[terms]{derived functor}

The main advantage of derived categories (compared to the 
more traditional approach)
is that the derived functors in this formalism 
can be composed, yielding what is known as {\em the Grothendieck
spectral sequence}.\index[terms]{spectral sequence!Grothendieck} Let $F:\; \cac_1 \arrow \cac_2$
and $G:\; \cac_2 \arrow \cac_3$ be two (say) right
exact functors, and $D^{+}(\cac_i)$ the corresponding 
derived categories, bounded from the left.
Consider an object $X\in \Ob(\cac_1)$.
Replacing $X$ by its injective resolution, we
obtain a complex $I^*\in \Ob(D^{+}(\cac_1))$ quasi-isomorphic to  
the complex $0 \arrow X \arrow 0$.
Then $RF(X) = F(I^*)$. To apply $RG$ to this complex, we need
to replace the complex $F(I^*)$ by a quasi-isomorphic complex
of injective objects from $\cac_2$. This is done in traditional
fashion (\cite{_Cartan_Eilenberg:Derived_}), 
by taking an injective resolution of each $F(I^k)$ 
in such a way that the following diagram 

\begin{equation}\label{_resolution_bicomplex_Equation_}
\begin{minipage}{0.85\linewidth}	\xymatrix{
&0\ar[d]&0\ar[d]&0\ar[d]\\
\ar[r]&F(I^k)\ar[r]^d\ar[d]^{\delta}&F(I^{k+1})\ar[r]^d\ar[d]^{\delta}&
F(I^{k+2})\ar[r]^d\ar[d]^{\delta}&\\
\ar[r]&J^{0,k}\ar[r]^d\ar[d]^{\delta}&J^{0,k+1}\ar[r]^d\ar[d]^{\delta}&J^{0,k+2}\ar[r]^d\ar[d]^{\delta}&\\
\ar[r]&J^{1,k}\ar[r]^d\ar[d]^{\delta}&J^{1,k+1}\ar[r]^d\ar[d]^{\delta}&J^{1,k+2}\ar[r]^d\ar[d]^{\delta}&\\
&&&
}
\end{minipage}
\end{equation}
has
exact columns, satisfies $d^2=\delta^2=0$, and commutes. 
According to \cite{_Cartan_Eilenberg:Derived_} (see also
\cite[Ch. 20, Lemma 9.5]{_Lang:Algebra_}),
this resolution can be chosen in such a way that
the cohomology of each row gives an injective 
resolution for the cohomology of the top row.

A commutative diagram with rows 
$...\stackrel d \arrow J^{i,j} \stackrel d \arrow J^{i+1,j}\stackrel d\arrow ...$
and the columns 
$...\stackrel \delta \arrow J^{i,j} \stackrel \delta \arrow J^{i,j+1}\stackrel \delta\arrow ...$
with $d^2=0$ and $\delta^2=0$ is called {\bf a bicomplex},
or {\bf a double complex}. The diagram
\eqref{_resolution_bicomplex_Equation_}
is an example of a bicomplex. \index[terms]{bicomplex (double complex)}

{\bf The total complex of the
  bicomplex} $(J^{*,*}, d, \delta)$
(sometimes called {\bf the convolution of the bicomplex})
is the complex \index[terms]{bicomplex (double complex)!total complex (convolution) of}
\[ ... \stackrel {d-
  \delta}\arrow \bigoplus_{i+j=p} J^{i, j}\stackrel {d+
  \delta}\arrow \bigoplus_{i+j=p+1} J^{i, j}\stackrel {d-
  \delta}\arrow \bigoplus_{i+j=p+2} J^{i, j}\stackrel {d+
  \delta}\arrow...
\]
its differential squares to zero because $d$ and $\delta$
commute. It is not hard to see that the total complex
of the bicomplex

\begin{equation}\label{_resolution_bicomplex_without_first_row_Equation_}
\begin{minipage}{0.85\linewidth}	\xymatrix{
		&0\ar[d]&0\ar[d]&0\ar[d]\\
		\ar[r]&J^{0,k}\ar[r]^d\ar[d]^{\delta}&J^{0,k+1}\ar[r]^d\ar[d]^{\delta}&J^{0,k+2}\ar[r]^d\ar[d]^{\delta}&\\
		\ar[r]&J^{1,k}\ar[r]^d\ar[d]^{\delta}&J^{1,k+1}\ar[r]^d\ar[d]^{\delta}&J^{1,k+2}\ar[r]^d\ar[d]^{\delta}&\\
		&&&
	}
\end{minipage}
\end{equation}
obtained from \eqref{_resolution_bicomplex_Equation_}
is quasi-isomorphic to the complex 
$$... \arrow F(I^k)  \arrow F(I^{k+1}) \arrow F(I^{k+2})
\arrow ...$$ 
We consider this total complex as an
object in $D^+(\cac_2)$; by definition, it is equal to
$RF(X)$. Then $RG(RF(X))= R(GF)(X)$ is the object of $D^+(\cac_3)$
obtained by applying $G$ to the total complex of
\eqref{_resolution_bicomplex_without_first_row_Equation_}.
In other words, $R(GF)(X)$ is the
total complex of the double complex
\begin{equation}\label{_G_of_bicomplex_Equation_}
\begin{minipage}{0.85\linewidth}	\xymatrix{
		&0\ar[d]&0\ar[d]&0\ar[d]\\
		\ar[r]&G(J^{0,k})\ar[r]^d\ar[d]^{\delta}&G(J^{0,k+1})\ar[r]^d\ar[d]^{\delta}&G(J^{0,k+2})\ar[r]^d\ar[d]^{\delta}&\\
		\ar[r]&G(J^{1,k})\ar[r]^d\ar[d]^{\delta}&G(J^{1,k+1})\ar[r]^d\ar[d]^{\delta}&G(J^{1,k+2})\ar[r]^d\ar[d]^{\delta}&\\
		&&&
	}
\end{minipage}
\end{equation}
Given a bicomplex $(K^{p, q}, d, \delta)$, one can
associate the standard spectral sequence to this bicomplex\index[terms]{spectral sequence!of a bicomplex}
(see e.g.   \cite{_Vakil:Friend_or_Foe_}). Its $E_2^{p,q}$-term
is $H^{p}(H^q(K^{*,*}, d), \delta)$, and it converges
to the cohomology of the total complex of $(K^{*,*}, d,\delta)$.

{\bf The \index[persons]{Grothendieck, A.} Grothendieck spectral sequence} is the
spectral sequence associated with the double complex
\eqref{_G_of_bicomplex_Equation_}. The \index[terms]{spectral sequence!Grothendieck}
columns of \eqref{_G_of_bicomplex_Equation_}
are injective resolutions of  $H^p(F(I^*))= R^pF(X)$,
hence the $E_2$-term of the spectral sequence associated
with the bicomplex \eqref{_G_of_bicomplex_Equation_}
is equal to $R^pG(H^p(F(I^*))= R^pG(R^pF(X))$.\index[terms]{spectral sequence!Grothendieck}

The \index[persons]{Grothendieck, A.} Grothendieck spectral sequence has a multitude of
applications; in this book, we are going to use 
the following two.

The Grothendieck spectral sequence  allows one to
compute the cohomology of a sheaf using derived direct
images. Consider a continuous map $X \stackrel f \arrow Y$. The
direct image functor $f_*:\; \Sh(X) \arrow \Sh(Y)$ 
takes a sheaf ${\cal F}$ to the sheaf $f_*{\cal F}$
that has the section space $f_*{\cal F}(W):={\cal  F}(\f^{-1}(W))$
 for any open subset $W\subset Y$. The functor $f_*$
is clearly left exact, but it is not, generally speaking,
right exact. When $Y$ is a point,
$f_*{\cal F}= H^0(X, {\cal F})$, and this
functor is not right exact (indeed, 
the derived functor of $H^0$ is the sheaf 
cohomology).

The derived direct image $R^kf_*({\cal F})$ can be
expressed as the sheafification of the presheaf
sending $W\subset Y$ to $H^k(f^{-1}(W), {\cal F})$
(\cite[Theorem 3, Section III.8]{_Gelfand_Manin_}).
This result is not hard to prove if one uses the
\index[persons]{Godement, R.} Godement resolution, that is  known to be injective\index[terms]{resolution!Godement}
(\cite{_Grothendieck:Tohoku_}). 

Let $X \stackrel f \arrow Y \stackrel g \arrow Z$
be a composition of continuous maps. 
The $E_2^{p,q}$-term of the Grothendieck's spectral sequence 
is $R^pg_*(R^qf_*{\cal F}))$, and it converges to
$R^*(gf)_*({\cal F})$. This is especially significant
when $Z$ is a point, because in this case,
$R^k(gf)_*({\cal F})= H^k(X, {\cal F})$, and 
the \index[persons]{Grothendieck, A.} Grothendieck spectral sequence converges to the
cohomology. We use this result to compute
the cohomology of vector bundles on Hopf manifolds.

For another application,
we could compute the cohomology
of $G$-equivariant sheaves. \index[terms]{sheaf!$G$-equivariant}
Consider a sheaf ${\cal F}$ over
a space $M$, and let $\tilde M\stackrel \pi \arrow M$ be a
Galois cover of $M$, such that\index[terms]{cover!Galois}
$M = \tilde M/G$, with $G$ acting
on $\tilde M$ freely.
It is not hard to see that the category
of sheaves on $M$ is equivalent to the
category of $G$-equivariant sheaves
on $\tilde M$ (\ref{_covering_equivariant_Theorem_}).
The equivalence of categories $\Psi:\; \Sh(M)\arrow \Sh^G(\tilde M)$
is given by the pullback functor $\pi^*$
with the additional equivariant structure.

The group $G$ naturally
acts on $H^0(\tilde M, {\cal F})$, for any $G$-equivariant sheaf ${\cal F}$.
From the proof of \ref{_covering_equivariant_Theorem_} 
it is apparent that $H^0(\tilde M,\Psi({\cal F}))=H^0(\tilde M, {\cal F})^G$,
where $(\cdot)^G$ denotes the space of $G$-invariants.
This means that the cohomology functor $H^0:\; \Sh(M)\arrow \Vec$
is a composition of two left exact functors, the zero cohomology functor
$H^0(\tilde M, \Psi(\dash)):\; \Sh(M)\arrow \Rep(G)$ to the category of $G$-modules,
and the functor of $G$-invariants $\Inv:\; \Rep(G)\arrow \Vec$.
The corresponding Grothendieck spectral sequence \index[terms]{spectral sequence!Grothendieck}
\begin{equation}\label{_Grothendieck_G_inv_Equation_}
E_2^{p,q}= H^p(G, H^q(\tilde M, \pi^*{\cal F}))\Rightarrow H^{p+q}(M, {\cal F})
\end{equation}
relates the cohomology of ${\cal F}$ with the
higher derived functor $R^i\Inv$ of $\Inv:\; \Rep(G)\arrow \Vec$,
identified with the functor $H^i(G, \dash)$ of group cohomology\index[terms]{cohomology!group}
(\ref{_group_coho_Definition_}), and\index[terms]{derived functor!higher}
the functor $H^q(\tilde M, \dash)= R^qH^0(\pi^*(\dash))$,
interpreted as a higher derived functor of $H^0(\tilde M, \Psi(\dash))$.

This spectral sequence
is used to compute the cohomology of any sheaf ${\cal F}$ on
$M$ in terms of the cohomology of its pullback $\pi^*{\cal F}$,
considered to be  a $G$-equivariant sheaf on $\tilde M$.

\section[Equivariant sheaves, local systems and cohomology]{Equivariant sheaves, local systems and\\ cohomology}

In this section, we introduce 
several important concepts in sheaf theory, and
state some preliminary results
(mostly well-known) that are  used further on
in the discussion of cohomology of vector bundles
on Hopf manifolds.
 
\subsection{Equivariant sheaves and equivariant objects}

To define the equivariant sheaves,
we start with an abstract notion of equivariant object in a category.
The model example for us is a group acting on a topological 
space $M$, and the induced action
on the category of sheaves on $M$.
Eventually, we use this notion to describe a class of sheaves on a circle.
For an alternative approach to equivariant sheaves, see 
\cite{_Bernstein_Lunts_}. The following definition is
found, for instance, in \cite{_Burciu_}.

\hfill

\definition\label{_equiv_object_Definition_}
Let $\cac$ be a category, and $G$ a group acting on $\cac$ by endofunctors,
that is, by functors from $\cac$ to itself.  {\bf A $G$-equivariant structure
on an object $X$ of $\cac$} is a family of isomorphisms
$R_g:\; g(X) \arrow X$ that is  associative in the
following sense (see also \ref{_equivariant_structure_Definition_}):
\[ g_2^*(R_{g_1}) \circ R_{g_2}= R_{g_1 g_2},\] where
$R_{g_1 g_2}$ is considered to be  an isomorphism
${(g_1g_2)}^*X \arrow X$, and $g_2^*(R_{g_1}) \circ R_{g_2}$
as a composition of the isomorphism
$g_2^*(R_{g_1}):\; g_2^*g_1^* (X)\arrow g_2^*(X)$
and $R_{g_2}:\; g_2^*(X) \arrow X$.

\hfill

\remark
Note that in this definition we consider $G$ as an abstract group.
In many applications, one needs to consider $G$ as a group object
in $\cac$, and then a $G$-equivariant structure is a
morphism $G\times X\arrow X$, where $\times$ denotes the
fibred product in $\cac$. This is the definition used, for example,
in \cite{_Vistoli_}.

\hfill

We define the equivariant sheaves as equivariant objects in the category of sheaves.
For vector bundles (considered to be  sheaves of modules over the structure sheaves),\index[terms]{sheaf!$G$-equivariant}
this gives precisely \ref{_equivariant_structure_Definition_} given above.

\hfill

\definition\label{_G_equiv_sheaf_Definition_}
Let $G$ be a group acting on a topological space $M$
by homeomorphisms. Then $G$ acts on the category of sheaves
on $M$ by endofunctors. {\bf A $G$-equivariant sheaf}
on $M$ is a $G$-equivariant object in the category of sheaves.

\hfill

\theorem\label{_covering_equivariant_Theorem_}
Let $M$ be a locally connected, locally simply connected
topological space, and $\pi:\; M_1 \arrow M$ a Galois cover,
that is, a covering equipped with a free action of the
group $G$ such that $M = M_1/G$ (\ref{_Galois_cover_definition_}).
Then the category of sheaves on $M$ is equivalent
to the category of $G$-equivariant sheaves\index[terms]{cover!Galois}
on $M_1$.

\hfill

\proof
Let ${\cal F}$ be a sheaf on $M$. Then the sheaf
$\pi^*({\cal F})$ is $G$-equivariant; indeed, 
$g\circ\pi= \pi$ for any $g\in G$, and hence 
$\pi^*({\cal F})= g^*(\pi^*({\cal F}))$,
that defines a $G$-equivariant structure on $\pi^*({\cal F})$.

Conversely, let ${\cal F}_1$ be a $G$-equivariant
sheaf on $M_1$,  $U\subset M$, and $U_1:=\pi^{-1}(U)$.
The equivariant structure on ${\cal F}_1$ defines
a $G$-action on the space of sections
${\cal F}_1(U_1)$. Indeed,
$g^*({\cal F}_1)(U_1)$ is by definition isomorphic to ${\cal F}_1(U_1)$,
and the isomorphisms $R_g:\; g^*({\cal F}_1)(U_1)\arrow {\cal F}_1(U_1)$
can be interpreted as self-maps on this space indexed by $g\in G$.
The equivariance relation $R_{g_1g_2}= g_2^*(R_{g_1}) R_{g_2}$
implies that the composition of these self-maps 
is compatible with the multiplication in $G$.

Consider the sheaf ${\cal F}_1^G$ on $M$ with the
space of sections ${\cal F}_1^G(U)$ equal to the space
of $G$-invariant sections of ${\cal F}_1(\pi^{-1}(U))$.
We claim that this functor from the category
of $G$-equivariant sheaves on $M_1$ to the
category of sheaves on $M$ is inverse to the
pullback functor defined above.

Clearly, any section of ${\cal F}$ on $U$ can be interpreted as
a $G$-invariant section of $\pi^*({\cal F})$ on $\pi^{-1}(U)$,
and any $G$-invariant section of $\pi^*({\cal F})$ comes from ${\cal F}(U)$.
Therefore, the composition of $\pi^*$ and $(\cdot )^G$
is equivalent to identity. To prove that
$\pi^*$ and $(\cdot )^G$ are inverse,
it remains to show that the functor
${\cal F}_1 \arrow \pi^*({\cal F}_1^G)$
is naturally isomorphic to identity 
(\cite[Section IV.4]{_Mac_Lane:Categories_}).

Since a $G$-invariant section of
${\cal F}_1$ on $U_1=\pi^{-1}(U)$ is a section of ${\cal F}_1(U_1)$,
one has a sheaf morphism $\Psi:\; \pi^*({\cal F}_1^G)\arrow {\cal F}_1$,
defining a natural transform of functors $(\pi^*(\cdot))^G \arrow \Id$.
It remains to show that $\Psi$ is always an isomorphism.

Let $U\subset M$ an open set
such that $\pi^{-1}(M)$ is a disjoint union
of spaces $U_g$, indexed by $g\in G$, and homeomorphic to $U$.
To prove that $\pi^*(({\cal F}_1)^G)= {\cal F}_1$
it would suffice to check this on $\pi^{-1}(U)$.

Choosing the connected component $U_e$ arbitrarily, the rest of the
components are determined by $U_g= g(U_e)$.
Then ${\cal F}_1(\pi^{-1}(U))= \prod_{g\in G}({\cal F}_1(U_g))$.
Since ${\cal F}_1$ is $G$-equivariant, the action of $G$
on the set $\{U_g\}$ defines isomorphisms 
\[ R_g:{\cal F}_1(U_e)= g^*({\cal F}_1)(U_g)\arrow {\cal F}_1(U_g).
\]
These isomorphisms define an action of the group $G$ on $\prod_{g\in G}({\cal F}_1(U_g))$.
Since $\prod_{g\in G}({\cal F}_1(U_g))$
is the product of $|G|$ copies of the same set ${\cal F}_1(U_e)$, 
one has 
\[ \left(\prod_{g\in G}({\cal F}_1(U_g))\right)^{\!\!\!G}= {\cal F}_1(U_e).\]
It follows that the natural $G$-invariant 
morphism $\pi^*({\cal F}_1^G)(U_e) \arrow {\cal F}_1(U_e)$
is the identity, and hence  $\Psi:\; \pi^*({\cal F}_1^G)\arrow {\cal F}_1$
is the identity transformation, and the functors $\pi^*$ and
${\cal F}_1 \mapsto ({\cal F}_1)^G$ are inverse.
\endproof

\subsection{Group cohomology, local systems and $\Ext$ groups}

\definition\label{_group_coho_Definition_}
Let $V$ be a representation of a discrete group $G$,
and $K(G, 1)$ the corresponding \index[persons]{Eilenberg, S.} Eilenberg-MacLane space,
that is, a space with $\pi_1(K(G, 1))=G$ and\index[terms]{space!Eilenberg-MacLane}
contractible universal covering. {\bf The group
cohomology} $H^i(G,V)$ is defined as
$H^i(G,V):=H^i(K(G, 1), {\Bbb  V})$, where
${\Bbb V}$ is the local system associated with 
the representation $V$ through the Riemann--Hilbert correspondence
(Section \ref{_Riemann_Hilbert_Section_}).\index[terms]{correspondence!Riemann--Hilbert}

\hfill

There is another definition of the group cohomology.\index[terms]{cohomology!group}
Let $\R[G]$ be the group algebra of a group $G$.\index[terms]{group algebra}
It is a vector space with the basis set-theoretically 
identified with $G$, and the product which
is defined on the basis vectors by $g_i \cdot g_j = g_i g_j$.
Representations of $G$ can be interpreted as $\R[G]$-modules.
Indeed, the category of $G$-representations \index[terms]{category!of $G$-representations}
coincides with the category of $\R[G]$-modules.

\hfill

\proposition\label{_group_cohomology_local_systems_Proposition_}
Let $G$ be a group, and $\R$ the trivial $G$-representation, 
considered to be  an $\R[G]$-module.
Then $H^i(G, V)=\Ext^i_{\R[G]}(\R, V)$
for any $G$-repre\-sen\-ta\-tion $V$, where $\Ext$
is taken in the category of $\R[G]$-modules.

\hfill

\pstep
The space $\Ext^i_{\R[G]}(V_1, V_2)$
is defined using a projective (e.g.   free) resolution
\footnote{By definition, a module over a ring is called
{\bf projective} if it is a direct summand of a free module.}
$... \arrow P_2 \arrow P_1 \arrow V_1\arrow 0$ 
as the cohomology of the complex\index[terms]{resolution!free}
\begin{equation}\label{_projective_res_Hom_Equation_}
... \longleftarrow \Hom_{\R[G]}(P_3, V_2) \longleftarrow \Hom_{\R[G]}(P_2, V_2)\longleftarrow 
\Hom_{\R[G]}(P_1, V_2) 
\end{equation}
Using the Riemann--Hilbert correspondence,
we interpret each $\R[G]$-module $V$ as a local system ${\Bbb V}$
on $K(G, 1)$. Then the sequence \eqref{_projective_res_Hom_Equation_} 
gives an exact sequence of local systems
\[
... \longleftarrow \intHom({\Bbb P}_2, {\Bbb V}_2)\longleftarrow 
\intHom({\Bbb P}_1, {\Bbb V}_2)\longleftarrow 
\intHom({\Bbb V}_1, {\Bbb V}_2)  
\longleftarrow  0.
\]
where $\intHom({\Bbb A}, {\Bbb B})$ denotes the interior $\intHom$ local system,
that is, the local system with the fiber\index[terms]{local system!interior Hom}
$\Hom\left({\Bbb A}\restrict x , {\Bbb B}\restrict x \right)$ at each $x\in K(G, 1)$.
The corresponding complex of global sections
\[
... \longleftarrow H^0(\intHom({\Bbb P}_3, {\Bbb V}_2))
\longleftarrow H^0(\intHom({\Bbb P}_2, {\Bbb V}_2))\longleftarrow 
H^0(\intHom({\Bbb P}_1,{\Bbb V}_2)) 
\]
is identified with \eqref{_projective_res_Hom_Equation_},
because $\Hom_{\R[G]}(A, B)=H^0(\intHom({\Bbb A}, {\Bbb B}))$.
This allows one to compute the $\Ext$-groups in terms
of the cohomology of local systems.

\hfill 

{\bf Step 2:}
Let $... \arrow P_2 \arrow P_1 \arrow \R\arrow 0$ 
be a free resolution of the trivial representation.\index[terms]{resolution!free}
Then $\Ext^i_{\R[G]}(\R, V)$ is the $i$-th cohomology of the 
complex
\begin{equation}\label{_H^0_exact_free_Equation_}
... \longleftarrow H^0(\intHom({\Bbb P}_2, {\Bbb V}))\longleftarrow 
H^0(\intHom({\Bbb P}_1,{\Bbb V})) \longleftarrow 0
\end{equation}
To prove \ref{_group_cohomology_local_systems_Proposition_}
it remains to show that the groups $H^i(G, V)$ 
are isomorphic to the cohomology of \eqref{_H^0_exact_free_Equation_}.

Consider the exact sequence of local systems
\[    
... \longleftarrow \intHom({\Bbb P}_2, {\Bbb V})\longleftarrow 
\intHom({\Bbb P}_1,{\Bbb V}) \longleftarrow {\Bbb V} \longleftarrow 0
\]
If we prove that the local systems $\intHom({\Bbb P}_i,{\Bbb V})$
are acyclic, that is, the cohomology $H^k(\intHom({\Bbb P}_i,{\Bbb V}))$
vanishes for all $k>0$, the cohomology of the complex
\eqref{_H^0_exact_free_Equation_} can be identified with $H^*(G,V)$
using the standard spectral sequence argument.
Indeed, for any exact sequence
$... \arrow F_2 \arrow F_1 \arrow {\Bbb V}\arrow 0$\index[terms]{local system!cohomology of}
of sheaves, there is a spectral sequence with $E_2$-term
$H^i(F_j)$ converging to $H^{i+j}({\Bbb V})$; when $F_i$ are acyclic,
only the $H^0(F_i)$-term of this spectral sequence survives,
hence it degenerates in $E_2$ and the cohomology groups of the complex
\[ ... \arrow H^0(F_i) \arrow H^0(F_{i-1}) \arrow ...
\]
are equal to $H^*(V)$ (\cite[Section 2.4.1]{_El_Zein_Tu_}).

\hfill

{\bf Step 3:}
To prove \ref{_group_cohomology_local_systems_Proposition_},
it remains to show that the local systems $\intHom({\Bbb P}_i, {\Bbb V})$
are acyclic. We are going to show that $H^{>0}(K(G, 1), W)=0$
for any local system $W= \intHom(P, U)$, where $P$
is a local system obtained from a free $\R[G]$-module, and any local system $U$. 

Let $p:\; \tilde M \arrow M$ be a covering, and $F$ a sheaf on $\tilde M$. Then 
$H^i(M, p_* F) = H^i(\tilde M, F)$ because the \v Cech cocycles of
$p_* F$ are in bijective correspondence with the \v Cech cocycles of $F$
for any covering of $\tilde M$ obtained from a sufficiently small covering 
$\{U_i\}$ of $M$ by the preimages of $U_i$. \index[terms]{cohomology!\v Cech}

Let $\pi:\; \tilde K(G, 1) \arrow K(G, 1)$
be the universal covering map.  Since $\tilde K(G, 1)$ is contractible,
to prove that $H^{>0}(K(G, 1), W)=0$ it would suffice
to show that $W= \pi_* \tilde W$, for some sheaf $\tilde W$ on $\tilde K(G, 1)$.
Assume that $P$ is associated with a rank 1 free module.
Then the local system $W= \intHom(P, U)$ is isomorphic to
$\pi_* \pi^* U$. Indeed, for any contractible open set 
$A\subset K(G, 1)$, the space of sections $P(A)$ is identified
with $\R[G]$, and the space of sections $\intHom(P, U)(A)$ is 
a map associating to any $g\in G$ a section of $U$ over $A$.
This is the same as a section of $\pi^* U$ over $\pi^{-1}(A)$,
because $\pi^{-1}(A)$ is the disjoint union of $G$ copies of $A$,
with the equivariant $G$ action shifting these $G$ copies
by left shifts. 

We completed the proof of
\ref{_group_cohomology_local_systems_Proposition_}.
\endproof

\hfill

\example\label{_Laurent_group_algebra_Z_Example_}
Let $G=\Z$. Then the ring $\R[\Z]$ is isomorphic
to the ring $\R[t, t^{-1}]$ of Laurent polynomials,
with $n\in \Z$ mapped to $t^n$. In this case,\index[terms]{Laurent polynomial}
$\R[\Z]$-modules are the same as coherent
sheaves\footnote{By definition, coherent sheaves
on an affine variety $X$ are the same
as finitely generated $R$-modules, for $R=H^0(X, \calo_X)$.} 
on the affine variety $\R \backslash 0= \Spec(\R[t, t^{-1}])$.
The trivial representation $\R$ of $\Z$ corresponds
to the skyscraper sheaf $\frac{\R[t, t^{-1}]}{\langle t-1  \rangle}$.
\index[terms]{sheaf!skyscraper}\index[terms]{sheaf!coherent}

\hfill

The following corollary reiterates the argument in 
\ref{_group_cohomology_local_systems_Proposition_}, Step 3.

\hfill

\corollary\label{_pi_*_pi^*_represented_by_free_Z_rep_Corollary_}
Let $G$ be a discrete group and
$\pi:\; \tilde K(G, 1) \arrow K(G, 1)$ the universal covering map.
Consider the local system $P$ on $K(G, 1)$ associated with the
free representation $\R[G]$ of $G$ (we can identify $P$ with
$\pi_* \R_{K(G, 1)}$, the pushforward of the constant sheaf). 
Then for any sheaf ${\cal F}$ on $K(G, 1)$, one has
$\pi_* \pi^* {\cal F}= \intHom(P, {\cal F})$,
where $\intHom$ denotes the internal Hom functor
(\cite[Section II.1.7]{gode}).\footnote{Given two sheaves
${\cal F}$ and ${\cal G}$ of vector spaces over 
a field $k$, the internal $\intHom({\cal F},{\cal G})$
is the sheaf associated with the presheaf
$\intHom({\cal F},{\cal G})(U)=\Hom({\cal F}\restrict U,{\cal G}\restrict U)$,
where $\Hom({\cal F}\restrict U,{\cal G}\restrict U)$ is the abelian
group of morphisms in the category of sheaves.}

\hfill

\proof
For a sufficiently small open set $U\subset K(G, 1)$,
the space of sections $\pi^*({\cal F})(\pi^{-1}(U))$
is the product of $G$ copies of ${\cal F}(U)$;
this space is naturally identified with 
$\intHom(P\restrict U, {\cal F}\restrict U)$.
To see that this isomorphism is functorial,
we notice that $\pi^*$ is right adjoint to $\pi_*$ and obtain  
\[
\intHom(P, {\cal F})= \intHom(\pi_* \R_{\tilde K(G,1)}, {\cal F})=
\pi_* \intHom(\R_{\tilde K(G,1)}, \pi^*{\cal F})=
\pi_*\pi^*{\cal F}. \vspace{-1em}
\]
\endproof

\section{Group cohomology of $\Z$}\index[terms]{cohomology!group}

The following result is well-known.

\hfill

\proposition
Let $V$ be a representation of $\Z$.
Then $H^0(\Z, V)$ is the space of $\Z$-invariants,
$H^i(\Z, V)=0$ for all $i>1$, and
$H^1(\Z, V)$ is isomorphic to
the quotient ${V}/{C}$
where $C$ is generated by all vectors $g(v)-v$ where
$v$ runs through $V$ and $g$ through $\Z$.\footnote{This
quotient can be defined for any group representation;
it is called {\bf the space of coinvariants}.}\index[terms]{space!of invariants/coinvariants}

\hfill

\proof
By definition,
$H^i(\Z, V)$ are the cohomology groups of a
local system on $S^1=K(\Z, 1)$, and hence  
$H^i(\Z, V)=0$ for $i>1$, and
$H^0(\Z, V)$ is the space of $\Z$-invariants.
It remains only to compute $H^1(\Z, V)$.

Using  
\ref{_group_cohomology_local_systems_Proposition_}, we 
interpret the group  $H^1(\Z, V)$
as $\Ext^1_{\R[\Z]}(\R, V)$.
Note that $\R[\Z]= \R[t, t^{-1}]$
(\ref{_Laurent_group_algebra_Z_Example_}),
and hence the trivial representation has the following free
resolution
\begin{equation}\label{_reso_trivi_Z-rep_Equation_}
0 \arrow \R[t, t^{-1}]\xlongrightarrow{x \mapsto (t-1) x} \R[t,
  t^{-1}]\arrow  \frac{\R[t, t^{-1}]}{\langle t-1
  \rangle}\arrow 0.
\end{equation}
Let $W$ be a local system over $S^1$ associated with 
the rank 1 free $\Z$-module, with $t$ denoting the generator
of the free action of $\Z$ on $W$. Applying the Riemann--Hilbert
correspondence to the exact sequence\index[terms]{correspondence!Riemann--Hilbert}
\eqref{_reso_trivi_Z-rep_Equation_}, we obtain an exact sequence of local systems
\begin{equation}\label{_resolu_for_R_S^1_Equation_}
0 \arrow W \xlongrightarrow{x \mapsto (t-1) x}  W \arrow
\R_{S^1}\arrow 0,
\end{equation}
where $\R_{S^1}$ denotes the trivial local system.
Let ${\Bbb V}$ be the local system on $S^1$
associated with the $\Z$-representation $V$.
Applying the functor $\intHom(\cdot, {\Bbb V})$, we obtain
an exact sequence of local systems
\[
0 \longleftarrow \intHom(W,{\Bbb V}) \longleftarrow  
\intHom(W,{\Bbb V}) 
\longleftarrow {\Bbb V} \longleftarrow 0,
\]
because $\intHom(\R_{S^1}, {\Bbb V})={\Bbb V}$.
By Step 3 in the proof of \ref{_group_cohomology_local_systems_Proposition_},
 the local systems $\intHom(W,{\Bbb V})$ are acyclic
whenever $W$ is associated with a free representation.
This gives a long exact sequence of cohomology
\begin{multline}\label{_cohomo_Z_long_exact_Equation_}
0 \longleftarrow H^1({\Bbb V}) \longleftarrow  H^0(\intHom(W,{\Bbb V})) 
\\
\longleftarrow
H^0(\intHom(W,{\Bbb V})) \longleftarrow H^0({\Bbb V})\longleftarrow 0.
\end{multline}
However, 
\[ H^0(\intHom(W,{\Bbb V}))= H^0(\pi_*\pi^* {\Bbb V})
\]
(\ref{_pi_*_pi^*_represented_by_free_Z_rep_Corollary_}), 
and $H^0(\pi_*\pi^* {\Bbb V})=V$.
Then  \eqref{_cohomo_Z_long_exact_Equation_} gives
$H^1(\Z, V) = H^1({\Bbb V})= 
\coker\left[ V \xlongrightarrow{x \mapsto (t-1) x} V\right]$.
\endproof

\section{Directed sheaves and cohomology of $\C^n\backslash 0$}
\label{_cohomo_C^n_without_0_Section_}

\subsection{Directed sheaves: definition and examples}\index[terms]{sheaf!directed}
 
We start from the following general result, that is  
a version of Hartogs theorem, and follows immediately
from the \index[persons]{Grauert, H.} Grauert's solution of the Levi problem
 (\ref{_Grauert_Levi_Theorem_}), or, \index[terms]{theorem!Hartogs}
for example, from the  Rossi and  Andreotti--Siu 
theorem (\ref{asr}).\index[terms]{theorem!Rossi, Andreotti--Siu}

\hfill

\theorem\label{_Hartogs_Stein_annulus_Theorem_}
Let $\phi:\; M \arrow [a, \infty[$
be a proper strictly plurisubharmonic function, bounded from below,
and the annulus $M_0:=\phi^{-1}(]b, c[)$,
with $]b, c[\subset [a, \infty[$.
Then any holomorphic function on $M_0$ can be extended
to $\phi^{-1}(]-\infty, c[)$.
\index[terms]{function!plurisubharmonic!(strictly) plurisubharmonic}
\endproof

\hfill

\definition \index[terms]{sheaf!directed} 
Let ${\cal F}$ be a sheaf on  $\R$ or on $]-\infty, a]$. 
We say that ${\cal F}$ is   {\bf directed} 
if ${\cal F}(]c_1, b[)\arrow {\cal F}(]c, b[)$
is an isomorphism for any $c\in ]c_1, b[$, and the restriction
map ${\cal F}(]c, b[)\arrow {\cal F}(]c_1, b_1[)$
is always injective. 

\hfill

\remark
Clearly, a sheaf on $\R$ 
is locally constant if for any embedded intervals \index[terms]{sheaf!locally constant}
$]a, b[ \subset ]a_1, b_1[\subset \R$,
the restriction map ${\cal F}(]a_1, b_1[)\arrow {\cal F}(]a, b[)$
is an isomorphism. All locally constant sheaves 
are directed, but not all directed sheaves are
locally constant.

\hfill

\example\label{_pushforward_directed_Example_}
Let $M$ be an LCK manifold with a proper potential \index[terms]{manifold!LCK!with potential}
(say, a Hopf manifold, or, more generally, a Vaisman manifold), $\tilde M$\index[terms]{manifold!Vaisman}
its K\"ahler $\Z$-cover, and 
$\phi:\; \tilde M \arrow \R^{>0}$ an automorphic
K\"ahler potential.
Consider the sheaf $\phi_* \calo_M$ on $\R^{>0}$ (we identify $\R^{>0}$ with $\R$
using any oriented diffeomorphism).
This sheaf is directed because
every holomorphic function on an 
annulus $\phi^{-1}(]a, b[)$
is naturally extended to $\phi^{-1}[0, e^b[)\subset \tilde M$
by Hartogs theorem (\ref{_Hartogs_Stein_annulus_Theorem_}; the
function $\phi$ is by definition plurisubharmonic).\index[terms]{theorem!Hartogs}

\hfill

\lemma\label{_cohomo_directed_zero_Lemma_}
Let ${\cal F}$ be a directed sheaf on $]-\infty, 0]$.
Then  $H^1(]-\infty, 0], {\cal F})= 0$.

\hfill

\proof
We use the \v Cech cohomology. Let $\{U_i\} $\index[terms]{cohomology!\v Cech}
be a cover of $\R$ by open sets. We can always drop
a set $U_i$ that belongs to some other $U_j$, and
split $U_i$ to its connected components. Therefore
we may assume that all $U_i$ are connected, and the
only non-trivial intersections are $U_i \cap U_{i+1}= W_i$. 
A 1-cocycle of ${\cal F}$ is represented by
a collection of sections  $\{\phi_i\in {\cal F}(W_i)\}$.
We assume that the segments $U_i$ are ordered in such a way that
$U_i$ and $W_i$ have common right end, and
the restriction ${\cal F}(U_i) \arrow {\cal F}(W_i)$
is an isomorphism because ${\cal F}$ is directed.

\centerline{\includegraphics[width=1\linewidth]{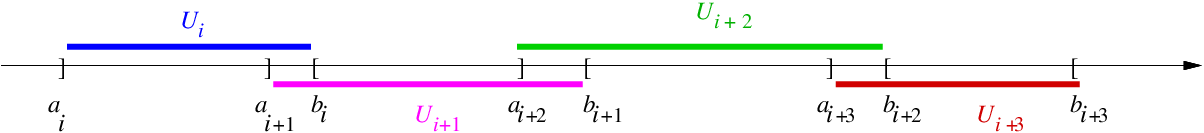}}

We number $\{U_i\}$ in such a way that $i\leq 0$,
and $U_0$ has right end in 0.

Consider a \v Cech 0-cocycle represented
by $\{f_i \in {\cal F}(U_i)\}$,
where all $f_i=0$ except $f_l= F$.
The coboundary $\6$ of this cocycle 
takes $W_i$ to $F$ if $i=l+1$, to 
$-F$ if $i=l$, and to 0 otherwise.

Now let  $\Phi$ be a \v Cech 1-cocycle,
associating $\phi_i\in {\cal F}(W_i)$ to $W_i$. The cocycle
condition is trivial, because all triple intersections
are empty. To prove the vanishing
of cohomology, it would suffice to find 
a 0-cocycle $\Psi:=\{f_i \in {\cal F}(U_i)\}$
such that $\6(\Psi)=\Phi$. We construct $\Psi$
by associating with each $U_i$ a section 
$\Psi_i\in {\cal F}(U_i)$.

Let $w_i\in \R$ be the right end of 
the interval $W_i$. 
Denote by $\tilde \phi_i$ the section of ${\cal F}(]-\infty, w_i[)$
which gives $\phi_i$ after restriction to $W_i$. This section
exists, because ${\cal F}$ is directed. 
Let $\Psi_n\in {\cal F}(U_n)$, $n\leq 0$ be the restriction
of $\sum_{i=n+1}^0 \tilde \phi_i$ to $U_n$. The difference between
$\sum_{i=n+1}^0 \tilde \phi_i$ and $\sum_{i={n}}^0 \tilde \phi_i$
is $\tilde \phi_{n}$, and hence  the restriction of 
$\Psi_n-\Psi_{n-1}$ to $W_n$ is equal to $\phi_{n}$.
Therefore, $\6(\Psi)=\Phi$.
\endproof

\subsection{Serre duality with compact supports}
\label{_Serre_s_duality_Section_}\index[terms]{duality!Serre!compact supports}

To go on, we need to recall the
Serre duality theorem on non-compact complex analytic manifolds.
Like the Poincar\'e duality theorem, it is most easily stated
and proven as a statement of duality between the usual
cohomology and the cohomology with compact support.
Indeed, it can be proven by the same argument as\index[terms]{duality!Poincar\'e}
the Poincar\'e duality theorem, using the Mayer--Vietoris
exact sequence and the \index[persons]{Serre, J.-P.} Serre duality for open balls (or polydisks),
that is  proven directly by a computation.\index[terms]{exact sequence!Mayer--Vietoris}

\hfill

We formulate the Serre  duality for vector bundles on
smooth manifolds; it is possible to state and prove it
for coherent sheaves on complex analytic spaces 
(\cite{_Ruppenthal_Samuelsson_Wulcan_,_Ramis_Ruget_}),
but the statement is more complicated and involves the
derived functors and derived categories. 

\hfill

Usually, one states the Serre  duality as a perfect pairing between
$H^i_c(M, B)$ and $H^{n-i}(B^*\otimes  K_M))$, where\index[terms]{pairing!perfect}
$H^i_c(M, B)$ denotes the cohomology with compact support.
This works very well when the groups $H^{n-i}(B^*\otimes  K_M))$
are finite-dimensional, e.g. when $M$ is compact.
However, to prove \index[terms]{duality!Serre} Serre duality (and for many other applications),
one wants to work with noncompact manifolds and infinite-dimensional
cohomology spaces. In this case,
one would prefer to put a topology on $H^{n-i}(B^*\otimes  K_M))$
such that $H^i_c(M, B)$ is its topological dual.
This is precisely how the Serre duality was stated
(and proven) in \cite{_Golovin:Duality_}.
We equip the \v Cech cohomology groups with
the topology that is  induced from the
\v Cech cochains with coefficients in $B$.

Let ${\cal F}$ be the sheaf of sections of a holomorphic 
vector bundle. To define the topology on \v Cech cochains
with coefficients in a sheaf ${\cal F}$,\index[terms]{topology!on \v Cech cochains}
we fix a sufficiently small Stein covering\index[terms]{cover!Stein}
$\{U_i\}$ on $M$. To avoid taking the limits,
we take $U_i$ such that all non-empty intersections
of $U_i$ are Stein; every covering admits such
a refinement, that is  not very hard to see,
if we choose an appropriate metric and make
sure that all $U_i$ are convex.
Then the cohomology of the sheaf ${\cal F}$
can be identified with the \v Cech
cohomology of $({\cal F}, \{U_i\})$.

We consider the space ${\cal F}(U)$
as a topological vector space, with open-compact (that
is, $C^0$-) topology\index[terms]{topology!$C^0$}. By Montel theorem, \index[terms]{theorem!Montel}
${\cal F}(U)$ is a Fr\'echet space (\ref{_Banach_bounded_holo_Theorem_}).
We consider ${\cal F}$ as a sheaf of topological vector spaces.\index[terms]{topological vector space!Fr\'echet}
 The topology on the space of ${\cal F}$-valued
 \v Cech cochains can be defined
for any sheaf ${\cal F}$ of topological vector spaces, as follows.

The space $\prod_{i_1, i_2, ..., i_k}{\cal F}(U_{i_1} \cap U_{i_2}\cap ...\cap U_{i_k})$
of $k$-cochains is equipped with the strongest topology such that 
all maps to ${\cal F}(U_{i_1} \cap U_{i_2}\cap ...\cap U_{i_k})$ are continuous.
Since the coboundary maps are continuous, 
the topology on the chain space induces the
topology on the relevant cohomology spaces.

Recall that a topological vector space $V$ is called 
{\bf Fr\'echet--Schwartz} if it is Fr\'echet and \index[terms]{topological vector space!Fr\'echet--Schwartz}
for every continuous map $A\in \Hom(V, W)$ to a normed 
topological vector space $W$ there exists 
an open subset $U\subset V$ such that $A(U)$ is precompact.
The notion of a Schwartz space \index[terms]{topological vector space!Schwartz}
is due to \index[persons]{Grothendieck, A.} Grothendieck (\cite{_Grothendieck:F_}), and this
characterization is \cite[Lemma 24.17]{_Meise_Vogt_}. 

Given a topological vector space $V$, 
let $V^*$ denote the dual vector space
(that is, the space of continuous linear functionals on $V$)
equipped with the topology of uniform convergence on
bounded subsets\footnote{For the definition of 
bounded subsets, see \ref{_bounded_set_Definition_}.} 
(this topology is called {\bf the strong topology}).\index[terms]{topology!strong}
A locally convex topological vector space\index[terms]{topological vector space!locally convex}\index[terms]{topological vector space!locally convex!reflexive}
is called {\bf reflexive} if the double dual
map $V \arrow (V^*)^*$ is an isomorphism.

A Fr\'echet-Schwartz space is reflexive
(\cite[Remark 24.24]{_Meise_Vogt_}).

V. D. \index[persons]{Golovin, V. D.} Golovin proved that 
the topology on the cohomology spaces
is Fr\'echet-Schwartz, and the  Serre duality is
the (strong) duality of the corresponding 
topological spaces.

\hfill

\theorem\label{_Serre_dualizing_sheaf_Theorem_}
(\cite{_Golovin:Duality_})
Let $M$ be an $n$-dimensional complex manifold.
Then the multiplication in
cohomology defines a perfect pairing\index[terms]{pairing!perfect}
\[
H^i_c(M, B) \times H^{n-i}(B^*\otimes  K_M))\arrow
H^n(K_M)=\C
\]
for any vector bundle $B$ on $M$.
Moreover, this duality identifies
the topological vector space $H^i(M, B)$ and its strong
dual $H^{n-i}(B^*\otimes K_M))^*$.\footnote{This space
is Fr\'echet-Schwartz, and hence  reflexive, and the isomorphism
$H^i(M, B)\cong H^{n-i}(B^*\otimes K_M))^*$ implies
$H^i(M, B)^*\cong H^{n-i}(B^*\otimes K_M))$.}
\endproof

%
%
%
%
%

\subsection{Cohomology of $\C^n\backslash 0$}

\proposition\label{_cohomo_of_C^n_withput_0_Proposition_}
The cohomology $H^i(\C^n\backslash 0, \calo_{\C^n})$
vanishes for $i \neq 0, n-1$.

\hfill

\pstep
Let $]u, v[\subset \R$ be an interval, and $\sigma:\; \C^n \arrow \R$
denote the function $z\arrow \log|z|$.
Then $R^i\sigma_* \calo_{\C^n}\restrict{]u, v[}$
is the sheafification of the presheaf given by $H^i(\sigma^{-1}(]u,v[), \calo_{\C^n})$
(\cite[Theorem 3, Section III.8]{_Gelfand_Manin_}).
Let $A_{u,v}$ denote the annulus $\sigma^{-1}(]u,v[)\subset \C^n$.
Denote by $D_t$ the set 
$D_t :=\{z\in \C^n\ \ |\ \ \sigma(z) < t\}$;
then $A_{u,v}$ is the interior of $D_{v}\backslash D_{u}$.

Consider the exact sequence
\begin{multline}\label{_long_exa_annulus_Equation_}
...\arrow H^i(A_{u,v}, \calo_{\C^n}) \arrow H^{i+1}_c(D_{u}, \calo_{\C^n})
\arrow H^{i+1}(D_{v}, \calo_{\C^n})\arrow\\ \arrow H^{i+1}(A_{u,v}, \calo_{\C^n}) \arrow ...
\end{multline}
where $H^*_c(Z, \dash)$ denotes the cohomology with compact support in
$Z$. By  Serre  duality\index[terms]{duality!Serre} with compact support
(\cite{_Golovin:Duality_}; see also Section
\ref{_Serre_s_duality_Section_}), there exists a duality of topological
vector spaces 
\begin{equation}\label{_Serre_s_duality_Equation_}
H^i_c(D_{u}, \calo_{\C^n})\times H^{n-i}(D_{u}, \calo_{\C^n}\otimes K_{D_{u}})
\arrow H^n_c(K_{D_{u}})=\C,
\end{equation}
where $K_{D_{u}}=\calo_{\C^n}$ is the canonical bundle.
Since $D_{u}$ is Stein,\index[terms]{manifold!Stein} the cohomology 
$H^{n-i}(D_{u}, \calo_{\C^n})$
vanishes for all $i\neq n$, and hence 
$H^i_c(D_{u},\calo_{\C^n})=0$ for all $i\neq n$.
Similarly, $H^{i}(D_{v}, \calo_{\C^n})=0$
for all $i\neq 0$.
We obtain that $R^q\sigma_* \calo_{\C^n}$
is equal to 0 for $q\neq 0, n-1$.

This computation expresses the cohomology groups 
$H^i(A_{u,v}, \calo_{\C^n})$ in terms of $H^0(D_{v}, \calo_{\C^n})$
and $H^n_c(D_{v}, \calo_{\C^n})$. The same computation actually
works when $v = \infty$; this gives the cohomology
of the complement $\C^n\backslash \bar D_{u}$.

\hfill

{\bf Step 2:}
We obtained that almost all terms
in the exact sequence \eqref{_long_exa_annulus_Equation_}
vanish, giving only two non-zero terms,
\[ H^{0}(D_{v},  \calo_{\C^n})= H^{0}(A_{u,v},  \calo_{\C^n}),
\]
and 
\[ H^{n}_c(D_{u},  \calo_{\C^n})= H^{n-1}(A_{u,v},  \calo_{\C^n}).
\]
Recall that the Grothendieck spectral sequence converges to  \index[terms]{spectral sequence!Grothendieck}
$H^{p+q}(D_{u}\backslash 0,  \calo_{\C^n})$
with $E_2$-term $H^p(]-\infty, 0], R^q\sigma_* \calo_{\C^n})$.
By Step 1, $R^q\sigma_* \calo_{\C^n}=0$ 
for $q\neq 0, n-1$.

The sheaf $\sigma_* \calo_{\C^n}$ is directed
(\ref{_pushforward_directed_Example_}). 
Since the cohomology $H^1(]-\infty, 0], {\cal F})$ 
of any directed sheaf ${\cal F}$ vanishes (\ref{_cohomo_directed_zero_Lemma_}),
we obtain that the term $H^1(]-\infty, 0], \sigma_*    \calo_{\C^n})=0$
does not contribute to the \index[persons]{Grothendieck, A.} Grothendieck spectral sequence.
Instead, the following three terms could in theory contribute:
$H^0(]-\infty, 0],\sigma_* \calo_{\C^n})$, 
$H^0(]-\infty, 0],R^{n-1}\sigma_* \calo_{\C^n})$
and $H^1(]-\infty, 0],R^{n-1}\sigma_* \calo_{\C^n})$.
However, the last term would add up to $n$-th
cohomology, and on a connected non-compact complex $n$-manifold,
the $n$-th cohomology of a vector bundle vanishes.
Indeed, by Serre duality, the $n$-th cohomology is dual
to 0-th cohomology with compact support, but a 
holomorphic bundle cannot have a non-zero holomorphic section
with compact support by the analytic continuation principle.\index[terms]{analytic continuation}

This implies that $H^{q}(D_{u}\backslash 0,  \calo_{\C^n})$
is equal to $H^0(]-\infty, 0],R^q\sigma_* \calo_{\C^n})$;
these groups vanish for $q\neq 0, n-1$ by Step 1.

To finish the proof, we  shall express  
$H^i(\C^n\backslash 0,  \calo_{\C^n})$ through 
 $H^{i}(D_{u}\backslash 0, \calo_{\C^n})$, that we have
already computed.

\hfill

{\bf Step 3:} We use
the Mayer--Vietoris exact sequence by decomposing \index[terms]{exact sequence!Mayer--Vietoris}
$\C^n\backslash 0$ onto a union of $D_{u}\backslash 0$
and an annulus $A_{x, \infty}$, with $x< u$.
The Mayer--Vietoris gives
\begin{multline}
...
\arrow H^i(\C^n\backslash 0, \calo_{\C^n})\arrow
H^i(A_{x, \infty}, \calo_{\C^n})\oplus H^i(D_{u}\backslash 0, \calo_{\C^n})
\arrow\\
\arrow H^i(A_{x, u}, \calo_{\C^n})\arrow  H^{i+1}(\C^n\backslash 0, \calo_{\C^n})\arrow ...
\end{multline}
By Step 2, all the cohomology terms in this sequence  vanish, except
for $i=0$ and $i=n-1$. For $i=0$, the sequence
\begin{multline}
0\arrow \H^0(\C^n\backslash 0, \calo_{\C^n})\arrow
H^0(A_{x, \infty}, \calo_{\C^n})\oplus H^0(D_{u}\backslash 0, \calo_{\C^n})
\arrow \\
\arrow H^0(A_{x, u}, \calo_{\C^n})\arrow 0
\end{multline}
is exact, because two functions which agree on the
intersection can be glued.
Therefore, $H^1(\C^n\backslash 0, \calo_{\C^n})=0$.
The same exact sequence together with the
computation in Step 2 implies that $H^i(\C^n\backslash 0, \calo_{\C^n})=0$
for all $i$ except $i=0, n-1, n$. The $n$-th cohomology
vanish by Serre duality (Step 2), and hence 
$H^i(\C^n\backslash 0, \calo_{\C^n})=0$
for all $i$ except $i=0, n-1$.
\endproof

\section[Contractions define  compact operators on holomorphic functions]{Contractions define  compact operators \\ on holomorphic functions}\index[terms]{operator!compact}

Let $\gamma:\; \C^n \arrow \C^n$ be a contraction map.
We consider $H^0(\C^n, \calo_{\C^n})$
as a topological vector space, with the open-compact topology,
also known as $C^0$-topology\index[terms]{topology!$C^0$} (\ref{_compact_open_topology_Definition_}).
Let $K\subset \C^n$ be a compact, and $\|f \|_K:=\sup_{x\in K} |f(x)|$. 
Then the base of topology in $H^0(\C^n, \calo_{\C^n})$ is given by open 
balls in the norm $\|\cdot \|_K$,
\[ B_K(C):= \{f\in H^0(\C^n, \calo_{\C^n})\ \ |\ \ \|f\|_K<C\}
\] 
for each compact $K$ and each $C\in \R^{>0}$. This sequence of norms
is complete, by Montel theorem, and defines a structure of a Fr\'echet
space on $H^0(\C^n, \calo_{\C^n})$.\index[terms]{theorem!Montel}\index[terms]{topological vector space!Fr\'echet}

Let $\gamma:\; \C^n \arrow \C^n$ be a contraction map.
In \ref{_contra_compact_Theorem_}, we proved that $\gamma$
is compact with respect to the completion of\index[terms]{operator!compact}
$H^0(\C^n, \calo_{\C^n})$ in each norm $\|\cdot \|_K$.
In this subsection, we are going to show that 
$\gamma^*:\; H^0(\C^n, \calo_{\C^n}^k)\arrow H^0(\C^n, \calo_{\C^n}^k)$ 
is compact in the $C^0$-topology.\index[terms]{topology!$C^0$} So far, the compact operators\index[terms]{operator!compact}
were defined only for Banach spaces\index[terms]{space!Banach}. Now we are going to define 
compact maps for arbitrary locally convex topological vector spaces.

A {\bf bounded set} of a locally convex topological
vector space $V$ is a subset $A\subset V$ such that
for any open neighbourhood $U\ni 0$, there exists
a number $k\neq 0$ such that $kA \subset U$.
A {\bf compact morphism} of locally convex topological
vector spaces is a continuous linear map taking every 
bounded set to a precompact set.

\hfill

\proposition\label{_contraction_compact_Frechet_Proposition_}
Let $\gamma:\; \C^n \arrow \C^n$ be a contraction map.
Then the map $\gamma^*:\; H^0(\C^n, \calo_{\C^n})\arrow H^0(\C^n, \calo_{\C^n})$
is compact.

\hfill

\proof
Since $H^0(\C^n, \calo_{\C^n})$ is separable for the Fr\'echet and\index[terms]{topology!Fr\'echet}
the $\|\cdot\|_K$-to\-po\-lo\-gies, it has a countable base, and
the compactness is equivalent to sequential compactness.
Therefore, it would suffice to prove that for
any bounded sequence $\{f_i\}\subset H^0(\C^n, \calo_{\C^n})$,
the sequence $\{\gamma^* f_i\}$ has a converging subsequence.

By Montel theorem, the 
completion of $H^0(\C^n, \calo_{\C^n})$ with respect to
$\|\cdot\|_K$ is the space $H^0(K,\calo_{\C^n})$ of germs of 
holomorphic functions in $K$; when $K$ is the closure of an open set $K_0$, 
the space $H^0(K,\calo_{\C^n})$ is the space of holomorphic functions
on $K_0$ that are  continuous up to the boundary.
Since $\gamma^*$ is compact in $\|\cdot\|_K$-topology
(\ref{_contra_compact_Theorem_}), the sequence
$\{\gamma^* f_i\}$ has a subsequence which converges in
$\|\cdot\|_K$ to a function $f\in H^0(K,\calo_{\C^n})$.
Choose an exhausting sequence of compact subsets $K_i\subset \C^n$;
for example, we can take for $K_i$ the closed ball of radius
$i\in \Z^{>0}$. Using the diagonal method, we find a subsequence
of $\{f_i\}$ which converges in each of $\|\cdot\|_{K_i}$-topologies,
hence, converges on all compacts in $\C^n$. Such a sequence clearly
converges in the Fr\'echet topology\index[terms]{topology!Fr\'echet} on $H^0(\C^n,\calo_{\C^n})$.
\endproof

\hfill

In the next subsection, we will apply this
result to explore the Fredholm-like properties of the
operator $\Id+ \gamma^*$. On Banach spaces\index[terms]{space!Banach}, compact
operators admit a spectral decomposition; this is the
content of Riesz--Schauder theory, \ref{_Riesz_Schauder_main_Theorem_}.
Even if there is no spectral decomposition
for compact operators\index[terms]{operator!compact} on general locally 
convex topological vector spaces, a 
part of it can be retained.

\hfill

\theorem \label{_decompo_compact_Theorem_}
Let $\kappa:\; W \arrow W$ be a compact operator\index[terms]{operator!compact}
on a locally convex topological vector space, and
and $A_\lambda= \lambda \Id + \kappa$ for some $\lambda\in \C^*$. 
Then $W$ admits a decomposition
$W= W_1 \oplus W_2$, where $W_1, W_2$ are closed
subspaces preserved by $\kappa$, the space $W_1$ is finite-dimensional,
$\kappa\restrict{W_1}$ is nilpotent, and $A_\lambda\restrict {W_2}$ 
is invertible.

\hfill

{\bf Proof:}
This is \cite[page 148, Theorem 1]{_Robertson_},
applied to $u=\Id$, $v=A_\lambda$ and $t=-\kappa$.
\endproof

\section[Mall's theorem on cohomology of vector bundles]{Mall's theorem on cohomology of\\ vector bundles}
\label{_Hopf_cohomo_Section_}


Let $H$ be a Hopf manifold, that is, a quotient
$\frac{\C^n\backslash 0}{\Z}$, $n \geq 2$, where $\Z$ acts on 
$\C^n$ by holomorphic contractions. Throughout this book,
we mostly used {\bf the linear Hopf manifold}, that is 
the quotient of $\C^n\backslash 0$ by a linear contraction.
In this chapter, we can relax this condition by allowing
a more general action.\index[terms]{manifold!Hopf}

\hfill

\claim
Let $\pi:\;  \C^n\backslash 0\arrow H$ be the
universal cover of a Hopf manifold, and 
$B= (TH)^{\otimes l} \otimes (T^*H)^{\otimes m}$ a tensor bundle on $H$.
Denote by $j:\; \C^n\backslash 0\hookrightarrow \C^n$
the standard embedding map. Then $j_* \pi^* B$ is a trivial
bundle on $\C^n$.

\hfill

\proof
Clearly, $\pi^* B$ is a tensor bundle on $\C^n \backslash 0$,
hence it is a trivial vector bundle, 
$\pi^* B= \calo_{\C^n \backslash 0}^{n(l+m)}$. Indeed,
$\pi^*B= (TX)^{\otimes l} \otimes (T^*X)^{\otimes m}$,
where $X=\C^n \backslash 0$, and this bundle is trivial.
However, $j_*\calo_{\C^n \backslash 0}=\calo_{\C^n}$
by Hartogs theorem.
\endproof

\hfill

The main result of this chapter is a theorem of D. \index[persons]{Mall, D.} Mall
(\cite{_Mall:Contractions_}); we need it only for the tangent
bundle, but it is naturally stated and proven in a more
general situation.

\hfill

Let $\pi:\;  \C^n\backslash 0\arrow H$ be the
universal cover of a Hopf manifold.
By \ref{_covering_equivariant_Theorem_}, 
the pullback of every sheaf on $H$\index[terms]{sheaf!equivariant}
is a $\Z$-equivariant sheaf on $\C^n\backslash 0$,
and, conversely, every $\Z$-equivariant sheaf
on $\C^n\backslash 0$ descends to a sheaf
on $H$. In the sequel, we shall consider
the pullback $\pi^*(F)$ of a sheaf
on $H$ as a $\Z$-equivariant sheaf on $\C^n\backslash 0$.
Consider the standard embedding map
$j:\; \C^n\backslash 0\hookrightarrow \C^n$.
Then the sheaf $j_* \pi^*(F)$
is also $\Z$-equivariant.

\hfill

\theorem\label{_Mall_cohomology_Theorem_}
\index[terms]{theorem!Mall}
Let $\pi:\;  \C^n\backslash 0\arrow H$,  $n \geq 3$, be the
universal cover of a Hopf manifold,  
$j:\; \C^n\backslash 0\hookrightarrow \C^n$
the standard embedding map, and 
$B$ a holomorphic vector bundle over $H$
such that $j_* \pi^* B$ is a locally trivial 
coherent sheaf on $\C^n$.%
\footnote{By \index[persons]{Oka, K.} Oka-Grauert 
homotopy principle, \cite[Theorem 5.3.1]{_Forstneric:Oka_book_},
any holomorphic vector bundle on $\C^n$ is trivial;\index[terms]{theorem!Oka-Grauert homotopy principle}
thus, instead of local triviality, we could assume 
that $j_* \pi^* B$ is a trivial vector bundle.}
Then $\dim H^0(H, B)=\dim H^1(H,B)$, and this group
is equal to the space of $\Z$-invariant sections
of $j_* \pi^* B$. Moreover, $H^i(H, B) =0$ for all
$i$ such that $1 < i < n-1$.

\hfill

\proof
By \ref{_cohomo_of_C^n_withput_0_Proposition_}, 
$H^i(\C^n \backslash 0, \calo_{\C^n \backslash 0})= 0$
for all $i\neq 0, n-1$. We use the \index[persons]{Grothendieck, A.} Grothendieck's spectral sequence
\eqref{_Grothendieck_G_inv_Equation_}\index[terms]{spectral sequence!Grothendieck}
\[ 
E_2^{p,q}= H^p(G, H^q(\tilde M, \pi^*{\cal F}))\Rightarrow H^{p+q}(M, {\cal F}),
\]
where $G=\Z$, $M=H$, ${\cal F}=B$, and $\tilde M=\C^n\backslash 0$:
\begin{equation}\label{_Grothendieck_for_Z_on_C^n_o_Equation_} 
E_2^{p,q}= H^p(\Z, H^q(\C^n\backslash 0, \pi^*B))\Rightarrow H^{p+q}(H, B).
\end{equation}
Since $\pi^*B$ is trivial, we have
 $H^q(\C^n\backslash 0, \pi^* B)=0$ 
for all $1\neq 0, n-1$ 
(\ref{_cohomo_of_C^n_withput_0_Proposition_}). However, 
$H^i(G, H^q(\C^n\backslash 0, \pi^* B))$
is the space of invariants (for $i=0$)
and coinvariants (for $i=1$) of the natural\index[terms]{space!of invariants/coinvariants}
$\Z$-action in $H^q(\C^n\backslash 0, \pi^* B))$.
This implies that the \index[persons]{Grothendieck, A.} Grothendieck spectral
sequence \eqref{_Grothendieck_for_Z_on_C^n_o_Equation_} 
degenerates, and the cohomology $H^i(H, B)$ of $B$
is non-zero only for $i=0,1,n-1, n$:
\begin{equation}\label{_cohomo_coinva_Mall_Equation_}
	\begin{split}
H^0(H, B)&=H^0(\Z,H^0(\C^n\backslash 0, \pi^* B)),\\ 
H^1(H, B)&=H^1(\Z,H^0(\C^n\backslash 0, \pi^* B)),\\
H^{n-1}(H, B)&=H^0(\Z,H^{n-1}(\C^n\backslash 0, \pi^* B)),\\ 
H^n(H, B)&=H^1(\Z,H^{n-1}(\C^n\backslash 0, \pi^* B)).
	\end{split}
\end{equation}
To prove \ref{_Mall_cohomology_Theorem_},
it remains only to show that 
\begin{equation}\label{_inva_coinva_Equation_}
\dim H^0(\Z,H^0(\C^n\backslash 0, \pi^* B))=\dim H^1(\Z,H^0(\C^n\backslash 0, \pi^* B)).
\end{equation}
Notice that $H^0(\C^n\backslash 0, \pi^* B)=H^0(\C^n, j_* \pi^*B)$,
where $j_* \pi^*B$ is a trivial bundle. Then 
\eqref{_inva_coinva_Equation_} is implied by the following lemma,
which finishes the proof.

\hfill

\lemma
Let $\gamma$ be an invertible holomorphic contraction of $\C^n$,
and $B$ a vector bundle on $\C^n$ that is  $\Z$-equivariant
with respect to the action of $\Z$ generated by $\gamma$ on $\C^n$.
Let $H^0(\C^n, B)^\Z$ denote the space of $\Z$-invariants on $H^0(\C^n, B)$,
and $H^0(\C^n, B)_\Z$ the space of $\Z$-coinvariants. Then the spaces
$H^0(\C^n, B)^\Z$ and $H^0(\C^n, B)_\Z$ are finite-dimensional and
have the same rank.

\hfill

\proof
A $\Z$-equivariant structure is an isomorphism
$\gamma^* B \stackrel \phi \arrow B$. Since 
any holomorphic bundle on $\C^n$ is trivial 
(\cite[Theorem 5.3.1]{_Forstneric:Oka_book_}),
we may assume that $B = \calo_{\C^n}^k$. This bundle
has a natural $\Z$-equivariant structure, that we
denote by $\phi_0$. Then $\phi_0 = A \circ \phi $, where
$A$ is an automorphism of $B$. 
By \ref{_contraction_compact_Frechet_Proposition_},
$\phi_0$ is a compact operator\index[terms]{operator!compact} on the space $H^0(\C^n, \calo_{\C^n}^k)$
equipped with the natural Fr\'echet topology\index[terms]{topology!Fr\'echet} (the $C^0$-topology).\index[terms]{topology!$C^0$}
Since $A$ is  $C^0$-continuous, $\phi$ is also compact.
Now, \ref{_decompo_compact_Theorem_}
implies that $H^0(\C^n, \calo_{\C^n}^k)= W_1 \oplus W_2$,
with $\Id - \phi$ invertible on $W_2$, $W_1$ finite-dimensional,
and $\phi\restrict{W_2}$ nilpotent. Then 
$H^0(\C^n, B)^\Z=\ker \phi\restrict {W_2}$ and
$H^0(\C^n, B)_\Z=\coker \phi\restrict {W_2}$;
the kernel and the cokernel of a nilpotent endomorphism of a
finite-dimensional space have the same dimension.
\endproof

\hfill

\remark
We did not state \ref{_Mall_cohomology_Theorem_}
when $M$ is a Hopf surface. However, from
its proof it follows that \index[terms]{surface!Hopf}
$$\dim H^1(H, B)= \dim H^0(H,B) + \dim H^2(H, B).$$
Indeed, there are two components in the first
cohomology group $H^1(H, B)$, separately obtained from 
$\Z$-invariant sections of $H^0(\C^n, \pi^* B)$ and 
from $\Z$-invariant sections of $H^0(\C^n, \pi^* (B^* \otimes K_H))$;
the first has the same dimension as $H^0(H,B) $, and the
second the same dimension as $\dim H^2(H, B)$.

\section{Exercises}

\begin{enumerate}[label=\textbf{\thechapter.\arabic*}.,ref=\thechapter.\arabic{enumi}]

\item
Recall that {\bf the Lebesgue dimension}\index[terms]{dimension!Lebesgue}
of a topological space $M$ is the smallest number
such that any covering has a refinement $\{U_i\}$
such that each point of $M$ is contained in no more than
$n+1$ elements of $\{U_i\}$.
\begin{enumerate}
\item Let $M$ be an $n$-manifold.
Prove that its Lebesgue dimension is $n$.
\item Let $M$ be a space of Lebesgue dimension $n$, 
and ${\cal F}$ a sheaf on $M$. Prove that
$H^i(M, {\cal F})=0$ for all $i >  n$.
\end{enumerate}

\item
Let $X$ be a topological space, and 
$\psi:\; F \arrow X$ a continuous map.
For any open set $U\subset X$, a continuous
section of $\psi$ is a continuous map $U \stackrel u\arrow F$ such that
$u \circ \psi= \Id_U$.  Let $\Psi(U)$ denote the
set of all sections of $\psi$ over $U$. 
For a smaller open set $U'\subset U$,
the restriction of a section from $U$ to $U'$ defines 
the natural map $\Psi(U) \arrow \Psi(U')$.
Prove that $U \mapsto \Psi(U)$ defines a sheaf of sets over $X$.

\item\label{_etale_space_Exercise_}
Let ${\cal F}$ be a sheaf over a topological space $X$.
\begin{enumerate}
\item Prove that there exists a topological space $F$ equipped
with a map $\psi:\; F \arrow X$ with the following properties:
(a) $\psi$ is a local homeomorphism; (b) for each open set $U\subset X$,
the sheaf of continuous sections of $\psi$ over $U$ is equal to ${\cal F}$.

\item Prove that the pair $(F, \psi:\; F \arrow X)$ is defined
uniquely by the sheaf ${\cal F}$.
\end{enumerate}

\definition
The space $F$ defined in Exercise 
\ref{_etale_space_Exercise_} is called {\bf the \'etal\'e space of 
the sheaf ${\cal F}$}.\index[terms]{space!\'etal\'e}

\item Let $X$ be a topological space equipped with a continuous
action of a group $G$, and ${\cal F}$ a sheaf of vector spaces on $X$. 
Prove that a $G$-equivariant structure on the sheaf ${\cal F}$
is the same as a $G$-action on its \'etal\'e space
compatible with the projection to $X$ and with the
vector space structure on each space of sections.

\item\label{_ext_zero_Exercise_}
Let ${\cal F}$ be a sheaf over $M$.
We say that ${\cal F}$ {\bf has the analytic continuation property}\index[terms]{analytic continuation}
when the restriction map ${\cal F}(U) \arrow {\cal F}(U')$
is injective for any non-empty $U' \subset U$ and any connected $U\subset M$.
Prove that ${\cal F}$ has the analytic continuation property
if and only if its \'etal\'e space is Hausdorff.
\index[terms]{sheaf!with analytic continuation property}
\item
Let ${\cal F}$ be a sheaf over an open set $U\stackrel j \hookrightarrow M$.
Denote by $j_!{\cal F}$ the {\bf extension by zero} of ${\cal F}$,
that is, the sheaf given by $j_!{\cal F}(W) = {\cal F}(W)$ when $W\subset U$
and  $j_!{\cal F}(W) = 0$ when $W\not\subset U$.
\begin{enumerate}
\item Prove that $j_!{\cal F}$ has the 
analytic continuation property 
if ${\cal F}$ has.
\item Let $\bar U$ denote the closure of $U$.
 Prove that $j_*{\cal F}$ 
does not have the analytic continuation property
if $H^0(U, {\cal F})\neq 0$ and $\bar U\neq M$.
\item 
Find an example of a non-trivial sheaf ${\cal F}$ on
an open subset $U\subsetneq M$
such that its direct image $j_* {\cal F}$ 
has the analytic continuation property.

\item Let $Z:= M \backslash U$, and 
$Z\stackrel i \hookrightarrow M$  denote the
closed embedding. Prove that for any sheaf ${\cal F}$ on $M$ 
the following sequence of sheaves is exact:
\[
0 \arrow j_!j^*{\cal F}\arrow {\cal F} \arrow i_*i^*{\cal F} \arrow 0.
\]
\end{enumerate}

\item
Let $U=\R \backslash \{0, 1\}$. Denote by  
$\R_U$ the constant sheaf on $U$, and let $j_! \R_U$ be its
zero extension to $\R$. Consider the exact sequence\index[terms]{sheaf!zero extension}
\[ 0 \arrow j_! \R_U\arrow \R_\R \arrow K\arrow 0,\]
where $j_! \R_U\arrow \R_\R $ denotes the tautological embedding,
and $K$ its cokernel.
\begin{enumerate}
\item Prove that $j_* \R_U = \R$.
\item Prove that $K$ is a direct sum of two
skyscraper sheaves, concentrated in 0 and in 1.\index[terms]{sheaf!skyscraper}
\item Prove that 
\[ \dim H^1(\R,j_! \R_U)\geq \dim H^0(\R,K)- \dim H^0(\R,\R_\R) >0.
\]
\end{enumerate}

{\em Hint:} Use Exercise \ref{_ext_zero_Exercise_} (d).


\item
Let ${\cal F}$ be a sheaf with analytic continuation \index[terms]{analytic continuation}property over $\R$. 
Prove that its \'etal\'e space is a disjoint
union of copies of intervals of form $]a, b[$, with
$a, b \in \R \cup \{\pm \infty\}$, injectively projected to $\R$.

\item
Let ${\cal F}$ be a sheaf over $\R$. Prove that 
${\cal F}$ is directed if and only if  its \'etal\'e space is a disjoint
union of copies of intervals of form $]-\infty, a[$ and
$]-\infty, \infty[$, injectively projected to $\R$.

\item
Define {\bf an $\R$-filtration} on a vector space $W$
as a collection of subspaces $W_t\subset W$ indexed by
$t\in \R$ such that $W_t \supset W_{t'}$ whenever $t < t'$.
Prove that the category of directed sheaves on $\R$ is
equivalent to the category of $\R$-filtered spaces.

\item
Let $\phi:\; \C^* \arrow \R^{> 0}$ be the map $z \mapsto |z|^2$,
and ${\cal F}:= \phi_* \calo_{\C^*}$ be the pushforward of the sheaf of
holomorphic functions. 
\begin{enumerate}
\item Prove that ${\cal F}$ is not a directed sheaf on $\R^{> 0}$.
\item Prove that $H^1(\R^{> 0}, {\cal F})=0$.
\end{enumerate}

{\em Hint:} Use the \index[persons]{Grothendieck, A.} Grothendieck spectral sequence
with the $E^2$-term \\ $H^p(\R^{> 0}, R^q\phi_*({\cal F}))$ converging
to $H^*(\C^*, \calo_{\C^*})$.\index[terms]{spectral sequence!Grothendieck}

\item \label{_Serre-Hochschild_sequence_Exercise_}
Let $G$ be a group, $H\subset G$ 
be a normal subgroup, and $V$ a $G$-module.
\begin{enumerate}
\item
Prove that the group $G/H$ naturally 
acts on the cohomology groups $H^i(H, V)$.
\item Construct a spectral sequence whose $E_2$-term
$H^p(G/H, H^q(H, V))$ converges to $H^{p+q}(G, V)$.\footnote{This
spectral sequence is called {\bf the \index[terms]{spectral sequence!Serre--Hochschild} Serre--Hochschild spectral sequence}.}
\end{enumerate}

{\em Hint:} Represent the functor of $G$-invariants
as a composition of functors of $H$-invariants and $G/H$-invariants
and use the \index[persons]{Grothendieck, A.} Grothendieck spectral sequence.

\item Let $\cac$ be an abelian category such that 
any exact sequence in $\cac$  splits, that is, is decomposed to 
a direct sum of complexes $0\arrow Q^i  \tilde \arrow Q^{i+1} \arrow 0$.
Prove that any left exact functor and any right exact functor on $\cac$
is exact.

\item
Let $\cac$ be an abelian category such that 
any short exact sequence in $\cac$ splits. Prove that any
exact sequence in $\cac$ splits.\index[terms]{category!abelian}

\item 
Let $\cac$ be the category of representations of a finite group $G$
over a field $k$, with $\Char k$ coprime with $|G|$.
\begin{enumerate}
\item
Let $0 \arrow A \arrow B  \arrow C\arrow 0$ be a short exact sequence.
Using the averaging map $a \arrow \frac{1}{|G|} \sum_{g\in G} g(a)$,
prove that the corresponding sequence 
of $G$-invariants $0 \arrow A^G \arrow B^G \arrow C^G \arrow 0$
is exact.
\item Prove that any exact sequence in $\cac$ splits.
\item Prove that the higher cohomology of a finite group over $\Q$ vanishes.
\end{enumerate}

{\em Hint:} To produce a section of $B \stackrel u \arrow C\arrow 0$
in the category of $G$-re\-pre\-sen\-tations,
show that the sequence $\Hom(C, B)^G \arrow \Hom(C, C)^G\arrow 0$
is exact, and take $u \in \Hom(C, B)^G$ that maps to identity. 

\item
Prove that any holomorphic line bundle on $\C^2 \backslash 0$
is trivial. 

{\em Hint:}  Use the exponential exact sequence.

\item
Let $U_i\subset \C^n$ be the complement
of the coordinate plane $t_i=0$, 
$W_k:= \bigcap_{i\neq k} U_i$, and 
$d:\; \bigoplus_k H^0(W_k, \calo_{\C^n}) \arrow H^0(\bigcap U_i,\calo_{\C^n})$
be the \v Cech differential. 
\begin{enumerate}
\item 
Prove that $\coker d$ is the quotient of the space of
all holomorphic functions on $\bigcap U_i$
by the space of holomorphic functions on $\bigcap U_i$
that can be extended to at least one of $W_i$.

\item 
Prove that $\coker d= H^{n-1}(\C^n \backslash 0, \calo_{\C^n})$.

\item
Let $H^{n-1}_{alg}(\C^n \backslash 0, \calo_{\C^n})$
denote the cohomology of the sheaf of regular
algebraic functions on the space $\C^n \backslash 0$ considered
as an algebraic manifold. Prove that 
$H^{n-1}_{alg}(\C^n \backslash 0, \calo_{\C^n})$
is the space \[  t_1^{-1}t_2^{-1}... t_n^{-1}\C[t_1^{-1},
  ..., t_n^{-1}]\] of principal parts of all Laurent polynomials
on $\C^n$.\index[terms]{Laurent polynomial!principal part of}
\end{enumerate}

\end{enumerate}


\chapter{Mall bundles and flat connections on Hopf manifolds}
\label{_Mall_bundles_Chapter_}

\epigraph{\it [Such] knowledge [is] too wonderful for me; it is high, I cannot [attain] unto it.}{\sc\scriptsize Psalm 139,\ \ King James' Bible}

\section{Introduction}


\subsection{Mall bundles and coherent sheaves}

In the previous chapter, we computed the cohomology
of a holomorphic vector bundle on a Hopf manifold
such that its pullback to $\C^n \backslash 0$ can be
extended to a holomorphic vector bundle on $\C^n$.
We call such vector bundles {\bf the \index[persons]{Mall, D.} Mall bundles}
(\ref{_Mall_bundle_Definition_}).\index[terms]{bundle!vector bundle!Mall}
Note that any holomorphic vector bundle on $\C^n$
is trivial, as follows from the \index[persons]{Oka, K.} Oka-Grauert homotopy principle
(\cite[Theorem 5.3.1]{_Forstneric:Oka_book_}).\index[terms]{theorem!Oka-Grauert homotopy principle}

In this chapter we explore the geometric properties of the
Mall bundles. 

The arguments used to compute the cohomology
of \index[persons]{Mall, D.} Mall bundles in Chapter \ref{_cohomo_on_Hopf_Theorem_}
can be used on any locally conformally K\"ahler manifold
with potential\index[terms]{manifold!LCK!with potential}, and, indeed, not only on Mall bundles.
However, this would require a deeper knowledge of algebraic
geometry and theory of coherent sheaves, and we decided 
to skip this part. For the readers who want to pursue
this direction, we state the main results from the
theory of coherent sheaves and the sheaf extensions.

A coherent sheaf is a locally finitely generated, and \index[terms]{sheaf!coherent}
locally finitely presented sheaf\footnote{A sheaf is ``locally finitely presented''
if it is locally isomorphic to a quotient of a finitely generated free sheaf of $\calo_M$-modules by
another finitely generated free sheaf.} of modules over holomorphic
functions. Any holomorphic vector bundle is a coherent
sheaf, but not vice versa: there are more coherent sheaves
than holomorphic bundles. We are interested in coherent
sheaves obtained as extensions of vector bundles.

The first natural notion which occurs in this context
is the notion of {\bf a normal sheaf}. A coherent sheaf
$\caf$ on a variety $M$ is {\bf normal} if for any \index[terms]{sheaf!normal}
complex subvariety $Z\subset M$, $\codim Z > 1$,
the natural map $\caf\arrow j_* j^* \caf$ is an isomorphism,
where $j:\; M \backslash Z \arrow M$ is a natural
embedding, and $j_*$ the sheaf pushforward (direct image).
One usually defines normal sheaves on $M$ in assumption
that $M$ is normal; normality of $M$ is
equivalent to the normality of the sheaf $\calo_M$.
By \ref{_reflexife_normal_Theorem_}, 
\index[terms]{sheaf!reflexive}
normal sheaves are the same as {\bf reflexive
  sheaves}, that is, sheaves $F$ such that the
natural morphism $\caf \arrow \caf^{**}$ from $F$ to its double dual is
an isomorphism.

In algebraic category, any coherent sheaf on 
$M \backslash Z$ can be extended to a coherent
sheaf on $M$; this statement is obtained directly
from the definitions. In complex geometry, this
statement is false, indeed, there exists
a holomorphic vector bundle on $\C^2 \backslash 0$
that cannot be extended to a coherent sheaf on $\C^2$
(\cite{_Ballico:bundles_on_C^2_}).

The subject of coherent sheaf extension was very
popular in late 1960s and early 1970s, indeed,
there are two books by Y.-T. \index[persons]{Siu, Y.-T.} Siu and  Siu-Trautmann dedicated
to this question (\cite{_Siu_techniques_,_Siu_gap_}).
The main result that we use in this chapter is \index[persons]{Siu, Y.-T.} Siu's theorem
(\ref{_Siu_extension_over_Z_Theorem_}) which allows
one to extend a vector bundle on $M \backslash Z$
to a coherent sheaf on $M$, if $M$ is smooth and
$\codim Z > 2$. This works nicely for Hopf manifolds:
from this theorem it follows that the pullback of a holomorphic bundle from
Hopf to $\C^n \backslash 0$ can be extended to 
a coherent sheaf on $\C^n$. If one needs to generalize 
this argument to arbitrary locally conformally K\"ahler
manifolds with potential\index[terms]{manifold!LCK!with potential}, one should use a stronger 
statement (\ref{_extension_over_a_point_Theorem_}), 
that is  obtained using the \index[persons]{Andreotti, A.} Andreotti-\index[persons]{Siu, Y.-T.}Siu extension
theorem. Let $F$ be a reflexive sheaf on $M \backslash Z$,
where $M$ is a normal variety and $\codim Z > 2$.
\ref{_extension_over_a_point_Theorem_} \index[terms]{theorem!Andreotti--Siu extension}
claims that its sheaf extension $j_* \caf$ is coherent and
reflexive, and the map $\caf \mapsto j_* \caf$
produces an equivalence between the category
$\Reff(M\backslash Z)$ of
reflexive sheaves on $M \backslash Z$ and $\Reff(M)$.

We do not pursue this direction further, restricting 
ourselves to \index[persons]{Mall, D.} Mall bundles on Hopf manifolds.

Let $\gamma$ be an invertible holomorphic contraction of\index[terms]{contraction!holomorphic}
$\C^n$ centred in 0, and $H= \frac{\C^n \backslash  0}{\langle \gamma\rangle}$
the corresponding Hopf manifold. 
Since the pullback of a Mall bundle to $\C^n \backslash 0$
is a trivial bundle, the category of Mall bundles
is equivalent to the category of  $\gamma$-equivariant
 bundles on $\C^n$.

Let $B$ be a complex vector bundle on a complex manifold $M$,
and $\nabla$ a flat connection. The Hodge component $\bar\6:=\nabla^{0,1}$
is a holomorphic structure operator on $B$. By \index[persons]{Koszul, J.-L.} Koszul-\index[persons]{Malgrange, B.}Malgrange theorem (\ref{_Koszul--Malgrange_Theorem_}), \index[terms]{theorem!Koszul--Malgrange}
the sheaf ${\cal B}:=\ker\bar\6$ is a holomorphic vector bundle on
$M$, with ${\cal B} \otimes_{\calo_M} C^\infty M = B$.
In this situation, we say that the flat\index[terms]{connection!flat}
connection $\nabla$ {\bf is compatible with the
  holomorphic structure on $B$} (\ref{_compatible_with_holo_Definition_}).

There are many examples of \index[persons]{Mall, D.} Mall bundles arising from the
geometry of Hopf manifolds. All tensor bundles, all line
bundles, and all extensions of Mall bundles are also Mall.
In many of those examples, the equivariant action of $\gamma$
on $B$ preserves a flat connection on $B$. In other words,
these Mall bundles are obtained from flat bundles
by taking the $(0,1)$-part of the connection.
It turns out that this situation is quite general, and 
an arbitrary \index[persons]{Mall, D.} Mall bundle admits a compatible flat
connection when the so-called ``non-resonance'' condition
is satisfied.\index[terms]{resonance}

This can be explained as follows. Recall that 
{\bf a holomorphic connection} (\ref{_holo_conne_Definition_}) on a holomorphic vector\index[terms]{connection!holomorphic}
bundle $B$ is a holomorphic differential operator
$\nabla:\; B \arrow B \otimes \Omega^1 M$ satisfying the
Leibniz rule, $\nabla(fb) = f \nabla(b) + df \otimes b$. 
We want to construct a holomorphic connection on a Mall
bundle on a Hopf manifold; this is equivalent to having
a $\gamma$-equivariant connection on its pullback to
$\C^n\backslash 0$ considered to be  a $\gamma$-equivariant
vector bundle.

Consider
the space ${\cal A}$ of all holomorphic connections on a trivial
$\gamma$-equi\-va\-riant holomorphic vector
bundle $R$ on $\C^n$. This is an affine  space
modeled on the vector space $H^0(\Omega^1 \C^n\otimes_{\calo_{\C^n}} \End(R))$,
and the equivariant action defines an affine endomorphism
of ${\cal A}$. The linearization $\rho$ of this action is a compact
endomorphism of $H^0(\Omega^1 \C^n\otimes_{\calo_{\C^n}} \End(R))$
considered to be  a topological vector space with $C^0$-topology\index[terms]{topology!$C^0$}
(\ref{_eigenvalues_diff_forms_contraction_Lemma_}).
If $\rho$ has all eigenvalues of absolute value $< 1$,
Banach fixed point theorem would imply that $R$
admits a $\gamma$-equivariant holomorphic connection.
In fact, it would suffice to check that
all eigenvalues $\lambda_i$ of $\rho$ 
are not equal to 1.

The ``resonance'' is a property of the
eigenvalues of the $\gamma$-equivariant action;
a $\gamma$-equivariant vector bundle $R$
on $\C^n$ {\bf has resonance} when the eigenvalues
of $D\gamma$ on the fibre $R \restrict 0$
are $\beta_1, ..., \beta_m$, the eigenvalues
of $D\gamma$ on $T_0 \C^n$ are $\alpha_1, ..., \alpha_n$,
and there exists a relation of the form\index[terms]{resonance!of a bundle}
$\beta_p = \beta_q \prod_{i=1}^n \alpha_i^{k_i}$, with
all $k_i$ non-negative integers, and $\sum_i k_i > 0$,
for some $\beta_p, \beta_q$, that are  not necessarily distinct.

Let $B$ be a vector bundle on a Hopf manifold, and
$R$ the extension of its pullback to $\C^n$. We prove that
$R$ has no resonance if and only if $H^0(H,
\Omega^1_H\otimes \End(B))=0$ (\ref{_B_on_Hopf_resonant_via_sections_Corollary_}).
 We also prove that any non-resonant holomorphic vector bundle
on a Hopf manifold admits a flat connection compatible
with the holomorphic structure (\ref{_Mall_flat_connection_Theorem_}).

\subsection{Flat affine structures and the development map}

The notion of resonance is classical, and harks back to \index[persons]{Poincar\'e, H.} Poincar\'e,
Latt\`es and \index[persons]{Dulac, H.} Dulac, who discovered the resonance while working on
the normal forms of ordinary differential equations.
In the modern language, they were looking at the normal form
of a real analytic or complex analytic vector field that has a simple
zero at a given point.

The ``normal form'' is a classical notion, that is  roughly equivalent,
in modern language, to ``the moduli space'', but includes a more 
explicit description in terms of coordinates. The original notion
of a normal form comes from the linear algebra: given a tensor
or collection of tensors of a given type, we are looking for
a basis in which this tensor is written in canonical form,
in such a way that equivalent tensors correspond to equivalent
expression in coordinates. Here, ``equivalent tensors'' means
tensors related by the $\GL(n)$-action.

For example, a pair of real quadratic
forms, one of them positive definite, has a normal form in
which the first one is represented by an identity
matrix, and the second is diagonal. An antisymmetric 2-form
is represented by a block matrix with all blocks zero or
$\small \begin{pmatrix} 0 & 1 \\ -1 & 0 \end{pmatrix}$, and so on.
If ${\cal W}$ denotes an appropriate set of tensors (antisymmetric
2-forms, pairs of quadratic forms and so on), {\em the normal
form} is a section from $\frac{{\cal W}}{\GL(n)}$ to ${\cal W}$.

For geometric objects on manifolds, the normal forms
are defined in a similar way, with the diffeomorphism
group playing the role of $\GL(n)$.

Let $S$ be a geometric
object on a manifold, such as a complex structure, a symplectic
structure, or a vector field; finding its normal form means
finding a local coordinate system in which this object is written
in a unique way, up to the diffeomorphism action. For example,
the Darboux theorem claims that a symplectic form admits
local coordinates $p_i, q_i$ in which it is written as $\sum dp_i \wedge dq_i$.
Another example is Morse lemma, which claims that, given a function $f\in C^\infty(M)$\index[terms]{theorem!Morse Lemma}
with a non-degenerate singularity in a given point, there exists a coordinate system
$x_1,..., x_n$ around this point, such that the function is written
as $f=a + \sum_{i=1}^n (\pm x_i^2)$. There are situations when the
normal form does not make sense. For example, the Riemannian
structures on a manifold do not have normal form, because (generally speaking)
the diffeomorphism group acts on the Riemannian structures freely.

For vector fields without zeroes, the normal form is  simple:
in an appropriate coordinate system, this vector field takes
the form $\frac d{dx_1}$; this is called ``straightening of a vector field''.
This result follows directly from the Peano and Picard theorems on
existence of solutions of ODE.\index[terms]{theorem!Peano--Picard}

The normal form theorem for a vector field with a simple zero 
is known as {\em  Poincar\'e-Dulac theorem},\index[terms]{theorem!Poincar\'e-Dulac} \cite{_Lattes_,_Dulac_,_Arnold:ODE+_}. 
We give a general outline of this theory, following \cite{_Enc_Math:Poincare_Dulac_}.

Let $\dot x = A(x) + u(x)$ be a formal
differential equation, where $x(t)\in \C^n$ is a time-dependent
point in $\C^n$, $A$ a non-degenerate linear operator, and $u(x)$ 
a Taylor series starting from the second order terms. It is said that {\bf $A$ has a resonance} if
there is a relation of the form $\lambda_i = \sum_{j=0}^n m_j\lambda_j$,
where $m_j \in \Z^{\geq 0}$ and $\sum_{j=0}^n m_j\geq 2$.\footnote{Later in this
chapter, we redefine this notion in such a way that this additive
relation becomes multiplicative, $\lambda_i = \prod_{j=0}^n \lambda_j^{m_j}$;
this is done because we work with holomorphic contractions and not with the vector fields.}

If $A$ has no resonance, then the normal form of 
this vector field is very simple: in appropriate coordinates
$y_1, ..., y_n$ it can be written as $\dot y = A(y)$.
If $A$ has a resonance, the vector field has a normal form, 
that is  written in a coordinate system  $y=(y_1, ..., y_n)$ as follows. 
Choose $y_i$ in such a way that $A$ is upper triangular
in the basis $\frac d{dy_i}$, and the diagonal
terms corresponding to $\frac d{dy_i}$ are $\lambda_i$.
Then the equation $\dot x = A(x) + w(x)$
has normal form $\dot y = A(y)+ \sum e_i w_i(y)$, where $e_i= \frac d{dy_i}$ is the
coordinate vector field, and $\sum e_i w_i(y)$ is a Taylor series
obtained as a sum of {\em resonant monomials}. A coordinate
monomial $e_i \prod_{j=1}^n y_j^{m_j}$ is {\bf resonant}
if $\lambda_i = \sum m_j \lambda_j$.

In general, it is hard to achieve convergence for these formal sums,
even when the differential equation is analytic. However, if $e^A$
is a contraction, the convergence is automatic, because the number
of resonant monomials is finite. Indeed, $e^A$ is a contraction
if and only if $\Re \lambda_i <0$ for all $i$, and 
the equation $\Re \lambda_i = \sum m_j \Re \lambda_j$, $m_i \in \Z^{\geq 0}$,
implies that $m_i \leq \max_{j,l}\frac {\Re \lambda_j}{\Re \lambda_l}$.

A similar result is true for germs of 
biholomorphic contractions, due to S. \index[persons]{Sternberg, S.} Sternberg (\cite{_Sternberg_contraction_}).
However, in this case, one should replace the linear resonance
by multiplicative, $\lambda_i = \prod_{j=0}^n \lambda_j^{m_j}$,
as in \ref{_resonant_matrix_Definition_}.

An invertible holomorphic contraction gives rise to a Hopf manifold,
and \index[persons]{Sternberg, S.} Sternberg's theorem can be interpreted as a structure theorem
about Hopf manifolds; this is how \index[persons]{Kodaira, K.} Kodaira used it in \cite{_Kodaira_Structure_II_}.\index[terms]{theorem!Sternberg}

In this chapter, we use the flat connection inherent on\index[terms]{bundle!vector bundle!Mall}
Mall bundles to give a new proof of the non-resonant part of the
 Poincar\'e-Dulac theorem. Let $\gamma$ be a germ of an invertible 
biholomorphic contraction of $\C^n$ with center in 0. We say that
$\gamma$ {\bf is non-resonant} if the differential $D_0 \gamma\in\End(\C^n)$
is a non-resonant matrix. Let $H$ be the corresponding Hopf manifold;
then the tangent bundle $TH$ is non-resonant, that is  equivalent
to $H^0(H, \Omega^1H \otimes \End(TH))=0$ (\ref{_B_on_Hopf_resonant_via_sections_Corollary_}).
This immediately implies that the flat connection in $TH$, given by
\ref{_Mall_flat_connection_Theorem_}, is torsion-free.\index[terms]{connection!torsion-free}

To prove the Poincar\'e-Dulac linearization theorem, we need to
find the coordinates on $\C^n$ in which $\gamma$ is linear.
To produce the flat coordinates, we use the developing map
defined within the framework of flat affine geometry\index[terms]{Cartan geometry}
(or, more generally, in \index[persons]{Cartan, E.} Cartan geometries).\index[terms]{map!developing}

Let's start from the definition of a Cartan geometry, to put flat affine
structures in a general context; however, the flat affine structures are
way more elementary, and the reader who is not interested in
Cartan geometries can skip this part.\index[terms]{structure!flat affine}

We follow \cite{_Alekseevsky_Michor:Cartan_}.
Let $\goth g \subset \goth h$ be Lie algebras, 
$G\subset H$ the corresponding Lie groups, and
$P\arrow M$ a principal $G$-bundle. {\bf A \index[persons]{Cartan, E.} Cartan connection
of type $\goth h/\goth g$} on a principal $G$-bundle $P$ \index[terms]{connection!Cartan}
is a trivialization of the tangent bundle $TP$ produced by a map 
$\kappa:\; \goth h\arrow TP$ that satisfies $[\kappa(x), \kappa(y)]= \kappa([x,y])$
when $x\in \goth g$, $y \in \goth h$: it is compatible with the
commutator if one of the terms belongs to $\goth g \subset \goth h$. 
Generally, $\kappa$ is not a Lie algebra homomorphism. 

The inverse map $\kappa^{-1}$ can be understood
as an $\goth h$-valued 1-form $\rho$ on $P$; it is also called
the \index[persons]{Cartan, E.} Cartan connection. {\bf The curvature} of the
Cartan connection $\rho$ is $\Theta_\rho:= d\rho + \frac 1 2 [\rho, \rho]$.
A Cartan connection is {\bf flat} if its curvature vanishes.
This is equivalent to $\kappa:\; \goth h\arrow TP$ being
a Lie algebra homomorphism (\cite[\S 4.1]{_Alekseevsky_Michor:Cartan_}).

{\bf The developing map} of a flat Cartan connection
$\rho$ is a map $\nu:\; P \arrow H$ such that
its derivative $d\nu\in T^* P \otimes {\goth h}$ is equal to $\rho$.
Locally, the developing map is uniquely determined by this relation,
because the local action of $H$ on $P$ is free and transitive. Indeed,
if $\rho$ is flat, we have a Lie algebra homomorphism $\kappa:\; \goth h \arrow TP$
that trivializes $TP$.
This also defines a canonical local diffeomorphism
$P/G\arrow H/G$, also called {\bf the developing map}.

\index[persons]{Cartan, E.} Cartan geometries are abstract and powerful, but
we are interested only in flat affine geometries, 
that are  a tiny special case of Cartan geometries.
In fact it is easier to study the flat affine manifolds
without referring to the Cartan connections.

A manifold $M$ is called {\bf affine}, or {\bf flat affine},
if it is equipped with an atlas of open sets, identified with
open subsets in $\R^n$, with the transition functions
affine. This is equivalent to having a torsion-free\index[terms]{connection!torsion-free}
flat affine connection on $M$\index[terms]{manifold!flat affine}
(\ref{_flat_affine_via_connection_Proposition_}).
The study of compact flat affine manifolds is 
ongoing, with many conjectures still open.
We refer to \cite{_Abels:survey_} for more details
and open questions.

For a flat affine manifold, the developing map 
can be defined directly, as follows. 
Let $M$ be a flat affine manifold, and
$P$ the principal $\GL(n,\R)$-bundle on $M$,
obtained as the bundle of frames in $TM$.
In this case, the \index[persons]{Cartan, E.} Cartan geometry is associated
with a pair $G\subset H$ of Lie groups, where
$G=\GL(n,\R)$ and $H=\Aff(n)$ is the group of
affine transforms of $\R^n$. The Cartan
connection is associated with the map 
$\kappa:\; \goth h \arrow TP$ constructed
as follows. For each frame $p\in P$ in $T_x M$,
$p= (e_1, ..., e_n)$ we can introduce local affine coordinates 
$z_1, ..., z_n$ on $M$ centred in $x$, with $e_i= \frac{d}{dz_i}$.
Using these coordinates, we define a local action of the group
$\Aff(n)$ on $P$, that is  by construction free and transitive.
This identifies a neighbourhood of $p\in P$
with an open set in $\Aff(n)$. This identification
takes the Lie algebra $\goth h$ of
left-invariant vector fields on $\Aff(n)$ 
to a subalgebra in $TP$. 

The developing map takes $P$ to $\Aff(n)$
by integrating these vector fields;
the corresponding map $M\stackrel \dev \arrow \frac{\Aff(n)}{\GL(n,\R)}=\R^n$
 preserves the flat affine structure.

Each flat affine manifold
is equipped with the natural flat, torsion-free connection $\nabla$.\index[terms]{connection!torsion-free}\index[terms]{connection!flat}
Using this connection, the developing map
can be defined as follows. Assume that $(M, \nabla)$
is a simply connected flat affine manifold.
Let $\theta_1, ..., \theta_n$ be parallel 1-forms
that trivialize the bundle $T^*M$. Since $\nabla$ is
torsion-free, all the forms $\theta_i$ are closed;
however, $H^1(M)=0$, which implies that $\theta_i$ are
exact, $\theta_i= dz_i$. Then the map $m \mapsto (z_1(m), ..., z_n(m))$
is identified with the developing map defined in terms
of \index[persons]{Cartan, E.} Cartan geometries, because the parallel transport along\index[terms]{parallel transport}
$u \in \goth h$ adds a constant to the coordinates, and hence 
the map $\dev:\; M \arrow \R^n$ is $\Aff(n)$-equivariant.

Further in this chapter, we define the developing map
in terms of the geodesics, as the inverse of the {\em exponential map}.
In this context, the exponential map \index[terms]{map!exponential}
takes a tangent vector to the point at time 1 on the geodesic 
tangent to this vector in time 0. This definition is equivalent to the one
given above, as indicated in the proof of \ref{_dev_for_complete_affine_Theorem_}.
A flat affine structure is {\bf complete}\index[terms]{structure!flat affine!complete}
when the geodesic equation $\nabla_{\dot \gamma_t} \dot \gamma_t$ can be solved
for all $t\in \R$ and all initial conditions $\gamma_0\in M$, $\dot \gamma_0\in T_{\gamma_0}M$.

The completeness condition is tricky and counter-intuitive;
indeed, even a compact flat affine manifold is not necessarily complete.
A textbook example of a non-complete flat affine manifold is a real linear Hopf
manifold $H$, obtained as a quotient of $\R^n \backslash 0$
by a linear contraction. This manifold is compact, but 
its universal cover is $\R^n \backslash 0$, and the 
developing map is an open embedding
$\R^n \backslash 0\hookrightarrow \R^n$. To obtain a non-complete geodesic,
one needs to start from a geodesic in $\R^n$ passing through
0; its image in $H$ is manifestly non-complete.

Let $M$ be a complete, simply connected affine manifold.
Then the developing map $\dev:\; M \arrow \R^n$
is an isomorphism of affine manifolds. This is a classical
result by \index[persons]{Auslander, L.} Auslander-\index[persons]{Markus, L.}Markus (\cite{_Auslander_Markus:holonomy_}) that we 
prove in \ref{_dev_for_complete_affine_Theorem_}.

For our present purposes, we need a variation
of this result, which ultimately implies the non-resonant
case of the  Poincar\'e-Dulac theorem. From \ref{_Mall_flat_connection_Theorem_},
it follows that any non-resonant 
Hopf manifold $M= \frac{\C^n \backslash 0}{\langle A \rangle}$ is equipped with a unique\index[terms]{connection!torsion-free}\index[terms]{connection!flat}
torsion-free flat affine connection compatible with the
complex structure. However, it is not complete, as we explained above.
We prove that this flat affine connection lifted to the universal 
covering $\C^n \backslash 0$ of $M$ can be extended to 0,
resulting in a complete flat affine structure on $\C^n$.
The corresponding developing map puts flat affine coordinates
on $\C^n$, and the contraction $A$ is affine in these coordinates, 
hence linear. This gives a new proof of the non-resonant
case of the Poincar\'e-Dulac theorem.

Part of this chapter reproduces results in \cite{ov_mall}.

\section{Coherent sheaves}

\subsection{Normal sheaves and reflexive sheaves}

We start by recalling some standard definitions and results of
complex geometry. In algebraic geometry\index[terms]{geometry!algebraic}, {\bf a coherent sheaf}
over an algebraic variety is a locally 
finitely generated sheaf\index[terms]{sheaf!coherent}
of $\calo_X$-modules, where $\calo_X$ is the sheaf of regular functions.
In complex geometry, {\bf a coherent sheaf} over a complex space $X$ 
is a locally finitely generated, locally finitely presented sheaf
of $\calo_X$-modules. \index[persons]{Oka, K.} Oka coherence theorem (a highly
non-trivial result) states that the kernel of
any surjective map of sheaves of $\calo_X$-modules
$\calo_X^n \arrow \caf$ is coherent, if $\caf$ is coherent.
Again, this result in algebraic geometry\index[terms]{geometry!algebraic} is trivial,
because the ring $H^0(U, \calo_U)$ of sections of $\calo_X$ 
over $U$ is Noetherian for any algebraic subvariety $U\subset X$.
In complex geometry, $H^0(\calo_U)$ is no longer Noetherian,
which makes Oka's coherency theorem non-trivial.\index[terms]{ring!Noetherian}\index[terms]{theorem!Oka's coherence}

However, here the difference between
coherent sheaves in complex geometry and algebraic
geometry ends: most results about coherent sheaves can be translated
from complex geometry to algebraic geometry\index[terms]{geometry!algebraic} and back.

An important exception from this rule is the 
existence of locally free resolutions: it is true
over projective manifolds, but for general complex
manifolds a locally free resolution often does not exist. 
For example, for a general complex torus $T$
of dimension $>2$, there is no locally free
resolution for a skyscraper sheaf $\calo_T/{\goth m}_x$
(\cite{_V:gen_tori_}).\index[terms]{sheaf!skyscraper}\index[terms]{resolution!locally free}

Given a coherent sheaf $\caf$ over $X$, let
$\caf^*:= \Hom(\caf, \calo_X)$ be {\bf the dual sheaf},\index[terms]{sheaf!reflexive}
that is, the sheaf of module homomorphisms to the ring\index[terms]{sheaf!dual}
of regular functions. The natural morphism of sheaves
$\caf\arrow \caf^{**}$ does not need to be an isomorphism:
for example, its kernel contains the torsion of $\caf$.
A coherent sheaf is called {\bf reflexive}
if the natural map $\caf\arrow \caf^{**}$ is an isomorphism.

For an introduction to the reflexive sheaves, 
see \cite{_oss_}; we just state some results that are 
relevant to this chapter. First of all, notice that
the sheaf $\caf^*$ is already reflexive
(\cite{_oss_}, Ch. II, Lemma 1.1.8).
Moreover, the singular set of a reflexive sheaf
over a normal variety has codimension $\geq 3$;
in particular, a reflexive sheaf over a complex
surface is locally free (\cite{_oss_}, Ch. II, 1.1.10).
Also, a reflexive sheaf of rank 1  over a smooth manifold
is locally free (\cite{_oss_},  Ch. II, Lemma 1.1.15).

The most important property of reflexive sheaves is {\bf normality}.
Let $Z\subset X$ be a subvariety of a normal 
complex variety, $\codim_X Z\geq 2$.\index[terms]{sheaf!coherent!normal}
Consider the open embedding map
$j:\; X \backslash Z \hookrightarrow X$, and let $j^*$ and $j_*$ be the
sheaf pullback and pushforward. For any sheaf $\caf$, there exists
a natural sheaf morphism  $\caf\arrow j_* j^* \caf$ taking a section
of $\caf$ to its restriction to $X \backslash Z$. A coherent sheaf
$\caf$ is {\bf normal} if 
the natural map $F\arrow j_* j^* \caf$ is an isomorphism, for any subvariety $Z\subset X$ of codimension $\geq 2$.

\hfill

The following theorem can be understood as 
a sheaf version of the Hartogs extension
theorem. \index[terms]{theorem!Hartogs}

\hfill

\theorem\label{_reflexife_normal_Theorem_}
A coherent sheaf $\caf$ over a normal variety 
is normal if and only if it is reflexive.\index[terms]{variety!normal}

\hfill

\proof
In algebraic category, this 
result is stated in \cite[Ch. II, Lemma 1.1.12]{_oss_}
for smooth $X$; however, the same proof works for normal varieties as well.
In analytic category, this is \cite[Proposition 7]{_Serre:Prolongement_}.
\endproof

\hfill

\corollary\label{_Hartogs_for_sheaves_Corollary_}
Let $\caf$ be a reflexive sheaf over a normal complex variety $X$,
$Z\subset X$ an algebraic subvariety, $\codim_X Z\geq 2$, 
and $f$ a section of $\caf$ over $X \backslash Z$. Then $f$
can be extended to a section of $\caf$ over $X$.
\endproof

%
%
%


\subsection{Extension of coherent sheaves on complex varieties}
\label{_extension_sheaves_Subsection_}

Here we present several results about extension
of coherent sheaves and their sections defined
on $M \backslash K$, where $K$ is a Stein domain 
with strictly pseudoconvex boundary. The work in this
direction started from the paper of \index[persons]{Andreotti, A.} Andreotti and \index[persons]{Grauert, H.} Grauert
\cite{_Andreotti_Grauert_}
and culminated in 1970-ies with the paper of \index[persons]{Andreotti, A.} Andreotti and \index[persons]{Siu, Y.-T.} Siu
and several papers and books of Siu and others 
(see 
\cite{andreotti_siu,_Siu_gap_,_Siu_techniques_,_Banica_Stanasila:French_} and
the references therein).

\hfill

We need to be able to extend coherent sheaves\index[terms]{sheaf!coherent}
over a Stein subset with strictly pseudoconvex boundary.

\hfill

\theorem \label{_extension_coherent_Theorem_}
{(\cite[Proposition 6.1]{andreotti_siu})}
Let $X$ be a compact complex manifold with boundary, consisting of
two components, strictly pseudoconvex and strictly\index[terms]{completion!Stein}
pseudoconcave. Suppose that $X$ admits a smooth strictly plurisubharmonic
function that is  constant on the components of the boundary. 
Denote by $\hat X$ its Stein completion.
Assume that $\dim X \geq 3$.
Then any coherent sheaf on $X$ can be extended to 
a coherent sheaf on $\hat X$.\index[terms]{theorem!Andreotti--Siu extension} 
\endproof

\hfill

\remark \label{_AS_sheaf_extension_unique_Remark_}
The extension of a sheaf $\caf$ from $X$ to $\hat X$
is not unique, because two non-reflexive sheaves can be isomorphic
outside of a finite set. However, when $\caf$ is reflexive, this
extension is unique up to isomorphism. Indeed, suppose that
$F_1$ and $F_2$ are two extensions of $\caf$ to $\hat X$.
Consider the identity map as a sheaf homomorphism
 $\nu \in \Hom(\caf_1, \caf_2)\restrict X$. By 
\cite[Proposition 2]{_Vajaitu_}
and \cite[pages 356-357]{_Banica_Stanasila:French_}, $\nu$
is extended to a map $\hat \nu \in \Hom(\caf_1, \caf_2)$
that is  an isomorphism outside of a closed
complex analytic (and hence, finite) subset
in the Stein manifold $\hat X \backslash X$.
Since $\caf_1, \caf_2$ are reflexive, they are normal
(\ref{_reflexife_normal_Theorem_}), 
and $\nu$ is an isomorphism everywhere.

\hfill

For our purposes, the full strength of 
\ref{_extension_coherent_Theorem_}
is often not necessary, and an earlier theorem of \index[persons]{Siu, Y.-T.} Siu
would usually suffice.

\hfill

\theorem\label{_Siu_extension_over_Z_Theorem_}
{(\cite[Main Theorem]{_Siu:Extension_})}
Let $M$ be a complex manifold,\index[terms]{theorem!Siu extension}
$Z\subset M$ a complex subvariety of codimension $\geq 3$,
and $B$ a holomorphic bundle on $M \backslash Z$.
Then $B$ can be extended to a coherent sheaf on $M$.
\endproof

\hfill

In \ref{_Siu_extension_over_Z_Theorem_}, the condition $\codim Z \geq 3$ is essential. 
Indeed, in \cite{_Ballico:bundles_on_C^2_} 
E. \index[persons]{Ballico, E.} Ballico has constructed
a holomorphic vector bundle on $\C^2\backslash 0$
without holomorphic rank 1 subsheaves; its sheaf extension
to $\C^2$ would be obviously non-coherent.

\hfill

The result that can be applied to 
reflexive sheaves on open algebraic cones\index[terms]{cone!algebraic} is the following.

\hfill

\theorem \label{_extension_over_a_point_Theorem_}
Let $M$ be a normal complex variety with isolated singularities,
$\dim_\C M \geq 3$, $x\in M$ a point, 
and $M_0:= M \backslash \{x\}\stackrel j\hookrightarrow M$.
Consider a reflexive sheaf ${\cal F}$ on $M_0$.
Then $j_*{\cal F}$ is a reflexive coherent sheaf.
Moreover, $j_*{\cal F}$ is the unique coherent
reflexive sheaf on $M$ that is  isomorphic to ${\cal F}$
on $M_0$.

\hfill

\pstep
Uniqueness of a reflexive extension
follows from \ref{_reflexife_normal_Theorem_};
indeed, a normal sheaf on $M$ is determined by
its restriction to any complement $M\backslash Z$,
when $Z$ is a complex subvariety of codimension $\geq 2$.

\hfill

{\bf Step 2:}
The existence of coherent extension of ${\cal F}$ follows from
\ref{_extension_coherent_Theorem_}, but to apply it,
we need to find an appropriate open neighbourhood of $x\in M$ 
with strictly pseudoconvex boundary. 

We show that $x\in M$ has a base of neighbourhoods
$\{U_i\}$ with smooth and strictly pseudoconvex boundaries.
Let $U$ be a neighbourhood of $x$ equipped with a complex-analytic embedding
$\phi:\; U \arrow \C^n$, mapping $x$ to 0; such an embedding exists, because
some neighbourhood of $x$ is Stein. Denote by $B_\epsilon$
an open ball of radius $\epsilon$ centred in 0.
Since $\6 B_\epsilon$ is strictly pseudoconvex, the same
is true for $\phi^{-1}(\6 B_\epsilon)$ (a complex analytic
map is compatible with the Levi form\index[terms]{form!Levi} on any codimension 1 
real submanifold).

\hfill

{\bf Step 3:} 
Let $U_i:= \phi^{-1}(B_{1/i})$, where $\phi$ is defined as above.
Then ${\cal F}\restrict M\backslash U_i$ can be extended
to a coherent sheaf on $M$ by Andreotti--Siu (\ref{_extension_coherent_Theorem_}).
By \ref{_AS_sheaf_extension_unique_Remark_}, this\index[terms]{theorem!Andreotti--Siu extension}
extension is unique, and hence  it is independent on  the
choice of an open set $U_i$. This implies that the pushforward
$j_* {\cal F}$ is a subsheaf of the coherent sheaf
$\underline {\cal F}$.
Replacing $\underline {\cal F}$ with its reflexization
$\underline {\cal F}^{**}$,
we may also assume that $\underline {\cal F}$ is
reflexive; indeed, ${\cal F}$ is reflexive, and hence 
the reflexization of its extension is equal to 
${\cal F}$ on $M_0$. However, a reflexive sheaf
is normal, and hence  \index[terms]{sheaf!normal}\index[terms]{sheaf!reflexive}
$\underline {\cal F}=j_* \left (\underline {\cal F}\restrict {M_0}\right)=
j_* \underline {\cal F}$.
\endproof

\hfill

For a broader discussion about the extension theorems for coherent sheaves,
see \cite{_Serre:Prolongement_,_Siu_gap_,_Siu_techniques_,_Banica_Stanasila:French_,_Bando_Siu_,_Vajaitu_,_Ovrelid_Vassiliadou_}.

\section{Dolbeault cohomology of Hopf manifolds}

\subsection{Degree of a line bundle}

\hfill

\definition
Let $M$ be a compact complex manifold, $\dim_\C M=n$, 
and $\omega^{n-1} \in \Lambda^{n-1, n-1}(M)$ a
Gauduchon form, that is, the $n-1$-th power of a 
Gauduchon (1,1)-form (\ref{_Gauduchon_definition_}). 
Consider a holomorphic line bundle $L$\index[terms]{form!Gauduchon}
on $M$; choose any Hermitian form\index[terms]{form!Hermitian} $h$ on $L$, and let
$\Theta_h\in \Lambda^{1,1}(M)$ be the curvature
of its  Chern connection.\index[terms]{connection!Chern} The {\bf degree} of
$L$ is the number $\1 \int_M \Theta_h \wedge \omega^{n-1}$.
\index[terms]{bundle!line!holomorphic!degree of}

\hfill

\remark
The degree is independent on  the choice of a Hermitian metric on $L$.
Indeed, for any Hermitian metrics $h_1, h_2$, with 
$h_1 = p h_2$, where $p$ is a positive function, 
and we have $\Theta_{h_1} = \Theta_{h_2} - \1 dd^c \log p$;
however, $\int_M dd^c \log p \wedge \omega^{n-1}=0$
because $dd^c \omega^{n-1}=0$ for any Gauduchon form.\index[terms]{form!Gauduchon}

%
%
%
%
%

\subsection{Computation of $H^{0,p}(H)$ for a Hopf manifold}

Dolbeault cohomology of Hopf manifolds is a classical subject,
but we could not find the computation for the general case.
For the classical Hopf manifold (a quotient of $\C^n \backslash 0$
by a constant times identity), an answer is given in \index[terms]{manifold!Hopf!classical}
\cite{_Ise_}. More general Hopf surfaces were defined and classified by \index[persons]{Kodaira, K.} Kodaira
(\cite{_Kodaira_Structure_II_}); he \index[terms]{surface!Hopf}
computed some of their cohomology in \cite{_Kodaira_Structure_III_}.
For a diagonal Hopf manifold, the Dolbeault cohomology
was computed by D. \index[persons]{Mall, D.} Mall\index[terms]{manifold!Hopf!diagonal}
(\cite{_Mall_}). For a reference to other special cases of this
theorem, see \cite{_Ramani_Sankaran_,_MO:Hopf_cohomology_}.

\hfill

Using Mall's theorem, one can easily compute the $(0,*)$-part of the 
Dolbeault cohomology of a Hopf $n$-manifold $H$. 
Indeed, $H^{0,1}(H)= H^1(\calo_H)$ has\index[terms]{theorem!Mall}
the same rank as $H^0(\calo_H)=\C$, and $H^{0,i}(H)=0$
for $1 < i < n-1$. By Serre duality, $H^{0,n-i}(H)= H^i(K_H)^*$,
where $K$ is the canonical bundle, and hence  to prove that\index[terms]{bundle!vector bundle!canonical}\index[terms]{duality!Serre}
$H^{0,n-1}(H)= H^1(K_H)^*$ and $H^{0,n}(H)= H^0(K_H)^*$
vanishes, it would suffice to show that the canonical
bundle on a Hopf manifold has no holomorphic sections.

\hfill

\theorem\label{_Dolbeault_for_Hopf_Theorem_}
Let $H= (\C^n \backslash 0)/\Z$ be a Hopf manifold,
that is, a quotient of $\C^n \backslash 0$ by a holomorphic contraction,
Then $H^i(\calo_H)=0$ unless $i = 0,1$, and
$\rk H^1(\calo_H)=\rk H^0(\calo_H)=1$.

\hfill

\proof
By \index[persons]{Mall, D.} Mall's theorem \index[terms]{theorem!Mall}
(\ref{_Mall_cohomology_Theorem_}), $\rk H^1(H,\calo_H)=\rk H^0(H,\calo_H)$,
and $H^i(H,\calo_H)=0$ for $1 <i < n-1$.
However, $\rk H^0(H,\calo_H)=1$ because $\calo_H$ is a trivial line bundle.
To finish the proof, it remains only to show that
$H^{n-1}(H,\calo_H)$ and $H^{n}(H,\calo_H)$ vanish.
By Serre duality, these two spaces are dual to $H^0(H, K_H)$ and
$H^1(H, K_H)$, that have the same rank by Mall's theorem again. 
It remains only to prove that
$H^0(H, K_H)=0$ for any Hopf manifold.\index[terms]{duality!Serre}

Suppose that $\eta$ is a non-zero element in $H^0(H, K_H)=0$;
we consider $\eta$ as a holomorphic volume form. Then
$\mu:\; =\eta\wedge \bar \eta$ is a measure on $H$, that is 
strictly positive outside of the zero divisor of $\eta$.

Consider the measure $\pi^* \mu:=\pi^* \eta\wedge \pi^* \bar\eta$ 
on $\C^n \backslash 0$. Since $\pi^* \eta$ is $\Z$-invariant, 
the measure $\pi^* \mu$ is $\Z$-invariant as well.
The canonical bundle of $\C^n$ is trivial, and hence ,
by Hartogs theorem, $\pi^*\eta$ can be extended to
a holomorphic section $j_* \pi^*\eta$ of $K_{\C^n}$.
Denote by $j_* \pi^* \mu$ the measure
$j_*\pi^* \mu:=j_*\pi^* \eta\wedge j_*\pi^* \bar\eta$.
This measure is finite on compacts, and is preserved
by the contraction $\gamma:\; \C^n \arrow \C^n$. 
This is impossible, unless $j_*\pi^* \mu=0$, because any bounded set
is mapped inside a given compact neighbourhood
of 0 by a sufficiently big power of $\gamma$.
This implies that $\eta=0$, and hence  $H^0(H, K_H)=0$.
\endproof

\subsection{Holomorphic differential forms on Hopf manifolds}

The results in this subsection generalize the vanishing 
$H^0(H, K_H)=0$ given in the proof of \ref{_Dolbeault_for_Hopf_Theorem_}.
It turns out that all holomorphic differential forms on Hopf
manifolds vanish. 

\hfill

In \cite{_Ise_} the vanishing of differential forms was proven
for the classical Hopf manifold $\frac{\C^n \backslash 0}{\lambda \Id}$;
we could not find other results in the literature, though the question
seems to be elementary and classical.


%
%
%

\hfill

We start from the following lemma.

\hfill

\lemma \label{_eigenvalues_diff_forms_contraction_Lemma_}
Let $\gamma:\; \C^n \arrow \C^n$ be an invertible holomorphic contraction
centred in 0, and $D \subset \C^n$  an open set such that
$\gamma(D)$ is precompact in $D$. Choose an Hermitian
metric on $\C^n$, and define the norm on the space
$H^0_b(D, \Omega^1D)$ of bounded holomorphic $1$-forms 
as $\|\eta\|:= \sup_{x\in D} |\eta_x|$.\footnote{%
By Montel theorem, $H^0_b(D, \Omega^1D)$
with this norm is a Banach space\index[terms]{space!Banach}.}  
Then the operator \index[terms]{theorem!Montel}
$\gamma^*:\; H^0_b(D, \Omega^1D)\arrow  H^0_b(D, \Omega^1D)$
is compact, and all its eigenvalues are smaller than 1 in absolute value.

\hfill

\remark
This statement can be formulated for any tensor bundle
on $\C^n$, but it is true only on differential forms and their tensor powers.
The $\gamma^*$-action on vector fields (or tensor fields containing
vector fields as a multiplier) is compact by \ref{_contra_compact_Theorem_}, 
but the eigenvalues are not always smaller than 1 in absolute value.
Indeed, by \ref{_S^1_action_exists_Corollary_} and \ref{_gene_Hopf_is_LCK_Corollary_},
there exist non-zero holomorphic vector fields on all Hopf manifolds,
which implies that the action of $\gamma^*$ on holomorphic vector fields can be 1. 
 
\hfill

{\bf Proof of \ref{_eigenvalues_diff_forms_contraction_Lemma_}:}
The operator 
$\gamma^*: H^0_b(D, \Omega^1D)\arrow  H^0_b(D, \Omega^1D)$
is compact by an argument that follows from Montel theorem\index[terms]{theorem!Montel}
(\ref{_contra_compact_Theorem_}). To obtain the
eigenvalue estimate, we use the same argument which proves
\ref{_gamma_on_T^*_contraction_Corollary_}.

The differential of $\gamma$
acts with all eigenvalues $|a_i|<1$ on $T_0\C^n$,
by \ref{_contra_compact_Theorem_}. 
Then $\gamma^*$ (the pullback of differential forms)
acts on $T^*_0\C^n$ with all eigenvalues $|a_i|<1$.
Consider the Taylor expansion of a function $f$ in 0.
The chain rule and the estimate of the eigenvalues
of $\gamma^*$ on $T^*_0\C^n$ imply that $\gamma^*$ acts
on the non-constant Taylor coefficients of $f$ with all eigenvalues
$|\cdot|<1$.

Let $a$ be an eigenvalue of the
compact operator \index[terms]{operator!compact}
$\gamma^*: H^0_b(\Omega^1D) \arrow H^0_b(\Omega^1D)$.
It remains to show that $|a|<1$.

The eigenvalues of $\gamma^*$ on the vector space
$\Omega^1\C^n\restrict 0$ are powers of its action on
the cotangent bundle, that are   all smaller than 1 in absolute value,
because $\gamma$ is a contraction. This implies that the eigenvalues
of $\gamma^*$ on the degree 0 term of the Taylor
expansion of $\eta$ are also smaller than 1 in absolute value.
Summarizing the above estimates, we obtain that
$\lim_n (\gamma^*)^n \eta=0$, and hence  $\gamma^*$
acts on differential forms with all eigenvalues
smaller than 1 in absolute value.
\endproof

\hfill

\proposition\label{_holo_forms_Hopf_vanish_Proposition_}
Let $H$ be a Hopf manifold, and $B$ a tensor power of $\Omega^1H$,
$B= (\Omega^1H)^{\otimes l}$. Then all holomorphic sections of $B$ vanish.
 
\hfill

\pstep
The universal cover of the Hopf manifold is $\C^n \backslash 0$
equipped with the free and properly discontinuous holomorphic $\Z$-action.\index[terms]{action!$\Z$-}
Choose precompact fundamental domains  $M_i$ of the $\Z$-action,
with the generator of $\Z$ mapping $M_i$ to $M_{i+1}$.
 Then we may assume that
$\tilde M_0 :=\{0\}\cup \bigcup_{i\leq 0} M_i$
is a bounded neighbourhood of 0 in $\C^n$.

Let $\gamma\in \Aut(\C^n)$ be the generator of the $\Z$-action
acting on $\C^n$ as a contraction.
By  \ref{_eigenvalues_diff_forms_contraction_Lemma_}, 
the action of $\gamma^*$ on $B$ defines a compact operator \index[terms]{operator!compact}
\[ H^0_b(\tilde M_0, B) \arrow H^0_b(\tilde M_0, B),\]
with all eigenvalues smaller than 1.

Given a holomorphic differential form $\alpha$ on $\C^n$,
we restrict it to $\tilde M_0$ and observe that it is bounded
because $\tilde M_0$ is precompact.
Using the compactness of $\gamma^*$-action and the estimate of its eigenvalues
(\ref{_eigenvalues_diff_forms_contraction_Lemma_}), we show that
the norm of  $(\gamma^n)^*\alpha$ on $\tilde M_0$ converges to 0.

\hfill

{\bf Step 2:} 
Let $\pi:\; \C^n\backslash 0 \arrow H$ be the universal cover of a Hopf
manifold, and $\beta\in H^0(H, B)$ a section of $B$.
The sheaf $(\Omega^1 \C^n)^{\otimes l}$ is normal by Hartogs theorem.
Therefore, the section 
$\pi^*\beta$ can be holomorphically extended to 
0, and this extension is
$\gamma^*$-invariant. By the argument in Step 1,
$\lim_n (\gamma^n)^*\alpha=0$
for any holomorphic form $\alpha$; applying
this to $\pi^*\beta$, we conclude that $\beta=0$.
\endproof

\section{Mall bundles on Hopf manifolds}

\subsection{Mall bundles: definition and examples}

We define {\bf \index[persons]{Mall, D.} Mall bundles} on a Hopf manifold
as bundles that satisfy the assumptions of 
\ref{_Mall_cohomology_Theorem_}.

\hfill

\definition\label{_Mall_bundle_Definition_}
Let $\pi:\;  \C^n\backslash 0\arrow H$ be the
universal cover of a Hopf manifold,  
$j:\; \C^n\backslash 0\hookrightarrow \C^n$
the standard embedding map, and 
$B$ a holomorphic vector bundle over $H$
such that $j_* \pi^* B$ is a locally trivial\index[terms]{bundle!vector bundle!Mall} 
coherent sheaf on $\C^n$, that is, a holomorphic vector bundle. 
Then $B$ is called {\bf a \index[persons]{Mall, D.} Mall bundle}.

\hfill

Since $j_*$ commutes with direct sums
and tensor products, the following observation is clear.

\hfill

\claim\label{_tensor_direct_component_Mall_Claim_}
A tensor product of Mall bundles, a direct sum
of \index[persons]{Mall, D.} Mall bundles, and any direct sum component
of a Mall bundle is again a Mall bundle.

\hfill

\proof For direct sums and tensor products the statement
is obvious, because these operations commute with the
functor $j_* \pi^*$. For a direct sum component, 
consider a coherent sheaf ${\cal F}$ on a complex manifold $M$, and let
${\cal F}|_x:= {\cal F}\otimes_{\calo_M} (\calo_M/{\goth m}_x)$.
This is a finite-dimensional space; if ${\cal F}$ is a vector bundle,
 ${\cal F}|_x$ is the fibre of ${\cal F}$ in $x$.
By Nakayama lemma, a coherent sheaf\index[terms]{lemma!Nakayama}
is generated by any collection $\{s_i\}$ of sections which
generates ${\cal F}|_x={\cal F}\otimes_{\calo_M} (\calo_M/{\goth m}_x)$ for all $x\in M$.
Indeed, let ${\cal F}_x$ be the stalk of ${\cal F}$ in $x$.
Then ${\cal F}_x$ is a finitely--generated module
over a local ring $\calo_{M, x}$ of germs of holomorphic functions,
that is  Noetherian by \index[persons]{Lasker, E.} Lasker theorem (\cite[Chapter II, Theorem B.9]{_Gunning_Rossi_}).\index[terms]{theorem!Lasker} 
Nakayama lemma tells that a finitely--generated module $A$ over a Noetherian local
ring is generated by any set of elements which
generate $A$ modulo the maximal ideal.
Therefore, every stalk ${\cal F}_x$ of ${\cal F}$  is generated
by the images of $s_i$ in  ${\cal F}_x$; by definition,
this implies that ${\cal F}$ is generated by $\{s_i\}$.

We have shown that ${\cal F}$ is a vector bundle
whenever $\rk {\cal F}_x$ is constant in $x$.

Given a direct sum decomposition
${\cal F} = {\cal F}' \oplus {\cal F}''$,
we immediately obtain $\rk {\cal F}_x = \rk {\cal F}'_x + \rk {\cal F}''_x$.
If $\rk {\cal F}|_x$ is constant in $x$, this implies that 
$ \rk {\cal F}'_x =\const$ and $\rk {\cal F}''_x=\const$
because $\rk {\cal F}'_x$ is upper semicontinuous as
a function of $x$ by Nakayama lemma. This
implies that for any direct sum decomposition
${\cal F} = {\cal F}' \oplus {\cal F}''$
of a vector bundle onto a direct sum of coherent
sheaves, the summands are also vector bundles.
\endproof

\hfill

For the next proposition, we need the following claim,
that is  almost trivial.

\hfill

\claim\label{_third_Mall_Claim_}
Let $0 \arrow A \arrow B \arrow C\arrow 0$
be an exact sequence of holomorphic bundles
on a Hopf manifold. Assume that $C$ and another of these bundles,
$B$ or $A$, are \index[persons]{Mall, D.} Mall. Then the third one is also
Mall.\index[terms]{bundle!vector bundle!Mall} 

\hfill

{\bf Proof:}
The functor $\pi^*$ is exact, and $j_*$ is left exact.
Therefore, applying the functor $j_*\pi^*$ to this sequence,
we obtain the following long exact sequence
\begin{equation}\label{_exact_exte_Equation_}
0 \arrow j_*\pi^*A \stackrel \alpha \arrow j_*\pi^*B \stackrel \beta \arrow j_*\pi^*C\arrow 
R^1 j_* \pi^* A^* \arrow ...
\end{equation}
where $R^ij_*$ is the higher derived pushforward functor
(Section \ref{_Derived_functors_Subsection_}). 
If $C$ is \index[persons]{Mall, D.} Mall, then the sheaf  $j_*\pi^*C$ is locally free. Therefore,
the map $\beta$ has a section $\beta':\;j_*\pi^*C\arrow j_*\pi^*B $. This gives 
a subsheaf $\alpha(j_*\pi^*A)\oplus  \beta'(j_*\pi^*C) \subset  j_*\pi^*B$,
that is  by construction normal and coincides with
$j_*\pi^*B$ outside of 0. By normality, it has
to coincide with $j_*\pi^*B$, thus we have a direct sum decomposition
$j_*\pi^*A\oplus  \beta'(j_*\pi^*C) = j_*\pi^*B$.
Therefore, one of the sheaves $j_*\pi^*A$ and $j_*\pi^*B$
is locally free when the other one is locally free
(the sheaf $j_*\pi^*C=\beta'(j_*\pi^*C)$ is locally free by assumption).
\endproof

\hfill

We now give a list of examples of Mall bundles. 

\hfill

\proposition\label{_Mall_exa_Proposition_}
Let $H$ be a Hopf manifold. 
Then any line bundle on $H$, any tensor bundle 
$TH^{\otimes p}\otimes_{\calo_H} T^*H^{\otimes q}$,
their tensor products and direct sum components are Mall.

\hfill

\proof
Tensor products and direct sum components were treated 
in \ref{_tensor_direct_component_Mall_Claim_}.
Let $L$ be a line bundle on a Hopf manifold,
and $j_*\pi^*L$ the corresponding sheaf on $\C^n$
(\ref{_Mall_bundle_Definition_}).
By \ref{_Siu_extension_over_Z_Theorem_}, $j_*\pi^*L$ is 
a normal (and, therefore, reflexive) coherent sheaf of
rank 1, and hence  it is locally free
(\cite{_oss_},  Ch. II, Lemma 1.1.15).
This implies that $L$ is \index[persons]{Mall, D.} Mall.

The tangent bundle $TH$ is Mall because
$\pi^* TH$ is the tangent bundle $T(\C^n\backslash 0)$;
since the tangent bundle $T\C^n$ is normal, (\ref{_reflexife_normal_Theorem_}), 
$j_* \pi^* TH= T\C^n$.
\endproof

\subsection{The Euler exact sequence and 
an example of a non-Mall bundle on a classical Hopf manifold}\index[terms]{exact sequence!Euler}

An example of a non-Mall bundle can be constructed as follows.
Consider the Euler exact sequence\index[terms]{exact sequence!Euler} over $\C P^{n-1}$
\[
0 \arrow \calo_{\C P^{n-1}} 
\stackrel \tau \arrow \calo_{\C P^{n-1}}(1)^n \arrow T\C P^{n-1} \arrow 0,
\]
where the image of $\tau$ is 
obtained by taking a line in 
$\calo_{\C P^{n-1}}(1)^n\restrict{\cal x}$ 
generated by $(x_1, ..., x_n)$ in every point ${\cal x}=[x_1:x_2:...:x_n]$ in
$\C P^{n-1}$. Here we understand $H^0(\C P^{n-1},\calo_{\C P^{n-1}}(1))$ as
the space of homogeneous polynomials in $x_1, ..., x_n$ of degree 1, and 
$(x_1, ..., x_n)$ is understood as a collection of $n$ homogeneous 
polynomials. The quotient $\calo_{\C P^{n-1}}(1)^n /\im \tau$
is identified with the tangent bundle using the following chain
of identifications.

Let $\phi:\; \C^n \backslash 0 \arrow \C P^{n-1}$ be the
standard projection defining the complex projective space. Since $H^0(\C P^{n-1},\calo_{\C P^{n-1}}(1))$ is
the space of homogeneous polynomials of degree 1, and $T\C^n$ is trivial,
the pullback bundle 
\begin{multline} \phi^*(\calo_{\C P^{n-1}}(1)^n)
=\phi^*((\calo_{\C P^{n-1}})^n \otimes \calo_{\C P^{n-1}}(1))\\
= (\calo_{\C^n})^n\otimes \phi^*\calo_{\C P^{n-1}}(1)=
T\C^n\otimes \phi^*\calo_{\C P^{n-1}}(1) 
\end{multline}
can be identified with the bundle of tangent vector fields
of homogeneity 1. The Euler vector field, associating\index[terms]{vector field!Euler}
the vector $(x_1, ..., x_n)$ to a point $(x_1, ..., x_n)\in \C^n \backslash 0$,
is homogeneous of degree 1; this gives an embedding
$$e:\; \calo_{\C^n \backslash 0} \hookrightarrow \phi^*(\calo_{\C P^{n-1}}(1)^n)$$
and the quotient $\frac{\phi^*(\calo_{\C P^{n-1}}(1)^n)}{\im e}$ is
projected bijectively to $T\C P^{n-1}$.
We have obtained the exact sequence
\begin{equation}\label{_pullback_of_Euler_equation_}
0 \arrow \phi^*\calo_{\C P^{n-1}}  
\arrow \phi^*\calo_{\C P^{n-1}}(1)^n \arrow \phi^*T\C P^{n-1} \arrow 0
\end{equation}
that is  identified with the pullback of the Euler 
exact sequence.\index[terms]{exact sequence!Euler}

Now, consider the exact sequence \eqref{_pullback_of_Euler_equation_},
and let $\rho:\; \C^* \arrow \Aut(\C^n \backslash 0)$
denote the homothety action. By construction, the
sequence \eqref{_pullback_of_Euler_equation_} is $\C^*$-equivariant.
Applying \ref{_covering_equivariant_Theorem_}, we obtain that
the exact sequence \eqref{_pullback_of_Euler_equation_}
is the pullback of the exact sequence \eqref{_Euler_on_Hopf_Equation_} 
of bundles on a classical 
Hopf manifold $H = \frac{\C^n \backslash 0}{\lambda \cdot \Id}$,
where $\lambda$ is a complex number, $|\lambda|>1$:
\begin{equation}\label{_Euler_on_Hopf_Equation_} 
0 \arrow \calo_H \arrow \calo_H(1)^n \arrow \phi_H^*T\C P^{n-1}\arrow 0,
\end{equation}
where $\phi_H:\; H \arrow \C P^{n-1}$ is the standard projection
of the classical Hopf manifold to $\C P^{n-1}$.

\hfill

\claim
In the above assumptions, $\phi_H^*T\C P^{n-1}$ is not a Mall bundle
when $n >2$.

\hfill

\proof
Let $j:\;\C^n \backslash 0\arrow \C^n$ be the natural embedding. 
Consider the exact sequence 
\[
0 \arrow j_*\calo_H \arrow j_*\calo_H(1)^n \arrow j_* \phi_H^*T\C P^{n-1}
\arrow R^1j_*\calo_H  \arrow ...,
\]
where $R^1j_*$ is the higher derived pushforward functor
(Section \ref{_Derived_functors_Subsection_}). 
A standard argument implies that the derived functor $R^1j_*\calo_H=0$
(Exercise \ref{_R^1j_*_vanish_Exercise_}). Let $k_0:= \calo_{\C^n}/{\goth m}_0$ 
denote the skyscraper sheaf in 0. 
If $j_* \phi_H^*T\C P^{n-1}$ is
a vector bundle, the group $\Tor^1(k_0, \phi_H^*T\C P^{n-1})$
would vanish,\footnote{By definition, $\Tor^1(k_0, \shortdash)$ is the left higher derived
functor for the tensor product $k_0\otimes_{\calo_{\C^n}} \shortdash$.}
and we would have an exact sequence
\[
0 \arrow k_0\otimes_{\calo_{\C^n}}j_*\calo_H \arrow k_0\otimes_{\calo_{\C^n}}j_*\calo_H(1)^n \arrow k_0\otimes_{\calo_{\C^n}}j_* \phi_H^*T\C P^{n-1}\arrow 0.
\]
However, the map 
$k_0\otimes_{\calo_{\C^n}}j_*\calo_H \arrow k_0\otimes_{\calo_{\C^n}}j_*\calo_H(1)^n$
vanishes because the Euler vector field vanishes in 0, and hence 
$j_* \phi_H^*T\C P^{n-1}$ cannot be locally free.
\endproof


\section{Resonance in Mall bundles}


\subsection{Resonant matrices}

\definition\label{_resonant_matrix_Definition_}
Let $A\in \GL(n, \C)$ be a matrix with eigenvalues $\alpha_1,..., \alpha_n$.
This matrix is called {\bf resonant} if there exists
a relation $\alpha= \prod_{i=1}^n \alpha_i^{k_i}$,
with $\alpha$ an eigenvalue of $A$ and 
$k_i \in \Z^{\geq 0}$, $\sum_i k_i \geq 2$, and {\bf non-resonant}
otherwise. \index[terms]{resonance}

\hfill

The reason we need this definition is the following elementary lemma.

\hfill

\lemma\label{_inva_tensors_vanish_no_resonant_Lemma_}
Let $\gamma:\; \C^n \arrow \C^n$ be a germ of a holomorphic diffeomorphism in 0.
Let $B= T\C^n\otimes (T^* \C^n)^{\otimes k}$, $k >1$.
Diffeomorphisms of $\C^n$ can be naturally extended to $T\C^n$ and to its
tensor powers. Let $B_0$ denote the space of germs of sections of $B$ in 0;
clearly, $\gamma$ induces a natural automorphism of $B_0$. 
Assume that the differential $A:=D_0\gamma$ of $\gamma$ in 0 is non-resonant. 
Then any $\gamma$-invariant germ $v\in B_0$ vanishes.

\hfill

\proof
Let $t_1, ..., t_n$ be coordinates on $\C^n$.
Consider the Taylor series decomposition for 
$v$:
\[
v = \sum_i v_i  P_i(t_1, ..., t_n)
\]
where $P_i$ is a homogeneous polynomial of degree $i$
and $v_i \in V \otimes (V^*)^{\otimes k}$, where $V= T_0 \C^n$.
Let $l$ be the smallest integer such that $P_l \neq 0$.
The corresponding Taylor term can be considered to be 
an element of $V \otimes (V^*)^{\otimes k}\otimes \Sym^l(V^*)$.
By the chain rule, this tensor is coordinate-independent, and hence 
it is also $A$-invariant. Let $\alpha_1, ..., \alpha_n$
be the eigenvalues of $A$ on $T_0 \C^n$. Then 
$\alpha_1^{-1}, ..., \alpha_n^{-1}$ are the 
eigenvalues of $A$ on $T_0^*\C^n$. This is clear,
because $A$ preserves the pairing between $T\C^n$ and 
$T^*\C^n$: the differential $D\gamma$ acts covariantly
on vector field and contravariantly on 1-forms.

We obtain that the eigenvalues
of $A$ on $V \otimes (V^*)^{\otimes k}\otimes \Sym^l(V^*)$ are products
of the form $\alpha_l\prod_{i=1}^k \alpha_{j_i}^{-1}\prod_{i=1}^l \alpha_{s_i}^{-1}$.
Unless $A$ acts on $V \otimes (V^*)^{\otimes k}\otimes \Sym^l(V^*)$ 
with an eigenvalue 1, there would be no $A$-invariant vectors. Therefore
$v=0$ unless $\alpha_l = \prod_{i=1}^n \alpha_i^{k_i}$, where
$\sum_i k_i \geq 2$, for some eigenvalue $\alpha_l$.
\endproof

\subsection{Resonant equivariant bundles}

\definition
Let $\gamma$ be an invertible holomorphic contraction on $\C^n$,
centred in 0, and $B$ a $\gamma$-equivariant holomorphic
vector bundle on $\C^n$. Let $\alpha_1, ..., \alpha_n$
be the eigenvalues of the differential $A=D_0\gamma$ of
$\gamma$ in 0, and $\beta_1, ..., \beta_m$ the 
eigenvalues of the differential of the equivariant action of $\gamma$ 
on the fibre $B\restrict 0$. We say that $B$ is {\bf non-resonant}
if there is no multiplicative relation of the form $\beta_i = \beta_j \prod_{l=1}^k\alpha_{i_l}$
for some integer $k \geq 1$, where $\{\alpha_{i_1}, ...,\alpha_{i_k}\}$ is any collection of 
$k$ eigenvalues, possibly repeating, and $\beta_i, \beta_j\in \{\beta_1, ..., \beta_m\}$
any two eigenvalues, possibly the same. \index[terms]{resonance!of a bundle}

\hfill

\remark
A linear map with eigenvalues $\alpha_1, ..., \alpha_n$ is
resonant if one of the eigenvalues is a product of two or more
eigenvalues. The data associated with a $\gamma$-equivariant
bundle consists in two linear operators, the differential $D\gamma\restrict {T_0\C^n}$
with the eigenvalues $\alpha_i$
and the differential  of the equivariant action of $\gamma$ 
on the fibre $B\restrict 0$, with the eigenvalues $\beta_j$.
It has resonance when one of $\beta_j$ is a product of another
and one or more $\alpha_i$'s. The resonance in the bundle
$T\C^n$ with the natural $\gamma$-equivariant structure
is the same as the resonance in the linear operator
$D\gamma\restrict {T_0\C^n}$.

In other words,
if $B$ is $T\C^n$ with the standard $\gamma$-equivariant structure,
the relation $\beta_i = \beta_j \prod_{l=1}^k\alpha_{i_l}$
becomes $\alpha_i = \alpha_j \prod_{l=1}^k\alpha_{i_l}$, $k \geq 1$,
that is, $B=T\C^n$ is non-resonant if and only if
the differential $A=D_0\gamma$ is non-resonant.

\hfill

Resonant automorphisms can be characterized in terms 
of invariant 1-forms with coefficients in 
endomorphisms of $B$, as follows.

\hfill

\theorem\label{_resonance_in_equiv_bundles_Theorem_}\index[terms]{contraction!holomorphic}
Let $\gamma:\; \C^n \arrow \C^n$ be an invertible holomorphic
contraction centred in 0, and $B$ a $\gamma$-equivariant vector bundle.
Then $B$ is resonant if and only if there exists a non-zero
$\gamma$-invariant section of the bundle $\Omega^1\C^n\otimes \End(B)$.

\hfill

\pstep
Let $R$ be a $\gamma$-invariant section of $\Omega^1\C^n\otimes \End(B)$. 
We are going to prove that $B$ is resonant.
By \cite[Theorem 5.3.1]{_Forstneric:Oka_book_}, $B$ is trivial.
Choose a basis $b_1,..., b_m$ of $B$. Let $b_{ij}\in \End(B)$ be the
corresponding elementary matrices, defining a basis in $\End(B)$.
We write $R= \sum_{i,j,l} f_{ijl} dz_l \otimes b_{ij}$, where $z_i$ are coordinates in $\C^n$, and 
$f_{ijl}\in \calo_{\C^n}$ a function.
We write each $f_{ijl}$ as Taylor series, $f_{ijl} = \sum_s P_s^{ijl}$, 
where $P_s^{ijl}(z_1, ..., z_n)$ is 
a homogeneous polynomial of coordinate functions $z_i$
of degree $s$. Let $d$ be the smallest number
for which not all $P_d^{ijl}$ vanish. Using the chain rule again, we obtain that 
$\sum_{i,j} P_d ^{ijl} dz_l \otimes b_{ij}$ is $D_0\gamma$-invariant. 
This is a $D_0\gamma$-invariant vector in the space
$W^* \otimes \End(B|_0) \otimes \Sym^d(W^*)$, where $W=T _0\C^n$, with the
 action of $D_0\gamma$ which comes from the tensor product.

The eigenvalues of $D_0\gamma$ on the space
$W^* \otimes \End(B|_0) \otimes \Sym^d(W^*)$ are eigenvalues of
$D_0\gamma$ on $W^*$ times the eigenvalues on $\End(B|_0)$
times a product of $d$ eigenvalues of
$D_0\gamma$ on $W^*$, that is, numbers of form 
$\beta_u \beta_v^{-1} \alpha_l^{-1}$. 

Since
 $\sum_{i,j} P_d ^{ijl} dz_l \otimes b_{ij}\in W^* \otimes \End(B|_0) \otimes \Sym^d(W^*)$ 
is $D_0\gamma$-invariant, one of these numbers is  1. 
This gives a relation
 $\beta_u  = \beta_v\alpha_l\prod_{l=1}^d\alpha_{i_l}$.

\hfill

{\bf Step 2:} Let $B$ be a resonant bundle;
we are going to produce an invariant section of 
$\Omega^1\C^n\otimes \End(B)$. 
We follow the standard scheme
(\ref{_contra_compact_Theorem_},
\ref{_Stein_by_contract_to_linear_Hopf_Theorem_}),
using the Riesz--Schauder theorem\index[terms]{theorem!Riesz--Schauder}
and compactness of the action of holomorphic contractions on
holomorphic functions; this time, the contraction
acts on the sections of an equivariant bundle.
Let $M_0\subset \C^n$
be a subset that satisfies $\gamma(M_0) \Subset M_0$.
We equip $B$ with a Hermitian metric, and notice that
the space $H^0_b(M_0, \Omega^1M_0\otimes \End(B))$ with sup-norm
is Banach, by Montel theorem. Denote by \index[terms]{theorem!Montel}
$V$ the bundle $\Omega^1M_0\otimes \End(B)$.
By the standard argument (\ref{_contra_compact_Theorem_}), the operator
$\gamma^*:\; H^0_b(M_0, V)\arrow H^0_b(M_0, V)$ is compact.

Consider the filtration on $H^0_b(M_0, V)$
by the powers of the maximal ideal ${\goth m}_0$ of zero,
\[
H^0_b(M_0, V) \supset
H^0_b(M_0, {\goth m}_0V)\supset 
H^0_b(M_0, {\goth m}_0^2V)\supset ...
\]
The finite-dimensional space 
$\frac{H^0_b(M_0, V)}{H^0_b(M_0,  {\goth m}^k_0V)}$ is interpreted as
the space of $(k-1)$-jets of the sections of $V$ in 0.

Using the integral Cauchy formula,  any 
derivative of a function in a point 
can be expressed through a certain integral
of this function.
Therefore, the restriction map 
\[ H^0_b(M_0, V)\arrow
\frac{H^0_b(M_0, V)}{H^0_b(M_0,
  {\goth m}^k_0V)}=  
V\restrict 0,
\]
taking a section of $V$ to its $(k-1)$-jet, is also
continuous in the sup-norm. This implies also 
that the subspaces $H^0_b(M_0, {\goth m}_0^kV)\subset H^0_b(M_0, V)$
of sections with vanishing $(k-1)$-jet are closed.

The map $\gamma^*:\; H^0_b(M_0, V)\arrow H^0_b(M_0, V)$
preserves this filtration. Therefore, the eigenvalues of 
$\gamma^*$ on $H^0_b(M_0, V)$ and on the associated graded
space $\bigoplus_i \frac{H^0_b(M_0, {\goth m}^{i-1}V)}{H^0_b(M_0,
  {\goth m}^i_0V)}$ are equal.
However, the space 
$\frac{H^0_b(M_0, {\goth m}^{d-1}V)}{H^0_b(M_0,
  {\goth m}^d_0V)}$ is naturally identified with
$W^* \otimes_\C \Sym^d W^* \otimes_\C \End(B|_0)$,
where $W= T_0 \C^n$ (Step 1). In Step 1, we identified
this space with the space of $(d-1)$-th Taylor coefficients of
a section of $V$.

The eigenvalues of
$\gamma^*$ on this space are 
$\prod_{i=1}^{d+1}\alpha_i^{-1} \beta_1 \beta_2^{-1}$,
where $\beta_1, \beta_2$  denotes some eigenvalues of
the equivariant action of $\gamma^*$ on $B|_0$, possibly equal,  and
$\alpha_1, ..., \alpha_{d+1}$ is a collection of some
eigenvalues of $D_0\gamma$ on $W$, also possibly equal.
The existence of resonance on $B$ implies that
$\prod_{i=1}^{d+1} \alpha_i^{-1} \beta_1 \beta_2^{-1}=1$
for an appropriate choice of the eigenvalues
$\alpha_1, ..., \alpha_{d+1}, \beta_1, \beta_2$,
and this is equivalent to $\gamma^*$ having eigenvalue 1
on some $\gamma^*$-invariant quotient of $H^0_B(M_0, V)$, 
and hence on $H^0_B(M_0, V)$ too.
\endproof

\hfill

We can rewrite this theorem as a result about \index[persons]{Mall, D.} Mall bundles.

\hfill

\corollary\label{_B_on_Hopf_resonant_via_sections_Corollary_}
Let $B_H$ be a Mall bundle on a Hopf manifold
$H= \frac{\C^n \backslash  0}{\langle \gamma \rangle}$,
and $B:=j_* \pi^* B_H$ the corresponding
$\gamma$-equivariant bundle on $\C^n$.
Then $B$ is resonant if and only if
$H^0(H, \Omega^1 H\otimes_{\calo_H}(\End B_H))\neq 0$.
\endproof

\hfill

When $B= \Omega^1 \C^n$, the equivariant bundle $B$
is resonant if and only if the action of $D_0\gamma^*$
on $T_0 \C^n$ is resonant:

\hfill

\corollary\label{_resonant_matrix_tangent_bundle_Corollary_}
Let $\gamma:\; \C^n \arrow \C^n$ be an invertible
contraction centred in zero, $W:= T_0\C^n$, and $A=D_0\gamma\in \End W$ 
its differential in zero. Then the following are equivalent.
\begin{description}
\item[(i)] The matrix $A$ is resonant.
\item[(ii)] There exists a non-zero $\gamma^*$-invariant section of
$\Omega^1(\C^n) \otimes \End(T\C^n)$.
\item[(iii)] The bundle $T\C^n$ with the natural 
equivariant structure induced by the action of $\gamma$ is
resonant.
\end{description}

\proof
From the definition it is clear that $A$ is resonant if
and only if the $\gamma$-equivariant bundle $B=\Omega^1 \C^n$
is resonant. By \ref{_resonance_in_equiv_bundles_Theorem_},
this is equivalent to the existence of non-zero $\gamma^*$-invariant sections of
$\Omega^1(\C^n) \otimes \End(T\C^n)$.
\endproof

\hfill

Consider an invertible
contraction $\gamma:\; \C^n \arrow \C^n$ centred in zero,
and let $B$ be a $\gamma$-equivariant vector bundle.
Let $R= (\Omega^1 \C^n)^{\otimes k}$ be the bundle of
$k$-multilinear forms on $\C^n$. We consider
$R$ as a $\gamma$-equivariant vector bundle as well.
Then the set of eigenvalues of the $\gamma^*$-action on
$H^0_b(M_0, R \otimes \End B)$ is equal to 
\[ 
 \beta_1\beta_2^{-1}\left\{\prod_{i=1}^{d+k} \alpha_i^{-1}\right\}.
\]
Here $d\in \Z^{\geq 0}$, $\alpha_1,..., \alpha_{d+k}$
is any collection of eigenvalues 
of $D_0\gamma^*$ on $W=T_0\C^n$, and $\beta_1, \beta_2$
some eigenvalues of $\gamma^*$ on $B|_0$.
This is proven by the same argument as proves
\ref{_resonance_in_equiv_bundles_Theorem_}, Step 2.
We obtained the following corollary.

\hfill

\corollary\label{_sections_tensor_w_coeff_in_End_B_Corollary_}
Let $\gamma:\; \C^n \arrow \C^n$ be an invertible
contraction of $\C^n$ centred in zero,
and $B$ a $\gamma$-equivariant vector bundle.
Assume that $B$ is non-resonant. Consider
$R= (\Omega^1 \C^n)^{\otimes k}$, $k \geq 1$, 
as a $\gamma$-equivariant bundle.
Then the space of $\gamma^*$-invariant sections
of $R\otimes_{\calo_{\C^n}} B$ is empty.
\endproof


\subsection{Holomorphic connections on vector bundles}


Further on, we need the notion of 
a holomorphic connection.

\hfill

\definition\label{_holo_conne_Definition_}
Let $B$ be a holomorphic vector bundle on 
a complex manifold, and $\nabla:\; B \arrow B\otimes_{\calo_M} \Omega^1 M$
a differential operator that satisfies $\nabla(fb) = df \otimes b + f \nabla b$
for any locally defined holomorphic function $f$ and any local section $b$ of $B$.
Then $\nabla$ is called {\bf a holomorphic connection}.\index[terms]{connection!holomorphic}

\hfill

\index[persons]{Atiyah, M. F.}
The notion of a holomorphic connection was introduced by
M. Atiyah in \cite{_Atiyah:hol_bundles_}; for more results and
references, see
\cite{_Biswas_,_Biswas_Dumitrescu:holomorphic_affine_,%
_Dumitrescu_Biswas:holomorphic_Riemannian_}.

\hfill

Every flat connection is holomorphic with respect
to the holomorphic structure induced by this connection,
but there are more holomorphic connections than there
are flat connections.  Indeed, holomorphic connections
can be realized as objects of differential geometry,
as follows.\index[terms]{connection!flat}

\hfill

Let $(B, \nabla)$ be a complex vector bundle with connection on a 
complex manifold, and $\bar\6= \nabla^{0,1}$ the
corresponding $\bar\6$-operator.
\index[terms]{structure!holomorphic}
By Koszul-\index[persons]{Malgrange, B.}Malgrange theorem (\ref{_Koszul--Malgrange_Theorem_}),
$\bar\6$ defines a holomorphic structure on $B$ if and only
if $\bar\6^2=0$, or, equivalently, when \index[terms]{theorem!Koszul--Malgrange}
$\nabla^2 \in [\Lambda^{2,0}(M)\oplus \Lambda^{1,1}(M)]\otimes \End(B)$.

\hfill

\proposition\label{_holo_conne_curvature_Proposition_}
Let $(B, \nabla)$ be a complex vector bundle with connection on a 
complex manifold, and $\bar\6= \nabla^{0,1}$ the corresponding $\bar\6$-operator.
Assume that $\bar\6^2=0$, and let ${\cal B}=\ker \bar\6$ be 
the holomorphic vector bundle obtained from $\bar\6$ 
using \ref{_Koszul--Malgrange_Theorem_}. Then the
following assertions are equivalent.
\begin{description}
\item[(i)] The operator $\nabla^{1,0}$ is a holomorphic
connection operator on ${\cal B}$.
\item[(ii)] $\nabla^2 \in \Lambda^{2,0}(M)\otimes \End(B)$.
\end{description}
\proof
If $\nabla^{1,0}$ is a holomorphic
connection operator, it commutes with $\bar\6$, and hence 
the $(1,1)$-part of the curvature of $\nabla$ vanishes.
Conversely, if $\nabla^2 \in \Lambda^{2,0}(M)\otimes \End(B)$,
this implies that $\nabla^{1,0}$ commutes with $\bar\6$,
hence $\nabla^{1,0}$ maps holomorphic sections of $B$
to the holomorphic sections of $B \otimes \Lambda^{1,0}(M)$.
This implies that $\nabla^{1,0}$ is a holomorphic connection.
\endproof

\hfill

\remark
A holomorphic vector bundle $B$\index[terms]{class!Chern}
has Chern classes represented by closed forms
of type $(p,p)$, because it admits  Chern connections\index[terms]{connection!Chern}. However, 
the Chern classes of a bundle equipped with a holomorphic
connection are represented by holomorphic forms of type $(2p,0)$.
On a K\"ahler manifold (or any other manifold admitting
the $(p,q)$-decomposition in cohomology), this is impossible,
unless $c_i(B)=0$ for all $i$. This is why holomorphic
connections rarely occur in K\"ahler geometry.\index[terms]{geometry!K\"ahler} However,
on non-K\"ahler manifolds they do not seem that rare.

\hfill

The Picard group $\Pic(M)$ of a complex manifold  $M$ (that is, the group
of line bundles, with the group operation defined by the
tensor multiplication) is naturally identified with  
$H^1(M, \calo_M^*)$. Similarly, one could identify
the group of 1-dimensional local systems with $H^1(M, \C_M^*)$.
Here $\C^*_M$
denotes the constant sheaf with the space of sections
$H^0(U, \C^*_M)=\C^*$ for each connected
open set $U\subset M$.  The natural map
$H^1(M, \C_M^*)\arrow H^1(M, \calo_M^*)$
can be interpreted as a forgetful map,
taking a flat line bundle $(L, \nabla)$
to the holomorphic line bundle $(L, \nabla^{0,1})$.

\hfill

For K\"ahler manifolds, the map 
$H^1(M, \C_M^*)\arrow H^1(M, \calo_M^*)$ is never an isomorphism,
unless $b_1(M)=0$. Indeed, from the exponential exact sequence
it follows that $\dim_\C H^1(M, \C_M^*)= b_1(M)$ and
$\dim_\C H^1(M, \calo_M^*)= \frac 1 2 b_1(M)$.

However, on a Hopf manifold, these two groups are isomorphic,
that is, every holomorphic line bundle admits a unique flat connection.

\hfill

\proposition\label{_Picard_Hopf_Proposition_}
(\cite[page 57]{_Kodaira_Structure_III_}) \\ Let $H$ be a Hopf manifold,
and $H^1(H, \C_H^*)$ the cohomology with coefficients in
the constant sheaf $\C_H^*$.
Then the natural map $H^1(H, \C_H^*)\arrow H^1(H, \calo_H^*)$
to the Picard group is an isomorphism.\index[terms]{group!Picard}

\hfill

\pstep
Consider the exponential exact sequences
\[
0 \arrow \Z_H \arrow \calo_H \arrow \calo_H^* \arrow 0,
\]
and 
\[
0 \arrow \Z_H \arrow \C_H \arrow \C_H^* \arrow 0.
\]
The corresponding long exact sequences of cohomology give
\[
0 \arrow H^1(\Z_H) \arrow H^1(\calo_H)\arrow H^1(\calo_H^*)\arrow H^2(\Z_H)=0,
\]
and 
\[
0 \arrow H^1(\Z_H) \arrow H^1(\C_H)\arrow H^1(\C_H^*)\arrow H^2(\Z_H)=0.
\]
It remains to show that the natural map 
$\nu:\; H^1(\C_H)\arrow H^1(\calo_H)$
is an isomorphism; however, both groups are equal to $\C$
(\ref{_Dolbeault_for_Hopf_Theorem_}), and hence  it would
suffice to show that $\nu \neq 0$.
The relevant $E_2$-term of the Dolbeault spectral sequence
is $(H^{0,1}(H)=\C; H^{1,0}(H)=0)$\index[terms]{spectral sequence!Dolbeault}
(\ref{_Dolbeault_for_Hopf_Theorem_}, \ref{_holo_forms_Hopf_vanish_Proposition_}).
The sum of dimensions of these spaces is 1. However, $b_1(H)=1$,
hence the higher differentials of this spectral sequence  vanish
on $E_2^{1,0}+ E_2^{0,1}$ and it degenerates in the $E_2^{1,0}+ E_2^{0,1}$-term.
This implies that the first de Rham cohomology\index[terms]{cohomology!de Rham} $H^1(H)$
is equal to the first Dolbeault cohomology\index[terms]{cohomology!Dolbeault}
$H^{0,1}(H)\oplus H^{1,0}(H)$.

\hfill

{\bf Step 2:}
The standard map
$\nu:\; H^1(\C_H)\arrow H^1(\calo_H)$ is
the natural map from the $E_\infty^{0,1}$-term of this
spectral sequence to $E_2^{0,1}$, that is  an
isomorphism because the spectral sequence degenerates.
\endproof


\subsection{The flat connection on a non-resonant Mall bundle}


\definition
Let $B_H$ be a Mall bundle on a Hopf manifold
$H= \frac{\C^n \backslash  0}{\langle \gamma \rangle}$,
and $B:=j_* \pi^* B_H$ the corresponding
$\gamma$-equivariant bundle on $\C^n$.
We call $B_H$ {\bf resonant}
if $H^0(H, \Omega^1_H\otimes \End B_H)\neq 0$,
or, equivalently, when $B=j_* \pi^* B_H$
is a resonant $\gamma$-equivariant bundle on $\C^n$
(\ref{_B_on_Hopf_resonant_via_sections_Corollary_}).

\hfill

\remark\label{_conne_unique_Remark_}
Let $B_H$ be a non-resonant Mall bundle on a Hopf manifold.
Since the difference of two connections is
a holomorphic 1-form with coefficients in endomorphisms,
and $H^0(H, \Omega^1_H\otimes \End B_H)=0$,
a holomorphic connection on $B_H$ is unique, if it exists.

\hfill

\theorem\label{_Mall_flat_connection_Theorem_}
Let $B$ be a holomorphic vector bundle over a Hopf manifold $H$,
$\dim_\C H \geq 3$. \index[terms]{bundle!vector bundle!Mall} 
Assume that $B$ admits a flat connection compatible
with the holomorphic structure. Then $B$ is Mall.
Conversely, any non-resonant \index[persons]{Mall, D.} Mall bundle on $H$
admits a flat connection $\nabla$.\footnote{All flat connections
are holomorphic (\ref{_holo_conne_curvature_Proposition_}). 
By \ref{_conne_unique_Remark_}, $B$
admits a unique holomorphic connection, and hence  $\nabla$ is unique.}

\hfill

\pstep 
Let $(B, \nabla)$ be a holomorphic bundle with a flat connection on $H$,
and $\pi:\; \C^n \backslash 0 \arrow H$ the
universal cover. Since $\pi_1(\C^n \backslash 0)=0$,
the flat bundle $\pi^* B$ is trivial, and hence  $B$ is Mall.

\hfill

{\bf Step 2:} 
Let $B$ be a non-resonant Mall bundle on a Hopf manifold $H$.
We are going to show that $B$ admits a holomorphic connection.
Locally, a holomorphic connection always exists. The
space of holomorphic connections on $B\restrict U$
is an affine space, with linearization $H^0(U, \Omega^1 U
\otimes_{\calo_U} \End (B))$; indeed, the difference of two
holomorphic connections is 1-form with coefficients in 
endomorphisms. Choose an open cover $\{U_i\}$ of $H$
such that all restrictions $B\restrict{U_i}$ are trivial,
and let $\nabla_i$ denote a connection in $B\restrict{U_i}$.
The connections $\{\nabla_i\}$, generally speaking, 
do not agree on the intersections $U_{ij}:=U_i \cap U_j$,
giving a section $\alpha_{ij}\in (U_{ij}, \Omega^1 U_{ij}
\otimes_{\calo_{U_{ij}}} \End (B))$. Clearly, $\alpha_{ij}$
satisfy the cocycle condition, giving an element 
in $H^1(H, \Omega^1 H\otimes_{\calo_H} \End (B))$. This cocycle
is \v Cech cohomologous to zero if and only if one
can find a connection in each $B\restrict {U_i}$
which agrees on the pairwise intersections.

This implies that the obstruction to the global existence 
of a holomorphic connection belongs to 
$H^1(H, \Omega^1 H\otimes_{\calo_H} \End (B))$. 
Since $B$ is Mall, we have
\[ \rk H^1(H, \Omega^1 H\otimes_{\calo_H} \End (B))= 
\rk H^0(H, \Omega^1 H\otimes_{\calo_H} \End (B))
\]
(\ref{_Mall_cohomology_Theorem_}). 
When $B$ has no resonance, 
$\rk H^0(H, \Omega^1 H\otimes_{\calo_H} \End (B))=0$
hence $B$ admits a connection.

\hfill

{\bf Step 3:} To finish
\ref{_Mall_flat_connection_Theorem_},
it remains to show that any holomorphic connection
$\nabla$ on a non-resonant Mall bundle $B$ is flat. 
However, its curvature is a holomorphic 2-form, 
and the space of holomorphic 2-forms with
coefficients in $\End(B)$ vanishes by 
\ref{_sections_tensor_w_coeff_in_End_B_Corollary_}.
\endproof


\section{Flat connections on Hopf manifolds}
\label{_flat_Hopf_Section_}

\subsection{Developing map for flat affine manifolds}

For an introduction to flat affine manifolds, see\index[terms]{manifold!flat affine}
\cite{_Abels:survey_}, \cite{_Shima:book_}, \cite{_Goldman:structures_book_} and the references therein.
Recall that {\bf an affine function} on a vector space is a linear function
plus constant.\index[terms]{map!developing}\index[terms]{function!affine}

\hfill

\definition\label{_flat_affine_Definition_}
Let $M$ be a manifold equipped with a sheaf ${\cal F}\subset \C^\infty M$.
We say that ${\cal F}$ {\bf defines a flat affine structure on $M$}
if for each $x\in M$ there exists a neighbourhood diffeomorphic
to the ball $B \subset \R^n$ such that in this coordinate patch, the sheaf ${\cal F\restrict{B}}$ is the sheaf of affine functions on $B$. In other words, $M$ is a flat
affine manifold if there exists an open cover $\{ U_i\}$
with all $U_i$ diffeomorphic to an open ball in $\R^n$ and
all transition maps are affine. 
Such a cover is called {\bf an affine atlas} of $M$.
The sheaf ${\cal F}$ is called {\bf the sheaf of affine functions}
on the flat affine manifold $M$.\index[terms]{structure!affine}\index[terms]{structure!flat affine}

\hfill

Flat affine structures can be equivalently described in terms
of torsion-free, flat connections.\index[terms]{connection!torsion-free}\index[terms]{connection!flat}

\hfill

\proposition\label{_flat_affine_via_connection_Proposition_}
Let $(M, {\cal F})$ be a flat affine manifold. Then
$M$ admits a unique torsion-free, flat connection $\nabla$ such that
the sections $f$ of ${\cal F}$ satisfy $\nabla( df)=0$. Conversely,
if $\nabla$ is a torsion-free flat connection on $M$,
then the sheaf $\{ f \in C^\infty M\ \ |\ \ \nabla(df)=0\}$
defines a flat affine structure.

\hfill

\proof
Let $(M, {\cal F})$ be a flat affine manifold, and
$\{ U_i\}$ its affine atlas. Each $U_i$ admits 
a connection $\nabla:\; TM \arrow TM \otimes \Lambda^1 M$ with 
\begin{equation}\label{_standard_flat_connection_Equation_} 
\nabla\left(\sum_j f_j \frac{\6}{\6x_j}\right) = \sum_j df_j \otimes \frac{\6}{\6x_j},
\end{equation}
where $x_j$ are coordinate functions; this connection is clearly
flat and torsion-free,\index[terms]{connection!torsion-free}\index[terms]{connection!flat} and $\ker (\nabla d)$ is the sheaf of affine
functions on $U_i$. We call \eqref{_standard_flat_connection_Equation_}
{\bf the standard flat connection} on $U\subset \R^n$.

Conversely, let $\nabla$ be a torsion-free affine connection on $M$.
Locally in an open subset $U\subset M$, 
the bundle $\Lambda^1M$ admits a basis $\theta_1, ..., \theta_n$ 
of $\nabla$-parallel sections,\index[terms]{section!parallel}
that are  closed because $\nabla$ is torsion-free (here we use
the relation $d\theta = \Alt(\nabla \theta)$, which holds
for any torsion-free connection).  This implies that
$\theta_i = dx_i$ whenever $U$ is simply connected.
The functions $x_i$ give a coordinate system on $U$,
because $dx_i$ are linearly independent, and in this
coordinate system the affine functions are those
that satisfy $\nabla (df)=0$.
\endproof

\hfill

Further on, we will also call a torsion-free flat connection\index[terms]{connection!torsion-free}\index[terms]{connection!flat}
{\bf an affine structure} and a pair $(M, \nabla)$ 
{\bf an affine manifold}, or {\bf a flat affine manifold}.
\index[terms]{structure!affine}\index[terms]{manifold!affine}

\hfill

\definition\label{_flat_affine_complete_Definition_}
Let $\nabla$ be a connection on $TM$, where $M$ is a smooth manifold. 
We say that $(M, \nabla)$ is {\bf complete} if for any
$x\in M$ and any $v\in T_x M$, there exists a solution
$\gamma_v:\; \R \arrow M$ of the geodesic equation 
$\nabla_{\gamma_v'(t)}\gamma_v'(t)=0$ for all $t\in ]-\infty, \infty[$. 
{\bf The exponential map} is the map $\exp_x:\; T_x M \arrow M$
taking $v \in T_x M$ to $\gamma_v(1)$. Clearly, $\exp_x$ is a
diffeomorphism in a neighbourhood $U$ of $0\in T_x M$.
{\bf The developing map} $\dev:\; \exp_x(U) \arrow U$
is the inverse of $\exp_x$; in general, it is
defined only in a neighbourhood of $x\in M$.

\hfill

The following classical theorem, due to \index[persons]{Auslander, L.} Auslander and
\index[persons]{Markus, L.} Markus, is well known.

\hfill

\theorem\label{_dev_for_complete_affine_Theorem_} 
(\cite{_Auslander_Markus:holonomy_}, \cite{_Goldman:structures_book_}) 
Let $(M, \nabla)$ be a simply connected, connected, and  flat affine manifold.
Then the developing map can be extended to an affine map
$\dev:\; M \arrow T_x M$, also called {\bf the developing map}.
If, in addition, $M$ is complete, the developing map is 
a diffeomorphism.

\hfill

\proof
Since $M$ is simply connected, the cotangent
bundle $\Lambda^1 M$ is trivialized by parallel (and hence, closed) 
sections $\theta_1, ..., \theta_n$. Since $H^1(M)=0$,
there exists 
functions  $x_1, ..., x_n$  on $M$ that satisfy $dx_i = \theta_i$. 
If we consider $x_i$ as coordinate functions, the geodesics
are given by affine maps $t \arrow \sum a_i t x_i + b_i$. Therefore,
the map $m \mapsto (x_1(m),..., x_n(m))$ defines an affine map
 $\dev:\; M \arrow \R^n=T_x M$, that is  inverse to the exponential map in a 
neighbourhood of $x$. Here we identify $T_x M$ with $\R^n$
by taking the map $v\arrow (\theta_1(v),..., \theta_n(v))$.

Suppose now that $M$ is complete.
By construction, for any $v\in T_x M$ we have
$\dev(\exp_x(v))=v$. To prove that $\dev$ is inverse to $\exp_x$,
it remains only to show that $\exp_x$ is surjective. 
Clearly, $\exp_x$ is an open map. To prove that
$\im \exp_x=M$, it remains to show that the image $\im\exp_x$ is closed.

By absurd, choose a point $y \in M \backslash \im \exp_x$ 
which lies in the closure of $\im \exp_x$,
and let $\{z_i\}$ be a sequence of points in $\im \exp_x$
converging to $y$. Since $\dev$ is continuous,
we have $\lim_i \dev(z_i) = \dev(y)$,
hence 
\[ \lim_{i\to \infty} \exp_x(\dev(z_i)) = \lim_{i\to \infty} z_i
=\exp_x(\dev(y)).
\]
We obtain $y\in \im \exp_x$, a contradiction.
\endproof

\subsection{Flat affine connections on a Hopf manifold}

The following result will be used further on in 
our proof of the classical  Poincar\'e linearization
theorem.

\hfill

\theorem\label{_flat_Hopf_is_linear_Theorem_}
Let $H= \C^n \backslash 0/\Z$ be a Hopf manifold.
Assume that $TH$ admits a torsion-free, flat connection\index[terms]{connection!torsion-free}\index[terms]{connection!flat}
$\nabla$ compatible with the holomorphic structure. Then
$H$ is isomorphic to a linear Hopf manifold \index[terms]{structure!holomorphic}
$\C^n \backslash 0/\langle A \rangle$, where $A\in \GL(n, \C)$ is a
linear endomorphism.

\hfill

\pstep
Let $\pi:\; \C^n \backslash 0\arrow H$ be the covering
map, and $j:\; \C^n \backslash 0\hookrightarrow \C^n$ 
the tautological open embedding map. Let $B$ be a tensor
bundle on $H$. By \ref{_Mall_exa_Proposition_}, $B$ is a Mall 
bundle, and hence   $j_* \pi^* B$ is a bundle on $\C^n$.
By \cite[Theorem 5.3.1]{_Forstneric:Oka_book_}, 
any vector bundle on $\C^n$ is trivial.
Let $\nabla_0$ be the trivial connection on $T\C^n$.
Then the connection form $\nabla_0 -\pi^* \nabla\in \Omega^1 \C^n
\otimes_{\calo_{\C^n}} \End (T\C^n)$ is a section of 
the trivial bundle $\Omega^1 \C^n\otimes_{\calo_{\C^n}} \End (T\C^n)$
over $\C^n \backslash 0$, and hence  it can be extended to a
holomorphic connection on $\C^n$. This connection
is also torsion-free, flat and holomorphic, by continuity
of the corresponding tensors. \index[terms]{connection!torsion-free}\index[terms]{connection!flat}\index[terms]{connection!holomorphic}

To avoid further confusion, we denote
$\C^n$ by $M$, and this connection by $\nabla_M$.
To finish the proof
of \ref{_flat_Hopf_is_linear_Theorem_}, it would
suffice to show that the flat affine manifold
$(M,\nabla_M)$ is isomorphic to $\C^n$ with the standard flat
connection. Since the $\Z$-action on $M$ preserves
the affine structure and fixes 0, it is defined by
a linear endomorphism, and $(M \backslash 0)/\Z$
is a linear Hopf manifold.

We conclude that it remains to show that the developing map
$\dev: M \arrow \C^n$ associated with $\nabla_M$ is an
isomorphism. This would follow from
\ref{_dev_for_complete_affine_Theorem_}
if we could prove that $\nabla_M$ is complete,
but we have no direct control over the connection
form of $\nabla_M$, and hence  the completeness is not obvious.

\hfill

{\bf Step 2:}
Let $\exp_0:\; U \arrow M$ 
be the exponential map, with $U\subset T_0M$
being a maximal neighbourhood of $0$ where $\exp_0$ is defined.

It remains to show that
$U=T_0M$; this would imply
that $(M, \nabla_M)$ is complete and the developing
map is an isomorphism. 

Denote by $A$ the generator
of the $\Z$-action on $M$ contracting $M$ to 0,
and let $A_0\in \End(T_0M)$ be its differential.
By construction, $\nabla_M$ is $A$-invariant. 
Therefore, $A$ preserves the geodesics of $(M, \nabla_M)$,
and we have a commutative diagram 
\[
\begin{CD}
U@>{\exp_0}>> M=\C^n \\
@V{A_0}VV  @VV{A} V \\
U @>{\exp_0}>> M= \C^n
\end{CD}
\]
This implies that $U\subset T_0 M$ is a $\Z$-invariant
neighbourhood of $0\in T_0 M$, where the action of $\Z=\langle A_0\rangle$ is
generated by $A_0$. Since 
$\bigcup_n A_0^{-n}(V)= T_0 M$ for any neighbourhood 
$V$ of 0, any $\langle A_0\rangle$-invariant neighbourhood
is equal to $T_0 M$. This implies that $(M, \nabla)$ is
complete, and the developing map is an isomorphism.
\endproof

\subsection[A new proof of Poincar\'e 
theorem about linearization of non-resonant contractions]{A new proof of Poincar\'e 
theorem about\\ linearization of non-resonant contractions}
\index[terms]{theorem!Poincar\'e linearization}\index[terms]{contraction!(non-)resonant}

The Poincar\'e-Dulac theorem (\cite{_Arnold:ODE+_})
gives a normal form of a 
smooth (or analytic) contraction; its non-resonant case
is sometimes called {\em the  Poincar\'e theorem}. It
proves that a contraction (or a germ of a contraction),
that is  non-resonant, becomes linear after an appropriate
coordinate change. We give a new\index[terms]{theorem!Poincar\'e-Dulac}
proof of this theorem based on complex geometry.
Note that the assumption $n \geq 3$ below is unnecessary;
we leave the case $n=2$ for the reader as an exercise.

\hfill

\theorem\label{_Poincare_new_proof_Theorem_}
Let $\gamma$ be invertible holomorphic contraction of $\C^n$ centred in 0, 
$n \geq 3$. Assume that the differential
$D_0\gamma\in \GL(T_0 \C^n)$ is non-resonant. Then there exists
a holomorphic diffeomorphism $U:\; \C^n \arrow \C^n$ such that
$U \gamma U^{-1}$ is linear.\footnote{The same statement can
be stated for a germ of a holomorphic diffeomorphism; the proof will be essentially
the same. We leave the required changes to the reader.}

\hfill

\proof
Let $H:= \C^n \backslash 0 /\langle \gamma \rangle$ be the Hopf
manifold associated with $\gamma$, and 
$\pi:\; \C^n \backslash 0\arrow H$ the universal covering map.
By \ref{_Mall_exa_Proposition_}, the tangent
bundle $TH$ is Mall. By \ref{_resonant_matrix_tangent_bundle_Corollary_}, 
it is non-resonant.
By \ref{_Mall_flat_connection_Theorem_},
$TH$ admits a flat 
holomorphic connection $\nabla$.  Since $\Omega^2 H \otimes TH$ is a direct sum
component of $\Omega^1 H \otimes \End(TH)$,
and $TH$ is non-resonant, 
\[ H^0(H,\Omega^2 H \otimes TH)
\subset H^0(H,\Omega^1 H \otimes \End(TH))=0,
\] 
(\ref{_resonant_matrix_tangent_bundle_Corollary_}), and hence  $\nabla$ is torsion-free.\index[terms]{connection!torsion-free}

By \ref{_flat_Hopf_is_linear_Theorem_},
the universal cover $\C^n \backslash 0$ of $H$ admits flat coordinates
such that $\gamma$ is linear in these coordinates. 
This proves \ref{_Poincare_new_proof_Theorem_}.
\endproof

\section[Harmonic forms on Hopf manifolds with coefficients in a bundle]{Harmonic forms on Hopf manifolds\\ with coefficients in a bundle}

In this section, we compute the multiplication of cohomology
of holomorphic bundles on Hopf manifolds in some
particular examples, and express the cohomology in terms
of the harmonic forms. This is used in Chapter \ref{_Kuranishi_Chapter_}
in our treatment of the Kuranishi space.\index[terms]{space!Kuranishi}

\subsection[The Hodge $*$ operator and cohomology of holomorphic bundles]{The Hodge $*$ operator and cohomology\\ of holomorphic bundles}

Recall that the cohomology of a holomorphic vector bundle\index[terms]{operator!Hodge $*$}\index[terms]{resolution!Dolbeault}
$B$ on a complex manifold $M$ can be computed using the Dolbeault resolution
\[
0 \arrow B \stackrel {\bar\6}\arrow B \otimes \Lambda^{0,1}(M) 
\stackrel {\bar\6}\arrow B \otimes \Lambda^{0,2}(M)
\stackrel {\bar\6} \arrow ...
\]
where $B$ and $B \otimes \Lambda^{0,i}(M)$ denote the space of
smooth sections of the bundles $B$ and $B \otimes \Lambda^{0,i}(M)$.
We equip $B$ and $M$ with a Hermitian metric.
Then $H^i(M, B)$ is identified with the kernel of the corresponding
Laplacian $\Delta_{\bar\6}:= \bar\6\bar\6^* + \bar\6^* \bar\6$.
The operator $\bar\6^*$ can be expressed using the 
appropriate Hodge $*$ operator, as follows.

The metric on $B$ and the Hermitian form\index[terms]{form!Hermitian} on $M$ can be
combined to give an Hermitian scalar product
$g$ on $B \otimes \Lambda^{p, q}(M)$. This
construction works for any tensor bundle, or 
the bundle of tensors
with coefficients in $B$: the scalar
product is extended to the tensor product
from its multipliers using
$(a\otimes b, a_1 \otimes b_1) = (a, a_1) (b, b_1)$.
We can consider $g$ as a map 
\begin{equation}\label{_Hermitian_product_on_forms_Equation_}
g:\; B \otimes \Lambda^{p, q}(M)\arrow \bar B^* \otimes \Lambda^{q, p}(M)^*
\end{equation}
(the complex conjugation is used to denote that $g$ is Hermitian,
hence complex anti-linear).

Another product can be defined using the exterior multiplication
of differential forms. Let $\Vol$ be the Hermitian volume form on $M$.
Denote by 
\[ 
u:\; B \otimes  B^* \otimes 
\Lambda^{p, q}(M) \arrow \Lambda^{p, q}(M)
\]
the map pairing a section of $B$ with a section
of $B^*$.
Let $\eta\in B \otimes \Lambda^{p, q}(M), 
\eta' \in B^* \Lambda^{n-p, n-q}(M)$, and
$h(\eta, \eta'):= \frac{u(\eta\wedge \eta')}{\Vol}$.
This map defines a non-degenerate $C^\infty M$-valued product
on differential forms, that we also interpret as an
isomorphism
\begin{equation}\label{_Poincare_product_on_forms_Equation_}
h:\; B \otimes \Lambda^{p, q}(M)\arrow B^* \otimes \Lambda^{n-p, n-q}(M)^*
\end{equation}
Combining the dualities from \eqref{_Hermitian_product_on_forms_Equation_}
and \eqref{_Poincare_product_on_forms_Equation_}, we obtain the Hodge $*$ operator
\[
*:= g h^{-1}:\; B \otimes \Lambda^{p, q}(M)\arrow
\bar B \otimes \Lambda^{n-q, n-p}(M).
\]
The following claim seems to be well known.

\hfill

\claim\label{_bar_6^*_via_*_Claim_}
The operator $\bar \6^*:\; B \otimes \Lambda^{0,p}(M) \arrow  
B \otimes \Lambda^{0,p-1}(M)$ can be written
as $\bar\6^*= - * \6 *$, where $*$ is the Hodge star operator,
and $\6:\; \bar B \otimes \Lambda^{p,q}(M)\arrow \bar B \otimes \Lambda^{p+1,q}(M)$
denotes the complex conjugate of $\bar\6$.

\proof
\cite[Chapter 1, Section 2, page 152]{griha}.
\endproof

\subsection[Multiplication in cohomology of holomorphic bundles on Vaisman-type Hopf manifolds]{Multiplication in cohomology of holomorphic\\ bundles on Vaisman-type Hopf manifolds}

Let $B$ be a holomorphic bundle on a Hopf manifold $H$,
satisfying the assumptions of \index[persons]{Mall, D.} Mall's theorem (\ref{_Mall_cohomology_Theorem_}).
The spaces $H^0(H, B)$ and $H^1(H, B)$ have the same
dimension, and hence  they are isomorphic, but this
isomorphism is somewhat mysterious. We were able to find its
geometric interpretation only for a special kind of Hopf manifolds,
namely for those that admit a Vaisman structure.

\hfill

\theorem\label{_multi_cohomo_of_B_Theorem_} 
Let $(H, \omega, \theta)$ be a Hopf manifold
equipped with an LCK structure,\index[terms]{structure!LCK} and $B$ a Mall bundle on $M$.
Let $[\theta^{0,1}]\in H^{0,1}(H) = H^1(H,\calo_H)$ denote
the Dolbeault class of the $(0,1)$-part of the Lee form.\index[terms]{form!Lee}
Consider the map $L_\theta:\; H^0(H, B) \arrow H^1(H, B)$
taking $\alpha \in H^0(H, B)$ to the image of 
$[\theta^{0,1}]\wedge \alpha$ under the multiplication
map $H^1(H, \calo_ H) \otimes H^0(H, B)\arrow H^1(H, B)$.
Assume that $H$ is Vaisman. Then $L_\theta$ is an
isomorphism.

\hfill

\pstep
Let $\alpha_0\in H^0(H,B)$ be a non-zero holomorphic section,
and $\alpha := \theta^{0,1}\otimes \alpha_0$ the corresponding
$\bar\6$-closed $B$-valued 1-form. Recall that  $\dim H^0(H, B)= \dim H^1(H, B)$. Therefore, in order 
to prove that the map $L_\theta:\; H^0(H, B) \arrow H^1(H, B)$
is an isomorphism, it would suffice to show that
$\alpha =\theta^{0,1}\otimes \alpha_0$ is not cohomologous to 0.
To prove this, it is enough to show that this form is harmonic.

{\bf Step 2:} 
Since $\bar\6(\alpha)=0$, harmonicity of $\alpha$ is
equivalent to $\bar\6^*(\alpha)=0$. By 
\ref{_bar_6^*_via_*_Claim_}, this is the same
as $\6(* \alpha)=0$. A trivial calculation
in coordinates implies that on every Vaisman manifold $(M, \omega, \theta)$,\index[terms]{manifold!Vaisman}
one has $* (\theta^{0,1})= -\1 (\omega_0)^{n-1}\wedge \theta^{0,1}$,
where $\omega_0$ is the transversal K\"ahler form\index[terms]{form!K\"ahler!transversal} (\ref{_Subva_Vaisman_Theorem_}). 
To see this, notice that $*(\theta) = I(\theta) \wedge \omega^{n-1}$
by \eqref{_*theta_Equation_}. This gives 
\begin{equation*}
	\begin{split}
*(\theta^{0,1})&= \frac 1 2 *(\theta- \1 I\theta)=
\frac 1 2  (I\theta\wedge\omega^{n-1} - \1 \theta \wedge\omega^{n-1})\\
&=
- \frac{\1} 2 (\theta \wedge\omega^{n-1} - \1  I\theta\wedge\omega^{n-1})=
- \1 \theta^{0,1} \wedge \omega^{n-1}.
\end{split}
\end{equation*}
However, $\omega= \omega_0 + \theta \wedge I\theta$ by 
\eqref{_omega_via_theta_Chapter_8_Equation_}, 
which implies $\theta^{0,1} \wedge \omega^{n-1}= \theta^{0,1} \wedge\omega_0^{n-1}$.
This gives $*(\alpha) = -\1\bar\alpha \otimes \theta^{0,1} \wedge \omega^{n-1}_0$.

{\bf Step 3:}
It remains to show that 
$\6(-\1\bar\alpha \otimes \theta^{0,1} \wedge \omega^{n-1}_0)=0$.
Since $\alpha$ is holomorphic, one has $\6(\bar\alpha)=0$; by Leibniz relation, 
it remains only to show that $\6(\theta^{0,1} \wedge \omega^{n-1}_0)=0$.
The form $\omega_0$ is closed (\ref{_Subva_Vaisman_Theorem_}),
and $d(\theta^{0,1}) = \frac {-\1}2\omega_0$, 
by \eqref{_omega_0_Equation_}. This gives
\[
\6(\theta^{0,1} \wedge \omega^{n-1}_0)= \frac {-\1}2\omega_0^n=0,
\]
since $\omega_0$ vanishes on the canonical foliation \index[terms]{foliation!canonical} $\Sigma$ 
(\ref{_Subva_Vaisman_Theorem_})
and $\dim_\C M=n$.
\endproof

\subsection[Appendix: cohomology of local systems on $S^1$
and the multiplication in cohomology of holomorphic bundles on Hopf manifolds]{Appendix: cohomology of local systems on $S^1$
  and\\ the multiplication in cohomology of holomorphic\\ bundles on Hopf manifolds}

Let $H$ be a Hopf manifold, and $B$ a holomorphic bundle on $H$.
We are going to express the multiplication map
$H^1(H, \calo_H) \otimes H^0(H, B)\arrow H^1(H, B)$
in explicit terms, in order to compute it 
for particular examples of Hopf manifolds.
Among other things, we show that the multiplication
map $L_{\theta}:\; H^0(H,B) \arrow  H^1(H,B)$ (\ref{_multi_cohomo_of_B_Theorem_})
is not always an isomorphism when $H$ is not Vaisman.

This subsection is more technical than the rest, and can
be safely skipped, because we are not going to use any of
the results proven here.

\hfill

\remark\label{_H^1_H^0,1_Hopf_Remark_}
Notice that $H^{1,0}(H)=0$ by \ref{_holo_forms_Hopf_vanish_Proposition_},
hence the natural map $H^1(H, \C)\stackrel \nu \arrow H^{0,1}(H)$ is
an isomorphism (\ref{_Picard_Hopf_Proposition_}, Step 2).

\hfill

Let $H$ be a Hopf manifold admitting an LCK potential \index[terms]{potential!LCK}
$\phi:\; \tilde H \arrow \R^{>0}$ on its $\Z$-cover
$\tilde H = \C^n \backslash 0$. Then the $\Z$-action
on $\tilde M$ takes $\log\phi$ to $\log\phi-\const$,
where $\const$ belongs to a discrete subgroup 
$\lambda\Z \subset \R$. Let
$\psi:\; H \arrow S^1$ be the map
taking $x\in H$ to $\log \phi(x) \mod \lambda\Z$.

\hfill

\lemma\label{_derived_pushforward_to_S^1_Lemma_}
In these assumptions, let $B$ be a Mall
bundle on $H$, $n=\dim H \geq 3$. 
Then the higher derived images $R^i\psi_*(B)$ 
vanish for all $0<i<n-1$, and $H^i(H, B)= H^i(S^1,
\psi_*(B))$ when $i < n$.

\hfill

\pstep
The Grothendieck spectral sequence\index[terms]{spectral sequence!Grothendieck} 
$E^{p,q}_2=H^p(R^q\psi_* B) \Rightarrow H^{p+q}(H,B)$
together with the vanishing of $R^i\psi_*(B)$ implies the isomorphisms
$H^i(H, B)= H^i(S^1,\psi_*(B))$ if $R^i\psi_*(B)=0$ when $0<i<n-1$.

\hfill

{\bf Step 2:} It remains to show that $R^i\psi_*(B)=0$ 
when $0<i<n-1$. As we have mentioned before,
the derived direct image $R^k\psi_*B$ can be
expressed as the sheafification of the presheaf
sending $W\subset S^1$ to $H^k(\psi^{-1}(W), B)$
(\cite[Theorem 3, Section III.8]{_Gelfand_Manin_}).
The cohomology of the annulus $\psi^{-1}(W)$ is
computed in \ref{_cohomo_of_C^n_withput_0_Proposition_},
Step 2. There we prove that $H^i(\psi^{-1}(W), B)=0$
when $0 <i < n-1$ and $B$ is a trivial bundle on the annulus.
However, the same proof works when $B$ is any Mall 
bundle, because $j_* \pi^* B$ is a trivial bundle
on $\C^n$, and hence  the restriction of $B$ to the annulus is trivial.
\endproof

\hfill

\proposition\label{_multiplica_by_theta_via_local_system_Proposition_}
Let $B$ be a Mall bundle on a Hopf manifold $H$ with LCK potential,\index[terms]{potential!LCK}
$[\theta^{0,1}]\in H^1(H, \calo_H)$ the generator
of the first Dolbeault cohomology group, and
$H^0(H, B) \xlongrightarrow{\wedge [\theta^{0,1}]} H^1(H, B)$
the multiplication map. Let $\psi:\; H\arrow S^1$
be the standard projection map (\ref{_derived_pushforward_to_S^1_Lemma_}).
Let $[\theta]\in H^1(H, \C)$ be the
corresponding generator of the cohomology,
and  $[\theta_0]\in H^1(S^1, \C)$ 
the generator that satisfies 
$\psi^*([\theta_0]) =[\theta]$. Denote by $\C_{S^1}$ the
constant sheaf on $S^1$.
Consider the tautological multiplication map
$\C_{S^1} \otimes \psi_* B \arrow \psi_*B$,
and let $H^0(\psi_* B) \xlongrightarrow{\wedge [\theta_0]} H^1(\psi_* B)$
be the corresponding multiplication in the cohomology.
Then the following diagram is commutative
(the vertical arrows are induced by the isomorphisms
$H^i(H,B)\stackrel \nu \arrow H^i(S^1,\psi_* B)$ for $i=0, 1$, 
\ref{_derived_pushforward_to_S^1_Lemma_})
\begin{equation}\label{_multi_commu_pushforward_Equation_}
\begin{CD}
H^0(H,B) @>{\wedge [\theta^{0,1}]}>> H^1(H,B)\\
@V\nu VV @VV\nu V\\
H^0(S^1,\psi_* B) @>{\wedge [\theta_0]}>> H^1(S^1,\psi_* B).
\end{CD}
\end{equation}
\proof
The tautological multiplication of sheaves
$\C_H \otimes B \arrow B$ induces 
multiplication in cohomology 
$H^1(H, \C_H)\times H^0(H, B) \arrow  H^0(H, B)$.
The natural embedding of $\C_H$ to $\calo_H$
takes the generator $[\theta]\in H^1(H, \C_H)= H^1(H, \C)$
to $[\theta^{0,1}]\in H^1(H, \calo_H)$
(\ref{_H^1_H^0,1_Hopf_Remark_}).
This makes the following diagram
commutative
\begin{equation}\label{_multi_commu_multi_Equation_}
\begin{CD}
H^0(H,B)\otimes H^1(H, \calo_H) @>>{\wedge [\theta^{0,1}]}> H^1(H,B)\\
@V\eta VV @VV{\Id}V\\
H^0(H,B)\otimes H^1(H, \C_H)  @>>{\wedge [\theta]}> H^1(H,B)
\end{CD}
\end{equation}
where $\eta$ is induced by the natural isomorphism
$H^1(H, \C_H)=H^1(H, \calo_H)$ (see \ref{_H^1_H^0,1_Hopf_Remark_}).
The multiplication 
$H^0(H,B)\otimes H^1(H, \C_H)  \stackrel{\wedge [\theta]}\arrow H^1(H,B)$
coincides with the multiplication
$H^0(S^1,\psi_* B)\otimes H^1(S^1, \C_{S^1})  \stackrel{\wedge  [\theta_0]}\arrow
  H^1(S^1,\psi_* B)$
under the isomorphisms 
$H^i(H,B)\stackrel \nu \arrow H^i(S^1,\phi_* B)$ 
given by \ref{_derived_pushforward_to_S^1_Lemma_}.
This identifies the diagrams
\eqref{_multi_commu_multi_Equation_}
and \eqref{_multi_commu_pushforward_Equation_}
and proves that the latter is commutative.
\endproof

\hfill

Let $H$ be a Hopf manifold and $B$ a Mall bundle on $H$,
and $\pi:\; \C^n \backslash 0 \arrow H$ the
standard $\Z$-covering. Then $\pi^* B$ is a
$\Z$-equivariant holomorphic bundle on $\C^n\backslash 0$
(\ref{_covering_equivariant_Theorem_}). Let $\gamma$ be
the generator of the $\Z$-action on $\C^n\backslash 0$.
Recall that a section $\alpha \in H^0(\C^n\backslash 0, \pi^* B)$
is called {\bf $\gamma$-finite} \index[terms]{function!$\gamma^*$-finite}if it is contained in 
a finitely-dimensional $\Z$-invariant subspace of 
$H^0(\C^n\backslash 0, \pi^* B)$. The following
result is a version of Mall's theorem
(\ref{_Mall_cohomology_Theorem_}), implicit in its proof.

\hfill

\theorem\label{_Mall_theorem_inva_coinva_Theorem_}
Let $B$ be a Mall bundle on a Hopf manifold $H$. As above, denote by
$W:=H^0(\C^n\backslash 0, \pi^* B)^\gamma$
the space of $\gamma$-finite sections, $n \geq 3$.
Then $H^0(H, B) = H^0(\Z, W)$ and  $H^1(H, B) = H^1(\Z, W)$,
where the $\Z$-action in $W$ is generated by $\gamma$.

\proof This is \eqref{_cohomo_coinva_Mall_Equation_}. \endproof

\hfill

Return now to the case when $H$ is equipped with an LCK
potential. Let ${\cal W}\subset \psi_* \pi^* B$ be the subsheaf
of all $\gamma$-finite sections\index[terms]{function!$\gamma^*$-finite} of $\psi_* \pi^* B$.
For any interval $I\subset S^1$ corresponding to an
annulus $\psi^{-1}(I) \subset \C^n\backslash 0$,
and any section of ${\cal W}$ on $I$,
the corresponding section of $H^0(\psi^{-1}(I), B)$
can be extended uniquely to the whole $\C^n\backslash 0$.
Then every local section of ${\cal W}$ can be uniquely 
extended to the universal cover $\R$ of $S^1$; this implies
that ${\cal W}$ is locally constant. By construction, the
$\gamma$-action on $H^0(\C^n\backslash 0, \pi^* B)^\gamma$
is equal to its monodromy. By \index[terms]{action!monodromy}
\ref{_Mall_theorem_inva_coinva_Theorem_}, 
$H^0(H, B) = H^0(S^1, {\cal W})$ and $H^1(H, B) = H^1(S^1,    {\cal W})$.
In \ref{_multi_cohomo_of_B_Theorem_}, 
we have shown that the multiplication
$L_\theta:\; H^0(H, B)\arrow H^1(H, B)$ is
an isomorphism when $B$ is a Mall bundle on
a Vaisman manifold. \index[terms]{manifold!Vaisman}

Using
\ref{_multiplica_by_theta_via_local_system_Proposition_}, 
we can compute this map
explicitly for any Hopf manifold
with an LCK potential.\index[terms]{potential!LCK}

\hfill

\proposition
Let $B$ be a Mall bundle over a Hopf manifold with LCK 
potential, and ${\cal W}$ the local system on $S^1$
obtained from the $\gamma$-finite sections\index[terms]{function!$\gamma^*$-finite} of $\pi^* B$
as above. Let $[\theta_0]\in H^1(S^1)$ be a generator
of $H^1(S^1)$. Consider the multiplication map
$H^0(S^1, {\cal W}) \stackrel {\wedge [\theta_0]}\arrow H^1(S^1, {\cal W})$
induced by the tautological multiplication $\C_{S^1}
\times {\cal W} \arrow {\cal W}$. 
Then, under the identifications
$H^0(H, B) = H^0(S^1, {\cal W})$ and $H^1(H, B) = H^1(S^1,{\cal W})$,
the multiplication
$H^0(S^1, {\cal W}) \stackrel {\wedge [\theta_0]}\arrow
H^1(S^1, {\cal W})$
corresponds to the multiplication
$L_\theta:\; H^0(H, B)\arrow H^1(H, B)$.

\proof Follows immediately from 
\ref{_multiplica_by_theta_via_local_system_Proposition_}.
\endproof

\hfill

The multiplication in cohomology of a local system on $S^1$ 
can be expressed explicitly as follows.

\hfill

\claim\label{_mult_on_coh_of_S^1_Claim_}
Let ${\Bbb V}$ be a finite-dimensional local system on $S^1$,
and $V$ the corresponding representation of $\Z$.
Consider the multiplication \[ L_{dt}:\; H^0(S^1, {\Bbb V})\arrow H^1(S^1, {\Bbb V})\]
taking $\alpha$ to $\alpha\wedge dt$,
where $dt$ is the unit 1-form interpreted as a generator of $H^1(S^1, \C_{S^1})$, 
and $\alpha$ a section of ${\Bbb V}$.\index[terms]{space!of invariants/coinvariants}
Then $H^0(S^1, {\Bbb V})= V^\Z$ is the space of $\Z$-invariants, 
$H^1(S^1, {\Bbb V})= V_\Z$  the space of $\Z$-coinvariants,
and $L_{dt}$ maps an invariant vector in $V$ to its image
in $V_\Z= \frac{V}{\sum_{\gamma\in\Z, l \in V} \gamma(l) -l}$.

\hfill

\proof
Let $\gamma$ be the generator of $\Z$-action.\index[terms]{action!$\Z$-} 
Partitioning $V$ onto a direct
sum of Jordan blocks, we can always assume that 
$V$ is just one Jordan block of the form
\[\small 
\begin{pmatrix}
\lambda & 1       & 0       & \cdots  & 0 \\
0       & \lambda & 1       & \cdots  & 0 \\
\vdots  & \vdots  & \vdots  & \ddots  & \vdots \\
0       & 0       & 0       & \lambda & 1       \\
0       & 0       & 0       & 0       & \lambda
\end{pmatrix}.
\]
If $\lambda\neq 1$, the invariants and the coinvariants are zero,
hence we can also assume that $\lambda=1$. 

\hfill

{\bf Step 2:} We are going to show that
the multiplication map $H^0(S^1, {\Bbb V})\stackrel {L_{dt}} \arrow H^1(S^1, {\Bbb V})$
is zero if the Jordan block is not diagonal, and an isomorphism otherwise;
this would easily imply the claim of \ref{_mult_on_coh_of_S^1_Claim_}.

We interpret the cohomology of ${\Bbb V}$
as the cohomology of the Morse--Novikov complex
\[ 0 \arrow C^\infty S^1 \otimes V \stackrel{d- \theta}\arrow 
\Lambda^1 S^1 \otimes V\arrow 0,
\]
where $\theta(l) = dt \otimes u(l)$ for all $l\in V$, with
$u=\log(\gamma)\in \End({\Bbb V} \otimes C^\infty S^1)$.
Since $V$ contains only one Jordan block, the spaces
$V^\Z= \ker(1-\gamma)$ and $V_\Z= \coker(1-\gamma)$
are 1-dimensional. Let $l_0 \in V$ be a $\gamma$-invariant
vector. The map $L_{dt}$ takes $l$ to $dt \otimes l$.
Since $d_\theta(l)=dl- dt \otimes u(l)$,
the ${\Bbb V}$-valued $dt \otimes l$ is $d_\theta$-exact if
$l\in \im \gamma$. This always happens unless 
$\gamma$ is diagonal; if $\gamma$ is diagonal,
it acts as identity, and in this case,
the multiplication $H^0(S^1) \otimes H^1(S^1) \arrow H^1(S^1)$
is an isomorphism because $H^0(S^1)$ is generated by the unity 
in the cohomology.
\endproof

\hfill

We can use \ref{_mult_on_coh_of_S^1_Claim_}
to show that \ref{_multi_cohomo_of_B_Theorem_} is false
for non-Vaisman Hopf manifolds.

\hfill

\proposition
Let $H=\C^n\backslash 0/\langle A \rangle$ be a linear Hopf manifold,
with the contraction matrix $A\in \End(\C^n)$ non-diagonalizable.
Then there exists a holomorphic line bundle $L$
on $H$ such that the multiplication map
\[ L_\theta:\; H^0(H, L) \arrow H^1(H, L)\]
(\ref{_multi_cohomo_of_B_Theorem_}) is not an isomorphism.

\hfill

\pstep
By \ref{_Picard_Hopf_Proposition_}, 
any holomorphic line bundle on $H$ admits a unique flat connection
compatible with the holomorphic structure.
Let $\chi:\; \Z \arrow \C^*$ be the 
character defining the corresponding local system.
Then the space $H^0(H, L)$ is identified with
the space of $\chi$-automorphic holomorphic
functions on $\C^n\backslash 0$ (\ref{_covering_equivariant_Theorem_}).
Let ${\Bbb V}$ be the local system on $S^1$
associated with the space of $A$-finite\index[terms]{section!$A$-finite}
sections $r\in H^0(\C^n\backslash 0, \pi^* L)$
as in \ref{_Mall_theorem_inva_coinva_Theorem_}.
By \ref{_Mall_theorem_inva_coinva_Theorem_},
$H^0(H, L)= H^0(S^1, {\Bbb V})$ and 
$H^1(H, L)= H^1(S^1, {\Bbb V})$.
Under these identifications, 
the multiplication $L_\theta:\; H^0(H, L) \arrow H^1(H, L)$
is equal to the multiplication 
$L_{dt}:\; H^0(S^1, {\Bbb V})\arrow H^1(S^1, {\Bbb V})$
(\ref{_multiplica_by_theta_via_local_system_Proposition_}).
By \ref{_mult_on_coh_of_S^1_Claim_},
to produce an example of a line bundle
where $L_\theta:\; H^0(H, L) \arrow H^1(H, L)$
is not an isomorphism, it would suffice to
find $L$ such that the Jordan blocks
of the local system ${\Bbb V}$ with the
eigenvalue 1 are non-diagonal.

\hfill

{\bf Step 2:}
Consider the covering map $\pi:\; \C^n \backslash 0 \arrow H$.
Since $\pi^*(L)$ is trivial, the space of $A$-finite \index[terms]{section!$A$-finite}
sections of $\pi^*(L)$ is identified with the space
of $A$-finite functions\index[terms]{function!$A$-finite} on $\C^n \backslash 0$. Since $A$ 
is linear, a function $f\in H^0(\C^n \backslash 0, \calo_{\C^n})$
is $A$-finite if and only if it is polynomial
(\ref{_cone_cover_for_LCK_pot_Theorem_}, Step 1).
The $\Z$-action on polynomial functions is induced
from the action of $A$ on $\C^n$; the corresponding
action of $\Z$ on $L$ is equal to $\lambda A$, 
where $\lambda=\chi(1)$ is the number associated
to the character  $\chi:\; \Z \arrow \C^*$.

Let $\alpha$ be an eigenvalue of $A$ which 
corresponds to some non-diagonal Jordan block,
 $\lambda:=(\alpha)^{-1}$, and $L$ be the line bundle
associated with the character $\chi$ such that
$\chi(1)=\lambda$. In Step 1 we use $L$ to produce the 
local system ${\Bbb V}$ on $S^1$; the corresponding
vector space is the space of all polynomial
functions on $\C^n$, and the monodromy\index[terms]{monodromy} is given
by $A\lambda\Id$. The multiplication
in the cohomology of this local system
is not an isomorphism, because the Jordan
blocks associated with the eigenvalue 1 are
non-diagonal (\ref{_mult_on_coh_of_S^1_Claim_}, Step 2).
\endproof



\section{Exercises}

\begin{enumerate}[label=\textbf{\thechapter.\arabic*}.,ref=\thechapter.\arabic{enumi}]

\item
Let $X$ be a Riemann surface, $Z\subset X$ a point,
and $j:\; X \backslash Z \arrow X$ an open embedding.
Prove that $j_* j^* \calo_X$ is not a coherent sheaf.

{\em Hint:} Show that $j_* j^* \calo_X$ is not 
locally finitely generated.\index[terms]{Riemann surface}

\item
Let $M$ be a complex variety,
$x\in M$ a point, and $\calo_{M,x}$
the ring of germs of holomorphic functions.
For any coherent sheaf $F$, denote by $F_x$
the stalk of $F$ in $x$; by definition, $F_x$
is a finitely--generated $\calo_{M,x}$-module.
A coherent sheaf $F$ on $M$
is called {\bf globally generated}
if the natural map \index[terms]{sheaf!globally generated}
$H^0(M, F) \otimes_{\calo_M} \calo_{M,x}\arrow F_x$
is surjective for all $x\in M$. 
\begin{enumerate}
\item  Suppose that a coherent sheaf $F$ on $\C^n \backslash 0$
can be extended to a coherent sheaf on $\C^n$.
Prove that $F$ is globally generated.
\item Suppose that $F$ is a globally generated sheaf
of $\calo_M$-modules on a complex manifold $M$. Prove that
there  exists a surjective sheaf morphism
$W\otimes_\C \calo_M\arrow F$,
where $W$ is a complex vector space, possibly
infinite-dimensional.
\item Let $M= \C^n \backslash 0$, $n \geq 2$, and
$F \subset W\otimes_\C \calo_M$ a subsheaf of a trivial
vector bundle. Prove that $F$ can be extended to 
a subsheaf $\tilde F$ of $W\otimes_\C \calo_{\C^n}$. 
Prove that $F$ is coherent if the space $W$ is finite-dimensional.

\item Let $F$ be a globally generated coherent sheaf on
  $\C^n \backslash 0$. Prove that there exists a
  neighbourhood $U$ of 0 in $\C^n$ such that the
  restriction of $F$ to $U\backslash 0$ is a quotient of
a sheaf $W\otimes_\C \calo_M$ by a coherent subsheaf,
with $W$ finite-dimensional.

{\em Hint:} Use the Noetherian property
for coherent sheaves.

\item Let $F$ be a globally generated coherent sheaf on
  $\C^n \backslash 0$, $n \geq 2$. Prove that it can be extended 
to a coherent sheaf on $\C^n$.
\end{enumerate}

 \item
Let $A\subset \C^n$ be an annulus
$A:= B_R\backslash B_r$. Prove that $H^1(A,  \calo)=0$,
where $\calo$ denotes the sheaf of holomorphic functions.

\item
Show that $H^1(\C^n \backslash 0, \calo)=0$.

{\em Hint:} Use the Grothendieck spectral sequence and
the previous exercise.\index[terms]{spectral sequence!Grothendieck}

\item
Prove that any holomorphic line bundle on $\C^n \backslash 0$
is trivial. 

{\em Hint:}  Use the exponential exact sequence and the previous exercise.

\end{enumerate}

\definition
A coherent sheaf over a complex manifold\index[terms]{sheaf!filtrable}
is called {\bf filtrable} if it has a filtration
by coherent subsheaves, with all subquotients of rank $\leq 1$.%
\footnote{We say that a coherent sheaf ${\cal F}$ over $M$ {\bf has rank $k$}
if it is a vector bundle of rank $k$ outside of its singular locus $Z\subset M$.}
 In dimension $>2$, all coherent sheaves over a diagonal Hopf manifold
are filtrable (\cite{_V:filtrable_}), and the lift of a holomorphic
bundle from Hopf manifold $H$ to its universal cover $\C^n \backslash 0$ 
can be extended to $\C^n$ (\ref{_extension_coherent_Theorem_}) 
when $n \geq 3$. In dimension 2,
not all reflexive sheaves are filtrable (\cite{_Brinzanescu_Moraru:rank2_}), and
the Andreotti-\index[persons]{Siu, Y.-T.}Siu extension theorem also does not work.\index[terms]{theorem!Andreotti--Siu extension} 
We call a reflexive coherent sheaf ${\cal F}$ \index[terms]{sheaf!reflexive}
on $H$ {\bf extensible} if its pullback to $\C^n \backslash 0$
can be extended to a coherent sheaf on $\C^n$.

\begin{enumerate}[label=\textbf{\thechapter.\arabic*}.,ref=\thechapter.\arabic{enumi}]
\setcounter{enumi}{5}

\item
Let ${\cal F}$ be a rank 1 coherent sheaf over a Hopf manifold. Prove
that it is extensible.

{\em Hint:} Use the previous exercise.

\item
Let $B$ be a filtrable vector bundle over a Hopf manifold.
Prove that it is extensible.

\item
Let $Z\subset M$ be a complex subvariety,
and $F$ a sheaf on $M$. Consider the functor
$F \mapsto F_Z$ taking $F$ to its subsheaf consisting
of all sections with support in $Z$
(that is, sections $f$ such that the restriction of
$f$ to $M \backslash Z$ vanishes). 
\begin{enumerate}
\item Prove that the functor $F \arrow F_Z$ is left exact.
Let ${\cal H}^i_Z(M,F)$ denote the corresponding 
derived functor, called {\bf local cohomology with
support in $Z$}. \index[terms]{cohomology!local}
\item Let $j:\; M \backslash Z \arrow M$ denote the 
tautological embedding.
Prove that the sequence 
\[
0 \arrow {\cal H}^0_Z(M,F) \arrow F \arrow j_* j^* F \arrow 0
\]
is exact when $F$ is flasque. Given a flasque resolution
$$0 \arrow F \arrow I_0 \arrow I_1 \arrow I_2 \arrow ...$$
prove that 
$$0 \arrow {\cal H}^0_Z(M,I_*) \arrow I_* \arrow
j_* j^* I_* \arrow 0$$\index[terms]{resolution!flasque}
 is an exact sequence of complexes.

\item
Let $0 \arrow A_* \arrow B_* \arrow C_* \arrow 0$
be an exact sequence of complexes of groups or sheaves. Prove
that there exists a long exact sequence of cohomology of these complexes
\begin{multline}
	 ... \arrow  H^i(C_*) \arrow H^{i+1}(A_*) \arrow  H^{i+1}(B_*)\arrow\\  
\arrow  H^{i+1}(C_*) \arrow H^{i+2}(A_*) \arrow ...
\end{multline}

\item
Construct an exact sequence
\[
0 \arrow {\cal H}^0_Z(M,F) \arrow F \arrow j_* j^* F \arrow
{\cal H}^1_Z(M,F) \arrow 0,
\]
and an isomorphism $R^ij_* j^*F \cong {\cal H}^{i+1}_Z(M,F)$
for all $i >0$ and all coherent sheaves $F$ on $M$.
\end{enumerate}

\item\label{_R^1j_*_vanish_Exercise_}
Let 
$0 \arrow A \arrow B \arrow C\arrow 0$ be an exact sequence 
of coherent sheaves 
on $\C^n \backslash 0$, where $A$ is locally free.
\begin{enumerate}

\item Let $H^*_0(\C^n, F)$ denote the cohomology
of a coherent sheaf $F$ on $\C^n$ with support in 0,
and let $j:\; \C^n \backslash 0 \arrow \C^n$ 
denote the standard embedding.
Prove that $R^ij_* j^*(\C^n, F)= H^{i+1}_0(\C^n, F)$, for any $i > 0$,
and any coherent sheaf $F$ on $\C^n$.

\item
Let $F$ be a locally free sheaf on $\C^n$.
In \cite[\S 4, Example 3]{_Hartshorne_LMS41_}, 
it is shown that $H^i_0(\C^n, F)=0$ when $i< n$.
Use this result to show that $R^1j_* \pi^* F=0$ when $n>2$.
\end{enumerate}

\item
Let $0 \arrow F \arrow G \arrow H \arrow 0$ be an exact sequence
of coherent sheaves on a complex manifold $M$. 
\begin{enumerate}
\item Let $U\subset M$ be an open set, $Z\subset U$ a complex
analytic subvariety, $\codim Z \geq 2$, and $V:= U\backslash Z$.
Consider the restriction maps
\[ \begin{CD}
0 @>>> H^0(U,F) @>>> H^0(U,G) @>\mu >> H^0(U,H)\\
&& @VrVV @VrVV @VrVV \\
0 @>>> H^0(V,F) @>>> H^0(V,G) @>>> H^0(V,H)
\end{CD}
\]
Prove that there exists an exact sequence
$A \arrow B \arrow C$, where 
\[ A= \ker r\restrict {\im \mu}, \ \ C= \coker r\restrict{H^0(U,G)},\ \ 
B= \coker r\restrict{H^0(V,F)}.
\]
\item Assume that
$G$ is normal and $H$ is torsion-free.\index[terms]{group!torsion-free} Prove that $F$ is normal.
\end{enumerate}

\item Let $B$ be a holomorphic bundle on a complex surface.
Assume that $B$ admits a holomorphic connection. Prove that
the rational  Chern classes $c_1(B, \Q)$, $c_2(B, \Q)$ \index[terms]{class!Chern}
(that is, the images of the Chern classes in $H^*(M,\Q)$)  vanish.

\item
A holomorphic bundle $B$ on $M$ is called {\bf simple} if $H^0(M, \End(B))=\C\Id$.
Let $B$ be a simple holomorphic bundle on the complex projective space $\C P^n$. Prove that $H^0(\C P^n, \Omega^k_{\C P^n}\otimes_{\calo_{\C P^n}} \End(B))=0$
for any $k >0$.\index[terms]{bundle!vector bundle!holomorphic!simple}

\item
Let $B$ be a simple bundle on $\C P^n$.\index[terms]{connection!holomorphic}
Prove that any holomorphic connection on $B$ is flat.

{\em Hint:} Use the previous exercise.

\item
Let $W \subset \R^n$ be a convex, open $\R^{>0}$-invariant subset, $0\notin W$,
and $G_W$ the group of affine automorphisms of $W$.
\begin{enumerate}
\item Prove that the image of $W$ in $\R P^{n-1}$
admits a $G_W$-invariant metric structure.
\item  Suppose  $W$ is a quadratic cone, defined in coordinates
by 
\[ W:=\left \{(x_1, ...., x_n) \ \ |\ \ x_1^2 > \sum_{i=2}^n x_i ^2\right\}.
\]
Prove that $G_W = \SO(1, n-1)$.
\item Consider a discrete, cocompact subgroup $\Gamma \subset G_W= \SO(1,n-1)$,
and let ${\Bbb P} W$ be the projectivization of $W$ equipped with 
a $G_W$-invariant Riemannian metric.
Prove that ${\Bbb P} W/\Gamma$ is a compact Riemannian orbifold of constant sectional curvature.
Prove that all compact Riemannian orbifolds of constant negative  sectional curvature
are obtained in this way.
\item Consider the quotient $\frac{W}{\Gamma\times \Z}$, where $\Gamma \subset G_W= \SO(1,n-1)$
is a discrete, cocompact subgroup of $\SO(1,n-1)$, and the action of $\Z$ on $W$ is generated 
by multiplication by $\lambda \in \R^{>1}$. Prove that $\frac{W}{\Gamma\times \Z}$
is a compact, flat affine manifold. Prove that it is not a complete flat affine manifold.
\end{enumerate}

\item
Let $\gamma:\; \C^2 \arrow \C^2$ be an automorphism mapping
 	$(x,y)$ to $\left(\alpha x, \beta
          y+x^2\right)$.
Prove that all eigenvalues of $\gamma^*$ on $\calo_{\C^n}$
are equal to $\alpha^m \beta^n$, for $m, n \in \Z^{\geq 0}$. Prove
that all these eigenvalues are realized by eigenvectors.

\item \label{_non_linear_Hopf_Exercise_}
Let $\gamma:\; \C^2 \arrow \C^2$ be an automorphism mapping
 	$(x,y)$ to $\left(\frac 1 2 x, \frac 1 4
          y+x^2\right)$.
\begin{enumerate}
\item
 Prove that it is a
          contraction. Prove that the
 	corresponding Hopf surface \index[terms]{surface!Hopf!resonant}
 	$H$ is resonant.

\item Prove that the space 
 of solutions of the equation
$\gamma^* f= \frac 1 2 f$, $f\in \calo_{\C^2}$ is 1-dimensional,
and generated by the coordinate function $x$.

\item Prove that the action of $\gamma^*$ preserves the 
ideal $I\subset \calo_{\C^2}$ generated by $x^3$, $xy$ and $y^2$.
Prove that the eigenvalues $\alpha_i$ of $\gamma^*$ on $I$ satisfy
 $|\alpha_i| \leq 1/8$.

\item Prove that any solution of the equation
$\gamma^* f= \frac 1 4 f$, $f\in \calo_{\C^2}$ is proportional
to $x^2$ modulo $I$. Prove that $f=x^2$.

\item Prove that there is no coordinate system $z_1, z_2\in \calo_{\C^2}$
on $\C^2$ such that $\gamma^*$ acts on the space $\langle z_1, z_2\rangle$
linearly.

\item
Prove that the bundle $TH$ does not admit a torsion-free flat
holomorphic connection.\index[terms]{connection!torsion-free}\index[terms]{connection!flat}\index[terms]{connection!holomorphic}
\end{enumerate}

\item\label{_torsion-free_from_torsion_Exercise_}
Let $\nabla$ be a holomorphic connection $\nabla$ on 
the tangent bundle $TH$ to a complex manifold $H$. Prove that
there exists a torsion-free holomorphic connection
on $TH$.

{\em Hint:} Replace $\nabla$ by $\nabla -\frac 1 2 T$, where
$T\in \Omega^2 H \otimes TH$ is the torsion form.

\item
Let $H$ be a Hopf surface\index[terms]{surface!Hopf} associated with a contraction
$\gamma:\; \C^2 \arrow \C^2$ centred in 0. Suppose that
the differential $D_0 \gamma$ has eigenvalues $1/2$ and $1/4$.
Prove that $H^0(H, \Omega^2 H \otimes_{\calo_H}\End(TH))=0$.

\item\label{_no_conne_on_Hopf_Exercise_}
Let $H$ be the non-linear Hopf surface constructed in 
Exercise \ref{_non_linear_Hopf_Exercise_}.
Prove that the holomorphic vector bundle $TH$
does not admit a holomorphic connection.

{\em Hint:} Use the previous exercise and 
Exercise \ref{_torsion-free_from_torsion_Exercise_}.

\item \label{_non_linear_Hopf_Example_}
Let $\gamma:\; \C^2 \arrow \C^2$ be a polynomial mapping
 	$(x,y)$ to $\left(\frac 1 2 x, \frac 1 4
          y+x^2\right)$.
\begin{enumerate}

\item Prove that $\gamma^*(dx) =  \frac 1 2  dx$, and
$\gamma^*(dy)=  \frac 1 4 dy + \frac 1 2  x dx$.

\item Consider a connection $\nabla_0$ on $\C^2$ 
such that $\nabla_0(d/dx)=\nabla_0(d/dy)=0$.
Prove that $\gamma^*(\nabla_0)(d/dx)=0$
and $\gamma^*(\nabla_0)(d/dy)=\frac 1 2 x dx\otimes d/dy$

\item Prove that $\gamma^*(\nabla_0)=\nabla_0 + A$,
where $A \in \Omega^1(\C^2)\otimes \End(T\C^2)$
is written as $A = dx \otimes Y$, where $Y(\frac d{dx})=0$
and $Y(\frac d{dy})= \frac 1 2 x \frac{d}{dy}$.

\item Let $\nabla = \nabla_0 + B$, where $B\in \Lambda^1(\C^2)\otimes \End(T\C^2)$
is any $\End(T\C^2)$-valued holomorphic 1-form. Prove that
$\gamma^*(\nabla) = \gamma^*(\nabla_0 + B) = \nabla_0 + A + \gamma^* B = \nabla - B + A + \gamma^* B$.
Prove that $\gamma^*$-invariant connections are in bijective correspondence
with solutions $B$ of the equation $B- \gamma^* B = A$.

\item Prove that $H$ admits a holomorphic connection if and only if
the image $[A]$ of $A$ in the group $H^0(\Lambda^1(\C^2)\otimes \End(T\C^2))_\Z$
of $\gamma^*$-coinvariants is zero. Prove that  $[A]\neq 0$.

\end{enumerate}

{\em Hint:} Use Exercise \ref{_no_conne_on_Hopf_Exercise_}.

\item
Let $B$ be a non-resonant Mall bundle on a Hopf manifold.
Prove that $B$ admits a unique holomorphic connection.

\item
Let $M$ be a non-resonant linear Hopf manifold,
$M=\frac{\C^n \backslash 0}{\langle A\rangle}$.
Prove that all automorphisms of $M$ act linearly on its
universal covering $\C^n\backslash 0$.

{\em Hint:} Use the previous exercise. 

\item
Let $H$ be a linear Hopf manifold, $H = \frac{\C^n\backslash 0} {\langle A \rangle}$,
and $G_A \subset \GL(\C^n)$ the Zariski closure of the group generated by the
linear contraction $A$. Consider a subsheaf $S \subset TH$, and
let $\tilde S \subset T(\C^n\backslash 0)$ be its pullback. 
\begin{enumerate}
\item Prove that $\tilde S$ is $G_A$-invariant.
\item Assume that $A= \lambda \Id$, where $\lambda \in \C$ satisfies
$0< |\lambda| < 1$. Consider the action $\tilde \rho$ of $\C^*$ on $\C^n\backslash 0$ defined
as $\tilde \rho(t)(z) = tz$, and let $\rho$ be the corresponding action
of $\C^*$ on $H$. Prove that any coherent subsheaf $S\subset TH$ is
$\rho$-invariant.
\end{enumerate}

\item
Recall that {\bf a rank 1 foliation} on a complex manifold
$M$ is a rank 1  coherent subsheaf $S\subset TM$ which
is {\bf saturated}, that is, non-singular in codimension 2
and normal. A rank 1 foliation is called {\bf smooth} if
it is a line bundle. Let $H = \frac{\C^n\backslash 0} {\langle A \rangle}$
be the classical Hopf manifold, $A= \lambda \Id$,
and $S\subset TH$ a rank 1 foliation. Suppose that
for some $x\in H$, the leaf of  $S$ in $x$
is tangent to a fibre of the natural projection
$\pi:\; H \arrow {\C P^{n-1}}$.
Prove that this leaf of $S$ 
coincides with a fibre of $\pi$.

{\em Hint:} Use the previous exercise.

\item \label{_foliation_on_classical_Hopf_compact_leaf_Exercise_}
The Baum--Bott theorem \index[terms]{theorem!Baum--Bott}
implies that any rank  1 foliation on 
$\C P^n$ is singular (\cite[page 287, remark C]{_Baum_Bott_}).
Let $S\subset H$ be a smooth rank 1 foliation on a classical Hopf manifold $H$.
\begin{enumerate} 
\item
Using the $\rho$-invariance of $S$, prove that there exists
a foliation $S_0\subset T\C P^{n-1}$ such that 
\[ v \in S\subset T_x H\ \Leftrightarrow D\pi(v) \in S_0 \subset T\C P^{n-1}.
\]
\item
Let $x\in \C P^{n-1}$ be a point where
$\pi(S)$ has a singularity; such a point exists, again, by the Baum--Bott theorem.
Prove that any leaf of $S$ intersecting the fibre $\pi^{-1}(x)$
is tangent to this fiber.
\item Prove that any smooth rank 1 foliation on a classical Hopf manifold $S$
has a compact leaf.\footnote{We are grateful to Jorge Vit\'orio \index[persons]{Pereira, J. V.} Pereira for this argument.}
\end{enumerate}
{\em Hint:} Use the previous exercise.

\item\label{_splitting_foliation_Exercise_}
Let $M$ be a smooth manifold,
and $\Sigma\subset TM$ be a bundle tangent
to a free action of a group $G=\C$ on $TM$,
and $TM = \Sigma \oplus B$ a $G$-invariant
splitting. Prove that the Frobenius
form $\Phi:\; \Lambda^2 B \arrow \Sigma=C^\infty(M,\C)$
is closed and basic. Prove that 
its basic cohomology\index[terms]{cohomology!basic} class is independent
from the choice of the splitting $TM = \Sigma \oplus B$.

\item\label{_parallelizable_Exercise_}
Let $M$ be a Vaisman manifold\index[terms]{manifold!Vaisman}
with the tangent bundle $TM$ holomorphically 
trivial.\footnote{A complex manifold with trivial
tangent bundle is called {\bf complex parallelizable.}}
The following sequence of exercises aims to show
that this is impossible.
\begin{enumerate}
\item Let $\Sigma\subset TM$ be the holomorphic sub-bundle of $TM$ generated
by the Lee field.\index[terms]{Lee field} Prove that $TM$ splits, as a holomorphic
vector bundle, $TM= \Sigma \oplus B$.
\item Prove that the Frobenius form\index[terms]{form!Frobenius} $\Phi:\; \Lambda^2 B \arrow \Sigma$ 
of $B$ is holomorphic.
\item Prove that all holomorphic vector fields on $M$ commute
with the Lee field.\index[terms]{Lee field} You may use 
\ref{_holo_tensor_on_Vaisman_Lee_invariant_Theorem_}.
\item Prove that any holomorphic decomposition
$TM= \Sigma \oplus B$ is Lee invariant.
\item Let $\f:\; TM \arrow B$ be the projection along $\Sigma$.
Consider the Frobenius form\index[terms]{form!Frobenius}
$\Phi:\; \Lambda^2 B \arrow \Sigma=\calo_M$ 
as a 2-form on $B$. Prove that 
$\Phi_0:=\Phi\circ f\in \Lambda^{2,0}(M)$
is a $\Sigma$-basic form on $M$.
\item Using Exercise \ref{_splitting_foliation_Exercise_},
prove that $\Phi_0$ is cohomologous to 
the transversal K\"ahler form\index[terms]{form!K\"ahler!transversal} $\omega_0$
in the group $H^2_\kah(M)$ of $\Sigma$-basic cohomology.\index[terms]{cohomology!basic}
\item 
Prove that this is impossible, because $\Phi_0$
is a (2,0)-form, and $\omega_0$ is a (1,1)-form.
\end{enumerate}

\end{enumerate}


\chapter[Kuranishi and Teichm\"uller spaces for LCK manifolds]{Kuranishi and Teichm\"uller spaces for LCK manifolds with po\-ten\-tial}\index[terms]{manifold!LCK!with potential}
\label{_Kuranishi_Chapter_} \index[terms]{space!Kuranishi} \index[terms]{space!Teichm\"uller}

{\setlength\epigraphwidth{0.6\linewidth}
\epigraph{{\em \cyr Ponyal teper{\cprime} ya: nasha svoboda\\
Tol{\cprime}ko ottuda b{\cprime}yushchi{\u i} svet,\\
Lyudi i teni stoyat u vkhoda\\
V zoologicheski{\u i} sad planet.\\[10pt]}  {\em
\font\tenrm = cmssi17 at 10pt
\tenrm (I've grasped it at last: our freedom\\
is only a light pulsating from far -\\
people and shadows stand at the entrance\\
to the zoo of the wandering stars.)}}
{{\sc \scriptsize Nikolay Gumilev, ``The Lost Tram'',\\ translated by Boris Dralyuk}} 
}

\section{Introduction}

\subsection{Deformation spaces}

``A deformation of complex structure'' usually means a continuous
or a complex analytic family of complex varieties. For most of this
book, we are interested in deformations of complex structures
on a manifold or on an orbifold. In this case, the smooth structure
of the underlying manifold does not change, and one can study the
deformation of the complex structure tensor on the same underlying
manifold (or orbifold). This is the approach taken by \index[persons]{Kodaira, K.} Kodaira-\index[persons]{Spencer, D. C.}Spencer
and \index[persons]{Kuranishi, M.} Kuranishi in their research on deformations.

Even earlier, Teichm\"uller \index[persons]{Teichm\"uller, O.} explored the space of complex structures
on a Riemann surface. In dimension 1, complex structures are the same
as conformal structures on an oriented 2-dimensional real manifold,
and Teichm\"uller\index[persons]{Teichm\"uller, O.} invented a way of parametrising the isotopy classes
of conformal structures on a Riemann surface; this gives another
way of producing the moduli of deformations of complex structures.

{\bf A dilatation} of a \index[terms]{dilatation}
smooth diffeomorphism $f:\; S \arrow R$ of Riemann surfaces\index[terms]{Riemann surface} in a point $x\in S$
is the excentricity\footnote{The ratio of the larger axis to the smaller one.}
 of the ellipse in $T_{\f(x)} R$ obtained as the image
of a circle in $T_x S$. It is equal to 1 whenever $f$ is conformal in $x$.
Teichm\"uller\index[persons]{Teichm\"uller, O.} has shown that in each homotopy class of diffeomorphisms
$f:\; S \arrow R$ there exists a unique map, called 
{\bf the Teichm\"uller extremal map}, which minimizes the 
supremum of the dilatation (\cite{_Ahlfors:quasi-conformal_,_Nag:Teichmuller_}).
The {\bf Teichm\"uller distance} between $S$ and $R$ is the\index[terms]{map!Teichm\"uller extremal}\index[terms]{distance!Teichm\"uller}
logarithm of the dilatation of the Teichm\"uller extremal map.

The {\bf Teichm\"uller space} $\Teich$ of complex structures on a Riemann
surface $S$ is the space of all complex structures on $S$
up to an isotopy. The Teichm\"uller distance defines a metric on $\Teich$.

For the Teichm\"uller extremal map $f:\; S \arrow R$, the
direction of the axes of the ellipsoid $T_{\f(x)} R$ holomorphically
depends on $x$, defining a quadratic differential on $S$. This allows
one to give complex coordinates on the space of complex structures on $S$,
identifying it with an open subset in the space of quadratic differentials
$H^0(S, \Sym^2(\Omega^1 S))$ (\cite{_Nag:Teichmuller_}). 

This defines a complex structure on the
Teichm\"uller space of Riemann surfaces.\index[terms]{Riemann surface} This space is Stein by Bers and Ehrenpreis (\cite{_Bers_Ehrenpreis_}).
The Teichm\"uller distance gives a Finsler metric on $\Teich$,\index[terms]{metric!Finsler}
that is  equal to its \index[persons]{Kobayashi, S.} Kobayashi pseudometric\footnote{We define the Kobayashi pseudometric in \index[terms]{space!Teichm\"uller}
Exercise \ref{_Kobayashi_pseudometric_Exercise_}.}\index[terms]{theorem!Royden}
by \index[persons]{Royden, H.} Royden's theorem, \cite{_Royden:Teichmuller_Kobayashi_}. 

Teichm\"uller's\index[persons]{Teichm\"uller, O.} approach to the deformation theory was very successful
in dimension 1, but higher dimensions required different methods. The
breakthrough in this direction came with a series of papers by
Kodaira and \index[persons]{Spencer, D. C.} Spencer,\index[terms]{deformation theory} \cite{_Kodaira_Spencer:Deformations1-2_,_Kod-Spen-AnnMath-1960_}.
\index[persons]{Kodaira, K.} Kodaira and Spencer realized that the deformations of complex structures can be obtained
from the first cohomology of vector fields. A couple of years later, \index[persons]{Kuranishi, M.} Kuranishi
constructed the Kodaira-\index[persons]{Spencer, D. C.}Spencer deformation space explicitly
(\cite{_Douady:Bourbaki_,_Kuranishi:new_,_Kuranishi:note_}).
We give a modern version of this construction in Section 
\ref{_Kuranishi_Section_} and compute the Kuranishi space
for the Hopf manifolds explicitly in Section \ref{_Kuranishi_Hopf_Section_}.

Independently from Teichm\"uller\index[persons]{Teichm\"uller}, \index[persons]{Grothendieck, A.} Grothendieck has constructed
the algebraic (axiomatic) version of the theory of moduli spaces
(\cite{_Grothendieck:Teichmuller_}).\index[terms]{deformation functor}
He did it from the general principles, by defining
the deformation functor, and showing that this
functor is representable (at least, locally). 
To complete the proof, he introduced the\index[terms]{Hilbert scheme}
notion of a Hilbert scheme (Section \ref{_Hilbert_scheme_Section_}), 
and showed that it can be used to construct the deformation space.

The Teichm\"uller space\index[terms]{space!Teichm\"uller} is defined as the quotient of the
space of all complex structure tensors by the action\index[terms]{group!of isotopies}
of the group of isotopies. In \cite{_Meersseman:stacks_,_Meersseman:Kuranishi_},
L. \index[persons]{Meersseman, L.} Meersseman considered the natural quotient map from\index[terms]{map!Kuranishi to Teichm\"uller}
the Kuranishi to Teichm\"uller space. From the definition
it is clear that the Teichm\"uller space can be obtained
as a union of open sets, given as images of the  Kuranishi
to Teichm\"uller maps.\index[terms]{map!Kuranishi to Teichm\"uller} Whenever the Teichm\"uller space
admits some kind of geometric structure, such as the
complex structure, it should be obtained from the
local charts produced by the Kuranishi
to Teichm\"uller maps. We explicitly compute the  Kuranishi 
to Teichm\"uller map for general Hopf manifolds (\ref{_Teich_nr_for_Hopf_}).

\subsection{Deformations of Hopf surfaces: a short survey}\index[terms]{surface!Hopf}

The deformations of Hopf surfaces is an old subject, which was
already well studied in 1980-ies.\index[terms]{surface!Hopf} Deformations of Hopf surface
are controlled by the first cohomology of the tangent bundle.
By Serre's duality,\index[terms]{duality!Serre}
the second cohomology of the tangent bundle $TH$ is the same
as sections of $\Omega^1 H \otimes K_H$, and this group
vanishes (\ref{_holo_forms_Hopf_vanish_Proposition_}). 
Then, the Riemann--Roch formula implies that $H^1(TH)= H^0(TH)$, and
the latter group (identified with the Lie algebra of the
automorphism group of $H$) is not very hard to compute. 
However, $H^0(TH)$ is wildly different for different
classes of Hopf manifolds. By the early 1980-ies,\index[terms]{Riemann--Roch formula} 
the Hopf surfaces were classified together with
their automorphism groups (M. \index[persons]{Namba, M.} Namba, \cite{_Namba:Hopf_}).
The answer that Namba obtained depends on the classification
of Hopf surfaces obtained by \index[persons]{Kodaira, K.} Kodaira. \index[terms]{surface!Hopf}
J. \index[persons]{Wehler, J.} Wehler (\cite{_Wehler_}) used the results of \index[persons]{Namba, M.} Namba
and the classification of Hopf surfaces due to Kodaira
to compute the versal deformations of all classes of
Hopf manifolds. In \index[persons]{Wehler, J.} Wehler's notation, there are\index[terms]{deformation!versal}
5 classes of Hopf surfaces: IV, III, IIa, IIb, IIc.
Class IV corresponds to classical Hopf surfaces,
$H=\frac{\C^2 \backslash 0}{\langle \gamma\rangle}$,
where $\gamma$ is a scalar contraction. Class III is all
other elliptic Hopf surfaces. Class IIc is non-linear Hopf
surfaces, class IIb is non-diagonal linear Hopf surfaces, and
class IIa is diagonal, non-elliptic Hopf surfaces.
The group of automorphisms and the Kuranishi space is computed \index[terms]{surface!Hopf}
separately for each class.  His deformation space is essentially 
the same as the \index[terms]{space!Kuranishi} Kuranishi space we describe in
Section \ref{_Kuranishi_Section_}.

Contemporaneously with \index[persons]{Wehler, J.} Wehler, C. \index[persons]{Borcea, C.} Borcea published a\index[terms]{deformation!of Hopf manifolds}
paper \cite{_Borcea_} on deformations of Hopf manifolds, where
he described the versal deformations of diagonal Hopf manifolds,
essentially proving the same statement as we formulated 
in \ref{_Kuranishi_on_Hopf_Theorem_}.
He also proved that the Kuranishi deformation of a Hopf
manifold is unobstructed.
Later, V. \index[persons]{Palamodov, V. P.} Palamodov (\cite{_Palamodov:Hopf_}) 
re-proved the results of \index[persons]{Wehler, J.} Wehler, interpreting 
the versal deformations of Hopf manifolds as the moduli of holomorphic
contractions, and used the  Poincar\'e-Dulac theorem\index[terms]{theorem!Poincar\'e-Dulac}
to classify the holomorphic contractions of $\C^n$, obtaining
an explicit description of the versal deformation.


Finally, A. \index[persons]{Haefliger, A.} Haefliger proved that the group 
$H^2(TH)$ vanishes for Hopf manifolds (\cite{_Haefliger:Hopf_});
he used this result to show that the  Kuranishi deformation
of the Hopf manifold is unobstructed, and described the
local deformation space explicitly. His graduate student
D. Mall extended this vanishing result to a wide class
of holomorphic bundles, that we called {\em the \index[persons]{Mall, D.} Mall bundles}
(Chapter \ref{_Mall_bundles_Chapter_}).

\subsection{Teichm\"uller space of Hopf manifolds and ap\-pli\-ca\-tions to LCK geometry}\index[terms]{space!Teichm\"uller}

In algebraic geometry\index[terms]{geometry!algebraic}, the  moduli spaces were defined\index[terms]{moduli space}
by A. \index[persons]{Grothendieck, A.} Grothendieck axiomatically, as objects representing
the deformation functors. Existence of representing objects
in the category of schemes or the category of algebraic
spaces follows from results about representable functors.
In complex geometry, this approach does not
seem to work; instead we have the Teichm\"uller space,
defined as the quotient of the space of all complex
structure operators by the group of isotopies
(\ref{_Teich_Definition_}). \index[terms]{space!Teichm\"uller}

Teichm\"uller spaces are often non-Hausdorff.
We call two points $x, y\in M$ in a topological space
{\bf non-separable} if any neighbourhood of $x$ intersects
with any neighbourhood of $y$.\index[terms]{non-separable points}

Non-resonant Hopf manifolds are equipped with a canonical
torsion-free flat connection\index[terms]{connection!torsion-free}\index[terms]{connection!flat} (\ref{_flat_Hopf_is_linear_Theorem_}),
identifying their Teichm\"uller space\index[terms]{space!Teichm\"uller} with the quotient of the 
space of linear contractions of $\C^n$ up to the conjugation.
This quotient space is already non-Hausdorff \index[terms]{manifold!Hopf!non-resonant}
(\ref{_semisimple_operator_approx_Proposition_} 
and Subsection \ref{_Teichmuller_LCK_with_pot_Subsection_});
we show that all points in the Teichm\"uller space 
of Hopf manifolds are non-separable from points
which correspond to linear Hopf manifolds (\ref{_non_separable_Remark_}).

The main practical application of this chapter is a
result about deformations of LCK manifolds with potential\index[terms]{manifold!LCK!with potential}.
We prove that for any LCK manifold $M$ with potential, $\dim_\C M \geq 3$
the corresponding point in the Teichm\"uller space is
non-separable from a point which corresponds to a 
Vaisman-type complex structure. This is used in Chapter
\ref{_Lee_classes_Chapter_} to describe the space of cohomology classes
of Lee forms\index[terms]{form!Lee} on an LCK manifold with potential.\index[terms]{manifold!LCK!with potential}

The proof of this observation actually uses the third 
construction of the moduli spaces, namely the Hilbert scheme.\index[terms]{Hilbert scheme}
We start with an LCK manifold $M$ with a potential embedded to
a linear Hopf manifold $H = \frac{\C^n \backslash 0}{\langle \gamma \rangle}$
(by \ref{embedding}, any LCK manifold with 
potential can be embedded). The preimage $\tilde M$ of $M$ in
$\C^n \backslash 0$ is an open algebraic cone\index[terms]{cone!algebraic} (\ref{_cone_cover_for_LCK_pot_Theorem_}).
 Using the Jordan--Chevalley theorem
we decompose $\gamma$ onto a semisimple and a unipotent part,
$\gamma=\gamma_s \gamma_u$, and notice that the algebraic
cone $\tilde M$ is necessarily $\gamma_s$ and $\gamma_u$-invariant.
Then the quotient $M_s:=\frac{\tilde M}{\langle \gamma_s \rangle}$
is a Vaisman manifold (\ref{vaisman_embed}). We prove that the corresponding\index[terms]{manifold!Vaisman}
point in the Teichm\"uller space is non-separable 
from $M= \frac{\tilde M}{\langle \gamma \rangle}$.

The proof uses the Hilbert scheme of algebraic
cones, interpreted as projective orbifolds in\index[terms]{cone!algebraic} 
$\C P^{n-1}= {\Bbb P}(\C^n)$. By \ref{_semisimple_operator_approx_Proposition_}, 
the point in the Teichm\"uller space
of a Hopf manifold $\frac{\C^n \backslash 0}{\langle \gamma \rangle}$
is non-separable from $\frac{\C^n \backslash 0}{\langle \gamma_s \rangle}$:
there exists a sequence $r_i\in \GL(n, \C)$ such that 
$\lim_i \gamma^{r_i} = \gamma_s$, where $x\mapsto x^{r_i}$ denotes the
conjugation, $x^{r_i} = r_i x r_i^{-1}$. 
However, the cone $\tilde M\subset \C^n$
is not $r_i$-invariant. 

To prove that
$M$ and $M_s$ correspond to non-separable points in the Teichm\"uller space\index[terms]{space!Teichm\"uller}, we construct
a sequence of diffeomorphisms $\nu_i\in \Diff_0(M)$
such that $\lim_i \nu_i(I)= I_s$, where
$I$ is the complex structure tensor on $M$ and
$I_s$ the complex structure on $M_s$. Using the Hilbert 
schemes, we prove that the sequence $\{\tilde M^{r_i}\}$
of algebraic cones\index[terms]{cone!algebraic} in $\C^n \backslash 0$ conjugated by $r_i$
converges to $\tilde M$ in $C^\infty$-topology.\index[terms]{topology!$C^\infty$}
Then 
\[ \lim_i \frac{\tilde M^{r_i}}{\langle \gamma^{r_i} \rangle}=
\frac{\tilde M}{\langle \gamma_s \rangle}
\]
giving the non-separability of the points 
$\frac{\tilde M}{\langle \gamma_s \rangle}$
and $\frac{\tilde M}{\langle \gamma \rangle}$
in the Teichm\"uller space.

\section{The Kuranishi space}
\label{_Kuranishi_Section_} \index[terms]{space!Kuranishi}

\subsection{Nijenhuis--Schouten and Fr\"olicher--Nijenhuis brackets}
\label{_F-N_bracket_Subsection_}

In this subsection,\index[terms]{bracket!Nijenhuis--Schouten}\index[terms]{bracket!Fr\"olicher--Nijenhuis} \index[persons]{Nijenhuis, A.}
we attempt to give an introduction to various brackets on the
multivector (polyvector) fields and vector fields with coefficients in differential forms.
We follow \cite{_Kosmann_Schwarzbach:schouten_}.

\hfill

{\bf The  Nijenhuis--Schouten bracket} (\cite{_Michor:Schouten-Nijenhuis_}), sometimes also called simply 
{\bf the Schou\-ten bracket},\index[terms]{bracket!Schouten} is enormously\index[terms]{structure!Poisson}
useful for the study of Poisson structures. It is a bracket 
$\Lambda^p TM \times \Lambda^q TM\arrow \Lambda^{p+q-1} TM$ that can be
defined axiomatically as follows. First, it is
graded anticommutative and satisfies the graded Jacobi identity
(see Subsection \ref{_superalgebras_Subsection_} below).
Second, for any $\xi \in \Lambda^p TM$, the \index[terms]{graded Jacobi identity}
bracket $[\xi, \dash]$, considered to be  a map $\Lambda^p TM\arrow \Lambda^{p+q-1} TM$\index[terms]{graded derivation}
is a graded derivation (that is, satisfies the graded Leibniz rule, 
\ref{_graded_derivation_Definition_}). Third, the Nijenhuis--Schouten bracket\index[terms]{bracket!Nijenhuis-Schouten}
$\Lambda^1 TM \times \Lambda^1 TM\arrow \Lambda^{1} TM$ is equal to the commutator on $\Lambda^1 TM =TM$.
The graded Leibniz rule allows to extend the  Nijenhuis--Schouten bracket 
 from $TM$ to $\Lambda^pTM$.\index[terms]{graded Leibniz rule}

A bivector $\pi\in \Lambda^2 TM$ defines a Poisson structure on $M$
if and only if $[\pi, \pi]=0$. This is
in fact one of the first definitions of Poisson structures,\index[terms]{structure!Poisson} due to \index[persons]{Lichnerowicz, A.} Lichnerowicz (\cite{_Lichnerowicz:Poisson_}).

Another bracket, called Nijenhuis or Fr\"olicher--Nijenhuis bracket, acting as \index[terms]{bracket!Nijenhuis}
\[
\Lambda^p M\otimes  TM \times \Lambda^q M\otimes TM\arrow \Lambda^{p+q} M\otimes TM,
\]
first appeared in \cite{_Nijenhuis:concomitants_} and \cite{_Frolicher_Nijenhuis:vector-valued_};
it can be defined the same way using the graded  Leibniz rule.  We explain its definition
by identifying the vector-valued differential forms with derivations of the de Rham algebra.

The algebra of vector fields acts on differential forms by Lie derivatives.
This is an action which commutes with de Rham differential and satisfies
the Leibniz identity;\index[terms]{Leibniz identity} in fact, the Lie derivatives can be characterized
as derivations $\delta:\; \Lambda^* M\arrow \Lambda^* M$ 
of the de Rham algebra which commute with de Rham differential and preserve the grading
(Exercise \ref{_deriva_de_Rham_Exercise_}).

We extend this action to vector-valued differential forms, using the following observation
(Exercise \ref{_deriva_de_Rham_forms_Exercise_}). Let $\delta:\; \Lambda^* M\arrow \Lambda^{*+d} M$ 
be a first order differential operator which  commutes with the de Rham differential and satisfies
the Leibniz identity.\index[terms]{Leibniz identity} Then $\delta$ is given by a section of $\Lambda^d M\otimes TM$,
where $TM$ acts on $\Lambda^* M$ as the Lie derivative, 
and $\Lambda^d M$ acts on $\Lambda^* M$ as the multiplication by a $d$-form.

The Fr\"olicher-Schouten bracket on $\Lambda^* M\otimes  TM$
maps the derivations $\delta_1, \delta_2:\; \Lambda^* M\arrow \Lambda^{*+*} M$ to the supercommutator\index[terms]{supercommutator}
$\{\delta_1, \delta_2\}$ (Exercise \ref{_supercommutator_Exercise_}).

\index[persons]{Nijenhuis, A.} Nijenhuis defined this bracket in conjunction with his research on almost complex structures.
Consider an almost complex structure $I \in \End(TM)= \Lambda^1 M\otimes  TM$
as a 1-form with coefficients in $TM$. This 1-form defines an odd derivation
$\Lambda^* M\arrow \Lambda^{*+1} M$, which equals $[d, W_I]$,
where $W_I$ acts on $(p,q)$-forms as multiplication by $\1(p-q)$ 
(Exercise \ref{_almost_complex_str_derivation_Exercise_}).
The anticommutator of this derivation
with itself is  a section of $\Lambda^2 M\otimes  TM$,
that is  equal to the \index[persons]{Nijenhuis, A.} Nijenhuis tensor.\index[terms]{tensor!Nijenhuis}

\subsection{Kuranishi space: the definition} \index[terms]{space!Kuranishi}

We briefly introduce the  Kuranishi 
space\index[terms]{space!Kuranishi}; for more details and reference,
see \cite{_Catanese:Moduli_} and \cite{_Barannikov_Kontsevich_}.

\hfill

\definition
Let $M$ be a complex manifold, and $\bar\6:\; \Lambda^{0,i}(M)\arrow  \Lambda^{0,i+1}(M)$
the complex structure operator. Let $\gamma\in \Lambda^{0,1} \otimes T^{1,0}M$
be a (0,1)-form with coefficients in (1,0)-vector fields. As in Subsection \ref{_F-N_bracket_Subsection_},
we interpret $\gamma$ as a derivation of the de Rham algebra, $\gamma:\; \Lambda^{0,i}(M)\arrow  \Lambda^{0,i+1}(M)$.
{\bf The Maurer--Cartan equation}\index[terms]{Maurer--Cartan equation}
is $(\bar\6 + \gamma)^2=0$, where $\bar\6$ and $\gamma$  are 
considered to be  first order differential operators on differential forms. 

\hfill

\remark
By Newlander--Nirenberg theorem, the integrability of an almost complex structure is equivalent
to $\bar\6^2=0$ (Exercise \ref{_Nijenhuis_tensor_Exercise_}). Therefore, $(\bar\6 + \gamma)^2=0$
is equivalent to the integrability of the complex\index[terms]{theorem!Newlander--Nirenberg}
structure associated with $\bar\6 + \gamma$ considered
as a holomorphic structure operator. In other words,\index[terms]{structure!holomorphic}
every solution of the Maurer--Cartan equation\index[terms]{Maurer--Cartan equation} gives
a deformation of complex structures. 
The Maurer--Cartan equation 
 is rewritten as $\bar\6  \circ\gamma + \gamma \circ \bar \6  =- \gamma^2$.
Sometimes this is further simplified to $\bar\6(\gamma)= - \frac 1 2 \{\gamma, \gamma\}$,
where $\{\gamma, \gamma\}$ denotes the graded commutator. We will use another
convention, which simplifies the computations, by writing
$\bar\6(\gamma)= - \{\gamma, \gamma\}$; this assumes
$\bar\6(\gamma):= \{\bar\6, \gamma\}$.\index[terms]{graded commutator}

\hfill

Write $\gamma= \sum_{n=0}^\infty t^n \gamma_n$, where $t$ is a
formal parameter. The Maurer--Cartan equation becomes
\begin{equation}\label{_MC_n_Equation_}
\bar\6(\gamma_n) = - \sum_{l+k=n-1}\{\gamma_l, \gamma_k\}.
\end{equation}
As in Subsection \ref{_F-N_bracket_Subsection_},
we understand $\gamma_n:\; \Lambda^{0,*}(M) \arrow \Lambda^{0,*+1}(M)$
as a first order differential operator, as follows:
\[ \left(\sum \xi_i^{0,1} \otimes X_i^{1,0}\right)(\alpha)=
   \sum \xi_i^{0,1}\wedge\Lie_{X_i^{1,0}}(\alpha),
\]
for any $\alpha\in \Lambda^{0,p}(M)$.

The supercommutator $\{\gamma_l, \gamma_k\}$ is the 
Fr\"olicher-Nijenhuis bracket, defined in Subsection \ref{_F-N_bracket_Subsection_}.\index[terms]{bracket!Fr\"olicher-Nijenhuis}
The bracket $\{\cdot, \cdot\}$ can be understood as
the graded commutator\index[terms]{graded commutator} of 
$\Lambda^{0,*}(M)$-valued vector fields, considered to be  a first
order differential operator $\Lambda^{0,*}(M) \arrow \Lambda^{0,*+1}(M)$:\footnote{See Subsection \ref{_superalgebras_Subsection_} for definitions
and conventions about superalgebras and graded commutators.}
\begin{multline*}
\left\{ \sum \xi_i^{0,1} \otimes X_i^{1,0}, \sum \zeta_k^{0,1} \otimes Y_k^{1,0}\right\}
=\\ = \sum\xi_i^{0,1}\wedge \zeta_k^{0,1}\otimes [X_i^{1,0}, Y_k^{1,0}]\in 
\Lambda^{0,2}(M) \otimes T^{1,0}M.
\end{multline*}
Notice that this graded commutator is actually commutative, because
it is anticommutative on both tensor components.

By induction, we may assume that 
\begin{multline}\label{_Maurer--Cartan_inductive_Equation_}
\bar\6\left(\sum_{l+k=n-1}\{\gamma_l, \gamma_k\}\right)
=\\
= \sum_{l+k+m=n-2}\!\!\!\!\{\{\gamma_l, \gamma_k\}, \gamma_m\}
\, - \sum_{l+k+m=n-2}\!\!\!\!\{\gamma_l, \{\gamma_k, \gamma_m\}\}=0
\end{multline}
(the sum vanishes because $\{A,\{A, A\}\}=0$
for any $A$). Then \eqref{_MC_n_Equation_}
has a solution if and only if the form
$\sum_{l+k=n-1}\{\gamma_l, \gamma_k\}\in \Lambda^{0,2}(M)\otimes T^{1,0}(M)$
is $\bar \6$-exact, that is, if its cohomology class in $H^2(M, TM)$ (the 
cohomology of holomorphic vector fields) vanishes.

\hfill

We introduce a Hermitian metric on $M$, obtaining the harmonic
decomposition 
\[ 
\Lambda^{0,1}(M) \otimes T^{1,0}(M)= \im \bar\6 \oplus \im \bar\6^*
\oplus \ker \Delta_{\bar\6}
\]
where $\Delta_{\bar\6}=\{ \bar\6, \bar\6^*\}$
with $\ker \Delta_{\bar\6}= H^1(M, TM)$ (the cohomology of the sheaf of
holomorphic vector fields). Using this decomposition, we group
the sum $\gamma= \sum \gamma_i$ in such a way that
$\gamma_0 \in  \ker \Delta_{\bar\6}\oplus \im \bar\6$, and
choose the solutions $\gamma_i, i>0$, of the equation \eqref{_MC_n_Equation_}
in $\im \bar\6^*$. 

\hfill

\claim\label{_M_i_Kuranishi_Claim_}
Let $\gamma_0 \in \Lambda^{0,1}(M)\otimes T^{1,0}(M)$
be a $\bar\6$-closed $T^{1,0}(M)$-valued (0,1)-form.
Assume that the equation  \eqref{_MC_n_Equation_} has a solution
for $n < n_0$. Denote by $M_{n_0}(\gamma_0)$ the cohomology class
of $\sum_{l+k=n_0-1}\{\gamma_l, \gamma_k\}$. Then $M_{n_0}(\gamma_0)$
is uniquely determined by the cohomology class of $\gamma_0\in H^1(M, TM)$
(the cohomology of holomorphic vector fields).\footnote{In deformation theory,
the operations $M_i$ are often called {\bf the Massey products};\index[terms]{Massey product}
see \cite{_Babenko_Taimanov_}  for an explanation of how the Maurer--Cartan equation is related
to the Massey products.}

\hfill

\proof
Suppose that $\gamma_0'$ and $\gamma_0''$ are two choices
of $\gamma_0$ in the same cohomology class, with $\gamma_0' - \gamma_0''=\bar\6\alpha$.
Then 
\begin{multline*} 
\bar\6(\{\alpha, \gamma_0''+\bar\6\alpha\})=
\bar\6(\{\alpha, \gamma_0'\}) = \{\gamma_0' - \gamma_0'',\gamma_0'\}+
 \{\alpha, \bar\6\gamma_0'\}= \\= \{\gamma_0' - \gamma_0'',\gamma_0'\}=
\{\gamma_0', \gamma_0'\}-\{\gamma_0'', \gamma_0'\}.
\end{multline*}
This implies that $\{\gamma_0', \gamma_0'\}$ and  $\{\gamma_0'', \gamma_0'\}$
are equal up to the image of $\bar\6$. 
By the same reason, $\{\gamma_0', \gamma_0'\}$
and $\{\gamma_0'', \gamma_0''\}$ are also equal up to the image of $\bar\6$.

Using induction, we obtain that after passing from $\gamma_0$ to $\gamma_0'$, the
sum \[ \sum_{l+k=n_0-1}\{\gamma_l, \gamma_k\}\] is modified by the addition of a term
that belongs to the image of $\bar\6$, and hence  the cohomology class of
$M_{n_0}(\gamma_0)$ remains unchanged.
\endproof

\hfill

\definition\label{_Kuranishi_Definition_}
{\bf The Kuranishi space}\index[terms]{space!Kuranishi} of deformations of the complex structures is
the set of all $\gamma_0\in H^1(M, TM)$ such that
$M_{i}(\gamma_0)=0$ for all $i>0$.

\hfill

From \ref{_M_i_Kuranishi_Claim_}, we obtain that each
$\gamma_0$ in the  Kuranishi space produces a formal solution of the
Maurer--Cartan equation. However, we can choose $\gamma_i$ inductively
by taking $\gamma_i:= G_{\bar\6}\sum_{l+k=n_0-1}\{\gamma_l, \gamma_k\}$,
where $G_{\bar\6}$ is the Green operator that is  inverse to $\bar\6$
on $\bar\6$-exact forms. This choice is in a sense optimal, because
the operator $G_{\bar\6}$ is defined in such a way that\index[terms]{operator!Green}
$w:=G_{\bar\6}(v)$ is the form with minimal $L^2$-norm\index[terms]{norm!$L^2$} which
satisfies $\bar\6(w)=v$. Using the eigenvalue estimate for $G_{\bar\6}$
and some combinatorial estimates for the number of terms in
the recursive formulas \eqref{_Maurer--Cartan_inductive_Equation_},
it is not hard to show that this formal series converges
for $\gamma_0$ sufficiently small (see e.g. \cite[Section 7]{_V:Hyperholomorphic_}).\index[terms]{space!Kuranishi}
This solution of the Maurer--Cartan equation gives a deformation of the
$\bar\6$-operator in a neighbourhood of 0 in the Kuranishi space.
Usually, one considers the  Kuranishi space as a germ of
the formal space we defined above, and the Kuranishi deformation\index[terms]{deformation!Kuranishi}
as the corresponding germ of the deformation of complex structures.

\hfill

\remark
Recall that a small deformation $D$ of a complex structure
is called {\bf versal} if any other small deformation can be obtained as a
pullback of the deformation $D$. Using the Kodaira--Spencer theory,
Kuranishi \index[persons]{Kuranishi, M.} has shown that the Kuranishi deformation space\index[terms]{deformation!versal}
is versal.\index[terms]{Kodaira--Spencer theory}

\subsection{Kuranishi to Teichm\"uller map}\index[terms]{map!Kuranishi to Teichm\"uller}

Deformation theory takes much effort to develop fully.
Fortunately, the only thing we need from the deformation
theory is the approximation-type results, and the
Teichm\"uller space, that is  easier to define, \index[terms]{space!Teichm\"uller}
is sufficient for this purpose.

\medskip

Recall that the {\bf $C^k$-topology}\index[terms]{topology!$C^k$}
on the space of sections of a bundle $B$
is the topology of uniform convergence 
of $b$, $\nabla b$, $\nabla^2 b$, ...,
$\nabla^k b$ on compacts, for some
connection $\nabla$ on $B$. \index[terms]{topology!$C^k$}
The {\bf $C^\infty$-topology} \index[terms]{topology!$C^\infty$}
is the topology of uniform convergence
of {\em all} derivatives. In other words,
a set is open in the $C^\infty$-topology
if it is open in all $C^k$-topologies.\index[terms]{topology!$C^k$}

\medskip

\definition\label{_Teich_Definition_}
Let $\Comp$ be the set of all integrable complex
structures on $M$, equipped with the $C^\infty$-topology,
and $\Diff_0$ the group of isotopies, that is, the
connected component of the group of diffeomorphisms.
{\bf The Teichm\"uller space} of complex structures
on $M$ is the quotient $\Teich:=\frac{\Comp}{\Diff_0}$
equipped with the quotient topology.

\medskip

\definition\label{_Kur_to_Teich_Definition_}
Let $X$ be a complex manifold, and ${\cal X} \stackrel
\phi \arrow U$
its Kuranishi deformation, that is  fibred
over the Kuranishi space $U \subset H^1(X, TX)$.
{\bf The 
Kuranishi to Teichm\"uller map}\index[terms]{map!Kuranishi to Teichm\"uller}
takes each $u \in U$ to the point of the Teichm\"uller
space of $X$ corresponding to $\phi^{-1}(u)$.

\medskip

The Kuranishi to Teichm\"uller map was introduced by L. \index[persons]{Meersseman, L.} Meersseman
(\cite{_Meersseman:stacks_,_Meersseman:Kuranishi_}).
In the ideal world, the  Kuranishi to Teichm\"uller map
should be an open embedding, but in practice this is not
always so. For hyperk\"ahler and Calabi--Yau manifolds,
the Kuranishi to Teichm\"uller map is always\index[terms]{manifold!hyperk\"ahler}
an open embedding (\cite{_Catanese:Moduli_}),\index[terms]{manifold!Calabi--Yau}
which allows one to equip the Teichm\"uller space
with a complex structure. \index[terms]{space!Teichm\"uller}

\medskip

Originally, the Teichm\"uller spaces were defined
by Teichm\"uller\index[persons]{Teichm\"uller, O.} to classify the complex structures
on Riemann surfaces. This notion was also useful
to prove the global Torelli theorem for hyperk\"ahler\index[terms]{theorem!global Torelli}
manifolds (\cite{_V:Torelli_},
\cite{_Verbitsky:ergodic_}). In this situation
the Teichm\"uller space is particularly nice.
For hyperk\"ahler manifolds, as well as for all Calabi--Yau,\index[terms]{manifold!Calabi--Yau}
the Teichm\"uller space is a smooth 
complex manifold (\cite{_Catanese:Moduli_}).
However, even for the hyperk\"ahler manifolds\index[terms]{manifold!hyperk\"ahler} 
the Teichm\"uller space is not Hausdorff. For 
other complex manifolds, the Teichm\"uller 
space is even worse. In \cite{_Meersseman:stacks_},
it was shown that the Teichm\"uller stack
for the \index[persons]{Hirzebruch, F.} Hirzebruch surface does not have\index[terms]{surface!Hirzebruch}
an underlying complex variety, and in fact
no two points are separable.

In \cite{_Meersseman:Kuranishi_},
it was shown that, when $M$ is K\"ahler, 
in a neighbourhood of a general point in $\Teich$, 
the  Kuranishi to Teichm\"uller map is a homeomorphism.\index[terms]{space!Kuranishi} \index[terms]{map!Kuranishi to Teichm\"uller}

In the next section, we compute the 
Kuranishi to Teichm\"uller map for
deformations of a general Hopf manifold,
and show that it is locally a 
homeomorphism (\ref{_Kuranishi_to_Teich_Hopf_Theorem_}).

\section{The Kuranishi space for  Hopf manifolds}
\label{_Kuranishi_Hopf_Section_}

\subsection{Vanishing of $H^2(TH)$ for a Hopf manifold}

We start from the following preliminary result,
that is  essentially due to \index[persons]{Haefliger, A.} Haefliger, \cite{_Haefliger:Hopf_}.

\hfill

\theorem\label{_H^2_TH_vanishes_Theorem_}
Let $H$ be a Hopf manifold. Then $H^2(TH)=0$.

\hfill

\proof
When $\dim H > 3$, this result immediately
follows from Mall theorem. Indeed, $TH$ is Mall
(\ref{_Mall_exa_Proposition_}), and $H^p(H, TH)=0$ for $1 < p < \dim H -1$
by \ref{_Mall_cohomology_Theorem_}.\index[terms]{theorem!Mall}

When $\dim H=3$, \index[persons]{Mall, D.} Mall theorem implies that \[ \dim H^2(H, TH)=\dim H^3(H, TH).\]
Serre duality gives $H^3(H, TH)\cong H^0(H, \Omega^1 H \otimes K_H)^*$,
and the latter group vanishes by \ref{_holo_forms_Hopf_vanish_Proposition_}.

When $\dim H=2$, we have $H^2(H, TH)\cong H^0(H, \Omega^1 H \otimes K_H)^*$
by Serre duality, and this group vanishes again,  by \ref{_holo_forms_Hopf_vanish_Proposition_}.\index[terms]{duality!Serre}
\endproof

\subsection{Kuranishi space and linear vector fields} \index[terms]{space!Kuranishi}

As explained in Section \ref{_Kuranishi_Section_},
the Kuranishi space for a complex manifold $X$
is the subset of $H^1(TX)$ consisting of all classes
$\gamma_0 \in H^1(TX)$ such that $M_1(\gamma_0)= M_2(\gamma_0) = ... =0$
(\ref{_M_i_Kuranishi_Claim_}).\index[terms]{vector field!linear}
Here, $M_i:\; H^1(TX)\arrow H^2(TX)$ are the higher Massey \index[terms]{Massey product}
products that are  defined on the zero set of the previous Massey products 
\[ \{ \gamma_0 \in H^1(TX)\ \ \ |\ \ \ M_1(\gamma_0)= M_2(\gamma_0) = ...= M_{i-1}(\gamma_0)=0\};
\]
the Massey product $M_i(\dash)$ is a homogeneous polynomial map of degree $i$.

For a Hopf manifold $H$, one has $H^2(TH)=0$
by \ref{_H^2_TH_vanishes_Theorem_},
hence the obstruction vanishes. Moreover, $H^1(TH)= H^0(TH)$.
The following theorem is essentially due to A. \index[persons]{Haefliger, A.} Haefliger
(\cite[Theorem A.3]{_Haefliger:Hopf_}).

\hfill

\theorem\label{_Kuranishi_on_Hopf_Theorem_}
Let $H$ be a Hopf manifold. Then 
the  \index[terms]{space!Kuranishi} Kuranishi space of $M$ is unobstructed; we can identify
it with a neighbourhood of zero in $H^1(H, TH)$. Moreover,
the dimension of $H^1(H, TH)$ is equal to the dimension of $H^0(H, TH)$,
and the identification $H^1(H, TH)=H^0(H, TH)$ is canonical
when $H$ is Vaisman.\index[terms]{space!Kuranishi}

\hfill

\proof 
The Kuranishi space of $H$ is unobstructed, because
obstructions to the deformation belongs to $H^2(H, TH)$,
and this space vanishes by \ref{_H^2_TH_vanishes_Theorem_}.

By \ref{_multi_cohomo_of_B_Theorem_} , 
whenever $H$ is Vaisman, the multiplication in cohomology
$H^0(H, TH)\times H^1(H, \C_H)\stackrel \mu \arrow H^1(H, TH)$ 
defines an isomorphism $\mu(\dash, \theta):\; H^0(H, TH)\arrow H^1(H, TH)$, 
where $\theta$ is the generator of $H^1(H, \C_H)=\C$. This implies that
the identification $H^1(H, TH)=H^0(H, TH)$ is canonical
when $H$ is Vaisman.
\endproof

\hfill

For a non-resonant linear Hopf manifold, the space $H^0(TH)$ can be
computed explicitly.\index[terms]{manifold!Hopf!linear}\index[terms]{manifold!Hopf!non-resonant}

\hfill

\definition\label{_linear_v_field_Definition_}
Let $H$ be a linear Hopf manifold, $H= \C^n \backslash 0/\langle A\rangle$,
where $A\in \GL(\C^n)$ is a linear contraction. A {\bf
  linear (polynomial) vector field}
on $\C^n$ is a vector field written as $\sum_i a_i \frac{\6}{\6z_i}$,
where $a_i$ are linear (polynomial) functions on $\C^n$ and $\frac{\6}{\6z_i}$ are 
coordinate vector fields. A linear vector field $\xi$ is called
{\bf $A$-invariant} if $A_*(\xi)=\xi$, where $A_*$ denotes the
action of the diffeomorphism $A\in \GL(C^n)$ on $T\C^n$.
A holomorphic vector field $\zeta$ on $H$ is called
{\bf linear} if its lift $\pi^{-1}(\zeta)$ to 
$\C^n\backslash 0\stackrel\pi \arrow H$ is linear.
The vector field $\pi^{-1}(\zeta)$ is clearly $A$-invariant.

\hfill

\remark\label{_linear_v_field_Remark_}
Let $V=\C^n$. The space of linear vector fields on $V$
coincides with $\End(V)$. To see this identification geometrically,
we take  $B \in \End(V)$ a point $x\in V$ and associate a vector $B(x)\in T_x V$ 
with $x$. This gives a linear vector field $\zeta_B$, and all linear vector fields
are obtained in this way. It is not hard to see that $\zeta_B$
is $A$-invariant if and only if $B$ commutes with $A$.

\hfill

\proposition\label{_global_vector_on_Hopf_Proposition_}
Let $H$ be a linear Hopf manifold, $H= \C^n \backslash 0/\langle A\rangle$,
where $A\in \GL(\C^n)$ is a linear contraction. 
Then 
\begin{description}
\item[(i)] $H^0(TH)$ coincides with the space of all 
polynomial vector fields on $\C^n$ commuting with $A$.
\item[(ii)] Moreover, if $H$ is non-resonant, $H^0(TH)$
is the space of all linear vector fields
commuting with $A$, and hence  $H^0(TH)=\{u\in \End(\C^n)\ \ |\ \ [A,u]=0\}$.
\end{description}

\hfill

\pstep
Let $\zeta\in H^0(TH)$ be a vector field on $H$, and
$\pi^{-1}(\zeta)$ its lift to $\C^n\backslash 0\stackrel\pi \arrow H$.
The vector field $\pi^{-1}(\zeta)$ is by construction $A$-invariant,
and any $A$-invariant holomorphic vector field defines a vector field in 
$H^0(TH)$. It remains only to show that any $A$-invariant
holomorphic vector field $\xi$ on $\C^n\backslash 0$ is polynomial.

By Hartogs theorem, $\xi$ is extended to a holomorphic
vector field on $\C^n$. Let $\xi = \sum_{d,i} P_{d,i}(z_1, ..., z_n)\frac{\6}{\6z_i}$
be the Taylor series for $\xi$, where $P_{d,i}$ are homogeneous polynomials
of order $d$. For every $d$, the degree $d$
homogeneous component $\sum_{i} P_{d,i}(z_1, ..., z_n)\frac{\6}{\6z_i}$
is also $A$-invariant. Since $H$ is compact, $H^0(TH)$
is finite-dimensional, and hence  the space of  $A$-invariant
holomorphic vector fields $\xi$ on $\C^n\backslash 0$ is
finite-dimensional. Then there are only finitely many
$d$ for which $\sum_{i} P_{d,i}(z_1, ..., z_n)\frac{\6}{\6z_i}$
is non-zero, and $\xi$ is therefore, polynomial. This proves
\ref{_global_vector_on_Hopf_Proposition_} (i).

\hfill

{\bf Step 2:} Now, assume that $H$ is non-resonant.
Let $\pi:\; \C^n \backslash 0\arrow H$ be the universal covering map.
We say that a vector field $X\in H^0(H, TH)$ {\bf is linear}
if $\pi^{-1}(X)$ is a linear vector field on $\C^n \backslash 0$.
By \ref{_Mall_flat_connection_Theorem_},  
$H$ is equipped with a flat connection $\nabla$. Moreover, $\nabla$ is
constant in the  flat coordinates, constructed in \ref{_Poincare_new_proof_Theorem_}.
Let $X\in H^0(H, TH)$ be a holomorphic vector field, and
$\nabla X$ its derivative. Clearly, $\nabla(\nabla X)=0$
if and only if $X$ is linear. However, if $\nabla (\nabla X)\neq 0$,
it gives a section of $\Omega^1 H \otimes \Omega^1 H \otimes TH$,
that has no sections because $H$ is non-resonant (\ref{_resonant_matrix_tangent_bundle_Corollary_}). 
\endproof

\subsection{Kuranishi to Teichm\"uller map for Hopf manifolds} \index[terms]{map!Kuranishi to Teichm\"uller}


For a general Hopf manifold $H$, we can describe the
 Kuranishi to Teichm\"uller map in a neighbourhood of $H$.
Surprisingly, this map is an isomorphism. 

Let $H=\frac{\C^n \backslash 0}{\langle \gamma\rangle}$ be a Hopf manifold.
The contraction $\gamma$ is extended to $\C^n$ by Hartogs.
A coordinate change in $\C^n$ takes the differential $D_0\gamma$ 
to a conjugate operator. Therefore, its eigenvalues remain
the same. This implies that the isotopy group maps
a non-resonant complex structure on a Hopf manifold 
to a non-resonant complex structure. Further on,
we shall speak of the Teichm\"uller space for
non-resonant Hopf manifolds, using this observation.
This space can be easily described using the
flat coordinates constructed in \ref{_Poincare_new_proof_Theorem_}.

\hfill

\proposition\label{_Teich_nr_for_Hopf_}
Let $\Teich_{nr}$ be the Teichm\"uller space\index[terms]{manifold!Hopf!non-resonant}
of non-resonant complex structures on a Hopf manifold $H$, 
$\dim_\C H =n \geq 3$. Consider the action of $\GL(n,\C)$
on the space $W$ of non-resonant linear invertible contractions of $\C^n$
by conjugation. Let $X=\frac{W}{\GL(n,\C)}$ be the orbit space,
equipped with the quotient topology.
Each non-resonant Hopf manifold $H = \frac{\C^n \backslash 0}{\langle \gamma\rangle}$\index[terms]{connection!flat}
is equipped with a flat connection, that is  unique (\ref{_Mall_flat_connection_Theorem_}) and gives
flat coordinates on its universal covering (\ref{_Poincare_new_proof_Theorem_}).
The map $\gamma$ is linear in these coordinates, and the
corresponding matrix is defined uniquely up to conjugation.
This defines a map $\tau:\; \Teich_{nr}\arrow X$.
Then $\tau$ is a homeomorphism.

\hfill

\proof
The map $\tau$ is bijective by construction.
The flat connection on $TH$ defined in \ref{_Mall_flat_connection_Theorem_}
can be obtained using the Banach fixed point theorem applied to the
space of all connections on $T\C^n$. Therefore, it\index[terms]{theorem!Banach fixed point}
continuously depends on the contraction $\gamma$. This implies that
for any converging sequence $\{I_i\}$ of non-resonant complex structures on $H$,
the corresponding sequence of connections will also converge,
producing a convergent sequence of linear contractions.
Hence, $\tau$ is continuous.
The inverse map is clearly continuous,
since a convergent sequence of linear contractions
defines a convergent sequence of Hopf manifolds.\footnote{``Convergence of Hopf manifolds''
is understood as the convergence of the corresponding
points in the Teichm\"uller space, \ref{_Teich_Definition_}.}
\endproof

\hfill

\remark
Let $H = \frac{\C^n \backslash 0}{\langle \gamma\rangle}$
be a resonant Hopf manifold. It is not hard to see that
the set of resonant contractions is closed, nowhere dense in the
set of all linear contraction matrices (Exercise \ref{_resonant_algebraic_locally_Exercise_}).\index[terms]{contraction!(non-)resonant}
This implies that any neighbourhood of $\gamma$ in the
set of holomorphic contractions contains non-resonant
contractions. Therefore, the space $\Teich_{nr}$
is open and dense in the Teichm\"uller space $\Teich$ of all
Hopf manifolds.

\smallskip

\remark\label{_non_separable_Remark_}
Let $\Teich_{lin}$ be the Teichm\"uller space of linear Hopf manifolds,
that is, of Hopf manifolds that can be obtained from linear contractions.
We say that two points $x, x'$ in a topological
space are {\bf non-separable} if any open neighbourhood
of $x$ intersects any open neighbourhood
of $x'$. Let $\{\gamma_i\}$ be a sequence of linear
contractions converging to $I\in \Teich$.
Then the eigenvalues of $\gamma_i$ also converge,
hence $\{\gamma_i\}$ is convergent\index[terms]{non-separable points}
inside $\Teich_{lin}$. Therefore,
each point $x$ in the complement $\Teich\backslash \Teich_{lin}$
is non-separable from $x'$ in $\Teich_{lin}$.

\hfill

\theorem\label{_Kuranishi_to_Teich_Hopf_Theorem_}
Let $H = \frac{\C^n \backslash 0}{\langle \gamma\rangle}$
be a non-resonant Hopf manifold, with $\gamma$ linear,
diagonal, and with all eigenvalues distinct.
Then the Kuranishi to Teichm\"uller map is \index[terms]{map!Kuranishi to Teichm\"uller}
a homeomorphism in a neighbourhood of $H$.

\hfill

{\bf Proof:} Since $\gamma$ is diagonal, it is Vaisman.
By \ref{_Kuranishi_on_Hopf_Theorem_}, the Kuranishi space \index[terms]{space!Kuranishi}
is naturally identified with a neighbourhood of 0 in 
$H^0(H, TH)$. By \ref{_global_vector_on_Hopf_Proposition_},
$H^0(H, TH)$ can be identified with the
space $S_\gamma$ of all matrices $A$ commuting with $\gamma$.
Since $\gamma$ is diagonal with distinct eigenvalues,
$S_\gamma$ is the space of all matrices that are  diagonal
in the same basis. \index[terms]{space!Kuranishi}

Now, the non-resonant Teichm\"uller space $\Teich_{nr}$ is an open subset of the space of 
conjugation orbits of $\GL(n,\C)$ (\ref{_Teich_nr_for_Hopf_}).
To finish the proof of
\ref{_Kuranishi_to_Teich_Hopf_Theorem_},
it remains to show that any 
conjugation orbit in a neighbourhood of $\gamma$
has a unique representative commuting with $\gamma$.
The following lemma finishes the proof of \ref{_Kuranishi_to_Teich_Hopf_Theorem_}. 

\hfill

\lemma
Let $\gamma\in \GL(n,\C)$ be a matrix with distinct
eigenvalues. Then there exists a sufficiently small 
neighbourhood $U\subset \GL(n,\C)$ of $\gamma$ 
and a sufficiently small neighbourhood $W$ of the 
class $[\gamma]$ in the space $\frac{\Mat(n,\C)}{\GL(n,\C)}$ of conjugation orbits
such that for any $[A] \in W$, there exists
a unique representative $A\in \Mat(n,\C)$ 
which commutes with $\gamma$.

\hfill

{\bf Proof:}
This is a simple linear algebra lemma, that can be proven
as follows.

Let $W\ni [\gamma]$ be an open neighbourhood of $[\gamma]$
in the space $\frac{\Mat(n,\C)}{\GL(n,\C)}$ of conjugation
orbits that does not contain any matrices with repeating
eigenvalues. We consider $W$ as a conjugation-invariant set
of matrices.

 Since all eigenvalues of
$A\in W$ are distinct, the conjugate orbit
of $A$ is the set of all matrices with the same eigenvalues.
A matrix commutes with $\gamma$ if and only if it
is diagonal in the same basis. Therefore, for each $A\in W$
there exists $n!$ matrices in its conjugate orbit
which commutes with $\gamma$. These matrices are
related by permuting the eigenvalues. However,
only one of these $n!$ representatives
can be found in a sufficiently small
neighbourhood $U$ of $\gamma$.
\endproof

\section{Hilbert schemes}\index[terms]{Hilbert scheme}
\label{_Hilbert_scheme_Section_}

In this section we briefly recall the construction of Hilbert scheme and
Hilbert variety. We need it further on to prove the smooth local triviality
of the explicit deformation of an LCK manifold with potential to the \index[terms]{manifold!LCK!with potential}
Vaisman manifold obtained in\index[terms]{manifold!Vaisman} \ref{_Vaisman_limit_of_LCK_pot_Theorem_}.
 For more details on Hilbert schemes and an introduction to moduli spaces
in algebraic geometry\index[terms]{geometry!algebraic}, please see \cite{_Harris:first_course_,_Viehweg_}.

Let $X\subset\C P^n$ be a projective variety,\index[terms]{graded ideal}
and $J^* \subset \bigoplus_i H^0(\C P^n, \calo(i))$
its graded ideal. One can define $J^*$ as the ideal
of all homogeneous polynomials $P(t_0, ..., t_n)\in \C[t_0, ..., t_n]$
vanishing in $X$. The {\bf Hilbert function} of $X$ is
 $F_X:\; \Z^{>0} \arrow \Z$ taking\index[terms]{Hilbert function}
$i$ to the dimension of $H^0(X, \calo(i)|_ X)$.
It is not very hard to show that the\index[terms]{Hilbert polynomial}
Hilbert function $k \mapsto F_X(k)$ is polynomial as a function of $k$
for $k \gg 0$ (Exercise \ref{_Hilbert_polyn_Exercise_}).
This polynomial is called {\bf the Hilbert polynomial} of $X$.

Hilbert basis theorem implies that $J^*$ is 
finitely generated as a module over $\bigoplus_i H^0(\C P^n,
\calo(i))$. Then \index[terms]{theorem!Hilbert's basis}
$J^*$ is generated by $J^k$ for $k$ sufficiently big.
Let $r$ be the dimension of $J^{k}\subset H^0(\C P^n,
\calo(k))$.

\hfill

\definition
In these assumptions, {\bf the Hilbert scheme} of
$X$ is the set of all $r$-dimensional subspaces 
$W\subset H^0(\C P^n, \calo(k))$ such that
the Hilbert polynomial of the ideal\index[terms]{Hilbert scheme}
generated by $W$ is equal to $F_X(k)$, for
$k$ sufficiently big.

\hfill

It is possible to show that the Hilbert scheme $S$ of
$X\subset \C P^n$ is an algebraic variety, and,
moreover, there exists a universal fibration
${\cal X} \subset X \times \C P^n$ projecting to 
$s \in S$ with the fibre given by the 
ideal $J_s$ that is  generated by the corresponding
subspace $W_s\subset H^0(\C P^n, \calo(k))$
(\cite[Theorem 21.8]{_Harris:first_course_}).
We think of the Hilbert scheme as of a variety
which parametrizes the flat\footnote{A family ${\cal X} \arrow S$
of subvarieties in $\C P^n$ parametrized by $s\in S$ is called
{\bf flat} if the Hilbert polynomial of the fibres is constant for
all $s\in S$. These fibres are called {\bf deformations},
and the whole family {\bf a deformation} of any of its fibres.}
deformations of $X \subset \C P^n$.\index[terms]{deformation!flat}

\hfill

Further on, we shall use the following proposition.

\hfill

\proposition\label{_Hilbert_smooth_locally_trivial_Proposition_}
Let $X\subset \C P^n$ be a projective orbifold, and
$S$ its Hilbert scheme. Denote by $S_0$ the set
of points of $S$ which correspond to orbifold
subvarieties $X'\subset \C P^n$, satisfying $X' = g(X)$
for some $g\in \PGL(n+1, \C)$. Then
\begin{description}
\item[(i)] $S_0$ is a quasi-projective subvariety of $S$.
\item[(ii)] Let $\pi:\; {\cal X}_0 \arrow S_0$
be the restriction of the universal family to $S_0 \subset S$.
Then $\pi$ is a locally trivial fibration.
\end{description}
\proof
By construction, $S_0\subset S$ is an orbit of $G$ acting on $S$.
Let $H\subset G=\PGL(n+1, \C)$ be the stabilizer of a point
$s\in S_0$. Then $S_0= G/H$; such a space is always quasi-projective
(see e.g. \cite[Theorem 1.16]{_Brion:intro_}).
To prove \ref{_Hilbert_smooth_locally_trivial_Proposition_} (ii),
take the unit $e\in H$, let $U\ni e$ be a neighbourhood, and $R\subset U$
a submanifold of complementary dimension and transversal to $H\subset G$
(this submanifold is usually called {\bf a slice} 
corresponding to the subgroup $H\subset G$).
Then the action of $G$ on $S_0$ gives a diffeomorphism
between $R$ and a neighbourhood of any given point $s\in S_0$,
hence it produces a diffeomorphism of $R \times X'$ and
a neighbourhood of the fibre $X'$ in the family ${\cal X}_0 \arrow S_0$.
\endproof

\section{The space of complex structures on LCK manifolds
  with potential}\index[terms]{manifold!LCK!with potential}

\subsection{The conjugation orbit of a linear operator}

Here we prove an explicit form of
Exercise \ref{_conjugate_of_Jordan_cells_Exercise_},
stating  that any conjugation orbit in $\GL(n, \C)$ contains  in its closure a semisimple
operator.\index[terms]{orbit!conjugation}

\hfill

\proposition\label{_semisimple_operator_approx_Proposition_}
Let $p\in \GL(n,\C)$ be a linear operator, and
$p=su$ its Jordan decomposition, with $s$ semisimple,
$u$ unipotent, and $su=us$. Then there exists\index[terms]{Jordan decomposition}
a sequence $r_i\in \GL(n,\C)$ of operators
commuting with $s$ and satisfying
$\lim_{i\to\infty} r_i p r_i^{-1}=s$.

\hfill

\proof
Since any operator is a sum of Jordan cells,\index[terms]{Jordan cell}
it would suffice to prove
\ref{_semisimple_operator_approx_Proposition_}
when $p$ is a single $k\times k$ Jordan cell,
{\[\small 
p =\begin{pmatrix} \alpha & 1 & 0 & \ldots & 0\\
0 & \alpha & 1 & \ldots & 0\\
\vdots &\vdots &\vdots & \cdots & \vdots \\
0&0&0 & \ldots &1\\
0&0&0 & \ldots &\alpha
\end{pmatrix}.
\]}
In this case, $s= \alpha\Id$, and hence 
it commutes with everything. Take
{\[\small 
r_i =\begin{pmatrix} 1 & 0 & 0 & \ldots & 0\\
0 & \epsilon_i & 0 & \ldots & 0\\
\vdots &\vdots &\vdots & \cdots & \vdots \\
0&0&0 & \ldots &0\\
0&0&0 & \ldots &\epsilon_i^{k-1}
\end{pmatrix}.
\]}
Then 
{\[\small 
r_i p r_i^{-1} =\begin{pmatrix} \alpha & \epsilon_i & 0 & \ldots & 0\\
0 & \alpha & \epsilon_i & \ldots & 0\\
\vdots &\vdots &\vdots & \cdots & \vdots \\
0&0&0 & \ldots &\epsilon_i\\
0&0&0 & \ldots &\alpha
\end{pmatrix}.
\]}
Taking a sequence $\epsilon_i$ converging to 0,
we obtain $\lim r_i p r_i^{-1}=s$.
\endproof 

\subsection{Diffeomorphism orbits of LCK structures with
potential have Vaisman limit points}\index[terms]{structure!LCK}

The following result, used later in this chapter, 
gives a more precise statement of the approximation 
 \ref{def_lckpot2Vai}
that we use throughout this book.

\hfill

\theorem\label{_Vaisman_limit_of_LCK_pot_Theorem_}
{(\cite{ov_lee})} 
Let $(M, J)$ be an LCK manifold with potential such that  $\dim_\C M \geq 3$.\index[terms]{manifold!LCK!with potential}
Then there exists a Vaisman-type complex structure $(M, J_\infty)$
such that the point $[J_\infty]$ in the Teichm\"uller
space $\Teich(M)$ of complex structures on $M$ belongs to the 
closure of $[J]\in \Teich(M)$. In other words, there exists
a sequence of diffeomorphisms $\nu_i \in \Diff_0(M)$
such that $\lim_i \nu_i(J)=J_\infty$, where the limit is taken
with respect to the $C^\infty$-topology on 
the space $\Comp(M)$ of complex structures on $M$.

\hfill

\proof
Without restricting the generality, we may
assume that $(M,J)$ has LCK rank 1\index[terms]{rank!LCK}, and its
$\Z$-cover $\tilde M$ is K\"ahler (\ref{defor_improper_to_proper}).
Denote by $\tilde M_c$ the Stein completion of $\tilde M$\index[terms]{completion!Stein}
(\ref{potcon}).  
Fix an embedding of $(M,J)$ into the Hopf manifold 
$H= \frac{\C^n\backslash 0}{\langle A \rangle}$
such that $\tilde M_c$ is obtained as a closure of the
corresponding $\Z$-covering $\tilde M \subset \C^n\backslash 0$.

Let $A=us$ be the Jordan--Chevalley decomposition\index[terms]{Jordan--Chevalley decomposition}
for $A\in \GL(n, \C)$. By \ref{_semisimple_operator_approx_Proposition_},
there exists a sequence $A_i= u_i s= r_i A r_i^{-1}$ 
of operators conjugated to $A$ such that $u_i$ converges
to $\Id$. Denote by $(H, I_i)$ the Hopf manifold 
$(H, I_i):= \frac{\C^n\backslash 0}{\langle A_i \rangle}$.
Since $(H, I_i)$ are all naturally 
isomorphic to $H$, one obtains the embedding
$\phi_i:\; (M,J)\arrow (H, I_i)$.

Since the operators $A_i= u_i s$ converge 
to $s$, the sequence $I_i\in \Comp(H)$
converges to $I_\infty$, where
$(H, I_\infty):= \frac{\C^n\backslash 0}{\langle s \rangle}$.
Indeed, $(H, I_i)= (\C^n\backslash 0) / \langle A_i\rangle$,
and the sequence $\{A_i\in \GL(n, \C)\}$ converges to $s$.

Let
$\gamma$ be the generator of the monodromy\index[terms]{action!monodromy} acting on $\tilde M_c$, 
and $\phi:\; \tilde M_c\arrow \C^n$
the embedding making the following diagram
commutative
\[
\begin{CD}
\tilde M_c@>{\phi}>> \C^n\\
@V{\gamma}VV @V{A}VV\\
\tilde M_c@>{\phi}>> \C^n.
\end{CD}
\]
Consider the map $\phi_i := r_i \phi r_i^{-1}$
as an embedding from $\tilde M_c$ to $\C^n$
making the following diagram commutative
\[
\begin{CD}
\tilde M_c@>{\phi_i}>> \C^n\\
@V{\gamma}VV @V{A_i}VV\\
\tilde M_c@>{\phi_i}>> \C^n.
\end{CD}
\]
We use the same letter $\phi_i$ to denote
the embedding 
$(M, J) \hookrightarrow \frac{\C^n\backslash 0}{\langle A_i  \rangle}$ associated with $\phi_i$.
Since $\phi(\tilde M_c)$ is $s$-invariant, and $A_i$
converge to $s$, the
sequence $\phi_i\restrict{\tilde M_c}$ converges to 
$\phi\restrict {\tilde M_c}$, giving
an embedding 
$\frac{\tilde M}{\langle s\rangle} \hookrightarrow
(H,I_\infty)$. The  limit manifold 
$(M, J_\infty)=\frac{\tilde M}{\langle s\rangle}$ is of Vaisman type,
because it is embedded to a diagonal Hopf manifold.

The maps $\phi_i$ do not converge to $\phi_\infty$ smoothly,
because the sequence $\{r_i^{-1}\}$ is not bounded.
However, the sequence $\{(\phi_i(\tilde M_c), A_i)\}$
converges to $(\phi_\infty(\tilde M_c), A)$ as a sequence
of pairs 
\[ \text{(algebraic subvariety $Z\subset \C^n$,  an
	automorphism $\psi\in \Aut(Z)$);}
\] hence the
corresponding points in $\Teich$ also converge.
This is what we are going to show.

Consider $\tilde M_c$ as a closed algebraic cone \index[terms]{cone!algebraic}in $\C^n$, and\index[terms]{graded ring}
 $\calo_{\tilde M_c}$ as the graded ring of polynomial functions on $\tilde M_c$.

Let $S\subset \calo_{\C^n}$, $\dim S=m$, be a finite-dimensional space
generating the ideal of $\phi_1(\tilde M_c)$
(\ref{_cone_cover_for_LCK_pot_Theorem_}).
By Step 1 of the proof of 
\ref{_cone_cover_for_LCK_pot_Theorem_}, 
we may assume that all elements of $S$ are  
polynomials of degree less than $d$. Denote by $V\subset \calo_{C^n}$ the space of
polynomials of degree $\leq d$. Let $X\subset \Gr_m(V)$ be 
a subset of the Grassmannian\footnote{Interpreting $X$ as a piece of the
	relevant Hilbert scheme, we obtain that it is an
	algebraic subvariety in $\Gr_m(V)$.} of $m$-dimensional planes in $V$
consisting of all subspaces $W\subset V$ which generate
an ideal $J_W\subset \C[\C^n]$ with the following property:
the graded ring $\C[\C^n]/J_W$ is isomorphic to the graded 
ring $\calo_{\tilde M_c}$.  \index[terms]{Hilbert scheme}

The sequence
$\{\sigma_i :=\phi_i(\tilde M_c)\}$ corresponds to
points in $X$ converging to $\sigma_\infty:=\phi_\infty(\tilde M_c)$. 
This gives the $C^0$-convergence of the
submanifolds $\phi_i(\tilde M)$ to $\phi_\infty(\tilde M)$
in $\C^n \backslash 0$.
Indeed, consider the ``universal''
fibration over $X$ with the fibre over $W\in X$ 
being the algebraic cone\index[terms]{cone!algebraic} associated with the ideal
$J_W\subset \calo_{\C^n}$ generated by $W$. 
The associated open cone fibration has
smooth fibres. 

To finish the proof, we need to prove that
the manifolds $(M,J_i)=\frac{\phi_i(\tilde M)}{\langle u_i s\rangle}$
smoothly converge to 
$(M, J_\infty)=  \frac{\phi(\tilde M)}{\langle s\rangle}$.
This would follow if we prove that 
the corresponding cones in $\C^n$ converge smoothly in each
annulus $B_R \backslash B_r$  around 0 (we need to restrict to the annulus,
because the cone itself is singular around zero, and hence  it makes
no sense to speak of $C^\infty$-convergence unless we remove a
neighbourhood of the origin).
Then $(M, J_i)$ and $(M, J_\infty)$ are quotients
of the respective cones by $A_i$ and $A$-actions
respectively, and $A_i$ converges to $A$ in $\GL(n, \C)$.

However, the cones $\phi_i(\tilde M_c)$ are smooth in each annulus, and they
converge to $\tilde M_c$ in $C^0$-topology\index[terms]{topology!$C^0$} (or in the Hausdorff metric)
by construction. For smooth families of compact\index[terms]{metric!Hausdorff}
complex manifolds, the  $C^\infty$-convergence of their fibres is automatic.
To finish the proof, we replace the cone fibration 
over $X$
by the corresponding fibration of compact complex orbifolds,
which also converges to the central fibre. The fibres
of a locally trivial fibration of compact orbifolds
$C^\infty$-converge to the central fibre by Ehresmann's theorem.\index[terms]{theorem!Ehresmann's fibration}

We use the Hilbert scheme parametrizing algebraic 
cones in $\C^n$. Since $r_i$ commutes with $s$, all subvarieties 
$\phi_i(\tilde M_c)\subset \C^n$ are invariant with respect
to the action of $s$. By definition (\ref{_alge_cone_Definition_}),
$\phi_i(\tilde M_c)$ is a closed algebraic cone\index[terms]{cone!algebraic}. 
By \ref{_cone_is_total_space_Proposition_}, it is 
the cone of a projective variety $Z_i \subset \C P^{n-1}$.
In fact, this variety is an orbifold, because
$\tilde M_c$ can be identified with a cone over a quasi-regular \index[terms]{manifold!Sasaki!quasi-regular}
Sasakian manifold (\ref{_S^1_qr_orbifold_Proposition_}). 
Let $S$ be the Hilbert scheme of $Z_i$, and
$S_0\subset S$ a $\PGL(n, \C)$-orbit of $Z_i$.
Let ${\cal X}_0\arrow S_0$ denote the universal fibration over $S_0$,
with the fibre over $s\in S_0$ a projective orbifold associated with $s$
(Section \ref{_Hilbert_scheme_Section_}). By 
\ref{_Hilbert_smooth_locally_trivial_Proposition_},
the variety $S_0$ is quasi-projective, and the fibration
${\cal X}_0\arrow S_0$ is locally trivial in complex analytic category.
Therefore, the corresponding cone fibration parametrizing
the deformations of $\tilde M_c$ is also locally trivial.
The central fibre $\phi(\tilde M_c)= \lim_i \phi_i(\tilde M_c)$ also belongs to this
family. \ref{_Vaisman_limit_of_LCK_pot_Theorem_}
claims that the family  $\frac{\phi_i(\tilde M)}{\langle A_i\rangle}$ 
smoothly converges to $\frac{\tilde M}{\langle s\rangle}$;
we have just shown that the family of algebraic cones\index[terms]{cone!algebraic} 
$\phi_i(\tilde M)$ smoothly converges to $\phi(\tilde M)$.

To obtain the convergence of the corresponding
LCK manifolds, we notice that  $\lim_i A_i=s$,
hence $\lim_i \frac{\phi_i(\tilde M)}{\langle A_i \rangle} = 
\frac{\phi(\tilde M)}{\langle s \rangle}= (M, J_\infty)$.
\endproof

\hfill

\corollary
Let $(M,J)$, $\dim_\C M\geq 3$, be a compact complex manifold
admitting an LCK structure with potential\index[terms]{structure!LCK!with potential}, and
$J_\infty$ the Vaisman-type complex structure
on $M$ obtained as in \ref{_Vaisman_limit_of_LCK_pot_Theorem_}.
Then any Vaisman-type Lee form\index[terms]{form!Lee} on $(M, J_\infty)$
can be realized as the Lee form of an LCK
structure with potential on $(M, J)$.

\hfill

\proof Let $(M, I, \omega, \theta)$ be a Vaisman structure on $M$,
with $\omega= d^c \theta + \theta \wedge \theta^c$, and
$I_i$ a sequence of complex structures on $M$ converging
to $I$, such that all $(M,I_i)$ are isotopic to  $(M,J)$ 
as complex manifolds. Then the sequence
$\omega_i= I_i d (I_i)^{-1} (\theta) + \theta \wedge I_i(\theta)$
converges to 
\[ d^c \theta + \theta \wedge \theta^c=
I d I^{-1} (\theta) + \theta \wedge I\theta.
\]
Since positivity is an open condition, the (1,1)-form
$\omega_i$ is positive for $i$ sufficiently big.
Then $(M, I_i, \omega_i, \theta)$ is LCK with potential\index[terms]{manifold!LCK!with potential},
and $\theta$ its Lee form.\index[terms]{form!Lee} However, $I_i$ is mapped
to $J$ by an isotopy, that preserves the cohomology
class of $\theta$, and hence  $\theta$ is a Lee class\index[terms]{class!Lee} on $(M, J)$.
\endproof

\hfill

\remark
The same argument was used in \cite[\S 5.1]{_Apostolov_Dloussky_} 
to precisely describe the Lee cone of a non-Vaisman 
Hopf surface.\index[terms]{cone!Lee}\index[terms]{surface!Hopf}

\subsection[The Teichm\"uller space of LCK manifolds with potential]{The Teichm\"uller space of LCK manifolds\\ with potential}\index[terms]{manifold!LCK!with potential}
\label{_Teichmuller_LCK_with_pot_Subsection_}

The ``Teichm\"uller space of matrices up to conjugation''
can be defined as the quotient space $X:=\frac{\Mat(\C^n)}{\GL(\C^n)}$
equipped with the quotient topology. Let $X_{nrcntr}\subset X$ be the open subset of
$X$ consisting of non-resonant invertible contractions. It follows from 
\ref{_Mall_flat_connection_Theorem_} and \ref{_Poincare_new_proof_Theorem_} 
that this space is homeomorphic to the\index[terms]{manifold!Hopf!non-resonant}
Teichm\"uller space of non-resonant Hopf manifolds.\index[terms]{space!Teichm\"uller}

\ref{_semisimple_operator_approx_Proposition_}
can be interpreted as a statement about the non-Hausdorff points in the space
$X=\frac{\Mat(\C^n)}{\GL(\C^n)}$. Indeed, 
\ref{_semisimple_operator_approx_Proposition_} is equivalent
to ``for each non-semisimple matrix $A\in \Mat(\C^n)$, 
the corresponding point $[A]\in X$ is non-separable
from the semisimple part $A_{ss}$ of $A$''.\footnote{As in \ref{_non_separable_Remark_},
two  points $x, y\in X$ in a topological space are called
{\bf non-separable} if any open neighbourhood of $x$ intersects
an open neighbourhood of $y$.} \ref{_non_separable_Remark_}
interprets this observation as a statement about the
Teichm\"uller space of Hopf manifolds.

\ref{_Vaisman_limit_of_LCK_pot_Theorem_}
can be formulated using the same language, as follows.

\hfill

\corollary
Let $M$ be an LCK manifold with potential,\index[terms]{manifold!LCK!with potential}
and $\Teich$ the Teichm\"uller space of 
complex structures on $M$ admitting an LCK
metric with potential. Then each point
$x\in \Teich$ is non-separable from 
a point $y\in \Teich$ which corresponds
to a complex structure of Vaisman type.
\endproof

\section{Notes}
Although LCK structures\index[terms]{structure!LCK} are not stable to small deformations, Vaisman structures are preserved by some particular type of deformations. For example, one can perturb the Lee vector field of a compact Vaisman manifold\index[terms]{manifold!Vaisman} in the Lie algebra of the closure of the group generated by its flow, as in \ref{defovai}, to obtain a new Vaisman structure that is  quasi-regular.\index[terms]{manifold!Vaisman!quasi-regular} Other types of deformations of Vaisman structures are obtained by fixing the complex structure and the Lee vector field and changing the Lee form \index[terms]{form!Lee}(\cite{bel}).

Another type of deformation was introduced in \cite{_os:deform_}. Given a Vaisman structure $(M,I,g,\theta)$, $\dim_\C M=n$, we define  $(M,J_{t},g_{t})$   with $J_{0}=J$, $g_{0}=g$ and $t$ sufficiently small, in such a way that the canonical foliation\index[terms]{foliation!canonical} is not affected and the deformation concerns only the transverse orthogonal complement and the transverse K\"ahler geometry.\index[terms]{geometry!transverse K\"ahler}  Namely, let $\zeta_t$ be a family of basic 1-forms w.r.t the canonical foliation $\Sigma$ of $(M,I,g,\theta)$. Define 	
\begin{equation*}
	\begin{split}
	I_{t}&=I-\theta^\sharp\otimes \zeta _{t}-(I\theta^\sharp)\otimes (\zeta _{t}\circ I),\\
	g_{t}&=g+d\zeta _{t}\circ (I\otimes \mathrm{Id})  \label{g_t}
	+\zeta _{t}\otimes \theta ^{c}+\theta ^{c}\otimes \zeta _{t}
	+\zeta _{t}\otimes \zeta _{t}.
	\end{split}
\end{equation*}
 For $I_t$ to be integrable, it is necessary that $d\zeta_t$ be of type (1,1) w.r.t $I$. To assure that $g_t$ is positive definite, one needs 
 \begin{equation*}\label{estimare1}
 	\|\zeta_t\|<\frac 1{(2^{n+1}-1){n}^\frac{n-1}{2}}, \quad\|d\zeta_t\|<\frac 1{(2^{n+1}-1){n}^\frac{n-1}{2}}.
 \end{equation*}
To construct such a basic 1-form, let $\phi$ be a basic function, and $d^c_b=Id_bI^{-1}$, where $d_b$ is the basic exterior differential. Then, for $t$ small enough, $\zeta_t:=td^c_b\phi$ is a basic form with the required properties. One easily verifies that the $(I_t,g_t)$ is a Vaisman structure with $\theta_t=\theta$ and $I\theta_t=I\theta+d\zeta_t$.

 This type of deformation is related  to the deformations of second type on Sasakian manifolds, as defined in \cite{bel}.

\section{Exercises}

\subsection{Hilbert polynomials}

In the sequel, we 
assume that all graded rings $A^*$ are graded by $\Z^{\geq
  0}$. The component
$A^0$ is equal to the base field, and $A^*$ is generated by $A^1$
and $A^0$.

\index[terms]{graded ring}
\begin{enumerate}[label=\textbf{\thechapter.\arabic*}.,ref=\thechapter.\arabic{enumi}]

\item
Let $A^*$ be a finitely generated graded ring, and
\[
0 \arrow M_1^* \arrow M^*_2 \arrow M_3^*\arrow 0.
\]
an exact sequence of finitely generated graded
$A^*$-modules.  Define {\bf the Hilbert function}
of $M^*_i$ as $F_{M_i}(k):= \dim M^k_i$. \index[terms]{Hilbert function}
Assume that the Hilbert functions $F_{M_a}(k)$ and $F_{M_a}(k)$ are
polynomial for $k\gg 0$, and two distinct indices
$a, b \in \{1, 2, 3\}$. Prove that $F_{M_c}(k)$ 
is polynomial for $k\gg 0$, where $c$ is the third index.

\item\label{_cyclic_Hilbert_Exercise_}
Recall that an $A^*$-module is called {\bf cyclic}
if it is generated by one element. Suppose that
for any cyclic $A^*$-module,
the Hilbert function is $F_{M}(k):= \dim M^k$
is polynomial for $k\gg 0$. Prove that it is
polynomial for $k\gg 0$ for any finitely generated
graded $A^*$-module.

{\em Hint:} Use the previous exercise.

\item\label{_exact_regular_Hilbert_Exercise_}
Let $A^*$ be a finitely generated graded ring,
$M^*$ a torsion-free graded $A^*$-module, and $a\in A^1$ 
an element that is  not a zero divisor. Consider the
exact sequence
\[
0 \arrow M^* \xlongrightarrow{m \mapsto am} M^{*+1}
\arrow  \frac{M^{*+1}}{a M^*}\arrow 0.
\]
Denote the $A^*$-module $\frac{M^{*+1}}{a M^*}$ by $U$.
Suppose that the Hilbert function
$F_{U}(k)$ is polynomial for  $k\gg 0$.
Prove that the Hilbert function
$F_{M}(k)$ is polynomial for  $k\gg 0$.

\item
Define {\bf a regular sequence}\index[terms]{regular sequence}
of elements $a_1, ..., a_n$ in a ring $A^*$ as
a sequence such that for each $k$, the element $a_k$ is not a zero divisor in
$\frac{A^*}{(a_1, ..., a_{k-1})}$, where $(a_1, ..., a_{k-1})$ 
denotes the ideal generated by $a_1, ..., a_{k-1}$.\index[terms]{depth (of a ring)}
Define {\bf the depth} $\depth(A^*)$ of $A^*$ as the maximal length of a
regular sequence. Prove that the depth of a finitely--generated ring is finite.

{\em Hint:} One can use induction on the Krull dimension\index[terms]{dimension!Krull}
$\dim A^*$ to show that $\depth(A^*) \leq \dim(A^*)$.

\item 
Let $M^*$ be a finitely--generated graded $\C[t]$-module.
Prove that the Hilbert function
$F_{M}(k)$ is polynomial for  $k\gg 0$.

\item
Let $M^*$ be a finitely--generated graded 
$A^*$-module. Assume that the Hilbert function
$F_{U}(k)$ is polynomial for $k\gg 0$, for all finitely--generated
$B^*$-mo\-dules $U^*$, where $B^*$ is a  finitely--generated
ring with $\depth B^* < \depth A^*$. Prove that
the Hilbert function $F_{M}(k)$ is polynomial for $k\gg 0$.

{\em Hint:} Use Exercises \ref{_cyclic_Hilbert_Exercise_}
and \ref{_exact_regular_Hilbert_Exercise_}.

\item\label{_Hilbert_polyn_Exercise_}
Let $M^*$ be a finitely--generated graded 
$A^*$-module, where $A^*$ is a finitely--generated graded
ring.  Prove that
the Hilbert function $F_{M}(k)$ is polynomial for $k\gg 0$.

{\em Hint:} Use the previous exercise and induction on
$\depth(A^*)$.

\end{enumerate}

\subsection{Deformation theory}

\definition
A map $\delta:\; \Lambda^*(M) \arrow \Lambda^{*+l}(M)$
that satisfies $\delta(x\wedge y) = \delta(x) \wedge y +
(-1)^{\tilde x l}x\wedge \delta(y)$\index[terms]{graded derivation}
where $\tilde x$ is the parity\footnote{{\bf The parity} of an operator
	$\Lambda^*(M) \arrow \Lambda^{*+l}(M)$ is $0$ when $l$ is even, and $1$ when $l$ is odd.} of $x$ is called {\bf a graded derivation} of the de Rham algebra.
An operator 
$\delta:\; \Lambda^*(M) \arrow \Lambda^{*+l}(M)$
{\bf supercommutes with the de Rham differential}
if $\delta d = (-1)^{l} d \delta$.

\begin{enumerate}[label=\textbf{\thechapter.\arabic*}.,ref=\thechapter.\arabic{enumi}]
\setcounter{enumi}{7}

\item\label{_deriva_determined_uniquely_Exercise_}
Let $\delta:\; \Lambda^*(M) \arrow \Lambda^{*+l}(M)$
be a graded derivation of  the de Rham algebra
which supercommutes with the de Rham differential.
Prove that $\delta$ is determined by its values on $\Lambda^0(M)= C^\infty (M)$.

\item\label{_deriva_de_Rham_Exercise_}
Let $\delta:\; \Lambda^*(M) \arrow \Lambda^*(M)$ be
a derivation which commutes with the de Rham differential and preserves the grading.
Prove that there exists a vector field $X$ on $M$
such that $\delta$ is the Lie derivative, $\delta= \Lie_X$.

{\em Hint:} Use the previous exercise.

\item Prove that any $\Lambda^{l}(M)$-valued derivation
$\delta:\; C^\infty(M) \arrow \Lambda^{l}(M)$ can be
extended to a unique graded derivation 
$\delta:\; \Lambda^*(M) \arrow \Lambda^{*+l}(M)$.

\item\label{_deriva_de_Rham_forms_Exercise_}
Let $\eta \otimes v\in \Lambda^l (M) \otimes TM$. Define a map
 $C^\infty M\xlongrightarrow {\rho(\eta \otimes v)} \Lambda^l (M)$ taking a function $f$
to $\eta \cdot \Lie_v f$.
Let $\delta:\; C^\infty(M) \arrow \Lambda^{l}(M)$
be a derivation. Prove that there exists
a section  $u \in \Lambda^l (M) \otimes TM$
such that $\delta(f) = \rho(u)(f)$.

\item\label{_supercommutator_Exercise_}
Let $\delta_1, \delta_2$ be two graded derivations
of the de Rham algebra shifting the grading by $p$ and $q$.
{\bf The supercommutator} $\{\delta_1, \delta_2\}$
is defined as $\delta_1 \delta_2 - (-1)^{\tilde \delta_1 \tilde \delta_2} \delta_2 \delta_1$,\index[terms]{supercommutator}
where, as elsewhere, $\tilde \delta_i$ denotes the parity of $\delta_i$.
Prove that the supercommutator of two graded derivations is a graded derivation.

\item\label{_almost_complex_str_derivation_Exercise_}
Let $I$ be an almost complex structure. Consider $I$ as a section
of $\Lambda^1 (M) \otimes TM$, and let $\delta_I:\; \Lambda^* (M)\arrow \Lambda^{*+1} (M)$
be the corresponding derivation, defined as in Exercise \ref{_deriva_de_Rham_forms_Exercise_}. \index[terms]{operator!Weil}
Denote by $W_I$ {\bf the Weil operator} taking
a $(p, q)$-form $\eta$ to $\1^{p-q}\eta$.
\begin{enumerate}
\item Prove that $W_I$ is a derivation.
\item Prove that the commutator
$[W_I, d]=\{d, W_I\}$ is a derivation that is  equal to $d^c= Id I^{-1}$ on functions.
\item Prove the graded Jacobi identity: \[ \{x, \{y, z\}\}=
  \{\{x, y\}, z\}+ (-1)^{\tilde x \tilde y} \{y, \{x, z\}\}.\]
Using the graded Jacobi identity, prove that $\{d, \{d, W_I\}\}=0$. 
\item Prove that $\delta_I=[W_I, d]$.\index[terms]{graded Jacobi identity}
\end{enumerate}
{\em Hint:} Use Exercise \ref{_deriva_determined_uniquely_Exercise_}.

\item
Let $I\in \End(TM)$ be an almost complex structure on a manifold,
and $\delta_I:\; \Lambda^* (M)\arrow \Lambda^{*+1} (M)$ the corresponding odd
derivation (Exercise \ref{_almost_complex_str_derivation_Exercise_}). 
Consider the anticommutator $\{ \delta_I, \delta_I\}:\;  \Lambda^* (M)\arrow \Lambda^{*+2} (M)$.
\begin{enumerate}
\item Prove that $\{ \delta_I, \delta_I\}$ is a derivation which corresponds
to a $TM$-valued 2-form $\nu \in \Lambda^2 M \otimes TM$ as in Exercise \ref{_deriva_de_Rham_forms_Exercise_}.
\item Define {\bf the \index[persons]{Nijenhuis, A.} Nijenhuis tensor} of an almost complex manifold
as a $TM$-valued 2-form given by $N(x, y):= [Ix, Iy]- I[Ix, y]-I[x, Iy]-[x, y]$.
Prove that $N(x,y)$ is linear in $x, y$, and $N$ vanishes if and only if $I$ is integrable.\index[terms]{tensor!Nijenhuis}
\item Prove that the $TM$-valued 2-form $\nu$ defined above 
is equal to the Nijenhuis tensor $N$ of $I$.
\end{enumerate}


%
%


\item
Let $r_1, ... r_n$ be a collection of positive numbers.
Prove that there exists an open subset $U\subset (\R^{>0})^n$
containing $(r_1, ..., r_n)$, and a number $d\in \Z$
such that for $(s_1, ..., s_n)\in U$ and any relation of form
$s_l= \sum_{i=1}^n m_i s_i$, $m_i \in \Z^{\geq 0}$,
the numbers $m_i$ satisfy $\sum_{i=1}^n m_i \leq d$.

\item
Let $A \in \GL(n, \C)$ be a linear contraction.
Prove that there exists an open subset $U \ni A$
and a number $d\in \Z$ 
such that for any $B\in U$, with eigenvalues 
$b_1, ..., b_n$, and equality of form $b_l =\prod_{i=1}^n b_i^{m_i}$, $m_i \in \Z^{\geq 0}$
the numbers $m_i$ satisfy $\sum_{i=1}^n m_i \leq d$.

{\em Hint:} Apply the previous exercise to $s_i = - \log |b_i|$.

\item\label{_resonant_algebraic_locally_Exercise_}
Let $A \in \GL(n, \C)$ be a linear contraction.
Prove that there exists an open subset $U \ni A$
such that the set of resonant matrices in $U$
is a finite union of closed algebraic subvarieties
of positive codimension.

{\em Hint:} Use the previous exercise.

\item
Let $\Comp$ be the space of complex structures on $S^2$.
\begin{enumerate}
\item Prove that the group $\Diff(S^2)$ of diffeomorphisms acts on $\Comp$
transitively. Prove that $\Comp= \frac{\Diff(S^2)}{\PSL(2, \C)}$.
\item Let ${\goth R}$ be the space of Riemannian structures on $S^2$.
Prove that $\Comp$ is diffeomorphic to the union
of two copies of $\frac{\goth R}{\C^\infty_{>0} M}$,
where the group ${\C^\infty_{>0} M}$ of positive real
functions acts on ${\goth R}$ by conformal changes.
\item Prove that ${\goth R}$ is contractible.
\item Prove that $\Comp$ is homeomorphic to a union of two contractible topological spaces.
\item Prove that $\Diff(S^2)$ is homotopy equivalent
to the union of two copies of $\PSL(2, \C)$.
\end{enumerate}

\item
Let $T^{2n}$ be a compact oriented torus of dimension $2n$, 
and $\Comp$ the space of oriented complex structures of K\"ahler type on $T^{2n}$. Let $I \in \Comp$. 
Denote by $\Diff_0$ the group of isotopies of $T^{2n}$. Prove that:
\begin{enumerate}
\item  The Albanese map \index[terms]{map!Albanese}
\[
(T^{2n}, I)\arrow \frac{H^{1,0}(T^{2n})^*}{H^1(T_{2n},\Z)} =\frac{H_1(T^{2n}, \R)}{H_1(T^{2n},\Z)}
\]
 can be used to construct
a canonical torsion-free flat connection $\nabla$ on $(T^{2n}, I)$. Prove that\index[terms]{connection!torsion-free}\index[terms]{connection!flat}
$I$ is $\nabla$-invariant.
\item Each connected component of 
$\frac{\Comp}{\Diff_0}$ is diffeomorphic to $\frac{\GL^+(2n, \R)}{\GL(n,\C)}$ (here $\GL^+(2n, \R)$ is the connected component of the identity in $\GL(2n, \R)$).
\item  $\frac{\GL^+(2n, \R)}{\GL(n,\C)}$ is Hausdorff.
\end{enumerate}

\item
Let $\Teich=\frac{\GL^+(2n, \R)}{\GL(n,\C)}$ 
be the Teichm\"uller space of K\"ahler-type complex structures on 
a complex torus $T$, $\dim_\C T =n > 1$,
and  ${\cal M}:= \frac{\Teich}{\Diff}$, where $\Diff$ is the group
of diffeomorphisms.
\begin{enumerate}
\item Prove that ${\cal M}$ is homeomorphic to the double coset space
$${\cal M}\simeq \GL(2n, \Z) \backslash \GL^+(2n, \R)/\GL(n,\C).$$
\item Prove that ${\cal M}$ is non-Hausdorff.
\end{enumerate}

\item Recall that two points
$x, y\in M$ in a topological space are called {\bf non-separable}
if any open neighbourhoods of $x$ and $y$ intersect.
\begin{enumerate}
\item Prove that this relation is not always transitive.
\item Let $\sim$ be an equivalence relation generated 
by the non-separability. Prove that the quotient space
$M/\!\!\sim$ is not always Hausdorff.
\end{enumerate}

\end{enumerate}

\definition\label{_DGA_Definition_}
{\bf A differential graded algebra} is a graded commutative algebra $A^*$
(\ref{_graded_commu_Definition_}) equipped with a differential 
$d:\; A^*\arrow A^{*+1}$ that satisfies $d^2=0$ and $d(x\wedge y)= dx \wedge y + (-1)^{\tilde x} x \wedge dy$,\index[terms]{differential graded algebra}
where $\tilde x$ denotes the parity of $x$. A homomorphism of 
differential graded algebras is called {\bf a quasi-isomorphism}
if it induces an isomorphism on cohomology. \index[terms]{quasi-isomorphism} 
Two differential graded algebras $A^*, B^*$  are called {\bf quasi-isomorphic}
if there exists a chain of quasi-isomorphisms $A^* \arrow A_1^*$,
$A_2^* \arrow A_1^*$,  $A^*_2 \arrow A_3^*$, $A^*_4 \arrow A_3^*$, ...,
starting from $A$ and ending with $B$. A differential graded
algebra is called {\bf formal} if it is quasi-isomorphic to its cohomology algebra.\index[terms]{differential graded algebra!formal} A manifold is {\bf formal} if its de Rham algebra is formal.

\begin{enumerate}[label=\textbf{\thechapter.\arabic*}.,ref=\thechapter.\arabic{enumi}]
\setcounter{enumi}{21}

\item
Let $M$ be a compact K\"ahler manifold. 
\begin{enumerate}
\item
Denote by $\ker d^c$ the
differential graded algebra $(\ker d^c, d)$ of $d^c$-closed differential forms. Using
the $dd^c$-lemma, prove that the embedding $(\ker d^c, d)\hookrightarrow (\Lambda^*M, d)$
is a quasi-isomorphism.
\item Prove that the natural homomorphism
$(\ker d^c, d)\arrow (H^*_{d^c}(M), d)$ is a quasi-isomorphism, where
$(H^*_{d^c}(M), d)$ is the algebra of the cohomology of $d^c$ with the differential $d$.
\item Prove that $(H^*_{d^c}(M), d)$ is isomorphic to the cohomology algebra of $M$
with zero differential, and $M$ is formal.
\end{enumerate}

\item
Let $P$ be a manifold, and $\alpha, \beta, \gamma$ 
three closed forms that satisfy $[\alpha \wedge \beta]=0$
and $[\beta\wedge \gamma]=0$ in the cohomology algebra. Consider a differential form
$a$ such that $da= \alpha \wedge \beta$ and $b$ such that
$db=\beta\wedge \gamma$. Then the form 
$M_1(\alpha, \beta, \gamma):=a \wedge \gamma - (-1)^{\tilde \alpha} \alpha \wedge b$ is
closed.
\begin{enumerate}
\item
 Prove that the  cohomology class of $M_1$ is well-defined modulo
$\im L_\gamma + \im L_\alpha$, where $L_\gamma, L_\alpha\in \End(H^*(P))$ are the 
operators of exterior multiplication by $\gamma$ and $\alpha$.
\item We consider the Massey product $M_1$ as a map from the set 
\[ \{\alpha, \beta, \gamma \in H^*(P)\ \ \ |\ \ \ \alpha \wedge \beta = \beta\wedge \gamma=0\}
\]
to $\frac{H^*(P)}{\im L_\gamma + \im L_\alpha}$.
Prove that $M_1=0$ when $P$ is formal.\index[terms]{Massey product}
\end{enumerate}

\item\label{_Kodaira_de_Rham_Exercise_}
Let $M$ be a \index[terms]{surface!Kodaira} Kodaira surface. Prove that its de Rham algebra
is quasi-isomor\-phic to the algebra $A^*= (\Lambda^*(V), d)$, where $V$
is the 4-dimensional space generated by $x, y, z, t$, and
the only non-zero differential is $d(z) = x\wedge y$.\index[terms]{surface!Kodaira}

{\em Hint:} Construct a nilmanifold \index[terms]{nilmanifold}diffeomorphic to the Kodaira surface.\index[terms]{surface!Kodaira}
Use \ref{kod_nil} to compute the cohomology
of left-invariant differential forms on $M$.
Prove that this algebra is quasi-isomorphic to $\Lambda^*(M)$.

\item 
Let $M$ be a Kodaira surface, and $A^*= (\Lambda^*(V), d)$
the differential graded algebra constructed in the previous exercise.
Prove that $M_1(y, x, x)= - z \wedge x$, and show that 
the corresponding class in $\frac{H^*(M)}{\im L_x + \im L_y}$
is non-zero. Prove that the \index[terms]{surface!Kodaira} Kodaira surface is not formal.

\end{enumerate}


\chapter{The set of Lee classes on LCK manifolds with potential}\index[terms]{manifold!LCK!with potential}\index[terms]{class!Lee}
\label{_Lee_classes_Chapter_}

\epigraph{\it Qui sait si la plupart de ces pens\'ees prodigieuses sur lesquelles tant de grands hommes, et un infinit\'e de petits, ont p\^ali depuis des si\`ecles, ne sont pas des monstres psychologiques, -- des {\em Id\'ees Monstres} --,  enfant\'ees par l'exercise na\"if de nos facult\'es int\'erogeantes que nous appliquons un peu partout, -- sans nous aviser que nous ne devons raisonnablement questionner que ce qui peut v\'eritablement nous r\'epondre?\\[10pt]

Who knows whether those of the prodigious thoughts over which so many great men and an infinity of lesser ones have grown pale for centuries are not, after all, psychological monsters -- {\em Monster Ideas} -- born of the na\"ive exercise of our questioning faculties, that we apply to anything at all, never realizing that we may reasonably question only what could give us an answer?}{\sc \scriptsize  Paul Val\'ery, \ \  M. Teste (Pr\'eface)\\ (translated by Jackson Mathews)}


\section{Introduction}

The first cohomological invariant one encounters when
dealing with the LCK manifolds is {\bf the Lee class},
the cohomology class of the Lee form.\index[terms]{form!Lee} Let\index[terms]{class!Lee} 
$(M, \theta,\omega)$ be a compact LCK manifold,
and $[\theta] \in H^1(M, \R)$ its Lee class.
By Vaisman's theorem (\ref{vailcknotk}), $[\theta]=0$ if
and only if $M$ is of K\"ahler type.

In this chapter, we are trying to find more ways to
relate the Lee class to the geometry of an LCK manifold.

The first question one would ask is: what is the set 
of Lee classes on a given compact complex manifold? 
A similar question about the set of K\"ahler classes\index[terms]{class!K\"ahler}
yields the notion of the ``K\"ahler cone'', that is \index[terms]{cone!K\"ahler}
one of the most important geometric features of 
a K\"ahler manifold. The set of Lee classes is
not called ``the Lee cone'' because it is not a 
cone; indeed, as shown by  V. \index[persons]{Apostolov, V.} Apostolov and G. \index[persons]{Dloussky, G.} Dloussky (\cite{ad2}), and A. \index[persons]{Otiman, A.} Otiman\index[terms]{cone!Lee}
(\cite[Theorem 3.11]{oti2}), for an
\index[terms]{surface!Inoue} Inoue surface of class $S^0$, the set
of Lee classes is a point. 

For locally conformally symplectic structures,
the set of Lee classes\index[terms]{class!Lee} is better understood, due to
\index[persons]{Eliashberg, Y.} Eliashberg and \index[persons]{Murphy, E.} Murphy, who proved that 
for any almost complex manifold, and
any non-zero class $\alpha\in H^1(M, \Q)$, 
there exists $C>0$ such that $C\alpha$ is a
Lee class of an LCS structure\index[terms]{structure!LCS} (\cite[Theorem 1.11]{em}).

In a recent preprint \cite{_Bertelson_Meigniez_},
\index[persons]{Bertelson, M.} Bertelson and \index[persons]{Meigniez, G.} Meigniez have made this result much
stronger, showing that any non-zero class $\alpha\in
H^1(M, \R)$ can serve as the Lee class of an LCS structure.
Moreover, there exists an LCS form\index[terms]{form!LCS} $\omega$ in any connected component
of the set of non-degenerate 2-forms with any prescribed
non-zero Lee class $[\theta]$ and any 
prescribed Morse--Novikov cohomology\index[terms]{cohomology!Morse--Novikov}
class $[\omega]\in H^2_\theta(M)$.

For complex surfaces with $b_1(M)=1$, the set ${\goth L}$ of 
Lee classes of LCK structures\index[terms]{structure!LCK} was studied by\index[terms]{class!Lee}
\index[persons]{Apostolov, V.} Apostolov and \index[persons]{Dloussky, G.} Dloussky, who proved that
${\goth L}$ is either open or a point.

The first important advance in this direction 
is due to K. \index[persons]{Tsukada, K.} Tsukada, who proved that the set of Lee classes on 
Vaisman manifolds is a half-space (\cite[Theorem 5.1]{tsu}), using the harmonic\index[terms]{manifold!Vaisman}
decomposition of 1-forms on Vaisman manifolds, due to T. \index[persons]{Kashiwada, T.} Kashiwada.

In this chapter we give a new proof of Tsukada's theorem, and extend
it to LCK manifolds with potential.\index[terms]{manifold!LCK!with potential}

We use the following decomposition theorem, that is  valid
for all compact LCK manifolds with potential:
\begin{equation}\label{_H^1_decompo_Equation_}
H^1(M, \C) = H_d^{1,0}(M) \oplus \overline{H_d^{1,0}(M)} \oplus \langle \theta \rangle,
\end{equation}
where $H_d^{1,0}(M)$ is the space of all closed holomorphic 1-forms, and
$\theta$ the Lee form.\index[terms]{form!Lee} \index[persons]{Tsukada, K.} Tsukada proved this
for Vaisman manifolds using the commutation formulae for Laplacians,\index[terms]{manifold!Vaisman}
and using the harmonic decomposition for Vaisman manifolds; we gave
a new proof of the decomposition theorem for harmonic 1-forms in \ref{_harmo_deco_1-form_Proposition_}.\index[terms]{decomposition!harmonic}

I. \index[persons]{Vaisman, I.} Vaisman conjectured that $b_1(M)$ is odd-dimensional for
any compact LCK manifold (\ref{_Vaisman_conjecture_disproved_}); 
this famous conjecture was disproven by
Oeljeklaus and \index[persons]{Toma, M.} Toma in \cite{ot}. The decomposition\index[terms]{conjecture!Vaisman}
\eqref{_H^1_decompo_Equation_} would imply that $b_1(M)$ is odd, and hence 
the counterexample of \index[persons]{Oeljeklaus, K.} Oeljeklaus-\index[persons]{Toma, M.}Toma does not satisfy 
\eqref{_H^1_decompo_Equation_}. However, the natural map
\begin{equation}\label{_H^1_holom_map_Equation_}
H_d^{1,0}(M) \oplus \overline{H_d^{1,0}(M)} \oplus \langle \theta \rangle
\arrow H^1(M, \C)
\end{equation}
is always injective (\ref{_H^1_holo_LCK_Lemma_}).

For LCK manifolds with potential\index[terms]{manifold!LCK!with potential}, we deduce 
\eqref{_H^1_decompo_Equation_} from a deformation argument,
by showing that an LCK manifold with potential $M_1$ obtained as a deformation
of a Vaisman manifold $M_2$ satisfies \index[terms]{manifold!Vaisman}
$\dim H_d^{1,0}(M_1)\geq \dim H_d^{1,0}(M_2)$, that follows
from the semicontinuity of Hodge numbers 
(\cite[Theorem 2.3]{_Bell_Narasimhan_}).
Unless $\dim H_d^{1,0}(M_1)= \dim H_d^{1,0}(M_2)$, this would
imply that $\dim H_d^{1,0}(M_1) > \frac{b_1(M)-1}{2}$,
that is  impossible because the map \eqref{_H^1_holom_map_Equation_}
is injective by \ref{_H^1_holo_LCK_Lemma_}.

Notice that the equality $\dim H_d^{1,0}(M) = \frac{b_1(M)-1}{2}$
is valid for non-K\"ahler complex surfaces as well. 

The above decomposition \eqref{_H^1_decompo_Equation_}
is the cornerstone for the description of the set of Lee classes\index[terms]{class!Lee}
on an LCK manifold with potential.\index[terms]{manifold!LCK!with potential} Indeed, consider the linear map
$\mu:\; H^1(M, \R)\arrow \R$ that vanishes 
on the codimension 1 subspace 
$H_d^{1,0}(M) \oplus \overline{H_d^{1,0}(M)}\subset H^1(M, \R)$
and is positive on the Lee form.\index[terms]{form!Lee} We prove that
$\alpha \in H^1(M, \R)$ is a Lee class
if and only if $\mu(\alpha) >0$
(\ref{_Lee_cone_on_LCK-pot_Theorem_}).

The description of the set of Lee classes
on LCK manifold with potential uses the
embedding to a Hopf manifold, that is 
valid in complex dimension $\geq 3$.
Conjecturally, it is also true in 
dimension 2 (\ref{_sphe_implies_dim2_Theorem_}).

On Vaisman manifolds, this\index[terms]{manifold!Vaisman}
result is valid in all dimension including
dimension 2, because the Vaisman surfaces are
classified (\cite{bel}), and in this case,
\ref{_sphe_implies_dim2_Theorem_} can be applied.

\section{LCK metrics on Vaisman manifolds}

We start with the following preliminary result 
which might be of separate interest.

\hfill

\proposition\label{_Lee_form_on_Vaisman_is_Vaisman_Proposition_} { (\cite{ov_lee})} 
Let $(M,\theta, \omega)$ be an LCK structure\index[terms]{structure!LCK}
on a compact Vaisman manifold. Then $\theta$\index[terms]{manifold!Vaisman}
is cohomologous to a Lee form\index[terms]{form!Lee} of a Vaisman structure.\index[terms]{structure!Vaisman}

\hfill

\pstep
Let $X$ be the Lee field\index[terms]{Lee field} of a Vaisman structure $(M, \omega^V, \theta^V)$\index[terms]{structure!Vaisman}
on $M$, and $G$ the closure of the group generated by 
exponents of $X$ and $I(X)$. Since $X$ and
$IX$ are Killing\index[terms]{vector field!Killing} and commute, $G$ is a compact
commutative Lie group, and hence  it is isomorphic to a compact torus.
This group acts on $M$ by holomorphic isometries
with respect to the Vaisman metric.

Averaging $\theta$ with the $G$-action, we obtain
a $G$-invariant 1-form $\theta^G$, corresponding
to another LCK structure \index[terms]{structure!LCK}in the same conformal class,
and cohomologous with $\theta$.
Without restricting the generality, we may assume
from the beginning that the form $\theta$ is $G$-invariant. 

Now, the equation
$d\omega=\omega\wedge\theta$ is invariant under the
action of $G$, because $\theta$ is $G$-invariant;
in other words,
$d(g^*\omega)= g^*\omega\wedge \theta$, for all $g\in G$.
This implies that
$\omega$ averaged with $G$ gives a form
$\omega^G$ that satisfies
$d(\omega^G)= \omega^G\wedge \theta$.
We have constructed a $G$-invariant
LCK structure\index[terms]{structure!LCK} $(M, \omega^G, \theta)$. 
Replacing $\omega$ with $\omega^G$,
we may assume that $(M, \omega, \theta)$
is $G$-invariant. We are going to prove
that $(M, \omega, \theta)$ is Vaisman.

\hfill

{\bf Step 2:} 
Let $\tilde X$ and $I(\tilde X)$ be the lifts of $X, I(X)$
to the universal cover $(\tilde M, \tilde \omega)$,
considered to be  a K\"ahler manifold. These vector fields are holomorphic
and conformal because $X$ and $I(X)$ are
holomorphic and Killing\index[terms]{vector field!Killing}\index[terms]{vector field!holomorphic} on $(M, \omega)$. 
By \ref{kami_or}, to prove that $(M, \omega)$
is Vaisman it suffices to show that the Lie algebra 
$\langle \tilde X,I(\tilde X)\rangle$ acts on 
$(\tilde M, \tilde \omega)$ non-isometrically.

\hfill

{\bf Step 3:}
Suppose that $\tilde X$ and 
$I(\tilde X)$ act on $(\tilde M, \tilde \omega)$
by isometries. By \cite[\S 2.132]{besse}, whenever $\tilde X$ and 
$I(\tilde X)$ are holomorphic Killing vector fields 
on a K\"ahler manifold $(N, I, g)$, these vector fields
are parallel.\index[terms]{vector field!parallel} However, $\tilde X$ and $I(\tilde X)$
are invariant under the homothety action of $\pi_1(M)$ on
$(\tilde M, \tilde \omega)$, and hence  their length cannot be constant,
and they cannot be parallel. 
\endproof

\hfill

This argument actually brings the following useful corollary.

\hfill

\corollary
Let $(M, g, \theta)$ be a Vaisman manifold,\index[terms]{manifold!Vaisman}
and $X, I(X)$ its Lee\index[terms]{Lee field} and anti-Lee fields.
Consider another LCK form $\omega'$
on $M$, and assume that $\Lie_X(\omega') = \Lie_{I(X)}(\omega')=0$.
Then $\omega'$ is Vaisman.
\endproof

\section{Opposite Lee forms on LCK manifolds with potential}\index[terms]{manifold!LCK!with potential}\index[terms]{form!Lee}

As another preliminary result, we need
the following non-existence claim. For Vaisman manifolds,
it was obtained by K. \index[persons]{Tsukada, K.} Tsukada (\cite{tsu}).

\hfill

\proposition\label{_Lee_cannot_be_opposite_Proposition_} 
{(\cite{ov_lee})} 
Let $(M, \theta, \omega)$ and $(M, \theta_1, \omega_1)$
be two LCK structures\index[terms]{structure!LCK} on the same
compact complex manifold. Suppose that
$(M, \theta, \omega)$ is an LCK structure with potential.
Then $\theta+\theta_1$ 
cannot be cohomologous to 0.

\hfill

\pstep
If $[\theta]$ is the Lee class for an\index[terms]{class!Lee} 
LCK structure with potential\index[terms]{structure!LCK!with potential} on $M$, then
$a[\theta]$ is also a Lee class for one, for any
$a>1$. To see this, consider an LCK structure with
preferred gauge, $\omega=d^c \theta + \theta \wedge\theta^c$.
Consider the corresponding K\"ahler potential
$\phi$ on the K\"ahler cover $(\tilde M, \tilde \omega)$,
with $\theta = -d\log\phi$.
Then $\phi^a$ is also a K\"ahler potential
on $\tilde M$,
\[ 
dd^c \phi^{a} = \phi^{a-2} (a \cdot \phi dd^c \phi + a(a-1)
d\phi\wedge d^c\phi).
\]
Indeed, the first summand $a\phi^{a-1} dd^c \phi$ is Hermitian, because
$dd^c\phi$ is Hermitian, and the second summand
$a(a-1)d\phi\wedge d^c\phi$ is positive.
The function $\phi^a$ is automorphic,
hence it defines an LCK structure\index[terms]{structure!LCK!with potential} with
potential on $M$, and the corresponding
Lee form \index[terms]{form!Lee}is  $- d \log (\phi^a)=a\theta$.

\hfill

{\bf Step 2:}
Let $\omega, \omega_1$ be LCK forms, and
$\theta, \theta_1$ the corresponding
Lee forms. Suppose that $k \theta + l \theta_1=0$.
Then 
\[ d(\omega^k \wedge \omega_1^l)= 
d(\omega^k) \wedge \omega_1^l + \omega^k \wedge
d(\omega_1^l)= k\theta \wedge \omega^k \wedge \omega_1^l
+ l\theta_1 \wedge \omega^k \wedge \omega_1^l=0.
\]
This computation can be interpreted as follows.
Let $L$ be the weight bundle for $(M, \omega, \theta)$
and $L_1$ the weight bundle for $(M, \omega_1, \theta_1)$.
We interpret $\omega$, $\omega_1$  as closed $L$- and\index[terms]{bundle!weight}
$L_1$-valued forms (Section \ref{_L_valued_Kahler_Subsection_}).
Then $\omega^k$ is a closed $L^{\otimes k}$-valued form,
$\omega_1^l$ is a closed $L_1^{\otimes l}$-form,
and $\omega^k \wedge \omega_1^l$ is a closed
form with coefficients in the flat bundle
$L^{\otimes k}\otimes L_1^{\otimes l}$,
that is  trivial.

Return now to the situation described in 
the assumptions of
\ref{_Lee_cannot_be_opposite_Proposition_}.
Let $n=\dim_\C M$.
Using Step 1, we replace the LCK structure\index[terms]{structure!LCK!with potential}
$(\omega, \theta)$ by another LCK structure with potential
in such a way that $\theta$ is replaced by
$(n-1)\theta$. Then $(n-1)\theta_1 = -\theta$,
and the volume form
$\omega \wedge \omega_1^{n-1}$ is closed.
However, $\omega$ is actually an exact $L$-valued
form, because $\omega= d_\theta (\theta^c)$, and hence 
$\omega \wedge \omega_1^{n-1}$ is an exact
$L \otimes L_1^{\otimes (n-1)}$-valued form,
that is, an exact form.

We verify this with an explicit computation:
\begin{equation*}
\begin{split}
d(\theta^c \wedge \omega_1^{n-1})&=
d(\theta^c) \wedge \omega_1^{n-1}-
\theta^c \wedge d(\omega_1^{n-1})\\
&=(\omega - \theta\wedge \theta^c)\wedge \omega_1^{n-1}
- (n-1)\theta_1 \wedge \theta^c \wedge \omega_1^{n-1}\\
&=\omega \wedge \omega_1^{n-1} -
(\theta\wedge \theta^c+(n-1)\theta_1 \wedge
\theta^c)\wedge \omega_1^{n-1}=
\omega \wedge \omega_1^{n-1}.
\end{split}
\end{equation*}
We have shown that the positive volume form
$\omega \wedge \omega_1^{n-1}$ on $M$ is exact,
that is  impossible.
\endproof

 \section[Hodge decomposition of $H^1(M)$ on LCK manifolds with potential]{Hodge decomposition of $H^1(M)$ on LCK\\ manifolds with potential}\index[terms]{manifold!LCK!with potential}
 \index[terms]{Hodge!decomposition}

We now prove a decomposition result
(\ref{_LCK_pot_Hodge_decompo_Theorem_}) for the first
cohomology group of an LCK manifold with potential, as a
direct sum of closed holomorphic 1-forms, closed
anti-holomorphic 1-forms, and the real line generated by
the Lee form. We start (\ref{_H^1_holo_LCK_Lemma_}) by
showing that the Lee form,\index[terms]{form!Lee} the closed holomorphic and
anti-holomorphic 1-forms can indeed be seen as
1-cohomology classes.

\hfill

 \lemma\label{_H^1_holo_LCK_Lemma_}
 Let $(M, \theta, \omega)$ be a compact LCK manifold,
 and $H_d^{1,0}(M)$ denote the space of closed holomorphic 1-forms on $M$.
 Then the natural map 
 \[ H_d^{1,0}(M)\oplus \overline{H_d^{1,0}(M)} \oplus \langle \theta\rangle\arrow H^1(M,\C)
 \]
 is injective,
 where $\langle \theta\rangle$ is the subspace generated by  $\theta$.
 
 \hfill
 
 \proof
A closed holomorphic form $\alpha$ belongs to $\ker d \cap \ker d^c$.
Indeed, $\bar\6 \alpha=0$ together with $d\alpha=0$ implies $d^c\alpha=0$.
 Therefore, if a real form $\beta\in H^{1,0}_d(M)+ \overline{H^{1,0}_d(M)}$
 is exact, one has $\beta = d f$ and $dd^c f=0$, that is  impossible by the 
 maximum principle. However, if $\theta$ is cohomologous to a 
 sum of holomorphic and antiholomorphic forms, this easily leads
 to a contradiction with \ref{_theta_not_d^c_closed_Corollary_}.
 Indeed, suppose that $\theta=\beta+ df$, where $\beta$ is
 a real form that satisfies
 $d\beta=d^c\beta=0$. Making a conformal change, we obtain
 another LCK structure\index[terms]{structure!LCK} that has Lee form\index[terms]{form!Lee} equal to $\beta$.
 This is impossible, again by \ref{_theta_not_d^c_closed_Corollary_}.
 \endproof
 
 \hfill
 
 This lemma immediately brings the following
 
 \hfill
 
\corollary\label{_inequa_holo_LCK_Corollary_}
Let $M$ be a compact LCK manifold, and $H^{1,0}_d(M)$ denote 
the space of closed  holomorphic 1-forms on $M$.
Then $\dim H_d^{1,0}(M) \leq \frac{b_1(M)-1}{2}$.
\endproof

 \hfill
 
 \theorem\label{_LCK_pot_Hodge_decompo_Theorem_} 
{(\cite{ov_lee})} 
 Let $(M, \theta, \omega)$ be a compact LCK manifold with potential,\index[terms]{manifold!LCK!with potential}
 and $H_d^{1,0}(M)$ denote the space of all closed holomorphic 1-forms on $M$.
 Then 
 $$H^1(M,\C) = H_d^{1,0}(M)\oplus \overline{H_d^{1,0}(M)} \oplus \langle \theta\rangle.$$

 \proof
 For Vaisman manifolds,\index[terms]{manifold!Vaisman} this result is already proven 
 (\ref{_harmo_deco_1-form_Proposition_}). 
 Indeed, by \ref{_harmo_deco_1-form_Proposition_}
 the space of harmonic forms in $\Lambda^1(M,\C)$ is 
 identified with $\ker d \cap \ker d^c\oplus \langle \theta\rangle$.
 
 If $M$ is just an LCK manifold with potential\index[terms]{manifold!LCK!with potential}, the injectivity of the map 
 \[ H_d^{1,0}(M)\oplus \overline{H_d^{1,0}(M)} \oplus \langle \theta\rangle
 \arrow H^1(M,\C)
 \]
 follows from \ref{_H^1_holo_LCK_Lemma_}.
 To prove the surjectivity, it would suffice to show that
 $\dim_\C H_d^{1,0}(M)= \frac{b_1(M)-1}{2}$. We prove it by deforming
 $M$ to a Vaisman manifold\index[terms]{manifold!Vaisman} $M_0$ and showing that 
 $\dim H_d^{1,0}(M)= \dim H_d^{1,0}(M_0)$.
 
 As always, we first deform the given LCK metric on $M$ to an LCK metric
 of LCK rank 1 \index[terms]{rank!LCK}(\ref{defor_improper_to_proper}). This operation does not affect the complex structure
 on $M$, and hence  $\dim H_d^{1,0}(M)$ does not change, and it will suffice to 
 prove that $\dim_\C H_d^{1,0}(M)= \frac{b_1(M)-1}{2}$ when 
 $M$ is an LCK manifold with proper potential.\index[terms]{manifold!LCK!with potential}
 
 Let $\tilde M$ be the open algebraic cone\index[terms]{cone!algebraic} associated with $M$ as in 
 \ref{_cone_cover_for_LCK_pot_Theorem_}, and denote by $A:\; \tilde M \arrow \tilde M$
 the generator of the deck group. Applying the Jordan--Chevalley\index[terms]{Jordan--Chevalley decomposition}
 decomposition $A=SU$ as in \ref{def_lckpot2Vai}, we can deform
 $\tilde M/\langle A\rangle$ to the Vaisman manifold\index[terms]{manifold!Vaisman}
 $M_0:=\tilde M/\langle S\rangle$. To prove that 
 $\dim_\C H_d^{1,0}(M) = \frac{b_1(M)-1}{2}$ it would suffice to show that
 all holomorphic, $S$-invariant 1-forms on $\tilde M$ are
 also $U$-invariant. 
 
 Consider $U$ as an automorphism of $M_0$. This automorphism
 is homotopy equivalent to the identity because
 $U= e^{N}$, where $N$ commutes with $S$.
 Since $U$ is an unipotent element of the group
 of automorphisms of the algebraic cone\index[terms]{cone!algebraic} $\tilde M$,
 the action of $U_t:= e^{tN}$ preserves $\tilde M$ and 
 commutes with $S$, and hence  it is well-defined on $M_0$.
 This gives a homotopy of $U=U_1$ to $\Id=U_0$.
 
 Since $U$ is homotopy equivalent to the identity,
 it acts trivially on $H^1(M_0)$, and hence 
 all $S$-invariant holomorphic forms
 on $\tilde M$ are also $SU$-invariant.
 This implies that $\dim_\C H_d^{1,0}(M) \geq \dim_\C H_d^{1,0}(M_0)  =\frac{b_1(M_0)-1}{2}$.
 The inequality in this expression is, in fact, an equality by
 \ref{_inequa_holo_LCK_Corollary_}. We proved \ref{_LCK_pot_Hodge_decompo_Theorem_}.
 \endproof

\section{The set of Lee classes on Vaisman manifolds}\index[terms]{manifold!Vaisman}\index[terms]{class!Lee}

\definition
{\bf The Lee class} of an LCK manifold $(M, \theta, \omega)$
is the cohomology class $[\theta]\in H^1(M,\R)$.\index[terms]{class!Lee}

\hfill

\question\label{Lee_class_on_Vaisman_Question_} Let
$(M,I)$ be a compact complex manifold that admits Vaisman
metrics. Which elements in $H^1(M,\R)$ can be Lee classes
for a Vaisman metric? 

\hfill

\question\label{Lee_class_on_potential_Question_} Let
$(M,I)$ be a compact complex manifold that admits LCK
metrics with potential. Which elements in $H^1(M,\R)$ can
be Lee classes for an LCK metric with potential?\index[terms]{manifold!LCK!with potential}

\hfill

The answer to \ref{Lee_class_on_Vaisman_Question_} is
surprisingly simple. \ref{_holomo_on_Vaisman_basic_Proposition_} below was obtained by
K. \index[persons]{Tsukada, K.} Tsukada (\cite{tsu}), using the harmonic decomposition
of 1-forms (\ref{_harmo_deco_1-form_Proposition_}). We\index[terms]{decomposition!harmonic}
give a new proof (according to \cite{ov_lee}), based on
\ref{_LCK_pot_Hodge_decompo_Theorem_}. To begin with, we
need to prove that on a compact Vaisman manifold,\index[terms]{manifold!Vaisman} the
closed holomorphic 1-forms are basic with respect to the
canonical foliation.\index[terms]{foliation!canonical}\footnote{In the framework of harmonic
decomposition, this is \ref{_harmo_deco_1-form_Proposition_},
Step 2.}
Recall (\ref{_basic_form:
  definition_}) that a form $\eta$ on $(M,\Sigma)$ is
basic with respect to
$\Sigma\subset TM$ if and only if for any vector field
$X \in \Sigma$, one has $i_X(\eta) = \Lie_X(\eta)=0$,
where $i_X$ denotes the contraction with $X$. In particular, by the Cartan formula, a closed form $\eta$ is basic with respect to
$\Sigma\subset TM$ if and only if $i_X(\eta) = 0$.\index[terms]{form!basic}

\hfill

\proposition\label{_holomo_on_Vaisman_basic_Proposition_}
Let $M$ be a compact Vaisman manifold,\index[terms]{manifold!Vaisman} and $\eta$ a 
closed holomorphic 1-form on $M$. Then $\eta$ is basic
with respect to the canonical foliation\index[terms]{foliation!canonical} $\Sigma$ on $M$.

\hfill

\proof
Let $n=\dim_\C M$,
and $\omega_0=d\theta^c\in \Lambda^{1,1}(M)$ 
the transversal K\"ahler form\index[terms]{form!K\"ahler!transversal} defined 
on every Vaisman manifold\index[terms]{manifold!Vaisman} (\ref{_Subva_Vaisman_Theorem_}).
Since $\eta$ is closed and $\omega_0$ is exact,
one has $\int_M \omega_0^{n-1}\wedge \eta\wedge\bar \eta=0$.
However, $-\1 \eta\wedge\bar \eta$
is a semi-positive form, and $\omega_0$
is strictly positive in the directions transversal
to $\Sigma$. This implies that
$-\1 \omega_0^{n-1}\wedge \eta\wedge\bar \eta$
is a positive volume form in every point
$x\in M$ such that $\eta\restrict {T_x M}$ does not vanish
on $\Sigma\restrict{T_x M}$. Since 
$\int_M \omega_0^{n-1}\wedge \eta\wedge\bar \eta=0$,
it follows that $\eta \restrict\Sigma=0$ everywhere.
By \ref{_basic_for_closed-Remark_}, $\eta$ is basic.
\endproof

\hfill

\theorem \label{_Lee_cone_on_Vaisman_Theorem_}
Let $M$ be a compact Vaisman manifold,\index[terms]{manifold!Vaisman} and
\[ H^1(M)= H_d^{1,0}(M) \oplus \overline{H_d^{1,0}(M)} \oplus \langle \theta\rangle\]
the decomposition established in \ref{_LCK_pot_Hodge_decompo_Theorem_}.
Consider a 1-form $\mu\in H^1(M)^*$ vanishing on 
$H^{1,0}(M) \oplus \overline{H^{1,0}(M)} \subset H^1(M)$
and satisfying $\mu([\theta])>0$. Then 
a class $\alpha\in H^1(M,\R)$ is a 
Lee class for some LCK structure \index[terms]{structure!LCK}if and only if $\mu(\alpha) >0$.

\hfill

\pstep
We start by proving that any $\alpha\in H^1(M,\R)$ satisfying
$\mu(\alpha) >0$ can be realized as a Lee class. 

From \eqref{_omega_via_theta_Equation_},
we have $\omega=d^c\theta+\theta\wedge I\theta$.
By \ref{_Subva_Vaisman_Theorem_} (ii),  the form
$\omega_0:=d^c\theta$ is semi-positive: 
it vanishes on the canonical foliation \index[terms]{foliation!canonical}$\Sigma$ 
and is strictly positive in the transversal directions. 
Let $u\in H_d^{1,0}(M) \oplus \overline{H_d^{1,0}(M)}$.
Then
\[
d^c(\theta+ u) + (\theta+u) \wedge (\theta^c + u^c)
= \omega_0 + (\theta+u) \wedge (\theta^c + u^c)
\]
is the sum of two semi-positive forms. Indeed,
since $u$ is $d^c$-closed, $d^c(\theta+ u) =\omega_0$
is semi-positive; the form  $(\theta+u) \wedge (\theta^c + u^c)$
is semi-positive of rank 1 by definition.

By \ref{_holomo_on_Vaisman_basic_Proposition_}, $u$ is basic.
Since
$\theta+u$ is the sum of $\theta$ and a
basic form, the restriction of
$(\theta+u) \wedge (\theta^c + u^c)$ to $\Sigma$
satisfies
\[
(\theta+u) \wedge (\theta^c + u^c)\restrict \Sigma =
\theta \wedge \theta^c \restrict \Sigma.
\]
The 
sum $\omega_0 + (\theta+u) \wedge (\theta^c + u^c)$
is (strictly) positive on all tangent
vectors $x\notin \Sigma$ because $\omega_0$ is positive
on these vectors,\footnote{When we say ``a positive
	(1,1)-form $\alpha$ is positive on a vector $v$'',
	we mean that $\1\alpha(v, I(v))>0$; a form is
	Hermitian if it is positive on all non-zero vectors.}
and positive on $x\in \Sigma$
because $\theta \wedge \theta^c \restrict \Sigma$
is positive on such $x$.

\hfill

{\bf Step 2:} It remains to show that
none of the classes $\alpha$ with $\mu(\alpha)\leq 0$ can be
realized as a Lee class of an LCK structure. \index[terms]{structure!LCK}
By \ref{_Lee_form_on_Vaisman_is_Vaisman_Proposition_},
any Lee class on $M$ is the Lee class of a Vaisman metric.
If $\mu(\alpha)=0$, we can represent $\alpha$ by
a $d, d^c$-closed form $\alpha_0$. This is impossible by
\ref{_theta_not_d^c_closed_Corollary_}. 

\hfill

{\bf Step 3:} It remains to show that
there are no LCK classes with 
$\mu(\alpha)< 0$. Suppose that 
such a class exists; by \ref{_Lee_form_on_Vaisman_is_Vaisman_Proposition_},
it is the Lee class of a Vaisman manifold,\index[terms]{manifold!Vaisman} hence
it has an LCK potential.\index[terms]{potential!LCK} This is impossible,
as two Lee classes for LCK structures\index[terms]{structure!LCK} with
potential cannot sum to zero, by \ref{_Lee_cannot_be_opposite_Proposition_}.
\endproof

\hfill

\remark 
Let $(M,I)$ be a compact complex manifold of
Vaisman type, and ${\cal L}_v\subset H^1(M, \R)$ the set
of all Lee classes for all Vaisman structures on $(M,I)$. 
From \ref{_Lee_classes_for_Vaisman_open_Corollary_}
it follows that the set of all Lee classes
is open in $H^1(M, \R)$, and from 
\ref{_Lee_cone_on_Vaisman_Theorem_} we obtain
${\cal L}_v$ is actually a half-space.

\section{The set of Lee classes on LCK manifolds with potential}\index[terms]{manifold!LCK!with potential}

\theorem\label{_Lee_cone_on_LCK-pot_Theorem_} 
{(\cite{ov_lee})} 
Let $(M, \theta, \omega)$ be a compact LCK manifold with potential,
and  $\mu:\; H^1(M, \R)\arrow \R$ a non-zero linear map vanishing
on the space $H_d^{1,0}(M)\oplus \overline{H_d^{1,0}(M)}$ which
has codimension 1 by \ref{_LCK_pot_Hodge_decompo_Theorem_}.
Assume that $\mu(\theta) >0$. Then $\xi\in H^1(M, \R)$ is the Lee class of an 
LCK structure with potential on $M$ if and only if $\mu(\xi)>0$.\index[terms]{structure!LCK!with potential}

\hfill

\proof
Let $(M, I_\infty)$ be a Vaisman manifold,\index[terms]{manifold!Vaisman}
and $\{I_i\}$ the sequence of complex structures
converging to $I_\infty$, such that all manifolds
$(M, I_k)$ are isomorphic to $(M,I)$
(\ref{_Vaisman_limit_of_LCK_pot_Theorem_}).
Given an LCK metric with preferred gauge
\[
\omega_\infty= d^c\theta_\infty + \theta_\infty \wedge
\theta^c_\infty =I_\infty d I_\infty ^{-1}\theta_\infty + \theta_\infty \wedge I_\infty(\theta_\infty)
\]
on $(M,I_\infty)$, the form 
$I_kdI_k^{-1}\theta_\infty + \theta_\infty \wedge I_k(\theta_\infty)$
remains strictly positive for almost all manifolds $(M, I_k)$,
because $\lim_k I_k = I_\infty$, and positivity is an open
condition. Therefore, $\theta_\infty$ is a Lee form\index[terms]{form!Lee} on $(M,I)$.
We already proved that the set ${\goth L}$ of Lee classes on $(M,I)$ 
contains the half-space $\{u\in H^1(M, \R)\ \ |\ \ \mu_0(u)>0\}$
for some linear map $\mu_0:\; H^1(M, \R)\arrow \R$.
By \ref{_Lee_cannot_be_opposite_Proposition_}, 
${\goth L}$ cannot be greater than a half-space.
It only remains to show that $\mu_0$ is proportional
to $\mu$. This will follow if we prove that
$\ker \mu=\ker \mu_0$. The space $\ker \mu_0$ is
the set of  classes $\alpha \in H^1(M, \R)$
such that neither $\alpha$ nor $-\alpha$ are Lee classes,
and $\ker \mu$ are classes represented by $d, d^c$-closed forms.
By \ref{_theta_not_d^c_closed_Corollary_},
a Lee class of an LCK manifold cannot be represented
by a $d^c$-closed form, which gives $\ker \mu=\ker \mu_0$.
\endproof

\hfill

As a corollary, we obtain a vanishing theorem for the Morse--Novikov cohomology
of a compact LCK manifold with potential\index[terms]{manifold!LCK!with potential}.

\hfill

\corollary\label{_MN_vanishes_Lee_class_Corollary_} {(\cite{ov_lee})} 
Let $(M, \theta, \omega)$, $\dim_\C M >2$, be a compact LCK manifold with potential,
 and $H^*_\theta(M)$ the Morse--Novikov cohomology
groups, that is, the cohomology of $(\Lambda^*(M), d_\theta)$.
Then $H^*_\theta(M)=0$.

\hfill

\proof
For Vaisman manifolds,\index[terms]{manifold!Vaisman} this follows from \cite{llmp}.
Indeed, in \cite{llmp} it is proven that $H^*_\theta(M)=0$ for 
all compact LCS \index[terms]{manifold!LCS}
manifolds admitting a Riemannian metric $g$ such that the
Lee form\index[terms]{form!Lee!parallel} is $g$-parallel.
Given an LCK manifold $(M, \theta, \omega)$ with potential\index[terms]{manifold!LCK!with potential},
by \ref{_Lee_cone_on_LCK-pot_Theorem_} there exists a 
Vaisman manifold $(M, \theta', \omega')$
with the same Lee class: $[\theta]=[\theta']$. 
Since $H^*_\theta(M)$ depends only on the cohomology
class of $\theta$, this implies $H^*_\theta(M)=0$.

For a direct proof of this result, not referring to \cite{llmp},
see \ref{mn_van}.
\endproof

\section{Notes}
\begin{enumerate}
	\item The first attempt to determine the space of 
          Lee forms\index[terms]{form!Lee} on a compact LCK manifold  appears in
          \cite{tsu}. In Theorem 5.1, \index[persons]{Tsukada, K.} Tsukada proves that
          on a compact Vaisman type manifold $M$, the Lee
          classes form a half-space in $H^1(M)$. His proof
          is based on the harmonic decomposition of
          1-forms on Vaisman manifolds.\index[terms]{manifold!Vaisman}

	\item The space of Lee classes on compact  LCK
          manifolds that do not admit metrics with LCK
          potential can be very small, as happens on the \index[terms]{surface!Inoue} Inoue surfaces $S^0$
          (see  Chapter \ref{inoue_lck}): in \cite[Theorem
            3.11]{oti2}, A. \index[persons]{Otiman, A.} Otiman proved that on these
          surfaces, the space of Lee classes is 
          a point.
\end{enumerate}

\section{Exercises}

\begin{enumerate}[label=\textbf{\thechapter.\arabic*}.,ref=\thechapter.\arabic{enumi}]


\item
Let $\omega_0$ be a semi-positive exact (1,1)-form 
on a compact complex manifold $M$, and $\Sigma=\ker \omega_0$;
assume that $\rk_\C \Sigma=1$. 
Consider a real 1-form $\eta$ such that $d\eta=d^c\eta=0$.
Prove that $\Sigma$ is involutive. 
Prove that $\eta$ is basic with respect to $\Sigma$.

\item
As usual, we denote the quaternion algebra by ${\Bbb H}$.
Let ${\Bbb H}^2 \backslash 0\arrow {\Bbb H}P^1=S^4$
be the standard projection. Identifying ${\Bbb H}^2$ and $\C^4$,
and considering $\C$ as a subalgebra in ${\Bbb H}$, we can 
factorize ${\Bbb H}^2 \backslash 0\arrow {\Bbb H}P^1=S^4$
through the map ${\Bbb H}^2\backslash 0 = \C^4 \backslash 0 \arrow \C P^3$.
This defines the Hopf map $\C P^3 \arrow S^4$.
Prove that all fibres of this map are complex lines
in $\C P^3$.\footnote{This map is called ``the twistor fibration''
for $S^4$.}

\item Let $(M,I)$ be a complex manifold, and 
$\Sigma\subset TM$ be a real sub-bundle of rank 2, satisfying
$I(\Sigma)=\Sigma$. Assume that $\Sigma$ is involutive.
Would this imply that $\Sigma$ is holomorphic?

{\em Hint:} Try to construct a counter-example,
using the previous exercise.

\item\label{_defo_Gauduchon_dd^c_exact_Exercise_}
Recall that a Hermitian form\index[terms]{form!Hermitian} $\eta$ on a complex $n$-manifold
is called {\bf Gauduchon} if $dd^c(\eta^{n-1})=0$.
Let $\eta$ be a Gauduchon form on $M$ that satisfies
$\eta^{n-1} \in \im(dd^c)$. Prove that a small complex 
deformation of $M$ admits a Gauduchon form\index[terms]{form!Gauduchon}
that satisfies $\eta^{n-1} \in \im(dd^c)$.

\item
Define {\bf the  Aeppli cohomology} of a complex manifold
as $\frac{\ker dd^c}{\im d + \im d^c}$.\index[terms]{cohomology!Aeppli}
Prove that the Aeppli cohomology is finite-dimensional.

{\em Hint:} Use the same argument as in \ref{_BC_f_dim_Theorem_}.

\item Let $\eta$ be a Gauduchon form on a 
compact complex manifold $M$. Assume that 
$M$ admits a semi-positive closed (1,1)-form.
Prove that the Aeppli class of $\eta^{n-1}$ is non-zero.

\item
Let $\eta$ be a Gauduchon form\index[terms]{form!Gauduchon} on $M$.
Prove that $\eta^{n-1} \notin \im(dd^c)$
\begin{enumerate}
\item when $M$ is a compact Vaisman manifold;\index[terms]{manifold!Vaisman}
\item when $M$ is an LCK manifold with potential\index[terms]{manifold!LCK!with potential}.
\end{enumerate}

{\em Hint:} For (a), use the previous exercise.
For (b), use Exercise \ref{_defo_Gauduchon_dd^c_exact_Exercise_}.

\item
Find two LCK structures\index[terms]{structure!LCK} $(\omega, \theta)$ and $(\omega', \theta')$
on the same compact (non-K\"ahler) complex manifold such that $\theta$ is cohomologous to $\theta'$,
but $\omega$ is not conformally  equivalent to $\omega'$.

\item
Let $(\omega, \theta)$ and $(\omega_1, \theta)$ be 
LCK structures on $M$. Prove that $\omega+\omega_1$ is an LCK form.

\item
Let $M$ be a complex manifold,
$[\theta]\in H^1(M, \R)$, and
${\cal C}_\theta$ the set of all
LCK structures with the Lee class $[\theta]$.
Prove that ${\cal C}_\theta$ is contractible.\index[terms]{contractible}

{\em Hint:} Use the previous exercise.

\item
Let $M$ be a compact complex manifold
admitting an LCK structure with potential.
Prove that the set of LCK structures on $M$
is contractible.\index[terms]{structure!LCK}\index[terms]{structure!LCK!with potential}

{\em Hint:} Use the previous exercise and \ref{_Lee_cone_on_LCK-pot_Theorem_}.

\item
Let $(M, \theta, \omega)$ be an LCK
manifold, and $\phi$ a positive automorphic
K\"ahler potential on its K\"ahler cover.
\begin{enumerate}
\item
Prove that
$dd^c(\log\phi)=\frac{dd^c\phi}{\phi}-\frac{d\phi\wedge d^c\phi}{\phi^2}$
and 
\[\frac{dd^c \phi^{a}}{a\phi^a} = 
\frac{dd^c\phi}{\phi}-(a-1)\frac{d\phi\wedge d^c\phi}{\phi^2}.
\]
Use this to show that
\[
\lim_{a\arrow 0} \frac{dd^c \phi^{a}}{a\phi^a} 
=\frac{dd^c\phi}{\phi}-\frac{d\phi\wedge d^c\phi}{\phi^2} = dd^c(\log\phi).
\]
\item
Assume that $\phi^a$ is an automorphic
K\"ahler potential for all $a>0$. Prove that 
$d^c \theta$ is a semi-positive (1,1)-form.
\end{enumerate}

\item
Let $M$ be a compact Vaisman manifold\index[terms]{manifold!Vaisman} equipped with a smooth elliptic
fibration $\pi:\; M \arrow X$. Prove that 
$H^1(M) = \pi^*(H^1(X)) \oplus \langle \theta \rangle$. 

\item
Let $(M, \theta, \omega)$ be a Vaisman manifold. Recall that it satisfies the equation
$\omega= d^c \theta + \theta \wedge \theta^c$. Prove that
$(M, a\theta, a d^c \theta + a^2\theta \wedge \theta^c)$
is Vaisman for any $a\in \R^{>0}$.

\item
Let $(M, \theta, \omega)$ be a compact $n$-dimensional LCK manifold with potential,\index[terms]{manifold!LCK!with potential}
and $(M, \theta_1, \omega_1)$ another LCK structure on the same\index[terms]{structure!LCK}
complex manifold. Assume that $\theta$ is cohomologous to $\lambda \theta_1$
for some $\lambda\in \R$. 
\begin{enumerate}
\item Suppose that $\theta + (n-1) \theta_1$ is cohomologous to 0,
that is, $\lambda = -\frac{1}{n-1}$.
Prove that the volume form $\omega\wedge \omega_1^{n-1}$ is exact,
giving a contradiction.
\item Suppose that $\lambda \leq -\frac{1}{n-1}$.
Prove that there exists an LCK structure with potential\index[terms]{structure!LCK!with potential}
 $(M, \theta', \omega')$ such that $\theta'+(n-1) \theta_1$ is exact.
\item
Prove that $\lambda >-\frac{1}{\dim_\C M-1}$.
\end{enumerate}

\item \label{_sum_Lee_even_dim_Exercise_}
Let $M$ be a compact complex manifold, and
 $(M, \theta, \omega)$ and $(M, \theta_1, \omega_1)$ LCK structures\index[terms]{structure!LCK}
such that $\theta+ \theta_1$ is exact. 
\begin{enumerate}
\item Prove that $\omega^k \wedge \omega_1^k$ is closed for any $k\in \Z^{>0}$.
\item Assume that $\omega$ is $d_\theta$-exact. Prove that
$\omega^k \wedge \omega_1^k$ is exact for any $k\in \Z^{>0}$.
\item Assume that $\dim_\C M$ is even, $\dim_\C M=2k$.
Prove that $d_\theta$-exactness of $\omega$ implies exactness
of the volume form $\omega^{k} \wedge \omega_1^{k}$,
hence $\omega$ cannot be $d_\theta$-exact.
\end{enumerate}

\item
Let $(M, \omega, \theta)$ be a locally conformally
symplectic (LCS) manifold. Assume that $\omega$ is\index[terms]{manifold!LCS}
$d_\theta$-exact, that is, 
$\omega = d\eta - \theta \wedge \eta$.
Prove that for any 1-form $\eta'$ that is  $C^1$-close\index[terms]{topology!$C^1$} to $\eta$
and any closed 1-form $\theta'$ that is  $C^1$-close to
$\theta$, the form $\omega' = d\eta' - \theta' \wedge \eta'$
is also locally conformally symplectic.

\item
{\bf The Lee class} of $(M, \omega,\theta)$ 
of an LCS manifold is the cohomology class of $\theta$. 
Prove that the set of all Lee classes of all LCS on a
compact manifold $M$ 
structures with $d_\theta$-exact form $\omega$ is open.

\item
Let $M$ be a compact manifold, $\dim_\R M = 4n$,
and $(\omega, \theta)$, $(\omega_1, \theta_1)$
LCS structures.\index[terms]{structure!LCS}
Assume that $\omega$ is $d_\theta$-exact.
Prove that $\theta+\theta_1$ cannot be exact.

{\em Hint:} Use the same argument as in 
Exercise \ref{_sum_Lee_even_dim_Exercise_}.

\item
Recall that a real $(p,p)$-form on an $n$-dimensional
complex vector space $V$ is called {\bf weakly positive}
if $\eta(X_1, I(X_1), X_1, I(X_2), ..., X_p, I(X_p))\geq 0$ 
for all tangent vectors $X_1, ..., X_p$, and
{\bf strongly positive} if $\eta$ is a linear combination
of products of positive (1,1)-forms with non-negative coefficients. Prove that the cone
of weakly positive $(p,p)$-forms is dual to the cone of strongly
positive $(n-p, n-p)$-forms with respect to the product
$\langle \eta_1 ,\eta_2\rangle:= \frac{\eta_1\wedge \eta_2}{\Vol V}$.

\item
Let $\alpha$ be a weakly positive $(p,p)$-form on a Vaisman $n$-manifold
such that $\alpha \wedge \omega_0^{n-p}=0$. 
\begin{enumerate}
\item Prove that
$\alpha$ vanishes on the canonical foliation\index[terms]{foliation!canonical} $\Sigma$.
\item Assume that $\alpha$ is in addition closed.
Prove that $\alpha$ is $\Sigma$-basic.
\end{enumerate}

{\em Hint:} Use the previous exercise.

\item \label{_holo_forms_on_Exercise_}
Let $\eta$ be a holomorphic $p$-form on a compact Vaisman $n$-manifold.
Prove that $\eta$ is $\Sigma$-basic if
\begin{enumerate}
\item $d\eta =0$ and $p < n$.
\item $p< n-1$.
\end{enumerate}

{\em Hint:} Use the previous exercise.

\item
Let $M= \frac{\Tot^\circ(L)}{\Z}$ be a Vaisman manifold\index[terms]{manifold!Vaisman}
constructed from the total space of all non-zero vectors in a
negative holomorphic line bundle $L$ over a projective orbifold $X$
(\ref{_Structure_of_quasi_regular_Vasman:Theorem_}).
Assume that the canonical bundle of $X$ is trivial. Prove
that the canonical bundle of $M$ is trivial. 

{\em Hint:} Use the adjunction formula.\index[terms]{adjunction formula}

\item
Let $(M, \omega,\theta)$ be an LCK manifold, and $G$ a finite group freely acting
on $M$. Suppose that $G$ preserves the Lee class of $M$.
\begin{enumerate}
\item 
Prove that there exists a $G$-invariant LCK metric conformally equivalent to
$\omega$.
\item Prove that the quotient $M/G$ admits an LCK
  structure.
\end{enumerate}

\item\label{_quotients_Vaisman_Extercise_}
Let $(M, \omega,\theta)$ be a compact LCK manifold with potential, \index[terms]{manifold!LCK!with potential}
and $G$ a finite group freely acting on $M$. 
\begin{enumerate}
\item Prove that the set of Lee classes of LCK
  structures with potential on $M$ is convex.
\item Prove that there exists an LCK metric with potential on $M$\index[terms]{metric!LCK!with potential}
  with $G$-invariant Lee class.
\item Prove that the quotient $M/G$ admits an LCK
  structure with potential.\index[terms]{structure!LCK!with potential}
\item 
Let $M$ be a compact Vaisman manifold,\index[terms]{manifold!Vaisman} and
$G$ a finite group freely acting on $M$. 
Prove that the quotient $M/G$ admits a Vaisman metric.
\end{enumerate}

\item\label{_taming_Lee_cone_Vaisman_Exercise_}
Let $(M,I)$ be a complex, or an almost complex, manifold.
A locally conformally symplectic form $\omega$ is {\bf taming} the almost
complex structure $I$ if $\omega(x, Ix)>0$ for any non-zero
tangent vector $x\in T_m M$.  The {\bf taming Lee cone}\index[terms]{cone!Lee!taming}
is the set of Lee classes of all LCS forms\index[terms]{form!LCS} $\omega$ taming $I$.
\begin{enumerate}
\item Show that $\omega$ is taming $I$ if and only
if its $(1,1)$-part $\omega^{1,1}$ is an Hermitian form.\index[terms]{form!Hermitian}
\item Let $\omega$ be a Hermitian LCS form 
on an almost complex manifold $M$, $\dim_\C M =n$,
and $\omega_1$ a taming LCS form, and $\theta$,
$\theta_1$ their Lee classes. Assume that 
$a\theta= - b \theta_1$, where $a, b$ are positive
integer numbers, and $a+b = n$. Prove that 
the Morse--Novikov classes\index[terms]{class!Morse--Novikov} $[\omega]\in H^2_\theta(M)$
and $[\omega_1]\in H^2_{\theta_1}(M)$ are non-zero.
\item Prove that the taming Lee cone on
an LCK manifold with potential\index[terms]{manifold!LCK!with potential} is equal to the Lee cone
of LCK structures with potential.\index[terms]{structure!LCK!with potential}
\end{enumerate}

{\em Hint:} Use the same
argument as in \ref{_Lee_cannot_be_opposite_Proposition_}
and then apply \ref{_Lee_cone_on_LCK-pot_Theorem_}.

\end{enumerate}


\chapter{Harmonic forms on Sasakian and Vaisman
  manifolds}\label{_Harmonic_forms_chapter_}

\epigraph{\it ``Everything you say is boring and
incomprehensible,'' she said, ``but that alone doesn't
make it true.''}{\sc \scriptsize Franz Kafka, ``Description of a Struggle''}


\section{Introduction}


The raison d'\^etre of this chapter is twofold.
First, we prove the harmonic decomposition on
Vaisman and Sasakian manifolds. Second, we explain
the supersymmetry which lies at the foundation of Hodge
theory, and apply this supersymmetry approach to
Sasakian manifolds, obtaining an interesting
Lie superalgebra acting on the differential forms.

Recall that a differential form on 
a manifold $M$ equipped with a foliation $F$
is called {\bf basic} if it is locally lifted from
the leaf space of $F$. A Sasakian manifold
is equipped with the Reeb foliation, that is 
actually transversally K\"ahler (that is,
its leaf space is equipped with a natural
K\"ahler structure). In this context,\index[terms]{foliation!transversally K\"ahler}
the harmonic decomposition expresses
the harmonic forms on the Sasakian
manifolds in terms of the basic
harmonic forms and the contact form
of the Sasakian manifold. This result
is very useful, because it allows to
define the Hodge decomposition on the
cohomology of the Sasakian manifold, and
use  other results about harmonic forms
on K\"ahler manifolds.\index[terms]{decomposition!Hodge}

A Vaisman manifold has \index[terms]{manifold!Vaisman}
two foliations, the Lee foliation
of real codimension one, and the 
canonical foliation,\index[terms]{foliation!canonical} that is  
holomorphic of complex dimension one.
The basic forms with respect to 
the canonical foliation (that is 
transversally K\"ahler, too) also
admit the Hodge decomposition.\index[terms]{decomposition!harmonic}
The harmonic decomposition allows 
interpreting the harmonic forms on a 
Vaisman manifold in terms of the
basic harmonic forms on this transversally K\"ahler
foliation, obtaining the Hodge
decomposition on cohomology.

The harmonic decomposition on Sasakian manifolds is 
due to Shun-ichi Ta\-chi\-bana \cite{_Tachibana_}, 1965.
\index[persons]{Tachibana, S.} Tachibana proved that on a compact Sasakian
manifold $Q$, $\dim_\R Q =2k +1$, every harmonic form in 
degree $\leq k$ is basic with respect to the 
Reeb foliation. The proof is so technical that
it is not even reproduced in the Boyer-Galicki
brilliant and comprehensive monograph \cite{bog} on Sasakian manifolds; they just
refer to Tachibana's result. We give a new proof of his theorem
in \ref{_harmo_Sasa_decompo_Theorem_}, 
and express the basic harmonic forms in terms
of the K\"ahler geometry\index[terms]{geometry!K\"ahler} of the Reeb 
leaf space. 

The same approach was used 15 years later 
by \index[persons]{Kashiwada, T.} Kashiwada and Kashiwada-\index[persons]{Sato, S.}Sato 
(\cite{KS}, \cite{kashiwada_kodai}) to study
the harmonic forms on Vaisman manifolds.\index[terms]{manifold!Vaisman} Later, Kashiwada's
results were re-proven by I. Vaisman,\index[persons]{Vaisman, I.}
using \index[persons]{Tachibana, S.} Tachibana's harmonic decomposition (\cite{va_gd}).
We give a new proof the results of Kashiwada
in \ref{_Vaisman_harmonic_forms_Theorem_}.

Many Hodge-theoretic results related to the Dolbeault cohomology\index[terms]{cohomology!Dolbeault} and the
related structures in the context of  this
harmonic decomposition were later obtained
by K. \index[persons]{Tsukada, K.} Tsukada (\cite{tsuk}).

Here we give new proofs of these classical results.
This chapter largely reproduces the contents 
of the paper \cite{_ov_super_sas_}.

\subsection{Basic cohomology and taut foliations}
\label{_basic_taut_Subsection_}

The approach we use in this chapter is more
topological. Let $M$ be a manifold equipped with
a foliation. The {\bf basic cohomology} is the cohomology
of the differential graded algebra of basic forms. Generally speaking,
it is hard to find geometric relevance in the basic\index[terms]{cohomology!basic}
cohomology; for example, the basic cohomology can be
even infinite-dimensional \cite{_Schwarz:foliation_}.
Everything changes when the foliation is Riemannian. 

Recall that a foliation $F \subset TM$ on a Riemannian
manifold $(M,g)$ is called {\bf a Riemannian foliation} if the 
locally defined map $\pi$ to the leaf space is a Riemannian
submersion, that is, there exists a Riemannian metric\index[terms]{foliation!Riemannian}
on the leaf space such that differential of $\pi$
restricted to $F^\bot$ is an isometry (Section 
\ref{_foliations_Notes_}). In this
case, the Riemannian form $g$ is called 
{\bf bundle-like}.\index[terms]{metric!bundle-like}

A $p$-form $\eta$ on a foliated manifold $(M, F)$
is called {\bf relatively closed} if \[ d\eta(X_1, ..., X_{p+1})=0\]
for any vector fields $X_1, ..., X_{p+1}$, with
$X_1, ..., X_{p+1}\in F$.\index[terms]{form!relatively closed}

The notion of a taut foliation was defined by D. \index[persons]{Sullivan, D.} Sullivan in
\cite{_Sullivan:taut_}. An involutive distribution $F \subset TM$
is {\bf taut} if there exists a metric on $M$ such that
all leaves of $F$ are minimal, that is, have vanishing 
mean curvature. \index[persons]{Rummler, H.} Rummler--\index[persons]{Sullivan, D.}Sullivan theorem \index[terms]{foliation!taut}
(\cite{_Rummler_,_Sullivan:taut_,_Haefliger:taut_} \index[terms]{theorem!Rummler--Sullivan}
states that $F$ is taut if and only if there
exists a relatively closed form on $M$ such that its restriction
to the leaves of $F$ is a volume form. In 
\cite{_Masa:duality_} it is shown that
$F$ is taut if and only if $H^q_b(M)$ is non-zero,
where $q= \codim F$. In this case, $H^q_b(M)$
is 1-dimensional, and the multiplication in basic cohomology defines 
a  Poincar\'e-type duality on $H^*_b(M)$.\index[terms]{duality!Poincar\'e}

Sasakian manifolds are clearly taut;
indeed, the contact form $\eta$ satisfies
$d\eta = \omega_0$, where $\omega_0$ is a
basic K\"ahler form\index[terms]{form!K\"ahler!basic} that vanishes on the
leaves of the Reeb foliation 
(\ref{_omega_0_trans_Kahler_Claim_}). In fact, any Riemannian foliation\index[terms]{foliation!taut}\index[terms]{foliation!Riemannian}
on a compact simply connected manifold
is taut (\cite{_Ghys:taut_}).

A cohomological obstruction to tautness is provided
by Jes\'us A. \'Alvarez L\'opez in \cite{_Alvarez:basic_}.
In this paper, \'Alvarez L\'opez\index[persons]{\'Alvarez L\'opez, J.} defines what
is now known as {\bf Alvarez class}, a closed
basic 1-form $\kappa_b$ that can be defined as the basic component
of the mean curvature form \index[terms]{form!curvature}of $F$ with respect to 
a bundle-like Riemannian structure. The cohomology
class $\kappa_b$ is independent on  the choice
of a Riemannian metric on $(M, F)$.
A Riemannian foliation $F$ is taut 
if and only if its \'Alvarez class vanishes.

Using the \'Alvarez 1-form $\kappa_b$, \index[persons]{Habib, G.} Habib and Richardson 
(\cite{_Habib_Richardson:basic_})
show that the space of basic harmonic forms
is isomorphic to the cohomology of ``twisted basic
differential'' $d+ \kappa\wedge$ on basic forms. 
When $F$ is taut, this result gives an isomorphism
between the basic harmonic forms and the basic cohomology.

Let $Q$ be a compact Sasakian manifold,\index[terms]{cohomology!basic}
and $H^*_b(Q)$ the basic cohomology associated with
the Reeb foliation. Since the Reeb foliation is
taut, Habib--Richardson theorem implies that
the basic cohomology of $Q$ is isomorphic to
the space of basic harmonic forms. We give an
independent proof of this result, that does
not use tautness and the \'Alvarez 1-form (\ref{_basic_harmo_Proposition_}).
Instead of the \'Alvarez 1-form, we use 
a collection of Killing fields\index[terms]{vector field!Killing}, tangent to the leaves,
such that the corresponding group action is transitive
on the leaves of the foliation.
This version of \index[persons]{Habib, G.} Habib--Richardson theorem\index[terms]{theorem!Habib--Richardson} 
is much weaker, but its proof is elementary,
and it can be applied to Vaisman manifolds
as well.\index[terms]{manifold!Vaisman}

\subsection{Hattori spectral sequence and Hattori
  differentials} \label{_filtered_spectral_sequence_Subsection_}

The notion of a spectral sequence is due to \index[persons]{Leray, J.} Leray and
\index[persons]{Serre, J.-P.} Serre, however, in this chapter we use a more abstract version of their
construction, called {\bf the spectral sequence of
a filtered complex}, that is  due to \index[persons]{Eilenberg, S.} Eilenberg and Moore
(\cite{_Eilenbrg_Moore:spectral_}). This is the point
where many topology textbooks start, as well, see for instance
\cite[Chapter 3]{_Fuks_Fomenko_}. A {\bf filtered  complex}\index[terms]{complex!filtered}
is a complex  $C^*=(C_k^*, d)$ equipped with a sequence of subcomplexes,
that is, with a sequence of $d$-invariant graded subgroups
$0=C_{-1}^*\subset C_0^*\subset C_1^* \subset ...\subset C_k^*$.
The corresponding spectral sequence is a sequence
of groups $E^{*,*}_r$, equipped with the differentials
${\goth d}_r:\; E^{p,q}_r\arrow E^{p-r,q+r+1}_r$ that satisfy
${\goth d}_r^2=0$ and $E_{r+1}^{*,*}= H^*(E^{*,*}_r, {\goth d}_r)$.
It starts from
$E^{p,q}_0=\frac{C_p^{p+q}}{C_{p-1}^{p+q}}$
and has as its limit $E_\infty^{p,q}=\frac{\im\big[(H^{p+q}(C_q)
  \to H^{p+q}(C)\big]}{\im\big[H^{p+q}(C_{q-1})
\to H^{p+q}(C)\big]}$. In \cite[Chapter 3]{_Fuks_Fomenko_}, 
this statement is called ``the \index[persons]{Leray, J.} Leray theorem''.\index[terms]{theorem!Leray}

Usually, the spectral sequence differentials
are denoted $d_r$; we use ${\goth d}_r$ 
to distinguish them from the \index[persons]{Hattori, A.} Hattori differentials
$d_i$ we use in the body of this chapter.\index[terms]{differential!Hattori}

To obtain this spectral sequence, we can use the following
trick. Assume that $C^*_i$ are complexes of vector spaces
over $\R$ or over $\C$, and introduce a compatible
Euclidean or Hermitian metric on these spaces.
Then $C^*=C_k^*$ can be decomposed into the orthogonal
direct sum $C_k^*= \bigoplus_{i=0}^k A_i^*$,
with $C_0^* = A_0^*$, $C_1^* =  A_0^*\oplus A_1^*$, ...,
$C_r^* = \bigoplus_{i=0}^r A_i^*$.
Let $d=d_0 + d_1 + ... + d_k$ be the
corresponding decomposition of the 
differential, with $d_r:\; A_i^* \arrow A_{i+r}^{*+1} $.
Then the partial sum $\delta_r:= \sum_{i=0}^r d_r$
satisfies $\delta_r^2=0$. The spectral
sequence differential ${\goth d}_r$
is obtained as a restriction of $\delta_r$
to the cohomology $E_r^{*,*}$ of ${\goth d}_{r-1}$.

The differential $d_r$ can be identified
with ${\goth d}_r$, after the following chain
of identifications. Using induction, we can assume
that $E_r^{*,*}$ is obtained by taking successive
cohomology of $\delta_i$, $i=0, 1, ..., r-1$.
Then $\delta_{r-1}$ vanishes on $E_r^{*,*}$,
hence $\delta_r= d_r + \delta_{r-1}$ restricted to
$E_r^{*,*}$ is equal to $d_r$. This means that
the page $E_r^{*,*}$ of the spectral sequence 
can be obtained by taking the cohomology
of $d_0$, then the cohomology of $d_1$, and
so on up to $d_{r-1}$. This identifies
the \index[persons]{Hattori, A.} Hattori differentials $d_r$ with the 
spectral sequence differentials ${\goth d}_r$.

\index[persons]{Leray, J.} Leray, \index[persons]{Serre, J.-P.} Serre and later \index[persons]{Eilenberg, S.} Eilenberg-Moore used
this construction to study the spectral sequence of a fibration
$\pi:\; E\arrow B$, often called {\bf the Leray-Serre spectral sequence.}
Let $B$ be a CW-space, $B_0$ a point, and $B_i$ the\index[terms]{spectral sequence!Leray-Serre}
$i$-skeleton of $B$. The homological Leray-Serre spectral sequence
is the spectral sequence of the cellular complex
of $C_*(E)$ filtered by $(C_i)_*(E)=C_*(\pi^{-1}(B_i))\subset C_*(E)$;
the cohomological Leray-Serre spectral sequence
is associated with the dual filtration on the dual
complex of cellular cochains. Instead of the
cellular chains and cochains, one could take
singular chains and cochains;  this would
result in a spectral sequence with the same groups $E^{p,q}_r$
for all $r \geq 2$. 

To do computations related to Vaisman and Sasakian
geometries, singular and cell complexes are not that
useful, because these geometries are defined
in terms of differential forms. This is why, instead
of the Leray-Serre spectral sequence, we
use the spectral sequence defined in
\cite{_Hattori_} by A. \index[persons]{Hattori, A.} Hattori (Subsection 
\ref{_Leray_Serre_Subsection_}).

Let $\pi:\; M \arrow B$ be a submersion of
smooth manifolds. The tangent bundle $TM$ 
contains the sub-bundle $T_\ver M$ of vertical tangent
vectors (that is, vectors tangent to the fibres of $\pi$).
Let $\Lambda_\hor^1(M)= \pi^* \Lambda^1 B$ be the space of
differential forms
vanishing on $T_\ver M$. Let $I_r^*\subset \Lambda^*(M)$
be the ideal in $\Lambda^*(M)$  generated
by $\Lambda_\hor^r(M)= \pi^* \Lambda^r B$,
and $C_p:= \bigoplus_{i=0}^p\Lambda^i(M)\oplus \bigoplus_{i=0}^n I_{i}^{p+i}$.
By construction, $C_p$ are differential
forms on $M$ which would vanish 
when contracted with $p+1$ vertical 
tangent vectors.

Clearly, $C_r$ is a differential ideal,
and $C_0 \subset C_1 \subset ... \subset C_k=
\Lambda^*(M)$, $k= \dim M - \dim B$
is a filtered complex. The spectral sequence 
associated with this filtered complex
is called {\bf the \index[persons]{Hattori, A.} Hattori spectral sequence}.
The corresponding direct sum decomposition\index[terms]{spectral sequence!Hattori}
is written explicitly in \eqref{_horis_vert_decompo_Equation_}:
\[
\Lambda^mM= \bigoplus_p \Lambda^p_\hor M \otimes \Lambda^{m-p}_\ver M, 
\]
with $\Lambda^p_\hor M \otimes \Lambda^{m-p}_\ver M$
being the orthogonal complement of $C_{m-p-1}$ in
$C_{m-p}$. 
The Hattori differentials $d_i$ are the components
of the de Rham differential which shift this grading 
by $i$, 
\[
d_i:\; \Lambda^p_\hor(M)\otimes \Lambda^{q}_\ver(M)\arrow 
\Lambda^{p+i}_\hor(M)\otimes \Lambda^{q+1-i}_\ver(M)
\]
The ${\goth d}_0$-differential of this spectral sequence 
takes $\Lambda^p_\hor(M)\otimes \Lambda^{q}_\ver(M)= E^{p,  q}_0$ 
to $\Lambda^p_\hor(M)\otimes\Lambda^{q+1}_\ver(M)$, and
the cohomology of $d_0={\goth d}_0$ is the fibrewise cohomology of the de
Rham algebra. This implies that the space $E_1^{p,q}$
is identified with $\Lambda^p_\hor(M)\otimes_{\R_B} R^q\pi_*\R_M$,
where $R^q\pi_*\R_M$ is the sheaf of cohomology of the
fibres of $\pi$. In other words, the
 $E_1^{*,*}$-term of this spectral sequence is identified with
the space of differential forms on the base with the coefficients
in the local system $R^q\pi_*\R_M$  of cohomology of the fibres.
This is similar to the Leray-Serre spectral
sequence, with the $E_1^{*,*}$-term identified
with singular cochains with coefficients
in $R^q\pi_*\R_M$. From this observation it
should be clear that the $E_i^{p,q}$-terms of the
\index[persons]{Hattori, A.} Hattori spectral sequence of a smooth fibration, 
are equal to the $E_i^{p,q}$-terms of the
Leray-Serre spectral sequence if $i >1$.

It makes sense to define the \index[persons]{Hattori, A.} Hattori spectral sequence 
in a more general context, when $M$ is a manifold equipped
with a foliation $F$. Locally $M$ admits
a smooth submersion to the leaf space of $F$.
The filtration
$C_0 \subset C_1 \subset ... \subset C_k=
\Lambda^*(U)$ is defined for
 each open set of $M$ where such a submersion
is defined. Since $C_r$ are subsheaves
on $\Lambda^*(M)$, this filtration gives
a filtration on $\Lambda^*(M)$, compatible
with the differential. We call the corresponding
spectral sequence ``the \index[persons]{Hattori, A.} Hattori spectral sequence''
as well.

This spectral sequence 
was re-invented independently on several instances after
\index[persons]{Hattori, A.} Hattori. In the book \cite[Section 1.6]{_Brylinski_} by J.-L. \index[persons]{Brylinski, J.-L.} Brylinski
it was described as ``\index[persons]{Cartan, H.}Cartan spectral sequence'', without 
reference. About the same time, it was described
in Vlad \index[persons]{Sergiescu, V.} Sergiescu's Ph. D. thesis  and in his subsequent papers
(\cite{_Sergiescu:thesis_,_El_Kacimi_Sergiescu_Hector_,_Sergiescu:Mexican_})
under the name ``Leray-Serre type spectral sequence''.
This work was quite influential, with a number of publications
citing Sergiescu's papers and his thesis 
(for example, \cite{_Alvares:finiteness_,_Alvares:duality_,_Domingues:finiteness_}). 
In \cite{_Alvares:finiteness_},
it was shown that all terms on $E_2$ page of this spectral
sequence are finite-dimensional when the foliation
admits a transversal Riemannian structure. In \cite{_Domingues:finiteness_} the same result
was proven for cohomology with coefficients in a local 
system. See also Subsection \ref{_Leray_Serre_Subsection_} 
for an alternative description of the Hattori spectral sequence.

As an application of this machinery, we use
the \index[persons]{Hattori, A.} Hattori spectral sequence to obtain the Griffiths
transversality relations (\ref{_Gauss_Manin_via_Hattori_Remark_}). 
However, the main application of the Hattori formalism
is an explicit expression for the Hattori differentials
on Sasakian manifolds,\index[terms]{condition!Griffiths transversality}
that are  used then to compute the cohomology.

Let $Q$ be a Sasakian manifold, and $R$ the Reeb
foliation, considered to be  a Riemannian foliation on $Q$.
We take an orthogonal decomposition of the
differential forms on $Q$ with respect to the
Hattori filtration, obtaining two terms,
\[ \Lambda^*(Q) = \Lambda^*_\hor(Q)\oplus \eta \wedge
\Lambda^*_\hor(Q).
\]
Under this decomposition, the de Rham differential decomposes
into three Hattori components, $d= d_0+d_1+d_2$.
In Subsection \ref{_LS_decompo_Sasakian_Subsection_}, we describe these components explicitly.
First, $d_0:\; \Lambda^p_\hor(Q)\arrow  \eta \wedge
\Lambda^p_\hor(Q)$ is the Lie derivative along the
Reeb field composed with the map $\alpha \mapsto \eta
\wedge\alpha$. Secondly,
$d_1:\; \Lambda^p_\hor(Q)\oplus \eta \wedge
\Lambda^p_\hor(Q) \arrow \Lambda^{p+1}_\hor(Q)\oplus \eta \wedge
\Lambda^{p+1}_\hor(Q)$ is the de Rham differential taken in
the direction orthogonal to the Reeb foliation.
Finally, $d_2:\;  \eta \wedge \Lambda^p_\hor(Q)\arrow \Lambda^{p+2}_\hor(Q)$
takes $\eta$ to the transversal K\"ahler form\index[terms]{form!K\"ahler!transversal}
$\omega_0\in \Lambda^{2}_\hor(Q)$ and
$\eta \wedge \alpha$ to $\omega_0\wedge \alpha$.

This description immediately implies an 
explicit formula for the co\-ho\-mo\-lo\-gy of $Q$.
Recall that {\bf the cone} of a morphism
of complexes $\phi: (A^*, d_A)\arrow (B^*, d_B)$
is the complex $\oplus_i C^i=A^{i-1} \oplus B^i$
with the differential $d_A - d_B + \phi$.
This gives a short exact sequence of complexes
\[ 
0 \arrow B^* \arrow C^* \arrow A^*[1]\arrow 0, 
\]
where $A^*[1]$ denotes the same complex shifted by 1.
The corresponding cohomology groups fit
into a long exact sequence
\begin{equation}\label{_cone_complex_long_exact_Equation_}
... \arrow H^i(B) \arrow H^i(C) \arrow H^{i+1}(A) \arrow H^{i+1}(B) \arrow ...
\end{equation}
Now, consider the $E^{*,*}_1$-term of the \index[persons]{Hattori, A.} Hattori spectral
sequence of a Sasakian manifold; it is identified with
the cohomology of $d_0= \Lie_\xi\circ L_\eta$,
where $L_\eta$ is the operation of multiplication by
$\eta$, and $\xi$ the Reeb field.
The cohomology of this operator is not hard to compute.
Since $\xi$ is Killing,\index[terms]{vector field!Killing} the closure $G$ of the group action
generated by $e^{t\xi}$ is compact. Therefore, $G$ is a
torus. Using the Fourier decomposition over $G$,
we decompose the space of forms on $M$
onto the weight components of $\xi$,
$\Lambda^*(M) = \hat \bigoplus_{u \in W} \Lambda^*(M)_u$
where $\hat \bigoplus$ is the completed direct sum,
$W\subset \C$ the set of weights of $\xi$, and
$\Lambda^*(M)_u= \{ \alpha \in \Lambda^*(M) \ \ | \Lie_\xi
\alpha=u\alpha\}$ (Exercise \ref{_Fourier_for_torus_on_de_Rham_Exercise_}).

Since $\Lie_\xi$ is surjective on $\hat \bigoplus_{u \neq 0} \Lambda^*(M)_u$,
the cohomology of $d_0= \Lie_\xi\circ L_\eta$ can be
identified with $\Lambda^*_\bas(Q)\oplus \eta \wedge\Lambda^*_\bas(Q)$, where $\Lambda^*_\bas(Q)$
denotes the space of basic forms on $Q$, that is, $\xi$-invariant
horizontal forms.  The action $d_1$ on
$\Lambda^*_\bas(Q)\oplus \eta \wedge\Lambda^*_\bas(Q)$
is identified with the basic differential (that is, the de Rham
differential restricted to the basic forms),
and $d_2(\eta \wedge \alpha) = \omega_0 \wedge \alpha$.
This implies that the complex
$(\Lambda^*_\bas(Q)\oplus \eta \wedge\Lambda^*_\bas(Q), d_1+d_2)$
can be identified with the cone of the morphism
$L_{\omega_0}:\; (\Lambda^*_\bas(Q), d)\arrow(\Lambda^{*+2}_\bas(Q), d)$
which takes a basic form and multiplies it by the transversal K\"ahler
form $\omega_0$. Then
\eqref{_cone_complex_long_exact_Equation_} 
gives the exact sequence
\begin{multline}\label{_Sasakian_cohomo_from_cone_intro_Equation_}
... \arrow H^i_\bas(Q) \stackrel {L_{\omega_0}}\arrow
H^{i+2}_\bas(Q) \arrow H^{i+2}(Q)  \arrow  \\ \arrow H^{i+1}_\bas(Q)  \stackrel {L_{\omega_0}}\arrow H^{i+3}_\bas(Q) \arrow ...
\end{multline}
The cohomology $H^*_\bas(Q)$ of $(\Lambda^*_\bas(Q), d)$ enjoys the same
properties as the cohomology of a K\"ahler manifold
(\ref{_Hodge_decompo_d_1_Sasa_Claim_}),
including the Lefschetz $\goth{sl}(2)$-action,
hence the natural map $L_{\omega_0}:\; H^i_\bas(Q)\arrow H^{i+2}_\bas(Q)$
is injective up to the middle degree, and surjective after
the middle. This implies that the arrow  $H^i_\bas(Q) \stackrel {L_{\omega_0}}\arrow
H^{i+2}_\bas(Q)$ is injective or surjective, depending on
$i$, and the exact sequence
\eqref{_Sasakian_cohomo_from_cone_intro_Equation_}
splits onto short exact sequences
\[
0 \arrow H^i_\bas(Q) \stackrel {L_{\omega_0}}\arrow
H^{i+2}_\bas(Q) \arrow H^{i+2}(Q)\arrow 0
\]
for $i\leq \frac{\dim_\R Q-1}{2}$, and
\[
0 \arrow H^{i+1}(Q)  \arrow H^{i}_\bas(Q)  \stackrel
{L_{\omega_0}}\arrow H^{i+2}_\bas(Q)  \arrow 0
\]
for $i\geq \frac{\dim_\R Q-1}{2}$.
The kernel of ${L_{\omega_0}}$ in the first exact sequence
is called {\bf the primitive part of the cohomology $H^*_\bas(Q)$},
the cokernel of ${L_{\omega_0}}$ in the second exact sequence
is called {\bf the coprimitive part of the cohomology $H^*_\bas(Q)$}.
We obtain that \eqref{_Sasakian_cohomo_from_cone_intro_Equation_}
is translated to isomorphisms between $H^*(Q)$ and
the primitive and coprimitive parts of $H^*_\bas(Q)$ (\ref{_Sasakian_cohomo_cone_Theorem_}).\index[terms]{manifold!Sasaki!quasi-regular}\index[terms]{manifold!Sasaki!regular}
When $Q$ is regular, or quasi-regular, the Reeb orbit
space is an orbifold, and $H^*_\bas(Q)$ is the cohomology of
this orbifold, that is  actually K\"ahler (and even projective).

These cohomology isomorphisms do not give the harmonic
decomposition yet. However, the harmonic decomposition\index[terms]{decomposition!harmonic}
easily follows from the dimension count if we realize each cohomology
class in $H^{i}(Q)$, $i\leq \frac{\dim_\R Q-1}{2}$, by a
primitive,\index[terms]{form!primitive} basic harmonic form, and 
 each cohomology class in $H^{i+1}(Q)$, $i\geq \frac{\dim_\R Q-1}{2}$,
by $\eta \wedge \alpha$ where $\alpha$ is a coprimitive,
basic harmonic form (\ref{_harmo_Sasa_decompo_Theorem_}). 

The same argument can be easily adjusted to obtain the
harmonic decomposition for Vaisman manifolds.\index[terms]{manifold!Vaisman}
First, we use an explicit form of \index[persons]{Hattori, A.} Hattori differentials
to express the cohomology in terms of the appropriate
cone of complexes (\ref{_Vaisman_coho_Theorem_}). 
Second, we realize the cohomology classes 
by basic harmonic forms and use
the dimension count (\ref{_Vaisman_harmonic_forms_Theorem_}).

\subsection[Supersymmetry and geometric structures on
  manifolds]{Supersymmetry and geometric structures\\ on
  manifolds}

The supersymmetry approach to Hodge theory 
is used in research literature to prove results
about the Hodge decomposition and harmonic
decomposition
(\cite{_Slupinski:Hodge_,_Verbitsky:HKT_,_Verbitsky:Vanishing_LCHK_,
_Cirici_Wilson:Dolbeault_,Ve2,_Verbitsky:G2_forma_,_CKT:cohomologies_,_Huang_Tan:parallel_}).
It was not used so far in basic textbooks on Hodge theory,
that is  a pity, because the supersymmetry simplifies
the deduction of standard identities and makes the
statements and the computations much simpler and cleaner.

We give a textbook treatment of K\"ahler identities
via supersymmetry in Section \ref{_superalgebra_Section_}.
We use this method, in particular, to prove the
K\"ahler identities for the basic cohomology of the
Vaisman and Sasakian manifolds, that are  transversally
K\"ahler by \ref{_omega_0_trans_Kahler_Claim_}.

The idea behind the superalgebra approach to the K\"ahler
identities is the following.
Consider the Lie superalgebra of endomorphisms of the
space $\Lambda^*(M)$ of differential forms on $M$
(Section \ref{_superalgebras_Subsection_}).
Then the K\"ahler identities can be understood
as identities in the Lie superalgebra\index[terms]{Lie superalgebra}
generated by the de Rham differential,\index[terms]{K\"ahler identities}
the Lefschetz $\goth{sl}(2)$-triple,\index[terms]{Lefschetz triple}
and the Weil operator $W_I$ acting
on $(p,q)$-forms as $\1(p-q)\Id$.
This algebra is, in fact, finite-dimensional,
and its relations are given by \ref{_kah_susy_Theorem_}.

A similar statement is true for \index[terms]{manifold!hyperk\"ahler}
hyperk\"ahler manifolds (Exercise
\ref{_supersymmetry_hyperkahler_Exercise_}).
For other geometric structures, one can define
an intrinsic finite-dimensional Lie superalgebra
containing an appropriate differential together
with its Hermitian adjoint operator (and hence,
containing the corresponding Laplacian).
If this Laplacian is central, the Lie superalgebra
actually acts on the cohomology.

One could think that this method could have been
used to produce the harmonic decomposition
(originally due to \index[persons]{Tachibana, S.} Tachibana) on Sasakian manifolds,
and its counterpart statement for
Vaisman manifolds,\index[terms]{manifold!Vaisman} due to \index[persons]{Kashiwada, T.} Kashiwada. However, this approach does not 
seem to work. Indeed, Tachibana and Kashiwada 
decomposition theorems work only for compact
manifolds, and the superalgebra action is local.

We were able to construct
a finitely generated, finitely presented algebra
of supersymmetry for Sasakian manifolds (\ref{_susy_Sasa_Theorem_}). The
full list of relations is not used in the proof of
\index[persons]{Tachibana, S.} Tachibana and \index[persons]{Kashiwada, T.} Kashiwada theorems, however, the
explicit description of the Hattori differentials
$d_0, d_1, d_2$ is heavily used in our proof of the harmonic
decomposition.

Unlike the K\"ahler and hyperk\"ahler supersymmetry
superalgebras, the Sasa\-kian supersymmetry algebra
is infinite-dimensional. However, it is finite-dimen\-sional
as an algebra over the ring $\R[\Lie_\xi]$ of polynomials, with the polynomial
action provided by the Lie derivative along the Reeb
field. 

The main difference between the Sasakian supersymmetry
algebra and the K\"ahler one, in addition to the
$\R[\Lie_\xi]$-action, is that it does not contain
the de Rham differential. Instead, it contains
the two \index[persons]{Hattori, A.} Hattori components $d_0$ and $d_1$ of the de Rham
differential, and the third component $d_2$
lies in its universal enveloping algebra.

Another major difference is that the anticommutator
of the transversal differential $d_1$ with itself is
non-zero, instead, we have $\{d_1, d_1\}=L_{\omega_0}
\circ \Lie_\xi$, which makes some of the K\"ahler
relations different and strange-looking.


\section[Lie superalgebras acting on the de Rham algebra]{Lie superalgebras acting on the de Rham\\ algebra}
\label{_superalgebra_Section_}


\subsection{Lie superalgebras and superderivations}
\label{_superalgebras_Subsection_}

In the following, all vector spaces and algebras are considered over $\R$. Let $A$ be a $\Z/2\Z$-graded vector space,
\[ A = A^\ev \oplus A^\odd.\] 
We define the structure of the {\bf Lie superalgebra} on $A$ axiomatically as follows.\index[terms]{Lie superalgebra}
We say that $a\in A$
is {\bf pure} if $a$ belongs to $A^\ev$ or $A^\odd$.
For a pure element $a\in A$, we
write $\tilde a =0$ if $a\in A^\ev$,
and $\tilde a =1$ if $a\in A^\odd$.
Consider a bilinear operator
\[ \{\cdot, \cdot\}:\; A \times A \arrow A,\]
called {\bf supercommutator}. Assume that\index[terms]{supercommutator} 
$\{\cdot, \cdot\}$ is {\bf graded an\-ti-\-com\-mu\-ta\-ti\-ve},
that is, satisfies
\[ \{ a, b\} = - (-1)^{\tilde a \tilde b}  \{ b,a \}
\]
for pure $a, b \in A$. Assume, moreover, that
$\{\cdot, \cdot\}$ is compatible with the grading:
the commutator $\{ a, b\}$ is even when both $a$, $b$ are even 
or odd, and odd if one of these elements is odd and another 
is even. We say that $(A, \{\cdot, \cdot\})$ is a {\bf Lie superalgebra} 
if the following identity (called {\bf the graded Jacobi
	identity}, or {\bf super Jacobi identity}) \index[terms]{graded Jacobi
	identity}
holds, for all pure elements $a, b, c\in A$: 
\begin{equation}
\{ a, \{ b, c\}\} = \{ \{ a, b\}, c\} + (-1)^{\tilde a \tilde b}\{b, \{ a,c \}\}.
\end{equation}
Up to a sign, this is the usual Jacobi identity.

Every reasonable property of Lie algebras has a 
natural analogue for Lie superalgebras, using the following  rule of thumb: 
every time one would exchange two elements $a$ and $b$,
one adds a multiplier $(-1)^{\tilde a \tilde b}$.

\hfill

\example 
Let $V= V^\ev\oplus V^\odd$ be a $\Z/2\Z$-graded vector space,
and $\End (V)$ its space of endomorphisms, equipped with the induced grading.
We define a supercommutator in $\End (V)$ by the formula:
\[ \{a, b\} = ab - (-1)^{\tilde a \tilde b}ba
\]
It is easy to check that $\left(\End (V), \{\cdot, \cdot\}\right)$
is a Lie superalgebra. 

\hfill

\remark
For a $\Z$-graded vector space $A$, one defines
$A^\ev$ as the direct sum of even components, and $A^\odd$ as
the direct sum of odd components. Then a $\Z$-graded
Lie superalgebra is given by a supercommutator on $A$
satisfying $\{A^p, A^q\}\subset A^{p+q}$ and satisfying the
graded Jacobi identity. In the sequel, all super Lie algebras
we consider are of this type. An endomorphism $u\in \End(A)$
is called {\bf even} if $u(A^\odd)\subset A^\odd$ and
$u(A^\ev)\subset A^\ev$, and {\bf odd} if
$u(A^\ev)\subset A^\odd$ and
$u(A^\odd)\subset A^\ev$. An endomorphism that is  either
odd or even is called {\bf pure}.

\hfill

\definition\label{_graded_commu_Definition_}
A graded algebra $A$ is called {\bf graded commutative}
if $\{a, b\} =0$ for all $a, b\in A$.

\hfill

The Grassmann algebra and de Rham algebra are clearly graded commutative.

\hfill

\definition\label{_graded_derivation_Definition_}
Let $\g$ be a graded commutative algebra.
A map $\delta:\; \g\arrow \g$ is called an {\bf even derivation}
if it is even 
and satisfies $\delta(xy)= \delta(x)y+ x\delta(y)$.
It is called an {\bf odd derivation} if it is odd and satisfies
$\delta(xy)= \delta(x)y+ (-1)^{\tilde x}x\delta(y)$.
It is called {\bf graded derivation} if 
it shifts the grading by $i$ and satisfies
$\delta(ab) = \delta(a) b + (-1)^{ij} a \delta(b)$,
for each $a \in A^j$. \index[terms]{graded derivation}

\hfill

\remark
The supercommutator of two derivations is again a derivation (Exercise \ref{_supercommutator_Exercise_}).
The\-refore, the derivations form a Lie superalgebra.

\subsection[Differential operators on graded commutative algebras]{Differential operators on graded commutative\\ algebras}
\label{_diffe_super_alge_Subsection_}

We need some algebraic results, that are  almost trivial,
and well-known for commutative algebras. We extend these statements
to graded commutative algebras; the proofs are the same
as in the commutative setting (see \cite{_co_}).

\hfill

\definition\label{_diff_ope_Definition_}
Let $A$ be a graded commutative algebra.
{\bf The algebra of differential operators} $\Diff(A)$ is an associative
subalgebra of $\End(A)$ generated by graded
derivations and $A$-linear self-maps.\index[terms]{operator!differential}
Let $\Diff^0(A)$ be the space of $A$-linear self-maps,
$\Diff^0(A)=A$, $\Diff^1(A)\supset \Diff^0(A)$
the subspace generated by $\Diff^0(A)$ and all graded derivations, and
$\Diff^i(A):= \Diff^{i-1}(A)\cdot \Diff^1(A)$.
This gives a multiplicative filtration
$\Diff^0(A)\subset \Diff^1(A) \subset \Diff^2(A)\subset \cdots$ 
The elements of $\Diff^i(A)$ are called 
{\bf differential operators of order $i$} on $A$.

\hfill

\claim\label{_diff_ope_commutator_Claim_}
Let $D\in \Diff^i(A), D'\in \Diff^j(A)$ be differential
operators on a graded commutative algebra. Then
$\{D, D'\} \in \Diff^{i+j-1}(A)$. 

\proof Since the 
commutator of two derivations is again a derivation, one has
$\{\Diff^1(A), \Diff^1(A)\}\subset \Diff^1(A)$.
Then we use induction in $i$ and the standard
commutator identities in the associative algebra. 
\endproof

\hfill

We shall apply this claim to geometric operations on
the de Rham algebra, obtaining differential operators of first order.
Note that the ``differential operator'' in the usual sense
is a different notion. For example, the  interior product operator $i_v$  of contraction 
with a vector field $v$ on a manifold $M$ is an odd derivation of the de Rham algebra, 
hence it is a first order differential operator; however,
$i_v$ is $C^\infty(M)$-linear, and thus it is not a differential operator in the usual sense.

\hfill

\claim\label{_first_order_vs_deriva_Claim_}
Let $D$ be a differential operator of first order on $A$.
Then $D(x)= D(1) x + \delta(x)$, where $\delta$ is a derivation.

\hfill

\proof This is seen by defining $\delta:=D-D(1)\cdot$, then observing that $\delta(1)=0$, which implies that $\delta\circ a-a\circ\delta=\delta(a)$ (we identify $a\in A$ with $a\in\End(A)$, $a(b)=ab$). We then have
$$\delta(ab)=(\delta\circ a)(b)=(a\circ\delta)(b)+\delta(a)(b)=a\delta(b)+\delta(a)b,$$
proving that $\delta$ is a derivation.
\endproof

\hfill

We use this formalism to compare first order 
differential operators on $A$ as follows.

\hfill

\claim\label{_deriva_dete_Claim_}
Let $\delta$ be a derivation on $A$. Then
$\delta$ is uniquely determined by the values it takes
on any set of multiplicative generators of $A$.
\endproof

\hfill

\claim\label{_first_or_dete_Claim_}
Let $D$ be a first order differential operator on $A$. Then
$D$ is uniquely determined by $D(1)$ and the values it takes
on any set of multiplicative generators of $A$.

\proof By \ref{_first_order_vs_deriva_Claim_}, 
$D-D(1)$ is a derivation, and hence  \ref{_deriva_dete_Claim_} implies 
\ref{_first_or_dete_Claim_}.
\endproof

\hfill

\corollary\label{_d^2=0_dete_Corollary_}
Let $d_1, d_2:\; \Lambda^*M \arrow \Lambda^{*+1}M $ be first order 
differential operators\footnote{Here, as elsewhere, ``differential operators
	on the de Rham algebra'' are understood in the algebraic sense, as in \ref{_diff_ope_Definition_}.} 
on  the de Rham algebra of a manifold $M$,
satisfying $d_1^2=d_2^2=0$, and $V\subset \Lambda^1M$ a subspace such that
the space $d_1(C^\infty (M)) + V$  generates $\Lambda^1M$
as a $C^\infty (M)$-module. Suppose that $d_1\restrict {C^\infty (M)}=d_2\restrict {C^\infty (M)}$
and $d_1\restrict {V}=d_2\restrict {V}$. Then $d_1=d_2$.

\proof Clearly, $d_1=d_2$ on $C^\infty (M)\oplus d_1(C^\infty (M))$.
Then $d_1=d_2$ on the set of multiplicative generators 
$C^\infty (M)+ d_1(C^\infty (M))+ V$, and \ref{_first_or_dete_Claim_}
implies \ref{_d^2=0_dete_Corollary_}. \endproof

\hfill

The following claim is used many times in the sequel.

\hfill

\claim\label{_d^2=0_super_Claim_}
Let $d$ be an odd element in a Lie superalgebra $\goth h$,
satisfying $\{d,d\}=0$. Then
$\{d, \{d, u\}\}=0$ for any $u\in {\goth h}$.

\proof By the super Jacobi identity,
\begin{equation}\label{_odd_squa_Jaco_Equation_}
\{d, \{d, u\}\} = - \{d, \{d, u\}\} + \{\{ d, d\}, u\} = - \{d, \{d, u\}\}.
\end{equation}
\endproof

This claim is a special case of the following: 

\hfill

\claim\label{_d^2neq0_super_Claim_}
Let $d\in \goth h$ be an odd element. Then
$2\{d, \{d, u\}\}= \{\{d,d\}, u\}$ for any $u\in {\goth h}$.

\proof Follows from \eqref{_odd_squa_Jaco_Equation_}.
\endproof

\subsection{Supersymmetry on K\"ahler manifolds}
\label{_susy_Kah_Subsection_}

One of the purposes of this paper is to obtain a natural superalgebra
acting on the de Rham algebra of a Sasakian manifold.
This is modeled on the superalgebra of a K\"ahler manifold,
generated by the de Rham differential, Lefschetz triple,
and other geometric operators. To make the analogy\index[terms]{supersymmetry algebra}
more clear, we recall the main results on the supersymmetry 
algebra of K\"ahler manifolds. We follow \cite[Section 1.3]{_Verbitsky:HKT_}.

\hfill

Let $(M, I, g, \omega)$ be a K\"ahler manifold. Consider $\Lambda^*M$ as a graded
vector space. The differentials $d, d^c:=  IdI^{-1}$
can be interpreted as odd elements in $\End(\Lambda^*M)$,
and the Hodge operators $L, \Lambda, H$ as even elements. 
To avoid confusion, we denote the supercommutator of odd
elements as $\{\cdot,\cdot\}$. In terms of the associative
algebra, $\{a,b\} = ab +ba$.
Let $d^*:=  \{\Lambda, d^c\}$, $(d^c)^*:= -\{\Lambda, d\}$.
The usual Kodaira relations can be stated as follows
\begin{equation}\label{_Kodaira_Equations_}
\begin{aligned} 
{} & \{ L, d^*\} = - d^c,\ \ \ \  \{ L, (d^c)^*\} = d, \ \ \ \ 
\{ d, (d^c)^*\} = \{ d^*, d^c\} =0, \\
{} & \{ d, d^c\} =  \{ d^*, (d^c)^*\} =0,\ \ \ \  
\{ d, d^*\} = \{ d^c, (d^c)^*\} = \Delta,
\end{aligned}
\end{equation}
where $\Delta$ is the Laplace operator,
commuting with $L, \Lambda, H$, and
$d$, $d^c$.

\hfill

\definition \label{_KdR_Definition_}
Let $(M, I, g, \omega)$ be a K\"ahler manifold.
Consider the Lie superalgebra $\ga \subset \End(\Lambda^*M)$
generated by the following operators:
\begin{enumerate}
	\item $d$, $d^*$, $\Delta$, constructed out of the Riemannian metric.
	
	\item $L(\alpha):= \omega\wedge \alpha$.
	
	\item $\Lambda(\alpha) := * L * \alpha$.  It is easily seen that $\Lambda= L^*$.
	
	\item The Weil operator $W\restrict{\Lambda^{p,q}M}=\1(p-q)$.\index[terms]{operator!Weil}
\end{enumerate}
This Lie superalgebra is called {\bf the algebra of supersymmetry
	of the K\"ahler manifold}.

\hfill

Using \ref{_kah_susy_Theorem_} below,
is easy to see that
$\goth a$ is in fact independent on  $M$.
This Lie superalgebra was studied from the physicists'
point of view in \cite{_FKS_}.

\hfill

\theorem\label{_kah_susy_Theorem_}
Let $M$ be a K\"ahler manifold, and 
$\goth a$ its  supersymmetry algebra
acting on $\Lambda^*M$. Then
$\goth a$ has dimension $(5|4)$
(that is, its odd part is 4-dimensional, and its even part
is 5-dimensional). The odd part is generated by
$d$, $d^c=IdI^{-1}$, $d^*$, $(d^c)^*$, 
while the even part\index[terms]{Lefschetz triple}
is generated by the Lefschetz triple
$L, \Lambda, H:=[L, \Lambda]$, the Weil operator $W$ and the
Laplacian $\Delta=\{d, d^*\}$.
Moreover, the Laplacian $\Delta$ is
central in $\goth a$, and hence  $\goth a$ also acts on the
cohomology of $M$. The following are the only non-zero
commutator relations in $\goth a$:
\begin{enumerate}
	\item $\gsl(2)$-relations in $\langle L, \Lambda, H\rangle$:
	\[
	[H, L]=2L, \ \  [H, \Lambda]=2\Lambda,\ \  [L, \Lambda]=H.
	\]
	For any operator $D$ of grading $k$, one has $[H, D]=kD$.
	
	\item The Weil operator acts as a complex structure on the odd part of $\goth a$:
	\[ [W, d]=d^c, \ \   [W, d^c]=-d, \ \ [W, d^*]=-(d^c)^*, \ \   [W, (d^c)^*]=d^*.\]
	\item The K\"ahler-Kodaira relations between the differentials and the Lefschetz operators:
	\begin{equation}\label{_Kodaira_rel_Equation_}
	[\Lambda, d] = (d^c)^*,\ \  [ L, d^*] = - d^c,\ \ [\Lambda, d^c] = - d^*, \ \ [ L, (d^c)^*] = d.
	\end{equation}
	\item Almost all odd elements supercommute, with the only exception
	\[
	\Delta=\{d, d^*\}=\{d^c, (d^c)^*\},
	\]
	and $\Delta$ is central. In other words, the odd elements of $\goth a$
	generate the odd Heisenberg superalgebra\footnote{For a definition of
odd Heisenberg algebra, see \ref{_odd_Heisenberg_Claim_}.}.
\end{enumerate}\index[terms]{Lie superalgebra!Heisenberg}

\proof
These relations are standard in algebraic geometry\index[terms]{geometry!algebraic} (see {e.g.} \cite{griha}), but probably the easiest way of proving them is using the results about Lie superalgebras
collected in Subsection \ref{_diffe_super_alge_Subsection_}.

\hfill

{\bf Proof of the Lefschetz $\gsl(2)$-relations:}
These relations would follow if we prove that
$H:=[L, \Lambda]$ acts on $p$-forms by multiplication
by $p-n$, where $n=\dim_\C M$. Since $L, \Lambda, H$ are
$C^\infty(M)$-linear, it would suffice to prove these
relations on a Hermitian vector space. 

Let $V$ be a real vector space equipped with
a scalar product, and fix an orthonormal basis $\{v_1, \ldots, v_{m}\}$.
Denote by $e_{v_i}:\; \Lambda^k V  \arrow \Lambda^{k+1} V$
the operator of multiplication, $e_{v_i}(\eta) = v_i \wedge \eta$.
Let $i_{v_i}:\; \Lambda^k V  \arrow \Lambda^{k-1} V$
be the operator of contraction with $v_i$.
The following claim is clear.

\hfill

\claim\label{_odd_Heisenberg_Claim_}
The operators $e_{v_i}$, $i_{v_i}$, $\Id$ form a basis of the {\bf 
	odd Heisenberg Lie superalgebra}, 
with the only non-trivial\index[terms]{Lie superalgebra!Heisenberg} 
supercommutator given by the formula $\{ e_{v_i}, i_{v_j}\} = \delta_{i,j}\Id$.
\endproof

\hfill

Now let  $V$ be an even-dimensional real vector space equipped with
a scalar product, and $\{x_1, \ldots, x_n, y_1, \ldots, y_{n}\}$ an orthonormal basis.
Consider the complex structure operator $I$ such that
$I(x_i)=y_i$, $I(y_i)=-x_i$. 
The fundamental symplectic form is given
by $\sum_i x_i \wedge y_i$, and hence 
\[
L= \sum_i e_{x_i} e_{y_i}, \ \ \ \Lambda= \sum_i i_{x_i} i_{y_i}.
\]
Clearly, for any odd elements $a, b, c, d$ such that
$\{a,b\}= \{a, d\} = \{b, c\}=\{c,d \}=0$, 
one has
$\{ab, cd\}= -\{a, c\}bd + ca \{b,d\}$. 
Then 
\begin{equation*}
\begin{split}
[L, \Lambda] &= 
\left [\sum_i e_{x_i} e_{y_i}, \sum i_{x_i} i_{y_i}\right]
=  \sum_{i=1}^{n} e_{y_i} i_{y_i} -\sum_{i=1}^{n}i_{x_i} e_{x_i}\\
& =
\sum_{i=1}^{n} e_{y_i} i_{y_i} -\sum_{i=1}^{n}(e_{x_i} i_{x_i} -1).
\end{split}
\end{equation*}
This term, applied to a monomial $\alpha$ of degree $d$,
would give $(d-n)\alpha$. This proves the Lefschetz $\gsl(2)$-relations.

\hfill

{\bf Proof of the relations between the Weil operator $W$ and the odd part of $\goth a$.}
Clearly, it is enough to prove $[W, d]=d^c$, the remaining relations follow by duality or
by complex conjugation. Writing the Hodge components of $d= d^{1,0}+d^{0,1}$, with
$d^{1,0}=\frac{d + \1 d^c}2$ and $d^{0,1}=\frac{d -\1 d^c}2$ ,
we obtain $[W, d]=\1 d^{1,0}-\1 d^{0,1}= d^c$. 

\hfill

{\bf Proof of the K\"ahler-Kodaira relations\index[terms]{relations!K\"ahler-Kodaira} between the Lefschetz $\gsl(2)$-o\-pe\-ra\-tors and the
	odd part of $\goth a$.} As before, it is enough to prove $[ L, d^*] = - d^c$,
the remaining K\"ahler-Kodaira relations follow by duality or
by complex conjugation. The operator $L$ is $\Lambda^*(M)$-linear, and hence 
it is a differential operator of order 0. The operator
$d^*$ can be written in a frame $\{v_i \}$ of $TM$ 
as 
\begin{equation} \label{_d^*_nabla_Equation_}
d^*(\eta)= \sum_i i_{v_i}\nabla_{v_i}\eta, 
\end{equation}
where $\nabla$ is the Levi--Civita connection of the metric $g$.
Since $\nabla_{v_i}$ and $i_{v_i}$ are
both derivations of the de Rham algebra, their product is
an order 2 differential operator (in the algebraic sense
as given by \ref{_diff_ope_Definition_}). From 
\ref{_diff_ope_commutator_Claim_} it follows that
$[ L, d^*]$ is a first order operator.

We shall prove $[ L, d^*] = - d^c$ by applying 
\ref{_d^2=0_dete_Corollary_}. First, let us show that
$[L,[L, d^*]]=0$. Clearly, $[\Lambda, d^*]= ([L, d])^*=0$,
and $[H, d^*]=-d^*$, and hence  $d^*$ is the lowest weight vector in
a weight 1 representation of $\goth{sl}(2)$. This gives
$[L,[L, d^*]]=0$. Now, the super Jacobi identity gives
\[
\{\{L, d^*\},\{L, d^*\}\}= \{L,\{d^*,\{L, d^*\}\}+ \{d^*,\{L,\{L, d^*\}\}.
\]
The first term in the RHS vanishes by \ref{_d^2=0_super_Claim_}, and 
the second term vanishes because $[L,[L, d^*]]=0$.
Then $([L, d^*])^2=0$. Clearly, $d^c(C^\infty (M))$ generates $\Lambda^1(M)$ over $C^\infty (M)$.
To deduce $[ L, d^*] = - d^c$ from  \ref{_d^2=0_dete_Corollary_}, it remains to show
that $[ L, d^*]\restrict{C^\infty (M)} = - d^c\restrict{C^\infty (M)}$.
This is clear for the following reason. For any function
$f\in C^\infty (M)$, one has $[ L, d^*] f= - d^*(f \omega)$. Writing
$d^*= \sum_i i_{v_i}\nabla_{v_i}$ as in \eqref{_d^*_nabla_Equation_},
and using $\nabla\omega=0$, we obtain
\[ 
d^*(f \omega)= \sum_i i_{v_i}\nabla_{v_i}(f\omega) = \sum_i \Lie_{v_i}(f) i_{v_i}(\omega)=
\sum_i - \Lie_{v_i}(f) I(v_i)= -d^cf.
\]
This finishes the proof of $[ L, d^*] = - d^c$.

\hfill

{\bf Proof of the commutator relations between the odd part of $\goth a$.}
We have already shown that $d^c=[W, d]$. Then $\{d, d^c\}= \{d,\{W, d\}\}=0$
by \ref{_d^2=0_super_Claim_}. Similarly, $d^*= - [\Lambda, d^c]$, giving
$\{d^c, d^*\}=0$. The relation $\{(d^c)^*, d\}=0$ is obtained by duality.
Finally, $\{d, d^*\}=\{d^c, (d^c)^*\}$ is obtained by
applying $[\Lambda, \cdot]$ to $\{d, d^c\}=0$. Using 
the K\"ahler-Kodaira relations \eqref{_Kodaira_rel_Equation_}, we obtain
\[
0=\{\Lambda, \{d, d^c\}\}= \{\{\Lambda,d\}, d^c\}+ \{d, \{\Lambda,d^c\}\}=
\{(d^c)^*,  d^c\}-\{d, d^*\},
\]
giving $\{d, d^*\}=\{d^c, (d^c)^*\}$.

We proved all the relations in the K\"ahler supersymmetry algebra
$\goth a$, finishing the proof of \ref{_kah_susy_Theorem_}.\index[terms]{supersymmetry algebra!K\"ahler}
\endproof

\hfill

In Section \ref{_Appendix_C_} we shall develop similar relations
for the superalgebra associated with a Sasakian manifold. 


\section{Hattori differentials on Sasakian manifolds}


\subsection{Hattori spectral sequence and associated differentials}
\label{_Leray_Serre_Subsection_}

Let $\pi:\; M \arrow B$ be a smooth fibration,
and $F_k\subset \Lambda^*M$ be the  ideal generated
by $\pi^*\Lambda^kB$. The \index[persons]{Hattori, A.} Hattori spectral sequence  
(\cite{_Hattori_}) is the spectral sequence 
associated with this filtration.\index[terms]{spectral sequence!Hattori}

The $E_1^{p,q}$-term of this sequence is $\Lambda^pB\otimes_\R R^q\pi_* \R_M$,
where $R^*\pi_* \R_M$ is the local system of cohomology of the fibres, and
the $E_2^{p,q}$-term is $H^p(B, R^q\pi_* \R_M)$. We explain its construction
below, in a more general context.

Let $M$ be a manifold equipped
with an integrable distribution $F\subset TM$.
The \index[persons]{Hattori, A.} Hattori spectral sequence is defined in this generality,
similar to the Hattori spectral sequence for fibred bundles.
This is done as follows.

\hfill

\definition
Let $M$ be a manifold, and $F\subset TM$ an
integrable distribution. A $k$-form $\eta\in \Lambda^* M$
is called {\bf basic} if for any vector field $v\in F$,
one has $\Lie_v\eta=0$ and $i_v(\eta)=0$.\index[terms]{form!basic}

\hfill

Let $\pi:\; M \arrow B$ be the projection of $M$ onto the leaf
space. Locally, such a projection is always defined. By \ref{_basic_forms_Frobenius_lifted_Theorem_}, 
a form is basic if and only if it is the pullback of a form on $B$.

\hfill

If $F$ is the tangent bundle of the fibres of a
fibration $\pi:\; M \arrow B$, then the space
of basic forms is $\pi^*\Lambda^kB$.
We are going to produce a spectral sequence which gives
the standard \index[persons]{Hattori, A.} Hattori spectral sequence when $F$ 
is tangent to the leaves of a fibration.

Define {\bf the Hattori filtration associated with $F$} as   
$F_n\subset F_{n-1}\subset \cdots$ 
by putting $F_k\subset \Lambda^*M$, where $F_k$ is the
ideal generated by the basic $k$-forms.

When $M$ is Riemannian and the metric is compatible with the foliation, 
this filtered bundle is decomposed 
into the direct sum of subquotients. This
gives a decomposition of the de Rham differential,
$d= d_0 + d_1 + d_2+\cdots + d_{r+1}$, where $r=\rk F$
with each successive piece associated with the differential
in $E_k^{p,q}$, as follows.

Using the metric, we split the cotangent bundle into orthogonal complements
as $\Lambda^1M= \Lambda^1_\hor M \oplus \Lambda^1_\ver M$,
where $\Lambda^1_\hor(M)$ is generated by basic 1-forms, and
$\Lambda^1_\ver M= F^*$ is its orthogonal complement.
Denote by $\Lambda^p_\hor M$, $\Lambda^q_\ver M$
the exterior powers of these bundles.

This gives the following splitting of the de Rham algebra of 
$M$:
\begin{equation}\label{_horis_vert_decompo_Equation_}
\Lambda^mM= \bigoplus_p \Lambda^p_\hor M \otimes \Lambda^{m-p}_\ver M, 
\end{equation}
with $F_{p}/F_{p-1}= \bigoplus_q \Lambda^p_\hor(M)\otimes \Lambda^{q}_\ver(M)$.

Consider the associated
decomposition of the de Rham differential,
$d= d_0 + d_1 + d_2+ \cdots + d_{r+1}$, where $r=\rk F$,
and 
\[ d_i:\; \Lambda^p_\hor(M)\otimes \Lambda^{q}_\ver(M)\arrow 
\Lambda^{p+i}_\hor(M)\otimes \Lambda^{q+1-i}_\ver(M).
\]
The terms $d_{i}$ vanish for $i> r+1$ because
$F$ is $r$-dimensional, and hence  for $i>r+1$ either $\Lambda^{q+1-i}_\ver(M)$
or $\Lambda^{q}_\ver(M)$ is 0.

These differentials are related to the differentials in
the \index[persons]{Hattori, A.} Hattori spectral sequence in the following way:
to find $E_1^{p,q}$, one takes the cohomology of $d_0$.
Then one restricts $d_1$ to $E_1^{p,q}$, and its cohomology
gives $E_2^{p,q}$, and so on. In this sense, the
Hattori differentials are indeed differentials
in the Hattori spectral sequence.

\hfill

\remark 
Each of the differentials $d_i$ is a derivation, because
the decomposition $\Lambda^*(M)= \bigoplus_{p,q}\Lambda^p_\hor(M)\otimes \Lambda^q_\ver(M)$
is multiplicative.

\hfill

For the next remark, we will need the following lemma.

\hfill

\lemma\label{_d_0_horizontal_commute_Lemma_}
Let $\pi:\; M \arrow B$ be a proper locally trivial fibration.
Consider a vector field $X\in TB$, and let
$\tilde X\in T_\hor M$ be its horizontal lift.
Then the contraction $i_{\tilde X}$ super-commutes with $d_0$.

\hfill

\proof
Since $d_0$ and $i_{\tilde X}$ are odd derivations,
their anticommutator $\{d_0, i_{\tilde X}\}$
is also a derivation. Clearly, this operator
vanishes on functions. To prove
\ref{_d_0_horizontal_commute_Lemma_} 
we use \ref{_deriva_dete_Claim_};
it would suffice to find a subspace
$V\subset \Lambda^1 M$ such that
$V$ and $C^\infty M$ generate $\Lambda^*M$
and $\{d_0, i_{\tilde X}\}(V)=0$.
Let $V$ be generated by closed basic
1-forms and fibrewise closed horizontal 1-forms.
Clearly, $C^\infty M\cdot V= \Lambda^1 M$.
On the other hand, $d_0$ vanishes on $V$
because all elements of $V$ are closed, and
$i_{\tilde X}$ vanishes on horizontal 1-forms.
To prove that $\{d_0, i_{\tilde X}\}(V)=0$
it remains to show that $d_0 i_{\tilde X}(\eta)=0$
for any closed basic 1-form $\eta= \pi^*\eta_0$. However,
$i_{\tilde X}\eta= \pi^*(i_X\eta_0)$, and $d_0$
is the vertical component of $d$, and hence  it vanishes
on $\pi^* C^\infty B$.
\endproof

\hfill

\remark\label{_Gauss_Manin_via_Hattori_Remark_}
\index[terms]{connection!Gauss--Manin}
The Gauss-\index[persons]{Manin, Yu. I.}Manin connection on the bundles of cohomology can
be expressed using the \index[persons]{Hattori, A.} Hattori differentials. This approach
can be used to prove the \index[persons]{Griffiths, P.} Griffiths transversality relations 
\eqref{_Griffiths_tranversality_Equation_}. \index[terms]{condition!Griffiths transversality}
Let $\pi:\; M \arrow B$ be a proper locally trivial fibration,
and $\eta$ a fibrewise closed $k$-form. The condition of being
fibrewise closed is equivalent to $d_0(\eta)=0$, where $d_0$ is the
Hattori differential.  This form gives a section of the fibrewise
cohomology bundle $R^k\pi_*\R_M\otimes_\R C^\infty B$.
Its Hattori differential $d_1(\eta)$ is $d_0$-closed,
because $\left(\sum_i d_i\right)^2=0$ and $d_0 \eta=0$.
This implies that $d_1\eta$ is fibrewise closed.
Consider a vector field $X$ on $B$, and let
$\tilde X\in T_\hor M$ be its horizontal lift.
By \ref{_d_0_horizontal_commute_Lemma_},
the differential $d_0$ commutes with $i_{\tilde X}$,
hence $i_{\tilde X}(d_1\eta)$ is fibrewise closed.
The Gauss--Manin connection $\nabla$
can be defined by setting $\nabla_X([\eta]):= [i_{\tilde X}(d_1\eta)]$,
where $[\cdot]$ denotes the fibrewise cohomology class.
To see that this definition is equivalent to the standard one,
we notice that this is indeed a connection, and it 
vanishes on any section of the fibrewise
cohomology bundle $R^k\pi_*\R_M\otimes_\R C^\infty B$
which comes from a closed form on $M$.
A connection that vanishes on the set
of horizontal sections of a flat connection $\nabla$
coincides with $\nabla$, and hence  the connection which
we defined coincides with the Gauss--Manin connection.

\subsection{Hattori differentials on Sasakian
	manifolds}
\label{_LS_decompo_Sasakian_Subsection_}

Now let  $Q$ be a Sasakian manifold, $\xi $ its Reeb field, normalized
in such a way that $|\xi |=1$, and
$R\subset TQ$ the 1-dimensional foliation generated by the Reeb field.
The corresponding \index[persons]{Hattori, A.} Hattori differentials are written as
$d=d_0 + d_1 + d_2$, because $R$ is 1-dimensional. 
Since $d^2=0$, one has $d_0^2=d_2^2=0$ and
$\{d_0, d_2\}=-\{d_1, d_1\}$. The differentials $d_0, d_2$
can be described explicitly as follows.

\hfill

\claim\label{_d_0_expli_Claim_}
Let $e_\xi :\; \Lambda^*(Q)\arrow \Lambda^{*+1}(Q)$ be the operator of multiplication
by the form $\xi ^{\,\flat}=\eta$ dual to $\xi $. Then $d_0= e_\xi  \Lie_\xi $.

\proof 
Locally, the sheaf $\Lambda^1(Q)$ is generated over $C^\infty (Q)$ by 
the basic 1-forms and $\xi ^{\,\flat}$. Clearly, the differential of a basic
1-form $\alpha$ belongs to $\Lambda^2_\bas(Q)$, and hence  $d_0(\alpha)=0$,
and $e_\xi  \Lie_\xi (\alpha)=0$. Also, 
$d_0(\xi ^{\,\flat})= e_\xi  \Lie_\xi (\xi ^{\,\flat})=0$. 
By \ref{_deriva_dete_Claim_}, to prove $d_0=e_\xi  \Lie_\xi $
it remains to show that these two operators are equal on 
$C^\infty (Q)$. However, on $C^\infty (Q)$, 
we have $d_0= e_\xi  \Lie_\xi $, because $d_0f$ is the orthogonal
projection of $df$ to $\Lambda^1_\ver(Q)=R^*$ generated
by $e_\xi $, for all $f\in C^{\infty}(Q)$.
\endproof

\hfill

Recall that the space of leaves of $R$ on a regular Sasakian manifold is equipped with \index[terms]{manifold!Sasaki!quasi-regular}
a complex and K\"ahler structure (\ref{regsas}). The corresponding 
K\"ahler structure can be described very explicitly.

\hfill

\claim\label{_omega_0_trans_Kahler_Claim_}
Let $Q$ be a  Sasakian manifold, and $\xi ^{\,\flat}$ its contact form\footnote{Usually, $\xi^{\, \flat}$ is denoted $\eta$, but in the context of this chapter we find the notation $\xi ^{\,\flat}$ more transparent.}.
Then $\omega_0:=d(\xi ^{\,\flat})$ is basic with respect to $R$,
and defines a transversally K\"ahler structure, that is, the 
K\"ahler structure on the space of leaves of $R$ (see \ref{_transvers_Kahler_Definition_} below).

\hfill

\proof Replacing $X$ by a smaller open subset, if necessary,
we may assume that the Reeb action is regular, that is,
the  space of orbits  $X=Q/\xi $ of the Reeb action is a manifold,
and the corresponding quotient map $Q \arrow X$ is smooth.
 Denote by $C(Q)$ the Riemannian cone\index[terms]{cone!Riemannian} of $Q$.
Then $X$ is a local quotient of $C(Q)$ by the holomorphic $\C$-action\index[terms]{action!$\C$-} generated
by $r=t\frac{d}{dt}$ and $\xi = I(r)$.
Therefore,
$X$ is a complex manifold (as  a quotient of a complex
manifold by a holomorphic action of a Lie group).
It is K\"ahler by \ref{_Sigma_defi_Claim_}. \endproof

\hfill

\proposition\label{_d_2_expli_Proposition_}
Let $L_{\omega_0}:\; \Lambda^*(Q)\arrow \Lambda^{*+2}(Q)$ be the operator of multiplication
by the transversal K\"ahler form\index[terms]{form!K\"ahler!transversal} $\omega_0=d(\xi ^{\,\flat})$, and
$i_\xi $ the contraction with the Reeb field. Then
$d_2 = L_{\omega_0} i_\xi $.

\hfill

\proof Clearly, the 
Hattori differentials 
$$d_i:\; \Lambda^p_\hor(Q)\times \Lambda^{q}_\ver(Q)\arrow 
\Lambda^{p+i}_\hor(Q)\otimes \Lambda^{q+1-i}_\ver(Q)$$
vanish on $\Lambda^0(Q)$ unless $i=0,1$. Therefore, the differentials
$d_2, d_3, ...$ are always $C^\infty(Q)$-linear.
By \ref{_deriva_dete_Claim_}, it only remains to show  that $d_2 = L_{\omega_0} i_\xi $ on 
some set of 1-forms generating $\Lambda^1(Q)$ over $C^\infty (Q)$.

Clearly, on $\Lambda^1_\hor(Q)$ the differential
$d_2$  should act as 
$$d_2:\;\Lambda^1_\hor(Q)\arrow \Lambda^{3}_\hor(Q)\otimes \Lambda^{-1}_\ver(Q),$$
hence $d_2\restrict{\Lambda^1_\hor(Q)}=0$.
In order to prove \ref{_d_2_expli_Proposition_} it remains to show that 
$d_2(\xi ^{\,\flat})=L_{\omega_0} i_\xi (\xi ^{\,\flat})$.
However, $d_2(\xi ^{\,\flat})$ is the $\Lambda^2_\hor(Q)$-part of  $d(\xi ^{\,\flat})$,
thus $d_2(\xi ^{\,\flat})=d(\xi ^{\,\flat})=\omega_0$, and $L_{\omega_0} i_\xi (\xi ^{\,\flat})=\omega_0$
because $i_\xi (\xi ^{\,\flat})=1$.
\endproof

\hfill

The Hattori differential $d_1$ is, heuristically speaking, the ``transversal component''
of the de Rham differential. Indeed, $d_1(\eta)=d(\eta)$ for any basic form $\eta$.
Since the leaf space of $R\subset TQ$ is equipped with a complex structure, it is natural to expect that
the Hodge components of $d_1$ have the same properties as the Hodge components
of the de Rham differential on a complex manifold.

\hfill

\claim\label{_Hodge_decompo_d_1_Sasa_Claim_}
Let $Q$ be a Sasakian manifold, and 
\[ \Lambda^m(Q)= \bigoplus_p \Lambda^p_\hor(Q)\otimes \Lambda^{m-p}_\ver(Q)\] the decomposition
associated with $R\subset TQ$ as in \eqref{_horis_vert_decompo_Equation_}.
Using the complex structure on the basic forms, consider the
Hodge decomposition \[ \Lambda^m_\hor(Q)=\bigoplus_p \Lambda^{m-p,p}_\hor(Q).\]
Then the differential \[ d_1:\; \Lambda^p_\hor(Q)\otimes \Lambda^{q}_\ver(Q)\arrow 
\Lambda^{p+1}_\hor(Q)\otimes \Lambda^{q}_\ver(Q)\] 
has two Hodge components $d_1^{0,1}$ and $d_1^{1,0}$.
Moreover, the differential $d_1^c:=d_1^{0,1}-d_1^{1,0}$
satisfies $d_1^c=I d_1 I^{-1}$, where $I$ acts as
$\1^{p-q}$ on $\Lambda^{p,q}_\hor(Q)\otimes \Lambda^m_\ver(Q)$.

\hfill

\proof
First, let us prove that $d_1$ has only two non-zero Hodge components.
A priori, $d_1$ could have several Hodge components,
$d_1= d^{-k, k+1}_1+ d^{-k+1, k}_1+ \cdots + d^{k, -k+1}_1+d^{-k,
	k+1}_1$. 
This is what happens with the Hodge components of de Rham
differential on an almost complex manifold.
All these components are clearly derivations.
However, 
\[ 
d(\Lambda^0(Q)) \subset \Lambda^{1}_\ver(Q)\oplus
\Lambda^{1,0}_\hor(Q) \oplus\Lambda^{0,1}_\hor(Q). 
\]
Therefore, only $d_1^{0,1}$ and $d_1^{1,0}$ are non-zero on
functions, the rest of the Hodge components are $C^\infty(Q)$-linear.
By \ref{_deriva_dete_Claim_}, it would suffice to show that
the other differentials vanish on a basis in $\Lambda^1(Q)$.

Since the space of leaves of $R$ is a complex manifold,
and $d=d_1$ on basic forms, we have $d_1=d_1^{0,1}+d_1^{1,0}$
on basic forms. Since $d(\xi ^{\,\flat})=\omega_0$, we find 
$d_1(\xi ^{\,\flat})=0$, which gives $d_1=d_1^{0,1}+d_1^{1,0}$
on $\xi ^{\,\flat}$. We proved the decomposition
$d_1=d_1^{0,1}+d_1^{1,0}$. The relation
$d_1^{0,1}-d_1^{1,0}=\1 I d_1 I^{-1}$ follows in the usual way,
because $\frac{d_1- \1 d_1^c}{2}$ has Hodge type $(1,0)$,
hence satisfies $d_1^{1,0}=\frac{d_1-\1 d_1^c}{2}$.
\endproof


\section{Transversally K\"ahler manifolds}


\definition\label{_transvers_Kahler_Definition_}
A manifold $M$ equipped with an integrable distribution
$F\subset TM$ is called a {\bf foliated manifold}. In the sequel, we shall
always assume that $F$ is orientable. 
Let $\omega_0\in \Lambda^2 (M)$ be a closed, basic 2-form 
on a foliated manifold $(M,F)$, vanishing on $F$ and
non-degenerate on $TM/F$. Let $g_0 \in \Sym^2(T^*M)$ 
be a basic bilinear symmetric form that is  positive definite on $TM/F$.
Since $g_0, \omega$ are basic, the operator $I:= \omega^{-1}_0 \circ g_0 :\; TM/F\arrow TM/F$
is well-defined on the leaf  space $L$ of $F$ (locally the leaf space always
exists by Frobenius theorem).  Assume that $I$ defines an integrable complex structure
on $L$ for any open set $U\subset M$ for which the leaf space is well-defined.
Then $(M, F, g, \omega)$ is called {\bf transversally K\"ahler}.
A vector field $v$ such that $\Lie_v(F)\subset F$, and
$\Lie_v I =0$ is called {\bf transversally holomorphic};
it is called {\bf transversally Killing}\index[terms]{vector field!Killing!transversally}\index[terms]{vector field!holomorphic!transversally} if, in addition,
$\Lie_v g_0 =0$.

\hfill

\definition
Let $F\subset TM$ be an integrable distribution and $\Lambda_\bas^*(M)$
the complex of basic forms. Its cohomology algebra is called the 
{\bf basic}, or {\bf transversal} cohomology of $M$. We denote the
basic cohomology by $H_\bas^*(M)$.\index[terms]{cohomology!basic}

\hfill

\remark
Note that $H_\bas^*(M)$ can be infinite-dimensional even
when $M$ is compact, \cite{_Schwarz:foliation_}.

\hfill

The main result of this section is the following theorem,
that is  a weaker form of the main theorem from
\cite{_Kacimi_}. Our result is less general,
but the proof is simple and self-contained.

\hfill

\theorem\label{_Transversal_Lefschetz_Theorem_}
Let $(M, F, \omega_0, g_0)$ be a compact, transversally K\"ahler manifold.
Assume, moreover, that:
\begin{itemize}
	\item[(*)] $M$ is equipped with a Riemannian metric 
	$g$ such that the restriction of $g$ to the orthogonal complement
	$F^\bot = TM/F$ coincides with $g_0$, and $F$ is generated by a
	collection of Killing vector fields\index[terms]{vector field!Killing} $v_1, ..., v_r$. 
	\item[(**)] There exists a closed differential form $\Phi$
	on $M$
	that vanishes on $F$ and gives a Riemannian volume
	form on $TM/F$.\footnote{The assumptions (*) and (**) hold
		for Sasakian manifolds (\ref{_Reeb_Sasakian_Theorem_}) and Vaisman manifolds\index[terms]{manifold!Vaisman} (\ref{_canon_foli_totally_geodesic_Remark_}).}
\end{itemize}
Then the basic cohomology\index[terms]{cohomology!basic} $H^*_\bas(M)$ of $M$ is finite-dimensional and admits
the Hodge decomposition and the Lefschetz $\gsl(2)$-action,
as in the K\"ahler case.

\hfill

\proof Consider the differential graded algebra of basic forms
$\Lambda^*_\bas(M)$.
This algebra is equipped with an action of the superalgebra of
K\"ahler supersymmetry $\goth a$ as in \ref{_kah_susy_Theorem_}. 
Indeed, define the Lefschetz $\gsl(2)$-action by taking
the Lefschetz triple $L_{\omega_0}, \Lambda_{\omega_0}:=*L_{\omega_0}*$
and $H_{\omega_0}:= [L_{\omega_0}, \Lambda_{\omega_0}]$, and
the transversal Weil operator $W$ acting in the standard
way on $\Lambda^*_\bas(M)$ and extended to $\Lambda^*(M)$
by acting trivially on $\Lambda^1_\bas(M)^\bot$
(or in any other way, it does not matter). 
Together with the de Rham differential
$d:\; \Lambda^*_\bas(M)\arrow \Lambda^*_\bas(M)$
these operators generate the Lie superalgebra
$\goth a\subset \End(\Lambda^*_\bas(M))$, which
is isomorphic to the superalgebra of
K\"ahler supersymmetry $\goth a$ (\ref{_kah_susy_Theorem_}), 
because $\Lambda^*_\bas(M)$ is locally identified
with the algebra of differential forms on the leaf space of $F$
that is  K\"ahler. 

\ref{_Transversal_Lefschetz_Theorem_} would follow if
we identify $\ker (\Delta_\bas\restrict {\Lambda^*_\bas(M)})$ with the space 
$H^*_\bas(M)$, where $\Delta_\bas\in \goth a$ is the transversal
Laplace operator, $\Delta_\bas=\{d, d^*_\bas\}$, where
$d^*_\bas$ denotes the $d^*$-operator on the leaf space.

We reduced  \ref{_Transversal_Lefschetz_Theorem_} to the
following result.

\hfill

\proposition\label{_basic_harmo_Proposition_}
Let $(M, F, g)$, $\rk F=r$, be a Riemannian foliated manifold,
$d=d_0+d_1+ \cdots+ d_{r+1}$ the Hattori decomposition of the differential,
and $\Delta_\bas:= \{d, d_\bas^*\}$ the basic Laplacian defined
on basic forms. Assume that:\index[terms]{basic Laplacian}
\begin{itemize}
	\item[(*)] $F$ is generated by a
	collection of Killing\index[terms]{vector field!Killing} 
	vector fields $v_1, v_2, \ldots$
	\item[(**)] The Riemannian volume form $\Phi\in\Lambda^r_\ver(M)$
	satisfies $d_1(\Phi)=0$. 
\end{itemize}
Then there exists a natural
isomorphism between the basic harmonic forms
and $H^*_\bas(M)$.

\hfill

We start from the following lemma.

\hfill

\lemma
In the assumptions of \ref{_basic_harmo_Proposition_}, let
$\alpha, \beta \in \Lambda^*_\bas(M)$ be two basic forms.
Then $g(d\alpha, \beta)= g(\alpha, d^*_\bas\beta)$,
where $d^*_\bas$ denotes the $d^*$-operator on the leaf space of $F$.

\hfill

\proof
This is where we use the assumption (**) of \ref{_basic_harmo_Proposition_}.
Let $d^*_h$ be the composition of $d^*$ with the orthogonal
projection to the horizontal part $\Lambda^*_\hor(M)$
(see Subsection \ref{_Leray_Serre_Subsection_}).
We only need  to show that 
\begin{equation}\label{_d^*_bas_Equation_}
g(d^*_h\alpha, \beta) = g(d^*_\bas\alpha, \beta)\ \ \text{for all}\ \ \alpha, \beta\in \Lambda^*_\bas(M). 
\end{equation}
By \ref{_kah_susy_Theorem_}, one has $d^*_\bas= - *_\bas d *_\bas$
where $*_\bas$ is the Hodge star operator on the leaf space.
Let $r=\rk F$ and $\Phi\in \Lambda^r F$ be the Riemannian volume form.
Using the assumption (**), we obtain that 
\begin{equation}\label{_d_1_Phi_via_grading_Equation_}
d\Phi\in \bigoplus_{i=0}^{r-1} \Lambda^*_\hor(M) \otimes \Lambda^i_\ver(M).
\end{equation}
Clearly, $*\alpha= \Phi \wedge *_\bas(\alpha)$, and
\[
*d*\alpha = *d( \Phi \wedge *_\bas(\alpha)) = d^*_\bas\alpha + *(d\Phi\wedge \alpha).
\]
By \eqref{_d_1_Phi_via_grading_Equation_},
the last term belongs to $\bigoplus_{i=1}^{r} \Lambda^*_\hor(M) \otimes \Lambda^i_\ver(M)$,
hence it is orthogonal to $\Lambda^*_\hor(M)$. This proves \eqref{_d^*_bas_Equation_}. \endproof

\hfill

Now we can prove  \ref{_basic_harmo_Proposition_}.
By \eqref{_d^*_bas_Equation_}, for any basic form $\alpha$ we have 
\begin{equation}\label{_delta_bas_dual_Equation_}
g(\Delta_\bas\alpha, \alpha)= (d\alpha, d\alpha)+ (d^*_\bas \alpha, d^*_\bas \alpha),
\end{equation}
hence a basic form belongs to $\ker \Delta_\bas$ if and only if it is 
closed and orthogonal to all exact basic forms. This gives
an embedding 
\begin{equation}\label{_harmo_to_bas_Equation_}
\ker \Delta_\bas\restrict {\Lambda^*_\bas(M)}\hookrightarrow H^*_\bas(M).
\end{equation}
It remains only to show that the map \eqref{_harmo_to_bas_Equation_} is surjective.
This is where we use the assumption (*) of \ref{_basic_harmo_Proposition_}.

Consider the \index[persons]{Hattori, A.} Hattori decomposition $d=d_0+d_1+\cdots+d_{r+1}$ (see 
Subsection \ref{_Leray_Serre_Subsection_}).
Let $v_1,\ldots, v_r\in F$ be the Killing,\index[terms]{vector field!Killing!transversally} transversally Killing vector fields,
postulated in (*), and $\Delta_s$ the ``split Laplacian'',
$\Delta_s:= \{d_1, d_1^*\}- \sum_i \Lie_{v_i}^2$. 
Clearly, $\Delta_s$ is\index[terms]{split Laplacian}
an elliptic, second order differential operator. 
By definition, $g(e_{v_i} i_{v_i}\alpha, \alpha)= g(i_{v_i}\alpha,i_{v_i}\alpha)$.
Since $v_i$ are Killing, one has $\Lie_{v_i}=-\Lie_{v_i}^*$. Therefore, 
$\Delta_s$ is self-adjoint and positive, with
\[
g(\Delta_s\alpha, \alpha) = g(d_1\alpha, d_1\alpha) +  
g(d^*_1\alpha, d^*_1\alpha) + \sum_i g(\Lie_{v_i}\alpha,\Lie_{v_i}\alpha). 
\]
We obtain that each $\alpha \in \ker \Delta_s$ satisfies
$i_{v_i} \alpha = \Lie_{v_i}\alpha=0$.

Consider the orthogonal projection map $\Pi_\hor:\; \Lambda^*(M)\arrow \Lambda^*_\hor(M)$.
Since the vector fields $v_i$ are Killing\index[terms]{vector field!Killing} and preserve $F$, 
the projection $\Pi_\hor$ commutes with $\Lie_{v_i}$. It is not hard to see
that $\Pi_\hor$ commutes with $d_1$ and $d_1^*$.
Indeed, $d_1$ is the part of $d$ that maps
$\Lambda^p_\hor(M)\otimes \Lambda^q_\ver(M)$
to $\Lambda^{p+1}_\hor(M)\otimes \Lambda^q_\ver(M)$,
and $d_1^*$ its adjoint. We obtain that 
$\Pi_\hor$ commutes with $\Delta_s$.

Since $\Delta_s$ is a positive, self-adjoint Fredholm operator,
its eigenvectors are dense in $\Lambda^*(M)$.
Since $[ \Pi_\hor, \Delta_s]=0$, the eigenvectors of 
$\Delta_s$ are dense in $\Lambda^*_\hor(M)$.

Let $G$ be the closure of the Lie group generated by
the action of $e^{\R v_i}$. Since each $e^{\R v_i}$ acts by isometries,
the group $G$ is compact. By construction, $G$ acts on $M$
preserving the foliation $F$, the metric and the transversal
K\"ahler structure, and hence  it commutes with the Laplacian.
Averaging on $G$, we obtain that
the eigenvectors of $\Delta_s$ are dense in the space
$\Lambda^*_\hor(M)^G= \Lambda^*_\bas(M)$ of all basic forms.

On basic forms, $d=d_1$, and hence  on $\Lambda^*_\bas(M)$ one has 
$[d, \Delta_s]=0$. Let $\Lambda^*_\bas(M)_\lambda$ be the eigenspace
of $\Delta_s\restrict {\Lambda^*_\bas}(M)$ corresponding to the eigenvalue $\lambda$.
For any closed $\alpha\in \Lambda^*_\bas(M)$
we have $\lambda\alpha = (dd^*+ d^*d)(\alpha)= dd^*\alpha$.
Therefore, any closed form in $\Lambda^*_\bas(M)_\lambda$ is
exact when $\lambda\neq 0$. We have shown that the only eigenspace of $\Delta_\bas$
which contributes to the basic cohomology\index[terms]{cohomology!basic} is 
$\Lambda^*_\bas(M)_0= \ker \Delta_\bas\restrict {\Lambda^*_\bas(M)}$.
This proves \eqref{_harmo_to_bas_Equation_} and ends the
proof of \ref{_basic_harmo_Proposition_} and \ref{_Transversal_Lefschetz_Theorem_}.
\endproof


\section[Basic cohomology and Hodge theory on Sa\-sa\-ki\-an manifolds]{Basic cohomology and Hodge theory on\\ Sa\-sa\-ki\-an manifolds}


\subsection[The cone of a morphism of complexes and cohomology of Sasakian manifolds]{The cone of a morphism of complexes and\\ cohomology of Sasakian manifolds}
\index[terms]{cone!of a morphism}
We recall first several notions in homological algebra, see e.g. 
\cite{_Gelfand_Manin_}. We have already given an introduction to the
derived functors in Section \ref{_Derived_functors_Subsection_}; 
it contains some of the same notions, with a
different spin.

\hfill

\definition
{\bf A complex} $(C_*, d)$ is a collection of vector spaces and homomorphisms
\[ 
\dots\stackrel d \arrow C_i  \stackrel d \arrow C_{i+1} \stackrel d \arrow \cdots
\]
(more generally, a collection of objects in an abelian category)
such that $d^2=0$. {\bf A morphism} of complexes is a collection of maps
$C_i \arrow C_i'$ from the vector spaces of 
a complex $(C_*, d)$ to the vector spaces of $(C'_*, d)$, commuting with the differential.
{\bf The cohomology groups} of a complex $(C_*, d)$  are the groups 
\[ H^i(C_*, d):= \frac{\ker d\restrict C_i}{\im d\restrict{C_{i-1}}}.\]
Clearly, any morphism induces a homomorphism in cohomology.
{\bf An exact sequence of complexes} is a sequence 
\[ 0\arrow A_* \arrow B_* \arrow C_*\arrow 0\] of morphisms 
of complexes such that the corresponding sequences
\[ 0\arrow A_i \arrow B_i \arrow C_i\arrow 0\]
are exact for all $i$.

\hfill

The following claim is very basic.

\hfill

\claim
Let $0\arrow A_* \arrow B_* \arrow C_*\arrow 0$
be an exact sequence of com\-plex\-es. Then
there is a natural long exact sequence of cohomology
\[
\cdots \rightarrow H^{i-1}(C_*, d_C) \rightarrow H^{i}(A_*, d_A)\rightarrow H^{i}(B_*, d_B)
\rightarrow H^{i}(C_*, d_C) \rightarrow \cdots \quad\text{\endproof}
\]

\definition\label{_cone_complexes_Definition_}
Let  $(C_*, d_C)\stackrel \phi \arrow (C'_*, d_{C'})$
be a morphism of complexes. Consider the complex
$C(\phi)_*$, with $C(\phi)_i=C_{i+1}\oplus C'_{i}$
and differential 
\[ d:=d_C+\phi-d_{C'}:\; C_i\oplus C'_{i-1}\arrow 
C_{i+1}\oplus C'_{i},
\]
or, explicitly:
\[
d(c_i,c'_{i-1})=\left(d_C(c_i), \phi(c_i)-d_{C'}(c_{i-1}')\right).
\]
The complex $(C(\phi)_*, d)$ is called\index[terms]{cone!of a morphism}
{\bf the cone of $\phi$}.\footnote{See also 
Subsection \ref{_Derived_functors_Subsection_}.}\index[terms]{derived functor}

\hfill

Denote by $C_*[1]$ the complex $C_*$ shifted by one,
with $C_i[1]=C_{i+1}$. 
The exact sequence of complexes
\[ 
0 \arrow C_*'\arrow C(\phi)_*\arrow C_*[1]\arrow 0
\]
gives the long exact sequence
\begin{multline}\label{_long_exact_from_cone_Equation_}
\cdots \arrow H^{i}(C)\stackrel \phi \arrow H^i(C') 
\arrow  H^i(C(\phi))\arrow\\ \arrow H^{i+1}(C) \stackrel \phi\arrow   H^{i+1}(C') \arrow \cdots
\end{multline}

\proposition\label{_Sasakian_vecr_inva_cone_Proposition_}
Let $Q$ be a Sasakian manifold, and $\xi $ the Reeb field.
Denote by $\Lambda^*_\xi (Q)$ the differential
graded algebra of $\Lie_\xi $-invariant forms,
and let $\Lambda^*_\bas(Q)\subset\Lambda^*_\xi (Q)$ be the algebra of 
basic forms. Denote by $\omega_0\in \Lambda^*_\bas(Q)$ the transversal 
K\"ahler form\index[terms]{form!K\"ahler!transversal} (\ref{_omega_0_trans_Kahler_Claim_}), and let 
$L_{\omega_0}:\; \Lambda^*_\bas(Q)\arrow \Lambda^{*+2}_\bas(Q)$
be the multiplication map. Then the complex $\Lambda^*_\xi (Q)$
is naturally identified with $C(L_{\omega_0})[-1]$,
where $C(L_{\omega_0})$ is the cone of the \index[terms]{cone!of a morphism}
morphism $L_{\omega_0}:\; \Lambda^*_\bas(Q)\arrow \Lambda^{*+2}_\bas(Q)$.

\hfill

\proof
Consider the Hattori decomposition of the differential in $\Lambda^*(Q)$, 
$d= d_0+ d_1 + d_2$ (Subsection \ref{_LS_decompo_Sasakian_Subsection_}).
By \ref{_d_0_expli_Claim_}, $d_0$ vanishes on $\Lambda^*_\xi (Q)$;
\ref{_d_2_expli_Proposition_} implies that $d_2=L_{\omega_0}i_\xi $.
Clearly, $\Lambda^*_\xi (Q)= \Lambda_\bas^*(Q)\oplus \xi ^{\,\flat} \wedge \Lambda_\bas^*(Q)$.
The operator $d_2\Lambda_{\omega_0}i_\xi $ acts trivially on 
$\Lambda_\bas^*(Q)$ and maps $\xi ^{\,\flat} \wedge \alpha$ to
$L_{\omega_0}(\alpha)$ for any $\alpha \in \Lambda_\bas^*(Q)$.
We obtain the following natural decomposition
\[
\Lambda^*_\xi (Q)= 
\Lambda_\bas^*(Q)\oplus \xi ^{\,\flat} \wedge \Lambda_\bas^*(Q) = \Lambda_\bas^*(Q)\oplus \Lambda_\bas^*(Q)[-1].
\]
\ref{_d_0_expli_Claim_} 
implies that $d_0$ vanishes on $\Lambda^*_\xi (Q)$, which implies $d=d_1+d_2$.
By \ref{_d_2_expli_Proposition_}, the differential $d_1+ d_2$ in 
$\Lambda^*_\xi (Q)$ gives the same differential as in 
$C(L_{\omega_0})[-1]=\Lambda_\bas^*(Q)\oplus \Lambda_\bas^*(Q)[-1]$:
\[
d(\alpha\oplus \xi ^{\,\flat} \wedge\beta)= d\alpha+ L_{\omega_0}(\beta)\oplus (- \xi ^{\,\flat} \wedge d \beta)
\]
for any $\alpha, \beta \in \Lambda_\bas^*(Q)$.
\endproof

\subsection{Harmonic form decomposition on Sasakian manifolds}

We give a new proof of the main result on harmonic forms
on compact Sasakian manifolds (see \cite[Proposition 7.4.13]{bog}, 
essentially based on a result by Tachi\-bana, \cite{_Tachibana_}).

\hfill

\theorem\label{_Sasakian_cohomo_cone_Theorem_}
Let $Q$ be a $2n+1$-dimensional compact Sasakian manifold,  $R$ the Reeb foliation, 
and $H^*_\bas(Q)$ the corresponding basic cohomology.
Consider the Lefschetz $\gsl(2)$-triple
$L_{\omega_0}, \Lambda_{\omega_0}, H_{\omega_0}$ acting
on the basic cohomology\index[terms]{cohomology!basic} (\ref{_Transversal_Lefschetz_Theorem_}).
Then:
\begin{equation*}
H^i(Q)= \begin{cases}
\ker L_{\omega_0} \restrict {H^i_\bas(Q)}, \qquad \text{for}\ i\geq n,\\[.2in]
\dfrac{H^i_\bas(Q)}{\im L_{\omega_0}}, \qquad\qquad \text{for}\ i< n.
\end{cases}
\end{equation*}


\hfill

\pstep
The cohomology of the algebra $\Lambda^*_\xi (Q)$ of $\Lie_\xi$-invariant
forms is equal to the cohomology of $\Lambda^*(Q)$. Indeed,
let $G$ be the closure of the action of $e^{t\xi }$.
Since $\xi $ is Killing,\index[terms]{vector field!Killing} $e^{t\xi }$ acts on $Q$ by isometries,
hence its closure $G$ is compact. Since $G$ is connected,
its action on cohomology is trivial. The averaging map
$\Av_G:\; \Lambda^*(Q)\arrow \Lambda^*_\xi (Q)$
induces an isomorphism on cohomology.

\hfill

{\bf Step 2:} Applying \eqref{_long_exact_from_cone_Equation_},
\ref{_Sasakian_vecr_inva_cone_Proposition_}, and taking into account the isomorphism 
$H^*(\Lambda^*_\xi (Q))= H^*(Q)$, we obtain the long exact sequence
\begin{multline}\label{_cone_Sasa_equation_}
\cdots \arrow H^{i-2}_\bas(Q)\xrightarrow{L_{\omega_0}}  H^{i}_\bas(Q) 
\arrow  H^{i}(Q)\arrow\\
\arrow H^{i-1}_\bas(Q) \xrightarrow{L_{\omega_0}}   H^{i+1}_\bas(Q) \arrow \cdots
\end{multline}
Since the Lefschetz triple $L_{\omega_0}, \Lambda_{\omega_0}, H_{\omega_0}$
induces an $\goth{sl}(2)$-action on $H^{i}_\bas(Q)$ (\ref{_Transversal_Lefschetz_Theorem_}),
the map $H^{i-2}_\bas(Q)\stackrel {L_{\omega_0}} \arrow H^{i}_\bas(Q)$
is injective for $i\geq n$ and surjective for $i<n$.
Therefore, the long exact sequence \eqref{_cone_Sasa_equation_}
gives the short exact sequences
\[
0\arrow H^{i-2}_\bas(Q)\stackrel {L_{\omega_0}} \arrow H^{i}_\bas(Q) 
\arrow  H^{i}(Q)\arrow 0
\]
for $i\geq n$ and 
\[
0 \arrow H^{i}(Q)\arrow H^{i-1}_\bas(Q) \stackrel {L_{\omega_0}} 
\arrow   H^{i+1}_\bas(Q)\arrow 0
\]
for $i< n$.
\endproof

\hfill

\theorem\label{_harmo_Sasa_decompo_Theorem_}
Let $Q$ be a $2n+1$-dimensional compact Sasakian manifold, and
${\cal H}^i_\bas(Q)$ the space of basic harmonic forms.
Let $\xi ^{\,\flat}$ be the contact form, dual to the Reeb field.
Denote by ${\cal H}^i$ the space of all $i$-forms $\alpha$ on $Q$ 
that satisfy
\begin{description}
	\item[for $i \leq n$:]
	$\alpha$ is basic harmonic (that is, belongs 
	to the kernel of the basic Laplacian) and satisfies
	$\Lambda_{\omega_0}(\alpha)=0$; 
	\item[for $i > n$:]
	$\alpha = \beta \wedge \xi ^{\,\flat}$
	where $\beta$ is basic harmonic and satisfies
	$L_{\omega_0}(\beta)=0$. 
\end{description}
Then all elements of ${\cal H}^*$ are harmonic, and, moreover,
all harmonic forms on $Q$ belong to ${\cal H}^*$.

\hfill

\pstep
We prove that all $\gamma\in {\cal H}^*$ are harmonic.
Let $*_\bas$ be the Hodge star operator on basic forms.
Then for any $\gamma \in \Lambda^*_\hor(Q)$ one has
$*(\gamma) = *_\bas(\gamma) \wedge \xi ^{\,\flat}$.
This implies that the  two classes of forms in
\ref{_harmo_Sasa_decompo_Theorem_} are exchanged by
$*$, and it suffices to prove that
all $\alpha\in \Lambda^i(Q)$ that are  basic harmonic
and satisfy  $\Lambda_{\omega_0}(\alpha)=0$ 
for $i\leq n$ are harmonic. Such a form $\alpha$ is  closed by \ref{_Transversal_Lefschetz_Theorem_}.
Since $\alpha$ is basic hence $\Lie_\xi $-invariant, 
it satisfies $d_0^*(\alpha)=0$ (\ref{_d_0_expli_Claim_}).
By \eqref{_delta_bas_dual_Equation_}, a
basic form $\alpha$ is basic harmonic if and only if
$d_1(\alpha)=d_1^*(\alpha)=0$. Finally, $d_2^*=\Lambda_{\omega_0}e_{\xi }$
(\ref{_d_2_expli_Proposition_}), and hence  $d_2^*(\alpha)=0$. This implies
$d^*\alpha= d_0^*+d_1^*+d_2^*(\alpha)=0$.

\hfill

{\bf Step 2:} 
We prove that all harmonic forms on $Q$ belong to ${\cal H}^*$.
By \ref{_Sasakian_cohomo_cone_Theorem_}, 
the dimension of $H^i(Q)$ is equal to the dimension of 
${\cal H}^i$, and hence  the embedding
${\cal H}^i \arrow H^i(Q)$ constructed in Step 1
is also surjective for all indices $i$.
\endproof


\section{Hodge theory on Vaisman manifolds}\index[terms]{manifold!Vaisman}


\subsection{Basic cohomology of Vaisman manifolds}


Let $M$ be a Vaisman manifold, $\theta^\sharp$ its Lie field,
and $\xi := I(\theta^\sharp)$. Using the local decomposition
of $M$ as a product of $\R$ and a Sasakian manifold (\ref{_Lee_field_holo_Proposition_}),
$\xi $ can be identified with the Reeb field on its Sasakian component.
The manifold $M$ is equipped with 2 remarkable foliations,
$\Sigma=\langle \xi , \theta^\sharp\rangle$ and 
${\cal L}=\langle \theta^\sharp\rangle$. Both of these
foliations satisfy the assumptions of \ref{_basic_harmo_Proposition_} (see also \ref{_canon_foli_totally_geodesic_Remark_} (iii)).

Consider the corresponding algebras of basic forms:
\begin{itemize}
	\item[(i)] $\Lambda^*_\sas(M)$ -- forms that are  basic with respect to ${\cal L}$,
	\item[(ii)] $\Lambda^*_\kah(M)$ -- forms that are  basic with respect to $\Sigma$.
\end{itemize}
Denote by $H^*_\sas(M)$, $H^*_\kah(M)$ the corresponding basic cohomology\index[terms]{cohomology!basic} 
algebras. 

\hfill

By \ref{_Transversal_Lefschetz_Theorem_}, the algebra $H^*_\kah(M)$
is equipped with the Lefschetz $\gsl(2)$-action and the Hodge decomposition
in the same way as for K\"ahler manifolds.

\hfill

The cohomology of Vaisman manifolds\index[terms]{manifold!Vaisman} can be expressed non-am\-bi\-gu\-ously in terms
of $H^*_\kah(M)$ and the Lefschetz $\gsl(2)$-action. 
The following theorem is the Vaisman analogue
of \ref{_Sasakian_cohomo_cone_Theorem_}.

\hfill

\theorem\label{_Vaisman_coho_Theorem_}
Let $M$ be a compact Vaisman manifold, $\dim_\R M=2n$, and $\theta$ its
Lee form.\index[terms]{form!Lee} Then 
\begin{align}\label{_coho_Sasakian_folia_Equation_}
H^i(M)& \cong H^i_\sas(M) \oplus \theta \wedge H^{i-1}_\sas(M), \\[.1in]
H^i_\sas(M)&=
\begin{cases} \label{_coho_Vaisman_folia_Equation_}
\displaystyle\ker L_{\omega_0} \restrict {H^i_\kah(M)}, & \text{for} \ \  i> n-1,\\[.2in]
\displaystyle\frac{H^i_\kah(M)}{\im L_{\omega_0}}, & \text{for}\ \   i\leq n-1
\end{cases} 
\end{align}
where $\{L_{\omega_0}, \Lambda_{\omega_0}, H_{\omega_0}\}$
is the Lefschetz $\gsl(2)$-triple acting on $H^i_\kah(M)$.

\hfill

\pstep We start by proving  \eqref{_coho_Sasakian_folia_Equation_}.
Let $G_{\theta^\sharp}$ be the closure of the
1-parametric group $e^{t\theta^\sharp}$ in the group
$\Iso(M)$ of isometries of $M$. Since
$\Iso(M)$ is compact, the group $G_{\theta^\sharp}$
is also compact. Clearly, averaging on 
$G_{\theta^\sharp}$ does not change the cohomology class
of a form. Therefore, the algebra of
$G_{\theta^\sharp}$-invariant forms has the same
cohomology as $\Lambda^*(M)$. 
Consider the decomposition
\eqref{_horis_vert_decompo_Equation_}
associated with the rank 1 foliation ${\cal L}$:
\begin{equation}\label{_horis_vert_2_decompo_Equation_}
\Lambda^m(M)= \bigoplus_p \Lambda^p_\hor(M)\otimes
\Lambda^{m-p}_\ver(M).
\end{equation}
The bundle $\Lambda^{1}_\ver(M)= {\cal L}^*$
is 1-dimensional and generated by $\theta$. Therefore,
\eqref{_horis_vert_2_decompo_Equation_}
gives
\[
\Lambda^m(M)= \Lambda^m_\hor(M)\oplus \theta \wedge \Lambda^{m-1}_\hor(M).
\]
Denote by $\Lambda^m(M)^{G_{\theta^\sharp}}$ the $G_{\theta^\sharp}$-invariant part of $\Lambda^m(M)$. 
Since $\theta$ is $G_{\theta^\sharp}$-invariant, and the
$G_{\theta^\sharp}$-invariant part of $\Lambda^m_\hor(M)$
is identified with $\Lambda^m_\sas(M)$, this gives
\[
\Lambda^m(M)^{G_{\theta^\sharp}}= 
\Lambda^m_\sas(M)\oplus \theta \wedge \Lambda^{m-1}_\sas(M).
\]
Taking cohomology, we obtain \eqref{_coho_Sasakian_folia_Equation_}.

\hfill

{\bf Step 2:} Now we shall prove
$H^i_\sas(M)= \ker L_{\omega_0} \restrict {H^i_\kah(M)}$
for $i\geq n-1$ and $H^i_\sas(M)= \dfrac{H^i_\kah(M)}{\im L_{\omega_0}}$
for $i< n-1$.

Let $\xi \in TM$ be the Reeb field defined as above, and
$G_\xi $ the closure of the 1-parametric group $e^{t \xi }$.
Then (as in Step 1 and in the proof of \ref{_Sasakian_cohomo_cone_Theorem_}),
the $G_\xi $-invariant part of $\Lambda^*_\sas(M)$ is written as
\begin{equation}\label{_Sasakian_inva_decompo_Equation_}
\Lambda^i_\sas(M)^{G_\xi }=\Lambda^i_\kah(M) \oplus \xi ^{\,\flat} \wedge \Lambda^{i-1}_\kah(M)
\end{equation}
with the differential acting (see \ref{_Sasakian_vecr_inva_cone_Proposition_})
as
\begin{equation}\label{_Sasakian_inva_decompo_differe_Equation_}
d(\alpha\oplus \xi ^{\,\flat} \wedge\beta)= d\alpha+ L_{\omega_0}(\beta)\oplus (- \xi ^{\,\flat} \wedge d \beta)
\end{equation}
for any $\alpha, \beta\in \Lambda^*_\kah(M)$.

Using the Cartan formula, we notice that $G_\xi $ acts trivially
on the cohomology of the complex $(\Lambda^*_\sas(M), d)$, and hence  the cohomology
of $(\Lambda^*_\sas(M)^{G_\xi }, d)$ is identified with $H^*_\sas(M)$.
From \eqref{_Sasakian_inva_decompo_Equation_} and
\eqref{_Sasakian_inva_decompo_differe_Equation_}, we obtain that
$\Lambda^i_\sas(M)^{G_\xi }$ is the cone of the morphism of complexes
$L_{\omega_0}:\; \Lambda^*_\kah(M)[-1] \arrow \Lambda^*_\kah(M)[1]$.
From the long exact sequence
\eqref{_Sasakian_inva_decompo_differe_Equation_}
we obtain a long exact sequence identical to \eqref{_cone_Sasa_equation_}:
\begin{multline}
\cdots \arrow H^{i-2}_\kah(M)\xrightarrow{L_{\omega_0}}  H^{i}_\kah(M) 
\arrow  H^{i}_\sas(M)\arrow \\
\arrow H^{i-1}_\kah(M) \xrightarrow{L_{\omega_0}}   H^{i+1}_\kah(M) 
\arrow\cdots
\end{multline}
Since $L_{\omega_0}:\; H^{i-2}_\kah(M)\arrow  H^{i}_\kah(M)$ is injective for $i\geq n-1$
and surjective for $i <n-1$, this long exact sequence breaks into short
exact sequences of the form
\[
0\arrow H^{i-2}_\kah(M)\stackrel {L_{\omega_0}} \arrow H^{i}_\kah(M) 
\arrow  H^{i}_\sas(M)\arrow 0
\]
for $i\geq n-1$ and 
\[
0 \arrow H^{i}(M)_\sas \arrow H^{i-1}_\kah(M) \stackrel {L_{\omega_0}} \arrow   H^{i+1}_\kah(M)\
\arrow 0
\]
for $i< n$. We finished the proof of 
\ref{_Vaisman_coho_Theorem_}. \endproof

\subsection{Harmonic forms on Vaisman manifolds}\index[terms]{manifold!Vaisman}

It turns out that (just like it happens in the Sasakian case)
the cohomology decomposition obtained in \ref{_Vaisman_coho_Theorem_}
gives a harmonic form decomposition. Together with the Hodge decomposition of the 
basic cohomology\index[terms]{cohomology!basic} of transversally K\"ahler structure\index[terms]{structure!K\"ahler!transversally} this allows us to represent
certain cohomology classes by forms of a given Hodge type. 
This theorem was obtained in \cite{KS}, see also \cite{va_gd} and \cite{tsu}.

\hfill

Recall that on a compact Vaisman manifold\index[terms]{manifold!Vaisman} $(M,I,g,\theta)$ with fundamental form $\omega$, the canonical foliation\index[terms]{foliation!canonical} $\Sigma$ is transversally K\"ahler (\ref{_Sigma_transv_K_}). 
We used $\xi$ to denote the vector field $g$-equivalent to $I\theta$, $\Lambda^*_\kah(M)$ the space of 
basic forms with respect to $\Sigma$, and 
$\Delta_\kah:\; \Lambda^*_\kah(M)\arrow \Lambda^*_\kah(M)$
the basic Laplacian (\ref{_basic_harmo_Proposition_}). From \ref{_Transversal_Lefschetz_Theorem_}
we obtain that the space ${\cal H}^i_\kah(M)$ 
of basic harmonic forms is equipped with the
Lefschetz $\gsl(2)$-action by the operators 
$L_{\omega_0}, \Lambda_{\omega_0}, H_{\omega_0}$. The main result of this section is:

\hfill

\theorem\label{_Vaisman_harmonic_forms_Theorem_}
Let $(M,I,g,\theta)$ be a compact Vaisman manifold\index[terms]{manifold!Vaisman} of complex dimension $n$, with fundamental form $\omega$, and canonical foliation $\Sigma$. \index[terms]{foliation!canonical}
Denote by ${\cal H}^i$ the space of all  basic $i$-forms 
$\alpha\in \Lambda^*_\kah(M)$
that satisfy:
\begin{description}
	\item[\hspace{.1in}for $i \leq n$:]
	$\alpha$ is basic harmonic (i.\,e.    $\Delta_\kah(\alpha)=0$) and satisfies
	$\Lambda_{\omega_0}(\alpha)=0$; 
	\item[\hspace{.1in}for $i > n$:]
	$\alpha = \beta \wedge I\theta$
	where $\beta$ is basic harmonic and satisfies
	$L_{\omega_0}(\beta)=0$. 
\end{description}
Then all elements of ${\cal H}^*\oplus \theta\wedge {\cal H}^*$ are harmonic and, moreover,
all harmonic forms on $M$ belong to ${\cal H}^*\oplus \theta\wedge {\cal H}^*$.

\hfill

\pstep
This statement is similar to \ref{_harmo_Sasa_decompo_Theorem_},
and the proof is essentially the same. We start by
proving that all $\gamma\in {\cal H}^*\oplus\theta\wedge {\cal H}^*$ are harmonic.
Since $\theta$ is parallel, and a product of a parallel form 
and a harmonic form is harmonic (\cite[Corollary 2.9]{_Verbitsky:G2_forma_};
see also the proof of \ref{_invariance_for_cohomology_Vaisman_Lemma_}), it suffices only to show
that all elements in ${\cal H}^*$ are harmonic:
\begin{equation}\label{_H^*_harmo_Equation_}
{\cal H}^*\subset \ker \Delta.
\end{equation}
Consider the foliation ${\cal L}=\langle \theta^\sharp\rangle$.
Since $M$ is locally a product of a Sasakian manifold
and a line (\ref{_Lee_field_holo_Proposition_}), 
all ${\cal L}$-basic, transversally harmonic forms are harmonic.
Therefore, \eqref{_H^*_harmo_Equation_} would follow if we prove
\begin{equation}\label{_H^*_transve_sasa_harmo_Equation_}
{\cal H}^*\subset \ker \Delta_\sas.
\end{equation}
where $\Delta_\sas$ is the basic Laplacian of 
the foliation ${\cal L}$.

It would suffice to prove \eqref{_H^*_transve_sasa_harmo_Equation_}
locally in $M$. However, locally  ${\cal L}$
admits a leaf space $Q$,
that is  Sasakian, and we can regard elements
from ${\cal H}^*$ as forms on $Q$ and $\Delta_\sas$ as the
usual Laplacian on $\Lambda^*(Q)$. 
Then \eqref{_H^*_transve_sasa_harmo_Equation_} follows from
\ref{_harmo_Sasa_decompo_Theorem_}, 
step 1.\footnote{The statement of \ref{_harmo_Sasa_decompo_Theorem_}
	is global, however, the proof of \ref{_harmo_Sasa_decompo_Theorem_}, 
	step 1 is local, and this statement is essentially identical to 
	\eqref{_H^*_transve_sasa_harmo_Equation_}.}

\hfill

{\bf Step 2:} We have shown that all elements of 
${\cal H}^*\oplus \theta\wedge {\cal H}^*$ are harmonic;
this gives a natural linear map
\[
\Psi:\; {\cal H}^*\oplus \theta\wedge {\cal H}^*\arrow H^*(M).
\]
To show that all harmonic form are obtained in this way,
it would suffice to prove that $\Psi$ is surjective.
This follows from \ref{_Vaisman_coho_Theorem_}
by dimension count.

 The cohomology of $M$ is isomorphic to
$H^i_\sas(M) \oplus \theta \wedge H^{i-1}_\sas(M)$ (by \ref{_Vaisman_coho_Theorem_}).
All forms in ${\cal H}^*$  are closed and belong to $\Lambda^*_\sas(M)$,
which gives a  map $\Psi_\sas:\; {\cal H}^*\arrow H^i_\sas(M)$. 
To prove that $\Psi$ is surjective, it remains to show that
$\Psi_\sas$ is surjective.

However, dimension of $H^i_\sas(M)$ is equal to
$\dim \ker L_{\omega_0}\restrict{H^i_\kah(M)}$ for $i \geq n-1$ and to 
$\dim \dfrac{H^i_\kah(M)}{\im L_{\omega_0}}$ for $i < n-1$.
The space ${\cal H}^i$ has the same dimension by 
the transversal Hodge decomposition
(\ref{_Transversal_Lefschetz_Theorem_}).
This finishes the proof of 
\ref{_Vaisman_harmonic_forms_Theorem_}.
\endproof

\hfill

\corollary {(\cite{KS}, \cite{va_gd})}\label{dec_har}
Let $M$ be a  compact Vaisman manifold\index[terms]{manifold!Vaisman} of complex
dimension $n$. Then a (real)  $p$-form $\eta$, $0\leq
p\leq n-1$, is harmonic ($\Delta\eta=0$) if and only if it
has a decomposition   \index[terms]{differential form!transversally harmonic}
\begin{equation}\label{desc}
\eta=\al+\theta\wedge\beta,
\end{equation}
with $\al$ and $\be$ $\Sigma$-basic, primitive\footnote{A form $\eta\in \Lambda^*_\kah(M)$ is
called {\bf primitive} if $\Lambda_{\omega_0} \eta=0$.} and transversally\index[terms]{form!primitive}
harmonic forms.
Moreover, this decomposition is unique. \endproof

\hfill

Let  $e_i:=\dim {\cal H}^i_\kah(M)$  the {\bf $\Sigma$-basic Betti
  numbers} (w. r. t.  the canonical foliation).\index[terms]{foliation!canonical} When
$\Sigma$ is regular, they are the Betti numbers of the
K\"ahler manifold $M/\Sigma$. This allows 
to relate the Betti numbers of $M$ with the $\Sigma$-basic
Betti numbers. When $\Sigma$ is not
regular, we have a similar result:
\index[terms]{Betti numbers!basic}\index[terms]{Betti numbers!rational}

\hfill

\corollary {(\cite{va_gd})}\label{betti_basic} On a
compact Vaisman manifold,\index[terms]{manifold!Vaisman} the $\Sigma$-basic Betti numbers $e_i$ are
finite and they are related to the  Betti numbers $b_i$
of $M$ by the relations:
\begin{equation*}
\begin{split}
b_i&=e_i+e_{i-1}-e_{i-2}-e_{i-3}, \qquad 0\leq i\leq n-1,\\
b_i&=e_{i-2}+e_{i-1}-e_{i}-e_{i+1}, \qquad n+1\leq i\leq 2n,\\
b_n&=2(e_{n-1}-e_{n-3}).
\end{split}
\end{equation*}
Using the Hodge decomposition on $\Sigma$-basic
cohomology, we obtain that $e_{2k}\neq 0$ and $e_{2k+1}$
are even (\ref{_Transversal_Lefschetz_Theorem_}). In particular:

\hfill

\corollary\label{_odd_first_betti_} A compact Vaisman manifold has odd first Betti number.\index[terms]{manifold!Vaisman}

%
%
%
%


\section[The supersymmetry algebra of a Sa\-sa\-kian ma\-ni\-fold]{The supersymmetry algebra of a Sa\-sa\-kian\\ ma\-ni\-fold}\label{_Appendix_C_}	


Let $Q$ be a Sasakian manifold.
In this section  we describe the Lie superalgebra $\goth q\subset \End(\Lambda^*Q)$ \index[terms]{supersymmetry algebra!Sasaki}
reminiscent of the supersymmetry algebra of a K\"ahler manifold 
(\ref{_kah_susy_Theorem_}). Unlike the K\"ahler supersymmetry algebra,
the algebra $\goth q$ is infinite-dimensional; however, it has a simple
and compact description, independent of the choice of $Q$.

\hfill

\remark Some of the relations we derive can be found in P.A.-Nagy's doctoral thesis, in a more general setting, for a compact Riemannian manifold endowed with a unitary vector field (see \cite{_Nagy_}). When the Sasakian manifold is also Einstein,
the superalgebra structure seems to be simpler.
Using this approach, J. \index[persons]{Schmude, J.} Schmude (\cite{_Schmude_})
obtained a closed formula for the de Rham Laplacian
in terms of the transversal Laplacian operator.

\hfill


We first recall the notations. Let $Q$ be a Sasakian manifold, $\xi $ its Reeb field,
and $R\subset TM$ the corresponding rank 1 distribution. 
Let 
$$\Lambda^*(Q):= \bigoplus\Lambda^{p,q}_\hor(Q)\otimes
\Lambda^m_\ver(Q)$$
be the corresponding decomposition
(Hodge and \index[persons]{Hattori, A.} Hattori) of the de Rham algebra, and 
$d=d_0+d_1+d_2$ the Hattori differentials (see Subsection
\ref{_LS_decompo_Sasakian_Subsection_}).
Denote by $W$ the Weil operator acting as multiplication
by  $\1(p-q)$ on
$\bigoplus\Lambda^{p,q}_\hor(Q)\otimes\Lambda^m_\ver(Q)$. 
Let $d_1, d_1^*, d_1^c$ and $(d_1^c)^*$ be the
differentials defined in \ref{_Hodge_decompo_d_1_Sasa_Claim_}.
Let $L_{\omega_0}, \Lambda_{\omega_0}:= L_{\omega_0}^*, H_{\omega_0}:= [L_{\omega_0}, \Lambda_{\omega_0}]$
be the Lefschetz $\goth{sl}(2)$-triple associated with the
transversal K\"ahler form \index[terms]{form!K\"ahler!transversal}$\omega_0$.
For an operator $A\subset \End(Q)$, 
commuting with $\Lie_\xi $, we denote by
$A(k)$ the composition $A\circ (\Lie_\xi )^k$.
Let $i_\xi $ be the
contraction with $\xi $ and $e_\xi $ the dual operator\footnote{Normally, we should write $e_{\xi ^{\, \flat}}$ or $e_\eta$, but in the context of this section we find the notation $e_\xi$ more transparent.}.

Denote by $\goth q$ the Lie superalgebra generated
by all operators 
\[ L_{\omega_0}(i),\ 
\Lambda_{\omega_0}(i),\ H_{\omega_0}(i), \ d_1(i),\  e_\xi (i),\ 
i_\xi (i),\  \Id(i), W(i)\ \text{for all $i\in \Z^{\geq 0}$}.
\] 
Since the vector field $\xi $ acts by Sasakian
isometries, all elements of $\goth q$ commutes with $\Lie_\xi $;
in other words,  $\Lie_\xi $ is central in $\goth q$.
We consider $\goth q$ as an $\R[t]$-module, with
$t$ mapping $A(i)$ to $A(i+1)$.
Then the Lie superalgebra $\goth q$ is a free 
$\R[t]$-module of rank $(6|6)$ (with 6 even and 6 odd
generators) over $\R[t]$. 
Its even generators are 
\[ L_{\omega_0},\ 
\Lambda_{\omega_0},\  H_{\omega_0},\  W,\  \Delta_1:= \{d_1,
d_1^*\},\  \Id,
\]
and the odd generators are 
$$d_1, \ d_1^*,\  d_1^c,\  (d_1^c)^*,\  e_\xi ,\  i_\xi .$$
Notice that $d_0=e_\xi (1)$ and $d_0^*=i_\xi (1)$ (\ref{_d_0_expli_Claim_}).

The main result of this section is:

\hfill

\theorem \label{_susy_Sasa_Theorem_}
The only non-zero commutator relations in $\goth q$
can be written as follows.
\begin{description}
	\item[(i)] The Lefschetz $\gsl(2)$-action: the even elements
	$L_{\omega_0},\Lambda_{\omega_0},H_{\omega_0}$ 
	satisfy the u\-su\-al $\gsl(2)$-relations (\ref{_kah_susy_Theorem_}) 
	and commute with 
	\[
	W, \ \Delta_1,\  \Delta_0, \ d_0,\  d_0^*.
	\]
	The operator
	$H_{\omega_0}$ acts as multiplication by 
	$p-n$ on $\Lambda^{p}_\hor(Q)\otimes\Lambda^m_\ver(Q)$,
	where $\dim_\R Q= 2n+1$.
	
	\item[(ii)] The Weil operator satisfies 
	\[
	[W, d]=d_1^c, \ [W, d^c_1]= - d_1,\ 
	[W, d_1^*]=-(d_1^c)^*,\  [W, (d^c_1)^*]= d_1.
	\]
	Also, $W$ commutes with the rest of the generators of $\goth  q$.
	
	\item[(iii)] The differentials $d_1, d_1^*, d_1^c, (d_1^c)^*$ have non-zero square:
	\[ 
	\{d_1, d_1\}=\{d_1^c, d_1^c\}=-L_{\omega_0}(1), \ \ 
	\{d_1^*, d_1^*\}=\{(d_1^c)^*, (d_1^c)^*\}=\Lambda_{\omega_0}(1).
	\]
	Moreover, $\{d_1, d_1^c\}=\{d_1^*, (d_1^c)^*\}=0$.
	\item[(iv)] 
	The usual \index[terms]{Kodaira identities} Kodaira identities still hold:
	\begin{equation}\label{_Kodaira_rel_Sasakian_Equation_}
	\begin{split}
	[\Lambda_{\omega_0}, d_1]& = (d_1^c)^*,\ \ \  \ \,
	\ \ [ L_{\omega_0}, d_1^*] = - d_1^c,\\
	[\Lambda_{\omega_0}, d_1^c]& = - d_1^*,
	\ \ \ \ [ L_{\omega_0}, (d^c_1)^*] = d_1.
	\end{split}
	\end{equation}
	
	\item[(v)] Unlike it happens in the K\"ahler case,
	the differentials $d_1^*, d_1^c$, etc. do not supercommute:
	\begin{equation}\label{_commutators_v_}
	\{d_1^*, d_1^c\}=\{d_1,
	(d_1^c)^*\}= -\frac 1 2 H_{\omega_0}(1).
	\end{equation}
	\item[(vi)] The only non-zero commutator  between
	$e_\xi , i_\xi , \Id$ is $\{e_\xi , i_\xi \}=\Id$.
	These elements commute with the rest of $\goth q$.
	
	\item[(vii)] The Laplacian $\Delta_1:= \{d_1, d_1^*\}$ 
	satisfies $\Delta_1= \{ d_1^*, (d_1^c)^*\}$
	and commutes with all
	even generators in $\goth q$ and with $e_\xi , i_\xi $.
	Its commutators with the other 4 odd generators are expressed as
	follows
	\begin{equation}\label{_Delta_1_Equation_}
	\begin{split}
	\{d_1,\Delta_1\}&= -\frac 12 d_1^c(1), \ \ \ \ \ \ \ \ \ \  \{d_1^c,\Delta_1\}= \frac 12
	d_1(1),\\
	\{d_1^*,\Delta_1\}&= - \frac 12 (d_1^c)^*(1),  
	\ \  \{(d_1^c)^*,\Delta_1\}= \frac 12 (d_1)^*(1).
	\end{split}
	\end{equation}
\end{description}

{\bf Proof of \ref{_susy_Sasa_Theorem_} (i):}
The $\goth{sl}(2)$-relations and the expression for 
$H_{\omega_0}$ are proven in Subsection
\ref{_susy_Kah_Subsection_}; the proof in the Sasakian
case is literally the same. Also, from the definition
it is clear that
$L_{\omega_0},\Lambda_{\omega_0},H_{\omega_0}$ 
commute with $W$, $i_\xi $ and $e_\xi $.
Since $d_0=e_\xi (1)$ and $d_0^*=i_\xi (1)$ 
(\ref{_d_0_expli_Claim_}), these operators
commute with the $\goth{sl}(2)$-action
and satisfy $\{d_0, d_0^*\}=\Lie^2_\xi $.
We postpone the commutator relation for $\Delta_1$ until
we proved \ref{_susy_Sasa_Theorem_} (vii).

\hfill

{\bf Proof of \ref{_susy_Sasa_Theorem_} (ii):}
Same as \ref{_kah_susy_Theorem_}, part 2.

\hfill

{\bf Proof of \ref{_susy_Sasa_Theorem_} (iii):}
Start from $\{d_1, d_1\}=-L_{\omega_0}(1)$.

Since $d^2=0$, one has $d_0^2=d_2^2=0$ and
$\{d_0, d_2\}=-\{d_1, d_1\}$. However,
$d_0= e_{\xi }(1)$ and $d_2=L_{\omega_0} i_\xi $,
hence $\{d_0, d_2\}=\{e_{\xi }, i_\xi \}
L_{\omega_0}(1)=L_{\omega_0} (1)$.
The squares of the rest of the differentials $d_1^c$, etc.,
are obtained by duality and complex conjugation.

To show that $\{d_1, d_1^c\}=0$, we use
$d_1^c=[W, d_1]$. By \ref{_d^2neq0_super_Claim_},
\[
\{d_1, d_1^c\}= \{ d_1, \{d_1, W \}\}= \frac 1 2 \{\{ d_1, d_1\}, W\}.
\]
Since $\{ d_1, d_1\}=-L_{\omega_0}(1)$ and $W$ commutes
with $L_{\omega_0}$, this gives $\{d_1, d_1^c\}=0$. The
equation $\{d_1^*, (d_1^c)^*\}=0$ is obtained by duality.
We finished \ref{_susy_Sasa_Theorem_} (iii).

\hfill

{\bf Proof of \ref{_susy_Sasa_Theorem_} (iv), the K\"ahler-Kodaira relations:}
Again, it suffices to prove 
$ [L_{\omega_0}, d_1^*] = - d_1^c$,
the rest is obtained by duality and complex conjugation.
We prove it by applying the same argument as used in
\ref{_kah_susy_Theorem_}.
As an operator on $\Lambda^*(Q)$, the 
commutator  $[L_{\omega_0}, d_1^*]$ has first order,
because $L_{\omega_0}$ is zero order, and $d_1^*$ is order
2 (\ref{_diff_ope_commutator_Claim_}).

The differential operators $[L_{\omega_0}, d_1^*]$
and $- d_1^c$ are equal on functions for the same
reasons as in \ref{_kah_susy_Theorem_}.

Clearly, $d_1^c(C^\infty (Q))$ generates
$\Lambda^1_\hor(Q)$.  Therefore, to prove  
$[L_{\omega_0}, d_1^*]=- d_1^c$ on 
$\Lambda^1_\hor(Q)$, we need only to show that
\begin{equation}\label{_commu_L_d_1^*_square_Equation_}
([L_{\omega_0}, d_1^*])^2= (- d_1^c)^2.
\end{equation}
This is implied by graded Jacobi identity, applied as follows.
First, we notice that $[\Lambda_{\omega_0}, d_1^*]=0$,
because $d_1(\omega_0)=0$. Then, the $\goth{sl}(2)$-representation
generated by the Lefschetz triple 
$\langle
L_{\omega_0},\Lambda_{\omega_0},H_{\omega_0}\rangle$
from $d_1^*$ has weight 1, and
$[L_{\omega_0},[L_{\omega_0}, d_1^*]]=0$.
Applying the graded Jacobi identity, and using
$[L_{\omega_0},[L_{\omega_0}, d_1^*]]=0$,
we obtain
\begin{equation*}
\begin{split}
\{[L_{\omega_0}, d_1^*], [L_{\omega_0}, d_1^*]\}&=
[L_{\omega_0},\{d_1^*, [L_{\omega_0}, d_1^*]\}]\\[.1in] 
&=\frac 12 [L_{\omega_0}(1), [L_{\omega_0}(1), \{d_1^*,d_1^*\}]]\quad \text{by \ref{_d^2neq0_super_Claim_}}\\
&=-\frac 1 2 [L_{\omega_0}(1), [L_{\omega_0}(1), \Lambda_{\omega_0}]]\quad \text{by \ref{_susy_Sasa_Theorem_} (iii)}\\
&=-L_{\omega_0}(1).
\end{split}
\end{equation*}

However, $(- d_1^c)^2=-L_{\omega_0}(1)$ by
\ref{_susy_Sasa_Theorem_} (iii), which proves
\eqref{_commu_L_d_1^*_square_Equation_}.

This implies that  
$ [L_{\omega_0}, d_1^*] = - d_1^c$
on $\Lambda^1_\hor(Q)$ .
The space $\Lambda^1(Q)$ is generated over
$C^\infty (Q)$ by $\Lambda^1_\hor(Q)$ and 
$V= \langle \xi ^{\,\flat} \rangle$.
Applying \ref{_d^2=0_dete_Corollary_}, we obtain that
$ [L_{\omega_0}, d_1^*] = - d_1^c$ if 
$ [L_{\omega_0}, d_1^*](\xi ^{\,\flat})=-d_1^c(\xi ^{\,\flat})$.
However, $d(\xi ^{\,\flat})=\omega_0$, and hence 
$d_1^c(\xi ^{\,\flat})=0$. Then, 
$ [L_{\omega_0}, d_1^*](\xi ^{\,\flat})=
d_1^*(\xi ^{\,\flat}\wedge \omega_0)$.
It is not hard to see that 
$*(\xi ^{\,\flat}\wedge \omega_0)=
\frac{1}{(n-1)!}\omega_0^{n-1}$, where $\dim_\R Q=2n+1.$
Therefore, 
\[ d^*_1(\xi ^{\,\flat}\wedge \omega_0) =\pm
*d_1*(\xi ^{\,\flat}\wedge \omega_0)=
\pm *d_1 \frac{1}{(n-1)!}\omega_0^{n-1}=0,
\]
because $d\omega_0^{n-1}=0$.
We proved that $ [L_{\omega_0},
d_1^*](\xi ^{\,\flat})=-d_1^c(\xi ^{\,\flat})$
and finished the proof of the K\"ahler-Kodaira relations.

\hfill

{\bf Proof of \ref{_susy_Sasa_Theorem_} (v), the commutator of $d_1^*, d^c$:}
The commutator $\{d_1^*, d^c\}$ is obtained from the
K\"ahler-Kodaira relations. Indeed, we have 
$\{d_1^*, d^c\}=\{d_1^*, \{d_1^*, L_{\omega_0}\}\}$.
Then \ref{_d^2neq0_super_Claim_} gives
\[ \{d_1^*, \{d_1^*, L_{\omega_0}\}= \frac 1 2 \{\{d_1^*,
d_1^*\}, L_{\omega_0}\}= \frac 1 2 [\Lambda_{\omega_0}(1),
L_{\omega_0}]=- \frac 1 2 H_{\omega_0}.
\]
The relation $\{d_1,
(d_1^c)^*\}=-\frac 1 2 H_{\omega_0}(1)$ is dual to
$\{d_1^*,d_1^c\}=-\frac 1 2 H_{\omega_0}(1)$.

\hfill

{\bf Proof of \ref{_susy_Sasa_Theorem_} (vi), the commutator of $e_\xi $ and $i_\xi $:}
The equation $\{e_\xi , i_\xi \}=1$ is standard. Vanishing of the commutators
between $e_\xi$, $i_\xi$ and  $L_{\omega_0}$, 
$\Lambda_{\omega_0}$, $H_{\omega_0}$, $W$ is standard linear algebra.
The only commutators for that we have to prove the vanishing
is between $e_\xi $, $i_\xi $ and $d_1, d_1^*, d_1^c$, $(d_1^c)^*$.
Using duality and complex conjugation, we reduce the vanishing
of these commutators to only two of them:
$\{e_\xi , d_1\}=0$ and $\{e_\xi , d_1^*\}=0$.
As $d^2=0$ and $\{d_0, d_1\}$ is the grading 1 part of
$d^2$, one has $\{d_0, d_1\}=0$. Since $d_0= e_\xi (1)$,
this also implies $\{e_\xi , d_1\}=0$. Twisting with $I$, we obtain
$\{e_\xi , d_1^c\}=0$. Applying the graded Jacobi identity to
$\{e_\xi , d_1^*\}= -\{ e_\xi , \{\Lambda_{\omega_0}, d_1^c\}\}$
(\ref{_susy_Sasa_Theorem_} (iv)) and using
$\{e_\xi , L_{\omega_0}\}=0$, we obtain
\[
\{e_\xi , d_1^*\}= -\{ e_\xi , \{\Lambda_{\omega_0}, d_1^c\}\}=
\{\{e_\xi , \Lambda_{\omega_0}\},d_1^c\}+ \{\Lambda_{\omega_0}, \{e_\xi , d_1^c\}\}=0.
\]
This finishes the proof of \ref{_susy_Sasa_Theorem_} (vi).

\hfill

{\bf Proof of   \ref{_susy_Sasa_Theorem_} (vii).}
The equation
\[ \{d_1, d_1^*\}=\{(d_1^c), (d_1^c)^*\}
\]
follows from $\{d_1, d_1^c\}=0$ (\ref{_susy_Sasa_Theorem_} (iii))
because 
\[
0=\{\Lambda, \{d_1, d_1^c\}\}= \{\{\Lambda,d_1\}, d_1^c\}+ \{d_1, \{\Lambda,d_1^c\}\}=
\{(d_1^c)^*,  d_1^c\}-\{d_1, d_1^*\}.
\]
This implies, in particular, that $[W,\{d_1, d_1^*\}]=0$.
The commutators between the Lefschetz operators and $\Delta_1$
are computed using the K\"ahler-Kodaira relations:
\[
\{L_{\omega_0}, \{d_1, d_1^*\}\}= \{\{L_{\omega_0}, d_1\}, d_1^*\}
+ \{d_1, \{L_{\omega_0},d_1^*\}\}=-\{d_1,d_1^c\}=0.
\]
We proved that $\Delta_1$ commutes with the even part of $\goth q$.
By duality and complex conjugation, to prove \eqref{_Delta_1_Equation_}
it would suffice to prove only one of these relations, say,
$\{d_1,\Delta_1\}= -\frac 12 d_1^c(1)$. This equation 
follows from \eqref{_Kodaira_rel_Sasakian_Equation_},  \eqref{_commutators_v_} (i.\,e.   (iv) and (v) of this theorem) and \ref{_d^2neq0_super_Claim_}:
\[
\{d_1,\{d_1, d_1^*\}\}=\frac 1 2 \{\{d_1, d_1\}, d_1^*\}=
-\frac 1 2 \{L_{\omega_0}(1), d_1^*\}=-\frac 1 2 (d_1^c)^*(1).
\]
We finished the proof of \ref{_susy_Sasa_Theorem_}. \endproof

\section{Notes} 

\begin{enumerate}

\item
Toyoko \index[persons]{Kashiwada, T.} Kashiwada proved the harmonic decomposition\index[terms]{decomposition!harmonic}
for Vaisman manifold\index[terms]{manifold!Vaisman} in \cite{kashiwada_kodai} and in
\cite{KS}, joint with Shizuko \index[persons]{Sato, S.} Sato.
This paper is notable because it was published in 
Ann. Fac. Sci. Univ. Nat. Za\"{\i}re (Kinshasa);
it was not widely available. Kashiwada is more
known for another of her results, \cite{kashiwada}: she proved
that any 3-Sasakian manifold is Einstein (\ref{_3_Sasaki_Einstein_}).

\item The connection between de Rham calculus on manifolds
with special geometry (such as K\"ahler and hyperk\"ahler)
and their supersymmetry appeared as early as in 1997 
(\cite{_FKS_}). It is well known that extra supersymmetries
of the $\sigma$-model force the target space to acquire
extra geometric structures: the $N=1$ supersymmetry 
implies K\"ahler structure on the target space,
the $N=2$ supersymmetry makes it hyperk\"ahler, and so on.
In \cite{_FKS_}, this supersymmetry was interpreted
in terms of the de Rham calculus on the target space.

In \cite{_Kim_Saberi_}, the connection between 
the supersymmetry and rational homotopy theory is further
expounded, with the constructions of rational homotopy
theory (such as \index[persons]{Sullivan, D.} Sullivan's minimal models) interpreted
in terms of quantum mechanics.

For Sasakian manifolds, this approach was pioneered
by \cite{_Tievsky_}, who used the transversal K\"ahler relations
to obtain results about rational homotopy of
Sasakian manifolds. This work was
applied to homotopy formality of Sasaki manifolds which
are not formal, but satisfy different weaker types of
formality  
(\cite{_Biswas_Fernandez_Tralle_}, \cite{_Papadima_Suciu_}, \cite{cape3}).
The same method was applied to study the geometry of Sasakian
nilmanifolds (\cite{cape1}).\index[terms]{nilmanifold}

\item Let $(M,I,\omega,\theta)$ be a compact Vaisman manifold, $\dim_\C M=n$.\index[terms]{manifold!Vaisman}  \ref{_Vaisman_harmonic_forms_Theorem_} implies the existence of isomorphisms
\[
H^k(M)\stackrel{\sim}{\arrow} H^{n-k}(M),\qquad 0\leq k\leq n-1.
\]
In \cite{cape2}, these  are called ``hard
Lefschetz'' (and their existence is named ``Hard Lefschetz
theorem''). Moreover, their action was constructed 
explicitly on the level of differential forms 
(\cite[Theorem 1.1]{cape2}). 

\item  
 For a Vaisman manifold,\index[terms]{manifold!Vaisman} one can  compute the
dimensions of the groups appearing in the second page of
the \index[persons]{Hattori, A.} Hattori spectral sequence of the canonical foliation\index[terms]{foliation!canonical} $\Sigma$. This leads  to several inequalities
between these invariants and the basic Betti numbers. Namely (\cite{os_spec}:
\[
\dim E_2^{u,1}\ge 2\, e_u, \quad \text{for} \quad 0\le u\le 2n-2.
\]
In particular, $\dim E_2^{0,1}\ge 2$. If $\Sigma$ is quasi-regular, then $\dim E_2^{u,1}=2\cdot e_u$.

These inequalities can be interpreted as an
obstruction for a Riemannian complex two-dimensional
foliation on a compact complex manifold to be the
canonical foliation of a Vaisman structure,
\end{enumerate}

\section{Exercises}
\begin{enumerate}[label=\textbf{\thechapter.\arabic*}.,ref=\thechapter.\arabic{enumi}]
	
	\item
	Let $A^*$ be a graded commutative
	algebra, and $D:\; A^* \arrow A^{*+i}$ be a linear map
	such that for all $x\in A$ there exists $N$ such that
	$D^N(x)=0$. Prove that $e^D:= 1 + D + \frac
	{D^2}{2} + \cdots + \frac{D^i}{i!} + \cdots $ is an automorphism
	of $A^*$ if and only if $D$ is a derivation.
	
	\item
	Prove that a supercommutator of graded derivations
	is also a graded derivation.

	\item
	Let $\tau:\; \Lambda^*(M) \arrow \Lambda^{*-1}(M)$
	be a derivation shifting grading by $-1$.
	Prove that there exists a vector field $v\in TM$ such that
	$\tau=i_v$, or find a counterexample.
	
	\item
	Let $(M,I)$ be an almost complex manifold,
	$d:\; \Lambda^*(M) \arrow \Lambda^{*+1}(M)$
	the de Rham differential, and $d=\bigoplus_{p+q=1}d^{p,q}$
	its Hodge decomposition. Prove that $d^{p,q}=0$ for all
	$p>2$. Prove that the operator $d^{2,-1}$ is
	$C^\infty(M)$-linear. 
	
	\item In these assumptions,
	prove that  $(d^{1,0})^2=0$ $\Leftrightarrow$ 
	the almost complex structure $I$ is integrable.
	
	\item 
	Let $\eta$ be a parallel differential form on a Riemannian manifold,
	and $\eta'$ a harmonic form. Prove that $\eta\wedge \eta'$
	is harmonic.

	\item Let $(M,I,g)$ be an almost complex Hermitian manifold,
	and $\omega$ its fundamental form.
	Assume that $d\omega=0$.\footnote{In this case, $(M,I,g)$
		is called ``almost K\"ahler''.}
	Prove that $\omega$ is harmonic.

\item Let $X$ be a holomorphic vector field on a complex
  manifold, and $X^c=I(X)$.
\begin{enumerate}
\item Prove that $\{d^c, i_X\}= - \Lie_{X^c}$.
\item
Prove that $\{\bar\6, i_X\}= \frac 1 2 (\Lie_X - \1
\Lie_{X^c})$, where $i_X$ is the contraction operator.\footnote{This result was proven in \cite{_Klemyatin_}.}
\end{enumerate}


\item
Let $X$ be a Killing\index[terms]{vector field!Killing} field on a compact Riemannian manifold $M$.
Prove that $M$ admits an action of a compact torus $G$ by isometries,
such that $X$ belongs to the Lie algebra of $G$.

\item\label{_Fourier_for_torus_on_de_Rham_Exercise_}
Let $G= (S^1)^r$ be a compact torus acting on a Riemannian manifold
$M$ by isometries. Prove that any differential form 
$\alpha$ on $M$ can be expressed as $\alpha= \sum_{j_1, ..., j_r\in \Z} \alpha_{j_1, ..., j_r}$,
with $t= \exp\left({\sum_{i=1}^r 2\pi \1 t_i}\right)$ acting on $\alpha_{j_1, ..., j_r}$ as on the
appropriate Fourier monomial,
\[ t(\alpha_{j_1, ..., j_r})= \exp\left(\sum_{i=1}^r 2\pi \1 j_i t_i\right)\alpha _{j_1, ..., j_r}.
\]
{\em Hint:} Use the Fourier transform over a torus.

\item
Let $M$ be a compact Riemannian manifold, $\theta\in \Lambda^1(M, \R)$
a closed 1-form on $M$, $d_\theta= d+ \theta\wedge$ the
twisted differential, and $X:=\theta^\sharp$ the dual
vector field. Assume that $X$ is Killing.\index[terms]{vector field!Killing}
\begin{enumerate}
\item Prove that $\{d_\theta, i_X\}= \Lie_X +
  |X|^2\Id$. 

\item Let $\alpha$ be a $d_\theta$-closed form.
Prove that $\alpha = \sum c_i \alpha_i$
where $c_i$ are scalars and $\alpha_i$ are 
$d_\theta$-closed forms that satisfy
$\Lie_X \alpha_i = w_i \alpha_i$, with $w_i$ constant.

\item Let $\alpha$ be a $d_\theta$-closed form
such that $\Lie_X \alpha = w \alpha$, where $w\in \C$ is constant.
Prove that $|X|^2\alpha$ is $d_\theta$-closed and
$d_\theta$-cohomologous to $-w\alpha$.

\item
 Let $\alpha$ be a $d_\theta$-closed form.
Prove that $|X|^2\alpha$ is $d_\theta$-closed.

\item Consider the action of $\Lie_X$ on
  $d_\theta$-cohomology $H^*_\theta(M, \R)$. Prove that for each
$\alpha \in H^*_\theta(M)$ such that $\Lie_X(\alpha) = w \alpha$, 
the constant $w$ is real.

\item 
Prove that each class $\alpha \in H^*_\theta(M, \R)$
can be represented by a form that satisfies $\Lie_X \alpha=0$.
\end{enumerate}

{\em Hint:} Use the previous exercise.

\item
Let $M$ be a Vaisman manifold,\index[terms]{manifold!Vaisman} $\xi$ the Reeb field, 
$d_\theta$ the twisted differential, $\6_\theta$ and
$\6_\theta^*$ its Dolbeault components, and $d_\theta^c:=I d_\theta I^{-1}$.
\begin{enumerate}
\item Prove that $\{i_\xi, d_\theta\}=\Lie_\xi$.
\item Prove that $\{i_{I\xi}, d_\theta\}=\Lie_{I\xi} +\Id$.
\item Prove that $\{i_{\xi}, d_\theta^c\}=-\Lie_{I\xi}+\Id$.
\item Prove that $\{i_{I\xi}, d_\theta^c\}=\Lie_{\xi}$.
\item Prove that $\{i_{\xi}, d_\theta^*\}=\Lambda_{\omega_0}$.
\item Prove that 
\[ \{i_{\xi},\bar\6_\theta\}=\frac1 2(\1\Id +
  \Lie_\xi +\1\Lie_{I\xi}).
\] 
\item 
Prove that 
\[ \{i_{I\xi},\bar\6_\theta\}=\frac1 2(\Id 
  +\Lie_{I\xi} +\1\Lie_{\xi}).
\]
\end{enumerate}

\item\label{_supersymmetry_hyperkahler_Exercise_}
A Riemannian manifold $M$ is called {\bf hyperk\"ahler}
if it is equipped with an action of the quaternion algebra in $TM$,\index[terms]{algebra!quaternion}
orthogonal and parallel with respect to the Levi--Civita connection.
Let $(M, I, J, K, g)$ be a hyperk\"ahler manifold.
\begin{enumerate}
\item Let $L=aI +b J + cK$ be a quaternion,
$L^2=-1$. Prove that $(M,L,g)$ is K\"ahler.
\item Let $d_I= I d I^{-1}, d_J= J dJ^{-1}, d_K=K d K^{-1}$ be twisted 
differentials, and ${\goth a}$ the Lie algebra generated by 
Lefschetz $\goth{sl}(2)$-triples for $I,J,K$.
Prove that the action of $\goth a$ on $\End(\Lambda^*(M))$ preserves
the 8-dimensional space $\langle d, d_I, d_J, d_K, d^*,  d_I^*, d_J^*, d_K^*\rangle$.

\item
Let $L_I, L_J, L_K$ and $\Lambda_I, \Lambda_J, \Lambda_K$ be the
generators of the Lie algebra $\goth a$. Prove that $L_I, L_J, L_K$
are 0-th order differential operators on the algebra $\Lambda^*(M)$
and  $\Lambda_I, \Lambda_J, \Lambda_K$ are 2-th order differential 
operators on $\Lambda^*(M)$. Prove that the commutators
$[L_u, \Lambda_v]$ are derivations of $\Lambda^*(M)$
for any $u, v= I, J,K$.

\item Using \ref{_deriva_dete_Claim_}, 
prove that $[L_I, \Lambda_J] = W_K$, $[L_J, \Lambda_K] = W_I$,
and $[L_K, \Lambda_I] = W_J$, where $W_I$, $W_J$, $W_K$
are the Weil operators (\ref{_KdR_Definition_}). 

\item Prove that the algebra $\goth a$ is 10-dimensional.

\item
Prove that any two generators of $\langle d, d_I, d_J, d_K, d^*,  d_I^*, d_J^*, d_K^*\rangle$
supercommute, except $\{d, d^*\}=\{d_I, d^*_I\}=\{d_J, d^*_J\}=\{d_K, d^*_K\}=\Delta$.
Prove that the metric $x, y \mapsto \frac{\{x, y\}}{\Delta}$
is non-degenerate of signature $(4,4)$.

\item Prove that $\goth a$ acts on
$\langle d, d_I, d_J, d_K, d^*,  d_I^*, d_J^*, d_K^*\rangle$
by isometries, preserving the metric $x, y \mapsto \frac{\{x, y\}}{\Delta}$.

\item
Let $V=\R^4=\Bbb H$ be a 1-dimensional quaternionic vector space.
Consider an algebra $H$ of operators commuting with the 
action of ${\Bbb H}$ on $V$. Prove that $H$ is isomorphic 
to the quaternionic algebra (one of these actions can be interpreted
as the left action by quaternions, the other one as the right action).

\item
Prove that $W_I, W_J, W_K$ generate
a quaternionic action on the space
$Q:=\langle d, d_I, d_J, d_K, d^*,  d_I^*, d_J^*, d_K^*\rangle$.
Let $H$ be an algebra of operators acting
on $\langle d, d_I, d_J, d_K\rangle$ and commuting with this
quaternionic action. Prove that $H$ is isometric to quaternions.

\item
Prove that $\goth a$ acts on $Q$ by $H$-invariant isometries, defining
a homomorphism from $\goth a$ to the Lie algebra $\goth{u}({\Bbb H}, 1,1)= \goth{sp}(1,1)$
of $H$-linear isometries. Prove that this map is actually an isomorphism.

\item
Let $M$ be a compact hyperk\"ahler torus, $M= (S^1)^4$.
Consider the action of $\goth a$ on the 5-dimensional space
$P$ generated by $1, \omega_I, \omega_J, \omega_K$, and $\Vol_M$.
Prove that this action preserves the metric
$x, y \arrow \int_M x\wedge y$ on $P$, defining
a homomorphism $\rho:\; \goth a \arrow  \goth{so}(P)$.

\item Prove that the bilinear symmetric form
$x, y \arrow \int_M x\wedge y$ on $P$ has signature $(1,4)$.
Prove that the homomorphism $\rho:\; \goth a\arrow \goth{so}(1,4)$
is an isomorphism. Use this observation to deduce that
$\goth{so}(1,4)$ is isomorphic to $\goth{sp}(1,1)$.

\item
Let $\goth g$ be a Lie superalgebra generated by
$\goth a$ and $Q$.
Prove that $\goth g$ has dimension $11|8$ (even 11, odd 8),
and is spanned by $\goth a, Q$ and $\Delta$, that is  central in $\goth g$.
\end{enumerate}

\item\label{_kernel_pullback_H^2_Exercise_}
Let $M= \frac{\Tot^\circ(L)}{\langle \lambda\rangle}$
be a Vaisman manifold\index[terms]{manifold!Vaisman} obtained as a quotient of a 
total space $\Tot^\circ(L)$ of the space 
of non-zero vectors of a negative line bundle on a 
projective orbifold $X$
(\ref{_Structure_of_quasi_regular_Vasman:Theorem_}).
Denote by $\tau:\; M \arrow X$ the corresponding elliptic
fibration. Prove that the kernel of the corresponding
map $\tau^*:\; H^2(X) \arrow H^2(M)$ is generated by $c_1(L)$.

\end{enumerate}


\chapter{Dolbeault cohomology of LCK manifolds with potential}
\label{_Dolbeault_Vaisman_Chapter_}\index[terms]{manifold!LCK!with potential}\index[terms]{cohomology!Dolbeault}

{\setlength\epigraphwidth{0.8\linewidth}
\epigraph{\it The attempt at art is in itself a virtue. }
{\sc \scriptsize David Burlyuk, ``From Now On, I Refuse to Speak Ill Even of the Work of Fools''}}


\section{Introduction}

In this chapter, we compute the Dolbeault cohomology\index[terms]{cohomology!Dolbeault}
of Vaisman manifolds.\index[terms]{manifold!Vaisman} We follow the same scheme as
in Chapter \ref{_Harmonic_forms_chapter_}: we try to represent
the cohomology classes by Lee- and anti-Lee-invariant forms,
and compute the cohomology of the algebra of the 
Lee- and anti-Lee-invariant forms.

For de Rham cohomology,\index[terms]{cohomology!de Rham} this part of the argument is
transparent. However, on non-K\"ahler manifolds, the
connected component of the automorphism group might act 
non-trivially on the Dolbeault cohomology.\index[terms]{cohomology!Dolbeault}

To circumvent this obstruction, we use a recent work of
N. \index[persons]{Klemyatin, N.} Klemyatin (\cite{_Klemyatin_}), who studied the action of the automorphism
group of a manifold on its Dolbeault cohomology. 

It is not hard to see that any path-connected group 
acts trivially on the de Rham cohomology.\index[terms]{cohomology!de Rham} Using the
Hodge to de Rham isomorphism, it also follows 
that any path-connected group of complex automorphisms
of a K\"ahler manifold acts trivially on its
Dolbeault cohomology. \index[terms]{cohomology!Dolbeault}

For Vaisman manifolds,\index[terms]{manifold!Vaisman} we eventually obtain that
the Dolbeault cohomology is isomorphic to the de Rham
cohomology (\ref{_E_1_degenerates_for_Vaisman_Remark_};
the same result is true for the LCK manifolds with potential,
\ref{_LCK_pot_degenerate_Corollary_}),\index[terms]{manifold!LCK!with potential} 
hence the same argument can be applied. 
However, we do not follow this route, because we use
the triviality of the action of the (anti-)Lee field\index[terms]{Lee field}\index[terms]{Lee field!anti-}
on the Dolbeault cohomology to compute the Dolbeault cohomology\index[terms]{cohomology!Dolbeault}
and to prove that Dolbeault is isomorphic to de Rham.

For all compact K\"ahler and many compact non-K\"ahler manifolds,
the Dolbeault cohomology is isomorphic to de Rham, and hence 
the connected component $\Aut_0(M)$ of the complex
automorphism group acts trivially on $H^{p,q}_{\bar\6}(M)$.
In complex dimension 2, Dolbeault is always isomorphic to de Rham
(\ref{_E_2_degenerates_Corollary_}). 

F. \index[persons]{Lescure, F.} Lescure was
the first one to notice that this action can be non-trivial
(\cite{_Lescure:93_}). Later, D. \index[persons]{Akhiezer, D.} Akhiezer found
a general way of constructing examples of manifolds with
non-trivial action of $\Aut_0(M)$ on the Dolbeault cohomology
(\cite{_Akhiezer:group_}).\index[terms]{cohomology!Dolbeault} 

Let $G$ be a complex Lie group, and $\Gamma\subset G$
a cocompact lattice.\index[terms]{lattice!cocompact} \index[persons]{Akhiezer, D.} Akhiezer has shown that 
\[ 
H^{p,q}(G/\Gamma)\cong \Lambda^p(\goth g) \otimes_\C H^{q}(G/\Gamma, \C)
\]
whenever $G$ is reductive. This isomorphism is compatible with the $G$-action,
where the action of $G$ on $\Lambda^p(\goth g)$ 
is induced by the adjoint representation, and the action on $H^{q}(G/\Gamma, \C)$
is trivial.  In this case, the $G$-action on $H^{p,q}(G/\Gamma)$ is non-trivial
for $p \neq 0, n$, and this action on $H^{0,q}(G/\Gamma)$ is trivial.
\index[persons]{Akhiezer, D.} Akhiezer also computes the $G$-action on $G/\Gamma$ for some complex
solvable groups, and shows that, in this case, the action
of $G$ on $H^{0,q}(G/\Gamma)$ can  be non-trivial as well.

Let $\rho$ be a holomorphic flow, that is, an action of $\C$
on a complex manifold $M$. Suppose that the compactification
of the image of $\rho$ is a compact torus, equivalently, that
there exists a Hermitian metric preserved by $\rho$.
\index[persons]{Klemyatin, N.} Klemyatin has shown that, in this situation, the action 
induced by $\rho$ on $H^{p,q}(M)$ is trivial. 
In \ref{_G_acts_on_Dolbeault_trivially_Theorem_}, we give a proof 
of this result. We use a holomorphic version of the Cartan formula,
$\{d^c, i_X\}= - \Lie_{X^c}$, true for any holomorphic vector field $X$ on $M$.
In comparison with the weight decomposition induced
by the torus action, this implies that the torus acts trivially on $H^{p,q}_{\bar\6}(M)$.

From this result, we easily deduce the
explicit expression for the Dolbeault cohomology\index[terms]{cohomology!Dolbeault} (\ref{_Vaisman_Dolbeault_coho_Theorem_})
and the formula for the harmonic decomposition on Vaisman manifolds\index[terms]{manifold!Vaisman}
(\ref{_Dolbeault_harmonic_on_Veisman_}), originally due to K. \index[persons]{Tsukada, K.} Tsukada (\cite{tsuk}).

Recall that any Vaisman manifold is equipped
with a 1-dimensional holomorphic totally geodesic foliation called
{\bf the canonical foliation}\index[terms]{foliation!canonical} (\ref{_Subva_Vaisman_Theorem_}). This
is the foliation generated by the Lee and the anti-Lee field\index[terms]{Lee field}\index[terms]{Lee field!anti-}
(that are  both holomorphic, Killing\index[terms]{vector field!Killing}\index[terms]{vector field!holomorphic} and of constant length,
also by \ref{_canon_foli_totally_geodesic_Remark_}).
A Vaisman manifold is
not K\"ahler; however, this foliation is equipped with
a transversal K\"ahler structure. In other words,
its leaf space is K\"ahler, on charts where the
leaf space is defined, and the transition maps
between such charts induce holomorphic isometries
of the local leaf spaces.

A form on a foliated manifold is called {\bf basic}
if it is locally obtained as the pullback of a
form on the leaf space. The space of basic forms
is preserved by the de Rham differential, and
the corresponding cohomology space is called
{\bf the space of basic cohomology}.\index[terms]{cohomology!basic} It is
an important invariant of a foliated manifold.

In Chapter \ref{_Harmonic_forms_chapter_},
we discussed the transversally K\"ahler structure
associated with the canonical foliation on a\index[terms]{foliation!canonical} 
Vaisman manifold.\index[terms]{manifold!Vaisman} We proved that the 
basic cohomology of a compact Vaisman
manifold has the same properties as on a
compact K\"ahler manifold: the \index[terms]{duality!Poincar\'e} Poincar\'e
duality, the Hodge decomposition, and the
Lefschetz $\goth{sl}(2)$-action.

In this chapter
we express the Dolbeault cohomology 
of a Vaisman manifold in terms of the
transversal Dolbeault cohomology \index[terms]{cohomology!Dolbeault}of the
canonical foliation.\index[terms]{foliation!canonical} 

It is interesting (and highly important, too)
that the total dimension of Dolbeault cohomology
is equal to the corresponding  Betti number:
\[ \bigoplus_{p+q=r} \dim H^{p,q}(M) = \dim H^r(M).\]
This is equivalent to the degeneration of the
Fr\"olicher (a. k. a. Hodge to de Rham)
spectral sequence (\ref{_E_1_degenerates_for_Vaisman_Remark_}).

In general, the Hodge numbers (that is, the
dimension of the Dolbeault cohomology groups) are
not constant in families of complex manifolds.
However, the Hodge numbers are semicontinuous
as functions of the complex structure on a complex manifold
(\cite[Theorem 2.3]{_Bell_Narasimhan_}).
On a K\"ahler manifold, the sum of the Hodge numbers is
equal to the corresponding Betti number. Therefore,
for any smooth family of K\"ahler manifolds, the
Hodge numbers stay constant.

It turns out that the same is true for LCK manifolds
with potential.\index[terms]{manifold!LCK!with potential}

The behaviour of Hodge numbers in families of complex manifolds
is already a classical subject. Originally, it was suggested as a 
question by S. \index[persons]{Iitaka, S.} Iitaka (\cite{_Iitaka:genus_I_}), who noticed that the Hodge numbers
are semicontinuous, and asked if they stay constant in families of complex manifolds. 
Very soon, I. \index[persons]{Nakamura, I.} Nakamura found an example
of a family of complex parallelizable manifolds
with non-constant Hodge numbers.

From the degeneration of the Fr\"olicher spectral sequence,
it is apparent that the jumping phenomenon cannot occur for
the deformations of Vaisman manifolds: the Hodge numbers
of Vaisman manifolds stay constant.\index[terms]{manifold!Vaisman}

The Dolbeault cohomology\index[terms]{cohomology!Dolbeault} of a Vaisman manifold
can be expressed through the basic cohomology\index[terms]{cohomology!basic} of the
canonical foliation,\index[terms]{foliation!canonical} and, conversely, the Hodge numbers
of the basic cohomology can be expressed through 
the Dolbeault cohomology of the Vaisman manifold
(\ref{_Vaisman_Dolbeault_coho_Theorem_}).
This implies, in particular, that the Hodge numbers
of the canonical foliation stay constant.

This result is not as obvious as it might seem,
because the topology of the canonical foliation
can change drastically, even under smooth changes
of the Vaisman structure.

A similar question about the Reeb foliation on
a Sasakian manifold $S$ was asked in \cite{_GNT:semicontinuous_}.
The Reeb foliation is also transversally K\"ahler (\ref{_omega_0_trans_Kahler_Claim_}).
Indeed, the leaf space of the Reeb foliation is identified
with the leaf space of the canonical foliation\index[terms]{foliation!canonical}
of the standard Vaisman structure on $S\times S^1$
(\ref{halfstr}, \ref{_qreg_Vaisman_via_Sasakian_Theorem_}).

In 2021, P. \index[persons]{Ra\'zny, P.} Ra\'zny proved that the 
basic Hodge numbers of the Reeb foliation stay constant in 
smooth families of Sasakian manifolds (\cite{_Razny:Hodge_numbers_}).
His argument was based on \cite[Theorem 7.4.14]{bog},
that is  a consequence of \index[persons]{Tachibana, S.} Tachibana's theorem (\cite{_Tachibana_})
claiming that the transversal Betti numbers stay constant
in smooth families. Then he uses the transversal $dd^c$-lemma
to obtain that the transversal Hodge numbers add up to
transversal Betti numbers. From the semicontinuity of
transversal Hodge numbers (\cite{_GNT:semicontinuous_}),
 Ra\'zny obtains that the transversal Hodge numbers stay constant
as well.

We give a direct proof of \index[persons]{Ra\'zny, P.} Ra\'zny's theorem, which
also includes transversal Hodge numbers of Vaisman manifolds,
even those not related to compact Sasakian ones.
The argument is based on Vaisman geometry.\index[terms]{manifold!Vaisman}
We express the de Rham cohomology\index[terms]{cohomology!de Rham} of Vaisman manifolds
in terms of the basic cohomology;\index[terms]{cohomology!basic} from this expression
it is apparent that the sum of basic Betti numbers
stays constant in the family. Using the harmonic
decomposition on basic cohomology, we obtain that
the sum of basic Dolbeault numbers is equal to the
basic Betti numbers, and hence  it also stays constant
in families. Then the semicontinuity of basic
Dolbeault numbers (\cite{_GNT:semicontinuous_})
implies that it stays constant.

Any Vaisman manifold $M$ of LCK rank 1 \index[terms]{rank!LCK}
admits a K\"ahler $\Z$-cover \index[terms]{manifold!Vaisman}
that is  holomorphically isometric to
the K\"ahler cone over a Sasakian manifold $S$
(\ref{str_vai}).
For such Vaisman manifolds, the local leaf space
of the canonical foliation\index[terms]{foliation!canonical} can be identified
with the local leaf space of the Reeb foliation on
$S$, as mentioned earlier in this introduction
(\ref{halfstr}, \ref{_qreg_Vaisman_via_Sasakian_Theorem_}). 
For LCK rank $>1$\index[terms]{rank!LCK}, the leaf space of the
canonical foliation is not related to 
a foliation on Sasakian manifolds.

Using the deformation theorem proven in
Chapter \ref{_Kuranishi_Chapter_},
the results about the cohomology can be extended
also to LCK manifolds with potential.\index[terms]{manifold!LCK!with potential}
Let $(M,I_0)$ be a compact LCK manifold with potential.
In \ref{_Vaisman_limit_of_LCK_pot_Theorem_}, we construct
a family of diffeomorphisms $\nu_i\in \Diff(M)$, $i=1,2,3,4, ...$
such that $\lim_i \nu_i(I_0)=I$, where $I$
is a Vaisman-type complex structure on $M$.
Using the semicontinuity of the Hodge numbers
(\cite[Theorem 2.3]{_Bell_Narasimhan_}), we obtain
that $H^{p,q}(M, I_0)\leq H^{p,q}(M,I)$.
However, the Vaisman manifolds \index[terms]{manifold!Vaisman}
satisfy $\sum_{p+q=r}\dim H^{p,q}(M,I)= \dim H^r(M)$
(\ref{_E_1_degenerates_for_Vaisman_Remark_}),
hence $\sum_{p+q=r}\dim H^{p,q}(M,I)\leq \dim H^r(M)$;
this can happen only if $\sum_{p+q=r}\dim H^{p,q}(M,I)= \dim H^r(M)$,
because the Fr\"olicher spectral sequence\index[terms]{spectral sequence!Fr\"olicher} converges to the
de Rham cohomology.\index[terms]{cohomology!de Rham} This implies, in particular,
that the Fr\"olicher spectral sequence degenerates
for LCK manifolds with potential, and the
Hodge numbers are constant in smooth families
of LCK manifolds with potential\index[terms]{manifold!LCK!with potential} (\ref{_LCK_pot_degenerate_Corollary_},
\ref{_Hodge_constant_LCK+pot_Theorem_}).

\section[Weights of a torus action on the de Rham algebra]{Weights of a torus action on the de Rham\\ algebra}

\remark
It is well known (see Exercise \ref{_Fourier_for_torus_on_de_Rham_Exercise_})
that for all 1-dimensional representations $\rho$
of $T^n=\U(1)^n$ there are numbers $p_1, ..., p_n$ such that
 \[ \rho(t_1, ..., t_n)(x)= \exp\left(2\pi\1\sum_{i=1}^np_it_i\right) x.\]
Moreover, any finite-dimensional representation is a direct
sum of one-di\-men\-si\-o\-nal representations.

\hfill

\definition
Let $T^n=\U(1)^n$ be a compact torus acting on a compact manifold $M$,
and $\zeta\in \Lambda^*(M)$ a differential form.
We say that $\zeta$  is {\bf pure of weight $p_1,..., p_n$} if for all $t=(t_1, ..., t_n)\in T$,
the pullback $t^*\zeta$ is proportional to $\zeta$, with the
coefficient of proportionality written as
$t^*(\zeta)= \exp(2\pi\1\sum_{i=1}^np_it_i) \zeta$.

\hfill

\definition
Let $V$ be a Hermitian space (possibly infinite-dimensional)
en\-dowed with a unitary action of $T^n$,
and $V_\alpha\subset V$ weight $\alpha$ representations, $\alpha\in \Z^n$.
The direct sum $\bigoplus_{\alpha\in \Z^n} V_\alpha$ is called {\bf 
the weight decomposition} for $V$ if it is dense in $V$.\index[terms]{weight decomposition}

\hfill

\claim\label{_Fourier_in_L^2_Claim_}
Let $W$ be a Hermitian vector space. Then
 the Fourier series provide the weight decomposition
on the space of maps $L^2(T^n, W)$.\index[terms]{Fourier series}

\proof
Exercise \ref{_maps_from_torus_Exercise_}.
\endproof

\hfill

\claim\label{_weight_in_subspaces_Claim_}
Let $W$ be a Hermitian vector space equipped with a
unitary $T^n$-action. Assume that $W$ admits a weight
decomposition. Then any closed $T^n$-invariant subspace $W_1\subset W$
also admits a weight decomposition.

\hfill

\proof
Take the weight decomposition of $x\in W_1$
within $W$ and project orthogonally to $W_1$.
Since the orthogonal projection commutes with 
the $T^n$-action, it will give a weight
decomposition in $W_1$.
\endproof

\hfill

\claim
Let $W$ be a Hermitian representation of $T^n$.
Then its completion $\hat W$ admits a weight decomposition
$V = \widehat{\bigoplus_{\alpha\in \Z^n} W_\alpha}$.

\hfill

\proof
We realize the completion $\hat W$ as the space of constant maps, 
 $\hat W \subset L^2(T^n, W)$, clearly $T^n$-invariant.
Then we use \ref{_Fourier_in_L^2_Claim_} to
obtain the weight decomposition of $L^2(T^n, W)$.
By \ref{_weight_in_subspaces_Claim_}, this gives also
a weight decomposition in $\hat W$.
\endproof

\hfill

\remark
Let $T^n$ act on a Riemannian manifold $M$ by isometries, and
$\alpha= \sum \alpha_{p_1, ..., p_k}$ be the weight decomposition of a differential form.
Each form $\alpha_{p_1, ..., p_k}$ 
is obtained by averaging
\[ 
e^{2\pi\1\sum_{i=1}^np_it_i}\alpha
= \Av_{{}_{T^n}} e^{2\pi\1\sum_{i=1}^n-p_it_i}\alpha,
\]
and hence this form is smooth. Therefore,
the space $\Lambda^*(M)$ admits the weight decomposition.

\hfill

\remark
Let $M$ be a manifold equipped with a $T^n$-action,
and \[ \Lambda^*(M)= 
\hat\bigoplus_{\alpha \in \Z^n} \Lambda^*(M)_{p_1, ..., p_k}
\]
be the weight decomposition on the differential forms.
Then the de Rham differential preserves each term
$\Lambda^*(M)_{p_1, ..., p_k}$. Indeed, $d$ commutes
with the action of the Lie algebra of $T^n$,
and $\Lambda^*(M)_{p_1, ..., p_k}$ are its
eigenspaces.

\hfill

\theorem
Let $M$ be a smooth manifold, and $T^n$ a torus
acting on $M$ by diffeomorphisms. Denote
by $\Lambda^*(M)^{T^n}$ the complex of 
$T^n$-invariant differential forms.
Then the natural embedding
$\Lambda^*(M)^{T^n}\hookrightarrow \Lambda^*(M)$
induces an isomorphism on the de Rham cohomology.\index[terms]{cohomology!de Rham}

\hfill

\pstep
Let $\alpha\in \Lambda^*(M)$ be a form
and $\alpha= \sum \alpha_{p_1, ..., p_n}$ its weight
decomposition, with $\alpha_{p_1, ..., p_n}\in \Lambda^*_{p_1, ..., p_n}(M)$ 
a form of weight $p_1, ... , p_n$.
 Since the $T^n$-action
commutes with the de Rham differential, 
these forms are closed when $\alpha$ is closed.

\hfill

{\bf Step 2:}
Consider the standard 
generators $r_1, ..., r_n$  of the Lie algebra of the torus $T^n$
rescaled in such a way that 
$\Lie_{r_k}(\exp(2\pi\1\sum_{i=1}^np_it_i))=\1 p_k$,
and $i_{r_k}:\; \Lambda^i(M)\arrow \Lambda^{i-1}(M)$
the interior product.
Using the Cartan formula\index[terms]{Cartan formula} $\Lie_{r_k}= \{d, i_{r_k}\}$, we have
$p_k\alpha_{p_1, ..., p_n}= d (i_{r_k}\alpha_{p_1, ..., p_n})$
whenever $\alpha_{p_1, ..., p_n}$ is closed.
Therefore,  all terms in the weight 
decomposition $\alpha= \sum \alpha_{p_1, ..., p_n}$ 
are exact except $\alpha_{0,0, ...,0}$.

\hfill

{\bf Step 3:}
In the direct sum decomposition of the de Rham complex
\[
\Lambda^*(M)= \Lambda^*(M)^{T^n}\oplus 
\widehat{\bigoplus\limits_{(p_1,  ..., p_k)\neq (0,0,..., 0)}}\Lambda^*_{p_1, ..., p_k}(M),
\]
the second component has trivial cohomology, because
the operator $\Lie_{r_k}$ is invertible on 
$\bigoplus_{p_k\neq 0}\Lambda^*_{p_1, ..., p_n}(M)$.
\endproof

\section{Dolbeault cohomology on  manifolds  with a group action}
\index[terms]{cohomology!Dolbeault} 

In this section, we reproduce a result by N. \index[persons]{Klemyatin, N.} Klemyatin \cite{_Klemyatin_}.
For K\"ahler manifolds, a connected group of automorphisms
acts trivially on the Dolbeault cohomology. On a general compact complex
manifold, this is false (\cite{_Akhiezer:group_}). Klemyatin found that
a complex Lie group acting on a Hermitian compact manifold by isometries
acts trivially on its Dolbeault cohomology.

\hfill

We start with the following curious exact sequence \index[terms]{cohomology!Bott--Chern}
relating the Bott--Chern and Dolbeault cohomology of a complex manifold.\index[terms]{cohomology!Dolbeault}

\hfill

\lemma\label{_Dolbeault_via_BC_Lemma_}
Let $M$ be a complex manifold. Consider 
the tautological map $\tau:\;H^{p,q}_{BC}(M)\stackrel\tau\arrow
H^{p,q}_{\bar\6}(M)$  from the Bott--Chern cohomology to the
Dolbeault cohomology \index[terms]{cohomology!Bott--Chern}\index[terms]{cohomology!Dolbeault}taking 
the Dolbeault class of a $\6$-closed and $\bar\6$-closed form.
Then the following sequence is exact:
\begin{equation}\label{_Dolbeault_via_BC_Equation_}
H^{p-1,q}_{BC}(M)\stackrel\tau\arrow
H^{p-1,q}_{\bar\6}(M) \stackrel{\6}{\arrow} H^{p,q}_{BC}(M).
\end{equation}
\proof
We take a $\bar\6$-closed $(p-1,q)$-form
$\alpha$, and assume that $\6\alpha$ is Bott--Chern exact,
$\6\alpha= \6\bar\6\beta$. Then $\alpha - \bar\6\beta$
is $\6$-closed and $\bar\6$-closed, and Dolbeault
cohomologous to $\alpha$, which gives 
$[\alpha] = [\tau(\alpha - \bar\6\beta)]$.
\endproof

\hfill

Now we prove the 
Dolbeault analogue of the Cartan formula.

\hfill

\lemma\label{_d^c_Lie_Lemma_}
Let $X$ be a holomorphic vector field, and $X^c=I(X)$. Then $\{d^c, i_X\}= - \Lie_{X^c}$.

\hfill

\proof Using $\{IdI^{-1}, i_X\}= I\{d,  I^{-1}i_XI\}I^{-1}$,
we obtain \[ \{d^c, i_X\}= - I \{d, i_{X^c}\}I^{-1}= I\Lie_{X^c}I^{-1}.\]
However, $X^c$ is holomorphic, and hence  $I\Lie_{X^c}I^{-1}=
\Lie_{X^c}$.
\endproof

\hfill

\proposition\label{_Dolbeault_commu_Proposition_}
Let $X$ be a holomorphic vector field, and $X^c=I(X)$.
Then $\{\bar\6, i_X\}= \frac 1 2 (\Lie_X - \1 \Lie_{X^c})$.

\hfill

\proof
$\bar\6= \frac 1 2 (d +\1 d^c)$, and hence 
\[
\{\bar\6, i_X\}= \tfrac 1 2 (\Lie_X + \1 \{d^c, i_X\})= \tfrac 1 2  
(\Lie_X - \1 \Lie_{X^c}). \ \ \ \endproof
\]

\remark \label{_Dolbeault_weight_commutes_Remark_}
Let $M$ be a complex manifold
equipped with a holomorphic action of the torus 
$T^n$. Then the action of $T^n$ commutes with $d$ and $d^c$.
Therefore, the operators $d, d^c$ preserve 
the eigenspaces of the corresponding Lie algebra.
These eigenspaces are components of the weight decomposition.
This implies that the Dolbeault differential
$\bar\6$ preserves the weight decomposition.\index[terms]{weight decomposition}

\hfill

\remark \label{_group_generated_X_X^c_Remark_}
Let $X$ be a holomorphic Killing vector field\index[terms]{vector field!Killing}\index[terms]{vector field!holomorphic}
on a compact Hermitian manifold $M$, and $X^c:=I(X)$.
Assume that $X^c$ is also Killing. Then the
flows generated by $X$ and $X^c$ act on $M$
by isometries, and hence  the closure of these 
flows is a compact Lie group (\cite{_Myers_Steenrod_,_Kobayashi_Transformations_}).
It is commutative, because it contains a dense subgroup
$\langle e^{tX}, e^{tX^c}\rangle$ that is  commutative 
(\ref{_holo_field_commute_Remark_}).

\hfill

\definition
Let $X$ be a holomorphic Killing vector field\index[terms]{vector field!Killing}\index[terms]{vector field!holomorphic}
on a compact Hermitian manifold $M$, and $X^c:=I(X)$.
We say that a differential form $\alpha$ {\bf has weight $(a, b)$
with respect to $X$ and $X^c$}, if $\Lie_X\alpha = a\alpha$
and $\Lie_{X^c} \alpha = b\alpha$.

\hfill

\claim\label{_weight_deco_holo_vector_field_Claim_}
Let $X$ be a holomorphic Killing vector field
on a compact Hermitian manifold $M$. Then 
each differential form $\alpha$ can be obtained
as a sum $\alpha = \sum_{i=0}^\infty \alpha_{(a_i,b_i)}$,
where $\alpha_{(a_i,b_i)}$ are forms of weight
$(a_i,b_i)$ with respect to $X, X^c$, and
$\{(a_i, b_i)\}$ is a sequence in $(\1 \R)^2$.

\hfill

\proof
Let $T^n$ be the closure of the group generated
by $e^{tX}, e^{tX^c}$ (\ref{_group_generated_X_X^c_Remark_}),
and $\alpha= \sum_p \alpha_p$ the corresponding weight
decomposition, with $p\in \Z^n$. Since $T^n$
multiplies $\alpha_p$ by a unitary complex number,
we have $\Lie_X\alpha_p = a\alpha_p$
and $\Lie_{X^c} \alpha_p = b\alpha_p$,
for some $a, b\in \1 \R$ depending on $p$.
\endproof

\hfill

Now we can prove \index[persons]{Klemyatin, N.} Klemyatin's theorem.\index[terms]{theorem!Klemyatin}

\hfill

\theorem\label{_G_acts_on_Dolbeault_trivially_Theorem_}
Let $X$ be a holomorphic vector field on 
a compact complex manifold $M$ equipped with a Riemannian structure
(not necessarily related to the complex structure).
Assume that $X$ is Killing\index[terms]{vector field!Killing} and $I(X)$ is also Killing.
Then the diffeomorphism flow generated by $X$ acts
on the Bott--Chern and Dolbeault cohomology\index[terms]{cohomology!Bott--Chern}\index[terms]{cohomology!Dolbeault} of $X$ trivially.

\hfill

\pstep
We start by proving that $e^{tX}$ acts trivially on
the Bott--Chern cohomology.
Let $\alpha$ be a $\bar\6$- and $\6$-closed form, and
$\alpha = \sum\alpha_{(a_i, b_i)}$ be its weight decomposition,
\ref{_weight_deco_holo_vector_field_Claim_}. By 
\ref{_Dolbeault_weight_commutes_Remark_}, the forms $\alpha_{(a_i, b_i)}$
are $\bar\6$- and $\6$-closed. To prove that $e^{tX}$
acts on $H^{p,q}_{BC}(M)$ trivially, it remains to show that 
$\alpha_{(a_i, b_i)}$ is $\6\bar\6$-exact when $a_i\neq 0$ or $b_i\neq 0$.

\hfill

{\bf Step 2:} Let $\alpha$ be a $\bar\6$- and $\6$-closed form of
weight $(a,b)$ with respect to $X$ and $X^c$.
\ref{_Dolbeault_commu_Proposition_} implies
that $\{\bar\6, i_X\}= \frac 1 2  (\Lie_X - \1 \Lie_{X^c})$ and 
$\{\6, i_X\}= \frac 1 2 (\Lie_X + \1 \Lie_{X^c})$.
Therefore,
\[
\{\bar\6, i_X\}(\alpha) = (a - \1 b)\alpha, \ \ \{\6, i_X\}(\alpha) = (a+ \1 b)\alpha,
\]
hence
\begin{equation*}
\begin{split}
\{\6\bar\6, i_X\}(\alpha)&=
(a^2 + b^2) \alpha-\6 (a - \1 b)\alpha-(a+ \1 b) \bar\6\alpha+ i_X \6\bar\6 \alpha\\
&=(a^2 + b^2) \alpha.
\end{split}
\end{equation*}
(the last three terms in the middle vanish, because $\alpha$ is $\bar\6$- and $\6$-closed).
Unless $a^2 + b^2=0$, this gives $(a^2 + b^2)^{-2}\6\bar\6 i_X\alpha=  \alpha$,
which implies that $\alpha$ is $\6\bar\6$-exact.
The numbers $a$ and $b$ are imaginary (\ref{_weight_deco_holo_vector_field_Claim_}), and hence 
$a^2 + b^2=0$ implies $a=b=0$. Therefore, the terms of non-zero weight do not
contribute to the Bott--Chern cohomology, and the action of $e^{tX}$ on 
$H^{p,q}_{BC}(M)$ is trivial.

\hfill

{\bf Step 3:} It remains to prove that the action of $e^{tX}$ on 
the Dolbeault cohomology\index[terms]{cohomology!Dolbeault} $H^{p,q}_{\bar\6}(M)$ is trivial.
For this purpose, \index[persons]{Klemyatin, N.} Klemyatin uses an exact sequence
developed by \index[persons]{Angella, D.} Angella and \index[persons]{Tomassini, A.} Tomassini in \cite{_Angella_Tomassini_}.
We use the exact sequence \eqref{_Dolbeault_via_BC_Equation_}
instead:
\[
H^{p-1,q}_{BC}(M)\stackrel\tau\arrow
H^{p-1,q}_{\bar\6}(M) \stackrel{\6}{\arrow} H^{p,q}_{BC}(M)
\]
Since $e^{tX}$  acts trivially on the leftmost and the rightmost
terms it also acts trivially on $H^{p-1,q}_{\bar\6}(M)$.
\endproof

\section{Dolbeault cohomology of Vaisman manifolds}\index[terms]{manifold!Vaisman}
\index[terms]{cohomology!Dolbeault}

In this section, we prove the harmonic decomposition theorem
for the Dolbeault cohomology of Vaisman manifolds.  It was originally
obtained by K. \index[persons]{Tsukada, K.} Tsukada, \cite{tsuk}, but our proof is different.
First, we compute the cohomology groups, following 
\index[persons]{Klemyatin, N.}
N. Klemyatin \cite{_Klemyatin_}. Then we compare
the Dolbeault cohomology to basic harmonic
forms, showing that they are isomorphic when
appropriate, by dimensional reasons. This
is the same scheme as used for our proof of 
the harmonic decomposition of the de Rham cohomology\index[terms]{cohomology!de Rham}
of Vaisman manifolds, originally due
to \index[persons]{Kashiwada, T.} Kashiwada (\ref{_Vaisman_harmonic_forms_Theorem_}).\index[terms]{manifold!Vaisman}

\subsection[Basic and Dolbeault cohomologies of Vaisman ma\-ni\-folds]{Basic and Dolbeault cohomologies of Vaisman\\ ma\-ni\-folds}\index[terms]{manifold!Vaisman}\index[terms]{cohomology!Dolbeault}

In this section, we compute the Dolbeault cohomology of Vaisman manifolds.
The first step reduces the Dolbeault complex to a complex constructed from the
basic forms. Recall that the Lee field\index[terms]{Lee field} $X$ of a Vaisman manifold $M$,
and its anti-Lee field $X^c:= I(X)$ are holomorphic and Killing \index[terms]{Lee field!anti-}\index[terms]{vector field!Killing}\index[terms]{vector field!holomorphic}
(\ref{_canon_foli_totally_geodesic_Remark_}, \ref{prop_lee}),
hence \ref{_G_acts_on_Dolbeault_trivially_Theorem_} can be applied
to a Vaisman manifold.\index[terms]{manifold!Vaisman}

\hfill

\claim\label{_Dolbeault_quasiiso_basic_Theorem_}
Let $(M, \omega, \theta)$ be a Vaisman manifold, $X$ its Lee field,\index[terms]{Lee field} 
$X^c:= I(X)$ the anti-Lee field,\index[terms]{Lee field!anti-} and $\Sigma=\langle X, X^c\rangle$
the canonical foliation.\index[terms]{foliation!canonical} Denote by $G\subset \Aut(M)$ the group generated
by the Lee and the anti-Lee flow.\index[terms]{Lee flow}\index[terms]{Lee flow!anti-} Let $\Lambda^*(M)^G$ be the space of $G$-invariant differential
forms, and $\Lambda^*_\Sigma(M)$ the space of $\Sigma$-basic differential forms. Then 
\begin{multline}\label{_G_inva_decomposition_Equation_}
\Lambda^*(M)^G = \Lambda^*_\Sigma(M) \oplus \bigg(\Lambda^*_\Sigma(M)\wedge \theta^{1,0}  \bigg)
\oplus  \\
\oplus
\bigg(\Lambda^*_\Sigma(M) \wedge \theta^{0,1}\bigg) 
\oplus \bigg(\Lambda^*_\Sigma(M) \wedge \theta^{1,0}\wedge \theta^{0,1}\bigg).
\end{multline}
\proof Clearly, basic forms are $G$-invariant, and the forms
$\theta^{1,0}$  and $\theta^{0,1}$ are also $G$-invariant. Therefore,
all forms on the right-hand side of \eqref{_G_inva_decomposition_Equation_}
are $G$-invariant. To obtain the converse inclusion, we choose a local
$G$-invariant decomposition $M= M_0 \times U$, where $U\subset \C$
is the disk, and the fibres of the projection $M\arrow M_0$ are
the leaves of $\Sigma$. We assume that $G$ acts on $U$ (locally),
and its action on $M_0$ is trivial.
Choose a $G$-invariant basis $\alpha, \beta\in \Lambda^1(U)$
on $U$. Then the space of $G$-invariant forms on $M_0 \times U$
is multiplicatively generated by $\Lambda^*(M_0)$ and
$\alpha, \beta$, which gives
\begin{multline*}
\Lambda^*(M_0\times U)^G=\Lambda^*(M_0) \oplus \bigg(\Lambda^*(M_0)\wedge \alpha \bigg)
\oplus\\ \oplus  \bigg(\Lambda^*(M_0) \wedge \beta\bigg)
\oplus \bigg(\Lambda^*(M_0) \wedge \alpha\wedge \beta\bigg).
\end{multline*}
Then \eqref{_G_inva_decomposition_Equation_}
follows by local patching.
\endproof

\hfill

Since the Dolbeault classes on any Vaisman manifold
are represented by $G$-invariant forms\index[terms]{manifold!Vaisman} (\ref{_G_acts_on_Dolbeault_trivially_Theorem_}), 
this immediately brings the following corollary.

\hfill

\corollary
Let $(M, \omega, \theta)$ be a compact Vaisman manifold,  $X$ its Lee field,\index[terms]{Lee field} 
$X^c:= I(X)$ the anti-Lee field,\index[terms]{Lee field!anti-} and $\Sigma=\langle X, X^c\rangle$
the canonical foliation.\index[terms]{foliation!canonical}
Then any $\bar\6$-closed form is $\bar\6$-cohomologous to a form
\[
\alpha =  \alpha_1 + \alpha_2\wedge \theta^{1,0}  
+\alpha_3 \wedge \theta^{0,1} 
\oplus  \alpha_4\wedge \theta^{1,0}\wedge \theta^{0,1},
\]
where all $\alpha_i$ are basic.
\endproof

\hfill

The form $\theta^{0,1}$ is $\bar\6$-closed, because  $\bar\6\theta^{0,1}$ is equal
to the $(2,0)$-part of $d\theta$, that vanishes.
Therefore, the operator $\bar \6$ acts on 
the four terms of \eqref{_G_inva_decomposition_Equation_}
\[
 \Lambda^*_\Sigma(M) \oplus \Lambda^*_\Sigma(M)\wedge \theta^{1,0}  
\oplus \Lambda^*_\Sigma(M) \wedge \theta^{0,1} 
\oplus \Lambda^*_\Sigma(M) \wedge \theta^{1,0}\wedge \theta^{0,1},
\]
preserving the terms $\Lambda^*_\Sigma(M) \wedge \theta^{0,1}$ and
$\Lambda^*_\Sigma(M)$ 
and taking $\alpha_2\wedge \theta^{1,0} \in\Lambda^*_\Sigma(M) \wedge \theta^{1,0}$
to 
\[ \bar\6\alpha_2  \wedge \theta^{1,0} + (-1)^{\tilde\alpha_2}\alpha_2 \wedge \bar\6(\theta^{1,0})=
\bar\6\alpha_2  \wedge \theta^{1,0} -(-1)^{\tilde\alpha_2}\1\alpha_2 \wedge \omega_0.
\]
because $\theta^{1,0}= \theta - \1 \theta^c$, $\bar\6=\frac{d+ \1 d^c}2$,
and $d\theta^c= -\omega_0$, $d^c\theta= \omega_0$, which gives 
\begin{equation}\label{_bar_6_theta_1_0_Equation_}
\bar\6(\theta^{1,0})=\1\omega_0
\end{equation}
Similarly, we have
\[
\bar\6(\alpha_4  \wedge \theta^{1,0}\wedge \theta^{0,1})=
\bar\6\alpha_4  \wedge \theta^{1,0}\wedge \theta^{0,1} +(-1)^{\tilde\alpha_4}\1\alpha_4 \wedge \omega_0\wedge \theta^{0,1}.
\]
\proposition
Let $(M, \omega, \theta)$ be a compact Vaisman manifold,  $X$ its Lee field,\index[terms]{Lee field} \index[terms]{manifold!Vaisman}
$X^c:= I(X)$ the anti-Lee field\index[terms]{Lee field!anti-}, and $\Sigma=\langle X, X^c\rangle$
the canonical foliation.\index[terms]{foliation!canonical} Denote by $G\subset \Aut(M)$ the group generated
by the Lee and the anti-Lee flow.\index[terms]{Lee flow}\index[terms]{Lee flow!anti-} Let $\Lambda^*(M)^G$ be the space of $G$-invariant differential
forms, and $\Lambda^*_\Sigma(M)$ the space of $\Sigma$-basic differential forms.
Then the Dolbeault complex $(\Lambda^*(M), \bar\6)$ is quasi-isomorphic
to the cone $C(f)$, where $f$ is the morphism of complexes
\begin{equation}\label{_f_acts_on_complex_Definition_}
\Lambda^*_\Sigma(M) \wedge \theta^{0,1}\oplus 
\Lambda^*_\Sigma(M) \wedge \theta^{1,0}\wedge \theta^{0,1}\stackrel f
\arrow \Lambda^*_\Sigma(M)   
\oplus \Lambda^*_\Sigma(M) \wedge \theta^{0,1} 
\end{equation}
taking  $\alpha_2\wedge \theta^{1,0}$
to  $-(-1)^{\tilde\alpha_2}\1\alpha_2 \wedge \omega_0$
and $\alpha_4  \wedge \theta^{1,0}\wedge \theta^{0,1}$
to $-(-1)^{\tilde\alpha_4}\1\alpha_4 \wedge \omega_0\wedge \theta^{0,1}.$

\hfill

\proof
We use the same argument as in \ref{_Sasakian_vecr_inva_cone_Proposition_}.
By \ref{_Dolbeault_quasiiso_basic_Theorem_}, the Dolbeault complex is equivalent
to $(\Lambda^*(M)^G, \bar\6)$, that is  split onto 4 terms according to
\eqref{_G_inva_decomposition_Equation_}. 
Now, the Dolbeault differential $\bar\6$ preserves
2 of these terms and acts on the second and the fourth as $\bar\6+f$
as shown above, and hence  the complex $(\Lambda^*(M)^G, \bar\6)$ 
is equal to the cone of $f$.\index[terms]{cone!of a morphism}
\endproof

\hfill

\lemma\label{_A_B_decomposition_Lemma_}
Let $(M, \omega,\theta)$ be a compact Vaisman manifold, $\dim_\R M=2n$.\index[terms]{manifold!Vaisman}
Denote by $L_{\theta^{0,1}}:\; H^{p,q}_{\bar\6}(M)\arrow H^{p,q+1}_{\bar\6}(M)$ the multiplication by
$\theta^{0,1}$. Let $A^{p,q}\subset H^{p,q}_{\bar\6}(M)$
be the part of the cohomology associated with $\Lambda^*_\Sigma(M)   
\oplus \Lambda^*_\Sigma(M) \wedge \theta^{1,0}$, and
$B^{p,q}$ the part associated with
$\Lambda^*_\Sigma(M) \wedge \theta^{0,1}  
\oplus \Lambda^*_\Sigma(M) \wedge \theta^{1,0}\wedge \theta^{0,1}$
in the direct sum decomposition \eqref{_G_inva_decomposition_Equation_}. 
Then the Dolbeault cohomology\index[terms]{cohomology!Dolbeault} of $M$
satisfies $H^{p,q}_{\bar\6}= A^{p,q}\oplus B^{p,q}$,
and $L_{\theta^{0,1}}$ induces an isomorphism 
between $(A^{p,q-1}, \bar\6)$ and $(B^{p,q},\bar\6)$.

\hfill

\proof
By  \ref{_G_acts_on_Dolbeault_trivially_Theorem_},
the Dolbeault cohomology of $M$ is isomorphic to the complex $(\Lambda^*(M)^G,\bar\6)$
of $G$-invariant forms.
In the complex \eqref{_G_inva_decomposition_Equation_},
the terms $A^*$ and $B^*$ are preserved by the Dolbeault differential $\bar\6$,
and $L_{\theta^{0,1}}$ commutes with $\bar\6$ and induces an isomorphism
of complexes $\bigg[\Lambda^*_\Sigma(M)   
\oplus \Lambda^*_\Sigma(M) \wedge \theta^{1,0}, \bar\6\bigg]$, and 
$\bigg[\Lambda^*_\Sigma(M)\wedge \theta^{0,1}   \oplus \Lambda^*_\Sigma(M) \wedge \theta^{1,0}\wedge \theta^{0,1}, \bar\6\bigg]$.
\endproof

\hfill

We are in position to prove the Dolbeault analogue
of \ref{_Vaisman_coho_Theorem_}.

\hfill

\theorem\label{_Vaisman_Dolbeault_coho_Theorem_}
 Let $(M, \omega,\theta)$ be a compact Vaisman manifold, $\dim_\R M=2n$.\index[terms]{manifold!Vaisman}
Denote by $H^*_\kah(M)$ the basic cohomology\index[terms]{cohomology!basic} associated with the
canonical foliation\index[terms]{foliation!canonical} $\Sigma$, and let $A^{p,q}$ be the group defined
in \ref{_A_B_decomposition_Lemma_}.\footnote{By 
\ref{_Transversal_Lefschetz_Theorem_}, the group $H^*_\kah(M)$
is equipped with the same structures as the cohomology of a 
compact K\"ahler manifold.} Then
\begin{equation}
\begin{aligned}\label{_coho_Dolbeault_Vaisman_folia_Equation_}
H^{p,q}_{\bar\6}(M)& \cong H^*(A^{p,q},\bar\6) \oplus  H^*(A^{p, q-1}, \bar\6)\wedge [\theta^{0,1}], \text{\ \ and}\\[.1in]
H^*(A^{p,q}, \bar\6)&=
\begin{cases} 
\displaystyle L_{\theta^{1,0}}\left(\ker L_{\omega_0} \restrict { H^{p,q}_\kah(M)}\right), & \text{for} \ \  p+q> n-1,\\[.2in]
\displaystyle\frac{H^{p,q}_\kah(M)}{\im L_{\omega_0}}, & \text{for}\ \   p+q\leq n-1
\end{cases} 
\end{aligned}
\end{equation}
where $L_{\omega_0}$ is the operator in the 
Lefschetz $\gsl(2)$-triple $\langle L_{\omega_0}, \Lambda_{\omega_0}, H_{\omega_0}\rangle$
 acting on $H^*_\kah(M)$.

\hfill

\proof
The statement of this theorem closely mimics
\ref{_Vaisman_coho_Theorem_}, and the proof is
entirely the same. First, we note that
the Dolbeault complex of $M$ is
quasi-isomorphic to the complex $(\Lambda^*(M)^G,\bar\6)$
 of $G$-invariant forms (\ref{_G_acts_on_Dolbeault_trivially_Theorem_}).
Then we decompose $(\Lambda^*(M)^G,\bar\6)$ onto a sum of
$(A^*, \bar\6)$ and $(B^*, \bar\6)$ as in
\ref{_A_B_decomposition_Lemma_}; this
gives the first equation of \eqref{_coho_Dolbeault_Vaisman_folia_Equation_}.
The second equation follows because the complex
$(A^*, \bar\6)$ is obtained as the
cone of the morphism
\[ f:\; (\Lambda^*_\kah(M) \wedge \theta^{1,0}, \bar\6)\arrow (\Lambda^*_\kah(M), \bar\6)\]
giving the long exact sequence
\begin{multline} ... \arrow H^{p-1,q-1}_{\kah}(M)\stackrel {L_{\omega_0}} \arrow H^{p,q}_{\kah}(M) \arrow A^{p,q} \\
\arrow H^{p,q}_{\kah}(M) \stackrel {L_{\omega_0}} \arrow H^{p+1,q+1}_{\kah}(M)\arrow ...,
\end{multline}
which gives \eqref{_coho_Dolbeault_Vaisman_folia_Equation_}
in the same way as in \ref{_Vaisman_coho_Theorem_}, because
$L_{\omega_0}$ is injective on basic $k$-cohomology with $k< n-1$ and surjective
on basic $k$-cohomology with $k\geq n-1$, as follows from the Lefschetz $\goth{sl}(2)$-action
on the basic cohomology\index[terms]{cohomology!basic} (\ref{_Transversal_Lefschetz_Theorem_}).
\endproof

\hfill

\remark\label{_E_1_degenerates_for_Vaisman_Remark_}
From \ref{_Vaisman_Dolbeault_coho_Theorem_}
and \ref{_Vaisman_Dolbeault_coho_Theorem_} it is apparent that
the Dolbeault cohomology \index[terms]{cohomology!Dolbeault}of a Vaisman manifold
has the same dimension as the de Rham cohomology.\index[terms]{cohomology!de Rham}
Therefore, the Fr\"olicher spectral sequence for\index[terms]{spectral sequence!Fr\"olicher}
a Vaisman manifold degenerates on the page $E_1^{*,*}$.\index[terms]{manifold!Vaisman}

\hfill

From the degeneracy of the Fr\"olicher spectral sequence
and \ref{_Vaisman_Dolbeault_coho_Theorem_}, we can deduce a
result which answers (among other things) the question asked in 
\cite{_GNT:semicontinuous_}.

\hfill

\question
Are the basic Hodge numbers\index[terms]{Hodge numbers} of Sasakian manifolds invariant
under general smooth deformations of the Sasakian structure?

\hfill

Since the paper \cite{_GNT:semicontinuous_} appeared, this
question was answered in affirmative by P. \index[persons]{Ra\'zny, P.} Ra\'zny in \cite{_Razny:Hodge_numbers_}.
We give another proof of Ra\'zny's result below.

\hfill

Since the basic Hodge numbers of a Sasakian manifold $S$
are the same as the numbers $\dim H^{p,q}_\kah(M)$ of the Vaisman manifold $M$\index[terms]{manifold!Vaisman}
associated with $S$, and any Sasakian deformation of $S$
extends to a smooth deformation of the Vaisman structure on 
$M$, the affirmative answer to this question is implied
by the following general result.

\hfill

\theorem\label{_basic_constant_Theorem_}
Let $(M,I_t, \omega_t, \theta_t)$ be a smooth deformation
of a Vaisman structure, depending on 
$t\in [0,1]$. Then the Hodge numbers $\dim H^{p,q}(M, I_t)$ and the basic Hodge numbers
$\dim H^{p,q}_\kah(M, I_t)$ are constant in $t$, for all $p,q$.

\hfill

\proof
By \cite[Theorem 2.3]{_Bell_Narasimhan_}, 
the Hodge numbers $\dim H^{p,q}(M, I_t)$ are semicontinuous in smooth families.
Then, from \ref{_Vaisman_Dolbeault_coho_Theorem_}, it follows that
the numbers $\dim H^{p,q}_\kah(M, I_t)$ are semicontinuous
(this also follows from \cite{_GNT:semicontinuous_}),
that is, they can only jump in closed subsets of the 
parameter set. On the other hand, the sum
$\sum_{p,q}\dim H^{p,q}(M, I_t)= \sum_r\dim H^r(M)$ is constant
by \ref{_Vaisman_Dolbeault_coho_Theorem_}, and hence 
the jumping cannot actually occur. Therefore, $\dim H^{p,q}_\kah(M, I_t)$ is constant in $t$.
\endproof

\subsection[Harmonic decomposition for the Dolbeault cohomology]{Harmonic decomposition for the Dolbeault\\ cohomology}\index[terms]{cohomology!Dolbeault}

Using the same logic as in \ref{_Vaisman_harmonic_forms_Theorem_}, 
we prove \index[persons]{Tsukada, K.} Tsukada's theorem about the
harmonic decomposition of the Dolbeault cohomology.
Recall that a cohomology class
$[\eta]\in H^{k}_{\kah}(M)$
is called {\bf primitive}\index[terms]{form!primitive}
if $\Lambda_{\omega_0} \eta=0$,
for $k \leq \dim_\C M -1$,
and {\bf coprimitive} if
$L_{\omega_0} \eta=0$, for
$k \geq \dim_\C M -1$, where
$L_{\omega_0}$ and $\Lambda_{\omega_0}$
are the transversal Lefschetz operators
acting on basic forms and on basic cohomology.\index[terms]{cohomology!basic}

For the next theorem, note that the Vaisman manifold
is not K\"ahler, and the Dolbeault Laplacian is not proportional\index[terms]{manifold!Vaisman}
to the de Rham Laplacian, as it is in the K\"ahler case.

\hfill

\theorem\label{_Dolbeault_harmonic_on_Veisman_} 
Let $(M, \omega, \theta)$ be a compact Vaisman manifold, 
and $\Lambda^k_{\kah}(M)$ the basic forms with respect 
to the canonical foliation.\index[terms]{foliation!canonical} 
Then
\begin{description}
\item[(i)] Let $\eta\in  \Lambda^k_{\kah}(M)$
be a primitive\index[terms]{form!primitive} basic harmonic form, $k \leq \dim_\C M -1$.
Then the forms $\eta$ and $\eta \wedge \theta^{0,1}$
are Dolbeault harmonic.
\item[(ii)] Let $\alpha\in  \Lambda^k_{\kah}(M)$
be a coprimitive basic harmonic form, $k \geq \dim_\C M -1$.
Then $\alpha\wedge \theta^{1,0}$ and  
$\alpha\wedge \theta^{1,0}\wedge \theta^{0,1}$ 
are Dolbeault  harmonic.
\item[(iii)] The space of Dolbeault harmonic forms is generated
by the Dolbeault harmonic forms listed in (i)-(ii).
\end{description}
\proof
The forms of (i)-(ii) are by construction  $\bar\6$-closed.
Indeed, a basic harmonic form is  $\bar\6$-closed,
$\theta^{0,1}$ is also $\bar\6$-closed, and
$\theta^{1,0}$ satisfies $\bar\6\theta^{1,0}= \1\omega_0$
by \eqref{_bar_6_theta_1_0_Equation_}.
Since $\alpha \wedge \omega_0=0$,
the forms $\alpha\wedge \theta^{1,0}$ and  
$\alpha\wedge \theta^{1,0}\wedge \theta^{0,1}$ 
are  $\bar\6$-closed.
Also the forms of (i) and (ii) are exchanged by
the composition of the Hodge star operator
and the complex conjugation, and hence  $\bar\6$-closedness
implies also that they are $\bar\6^*$-closed
(\ref{_bar_6^*_via_*_Claim_}). 

This gives an embedding from 
the primitive part of $H^k_\kah(M)\oplus H^k_\kah(M)\wedge \theta^{0,1}$, 
$k \leq \dim_\C M -1$ and from the coprimitive part of
$H^k_\kah(M)\wedge \theta^{1,0}\oplus H^k_\kah(M)\wedge \theta^{1,0}\wedge \theta^{1,0}$
$k \geq \dim_\C M -1$ to $H^*(M)$. From \ref{_Vaisman_Dolbeault_coho_Theorem_}
we obtain that the dimensions of these spaces are equal,
hence this map is an isomorphism. Therefore, 
all Dolbeault  harmonic forms on $M$ are obtained from
(i) and (ii).
\endproof

\section[Dolbeault cohomology of LCK manifolds with potential]{Dolbeault cohomology of LCK manifolds\\ with potential}\index[terms]{cohomology!Dolbeault}\index[terms]{manifold!LCK!with potential}

Previously we proved the degeneration of the Fr\"olicher 
spectral sequence for Vaisman manifolds\index[terms]{spectral sequence!Fr\"olicher}\index[terms]{manifold!Vaisman} (\ref{_E_1_degenerates_for_Vaisman_Remark_}).
Unlike the standard proof of the degeneration for the K\"ahler
manifolds, our proof does not use the harmonic decomposition.
However, the argument is very close to the proof
of the existence of the harmonic decomposition anyway.
It is interesting that the  Fr\"olicher 
spectral sequence degenerates even for the LCK manifolds
with potential,\index[terms]{manifold!LCK!with potential} where no harmonic decomposition
is apparent. 

We start with the following general observation.

\hfill

\proposition\label{_degenerate_Frolicher_limits_Proposition_}
Let $I\in \End(TM)$ and $\{I_j\in \End(TM), j=0, 1, 2, ...\}$ be a collection of complex structures
on a compact manifold $M$. Assume that there exists
a collection of diffeomorphisms $\nu_j\in \Diff(M)$
such that $\nu_j(I_0)=I_j$, and $\lim_j I_j = I$,
where the limit is taken in the $C^\infty$-topology.\index[terms]{topology!$C^\infty$}
Assume that the Fr\"olicher spectral sequence
on $(M,I)$ degenerates in the $E_1^{*,*}$-page, that is, 
$\sum_{p+q=r} \dim H^{p,q}_{\bar\6}(M,I)= \dim H^r(M)$ for all $r$.
Then $E_1^{*,*}$-page degenerates for $(M,I_0)$ (and hence, for all $(M,I_j)$).

\hfill

\proof
 For any $p,q$, 
the numbers $\dim H^{p,q}_{\bar\6}(M,J)$ are semicontinuous
as functions of a complex structure $J$, that is, they can only jump in the special
points (see  \cite[Theorem 2.3]{_Bell_Narasimhan_}). Then $\dim H^{p,q}_{\bar\6}(M,I)\geq\dim H^{p,q}_{\bar\6}(M,I_0)$.
On the other hand, $\sum_{p+q=r} \dim H^{p,q}_{\bar\6}(M,I)= \dim H^r(M)$
and $\dim H^{p,q}_{\bar\6}(M,I_0)\geq \dim H^r(M)$ because
the limit of the Fr\"olicher spectral sequence for $I_0$ is the de Rham
cohomology. Comparing these two inequalities, we obtain
\ref{_degenerate_Frolicher_limits_Proposition_}.
\endproof

\hfill

The relevant result on LCK manifolds with potential \index[terms]{manifold!LCK!with potential}follows immediately.

\hfill

\corollary\label{_LCK_pot_degenerate_Corollary_}
Let $(M, I_0, \omega, \theta)$ be a compact LCK manifold with potential.
Then the Fr\"olicher spectral sequence\index[terms]{spectral sequence!Fr\"olicher}
on $(M,I_0)$ degenerates in the $E_1^{*,*}$-page, that is, 
$\sum_{p+q=r} \dim H^{p,q}_{\bar\6}(M,I_0)= \dim H^r(M)$ for all $r$.

\hfill

\proof
When $M$ is a complex surface, \ref{_LCK_pot_degenerate_Corollary_}
follows from \ref{_E_2_degenerates_Corollary_}. 
Otherwise, \ref{_Vaisman_limit_of_LCK_pot_Theorem_} can be applied.
Using \ref{_Vaisman_limit_of_LCK_pot_Theorem_}, 
we find a sequence of diffeomorphisms
$\nu_j\in \Diff(M)$, such that $\lim_j I_j = I$,
where $I$ is a Vaisman-type complex structure on $M$
and $I_j= \nu_j(I_0)$. Then \ref{_LCK_pot_degenerate_Corollary_}
follows from \ref{_degenerate_Frolicher_limits_Proposition_}
and \ref{_E_1_degenerates_for_Vaisman_Remark_}.
\endproof

\hfill

We also obtain that the Hodge numbers of LCK manifold with 
potential stay constant in families, generalizing\index[terms]{Hodge numbers} \ref{_basic_constant_Theorem_}.

\hfill

\theorem \label{_Hodge_constant_LCK+pot_Theorem_}
Let $(M,I_t, \omega_t, \theta_t)$ be a  deformation
of an LCK manifold with potential,\index[terms]{manifold!LCK!with potential} smoothly depending on the parameter
$t\in [0,1]$. Then the Hodge numbers
$\dim H^{p,q}(M, I_t)$ are constant in $t$.

\proof
By \cite[Theorem 2.3]{_Bell_Narasimhan_}, 
the Hodge numbers $\dim H^{p,q}(M, I_t)$ are semicontinuous in smooth families,
that is, they can only jump in closed subsets of the 
parameter set. However, the sum
$\sum_{p,q}\dim H^{p,q}(M, I_t)= \sum_r\dim H^r(M)$ is constant
by \ref{_LCK_pot_degenerate_Corollary_}, and hence 
the jumping cannot actually occur. This implies
that $\dim H^{p,q}(M, I_t)$ is constant in $t$.
\endproof

\section{Exercises}

\subsection{Isometry groups}

\definition
\index[terms]{space!metric} {\bf A geodesic} of a metric space $M$ is a local isometry $[a, b] \arrow M$
from an interval with the standard metric.\index[terms]{space!metric}

\begin{enumerate}[label=\textbf{\thechapter.\arabic*}.,ref=\thechapter.\arabic{enumi}]
\item
Let $M$ be a metric space.
Define {\bf the length} $L_d(\gamma)$ of a continuous map
$[a, b] \stackrel \gamma \arrow M$
as the upper limit of $\sum d(\gamma(a_i), \gamma(a_{i+1}))$
taken over all partitions $[a_0, a_1], ...., [a_{n-1}, a_n]$ of
$[a, b]$, with $a=a_0, b=a_n$.\index[terms]{intrinsic/local length}
\begin{enumerate}
\item
Suppose that every two points of $M$ can be connected by
a path of finite length, and define {\bf the path metric}\index[terms]{metric!path}
$M \times M \stackrel {d_1}\arrow \R^{\geq 0}$ in such a way that $d_1(x, y)$ is equal to the infimum of $L_d(\gamma)$
for all $\gamma$ connecting $x$ to $y$. Prove that $d_1$ is a metric
and $d_1 \geq d$.
\item
Let $d_1$ be the path metric constructed above.
Prove that $L_{d_1}(\gamma)= L_d(\gamma)$ for any 
path $[a, b] \stackrel \gamma \arrow M$.
\end{enumerate}

\item
Let $(M,d)$ be a metric space.\index[terms]{space!metric} We say that 
$d$ {\bf is an intrinsic metric} if $d(x,y)$
is  equal to the infimum of $L_d(\gamma)$
for all $\gamma$ connecting $x$ to $y$.\index[terms]{metric!intrinsic}
Suppose that $(M,d)$ is complete and locally compact,
and has intrinsic metric. Prove that every two
points of $M$ can be connected with a geodesic.

{\em Hint:} Use the Arzel\`a--Ascoli theorem.\index[terms]{theorem!Arzel\`a--Ascoli}

\item
Let $(M, g)$ be a Riemannian manifold, 
and $d$ the Riemannian distance function
associated with $g$. Prove that $d$ is an
intrinsic metric.

\item
Let $(M,d)$ be a metric space,\index[terms]{space!metric} and $\{U_i\}$ an open covering.
Consider the set ${\cal S}_{\{U_i\}}$ of all metrics
$d'$ on $M$ such that $d'\restrict{U_i}\leq d\restrict{U_i}$
for every element of the cover $\{U_i\}$. Let
 $d(\{U_i\})$ be the supremum of all metrics
in the family ${\cal S}_{\{U_i\}}$.
A metric is called {\bf local}\index[terms]{metric!local}
if $d(\{U_i\})=d$ for any covering $\{U_i\}$.
\begin{enumerate}
\item Prove that any intrinsic metric is local.
\item Let $(M,d)$ be a complete metric space.
Prove that $d$ is intrinsic if $d$ is local.
\end{enumerate}

\item
We say that two geodesics $\gamma(t), \gamma'(t)$ in a metric space\index[terms]{space!metric} starting from $m_i\in M_i$
{\bf have the same direction} if $\lim_{t\to 0}\frac{d(\gamma_i(t), \gamma_i'(t))}{t}=0$.
{\bf A direction} is an equivalence class of geodesics having the same direction.
Prove that a geodesic on a Riemannian manifold $M$
is uniquely determined by its direction, and $T_{m}(M)$ is the set of directions.


\item
Let $\gamma_1:\; [0,a]\arrow M$, $\gamma_2:\; [0,b]\arrow M$
be geodesics in a Riemannian manifold $M$, $\gamma_1(0)=\gamma_2(0)=m$, and
$\bar \gamma_1, \bar\gamma_2\in T_m M$ be their directions.
Define $d(\bar \gamma_1, \bar\gamma_2):= \lim_{t\to 0} \frac{d(\gamma_1(t), \gamma_2(t))}{t}$.
Prove that this limit is well-defined, and induces the standard Euclidean
metric on the space of directions, identified with $T_mM$.

\item
Let $u:\; M_1 \arrow M_2$ be an isometry of Riemannian manifolds
(here the Riemannian manifolds are 
understood as metric spaces;\index[terms]{space!metric} the isometry is not a priori assumed to be differentiable).
\begin{enumerate}
\item 
 Prove that $u$ maps geodesics to geodesics.
\item
Let $m_1\in M_1$ and $m_2:= u(m_1)$.
Prove that $u$ induces an isometry $D:\; T_{m_1}M_1 \arrow T_{m_2}(M_2)$ 
on the space of directions.
\item
Prove that any isometry of Euclidean spaces is linear.
\item Prove that $u$ is differentiable in $m_1$, and $D$ is its differential.
\end{enumerate}

{\em Hint:} Use the previous exercise.

\item
Let $M$ be a compact metric space,\index[terms]{space!metric} and $G$ its isometry group.
Define the metric on $G$ as $d(g_1, g_2):= \sup_{m\in M} d(g_1(m), g_2(m))$.
Prove that $d$ is a metric.
Prove that $G$ is compact in the topology induced by this metric.

\item
Let $M_1, M_2$ be connected Riemannian manifolds, and
$u, u':\; M_1 \arrow M_2$ isometries that satisfy
$u(m)= u'(m)$ for some $m\in M$. Suppose that the
differential of $u$ in $m$ is equal to the differential
of $u'$ in $m$. Prove that $u=u'$.

\item\label{_maps_from_torus_Exercise_}
Let $W$ be a Hermitian vector space. Consider the
space $L^2(T^n, W)$ as a unitary $T^n$-representation. Prove that 
 the Fourier series provides the weight decomposition
on the space of maps $L^2(T^n, W)$.\index[terms]{Fourier series}

\item 
Let $A$ be an orthogonal automorphism of a Hilbert space $H$.
Assume that $H$ has $A$-invariant vectors.
Prove that $\lim_{n\to \infty} \frac 1 n \sum_{i=1}^n A^i(x)=0$ for any $x\in H$.

\item Let $(M,I,h)$ be a compact Hermitian manifold, and $X$
a holomorphic vector field that acts on $M$ by isometries.
Find an example of $(M,I, h, X)$ such that the closure of
the Lie group generated by the flows of $X, I(X)$ is not
compact.

\item Let $M$ be a compact complex manifold, and
 $G$ the connected component of its group of
holomorphic automorphisms. Prove that $G$ is closed
in the group of self-homeomorphisms of $M$ equipped
with the $C^0$-topology.\index[terms]{topology!$C^0$}

{\em Hint:} Use Montel theorem.\index[terms]{theorem!Montel}

\item
Let $M$ be a Vaisman surface that is  not of class VII,
and $G$ the connected component of its group of
holomorphic automorphisms.
Prove that $G$ is compact. 

\item Let $M$ be a linear Hopf surface, and $G$  its group of
holomorphic automorphisms. Prove that $G$ is not compact.\index[terms]{surface!Hopf!linear}

\item Let $M$ be a compact Vaisman manifold.\index[terms]{manifold!Vaisman}
Prove that the space of Dolbeault harmonic forms
does not coincide with the space of de Rham harmonic
forms.

\subsection{Aeppli and Dolbeault cohomologies of Vaisman
  ma\-ni\-folds}\index[terms]{cohomology!Dolbeault}

\item 
Define {\bf the  Aeppli cohomology} $H^{p,q}_{\Ae}(M)$ of a
compact complex $n$-di\-men\-sio\-nal manifold $M$
as $\frac{\ker dd^c}{\im d + \im d^c}$.\index[terms]{cohomology!Aeppli}
Consider the bilinear map \[H^{n-p, n-q}_{BC}(M) \times
H^{p,q}_{\Ae}(M)\arrow \C\]
taking the cohomology classes
$[\alpha], [\beta]$ represented by the forms
$\alpha$ and $\beta$ to $\int_M \alpha\wedge \beta$.
\begin{enumerate}
\item Prove that this pairing is independent on  the
choice of the representatives $\alpha$, $\beta$.
\item Prove that this pairing
is non-degenerate.
\end{enumerate}

\item
Let $M$ be a complex manifold, and
\[
H^{p-1,q}_{\Ae}(M) \stackrel \6 \arrow
H^{p,q}_{\bar\6}(M)\stackrel \tau\arrow H^{p,q}_{\Ae}(M)
\]
the sequence formed by $\6$ and the tautological
map $\tau$. Prove that it is exact.

\item
Let $M$ be a complex manifold, and
$\tau:\; H^{p,q}_{BC}\arrow H^{p+q}(M)$ the tautological map.
Prove that the following sequence
\begin{equation}\label{_BC_via_Dolbeault_and_dR_Equation_}
H^{p-1, q}_{\bar \6}(M) \oplus H^{p, q-1}_{\6}(M)\xlongrightarrow{\6 \oplus \bar\6}
H^{p,q}_{BC} \stackrel \tau \arrow H^{p+q}(M) 
\end{equation}
is exact.

\item \label{_BC_Vaisman_Exercise_}\index[terms]{manifold!Vaisman}
Let $M$ be a compact Vaisman manifold, $\dim_\C M=n$.
As usual, $H^*_\kah(M)$ denotes the basic cohomology\index[terms]{cohomology!basic} of
the canonical foliation \index[terms]{foliation!canonical}$\Sigma$, and $L_{\omega_0}$
the operator of multiplication by $\omega_0$, considered to be  a part of
the Lefschetz $\goth{sl}(2)$-triple  acting on $H^*_\kah(M)$.
\begin{enumerate}
\item
Let $G\subset \Iso(M)$ be the group generated 
by the Lee- and the anti-Lee flows.\index[terms]{Lee flow}\index[terms]{Lee flow!anti-}
Prove that all elements of $H^{p,q}_{BC}(M)$, 
can be represented by  $G$-invariant forms.

\item
Let $\alpha$ be a closed basic form.
Prove that $L_{\omega_0}^2(\alpha)$ is Bott--Chern exact.
\item Prove that any Bott--Chern class\index[terms]{class!Bott--Chern} can be represented 
as a linear combination 
\[ \alpha =  \alpha_1 + \alpha_2\wedge \theta^{1,0}  
+\alpha_3 \wedge \theta^{0,1} 
\oplus  \alpha_4\wedge \theta^{1,0}\wedge \theta^{0,1},
\]
where all $\alpha_i$ are 
basic.

\item Analyzing the 4-term sum given above,
prove that for  $p+q<n$ only the first of these 4 terms can be non-zero.
Show that all elements of $H^{p,q}_{BC}(M)$, $p+q<n$,
can be represented by  $\Sigma$-basic forms.

\item
Let $C^{p,q}\subset H^{p,q}_\kah(M)$ be the space of 
primitive basic cohomology\index[terms]{cohomology!basic} clas\-ses.
Prove that $H^{p,q}_{BC}(M)= 
C^{p,q}\oplus L_{\omega_0}(C^{p-1,q-1})$, for $p+q<n$.
\end{enumerate}

{\em Hint:} Use the exact sequence
\eqref{_BC_via_Dolbeault_and_dR_Equation_}.\footnote{See also \cite{_Istrati_Otiman:BC_}, where these groups
were computed using different methods.}

\item
Let $M$ be a Hopf $n$-manifold. Prove that
the  Aeppli cohomology group $H^{n-1, n-1}_{Ae}(M)$ is 1-dimensional
and generated by $\omega^{n-2} \wedge \theta \wedge \theta^c$,
where $(M,\omega, \theta)$ is an LCK structure with potential\index[terms]{structure!LCK!with potential}.

\item
Let $M$ be a Vaisman $n$ manifold, and $L$ the weight bundle
of any LCK structure. Prove that the curvature of 
$L$ is equal to $\omega_0:= d^c \theta$ (\ref{_Subva_Vaisman_Theorem_}), 
and it is a  semi-positive Hermitian form\index[terms]{form!Hermitian} with exactly 
$n-1$ positive eigenvalues, vanishing in the direction
of the canonical foliation.\index[terms]{foliation!canonical}

\item
Let $(M, \omega, \theta)$ be a Vaisman manifold, and $L$ its weight bundle.\index[terms]{manifold!Vaisman}
Prove that $\deg_\omega L>0$ (see Subsection
\ref{_BC_degree_Subsection_} for the
  definition of the degree).

\item
Let $(M, \omega, \theta)$ be a Vaisman manifold,
such that its K\"ahler cover\index[terms]{cover!K\"ahler} $\tilde M$ is Ricci-flat.
Prove that $\deg_\omega K_M <0$.

{\em Hint:} Prove that $K_M = L^{-\dim_\C M}$ and use the
previous exercise.

\item\label{_balanced_degree_Exercise_}
Recall that a Hermitian metric $\omega$ on a complex $n$-manifold
is called {\bf balanced} if $d(\omega^{n-1})=0$.
Prove that $\deg_\omega L=0$ for any topologically
trivial line bundle.

\item Let $M$ be a complex manifold admitting
an exact, closed, positive (1,1)-form.
Prove that $M$ cannot admit a balanced metric.

\item\label{_Vaisman_non_balanced_Exercise_}
Prove that a compact Vaisman manifold cannot admit a balanced metric.\index[terms]{manifold!Vaisman}

{\em Hint:} Use the previous exercise (see also \cite{_Angella_Otiman_}).

\item
Let $(M, \omega_0, I_0)$ be a compact LCK manifold with potential,\index[terms]{manifold!LCK!with potential}
and $(M, I_i)$ a sequence of complex manifolds isomorphic to
$(M,I_0)$ converging to a Vaisman manifold $(M, I, \omega)$
(\ref{_Vaisman_limit_of_LCK_pot_Theorem_}).
Consider an LCK metric $\omega_i$ on $(M, I_i)$
converging to  $\omega$, and let $L$ be the corresponding
weight bundle. Prove that $\deg_{\omega_i} L_i = \deg_\omega L$.

\item \label{_LCK+pot_not_balanced_Exercise_}
Prove that an LCK manifold with potential
cannot have a balanced metric.\index[terms]{metric!balanced}

{\bf Hint:} Use the previous exercise and 
Exercise \ref{_balanced_degree_Exercise_}.

\subsection[Aeppli cohomology and strongly Gauduchon
me\-trics]{Aeppli cohomology and strongly Gauduchon\\
  me\-trics}

\item\label{_strongly_Gauduchon_Exercise_}
Let $M$ be a complex manifold.
The Dolbeault differential
is mapping the \index[terms]{cohomology!Aeppli} Aeppli cohomology to the Dolbeault
cohomology,
\[ \6:\; H^{p,q}_\Ae(M) \arrow H^{p+1, q}_{\bar \6}(M).\]
A Gauduchon form\index[terms]{form!Gauduchon} $\omega$ is called {\bf strongly Gauduchon}
(\cite{_Popovici:SG_}, \cite[Chapter 4]{_Popovici:book_}) if $\6(\omega^{n-1})$ is Dolbeault exact.
\begin{enumerate}
\item Let $(M, \omega, I)$ be a compact Vaisman manifold,\index[terms]{manifold!Vaisman} and 
$\theta^{0,1}$ the (0,1)-part of the Lee form.\index[terms]{form!Lee}
Prove that $\int_M \6(\omega^{n-1})\wedge \theta^{0,1}>0$.
\item Prove that the integral 
$\int_M \6(\omega^{n-1})\wedge \theta^{0,1}$ only depends 
on the Dolbeault class of $\6(\omega^{n-1})$. Deduce that
$\omega$ is not strongly Gauduchon.\index[terms]{form!strongly Gauduchon}
\end{enumerate}

\item\label{_Gauduchon_cone_Exercise_}
Let $(M, \omega)$ be  a complex $n$-manifold with the Gauduchon
form $\omega$. The  Aeppli cohomology class of
$\omega^{n-1}$ is called  {\bf the Gauduchon class}\index[terms]{class!Gauduchon} of $(M, \omega)$.
{\bf The Gauduchon cone} (\cite{_Popovici:Gauduchon_})
is the set of all Gauduchon
classes in the Aeppli cohomology group $H^{n-1, n-1}_\Ae(M)$. Prove that the 
Gauduchon cone of $M$ is an open, convex, homothety invariant
subset in $H^{n-1, n-1}_\Ae(M)$.\index[terms]{cone!Gauduchon}

\item
Let $\eta$ be a $dd^c$-closed $(p,p)$-form 
on a complex manifold, and $Z$ a $p$-dimensi\-onal complex subvariety.
Prove that $\int_Z \eta$ depends only on the Aeppli
class of $\eta$.\index[terms]{class!Aeppli} 

\item
Let $M$ be a Vaisman manifold,\index[terms]{manifold!Vaisman} embedded to a Hopf manifold $H$, 
and $Z\subset M$ a divisor, obtained as an 
intersection $M \cap H_0$, where $H_0$
is a divisor in $H$. 
\begin{enumerate}
\item Prove that the fundamental class of $Z$
vanishes.
\item 
Let $[Z]$ denote the Bott--Chern class\index[terms]{class!Bott--Chern} of 
the integration current \index[terms]{current!of integration}of $Z$. Prove that
$[Z]$ is proportional to $\omega_0$.
\end{enumerate}

{\em Hint:} Prove first that the fundamental class
of $H_0$ vanishes. For the second assertion, 
compute the kernel of the tautological map
$H^{1,1}_{BC} \arrow H^2(M)$ and prove that it
is 1-dimensional (Exercise \ref{_BC_Vaisman_Exercise_}).

\item\label{_Gauduchon_cone_in_half-space_Exercise_}
Let $M$ be a compact complex manifold
and $\omega_0$ a positive, closed (1,1)-form. Prove that the
Gauduchon cone\index[terms]{cone!Gauduchon} of $M$ is contained in
a half-space 
$\{\eta \in H^{n-1, n-1}_\Ae(M)\ \ |\ \ \int_M \omega_0\wedge\eta>0\}$.


\item 
As usual, we denote 
the $\Sigma$-basic cohomology on a Vaisman manifold \index[terms]{manifold!Vaisman}by $H^{*}_\kah(M)$.
We call a basic cohomology\index[terms]{cohomology!basic} class $\eta \in H^{p, q}_\kah(M)$
{\bf coprimitive} if $\eta \wedge \omega_0=0$.
Denote by $C^{p,q}\subset H^{p, q}_\kah(M)$ the space
of coprimitive cohomology classes.
\begin{enumerate}
\item Prove that 
\[ H^{n-1, n-1}_\Ae(M)= C^{n-2,n-2}\wedge\theta \wedge
\theta^c \oplus \omega_0^{n-2} \wedge\theta \wedge
\theta^c.
\] 
\item
Prove that $H^{n, n-1}_{\bar\6}(M)= H_\kah^{2n-1}(M) \oplus \R\kappa$, 
where $\kappa$ is the Dolbeault cohomology\index[terms]{cohomology!Dolbeault}
class represented by $\omega_0^{n-1} \wedge \theta^{0,1}$.
\end{enumerate}

{\em Hint:} Use Exercise \ref{_BC_Vaisman_Exercise_}
and \ref{_Dolbeault_harmonic_on_Veisman_}.

\item 
Let $M$ be a Vaisman manifold.\index[terms]{manifold!Vaisman} Prove that the natural map 
$$\6:\; H^{n-1,n-1}_\Ae(M) \arrow H^{n, n-1}_{\bar\6}(M)$$
takes a class represented by 
$\alpha\wedge \theta^{1,0} \wedge \theta^{0,1}$ 
to $\alpha\wedge \omega_0 \wedge \theta^{0,1}$.
Prove that the image of this map 
is non-zero if and only if the class
$\alpha \wedge \omega_0\in H^{n-1,n-1}_\kah(M)$
is non-zero.

{\em Hint:} Use the previous exercise.

\item\label{_Vaisman_all_are_not_sG_Exercise_}
Let $M$ be a compact Vaisman manifold.
Prove that $M$ does not admit a strongly
Gauduchon metric.

{\em Hint:} Use the previous exercise together with
Exercise \ref{_Gauduchon_cone_in_half-space_Exercise_}.

\item
Let $M$ be a compact Vaisman manifold, $\dim_\C M =n$,
and $\Omega$ a holomorphic $(p,0)$-form, $p <n$.
Prove that $\Omega$ is basic.

{\em Hint:} Use \ref{_Dolbeault_harmonic_on_Veisman_}.

\item\label{_Vaisman_never_holo_symple_Exercise_}
Let $M$ be a compact Vaisman manifold, $\dim_\C M =2n, n > 1$,
and $\Omega$ a holomorphic $(2,0)$-form. Prove that
$\Omega^n=0$.

{\em Hint:} Use the previous exercise.

\end{enumerate}


\chapter{Calabi--Yau theorem for Vaisman manifolds}\index[terms]{manifold!Vaisman}

{\setlength\epigraphwidth{0.8\linewidth}
\epigraph{\it
I am fortunate to have broken out of that torture-chamber
of the Inquisition that is  academic art. I have arrived
at the surface and can take the dimension of a living
body. But I shall use the dimension from which I shall
create the new. I have released all the birds from the
eternal cage, and opened the gates to the animals in
zoological gardens. May they tear to pieces and devour the
remains of your art. And may the freed bear bathe his body
in the ice of the frozen north and not languish in the
aquarium of boiled water in the academic garden.}{\sc\scriptsize Kazimir Malevich,  ``From
Cubism and Futurism to Suprematism: The New Realism in
Painting'' (1915)}
}



\section{Introduction}


E. \index[persons]{Calabi, E.} Calabi
(\cite{_Calabi:ICM54_,_Calabi:vanishing_canonical_}) 
has noticed that a K\"ahler metric is uniquely
determined by its K\"ahler class\index[terms]{class!K\"ahler} and its volume form.
The \index[persons]{Calabi, E.} Calabi conjecture, proven by S.-T. \index[persons]{Yau, S.-T.} Yau some 23 years
later, claims that on any compact K\"ahler $n$-manifold $M$
there exists a unique K\"ahler metric $\omega$\index[terms]{conjecture!Calabi}
with a given volume form $V$ and a given K\"ahler class\index[terms]{class!K\"ahler}
$[\omega]$ if $\int_M V= \int_M[\omega]^n$.\index[terms]{conjecture!Calabi}
This statement is equivalent to the existence
and uniqueness of the solutions of the complex
Monge-Amp\`ere equation $(\omega +dd^c \phi)^n= V$.\index[terms]{Monge-Amp\`ere equation}

From this observation, \index[persons]{Calabi, E.} Calabi obtained that (conditional
on the Calabi conjecture)\index[terms]{conjecture!Calabi} every compact K\"ahler manifold
$M$ with $c_1(M)=0$ admits a unique Ricci-flat metric
in any given K\"ahler class;\index[terms]{class!K\"ahler} this Ricci-flat metric
is clearly Einstein. Finding and classifying the Einstein metrics 
(and, more generally, the ``extremal metrics'', defined\index[terms]{metric!extremal}
by \index[persons]{Calabi, E.} Calabi in \cite{_Calabi:extremal_} as a generalization of Einstein
and constant scalar curvature K\"ahler metrics) 
was one of the central subjects of the complex
algebraic geometry\index[terms]{geometry!algebraic} since the 1980-ies.

The Calabi--Yau theorem has two facets, equally important.
One is related to finding the Ricci-flat metric on manifolds
with trivial canonical bundle. The other one is not
related to canonical bundle in any way: it is 
a result which claims that a K\"ahler metric\index[terms]{theorem!Calabi--Yau}
is uniquely defined by its volume form and
the K\"ahler class.\index[terms]{class!K\"ahler} The existence of the Ricci-flat
metrics follows directly from this, more general, result.

We prove a version of the  Calabi--Yau theorem\index[terms]{theorem!Calabi--Yau} for Vaisman manifolds.\index[terms]{manifold!Vaisman}

Let $M$ be a compact smooth
manifold, and $F\subset TM$ a smooth
foliation. It is called {\bf transversally K\"ahler}
if the normal bundle $TM/F$ is equipped with\index[terms]{foliation!transversally K\"ahler}
a Hermitian structure (that is, a complex structure
and a Hermitian metric) that is  locally obtained
as the pullback of a K\"ahler structure on the leaf
space. Sasakian manifolds are prime examples of
transversally K\"ahler manifolds (the leaf space
of the Reeb foliation on a Sasakian manifolds
is K\"ahler).

A differential form on $M$ is called {\bf basic}
if it vanishes on $F$ and is locally obtained
as the pullback of a form on the leaf space.
The basic forms are preserved by de Rham differential,
and the cohomology of the basic forms is called
{\bf the basic cohomology}. \index[terms]{foliation!taut}

A foliation is {\bf taut} if the top basic cohomology\index[terms]{cohomology!basic}
is non-zero (this is not the usual definition, but a theorem of
\index[persons]{Habib, G.} Habib and Richardson, see 
Subsection \ref{_basic_taut_Subsection_}). 
For taut foliations, one has also
Poincar\'e duality\index[terms]{duality!Poincar\'e} on the basic cohomology, and the
identification between the basic cohomology and the basic
harmonic forms if a transversal Riemannian structure
is given.\index[terms]{cohomology!basic}

When $F$ is taut and transversally\index[terms]{foliation!transversally K\"ahler}\index[terms]{foliation!taut}
K\"ahler, the basic cohomology should satisfy all the nice properties
of the cohomology of the K\"ahler manifolds: the $dd^c$-lemma,
the Hodge decomposition, the Hodge structure\index[terms]{structure!Hodge}, Lefschetz
$\goth{sl}(2)$-action and so on. We prove it in the
situation when $F$ is trivialized by a group action
that preserves the transversally K\"ahler structure
(\ref{_Transversal_Lefschetz_Theorem_}); 
a similar result was proven much earlier
by A. \index[persons]{El Kacimi-Alaoui, A.} El Kacimi-Alaoui (\cite{_Kacimi_}).

In this situation, A. \index[persons]{El Kacimi-Alaoui, A.} El Kacimi-Alaoui
proves the transversal Calabi--Yau theorem\index[terms]{theorem!transversal Calabi--Yau}
showing that the transversal K\"ahler structures
are uniquely defined by the transversal volume
form and the transversal K\"ahler class.\index[terms]{class!K\"ahler!transversal} 
We give the uniqueness part of the proof
in \ref{_transversal_volume_defines_omega_0_Lemma_},
and for the existence, refer to \cite{_Kacimi_}.

We apply these results to Vaisman geometry.
Any Vaisman manifold comes equipped by\index[terms]{manifold!Vaisman} 
a transversally K\"ahler, holomorphic 
foliation $\Sigma$ (\ref{_Subva_Vaisman_Theorem_}).
The transversally K\"ahler structure of this
foliation depends on the Vaisman metric,
however, the transversally complex structure
depends only on the complex structure of
the Vaisman manifold. We construct a
correspondence between the set of
Vaisman metrics on $M$ and the
set of transversal K\"ahler structures.

The transversal Calabi--Yau theorem implies 
an important result about the Vaisman metrics
(\ref{_CY_Vaisman_Theorem_}), showing that
the Vaisman metric is defined uniquely, up to a
constant multiplier, by the volume and the Lee class
$[\theta]\in H^1(M)$. The space of possible\index[terms]{cone!Lee}
Lee classes on $M$ (its ``Lee cone'')\index[terms]{class!Lee}
is described in \ref{_Lee_cone_on_Vaisman_Theorem_}:
it is identified with a certain open half-space in $H^1(M,\R)$.
Then, similarly to the Calabi theorem parametrizing the K\"ahler\index[terms]{theorem!Calabi}
forms, the set of all Vaisman structures on $(M,I)$ is parametrized by
the cohomological data together with the volumes.

The Vaisman  Calabi--Yau theorem is deduced\index[terms]{theorem!Calabi--Yau}
directly from the transversal Ca\-la\-bi-Yau theorem,
because the transversal K\"ahler form\index[terms]{form!K\"ahler!transversal} of\index[terms]{theorem!transversal Calabi--Yau}
a Vaisman manifold,\index[terms]{manifold!Vaisman} together with its Lee
class, uniquely defines the Vaisman structure
(\ref{_Vaisman_from_omega_0_and_Lee_Lemma_}).
On the other hand, the transversal volume form
uniquely defines, and is uniquely defined,
by the volume form of a Vaisman manifold 
\ref{_Vaisman_volume_from_omega_0_Lemma_}.
This follows from a curious
observation, made by K. \index[persons]{Tsukada, K.} Tsukada
in \cite{tsuk}, who has proved
that the direction of the Lee
field of a Vaisman manifold is
determined by its complex structure.
This chapter is based on the paper \cite{ov_CY}.\index[terms]{manifold!Vaisman}

\section{The Lee field on a compact 
Vaisman manifold}\index[terms]{Lee field} 

In \cite{tsuk}, K. \index[persons]{Tsukada, K.} Tsukada has shown that the
direction of the Lee field\index[terms]{Lee field} is uniquely determined by the
complex structure of the Vaisman manifold. We give
another proof of this result here; for other proofs,
see \cite{_Madani_Moroianu_Pilca:holo_} and 
\cite{_OV:EW_}.\footnote{The proof in \cite{_OV:EW_}
is valid only for $b_1(M)=1$.}\index[terms]{manifold!Vaisman}

\hfill

\proposition\label{_Lee_field_direction_Proposition_}
Let $M$ be a compact complex manifold of Vaisman type,
and $\theta^\sharp$ the Lee field\index[terms]{Lee field} of a Vaisman structure
`$(\omega, \theta)$. Then $\theta^\sharp$ is determined by the
complex structure on $M$ uniquely up to a real multiplier.

\hfill

\pstep
Recall that $\theta^\sharp$ is 
holomorphic and Killing by\index[terms]{vector field!Killing}\index[terms]{vector field!holomorphic} \ref{_Lee_field_Killing_Corollary_}.
Denote by $\Sigma$ the canonical foliation \index[terms]{foliation!canonical}of $M$
(\ref{_Subva_Vaisman_Theorem_}).
Since $\Sigma$ has a non-degenerate section $\theta^\sharp$, it is
trivial as a holomorphic line bundle. Then $H^0(M, \Sigma)=\C$,
and the space of holomorphic vector fields tangent
to $\Sigma$ has real dimension 2. Since $\theta^\sharp$
is Killing, it acts conformally on the 
K\"ahler cover\index[terms]{cover!K\"ahler} $\tilde M$ of $M$.
However, a holomorphic conformal vector field
on a K\"ahler manifold multiplies the K\"ahler
form by a constant, see Subsection \ref{_LCK_chapter_Intro_}.
This gives a character $\sigma:\; \goth s \arrow \R$
on the Lie algebra of all holomorphic 
vector fields tangent to $\Sigma$, $X \mapsto \frac{\Lie_X\tilde \omega}{\tilde \omega}$, where $\tilde\omega$ denotes the K\"ahler form on $\tilde M$.
We claim that $\sigma$ is uniquely
determined by the cohomology class of $\theta$.
Indeed, let $f\in C^\infty\tilde M$ be a function
such that $df =\theta$, and $X\in \goth s$, 
$X= a\theta^\sharp + bI(\theta^\sharp)$. Then 
\[ \Lie_X f= \langle \theta, X\rangle = 
a  \langle \theta, \theta^\sharp\rangle = a,\]
because $\nabla \theta=0$, and $a$ is equal to $\sigma(X)$
since
$\sigma(X)\tilde \omega =\Lie_X\tilde \omega = 
a \Lie_{\theta^\sharp}\tilde \omega+ 
b \Lie_{I(\theta^\sharp)}\tilde \omega=a\tilde \omega$.

\hfill

If we change the Vaisman metric and
modify $\theta$ in its cohomology class
by replacing it with $\theta + du$, the
product $\langle \theta+du , X\rangle$
becomes $a+ \Lie_X u$; since it remains constant,
this implies that $\Lie_X u=\const$, and hence 
$\Lie_X u=0$.\footnote{On a compact smooth manifold, the equation
	$\Lie_X u=\const$ has no solutions unless $\Lie_X u=0$;
	indeed, $u$ has to reach maximum somewhere, but $\Lie_X u=0$
	in its maximum point.}
Therefore, $\sigma(X)$
is determined by the cohomology class of $\theta$.

\hfill

{\bf Step 2:} 
The anti-Lee field\index[terms]{Lee field!anti-} $I(\theta^\sharp)$
is distinguished by $\sigma(I(\theta^\sharp))=0$,
and the direction of $\theta^\sharp$ is determined by
$\sigma(\theta^\sharp)>0$. Therefore, to show that
$\theta^\sharp$ is independent on  the choice
of the Vaisman structure, it would suffice to
prove that $\sigma$ is independent.

Let $\theta_1, \theta_2$ be two Lee classes\index[terms]{class!Lee} of the
Vaisman metrics $\omega_1, \omega_2$ on $M$. 
Using the harmonic decomposition for 1-forms
(\ref{_harmo_deco_1-form_Proposition_}), 
we obtain that any harmonic 1-form on $(M, \omega_1)$
is proportional to $\theta_1 + \alpha$, where
$\alpha$ is $\Sigma$-basic. The set
of possible Lee classes for Vaisman structures 
is a half-space in $H^1(M,\R)$ (\ref{_Lee_cone_on_Vaisman_Theorem_}),
with the boundary represented by $\Sigma$-basic forms.\footnote{
This result is also proven by K. \index[persons]{Tsukada, K.} Tsukada, \cite{tsuk}.}
Therefore, we can always replace the Lee class\index[terms]{class!Lee} $[\theta]$ of a Vaisman
manifold by $[\const \theta]$, where $\const$ is positive.
Then, for some positive real constant $A$, the class
$[A\theta_1-\theta_2]$ is equal to $[\alpha]$, 
where $\alpha$ is a $\Sigma$-basic closed 1-form. 
Since the form $\alpha$ is basic, 
$\langle X, \alpha\rangle=0$ for any
$X$ tangent to $\Sigma$. Therefore, 
the character $\sigma_1:\; \goth s \arrow \R$ associated with
$[\theta_1]$ is proportional to the character
$\sigma_2:\; \goth s \arrow \R$  associated with
$[\theta_2]=[A\theta_1+\alpha]$.
\endproof


\section{The complex Monge-Amp\`ere equation}\index[terms]{Monge-Amp\`ere equation}


We start by introducing a version of a result of
\cite{_OV:EW_} which proves the uniqueness of the solution
of the Monge-Amp\`ere equation on Vaisman manifolds.\index[terms]{manifold!Vaisman}
Using a theorem of A. \index[persons]{El Kacimi-Alaoui, A.} El Kacimi-Alaoui (\cite{_Kacimi_}),
we show that the solution always exists. This result
is a complete analogue of the existence and uniqueness
of the solutions of the complex Monge-Amp\`ere equations
on a compact K\"ahler manifold proven by S.-T. \index[persons]{Yau, S.-T.} Yau.\index[terms]{Monge-Amp\`ere equation}

\hfill

\theorem\label{_CY_Vaisman_Theorem_} 
(\cite{ov_CY}) 
Let $(M,\omega,\theta)$ be a compact Vaisman manifold,\index[terms]{manifold!Vaisman}
and $V$ a volume form on $M$, satisfying $\int_M V = \int_M \omega^n$. 
Then there exists a unique Vaisman metric on $M$ with the same Lee class\index[terms]{class!Lee}
and the volume form $V$.

\hfill

This theorem 
is proven later in this section (\ref{_CY_Vaisman_by_volume_Theorem_}).
Recall that 
{\bf a trans\-ver\-sal volume form}, or {\bf basic volume form}, on a foliated manifold 
is a basic form $V\in \Lambda^k_F(M)$, $k= \codim F$, that defines a non-degenerate volume
locally on the leaf spaces of $F$. The following lemma
is a transversal form of  Calabi--Yau theorem, essentially
due to \cite[\S 3.5.5]{_Kacimi_}.\index[terms]{theorem!Calabi--Yau}

\hfill

\lemma\label{_transversal_volume_defines_omega_0_Lemma_}
Let $M$ be a compact Vaisman $n$-manifold,
and $\Sigma$ its canonical foliation. \index[terms]{foliation!canonical}
Then for any $\Sigma$-basic volume form $V_0$ that is 
cohomologous to an $(n-1)$-th power
of a transversal K\"ahler form $\eta_1$, there
exists a unique transversal K\"ahler form\index[terms]{form!K\"ahler!transversal} $\eta_2$
in the same basic cohomology \index[terms]{cohomology!basic}class 
such that $\eta_2^{n-1}=V_0$.

\hfill

\proof
We start by proving the uniqueness of a transversally
K\"ahler form with a given transversal volume.
Using the transversal $dd^c$-lemma 
(\ref{_Transversal_Lefschetz_Theorem_}), we obtain
$\eta_1=\eta_2 + dd^c f$, where $f$ is a $\Sigma$-basic
function (that is, a function that is  constant on the leaves of $\Sigma$).
Then $\eta_1^{n-1}-\eta_2^{n-1}= dd^c f\wedge P$, where
$P=\sum_{i=0}^{n-2} \eta_1^{i}\wedge \eta_2^{n-2-i}$.
Consider the operator
\[
f \xlongrightarrow{\ D_P\ } \frac{dd^c f\wedge P}{\eta_1^{n-1}}
\]
taking basic functions to basic functions. 
Let $U\subset M$ be a sufficiently small open set, and $X_U$
the leaf space of $\Sigma$ on $U$.
Clearly, $D_P:\; C^\infty(X_U) \arrow C^\infty(X_U)$ 
is a second order elliptic operator; by Hopf maximum
principle (\ref{hopf_theorem}), any non-constant $f \in \ker D_P$ cannot have a maximum.
However, any $\Sigma$-basic function on $M$ has a maximum
somewhere, because $M$ is compact. Therefore,
any $f\in \ker D_P$ is constant. This proves the
uniqueness of solutions. The existence of solutions
is obtained by repeating \index[persons]{Yau, S.-T.} Yau's argument in the transversal setup,
as done in \cite[\S 3.5.5 (iv)]{_Kacimi_}.
\endproof

\hfill

\lemma\label{_Vaisman_from_omega_0_and_Lee_Lemma_}
Let $M$ be a compact complex $n$-manifold of Vaisman type.
Then the Vaisman structure on $M$ is uniquely determined
by its transversal K\"ahler form\index[terms]{form!K\"ahler!transversal} $\omega_0$ (\ref{_Subva_Vaisman_Theorem_})
and the Lee class\index[terms]{class!Lee} $[\theta]\in H^1(M, \R)$.\footnote{The Lee
class $[\theta]$ is uniquely determined by the homothety character \index[terms]{homothety character}$\chi$
(\ref{_homothety_character_Definition_}) 
and determines it; or, equivalently, $[\theta]$ is uniquely determined by 
the weight local system (\ref{_weight_bun_min_Kahler_Definition_}).}

\hfill

\proof
By \eqref{_omega_via_theta_Chapter_8_Equation_}, 
we have $\omega=\omega_0 + \theta \wedge\theta^c$. 
Therefore, it would suffice
to show that the Lee form $\theta$ is uniquely determined 
by $\omega_0$ and the Lee class.\index[terms]{class!Lee}
Let $\theta$ and $\theta'$ be two Lee forms\index[terms]{form!Lee} of Vaisman manifolds,\index[terms]{manifold!Vaisman}
with the same transversal K\"ahler form\index[terms]{form!K\"ahler!transversal} $\omega_0$. 
Denote by $\eta$ the 1-form $\theta- \theta'$.
Since $\omega_0 = d^c \theta=  d^c \theta'$ 
\eqref{_omega_0_Equation_}
this would imply $d^c\eta= d\eta=0$. 
Such a 1-form cannot be exact, because 
if $\eta=df$, one has $dd^c f=0$;\index[terms]{maximum principle}
however, pluriharmonic functions are constant on any
compact manifold by the maximum principle. Therefore, $\theta$ 
cannot be cohomologous to $\theta'$.\index[terms]{function!pluriharmonic}
\endproof

\hfill

\lemma\label{_Vaisman_volume_from_omega_0_Lemma_}
Let $(M, \omega, \theta)$ be a Vaisman $n$-manifold,
and $\omega_0$ its transversal K\"ahler form (\ref{_Subva_Vaisman_Theorem_}).
Then $i_{I(\theta^\sharp)} i_{\theta^{\sharp}} \omega^n = n \omega_0^{n-1}$,
where $\theta^\sharp$ is the Lee field.\index[terms]{Lee field}

\hfill

\proof
By \eqref{_omega_via_theta_Chapter_8_Equation_}, 
we have $\omega=\omega_0 + \theta \wedge\theta^c$. 
This implies $\omega^n = \omega_0^{n-1} \wedge \theta \wedge\theta^c$.
The form $\omega_0$ is $\Sigma$-basic, and $|\theta^{\sharp}|=1$.
Then 
\[ 
i_{I(\theta^\sharp)} i_{\theta^{\sharp}} \omega^n = 
n i_{I(\theta^\sharp)}i_{\theta^{\sharp}} (\theta \wedge\theta^c)\wedge \omega_0^{n-1}
=n\omega_0^{n-1}.\ \ \  \endproof
\]

\theorem\label{_CY_Vaisman_by_volume_Theorem_} (\cite{ov_CY}) 
Let $(M,\omega,\theta)$ be a compact Vaisman $n$-manifold.
Then for each volume form $V'$ on $M$ such that
$\int_M\omega^n = \int_M V'$, there exists a unique 
Vaisman metric $\omega'$ with volume $V'$ and 
the same Lee class.\index[terms]{class!Lee}

\hfill

\proof
The transversally K\"ahler form\index[terms]{form!K\"ahler!transversal} $\omega_0$ 
associated with the Vaisman structure
is uniquely defined by its
transversal volume form $V_0=\omega_0^{n-1}$ 
and its basic cohomology\index[terms]{cohomology!basic} class
$[\omega_0]\in H^2_\Sigma(M)$, assuming that the transversal
cohomology classes of $\omega^{n-1}$ and $V$ are equal
(\ref{_transversal_volume_defines_omega_0_Lemma_}). However,
the class $[\omega_0]$ generates the kernel of the natural map
$H^2_\Sigma(M)\arrow H^2(M)$ (\ref{_Vaisman_harmonic_forms_Theorem_}).
This determines the basic cohomology class of $\omega_0$ up to a constant;
the constant is fixed if the transversal volume is fixed.
However, the transversal volume form $\omega_0^{n-1}$
is determined uniquely, up to a constant multiplier,
 by the volume form 
$V$, as follows from \ref{_Vaisman_volume_from_omega_0_Lemma_}.
We obtain that the volume form
of a Vaisman metric uniquely defines $\omega_0$, in such a way
that $\omega_0^{n-1}=\const i_{I(\theta^\sharp)} i_{\theta^{\sharp}} V$
and the Lee class\index[terms]{class!Lee} of the Vaisman structure associated with $\omega_0$
is proportional to $[\theta]$.

By \ref{_Lee_field_direction_Proposition_}, 
the Lee field\index[terms]{Lee field} $\theta^\sharp$ of a Vaisman structure is uniquely
(up to a constant) determined by the complex structure of $M$.
By \ref{_Vaisman_volume_from_omega_0_Lemma_}, 
$V_0= n i_{I(\theta^\sharp)}i_{\theta^\sharp}V$ uniquely (up to a constant) defines
the transversal volume form  $\omega_0^{n-1}$, and hence  $\omega_0$
exists and is uniquely defined, up to a constant.
Now, the Vaisman metric $\omega= \omega_0+ \theta \wedge \theta^c$
is uniquely defined by
$\omega_0$ and the Lee class\index[terms]{class!Lee} (\ref{_Vaisman_from_omega_0_and_Lee_Lemma_}),
hence the Vaisman metric $\omega$ is uniquely 
(up to a constant) defined by its volume form 
$V$ and the Lee class. The constant is also fixed,
because $V= \omega^n$.

This proves uniqueness of a Vaisman metric with prescribed volume.
To see that a metric with a prescribed volume
form $V'$ exists, we write the corresponding transversal
volume form $V_0':= n^{-1} i_{I(\theta^\sharp)}i_{\theta^\sharp}V'$
and solve the transversal Calabi--Yau equation 
(\ref{_transversal_volume_defines_omega_0_Lemma_}), arriving at
a transversal K\"ahler form $\omega_0'= d^c\theta+ dd^c f =d^c\theta'$ satisfying
$(\omega_0')^{n-1}=V_0'$. Then 
$\omega':=\omega_0' + \theta'\wedge (\theta')^c$ is an LCK form 
that is  invariant under the Lee field action; by 
\ref{kami_or},
any LCK metric $\omega'$ admitting a conformal holomorphic
flow, non-isometric on its K\"ahler cover, is Vaisman.
\endproof

\section{Exercises}

The exercises in this chapter are 
closer to the research problems; the
reader should not be discouraged.

\begin{enumerate}[label=\textbf{\thechapter.\arabic*}.,ref=\thechapter.\arabic{enumi}]

\item
Find a Vaisman Hopf surface such that a general
leaf of the canonical foliation is\index[terms]{foliation!canonical} dense.\index[terms]{surface!Hopf}

\item Let $M$ be a Vaisman manifold such that the general
leaf of the canonical foliation is dense. Prove that the
Vaisman structure on $M$ is unique in any Lee class.\index[terms]{manifold!Vaisman}\index[terms]{class!Lee}

\item \label{_adjunction_Exercise_}
(``The adjunction formula''.)
Let $M$ be a quasi-regular Vaisman manifold,\index[terms]{manifold!Vaisman!quasi-regular} and
$\tau:\; M \arrow X$ the projection to the leaf
space of the canonical foliation. Prove that
$K_M\cong \pi^* K_X$, where $K_M$, $K_X$ 
denotes the canonical bundle.

\item Let $M$ be a Vaisman manifold with trivial canonical
bundle. Prove that $M$ is a deformation of a quasi-regular
Vaisman manifold $M_1$ such that the leaf space of the
canonical foliation is a Calabi--Yau orbifold.\index[terms]{orbifold!Calabi--Yau}

{\em Hint:} Use the previous exercise.

\item\label{_linearizable_Exercise_}
Let $X$ be a projective manifold, $L$ an ample
line bundle on $X$, and $\omega$ its curvature, that is  
a K\"ahler form. \index[terms]{form!K\"ahler}Consider a holomorphic vector field
$v$ such that $\Lie_v \omega=0$. Prove that 
the following are equivalent (See also \cite[pp. 72-74, 187]{_Gauduchon:book_}):
\begin{description}
\item[(i)] The action of $v$ is Hamiltonian 
with respect to $\omega$.
\item[(ii)] The action of the group $e^{tv}$
can be holomorphically extended to an action on the total
space $\Tot(L)$ commuting with the standard $\C^*$-action.
\item[(iii)] The action of $e^{tv}$ on $\Pic(X)$
fixes $[L]\in \Pic(X)$.
\end{description}
An automorphism $s \in \Aut(X)$ is called {\bf linearizable}
if its action can be extended to an equivariant action on $L$,
that is, extended to $\Tot(L)$ holomorphically.

\item
Let $M$ be a quasi-regular Vaisman manifold,\index[terms]{manifold!Vaisman!quasi-regular} 
and $\tau:\; M \arrow X$ the projection to the leaf space.
By \ref{_Structure_of_quasi_regular_Vasman:Theorem_}, 
$M$ is isomorphic to a $\Z$-quotient
of the space $\Tot^\circ(L)$ of non-zero vectors in
an ample bundle $L$. Prove that the action of a 
holomorphic vector field  $v\in TX$
can be extended to $M$ if and only if it is 
linearizable.

\item
Let $X$ be an abelian manifold 
(i.e. a projective complex torus). Prove that $X$ does not
admit holomorphic, Hamiltonian vector fields.
Prove that the group $\Aut_L(X)$ of
all automorphisms of $X$ extending to an equivariant
action on $L$ (that is, the linearizable automorphisms) 
is finite.

\item
Let $M$ be a Vaisman manifold such that\index[terms]{manifold!Vaisman!quasi-regular}
its $\Z$-cover is K\"ahler and isomorphic
to the space $\Tot^\circ(L)$ of non-zero vectors in
an ample bundle $L$ on an abelian manifold.
Prove that $M$ is quasi-regular.

{\em Hint:} State and solve an orbifold version of the
previous exercise and use it.

\item\label{_no_linearizable_then_qreg_Exercise_}
Let $M$ be a Vaisman manifold such that
its $\Z$-cover is K\"ahler and isomorphic
to the space $\Tot^\circ(L)$ of non-zero vectors in
an ample bundle $L$ on a projective manifold $X$.
Assume that the group of $L$-linearizable
automorphisms of $X$ is finite.
Prove that $M$ is quasi-regular.

\item\label{_holo_parallel_Exercise_}
Let $X$ be a compact Ricci-flat K\"ahler orbifold.
Prove that all holomorphic vector fields on $X$
are parallel. \index[terms]{vector field!parallel}\index[terms]{vector field!holomorphic}

{\em Hint:} Use an orbifold version of 
\index[persons]{Bochner, S.} Bochner vanishing theorem.

\item\label{_parallel_not_lineara_Exercise_}
Let $v$ be a parallel vector field on 
a compact projective K\"ahler manifold $X$.
\begin{enumerate}
\item Prove that the Albanese manifold\index[terms]{manifold!Albanese}
$\Alb(X)$ is non-trivial.
\item Prove that the action of $e^{tv}$
on $\Alb(X)$ is non-trivial.
\item Prove that the action of $e^{tv}$ on $X$
is not linearizable with respect to any
ample line bundle.
\end{enumerate}

{\em Hint:} Use the isomorphism $\Alb(X) = \Pic(X)^*_0$
and Exercise \ref{_linearizable_Exercise_}.

\item\label{_Ricci_flat_linearizable_Exercise_}
Let $X$ be a compact Ricci-flat projective orbifold.
Prove that the group of linearizable holomorphic
isometries of $X$ is finite.

{\em Hint:} Use Exercise \ref{_holo_parallel_Exercise_}
and Exercise \ref{_parallel_not_lineara_Exercise_}.

\item
Let $M$ be a quasi-regular 
Vaisman manifold,\index[terms]{manifold!Vaisman!quasi-regular}
$\tau:\; M \arrow X$ the projection to the leaf space,
and $\pi:\; \tilde M \arrow M$ a K\"ahler $\Z$-covering.
Assume that $X$ is a Calabi--Yau orbifold.
Prove that any Vaisman $\Z$-quotient $\frac{\tilde M}{\Z}$
is quasi-regular.\index[terms]{orbifold!Calabi--Yau}

{\em Hint:} Exercise \ref{_Ricci_flat_linearizable_Exercise_} 
and  Exercise \ref{_no_linearizable_then_qreg_Exercise_}.

\item\label{_trivial_canonical_then_QR_Exercise_}
Let $M$ be a compact Vaisman manifold with trivial canonical
bundle. Prove that $M$ is quasi-regular.

{\em Hint:} Use the previous exercise and
Exercise \ref{_adjunction_Exercise_}.

\item
Let $M$ be a compact  Calabi--Yau orbifold with $c_2(M)=0$.
Prove that the Ricci-flat metric on $M$ is flat.
Use the  Bogomolov-L\"ubke inequality 
(\cite[\S 5]{_Lubke:Chernklassen_}).\index[terms]{inequality!Bogomolov-L\"ubke}

\item\label{_parallelizable_Vaisman_Exercise_}
Let $M$ be a compact Vaisman manifold\index[terms]{manifold!Vaisman} with trivial canonical
bundle, that is  quasi-regular by Exercise 
\ref{_trivial_canonical_then_QR_Exercise_}.
Denote by $X$ the leaf space of the canonical
foliation on $M$, considered to be  a projective orbifold.
Assume that the tangent bundle $TM$ is holomorphically
trivial. Prove that $X$ is a K\"ahler orbifold with
vanishing curvature. \footnote{By Bieberbach's solution
of Hilbert's 18 problem, any compact flat Riemannian
orbifold is a finite quotient of a torus (\cite{_Buser:Bieberbach_}).
Therefore, all Vaisman manifolds with trivial canonical bundle and holomorphically trivial
tangent bundle  are finite quotients of
Heisenberg-type Vaisman nilmanifolds,\index[terms]{nilmanifold} \ref{kod_nil}.}

{\em Hint:} Use the previous exercise. 

\item
Let $M$ be a compact Vaisman manifold and $\tilde M$
its K\"ahler $\Z$-cover. Assume that $\tilde M$
is Ricci-flat (such a Vaisman manifold is called
{\bf Einstein--Weyl}, see Section \ref{eiwey}).\index[terms]{manifold!Vaisman}\index[terms]{manifold!Einstein--Weyl}
Prove that the canonical bundle $K_M$ is non-trivial.

\end{enumerate}


\chapter{Holomorphic tensor fields on LCK manifolds with potential}\index[terms]{manifold!LCK!with potential}
\label{_tensor_invariant_Chapter_}


{\setlength\epigraphwidth{0.7\linewidth}
\epigraph{\it The desire to kill the old ones. The old ones are rotten. The young are mediocre and stupid. The young have always been insufferable to me, especially when I was young. The masters of thought courting the youngsters, walking after the client\`ele, is one of the most humiliating things I know. What lack of dignity, what cowardice, what nonsense!}{\sc\scriptsize Eug\`ene Ionesco, \ \ Fragments of a Journal}
 }


\section{Introduction}


Every Vaisman manifold\index[terms]{manifold!Vaisman} $M$ is equipped with an action of
the Lee\index[terms]{Lee field} and anti-Lee fields\index[terms]{Lee field!anti-}, that are  Killing and holomorphic\index[terms]{vector field!Killing}\index[terms]{vector field!holomorphic}
(\ref{_Lee_field_Killing_Corollary_}).
In \cite{_Verbitsky:Sta_Elli_}, 
it was shown that any stable holomorphic bundle on $M$
is equivariant with respect to this action. 
From \ref{_Subva_Vaisman_Theorem_}, it 
follows that any complex subvariety is invariant
under this action as well.

Let $G$ be the group of holomorphic isometries
of a Vaisman manifold\index[terms]{manifold!Vaisman} $M$ obtained as the closure of the
group generated by the Lee and the anti-Lee flow.\index[terms]{Lee flow}\index[terms]{Lee flow!anti-}
In this chapter, we prove that any
holomorphic tensor field on $M$ is
$G$-invariant. 

For holomorphic vector fields and differential forms this result was obtained
by K. \index[persons]{Tsukada, K.} Tsukada in \cite[Theorem 3.3, Theorem 4.2]{tsuk}.

We prove a similar result for LCK manifolds with potential.
Let $M$ be an LCK manifold with potential,\index[terms]{manifold!LCK!with potential} and $\tilde M$
its K\"ahler $\Z$-cover. We assume that $\tilde M$ is an open 
algebraic cone;\index[terms]{cone!algebraic} this is true automatically when $\dim_\C M\geq 3$
(\ref{_cone_cover_for_LCK_pot_Theorem_}).
When $\dim_\C M=2$, this would follow if we assume 
the GSS conjecture\index[terms]{conjecture!GSS}
(\ref{_sphe_implies_dim2_Theorem_}). Since $\tilde M$ is algebraic,
it makes sense to speak of the Zariski closure of the
$\Z$-action on $\tilde M$. We consider the Zariski 
closure\index[terms]{Zariski closure}
${\cal G}$ of the $\Z$-action as the smallest algebraic
subgroup of $\Aut(\tilde M)$ containing this $\Z$-action;
this group is clearly commutative and positive-dimensional. Since 
the action of ${\cal G}$ commutes with the action of the deck
transform group $\Z\subset {\cal G}$ on $\tilde M$, 
the quotient ${\cal G}/\Z$ naturally acts on $M$.

We prove that any holomorphic tensor field on 
an LCK manifold with potential\index[terms]{manifold!LCK!with potential} is ${\cal G}$-invariant
(\ref{_LCK_pot_Zariski_closure_Proposition_}).
When the $\Z$-action on $\tilde M$ admits a logarithm,
this also implies that any holomorphic tensor
field is invariant with respect to the logarithmic
vector field action. \index[terms]{action!logarithmic}

However, the relation between 
the logarithm field and the
Lee field on a Vaisman manifold\index[terms]{manifold!Vaisman}
was not very clear so far. \ref{_Lee_field_direction_Proposition_}
implies that the Lee field\index[terms]{Lee field} is determined uniquely, up to a 
complex multiplier, by the complex structure.
In this chapter, we describe the Lee field explicitly
in terms of the holomorphic embedding to a Hopf manifold.
From \ref{_Subva_Vaisman_Theorem_} it follows that the
Lee field\index[terms]{Lee field} commutes with holomorphic immersions of compact
Vaisman manifolds, and hence  to determine the Lee field \index[terms]{Lee field}we
need only to write it down on a Hopf manifold.

Let $H = \frac{\C^n \backslash 0}{\langle A \rangle}$
be a diagonal Hopf manifold, and ${\cal G}$ the
Zariski closure of ${\langle A \rangle}$.
Let $\alpha_i$ be the eigenvalues of $A$,
and $A_1$ an operator that has the eigenvalues $|\alpha_i|$
in the same basis. In  \ref{_Zariski_contains_real_part_Theorem_},
we show that the linear field $\log \!A_1$ belongs to the
Zariski closure of ${\langle A \rangle}$, and 
in \ref{_holo_tensor_on_Vaisman_Lee_invariant_Theorem_},
Step 3, we prove that $\log A_1$ can be obtained
as the Lee field \index[terms]{Lee field}of a Vaisman structure on $H$.

This implies that the Lee and anti-Lee field\index[terms]{Lee field!anti-}\index[terms]{Lee field} of a Vaisman 
manifold belong to the Lie algebra of ${\cal G}$
(the Zariski closure of the $\Z$-action).
Therefore, any holomorphic tensor field is Lee and anti-Lee
invariant.


\section{Holomorphic tensors on LCK manifolds with potential}\index[terms]{manifold!LCK!with potential}


In this section, we prove the following useful result.
A weaker version of this proposition appeared
in \cite{ov_imrn_10}.

\hfill

\proposition\label{_LCK_pot_Zariski_closure_Proposition_}
Let $M$ be a compact LCK manifold with potential,\index[terms]{manifold!LCK!with potential}
$\tilde M$ its K\"ahler $\Z$-cover, considered to be 
an open algebraic cone,\index[terms]{cone!algebraic} and 
$\Phi\in H^0(M, B)$  be a holomorphic tensor field on $M$,
where  $B=(\Omega^1 M)^{\otimes k}\otimes TM^{\otimes l}$. Denote by  $\tilde \Phi$
its lift to $\tilde M$, and let ${\cal G}$ be the Zariski closure
of the $\Z$-action on $\tilde M$. Then $\tilde\Phi$ is
${\cal G}$-invariant.\index[terms]{Zariski closure}

\hfill

\pstep
Consider a $\Z$-action on a finite-dimensional space $W$.
Then any $\Z$-invariant vector $w\in W$ is invariant under the
Zariski closure of $\Z$. This is why the statement
of \ref{_LCK_pot_Zariski_closure_Proposition_}
is not surprising. However, the space of 
tensor fields on $\tilde M$ is not
finite-dimensional, which makes \ref{_LCK_pot_Zariski_closure_Proposition_}
non-trivial. Denote by ${\goth m}$ the maximal ideal
of the origin in the closed algebraic cone\index[terms]{cone!algebraic} $\tilde M_c$.
We think of the quotients $\calo_{\tilde M_c}/ {\goth m}^k$
as of the spaces of jets of holomorphic (or algebraic)
functions. As in Subsection \ref{_normality_of_cone_Subsection_},
we choose $\tilde M_c$ normal, so that the $\Z$-action
is extended to $\tilde M_c$. Then any
$\Z$-invariant jet $u \in \calo_{\tilde M_c}/ {\goth m}^k$
is also ${\cal G}$-invariant.

\hfill

{\bf Step 2:} From now on, we use the same letter $\Phi$
to denote the lift of $\Phi$ to $\tilde M$.
Consider $\Phi$ as a
$\Z$-invariant section of the appropriate
tensor bundle $B=(\Omega^1 M)^{\otimes k}\otimes TM^{\otimes l}$,
that is  by construction $\Z$-equivariant.
Using \ref{_extension_over_a_point_Theorem_}, we extend
$B$ to a reflexive coherent sheaf $B_c$ on $\tilde M_c$.
Since reflexive sheaves are normal \index[terms]{sheaf!normal}
(\ref{_reflexife_normal_Theorem_}), the section $\Phi$
admits a holomorphic extension\index[terms]{sheaf!reflexive}
to $\tilde M_c$, denoted as  $\Phi_c\in H^0(\tilde M_c, B_c)$.
Consider the space $J_B^k:=\frac{B_c}{{\goth m}^{k+1} B_c}$
of $k$-jets of $B_c$. The $\Z$-action on $H^0(\tilde M_c, B_c)$
preserves 
$H^0(\tilde M_c, {\goth m}^k B_c)\subset H^0(\tilde M_c, B_c)$.
Let $\Phi^k_c \in J_B^k$ be the $k$-jet of $\tilde\Phi_c$.
Since $\Phi_c$ is $\Z$-invariant, and the space of $k$-jets
is finite-dimensional, the $k$-jet $\Phi^k_c$
is ${\cal G}$-invariant. Let \index[terms]{completion!adic}
$\hat B_c:= \lim\limits_{\leftarrow} J_B^k$
be the ${\goth m}$-adic completion of $B_c$, and
$\hat \Phi_c$ the image of $\Phi_c$ in this
completion. Since all $k$-jets $\Phi^k_c$ are ${\cal G}$-invariant,
the completion $\hat \Phi_c$ is also ${\cal G}$-invariant.
Since the ring of germs of holomorphic
functions on $\tilde M_c$ in $c$ is Noetherian 
(\cite[Chapter II, Theorem B.9]{_Gunning_Rossi_}), and $B_c$, 
being reflexive, is torsion-free, the\index[terms]{ring!torsion-free}\index[terms]{ring!Noetherian}
natural map $H^0(\tilde M_c, B_c)\arrow \hat B_c$
is injective (\cite[Theorem 10.17]{_Atiyah_MacDonald_}). 
Hence, $\Phi\in  H^0(\tilde M_c, B_c)$
is also ${\cal G}$-invariant.
%
\endproof


\section{Zariski closures and the Chevalley theorem}


We recall the following famous theorem due to \index[persons]{Chevalley, C.} Chevalley.

\hfill

Let $G\subset \GL(V)$ be an algebraic
group, and $W= V^{\otimes k} \otimes (V^*)^{\otimes l}$
a tensor representation of $G$.  A $G$-invariant
vector $v \in W$ is called {\bf a tensor invariant of $G$}.
A point $x\in {\Bbb P} W$ is called
{\bf a projective tensor invariant of $G$}
if it is $G$-invariant. 

\hfill

\theorem\index[terms]{theorem!Chevalley}
{\bf (Chevalley theorem,} \cite{_Morris:Ratner_}).\\
An algebraic\index[terms]{group!algebraic}
group is uniquely determined by the set
of its projective tensor invariants.
\endproof

\hfill

When an algebraic group is {\em reductive} (over $\C$, a group
is reductive if and only if it has a compact real form),
a stronger version of  \index[persons]{Chevalley, C.} Chevalley theorem is available.

\hfill

\theorem\label{_Chevalley_reductive_Theorem_}
A reductive algebraic\index[terms]{group!algebraic!reductive}
group is uniquely determined by the set
of its tensor invariants.

\proof 
\cite[Proposition 3.1 (c)]{_Deligne:Hodge_cycles_}.
\endproof

\hfill

We are going to prove the following result, used further
on to determine the Lee field\index[terms]{Lee field} action on a Vaisman manifold.\index[terms]{manifold!Vaisman}

\hfill

\theorem\label{_Zariski_contains_real_part_Theorem_}
Let $A \in \GL(n, \C)$ be a diagonal linear operator
with the eigenvalues $\alpha_1, ..., \alpha_n$,
and $A_1$ an operator that is  diagonal in the
same basis, with the eigenvalues $|\alpha_1|,..., |\alpha_n|$.
Denote by ${\cal G}$ the Zariski closure of the
group $\langle A \rangle$ in $\GL(n, \C)$.
Then ${\cal G}$ contains $A_1$.

\hfill

\proof
Since the operator $A$ is diagonalizable, 
the algebraic closure of $\langle A \rangle$ 
is a commutative group that contains only
semisimple elements. By 
\cite[Proposition 1.5]{_Borel_Tits:Groupes_Reductifs_},
its connected component of unity
is $(\C^*)^d$, and hence  it is reductive.
Then \ref{_Chevalley_reductive_Theorem_}
implies that ${\cal G}$ is the set of all
$g\in \GL(n, \C)$ that preserve all
$A$-invariant vectors 
$w\in W= V^{\otimes k} \otimes (V^*)^{\otimes l}$.
The eigenvalues of $A$ on $W$
are products of $k$ instances of $\alpha_i$
and $l$ instances of $\alpha_i^{-1}$.
The eigenvectors $w \in W$ 
are $
w = \bigotimes_{i=1}^k z_{m_i} \otimes \bigotimes_{i=1}^l \zeta_{n_i}$, where
$z_{n_i}$ are the eigenvectors in $V$ with the
eigenvalues $\alpha_{n_i}$
and $\zeta_{m_i}$ are the eigenvectors in $V^*$ with the
eigenvalues $\alpha_{m_i}^{-1}$. Then
$A(w) = \prod_{i=1}^k \alpha_{n_i} \prod_{i=1}^l \alpha^{-1}_{m_i}$.
The eigenvector $w$ is $A$-invariant
if and only if $\prod _{i=1}^k \alpha_{n_i} \prod_{i=1}^l \alpha^{-1}_{m_i}=1$.
This implies $\prod _{i=1}^k |\alpha_{n_i}| \prod_{i=1}^l |\alpha^{-1}_{m_i}|=1$,
hence $w$ is $A_1$-invariant as well.
\endproof


\section{Holomorphic tensors on Vaisman manifolds}\index[terms]{manifold!Vaisman}


Let $M$ be a Vaisman manifold, and $\theta^\sharp$ its Lee 
field. Then $\theta^\sharp$ and $I(\theta^\sharp)$
are holomorphic Killing fields\index[terms]{vector field!Killing}\index[terms]{vector field!holomorphic} (\ref{_canon_foli_totally_geodesic_Remark_}).
The closure $G$ of the group 
$e^{t_1 \theta^\sharp+t_2 I(\theta^\sharp)}$ in the group of isometries of 
$M$ is compact and  commutative, and hence  it is a compact torus
(\ref{_invariance_for_cohomology_Vaisman_Lemma_}).

\hfill

\theorem\label{_holo_tensor_on_Vaisman_Lee_invariant_Theorem_}
Let $M$ be a compact Vaisman manifold, \index[terms]{manifold!Vaisman}
$\Phi\in H^0(M, B)$ 
a holomorphic tensor, where $B= (\Omega^1 M)^{\otimes k}\otimes TM^{\otimes l}$,
and $G$ the smallest closed Lie group containing the
flows generated by the Lee\index[terms]{Lee field} and the anti-Lee fields.\index[terms]{Lee field!anti-}
Then $\Phi$ is $G$-invariant. 

\hfill

\pstep
Choose a Vaisman metric on $M$ with LCK rank\index[terms]{rank!LCK} 1, and
let $\tilde M$ be the corresponding 
K\"ahler $\Z$-cover of $M$, that is  considered
as an open algebraic cone\index[terms]{cone!algebraic}. Let ${\cal G}\subset \Aut(\tilde M)$ be 
the Zariski closure of the $\Z$-action.\index[terms]{action!$\Z$-} Then $\Phi$ is
${\cal G}$-invariant by \ref{_LCK_pot_Zariski_closure_Proposition_}.
To finish the proof of \ref{_holo_tensor_on_Vaisman_Lee_invariant_Theorem_},
it would suffice to show that $G \subset {\cal G}$;
since $ {\cal G}$ is closed, this would follow if
we prove that the Lee field\index[terms]{Lee field} $\theta^\sharp$
is tangent to ${\cal G}$.

\hfill

{\bf Step 2:} 
Let $M \arrow H$ be a holomorphic embedding of $M$ to a
diagonal Hopf manifold. By \ref{_Subva_Vaisman_Theorem_}, this embedding
commutes with the group generated by the
Lee and the anti-Lee flow,
\index[terms]{Lee flow!anti-} hence it would suffice
to show that $\theta^\sharp$ is tangent to ${\cal G}$ when $M=H$.
It remains to prove  \ref{_holo_tensor_on_Vaisman_Lee_invariant_Theorem_}
assuming that $M$ is a Hopf manifold. For this purpose,
we compute the groups $G$ and ${\cal G}$ for 
a diagonal Hopf manifold. This can be done explicitly
using the expression for the Vaisman metric on a Hopf manifold
obtained in \ref{shell_char}.

\hfill

{\bf Step 3:} The rank 2 algebra generated by
the Lee and anti-Lee flows\index[terms]{Lee flow!anti-} on a Vaisman
manifold $M$ is uniquely determined by 
the complex structure on $M$ (\ref{_Subva_Vaisman_Theorem_}). 
We express the Lee flow on a Hopf manifold
using a Vaisman structure
we shall construct explicitly, and
prove that it belongs to $\Lie({\cal G})$.
This would imply that the original Lee flow 
also belongs to $\Lie({\cal G})$.

Let $H = \frac{\C^n \backslash 0}{{\langle A \rangle}}$,
where $A \in \GL(n, \C)$ is a diagonal linear contraction with eigenvalues
$\alpha_1, ..., \alpha_n$.  Let $A_1$ be the matrix that is 
diagonal in the same basis as $A$, and has eigenvalues
$|\alpha_i|$. Then $\vec r:=\log A_1$
is a holomorphic contraction vector field,
and  $e^{t\vec r}$ is an $A$-invariant
contraction flow. The sphere $S^{2n-1}$ 
is a pseudoconvex shell\index[terms]{pseudoconvex shell} compatible with this contraction flow.
Using \ref{shell_char}, we obtain a $e^{t\log A_1}$-automorphic
plurisubharmonic function $\phi$. We are going to prove
that $\phi$ is $A$-automorphic.

Consider the diagonal $\U(1)^n$-action on $\C^n$.
The plurisubharmonic function $\phi$ is $\U(1)^n$-invariant
because $\vec r$ and $S^{2n-1}$ are $\U(1)^n$-invariant.
Since $A\in \U(1)^n \cdot A_1$,
and $\phi$ is $A_1$-automorphic, it is $A$-automorphic.
Therefore, $\frac{dd^c\phi}{\phi}$ defines a Vaisman metric 
on $H$. 
By definition, the Lee flow\index[terms]{Lee flow} maps the level sets
of $\phi$ to the level sets of $\phi$, with 
$\phi(e^{t\theta^\sharp}z)= t\phi(z)$.
By construction of $\phi$,
the same is true for $e^{-t\log A_1}$, that is,
$\phi(e^{-t\log A_1}z)= t\phi(z)$.
Therefore, $\theta^\sharp = -\log A_1$, where $\log A_1$
is the diagonal matrix with the eigenvalues $\log |\alpha_i|$
in the same basis. This is the Lee flow for the
Vaisman metric $\frac{dd^c\phi}{\phi}$. 

The vector field
$\log A_1$ belongs to $\Lie({\cal G})$,
because any power of $A_1$ fixes all
tensor invariants of ${\cal G}$
by \ref{_Zariski_contains_real_part_Theorem_},
hence the corresponding 1-parametric
subgroup also belongs to $\Lie({\cal G})$;
this gives $\log A_1\in \Lie({\cal G})$.
\endproof

\section{Exercises}

\begin{enumerate}%
[label=\textbf{\thechapter.\arabic*}.,ref=\thechapter.\arabic{enumi}]

\item
Let $V$ be a holomorphic vector field on a Vaisman manifold.\index[terms]{manifold!Vaisman}
Prove that $[V, \theta^\sharp]=0$, where $\theta^\sharp$ is the
Lee field.\index[terms]{Lee field}

\item \label{_decompo_holo_on_Vaisman_Exercise_}
Prove that any real holomorphic
vector field $\xi$ on a Vaisman manifold 
can be written as a sum $\xi= \xi_0 + a\theta^\sharp + b I (\theta^\sharp)$,
where $a, b \in \R$, and $\xi_0$ is a $\Sigma$-basic
holomorphic vector field.

\item Prove that all holomorphic differential forms on
a compact Vaisman manifolds are basic with respect to the
canonical foliation.\index[terms]{foliation!canonical}

\item
Let $H= \frac{\C^2 \backslash 0}{\lambda \Id}$
be a classical Hopf manifold. Prove that 
the Lee field\index[terms]{Lee field} is proportional to $z_1 \frac d{dz_1} + z_2 \frac d{dz_2}$
with a positive real coefficient.

\item
Find an example of a quasi-regular Vaisman manifold\index[terms]{manifold!Vaisman!quasi-regular}
such that all orbits of the Lee field are non-compact.

{\em Hint:} Use the previous exercise.

\item 
Let $M$ be a quasi-regular Vaisman manifold of LCK rank 1\index[terms]{rank!LCK},
and $\tilde M$ its K\"ahler $\Z$-covering. 
Denote by $\gamma$ the generator of the 
deck transform group acting on $\tilde M$.
Find an example of a Vaisman manifold of LCK rank 1
such that $\langle \gamma\rangle$ acts freely on
the space of leaves of the canonical foliation.

\item
Let $M_1$, $M_2$ be Vaisman manifolds\index[terms]{manifold!Vaisman} of LCK rank 1,\index[terms]{rank!LCK}
such that their K\"ahler $\Z$-covers $\tilde M_1$, $\tilde M_2$ 
are holomorphically
isometric. 
\begin{enumerate}
\item Find an example of such manifolds
with $M_1$ not biholomorphic to $M_2$.
\item
Suppose that $\tilde M_1$ does not admit
a holomorphic function that is  automorphic with
respect to the deck transform $\Z$-action.\index[terms]{action!$\Z$-} Let
$X_1$, $X_2$ be the Lee fields\index[terms]{Lee field} of $M_1, M_2$
lifted to $\tilde M_1$, $\tilde M_2$.
Prove that the holomorphic isometry
$\phi:\; \tilde M_1 \arrow \tilde M_2$
takes $X_1$ to $X_2$.
\end{enumerate}

\item
Let $M$, $\dim_\C M=n$, be a Vaisman manifold\index[terms]{manifold!Vaisman}
such that its K\"ahler $\Z$-cover is a cone over 
a Sasakian manifold $S$.
Assume that the closure of a general Reeb orbit
in $S$ is $n$-dimensional. Prove that 
the union of all proper complex subvarieties
$Z \subset M$ is not equal to $M$.

The following exercise assumes a certain knowledge of the\index[terms]{GIT}
GIT (Geometric Invariant Theory), \cite{_Mumford_Fogarty_Kirwan:GIT_}.


\item
Let $M$ be a Vaisman manifold,\index[terms]{manifold!Vaisman} $\tilde M$
its K\"ahler $\Z$-cover, and ${\cal G}$ the Zariski closure
of the $\Z$-action on $\tilde M$, acting on $M$.
Since $\tilde M$ is an algebraic cone\index[terms]{cone!algebraic},
it can be obtained as the total space
$\Tot^\circ(L)$ of non-zero vectors in 
an ample line bundle over an orbifold
$X$. 
\begin{enumerate}
\item Prove that the group ${\cal G}$
naturally acts on $X$.
\item Prove that any ${\cal G}$-invariant meromorphic
function on $X$ can be extended to a 
meromorphic function on $M$.
\item Suppose ${\cal G}$ acts
on $X$ with all orbits of dimension
$<\dim X$. Prove that there
exists a non-constant ${\cal G}$-invariant
meromorphic function on $X$.
\end{enumerate}

\item
A Vaisman manifold\index[terms]{manifold!Vaisman} is called\index[terms]{manifold!toric}
{\bf toric} if it admits a holomorphic
$(\C^*)^n$-action with an open orbit.
Prove that a Vaisman manifold \index[terms]{orbit!open}
without global meromorphic functions
is toric.\footnote{See Section \ref{toric_lck}
for the definition  of toric Vaisman manifolds.}
\index[terms]{function!meromorphic}

{\em Hint:} Use the previous exercise.

\item
Let $A \in \End(\C^n)$ be a linear contraction, 
that is  not diagonal, and ${\cal G}_A$ the Zariski
closure of $\langle A \rangle$. Prove that
$\dim_\C {\cal G}_A> 1$.

\item
Let $\alpha_1, ..., \alpha_n$ be complex
numbers such that the numbers $\log |\alpha_i|$
are linearly independent over $\Q$, and
$A\in \End(\C^n)$  a diagonal endomorphism
with eigenvalues $\alpha_1, ..., \alpha_n$.
\begin{enumerate}
\item Prove that the Zariski closure ${\cal G}_A$ 
of $\langle A \rangle$
is $(\C^*)^n$. 

\item Construct an example of a diagonal $A\in \End(\C^n)$
with eigenvalues $\alpha_1, ..., \alpha_n$
such that the numbers $\log |\alpha_i|$
are linearly dependent over $\Q$, but
${\cal G}_A\cong (\C^*)^n$.
\end{enumerate}

\item Let $A \in \End(\C^n)$ be a linear contraction, 
and ${\cal G}_A$ the Zariski closure of $\langle A \rangle$. 
Find a contraction $A$ such that  ${\cal G}_A$ does not contain $\C^* \cdot \Id$.

{\em Hint:} Use \ref{_Chevalley_reductive_Theorem_}.

\item
Find a linear contraction
$A \in \End(\C^n)$ such that
${\cal G}_A$ is disconnected.

\item
Construct a Hopf manifold which 
admits a holomorphic connection with non-zero torsion.

\item
Let $M$ be a compact complex manifold.
Prove that a meromorphic function on $M$ is
uniquely defined (up to a constant multiplier) 
by its pole and zero divisors.

\item
Let $M$ be a compact Vaisman manifold,\index[terms]{manifold!Vaisman}
and $\mu$ a meromorphic function on $M$.
Prove that $\Lie_{\theta^\sharp}\mu = \const \mu$,
where $\theta^\sharp$ is the Lee field\index[terms]{Lee field} of $M$.

{\em Hint:} Use the previous exercise.

\item
Let $\mu$ be a meromorphic function on a complex manifold $M$,
and $D$ its pole divisor, considered to be  an effective
divisor on $M$. Denote by $\calo(D)$ the line
bundle whose sections are meromorphic functions
with pole divisor $D_1$ such that $D-D_1$ is effective.
Then $\mu$ is realized as a holomorphic section of
$\calo(D)$. Two divisors $D_1, D_2$ are called
{\bf rationally equivalent} if
$D_1$ is the pole divisor of a 
meromorphic function, and $D_2$ is its zero divisor.
Prove that ${\Bbb P}H^0(M, \calo(D))$
is identified with the set of all
divisors rationally equivalent to $D$.

\item\label{_mero_fraction_sections_Exercise_}
Prove that any global meromorphic function on $H$
can be obtained as $\frac f g$, where $f, g$ are
sections of the same holomorphic line bundle.

{\em Hint:} Use the previous exercise.

\item
Let $H= \frac{\C^n\backslash 0}{\langle A \rangle}$
be a Hopf manifold, and ${\cal G}_A\subset \GL(n,\C)$
the Zariski closure of $\langle A \rangle$.
\begin{enumerate}
\item
Let $L$ be a holomorphic line bundle on $H$.
Prove that there exists $\alpha\in \C^*$
and a holomorphic flat connection on $L$
with the monodromy action\index[terms]{action!monodromy} equal to multiplication
by $\alpha$. Prove that this correspondence
identifies the group $\Pic(H)$ of holomorphic line bundles
and $\C^*$.

\item
{\bf An automorphic line bundle} on $\C^n$
is a trivial line bundle on $\C^n$ equipped
with an equivariant action of $\langle A \rangle$
that multiplies a constant section by a number $\alpha$.
Prove that the action of $\langle A \rangle$
on $H^0(\C^n, 0)$ can be extended to the
action of ${\cal G}_A$ on $H^0(\C^n, 0)$.

\item Prove that any holomorphic line bundle $L$ on $H$ is 
equipped with a natural ${\cal G}_A$-equivariant
structure. Prove that any section of $L$ is
${\cal G}_A$-invariant.

\item
Using Exercise \ref{_mero_fraction_sections_Exercise_},
prove that any meromorphic function on $H$
is ${\cal G}_A$-invariant.
\end{enumerate}

\item\label{_orbit_space_and_a(H)_Exercise_}
Let $H= \frac{\C^n\backslash 0}{\langle A \rangle}$
be a Hopf manifold, and ${\cal G}_A\subset \GL(n,\C)$
the Zariski closure of $\langle A \rangle$.
\begin{enumerate}
\item 
Let $\mu$ be a meromorphic function 
on $\C^n$ that is  constant on orbits
of ${\cal G}_A$. Prove that
$\mu$ is the pullback of a meromorphic
function on $H$. 

\item
Prove that the field of global meromorphic
functions on $H$ is isomorphic
to the field of ${\cal G}_A$-invariant
meromorphic functions on $\C^n$.

\item
Assume that the group $\C^* \cdot \Id$ of scalar matrices
belongs to ${\cal G}_A$. Suppose that the general 
orbit of ${\cal G}_A$ on ${\Bbb P}\C^n$ has dimension $k$.
Prove that $a(H)= n-1-k$.
\end{enumerate}

\item\label{_zari_closure_diagonal_Exercise_}
Let $A$ be a diagonal contraction,
with the eigenvalues $a_1, a_1, ..., a_1, a_2, ...,$ $a_{n-k+1}$,
where $a_1$ is taken $k$ times, and the rest of the eigenvalues once,
written in a basis $e_1, ..., e_n$.
Assume that the numbers $\log |a_i|$ are linearly
independent over $\Q$. Prove that 
${\cal G}_A= (\C^*)^{n-k+1}$, acting 
on $\C^n$ with $\C^*\Id_k$ acting as $\C^*  \cdot\Id$ 
on $\langle e_1, ..., e_k\rangle$, with the rest
of $\C^*$ factors acting on each $e_i$, $i = k+1, ..., n$.

\item\label{_alge_dime_Hopf_Exercise_}
Prove that  $a(H)$ 
can take any values between 0 and $\dim(H)-1$,
when $H$ is a diagonal Hopf manifold.

{\em Hint:} Use 
Exercise \ref{_orbit_space_and_a(H)_Exercise_} and
Exercise \ref{_zari_closure_diagonal_Exercise_}.

\medskip

\definition
Let $M$ be a compact complex manifold.
The {\bf \index[terms]{dimension!Kodaira} Kodaira dimension} $\kappa(M)$  
is defined as $\kappa(M):=\limsup_n \frac {\log(\dim H^0(K_M^n))} {\log n}$.
If the function $n \mapsto \dim H^0(K_M^n)$
grows as a polynomial of degree $d$, the Kodaira dimension
of $M$ is $d$. It is equal to $-\infty$ if $H^0(K_M^n)=0$
for almost all $n>0$, and to 0 if $\dim H^0(K_M^n)$ is bounded, but non-zero
for infinitely many $n$.

\item\label{_Kodaira_dimension_Vaisman_}
Let $M$ be a Vaisman manifold, \index[terms]{manifold!Vaisman}
$G$ the closure of the group generated by the
Lee and the anti-Lee action, $\tilde M$ 
a K\"ahler $\Z$-covering of $M$, and $\tilde G$ the lift
of $G$ to $\tilde M$.
\begin{enumerate}
\item Prove that there exists $\gamma\in \tilde G$
such that the quotient $M' := \frac{\tilde M}{\langle \gamma\rangle}$
is Vaisman. Prove that $\gamma$ can be chosen in such a way that 
$M'$ is quasi-regular.\index[terms]{manifold!Vaisman!quasi-regular}
\item Prove that $\dim H^0(M',K_{M'}^d)\geq \dim H^0(M,K_{M}^d)$
for all $d\in \Z$.
\item Assume that $M'$ is quasi-regular, and 
let $X$ be the space of leaves of its canonical
foliation, considered to be  a projective orbifold.
Prove that $\kappa(M) \leq \kappa(X)= \kappa(M')$.
\end{enumerate}

{\em Hint:} Apply
\ref{_holo_tensor_on_Vaisman_Lee_invariant_Theorem_}
to sections of $K^d_M$ considered to be  holomorphic
tensors on $M$.

\end{enumerate}

{\setlength\epigraphwidth{0.9\linewidth}
\partepigraph{
\it \qquad Russian readers have grown used to emasculated words, and they
see in them algebraic symbols that mechanically solve the problems of
petty thinking, meanwhile everything alive, superconscious in the
word, everything which connects it to its wellsprings, its sources of
existence, goes unnoticed.\\ \smallskip

\qquad
Art can be concerned only with the living, it cannot bother with the
dead! \\ \smallskip

\qquad
Like warriors on a foggy morning, we have attacked by surprise our
idle enemies-and now they, to the amusement of the victors and of
the whole world, kick each other, pull each other's hair, and all they
can throw at us is dirt and abuse.\\ \smallskip

\qquad Look, you thick-lipped ones!\\ \smallskip

\qquad 
Before us there was no verbal art}
{\sc\scriptsize New Ways of the Word
(the language of the future,
death to Symbolism)
A. Kruchenykh}

\part{Topics in locally conformally K\"ahler geometry}

\removeepigraph
}

\chapter{Automorphism groups of LCK manifolds}\index[terms]{geometry!K\"ahler}

{\setlength\epigraphwidth{0.6\linewidth}
	\epigraph{\em The end is in the beginning and yet you go on.}{\sc\scriptsize Samuel Beckett, Endgame}}


 There are various kind of transformations that can be
 considered  on a given  LCK manifold --- diffeomorphisms
 compatible with one or several of the structures
 associated with the LCK one. In this chapter, we define
 them and study their main properties and discuss the
 inclusion relations among them.

 \smallskip
 
 Let $(M,I,\omega,\theta)$ be an LCK manifold and let $c$
 denote the conformal class of the given LCK metric,
 $\nabla$ (respectively $D$)  its Levi--Civita
 (respectively Weyl) connection. The groups we shall deal
 with are:
 
 \begin{itemize}
\item $\Aff(M,D)$, the group of affine transformations, i.e.  preserving the Weyl connection. 
\item ${\cal H}(M,I)$, the group of biholomorphisms with respect to $I$.
\item $\Conf(M,c)$, the group of conformal transformations. 
\item $\Aut(M)$:= ${\cal H}(M,I)\cap \Conf(M,c)$,
  the group of LCK automorphisms.
\end{itemize}

The groups $\Aff(M,D)$, ${\cal H}(M,I)$, $\Conf(M,c)$, $\Aut(M)$
 are closed Lie subgroups of the group $\Diff(M)$. 
When $M$ is compact, 
$\Aut(M)$ is compact too (\cite{kor}). Denote with  $\mathfrak{aff}(M,D)$, $\mathfrak{h}(M,I)$,
$\mathfrak{conf}(M,c)$, $\mathfrak{aut}(M)$ their
respective Lie algebras. Vector fields from
$\mathfrak{conf}(M,c)$
are called {\bf conformal vector fields}.\index[terms]{vector field!conformal}

\hfill

The Weyl connection is the Levi--Civita connection on the K\"ahler
cover of an LCK manifold, and the K\"ahler metric on the cover is defined
uniquely up to a constant multiplier. Then the Weyl connection is uniquely determined
 by the LCK structure.\index[terms]{structure!LCK} This implies that 
$\Aut(M)\subset\Aff(M,D)$.


\section{Infinitesimal automorphisms}
The results in this paragraph are important for their own
sake, and will be useful when discussing twisted 
Hamiltonian actions and toric LCK manifolds in Section
\ref{toric_lck}.

\hfill

\remark
Any LCK automorphism of an LCK manifold $(M,I, \theta)$ defines
a conformal holomorphic automorphism of its K\"ahler cover. However,
a conformal map of K\"ahler manifolds in $\dim_\C > 1$ is a homothety,
that is, multiplies the K\"ahler metric by a constant. This defines
a character on the group $\Aut(M)$. The corresponding
Lie algebra homomorphism is defined  in \ref{lema_mmp} below.

\hfill

\proposition\label{lema_mmp} (\cite{va_lcs}, \cite{_Madani_Moroianu_Pilca:toric_}) 
Let $(M, \omega, \theta)$ be an LCK manifold.
Then there exists a morphism of Lie algebras $\sigma:\mathfrak{aut}(M)\ra \R$ such that:
\begin{align}
d_\theta(i_X\omega)=\sigma(X)\omega\label{primo},\\
\Ll_Xg=(\theta(X)+\sigma(X))g\label{secondo}.
\end{align}
for each $X\in \mathfrak{aut}(M)$.

\hfill

\proof
Locally, we may assume that $X$ can be
integrated to a conformal automorphism $e^{tX}$ of $M$.
Let $\tilde M$ be its K\"ahler cover, and
$\tilde X$ the vector field $X$ lifted to $X$.
Since $e^{t\tilde X}$ is a holomorphic automorphism,
it maps the K\"ahler form\index[terms]{form!K\"ahler} $\tilde \omega$ to a K\"ahler form;
since it is conformal, it actually
maps $\tilde \omega$ to $f\tilde \omega$,
for some $f\in C^\infty \tilde M$.
Since $f\tilde \omega$ is K\"ahler, 
$d(f\tilde \omega)=0$, and hence  $f$ is constant
(see also \eqref{_conformal_then_homothety_Equation_}).
This implies that $e^{t\tilde X}$ is a homothety,
and $\Lie_{\tilde X} \tilde \omega= \lambda_X
\tilde\omega$, where $\lambda_X$ is a constant.
We put $\sigma(X):= \lambda_X$.
The Cartan formula\index[terms]{Cartan formula} for $\Lie_{\tilde X}$ gives
$d(i_{\tilde X} \tilde \omega) = \sigma(X)\tilde
\omega$, which implies \eqref{primo} and  \eqref{secondo}.
\endproof

%

\hfill

\remark The map $\sigma$ was introduced by Vaisman in \cite{va_lcs} in the context of LCS geometry, under the name of ``Lee homomorphism'', and denoted $l$. He used it to define LCS manifolds ``of the first kind'', a notion widely studied afterwards.\index[terms]{manifold!LCS!of the first kind}\index[terms]{Lee homomorphism}\index[terms]{geometry!LCS}

%

\hfill

\definition\label{special_autom} Let $(\tilde M,\tilde g)$
be the K\"ahler universal cover of $M$. An automorphism $f\in
\Aut(M)$ is called {\bf special} (\cite{_Madani_Moroianu_Pilca:toric_}) if it lifts to an
isometry of $(\tilde M,\tilde g)$. Let $\Aut_s(M)$
be the subgroup of special automorphisms, and
$\mathfrak{aut}_s(M)$  its Lie algebra. 

\hfill

\proposition (\cite{_Madani_Moroianu_Pilca:toric_})\label{aut_char} Let $(M,I, c)$ be a compact LCK manifold. Let $g_0\in c$ be the Gauduchon metric,\index[terms]{metric!Gauduchon} and $\theta_0$ the corresponding Lee form\index[terms]{form!Lee}. Then:

(i) Every $X\in \mathfrak{aut}(M)$ is $g_0$-Killing.

(ii) $\mathfrak{aut}_s(M)\subseteq\ker\theta_0$.

(iii) The lift $\tilde X$ (resp. $\hat X$) of an $X\in \mathfrak{aut}_s(M)$ to the universal cover $\tilde M$ (resp. minimal cover $\hat M$) is Killing\index[terms]{vector field!Killing} with respect to the K\"ahler metric $\tilde g$ (resp. $\hat g$). \index[terms]{cover!minimal}

\hfill

\proof Let $\f\in\Aut(M)$. Then $[\f^*g_0]= c$ and
$\f^*g_0$ is a Gauduchon metric, too. By uniqueness of the
Gauduchon metric of a conformal class on a compact
manifold, $\f^*g_0=ag_0$, with $a\in\R^{>0}$; hence
$\f$ is a homothety of scale factor $a$. As $M$ is
compact, $a=1$ and $\f$ is an isometry of $g_0$. In
particular, the flow of $X$ will consist in isometries of
$g_0$, thus $X$ is $g_0$-Killing,\index[terms]{vector field!Killing} proving (i). Now
\eqref{secondo} proves (ii), while  (iii) follows directly
from \ref{special_autom}. \endproof

\section[Lifting a transformation group  to a K\"ahler co\-ver]{Lifting a transformation group  to a K\"ahler\\ co\-ver}

Let $\Gamma$ be a subgroup of $\pi_1(M)$,
$\ker(\chi)\subseteq \Gamma\subseteq \pi_1(M)$,  where
$\chi:\pi_1(M)\ra\R^{>0}$ is the character in
\ref{_homothety_character_Definition_}, and
$M_\Gamma$ the corresponding cover
(\ref{_Galois_covers_Theorem_}).
Then  $M_\Gamma$ is K\"ahler.

Denote by $\tilde M$ the universal cover and $\hat M$ the
minimal cover, that is, the cover corresponding to
$\ker(\chi)$. Since $\ker(\chi)$ is a normal subgroup
of $\pi_1(M)$, the cover $\hat M$ is a Galois cover,\index[terms]{cover!Galois}
and its automorphism group over $M$ is $\pi_1(M)/\ker(\chi)$.


Clearly, any diffeomorphism of $M$ can be lifted to
$\tilde M$,  but not necessarily to all $M_\Gamma$. 
Instead, the minimal cover plays a special role:

\hfill

\proposition (\cite{gopp})\label{lift_to_min} Any
$f\in\Aut(M)$ lifts to a
$\pi_1(M)/\ker(\chi)$--equi\-va\-riant holomorphic
homothety $\hat f$ of $\hat M$. 

\hfill

\proof
The subgroup $\ker(\chi)\subset \pi_1(M)$ 
is uniquely determined by the LCK structure\index[terms]{structure!LCK}; hence 
any automorphism $f\in\Aut(M)$ preserves 
$\ker(\chi)\subset \pi_1(M)$.

The homotopy lifting lemma allows to lift $f$ to a
$\pi_1(M)$-equivariant automorphism of $\tilde M$.
However, $\hat M = \frac{\tilde M}{\ker(\chi)}$.
Since the action of $f$ on $\pi_1(M)$ preserves
$\ker(\chi)$, its lift to $\tilde M$\index[terms]{lemma!homotopy lifting}
commutes with the $\ker(\chi)$-action,
hence induces an automorphism of $\hat M = \frac{\tilde M}{\ker(\chi)}$.
\endproof

\hfill

\remark \label{homo_con} Now suppose $M$ is a rank 1  compact Vaisman manifold.\index[terms]{manifold!Vaisman} Then $\tilde M$ is a cone $C(W)$ over a compact Sasakian manifold $W$. Any $f\in \Aut(M)$ can then be lifted to a biholomorphic homothety $\tilde f$ of $C(W)$ which, by Exercise \ref{iso_con}, is of the form $(\psi_f, a\cdot \id)$, with $\psi_f\in\Iso(W)$, and $a=\chi(\tilde f)$. In particular:

\hfill

\corollary (\cite{gopp}) \label{isom_con} For a compact Vaisman manifold
$M$ of rank 1 , with associated compact Sasakian manifold
$W$, then $\Iso(\hat M)\simeq \Iso(W)$.
\endproof

\section{Affine vector fields}

The following are local results. 

\hfill

\lemma  \label{_Aff_then_holo_Lemma_}
{(\cite{mo})} 
 Let $M$ be an LCK manifold such that the holonomy
 group\index[terms]{group!holonomy} of the Weyl connection
 $\Hol_0(D)$ is irreducible. Then:
 \begin{description}
 \item[(i)] Any $f\in  \Aff(M,D)$ is conformal, and
 
 \item[(ii)] If $M$ is not locally conformally\index[terms]{manifold!LCHK}
hyperk\"ahler (LCHK)\footnote{See Section \ref{lchk} for details about LCHK manifolds.}, then  
any transformation $f\in  \Aff(M,D)$ is holomorphic or
anti-holomorphic. 

\item[(iii)] If $M$ is LCHK, then there exists a
  quaternion 
$h\in \Sp(1)=\U({\Bbb H}, 1)$ such that $f^*(L)=hLh^{-1}$, where
$L$ is any of  the complex structures $I,J,K$.
\end{description}

\proof We prove (ii) and (iii), the proof of (i) is very
similar. This proof is different from the one given in
\cite{mo}.

Let $F:\; M \arrow M$ be an affine map, and $dF:\; T_xM
\arrow T_{F(x)}M$ its
differential. Then $dF$ normalizes the holonomy group of
the Weyl connection, that is: $dF(\Hol_x(M))=\Hol_{F(x)}M$.
Denote the group of conformal isometries of $W=\R^n$ by $\R^+\OO(W)$.
Identifying $W=T_x M$ and $T_{F(x)}M$ in a
$\Hol(M)$-invariant way, we may consider $dF$ as
a map $A\in \R^+\OO(W)$. Then $A$ normalizes the group $\Hol(M)\subset
\GL(W)$, that is, satisfies $A \Hol(M) A^{-1}=\Hol(M)$.

To advance the proof we need to classify the normalizers\index[terms]{normalizer}
for all possible holonomy groups\index[terms]{group!holonomy}, that are  $\U(n)$,
$\SU(n)$ and $\Sp(n)$ by \index[persons]{Berger, M.} Berger's classification
of irreducible holonomies. All the automorphisms of $\U(n)$ and
$\SU(n)$ are interior or interior composed with the complex
conjugation (\cite{besse}). Denote by $\CU(n)$ the group generated by 
$\U(n)$ and the complex conjugation; it consists of two
connected components, both isomorphic to $\U(n)$, and
acts on $\U(n)$ by automorphisms.

The automorphisms of $\Sp(n)$ are interior ones 
composed with the action by a qua\-ter\-nion $h\in \Sp(1)$ (\cite{besse}).

Any endomorphism of a vector space
that acts trivially on $\U(n)$, $\SU(n)$, $\Sp(n)$ 
belongs to $\R^+$, because the
adjoint representation of these groups is irreducible.
Therefore, the normalizer group of $\U(n)$ and $\SU(n)$
is $\R^+\CU(n)$ and the normalizer of $\Sp(n)$ is 
$\R^+ \Sp(1) \Sp(n)$.\index[terms]{normalizer}

In the case when $\Hol(M)=\U(n)$ or $\SU(n)$, 
one has $dF\in \R^+ \U(W)$, and in this case, $F$ is
holomorphic, or $dF$ is a composition of $A\in \R^+ \U(W)$
and a complex conjugation, in which case $F$ is
antiholomorphic; this proves (ii).

If $\Hol(M)=\Sp(n)$ ,
one has $dF\in \R^+ \Sp(1)\Sp(n)$; this proves (iii).
\endproof

\hfill

\remark
The assumption of Weyl  irreducibility
(\ref{_Aff_then_holo_Lemma_}) is not very restrictive.
Indeed, a Vaisman manifold\index[terms]{manifold!Vaisman} is Weyl irreducible
unless it is Hopf (this follows from \index[persons]{Gallot, S.} Gallot's theorem,
\ref{rem_mo}).

%
%

\hfill

\corollary \label{af=aut}
 Let $M$ be an LCK manifold that is  not LCHK and such that
$\Hol_0(D)$ is irreducible. Then  $\mathfrak{aff}(M,D)=\mathfrak{aut}(M)$. 
\endproof

\hfill

\remark \label{rem_mo}
  We still do not know whether the assumption of irreducibility of
  the Weyl connection can be relaxed or not, at least on compact
  $M$. But we  know that \emph{a compact Vaisman manifold,\index[terms]{manifold!Vaisman} which
    is not a (diagonal) Hopf manifold, is Weyl-irreducible}. Indeed,
  by the \ref{str_vai}, $M$ is a mapping torus\index[terms]{mapping torus} of a
  Sasakian isometry and its universal cover is a K\"ahlerian cone over\index[terms]{isometry!Sasakian}
  the compact, and thus  complete, Sasakian fibre. Then, as $D$ is,
  locally, the Levi--Civita connection of the local K\"ahler metrics,
  if $D$ is reducible, the K\"ahler metric of the covering cone is
  reducible. Now, by  \cite[Proposition 3.1]{gal}, a cone over a
  complete manifold is reducible if and only if it is flat. As  a flat
  cone is the cone over a sphere,  $M$ is a Hopf manifold.  

\hfill

\remark Let $W$ be a 3-Sasakian manifold\footnote{See
    \cite{bog} for the definition of 3-Sasakian manifolds.}. Then $M:=W\times S^1$ is LCHK
(for $W=S^{2n-1}$, $M$ is a Hopf manifold), equipped with the three Killing\index[terms]{vector field!Killing}
fields on $W$ generating the $\SU(2)$
action. The corresponding action is affine with respect to
the Weyl connection on $M$. Each of them is holomorphic
with respect to one and only one of the three complex
structures (see {e.g.} \cite[Chapter 11]{do}).

\section{Conformal vector fields on compact LCK ma\-ni\-folds}

The following theorem is well known.

\hfill

\theorem  \label{_Lichnerowicz_conformal_Theorem_}
{(\index[persons]{Lichnerowicz, A.}Lichnerowicz, \cite[\S 90]{li})}  \\
Let $(M,I,g)$ be a compact K\"ahler manifold. Then any $[g]$-conformal vector field $\xi$ is Killing.\index[terms]{vector field!Killing}\index[terms]{theorem!Lichnerowicz}

\proof Follows immediately from the maximum principle
(\ref{hopf_theorem}) and
\ref{_conf_on_K_harm_div_} below.
\endproof

%
%

\hfill

Recall that on a compact K\"ahler manifold $(X,I,h,F)$, a Killing field $\xi$ is holomorphic.\index[terms]{vector field!Killing}\index[terms]{vector field!holomorphic} Indeed, being $h$-parallel, the K\"ahler form\index[terms]{form!K\"ahler} $F$ is harmonic. Since the Lie derivative commutes with the Laplacian, we find that also $\Lie_\xi F$ is harmonic. Note that by the Cartan formula, $\Lie_\xi F$ is exact; hence, by Hodge theory, it vanishes. Now $\Lie_\xi F=0$ together with $\Lie_\xi h=0$ imply $\Lie_\xi I=0$, which proves the statement.

\index[persons]{Lichnerowicz, A.} Lichnerowicz's theorem then implies that on a compact K\"ahler manifold  a conformal vector field $\xi$ is holomorphic.

This result was extended to the LCK context by \index[persons]{Moroianu, A.} Moroianu and \index[persons]{Pilca, M.} Pilca in the following form:

\hfill

\theorem {(\cite{_Moro_Pilca_})}
Let $(M,I,\omega,\theta)$ be a compact LCK manifold of complex dimension $n$. Then:
\begin{description}
	\item[(i)] Any $[g]$-conformal vector field is Killing\index[terms]{vector field!Killing}\index[terms]{vector field!conformal} with respect to the Gauduchon metric $g_0\in[g]$. 
	\item[(ii)] If $(M,I,\omega,\theta)$ is neither Hopf nor locally conformally hyperk\"ahler, any $[g]$-conformal vector field is holomorphic.\index[terms]{vector field!holomorphic}
\end{description}

The proof uses methods of conformal and Riemannian geometries. We shall explain the main steps without entering the technical details.

For simplicity, let $c$ denote the conformal class $[g]$. Then, on the K\"ahler universal cover $(\tilde M,I,\tilde g,\tilde\omega)\stackrel{\pi}{\arrow} M$, we have $\tilde g=e^{-\phi}\pi^*g\in\tilde c$, where   $\pi^*\theta=d\phi$. 

Let $\xi$ be a $c$-conformal vector field on $M$ and $\tilde\xi$ its lift to $\tilde M$. Clearly $\tilde \xi$ is $\tilde c$-conformal. We are thus led to understand conformal vector fields on K\"ahler manifolds (not necessarily compact). 

\hfill

\claim\label{_conf_on_K_harm_div_}
Let $(N,I,h)$ be a K\"ahler manifold of complex dimension
at least 2 and $V$ an $[h]$-conformal vector field. Then
$\div_h V$ is a harmonic function with respect to $h$.   

\hfill

\proof
Let $\alpha:= V^\flat$ be the 1-form dual to $V$.
By \ref{_Killing_via_nabla_Proposition_},
\[ \Lie_V(h)(X,Y)= h(\nabla_X V, Y) +  h(\nabla_Y V, X),\]
and $(\nabla_X V)^\flat= \nabla_X \alpha$.
Therefore,
\begin{equation}\label{_Lie_of_metric_Equation_} 
\Lie_V(h)= \Sym(\nabla(\alpha)).
\end{equation}

Since $V$ is conformal, $\Lie_V h= \mu h$,
where $\mu$ is a function. By \eqref{_Lie_of_metric_Equation_},
$\nabla(\alpha) = \Lie_V h + d\alpha= \mu h + d\alpha$.
Since $\nabla(I)=0$, we have 
$\nabla(I\alpha) = \Id\otimes I(\nabla(\alpha))\in
\Lambda^1N\otimes \Lambda^1 N$. Then 
\begin{multline*}
 (d^c(\alpha))^{1,1}= I^{-1}(\Alt(\nabla(I\alpha)))^{1,1}= -(\Alt(\nabla(I\alpha)))^{1,1}= \\
= \Alt(\Id\otimes I(\nabla(\alpha)))^{1,1}= \mu\omega - \1
(d\alpha)^{1,1}.
\end{multline*}
(the last equation is directly implied by $\nabla(\alpha)=\mu h + d\alpha$).
This gives 

\begin{equation}\label{_dd^c_to_three_components_Equation_}
dd^c((d^c(\alpha))^{1,1}+ \1dd^c(d\alpha)^{1,1}=
dd^c\mu \wedge \omega.
\end{equation} Since $dd^c$ maps (1,1)-forms to (2,2)-forms,
one has $dd^c \Theta^{1,1}= (dd^c \Theta)^{2,2}$ for any
2-form $\Theta$. Then
\eqref{_dd^c_to_three_components_Equation_}
gives
\[
dd^c\mu \wedge \omega= (dd^cd^c\alpha)^{2,2} + \1(dd^cd\alpha)^{2,2}=0.
\]
On the other hand, by the K\"ahler identities (\ref{_kah_susy_Theorem_}) we derive  
$dd^c\mu\wedge \omega^{\dim_\C N -1}= \Delta\mu\Vol_h$;  
hence $dd^c\mu\wedge \omega=0$ gives $\Delta\mu=0$
if $\dim_\C N >1$.
\endproof

\hfill

\ref{_conf_on_K_harm_div_}
implies an LCK analogue of the Lichnerowicz
theorem (\ref{_Lichnerowicz_conformal_Theorem_}).
Originally it was proven in \cite{_Moro_Pilca_}
using a different approach.\index[terms]{theorem!Lichnerowicz} 

\hfill

\corollary
(\cite[Theorem 5.1]{_Moro_Pilca_})\label{_Conformal_is_Killing_Corollary_}
Let $\xi$ be a conformal vector field on a compact
LCK manifold $(M, I, g, \theta)$. Then $\xi$ is Killing with respect to\index[terms]{vector field!Killing}
the Gauduchon metric.\index[terms]{metric!Gauduchon}

\hfill

\proof Let $\tilde M\stackrel \pi\arrow M$ be a K\"ahler
cover of $M$, and $\tilde g$ its Riemannian bilinear symmetric form. 
Slightly abusing the notation, 
use the same letter $\xi$ to denote the lift of $\xi$
to $\tilde M$. Let $\mu$ be the function on $\tilde M$
that satisfies $\Lie_\xi \tilde g = \mu\tilde g$, and
$\psi$ a function on $\tilde M$ such that $\pi^* g =
\psi \tilde g$.
Then 
\begin{multline*}
\Lie_\xi(\pi^* g) = \Lie_\xi(\psi \tilde g) =
\psi \Lie_\xi \tilde g + \Lie_\xi(\psi) \tilde g =\\=
\psi \mu \tilde g + \frac{\Lie_\xi(\psi)}{\psi} \pi^* g=
\mu \pi^* g + \Lie_\xi(\log \psi) \pi^* g.
\end{multline*}
Since the deck transform map multiplies $\psi$ by a
constant, the 1-form $d\log \psi= -\theta$ is the pullback
of a 1-form on $M$. This gives $\pi^* \Lie_\xi g =
\pi^*(\langle \xi, \theta\rangle g) +  \mu \pi^* g$,
thus $\mu$ is invariant under the deck transform,
giving $\mu = \pi^* \mu_1$. Consider the 
differential operator $\Delta_1:\; C^\infty M \arrow C^\infty M$
taking $f\in  C^\infty M$ to $\psi\Delta_{\tilde g} \pi^*f$.
Since the deck transform multiplies $\tilde g$ and $\psi^{-1}$
by the same constant, the function $\psi\Delta_{\tilde g} \pi^*f$
is deck invariant, and hence  it is the pullback of a function $\Delta_1(f)$ on $M$.
Since the symbol of $\Delta_{\tilde g}$ is equal to $\tilde g$,
the symbol of $\Delta_1$ is equal $g$. Therefore, it is
is elliptic. Since $\Delta_1$ vanishes on constants, it
satisfies the assumptions of the maximum principle (\ref{hopf_theorem}). 
Then \ref{hopf_theorem} implies that any function in
$\ker \Delta_1$ is constant. By \ref{_conf_on_K_harm_div_},
$\mu_1 \in \ker \Delta_1$, thus  it is constant.

Since $u_t:=e^{t\xi}$ acts on $\tilde M$ by homotheties,
this diffeomorphism preserves 
 the Weyl connection $D$ on $M$ induced by the
Levi--Civita connection on $\tilde M$.

Let $\theta_t$ be the Lee form\index[terms]{form!Lee} of the metric
$g_t:= u_t^*g$, conformal to $g$, and $d^*_{g_t}$ be the codifferential
associated with $g_t$. Since $D(g_t)= g_t \otimes \theta_t$,
and $u_t^*(D(g)) = D(g_t)$, this implies $u_t^*\theta= \theta_t$.
Since $u_t:\; (M, g_t) \arrow (M,g)$ 
is an isometry, $d^*_{g_t}(\theta_t)= 0$ whenever $d^*\theta=0$.
We obtain that $g_t$ is a Gauduchon metric on $(M, I)$
if $g$ is Gauduchon; since the Gauduchon metric is unique
up to a constant multiplier, this implies that $g=g_t$.
This constant multiplier is 1, because the Riemannian volume
is diffeomorphism invariant.
Therefore, $\xi$ is Killing.\index[terms]{vector field!Killing}
\endproof

\hfill

Assume now that $M$ is neither a Hopf manifold, nor
locally conformally hyperk\"ahler. Then its universal
cover is neither flat, nor Calabi--Yau (with respect to its
K\"ahler structure). It is immediate that the vector field
$\xi$ is holomorphic if and only if its lift $\tilde\xi$
is holomorphic. Hence, we are led to consider a conformal
vector field on the non-compact K\"ahler manifold $(\tilde
M,  \tilde g, \tilde\omega)$. 

Recall that $\tilde\xi$ is $\tilde g$-homothetic, and hence  it
is affine with respect to the Levi--Civita connection
$\nabla^{\tilde g}$. Two cases occur: $\Hol(\tilde g)$ can be irreducible, or reducible.

If $\Hol(\tilde g)$ is irreducible,
 $\tilde\xi$ (\ref{af=aut}) is holomorphic,
finishing the proof.

The second case, $\Hol(\tilde g)$ reducible is much more
involved. We need a  result of \index[persons]{Kourganoff, M.} Kourganoff. Before stating
it, recall that a {\bf similarity structure} on a
connected manifold $X$ is a metric $h$ such that
$\pi_1(M)$ acts on the universal cover $(\tilde X, \tilde
h)$ by homotheties with strictly positive scale factor
(see \cite{_Kourganoff_,_Vaisman_Reischer_}). LCK
geometry is thus a special case of similarity geometry and
the following result can be applied:

\hfill

\theorem {(\cite[Theorem 1.5]{_Kourganoff_})}  Consider a compact manifold $M$ endowed with a similarity structure,
and its universal cover $\tilde M$ equipped with the corresponding Riemannian structure $\tilde g$.
Assume that $M$ is not globally Riemannian, i.e.  $\pi_1(M)$ is not a subgroup of $\Iso(\tilde M)$.
Then we are in exactly one of the following situations:
\begin{description}
\item[(i)] $\tilde M$ is flat.
\item[(ii)] $\tilde M$ has irreducible holonomy and $\dim(\tilde M) \geq 2$.
\item[(iii)] $\tilde M=\R^q\times N$, where $q\geq 1$, $\R^q$ is the Euclidean space, and $N$ is a non-flat,
non-complete Riemannian manifold that has irreducible holonomy. \endproof
\end{description} 

Hence, in our case, we are in the third situation:
$(\tilde M, \tilde g)=(\R^q, g_{\text{{\tiny \sf
      flat}}})\times (N,g^N)$. Now, a very detailed
analysis of homothetic flows on such a product manifold
proves that $\tilde\xi$ has to be constant on the flat
factor, and thus $\tilde \xi$ is tangent to $N$
(\cite[Proposition 6.1]{_Moro_Pilca_}). Call $\zeta$ the
restriction of $\tilde\xi$ to $N$.

If, by absurd, $\zeta$ is not holomorphic, again by
\ref{af=aut} we conclude that $(N,g^N)$ is hyperk\"ahler,
and hence Ricci flat, thus $(\tilde M,I,\tilde g)$ is
 Calabi--Yau. This means that $(M,I,[g])$ is Einstein--Weyl,\index[terms]{manifold!Einstein--Weyl}
thus Vaisman with respect to the Gauduchon metric
$g_0\in[g]$ (\ref{ewcy}), in particular the vector field
$\theta_0^\sharp$ that is  $g_0$-equivalent to
$\theta_0$  is Killing and\index[terms]{vector field!Killing} its lift
$\tilde\theta_0^\sharp$ is homothetic on $\tilde M$.   
\cite[Proposition 6.1]{_Moro_Pilca_} implies now
that $\tilde\theta_0^\sharp$ is constant on $\R^q$, thus 
the homothety factor vanishes. Since the homothety
factor is $2\Vert\tilde\theta_0^\sharp\Vert_{\tilde g}$,
this yields that $\tilde\theta_0^\sharp=0$,
contradiction. \endproof

\hfill

\remark  \ref{_Conformal_is_Killing_Corollary_}
superseeded a previous attempt in the particular case of
Vaisman manifolds.\index[terms]{manifold!Vaisman} In \cite{mo} it was
proven that every conformal vector field on the mapping
torus over a compact Riemannian manifold is Killing.\index[terms]{vector field!Killing}\index[terms]{vector field!conformal} This
showed that on compact Vaisman manifolds of LCK rank 1 \index[terms]{rank!LCK},
that are  neither Hopf, nor locally conformally
hyperk\"ahler, conformal fields are Killing. The
restriction on the LCK rank was then eliminated in
\cite[Proposition 2.1]{_gauduchon_moroianu_} where a more
general result is proven: {\em  Let $(M,g)$ be a connected
	compact oriented Riemannian manifold of dimension $n\geq
	2$, carrying a non-trivial parallel vector field
	$T$. Let $\xi$ be any conformal\index[terms]{vector field!parallel}
	Killing vector field on $(M, [g])$. Then, $\xi$ is Killing
	with respect to g. Moreover, it commutes with $T$ and the inner product $g(\xi, T)$ is
	constant.}  

\section{Holomorphic Lee field}

Holomorphy of the Lee field\index[terms]{Lee field} is not enough  to ensure its parallelism, not even on compact manifolds. To explain this, we show how we can modify the Vaisman structure of a compact complex manifold to obtain a non-Vaisman LCK one with the same Lee field.\index[terms]{Lee field}

\hfill

\lemma (\cite{_moroianu_moroianu_ornea_}) \label{lema}  
Let $(M,g,I,\omega, \theta)$ be a compact Vaisman manifold\index[terms]{manifold!Vaisman} with $|\theta|=1$.
Let $f\in C^\infty M$ be a non-constant function 
with gradient collinear to $\ts$ (i.\, e. $df\wedge\theta=0$)
and such that $f>-1$. Define 
$\hat\omega:=\omega+f\theta\wedge I\theta$, and let $\hat
g$ be the corresponding metric. Then $(M,\hat
g,I,\hat\omega)$ is LCK, with Lee form\index[terms]{form!Lee}
$\hat\theta=(1+f)\theta$ and Lee  field\index[terms]{Lee field}
$\hat\ts=\ts$.

\hfill

\proof 
The condition $f>-1$ ensures that $\hat g$ is positive-definite. Using \eqref{_omega_via_theta_Chapter_8_Equation_} we have:
\begin{align*}
	d\hat\omega&=\theta\wedge\omega+df\wedge\theta\wedge I\theta-f\theta\wedge dI\theta\\
	&=\theta\wedge\omega+f\theta\wedge\omega=(1+f)\theta\wedge\hat\omega.
\end{align*}
Hence $(M,\hat g,J,\hat\omega)$ is LCK, with Lee form\index[terms]{form!Lee} $\hat\theta=(1+f)\theta$. Then:
\[i_{\ts}\hat\omega=I\theta+fI\theta=(1+f)I\theta= I\hat\theta=i_{\hat\ts} \hat\omega.
\]
Therefore, $\hat\theta^\sharp=\ts$ as claimed. 
\endproof  

\hfill

\remark (\cite{_moroianu_moroianu_ornea_})\label{_Example_of_good_f_Remark_} The metric $\hat g$ is not Vaisman since the $\hat g$-norm of $\ts$ is $1+f$. An example of a Vaisman manifold\index[terms]{manifold!Vaisman} and a function $f$
satisfying the hypotheses of \ref{lema} can be
constructed as follows.
Let $M$ be the diagonal Hopf manifold $S^1\times
S^{2n-1}\simeq (\mathbb{C}^n\setminus 0)/\Gamma$,
where $\Gamma=\langle z\mapsto az\rangle$, $a>1$. 
For an appropriate map $\pi:\; M \arrow S^1$,
we have $\theta = \pi^* dt$, where $t$ is the
parameter on $S^1$. Let 
$f\in C^\infty_\R M$, $f>-1$, be obtained as  pullback of 
a function on the circle. Then $df$ is proportional
to $\theta$, and it satisfies the assumptions of
\ref{lema}.

\hfill 

\remark (\cite{_moroianu_moroianu_ornea_})
Given any nowhere vanishing closed form $\theta$ on a connected smooth manifold $M$, there exists a non-constant function $f$ such that $df\wedge\theta=0$  if and only if the line spanned by the cohomology class $[\theta]\in H^1(M,\mathbb{R})$ contains an integer class.

\hfill

\remark If, in the example presented in
\ref{_Example_of_good_f_Remark_}, the LCK rank\index[terms]{rank!LCK} of $M$ is 1
and $\theta^\sharp$ is quasi-regular, then $M$ is\index[terms]{manifold!Vaisman!quasi-regular}
isometric with $N\times S^1$, where $N$ is a Sasakian
orbifold.\index[terms]{orbifold!Sasakian} Then any non-constant
function on $S^1$, bounded below by $-1$, induces on $M$ a
function $f$ satisfying the assumptions of \ref{lema}. One can use this
observation to show that, in this case, the LCK metric
corresponding to $\hat\omega$ has positive potential. The
proof consists in a careful analysis of the differential
equation satisfied by $f$, see \cite{nico2}.

\hfill

We proved that an LCK manifold with holomorphic Lee field\index[terms]{Lee field} is not necessarily Vaisman. However, adding  mild conditions to the holomorphy of the
Lee field, the Vaisman condition is still assured. 

\hfill

\theorem  (\cite{_moroianu_moroianu_ornea_}) \label{_Holo_Lee_field_main_Theorem_} Let $(M,I,g,\theta)$ be a compact LCK manifold with holomorphic Lee field\index[terms]{Lee field} $\theta^\sharp$. Suppose that one of the following conditions is satisfied:
\begin{description}
	\item[(i)] The norm of the Lee form\index[terms]{form!Lee} $\theta$ is constant, {\em or}
	\item[(ii)] The metric $g$ is Gauduchon.\index[terms]{metric!Gauduchon}
\end{description}
\noindent Then $(M,I,g)$ is Vaisman.

\hfill

\pstep Suppose the Lee field\index[terms]{Lee field} is holomorphic. Then we have the following local formula (compactness is not used at this stage) relating the Laplacian of the squared length of the Lee form, the codifferential of the Lee form\index[terms]{form!Lee} and its covariant derivative: 
\begin{equation}\label{cinci}
	\Delta|\theta|^2+\ts(|\theta|^2)+|\theta|^2d^*\theta+2|\nabla \theta|^2-\ts(d^*\theta)=0.
\end{equation}
The proof of this equation is technical and exceeds the scope of this book.

\hfill

{\bf Step 2:} From now on, we assume  $M$ is compact. 

Suppose the length $|\theta|$ is constant, let $x_0$ be a point where $d^*\theta$ reaches its maximum. As $\int_M d^*\theta \vol_g=0$, we see that $(d^*\theta)(x_0)\geq 0$.

On the other hand, the first and second terms in \eqref{cinci} vanish because of the assumption $|\theta|=\const$, and $\theta^\sharp(d^*\theta)(x_0)=0$ since $x_0$ is an extremum of $d^*\theta$. Then \eqref{cinci} implies:
$$|\theta|^2(x_0)(d^*\theta)(x_0)+2|\nabla \theta|^2(x_0)=0.$$
Both terms are positive, consequently  
$(d^*\theta)(x_0)=0$, therefore, $d^*\theta\leq 0$. Then 
$d^*\theta$ vanishes identically on $M$, since $\int_Md^*\theta=0$. Finally,
\eqref{cinci} reduces to $\nabla\theta=0$, and  \ref{_Holo_Lee_field_main_Theorem_} (i) is proven.

Now suppose the codifferential $d^*\theta$ vanishes.  
We integrate \eqref{cinci} on $M$:
\begin{equation*}
	\begin{split}
		0&=\int_M (\Delta|\theta|^2+\ts(|\theta|^2)+2|\nabla \theta|^2)\vol_g\\
		&=\int_M (|\theta|^2d^*\theta +2|\nabla \theta|^2)\vol_g=2\int_M |\nabla \theta|^2\vol_g,
	\end{split}
\end{equation*}
which again implies $\nabla\theta=0$ and proves
\ref{_Holo_Lee_field_main_Theorem_} (ii). \endproof

\hfill

\corollary (\cite{ov_surf_in_lck_pot}, \cite{nico2}) Let
$(M,I,g,\theta)$ be a compact LCK 
manifold with potential.\index[terms]{manifold!LCK!with potential} If $\theta^\sharp$ is
holomorphic, then $g$ is Vaisman.

\proof Indeed, $\theta^\sharp$ holomorphic implies
$I\theta^\sharp$ holomorphic, and then Cartan formula\index[terms]{Cartan formula}
easily gives $d_\theta(|\theta^\sharp|^2-1)=0$ which
implies $|\theta^\sharp|=\const$. Now \ref{_Holo_Lee_field_main_Theorem_}
applies. 
\endproof

\chapter{Twisted Hamiltonian actions and LCK  reduction}\label{reduction}\index[terms]{action!twisted Hamiltonian}

\epigraph{
{\cyr Ya slovo pozabyl, chto ya khotel skazat{\cprime}.\\
Slepaya lastochka v chertog tene{\u i} vernet{\cydot}sya\\
Na kryl{\cprime}yakh srezannykh, s prozrachnymi igrat{\cprime}.\\
V bespamyat{\cydot}stve nochnaya pesn{\cprime} poet{\cydot}sya.}\\ \medskip
\em
\font\tenit = cmssi17 at 10pt \tenit
(It is as if I have forgotten what I wanted to say.\\
The blind swallow will return to the hall of shades,\\
On cut wings, to play with the transparent ones.\\
In unconsciousness the night song is sung.)}{\sc\scriptsize Osip Mandelstam, transl. by Natalie Staples}

In this chapter we extend the symplectic reduction
\index[terms]{reduction!symplectic}associated with a Hamiltonian
action of a Lie group to LCK geometry. As we shall
consider holomorphic and conformal
actions, the relevant notions will
then be  twisted  - formally, in the symplectic
definitions, we replace the operator $d$ with
$d_\theta$. We shall first describe the twisted
Hamiltonian actions, then define the LCK momentum map
\index[terms]{momentum map} and build the LCK quotient, and
finally show that compact toric LCK manifolds are
Vaisman.\index[terms]{manifold!toric}

\section{Twisted Hamiltonian actions}

The natural context for the notion of a twisted Hamiltonian vector field is LCS geometry\index[terms]{geometry!LCS} (and indeed, the definition first appears in \cite{va_lcs}, with motivation coming from theoretical mechanics), but we shall restrict to the LCK case.\index[terms]{vector field!twisted Hamiltonian}\index[terms]{action!twisted Hamiltonian}

Let $(M,I,g,\omega,\theta)$ be an LCK manifold\index[terms]{manifold!LCK}. 

\hfill

\definition A vector field $X$ is {\bf twisted
  Hamiltonian} if there exists a function $f_X\in C^{\infty}(M)$ such
that $i_X\omega=d_\theta f_X$ (equivalently,
$X=(df_X)^{\sharp_\omega}$, i.e.  $X$ is the $\omega$ dual
of $d_\theta f_X$). The set of all twisted Hamiltonian
vector fields is denoted $\ham_\theta(M)$. An action of a
group $G$ is called twisted Hamiltonian if
$\mathfrak{g}\subseteq \ham_\theta(M)$, after
identification of the elements of $\mathfrak{g}$ with the
corresponding fundamental vector fields on $M$.  
Equivalently, a vector field $X$ is twisted Hamiltonian if
the lift of $X$ to the K\"ahler cover\index[terms]{cover!K\"ahler} $\tilde M$
is Hamiltonian (see \ref{conf_inv} (ii)).

\hfill

\example\label{_Lee_twisted_Ham_Example_}
Let $(M, \omega, \theta)$ 
be an LCK manifold, and $V:= I(\theta^\sharp)$ the anti-Lee
field. Then $V$ is twisted Hamiltonian. Indeed,
$V$ is the symplectic dual to the Lee form,\index[terms]{form!Lee}
that is  $d_\theta$-exact, $\theta = - d_\theta(1)$.

\hfill

\remark\label{conf_inv} 
(i)  Let $g'=e^\f g$ (so that $\theta'=\theta+d\f$). Then:
$$d_\theta f=e^{-\f}d_{\theta'}(e^\f f),\ \ \text{and}\ \  (d_\theta f)^{\sharp_\omega}=\big(d_{\theta'}(e^\f f)\big)^{\sharp_{\omega'}},$$ hence being twisted Hamiltonian is a conformal notion. \index[terms]{action!twisted Hamiltonian}

If one defines  the {\bf twisted Poisson bracket}\index[terms]{bracket!Poisson} by
$\{f,h\}:=\omega(X_f,X_h)$, then it is easily seen that
$X_{\{f,h\}}=[X_f,X_h]$, thus $\ham_\theta(M)$ is a Lie
subalgebra of the Lie algebra $TM$ of vector fields.

(ii) Moreover, multiplication with $e^\f$ is an isomorphism between the Lie algebras $(C^{\infty}(M), \{,\})$ and $(C^{\infty}(M), \{,\}')$: indeed $\{e^\f f, e^\f h\}'=e^\f\{f,h\}$. This isomorphism commutes with the corresponding isomorphisms of taking the $\omega-$ and $\omega'-$ duals. In particular, if $M$ is GCK, we can choose $\theta'=0$, and thus {\em twisted Hamiltonian vector fields coincide with the usual Hamiltonian vector fields of the global K\"ahler metric}. Hence, using also  \ref{lift_to_min}, one can prove:\index[terms]{vector field!twisted Hamiltonian}

\hfill 

\proposition (\cite{oti1}) Let $G\subset \Aut(M)$ and let $\tilde G_0$ be the connected component of its universal cover. Then the action of $G$ is twisted Hamiltonian if and only if the action of $\tilde G_0$ on $\hat M$ is Hamiltonian with respect to the K\"ahler metric $\hat g$.\index[terms]{action!twisted Hamiltonian}

\hfill



\example (\cite{gop}, \cite{oti1})\label{ex_toric} Consider the weighted Sasaki sphere $S^{2n-1}_a$, where $a:=(a_1,\ldots a_n)$ with $a_i$ real, satisfying $0<a_1\leq a_2\leq\cdots a_n$, endowed with the contact form\index[terms]{form!contact} 
$$\eta_a=\frac{1}{\sum a_i|z_i|^2}\sum(x_idy_i-y_idx_i), \quad z_i=x_i+\1 y_i.$$
On  $S^1\times S^{2n-1}_a$, we have the Vaisman
metric corresponding to the LCK form
$\omega_a:=d_\theta\eta_a$ that is  associated with 
the K\"ahler form\index[terms]{form!K\"ahler} $d(e^t\eta_a)$ on the K\"ahler cone over
the weighted sphere.\index[terms]{weighted sphere} 
Let $T^k$, $1\leq k\leq n$, act on
$S^1\times S^{2n-1}_a$ by:
$$(u_1,\ldots, u_k)\cdot(s,(z_1,\ldots,z_{n}))=(s,
u_1z_1,\ldots,u_kz_k,z_{k+1},\ldots,z_{n})), \ u_i\in \U(1)
$$
This action preserves $|z_i|^2$ and
$d^c|z_i|^2=x_idy_i-y_idx_i$; then  it preserves $\omega_a$.
The Lee form \index[terms]{form!Lee}$\theta$ of $M=S^1\times S^{2n-1}_a$ is obtained
as the lift of the length element $\hat\theta\in
\Lambda^1(S^1)$ to $S^1\times S^{2n-1}_a$. Therefore,
$\theta$ vanishes on the fundamental vector
fields, and hence  the action is twisted Hamiltonian. \index[terms]{action!twisted Hamiltonian}

\subsection{The LCK momentum map} \index[terms]{momentum map}

Let $G$ be a connected Lie group acting on
$(M,I,[g],[\theta])$ in a twisted Hamiltonian way.
To arrive at an LCK structure\index[terms]{structure!LCK} $(M,I,g,\theta)$, one
needs to fix a choice of $g$ in its conformal class.
This also fixes the fundamental form $\omega$
and the Lee form\index[terms]{form!Lee} $\theta$.
 For any pair $(\omega, \theta)$, any  fundamental field
$X\in\mathfrak{g}$ satisfies $i_X \omega=d_\theta
f_X$, and produces a momentum map
$$\mu^\omega:M\ra \mathfrak{g}^*,\quad \langle \mu^\omega(x), X\rangle=f_X(x).$$
Equivalently, the momentum map takes a Hamiltonian vector field $X \in \g$ to its
Hamiltonian function, producing a $C^\infty(M)$-valued functional on $\g$.

\hfill

\remark  For $\omega'=e^\f\omega$, one has
$\mu^{\omega'}=e^\f\mu^\omega$. Therefore, we can speak
about the momentum map $\mu$ associated with the conformal
action of $G$, considering it as a family of momentum
mappings conformally depending on a conformal factor of the LCS form. \index[terms]{form!LCS}
In particular, the level set $\mu^{-1}(0)$ is
well-defined.

\hfill

\remark
The momentum map is not necessarily equivariant with
respect to the co-adjoint action,\index[terms]{action!co-adjoint}
because the co-adjoint action is not conformally
invariant. Still, one can modify the co-adjoint orbit such
that the momentum map become equivariant, see \cite{hr2}.

\hfill

\remark In the symplectic case, the momentum map is not
unique,\index[terms]{momentum map} the parameter space being
$H^1(\mathfrak{g})$. Here instead, the parameter space is
$H^1(\mathfrak{g}, N)$, where
$N:=\ker(d_\theta:C^\infty(M)\ra\Omega^1(M))$. From
Exercise  \ref{dtheta_inj}, if $M$ is not GCS, then
$N\equiv 0$, and hence the momentum map associated with a
twisted Hamiltonian action is unique.\index[terms]{action!twisted Hamiltonian}

\hfill

\remark By definition, if $X\in \ham_\theta(M)$
then the lift of $X$ to the symplectic cover
of $M$ is a Hamiltonian symplectomorphism.  
Therefore, $\ham_\theta(M)\subseteq\aut_s(M)$ (see
\ref{special_autom}). Conversely, if
$\omega=d_\theta\eta$, and $G$ is a compact Lie group
that acts conformally and $\mathfrak{g}\subset
\aut_s(M)$, the action is twisted Hamiltonian,
\cite{nico1}. Indeed, by averaging on $G$, $\omega$ and
$\theta$ can be taken $G$-invariant. Define
$\eta^G:=\int_Ga^*\eta d\vol_G$ (with $a\in G$). Then
$d_\theta\eta^G=\omega$. Then the momentum map is given by
$-\eta^G$.

\subsection{LCK reduction at $0$} The analogue of the K\"ahler reduction now reads:\index[terms]{reduction!LCK}

\hfill

\theorem (\cite{gop})\label{red_th} Let $(M,I,[g])$ be an
LCK manifold, and $G\subseteq \Aut(M)$ acting in a twisted
Hamiltonian way, with momentum map $\mu$. Suppose\index[terms]{momentum map}
$0\in\mathfrak{g}$ is a regular value for $\mu$, and let
$\iota:\mu^{-1}(0)\hookrightarrow M$ be the canonical
inclusion. Suppose $G$ acts freely and properly on the
level set $\mu^{-1}(0)$, and let $\pi:\mu^{-1}(0)\ra
\mu^{-1}(0)/G$ be the natural projection. Then there
exists an LCK structure\index[terms]{structure!LCK} $(I_r, [g_r])$ on
$M/\!\!/G:=\mu^{-1}(0)/G$ such that $\pi^*g_r\in
[\iota^*g]$.

\hfill

There are at least three ways in which one may prove this result:

(i) Use the LCS reduction\index[terms]{reduction!LCS} theorem proved in \cite{hr1} and
show that the whole process is compatible with the complex
structure, as in the usual proofs for the K\"ahler
reduction.\index[terms]{reduction!K\"ahler} This is essentially done in \cite{gop}.

(ii) Translate the setting in the language of conformal geometry (then, for example, the momentum mapping is $\mu:M\ra\mathfrak{g}\otimes L^2$, where $L\ra M$ is the (real) weight bundle). This was sketched in an unpublished paper by O. \index[persons]{Biquard, O.} Biquard and P. \index[persons]{Gauduchon, P.} Gauduchon, and was detailed in the arxiv version of \cite{gop}.\index[terms]{bundle!weight}

(iii) Using the K\"ahler reduction of the universal cover. The essential point here is that, by \ref{conf_inv}, the twisted Hamiltonian action\index[terms]{action!twisted Hamiltonian} of $G$ on $M$ lifts to a Hamiltonian action on $\tilde M$, with respect to the K\"ahler metric $\tilde g$. One then relates the K\"ahler momentum mapping $\tilde \mu$ with the LCK one and proves that the LCK reduced space is covered by the K\"ahler reduced space:
$$\mu^{-1}(0)/G\simeq (\tilde\mu^{-1}(0)/G)/\pi_1(M),$$
and the isomorphism is indeed one of LCK structures\index[terms]{structure!LCK}.

\hfill

\remark (\cite{gop}) 
The LCK rank\index[terms]{rank!LCK} (\ref{lck_rank}) of $M/\!\!/G$ equals the LCK rank of $M$.

\hfill

In the particular case of compact Vaisman manifolds,\index[terms]{manifold!Vaisman} one
is interested in relating the LCK reduction to the
Sasakian reduction, as described in \cite{gro}.\index[terms]{reduction!Sasakian}
Let $M$ be a compact, Vaisman manifold of LCK rank 1\index[terms]{rank!LCK},
with  associated compact Sasakian manifold $W$.
By \ref{homo_con} and \ref{isom_con}, 
the twisted Hamiltonian action of $G$ on $M$ induces  an
action by Sasakian automorphisms on $W$. Moreover, it is
proven in \cite{pilca}, \cite{bgp} (see also \cite{gop})
that the momentum mapping of a twisted Hamiltonian action\index[terms]{action!twisted Hamiltonian}
on a Vaisman manifold is $\mu(X)=\theta(IX)$. This is, up
to sign, the Sasaki momentum mapping. Using also
the above Remark, one has:\index[terms]{momentum map}

\hfill

\theorem (\cite{gop}) Let $M$ be an LCK rank\index[terms]{rank!LCK} one compact
Vaisman manifold \index[terms]{manifold!Vaisman}with associated Sasakian manifold
$W$. Suppose $G$ acts in a twisted Hamiltonian way  on $M$
satisfying the hypothesis of \ref{red_th}, and lifts to an
action of $\hat G$ on the minimal cover $\hat M$. Then the
LCK reduced space is Vaisman, of rank 1 ,  with
associated  Sasakian manifold $W/\!\!/G$.

\hfill

\remark Left aside the theoretical importance of the reduction procedure, there was a hope that this way one could obtain new examples of compact LCK manifolds. As far as we know, no new examples were found up to now. Several worked examples are given in \cite{gop}, reducing Vaisman Hopf manifolds, but they all fit into the above scheme and lead to Vaisman manifolds\index[terms]{manifold!Vaisman} of type $S^1\times$ a Sasakian space.

\hfill

\remark   It is completely non-trivial to extend the LCS and  LCK reductions to non-zero regular values of the momentum mapping. This is because the co-adjoint\index[terms]{reduction!LCK}\index[terms]{reduction!LCS} orbit\index[terms]{orbit!co-adjoint} does not have a canonical LCS (or LCK) structure, and thus the shifting trick, see \cite{ms}, cannot be performed. There are various attempts to define non-zero LCS and LCK reduction, see, for example, \cite{miron}. \index[terms]{momentum map}

\hfill

\remark In \cite{david}, a notion of {\bf Sasaki-Weyl}
manifold is introduced, in the context of CR-Weyl
geometry. The cone of such a manifold is LCK. Moreover, a
reduction scheme is developed for Sasaki-Weyl manifolds,
that is  compatible with the LCK reduction of their LCK
cone.

\section[Complex Lie group acting by holomorphic isometries]{Complex Lie group acting by holomorphic\\ isometries}

In \cite{_Klemyatin_}, N. \index[persons]{Klemyatin, N.} Klemyatin has studied 
complex manifolds equipped with a holomorphic
Killing \index[terms]{vector field!Killing}\index[terms]{vector field!holomorphic}vector field $X$ such that $X$ and $IX$ 
are both tangent to a compact torus action.
He proved that in this situation the 
diffeomorphism flow generated by $X$ acts
trivially on the Dolbeault cohomology.\index[terms]{cohomology!Dolbeault} We
used some of his computations in Chapter \ref{_Dolbeault_Vaisman_Chapter_}
to prove the Hodge decomposition theorem 
for Vaisman manifolds.\index[terms]{manifold!Vaisman}

It turns out that LCK manifolds that admit such 
vector fields are automatically Vaisman. This result
was obtained by N. \index[persons]{Istrati, N.} Istrati in 2018, even earlier than
\index[persons]{Klemyatin, N.} Klemyatin's paper appeared. 

\hfill

The following theorem generalizes \ref{kami_or}.
We present a new proof of this important result.

\hfill

\theorem {(\cite{nico2})}\label{_Killing_holo_Nicolina_Theorem_} 
\\ Let
$(M,I,\omega,\theta)$ be a compact LCK manifold, not GCK.
Consider a compact torus $T$ acting on $(M,I)$ by
biholomorphic diffeomorphisms.  Let $\goth t\subset TM$
be the Lie algebra of vector fields tangent
to this action. Assume that $I(\goth t) \cap \goth t\neq 0$.
 Then $(M,I)$ is of Vaisman type.

\hfill

\pstep 
For any compact,  connected Lie group $K$ acting
on an LCK manifold by biholomorphic diffeomorphisms,
there exists another LCK structure\index[terms]{structure!LCK} that is  $K$-invariant
(\ref{_LCK_averaging_Lemma_}).
This LCK structure is obtained by a double averaging
argument. One first averages the conformal multiplier
to have a Lee form\index[terms]{form!Lee} that is  $K$-invariant. Then one
averages the LCK form $\omega$ which already has
a $K$-invariant Lee form.


\hfill

{\bf Step 2:}
Choose $\xi\in \Lie(T)\setminus\{0\}$ 
with $I\xi\in \Lie(T)$. Since $\xi$ and $I\xi$
act on $M$ by isometries, their lifts act on
the K\"ahler cover $\tilde M$ by conformal
biholomorphisms. However, any conformal
biholomorphic automorphism of a K\"ahler
manifold of $\dim_\C >1$ is a homothety 
(Section \ref{_LCK_chapter_Intro_}). 
If the lift of $\xi$ to $\tilde M$ is
not an isometry,  we apply \ref{kami_or} 
and obtain that $(M, I)$ is Vaisman type. 
It remains to finish the proof of 
\ref{_Killing_holo_Nicolina_Theorem_} by absurd,
assuming that the lifts $\tilde \xi$ and $I(\tilde \xi)$ of $\xi$ and $I\xi$
to $\tilde M$ are holomorphic and Killing.\index[terms]{vector field!Killing}\index[terms]{vector field!holomorphic}

\hfill

{\bf Step 3:}  
Let $X$ be a Killing holomorphic vector field
on a K\"ahler manifold, such that $I(X)$ is also
Killing. By Exercise
\ref{_Killing_holomorphic_is_parallel_Exercise_},
$X$ is parallel\index[terms]{vector field!parallel} with respect to the Levi--Civita 
connection. Applying this to the vector field
$\tilde \xi$, we obtain that it is 
parallel, therefore, it has constant length.
On the other hand, $\tilde \xi$ is 
invariant under the homothety action; in particular, 
it cannot have constant length.
We arrived at the contradiction,
and finished the proof of
\ref{_Killing_holo_Nicolina_Theorem_}.
\endproof

\hfill

When $M$ is Vaisman, and $\goth t\subset TM$ is the algebra
of holomorphic Killing vector fields, the intersection
$\goth t\cap I(\goth t)$ can be described explicitly.
The following result is also due to \cite{nico2};
the proof is new.

\hfill

\proposition\label{_Killing_I(X)_not_killing_Proposition_}
(\cite[Proposition 4.1]{nico2})\\
Let $X$ be a holomorphic Killing\index[terms]{vector field!Killing}\index[terms]{vector field!holomorphic} vector field
on a Vaisman manifold,\index[terms]{manifold!Vaisman} such that $I(X)$ is also
Killing. Then $X$ is tangent to the canonical foliation.\index[terms]{foliation!canonical}

\hfill

\pstep
We say that a holomorphic tensor over a Vaisman manifold
is ``Lee invariant'' if it is invariant with respect
to the action of $e^{t\theta^\sharp}$, where
$\theta^\sharp$ is the Lee field.\index[terms]{Lee field} 
From \ref{_holo_tensor_on_Vaisman_Lee_invariant_Theorem_} 
it follows that 
all holomorphic tensors, in particular,
all holomorphic vector fields on a compact Vaisman
manifolds are Lee invariant. Therefore, 
any holomorphic vector field commutes
with the Lee field.\index[terms]{Lee field} This allows one to 
decompose any holomorphic vector field
$X$ onto $\Sigma$-tangent part and the $\Sigma$-basic part.
Let ${\goth t}_0$ be the Lie algebra of $\Sigma$-basic
holomorphic Killing vector fields.\index[terms]{vector field!Killing}\index[terms]{vector field!holomorphic} To finish the
proof, it remains to show that ${\goth t}_0\cap I({\goth t}_0)=0$.

\hfill

{\bf Step 2:}  
Let $X$ be a holomorphic basic vector field,
and $\tilde X$ its lift to the K\"ahler 
cover $\tilde M$ of $M$\!, locally isometric to a cone
over a Sasakian manifold, $\tilde M \cong C(S)= S \times \R^{>0}$. 
Since $\tilde X$ is orthogonal to the Lee field\index[terms]{Lee field} $t \frac d{dt}$, 
it is tangent to the level set of $t$. Since $dd^c (t^2)$ is the
K\"ahler metric on $\tilde M$, the vector field $\tilde X$
acts on $\tilde M$ by holomorphic isometries.
Applying Exercise
\ref{_Killing_holomorphic_is_parallel_Exercise_}
again, we find that $\tilde X$ is parallel on $\tilde M$
whenever both $X$ and $I(X)$ belong to ${\goth t}_0$.
The same argument as in \ref{_Killing_holo_Nicolina_Theorem_},
Step 3, implies that this is impossible: a parallel vector
field has\index[terms]{vector field!parallel} constant length, but $\tilde X$ is invariant
under the homothety action; this means that  it cannot have constant length.
We proved that  ${\goth t}_0\cap I({\goth t}_0)=0$.
\endproof

\hfill

\proposition\label{_Killing_Hamiltonian_Proposition_} 
Let $M$ be a Vaisman manifold,\index[terms]{manifold!Vaisman} $\Sigma$ its
canonical foliation,\index[terms]{foliation!canonical} and $V\in TM$ a $\Sigma$-basic
Killing holomorphic vector field.\index[terms]{vector field!Killing}\index[terms]{vector field!holomorphic} Denote by $\tilde V$
its lift to a K\"ahler cover $(\tilde M, \tilde \omega)$ of $M$.
Then $\tilde V$ is a Hamiltonian vector field
on $\tilde M$.\index[terms]{vector field!Hamiltonian}

\hfill

\proof
Let $Y:= I(\tilde V)$, and 
$\phi=\tilde \omega(\theta^\sharp, I\theta^\sharp)$ be the automorphic
K\"ahler potential. Since $\tilde V$ is holomorphic and Killing
on $\tilde M$ (\ref{_Killing_I(X)_not_killing_Proposition_}, Step 2), we have
$\Lie_{\tilde V}d^c\phi=0$, and
\[ 
i_{\tilde V} \tilde \omega=  i_{\tilde V}dd^c\phi=
\Lie_{\tilde V}d^c\phi- d (i_{\tilde V} d^c\phi)
= d \langle Y, d\phi\rangle= d(\Lie_{Y}\phi).
\]
which implies that $\Lie_{Y}\phi$ is the Hamiltonian for $\tilde V$.
\endproof

\hfill

\ref{_Killing_Hamiltonian_Proposition_} has the following
curious corollary.

\hfill

\corollary\label{_Sasa_Reeb_isotro_Corollary_}
Let $S$ be a compact Sasakian manifold, $\dim_\R S = 2n-1$,
and $G \subset \Iso(S)$ the closure of the subgroup generated
by the Reeb flow in the group of isometries of $S$.
Then $\dim_\R G \leq n$, and the equality can be realized
only when $S$ is a toric Sasakian 
manifold.\footnote{For a definition of toric Sasakian
manifold, see \ref{_toric_Sasakian_Definition_}.}

\hfill

\proof
The group $G$ is compact, because $\Iso(M)$ is compact, and commutative,
because it is a closure of a commutative group.
Consider the action of $G$ on the corresponding Vaisman manifold\index[terms]{manifold!Vaisman}
$M:=C(S)/\Z$. Then all vector fields on $M$ induced by the  
$\Lie(G)$-action are holomorphic and Killing. \index[terms]{vector field!Killing}\index[terms]{vector field!holomorphic}
By \ref{_Killing_Hamiltonian_Proposition_},
$\Lie(G)$ is a commutative algebra of Hamiltonian
vector fields on $\tilde M = C(S)$.
It is not hard to see that the Poisson bracket\index[terms]{bracket!Poisson} of two functions
$f$ and $g$ is $\tilde \omega(X_f, X_g)$, where $X_f, X_g$ are
their Hamiltonian vector fields. Since all vector fields in
$\Lie(G)$ commute, one has $\tilde \omega(X, Y)=0$ for all
$X, Y \in \Lie(G)$. This implies that any orbit
of $G$ is isotropic in $\tilde M$; hence 
has dimension at most $n$.\footnote{In classical
mechanics, the isotropicity of an orbit
of a group generated by commuting Hamiltonians
is called ``the Liouville theorem'', see 
\cite[Corollary 16.12]{_Libermann_Marle_}.}
\endproof

\hfill

\remark Let $(M,I, \omega, \theta)$ be a Vaisman
manifold, and $T$ a commutative group 
of holomorphic isometries of $M$.  \index[terms]{metric!Vaisman}
Denote by $\goth s\subset TM$ the intersection $\Lie(T)\cap I(\Lie(T))$.
By \ref{_Killing_I(X)_not_killing_Proposition_},
$\Lie(T)\cap I(\Lie(T))=\spp\{\theta^\sharp,I\theta^\sharp\}$. This
immediately
implies the following criterion for non-existence of LCK
metrics:

\hfill

\corollary {(\cite{nico2})} Let $M$ be a complex, compact
manifold such that its group of biholomorphisms contains a
torus $T$ with $\dim(\Lie(T)\cap I(\Lie(T)))>2$. Then $M$ admits no LCK metric.

\hfill

\corollary {(\cite{kami})} Let $M$ be a compact, strictly 
LCK  manifold endowed with a holomorphic action
 of an elliptic curve $T^1_\C$. 
Consider the natural homomorphism  
$H_1(O)\arrow H_1(M,\Z)$, where $O$ is an orbit of $T^1_\C$.
Then it induces a homomorphism $H_1(O)=\Z^2\ra H_1(M,\Z)$ of rank 1.

\hfill

\proof  
By \ref{_Killing_I(X)_not_killing_Proposition_},
 $M$ is of Vaisman type, and the action of $T^1_\C$
preserves the canonical foliation.\index[terms]{foliation!canonical} This implies
that $M$ is a quasi-regular Vaisman manifold,\index[terms]{manifold!Vaisman!quasi-regular}
and any $T^1_\C$-orbit is a fibre of the
elliptic fibration $M \arrow X$.

\hfill

\proposition {(\cite{ov_imrn_10})}\label{fund_group}
\index[terms]{fundamental group!of a Vaisman manifold}
Let $M$ be a compact, quasi-regular Vaisman manifold,
and $\pi:\; M \arrow X$ be the elliptic fibration.
Let $F$ be a general fibre of $\pi$.
Then the following sequence of the
fundamental groups\index[terms]{fundamental group} is exact
	\[
	\pi_1(F) \stackrel \tau \arrow \pi_1(M) \arrow \pi_1(X) \arrow 0,
	\]
where $\pi_1(X)$ is the fundamental group of $X$, considered
as a complex variety. Moreover, the kernel $\tau$ has rank 1.

\hfill

\proof The long exact homotopy sequence then gives
\begin{equation}\label{ex_seq_homo}
\pi_2(X) \stackrel \delta
\arrow \pi_1(T^2) \arrow \pi_1(M) \arrow \pi_1(X) \arrow 0
\end{equation}
The boundary operator $\delta$ can be described
as  follows. Let $\gamma:\; \Z^2 \arrow H^2(X)$
be the map  representing the  Chern
classes of the corresponding $S^1\times S^1$-fibration.
We may interpret this map as a differential
of the corresponding \index[persons]{Leray, J.} Leray spectral sequence, 
which\index[terms]{spectral sequence!Leray}
gives us an exact sequence
\[
0 \arrow H^1(X) \arrow H^1(M) \arrow
H^1(T^2)\stackrel \gamma\arrow H^2(X).
\]
Dualizing and using the Hurewicz theorem,
we obtain that the  boundary map $\pi_2(X) \stackrel \delta
\arrow \pi_1(T^2)$ is obtained as a
composition of $\gamma^*$ and the
Hurewicz homomorphism $\pi_2(X) \arrow H^2(X)$.
The Chern classes of the $S^1\times S^1$-fibration
are easy to compute: one of them
is trivial (because $M$ is fibred over a circle),
and the other one is non-trivial,\index[terms]{class!Chern}
because $M$ is non-K\"ahler, and the total space of
an isotrivial elliptic fibration with trivial Chern classes is K\"ahler
(\ref{_Blanchard_Theorem_}).
Therefore, the image of $\delta$ has rank $1$ in $\pi_1(T^2)$.
\endproof


\section{LCK manifolds admitting a torus action with an open orbit}


The following elementary lemma is used further on in
this section.

\hfill

\lemma\label{_open_orbit_Lemma_}
Let $M$ be a connected complex $n$-manifold,
and $G$ an $n$-di\-men\-si\-o\-nal complex Lie
group acting on $M$ with an open orbit.
Then this open orbit is unique.

\hfill

\proof
Let $x_1, ..., x_n$ be a basis in the
Lie algebra of $G$, and $\xi:=x_1 \wedge ... \wedge x_n\in H^0(\Lambda^n TM)$
the corresponding section of the anticanonical bundle.
Clearly, a point $m\in M$ belongs to an open orbit
if and only if all $x_i$ are linearly independent in $m$,
equivalently, if $\xi\restrict m\neq 0$.
Since the zero set $Z$ of $\xi$ is complex
analytic, the union $M\backslash Z$ of all open orbits is 
connected. A connected set cannot be represented
as a non-trivial disjoint union of its open subsets,
and hence there is just one open orbit.
\endproof

\hfill

The following theorem is used further on to develop
the theory of toric LCK manifolds. Originally, it was proven
in a different setting (assuming a twisted Hamiltonian
action of a compact torus of maximal dimension)
by N. \index[persons]{Istrati, N.} Istrati, \cite{nico1}. We prove essentially the same
statement using different methods.

\hfill

\theorem\label{_Toric_Vaisman_Nicolina_Theorem_} 
(\cite{nico1}) 
Let $M$ be a compact complex $n$-manifold admitting an LCK structure,\index[terms]{structure!LCK} and
 $T_\alg = (\C^*)^n$ an algebraic torus acting on $M$
holomorphically with an open orbit. Assume that the action of $T_\alg$ is lifted
to the minimal K\"ahler cover $\tilde M$. Then $M$ 
is of Vaisman type. 

\hfill

\pstep
Denote the deck transform group of $\tilde M$ by $\Gamma$.
Then $M = \tilde M/\Gamma$.
Let $\tilde U\subset \tilde M$ be
the open orbit\index[terms]{orbit!open}\index[terms]{orbit!dense} 
of $T_\alg$; it is unique by \ref{_open_orbit_Lemma_}.  
Since $\tilde U$ is biholomorphic to
a quotient of $T_\alg$, the automorphisms of $\tilde U$ commuting
with the  $T_\alg$-action belong to $T_\alg$. 
Therefore, the action of $\Gamma$ on $\tilde U$
factorizes through $T_\alg$, and, moreover,
$\Gamma$ is identified with the kernel of the natural map 
$T_\alg \arrow \Aut(U)$, where $U$ is the image 
of $\tilde U$ in $M$.

Now pick $a\in \Gamma=\Z^k$ that acts by
non-isometric homotheties on $\tilde M$ and consider its logarithm $X$
in the Lie algebra of $T_\alg$. The flow
generated by $X$ acts on $M$ by holomorphic automorphisms,
and the corresponding group is isomorphic
to $S^1$ because $a$ acts on $M$ trivially. 

Let $T^n\subset (\C^*)^n$ be the maximal compact
subgroup. Choose an LCK structure on $M$ which
is $T^n$-invariant. Such an LCK structure\index[terms]{structure!LCK} can
be constructed using the double averaging
technique, \ref{_LCK_averaging_Lemma_}.

Averaging this LCK structure with this
$S^1$-action,\index[terms]{action!$S^1$-} we can assume that
the LCK structure\index[terms]{structure!LCK} $(g, \theta)$ 
on $M$ is also $S^1$-invariant. Since this
$S^1$-action comes from $T_\alg$, it 
commutes with the $T^n$-action,
hence $(g, \theta)$ is invariant
with respect to the group $G$ generated by
$S^1$ and $T^n$. 

\hfill

{\bf Step 2:} We prove that $G\cong T^{n+1}$.
Since $T^n$ is a compact group which 
acts on $M$ by isometries and on $\tilde M$
by homotheties, it acts on $\tilde M$ by
isometries (a compact group cannot act by 
non-isometric homotheties).
Therefore, the compact torus $G$ is $n+1$-dimensional.

\hfill

{\bf Step 3:}
Denote by ${\goth t}_\alg$ the Lie algebra of $T_\alg$,
and let ${\goth t}^{n+1}\subset {\goth t}_\alg$ be the Lie
algebra of the $T^{n+1}$-action obtained in Step 2. 

Since ${\goth t}_\alg$
is a complex $n$-dimensional Lie algebra, the real
$(n+1)$-dimensional subspace ${\goth t}^{n+1}\subset {\goth  t}_\alg$
contains a 1-dimensional complex subspace ${\goth u}= \langle X, I(X)\rangle$.
An LCK manifold admitting two holomorphic\index[terms]{vector field!Killing}\index[terms]{vector field!holomorphic}
Killing vector fields $X$ and $I(X)$ is Vaisman
by \index[persons]{Istrati, N.} Istrati's theorem (\ref{_Killing_holo_Nicolina_Theorem_}).
\endproof


\section{Toric LCK manifolds}\label{toric_lck}
\index[terms]{manifold!toric}


Consider now a compact torus $T^n$ acting\index[terms]{action!twisted Hamiltonian} effectively, in
a twisted Hamiltonian way, on the connected, compact LCK
manifold $(M,I,g, \theta)$, with $\dim_\C M=n$. We say that
$M$ is a {\bf toric LCK manifold}. Denote
by $\mathfrak{t}$ the Lie algebra of $T^n$. Note that, as
$T^n$ is compact, the averaging procedure guarantees the
existence of an LCK form in the conformal class with
respect to which $T^n$ acts by holomorphic  isometries.

\hfill

\remark
In algebraic geometry\index[terms]{geometry!algebraic}, a ``toric manifold'' is a
projective manifold equipped with an action of an
``algebraic torus'' $T_\alg=(\C^*)^n$ that has an open orbit.
In symplectic geometry,\index[terms]{geometry!symplectic} a ``toric manifold''
is a symplectic manifold of real dimension $2n$ equipped with an effective
Hamiltonian action of a compact torus $T^n$ of dimension $n$. To pass
from the symplectic definition to the algebraic-geometric,
one replaces the compact torus $T^n$ by its
complexification $T_\alg$. 

\hfill

In fact, these two definitions agree in LCK geometry 
just as well as in the K\"ahler geometry.\index[terms]{geometry!K\"ahler}\index[terms]{geometry!LCK}

\hfill

\theorem
Let $M$ be a compact complex manifold of LCK type, $\dim_\C M=n$.
Then the following are equivalent
\begin{description} 
\item[(i)] $M$ admits a holomorphic action 
of the Lie group $(\C^*)^n$ with an open orbit that can be lifted to a 
 K\"ahler cover of $M$.
\item[(ii)]  $M$ admits an LCK structure\index[terms]{structure!LCK} with an effective twisted
Hamiltonian action of a compact torus $T^n$.\index[terms]{action!twisted Hamiltonian}
\end{description}

\pstep We deduce (i) from (ii).
Since the $T^n$-action is twisted Hamiltonian,
the space generated by $\Lie(T^n)$ in $TM$ is 
isotropic with respect to the Hermitian form\index[terms]{form!Hermitian} $\omega$
by Liouville theorem\index[terms]{theorem!Liouville}
(see the proof of \ref{_Sasa_Reeb_isotro_Corollary_}).
Then $\Lie(T^n)\cap I(\Lie(T^n))=0$, thus  
$\Lie(T^n)+ I(\Lie(T^n))$ generate $TM$ at any point where the
action of $\Lie(T^n)$ is effective. Therefore
the corresponding complex Lie group acts with
an open orbit. 

\hfill

{\bf Step 2:} It remains to prove the implication
(i) $\Rightarrow$ (ii): starting from a compact
LCK manifold $M$ admitting an action
of $(\C^*)^n$ with an open orbit, prove that
the action of the maximal compact subgroup $T^n \subset (\C^*)^n$
is twisted Hamiltonian for an appropriate LCK metric.
By another theorem of \index[persons]{Istrati, N.} Istrati
(\ref{_Toric_Vaisman_Nicolina_Theorem_}),
$M$ is Vaisman, because it admits an action of $(\C^*)^n$
with an open orbit. 

Using the double averaging technique
(\ref{_LCK_averaging_Lemma_}), we may assume
that the maximal torus $T^n\subset (\C^*)^n$ acts on $\tilde M$ and $M$ by holomorphic
isometries. It remains to prove that this
action is twisted Hamiltonian. 
Any holomorphic
Killing vector field on a Vaisman manifold \index[terms]{manifold!Vaisman}$M$ 
which lifts to a Killing field on its K\"ahler
cover is Hamiltonian by\index[terms]{vector field!Killing}\index[terms]{vector field!holomorphic} \ref{_Killing_Hamiltonian_Proposition_}.
This implies that the action of $T^n$ is twisted Hamiltonian. 
\endproof

\hfill

\example Consider the weighted  Hopf manifold in
\ref{ex_toric}, associated with the weighted Sasakian sphere\index[terms]{weighted sphere}
$S^{2n-1}_a$. Let the torus $(\C^*)^n$ act in the 
with the same weights on $\C^n\backslash 0$, and
and induce this action on
the corresponding Vaisman manifold \index[terms]{manifold!Vaisman}
$M = C(S^{2n-1}_a)/\Z$. We obtain a toric LCK
manifold.  Note that the weighted Hopf manifold is
Vaisman, that can be seen now from 
\ref{_Toric_Vaisman_Nicolina_Theorem_}.

\hfill

\proposition \label{_toric_LCK_surfaces_Proposition_}
(\cite{_Madani_Moroianu_Pilca:toric_})
Let $M$ be a compact complex surface admitting an LCK
toric structure. Then $M$ is a finite quotient of a
diagonal Hopf surface.\index[terms]{manifold!toric}\index[terms]{surface!Hopf}

\hfill

\pstep
This proof 
is different from the one given in \cite{_Madani_Moroianu_Pilca:toric_}.
As shown by F. Belgun (\ref{_Vaisman_is_Hopf_or_elli_Proposition_}),
a complex surface $M$ is Vaisman if it is Hopf or principal
elliptic. In the latter case, the canonical foliation\index[terms]{foliation!canonical} of
$M$ has compact leaves. Consider the map $\pi:\; M\arrow S$ to
the leaf space. The group $\Aut(M)$ of complex automorphisms of a
toric surface $M$ is at least 2-dimensional. However, the group
of automorphisms of the fibres of $\pi$ is 1-dimensional. Therefore,
$\Aut(M)$ has to act on $S$ with at least 1-dimensional
orbits. This immediately rules out the case when the curve
$S$ has genus $> 1$, because $\Aut(S)$ is finite.

\hfill

{\bf Step 2:}
When $S$ is an elliptic curve and
$M$ is a  Kodaira surface, the group $\Aut(M)$\index[terms]{surface!Kodaira}
is 1-dimensional, that is  easy to see from the
following argument. As shown in Chapter \ref{comp_surf},
a primary Kodaira surface $M$ is obtained as the 
quotient $\Tot(L^\circ)/\Z$ of the total space  of non-zero
vectors in an ample bundle $L$ over a torus $S$ by $\Z$ acting by homotheties.
Then any automorphism of $M$ preserves the $\C^*$-bundle\index[terms]{bundle!line!ample}
$L^\circ$; however, the bundle $L$ can be recovered from $L^\circ$ 
up to taking powers. Since the group of automorphisms of
$S$ acts effectively on $\Pic(S)$, it cannot preserve the
 $\C^*$-bundle $L^\circ$. Therefore, the group of
automorphisms of a primary \index[terms]{surface!Kodaira} Kodaira surface is 1-dimensional.
A secondary Kodaira surface is a finite quotient
of a primary one, and the same proof works for it as well.

For an alternative argument showing that Kodaira surfaces
are not LCK toric, see \ref{_b_1_LCK_Toric_Theorem_}.\index[terms]{surface!Kodaira}

\hfill

{\bf Step 3:}
We have shown that any toric LCK surface $M$ is Hopf;
taking a finite covering, we may assume that $M$ is
a primary Hopf surface. Using \index[terms]{surface!Hopf!primary}
\ref{_Hopf_surface_Vaisman_Theorem_} below,
we obtain that $M$ is a diagonal Hopf surface. \index[terms]{surface!Hopf!diagonal}
\endproof

%

\hfill

\definition\label{_toric_Sasakian_Definition_}
Let $G$ be a torus acting on a Sasakian manifold $S$\index[terms]{manifold!Sasaki}
by Sasakian automorphisms, $C(S)$ its K\"ahler
cone, and $R\subset \Aut(C(S))$ the group generated by the
Reeb field. Denote by $G_1\subset  \Aut(C(S))$ 
the group generated by $G$ and $R$. The Sasakian manifold
$S$ is called {\bf toric} if the action of $G_1$ on $C(S)$\index[terms]{manifold!toric}
is toric.

\hfill

A Sasakian manifold is toric if and only if the Vaisman
manifold obtained as a $\Z$-quotient of its cone is toric,
as the following theorem implies.

\hfill

\theorem (\cite{pilca}) Let $M$ be a compact toric Vaisman
manifold of rank 1 . Then the associated compact Sasakian
manifold $S$ is toric. Conversely, let $S$ be a toric
Sasakian manifold, and $M:=C(S)/\Gamma$ a compact Vaisman
manifold associated with it; if  the torus action, extended
trivially to the cone, commutes with the action of
$\Gamma$, then $M$ is toric.
\endproof

\hfill

A K\"ahler toric manifold is simply connected (\cite[\S 3.2]{_Fulton:Toric_}).
Similarly, a Vaisman toric manifold satisfies
$b_1(M)=1$.

\hfill

\theorem (\cite{_Madani_Moroianu_Pilca:toric_}) \label{_b_1_LCK_Toric_Theorem_}
Let $M$ be a compact, toric, Vaisman manifold.\index[terms]{manifold!Vaisman} Then:\index[terms]{manifold!toric}
\begin{description}
\item[(i)]  $b_1(M)=1$. 
\item[(ii)]  The  Kodaira dimension $\kod \ M=-\infty$.\index[terms]{dimension!Kodaira}
\end{description}

\proof  We present a different proof from the original
one.  We prove (ii) for compact complex manifolds
admitting a holomorphic action  of the complex Lie group $G$ 
such that $G$ has an open orbit in $M$, and at least
one vector field contained in its Lie algebra ${\goth g}$
vanishes somewhere.

When $M$ is an LCK toric manifold, the complex Lie group
$G$ is obtained by complexifying the compact torus action;
it has an open orbit by dimension count. In this case, the 
vanishing of some $x\in {\goth g}=\Lie(G)$ is clear, 
because otherwise $M=G$ is a compact complex torus, 
and then $M$ is K\"ahler.

To prove (ii), consider a point $x\in M$ such that
the Lie algebra ${\goth g}$ of $G$ generates the
tangent space $T_x M$. Take the vector fields
$\xi_1, ..., \xi_n\in {\goth g}\subset H^0(TM)$
such that $\xi_1, ..., \xi_n$ give a basis in $T_x M$,
but one of $\xi_i$ vanishes somewhere.
Since the vector fields $\xi_i$ are 
holomorphic and linearly independent in $x$,
 the exterior product $\Xi:=\xi_1\wedge ...\wedge \xi_n$ 
is a non-zero section of the anticanonical bundle\index[terms]{bundle!vector bundle!anticanonical}
$K_M^{-1}$. 
If some power $K_M^{\otimes k}$ of the
canonical bundle has a section, we pair this
section with $\Xi^k$, obtaining a non-zero holomorphic
function on $M$. This function is non-constant
because some of $\xi_i$ vanish somewhere on $M$. This gives a contradiction.

To prove (i), we replace the Vaisman metric by a Vaisman
metric of LCK rank 1\index[terms]{rank!LCK}, retaining the complex structure (\ref{defor_improper_to_proper}).
Then the minimal K\"ahler cover $\hat M$ of $M$ is an algebraic
cone over a projective orbifold $X$, giving $\hat M =C(X)$. The group $G$ acts on \index[terms]{cone!algebraic}
$\hat M$ preserving the Lee field,\index[terms]{Lee field} and hence it
commutes with the standard $\C^*$-action on $C(X)$; indeed, the
$\C^*$-action belongs to the closure of the group
generated by the Lee field.\index[terms]{Lee field} Therefore,
$X$ is a toric K\"ahler orbifold, and $\pi_1(X)=0$ by 
\cite[\S 3.2]{_Fulton:Toric_}.


Writing the exact sequence
\[
0 \arrow H^1(X) \arrow H^1(M) \arrow
H^1(T^2)\stackrel \gamma\arrow H^2(X).
\]
as in \ref{fund_group}, we observe that the rank of $\gamma$ is 1.
Therefore, $b_1(X)=0$ implies that
$b_1(M)=1$.
\endproof

\hfill

An interesting relation between toric LCK manifolds and existence of LCK potentials\index[terms]{potential!LCK} was recently proven:

\hfill

\theorem {(\cite{nico3})}  Let $(M,J,g,\theta)$ be a  compact toric LCK manifold.\index[terms]{manifold!LCK} Denote by $G$ the  torus acting on $M$. Then any $G$-invariant LCK metric in the conformal class $[g]$ admits a unique $G$-invariant positive LCK potential.

\hfill

\remark An LCS analogue of  Guillemin--Sternberg, and  Atiyah theorem\index[terms]{theorem!Guillemin--Sternberg, Atiyah} on the convexity of the image of the momentum map of a toric symplectic manifold \index[terms]{manifold!toric} was recently proven in \cite{bgp}, but an analogue of \index[persons]{Delzant, T.} Delzant's construction\index[terms]{Delzant's construction} is still missing.\index[terms]{theorem!Atiyah, Guillemin--Sternberg}


\chapter{Elliptic curves on Vaisman manifolds}\label{_Elliptic_curves:Chapter_}\index[terms]{curve!elliptic}\index[terms]{manifold!Vaisman}

\epigraph{\it \hspace{.1in} Nous enjambons la robe de la Reine, toute en dentelle,
avec deux bandes de couleur bise (ah! que l'acide
corps de femme sait tacher une robe \`a l'endroit de l'aiselle!).

\hspace{.1in} Nous enjambons la robe de Sa fille, toute en dentelle, avec deux bandes de couleur vive (ah! que la langue du l\'ezard sait cueillir les fourmis \`a l'endroit de l'aiselle!).

\hspace{.1in} Et peut-\^etre le jour ne s'\'ecoule-t-il point qu'un m\^eme homme a br\^ull\'e pour une femme et pour sa fille.}{\sc\scriptsize St.-John Perse, \ \ Anabasis}

\section{Counting elliptic curves}

All complex curves on compact Vaisman 
manifolds are elliptic, as shown by the next
result.

\hfill

\theorem Let $C$ be a complex curve on a compact Vaisman manifold $(V, \omega)$. Then $C$
is a leaf of the canonical foliation\index[terms]{foliation!canonical} $\Sigma$. In particular, $C$ is an elliptic curve.\index[terms]{curve!elliptic}

\hfill

\proof Obviously $\int_C\omega=0$, by Stokes'  theorem, thus  $C$
is tangent to $\Sigma=\ker\omega$. But all compact\index[terms]{theorem!Stokes}
leaves of $\Sigma$ are elliptic since the tangent bundle
$T\Sigma$ is trivial by construction. \endproof

\hfill

\remark Since  a compact Vaisman manifold $V$ can be holomorphically\index[terms]{manifold!Vaisman}
embedded in a diagonal Hopf manifold $H$ (see
Section \ref{embedding_vaisman_hopf}), intersecting $V$ with two complementary
flags of Hopf submanifolds (which exist due to a result of Ma. \index[persons]{Kato, Ma.} Kato, 
\cite{kato2}), we see that $V$ contains at least two
elliptic curves. In fact, there are many more elliptic
curves, as shown below.

\hfill

Recall first that, if one is interested in complex
geometry only every compact
Vaisman manifold admits a Vaisman metric of LCK rank\index[terms]{rank!LCK} one (\ref{_Vaisman_defo_transve_Proposition_}).
In other words, there exists a Vaisman structure 
with K\"ahler cover and the deck transform group $\Z$. Moreover,
every compact Vaisman manifold\index[terms]{manifold!Vaisman!quasi-regular} can be deformed to a
quasi-regular one (\ref{defovai}). With these
preparations, we can state the main result of this
section:

\hfill

\theorem\label{_Vaisman_number_Theorem_} {(\cite{_ov_closed_reeb_})}
Let $V$ be a Vaisman manifold, $V'$
its quasi-regular deformation, and $X:=V'/\Sigma_{V'}$ 
the corresponding K\"ahler orbifold. Assume that
$V$ admits only finitely many elliptic curves.
Then the number of elliptic curves on $V$ 
is equal to $\sum_{i=0}^{\dim_\R X} b_i(X)$.

\hfill

\proof Let $\theta^\sharp$ be the Lee field\index[terms]{Lee field} of $V$, and $G$ the $C^0$-closure\index[terms]{topology!$C^0$} of the flow of $\theta^\sharp$
and $I(\theta^\sharp)$ in the
group of diffeomorphisms of $V$. The key remark is that
since the canonical foliation\index[terms]{foliation!canonical} is unique, $G$ only depends
on the complex structure of $V$. Let  $\tilde G$ be the
connected Lie group obtained as the lift of $G$ \index[terms]{foliation!canonical}
to the $\Z$-cover \index[terms]{cover!K\"ahler $\Z$-}$\tilde V$. From \ref{_Vaisman_Lee_action_contains_monodromy_Corollary_} 
it follows that  $\Gamma\cap \tilde G$ has finite index in $\Gamma$.

Clearly, the elliptic curves on a Vaisman manifold $V$ 
coincide with the closed 2-dimensional orbits of $G$.
We are going to count the number of such orbits.\index[terms]{manifold!Vaisman}

Now, the number of 2-dimensional orbits of $G= \tilde
G/(\Gamma\cap \tilde G)$ on $V$,
of $\tilde G$ on $\tilde V=\tilde V'$ and of $G':= \tilde
G/(\Gamma\cap \tilde G)$
on $V'$ is the same. However, each $k$-dimensional
orbit of $G'$ on $V'$ is elliptically fibred
over a $(k-2)$-dimensional orbit of $G'$ on
the leaf space $X:=V'/\Sigma_{V'}$ of the
corresponding canonical foliation.\index[terms]{foliation!canonical} Therefore,
the number of elliptic curves on $V$
is the same as the number of fixed points
of $G'$ acting on the orbifold $X$. This number can
be found using the\index[terms]{theorem!Bia\l ynBaum--Bott}
\index[persons]{Bia\l ynBaum--Bott, A.} Bia\l ynBaum--Bott's theorem (\cite{_Birula_}) counting fixed points of
algebraic groups and its extension to orbifolds by \index[persons]{Fontanari, C.} Fontanari:

\hfill

\theorem {(\cite{_Fontanari_})} Let $\zeta$ be a holomorphic
vector field on a compact projective orbifold $X$ of
complex dimension $k$. Then the number of zeros of $\zeta$ is
equal to $\sum_{i=0}^{2k}b_i(X)$, the sum of all Betti numbers of $X$.  \endproof

\hfill

Using the Lefschetz $\goth{sl}(2)$-action on the cohomology
of $X$, we obtain a convenient lower bound.

\hfill

\corollary\label{_Vaisman_bound_from_Lfsch_Corollary_}
Let $V$ be a Vaisman manifold that contains \index[terms]{manifold!Vaisman}
$r$ distinct elliptic curves. Then $r\geq \dim_\C V$.
\endproof


\section{Application to Sasaki manifolds: closed Reeb
  orbits}\label{_Closed_Reeb_orbits:Section_}


\subsection{Boothby--Wang theorem for Besse contact ma\-ni\-folds}

A contact manifold with a fixed contact form\index[terms]{form!contact} 
is called {\bf regular} if its Reeb orbits are 
the orbits of a free $S^1$-action.\index[terms]{action!$S^1$-} The\index[terms]{manifold!contact!regular}  Boothby--Wang theorem\index[terms]{theorem!Boothby--Wang}
gives an explicit description of regular contact manifolds
in terms of prequantized symplectic structures, which
was the antecedent of the structure theorem for Sasakian
manifolds (\ref{regsas}). Recall that {\bf
  prequantization} of a symplectic manifold\index[terms]{prequantization}
$(X,\omega)$ with integer cohomology class of $\omega$
is a choice of complex Hermitian line bundle $L$ with
connection, such that the curvature of $L$ 
is equal to $-\1 \omega$. Let $S$ be the
space of unit vectors in $L$.\index[terms]{connection!Ehresmann}
Then the Frobenius form\index[terms]{form!Frobenius} of the horizontal bundle $T_\hor S$
of the corresponding Ehresmann connection is
equal to its curvature 
$\omega$; hence  $T_\hor S\subset TS$ is a contact 
distribution on $S$.

\hfill

\theorem \index[terms]{theorem!Boothby--Wang}
(\index[persons]{Boothby,  W. M.}Boothby-\index[persons]{Wang, H. C.}Wang, \cite{_Boothby_Wang_})\\
Let $S$ be a compact regular contact manifold,\index[terms]{manifold!contact!regular}
and $X$ the leaf space of the Reeb flow. Then
$X$ is equipped with a prequantization $(L, \nabla)$, in such
a way that the space of unit vectors of $L$
is naturally contactomorphic to $S$.\index[terms]{Boothby--Wang fibration}
\endproof

\hfill

\definition
A contact manifold $S$ with a fixed contact form \index[terms]{form!contact}is called
a {\bf \index[persons]{Besse, A.} Besse contact manifold} if its Reeb orbits\index[terms]{manifold!contact!Besse}
are orbits of an $S^1$-action\index[terms]{action!$S^1$-}; the corresponding\index[terms]{form!contact!Besse}
contact form is called {\bf the Besse contact form}.

\hfill

Clearly, the Reeb orbits of a Besse contact manifold
are compact. The converse assertion is also true:\index[terms]{theorem!Sullivan}
by \index[persons]{Sullivan, D.} Sullivan \cite{_Sullivan:geodesics_}, for any contact
manifold $S$ with a fixed contact form, there exists a Riemannian metric
such that the Reeb orbits of $S$ are geodesic. By \index[persons]{Wadsley, A. W.} Wadsley
(\cite{_Wadsley_}), this implies that any contact\index[terms]{theorem!Wadsley}
manifold with compact Reeb orbits admits an $S^1$-action\index[terms]{action!$S^1$-}
tangent to the Reeb orbits. For Sasakian manifolds, this 
theorem is proven in \ref{_Reeb_action_factorizes_Corollary_}.

In \cite{_Kegel_Lange_}, \index[persons]{Kegel, M.} Kegel and Lange prove a
version of the Boothby--Wang theorem for \index[persons]{Besse, A.} Besse manifolds.\index[terms]{manifold!Besse}
Like its Sasakian version
\ref{_quasireg_Sasakian_orbibundles_Theorem_},
\index[persons]{Kegel, M.} Kegel and \index[persons]{Lange, C.} Lange claim that any Besse contact manifold
is obtained as the space of unit vectors
in a complex Hermitian line $L$ bundle over 
a prequantized symplectic orbifold $(X, \omega, L, \nabla)$.

We expect that the methods used below to count
the closed Reeb orbits on Sasakian manifolds 
can be naturally extended to contact manifolds 
that admit a Besse contact form.

\subsection{Weinstein conjecture for Sasakian manifolds}

Earlier, we were counting closed
leaves of the canonical foliation\index[terms]{foliation!canonical}
using the action of the group generated
by the Lee\index[terms]{Lee field} and the anti-Lee field.\index[terms]{Lee field!anti-}
The same method can be applied to study the Reeb
dynamics on Sasakian manifolds. The problem is related to
the \index[persons]{Weinstein, A.} Weinstein conjecture in contact geometry\index[terms]{geometry!contact}:\index[terms]{conjecture!Weinstein}

\hfill

\noindent{\bf Weinstein Conjecture:} On any closed contact
manifold $(N,\eta)$ the Reeb field has at least one closed
orbit.

\hfill

For history and background on Weinstein conjecture, 
see the survey \cite{_Pasquotto:W_conj_}.

Indeed, with our method, we obtain the following
confirmation of \index[persons]{Weinstein, A.} Weinstein conjecture: 

\hfill

\theorem\label{_Main_} {(\cite{_Rukimbira_},
  \cite{_ov_closed_reeb_})} Let $M$ be a compact Sasakian
manifold of dimension $2n+1$. Then its Reeb field has at
least $n+1$ closed orbits. 

\hfill

This result has been previously obtained by P. \index[persons]{Rukimbira, P.} Rukimbira,
\cite{_Rukimbira_}, for the more general case of K-contact
manifolds using topological tools.

For our proof, let $R$ be the Reeb field of the Sasakian
structure, $(B,I)$  the  CR-structure \index[terms]{structure!CR}subjacent to the
contact one and $G$ the closure of the group generated by
the flow of $R$ in the group of CR-automorphisms
$\Aut(M,B,I)$. As above, there exists a one-to-one
correspondence between 1-dimensional orbits of $G$ and
closed orbits of the Reeb field $R$.

Now remember that any Sasakian Reeb field can be
approximated by a quasi-regular one (see
\ref{_defo_qr_sas_} and \cite{ru}). Let this one be $R'$
and let $X$ the associated projective
orbifold. \index[terms]{vector field!Reeb}

\hfill

\remark The group $G$ acts on $X$ by holomorphic
isometries. Moreover, there is a one-to-one correspondence
between 1-dimensional orbits of the action of $G$ on $M$
and fixed points of the action of $G$ on $X$.

\hfill

We thus reduced the problem of counting closed orbits of
the Reeb field on $M$ to counting fixed points of a group
acting by holomorphic isometries on a projective
orbifold. This, in turn, is equal to the number of zeros
of a generic vector field $\zeta\in\Lie(G)$. The result then follows  again from 
\cite{_Fontanari_}.

\hfill

Observe that by Lefschetz theorem,\index[terms]{theorem!Lefschetz} the sum of all Betti
numbers is at least $\dim_\C X+1$, and in our case
$\dim_\C X=n$, which completes the proof of
\ref{_Main_}. \endproof

\hfill

\remark Our proof does not explicitly use the contact
or Sasakian geometry,\index[terms]{geometry!Sasaki} but transfers the problem to the
framework of complex geometry.

\chapter{Submersions and bimeromorphic maps of LCK manifolds}\label{blowup}

\epigraph{\it
Wu Shi said: ``Dao can't be heard about, as what's heard about is contradictory. Dao can't be observed, as what's observed is contradictory. Dao can't be spoken about, as what's spoken about is contradictory. If you recognize its shape, that shape isn't its shape. Dao can't be regarded as something that can be named.''}
{\sc\scriptsize Chuang Tz\v u Chapter 22: Knowledge Wanders North}

The existence of holomorphic maps from LCK manifolds is\index[terms]{map!bimeromorphic}
rather restricted, see also \cite[Chapter 10]{do}. In this
chapter we shall be concerned with  holomorphic
submersions and study 
non-K\"ahler LCK manifolds that admit non-trivial holomorphic
submersions to a complex manifold. The main result relies
on a simple topological criterion that we prove below. The
same criterion is useful in discussing the existence of
LCK metrics on blow-ups of LCK manifolds, and this is the
reason of including the results on birational maps of LCK
manifolds in this chapter. It is known that a blow-up of a
K\"ahler manifold is K\"ahler. However,  while the blow-up
at points of an LCK manifold remains LCK,  the blow-up
along submanifolds is more complicated. Only a special
class of submanifolds can be blown-up to an LCK manifold. \index[terms]{blow-up}

\section{A topological criterion}\label{_Topological_Criterion_}

The following result says that, under some natural
topological restrictions, existence of a holomorphic map
from a compact LCK manifold forces the source manifold to
be GCK.

\hfill

\proposition  {(\cite{ovv1}, \cite{opv})}\label{fibr} 
Let $M$ be an LCK manifold, let $B$ be a path-connected topological space and
let $\pi: M \rightarrow B$ be a continuous, proper map. Assume that either:
\begin{description}
\item[(i)] $B$ is an irreducible complex variety, and
  $\pi$ is holomorphic, or 
\item[(ii)] $\pi$ is a locally trivial fibration with fibres complex subvarieties
of $M$.
\end{description}
Suppose that the map
$$\pi^*:H^1(B)\rightarrow H^1(M)$$
is an isomorphism and the fibres of $\pi$ are positive 
dimensional. Then the LCK structure\index[terms]{structure!LCK} of $M$ is GCK.

\hfill

\proof  Let $\theta$ be the Lee form\index[terms]{form!Lee} of $M$, and 
$\hat M$ be the minimal GCK covering of $M$ on which $\theta$ is exact. 
Since $H^1(B)\cong H^1(M)$, there exists a covering $\hat B \arrow B$
such that the following diagram is commutative, {and the  fibres  of $\hat \pi$ are compact}:
\[
\begin{CD}
\hat M @>>> M\\
@V{\hat \pi}VV @VV{\pi}V\\
\hat B @>>> B
\end{CD}
\]
Let $\hat B_0\subset \hat B$
be the set of regular values of $\hat \pi$,
and let $F_b:= \hat \pi^{-1}(b)$ be the
regular fibres of $\hat \pi$, $\dim_\C F_b =k$. Since
$B$ is irreducible, the set 
$B_0$ is connected, all $F_b$ represent
the same homology class in $H_{2k}(\hat M)$.

Denote the K\"ahler form \index[terms]{form!K\"ahler}of $\hat M$ by $\hat \omega$. It is conformally equivalent to $\pi^*\omega$.

Since all $F_b$ represent the same homology class, the
Riemannian volume
\[
\Vol_{\hat \omega}(F_b) := \int_{F_b}\hat\omega^k
\]
is independent of $b\in B_0$. This gives 
\[
\Vol_{\hat \omega}(F_b)= \int_{F_b}\hat\omega^k
= \int_{F_{\gamma^{-1}(b)}}\gamma^*\hat\omega^k =
\int_{F_{\gamma^{-1}(b)}}\chi(\gamma)^k \hat\omega^k=
\chi(\gamma)^k\Vol_{\hat \omega}(F_b),
\]
where $\chi:\pi_1(M)\rightarrow \R^{>0}$ is the homothety character\index[terms]{homothety character}
(\ref{_homothety_character_Definition_}). This proves that the constant $\chi_\gamma$ is equal to 1 for
all $\gamma \in \Gamma$. Therefore, $\hat\omega$
is $\Gamma$-invariant, and $M$ is globally conformally
K\"ahler. \endproof

\hfill

The same argument also proves the following
observation.

\hfill

\corollary\label{_fibres_chi_non-trivial_Corollary_}  
Let $\pi:\; M \arrow B$ be a proper surjective
holomorphic map, with the fibres of positive dimension,
and base $B$ an irreducible complex variety.
Assume that $M$ is LCK. Let $\chi:\; \pi_1(M) \arrow\R^{>0}$
be the homothety character,\index[terms]{homothety character} and $F$ a general fibre of $\pi$.
Then either $M$ is GCK or the restriction
$\chi\restrict{\pi_1(F)}$ is non-trivial.

\hfill

\proof
Assume that $M$ is not GCK, but $\chi\restrict{\pi_1(F)}$
is trivial. Then the minimal K\"ahler cover $\hat M$ of $M$
is projected to a cover $\hat B$ of $B$ with compact
fibres that are  isomorphic to the fibres of $\pi$.
We denote this projection by $\hat \pi:\; \hat M \arrow \hat B$.
The general fibres of $\hat \pi$ are homologous, and hence 
they have the same K\"ahler volume. This is impossible,
because the deck transform group acts on these fibres
and multiplies the K\"ahler form\index[terms]{form!K\"ahler} of $\hat M$ by a 
non-unit constant.
\endproof

\hfill

The following statement is proven by the same method.

\hfill

\theorem\label{_LCK_no_families_subvarieties_Theorem_}
Let $M$ be a compact LCK manifold, and $S\subset M$
a complex subvariety. Assume that the Barlet\index[terms]{space!Barlet}
space ${\goth B}$ of deformations of $S$ inside $M$ is compact,
and the union of all $S_1\in {\goth B}$ contains\index[terms]{deformation}
an open subset of $M$. Then the restriction of the
weight bundle to $S$ is non-trivial; in particular,
$S$ cannot be simply connected.\index[terms]{bundle!weight}

\hfill

\pstep
Let $\tilde M \stackrel \pi \arrow M$ be the smallest K\"ahler cover\index[terms]{cover!K\"ahler}.
Assume that the restriction of the weight bundle to $S$ is trivial. 
Let $\chi:\; \pi_1(M)\arrow \R^{>0}$ be the homothety
character associated with the weight bundle. Then
an element $[\gamma]\in \pi_1(M)$ is lifted to a closed
loop in $\tilde M$ if and only if $[\gamma]\in \ker \chi$.
Since $\chi\restrict {\pi_1(S)}$ is trivial,
any closed loop in $S$ is lifted to a closed loop
in $\tilde M$. Therefore, 
$\pi^{-1}(S)$ is a union of disjoint copies of $S$.

\hfill

{\bf Step 2:}
Let ${\goth B}_m$ be the incidence variety
associated with ${\goth B}$, that is the variety
of all pairs $S' \in {\goth B}$, $m \in S'$
(in Section \ref{_Barlet_Section_}, we call it 
{\em the marked \index[persons]{Barlet, D.} Barlet space}). 

The Barlet family ${\cal B}$ is compact,
hence the forgetful map $\tau$ from ${\goth B}_m$ to $M$
taking $(S', m)$ to $m\in M$ is proper.
By  Remmert proper mapping theorem  (\cite[Theorem V.C.5]{_Gunning_Rossi_}, \cite[\S II.8.2]{demailly}),\index[terms]{theorem!Remmert's proper mapping} 
the image $\tau({\goth B}_m)$ is a complex subvariety in $M$;
since $\tau({\cal B}_m)$ contains an open subset, its image is $M$.
Therefore, every point $m\in M$ is contained in a deformation of $S$.

Let $S_1$ and $S_2\subset \tilde M$ be two connected components
of the set $\pi^{-1}(S)$, both isomorphic to $S$. 
Since $\pi$ is the minimal K\"ahler cover, 
the deck transform of $\tilde M$ mapping $S_1$
to $S_2$ is a non-trivial homothety.
Let $\tilde\gamma:\; [0, 1] \arrow \tilde M$ be a path-connecting a point $x_1\in S_1$ 
to $x_2\in S_2$, with $\pi(x_1)=\pi(x_2)=x$, and $\gamma:\; [0,1]\arrow M$
the corresponding closed loop. Denote by ${\cal S}_t\subset {\cal B}$
the set of all $S'\in {\goth B}$ containing $\gamma(t)$. Since
the family ${\cal S}_t$ is continuous
and each ${\cal S}_t$ is compact, a finite power of the loop $\gamma$ 
is lifted to a loop in ${\cal B}$. This gives a  continuous
family of subvarieties of $\tilde M$ connecting $S_1$ to $S_2$.
However, the volume of $S_1$ is not equal to the volume of $S_2$,
because the deck transform of $\tilde M$ mapping
$S_1$ to $S_2$ is a non-trivial homothety.
This gives a contradiction, since a deformation
of a complex subvariety in a K\"ahler manifold
has constant volume by Stokes' theorem.
\endproof

\hfill

\definition\label{_uniruled_Definition_}\index[terms]{manifold!uniruled}
A complex manifold $M$ is called {\bf uniruled}
if there is a complex manifold $N$ and a surjective
map $N \times \C P^1\stackrel u \arrow M$ such that 
the restriction of $u$ to $\{n\}\times \C P^1$
is non-trivial for any $n\in N$. In other words,
$M$ is uniruled if it is covered by a family
of rational curves.

\hfill

\corollary\label{_non-uniruled_Corollary_} 
Let $M$ be a compact uniruled LCK manifold.
Assume that the family $N$ of the rational curves
in \ref{_uniruled_Definition_} can be chosen compact. 
Then $M$ is globally conformally K\"ahler.\index[terms]{manifold!uniruled}

\proof
Since $\pi_1(\C P^1)=0$, the weight bundle restricted to $\C P^1$
is trivial. Therefore, \ref{_LCK_no_families_subvarieties_Theorem_}
can be applied. \endproof

\hfill

As another  application, consider the twistor spaces. Twistor
spaces are complex manifolds associated with
even-dimensional differentiable manifolds. There are
several versions:
\begin{description}
\item[(i)]  twistor\index[terms]{twistor space}
spaces of half-conformally flat 4-dimensional Riemannian
manifolds (\cite{_AHS:self-duality_}),
\item[(ii)] twistor spaces of quaternionic-K\"ahler,
  hyperk\"ahler or hypercomplex manifolds (\cite{_Salamon_,besse}), and 
\item[(iii)] Riemannian twistor spaces of conformally flat
  manifolds (\cite{_Slupinski:twistors_,gau_weyl}).
\end{description}

Note that any of the above twistor spaces  is a locally trivial fibration
$Z \arrow B$, with complex analytic fibres that are 
simply connected K\"ahler spaces, and hence
\ref{fibr} can be applied and gives:

\hfill

\corollary 
If $M$ is one of the above twistor space, and if  $M$ admits an LCK
metric, then this metric is actually GCK.
\endproof

\hfill

\remark 
It was known from \cite{gau_weyl, oleg} that the natural metrics (with
respect to the twistor submersion) cannot be LCK. This result refers to any
metric on the twistor space, not necessarily related to the twistor
submersion.  On the other hand, as shown by \index[persons]{Hitchin, N.} Hitchin, the twistor space
of a compact 4-dimensional manifold is not of
K\"ahler type, unless it is
biholomorphic to $\C P^3$ or to the flag variety $F_2$,  \cite{hitchin}.


\section{Holomorphic submersions}

\subsection{LCK metrics on fibrations}


The problem of existence of holomorphic submersions from compact LCK manifolds is not solved in full generality. Examples exist: if the canonical foliation\index[terms]{foliation!canonical} $\Sigma$ of a compact Vaisman  manifold $M$ is quasi-regular,\index[terms]{manifold!Vaisman!quasi-regular} then $M$ admits a holomorphic submersion (that is  moreover a principal bundle map) over a projective orbifold.  Here we prove:\index[terms]{bundle!principal}

\hfill

\proposition {(\cite{opv})}\label{holsub} Let $M$ be a compact complex manifold that admits a holomorphic submersion  $\pi:M\ra B$ with  fibres of complex dimension at least 2. Assume one of the fibres of $\pi$ is of K\"ahler type. If $M$ has an LCK metric $g$, then $g$ is GCK.

\hfill

 

\proof Let $F_0$ be a fibre of K\"ahler type and  $F$  any fibre of $\pi$. Also, let $i_0$ and $i$ be the respective immersions of the fibres in $M$.  As $B$ is path-connected, from Ehresmann's fibration theorem it follows that $F$ and $F_0$ have the same homotopy type, and $\theta\restrict{F_0}$ is exact if and only if $\theta\restrict{F}$ is exact.

By Vaisman's  Theorem (\ref{vailcknotk}) the restriction of the LCK metric $g$ to the fibre is actually GCK if $\dim_\C F_0\geq 2$. Hence, if $F_0$ has dimension at least $2$,  it follows that $i_0^*([\theta])=0.$ From the exact sequence:
$$
\begin{CD}
0 @>>> H^1(B) @ >\pi^*>> H^1(M) @>i^*>> H^1(F),
\end{CD}
$$
 we see $[\theta]$ is a pullback; then \ref{fibr} implies that  $g$ is GCK.   \endproof

\hfill

\corollary\label{no_kaehler_fibres}
 A compact, connected
LCK manifold, that is  not GCK, cannot be the total space
of a locally trivial fibration with K\"ahler fibres of
dimension $\geq 2$.



\hfill

\remark 
If the fibres are genus 0 curves, then the above proof applies \textit{ad litteram}, since in this case, $H^1(F)=0$.

\hfill

\proposition
Let $M$ be a compact LCK manifold equipped with a smooth
elliptic fibration $M \arrow B$. Then $M$ is a 
finite quotient of a Vaisman manifold.\index[terms]{manifold!Vaisman}

\hfill

\proof
Since the $j-$invariant is a holomorphic map to $\C$ and $B$
is compact, we see that all the fibres are isomorphic; 
hence, by Fischer-Grauert theorem \cite{fg}, the map $\pi$
is a locally trivial fibration. Locally trivial fibrations
are classified by the first cohomology of the base with coefficients in 
the automorphism group of the fibre. Replacing $M$
by a finite cover, we may reduce the structure group
of this fibration to the group of translations of the
elliptic curve. Replacing $M$ by its finite cover if
necessary, we can assume that $M$ is a principal elliptic
bundle, that is, the fibre $F$  (that is  identified with 
a complex torus $T^1_{\CC}$) acts holomorphically on $M$. 
Replacing the LCK metric on $M$ by a $T^1_{\CC}$-invariant
metric as in Subsection \ref{_averaging_subsection_},
we may assume that $T^1_{\CC}$ acts on $M$ by isometries.

Let $\hat M$ be the minimal K\"ahler cover of $M$,
and $\chi:\; \pi_1(M) \arrow \R^{>0}$ the homothety
character. If $\chi\restrict {\pi_1(F)}$ is trivial,
$M$ is GCK by \ref{_fibres_chi_non-trivial_Corollary_}.
Therefore, this restriction is non-trivial, and
the preimage of $F$ in $\hat M$ is non-compact.
This implies that at least one of the $S^1$-actions \index[terms]{action!$S^1$-}
on $M$ associated with the principal $T^1_{\CC}$-structure
is lifted to a proper action of $\R$ on $\hat M$.
Let $\sigma\in \R \subset \Aut(\hat M)$ be a non-trivial
automorphism obtained from this flow, that acts on $M$
trivially. Then $\sigma$ belongs to the deck transform
group $\Gamma$ of $\hat M$ over $M$. Since $\hat M$ is the
minimal
K\"ahler cover, all elements of $\Gamma$ act on $\hat M$
non-isometrically. We obtain that the principal torus
action on $M$ is lifted to a non-isometric
conformal holomorphic action on $\hat M$.

Now, it follows from \ref{kami_or} that $M$ is Vaisman.
\endproof


\subsection{LCK metrics on products}\label{noprod}
The above result has implications on the existence of LCK metrics on products of compact LCK  manifolds.

It is clear that the product of two LCK metrics is not LCK, but this does not prohibit the existence of an LCK metric on a product of two LCK manifolds. Still, by the above result, two cases can be excluded:

\hfill

\corollary  {(\cite{tsu}, \cite{opv})}\label{no_lck_on_prod}
Let $M_1, M_2$ be compact quasi-regular Vaisman
manifolds.\index[terms]{manifold!Vaisman!quasi-regular} Then the complex manifold $M_1\times M_2$
carries no LCK metric.

\hfill

\proof Let the product $M_1\times M_2$ have an LCK structure\index[terms]{structure!LCK}. Since $M_1$ and $ M_2$ are regular compact Vaisman manifolds,\index[terms]{manifold!Vaisman} they are total spaces of  holomorphic submersions $\pi_i:M_i\ra B_i, i=1, 2$ onto (compact Hodge) manifolds $B_1, B_2$ with fibres elliptic curves $F_1, F_2$. Then 
$$\pi:M_1\times M_2\ra B_1\times B_2, \qquad \pi(x_1, x_2)=(\pi_1(x_1), \pi_2(x_2))$$
 is a holomorphic submersion with typical fibre $F_1\times F_2$ that is  a $2$-dimensional torus thus  of K\"ahler type. As the first Betti number of a compact Vaisman manifold is odd,  $b_1(M_1\times M_2) $ is even and we see that $M_1\times M_2$ is not biholomorphic to a Vaisman manifold.\index[terms]{manifold!Vaisman} Then, from the \ref{holsub} it follows that $M_1\times M_2$ is of K\"ahler type. This forces $M_1, M_2$ to be of K\"ahler type as well, a contradiction. \endproof

\hfill

Moreover, a direct application of \ref{holsub} proves:

\hfill

\corollary {(\cite{opv})}\label{homvai}
 The product of a compact LCK non-K\"ahler manifold with a compact K\"ahler manifold of complex dimension at least $2$ admits no  LCK metric.

\hfill

Instead, when one of the factors is a curve, we have:

\hfill

\theorem (\cite{nico2}) The product of a compact complex curve with a compact complex manifold $M$ admits an LCK metric only if $M$ is LCK with potential.\index[terms]{manifold!LCK!with potential}

\hfill

\proof Let $C$ denote the compact curve, and $P=C\times M$
the product  manifold, equipped with an LCK structure\index[terms]{structure!LCK}
$(I,\omega,\theta)$. Let $\pi_M, \pi_C$ be the natural
projection maps $P\arrow M, C$. We decompose
$\theta$ as
$\theta = \theta_C+\theta_M$, with $\theta_M \in
\pi_M^*\Lambda^1(M)$ and $\theta_M \in
\pi_M^*\Lambda^1(M)$. After a conformal change of LCK metric, we
may suppose that $\theta_C=\pi_C^* \hat\theta_C$, where
$\hat \theta_C$ is a closed 1-form on $C$ and
that it is the real part of
a holomorphic one-form (and hence
$d^c\theta_C=dI\theta_C=0$). 

We decompose $\omega$ according to the grading induced by
the identification 
\[ \Lambda^2(P)=\bigg(\pi_C^*\Lambda^2(C)\bigg) \oplus
\bigg(\pi_C^*\Lambda^1(C)\otimes_{C^\infty (P)}
\pi_M^*\Lambda^1(M)\bigg) \oplus
\bigg(\pi_M^*\Lambda^2(M) \bigg),
\] 
since $\omega=\omega_C+\omega'+\omega_M$. Now, identifying
also $d=d_C+d_M$, the equation
$d\omega=\theta\wedge\omega$ gives:
\begin{equation}\label{_d_c_omega'_+_d_m_omega_C_Equation_}
d_C\omega'+d_M\omega_C=\theta_M\wedge\omega_C+\theta_C\wedge\omega'.
\end{equation}
Indeed,
\eqref{_d_c_omega'_+_d_m_omega_C_Equation_}
is the part of $d\omega=\theta\wedge \omega$
that belongs to the component \[ \pi_c^*\Lambda^2(C)\otimes_{C^\infty (P)}
\pi_M^*\Lambda^1(M)
\] in the decomposition
\begin{align*} \Lambda^3(P)=&\bigg(\pi_C^*\Lambda^3(C)\bigg) \oplus
\bigg(\pi_M^*\Lambda^3(M) \bigg)\oplus
\bigg(\pi_C^*\Lambda^2(C)\otimes_{C^\infty (P)}
\pi_M^*\Lambda^1(M)\bigg)\\ \oplus &
\bigg(\pi_C^*\Lambda^1(C)\otimes_{C^\infty (P)}
\pi_M^*\Lambda^2(M)\bigg).
\end{align*}
Let ${\goth p}:\;\Lambda^{p,q}(P)\arrow \Lambda^{p-1,q-1}(M)$ 
be the push-forward of forms (by fibrewise integration):
$$\goth p(\al)_m:=\int_{C\times\{m\}}\al.$$
Let $f:=\goth p(\omega_C)$ that is  a strictly positive smooth function on $M$. Then 
  the push-forward of \eqref{_d_c_omega'_+_d_m_omega_C_Equation_} to $M$ reads:
$$d_{\theta_M}f=\goth p(\theta_C\wedge\omega'),$$
where $d_{\theta_M}(\beta)= d_M - \theta \wedge \beta$ is
the Morse--Novikov differential on $M$.
Applying $d^c$
and taking care of the components on each factor, we obtain:
$$d^cd_{\theta_M}f-I\theta_M\wedge d_{\theta_M}f=-\goth p(\omega_M\wedge\theta_C\wedge I\theta_C).$$
Let $\al:=\omega_M\wedge\theta_C\wedge I\theta_C$, and $d^c_{\theta_M}:= Id_{\theta_M}I^{-1}$. Then we
have $$\bar\omega_M:=\goth
p(\al)=d_{\theta_M}d^c_{\theta_M}f,$$ and thus 
$\bar\omega_M$ is an LCK form with potential\index[terms]{form!LCK!with potential}, with Lee
form $\bar\theta_M:=\theta_M-d\log f$, provided $\goth p(\al)$ is strictly positive. This follows because $\al$ is semipositive on $P$ and strictly positive on $U\times M$, where $U:=\{x\in C\; ;\; \theta_C\restrict{x}\neq 0\}$. \endproof

\section{Blow-up at points}\label{_Blow_up_at_points_section_}
\index[terms]{blow-up}

It is a classical fact that blowing-up a K\"ahler manifold at a finite number of points  produces a new K\"ahler manifold, \cite{griha}. LCK manifolds behave the same. The result was announced in \cite{tric} and completely proven in \cite{vu1}. 

\hfill

Let $(M,I,\omega,\theta)$ be an LCK manifold and let $(\check{M}, I)$ be its blow-up at a point $P\in M$, considered to be  a complex manifold.  Let $E:=\pi^{-1}(P)$ be the exceptional divisor.\index[terms]{divisor!exceptional}

\hfill

\lemma {(\cite{tric})}\label{tricc}
There exists an open
neighbourhood $U$ of $E$ and a smooth function $f:\check{M}\arrow \R$ such that $\omega':=
e^f\pi^*\omega$ satisfies $d\omega'=\theta'\wedge \omega'$ and such that
$\theta'\restrict{U}\equiv 0$ for some 1-form $\theta'$ on $\check{M}$.
\endproof

\hfill

\theorem\label{_blow_up_is_LCK_Theorem_} 
{(\cite{tric,vu1})}
The blow-up at a point of an LCK manifold is LCK.

\hfill

\proof In the above notations, we observe that
$\pi^*\omega$  is a  $(1,1)-$ form on $\check{M}$, which 
 is positive definite outside of the exceptional divisor and satisfies 
 $d\pi^*\omega=\pi^*\theta\wedge \pi^*\omega$.
 
 We want to glue it with a metric defined around the
 exceptional divisor, such that the result is
an LCK metric. We cannot use the K\"ahler metric of the divisor, as in the K\"ahler case. Instead, we shall work with the curvature form\index[terms]{form!curvature}  of a line bundle on $E$. The construction is standard:
 
\hfill

\claim  {(\cite[pp 185-187]{griha})} \label{gri_ha}
 Let $\calo_{\check{M}}(E)$ be the holomorphic line bundle associated with $E$. There exists a Hermitian metric on $\calo_{\check{M}}(E)$ such that the  curvature  $\Omega_E$ of its
Chern connection\index[terms]{connection!Chern}  satisfies the conditions:
\begin{description}
\item[(i)] $\Omega_E(v, Iv)<0$ for every non-vanishing vector  $v\in T_Q(E)$ and
for every $Q\in E$ (negative  definite along $E$),
\item[(ii)]  $\Omega_E(v,
Iv)\leq 0$ for any $Q\in E$ and any $v\in T_Q(\check{M})$ (negative semidefinite at points of $E$), and 
\item[(iii)] $\Omega_E=0$ outside $U$.  \endproof
\end{description}
 
 Notice that $\Omega_E$ is closed.
We now combine $\omega'$ and $\Omega_E$. Let $N\in \R^{>0}$ and consider the (1,1)-form 
$$\widehat\omega:=N\omega'-\Omega_E.$$
The form $\omega'$  obviously satisfies:
$$d\omega'=\theta'\wedge N\omega'.$$
On the other hand $\theta'\wedge\Omega_E=0$, as $\theta'$ and $\Omega_E$ have disjoint supports, and hence one can write:
$$d\widehat\omega=\theta'\wedge N\widehat\omega-\theta'\wedge\Omega_E=\theta'\wedge\widehat\omega.$$
Then $\widehat\omega$ is  $d_{\theta'}$-closed, and to show it is LCK it remains to prove the following:

\hfill

\begin{claim}
There exists $N\in \R^{>0}$ such that $\widehat\omega$ is positive definite.
\end{claim}

\hfill

We need to check the positivity in three situations:
\begin{enumerate}
\item Outside $U$: here  $\Omega_E=0$ and $N\omega'$ is
positive definite for any $N>0$.
\item Along $E$: here $\omega'$ vanishes when applied to
vectors tangent to $E$ and $\Omega_E\geq 0$ (by
1. above). In general, we decompose a vector field  along
$E$ in tangent and transversal parts and apply the above reasoning.
\item Finally, on $U$: we observe that for each $x\in U$
there exists some $n_x\in \R^{>0}$ such that $N\omega'-\Omega_E$ is positive definite at $x$
for all $N\geq n_x$, and hence, by continuity, also positive definite on a neighbourhood of $x$. As 
$U$ is relatively compact, we can cover it by finitely many such neighbourhoods, and
take $N$ as  the maximum of the numbers $n_x$. \endproof
\end{enumerate}

\section{Blow-up along submanifolds}\label{_Blow_up_on_submanifolds_section_}

The blow-up of points was discussed in Section \ref{_Blow_up_at_points_section_}. Here we treat the more delicate problem of blowing up (complex) submanifolds.

Prior to the discussion on LCK manifolds, we present a general result\index[terms]{blow-up} concerning the blow-up of K\"ahler manifolds along submanifolds, that we shall use:

\hfill

\lemma \label{paun} Let $(U, g)$
 be a K\"ahler  manifold,  $Y\subset U$   a compact submanifold,
and  $c : \tilde U\rightarrow U$  the blow-up of $U$ along $Y$. Then, for any open neighbourhood $V$ of $Y$,
there exists a K\"ahler metric $\tilde g$ on $\tilde U$ such that
$\tilde g\restrict {\tilde U \setminus c^{-1}(V)}
= c^*(g\restrict{U\setminus V})$.

\hfill

\proof  (communicated by Mihai P\u aun; see also
\cite{vu1} or \cite[Chapter 3]{vois}). 
The argument is essentially already used in \ref{gri_ha}. 
If we find a (non-singular) metric on $\calo_{\tilde U} (-D)$  (where $D$ is the exceptional divisor\index[terms]{divisor!exceptional}
of the blow-up) such that:
\begin{description}
\item[(a)] its curvature is zero outside $c^{-1}(V)$, and

\item[(b)] its curvature is strictly positive at every point of $D$ and in any direction tangent to $D$,
\end{description}
then everything follows, as the curvature
of this metric plus a sufficiently large multiple of $c^*(g)$ will be positive
definite on $\tilde U$.

To finish the proof, we note that the existence of a metric $h$ with property
(b) is clear, because the restriction of $\calo _{\tilde U} (-D)$ to $D$ is relatively ample.

Now let $\al$ be the  curvature of $h$. Then $\al -\1 \6\bar\6 \tau =-[D]$ for some function $\tau$ , with at most
logarithmic poles along $D$, bounded from above, and nonsingular on $\tilde U \setminus D$. Consider the
function $\tau_0:= \max(\tau,-C)$, where $C$ is some positive constant, big enough such that on
$\tilde U\setminus c^{-1}(V)$ we have $\tau >-C$. Clearly, on a (possibly smaller) neighbourhood of $D$ we shall have $\tau_0 =-C$, such that the new metric $e^{-\tau_0}h$ on $\calo_{\tilde U} (-D)$ also satisfies (a).
 \endproof

\hfill

Before stating the results, we need to  introduce the following:

\hfill

\definition \label{_IGCK_Definition_}
Let $(M,I,g,\theta)$ be an LCK manifold and let $j:Y\hookrightarrow M$ be a complex subvariety. We say that $Y$ is of {\bf  induced globally
conformally K\"ahler type} (IGCK) if 
the cohomology class $j^*[\theta]$ vanishes, where $[\theta]$  denotes
the cohomology class of the Lee form\index[terms]{form!Lee} on $M$.

\hfill

\remark
A subvariety $Y\subset M$ is IGCK if and only if the weight bundle
$(L, \nabla)$ restricted to $Y$ has trivial monodromy. Indeed, the
monodromy of $\nabla$ on a path $\gamma\subset M$ is $\int_\gamma \theta$.

\hfill

\remark \label{igck}
Note that any IGCK submanifold admits a K\"ahler metric. On the other hand, by \ref{vailcknotk}, on   compact and smooth $Y$ that are  not  curves, the IGCK condition is equivalent to $Y$ being K\"ahler; it follows that this condition is more interesting on curves.

\hfill

\remark \label{notigck}
There exist LCK manifolds with (K\"ahler) curves that are  not IGCK. 
 For instance, if $M$ is a regular Vaisman manifold,\index[terms]{manifold!Vaisman!regular}
and if $Y$ is a fibre of its elliptic fibration, then $Y$ is not 
IGCK, because the minimal K\"ahler cover $\hat M$ of $M$ has K\"ahler potential,
hence $\hat M$ has no compact complex subvarieties; however, the 
pullback of an IGCK variety $Y$ to the minimal K\"ahler cover
is a collection of several copies of $Y$.

\hfill

As the fibres of a blow-up along a submanifold are  positive-dimensional and compact, \ref{fibr} directly implies the following:

\hfill

\theorem {(\cite{ovv1})}\label{bdown}
Let $M$ be a complex variety, and let $\check M\rightarrow  M$ be  the blow-up of a compact subvariety
$Y\subset M$. Assume that $\check M$ is smooth and admits an LCK metric. Then the blow-up
divisor $\tilde Y\subset \check M$ is an IGCK
subvariety.

\hfill

\corollary \label{puf}
Let $M$ be an LCK manifold, and  $Y\subset M$ a smooth compact subvariety,
such that the blow-up of $M$ along $Y$ admits an LCK metric. Then $Y$ is a IGCK subvariety.

\hfill

\proof Indeed, the variety $\tilde Y$  is of K\"ahler
type, because it is IGCK. Then, if  $Y$ is smooth, it  is K\"ahler, as shown in \cite[Th\'eor\`eme II.2]{_Blanchard_} and \ref{igck} finishes the proof. \endproof

\hfill

\remark It is possible (and not very hard) to state and
prove \ref{puf} without assuming that $Y$ is smooth.

\hfill

We now prove the main result of this section:

\hfill

\theorem {(\cite{ovv1})}\label{main_blow}
Let $M$ be an LCK manifold, $Y\subset M$ a smooth complex IGCK subvariety,
and  $c:\check M\rightarrow M$  the blow-up of $M$ along $Y$. Then $\check M$ is LCK.\index[terms]{blow-up}

\hfill

\proof
Let $g$ be an LCK metric on M and  $\theta$  its Lee form.\index[terms]{form!Lee} Since $Y$ is IGCK, $\theta\restrict Y$  is exact.

Let  $U$ be a neighbourhood of $Y$ such that the inclusion $Y\hookrightarrow U$ induces an isomorphism $H^1(Y)\simeq H^1(U)$ 
in the first de Rham cohomology.\index[terms]{cohomology!de Rham} Then $\theta\restrict U$ is  exact too, and hence, after possibly conformally rescaling
$g$, we may assume that $\theta\restrict U =0$, implying that $g\restrict U$ is K\"ahler. In particular, $\supp(\theta)\cap  U =\emptyset$. We now 
choose a smaller neighbourhood $V$ of $Y$ and apply
\ref{paun}. This yields a K\"ahler metric $\tilde g$ on 
the blow-up of $U$ which equals $c^*(g)$ outside $c^{-1}(V)$, so that it glues to $c^*(g)$ and gives  an LCK metric on $\check M$.
\endproof

%
\hfill

\remark {(\cite[Main Theorem 1.1]{yz})} \\
The Morse--Novikov cohomology groups of $(M,\theta)$ and $(\check M, \check\theta)$ are related (through the de Rham cohomology \index[terms]{cohomology!de Rham}of $Y$) by the formula:
$$H^k_\theta(M)\oplus\left(\bigoplus_{i=0}^{r-2}H^{k-2i-2}_{dR}(Y)\right)\simeq H^k_{\check\theta}(\check M).$$
The proof uses a version of Mayer--Vietoris sequence for Morse--Novikov cohomology.\index[terms]{exact sequence!Mayer--Vietoris}\index[terms]{cohomology!Morse--Novikov}

\hfill

\corollary {(\cite{ovv1})}
The blow-up of a compact Vaisman manifold\index[terms]{manifold!Vaisman} along a compact complex
submanifold $Y$ of dimension at least 1 cannot have an LCK metric.
\index[terms]{blow-up}

\hfill

\proof A subvariety of a Vaisman manifold is tangent to
its  canonical foliation\index[terms]{foliation!canonical} $\Sigma$ (\ref{_Subva_Vaisman_Theorem_} (iii)). Therefore, it is either
a Vaisman manifold (\ref{_Subva_Vaisman_Theorem_} (iv)) or an elliptic curve $T^2$ which\index[terms]{manifold!Vaisman}
coincides with the leaf of $\Sigma$. In the first case
$\theta$ is not cohomologous to 0 because an LCK manifold
is not K\"ahler (Vaisman's  \ref{vailcknotk}). In the
second case, $\theta$ is constant on $T^2$; hence  it
is also non-exact. In both cases, $Y$ is not IGCK.
\endproof

\hfill

\remark
The blow-up at points and along submanifolds of LCS manifolds\index[terms]{manifold!LCS} was treated with similar techniques. It is shown in \cite{cy} that the blow-up at points preserves the LCS class, while in \cite{yyz} it is proven that the blow-up along symplectic submanifolds of an LCS manifolds always carries LCS forms.\index[terms]{form!LCS}
\index[terms]{blow-up}

\section{Weak LCK structures}\index[terms]{structure!weak LCK}

The LCK definition can be weakened such that the blow-up along submanifolds of strictly positive dimension  preserve this weaker condition.

\hfill

\definition {(\cite{_Angella_Parton_Vuletescu_})} 
A {\bf weak locally conformally K\"ahler} (WLCK) structure
on a complex manifold $(M,I)$ is given by a (1, 1)-form $\omega$ and a real 1-form
$\theta$ such that:
\begin{itemize}\index[terms]{structure!weak LCK}
	\item $d\omega=\theta\wedge\omega$ and $d\theta=0$;
\item $\omega$ is strictly positive definite outside a proper analytic subset.
\end{itemize}

\remark A \index[persons]{Calabi, E.} Calabi-Eckmann manifold $M=S^{2n-1}\times
S^{2m-1}$, $n,m\geq 2$,  does not admit any WLCK
structures. Indeed, if $\omega$ is WLCK on $M$, then
$\theta$ is exact and, after a conformal change, we
can assume $d\omega=0$. Then $\omega$ is exact
(since $b_2(M)=0$) and thus $\int_M\omega^{m+n-1}=0$. But
$\omega$ is strictly positive outside the bad locus which
has Lebesgue measure 0, and thus  
$\int_M\omega^{m+n-1}>0$,
 contradiction.\index[terms]{manifold!Calabi-Eckmann}

%

\hfill

\proposition {(\cite{_Angella_Parton_Vuletescu_})}\label{_blow_up_of_weak_is_weak_}
Let $\sigma: \hat M\arrow M$ be a proper modification (e.g. a composition
of blow-ups). If $M$ admits a WLCK metric, then $\hat M$ admits a WLCK metric. \endproof

\hfill

\ref{vailcknotk}  can be extended to the weak LCK structures as follows.\index[terms]{structure!weak LCK}

\hfill

\proposition {(\cite{_Angella_Parton_Vuletescu_})}
Let $M$ be a compact complex manifold, $\dim_\C M\geq 2$, endowed
with a WLCK structure $(\omega,\theta)$. If $M$ admits a K\"ahler metric, then $\theta$ is exact. \endproof

\hfill

\corollary {(\cite{_Angella_Parton_Vuletescu_})}\label{_Vaisman_for_weak_}
Let $M$ be an LCK manifold, $Z\subset M$ a compact complex subspace
(with reduced structure),  $\dim_\C Z\geq 2$. If the desingularization $\hat Z$ is of K\"ahler
type, then the restriction $\theta\restrict Z$ is exact.

\hfill

\proof Consider the blow-up $\hat M$ of $M$ along the singularities of $Z$. By \ref{_blow_up_of_weak_is_weak_},  $\hat M$ is WLCK. Since $\hat Z\subset \hat M$ and $\mathrm{\sf sing}(Z)\neq Z$ (because $Z$ is reduced) we see that $\hat Z$ inherits
a WLCK structure,\index[terms]{structure!weak LCK} and \ref{_Vaisman_for_weak_} applies. \endproof

\hfill

\remark   The notion of WLCK structure was used in
\cite{_Angella_Parton_Vuletescu_} to give a bimeromorphic
classification of compact LCK manifolds $X$ of dimension 3
and algebraic 
dimension\footnote{Recall that the {\bf algebraic dimension}
  $a(X)$ of a compact complex manifold $X$ is the  transcendence
  degree of the field of meromorphic functions on $X$. It
  is known that $0\leq a(X)\leq \dim X$.} $a(X)=2$. Before
stating this result, recall: \index[terms]{dimension!algebraic}

\hfill

\conjecture {(The Strong
		Factorization Conjecture,\index[terms]{conjecture!Strong Factorization} \cite{_Hironaka_resolution_})} 
Any  bimeromorphism 
$\xymatrix{
	\psi \colon X \ar@{<-->}[r] & X^* }$ between compact complex manifolds can be factorized as
{\scriptsize \begin{equation*}\label{eq:alg-red-assumption}
	\xymatrix{
		& \hat X \ar[ld]_{\sigma} \ar[rd]^{c} & \\
		X \ar@{<-->}[rr]_{\psi} & & X^* 
	}
\end{equation*}}
where $\sigma$ and $c$ are (compositions of) blow-ups with smooth centers.

\hfill

\remark Let $X$ be a compact complex manifold. Recall that there exists a projective smooth manifold $B$ and a meromorphic map $f:X\arrow B$ such that
$f^*:\C(B)\arrow{\mathcal M}(X)$ is an isomorphism, where $\C(B)$ is the $\C-$algebra of rational functions on $B$. Usually, $B$ is called an {\bf algebraic reduction} of $X$. Notice that $f$ and $B$ are only defined up to bimeromorphisms. Moreover, if $a(X)=\dim(X)-1$, then the general fibres of the algebraic reduction are (smooth) elliptic curves. For example, if $X$ is a Hopf surface,\index[terms]{surface!Hopf} then the map $f:X\arrow \C P^1$, given by $f(z_1,z_2)=[z_1, z_2]$ is the algebraic reduction of $X$. In fact, this map gives $X$ a structure of elliptic principal bundle over $\C P^1$. 

\hfill

\theorem {(\cite{_Angella_Parton_Vuletescu_})}
Let $(M,I,\omega,\theta)$ be a compact LCK manifold with algebraic dimension $a(M)=2$. Assume the Strong Factorization Conjecture, and let $F: M^*\arrow B$ be the algebraic reduction of $M^*$:

\vspace{-3mm}

{\scriptsize \begin{equation*}
	\xymatrix{
		& \hat M \ar[ld]_{\sigma} \ar[rd]^{c} \\
		M \ar@{<-->}[rr]_{\psi} & & M^* \ar[d]	_{f} \\
		& & B .
	}
\end{equation*}
}

\noindent Then all the smooth one-dimensional fibres of $f$ are isomorphic. If, moreover, (a) the morphism $f$ is flat, (b) $M^*$ can be chosen such that $K_{M^*}$ be nef, and (c)  the singular locus  of the fibration $f$ is a simple normal crossing divisor, then $f$ is a quasi-bundle.\footnote{Recall that a fibration is a quasi-bundle (or isotrivial) if its singular fibres  are all multiples of the general fiber, \cite{_Brinzanescu:bundles_}.} Finally, if $\hat M$ is LCK, then there exists a bimeromorphic morphism $\varphi:M\ra X^*$
which  is a proper modification of $M^*$ at points. \endproof\index[terms]{modification}\index[terms]{divisor!normal crossing}
\index[terms]{quasi-bundle}\index[terms]{bundle!vector bundle!nef}\index[terms]{reduction!algebraic}


\section{Moishezon manifolds are not LCK}


Recall that the algebraic dimension $a(M)$ of a compact complex manifold $M$
satisfies $a(M) \leq \dim_\C (M)$, and the equality is realized
if and only if $M$ is bimeromorphic to a projective manifold (\cite{moi}).
Compact complex manifolds that satisfy $a(M)=\dim_\C M$ are called {\bf Moishezon}.
This is a subclass of the Fujiki class C manifolds, that is, compact
complex manifolds that are  bimeromorphic to K\"ahler.
As shown in \cite[Corollary 5.23]{_DGMS:Formality_}, the $dd^c$-lemma is 
satisfied for the Moishezon manifolds and the Fujiki class C manifolds.
We extend Vaisman's \ref{vailcknotk} to the Moishezon 
and the Fujiki class C manifolds.

\hfill

The following theorem is probably due to \cite{au} (see
also \ref{vai_gen_au}).

\hfill

\theorem\label{_dd^c_then_non-LCK_Theorem_}
Let $M$ be a compact complex manifold satisfying the
$dd^c$-lemma, for instance, a manifold of Fujiki class C.
Then all LCK structures\index[terms]{structure!LCK} on $M$ are globally conformally K\"ahler.

\hfill

\proof
The proof repeats the proof of \ref{vailcknotk},
Let $\theta$ be the Lee form\index[terms]{form!Lee} of an LCK structure on a compact Fujiki class C manifold $M$.
Since $M$ satisfies the $dd^c$-lemma, $\theta$ is cohomologous to a real part of
a holomorphic 1-form. After a conformal change, we can assume that 
the Lee form of $M$ is $d^c$-closed. Then
$dd^c(\omega^{n-1})= (n-1)^2 \omega\wedge \theta \wedge \theta^c$.
This form is positive, and, unless $\theta=0$, it
satisfies $\int_M dd^c(\omega^{n-1})>0$, that is 
impossible by Stoke's theorem. Therefore, $\theta=0$,
and $M$ is globally conformally K\"ahler.
\endproof


\section{LCK currents and Fujiki LCK class}\index[terms]{Fujiki LCK class}


\subsection{LCK manifolds in terms of currents}

Let $M$ be a smooth manifold, $\dim_\R M=n$.
Consider the space $\Lambda^{n-r}_c(M,L)$ 
of compactly supported smooth differential
forms with coefficients in a line bundle $L$ and $C^\infty$-topology,\index[terms]{topology!$C^\infty$}
that is, the Fr\'echet topology\index[terms]{topology!Fr\'echet} of convergence of all derivatives.
The {\bf space of $r$-currents with coefficients in $L$}
is the space of functionals on $\Lambda^{n-r}_c(M,L^*)$
that are  continuous in each of the $C^i$-topologies. The space of
currents is equipped with the weak-$*$ topology (a sequence of currents
converges if it converges on each test-form).

The $L$-valued differential 
$r$-forms can be understood as $L$-valued currents:
given an $L$-valued form $T\in \Lambda^{r}(M,L^*)$, 
and an $L^*$-valued compactly supported form
$\alpha\in \Lambda^{n-r}_c(M,L^*)$, we define 
$\langle T, \alpha\rangle:= \int_M T \wedge \alpha$,
where $T\wedge \alpha$ is understood as a $L^*\otimes L$-valued
compactly supported $n$-form, that is a volume form on $M$.

Let $d_\theta:= d - \theta\wedge$ be the $L$-valued de Rham differential,
called {\em the twisted de Rham differential} elsewhere.
We define the twisted de Rham differential on the space
of currents writing 
\[ 
 \langle d_\theta T, \alpha\rangle:= -(-1)^{\tilde T}
\langle T, d_{-\theta} \alpha\rangle.
\]
This is compatible with the embedding of $r$-forms to $r$-currents,
because
\begin{equation*} 
\begin{split}
\langle d_\theta T, \alpha\rangle&= \int_M d_\theta T \wedge \alpha
= - \int_M \theta \wedge T\wedge \alpha - (-1)^{\tilde T}\int_M T \wedge d\alpha\\
&=- (-1)^{\tilde T}\int_M T \wedge d_{-\theta}\alpha,
\end{split}
\end{equation*}
by Stokes' theorem.

\hfill

Further on, we shall trivialize $L$, and consider the differential $d_\theta$
as an operator acting on the usual differential forms and currents.
However, it helps to think of these forms as of $L$-valued or $L^*$-valued forms.

\hfill

\ref{_HL_for_L_valued_Theorem_} below generalizes the intrinsic characterization of K\"ahler manifolds in \cite{hl}. We first state the Hahn-Banach extension theorem
in the form that will  be used in the proof.

\hfill

\theorem \label{_Hahn-Banach_extension_Theorem_}
Let $V$ be a locally convex topological vector space,
$A \subset V$ a closed subspace, $W\subset V$ an open 
convex subset, and $h$ a continuous linear functional on $A$
that is  positive on $W \cap A$. Then $h$
can be extended to a continuous linear functional on $V$ 
that is  positive on $W$.

\proof \cite{_Bourbaki:TVS_}. \endproof

\hfill

\theorem \label{_HL_for_L_valued_Theorem_} 
{(\cite{oti})} 
Let $M$ be a compact, complex manifold, $\dim_\C M\geq 2$, and $\theta$ a closed 1-form on
$M$. Then $M$ admits an LCK metric with Lee form \index[terms]{form!Lee}$\theta$
if and only if it does not admit nontrivial positive
currents that are  $(n-1,n-1)$-components of
$d_{-\theta}$-boundaries.
  
\hfill

\proof  Assume there exist no positive currents that are 
$(n-1,n-1)$ components of $d_{-\theta}$-boundaries. We shall use
the Hahn-Banach theorem to construct an LCK metric on $M$.

Let $V$ be the space of $(1,1)$-forms on $M$ 
with the $C^\infty$-topology.\index[terms]{topology!$C^\infty$}
Let $A$ be the subspace of
$d_{\theta}$-closed $(1,1)$ forms. The subspace $A\subset V$
is clearly closed. Now, let $W$ be the subset of strictly
positive $(1,1)$-forms, that is  convex. If, by absurd,
there is no LCK form on $M$, then $A\cap W=\emptyset$ and the
Hahn-Banach extension theorem applies producing a
functional that vanishes on $A$ and positive on $W$. By definition, this
functional is an $(n-1,n-1)$-current; call it $T$.

Since $T$ vanishes on $A$, we have $T(\eta)=0$ for all
$d_{\theta}$-closed $\eta$. Then
$T$ is $d_{-\theta}$-closed, and then $d_{-\theta}$-exact. Indeed, 
$\langle T, d_{\theta} \alpha \rangle=0$ for all $\alpha$,
which implies that $\langle d_{-\theta} T, \alpha \rangle=0$,
hence $d_{-\theta} T=0$. Assume that $T$ is not $d_{-\theta}$-exact,
and let $[\alpha]\in H^{n-2}_{d_{-\theta}}(M)= H^2(M,L^*)$ be its $d_{-\theta}$-cohomology class.
By  Poincar\'e duality between $H^2(M,L)$ and $H^{n-2}(M,L^*)$,
there exists a $d_\theta$-closed 2-form $\alpha$ such that
$\int_M T \wedge \alpha\neq 0$; this contradicts the
assumption $T\restrict A=0$. Therefore, $T$ is exact.
The converse is immediate.
\endproof
 
\subsection{An analogue of Fujiki class C}
\label{_LCK_class_C_Subsection_}
\index[terms]{Fujiki class C}

\definition 
Let $M$ be a compact complex manifold, $\theta$ a closed real one-form on $M$,
and $\Xi$ a positive, real (1,1)-current satisfying $d\Xi
= \theta\wedge \Xi$. Assume that $\Xi\geq \omega$ for some Hermitian
form $\omega$ on $M$. Then $\Xi$ is called an {\bf LCK current}.

\hfill

\remark \label{fibr_to_currents}
As shown in \cite{shi}, \ref{fibr} still holds for LCK currents. The idea is to decompose an LCK current into regular and singular parts and to use \index[persons]{Demailly, J.-P.} Demailly's regularization to deform the current such as to have singularities only in codimension 2, that will  not affect integration of the current, then rewrite the proof of \ref{fibr}. The following generalizes \ref{main_blow}:

\hfill

\theorem {(\cite{shi})} Let $X$ be a compact complex manifold, $\dim X \geq 3$, and $D\subset X$
a submanifold of $\codim D \geq 2$. Let $f: \tilde X \rightarrow X$ be a blow-up of $X$ along $D$,\index[terms]{blow-up}
and $E$ the exceptional set. Suppose that there exists an LCK current $\tilde T$ on $\tilde X$.
Let $\theta$  be the Lee form\index[terms]{form!Lee} of $\tilde T$. Then $\theta\restrict{E}=0\in H^1(E,\R)$.
In particular, the blow-up of a compact Vaisman manifold along a submanifold of codimension greater than 2 does not admit any LCK current.\index[terms]{manifold!Vaisman}

\hfill

\remark 
A compact complex variety $X$ is said to belong to {\bf Fujiki class C} if $X$ is bimeromorphic to a\index[terms]{map!bimeromorphic}
K\"ahler manifold. The Fujiki class C manifolds are closed under many natural operations,
such as taking a subvariety, or the moduli of subvarieties, and play an important role in
K\"ahler geometry.\index[terms]{geometry!K\"ahler}

By analogy, one can say that  a compact complex variety is 
 {\bf LCK class C}, or {\bf of LCK Fujiki class C}, if it is
bimeromorphic to an LCK manifold. A sufficient condition
for a manifold to be of LCK class C is the following:

\hfill

\theorem {(\cite{shi})} Let $X$ be a compact complex manifold with an LCK current $T$, and $\theta$ the Lee form.\index[terms]{form!Lee} We assume that
there exists an open neighbourhood $V$ of the singular locus $B$ of $T$ such that $\theta$ is exact on $V$. Then $X$ admits a modification $\tilde X \rightarrow  X$  from an LCK manifold $\tilde X$.

\hfill

\remark \label{_DP_Theorem_}
In \cite{dp}, \index[persons]{Demailly, J.-P.} Demailly and P\u aun  proved that a compact complex manifold  belongs to Fujiki 
class C if and only if it admits a positive K\"ahler current (recall that a (1,1) current $T$ is positive if $T(\al)\geq 0$ for all $(n-1,n-1)$-forms $\al$ and K\"ahler if it is closed and $T-\omega\geq 0$ for some Hermitian form\index[terms]{form!Hermitian} $\omega$). On the contrary, in LCK setting, by using \index[persons]{Hironaka, H.} Hironaka's example of non-projective Moishezon manifolds,\index[terms]{manifold!Moishezon} we have:

\hfill

\proposition  {(\cite{shi})}
Let $X$ be an LCK manifold. Suppose that $\dim_\C X \geq 3$. Then
there exists an LCK class C manifold that is bimeromorphic to $X$ but not an LCK
manifold.

\hfill

\remark All Fujiki class C manifolds are balanced\index[terms]{Fujiki class C} (\cite{_Alessandrini_Bassanelli:bimero_}).\index[terms]{metric!balanced}
In \cite{shi} it was shown
that all compact LCK class C manifolds are locally
conformally balanced\index[terms]{metric!locally conformally balanced}
(see Section \ref{lcb} for the definition of locally
conformally balanced manifolds).


\chapter{Bott--Chern cohomology of LCK manifolds with
  potential}

\epigraph{\it PUPIL: I can count to...  to infinity.
	 
	 PROFESSOR: That's not possible, miss. 
	 
	 PUPIL: Well then, let's say to sixteen. 
	 
	 PROFESSOR: That is enough. One must know one's limits. 
Count.}{\sc\scriptsize Eug\`ene Ionesco, \ \  The Lesson}


In this chapter, we give further information about the
Bott--Chern cohomology of LCK manifolds with potential\index[terms]{manifold!LCK!with potential}. The
main result will be a generic vanishing
theorem.\index[terms]{cohomology!Bott--Chern}\index[terms]{cohomology!Dolbeault}

  \section{Bott--Chern {\em versus} Dolbeault cohomology}\index[terms]{cohomology!Bott--Chern}\index[terms]{cohomology!Dolbeault}

The Bott--Chern cohomology were defined in Section \ref{_Bott_Chern_Section_},
The weighted Bott--Chern cohomology (Bott--Chern cohomology
with coefficients in a local system) were defined in Section \ref{twisted_BC}.

 There are natural (and functorial) maps from the
Bott--Chern cohomology to the Dolbeault cohomology
$H^*(\Lambda^{*,*}(M), \bar\6)$ and to the de Rham
cohomology, but no morphisms between de Rham and
Dolbeault cohomology.

In particular, the following generalization of \ref{vailcknotk} is available:

\hfill

\theorem {(\cite{au})}\label{vai_gen_au} 
Let $M$  be a compact complex manifold of complex dimension $n$ such that the natural map $H^{n-1,n}_{BC}(M)\rightarrow H^{n-1, n}_{\bar\6}(M)$ is injective. Then any LCK structure\index[terms]{structure!LCK} on $M$ is GCK.

\hfill

For Vaisman manifolds,\index[terms]{manifold!Vaisman} the Bott--Chern cohomology
can be determined using the transversal cohomology 
(\cite{_Istrati_Otiman:BC_}; see also Chapter \ref{_Harmonic_forms_chapter_}).
 
\hfill

\theorem \label{_BC_kernel_Proposition_}
Let $(M, I, \omega, \theta)$ be a compact Vaisman manifold, and let 
\[ \psi:\; H^{1,1}_{BC}(M)\arrow H^2(M)\]
be the tautological map.\index[terms]{map!tautologic}
Then $\ker \psi$ is 1-dimensional and it is generated by the semi-positive (1,1)-form
$\omega_0:=d^c\theta$.

\hfill

The proof is preceded by a sequence of lemmas. To begin, we need the following:

\hfill

\definition 
Let $\omega$ be a Gauduchon metric on a compact complex
manifold. Then the functional $\alpha \mapsto \int_M \alpha \wedge \omega^{n-1}$ 
vanishes on $dd^c$-exact 2-forms. Therefore, 
the map $\deg:\; H^{1,1}_{BC}(M)\arrow \C$,
taking  $[\alpha]$ to $\int_M \alpha \wedge \omega^{n-1}$,
is well-defined:
\begin{equation}\label{_omega_Gaud_integra_Equation_}
\int_M dd^c f \wedge \omega^{n-1} = -\int_M  f \wedge dd^c\omega^{n-1}=0.
\end{equation}
 It is called the {\bf degree map} 
(see also Subsection \ref{_BC_degree_Subsection_}).

\hfill

 Define now a second order differential operator on functions by:
 $$D(f):= \frac{dd^c f\wedge \omega^{n-1}}{\omega^n}.$$
Since $D$ has the same symbol as the Laplacian, 
 it has zero index and is elliptic.
 
 \hfill
 
 \remark 
 Note that $\ker D$  only contains constants by
the Hopf maximum principle (\ref{hopf_theorem}).
\index[terms]{maximum principle} Therefore, 
$\coker D$ is 1-dimensional. By
\eqref{_omega_Gaud_integra_Equation_},
 $\im D$ is the space of functions $g$ such that
$\int_M g\omega^n=0$.

%

\hfill
 
\lemma \label{_primitive_represe_degree0_Lemma_}
Let $[\alpha]\in H^{1,1}_{BC}(M)$ be a degree 0 
Bott--Chern cohomology class.\index[terms]{cohomology!Bott--Chern} Then $[\alpha]$ can be
represented by a closed primitive (1,1)-form.\index[terms]{form!primitive}

 \hfill
 
\proof Let $\alpha_1$ be a (1,1)-form representing 
$[\alpha]$. Then $\int_M \alpha_1 \wedge \omega^{n-1}=0$,
and hence $\alpha_1 \wedge  \omega^{n-1}= dd^c f \wedge
\omega^{n-1}$ for some $f\in C^\infty(M)$. 
Then $\alpha:= \alpha_1 - dd^c f$ is a primitive representative of $[\al]$.
\endproof

\hfill

The next step is essentially \cite[Proposition 4.2]{_Verbitsky:Sta_Elli_}, and its proof
is entirely similar:

\hfill

\lemma \label{_primitive_transversal_Lemma_}
Let $\eta\in \Lambda^{1,1}(M)$ be a primitive, closed
(1,1)-form on a compact Vaisman\index[terms]{form!primitive} manifold.\index[terms]{manifold!Vaisman} Then $\eta$ is
basic with respect to the canonical foliation.\index[terms]{foliation!canonical}

\hfill

\proof Choose an orthonormal basis $z_i$ in $\Lambda^{1,0}(M)$
in such a way that $\omega = -\1 \sum_{i=0}^{n-1} z_i
\wedge \bar z_i$ and $\omega_0 = -\1 \sum_{i=1}^{n-1} z_i
\wedge \bar z_i$. Clearly $z_0=\theta$.  It suffices to prove \ref{_primitive_transversal_Lemma_}
for real (1,1)-forms; then  we may write
$\1\eta = \sum_i a_i z_i \wedge \bar z_i + 
\sum_{i\neq j} b_{ij} z_i \wedge \bar z_j$,
where $a_i$ are real numbers and $b_{ij}= \bar b_{ji}$. 

Since $\eta$ is exact, one has 
$\int_M \eta\wedge \eta\wedge \omega_0^{n-2}=0$.
However, at each point
\begin{equation}\label{_square_primitive_explicit_Equation_}
\frac{\eta\wedge \eta\wedge \omega_0^{n-2}}{\omega^n}=
a_0\sum_{i=1}^{n-1} a_i - \sum b_{0i} b_{i0} = 
a_0\sum_{i=1}^{n-1} a_i - \sum |b_{0i}|^2.
\end{equation}
Since $\eta$ is primitive,\index[terms]{form!primitive} $\sum_{i=0}^{n-1} a_i =0$,
hence $a_0\sum_{i=1}^{n-1} a_i=-a_0^2$. Then 
\eqref{_square_primitive_explicit_Equation_}
becomes 
\[ 
 \frac{\eta\wedge \eta\wedge \omega_0^{n-2}}{\omega^n}
 =-a_0^2- \sum |b_{0i}|^2 \leq 0,
\]
with equality reached only when $a_0=0$ and $b_{0i}=0$ for all $i$.
However, the form $\eta\wedge \eta\wedge \omega_0^{n-2}$
is exact, and thus  the integral $\eta\wedge \eta\wedge
\omega_0^{n-2}$ cannot be positive, giving
$\eta\wedge \eta\wedge \omega_0^{n-2}=0$ and
$a_0=b_{0i}=0$. The latter equality means that $\eta$ is basic.
\endproof

\hfill

\corollary \label{_BC_primitive_Corollary_}
Let $[\alpha]\in H^{1,1}_{BC}(M)$ be a Bott--Chern class \index[terms]{class!Bott--Chern}
 on a compact Vaisman manifold.\index[terms]{manifold!Vaisman} Then $[\alpha]$ can be represented by
a transversal form.

\hfill

\proof See also \cite[Theorem 2.1]{tsu}. The form $\omega_0$ is positive, and hence  it 
has positive degree. Then $[\alpha_1]:=[\alpha]-c[\omega_0]$ has
degree 0, for appropriate $c$. \ref{_primitive_represe_degree0_Lemma_}
and \ref{_primitive_transversal_Lemma_}  
imply that any degree 0 Bott--Chern (1,1)-class $[\alpha_1]$ can
 be
represented by a primitive,\index[terms]{form!primitive} basic form $\alpha_1$.\index[terms]{form!basic}
As $\omega_0$ is closed and $i_{\theta^\sharp}\omega_0=i_{I(\theta^\sharp)}\omega_0=0$,  the form $\omega_0$ is also basic. 
Then $\alpha_1 +c \omega_0$ is a basic form 
representing $[\alpha]$. \endproof





\hfill

We can finally give the

\hfill

\noindent {\bf Proof of  \ref{_BC_kernel_Proposition_}:} Clearly,
$[\omega_0]_{BC}\in\ker\psi$: $\omega_0$ is
exact, but has positive degree; thus $[\omega]$ is  non-zero in 
$H^{1,1}_{BC}(M)$. To prove
\ref{_BC_kernel_Proposition_},
it is enough  to show that the kernel of the map 
$\psi:\; H^{1,1}_{BC}(M)\arrow H^2(M)$
is generated by $\omega_0$.


Let $\eta\in H^{1,1}_{BC}(M)$ be a Bott--Chern cohomology
class vanishing in de Rham cohomology.\index[terms]{cohomology!de Rham}\index[terms]{cohomology!Bott--Chern}  
By \ref{_BC_primitive_Corollary_} we can represent $\eta$ by a 
basic harmonic form. By \ref{_Vaisman_harmonic_forms_Theorem_}, the kernel  of the tautologic
map\index[terms]{map!tautologic} $H^*_\Sigma(M)\arrow H^2(M)$
from the basic cohomology\index[terms]{cohomology!basic} to the
de Rham cohomology is multiplicatively generated by
$\omega_0$. This means that  
$\ker \psi$ is  generated by $\omega_0$. 
 \endproof

  \section[Generic vanishing of Bott--Chern co\-ho\-mo\-lo\-gy]{Generic vanishing of Bott--Chern\\ co\-ho\-mo\-lo\-gy}\index[terms]{cohomology!Bott--Chern}
  
  The main result in this chapter is:
  
  \hfill
  
  \theorem {(\cite{ovv2})} \label{_BC_vanishing_Corollary_}
Let $M$ be an LCK manifold with proper potential,\index[terms]{manifold!LCK!with potential}
 $\alpha\in \C$ and $L_\al$ the flat line bundle
corresponding to $\al\cdot\theta$. Then
$H^{p,q}_{BC}(M,L_\al)=0$ for all
$\al\in \C$ but a  countable subset.

\hfill

 \remark $H^{p,q}_{BC}(M,L_\alpha)=0$
 implies the $d_{\al\theta}d_{\al\theta}^c$-lemma at level $(p,q)$, thus 
 our result 
says that, generically,  a compact LCK manifold with proper 
potential satisfies the $d_{\al\theta}d_{\al\theta}^c$-lemma for all $(p,q)$. On the other hand, \index[persons]{Goto, R.} Goto proved, \cite{goto}, that the $d_{\al\theta}d_{\al\theta}^c$-lemma does not hold in general on compact LCK manifolds.

\hfill

For the proof of \ref{_BC_vanishing_Corollary_}, we need
the following vanishing theorem for the  Dolbeault
cohomology. Its proof takes some background knowledge
in algebraic geometry\index[terms]{geometry!algebraic}, which makes it too technical
for this book. The reader is directed to \cite{ovv2} and to Chapter
\ref{_cohomo_on_Hopf_Theorem_}, where many similar results are proven.

\hfill

\theorem {(\cite{ovv2})} \label{mainresult} For any $q\in \Z^{>0}$,
$H^q(M, \Omega^p_M\otimes L_\alpha)=0$  for all $\alpha\in \C$ but a countable subset.

\hfill

\remark
For  some Hopf manifolds, Ise \cite{_Ise_} and \index[persons]{Mall, D.} Mall \cite{_Mall_}  obtained stron\-ger  vanishing results.
In these cases, the set of exceptions is made explicit. 

\hfill

We give a proof of \ref{_BC_vanishing_Corollary_} using
\ref{mainresult}.

\hfill

\claim Let be an LCK manifold with a proper potential.\index[terms]{manifold!LCK!with potential} Then
the following sequence is exact  for all
$\al\in \C$ but a countable subset:
\begin{multline}\label{exgen}
 H^{q-1}_{\bar \6}(\Omega^p_M\otimes L_\al)\oplus
 \overline{H^{p-1}_{\bar \6}(\Omega^q_M\otimes
   L_\al)}\xlongrightarrow{\6_\theta+\bar\6_\theta}
\\ \xlongrightarrow{\6_\theta+\bar\6_\theta}
H^{p,q}_{BC}(M,L_\al)\stackrel{\nu}{\longrightarrow}H^{p+q}(M,L_\al(\C))
\end{multline}
where $\nu$ is the tautological map\index[terms]{map!tautologic}, $\6_\theta=\6-\theta^{1,0}$ and $\bar\6_\theta=\bar\6-\theta^{0,1}$.

\hfill

\proof  We prove that
$\im(\6_\theta+\bar\6_\theta)=\ker \nu$. Let $\eta$ be a
$(p,q)$-form with values in $L_\al$ with vanishing cohomology class
in the cohomology of the local system $L_\al(\C)$.\index[terms]{local system} Then
$\eta=d_\theta\be$.

\hfill

{\bf Step 1.} Suppose that $\beta$ has only two Hodge components,
$\beta=\be^{p,q-1}+\be^{p-1,q}$.
Then $\eta$ decomposes as
$\eta=\bar\6_\theta \be^{p,q-1}+\6_\theta\be^{p-1,q}$.
On the other hand, as $\eta$ is of bidegree $(p,q)$, we have
$\6_\theta \be^{p,q-1}=0$ and $\bar\6_\theta
\be^{p-1,q}=0$, and hence $\be^{p,q-1}$ and $\be^{p-1,q}$
produce the cohomology classes in
$$[\be^{p-1,q}]\in H^{q-1}_{\bar \6}(\Omega^p_M\otimes L_\al), \ \text{ and} \ \ 
[\be^{p,q-1}]\in \overline{H^{p-1}_{\bar    \6}(\Omega^q_M\otimes L_\al)}.$$
Then
$[\eta]_{BC}=\6_\theta
[\be^{p-1,q}]+\bar\6_\theta[\be^{p,q-1}]$.\\

{\bf Step 2. Reduction to Step 1.}  Recall that $H^{p,q}(L_\alpha)=0$ for all $p
,q$ (\ref{mainresult}).
We use induction on the number of Hodge components.
Take the outermost Hodge component of
$\beta$, say, $\beta^{p-d-1, q+d}$, for $d> 0$. 
Then  $\bar\6_\theta(\beta^{p-d-1, q+d})=0$, and then, by vanishing
of the Dolbeault cohomology\index[terms]{cohomology!Dolbeault} group $H^{p-d-1,q+d}(L_\alpha)$,
we have  $\beta^{p-d-1, q+d}=\bar\6_\theta(\gamma)$, where
$\gamma\in \Lambda^{p-d-1, q-1+d}(M, L_\alpha)$ is an
$L_\alpha$-valued $(p-d-1, q-1+d)$-form.
Now if we replace
$\beta$ by $\beta-d_\theta\gamma$, we obtain another form $\beta'$ such that
$\eta=d_\theta\beta'$, and $\beta'$ has a smaller number
of Hodge components. 
\endproof

\hfill

\noindent By \ref{def_lckpot2Vai}, compact LCK manifolds with
potential are diffeomorphic to Vaisman
manifolds. Then, by \ref{mn_van}, the cohomology of the
local system $L_\al(\C)$ vanishes identically. 
\ref{_BC_vanishing_Corollary_} follows now from 
\ref{mainresult}. \endproof

\chapter{Hopf surfaces in LCK manifolds with potential}\index[terms]{surface!Hopf}\index[terms]{manifold!LCK!with potential}

{\setlength\epigraphwidth{0.7\linewidth}
\epigraph{
\it All pictures that's painted with sense or with thought\\
Are painted by madmen, as sure as a groat;\\
For the greater the fool, in the Art the more blest,\\
And when they are drunk they always paint best.\\
They never can Raphael it, Fuseli it, nor Blake it:\\
If they can't see an outline, pray how can they make it?\\
All men have drawn outlines whenever they saw them;\\
Madmen see outlines, and therefore, they draw them.}
{\sc\scriptsize ``Life of William Blake, Volume 2'' by William Blake  (1880)}
}



\section{Diagonal and non-diagonal Hopf surfaces}


\subsection{Complex curves in non-diagonal Hopf surfaces}

Existence of subvarieties in Hopf manifolds and LCK manifolds with potential\index[terms]{manifold!LCK!with potential}
was a subject of several papers, starting from the work of Ma. \index[persons]{Kato, Ma.} Kato in 1970-ies
(\cite{kato2,_Kato:submanifolds_,ov_surf_in_lck_pot}). In 
\cite{_ov_closed_reeb_} we count the elliptic curves on Vaisman manifolds.\index[terms]{manifold!Vaisman}
Here we prove uniqueness of an elliptic curve on non-Vaisman
Hopf surface; this result is due to I.\index[terms]{surface!Hopf} \index[persons]{Nakamura, I.} Nakamura
(\cite[Theorem 5.2]{_Nakamura:curves_}).

\hfill

\theorem\label{_inva_divisor_unique_Theorem_}
Let $A\in \GL(V)$ be a non-diagonalizable linear 
contraction map, where $V=\C^2$; we choose the basis
such that $A = \begin{pmatrix} \alpha & 1\\ 0  &\alpha
\end{pmatrix}$. Let $\gamma$ denote the holomorphic diffeomorphism of $V$
induced by $A$. Then:
\begin{description}
\item[(i)]
All eigenspaces of the $\gamma^*$-action on  $H^0(V, \calo_V)$
are 1-dimensional.
\item[(ii)] We use the same letter for the action of
  $\gamma$ on the dual space $V^*$. Let
$z\in V^*$ be an eigenvector of $\gamma$.\footnote{This
  eigenvector is unique up to a constant multiplier,
  because $\gamma$ is linear and has a unique Jordan
  cell.} We consider $z$ as a linear function on $V$.
Let $U_\lambda$ be the eigenspace of the
  $\gamma^*$-action on  $H^0(V, \calo_V)$
associated with an eigenvalue $\lambda$. Then
$U_\lambda=0$ if $\lambda$ is not a power of 
$\alpha$, and $U_\lambda= \langle z^k\rangle$ when
$\lambda = \alpha^k$, $k=0,1, 2, ...$.

\item[(iii)] Let $Z$ be an irreducible 
$\gamma$-invariant divisor in $V$. Then
$Z$ is the zero set of $z$ considered to be  a function on $V$.
\end{description}
\pstep
Since (i) is clearly implied by (ii), we start the proof
from (ii). Consider the action of $A=D\gamma$ on
$\Sym^k V^*$. We are going to show that the only eigenspace of $A$
on $\Sym^k V^*$  is 1-dimensional, that
is, the Jordan decomposition of the action of $A$
on $\Sym^k V^*$ has a unique Jordan cell.

We consider the action of 
$\log A - \log \alpha={\scriptsize \begin{pmatrix} 0 & 1\\ 0  & 0\end{pmatrix}}$ on $V^*$ as an element
$\goth f={\scriptsize \begin{pmatrix} 0 & 1\\ 0  & 0\end{pmatrix}}$ of the Lie algebra 
$\goth{sl}(2, \C)$ acting on $V^*$ as
on its fundamental representation.
All eigenvalues of $A$ on $\Sym^k V^*$ 
are equal to $\alpha^k$.
Therefore, the eigenspace of $A$ is identified with the
kernel of the action of $\log A- \log (\alpha^k)$. 
Extending this Lie algebra representation to $\Sym^k V^*$, we obtain that
$\log A- \log (\alpha^k)={\goth f}$ on $\Sym^k V^*$.
By the standard arguments based on the Clebsch--Gordan\index[terms]{theorem!Clebsch--Gordan}
theorem (\cite{_Serre:Lie_}), $\Sym^k V^*$ is a weight $k+1$ irreducible
representation of $\goth{sl}(2, \C)$. Since $\ker {\goth  f}$ 
is the space of highest weight vectors and $\Sym^k V^*$ is
irreducible, this kernel is 1-dimensional.

Let $f\in H^0(V, \calo_V)$ be an eigenvector of 
$\gamma^*$, and $f= \sum_{i=0}^\infty P_i$ the Taylor series for $f$,
where $P_i \in \Sym^k V^*$ are homogeneous polynomials of degree $i$.
Then all $P_i$
are eigenvectors of $A$ acting on $\Sym^k V^*$, and hence  the homogeneous
Taylor components $P_i$ of $f$ are proportional to $z^k$. 
This implies that any eigenvector of the $\gamma^*$-action on 
$H^0(V, \calo_V)$ is equal to $z^k$, and its eigenvalue is $\alpha^k$.

\hfill 

{\bf Step 2:} It remains to prove 
\ref{_inva_divisor_unique_Theorem_} (iii).
Let $I_Z\subset H^0(V, \calo_V)$ be the
ideal of all functions on $V$ vanishing in $Z$.
Since the action of $\gamma^*$ on $H^0(V, \calo_V)$
is compact, the same is true for $I_Z$ (\ref{_Stein_by_contract_to_linear_Hopf_Theorem_}).
 Let $f\in I_Z$ be an eigenvector of the
 $\gamma^*$-action. Then $f= z^k$, and therefore
the zero set of $f$ is equal to the
line $V_z:=\{ v\in V \ \ | z(v)=0\}$.
However, $Z$ is irreducible, and thus  
it is equal to any divisor contained in $Z$.
This implies $V_z=Z$.
\endproof

\subsection{Gauduchon metrics on LCK manifolds with potential}\index[terms]{manifold!LCK!with potential}\index[terms]{metric!Gauduchon}

Let $M$ be an LCK manifold, and  $\psi$ a potential on its cover $\tilde M$.
In general, the norm  of the Lee form\index[terms]{form!Lee} $d\psi$ with respect to the LCK metric is
not constant. The constancy of the norm of the Lee form 
is equivalent to the LCK metric being Gauduchon (see \index[terms]{metric!Gauduchon}
\ref{_Gauduchon_Lee_length_Proposition_} below) and
Vaisman, as shown in Section \ref{pluri_cond}. In fact we have:
 
 \hfill
 
 \proposition\label{_Gauduchon_Lee_length_Proposition_}
 Let $(M, \omega, \theta)$ be a compact LCK manifold with   
 potential and preferred gauge, 
$\omega=d^c\theta + \theta\wedge \theta^c$.  Then the LCK form $\omega$ 
 is Gauduchon if and only if $|\theta|=\const$.
 
 \hfill
 
 {\bf Proof:} The Hermitian form\index[terms]{form!Hermitian} $\omega$ is Gauduchon\index[terms]{form!Gauduchon} if and only if 
 $dd^c\omega^{n-1}=0$.
 
We derive
 \[
 dd^c\omega^{n-1}=(n-1)^2\omega^{n-1}\wedge\theta\wedge\theta^c+(n-1)\omega^{n-1}\wedge
 d\theta^c.
 \]
 On the other hand,
 $$\omega^{n-1}\wedge\theta\wedge\theta^c=\frac 1n |\theta|^2\omega^n,$$
 and $\omega=d^c\theta + \theta\wedge \theta^c$ gives $d\theta^c = -\omega + \theta\wedge \theta^c$, and hence 
\[d\theta^c\wedge\omega^{n-1}=-\omega\wedge\omega^{n-1}+
 \theta\wedge\theta^c\wedge\omega^{n-1}=\left(\frac 1n |\theta|^2-1\right)\omega^n.
\]
 All in all we get:
 \[ dd^c\omega^{n-1}=\frac{(n-1)^2}{n}|\theta|^2\omega^n+(n-1)\left(\frac
 1n|\theta|^2-1\right)\omega^n=(n-1)\big(|\theta|^2-1\big)\omega^n.
 \]
 Then $dd^c\omega^{n-1}=0$ if and only if $|\theta|=1$. This finishes
 the proof.   \endproof
 
 \hfill
 
 \remark \label{_sign_1-theta^2_Remark_}
The eigenvalues of $\omega_0:=d^c\theta$
 are $1$ (with multiplicity $n-1$) and $1-|\theta|^2$. As
 $\omega_0=-I d\theta^c$ is exact, by Stokes' theorem its top power
 cannot be sign-definite whenever $M$ is compact. Two
 possibilities occur:
 \begin{enumerate}
 	\item $|\theta|$ is non-constant, and then  $1-|\theta|^2$ has to change sign on $M$;
 	\item $|\theta|=\const$, and then $|\theta|=1$ and $\omega_0$ is semi-positive.
 \end{enumerate}
 Hence we have:
 
 \hfill
 
\corollary\label{_preferred_gauge_constant_theta_Corollary_} 
Let $(M,\omega,\theta)$ be a compact LCK manifold with
potential with preferred gauge, $\omega=d^c\theta + \theta\wedge \theta^c$. 
Then the following are equivalent:
\begin{description}
\item[(i)] 
the LCK metric is Gauduchon \index[terms]{metric!Gauduchon}
\item[(ii)]
the form $\omega_0=d^c\theta$ is semi-positive
\item[(iii)] $|\theta|=\const$.
\item[(iv)] $|\theta|=1$.
\item[(v)] $(M,\omega,\theta)$ is Vaisman.
\end{description}
\proof
The equivalence of (ii), (iii) and (iv) follows from 
\ref{_sign_1-theta^2_Remark_}, and the equivalence of 
(i) and (iii) is 
\ref{_Gauduchon_Lee_length_Proposition_}.
Finally, the equivalence of (v) and (ii) follows
from \ref{main_0}  below. Note that \ref{main_0} is non-trivial
and uses many results from the rest of this chapter in its proof;
however, we do not use \ref{_preferred_gauge_constant_theta_Corollary_} (v)
before we prove \ref{main_0}.
\endproof 

 \subsection{Complex surfaces of K\"ahler rank 1}
\label{_Kahler_rank_one_Subsection_}\index[terms]{rank!K\"ahler}
 
 \definition {(\cite{hl,ct,_Brunella:Inoue_})}\\ A compact complex surface 
is {\bf of K\"ahler rank 1}\index[terms]{rank!K\"ahler} if  it is not K\"ahler but 
 admits a non-trivial closed semipositive (1,1)-form.
 
 \hfill
 
\lemma A  compact  LCK surface $(M, \theta, I)$ with potential\index[terms]{surface!LCK!with potential} and
a semi-positive form  has K\"ahler rank 1\index[terms]{rank!K\"ahler}. 
  
\hfill

 {\bf Proof:} We need only to show that any LCK surface
$M$ with potential cannot be K\"ahler.  In dimension $>2$,
we would have used  \ref{def_lckpot2Vai} to show that $M$
is Vaisman, and this would imply that
$b_1(M)$ is odd (\ref{_odd_first_betti_}).
For dimension 2 this argument does not work, because
our proof of \ref{def_lckpot2Vai} in dimension 2 depends on the validity of the\index[terms]{conjecture!GSS}
GSS conjecture (\ref{_sphe_implies_dim2_Theorem_}).
However, Vaisman's theorem (\ref{vailcknotk})
is still valid, and it implies that $\theta$ is exact:
$\theta= df$. Taking the preferred gauge, 
$\omega=d^c\theta + \theta\wedge \theta^c$, we can assume
that $\omega= d^c df + df\wedge d^c f$.
Let $x\in M$ be the point where $f$ has maximum.
Then $dd^c f\restrict {T_xM}$ is a negative (1,1) form,
hence $\omega=d^c df + df\wedge d^c f$ has at most one
positive eigenvalue at $x$.
Therefore, $\omega$ cannot be Hermitian.
 \endproof
 
 \hfill
 
 Recall that a {\bf Hopf surface} \index[terms]{surface!Hopf!diagonal}
 $H$ is a quotient of $\C^2 \backslash 0$
 by a holomorphic contraction. 
 A Hopf surface is {\bf diagonal} if this 
 contraction is expressed by a 
 diagonalizable linear automorphism.
 
 \hfill
 
 Compact surfaces of K\"ahler rank 1\index[terms]{rank!K\"ahler} have been 
 classified in \cite{ct} and \cite{_Brunella:Inoue_}.\footnote{%
For more details on the classification of non-K\"ahler complex surfaces,
see Chapter \ref{comp_surf}.}
 They can be:
 \begin{enumerate}
 	\item Non-K\"ahler elliptic fibrations; 
 	\item Diagonal Hopf surfaces and their blow-ups; 
 	\item Inoue surfaces and their\index[terms]{surface!Hopf!diagonal} blow-ups.\index[terms]{blow-up}\index[terms]{surface!Inoue}
 \end{enumerate}

We want to know which surfaces in this list can admit a Vaisman metric.

 The LCK Inoue surfaces cannot have  LCK metrics with potential, \index[terms]{manifold!LCK!with potential}
as shown in \ref{_no_exact_on_inoue_} (see also  \ref{nopot}). 

 A cover of a blow-up of any complex manifold cannot admit
 plurisubharmonic functions. Indeed, by the homotopy lifting
 lemma, the projective spaces  contained in the
 blow-up lift to the cover. Thus, blow-ups cannot have
 global potential.

 We are left with non-K\"ahler elliptic fibrations and
 diagonal Hopf surfaces, known to admit Vaisman
 metrics, as shown in Chapter \index[terms]{surface!Hopf!diagonal}
\ref{comp_surf}; see also \cite{bel}. We then have the following:
 
 \hfill
 
\proposition\label{surfvai} 
(\cite{ov_surf_in_lck_pot}) \\
Let $M$ be a compact LCK surface with potential \index[terms]{surface!LCK!with potential}
that has K\"ahler rank 1\index[terms]{rank!K\"ahler}. Then $M$ is a diagonal Hopf surface
or an elliptic surface.\index[terms]{surface!elliptic} In particular, $M$ is Vaisman.
 \endproof

 \hfill
 
 For further use, it is convenient to list 
 all criteria used to distinguish the Vaisman
 Hopf surfaces from the non-Vaisman ones.
 
 \hfill
 
 \theorem\label{_Hopf_surface_Vaisman_Theorem_} 
(\cite{ov_surf_in_lck_pot})
 Let $M$ be a linear Hopf surface. Then the following
 are equivalent.\index[terms]{surface!Hopf!linear}
 \begin{description}
 	\item[(i)] $M$ is Vaisman.
 	\item[(ii)] $M$ is diagonal.
 	\item[(iii)] $M$ has K\"ahler rank 1\index[terms]{rank!K\"ahler}.
 	\item[(iv)] $M$ contains at least two distinct elliptic curves.
 \end{description}
 
 {\bf Proof:} The conditions (i) and (iii) are equivalent by 
\ref{surfvai}. A diagonal Hopf surface is Vaisman by \ref{semihopf}.
A non-diagonal Hopf surface cannot have K\"ahler 
rank 1  by \ref{surfvai}. This gives the equivalence
(ii) $\Leftrightarrow$ (i). \index[terms]{surface!Hopf!non-diagonal}

The equivalence of (iv) and (ii) is shown in
(\cite[Theorem 5.2]{_Nakamura:curves_}). It follows
immediately from \ref{_inva_divisor_unique_Theorem_}.
\endproof

 
 \section{Surfaces in compact LCK manifolds with potential}\index[terms]{manifold!LCK!with potential}
 

 \subsection{Algebraic groups}

 We now need to recollect several facts about algebraic groups (see also Subsection \ref{alg_gr_jc}). 
 Fix an $n$-dimensional complex vector space $V$.
 
 \hfill

 \lemma 
 Let $A\in \GL(V)$ be a linear operator,  and
 $\langle A \rangle$ the group generated by $A$. Denote by
 $G$ the Zariski closure of $\langle A \rangle$ in $\GL(V)$.
 Then, for any $v\in V$, the Zariski closure $Z_v$
 of the orbit $\langle A\rangle\cdot v$ is equal to the usual
 closure of $G\cdot v$.\index[terms]{Zariski closure}
 
 \hfill
 
 {\bf Proof:} Clearly, $Z_v$ is $G$-invariant. Indeed,
 its normalizer $N(Z_v)$ in $\GL(V)$ is an algebraic group
 containing $\langle A\rangle$,\index[terms]{normalizer}\index[terms]{group!algebraic}
  hence $N(Z_v)$ contains $G$. The converse is also
 true: since  $\langle A\rangle$ normalizes  $\langle A\rangle\cdot v$,
 its Zariski closure $G$ normalizes the Zariski closure
 $Z_v$ of the orbit. Therefore, the orbit  $G\cdot v$
 is contained in $Z_v$. Since $G\cdot v$ is a constructible
 set\index[terms]{constructible set}, its Zariski closure
 coincides with its usual closure, 
\cite[Expos\'e XII, Proposition 2.2]{_SGA1:new_}.
 This gives $\overline{G\cdot v}\subset Z_v$. As $\overline{G\cdot v}$ is
 algebraic and contains $\langle A\rangle\cdot v$, the inclusion $Z_v\subset \overline{G\cdot v}$ is also true.
 \endproof
 
 \hfill
 
 The reason we take the Zariski closure is 
explained in the following (see also \cite[Theorem 2.1]{ov_imrn_10}):
 
 \hfill
 
 \claim Let $I\subset \C[z_1, ..., z_n]$
 be an ideal that is  invariant with respect
 to  an automorphism $A$ of the space
 $\langle z_1, \ldots, z_n\rangle$ acting
 on the polynomial ring. Then $I$
 is invariant with respect to the
 Zariski closure $G$ of $\langle A\rangle$.
 
 \hfill
 
 {\bf Proof:} First, we show that the 0-adic completion $\hat I$ of
 $I$ is $G$-invariant in the 0-adic completion
 of the polynomial ring $\C[z_1, ..., z_n]$, that is  the
 ring of formal power series $\C[[z_1, ..., z_n]]$. However, any
 $A$-invariant subspace in a finite-dimensional
 space is $G$-invariant by definition of $G$, and
 the ideal $\hat I$ is obtained as an inverse
 limit of finite-dimensional subspaces of
 finite quotients of the polynomial ring.
 Therefore, $\hat I$ is $G$-invariant.
 The ideal $I$ is $G_A$-invariant, 
 because
 $I=\hat I \cap \C[z_1, ..., z_n]$. \endproof
 
 \hfill

 Using the Jordan--Chevalley decomposition (see \ref{jcdec}) we have:\index[terms]{Jordan--Chevalley decomposition}
 
 \hfill 
 
 \corollary\label{_algebraic_action_on_LCK_with_pot_Corollary_}
 Let $M$ be a subvariety of a linear Hopf manifold
 $H=(V\backslash 0)/A$, $\tilde M \subset V\backslash 0$
 its $\Z$-covering\index[terms]{cover!K\"ahler $\Z$-}, and $G$ the Zariski closure of \index[terms]{Zariski closure}
 $\langle A \rangle$ in $\GL(V)$. Then $\tilde M$ is $G$-invariant. 
Moreover, the connected component of identity
of $G$ is a product of $G_s:=(\C^*)^k$ and a commutative unipotent group
 $G_u$ commuting with $G_s$, and both of these groups
 preserve $\tilde M$.
 
 \hfill
 
 {\bf Proof:}
 Let $X$ be the closure of $\tilde M$ in $\C^N$.
 The ideal $I_X$ of $X$ is generated 
 by polynomials, as shown in \ref{_cone_cover_for_LCK_pot_Theorem_}, Step 1. 
In each degree, the space $I_d$ of graded polynomials of degree $d$ that belong
to $I$ is finite-dimensional, and $I_X$ is generated by $\bigoplus_d I_d$. 
By definition of the Zariski closure, each space $I_d$ is $G$-invariant; hence 
$I_X$ is also $G$-invariant. This implies that the set
 $\tilde M\subset \C^N$ is $G$-invariant (see also \ref{def_lckpot2Vai}, Step 1).
 
 The last assertion of 
 \ref{_algebraic_action_on_LCK_with_pot_Corollary_}
 is implied by 
\cite[Proposition 1.5]{_Borel_Tits:Groupes_Reductifs_},
which claims that any  commutative algebraic
group consisting of semisimple elements is isomorphic to $(\C^*)^k$.
 \endproof

\hfill

The following example is sometimes very useful; we give it
to illustrate
\ref{_algebraic_action_on_LCK_with_pot_Corollary_}.

\hfill

\claim\label{_Zariski_closure_of_Jordan_cell_Claim_}
Let $A\in \GL(n,\C)$ be a contraction that has only one Jordan cell
of maximal rank, 
\[A= \small \begin{pmatrix} \alpha & 1 & 0 & \ldots & 0& 0\\
0 & \alpha & 1 & \ldots & 0& 0\\
\vdots &\vdots &\vdots & \cdots & \vdots \\
0&0&0 & \ldots & \alpha &1\\
0&0&0 & \ldots & 0  &\alpha
\end{pmatrix}.
\]
Then the Zariski closure of $\langle A\rangle$ is
a 2-dimensional Lie group generated by the homotheties
$\C^* \Id$ and the unipotent group $e^{\C A_0}$,
where
\[A_0= \small\begin{pmatrix} 0 & 1 & 0 & \ldots & 0& 0\\
0 & 0 & 1 & \ldots & 0& 0\\
\vdots &\vdots &\vdots & \cdots & \vdots \\
0&0&0 & \ldots & 0 &1\\
0&0&0 & \ldots & 0  &0
\end{pmatrix}.
\]
is the nilpotent part of $\log A$.

\hfill

\proof
The decomposition $A= (\alpha \Id) e^{A_0}$
coincides with the Jordan-Che\-val\-ley decomposition of
$A$; indeed, the operator $\alpha \Id$ is semisimple,
$e^{A_0}$ is unipotent, and they commute.
Therefore, the Zariski closure of $\langle A \rangle$
is the product of the Zariski closure of $\langle\alpha \Id\rangle$
and the Zariski closure of $\langle e^{A_0}\rangle$.
On the other hand, the group $\C^* \Id$ is 1-dimensional
and contains $\langle\alpha \Id\rangle$, the group
$e^{\C A_0}$ is 1-dimensional and contains 
$\langle e^{A_0}\rangle$; these groups are 
Zariski closures of $\langle\alpha \Id\rangle$
and $\langle e^{A_0}\rangle$, because they are
connected and have the minimal possible dimension.
\endproof

\subsection{Orbits of algebraic groups in Hopf manifolds}

 \lemma\label{_surface_exists_Lemma_}
 Let $M$ be a submanifold of a linear Hopf manifold
 $H=(V\backslash 0)/A$, $\dim_\C M \geq 3$, and $G=G_sG_u$ the Zariski
 closure of $\langle A \rangle$ with its Jordan--Chevalley 
 decomposition. Then $M$ contains a surface $M_0$, possibly singular, 
 with $G_u$ acting non-trivially on its $\Z$-covering\index[terms]{cover!K\"ahler $\Z$-} 
 $\tilde M_0\subset V$.
 
 \hfill
 
 {\bf Proof:} Another form of this statement 
 is proven by Masahide \index[persons]{Kato, Ma.} Kato (\cite{kato2}).
 
 We shall use induction on dimension of $M$.
 It would suffice
 to find a subvariety $M_1 \subset M$ of codimension 1
 such that $G_u$ acts non-trivially on its $\Z$-covering \index[terms]{cover!K\"ahler $\Z$-}
 $\tilde M_1\subset \C^n \backslash 0$. Replacing 
 $V$ by the smallest $A$-invariant subspace containing
 $\tilde M$, we may assume that the intersection
 $\tilde M \cap V_1$ is not equal to $V_1$ for each proper subspace
 $V_1\subsetneq V$. Now take a codimension 1 
 subspace $V_1\subsetneq V$ that is  $A$-invariant and 
 such that $G_u$ acts on $V_1$ non-trivially
 (equivalently, such that $A$ acts on $V_1$
 non-diagonally). Using the Jordan decomposition\index[terms]{Jordan decomposition}
 of $A$, such $V_1$ is easy to construct.
 Then $\tilde M_1':= V_1 \cap \tilde M$ 
 gives a subvariety of $M$ of codimension 1
 and with non-trivial action of $G_u$.
 \endproof
 
 \hfill
 
 The same argument gives the following corollary,
 also parallel to a theorem by Ma. \index[persons]{Kato, Ma.} Kato.
 
 \hfill
 
 \corollary\label{flag} 
(\cite{ov_surf_in_lck_pot}) 
 Let $M$ be a compact LCK subvariety of a Hopf manifold,
in particular, an LCK manifold with potential\index[terms]{manifold!LCK!with potential}.
 Then $M$ has a flag of embedded subvarieties
 $M\supset M_1\supset M_2\supset \cdots\supset M_{\dim M-1}$
 with $\codim M_i =i$.
 \endproof
 
 \hfill
 
 
 Recall that \index[terms]{manifold!Oeljeklaus--Toma (OT)} Oeljeklaus--Toma manifolds  do not admit complex curves
 (see \ref{_nocurves_in_OT_Theorem_}). Then \ref{flag} gives a new proof for:
 
 \hfill
 
 \corollary\label{nopot} 
(\cite{oti2}) 
The Oeljeklaus--Toma
 manifolds, in particular, \index[terms]{surface!Inoue} Inoue surfaces,
 cannot admit LCK structures with potential\index[terms]{structure!LCK!with potential}.
 \endproof
 
 \hfill

 \proposition\label{_surface_non-diagonal_Proposition_}
 Let $M_1\subset H=(V\backslash 0)/\langle A\rangle$
 be a surface in a linear Hopf manifold, possibly singular,
 and $G=G_sG_u$ the Zariski closure of $\langle A \rangle$ 
 with its Jordan--Chevalley \index[terms]{Zariski closure}
 decomposition. Assume that $G_u$ acts on
 $M_1$ non-trivially. Then the normalization $M$ of \index[terms]{Jordan--Chevalley decomposition}\index[terms]{normalization}
 $M_1$ is a non-diagonal (in particular, non-Vaisman) Hopf surface.\index[terms]{surface!Hopf!non-diagonal}
 
 \hfill
 
 {\bf Step 1:}  We prove that $M_1$ is smooth
and equipped with an effective action of a 2-dimensional Lie group.
Replacing $G$ by its quotient
 by the subgroup acting trivially on the $\Z$-cover\index[terms]{cover!K\"ahler $\Z$-}
 $\tilde M_1$ if necessary, we may assume that
 $G$ acts faithfully on a general orbit in $\tilde M_1$. 
 Then $G$ is at most 2-dimensional. However, it cannot be
 1-dimensional because $G_s$ contains contractions
 (hence cannot be 0-dimensional) and $G_u$ acts
 non-trivially. Therefore, the connected component
of identity in $G_s$ is isomorphic to  $\C^*$ and $G_u\simeq\C$
(\ref{_algebraic_action_on_LCK_with_pot_Corollary_}).
 
 Since $G_s$ acts by contractions, the quotient
 $S:= \tilde M_1/G_s$ is a compact curve, equipped with the
 $G_u$-action that has a dense orbit.\index[terms]{orbit!dense} The group
 $G_u$ can act non-trivially only on a genus 0 curve,
 and there is a unique open orbit $O$ of $G_u$, with 
 $S \backslash O$ being one point. 
 
  Since the singular
 set of the normal surface $M$ is $G_s$-invariant, it has dimension at least 1, 
 and since $M$ is normal, it is non-singular in codimension 1,
 thus smooth. 

\hfill

{\bf Step 2:} We prove that $M$ cannot be Vaisman.
In Step 1, we proved $S:= \tilde M_1/G_s$ is $\C P^1$, and the unipotent group $G_u$
acts on $S$ with an open orbit $O=\C$.
This implies, in particular, that $M$ has only one curve:
otherwise it would have been elliptic, that is  impossible,
because all subvarieties of $M$ are $G$-invariant 
(\ref{_algebraic_action_on_LCK_with_pot_Corollary_}), and $G$
acts on $M$ with an open orbit.
However, a Vaisman surface has at least two elliptic curves  by
\ref{flag}. \index[terms]{curve!elliptic}

\hfill

{\bf Step 3:} Now we can prove that $M$ is a Hopf surface.\index[terms]{surface!Hopf}

Since $\tilde M_1$ is $G_s$-invariant, and $G_s$ acts
on this variety (that is  actually an open algebraic cone,\index[terms]{cone!algebraic} 
\ref{_alge_cone_Definition_}) the quotient 
$\tilde M_1/\langle A_s\rangle$ is a subvariety in
a Vaisman manifold.\index[terms]{manifold!Vaisman} Restricting the form $\omega_0$
to $\tilde M_1/\langle A_s\rangle$, we obtain that
this variety has K\"ahler rank 1\index[terms]{rank!K\"ahler}. Therefore,
$M$ is a deformation of a surface $M'$ that has 
K\"ahler rank 1.\footnote{We cannot claim that $M'$ 
is Vaisman immediately, because it is not a subvariety
of a diagonal Hopf manifold; by construction, $M'$
is a normalization of a subvariety of a diagonal Hopf.}\index[terms]{normalization}
Since $M'$ is the normalization of a subvariety in a Hopf manifold,
it is minimal; indeed, Hopf manifolds cannot contain rational curves.
In Subsection \ref{_Kahler_rank_one_Subsection_} it was shown that a surface of K\"ahler rank 1\index[terms]{rank!K\"ahler}
is Inoue or Vaisman if it is minimal. However, the Inoue surfaces
do not have complex curves, and subvarieties of
a Hopf manifold always have curves by \ref{flag},
hence $M'$ is Vaisman. 

By \ref{_Vaisman_is_Hopf_or_elli_Proposition_}, $M'$ is
an elliptic surface or a Hopf surface (or both).\index[terms]{surface!Hopf}\index[terms]{surface!elliptic}
However, any surface that is  not of class VII is elliptic,
and all its deformations are elliptic, because class VII
is closed under deformations. This is impossible,
because $M$ has only one curve. Therefore,
$M$ is a deformation of a Hopf surface.
Now, the classification of non-K\"ahler
surfaces with $b_2=0$ (Section \ref{_embedding_surfaces_pot_})
implies that a  deformation of a Hopf surface is again Hopf.
\endproof

\hfill

\remark
For most of this book, we considered linear Hopf manifolds,
that is, quotients of $\C^n \backslash 0$ by a cyclic
group generated by a linear contraction. In Chapter 
\ref{_non-linear_Hopf_Chapter_} and Chapter \ref{_Mall_bundles_Chapter_}, 
we consider a more general case, when $A$ is a holomorphic contraction,
not necessarily linear. In \ref{_surface_non-diagonal_Proposition_},
we have shown that the normalization $M$ of a certain surface inside a Hopf manifold\index[terms]{normalization}
is always Hopf; however, we do not claim that $M$ is linear, and in fact
it can be non-linear.

\subsection{Hopf surfaces in LCK manifolds with potential}\index[terms]{manifold!LCK!with potential}

\theorem (\cite{ov_surf_in_lck_pot})\label{main_surf_in_pot} 
A non-Vaisman compact LCK manifold $M$ with potential,  $\dim_\C M\geq 3$,
contains a surface with normalization biholomorphic to a non-diagonal 
Hopf surface.\index[terms]{surface!Hopf!non-diagonal}

\hfill

\proof
Let $M$ be a compact LCK manifold with potential\index[terms]{manifold!LCK!with potential}. Then  $M$ is holomorphically embedded
 into a Hopf manifold $\C^N\backslash 0 /\langle A\rangle$,
 where $A\in \GL(N,\C)$  is a linear operator, see  \ref{embedding}.
 Applying \ref{_surface_exists_Lemma_} 
 and \ref{_surface_non-diagonal_Proposition_},
 we find a surface in $M$ biholomorphic to a  non-diagonal Hopf surface .
 \endproof

 \hfill

 As an application, we now prove the following 
characterization of Vaisman manifolds:\index[terms]{manifold!Vaisman}
 
 \hfill

 \corollary\label{main_0} 
 Let $(M, \omega, I)$ be a compact LCK manifold with
 potential. Assume that the $(1,1)$-form
 $\omega_0:=d\theta^c$ is semipositive.
Then the LCK metric $\omega$ is Vaisman.  
 
 \hfill
 
 {\bf  Proof:} 
If $\dim_\C M=2$, this is just \ref{surfvai}.
 
 If $\dim_\C M\geq 3$, we first prove that $M$ is of
Vaisman type. Otherwise $M$ contains a surface whose 
 normalization is a non-diagonal Hopf surface\index[terms]{normalization}
 $H$ (\ref{main_surf_in_pot}). Then $\omega_0$ restricts to a semi-positive\index[terms]{surface!Hopf!non-diagonal}
 definite (1,1) form on $H$. By \ref{surfvai}, $H$ is
 Vaisman, and thus diagonal, contradiction.

It remains to show that the LCK metric of $(M,I)$ is, indeed, 
Vaisman.  Since $\omega_0$ is
semipositive and closed, it vanishes on the 
canonical foliation\index[terms]{foliation!canonical} $\Sigma$ on $(M,I)$. By 
\ref{_harmo_deco_1-form_Proposition_}, Step 2, 
$\omega_0$  is basic with respect to $\Sigma$. 

Let $G$ be the Lie group generated by\index[terms]{form!basic}
the Lee field\index[terms]{Lee field} $X$ and the anti-Lee field\index[terms]{Lee field!anti-} $X^c$ of a Vaisman
structure on $M$. Since $\omega_0$ is basic, it is $G$-invariant.
Since $d\theta$, $d^c\theta$, $d\theta^c$ and $d^c\theta^c$
 are $G$-invariant,
the action of the Lie algebra $\Lie(G)$ of $G$ 
takes $\theta$ and $\theta^c$
to forms that are  $d$ and $d^c$-closed.
The (1,0)-part of a form that is  $d$ and $d^c$-closed
is clearly holomorphic. 

However, the Cartan formula implies
that the image of the action of $\Lie(G)$ on closed forms is exact.
On a compact manifold, all holomorphic
exact 1-forms vanish, by Hopf maximum principle\index[terms]{maximum principle}
(\ref{_holo_exact_1-form_Lemma_}). Therefore,
the form $\theta$ is $G$-invariant.
The Lee action on the K\"ahler cover $(\tilde M, \tilde \omega)$ of $(M, \omega)$
cannot fix a K\"ahler potential, because it acts on 
$\tilde M$ by contractions. Therefore, the 
lift of the $G$-action acts on $\tilde \omega$ by non-isometric
homotheties. Applying \ref{kami_or}, we obtain that
$(M, \omega, \theta)$ is conformal to a Vaisman 
manifold. However, by \ref{_sign_1-theta^2_Remark_}, 
$\omega$ is Gauduchon,
and in each conformal class of a Vaisman metric,
the Gauduchon metric is Vaisman.
 \endproof


\section{The pluricanonical condition}\label{pluri_cond}


The pluricanonical condition was introduced in \cite{ko}
as a weakening of the Vaisman condition. 

\hfill

\definition 
An LCK manifold $(M,J,g,\theta)$ satisfying $(\nabla
\theta)^{1,1}=0$ is called {\bf pluricanonical}.\index[terms]{manifold!pluricanonical}

\hfill

Using \eqref{dthetac} and the fact that $d$ is the
anti-symmetrization\index[terms]{anti-symmetrization} of
$\nabla$, one proves (see also \cite{ov_imrn_10,_moroianu:pluricanonical_}):

\hfill

\claim\label{_pluricanonical_LCK-pot_Claim_}
Let $(M, I, \omega, \theta)$ be an LCK manifold.
Then the following are equivalent:
\begin{description}
\item[(i)] M is pluricanonical. 

\item[(ii)] The following equation is satisfied:
\begin{equation}\label{thetac_pluri}
d\theta^c=-\theta\wedge\theta^c+|\theta|^2\omega,
\end{equation}
and $|\theta|=\const$.
\item[(iii)] $\Lie_{\theta^\sharp}\omega=0$.
\end{description}
\proof
The 2-form $\nabla\theta$ is symmetric, because
$\theta$ is closed, and the Levi--Civita connection is torsion-free.\index[terms]{connection!torsion-free}
The pluricanonical condition is equivalent to
$$(\nabla_{IX}\theta)(IY)+(\nabla_X\theta)(Y)=0.$$
Changing $Y$ into $IY$ and using 
$\nabla\theta \in \Sym^2T^* M$,
this gives
$$(\nabla_Y\theta)(IX)-(\nabla_X\theta)(IY)=0.$$
The formula \eqref{dthetac} gives
\begin{equation}
d\theta^c(X,Y)=\left(-|\theta|^2\omega+
\theta\wedge\theta^c\right)(X,Y)+\frac 12\left
((\nabla_Y\theta)(IX)-(\nabla_X\theta)(IY)\right).
\end{equation}
Then the pluricanonical condition 
$(\nabla_Y\theta)(IX)-(\nabla_X\theta)(IY)=0$
is equivalent to
\begin{equation}\label{pluri_1}
d\theta^c=-\theta\wedge\theta^c+|\theta|^2\omega.
\end{equation}
We now take the exterior derivative of the above,
 using $d\omega=\theta\wedge\omega$ and 
$d\theta^c=-\theta\wedge\theta^c+|\theta|^2\omega$, obtaining:
$$0=-\theta\wedge d\theta^c-d|\theta|^2\wedge\omega-
|\theta|^2\theta\wedge\omega=-d|\theta|^2\wedge\omega,
$$
which implies $|\theta|=\const$,
because the exterior
multiplication by $\omega$ is injective on 1-forms.

We proved the equivalence of (i) and (ii).
Now, (iii) is equivalent to 
$0=\Lie_{\theta^\sharp}\omega= d\theta^c+ i_{\theta^\sharp}(d\omega)$;
since $i_{\theta^\sharp}(d\omega)=-\theta\wedge\theta^c+|\theta|^2\omega$,
this implies that $\Lie_{\theta^\sharp}\omega=0$ is equivalent to
\eqref{thetac_pluri}.
\endproof

\hfill

In \cite{_moroianu:pluricanonical_}, A. \index[persons]{Moroianu, A.} Moroianu and \index[persons]{Moroianu, S.} S. Moroianu have
proven that any compact pluricanonical manifold is Vaisman.
Here we give another proof of this statement.

\hfill

\theorem {(\cite{_moroianu:pluricanonical_})}\label{_pluri_vais_} 
A compact pluricanonical manifold is Vaisman.
  
\hfill

\proof 
By \ref{_pluricanonical_LCK-pot_Claim_}, 
$d\theta^c=-\theta\wedge\theta^c+|\theta|^2\omega,$
and $|\theta|=\const$. Rescaling $\omega$ by a constant
multiplier, we can assume that $|\theta|^2=1$.
The equation 
$d\theta^c=-\theta\wedge\theta^c+\omega$
is a special form of the ``LCK with potential'' condition;
in \ref{_preferred_gauge_Definition_}, we called it 
``the LCK potential\index[terms]{potential!LCK} with preferred conformal gauge''.\index[terms]{manifold!LCK!with potential and preferred gauge} 
Let $\tilde M\stackrel \pi\arrow M$ be a K\"ahler cover\index[terms]{cover!K\"ahler} of $M$,
such that $\pi^*\theta$ is exact. Then
$\pi^* \theta = d \psi$. Consider the function
$\phi= e ^{-\psi}$. The equation 
$d\theta^c=-\theta\wedge\theta^c+\omega$ easily implies
that $\phi$ is a positive,
automorphic K\"ahler potential on $\tilde M$, as 
we are going to show now.


Clearly, $d^c\phi= d^c(e^{-\psi})=-\phi d^c\psi= - \phi \pi^*\theta^c$, 
and $d\phi=d(e^{-\psi})=-\phi d\psi= - \phi \pi^*\theta$, which 
gives
\begin{equation*}
	\begin{split}
dd^c \phi& = - d (\phi \pi^*\theta^c)= - \phi^2 \pi^*\theta \wedge \pi^*\theta^c
- \phi d(\pi^*\theta^c)\\
&= - \phi^2 \pi^*\theta \wedge \pi^*\theta^c+
\phi^2\pi^*\theta\wedge\pi^*\theta^c + \phi \pi^*\omega= \phi\pi^*\omega
	\end{split}
\end{equation*}
using $d\theta^c=-\theta\wedge\theta^c+\omega$.

Now, $M$ is Vaisman by \ref{_preferred_gauge_constant_theta_Corollary_}.
\endproof  

\hfill

\remark 
For the time being, there are no examples of non-compact
pluricanonical manifolds that are  not Vaisman. It is not
known if there exists a purely local proof of the above
result.


\chapter{Riemannian geometry of LCK manifolds}


\epigraph{\it To those who deny the Word, to those who
deny the Jeweled Veil and the Visage, runs an imprecation
from the Rosa Secreta, I vow a wondrous Hell, for each
person who so denies shall reign over 999 empires of fire,
and in each empire shall be 999 mountains of fire, and
upon each mountain there shall be 999 towers of fire, and
each tower shall have 999 stories of fire, and each story
shall have 999 beds of fire, and in each bed shall that
person be, and 999 kinds of fire, each with its own face
and voice, shall torture that person throughout eternity.}
{\sc\scriptsize ``Hakim, the Masked
Dyer of Merv'', by Jorge Luis Borges}

\section{Existence of parallel vector fields}\index[terms]{vector field!parallel}
The aim of this section is to prove the following:

\hfill

\theorem {(\cite{_moroianu:parallel_},
    \cite{_Madani_Moroianu_Pilca:Einstein_})}\label{parafield} \\
Let $(M,I,g,\theta)$ be a compact, connected LCK manifold of complex dimension $n\geq 2$, of non-K\"ahler type. If $M$  admits a non-trivial parallel vector field $X$, then $g$ is Vaisman and $\theta^\sharp=X$.\index[terms]{vector field!parallel}

\hfill

The rest of this section will be devoted to the proof of
this theorem (that is  rather technical, we shall indicate
the main points of the authors' computations). Without
loss of generality, we may assume $X$ has unit length. 

By a rather involved computation, one first proves that the Lee field\index[terms]{Lee field} belongs to the distribution generated by $X$ and $IX$:
$$\theta^\sharp=aX+b(IX),\quad a,b\in C^\infty M.$$
This is done by first assuming that
$\theta=aX^\flat+b(IX)^\flat+\theta_0$, 
with $\theta_0^\sharp$ orthogonal to 
$\langle X, I(X)\rangle$,
and showing that the codifferential of $\theta_0$ satisfies
$$d^*\theta_0=c (\Lie_X(a)-|\theta_0|^2), \quad c\in\R^{>0}.$$
Since $\int_M d^*\theta_0 \,d\!\Vol=0$, we have $\int_M\Lie_X(a)d\Vol= \int_M|\theta_0|^2d\Vol.$ However,  
\[ \int_M \Lie_X(a)d\!\Vol =-\int_M a \Lie_X(d\!\Vol)=0
\]
(since $X$ is parallel,\index[terms]{vector field!parallel} its divergence $\Lie_X(d\!\Vol)$
vanishes). This proves $\theta_0=0$. The next step is:

\hfill

\claim\label{_ab=0_Claim_} 
$a=\const$, and $ab=0.$

\hfill

\proof
Indeed, since $d\theta=0$ and $X^\flat$ is parallel, the
derivative of 
$\theta=aX^\flat
+b(IX)^\flat$ gives $da\wedge X^\flat+db\wedge (IX)^\flat+b\,d((IX)^\flat)=0$.

On the other hand, \eqref{DJ} together with the parallelism of $X^\flat$ implies  
\begin{equation}\label{_d_I_X_flat_Equation_}
d((IX)^\flat)=a(X^\flat\wedge (IX)^\flat-\omega),
\end{equation}
 and hence
\begin{equation}\label{_da_db_ab_Equation_}
	da\wedge X^\flat+db\wedge (IX)^\flat+ab(X^\flat\wedge (IX)^\flat-\omega)=0.
\end{equation}
Fix a point $x\in M$ and evaluate the above equality at $x$. Let $v\in T_xM$ be orthogonal to $X_x$ and $(IX)_x$ (this is possible because $\dim_\R M\geq 4$). Then the scalar product of the left-hand side of \eqref{_da_db_ab_Equation_} with $v^\flat\wedge(Iv)^\flat$ gives $ab=0$ at $x$. Since $x$ was chosen arbitrarily, $ab=0$ on $M$. Now \eqref{_da_db_ab_Equation_} reduces to:
\begin{equation}\label{_daX+dbIX_}
	da\wedge X^\flat+db\wedge (IX)^\flat=0.
\end{equation}
We differentiate \eqref{_d_I_X_flat_Equation_}, taking into account that $d(X^\flat)=0$ and $ab=0$, and obtain:
\begin{equation}
	\begin{split}
		 0 &=da\wedge (X^\flat\wedge (IX)^\flat-\omega)+a(-X^\flat\wedge d((IX)^\flat)-\theta\wedge\omega)\\
		 &=da\wedge (X^\flat\wedge (IX)^\flat-\omega)\\
		 &+a \big(-X^\flat\wedge a(X^\flat\wedge (IX)^\flat-\omega)-(aX^\flat+b(IX)^\flat)\wedge\omega\big)\\
		 &=da\wedge (X^\flat\wedge (IX)^\flat-\omega).
	\end{split}
\end{equation}
Together with 	\eqref{_daX+dbIX_}, this gives
$-da\wedge\omega=0$, implying $da=0$, because $\omega$ is
non-degenerate.	 Since $M$ is connected, we have
$a=\const$. This finishes the proof of 
\ref{_ab=0_Claim_}.
\endproof

\hfill

Since $a=\const$ and $ab=0$,
only two possibilities arise. 

Either  $a\neq 0$. Then $b=0$  and $\theta^\sharp$ is a
constant multiple of $X$, thus parallel, and hence  $M$ is
Vaisman.\index[terms]{vector field!parallel}

Or $a= 0$, 
and then $\theta=b(IX)^\flat$. 
Using again \eqref{DJ} and \eqref{domega} we can prove that
\begin{equation}\label{ljvg}
\Lie_{IX}g=b(g-X^\flat\otimes X^\flat-(IX)^\flat\otimes (IX)^\flat).
\end{equation}
Observe now that $[X,IX]=0$, 
because 
\[ [X,IX]= \nabla_X IX - \nabla_{IX} X =\nabla_X(I)(X),
\] since
$X$ is parallel\index[terms]{vector field!parallel}. Then,
by \eqref{DJ}, we have
\[
[X,IX]= \nabla_X(I)(X) =
\frac 12 \left(\theta(IX)X-\theta(X)IX+g(X,X)I\theta^\sharp\right).
\]
Since $\theta=b(IX)^\flat$, we have $\theta(X)=0$.
Using $|X|=1$ and $\theta^\sharp=bIX$, we obtain
\[
[X,IX]= \frac 12 \left(\theta(IX)X + I\theta^\sharp\right)=
\frac 12 \left(b X +  I(bIX) \right)=0.
\]
Therefore, the distribution
$\cav:=\langle X, IX\rangle$ is a foliation. Let
$\cad:=\cav^\perp$. By \eqref{_daX+dbIX_}, $db\wedge(IX)^\flat=0$, 
which gives
$$0=d\theta=db\wedge(IX)^\flat+bd(IX)^\flat=bd(IX)^\flat,$$
hence $d(IX)^\flat=0$. Since also  $dX^\flat=0$, because $X$ is parallel,\index[terms]{vector field!parallel} $\cad$ is a foliation. Moreover, $[X,\cad]\subset\cad$, 
because $X$ is parallel, and
$[IX,\cad]\subset\cad$ by another application of 
\eqref{DJ}. Then, by Frobenius theorem,\index[terms]{theorem!Frobenius} one finds local coordinates $(s,t,x)$ on $M$ such that $X=\6/\6 s$, $IX=\6 /\6 t$, and $g=ds^2+dt^2+h(s,t)$ were $h(s,t)$ is  a metric on the leaves of $\cad$. The complex structure $I$ restricts to the leaves of $\cad$ and the corresponding Hermitian form\index[terms]{form!Hermitian} is precisely the restriction of $\omega$. But then $d\omega=\theta\wedge\omega$ vanishes on $\cad$ since $\theta\restrict\cad=0$ (because it is a multiple of $(IX)^\flat$). This proves that $h(s,t)$ is a K\"ahler metric.

In fact, since $X$ is parallel\index[terms]{vector field!parallel}, its flow preserves $g$ and then $h$ only depends on $t$. One then shows that  also $b$ only depends on $t$: $b=b(t)$. Now, from \eqref{ljvg} we find
$$\frac{\6 h}{\6 t}=2b(t)h,$$
therefore, $h(t)=e^{c(t)}h(0)$, with $c(t)=\int_0^tb(u)du$. 

By passing to the universal cover we obtain:

\hfill

\claim
Let $X$ be a parallel\index[terms]{vector field!parallel} vector field on an LCK manifold
satisfying $\theta=b(IX)^\flat$.
Then the universal cover of $M$ is biholomorphic with
$\R^2\times (N,g_N)$, where $N$ is a simply connected
K\"ahler manifold and $g_N=h(0)$.  \endproof

\hfill

\remark 
The K\"ahler form\index[terms]{form!K\"ahler} $\tilde \omega=ds\wedge
dt+e^{c(t)}\omega_N$ of $\R^2\times N$ satisfies
$d\tilde\omega=dc\wedge \tilde\omega$,thus being 
globally conformally K\"ahler, with Lee form\index[terms]{form!Lee}
$\tilde\theta=dc$.

\hfill

To finish the proof of \ref{parafield}, it now suffices to
show that the GCK metric of the covering (or a multiple of
it) descends to $M$, yielding a contradiction. From the
above remark, all we need to check is that the real function $c$ is
invariant with respect to the deck group $\Gamma$ of the
covering. As the next result shows, the K\"ahlerianity of
$g_N$ does not play any role here as one sees from the
following more general:

\hfill

\lemma \label{_action_on_R^2_times_N_Lemma_}
Let $(N,g_N)$ be complete and simply connected Riemannian
manifold, $c:\R\rightarrow \R$ a smooth function, and
$\Gamma$ a cocompact  group of isometries acting totally
discontinuously on the
Riemannian product $(\bar N, \bar g)=(\R^2\times N,
ds^2+dt^2+e^{c(t)}g_N)$. If $\Gamma$ preserves the vector
fields $\frac{d}{ds}$ and $\frac{d}{dt}$,  then $c$ is
$\Gamma$-invariant.

\hfill

\proof 
Since $\Gamma$ preserves the vector fields $\frac{d}{ds}$ and $\frac{d}{dt}$, its elements can be written 
$$\ga(s,t,x)=(s+s_\ga,t+t_\ga, \psi_\ga(s)), \quad \psi_\ga\in\Diff(N), s_\gamma\in\R, t_\gamma\in\R.$$
This is an isometry of $\bar g$ if and only if\index[terms]{isometry}
$$e^{c(t)}g_N=e^{c(t+t_\ga)}\psi_\ga^*g_N,$$
and hence $\psi_\gamma$ is a homothety of $(N,g_N)$ of ratio 
\begin{equation}\label{_rho_gamma_Equation_}
	\rho_\ga:=e^{c(t)-c(t+t_\ga)}.
\end{equation}
For later use, note that for all
$\gamma_0\in \Gamma$, one has
\begin{equation}\label{nc}
c(nt_{\ga_0})=c(0)-n\log (\rho_{\ga_0}),\quad n\in\N.
\end{equation}
Indeed, \eqref{_rho_gamma_Equation_} implies
that the map $\gamma \arrow c(t) - c(t+t_\gamma)$
defines a group homomorphism $\Gamma \arrow (\R,+)$ for
all $t\in \R$, while the map $\gamma \arrow t_\gamma$ is also a
homomorphism. Then $c(t+t_{n\gamma})- c(t)= -n \log(\rho_{\gamma})$.
Setting $t=0$ and $\gamma=\gamma_0$, we obtain
$c(nt_{\gamma_0})= c(0) - n \log(\rho_{\gamma_0})$.

Return to the proof of \ref{_action_on_R^2_times_N_Lemma_}.
If, by absurd, $c$ is not $\Gamma$-invariant, we may
choose $\ga_0$ such that $\rho_{\ga_0}<1$. This implies
that $\psi_{\ga_0}$ is a contraction on $(N,g_N)$,
considered to be  a metric space;\index[terms]{space!metric} therefore
$\psi_{\ga_0}$ has a unique fixed point, denoted by
$x_0$. Moreover, we have $\lim_k\psi_{\ga_0}^k(x)=x_0$.

It follows that, for any $\ga\in\Gamma$,  the sequence
$$y_k:=(\ga_0^k\circ\ga\circ\ga_0^{-k})(0,0,x_0)=(s_\ga,t_\ga,\psi_{\ga_0}^k(\psi_\ga(x_0)))$$ 
converges to
$$y_0:=(s_\ga,t_\ga, x_0).$$
However, the action is totally
discontinuous, so
$y_k=y_0$ for $k$ large enough. Acting on both sides of
the equation
 $x_0= \psi_{\ga_0}^k(\psi_\ga(x_0))$ by
$\psi_{\ga_0}^{-k}$,
we obtain $x_0=\psi_\ga(x_0)$; 
$x_0$ is then a common
fixed point for all $\psi_\ga$.

Now let  $f:\bar N\rightarrow \RR^{>0}$,
$f(s,t,x)=e^{c(t)}d_N(x,x_0)$, where $d_N$
is the metric on $N$ induced by the Riemann
structure $g_N$. By \eqref{nc}, $c$ is surjective on $\R$,
and hence $f$ is surjective on $\R^{>0}$.

From the definition of $\rho_\ga$, \eqref{_rho_gamma_Equation_}, we  derive:
$$(\ga^*f)(s,t,x)=f(s,t,x).$$
This implies that $f$ is $\Gamma$-invariant, thus induces a   map $\tilde f:\frac{\R^2\times N}{\Gamma}\arrow \R$ that is  continuous and surjective, because $f$ is surjective. This contradicts the fact that $\Gamma$ acts co-compactly on $\R^2\times N$.
\endproof

\section{LCK metrics are not Einstein}

One of the most important results concerning the Riemannian
geometry of LCK metrics is the following:

\hfill

\theorem \label{_LCK_Einstein_then_Kahler_Theorem_}
{(\cite{_Madani_Moroianu_Pilca:Einstein_})} Let 
$(M,I,g,\omega,\theta)$ be
a compact $n$-dimensional LCK manifold, $n\geq 2$,
that is  not K\"ahler. Then $g$ cannot be Einstein.

\hfill

\proof  By \index[persons]{Myers, S. B.} Myers' theorem,\index[terms]{theorem!Myers} 
it is clear that an  LCK metric on $M$
cannot be Einstein if $M$ is compact and
the Einstein constant is positive. Indeed, by Myers', the fundamental
group\index[terms]{fundamental group} of such a manifold is finite.
However, $\pi_1(M)$  has to admit a non-trivial character
$\chi:\pi_1(M)\ra\R^{>0}$ (``homothety character'',\index[terms]{homothety character} see Subsection
\ref{carac}). By contrast, the proof of 
\ref{_LCK_Einstein_then_Kahler_Theorem_} for non-positive
Einstein constant $\lambda$ is completely non-trivial.

By a long and very ingenious computation, 
\index[persons]{Madani, F.} Madani, \index[persons]{Moroianu, A.} Moroianu and \index[persons]{Pilca, M.} Pilca prove that for some constant
$c\in\R^{>0}$, the smooth function
$f=d^*\theta+c|\theta|^2$
 satisfies the equality:
\begin{equation}\label{_Differential_of_f:Equation_}
	df=(2\la-c'f+\phi(n)\vert\theta\vert^2)\theta,
\end{equation}
where $c'\in\R^{>0}$ and $\phi$ is a function
$\phi:\; \Z^{\geq2} \to \R^{\leq 0}$
(the precise values of $c, c', \phi$ 
are not important: what matters are their signs). 
Recall that for any 1-form $\eta$, the codifferential is 
$d^*\eta=-\tr(\nabla\eta)=-i_{e_j}\nabla_{e_j}\eta$, where $\{e_j\}$ is any local orthonormal base of vector fields. Then the Laplacian of $f$ is computed as
\begin{equation}\label{_Laplacian_of_f:Equation_}
	\Delta
        f=d^*df=(2\la-c'f+\phi(n)\vert\theta\vert^2)d^*\theta
        -
c'\Lie_{\theta^\sharp}(f)-\phi(n)\Lie_{\theta^\sharp}(\vert\theta\vert^2).
\end{equation}
Suppose now, by absurd, that 
$\theta$ is not identically zero. From the definition of
$f$, $f=d^*\theta+c|\theta|^2$, we obtain $\int_M
f=c\int_M\vert\theta\vert^2d\!\Vol >0$. 
Since $M$ is compact, there exists a point $x_0\in M$ where $f$ attains
its maximum. We must have $f(x_0)>0$, because the integral of $f$ on $M$ is strictly positive.   Then
$(df)\restrict{x_0}=0$, and $(\Delta f)\restrict{x_0}\geq 0$, so that equation
\eqref{_Differential_of_f:Equation_}, evaluated in $x_0$, implies
$\theta\restrict{x_0}=0$, because $2\la-c'f+\phi(n)\vert\theta\vert^2<0$. We thus have
$$0<f(x_0)=(d^*\theta)(x_0)+c|\theta|^2(x_0)=(d^*\theta)(x_0).$$ 
We evaluate \eqref{_Laplacian_of_f:Equation_} at $x_0$ and obtain
$$(\Delta f)(x_0)=(2\la-c'f(x_0))(d^*\theta)(x_0)<0,$$
contradicting the choice of $x_0$ as a maximum point for $f$. \endproof

\hfill

Recall that a Riemannian manifold $(M,g)$ is {\bf locally symmetric}\index[terms]{manifold!locally symmetric} if its curvature tensor $R^\nabla$ is parallel with respect to the Levi--Civita connection$\nabla$ of $g$. A locally symmetric space is {\bf irreducible} if it is not the Riemannian product of two locally symmetric spaces. We\index[terms]{space!locally symmetric} have:\index[terms]{tensor!curvature!parallel}

\hfill

\proposition (\cite{_Madani_Moroianu_Pilca:Einstein_}) A compact, irreducible, and locally symmetric LCK manifold $(M, g, I, \theta)$
has vanishing Lee form,\index[terms]{form!Lee} i.e $(M,g,I)$ is K\"ahler.

\hfill

\proof 
Recall that a complete and simply connected locally symmetric space is symmetric (\cite[Theorem 6.3, vol. 2]{KoNFD2}. Then a compact irreducible symmetric spaces is a finite quotient of an irreducible symmetric space; hence it is Einstein (\cite[Corollary 8.7, vol. 2]{KoNFD2}. 

From the proof of \ref{_LCK_Einstein_then_Kahler_Theorem_}, it is clear that if the Einstein constant is non-positive, then $\theta=0$.

Instead, when the Einstein constant is positive,  from \ref{_LCK_Einstein_then_Kahler_Theorem_} we can only deduce that $(M,g,I,\theta)$ is GCK. Let then $\theta=d\phi$, $\phi\in C^\infty M$ and $g'=e^{-\phi}$ the corresponding global K\"ahler metric. Let $X$ be any $g$-Killing field\index[terms]{vector field!Killing} (it exists, because a Riemannian symmetric space is, in particular, Riemannian homogeneous). Then $X$ is $[g']$-conformal. By \ref{_Lichnerowicz_conformal_Theorem_}, $X$ is $g'$-Killing. Then $X$ is Killing for both metrics $g$ and $g'$ that are  conformal: $\Lie_X(g')=\Lie_Xg=0$. This gives $$\Lie_X(g')=\Lie_X(e^{-\phi}g)=(\Lie_X(e^{-\phi}))g+e^{-\phi}\Lie_X(g)=(\Lie_X(e^{-\phi}))g=0,$$
therefore, $\Lie_X(\phi)=0$. Since $X$ was arbitrarily chosen and $(M,g)$ is homogeneous, it follows that $\phi=\const$ and $\theta=d\phi=0$. \endproof

\section[The Vaisman condition in terms 
of the Bismut connection]{The Vaisman condition in terms 
of the\\ Bismut connection}\label{_Vaisman_via_Bismut_Section_}
 
Let $(M, I, \omega, g)$ be a
complex Hermitian manifold. Consider a 
connection $\nabla^b:\; TM \arrow TM \otimes \Lambda^1 M$ 
with the following properties
\begin{description}
\item[(i)] $\nabla^b(I)=\nabla^b(\omega)=0$.
\item[(ii)] Let $T^b\in \Hom(\Lambda^2TM, TM)$
be the torsion tensor of $\nabla^b$.
Then the 3-form $X,Y,Z \arrow g(T(X,Y), Z)$
is skew-symmetric.
\end{description}

J.-M. \index[persons]{Bismut, J.-M.} Bismut proved the following theorem:

\hfill

\theorem (\cite{_Bismut:connection_})
Let $(M, I, \omega, g)$ be a
complex Hermitian manifold.
Then the connection $\nabla^b$
satisfying the properties (i)-(ii)
exists and is unique.
\endproof

\hfill

This connection is called
{\bf the \index[persons]{Bismut, J.-M.} Bismut connection}.
It can be expressed in terms of the 
Levi--Civita connection $\nabla$ of $g$ as
\[
g(\nabla^b_XY,Z)=g(\nabla_XY,Z)+\frac 12 d\omega(IX,IY,IZ). 
\]
(Exercise \ref{_Bismut_torsion_Exercise_}).


%


On an LCK manifold, $d\omega= \omega\wedge\theta$,
which gives the following expression of the \index[persons]{Bismut, J.-M.} Bismut torsion:
\[
2T^b(X,Y)=\theta(IX)IY-\theta(IY)IX-\omega(X,Y)I\theta^\sharp.
\]

\proposition 
(\cite{_Andrada_Villacampa_}) Let
$(M,I,g,\eta)$ be an LCK manifold. Then $\nabla\theta=0$
if and only if $\nabla^b\theta=0$.
\endproof


\hfill

\proposition\label{_Bismut_Vaisman_AV_Proposition_} 
(\cite{_Andrada_Villacampa_}) An LCK manifold is Vaisman
if and only if the torsion 3-form of the  Bismut connection
is parallel with respect to the Bismut connection.
\endproof

\hfill

Let $(M,I,g)$ be a Hermitian manifold and $\nabla$ a
metric connection with curvature tensor $R^\nabla$. Recall
(\cite{_Angella_Otal_Ugarte_Villacampa_}) that $\nabla$
is called {\bf K\"ahler-like} if it satisfies the first
Bianchi identity, $\sum_{cycl}R^\nabla(x,y)z=0$, and is of
type $(1,1)$, i.e. $R^\nabla(Ix,Iy)=R^\nabla(x,y)$. In
\cite{_Yau_Zhao_Zheng_}, a Hermitian manifold whose \index[persons]{Bismut, J.-M.} Bismut
connection is K\"ahler-like is called {\bf \index[persons]{Strominger, A.} Strominger
  K\"ahler-like} (SKL for short). The following result is
proven in :

\hfill

\proposition (\cite{_Yau_Zhao_Zheng_})\\
Let $(M,I,g)$, $\dim_\C M\geq 3$,  be a
compact SKL manifold with $g$ not K\"ahler. Then $(M,I)$
cannot admit any Vaisman metric.

\hfill

Moreover,  Yau, Zhao,
Zheng conjecture that a compact SKL manifold of dimension $\geq 3$ and with $g$ not K\"ahler cannot admit LCK metrics.
\index[terms]{conjecture!Yau, Zhao and Zheng}
\section{Notes}

\subsection{Bismut connections}\label{_Bismut_connection_Note_}

The \index[persons]{Bismut, J.-M.} Bismut connection appeared implicitly in
\cite{_Strominger_Symmetries_}. This is why Yau, Zhao,
Zheng call it {\bf \index[persons]{Strominger, A.} Strominger connection} 
(\cite{_Yau_Zhao_Zheng_}).

\hfill

Let $(M,I, \omega, g)$ be complex Hermitian manifold, and $\nabla$ the Levi--Civita connection of $g$.  In \cite{_Gauduchon:connections_}, P. \index[persons]{Gauduchon, P.} Gauduchon constructed the following  
1-parameter family of Hermitian connections: 
\[
g(\nabla^t_XY,Z)=g(\nabla_XY,Z)-\frac{t-1}{4} d\omega(IX,IY,IZ) -\frac{t+1}{4}d\omega(IX,Y,Z). 
\]
The \index[persons]{Bismut, J.-M.} Bismut connection corresponds to $t=-1$ in this family.

\subsection{Curvature properties}

The Riemannian geometry of an LCK metric is rather restrictive. For example, a manifold with LCK metric of zero sectional
curvature is K\"ahler. This was proven in \cite{va_tr2} for compact manifolds (see also \cite[Corollary 3.2]{_Madani_Moroianu_Pilca:Einstein_}); however,  in 
\cite{bbh} it was observed that the result is local.  Also, if on a compact manifold, the Riemann curvature tensor of an LCK metric satisfies 
$R_{XYZV}=R_{XY(IZ)(IV)}$ (or, equivalently, the 
curvature operator $R_{XYZ}^V$ belongs to $\Lambda^2(M)\otimes \goth{u}(TM)$), 
the metric is K\"ahler (\cite{va_tr2}). For other curvature properties of the LCK metrics, see \cite{do}.

\subsection{Harmonic maps and distributions} 

Let $f: (M^m,g)\arrow (N^n,h)$ be a smooth map between Riemannian manifolds. Its differential $Df$ can be viewed as a section of $T^*M\otimes f^{-1}TN$. The {\bf Hilbert--Schmidt norm}\index[terms]{Hilbert--Schmidt norm} is defined on sections of $T^*M\otimes f^{-1}TN$ in the following way. Let $\{e_1,...,e_m\}$ be a local orthonormal frame on $M$ and put
\[
\Vert Df\Vert:=\tr_g(f^*h)=\sum_{i=1}^m h(Df(e_i), Df(e_i)).
\]
The {\bf energy density} $e(f)$ of the map $f$ is the smooth function $M\ni x\mapsto \Vert Df\Vert^2\in \R^{>0}$. Now let  $M$ be compact and define the  {\bf energy functional}\index[terms]{energy functional} on $C^\infty(M,N)$ by $E(f):=\int_M e(f)$. Then $f$ is called {\bf harmonic} if it is a critical point of the energy functional, \cite{_Eells_Sampson_}. \index[terms]{map!harmonic}

\index[persons]{Lichnerowicz, A.} Lichnerowicz proved (\cite{_Lichnerowicz:harmonic_}) that a holomorphic map between K\"ahler manifolds is harmonic. For LCK manifolds we have:\index[terms]{theorem!Lichnerowicz}

\hfill

\theorem {(\cite{iov}, \cite{ghe})} A holomorphic map
between compact LCK ma\-nifolds is harmonic if and only if
it is constant on the orbits of the Lee field.\index[terms]{Lee field}

%
	
\hfill

%
%
%

A $k$-dimensional distribution on a manifold $M$ can be
viewed as a smooth map from $M$ to the oriented
Grassmannian of $k$-planes in $TM$. A Riemannian metric on
$M$ induces a Riemannian metric on this Grassmannian,
called {\bf the \index[persons]{Sasaki, S.} Sasaki metric}, generalizing the Sasaki metric
on the tangent bundle, see \cite{ggv}. One can thus speak
about harmonic and minimal distributions. Vaisman manifolds offer natural examples\index[terms]{manifold!Vaisman} of\index[terms]{metric!Sasaki}\index[terms]{distribution!harmonic}\index[terms]{distribution!minimal}
this kind. One obtains:

\hfill

\theorem {(\cite{ova})} Let $(M,I,\omega,\theta)$ be a compact Vaisman manifold. Then:
\begin{description}
\item[(i)] $I\theta^\sharp$ is  harmonic and minimal. It
  is also harmonic as a map into the unit tangent bundle
  of $M$ with the \index[persons]{Sasaki, S.} Sasaki metric.
\item[(ii)] The canonical foliation\index[terms]{foliation!canonical} induces a harmonic map into $\Gr^{or}_2(M)$.
\end{description}

On the other hand, the fundamental form $\omega$ of an LCK
manifold $M$ has constant length (with respect to the LCK
metric). After appropriate normalization, each power $\omega^r$, $r\leq n=\dim_\C \index[terms]{normalization}
M$, can be viewed as a section of the unit sphere bundle in
$\Lambda^{2r}(M)$. 

\hfill

\theorem {(\cite{gms})} Let $(M,\omega,\theta)$ be a
compact  LCK manifold, $\dim_\R M=2n$. If the Lee form \index[terms]{form!Lee} is nowhere zero,
and $r<n$, where then $\omega^r$ is a harmonic
section of its corresponding sphere bundle if and only if
$2r = n$. On a $4n$-dimensional compact Vaisman manifold,
$\omega^n$ is harmonic as a map from $M$ to\index[terms]{manifold!Vaisman}
$\Tot(\Lambda^{2n}(M))$ equipped with the Sasaki metric.


\chapter{Einstein--Weyl manifolds and the Futaki invariant}\index[terms]{manifold!Einstein--Weyl}

\epigraph{\it The truth of life lies in the impulsiveness of matter. The mind of man has been poisoned by concepts. Do not ask him to be content, ask him only to be calm, to believe that he has found his place. But only the madman is really calm.}{Antonin Artaud}

\hfill

In this chapter, we discuss several works on LCK geometry
related to  the Ca\-la\-bi-Yau theorem: the
\index[terms]{Futaki invariant} Futaki invariant, which was (in the K\"ahler context)
the first known obstruction to the existence of the 
K\"ahler--Einstein metrics, and the Einstein--Weyl
condition, that is  the closest approximation\index[terms]{Futaki invariant}
of the Einstein condition in the LCK context
(if an LCK metric is Einstein, it is in fact
K\"ahler, as shown in 
\cite{_Madani_Moroianu_Pilca:Einstein_}; see 
also \ref{_LCK_Einstein_then_Kahler_Theorem_}).
As proven by P. \index[persons]{Gauduchon, P.} Gauduchon (\ref{ewcy}), any compact Einstein--Weyl\index[terms]{manifold!Einstein--Weyl} manifold
is Vaisman; all results of this chapter 
properly belong to the Vaisman geometry.


\section{The Einstein--Weyl condition}\label{eiwey}


The class of Einstein metrics\index[terms]{metric!Einstein} is not conformally invariant. In the
conformal world, the appropriate Einstein-like notion is
called Einstein--Weyl\index[terms]{manifold!Einstein--Weyl} and refers to the {\bf symmetrized
  Ricci tensor of the Weyl connection}.\index[terms]{tensor!Ricci!symmetrized}
Recall that a Weyl connection (Appendix \ref{_Weyl_conne_Appendix_})\index[terms]{connection!Weyl}
on a conformal manifold is a
connection that is  torsion-free \index[terms]{connection!torsion-free}and preserves the
conformal class. Clearly, any Levi--Civita connection
preserving a metric in this conformal class in Weyl;
however, there are Weyl connections that do not
preserve any metric.

\hfill

\definition \label{_EW_Definition_}
An {\bf Einstein--Weyl manifold}\index[terms]{manifold!Einstein--Weyl} is a Riemannian manifold locally conformally equivalent to an Einstein manifold.\index[terms]{manifold!Einstein}

\hfill

Let $D$ be the Weyl connection of an LCK manifold
$(M,I,g,\theta)$, see \ref{_Weyl_conn_definition_Remark_}.  Note that the Ricci\index[terms]{tensor!Ricci}
tensor of a Weyl connection is not symmetric, in general,
see e.g. \cite{gau_wein} (indeed, a Weyl connection is not 
necessarily a Levi--Civita connection of a metric). However, the Ricci tensor of the
Weyl connection of an LCK manifold is symmetric, because it is the Levi--Civita 
connection of a K\"ahler metric on a K\"ahler cover\index[terms]{cover!K\"ahler}.

\hfill

\definition\label{_Einst_Weyl_Definition_} 
A conformal manifold with Weyl connection is {\bf Einstein--Weyl}\index[terms]{manifold!Einstein--Weyl}
if the Ricci tensor of the Weyl connection is conformal to
the metric: 
\begin{equation}\label{_Einstein_Weyl_Equation_}
\Ric(D)=\lambda \omega, \ \  \lambda\in C^\infty (M).
\end{equation}
  An Einstein--Weyl LCK-manifold is also called
{\bf Hermitian Einstein--Weyl}. The function $\lambda$
is called {\bf the Weyl scalar curvature}, denoted by
$\Scal^D$.\index[terms]{manifold!LCK!Einstein--Weyl}

\hfill

\remark
Let $M$ be an LCK manifold, $\tilde M$ its K\"ahler cover,
and $D$ its Weyl connection. Then $D$ is the Levi--Civita
connection on $\tilde M$, and $\Ric(D)$ is proportional
to the metric tensor $\tilde g$ on $\tilde M$:
\[ \Ric(D)=\lambda \tilde g,\] with $\lambda=\const$. Since $\Ric(D)$
is defined on $M$, this means that, unless $\lambda=0$,
the LCK manifold $M$ is globally conformally K\"ahler.
Since we assume from the beginning that our LCK manifolds
are not globally conformally K\"ahler, for an
LCK manifold the equation \eqref{_Einstein_Weyl_Equation_}
becomes $\Ric(D)=0$.

\hfill



\theorem{(\cite{gau_wein})}\label{lchk_vai} 
On a compact Einstein--Weyl manifold of real dimension $m\geq 3$\index[terms]{manifold!Einstein--Weyl}
and with non-trivial, co-closed Lee form,\index[terms]{form!Lee} the Lee field\index[terms]{Lee field} is
Killing.\index[terms]{vector field!Killing} If, moreover,  the Lee form is closed, then  the
Ricci tensor\index[terms]{tensor!Ricci}  of the Weyl connection vanishes identically
and the Lee form is parallel.\index[terms]{form!Lee!parallel}

\hfill

\proof By direct computation, starting from the relation
between the Weyl connection $D$ and Levi--Civita connection
$\nabla$  (Appendix \ref{_Weyl_conne_Appendix_}), one proves the relation:
\begin{equation*}
\frac 1m d\Scal^D+2d^* \nabla\theta -
d^* d\theta+2\nabla_{\theta^\sharp}\theta +(m-3)d|\theta|^2=0,
\end{equation*}
where all traces and operators are considered with respect
to $g$, and $d^*$ is the dual to the de Rham differential.
After contracting with $\theta$ and then integrating on $M$, one finds
\begin{equation}\label{tod}
2\nabla\theta=d\theta,
\end{equation}
which says that $\nabla\theta$, and thus $\nabla\theta^\sharp$  is $g$-antisymmetric, which means that  $\theta^\sharp$ is Killing and,\index[terms]{vector field!Killing} if $\theta$ is closed, it is in fact parallel.\index[terms]{vector field!parallel}
As for the Ricci tensor,\index[terms]{tensor!Ricci} one starts with the equation:
$$\Delta\theta=\frac 2m \Scal^D\cdot\theta,$$
where $\Scal^D$ is the Weyl scalar curvature (\ref{_Einst_Weyl_Definition_}).
Contracting this equation and integrating on $M$, we obtain
$$\int_M|d\theta|^2 =\frac 2m\int_M\Scal^D\cdot |\theta|^2,$$
which, if the Lee form\index[terms]{form!Lee} is closed, implies $\Scal^D=0$. \endproof

\hfill

Since the Lee form\index[terms]{form!Lee} of a Gauduchon metric on an LCK\index[terms]{metric!Gauduchon}
manifold is harmonic, and the Weyl connection $D$ locally
coincides with the Levi--Civita connection of the local
K\"ahler metrics, we have:

\hfill

\corollary\label{ewcy} 
The Gauduchon metric\index[terms]{metric!Gauduchon} on a compact
Einstein--Weyl LCK manifold of complex dimension at least
$2$ is Vaisman.\index[terms]{manifold!LCK!Einstein--Weyl} The K\"ahler covers of a compact
Einstein--Weyl LCK manifold are Ricci-flat, that is,
have the holonomy in $\SU(n)$, the same as for the
 Calabi--Yau manifolds. In particular,
the K\"ahler cone of a compact Einstein--Weyl manifold is the Riemannian cone\index[terms]{cone!Riemannian} over an\index[terms]{manifold!Einstein--Weyl} Sasaki--Einstein\index[terms]{manifold!Sasaki--Einstein} manifold.\index[terms]{manifold!Calabi--Yau} 

\hfill

\remark Locally conformally hyperk\"ahler manifolds
(LCHK manifolds) are LCK Einstein--Weyl, see Chapter \ref{other}, \S \ref{lchk}.\index[terms]{manifold!LCHK}

\hfill

We end with the following result:

\hfill

\proposition {(\cite[Proposition 5.6]{_Verbitsky:Vanishing_LCHK_})}  Let $M$ be a compact Einstein--Weyl LCK manifold of complex dimension $n$, $K$ its canonical bundle and $L_\C$ its complexified weight bundle, both considered to be  Hermitian holomorphic bundles with the metric induced from $M$. Then $L_\C^n\simeq K^{-1}$.\index[terms]{manifold!Einstein--Weyl}\index[terms]{bundle!weight}

\hfill

\proof  By \ref{ewcy}, the universal cover $\tilde M$ of $M$ is Calabi--Yau, that is $\Hol(\tilde M)\subset \SU(n)$. On the other hand,  
$\Hol(M)\subset \SU(n)\cdot \R^{>0}$, since the Levi--Civita connection of the LCK metric preserves the fixed conformal class. By the  definition of the weight bundle $L$ as associated with the representation $\det(TM)^{\frac{1}{n}}$, $n=\dim_\C M$, we find that the holonomy of the flat bundle $L$ is the quotient $G:=\Hol(M)/\Hol(\tilde M)\subset \R^{>0}$. The action of $G$ can be made explicit by: 
$$(\al, l)\mapsto \al\cdot l, \quad \al\in G, l\in  L.$$

Note that, since $\Hol(\tilde M)\subset \SU(n)$, the action of  $\Hol(M)$  on $K=\det(\Lambda^{1,0}M)$ factors through $G$ as:
$$(\al,\eta)\mapsto \al^{-n}\eta, \quad \al\in G, \eta\in K,$$ 
relating the monodromy\index[terms]{action!monodromy} of $K$ with holonomy action. As $L^n$ and $K^{-1}$ are both flat, they have the same monodromy, thus being isomorphic. \endproof\index[terms]{group!holonomy}

\hfill

This result has the following vanishing consequence:

\hfill

\corollary
    {(\cite{_Verbitsky:Vanishing_LCHK_})}\label{EW_vanish}
    Compact Einstein--Weyl, non-K\"ahler,  LCK manifolds
    satisfy: $H^i(\calo_M)=0$ for $i>1$, and $\dim
    H^1(\calo_M)=1$.\index[terms]{manifold!Einstein--Weyl}

\hfill

Indeed the two properties hold, more generally, for all compact Vaisman manifolds\index[terms]{manifold!Vaisman} whose canonical bundle is a negative power of the weight bundle, \cite[Theorem 8.4]{_Verbitsky:Vanishing_LCHK_}).\index[terms]{bundle!weight}
 
\hfill 

As the Dolbeault spectral sequence converges to the de Rham cohomology,\index[terms]{cohomology!de Rham} we have:\index[terms]{spectral sequence!Dolbeault}
$$h^1(M)\leq \dim H^0(\Omega^1(M))+H^1(\Omega^0(M))=1.$$
However, $h^1(M)$ cannot be zero: were it zero, then the weight bundle of $M$ would have trivial monodromy\index[terms]{monodromy} and $M$ would be K\"ahler, contradiction. We then have:

\hfill

\corollary\label{betti_1}
{(\cite{kashiwada2}, \cite{_Verbitsky:Vanishing_LCHK_})} 
Let $M$ be a compact Einstein--Weyl LCK manifold. Assume
that $M$ is not K\"ahler. Then $b_1(M)=1$.\index[terms]{manifold!Einstein--Weyl}

\section{The Futaki invariant of Hermitian manifolds}

 The \index[persons]{Futaki, A.} Futaki  integral invariant is an essential tool in
 K\"ahler geometry.\index[terms]{geometry!K\"ahler} It is defined on the Lie algebra of
 holomorphic vector fields of a compact K\"ahler manifold
 and its vanishing is necessary for the existence of a
 K\"ahler--Einstein metric\index[terms]{metric!K\"ahler--Einstein}
 on a compact Fano manifold (\cite{besse}). A similar
 integral invariant can be defined on Hermitian
 manifolds. In this Section we consider the Futaki
 invariant in the LCK context
 and show it vanishes on Vaisman manifolds.\index[terms]{manifold!Vaisman}

Let $M$ be a compact complex manifold.
The {\bf Futaki invariant} is a character on its
group of biholomorphic automorphisms, defined by
A. \index[persons]{Futaki, A.} Futaki when $M$ is K\"ahler in \cite{_Futaki:kahler_} and
for the general complex manifolds in \cite{_Futaki:complex_}.
For K\"ahler manifolds, non-zero Futaki invariant is one of the
obstructions to the existence of metrics with constant
scalar curvature (and hence, the K\"ahler--Einstein metrics).
For K\"ahler manifolds, one can write the Futaki invariant as
\[ X \mapsto \frac1 {\int_M \omega^n}\int_M \Lie_X (S) \omega_M^{n}=
 -  \frac1 {\int_M \omega^n}\int_M  S \Lie_X(\omega_M^{n}),
\]
where $(M, \omega)$ is an $n$-dimensional complex manifold,
$X$ a holomorphic vector field on $M$, and $S$ the scalar 
curvature of $\omega$ (\cite[(4.12.1)]{_Gauduchon:book_}.
From this definition, it is clear that this invariant vanishes
on K\"ahler manifolds of constant scalar curvature. A few
years after \index[persons]{Futaki, A.} Futaki discovered the eponymous invariant, he found another
(and much simpler) form of the Futaki invariant, that is  
well-defined for all complex manifolds.\index[terms]{Futaki invariant}

\hfill

\definition \label{_Futaki_Ricci_Definition_}
 Let $(M,I)$ be a complex manifold.
A volume form $\Omega$ on $M$ can be interpreted as 
a Hermitian structure on its canonical bundle $K_M$,
with $|\alpha|^2:= \frac{\alpha\wedge\bar\alpha}{\Omega}$,
for any section $\alpha$ of $K_M$.
 {\bf The Ricci form}\index[terms]{form!Ricci}  $\rho_\Omega$ is half of the curvature
of the Chern connection\index[terms]{connection!Chern} associated with this metric,
expressed as in \eqref{_Chern_conne_curvature_dd^c_log_Equation_} by the formula \index[terms]{form!Ricci}
\begin{equation}\label{_ric_form_}
\rho_\Omega=-\1\6\bar\6 \log |\alpha|.
\end{equation}
where $\alpha$ is a local non-degenerate holomorphic section of $K_M$.
By \ref{_dd^c_log_independent_Claim_}, this expression is independent on  the
choice of $\alpha$. Traditionally, the Ricci form\index[terms]{form!Ricci} associated
with a volume form $\Omega$ is expressed by  $\rho_\Omega=-\1\6\bar\6 \log \Omega.$

{\bf The divergence} of a vector field $X$ is a function
defined as
$$\div (X):=\frac{d(i_X\Omega)}{\Omega}= \frac{\Lie_X\Omega}{\Omega}.$$
If $X$ is a vector field of type $(1,0)$, the formula for the divergence reduces to:
$$\div (X)\cdot\Omega=\6(i_X\Omega).$$
Usually, in this book, when we write a ``holomorphic vector
field'', we mean a real vector field $X$ that satisfies
$\Lie_X(I)=0$. However, the space of $(1,0)$-vectors on
$M$ is naturally equipped with a structure of a
holomorphic vector bundle. The projection
$T^{1,0}M \arrow TM$ along  $T^{0,1}M$ identifies
$TM$ with $T^{1,0}M$, and under this identification,
holomorphic sections of $T^{1,0}M$ correspond
to real vector fields that satisfy $\Lie_X(I)=0$
(\ref{_holo_1,0_real_Remark_}).

In this chapter, when we say ``a holomorphic vector
field'', we actually mean a holomorphic section of $T^{1,0}M$.

Now let  $M$ be a compact complex manifold, 
and $\mathfrak{h}(M)$ its Lie algebra of  holomorphic vector fields. 
For a volume form $\Omega$, the {\bf Futaki invariant} 
$f:\mathfrak{h}(M)\rightarrow \C$ is defined as:\index[terms]{Futaki invariant}
\begin{equation} \label{_Futaki_main_definition_Equation_}
f(X)=\int_M \div(X)\rho_\Omega^n.
\end{equation}
One can prove that for Fano manifolds the
definition reduces to \index[persons]{Futaki, A.} Futaki's original one (see
\cite{futaki} for details).

In this section, we prove that the Futaki 
invariant is independent on  the choice of the volume form.
In fact, the map $X \mapsto f(X)$ defines a character
on the Lie algebra of holomorphic vector fields, as shown in
\cite{_Futaki:complex_}. In this paper, Futaki
used a group version of the Futaki invariant,
defined as follows. Let $\sigma\in \Aut(M)$
be a holomorphic automorphism; then
\[ 
f(\sigma):= \int_M
\log\frac{\sigma^*\Omega}{\Omega}\sigma^*(\rho_\Omega^n)
\]
We do not use this group version of the \index[persons]{Futaki, A.} Futaki invariant.
It is not hard to see that the group version of the Futaki
invariant is compatible with the Lie algebra version;
we leave this as an exercise to the reader.

One can show further (\cite{futaki}) that the Futaki invariant can be written as:
\begin{equation}\label{fut2}
f(X)=-\int_MX\left(\frac{\rho_\Omega^n}{\Omega}\right)\Omega.
\end{equation}
The vanishing of $f$ is then an obstruction to the
existence of a volume form with constant 
$\frac{\rho_\Omega^n}{\Omega}$.

\hfill

\theorem {(\cite{fho})}\label{_fut_ind_}  
Let $M$ be a compact complex manifold, $\Omega$ a volume
form, $X$ a (1,0)-holomorphic vector field, and
$X \mapsto f(X)$ the corresponding \index[persons]{Futaki, A.} Futaki invariant.
Then $f(X)$ is independent on  the choice of $\Omega$.

\hfill

\proof  Let $\Omega_0, \Omega_1$ be volume forms 
on $M$. We connect $\Omega_0$ to $\Omega_1$ by a path
$\Omega_t= \phi^t \Omega_0$, where 
$\phi= \frac{\Omega_1}{\Omega_0}$ and $\phi^t= e^{t\log\phi}$.

Let $\div_t X:= \frac{\Lie_X \Omega_t}{\Omega_t}$ be the divergence function
 computed with respect to the volume form $\Omega_t$. Then \eqref{_Futaki_main_definition_Equation_}
 gives a 1-parameter family of Futaki invariants
$$f_t(X)=\int_M \div_t(X)\rho^n_{\Omega_t}.$$
To prove \ref{_fut_ind_} it is enough to show that
$\frac{d}{dt}f_t(X)=0$ for all   holomorphic vector fields $X$.

From \eqref{_ric_form_} we obtain:
$$\frac{d}{dt}(\div_tX)=X(\log\f),\qquad \frac{d}{dt}(\rho_{\Omega_t})=-\1\6\bar\6\log\f,$$
and hence
\begin{equation*}
	\begin{split}
		\frac{d}{dt}f_t(X)&=\int_M
                X(\log\f)\rho^n_{\Omega_t}-
\int_M\div_tX\1\6\bar\6(\log\f)\wedge n \rho^{n-1}_{\Omega_t}\\
		&=\int_M X(\log\f)\rho^n_{\Omega_t}+
		\int_M\bar\6 \left(\div_tX \wedge\6(\1\log\f)\wedge n\rho^{n-1}_{\Omega_t}\right)\\
		&-\int_M \bar\6(\div_tX)\wedge\6(\1\log\f)\wedge n\rho^{n-1}_{\Omega_t}.
	\end{split}
\end{equation*}
Now,  $$\int_M\bar\6 \left(\div_tX \wedge\6(\1\log\f)\wedge n\rho^{n-1}_{\Omega_t}\right)=0$$ 
by Stokes' theorem. We proceed as follows, using
\[ i_X (\rho^n_{\Omega_t}\wedge\6\log\phi)=
i_X(\rho^n_{\Omega_t}) \wedge\6\log\phi +
\rho^n_{\Omega_t} \cdot i_X (\6 \log\phi)
\]
for the last equation.
\begin{equation*}
	\begin{split}
		\frac{d}{dt}f_t(X)&=\int_M
                X(\log\f)\rho^n_{\Omega_t}-
\int_M \bar\6(\div_tX)\wedge\6(\1\log\f)\wedge n\rho^{n-1}_{\Omega_t}\\
		&=\int_M X(\log\f)\rho^n_{\Omega_t}-
\int_M(i_X\6\bar\6\log(\f^t))\wedge \6(\1\log\f)\wedge n\rho^{n-1}_{\Omega_t}\\
		&=\int_M
X(\log\f)\rho^n_{\Omega_t}+\int_M(i_X
\rho^n_{\Omega_t})\wedge\6\log\f\\ 
		&=\int_M i_X (\rho^n_{\Omega_t}\wedge\6\log\phi)=0,
	\end{split}
\end{equation*}
because $\deg(\rho^n_{\Omega_t}\wedge\6\log\f)>\dim M$.
\endproof

\hfill

%
%
%
%
%
%
%
%
%
%


\section{The Futaki invariant on LCK manifolds}


%
%
%

The main result of this chapter is the following theorem.

\hfill

\theorem \label{_Futaki_vanishes_Vaisman_Theorem_}
{(\cite{fho})}
 The Futaki invariant vanishes on compact
Vaisman manifolds.\index[terms]{Futaki invariant}\index[terms]{manifold!Vaisman}

\hfill

\pstep Let $(M,I,g)$ be a compact Vaisman manifold of complex dimension $n$. Recall that its universal cover is a cone $S\times \R^{>0}$ over a compact Sasakian manifold $S$ of real dimension $2n-1$, with Reeb field $\xi=I(r\6 r)$, where $r$ is the coordinate on $\R^{>0}$.\index[terms]{manifold!Sasaki} The restriction of the complex structure of the cone to $S$ is $I\restrict{S}=\nabla\xi$ (here $\nabla$ is the Levi--Civita connection of the \index[persons]{Sasaki, S.} Sasaki metric). It is easy to see  that the Ricci tensor\index[terms]{tensor!Ricci} of the cone metric vanishes on the Reeb field:
$\tilde\Ric(\xi, \xi)=0$, 
and thus the Ricci tensor degenerates on the orbits of the (1,0) vector field $X:=\xi-\1 I(\xi)$. Therefore, if $\tilde\Omega$ is the volume form associated with the K\"ahler metric of the cone, $\tilde\omega=dd^c r^2$, then $\rho_{\tilde\Omega}^n=0$.

\hfill

{\bf Step 2:} We now find a metric $g_0\in[g]$ such that
the Ricci form\index[terms]{form!Ricci} $\rho_0$ associated with its Riemannian
volume form $\Omega_0$
satisfies $\rho_0^n=0$. This will imply that the integrand
in \eqref{fut2} vanishes  identically.

Recall that  $\frac{1}{r^2}\tilde\omega=\frac{1}{r^2}dd^c
r^2$ is $\Gamma$-invariant and descends to $M$, producing
an LCK metric on $M$.  On the other hand, the Ricci form
of $\Omega_0$ is
$\rho_{\tilde\Omega}-\frac 1 2 dd^c \log r^2=
\rho_{\tilde\Omega}-\omega_0$, where $\omega_0$ is the
transverse K\"ahler form\index[terms]{form!K\"ahler} of the canonical foliation\index[terms]{foliation!canonical}
(\ref{_Subva_Vaisman_Theorem_}). This form degenerates on the orbit of
the flow of $X$. Hence, the $n$-th
power of the Ricci form\index[terms]{form!Ricci} of the Hermitian form\index[terms]{form!Hermitian}
$\frac{1}{r^2}\tilde\omega$ vanishes identically on
$M$. \endproof

\hfill

\remark Using a localization formula, \cite{futaki}, a
detailed computation is provided in \cite{fho} for the
\index[persons]{Futaki, A.} Futaki invariant on the blow-up at a point of a
Hopf surface. As expected, the result is non-zero,
confirming the known fact that blow-ups of LCK manifolds
cannot be Vaisman.\index[terms]{surface!Hopf}

\chapter{LCK structures on homogeneous manifolds}\index[terms]{manifold!homogeneous}\index[terms]{structure!LCK}

\epigraph{\it But if only you'd been with us on Sunday! We saw a red vineyard, completely red like red wine. In the distance it became yellow, and then a green sky with a sun, fields violet and sparkling yellow here and there after the rain in which the setting sun was reflected.}{\sc\scriptsize Vincent van Gogh, \ \ Letter to Theo van Gogh, Nov. 3, 1888}

\section{Introduction}

In this chapter, we discuss the existence of invariant LCK
structures on compact quotients of Lie groups. Two main
cases appear: the complex structure can be invariant or
not. In the first case, we have homogeneous LCK manifolds,
in the second\index[terms]{manifold!LCK!homogeneous} one, locally
homogeneous manifolds (we discuss them in the next
chapter). \index[terms]{manifold!locally homogeneous}


\definition
	An LCK manifold $(M,I,g)$ is {\bf homogeneous}\index[terms]{manifold!LCK!homogeneous} if it admits a  transitive and effective action of a   connected Lie group $G$  by holomorphic isometries of the LCK metric. In this case, $M=G/H$, with $H$  the stabilizer subgroup of $G$.

\hfill

\remark
	(i) A group $G$ as above also preserves $\omega$ and $\theta$. 
	
	(ii) Recall that the group $\Aut(M)$ of biholomorphic conformalities of a compact LCK manifold is compact. Indeed, any holomorphic conformality preserves the corresponding Gauduchon metric\index[terms]{metric!Gauduchon} that is  unique up to a constant. We can therefore, assume that $M$ is homogeneous under the group of holomorphic conformalities.

\hfill

The main result of this chapter is:

\hfill

\theorem \label{hom_vai} A compact homogeneous LCK manifold is Vaisman.\index[terms]{manifold!LCK!homogeneous}

\hfill

\remark This result was announced in \cite{HK1}, see also \cite{HK2} and \cite{guan} for a subsequent discussion. We  present a proof using the methods developed in this book, in particular centred on the notion of LCK structure\index[terms]{structure!LCK!with potential} with potential (see \cite{ov_parma}). This is different from the original proof in \cite{gmo}.

\section{Homogeneous LCK manifolds} 

The idea  for the proof of \ref{hom_vai} is to show that the homogeneity
implies  the existence of a holomorphic
circle action \index[terms]{action!$S^1$-} whose lift to the
universal cover $\tilde M$ does not
act by isometries of the K\"ahler structure on $\tilde M$. This will
imply that the manifold is LCK with potential\index[terms]{manifold!LCK!with potential}. We then
observe that the Lee form\index[terms]{form!Lee} of this LCK structure with
potential has constant length, which characterizes the
Vaisman metrics among the LCK metrics with potential.

\hfill

We start with the following lemma.

\hfill

\lemma\label{_S^1_and_theta_Lemma_}
Let $(M, \omega, \theta)$ be an LCK manifold, and 
$A$ a vector field acting on $M$ by holomorphic
isometries. Assume that 
the function $\theta(A)$ is not identically zero. Then $A$ does not act by isometries
when lifted to a K\"ahler cover of $M$.

\hfill

\proof
Let $\tilde M\stackrel \pi \arrow M$ be a K\"ahler cover of $M$.
Then  $\theta_1:= \pi^*\theta$ 
is exact: $\theta_1= d\phi$, and the corresponding K\"ahler form\index[terms]{form!K\"ahler}
on $\tilde M$ can be written as $\tilde\omega=e^{\phi}\pi^*\omega$.
Denote by $A_1$ the lift of $A$ to $\tilde M$. Then 
$\Lie_{A_1}(\tilde\omega) = (\Lie_{A_1} \phi) \tilde\omega$,
in other words, $A_1$ acts by isometries if and only if
$\Lie_{A_1} \phi=0$. However, the Cartan formula gives \index[terms]{Cartan formula}
\[ \Lie_{A_1} \phi =i_{A_1} d\phi=d\phi(A_1)=  \theta_1(A_1).
\]
Since $\theta(A)\neq 0$ somewhere on $M$, we have
$\theta_1(A_1)\neq 0$ somewhere on $\tilde M$.
\endproof

\hfill

Every homogeneous manifold $M$ can be obtained 
as $M=G/H$, where $G$ acts on $M$ transitively 
by automorphisms, and $H$ is the stabilizer of a point.
In our case, $M$ is LCK, and $G$ acts on $M$
by holomorphic LCK isometries. Since
the group of isometries of a compact manifold
is compact, we may freely assume that $G$ is
compact. This will be our running assumption
from now on.

\hfill

\proposition\label{_homo_admits_Vais_Proposition_}
Let $(M, \omega, \theta)$ be a compact, homogeneous LCK
manifold, $M=G/H$.
Then $M$ admits a homogeneous Vaisman metric with the same\index[terms]{metric!Vaisman}
Lee form.\index[terms]{form!Lee}\index[terms]{manifold!LCK!homogeneous}

\hfill

\pstep Let $\g=\Lie(G)$. Choose $A\in \g$
such that $\theta(A)\neq 0$ somewhere on $M$, and let
$T\subset G$ be the closure of the
one-parametric subgroup of $G$ generated
by $e^{tA}$. Since $T$ is a compact torus,
$A$ can be obtained as a limit of vector
fields $A_i$ such that the one-parametric subgroup of $G$ generated
by $e^{tA_i}$ is a circle. Then for some $i$ the vector
field $A_i$ satisfies $\theta(A_1)\neq 0$ somewhere on $M$.
By \ref{_S^1_and_theta_Lemma_}, $A_i$
it cannot act by isometries on the K\"ahler cover
of $\tilde M$. By \ref{_S^1_potential_Theorem_}, $M$ admits an LCK
metric with potential.\index[terms]{metric!LCK!with potential}

\hfill

{\bf Step 2:}
The fundamental form of an LCK  metric with potential is written as  
$\omega_\psi=d_\theta d_\theta^c(\psi)$, for a positive function $\psi$ on $M$. Averaging
on the compact group  $G$, we obtain a metric $g_1$ with the fundamental form 
\[ \omega_1=\Av_G(\omega_\psi)=d_\theta d_\theta^c(\Av_G(\psi)).
\]
Then $g_1$ is still LCK with potential, but the potential is now
$G$-invariant, and the metric is also $G$-invariant. Since $\theta$ is already $G$-invariant, the averaging process does not change it. Then $\theta$ is the Lee form\index[terms]{form!Lee} of $g_1$ too, and 
its norm  in this metric is constant:
$|\theta|_{\omega_1}=\const$. 
Then \ref{_preferred_gauge_constant_theta_Corollary_}  implies that $g_1$ is Vaisman.
\endproof

%
%
%

\hfill

\theorem\label{_LCK_homo_Vaisman_Theorem_}
Let $(M, g, I, \omega, \theta)$ be a compact, homogeneous LCK
manifold, $M=G/H$.
Then $(M, \omega, \theta)$ is Vaisman.\index[terms]{manifold!LCK!homogeneous}

\hfill

By \ref{_homo_admits_Vais_Proposition_},
$M$ admits a $G$-homogeneous Vaisman metric
$\omega_1$ with the same Lee form.\index[terms]{form!Lee} 
By \ref{_S^1_potential_Theorem_}, 
we can assume that $\omega$ is an LCK
metric with potential. Then \ref{_LCK_homo_Vaisman_Theorem_}
follows if we prove that
a $G$-homogeneous metric with potential\index[terms]{metric!LCK!with potential} is uniquely
determined, up to a constant multiplier, by its Lee form.\index[terms]{form!Lee}

\hfill

\proposition
Let $M=G/H$ be a $G$-homogeneous
complex manifold, and $\omega_1, \omega_2$ two
$G$-homogeneous LCK metrics with potential
and the same Lee form.\index[terms]{form!Lee} Then 
$\omega_1$ is proportional to $\omega_2$.

\hfill

\proof By \ref{_LCK_pot_via_d_theta_d^c_theta_Proposition_},
$\omega_i = d_\theta d^c_\theta(\phi_i)$.
Averaging $\phi_i$ with $G$ if necessary, 
we can assume that it is  $G$-invariant, 
hence constant. Then 
$\omega_i = d_\theta d^c_\theta(a_i)$,
where $a_i$ are constant functions. This says that  
these two forms are proportional.
\endproof

\section{Homogeneous Vaisman manifolds}\index[terms]{manifold!Vaisman!homogeneous}

\theorem \label{homo_reg} 
The canonical foliation\index[terms]{foliation!canonical} of a  compact  homogeneous Vaisman manifold $(M=G/H, I, g)$ is regular.\index[terms]{manifold!Vaisman!regular} 

\hfill

\proof We give a proof different from the original one in
\cite{va_torino}. Since $M$ is compact, $\Sigma$ has at
least one compact leaf (this was proven in \cite{tsu}, but
follows also from results of \index[persons]{Kato, Ma.} Kato (in \cite{kato2}) or
from \ref{main_surf_in_pot}). Then, by homogeneity, all leaves
are compact, and hence $\Sigma$  is
quasi-regular. Consider the elliptic fibration
$\pi:M\rightarrow M/\Sigma$. This map has at least one
smooth fibre. To see this, just consider $\pi$ with values
in the smooth part of $M/\Sigma$ and apply Sard's
lemma.\index[terms]{theorem!Sard} Then, again by
homogeneity, all fibres are smooth. \endproof

\hfill

\remark
On homogeneous Vaisman manifolds \index[terms]{manifold!Vaisman!homogeneous}
the foliation $\Sigma$ is regular, but the foliation
$\langle \theta^\sharp\rangle \subset \Sigma$ generated by the
Lee field \index[terms]{Lee field}is not necessarily regular.
Indeed, consider the Hopf surface $H=\frac{(\C^2\backslash 0)}{\langle A\rangle}$,\index[terms]{surface!Hopf}
where $A=\alpha \Id$ and $\alpha$ is a complex number such that
$\frac{\alpha}{|\alpha|}$ is not a root of unity. In this case, the Lee field \index[terms]{Lee field}
is a radial vector field on $\C^2$, and its trajectory 
applied to a 3-dimensional sphere gives a diffeomorphism
$v\arrow \frac{\alpha}{|\alpha|}$ that has infinite order.

\hfill

\claim \label{_compact_homo_Vaisman_b_1_Claim_}
A compact homogeneous Vaisman manifold $M$ has $b_1=1$.\index[terms]{manifold!Vaisman!homogeneous}

\hfill

\proof
Since the canonical foliation\index[terms]{foliation!canonical} $\Sigma$ is regular, $M$ is a $T^2$-fibration
over a homogeneous projective manifold $P$. By
\cite[\S2, Cor. 2 to Th. 1]{_Chevalley:group_varieties_}
(see also \cite{_Chin_Zhang:homogeneous_}), all
homogeneous projective manifolds are 
rational varieties. Using \cite[Chapter 4]{_Debarre_},\index[terms]{manifold!projective!homogeneous}
we obtain that $P$ is simply connected. Writing the exact sequence
\[
0 \arrow H^1(P) \arrow H^1(M) \arrow
H^1(T^2)\stackrel \gamma\arrow H^2(P).
\]
as in \ref{fund_group}, we observe that the rank of $\gamma$ is 1
(see the proof of \ref{fund_group}).
Therefore, $b_1(P)=0$ implies that
$b_1(M)=1$ (note that the same argument was 
used in the proof of \ref{_b_1_LCK_Toric_Theorem_}).
\endproof

\hfill

\remark One can construct homogeneous Vaisman manifolds as in
\cite{va_torino}, starting from a compact homogeneous
projective manifold $P=G/H$ such that the action of $G$ on
$P$ can be linearized (i.e.  $P$ admits an equivariant ample
line bundle). Note that the structure of compact
homogeneous K\"ahler manifolds\index[terms]{manifold!K\"ahler!homogeneous} is clarified in
\cite{mats}: up to biholomorphisms, they are products of
flag manifolds.

\hfill

In the smallest dimension, a case by case analysis, using 
the classification of Vaisman surfaces 
(\cite{bel}; see also Chapter \ref{comp_surf}), proves:

\hfill

\proposition {(\cite{HK2})} The only  homogeneous compact Vaisman surface is, up to biholomorphism, the classical Hopf surface $\C^2\setminus\{0\}/A$, with $A=\al\Id$, $\al\in\C\setminus\{0\}$.\index[terms]{surface!Hopf!classical}

\section{Notes}
Let $G$ be a compact, connected real Lie group acting by
holomorphic isometries on a connected Hermitian manifold
$(M,I,g)$. Then $M$ is said to be {\bf of cohomogeneity one} if
the orbit space $M/G$ is 1-dimensional, i.e.  homeomorphic to one of the following spaces: $\R$,  $S^1$, or $[0,1]$. Long time it was a common belief that also cohomogeneity one LCK manifolds should be Vaisman. However, \index[persons]{Angella, D.} Angella and \index[persons]{Pediconi, F.} Pediconi recently constructed (\cite{_angella_pediconi_}) an example of non-Vaisman cohomogeneity one LCK manifold on some special kind of \index[persons]{B\'erard-Bergery, L.} B\'erard-Bergery manifolds (see \cite{_bergery_}).\index[terms]{manifold!of cohomogeneity one}


\chapter{LCK structures on nilmanifolds and solvmanifolds}
\label{_nil_and_solv_}\index[terms]{structure!LCK}
\index[terms]{nilmanifold}\index[terms]{solvmanifold}

\epigraph{\it The sense of the world must lie outside the world. In the
	world everything is as it is, and everything happens as it
	does happen: in it no value exists -- and if it did exist, it
	would have no value.
	
	\hspace{.1 in}If there is any value that does have value, it must lie
	outside the whole sphere of what happens and is the
	case. For all that happens and is the case is accidental.
	
	\hspace{.1 in}What makes it non-accidental cannot lie within the
	world, since if it did it would itself be accidental.
	
	\hspace{.1 in}It must lie outside the world.}{\sc \scriptsize Ludwig Wittgenstein,\ \  Tractatus Logico-Philosophicus (6.41)}

\hfill

The subject of this chapter is locally invariant
LCK structures on nilmanifolds\index[terms]{nilmanifold} and solvmanifolds.\index[terms]{solvmanifold} The full
classification is not yet available, but we can give a
classification on nilmanifolds, reproving a result of
\index[persons]{Sawai, H.} Sawai. We rely on a vanishing result for the twisted
Dolbeault  cohomology of a nilpotent Lie algebra, as
proven in \cite{_ov_twisted_dolbeault_} (we reproduce this
proof below, see \ref{_twi_Dolbeault_Vanish_}). 

\section{Invariant geometric structures on Lie groups}
\label{_loc_inva_Section_}

Let $G$ be a Lie group, $\Lambda$ a discrete subgroup, and
$I$ a left-invariant complex structure on $G$. Consider
the quotient space $G/\Lambda$, where $\Lambda$ acts by
left translations. Since $I$ is left-invariant,
the manifold $G/\Lambda$ is equipped with a natural
complex structure. This construction is usually 
applied to solvmanifolds\index[terms]{solvmanifold} or nilmanifolds,\index[terms]{nilmanifold} but
it makes sense with any Lie group. When we need to refer to this
particular kind of complex structures on $G/\Lambda$, we call them 
{\bf locally $G$-invariant}.

\hfill

A left-invariant almost complex structure $I$ on $G$ is
determined by its restriction to the Lie algebra
$\g=T_e G$, giving the Hodge decomposition
$\g\otimes_\R \C= \g^{1,0}\oplus \g^{0,1}$. 
This structure is integrable
if and only if the commutator of $(1,0)$-vector fields
is again a $(1,0)$-vector field, that is  equivalent
to $[\g^{1,0},\g^{1,0}]\subset \g^{1,0}$.
This allows one to define a complex structure
on a Lie algebra.

\hfill

\definition\label{_complex_on_nil_alg_Definition_}
A complex structure on a Lie algebra $\g$ is a subalgebra
$\g^{1,0}\subset \g\otimes_\R \C$
that satisfies $\g^{1,0}\oplus
\overline{\g^{1,0}}=\g\otimes_\R \C$.

\hfill

Indeed, an almost complex structure operator
$I$ can be reconstructed from the decomposition
 $\g^{1,0}\oplus\overline{\g^{1,0}}=\g\otimes_\R \C$ by making it act
as $\1$ on $\g^{1,0}$ and $-\1$ on $\overline{\g^{1,0}}$.

\hfill

Similarly, one could define 
symplectic structures or LCK structures\index[terms]{structure!LCK} on a Lie algebra $\g$.
Recall that the Grassmann algebra\index[terms]{Grassmann algebra} $\Lambda^*(\g^*)$ is
equipped with a natural differential, called 
{\bf \index[persons]{Chevalley, C.} Chevalley differential}, that is  equal \index[terms]{Chevalley differential}
to the de Rham differential on left-invariant
differential forms on the Lie group if we identify
the space of such forms with $\Lambda^*(\g^*)$. 
In degree 1  it can be written explicitly as follows:
 if $\la\in\goth{g^*}$, then $d\la(x,y)=-\la([x,y])$.
This differential is extended to $k$-forms by
the Leibniz rule, and the condition $d^2=0$
is equivalent to the Jacobi identity.
On 2-forms, the differential is written as 
\[ d\be(x,y,z)=-\be([x,y
 ], z) -\be([y,z], x) + \be([z,x], y ),
\] for $\be\in\Lambda^2\goth{g}^*$.
The corresponding complex is called  
{\bf the \index[persons]{Chevalley, C.} Chevalley--Eilenberg complex}.\index[terms]{complex!Chevalley--Eilenberg}

Further on, we shall always interpret the elements
of $\Lambda^p(\g^*)$ as left-invariant differential
forms on the corresponding Lie group, and refer
to them as to ``$p$-forms'', with all the usual
terminology (``closed forms'', ``exact forms'')
as used for the elements of de Rham algebra.

\hfill

\definition
{\bf A symplectic structure}\index[terms]{structure!symplectic} on a Lie algebra $\g$ is a
non-degenerate, closed 2-form $\omega\in
\Lambda^2(\g^*)$. A {\bf K\"ahler structure}\index[terms]{structure!K\"ahler}
on a Lie algebra $\g$ is a complex structure $I$
and a Hermitian form\index[terms]{form!Hermitian} $h$ on $\g$ such that 
the fundamental 2-form 
$\omega(\cdot, \cdot):= h( I(\cdot),\cdot)$
is closed. 

\hfill

\remark
In \cite{bens_gor} 
(see also \cite{hase}) it was shown that any
nilpotent Lie algebra admitting a K\"ahler structure is
actually abelian.\index[terms]{Lie algebra!nilpotent!}

\hfill

\definition\label{_LCK_Lie_algebra_Definition_}
{\bf An LCK structure}\index[terms]{structure!LCK} on a Lie algebra is a
complex structure $I$
and a Hermitian form\index[terms]{form!Hermitian} $h$ on $\g$ such that 
the fundamental 2-form 
$\omega(\cdot, \cdot):= h(I(\cdot),\cdot)$
satisfies $d(\omega) = \theta \wedge \omega$, where
$\theta\in \Lambda^1(\g^*)$ is a closed 1-form.

\hfill

\definition
A {\bf nilmanifold}\index[terms]{nilmanifold} is a 
quotient of a simply connected nilpotent\index[terms]{Lie group!nilpotent}
Lie group by a cocompact  discrete
subgroup. Alternatively (\cite{mal}), one can define
nilmanifolds as compact manifolds which
admit a homogeneous action by a nilpotent
Lie group. A {\bf solvmanifold} \index[terms]{solvmanifold} 
is  a  quotient of a solvable
Lie group by a cocompact  discrete
subgroup. 

\hfill

\definition\label{_LCK_nilma_Definition_}
Let $(G, I)$ be a nilpotent (resp.  solvable) Lie group\index[terms]{Lie group!solvable}\index[terms]{Lie group!nilpotent}
equipped with a left-invariant complex structure.
For any cocompact  discrete
subgroup $\Gamma\subset G$, the (left) quotient
$G/\Gamma$ is equipped with a natural complex structure.
This quotient is called  {\bf a complex nilmanifold\index[terms]{nilmanifold}
  (resp. solvmanifold)}.\index[terms]{solvmanifold} Similarly, if $G$ is a Lie group
with a left-invariant LCK structure\index[terms]{structure!LCK}, the quotient
$G/\Gamma$ is called {\bf an LCK nilmanifold
  (resp. solvmanifold)}. This structure is clearly
locally homogeneous. \index[terms]{manifold!locally homogeneous}When we need to refer to this
particular kind of structures, we call them 
{\bf locally $G$-invariant}.

\hfill

\remark The existence of a cocompact  lattice\index[terms]{lattice!cocompact} already
imposes strong restrictions: the group should be
unimodular\index[terms]{Lie group!unimodular},
\cite[\S 1.1, Exercise 14b]{_Morris:Ratner_}. Nilpotent Lie groups are unimodular, but
not all solvable ones are. Moreover, unimodularity is not sufficient for the existence of a cocompact lattice. According to \cite{mal}, a nilpotent Lie group admits a cocompact lattice if and only if its Lie algebra admits a basis in which the structural constants are all rational.

\hfill

In this chapter
we present results concerning LCK structures\index[terms]{structure!LCK} on
nilmanifolds \index[terms]{nilmanifold}and solvmanifolds\index[terms]{solvmanifold}. A classification of the
solvable Lie algebras that admit LCK structures is still\index[terms]{Lie algebra!solvable}
missing in full generality, but nilmanifolds and certain particular cases
of solvmanifolds\index[terms]{solvmanifold} are covered. 

\hfill 

\remark  \label{_Malcev_comple_Remark_}
A torsion-free group $G$ is called {\bf \index[persons]{Maltsev, A.} Maltsev complete}\index[terms]{group!torsion-free}
if for any $n \in \Z^{>0}$, and any $g\in G$, there exists
a $g_n\in G$ satisfying $g_n^n=g$. In \cite{_Maltsev:nilpotent_}, \index[persons]{Maltsev, A.} Maltsev 
developed the theory of complete groups and used it to 
study the nilmanifolds\index[terms]{nilmanifold}. He has shown that
any torsion-free nilpotent group $G$ can be embedded to a
complete group $\hat G$ such that all elements $g\in \hat
G$ satisfy $g^n\in G$ for some $n \in \Z^{>0}$.
He proved that the group $\hat G$ is defined by $G$
uniquely and functorially; it is called {\bf the
Maltsev completion of $G$}. The \index[persons]{Maltsev, A.} Maltsev completion
\index[terms]{completion!Maltsev}
is always a rational algebraic nilpotent group; the corresponding
real algebraic group $\hat G_\R$ is a nilpotent Lie group that contains
$G$ as a lattice.\index[terms]{lattice} Using this construction, we obtain 
that any nilmanifold is uniquely determined by its fundamental
group, that is  a discrete nilpotent torsion-free group,
\index[terms]{group!torsion-free} and any such group
uniquely determines a nilmanifold.\index[terms]{nilmanifold}

\hfill

\remark Note that right translations are not holomorphic
with respect to left-invariant complex structures. This
means that  for a lattice $\Lambda$, the complex structure
induced on the manifold $G/\Lambda$ is not necessarily
invariant, yielding that $G/\Lambda$ is not a homogeneous LCK
manifold.\index[terms]{manifold!LCK!homogeneous}  This is why we
prefer use ``locally $G$-invariant'' referring to such structures.

\hfill

\remark It might happen that a K\"ahler manifold is
diffeomorphic to a complex manifold of non-K\"ahler type. For example,
any homogeneous complex structure on a 6-dimensional
torus $T^6$ is clearly of K\"ahler type. However, there exists
a non-K\"ahler complex structure on $T^6$. For example,
let $E$ be  an elliptic curve that is   a 2-fold ramified
covering of  $\C P^1$, and take the fibred product
$E\times_{\C P^1} \Tw(T^4)$,  where $\Tw(T^4)$ is the 
twistor space of a 4-torus. This fibred product is
actually a 6-torus, but with an inhomogeneous complex
structure.
Indeed, it has non-trivial canonical bundle, in fact its anticanonical bundle
is semi-positive, with many non-trivial sections (compare
with \cite[\S 5]{_Calabi:hk_}). This example suggests that
the same phenomenon could occur for LCK nilmanifolds.\index[terms]{nilmanifold} It is
unknown at the moment if there exists an LCK structure\index[terms]{structure!LCK} $(I,\omega, \theta)$
on a nilmanifold $M = G/\Gamma$ that is  not locally $G$-invariant,
for an appropriate $G$-action.

\section{Twisted Dolbeault cohomology on nilpotent Lie algebras}\index[terms]{Lie algebra!nilpotent}
\label{_twisted_Dolbeault_nilma_Section_}\index[terms]{cohomology!Dolbeault!twisted}

Let $\theta$ be a closed 1-form on a complex manifold.\index[terms]{cohomology!Dolbeault!twisted}
Then $\nabla -\theta$ defines a connection on the trivial bundle.
We denote the corresponding local system by $L$.
The {\bf twisted Dolbeault differentials} $\6_\theta:= (d_\theta)^{1,0}=\6-\theta^{1,0}$ and
$\bar\6_\theta:= (d_\theta)^{0,1}=\bar\6-\theta^{0,1}$
are Morse--Novikov counterparts to the usual Dolbeault differentials.
The cohomology of these differentials corresponds to the 
Dolbeault cohomology with coefficients in the 
holomorphic line bundle $L_\C$ obtained as complexification 
of the local system $L$.

Working with LCK nilmanifolds and nilpotent Lie algebras, \index[terms]{nilmanifold}
it is natural to consider their twisted Dolbeault cohomology. 
As shown by Console and \index[persons]{Fino, A.} Fino (\cite{_Console_Fino_}), the usual 
Dolbeault cohomology\index[terms]{cohomology!Dolbeault!twisted} of complex nilmanifolds\index[terms]{nilmanifold} is equal to the \index[terms]{cohomology!Dolbeault}
cohomology of the Dolbeault version of the corresponding 
\index[persons]{Chevalley, C.} Chevalley-Eilenberg complex on the Lie algebra when the complex structure is rational (in dimension up to 6 the rationality of the complex structure is not needed, see \cite{_Fino_Rollenske_Ruppenthal_}, but for higher dimensions the question is still open). This result
does not hold in the twisted Dolbeault cohomology.


Let $E_1, E_2$ be elliptic curves.
Recall that a  Kodaira surface $M$ is an elliptic Vaisman \index[terms]{surface!Kodaira}
 surface, obtained as the total space of a principal holomorphic
$E_1$-bundle  $\pi:\; M \arrow E_2$\index[terms]{bundle!principal} 
(Exercise \ref{_Kodaira_surface_Exercise_}).
Such bundles are classified by the first  Chern
class $c_1(\pi)$ of the fibration  (\ref{_c_1_toric_Definition_}). 
This class can be identified with the
$d_2$-differential of the corresponding \index[persons]{Leray, J.} Leray 
spectral sequence, mapping from $H^1(E_1)$ to $H^2(E_2)$ (\cite[\S 23.5]{_Fuks_Fomenko_}).
In this context, the restriction of $d_2$ to $H^1(\text{fiber})$ is usually
called {\bf the transgression map}.
\index[terms]{spectral sequence!Leray}
By \index[persons]{Blanchard, A.} Blanchard's  \ref{_Blanchard_Theorem_}, $M$ is non-K\"ahler if and \index[terms]{theorem!Blanchard}
only if $c_1(\pi)$ is non-zero. The \index[terms]{surface!Kodaira} Kodaira surface is an example
of a nilmanifold \index[terms]{nilmanifold}obtained from the Heisenberg group 
(see \ref{_Kodaira_as_nilmanifold_Example_} below).\index[terms]{group!Heisenberg}

\hfill

\example  {\bf (Console-Fino theorem with coefficients in a local system fails)} \label{_console_fino_ex_}\\
We are going to compute the weighted Dolbeault cohomology of the Kodaira surface explicitly.\index[terms]{cohomology!Dolbeault!weighted}
Let $\omega_{E_2}$ be a K\"ahler form \index[terms]{form!K\"ahler}on $E_2$. 
Then $\pi^*(\omega_{E_2})\in \im d_2$ is exact, giving $\pi^*(\omega_{E_2})= d\xi$.
We choose $\xi$ in such a way that $d(\theta)=0$, where $\theta:=I(\xi)$
(\ref{_H^1_odd_Theorem_} or \ref{_Hodge_chasing_Corollary_}). 
Consider a trivial complex line bundle $(L, \nabla_0)$ with the connection
defined by $\nabla_0 - \theta$; this line bundle is flat, and hence  equipped with 
a holomorphic structure operator $\bar\6:=\nabla_0^{0,1}- \theta^{0,1}$.
Using the standard Hermitian metric, we express its
Chern connection \index[terms]{connection!Chern}as $\nabla_0 +\theta^{1,0} - \theta^{0, 1}$,
thus the curvature of this bundle satisfies $\Theta_L= -\1 \pi^*(\omega_{E_2})$.
Rescaling $\theta$ and $\omega_{E_2}$ if necessarily, we can assume
that the cohomology class of $\omega_{E_2}$ is integer.
Then $L$ is the pullback of an ample bundle $B$ on $E_2$.\index[terms]{bundle!line!ample} 
Now, $H^0(M,L)= H^0(E_2, B)\neq 0$,
but this space can be interpreted as $H^0(\Lambda^{0,*}(M), \bar\6_\theta)$.
In this example the twisted Dolbeault cohomology of $M$ is  non-zero,
but the corresponding twisted Dolbeault cohomology\index[terms]{cohomology!Dolbeault!twisted} of the Lie algebra vanishes
(\ref{_twi_Dolbeault_Vanish_}).
This gives a counterexample to Console-Fino theorem for the twisted cohomology. 

\hfill

\definition
Let $\g$ be a Lie algebra and 
$\theta\in \g^*$ a closed 1-form and \[d_\theta:=d-\theta\in \End(\Lambda^*(\g^*)),\]
where $d$ is the \index[persons]{Chevalley, C.} Chevalley differential and $\theta$
denotes the operation of multiplication by $\theta$.
The twisted Dolbeault differentials on
$\Lambda^*(\g^*)$ are $\6_\theta:= (d_\theta)^{1,0}$ and
$\bar\6_\theta:= (d_\theta)^{0,1}$.

\hfill

The cohomology of $d_\theta$ is always zero 
for non-zero $\theta$ (\cite{_Alaniya_,mili}). This classical theorem
is due to J. \index[persons]{Dixmier, J.} Dixmier, \cite{_Dixmier_}. 
The main result of the present section is the following theorem,
that is  the Dolbeault version of \index[persons]{Dixmier, J.} Dixmier's and \index[persons]{Alaniya, L. A.} Alaniya's theorem.

\hfill

\theorem\label{_twi_Dolbeault_Vanish_}
{(\cite{_ov_twisted_dolbeault_})} 
Let $\g$ be a nilpotent Lie algebra with a complex structure 
(\ref{_complex_on_nil_alg_Definition_}), and 
$\theta\in \g^*$ a non-zero, closed real 1-form.
Then the cohomology of $(\Lambda^{0,*}(\g^*),\bar\6_\theta)$ vanishes.

\hfill

\proof
Denote the  central series of the Lie algebra
$\g^{0,1}$ by
\[W_0=  \g^{0,1},\ W_1=  [W_0, W_0], ..., 
W_k=[W_0, W_{k-1}].
\]
Let $A_k\subset (\g^{0,1})^*$
be the annihilator of $W_k$. For any 1-form $\lambda\in(\g^{0,1})^*$, one has 
$\bar\6(\lambda)(x,y) = \lambda([x,y])$, and then   
$\bar\6(A_k)\subset \Lambda^2(A_{k-1})$. Consider the filtration
$\Lambda^*(A_1)\subset \Lambda^*(A_2)\subset ...$ on $\Lambda\in(\g^{0,1})^*$.
Since $\bar\6(A_k)\subset \Lambda^2(A_{k-1})$, the operator $\bar\6$
shifts the filtration by 1: $\bar\6(\Lambda^*(A_k))\subset \Lambda^*(A_{k-1})$.

Consider the spectral sequence of the complex $(\Lambda^*(\g^{0,1})^*,\bar\6_\theta)$
filtered by $V_0\subset V_1\subset V_2 \subset ...$, where\index[terms]{spectral sequence}
$V_k := \Lambda^*(A_k)$ (Subsection \ref{_filtered_spectral_sequence_Subsection_}).
Since  $\bar\6(V_k) \subset V_{k-1}$, the operator $\bar\6_\theta$
acts on $\bigoplus_k V_k/V_{k-1}$ as multiplication by $\theta^{0,1}$. 
The corresponding associated graded complex, that is  $E_0^{*,*}$ of this
spectral sequence, is identified with
$\left(\bigoplus_k V_k/V_{k-1}, \theta^{0,1}\right)$. 
We identify $\bigoplus_k V_k/V_{k-1}$ with the
Grassmann algebra of $(\g^{01})^*$: \[ \bigoplus_k V_k/V_{k-1}=\Lambda^*(\g^{01})^*.\]
After this identification, the multiplication
\[ \theta^{0,1}:\;\bigoplus_k V_k/V_{k-1}
\arrow \bigoplus_k V_k/V_{k-1}
\]
becomes the Grassmann multiplication 
by a vector of $(\g^{01})^*$ associated with $\theta^{0,1}$.
The cohomology of multiplication by a non-zero vector always vanishes. Then the
$E_1^{*,*}$-page of the spectral sequence vanishes, which
implies vanishing of the cohomology $H^*(\Lambda^*(\g^{0,1}),\bar\6_\theta)$.
\endproof

\hfill

The following corollary is used in the classification of LCK
structures on nilpotent Lie algebras. 

\hfill

\corollary\label{_Hodge_chasing_Corollary_}
Let $(\goth{g},  I)$ be a $2n$-dimensional
nilpotent Lie algebra with complex structure, $\theta\in \Lambda^1(\g^*)$
a non-zero, closed real 1-form, and  $\omega\in \Lambda^{1,1}(\g^*)$
a $d_\theta$-closed (1,1)-form.
Then there exists a 1-form $\tau\in \Lambda^1(\g^*)$
such that  $\omega=d_\theta(\tau)$,
where $d_\theta (I\tau)=0$.

\hfill

\proof
Denote by $H^{1,1}_{d_\theta d^c_\theta}(\g^*)$ the twisted
Bott--Chern cohomology, that is, all $d_\theta$-closed (1,1)-forms\index[terms]{cohomology!Bott--Chern}
up to the image of $d_\theta d^c_\theta$.
The standard exact sequence \eqref{_exact_sequence_BC_MN_Equation_} 
\[ H^{0,1}_{\bar\6_\theta}(\g^*) \oplus H^{1,0}_{\6_\theta}(\g^*)
\stackrel d\arrow H^{1,1}_{d_\theta d^c_\theta}(\g^*) \stackrel \mu \arrow H^2_{d_\theta}(M)
\]
implies that
the kernel of $\mu$ vanishes when  
$H^{0,1}_{\bar\6_\theta}(\g^*)=H^{1,0}_{\6_\theta}(\g^*)=0$.
The last equation follows from  \ref{_twi_Dolbeault_Vanish_},
hence $\mu$ is injective. The LCK form $\omega$ belongs to the
kernel of $\mu$, because the $d_\theta$-cohomology of
$\Lambda^*(\g^*)$ vanishes by \index[persons]{Dixmier, J.} Dixmier and Alanya theorem 
(\ref{_twi_Dolbeault_Vanish_}).\index[terms]{theorem!Dixmier, Alanya}
Therefore, $\omega$ is twisted Bott--Chern exact, i.e.   
$\omega= d_\theta d_\theta^c f$ for some
$f\in \Lambda^0(\g^*)=\R$. Set $\tau:=d^c_\theta f$. Then $I\tau=d_\theta f$ and $d_\theta(I\tau)=0$ as stated. \endproof

\hfill

\remark
In \cite{_Dixmier_}, \index[persons]{Dixmier, J.} Dixmier proved that Morse--Novikov
cohomology of a nilmanifold\index[terms]{nilmanifold} always vanish. In \cite{_Chen:MN_},
X. Chen generalized this result by proving
that the Morse--Novikov cohomology (that is,
the cohomology with coefficient in a local system $L$)
of any manifold $M$ with nilpotent fundamental group\index[terms]{fundamental group}
vanishes, provided that $L$ is non-trivial,
and the action of $\pi_1(M)$ on the higher
homotopy groups of $M$ is nilpotent.

\section{LCK nilmanifolds}\label{_LCK_nilmanifolds_section_}\index[terms]{nilmanifold}

\example\label{vai_hei}\label{kod_nil}  
Let $H_{2n-1}$ be {\bf the Heisenberg group}, that is, the 
group of matrices of the form\index[terms]{group!Heisenberg}\index[terms]{Lie algebra!Heisenberg}\index[persons]{Heisenberg, W.}
\[ \small
\begin{pmatrix}
1&A&c\\
0&I_{2n-1}&B^t\\
0&0&1
\end{pmatrix}, \quad c\in\R,\quad A,B\in\R^{n-1}.
\]
Its Lie algebra $\goth{h}_{2n-1}$
can be defined in terms of generators and relations as follows:
it has a basis $\{X_i,Y_i,Z\}$, $i=1,\ldots,n-1$, 
with $[X_i,Y_i]= Z$, and the rest of the Lie brackets trivial.
Clearly, it is a nilpotent Lie algebra. After taking the
quotient of
the corresponding Lie group
by a cocompact discrete subgroup $\Lambda$, one obtains 
a nilmanifold\index[terms]{nilmanifold} called {\bf the Heisenberg nilmanifold}. It
has an invariant Sasakian structure given by requiring the
above basis to be orthonormal, and defining the contact
form as $\sum y_idx_i+dz$ (where $dx_i,dy_i, dz$ are the dual
invariant one-forms). Then $Z$ is its Reeb field.\index[terms]{manifold!Sasaki}

Now, the product $H_{2n-1}/\Lambda\times S^1$ is a Vaisman
nilmanifold whose universal cover is the product $
H_{2n-1}\times\R$, a Lie group with Lie algebra
$\goth{g}:= \goth{h}_{2n-1}\times\R$, $\R=\langle
T\rangle$,  in which the brackets of the basis
$\{X_i,Y_i,Z, T\}$ are the above, and $T$ is central. 
The linear LCK structure \index[terms]{structure!LCK}on $\goth{g}$ is given by
asking the basis to be orthonormal, and by defining  the
complex structure $IX_i=Y_i$, $IZ=-T$. The fundamental
form is $\omega=\sum (X_i^*\wedge Y_i^*)-Z^*\wedge T^*.$ A direct
computation shows that:
$d\omega=T^*\wedge \omega$; consequently, the  Lee form\index[terms]{form!Lee!parallel} is $T^*$
that can be seen directly to be $g$-parallel. 

One can verify that the above complex structure on
$\goth{g}$ satisfies the equation $[IX,IY]=[X,Y]$ for
all $X,Y\in\goth{g}$. Such complex structures on Lie
groups are called {\bf abelian}. \index[terms]{complex structure!abelian}

\hfill

\remark An equivalent definition
of abelian complex structures is given in
\cite{_BD:Abelian_}.  Let $I$ be a complex structure
on a Lie algebra $\g$, and $\g\otimes_\R \C=
\g^{1,0}\oplus \g^{0,1}$ its Hodge decomposition.
Then $I$ is abelian if and only if the Lie algebra
$\g^{1,0}$ is commutative.

\hfill

\example \label{_Kodaira_as_nilmanifold_Example_}
As in 
Exercise \ref{_Kodaira_de_Rham_Exercise_}, we may consider
the \index[terms]{surface!Kodaira} Kodaira surface as a nilmanifold\index[terms]{nilmanifold} associated with the Heisenberg group
as above. The corresponding Heisenberg group\index[terms]{group!Heisenberg} is given by 
\[ \small
\begin{pmatrix}
1&a&c\\
0&1&b\\
0&0&1
\end{pmatrix}, \quad a,b,c\in\R.
\]
The corresponding nilmanifold\index[terms]{nilmanifold} $N$ is a 
$S^1$-fibration over $T^2= S^1 \times S^1$,
with non-zero  Chern class.\index[terms]{class!Chern} Then $N \times S^1$
is fibred over an elliptic curve with the
fibre again identified with an elliptic 
curve, and a non-trivial Chern class.
This is a Kodaira surface,\index[terms]{surface!Kodaira} as explained in 
Section \ref{_twisted_Dolbeault_nilma_Section_}.

\hfill

\remark It is proven in \cite{aor1} that a compact LCK
solvmanifold with a\-be\-li\-an complex structure is in fact
isomorphic with the above Vaisman nilmanifold\index[terms]{nilmanifold} $S^1\times
H_{2n-1}/\Lambda$. This is a special case of a
conjecture stated by  L. \index[persons]{Ugarte, L.} Ugarte: \index[terms]{solvmanifold}

\hfill

\conjecture {(\cite{uga})}\label{conj_uga} 
Let $M$ be a differentiable nilmanifold\index[terms]{nilmanifold} admitting an LCK structure\index[terms]{conjecture!Ugarte}\index[terms]{structure!LCK}
(not necessarily locally $G$-invariant, see \ref{_LCK_nilma_Definition_}).
Then $M$ is conformally biholomorphic to a quotient $S^1\times H_{2n-1}/\Lambda$ 
with the Vaisman structure described above.

For LCK nilmanifolds, this conjecture was proven
by H. \index[persons]{Sawai, H.} Sawai\index[persons]{Sawai, H.}:

\hfill

\theorem {(\cite{saw1})}\label{inv_nil_vai} 
Let
$(M,I)$  be a complex nilmanifold\index[terms]{nilmanifold}, $M=G/\Lambda$. 
If $(M, I)$ admits an LCK structure, then it is 
conformally equivalent to Vaisman. 
Moreover, it is {\bf of Heisenberg type}, that is,
biholomorphic to a nilmanifold obtained\index[terms]{manifold!Vaisman!Heisenberg type}
as the quotient of the product $(H_{2n-1}\times\R, I)$.

\hfill

We give a new proof of this result, different from 
the one in \cite{saw1}.

\hfill

\proof (\cite{_ov_twisted_dolbeault_}) {\bf Step 1:} 
We  replace the LCK metric by a locally $G$-invariant
LCK metric. To this
purpose, we apply the averaging trick that\index[terms]{averaging trick}
Belgun\index[persons]{Belgun, F. A.} used for
the \index[terms]{surface!Inoue} Inoue surface $S^+$ (\cite[Proof of Theorem 7]{bel});
later, it was generalized
by \index[persons]{Fino, A.} Fino and \index[persons]{Grantcharov, G.} Grantcharov (\cite{_Fino_Gra_}). This approach
does not work, in general, for solvmanifolds\index[terms]{solvmanifold}, but it works 
well for all nilmanifolds. For LCK structures\index[terms]{structure!LCK} on nilmanifolds\index[terms]{nilmanifold}, this construction
is due to L. \index[persons]{Ugarte, L.} Ugarte:

\hfill

\theorem \label{bel_av} 
{(\cite[Proposition 34]{uga})}\\
Let $M=G/\Lambda$ be a complex nilmanifold 
admitting an LCK structure $(\omega,\theta)$. Then it
also admits a structure of LCK nilmanifold. In other
words, there exists a left-invariant LCK structure
on $G$ that induces an LCK structure\index[terms]{structure!LCK} on $M$.

\hfill

\proof As $G$ admits a cocompact  lattice,\index[terms]{lattice!cocompact} it is
unimodular; it thus admits a bi-invariant measure. Let
$d\mu$ be a bi-invariant volume element on $M$, and
suppose $\vol(M)=1$. Consider a left-invariant 1-form 
form $\check \theta$ on $\g =\Lie(G)$ such that the corresponding
1-form $\theta_0$ on $G/\Lambda$ is cohomologous to $\theta$.
Such a 1-form exists because 
$H^1(G/\Lambda, \R)= H^1(\Lambda^*(\g^*), d)$
(\cite{nomizu}).
Replacing $\omega$ by a conformally equivalent
LCK form $\omega_0$, we can assume that
$d\omega_0=\omega_0\wedge \theta_0$.
Let $D\subset G$ be the fundamental domain of the action
of $\Lambda$.
Given left-invariant vector fields $X, Y$ on $G$, the 2-form
$$\check \omega(X,Y):=\int_D \omega_0(X,Y)d\mu$$
defines an Hermitian structure on the Lie algebra $\g$ of $G$. Moreover, as 
\begin{multline*}
d(\check \omega)(X, Y, Z)\\=  
-\int_D \omega_0([X,Y],Z)d\mu - \int_D
\omega_0(X,[Y,Z])d\mu + \int_D \omega_0(Y,[X,Z])d\mu \\= - \int_D
d(\omega_0)(X,Y,Z) d\mu  = - \int_D (\theta_0\wedge\omega_0)
(X,Y,Z) d\mu= -(\check\theta\wedge\check \omega)(X,Y,Z),
\end{multline*}
it follows that $\check \omega$ is an LCK form on $\g$. \endproof

\hfill

{\bf Proof of \ref{inv_nil_vai}, Step 2:}
Now we prove that any LCK nilmanifold\index[terms]{nilmanifold} is Vaisman.
Using \ref{_Hodge_chasing_Corollary_}, we obtain that $\omega=d_\theta(\tau)$,
where  $d_\theta(I\tau)=0$.

Applying \index[persons]{Dixmier, J.} Dixmier and \index[persons]{Alaniya, L. A.} Alaniya's theorem again,
we obtain that $d_\theta(I\tau)=0$ implies $I\tau=d_\theta(v)$,
where $v\in \Lambda^0(\g^*)$ is a constant. Therefore,
$\omega= d_\theta d^c_\theta(\const)$, that is  the 
equation for the LCK manifold with potential\index[terms]{manifold!LCK!with potential}.\footnote%
{Note that this constant
should be positive by \ref{_strictly_negative_pots_Theorem_}.}
However, as shown in \ref{_preferred_gauge_constant_theta_Corollary_}, an LCK manifold with
potential and constant $|\theta|$ is Vaisman.

\hfill

{\bf   Step 3:} 
We show that the Lee field\index[terms]{Lee field}
$\theta^\sharp\in \g$ and the anti-Lee field\index[terms]{Lee field!anti-}
$I(\theta^\sharp)\in \g$ generate an ideal in $\g$. 
This observation also follows from 
\cite[Theorem 3.12]{_Fino_Gra_Ve_}.

Let $\Sigma=\langle \theta^\sharp,I(\theta^\sharp)\rangle$
be the canonical foliation\index[terms]{foliation!canonical} on the Vaisman manifold $M$,
and ${\goth s}\subset \g$ the corresponding subspace of
the Lie algebra of $G$. Since $\theta^\sharp$ and $I(\theta^\sharp)$\index[terms]{manifold!Vaisman}
are Killing\index[terms]{vector field!Killing} and have constant length, their trajectories
are geodesics (\ref{_canon_foli_totally_geodesic_Remark_}).
Let $S\subset \Iso(M)$ be the group of isometries of $M$
generated by exponentials of $s\in \goth s$. 
By construction, $S$ is a subgroup of $G$ that acts on
$M= G/\Lambda$. Therefore, it is invariant under conjugation with
$\Lambda$, and the Lie algebra $\goth s \subset \g$ of $S$ is invariant under the
adjoint action of $\Lambda$ on $\g$. 
Since $G$ is the Maltsev completion of $\Lambda$, its image in $\End(\g)$ \index[terms]{Maltsev completion}
coincides with the Zariski closure of the image of $G$ under the\index[terms]{Zariski closure}
adjoint action (see Property 1.5 for the \index[persons]{Maltsev, A.} Maltsev 
completion functor in \cite{_Grunewald_O'Halloran_}). Therefore, $\goth s$ is $G$-invariant.



\hfill

{\bf   Step 4:} 
Now we prove that for a Vaisman nilmanifold\index[terms]{nilmanifold} $M=G/\Lambda$,
the group $G$ is a product of $\R$ and the Heisenberg
group. Let $\goth s\subset \g$ be the ideal associated with the
canonical foliation\index[terms]{foliation!canonical} $\Sigma$ (Step 3). The leaf space of $\Sigma$
is K\"ahler because $M= G/\Lambda$ is Vaisman.
However, this leaf space is a nilmanifold
associated with the Lie algebra $\g/\goth s$.
A K\"ahler nilpotent Lie algebra is abelian
(\cite{bens_gor}), and hence  $\g$ is a central
extension of an abelian algebra. Since
$M$ is Vaisman, it is locally isometric to a product
of $\R$ with a Sasakian manifold;\index[terms]{manifold!Sasaki} 
the corresponding Lie algebra is then the product
of $\R$ and a 1-dimensional central extension\index[terms]{extension!central}
of $\R^{2n-2}$. This ends the proof.
\endproof

\hfill

The general case of \ref{conj_uga} was proven only for
nilmanifolds\index[terms]{nilmanifold} of Vaisman type, by
G. Bazzoni\index[persons]{Bazzoni, G.}:

\hfill

\theorem {(\cite{baz})} Let $M=G/\Lambda$ be a smooth
$2n$-dimensional nilmanifold admitting a Vaisman structure
(not necessarily locally $G$-invariant). 
Then  $G$ is isomorphic to $\R\times H_{2n-1}$.

\hfill

\proof The idea of the proof is to use \ref{str_vai},
consider the associated Sasakian manifold, and  prove
that $G$ is the Heisenberg group\index[terms]{group!Heisenberg} by finding its minimal
model.\index[terms]{minimal model}

 Let then $S$ be the compact Sasakian manifold associated
 to $M$ by the mapping torus \index[terms]{mapping torus}in the Structure
 Theorem. Then $S$ is
 aspherical\index[terms]{manifold!aspherical} (all $\pi_k(S)=0$
 for $k>1$) and its fundamental group\index[terms]{fundamental group} $\Gamma$ fits into the short exact sequence
$$0\ra \Gamma\ra\Lambda\ra\Z\ra 0.$$
Then $\Gamma$ is  nilpotent,\index[terms]{group!nilpotent} torsion-free and finitely
generated.\index[terms]{group!finitely generated}\index[terms]{group!torsion-free} Therefore, there exists  a connected and simply
connected nilpotent group $H$ containing $\Gamma$ as a
lattice\index[terms]{lattice!cocompact} and such that $P:=H/\Gamma$ is a compact
nilmanifold (\ref{_Malcev_comple_Remark_}), that is  also aspherical and
satisfies $\pi_1(P)=\pi_1(S)$; then $P$ and $S$ have the same
rational homotopy type, and thus the same minimal
model. We thus have a compact nilmanifold\index[terms]{nilmanifold} and a compact
Sasakian manifold with the same minimal model.

It is proven in \cite{cape1} that a $2n-1$-dimensional compact nilmanifold with Sa\-sa\-kian structure is diffeomorphic to the quotient of $H_{2n-1}$ by a lattice.\index[terms]{lattice!cocompact} By a clever analysis of that proof, \index[persons]{Bazzoni, G.} Bazzoni concludes that $H=H_{2n-1}$ and then $P=H_{2n-1}/\Gamma$.

Being a mapping torus over the circle, $M$  is $m$-covered (for some $m>0$) by the product $S\times S^1$ and the following commutative diagram exists, \cite{bmo}:
$$
\xymatrix{
S\ar[r]\ar[d]^{\id}&S\times S^1\ar[d]^{\pi}\ar[r]&S^1\ar[d]^m\\
S\ar[r]&M\ar[r]&S^1
}
$$	


This induces the following commutative diagrams of fundamental groups\index[terms]{fundamental group}:

$$
\xymatrix{
1\ar[r]&\Gamma\ar[r]\ar[d]^{\id}&\Gamma\times\Z\ar[r]\ar[d]^{\pi_*}&\Z\ar[r]\ar[d]^m&1\\
1\ar[r]&\Gamma\ar[r]&\Lambda\ar[r]&\Z\ar[r]&1	
}
$$


By \cite[Theorem II.2.11]{ragu} (see
also \ref{_Malcev_comple_Remark_}), there exists a
commutative diagram involving the simply connected Lie
groups with cocompact, discrete subgroups  entering  the above
diagram:

$$
\xymatrix{
	1\ar[r]&H_{2n-1}\ar[r]\ar[d]^{\id}&H_{2n-1}\times\R\ar[r]\ar[d]^{\psi}&\R\ar[r]\ar[d]^m&1\\
	1\ar[r]&H_{2n-1}\ar[r]&G\ar[r]&\R\ar[r]&1	
}
$$


Here the maps are continuous (hence smooth) homomorphisms. By the simply connectedness, $\psi$ is a group isomorphism and the proof is complete. \endproof	
	
\section{LCK solvmanifolds}\index[terms]{solvmanifold}

For complex solvmanifolds equipped with an LCK structure\index[terms]{structure!LCK} one\index[terms]{solvmanifold}
can sometimes apply the \index[persons]{Belgun, F. A.} Belgun-Fino-Grantcharov averaging trick 
(see the proof of \ref{bel_av}).\index[terms]{averaging trick}
This allows one to replace the LCK structure with a locally $G$-invariant LCK structure.
Note that the \index[terms]{surface!Inoue} Inoue surfaces and the OT-manifolds are solvmanifolds, their
complex structure is locally $G$-invariant, and the averaging procedure \index[terms]{manifold!Oeljeklaus--Toma (OT)}
can be applied, following \cite{bel}.
Clearly, on any LCK-solvmanifold, the length of the 
Lee field\index[terms]{Lee field} is constant. If the LCK structure has potential, 
then, by \ref{_Gauduchon_Lee_length_Proposition_}, it is Vaisman. 
This can be used to prove the following:

\hfill

\claim {(\cite[Main Theorem]{saw2}, \cite[Main Theorem 2]{saw3})}\\
 Let $(M=G/\Lambda, I, \omega, \theta)$ be an LCK solvmanifold, $\dim M\geq 4$. Assume that
the fundamental 2-form $\omega$ of the LCK structure can
be written $\omega=\theta\wedge\theta^c+d\theta^c$. Then
$M$ is an LCK solvmanifold\index[terms]{solvmanifold} of Vaisman type. \endproof

\hfill 

Vaisman structures on differentiable
solvmanifolds are not yet classified. 
They are better understood
when the group is {\bf completely solvable} \index[terms]{Lie group!completely solvable}(i.e. $\ad(X)$ has only real
eigenvalues, for all $X\in\goth{g}$). For completely\index[terms]{Lie algebra!completely solvable}
solvable algebras, each form is cohomologous with an
invariant one, and thus the structure can be assumed
invariant. Then, with techniques similar to those in
\ref{inv_nil_vai}, \index[persons]{Sawai, H.} Sawai proves:

\hfill

\theorem {(\cite{saw4})}\label{complete_solv} Let
$G/\Lambda$ be a completely solvable solvmanifold\index[terms]{solvmanifold}
with a left--invariant complex structure. If $(G/\Lambda,
I)$ admits a Vaisman structure, then $G$ is isomorphic
with $H_{2n-1}\times\R$.

\hfill

\remark  For a recent survey of the results concerning LCK solvmanifolds, see \cite{aor2}.
For a classification of locally homogeneous LCK metrics on
compact real 4-manifolds, see \cite{hk3}.\index[terms]{manifold!LCK!locally homogeneous}

\chapter{Explicit LCK metrics on Inoue surfaces}\label{inoue_lck}\index[terms]{surface!Inoue}

{\setlength\epigraphwidth{0.7\linewidth}
\epigraph{
\it \scriptsize THE BATTLE OF THE ANTS\\ \medskip That is not that is .\\
           The only Word is Silence.\\
           The only Meaning of that Word is not.\\
           Thoughts are false.\\
           Fatherhood is unity disguised as duality.\\
           Peace implies war.\\
           Power implies war.\\
           Harmony implies war.\\
           Victory implies war.\\
           Glory implies war.\\
           Foundation implies war.\\
           Alas!  for the Kingdom wherein all these are at war.}
	   {\sc\tiny\;\;\;\;\; The Book of Lies Which is
           also Falsely Called BREAKS. The Wanderings or
           Falsifications of the One Thought of Frater
           Perdurabo, which Thought is itself Untrue \hfil
           (1912)}
}

What we now call ``Inoue surfaces''\index[terms]{surface!Inoue} were introduced in \cite{inoue} as examples of
compact surfaces with no curves and with $b_1=1$,
$b_2=0$. There are three types, known as $S^0$, $S^+$ and $S^-$, all obtained as quotients of
$\H\times \C$ by a discrete solvable group of complex affine
transformations, and diffeomorphic to solvmanifolds.\index[terms]{solvmanifold}

\section{Inoue surfaces of class ${S^0}$} 

They are also called of
class $S_M$ (see below).  These surfaces were generalized to higher
dimension by \index[persons]{Oeljeklaus, K.} Oeljeklaus and \index[persons]{Toma, M.} Toma (see Chapter
\ref{OT_manifolds}). In Chapter \ref{OT_manifolds}
we presented an alternative definition of Inoue surfaces of class
${S_M}$ constructed using number fields of degree 3.

Let $M\in\SL_3(\mathbb{Z})$ with one real
eigenvalue $\alpha >1$ and two complex eigenvalues $\beta$
and $\bar{\beta}$. Let  $(a_1,a_2,a_3)^t$  be 
a real eigenvector of $\alpha$ and  $(b_1, b_2, b_3)^t$  
a complex eigenvector of $\beta$. Let $G_M$ be the group of affine transformations of $\C \times \H$ generated by the transformations:
\begin{equation*}
	\begin{split}
		(z, w)& \mapsto (\beta z, \alpha w),\\
		(z, w)& \mapsto (z+b_i, w+a_i).
	\end{split}
\end{equation*}
for all $i=1, 2, 3$.

As a complex manifold,  ${S}_M$ is  $(\C \times \H)/G_M$.


The LCK structure\index[terms]{structure!LCK} constructed by
Tricerri\index[persons]{Tricerri, F.} in \cite{tric}
is given as a $G_M$-invariant globally conformally
K\"ahler structure on $\C \times \H$. In the
coordinates $(z, w)$, the formulae for the metric and
the Lee form,\index[terms]{form!Lee} respectively, are:
\begin{equation}\label{lck_on_S0}
	\begin{split}
		g&=\frac{dw \otimes d\overline{w}}{w_2^2}+w_2dz \otimes d\bar{z},\\
		\theta& = \frac{dw_2}{w_2},
	\end{split}
\end{equation}
where $w_2=\Im(w)$. One can see that the metric is a
warped product of the flat metric with the metric of
constant negative curvature of $\H$. See also Chapter
\ref{OT_manifolds} for a generalization of this LCK structure.

\hfill

\remark It is interesting to note that this example can be
characterized as the only compact LCK surface with
non-flat Weyl reducible structure, as shown in
\cite{_Madani_Moroianu_Pilca:Weyl_}.

\section{Inoue surfaces of class $S^+$}\label{_Inoue_original_S^+_Section_}

The definition given in \cite{inoue} for the surfaces of type $S^+$ is the following.

Let $N=(n_{ij}) \in
\SL(2,\Z)$ be a  unimodular matrix with real eigenvalues
$\alpha$ and $1/\alpha$, with $\alpha >1$. Let $(a_1, a_2)^{t}, (b_1, b_2)^{t}$ be  real eigenvectors corresponding to $\alpha$ and $1/\alpha$. Choose    $p, q, r\in\Z$, with $r \neq 0$,  and $t\in\C$. Let $e_1, e_2$, and $c_1, c_2$ be defined respectively by:
\begin{equation*}
	\begin{split}
		e_i& = \tfrac{1}{2}n_{i1}(n_{i1}-1)a_1b_1 + \tfrac{1}{2}n_{i2}(n_{i2}-1)a_2b_2 + n_{i1}n_{i2}b_1a_2,\\ 
		(c_1, c_2)&=(c_1, c_2) \cdot N^t + (e_1, e_2) + \frac{b_1a_2-b_2a_1}{r}(p, q).
	\end{split}
\end{equation*}
Let $G^+_{N, p, q, r, t}$ be the group of analytic automorphisms of $\H\times \C$, with coordinates $(w,z)$, generated by the following transformations:
\begin{equation*}
	\begin{split}
		g_0(w,z)& = (\alpha w, z+ t)\\
		g_i(w,z) & =(w+a_i, z+b_iw+c_i), \qquad i=1, 2\\
		g_3(w,z) & =(w, z+ \frac{b_1a_2-b_2a_1}{r}).
	\end{split}
\end{equation*}
\index[persons]{Inoue, Ma.} Inoue shows that the action of $G^+_{N, p, q, r, t}$ on $\H\times \C$ is properly discontinuous and without fixed points. Then  ${S}^+_{N, p, q, r, t}$ is the quotient $(\H\times\C)/G^+_{N,p, q, r, t}$.

The transformations $g_0, g_1, g_2$ and $g_3$ satisfy the relations:
\begin{equation*}
	\begin{split}
		g_3g_i&= g_ig_3,  \ \ \ \ \qquad {\rm {for}} \qquad i=0, 1, 2,\\
		g_3^r& =g_1^{-1}g_2^{-1}g_1g_2\\
		g_0g_1g_0^{-1}& =g_1^{n_{11}}g_2^{n_{12}}g_3^{p},\\
		g_0g_2g_0^{-1}& =g_1^{n_{21}}g_2^{n_{22}}g_3^{q}.
	\end{split}
\end{equation*}

\remark (\cite{inoue}) As an abstract group, $G^+_{N, p, q, r, t}$ does not depend on $t$.

\hfill

These surfaces are solvmanifolds. Their solvmanifold\index[terms]{solvmanifold} structure was clarified by C. T. C. Wall. We shall briefly describe this construction following   \cite{_Wall_proceedings_,_Wall_Topology_}.\footnote{We are grateful to Victor \index[persons]{Vuletescu, V.} Vuletescu for many insightful conversations about the construction of the \index[terms]{surface!Inoue} Inoue surfaces.}
 
 \subsection{The solvable group $\Sol^4_1$}

 \subsubsection{The structure of complex Lie groups}
  
 Consider the following real 1-parameter family of sets of matrices:
{ 
 \begin{equation}
 	G_u:={\small\left\{\begin{pmatrix} 1 & a & {\scriptstyle a + u\1 \ln  \alpha}\\[1mm] 0 & \alpha& b\\[1mm] 0&0&1\end{pmatrix}\ ; \ a,b,c\in\R,\ \alpha\in\R^{>0}\right\}}.
 	\end{equation} }
The following claim is almost obvious:

\hfill

\claim\label{_transitive_HxC_Claim_}
For each $u\in\R$, the set $G_u$
is a solvable Lie group that acts transitively 
on $\H\times\C$ with trivial isotropy (that is, 
the stabilizer group of a point is trivial).

\hfill

\proof It is trivial to check that $G_u$ is a solvable group. To describe the action on $\H\times\C$, let $A_{a,b,c,\alpha}\in G_u$ and write the couple $(w,z)\in \H\times \C$ as a column vector $(z,w,1)^t$. Then 
$$\small A_{a,b,c,\alpha}\cdot \begin{pmatrix}z\\w\\1\end{pmatrix}:=
\begin{pmatrix}z+aw+b+uc\ln
  \alpha\\ \alpha w+c\\ 1\end{pmatrix},$$
i.e. $A_{a,b,c,\alpha}\cdot(w,z)=(\alpha w+c, z+aw+b+uc\ln \alpha)$. 
The action is clearly transitive; it is thus enough to compute the isotropy at $(0,\1)$. A simple computation shows that this is $\{\Id\}$. \endproof

\hfill

\corollary\label{_Isomorphism_with_C_times_H_Corollary_}
Fix $p\in \H \times\C$.
Then the map $G_u \stackrel {\rho_p} \arrow \H \times\C$,  
$g \mapsto g(p)$ (\ref{_transitive_HxC_Claim_}) is
bijective, and can be used to define the group structure on
$\H \times\C$. Moreover, the map $\rho_p$ induces a 
left-invariant complex structure on $G_u$, independent
from the choice of a reference point $p$.

\hfill

Let $\g_u$ be the Lie algebra of the group $G_u$. The following claim is immediate.

\hfill

\claim
The generators of $\g_u$ (as a real vector space) are:
\begin{equation*}
{\scriptstyle 
	Y:=\begin{pmatrix}0&1&0\\0&0&0\\0&0&0\end{pmatrix},\ 
	U:= \begin{pmatrix}0&0&0\\0&0&1\\0&0&0\end{pmatrix},\ 
	Z:=\begin{pmatrix}0&0&1\\0&0&0\\0&0&0\end{pmatrix},\ 
	T_u:= \begin{pmatrix}0&1&\1 u\\0&0&0\\0&0&0\end{pmatrix}.}
\end{equation*}
The element $Z$ is central. The non-zero brackets are:
\begin{equation}\label{_Non_zero_commutators_Equation_}
[Y,U]=Z,\ [Y,T_u]=Y, \ [U, T_u]=U.\qquad\qquad
\end{equation}
\hfill{\endproof}

\remark Note that through our identification of $G_u$ with $\H\times \C$, the vector space $\g_u$ is identified with $T_{(0,\1)}(\H\times \C)$. If we let $z:=(z_1,z_2)\in\C$ and $w:(w_1,w_2)\in\H$ (with $w_2>0$), then we can identify:
$$Y=\frac d{dz_2},\ U=\frac d{dw_1},\ Z=\frac d{dz_1}, \ 
T_u=\frac d{dw_2}+u\frac d{dz_2}.$$

The almost complex structure  $I_u$ on $\g_u$ induced by the above identification reads:
\begin{equation}\label{_action_of_I_u_Equation_}
	\begin{split}
		I_u(Y)&=-Z,\qquad\quad I_u(Z)=Y,\\
        I_u(U)&=T_u-uY,\ \ \ \ I_u(T_u)=-U-uZ.
    \end{split}
\end{equation}
One immediately checks that this almost complex structure is integrable. 

\hfill

\remark
All groups  $G_u$ are isomorphic. All complex manifolds $(G_u, I_u)$ are biholomorphic. However, the biholomorphism is not necessarily a group isomorphism.  The next result clarifies this observation.

\hfill

\theorem(\cite{_Wall_proceedings_,_Wall_Topology_})\label{_Main_Wall_Theorem_}
\begin{description}
	\item[(i)] The groups $(G_u, I_u)$ and $(G_s,I_s)$ are biholomorphically isomorphic for all $u,s\in\R\backslash\{0\}$.
	\item[(ii)] For $u\neq 0$, there exists no map $f:G_0\to G_u$ that is  simultaneously a group isomorphism and a biholomorphism.
\end{description}

\proof 
Let $\f:(\g_u,I_u) \to (\g_s,I_s)$ be an isomorphism of Lie algebras (it exists by \ref{_Isomorphism_with_C_times_H_Corollary_}).

The generator $Z$ is central in both $\g_u$ and $\g_s$, and hence there exists $\la\in\R, \la\neq 0$, such that $\f(Z)=\la Z$. 

Suppose  $\f$ is holomorphic. Then $\f(I_u(Z))=I_s(\f(Z))$ which, by 
\eqref{_action_of_I_u_Equation_}, implies 
\begin{equation}\label{_phi(Y)_Equation_}
	\f(Y)=\la Y.
\end{equation}
Now let  $\f(T_u)=\al Z+\be Y+\ga T_s+\delta U$. Note that $[\f(Y),\f(T_u)]=\f([Y,T_u])=\f(Y)$, because $\f$ is a morphism of Lie algebras and since $[Y,T_u]=Y$ by \eqref{_Non_zero_commutators_Equation_}. A short, direct computation now brings $\ga=1$ and $\delta=0$, therefore
\begin{equation}
	\f(T_u)=\al Z+\be Y + T_s.
\end{equation}
Similarly, we have $\f(U)=mZ+nY+pT_s+qU$. Applying $\f$ to the equality $[Y,U]=Z$ (\eqref{_Non_zero_commutators_Equation_}), brings
\begin{equation}
	\f(U)=mZ+nY+U.
\end{equation}
In the same manner, applying $\f$ to the equality $[T_u,U]=U$ brings
$$[\al Z+\be Y +T_s, mZ + nY + U]=mz + nY + U.$$
After a straightforward computation, using the linear independence of $Z, Y, U$, we obtain $n=0$ and $\be=m$. Therefore:
\begin{equation}\label{_phi(T_u)_and_phi(U)_Equation_}
	\f(T_u)=\al Z + \be Y +T_s,\qquad \f(U)=\be Z + U.
\end{equation}
We have $\f(I_u(U))=I_s(\f(U))$, and hence $\f(T_u-U)=I_s(\be Z+U)$. This implies 
$$\al Z + \be Y +T_s - t\la Y =\be Y + T_s -sY,$$
thus $\al=0$ and $\la u=s$. 

Moreover, the equation $\f(I_u(T_u))=I_s(\f(T_u))$ gives into $\f(-U-uZ)=I_s(\be Y+T_s)$ which puts no further restriction on $\be$ and implies again $\la u=s$.

In conclusion, if $u=0$ and $s\neq 0$, the equation $\la
u=s$ is absurd, implying  that  the isomorphism $\f$ cannot be
holomorphic. Whereas for $u,s\neq 0$, the equation $\la
u=s$ can be solved uniquely  for $\la$. Then  
\eqref{_phi(T_u)_and_phi(U)_Equation_} in which $\al=0$ tells that for each $\be\in \R$ an isomorphism $\f$ exists and, by the above proof, it is holomorphic.  
\endproof

\hfill

In view of Wall's \ref{_Main_Wall_Theorem_}, it is enough to discuss the groups $(G_0, I_0)$ and $(G_1,I_1)$. They will be denoted respectively by $\Sol^4_1$ and ${\Sol'}^4_1$.
We resume their definitions.

\subsubsection{The group $(\Sol^4_1, I_0)$}
\begin{equation}
	\Sol^4_1=\left\{\begin{pmatrix} 1 & a & c\\ 0 & \alpha& b\\ 0&0&1\end{pmatrix}\ ; \ a,b,c\in\R,\ \alpha\in\R^{>0}\right\},
\end{equation}
with Lie algebra ${\sol}^4_1$ generated by
\begin{equation}
	{\small 
		Y:=\begin{pmatrix}0&1&0\\0&0&0\\0&0&0\end{pmatrix},\ 
		U:= \begin{pmatrix}0&0&0\\0&0&1\\0&0&0\end{pmatrix},\ 
		Z:=\begin{pmatrix}0&0&1\\0&0&0\\0&0&0\end{pmatrix},\ 
		T:= \begin{pmatrix}0&0& 0\\0&1&0\\0&0&0\end{pmatrix},}
\end{equation}
with commutators
\begin{equation}\label{_sol_4_commutators_Equation_}
	Z\ \text{central},\ 	[Y,U]=Z,\ [Y,T]=Y,\ [T,U]=U,
\end{equation}
and complex structure acting as follows:
\begin{equation}
\begin{split}
	I_0(Y)&=-Z,\quad I_0(Z)=Y\\
	I_0(U)&=T,\qquad I_0(T)=-U.
\end{split}
\end{equation}		
	
\subsubsection{The group $({\Sol'}^4_1, I_1)$}
\begin{equation}
	{\Sol'}^4_1={\small\left\{\begin{pmatrix} 1 & a &
        {\scriptstyle c+\1 \ln \alpha}\\ 0 & \alpha& b\\ 0&0&1\end{pmatrix}\ ; \ a,b,c\in\R,\ \alpha\in\R^{>0}\right\}},
\end{equation}
with Lie algebra ${\sol'}^4_1$ generated by
\begin{equation}
	{\scriptstyle 
		Y:=\begin{pmatrix}0&1&0\\0&0&0\\0&0&0\end{pmatrix},\ 
		U:= \begin{pmatrix}0&0&0\\0&0&1\\0&0&0\end{pmatrix},\ 
		Z:=\begin{pmatrix}0&0&1\\0&0&0\\0&0&0\end{pmatrix},\ 
		T:= \begin{pmatrix}0&0& \1\\0&1&0\\0&0&0\end{pmatrix},}
\end{equation}
with same commutators as for $\sol^4_1$ (\ref{_sol_4_commutators_Equation_}) 
and complex structure acting as follows:
\begin{equation}
	\begin{split}
		I_1(Y)&=-Z,\qquad \quad I_1(Z)=Y\\
		I_1(U)&=T-Y,\qquad I_1(T)=-U-Z.
	\end{split}
\end{equation}		
	
\subsubsection{Non-existence of LCK metrics on ${\sol'}^4_1$}

\theorem (\cite{va_torino_87})\label{_No_LCK_on_Sol'_Theorem_}
The group ${\Sol'}^4_1$ admits no left-invariant LCK metric.

\hfill

\proof The following proof is essentially the one in \cite[Proposition 18]{bel}. It is enough to show that the Lie algebra
$({\sol'}^4_1, I_1)$ admits no LCK metric.\footnote{For
a definition of LCK Lie algebras, see 
\ref{_LCK_Lie_algebra_Definition_}.} For simplicity,
we drop the subscript $1$ and denote the complex structure
by $I$.

We will prove the statement by contradiction. Assume there exists an LCK structure\index[terms]{structure!LCK} $(\omega,\theta)$ on $({\sol'}^4_1, I)$. Recall that on a Lie algebra $\g$, the derivative of a 1-form reads $d\eta(A,B)=-\eta([A,B])$ for all elements of  $g$. In  particular, a closed 1-form vanishes on all commutators.  Since, by \eqref{_Non_zero_commutators_Equation_}, the generators $Y$, $Z$ and $U$ are commutators, we have
$$\theta(U)=\theta(Z)=\theta(Y)=0.$$
This forces $\theta$ to not vanish on the generator $T$,
and hence $T$ should be a multiple of the Lee
field.\footnote{In the above computation we only used the
  brackets in which $t$ does not appear. It follows that
  the space of Lee forms\index[terms]{form!Lee} on $({\sol'}^4_1, I_1)$ and on $({\sol}^4_1, I_1)$, if it is non-empty, is a single point. This was observed in \cite{_Apostolov_Dloussky_} for the more general case of LCS structures\index[terms]{structure!LCS} on compact complex surfaces.}

Apply the equality $d\omega=\theta\wedge \omega$ on $(T,Z, IZ)$, using \eqref{_Non_zero_commutators_Equation_} and \eqref{_action_of_I_u_Equation_}. We obtain
\begin{equation*}
	\begin{split}
		d\omega(T,Z, IZ)&=\omega([T,Z],IZ)+
\omega([Z,IZ],T)+\omega([IZ,T],Z)\\
		&=\omega([IZ,T],Z)=\omega([Y,T],Z)=\omega(Y,Z)\\
		&=\omega(IZ,Z)\neq 0;\\
		(\theta\wedge\omega)(T,Z, IZ)&=\theta(T)\omega(Z,IZ),
	\end{split}
\end{equation*}
and thus $\theta(T)=-1$.

Apply now $d\omega=\theta\wedge \omega$ on $(T,IT, IZ)$. With similar computations, using also $\theta(T)=-1$, we obtain
$$d\omega(T,IT,IZ)=\omega(Z,T)-\omega(IZ,Z),\quad (\theta \wedge\omega)(T,IT,IZ)=-\omega(T,Z).$$
Therefore
$$\omega(Z,T)-\omega(IZ,Z)=\omega(Z,T),$$
and hence $\omega(IZ,Z)=0$, that is  a contradiction. \endproof

\hfill

On the contrary:

\hfill

\theorem (\cite{tric,va_torino})\label{_LCK_metric_on_Sol'4_1_Theorem_}
The group ${\Sol}^4_1$ admits a left-invariant LCK metric.

\hfill

\proof The metric written by \index[persons]{Tricerri, F.} Tricerri uses the coordinates of $\H\times\C$ transferred on $G_0$ by the above discussed isomorphism. We prefer to give a coordinate-free expression of this metric, by writing it directly on the Lie algebra $\sol^4_1$.

Denote with $Y^\vee$, $U^\vee$, $T^\vee$, $Z^\vee$ the left-invariant 1-forms on ${\Sol'}^4_1$ dual to the vector fields $Y, U, T, Z$. Then one can use the commutator relations  \eqref{_sol_4_commutators_Equation_} to show that $$\omega:=Z^\vee\wedge Y^\vee+U^\vee\wedge T^\vee$$
 is an LCK form with Lee form\index[terms]{form!Lee} $\theta:=-T^\vee$. \endproof

\subsection{Cocompact lattices in $\Sol^4_1$ and ${\Sol'}^4_1$}\index[terms]{lattice!cocompact}

The complex structure on the solvable group defined
above can be used to obtain a clean and explicit
definition of \index[terms]{surface!Inoue} Inoue surface as a quotient space
with the induced locally invariant complex structure. 
In the sequel, we are going to show that this
definition is equivalent to the one given by Inoue.
We start by describing a cocompact
lattice in $\Sol^4_1$ and ${\Sol'}^4_1$.

Recall the definition of the Heisenberg group:\index[terms]{group!Heisenberg}
\[
H_3:=\small\left\{\begin{pmatrix}
	1&a&c\\
	0&1&b\\
	0&0&1\end{pmatrix}\ ;\ a,b,c\in\R\right\}.
\]
It is easy to see that the Heisenberg group fits into the following exact sequence
\[
1\arrow  H_3\arrow G_u\stackrel{p_1}\arrow {Q}_1\arrow 1,
\]
where ${Q}_1=\R^{>0}$ and $p_1$ acts by
\[
\small\begin{pmatrix} 1 & a & {\scriptstyle c + u\1 \ln
    \alpha}\\[1mm] 0 & \alpha& b\\[1mm] 0&0& 1\end{pmatrix}\stackrel{p_1}\arrow \alpha.
\]
This exact sequence is split, since $p_1$ has the section
\[
\alpha\mapsto\small \begin{pmatrix} 1 & 0 & u\1\ln \alpha\\ 0 & \alpha& 0\\ 0&0&1\end{pmatrix}.
\]
This means that the group $G_u$ is the semidirect product $G_u={Q}_1\ltimes H_3$.

Let $Z$ be the center of $H_3$ (and of $G_u$ as well). We also have the exact sequence
\begin{equation}\label{_H_3_exact_sequence_}
1\arrow Z\arrow H_3\stackrel{p_2}\arrow Q_2\arrow 1,
\end{equation}
where $Q_2=\R^2$ and $p_2$ acts by
\[
{\small\begin{pmatrix}
	1&a&c\\
	0&1&b\\
	0&0&1\end{pmatrix}} \arrow (a,b).
\]

In the next proposition we see that $G_u$ has a unique (up
to an isomorphism)  cocompact lattice.\index[terms]{lattice!cocompact}

\hfill

\proposition\label{_Lattice_in_Sol_Proposition_}
Let $\Lambda$ be a cocompact lattice \index[terms]{lattice!cocompact}in $G_u$. Then there exist the integers $p,q,r\in\Z$, $r\neq 0$,  and a matrix $N=(n_{ij})\in \SL(2,\Z)$ such that the group $\Lambda$ has the following finite presentation:
\[
\Lambda=\langle g_0,g_1,g_2,g_3\ ;\ w_{30},w_{31},w_{32}, w_{12}, w_{01}, w_{02}\rangle,
\]
where the relations $w_{ji}$ ($j=0,1,3$, and  $i=0,1,2$) are:
\begin{align}
	w_{3i}&=g_3g_ig_3^{-1}g_i^{-1},\quad i=0,1,2,\label{_w_3i_relation_}\\ 
	w_{12}&=g_1g_2g_1^{-1}g_2^{-1}g_3^{-r},\label{_w_12_relation_}\\
	w_{01}&=g_1^{n_{11}}g_2^{n_{12}}g_3^pg_0^{-1}g_1g_0,\label{_w_01_relation_}\\
	w_{02}&=g_1^{n_{21}}g_2^{n_{22}}g_3^pg_0^{-1}g_1g_0.\label{_w_02_relation_}
\end{align}

\proof Since $H_3$ is nilpotent and 3-dimensional, it is the maximal nilpotent subgroup of $G_u$.  By \index[persons]{Mostow,  G. D.} Mostow's
theorem (\cite[Theorem 3.3]{ragu}), it follows that
$\Lambda_3:=\Lambda\cap H_3$ is also cocompact in $H_3$.
Now \cite[Corollary 1, p. 31]{ragu} 
implies that $\Lambda_3\cap Z$
is cocompact in the center $Z$ of $H_3$. Therefore, as
$Z\simeq \R$ is 1-dimensional, $\Lambda_3\cap Z$ has rank
one;  it is thus generated by a  single element, call
it $g_3$. On the other hand, $\Lambda_3\cap Z$ is
cocompact in $Z$, thus its projection $p_2(\Lambda_3\cap
Z)$ will be cocompact in $Q_2\simeq \R^2$. This means that
$p_2(\Lambda_3\cap Z)$ is generated by   two elements, say
$p_2(g_1)$ and $p_2(g_2)$.

Therefore, $\Lambda_3\cap Z$ 
 is generated by by the elements $g_3, g_2, g_1$ and the relations \eqref{_w_3i_relation_} and \eqref{_w_12_relation_}
follow.

Let $\hat g$ denote the class $p_2(g)$ of an element of $H_3$ in $Q_2$ (see \eqref{_H_3_exact_sequence_}). Since $Q_2$ is the abelianization of $H_3$, we pass to additive notation. Then the relations \eqref{_w_01_relation_}, \eqref{_w_02_relation_} project in $Q_2$ into the relations:
\begin{equation*}
	\begin{split}
		n_{11}&\hat{g}_1+n_{12}\hat{g}_2=\widehat{\Ad_{g_0}(g_1)},\\
		n_{12}&\hat{g}_1+n_{22}\hat{g}_2=\widehat{\Ad_{g_0}(g_2)}.
	\end{split}
\end{equation*}
These relations are used to define $N=(n_{ij})\in \SL(2,\Z)$.
Since $\Ad_{g_0}$ induces an automorphism
$\widehat{\Ad_{g_0}}$ of $Q_2$, it follows that $N$
induces an automorphism of the lattice \index[terms]{lattice}$p_2(\Lambda_3)$.
The commutator in Heisenberg group\index[terms]{group!Heisenberg} defines a 
non-degenerate skew-symmetric form on $Q_2$,
taking values in $Z$. Therefore, an automorphism of
$H_3$ preserves the orientation on $Q_2$, and
hence $N\in\SL(2,\Z)$. 

Finally, $r\neq 0$ because $g_1$ and $g_2$ cannot commute
(otherwise, $H_3$ would be abelian).
 \endproof

\subsection[Equivalence with the Inoue's description of the surfaces of class $S^+$.]{Equivalence with the Inoue's description of the\\ surfaces of class $S^+$.} 

At this stage, we can finally sketch the equivalence
between the solvamanifold description and the original
description in \cite{inoue} (see also
\cite{_Hasegawa:solv_}).

Note first that for an Inoue surface $S^+_{N,p,q,r,t}$,
with $t=\eta+\1\zeta$  
(see Section \ref{_Inoue_original_S^+_Section_}), the transformation
$$\H\times\C\ni(w,z)\arrow (w, z+\zeta)$$
is a biholomorphism between $S^+_{N,p,q,r,t}$ and
$S^+_{N,p,q,r,\1\zeta}$. Therefore, 
we may always assume that 
the parameter $t$ of $S^+_{N,p,q,r,t}$ is purely
imaginary.

Consider the action of $G_u$ on $\H\times \C$
defined in \ref{_transitive_HxC_Claim_},
and let $g_0, g_1, g_2, g_3$ be the generators of 
a cocompact lattice\index[terms]{lattice!cocompact} $\Lambda\subset G_u$ 
defined in \ref{_Lattice_in_Sol_Proposition_}.

After a suitable change of coordinates on $\H\times \C$,
the action of $g_0$ can be written as
$$
  g_0(w,z)=(\alpha w, z+i\zeta).
$$
Then, since $g_1,g_2\in H_3$, they assume the form:
$$
g_i=\small\begin{pmatrix}1&a_i&c_i\\0&1&b_i\\0&0&1\end{pmatrix},
$$
in other words, they act as
$$
g_i(w,z)=(w+a_i,z+b_iw+c_i).
$$
On the other hand, $g_3\in Z$ (the central element) has the form
$$
g_3=\small\begin{pmatrix}1&0&c_3\\0&1&0\\0&0&1\end{pmatrix},
$$
i.e. it acts as
$$
g_3(w,z)=(w, z+c_3).
$$

Therefore, the relation \eqref{_w_12_relation_} implies
\[
c_3=\frac{a_1b_2-a_2b_1}{r},
\]
and 
$$
g_3(w,z)=(w, z+ \frac{a_1b_2-a_2b_1}{r}).
$$
Similar computations lead to all the relations among 
the generators in Section \ref{_Inoue_original_S^+_Section_}.

\subsection{The LCK metric on $S^+_{N,p,q,r,0}$}

In \cite{tric}, \index[persons]{Tricerri, F.} Tricerri wrote a $G^+_{N, p, q, r, t}$-invariant LCK metric  on ${S}^+_{N, p, q, r, t}$ for $t \in \R$. In view of our preceding discussion, we may assume $t=0$. 

Obviously, the  LCK metric that we wrote in \ref{_LCK_metric_on_Sol'4_1_Theorem_} is preserved by the action of the lattice,\index[terms]{lattice!cocompact} because it was written directly on the Lie algebra. It then descends to an LCK metric on  ${S}^+_{N, p, q, r, 0}$.

In the coordinates $(w,z)$ of $\H\times \C$, with $w=w_1+\1 w_2$, $w_2>0$, this metric reads: 
\begin{equation}\label{lck_on_S+}
\begin{split}
\omega& = \frac{1}{w_2^2}\left((w_2dz-z_2dw)\otimes (w_2d\bar z-z_2d\bar w)+dw\otimes d\bar w\right),\\
\theta&=\frac{dw_2}{w_2}.
\end{split}
\end{equation}

\subsection{Non-existence of LCK metrics on
  $S^+_{N,p,q,r,t}$, $t\neq 0$}

\theorem (\cite{bel})\label{_Existence_of invariant_LCK_on_Sol^4'_}
 If the \index[terms]{surface!Inoue} Inoue surface  ${\Sol'}^4_1/\Lambda$ admits an LCK metric, then it also admits an invariant LCK metric.

\hfill

\index[persons]{Belgun, F. A.} Belgun's proof relies on following important result.

\hfill

\lemma \label{_Bi_invariant_volume_Lemma_}
The group ${\Sol'}^4_1$ admits a bi-invariant volume form. 

\hfill

\proof Indeed, in the notations of \ref{_LCK_metric_on_Sol'4_1_Theorem_}, the form $V=Y^\vee\wedge Z^\vee\wedge T^\vee\wedge U^\vee$ is $\ad$-invariant. After rescaling, we may suppose that the volume of the Inoue surface $X:=S^+_{N,p,q,r,t}$ is 1. \endproof

\hfill

The existence of the bi-invariant volume form allows performing an ``average trick'' that associates to each tensor field on $X$ a left-invariant tensor field on ${\Sol'}^4_1$.

\hfill

\lemma\label{_Average_trick_Lemma}
Let $\eta\in\Lambda^mX$. Define $\eta_V\in\Lambda^m({\sol'}^4_1)$ by 
\begin{equation}\label{_Invariant_form_definition_}
	\eta_V(A_1,...,A_m):=\int_X \eta_x(A_{1x},...,A_{mx})V.
\end{equation}
Then $(d\eta)_V=d(\eta_V)$.

\hfill

\proof The argument uses the Cartan formula\index[terms]{Cartan formula} for the exterior derivative:
\begin{equation}
	\begin{split}
d\eta(A_0,..,A_m)&=\sum_i\Lie_{A_i}(\eta(A_0,...,\hat A_i,...,A_m)\\
&+\sum_{i<j}\eta([A_i,A_j],A_0,...,\hat A_i,...,\hat A_j,..., A_m).
	\end{split}
\end{equation}
The integrals of the terms in the first row are of type 
$\int_X (\Lie_A f)V$ where $A$ is any left-invariant vector field and $f\in C^\infty X$. By Cartan formula\index[terms]{Cartan formula} and by the invariance of $V$, we have
$$(\Lie_A f)V=\Lie_A (fV)= i_A (d(fV))+d(i_A (fV)).$$
Now $i_A(d(fV)=0$ by dimension reason. Then 
$$\int_X (\Lie_A f) V=\int_X d(i_A (fV))=0,$$
by Stokes. The conclusion follows. \endproof

\hfill

We can now give:

\hfill

{\bf Proof of \ref{_Existence_of invariant_LCK_on_Sol^4'_}:} 
Let $(\omega, I, \theta)$ be an LCK structure\index[terms]{structure!LCK} on $X$. 
Using \ref{_Average_trick_Lemma}, we can associate
to $\omega$ a left-invariant form $\omega_V$ on ${\Sol'}^4_1$. 
 
Since $H^1(X,\R)=1$, and the 1-form
$\theta_0=-T^\vee=\frac{dw_2}{w_2}$ is invariant and
descends to $X$, we see that, after a conformal change, we
can assume that the Lee form\index[terms]{form!Lee} of $\omega$ is
$\theta=k\theta_0$, with $k\in\R$. In particular, the Lee
form of $\omega$ is left-invariant and corresponds to a
closed 1-form on ${\sol'}^4_1$.

Then, by \ref{_Average_trick_Lemma} again,   
$$d(\omega_V)=(d\omega)_V=(k\theta_0\wedge\omega)_V=k\theta_0\wedge \omega_V,$$
and thus $\omega_V$ is an LCK form on ${\sol'}^4_1$, contradicting \ref{_No_LCK_on_Sol'_Theorem_}. \endproof

\hfill

Together with \ref{_No_LCK_on_Sol'_Theorem_}, this brings:

\hfill

\corollary (\cite{bel})\label{_no_lck_} The \index[terms]{surface!Inoue} Inoue surface $S^+_{N,p,q,r,z}={\Sol'}^4_1/\Gamma$, with $z\in\C\backslash\R$, admits no LCK metric.

\hfill

In particular, this proves:

\hfill

\theorem {(\cite{bel})} The LCK class is not stable at small deformations.

\hfill

\proof The surfaces  $S^+_{N,p,q,r,t}$ with $t\in\C\setminus \R$, which, by \ref{_no_lck_}, do not admit LCK metrics,  are complex deformations of the surfaces $S^+_{N,p,q,r,t}$ with $t\in\R$ that admit LCK metrics. \endproof

 \section{Inoue surfaces of class $S^-$.}
 
These are defined similarly to $S^+$, such that  ${S}^+_{N^2, p_1, q_2, r, 0}$ is a double unramified cover of ${S}^-_{N, p, q, r}$. The LCK structure \eqref{lck_on_S+} induces an LCK structure on $S^-$. \index[terms]{structure!LCK}
 
 \hfill
 
 We end this chapter with the following important result:
 
 \hfill
 
 \theorem {(\cite{oti2})}\label{_no_exact_on_inoue_} An
 LCK \index[terms]{surface!Inoue} Inoue surface cannot support a $d_\theta$-exact LCK
 structure $(M, \omega,\theta)$. In particular, no LCK
 structure on an Inoue surface can be LCK with potential.\index[terms]{structure!LCK!with potential}
 
 \hfill 
 
\proof The proof follows from the following facts:
\begin{itemize}
	\item[(i)] The only possible Lee class\index[terms]{class!Lee}
          $[\theta]\in H^1(M)$ on
          an Inoue surface are the one corresponding to
          the metrics found by Tricerri (see above):
          \cite{oti2}, \cite{ad2}. In particular,
$\lambda\theta$ is never a Lee form\index[terms]{form!Lee} for an LCK structure,
unless $\lambda=1$

	\item[(ii)] From \ref{_power_of_pote_Remark_}
(see also \cite[page 331]{ov_pams} and \cite{ad2}), it
          follows that if
$\theta$ is the Lee form of an LCK metric with potential,
          then  $t\theta$ is the Lee form\index[terms]{form!Lee} of an LCK metric
          with potential\index[terms]{metric!LCK!with potential} for all $t\geq 1$. \endproof
 \end{itemize}

	
\chapter{More on Oeljeklaus--Toma manifolds}\label{moreot}


\epigraph{\it If one looks at a thing with the intention of trying to discover what it means, one ends up no longer seeing the thing itself, but of thinking of the question that is raised.}{\sc\scriptsize Ren\'e Magritte}

\hfill

This short chapter is intended to give more insight on the
manifolds introduced in Chapter \ref{OT_manifolds} (with the
notations we preserve here) and to present several recent
results.

\section{Cohomology of OT-manifolds}
\label{_coho_OT_Section_}\index[terms]{manifold!Oeljeklaus--Toma (OT)}

Let $M=X(K,U)$ be an OT-manifold, with $K$ a number field
and $U\subset \calo_K^*$ a subgroup of the group of units.
By construction, $M$ is fibred over the torus $T_U :=U\otimes_\Z \R/U$ with the
torus fibre $T_K:=\calo_K\otimes_\Z\R/\calo_K$. This fibration is in fact
locally trivial and possesses a flat Ehresmann connection\index[terms]{connection!Ehresmann}. The corresponding
monodromy action\index[terms]{action!monodromy} of $U=\pi_1(T_U)$ on the torus $T_K$ is induced by the
multiplicative action of $U$ on $\calo_K$.

The corresponding Leray-Serre spectral sequence\index[terms]{spectral sequence!Lyndon-Hochschild-Serre}\footnote{
It is sometimes called ``the Lyn\-don-\-Hoch\-schild-\-Serre spectral sequence'' 
for the group cohomology; see Exercise
\ref{_Serre-Hochschild_sequence_Exercise_} for its
construction using the Grothendieck spectral sequence.}
has $E_2$-term \[ E_2^{*,*}=H^*(T_U, H^*(T_K)).\] 

The $\pi_1(T_U)$-action on $\calo_K\otimes_\Z\C$
is diagonalizable, and this representation splits onto a direct
sum of 1-dimensional representations corresponding to the
field decomposition $\calo_K\otimes_\Z\C=\bigoplus \C$,
with each component of $\bigoplus \C$ representing
the corresponding complex or real embedding of $K$.

A 1-dimensional representation $V$ of\,  $\Gamma=\Z^n$
has $H^*(\Gamma, V)= \Lambda^*(V)^\Gamma$,
where $\Lambda^*(V)^\Gamma$ denotes the space of $\Gamma$-invariant
vectors in the corresponding Grassmann algebra. \index[terms]{Grassmann algebra}
This is deduced using the K\"unneth formula from the analogous
result for $H^*(\Z, V)$ (\ref{_mult_on_coh_of_S^1_Claim_}). 
This gives $E_2^{*,*}=H^*(T_U) \otimes H^*(T_K)^U$.

Since the connection in the
fibration $M \arrow T_U$ is flat, the Leray-Serre
 spectral sequence degenerates in $E_2$.\index[terms]{spectral sequence!Leray-Serre}
Indeed, a monodromy invariant closed form on $T_K$ can be extended\index[terms]{monodromy}
to a closed form on $M$. 

We obtained the following theorem,
that is  due to \index[persons]{Istrati, N.} Istrati and \index[persons]{Otiman, A.} Otiman \cite{isot}.

\hfill

\theorem {(\cite[Theorem 3.1]{isot})}\label{OT_cohomology}
Let $M=X(K,U)$ be an OT-manifold, with $K$ a number field
and $U\subset \calo_K^*$ a subgroup of the group of units.
Then the de Rham cohomology\index[terms]{cohomology!de Rham} of $M$ is expressed
as  $H^*(M)= H^*(T_U, H^*(T_K))$,\index[terms]{manifold!Oeljeklaus--Toma (OT)}
where $T_U :=U\otimes_\Z \R/U$ and $H^*(T_K)^U$ is the space of all $U$-invariant cohomology classes
on the torus $T_K:=\calo_K\otimes_\Z\R/\calo_K$.
\endproof

\hfill

The first cohomology was computed in the original paper
of \index[persons]{Oeljeklaus, K.} Oeljeklaus and \index[persons]{Toma, M.} Toma.

\hfill

\theorem {(\cite{ot})}\label{topinv}
Let $X(K,U)=(\C^t \times {\H}^s)/\Gamma$. Then
$H^1(X(K, U))=U\otimes_\Z \R$. In particular, $b_1(X(K,U))=s$.

\proof
The action of $U$ on $\calo_K$ has no invariants,
hence $H^1(M)= H^1(T_U)\oplus H^1(T_K)^U= H^1(T_U)$.
\endproof

\hfill

In \cite{_Otiman_Toma:Cousin_}, a similar argument was used to prove that 
the Fr\"olicher spectral sequence degenerates for all OT-manifolds.\index[terms]{spectral sequence!Fr\"olicher}

\hfill

\theorem \label{_OT_degenerate_Theorem_}
(\cite[Theorem 4.5]{_Otiman_Toma:Cousin_})\\
Let $M$ be an OT-manifold. Then\index[terms]{manifold!Oeljeklaus--Toma (OT)}
$\sum_{p+q=r} \dim H^{p,q}_{\bar\6}(M)= \dim H^r(M)$.
\endproof

\section{LCK structures on general OT manifolds}\index[terms]{structure!LCK}

\proposition {(\cite{vu2}, \cite[Appendix by L. Battisti]{du})}\label{crit_lck}
 Let $M=X(K,U)$ be an OT-manifold
associated with a number field $K$ with $s$ real
embeddings and $2t$ complex embeddings and 
an admissible group of units $U$. Then it admits an LCK metric if and only if 
 \begin{equation}\label{critlck}
|\tau_1(u)|=\cdots =|\tau_t(u)|,\quad \text{for all}\,\, u\in U.
 \end{equation}
where $\tau_1, ..., \tau_t$ is a set of pairwise
non-conjugate complex embeddings of $K$.

\hfill

\proof Necessity of \eqref{critlck} is proven in
\cite[Proposition 2.9]{ot}. Let $\tilde M$ be the
$U$-cover of $M$, realized as a product of a torus
$T_K=\calo_K\otimes_\Z \R/\calo_K$ and $(\R^{>0})^s$. Since
$H^1(M)=U\otimes_\Z \R$ (\ref{topinv}), the Lee form\index[terms]{form!Lee} is
exact  on $\tilde M$. Let $\tilde \omega$
be a K\"ahler form\index[terms]{form!K\"ahler} on $\tilde M$ conformal to
the LCK form on $M$.
Averaging $\tilde\omega$ with the $T_K$-action, we obtain 
a closed $T_K$-invariant form $\tilde\omega_1$
that is  automorphic. The coefficients of
$\omega$ depend only on variables
$z_{t+1},..., z_{t+s}$ corresponding
to the real embeddings. The automorphy
of $\tilde\omega_1$ under the action of $U$ produces a
character $\chi:\; U\arrow \R^{>0}$ such that
\begin{equation}\label{_weight_on_OT_Equation_}
|\tau_1(u)|^2(dz_1\wedge d\bar z_1) = \chi(u)dz_1\wedge
d\bar z_1, \ \  ...,  \ \  |\tau_t(u)|^2(dz_t\wedge d\bar z_t) = \chi(u)dz_t\wedge
d\bar z_t.
\end{equation}
This proves that \eqref{critlck} necessarily holds for any OT-manifold admitting\index[terms]{manifold!Oeljeklaus--Toma (OT)}
an LCK structure.

To prove that \eqref{critlck} is also sufficient,
define the following real function on $\C^t\times\H^s$
(that is  a natural extension of the potential  $\f$ used
in the definition of the LCK metric, see \ref{_Existence_of_LCK_metric_on_OT_}):
\[
\phi(x_1,\ldots, x_t, \zeta_1, ..., \zeta_s)=
\sum_{j=1}^t|x_j|^2 + \prod_{i=1}^s
\Im(\zeta_i)^{\frac{-1}{t}}
\]
It is clear that $\f$ is a K\"ahler
potential. Then one examines the action of a $g=(u,k)\in
U\ltimes \calo_K$ on the K\"ahler form $\ri\6\bar\6 \f$
and finds that $g^*\ri\6\bar\6 \f=|\sigma_j(u)|^2\ri\6\bar\6
\f$. Hence, for the K\"ahler form\index[terms]{form!K\"ahler} to be acted on by
homotheties the condition \eqref{critlck} is sufficient. 
\endproof

\hfill

The following result was achieved in several steps, starting with \cite{ot}, passing through the papers \cite{vu2} and  \cite{du}, and finally by \cite{_Deaconu_Vuletescu_}. The proofs rely heavily on number theoretical results.

\hfill

\theorem \label{_OT_LCK_t=1_Theorem_}
Let $M=X(K,U)$ be an OT manifold associated with a
number field $K$ with $s$ real embeddings and $2t$ complex
embeddings. Then $M$ has an LCK metric if and only of
$t=1$.

\section{Cohomology of LCK OT-manifolds}\index[terms]{manifold!Oeljeklaus--Toma (OT)}


Before it was known that LCK OT-manifolds satisfy $t=1$
(\ref{_OT_LCK_t=1_Theorem_}), the task of computing the cohomology
of such a manifold was already of considerable importance and interest.
The problem was solved in \cite{isot}, using \ref{OT_cohomology} and \ref{crit_lck}.
Here we give an updated version of this result, relying on \ref{_OT_LCK_t=1_Theorem_}.

\hfill

In this section, we deal only with OT-manifolds associated with number fields\index[terms]{manifold!Oeljeklaus--Toma (OT)}
that have exactly 2 complex conjugate complex embeddings.
Let $\sigma_1, \ ...,\ \sigma_s$
be all real embeddings of $K$, and $\tau$ one of the two complex-conjugate complex embeddings. 
Elements of $U$ are units of the ring $\calo_K$ that are  positive in 
all real embeddings.

\hfill

\theorem\label{LCKbun}
Let $M=X(K,U)$ be an OT manifold of type $(s,1)$.
Then  its de Rham cohomology\index[terms]{cohomology!de Rham} algebra
$H^\bullet(X,\CC)$ is isomorphic to the graded algebra
over $\CC$ generated by: 
\begin{equation}\label{_generators_for_OT_LCK_Equation_}
dr_1,\ldots, dr_s, dz_1\wedge d\bar z_1\wedge dv_1\wedge ... \wedge dv_s,
\end{equation}
where $r_i$ is 
the absolute value of the complex
coordinate on $i$-th hyperbolic plane ${\Bbb H}_i$, and $v_i$ the real coordinate on ${\Bbb H}_i$.
In particular, its Betti numbers are:\index[terms]{Betti numbers!rational}
\begin{align*}
b_l=b_{2m-l}&=\binom sl \text{ for } 0\leq l\leq s,\\
b_l&=0 \text{ for } l=s+1,
\end{align*} 
where $m=\dim_\C M=s+1$.

\hfill

\proof 
The norm of $u\in U$ is given by
\[
N(u)=|\tau(u)|^2\prod_{j=1}^s\sigma_j(u).
\]
Since the norm of any $u\in U$ is positive, integer and invertible, we have $N(u)=1$, giving
\begin{equation}\label{_norms_LCK_Equation_}
|\tau(u)|^{-2}=\prod_{j=1}^s\sigma_j(u).
\end{equation}
Let  $z, \bar z$ be the coordinates corresponding to $\tau$.
Consider the monomial
\[ Q_{I,\alpha} = \alpha \wedge \prod_{I=k_1< ...< k_r} dv_{k_1}\wedge dv_{k_2} \wedge \cdots \in \Lambda^*(T_K),
\]
where $\alpha$ is a monomial on $dz$, $d\bar z$.
The monomial $Q_{I,\alpha}$ is $U$-invariant if and only if
\begin{equation}\label {_prod_sigma_trivial_Equation_}
\tau_\alpha(u) \prod_r
 \sigma_{k_n}(u) =1,
\end{equation}
where $\tau_\alpha$ is the number $\tau_\alpha=\frac{\tau(u)^*(\alpha)}{\alpha}$,
that is, $\tau_\alpha=\tau(u)$ for
$\alpha=dz$,  $\tau_\alpha=\bar\tau(u)$ for
$\alpha=d\bar z$, $\tau_\alpha=|\tau(u)|^2$ for
$\alpha=dz\wedge d\bar z$, and $\tau_\alpha=1$ for $\alpha=1$.

Let $Q_{I,\alpha}$ be a monomial that is  invariant under
the action of $U$ on $\C\times {\Bbb H}^s$. Then \eqref{_prod_sigma_trivial_Equation_} gives
a multiplicative relation between $\sigma_i(u)$ and $|\tau(u)|$. Using
\eqref{_prod_sigma_trivial_Equation_} and taking
the absolute value, we obtain that $u^*(Q_{I,\alpha})=Q_{I,\alpha}$
implies that
\begin{equation}\label{_product_sigma_i_via_monomial_Equation_}
1= |\tau(u)|^d\prod_{i=1}^r \sigma_{k_i}(u)=
\prod_{i=1}^s \sigma_i(u)^{-\frac d 2} \prod_{i=1}^r \sigma_{k_i}(u),
\end{equation}
for all $u\in U$ where $d$ is the degree  of the monomial $\alpha$.
By construction of an OT-manifold, the map
$\sigma_1(u), ..., \sigma_s(u)$ defines a lattice\index[terms]{lattice} in
$(\R^{>0})^s$; then  the matrix $\bigg[\log(\sigma_i(u_j))\bigg]$
is non-degenerate for any basis $\{u_i\}$ in $U=\Z^s$.
Therefore, a relation of form
\eqref{_product_sigma_i_via_monomial_Equation_} is impossible,
unless $dv_{k_1}\wedge dv_{k_2} \wedge \cdots dv_{k_r}=dv_1\wedge dv_2 \wedge ... \wedge dv_n$ and
$\alpha=dz\wedge d\bar z$.

Now, \ref{OT_cohomology}
implies that the de Rham algebra of $M$ is freely generated
by the $U$-invariant monomials on $dv_1, ..., dv_s$, $dz, d\bar z$ and
all monomials on $dr_1, ..., dr_s$, which gives 
\eqref{_generators_for_OT_LCK_Equation_}.
\endproof

\hfill

This has the following important consequence which
generalizes the corresponding result for \index[terms]{surface!Inoue} Inoue surfaces in
\cite{oti2}:

\hfill

\theorem {(\cite{isot})} \label{unic}
Let $X=X(K,U)$ be an OT manifold of type $(s,1)$. There
exists at most one Lee class\index[terms]{class!Lee} of an LCK metric on $X$,
namely the one represented by the $U\ltimes O_K$-invariant
form on $\HH^s\times\C$, 
\begin{equation}\label{theta_for_OT}
\theta=d\log \left(\prod_{k=1}^sv_k\right).
\end{equation}
\proof
From \ref{topinv} we obtain that $\theta\in H^1(T_U)$,
where $T_U$ is the torus considered in Section \ref{_coho_OT_Section_}.
Assume that $u\in U$ acts on $\HH^s\times\C$ as a homothety 
with the factor $\chi(u)$.
The automorphy conditions \eqref{_weight_on_OT_Equation_}
together with the formula \eqref{_norms_LCK_Equation_}
imply that $\chi(u)=\prod_{i=1}^s\sigma_i(u)$.
Then $\chi$ is determined uniquely. \endproof

 \hfill
 
The following corollary immediately follows.

\hfill
 
 \corollary {(\cite{isot}; see also \cite{goto})}  On
 any OT manifold of LCK type, the set of Lee classes\index[terms]{class!Lee} of 
LCK structures\index[terms]{structure!LCK} is not open in $H^1(M, \R)$. 
\endproof

\hfill

\remark
Using the notion of ``stability'' of \index[persons]{Kodaira, K.} Kodaira and \index[persons]{Spencer, D. C.} Spencer
\cite{_Kod-Spen-AnnMath-1960_}, this can be expressed by
saying that ``LCK structures are not
 stable under the small deformations of the Lee class''.
 
 \hfill
 
\remark
The same proof as in \ref{LCKbun} 
also gives a complete description of the
Morse--Novikov cohomology, \cite[Proposition 6.7]{isot}.

 
%

\hfill
 
\remark
A. \index[persons]{Otiman, A.} Otiman proved in \cite{oti2} that LCK OT manifolds 
 cannot have $d_\theta$-exact LCK forms (in
particular, they cannot have LCK metrics with
potential).

\section{LCK rank of OT manifolds}\index[terms]{rank!LCK}

OT manifolds provide examples of compact LCK manifolds
with LCK rank (see \ref{lck_rank}) strictly greater than 1.

\hfill

\theorem  {(\cite{pv})} \label{PV}
Let $X$ be an LCK OT-manifold associated with a number
field $K$ with 1 complex and $s$ real embeddings, and
$r$ its LCK rank. Assume that \index[terms]{manifold!Oeljeklaus--Toma (OT)}
$K$ is a quadratic extension of a totally real number
field. Then   $r=\frac{b_1(X)}{2}$.
Otherwise $r=b_1(X)$.


\hfill

\proof  
From \ref{OT_cohomology} it follows that
$H_1(X,\R)=U\otimes_\Z \R$.
Therefore, the homothety character \index[terms]{homothety character}$\chi$
can be considered to be  a homomorphism $U \arrow \R^{>0}$.
 For any $u\in U$ (seen as an element in
$\pi_1(X)$), the automorphy factor $\chi(u)$ is
$\vert\tau(u)\vert$, where $\tau$ is the only (up to
complex conjugation) complex embedding of $K$.  If $\rk X
\neq b_1(X)$, then at least one $u\in U$ must have
$\vert\tau(u)\vert=1$. This forces $u$ to be a {\bf
  reciprocal unit,} i. e.\  its minimal polynomial over
$\Q$ to be a reciprocal one. However, if $u$ is a reciprocal
unit, then the field $K'=\Q(u+\frac{1}{u})$ is a subfield
of $K$,  of relative degree $2$,  and it can be  shown
that $K'$ must  be totally real.
\endproof


\chapter{Locally conformally parallel and non-parallel structures}\label{other}\index[terms]{structure!locally conformally parallel}

{\setlength\epigraphwidth{0.75\linewidth}
\epigraph{\it
Let's go swimming in the womb of nothingness, where our
conversations can fit together like dove-tailed joints
without coming to a stressful end. Let's go do nothing!
Scattering our words like placidly sprinkling water!
Indifferent as pure essence! Tuned in and leisurely! Our
aspirations will have already become vacant, and as we
move on we won't know where we'll end up. }
{\sc\scriptsize Chuang Tz\v u Chapter 22: Knowledge Wanders North}}


This chapter is devoted to a brief description of
different extensions of the ideas and/or techniques from
LCS and LCK geometries\index[terms]{geometry!LCS}\index[terms]{geometry!LCK} to other settings. Many geometric
structures have locally conformal versions with
interesting geometry.

To clarify our choices for the topics presented below,
recall the \index[persons]{Berger, M.} Berger list of possible holonomies of\index[terms]{theorem!Berger list}
irreducible, non locally symmetric  Riemannian manifolds
(see \cite{besse}):

\hfill

\begin{center}

\begin{tabular}{|l|l|}
\hline
\multicolumn{2}{|c|}{\bf \color{red}Berger's list}\\[1mm]
\hline
\it \bf Holonomy  & \it \bf Geometry\\[1mm]
\hline
$\SO(n)$ acting on $\R^n$ & Riemannian manifolds\\[1mm]
\hline
$\U(n)$ acting on $\R^{2n}$ & K\"ahler manifolds\\[1mm]
\hline
$\SU(n)$ acting on $\R^{2n}$, $n>2$ & Calabi--Yau manifolds\\[1mm]
\hline
$\Sp(n)$ acting on $\R^{4n}$ & hyperk\"ahler manifolds\\[1mm]
\hline
$\Sp(n)\times \Sp(1)/\{\pm 1\}$ & 
quaternionic-K\"ahler\\[1mm] acting on $\R^{4n}$, $n>1$ &  manifolds\\[1mm]
\hline
$G_2$ acting on $\R^7$ & $G_2$-manifolds \\[1mm]
\hline
$\Spin(7)$ acting on $\R^8$ & $\Spin(7)$-manifolds\\[1mm]
\hline
\end{tabular}
\end{center}

Recall that the Riemannian cone\index[terms]{cone!Riemannian} over a Riemannian
manifold $(X, g)$ is the Riemannian manifold
$C(X)=(X\times \R^{>0}, t^2g+ dt\otimes dt)$, where $t$
is the usual parameter in $\R^{>0}$.
 We are interested in manifolds that are  locally
 conformally equivalent to manifolds in the \index[persons]{Berger, M.} Berger's
 list. In many significant cases, the universal cover of
 such a manifold is a Riemannian cone. For example,  the
 universal cover of a Vaisman manifold\index[terms]{manifold!Vaisman} is a cone over a
 Sasakian manifold\index[terms]{manifold!Sasaki} (\ref{halfstr}, Subsection \ref{str_th}). One can
 establish a dictionary between the  geometry of a
 manifold and the geometry of its Riemannian cone.\index[terms]{cone!Riemannian} We
 summarize it in the following table:

\hfill

\begin{center}\index[terms]{manifold!3-Sasakian}\index[terms]{manifold!Calabi--Yau}\index[terms]{manifold!hyperk\"ahler}\index[terms]{manifold!quaternion-K\"ahler}\index[terms]{manifold!$G_2$}\index[terms]{manifold!$\Spin(7)$}\index[terms]{group!holonomy!special}
	
\begin{tabular}{|l|l|l|}
\hline
\multicolumn{3}{|c|}{\bf \color{red} Riemannian cones with
special holonomy}\\[1mm]
\hline
\it\bf \color{blue} Holonomy of $C(M)$  &\it\bf\color{blue} Geometry of $C(M)$  & \it\bf\color{blue} Geometry of $M$\\[1mm]
\hline
$\SO(n)$ &  Riemannian & --- \\[1mm]
\hline
$\U(n)$ &  K\"ahler & Sasakian\\[1mm]
\hline
$\SU(n)$  & Calabi--Yau &
Sasaki--Einstein \\[1mm]
\hline
$\Sp(n)$ &  hyperk\"ahler & 3-Sasakian\\[1mm]
\hline
$\Sp(n)\times \Sp(1)$ & 
quaternion-K\"ahler & --- \\[1mm]
\hline
$G_2$  & $G_2$-manifolds & nearly K\"ahler \\[1mm]
\hline
$\Spin(7)$  & $\Spin(7)$-manifolds & nearly $G_2$-manifolds \\[1mm]
\hline
\end{tabular}
\end{center}

As regards the reduced 
holonomy group of an LCK metric on a compact manifold the following result is available:

\hfill

\theorem (\cite{_Madani_Moroianu_Pilca:Einstein_}) Let $(M,I,g,\theta)$ be a compact, strict LCK manifold (i.e. $\theta$ is not exact), with non-generic
holonomy group $\Hol(M,g)\subsetneq \SO(2n)$. Then $\Hol(M,g)\cong  \SO(2n-1)$ and $(M,I,g,\theta)$ is Vaisman. \endproof

\hfill

\remark Another striking result is \cite[Theorem 5.1]{_Madani_Moroianu_Pilca:Einstein_}. It states that if $(M,I,g)$ is a compact strict LCK manifold, there cannot exist any integrable complex structure $J$ on $M$ such that $(g,J)$ is K\"ahler.

\hfill

For later use, recall:

\hfill

\theorem {(\cite[Proposition 3.1]{gal})}\label{gallot} The Riemannian cone\index[terms]{cone!Riemannian} over a
  complete manifold has reducible holonomy  if and only if it is flat, in which case it is 
  a  cone over a sphere.  

  
\section{Locally conformally hyperk\"ahler
  structures} \label{lchk}


Let $M$ be a smooth manifold of real dimension $4k$ equipped with three operators $I_1, I_2, I_3 \in \End (TM)$
satisfying the quaternion relations 
\[ I_1^2 = I_2^2=I_3^2 = -1, I_1\circ I_2 = - I_2\circ I_1 = I_3.
\]
Assume that 
the operators $I_1, I_2, I_3$ induce integrable complex structures on
$M$. Then $M$ is called {\bf hypercomplex}\index[terms]{manifold!hypercomplex}. A 
Riemannian metric $g$ on a hypercomplex manifold is called 
{\bf quaternionic Hermitian} if it is  Hermitian with
respect to each of the three complex structures.
By a theorem 
of \index[persons]{Obata, M.} Obata (\cite{_Obata_}), a hypercomplex manifold admits
a unique torsion-free\index[terms]{connection!torsion-free} connection\index[terms]{connection!Obata} preserving $I_1, I_2, I_3$.
If the Obata connection preserves a quaternionic Hermitian metric
$g$, the manifold $(M,g)$ is called {\bf hyperk\"ahler}. This 
is equivalent to $(M, g, I_i)$ being K\"ahler for each $i=1,2,3$.

Let $\omega_i$ be the fundamental form associated with
$(g,I_i)$ and $\Theta:=\sum_{i=1}^3\omega_i\wedge\omega_i$.

\hfill

\definition A hyper-Hermitian manifold $(M, g, I_1, I_2,
I_3)$ is {\bf locally conformally 
hyperk\"ahler} (LCHK) if $g$ is locally conformal to a
hyperk\"ahler metric.

\hfill

\remark
Locally, $g$ is conformally equivalent to a 
hyperk\"ahler metric, $g= e^f g_0$.
Let $\theta:= d f$. Then $d\omega_i = \theta \wedge \omega_i$.

\hfill

\remark
If $\dim_{\Bbb H}M >1$, the LCHK condition 
immediately implies that
\begin{equation}\label{lcqk}
d\Theta=2\theta\wedge\Theta. 
\end{equation}
It is not hard to show that this condition,
in turn, implies the LCHK condition.

\hfill

\remark\label{_LCHK_Gauduchon_Remark_}
This implies that the Gauduchon metric 
on $(M, g, I_i)$ is the same for $i=1, 2, 3$.
Indeed, the metric is Gauduchon if and only if
the Lee form \index[terms]{form!Lee}satisfies $d^* \theta=0$ (Subsection \ref{_gauduchon_basics_Subsection_}),
and the Lee form on $(M, g, I_i)$ is the same for $i=1, 2, 3$.

\hfill 

\remark On a compact 4-dimensional manifold, all
hypercomplex structures admit compatible K\"ahler or
LCHK metrics (\cite{boy}). This is no longer true in
higher dimensions. In the following, we concentrate on real
dimension at least 8.

\hfill

\remark {\bf Almost quaternionic Hermitian structure}
on a real manifold of dimension $4k$ is reduction
of its structure group to $\Sp(k)\cdot \Sp(1)$.
A classification of almost quaternionic Hermitian
manifolds of dimension $4k>8$ into 64 classes was obtained
by \index[persons]{Swann, A.} Swann by using representation theory, precisely by
acting with $\Sp(k)\cdot \Sp(1)$ on the
space of tensors having the same symmetries as
$\nabla\Theta$ where $\nabla$ is the Levi--Civita
connection and $\Theta$ the fundamental 4-form. 
LCHK manifolds appear in this
classification and a number of tensorial conditions
satisfied by the three fundamental forms of a
hyperhermitian structure in order to be LCHK are
given. See \cite{swann}, \cite{cab},
 \cite{cab_sw}. 
 
\hfill

Hyperk\"ahler manifolds are Ricci-flat; consequently, from \ref{lchk_vai} we obtain
(see e.g. \cite{pps}, \cite{do}):\index[terms]{manifold!Ricci-flat}

\hfill

\theorem Let $M$ be a compact LCHK manifold, that is  not
globally conformally hyperk\"ahler. Then $M$
Vaisman for an appropriate choice of the metric in the
given conformal class.

\hfill

From now on, when working on compact LCHK manifolds, we
shall tacitly assume that the Lee form \index[terms]{form!Lee}is parallel
(i.e. each LCK structure\index[terms]{structure!LCK} $(g,I_i)$ is Vaisman).

The analogue of the canonical foliation\index[terms]{foliation!canonical}
$\Sigma=\spp\{\theta^\sharp, I\theta^\sharp\}$  on LCK
manifolds is the 4-dimensional foliation $\cad$ locally
generated by $\theta^\sharp, I_1\theta^\sharp,
I_2\theta^\sharp, I_3\theta^\sharp$. Its properties are
summarized in:

\hfill 

\theorem {(\cite{op})} Let $M$ be a compact LCHK manifold,
$\dim_\R M=4k$. Then the canonical foliation $\cad$ is a Riemannian totally geodesic foliation. Its
leaves, if compact, are hypercomplex  Hopf surfaces (primary or\index[terms]{surface!Hopf!hypercomplex}
non-primary).  Its orthogonal complement $\cad^\perp$ is not
integrable. Any submanifold $X\subset M$ such that
$TX\subset\cad^\perp$,\footnote{Such submanifolds are
  called ``integral manifolds'' of the distribution
  $\cad^\perp$.} is totally real with
respect to each $I_i$ and its maximal dimension is
$k-1$. 

\hfill

\remark\label{_3_Sasaki_Einstein_} 
Recall (\cite{bog}) that a {\bf 3-Sasakian}\index[terms]{manifold!3-Sasakian} ma\-ni\-fold is a
Riemannian ma\-ni\-fold $(S,h)$ such that its Riemannian cone\index[terms]{cone!Riemannian}
$C(S)$ is hyperk\"ahler (in particular, Ca\-la\-bi-Yau). Equivalently, $M$ has three
Killing fields $\xi_2,\xi_2, \xi_3$\index[terms]{vector field!Killing} which generate a
$\mathfrak{su}(2)$ action  and such that each $(g,\xi_i)$
is Sasakian. Note that 3-Sasakian manifolds are Einstein\index[terms]{manifold!Einstein}
(\cite{kashiwada}), i.e.  each of the three Sasakian structures are Sasakian-Einstein. Then we have:

\hfill

\proposition {(\cite{op})} The  foliation $\ker\theta$ is Riemannian and totally ge\-o\-de\-sic. Its leaves are 3-Sasakian.

\hfill

By \ref{betti_1}, a compact non-K\"ahler LCHK manifold has
$b_1=1$. Then the Structure theorem for Vaisman manifolds\index[terms]{manifold!Vaisman}
implies:

\hfill

\theorem \label{_LCHK_mapping_torus_Theorem_}
{(\cite{_Verbitsky:Vanishing_LCHK_})} Let $M$ be a compact  LCHK 
manifold that is  not hyperk\"ahler. Then $M$ is a mapping torus over the circle with 3-Sasakian fiber, that is there exists\index[terms]{mapping torus}
a 3-Sasakian manifold $S$ and a\index[terms]{manifold!3-Sasakian}
3-Sasakian automorphism $\phi:\; S \arrow S$
such that $M \simeq
 C(S)/\sim_{\phi, q}$,
for some $q\in \R, q >1$,
where $\sim_{\phi, q}$ is an equivalence relation
generated by $(x,t)\sim_{\phi, q} (\phi(x), qt)$.
Moreover, the manifold $S$ and the automorphism
$\phi$ are determined uniquely from the
LCHK structure on $M$, up to a rescaling
of a  metric on $S$.

\hfill

\corollary 
{(\cite{_Verbitsky:Vanishing_LCHK_})} 
Let $M$ be a compact LCHK manifold,\index[terms]{manifold!LCHK}
obtained as above 
{}from a 3-Sasakian manifold $S$ and a 3-Sasakian\index[terms]{manifold!3-Sasakian}
automorphism $\phi:\; S \arrow S$. Then 
$\cad$ is quasi-regular if and only if 
$\phi$ is a finite order automorphism.\index[terms]{manifold!LCHK!quasi-regular}

\hfill

Next, we give a characterization of compact LCHK
manifolds, in the spirit of \index[persons]{Beauville, A.} Beauville's characterization of compact
hyperk\"ahler manifolds. Recall that the latter are equivalent to
K\"ahler, endowed with a holomorphic symplectic form, \index[terms]{form!symplectic!holomorphic}
see \cite{beauville}.\index[terms]{theorem!Beauville}

\hfill

\theorem {(\cite{oo})}\label{crit_lchk}  Let $(M, I, g, \theta)$ be a compact locally conformally K\"ahler manifold of real dimension $2n$. Then $g$ is locally conformally hyperk\" ahler if and only if there exists a non-degenerate $(2, 0)$-form $\Omega$  such that:
\begin{description}
\item[(i)] $d\Omega=\theta \wedge \Omega$, and
\item[(ii)] $\Omega^{\frac{n}{2}} \wedge \overline{\Omega}^{\frac{n}{2}}=c \cdot d\vol_g,$ where $c \in \mathbb{R}_{+}$, and $d\vol_g$ is the volume form of $g$.
\end{description}

\remark The existence of a nondegenerate $(2, 0)$-form implies that the real dimension is a multiple of 4.

\hfill

\remark 
A complex manifold admitting a nondegenerate $(2, 0)$-form $\omega$ and a closed one-form $\theta \in \Lambda^1(M, \C)$ such that $d \omega = \theta \wedge \omega$, is called {\bf complex locally conformally symplectic} (CLCS). The Lee form\index[terms]{form!Lee} $\theta$ can be real or complex. CLCS manifolds    appeared in \cite[Section 5]{lic1}, as even-dimensional leaves of the natural generalized foliation of a complex Jacobi manifold (recall that real LCS structures also appear as leaves of real Jacobi manifolds).\index[terms]{manifold!CLCS} 
\index[terms]{manifold!LCS}

\hfill

The proof of \ref{crit_lchk} is technical. It relies on a kind of ping-pong from objects on the universal cover $\pi: \tilde M\ra M$  to objects on the compact $M$ in order to can apply a Weitzenb\"ock formula. First note that 
if $\pi^*\theta=df$, then $\tilde\Omega:=\pi^*\Omega$ is holomorphic. Using (iii), $\Vert\tilde\Omega\Vert_{\tilde g}=\const$, so $\tilde g$ is Ricci-flat,\index[terms]{manifold!Ricci-flat} thus $g$ is Einstein\index[terms]{metric!Gauduchon}-Weyl, and hence  Vaisman. Then one can work with the Gauduchon metric. It remains to show 
that $\tilde\Omega$ is $\tilde g$-parallel that is  done by translating this on $M$. \endproof  

\hfill

\remark\label{_qK_Remark_}
Recall that a Riemannian manifold is called {\bf
  quaternion-K\"ahler} if its Riemannian holonomy \index[terms]{manifold!quaternion-K\"ahler}
is a subgroup of $\Sp(n)\cdot \Sp(1)$ (\cite{_Salamon_}).
When the foliation $\cad$ is quasi-regular, one obtains the following relations among LHCK, 3-Sasakian and quaternion-K\"ahler\index[terms]{manifold!3-Sasakian} geometries:

\hfill

\theorem {(\cite{op})} Let $M$ be a compact quasi-regular LCHK manifold. Then we have the following\index[terms]{manifold!LCHK!quasi-regular}
commutative diagram of fibre bundles (and Riemannian submersions) in the
orbifold category:\index[terms]{submersion!Riemannian}


\begin{equation}\begin{minipage}{0.85\linewidth}\label{lchk_diag}
\xymatrix{
&M \ar[dl]_{T^1_\C} \ar[dd] \ar[dr]^{S^1}&\\
Z\ar[dr]_{S^2}&&S\ar[ll]_{\qquad\quad S^1}\ar[dl]^{S^3\!/\Gamma}\\
&P&
}
\end{minipage}
\end{equation}
Here $P$ is a quaternion-K\"ahler orbifold\index[terms]{orbifold!quaternion-K\"ahler} with positive
scalar curvature, and the fibres of the projection
$M\rightarrow P$ are (not necessarily primary) Hopf
surfaces $S^1\times S^3/G$; $Z$ is its K\"ahler--Einstein
twistor space.\index[terms]{manifold!K\"ahler--Einstein}\index[terms]{twistor space}
Conversely, given a quaternion-K\"ahler manifold\index[terms]{manifold!quaternion-K\"ahler} $P$ of positive scalar curvature,
there exists a commutative diagram as above with manifolds $M'$,  $Z'$,  $S'$ respectively,
LCHK, K\"ahler--Einstein and 3-Sasakian and with\index[terms]{manifold!3-Sasakian}
fibres as described with $G=\Gamma=\Z/2\Z$.

\hfill

\remark There are many examples of compact 3-Sasakian
manifolds, in any dimension and with arbitrary second
Betti number (\cite{_BGP:quotients_,_BGGR:7_manifolds_,bog}.
 Therefore, by
the above constructions we obtain a wealth of compact LCHK
manifolds. In particular, the quaternionic Hopf manifolds
$S^1\times S^{4n+3}$ are always LCHK.

\hfill 

\remark\label{_fundamental_3-Sasakian_foliation_Remark_}
The analogue of the Reeb foliation on a 3-Sasakian
manifold is the 3-dimensional totally geodesic foliation
generated by the 3 Reeb fields. 
In the sequel, we call this foliation {\bf the fundamental foliation}.
Its leaves are locally isometric to a 3-sphere and complete.
By Meyers' theorem, they are
compact. Therefore, a 3-Sasakian manifold is
quasi-regular, and all the terms in the\index[terms]{manifold!3-Sasaki!quasi-regular}
diagram \eqref{lchk_diag} are globally defined
orbifolds, whenever the 3-Sasakian orbifold leaf space $S$
(that is, the leaf space of the foliation
tangent to the Lee field\index[terms]{Lee field}) exists.
However, the existence of $S$ is not a priori given,
since the Lee field action might be irregular
even on a quaternionic Hopf manifold.

\hfill

\remark Since compact LCHK manifolds are Einstein--Weyl,\index[terms]{manifold!Einstein--Weyl}\index[terms]{manifold!LCHK}
\ref{EW_vanish} and \ref{betti_1} apply. In particular,
$b_1(M)=1$ (see also \cite{op}). Moreover, from the Gysin
sequence in the above diagram one obtains (\cite{op})
$$b_{2p}(M)=b_{2p+1}(M)=b_{2p}(N)-b_{2p+4}(N),\quad (0\leq
2p\leq 2k-2),$$
when the fundamental foliations on the LCHK manifold are regular.
In this case, $b_{2k}=0$. Moreover, for $k\geq 1$, (\cite{gs}):
\begin{description}
\item[(i)] $b_2(M)=b_3(M)=0$ unless $M=S^1\times \U(n+1)/(\U(n-1)\times \U(1))$;
\item[(ii)] if $b_4(M)=0$ and $\dim M=16$ or $20$, then $M$ is either $S^1\times S^{4k-1}$ or one of the two flat circle bundles over $\R P^{4k-1}$;
\item[(iii)] $\displaystyle\sum_{i=1}^{k-1}i(k-i+1)(k-2i+1)b_{2i}(M)=0.$
\end{description}

\remark A compact, homogeneous LCHK manifold is
regular (compare with \ref{homo_reg}). Therefore, all terms in
the diagram \eqref{lchk_diag}\index[terms]{manifold!LCHK!homogeneous} 
are smooth. Moreover, all manifolds involved in it will be
homogeneous, too. Using the classification of compact,
homogeneous 3-Sasakian
manifolds,\index[terms]{manifold!3-Sasakian!homogeneous}
\cite{bog}, we obtain:\index[terms]{manifold!3-Sasakian}
 
 \hfill
 
\proposition {(\cite{op})} A compact homogeneous  LCHK manifold is one of the following:
\begin{description}
\item[(i)] A product $G/H \times S^1$, where $G/H$ can be any of the following  homogeneous 3-Sasakian manifolds: $S^{4k-1}$, $\R P^{4k-1}$, $\SU(m)/\mathrm{\sf  S}(\U(m-2)\times\U(1))$, $m\geq 3$, $\SO(m)/(\SO(m-4)\times \Sp(1))$, $m\geq 7$, $G_2/\Sp(1)$, $F_4/\Sp(3)$, $E_6/\SU(6)$, $E_7/\Spin_{12}$, $E_8/E_7$.
\item[(ii)] The unique non-trivial $S^1$ bundle over $\R P^{4k-1}$.
\end{description}

\remark LCHK geometry \index[terms]{geometry!LCHK}is also related to HKT (hyperk\"ahler with torsion) geometry. It is proven in \cite{ops} that if $(g,I_1, I_2, I_3)$ is an LCHK structure with parallel Lee form\index[terms]{form!Lee!parallel} $\theta$, then the metric
$$g_0=g-\frac 12\big\{\theta\otimes\theta + I_1\theta\otimes I_1\theta + I_2\theta\otimes I_2\theta + I_3\theta\otimes I_3\theta\big\}$$
is HKT.\index[terms]{manifold!HKT} Conversely, among the HKT manifolds, the LCHK ones
are identified by the existence of certain symmetries 
(\cite[Theorem 19]{ops}).
The first examples of compact HKT-manifolds that are  neither homogeneous nor LCHK were obtained in \cite{_Verbitsky:HKT_}. \index[terms]{manifold!HKT}


Recall that a {\bf quaternionic-K\"ahler manifold}\index[terms]{manifold!quaternionic-K\"ahler}
 (\cite{_Salamon_,besse}) is a $4n$-dimensional 
Riemannian manifold with holonomy in $\Sp(1)\cdot \Sp(n)$.
The action of $\Sp(1)\cdot \Sp(n)= \frac{\Sp(1)\times
  \Sp(n)}{\pm 1}$ in $\R^{4n}={\Bbb H}^n$ is by
$\Sp(1) \subset {\Bbb H}$ from the right and
by $\Sp(n) = \U({\Bbb H}, n)$ from the left.
A quaternionic-K\"ahler manifold $P$ is Einstein
with Einstein constant $c$, that vanishes if and only
if $P$ is a finite quotient of a hyperk\"ahler
manifold by isometries. Such a manifold is called
{\bf locally hyperk\"ahler}. In the literature, the name
``quaternionic-K\"ahler'' is usually reserved for
manifolds with holonomy in $\Sp(1)\cdot \Sp(n)$
and non-zero Ricci curvature, however, in the
sequel we include the locally hyperk\"ahler case
as a special case of quaternionic-K\"ahler geometry.\index[terms]{geometry!quaternionic-K\"ahler}

\hfill

\remark
Let $P$ be a quaternionic K\"ahler orbifold of
positive Ricci curvature, and $U$ the principal
$\SO(3)=\Sp(1)/\pm 1$-bundle associated with the
$\Sp(1)$-factor of $\Sp(1)\cdot \Sp(n)$.
Then the total space of $U$ is 3-Sasakian.
Since the fundamental 3-dimensional
foliation on a 3-Sasakian manifold is
quasi-regular\index[terms]{manifold!3-Sasaki!quasi-regular}
(\ref{_fundamental_3-Sasakian_foliation_Remark_}), 
every 3-Sasakian manifold
is obtained from this construction.
If the Ricci curvature is negative,
$U$ is equipped with a pseudo-Riemannian
version of a 3-Sasakian geometry\index[terms]{geometry!Sasaki}, and
if it vanishes, then $U$ is locally 
isometric to a product $Q \times \R P^3$.

\hfill

\remark  {\bf A locally conformally quaternionic-K\"ahler}
(LCQK) ma\-ni\-fold can also be  defined, see \cite{op}, by
starting with a quaternion-Hermitian ma\-ni\-fold and asking
that the quaternion-Hermitian metric be locally conformal
with a qu\-a\-ter\-nion-K\"ahler metric. It can be seen that
this is equivalent to the same equation \eqref{lcqk} as
in the LCHK case (recall than on quaternionic-K\"ahler
manifolds, the four-form $\Theta$ is globally defined,
even though the two-forms 
$\omega_i$ are just local). Being locally conformally
quaternionic-K\"ahler implies that the metric is
Einstein--Weyl.\index[terms]{manifold!Einstein--Weyl} Moreover, the
scalar curvature of the Weyl connection should
vanish. This is a local fact, observed directly in
\cite[Theorem 2.1]{op}. An alternative argument is as
follows: on the universal cover the lifted metric $\tilde
g$ is quaternionic-K\"ahler, and hence  Einstein, with constant
scalar curvature $\lambda$. If $\lambda\neq 0$, all
homotheties in the deck group preserve the metric 
$\tilde g'=\pm\lambda\tilde{\Ric}$, and thus they are
isometries with respect to $g'$. This is a
contradiction. In conclusion, LCQK manifolds are locally
LCHK, that is the only examples are finite quotients of
LCHK manifolds.

\section{Locally conformally balanced structures}\label{lcb}

\definition A Hermitian metric on a com\-plex manifold $(M,I)$ of com\-plex dimension $n$ is called {\bf balanced}\footnote{Note that the notion we work with is different from the
	one used by S. \index[persons]{Donaldson, S. K.} Donaldson in \cite{don}.} if $d\omega^{n-1}=0$.\index[terms]{metric!balanced}

\hfill

\remark Clearly, on surfaces all balanced metrics are K\"ahler, and thus the notion is consistent for $n\geq 3$. 

\hfill

\remark Balanced metrics are also called semi-K\"ahler. They appear as class
$W_3$ in \cite{gh}. Balanced manifolds were intensively studied from the
geometric and analytic  viewpoint by Michelson,\index[terms]{manifold!semi-K\"ahler}
\index[persons]{Alessandrini, L.} Alessandrini-Bassanelli etc 
(\cite{_Michelson_,_Alessandrini_Bassanelli:bimero_,_AB:modifications_balanced_}, see also \cite[Chapter 4]{_Popovici:book_}). 
It is known that all
Riemannian twistor spaces admit balanced metrics,  and
the twistor spaces of hypercomplex manifolds are balanced
(\cite{_Tomberg_}). On the other hand, the
\index[persons]{Calabi, E.} Calabi-Eckmann manifolds are not balanced: indeed, a
complex manifold $M$ that contains a divisor homologuous
with $0$ cannot be
balanced. \index[terms]{manifold!Calabi-Eckmann}\index[terms]{divisor}


\hfill

\definition A Hermitian manifold that admits an open cover such that on each open set of the cover the restriction of the metric is conformal to a local balanced metric is called {\bf locally conformally balanced} (LCB for short).\index[terms]{manifold!locally conformally balanced}

If the cover has only one open set, the manifold is called {\bf globally conformally balanced} (GCB).

\hfill

\proposition (\cite{au}, \cite{lya})A complex manifold $(M,I)$ of complex dimension $n\geq 3$ is LCB if and only if its Lee form\index[terms]{form!Lee} $\theta:=-\frac{1}{n-1}Id^*\omega$ is closed.
$(M,I)$ is GCB if and only if $\theta$ is exact. \endproof

\hfill

\remark 
\begin{description}
	\item[(i)] The Lee form\index[terms]{form!Lee} of any Hermitian manifold satisfies  $d\omega^{n-1}=(n-1)\theta\wedge \omega^{n-1}$, thus LCK structures\index[terms]{structure!LCK}  are LCB. On the other hand, on surfaces the two notions coincide.
\item[(ii)] The LCB condition is conformally invariant.
\item[(iii)] Unlike the LCK case, if $(M_i, g_i,I_i)$, $i=1,2$, is\index[terms]{manifold!LCB}
LCB with (closed) Lee form\index[terms]{form!Lee} $\theta_i$, then the product
$M_1\times M_2$ is LCB with respect to the product
Hermitian structure. This gives 
a large class of examples of compact LCB
manifolds that are  not LCK: products of K\"ahler
manifolds of dimension $>1$ and LCK manifolds are LCB
manifolds that do not
admit an LCK structure\index[terms]{structure!LCK} (\ref{no_kaehler_fibres}).
\item[(iv)] An example of a compact complex manifold that
  cannot admit any LCB structure is obtained by
  H. \index[persons]{Shimobe, H.} Shimobe in \cite[Proposition 3.4]{shi2}. 
Let $X$ be an \index[terms]{surface!Inoue} Inoue surface that cannot
  admit LCK structures (see \ref{_no_lck_}), and
  $M:=X\times X$. Shimobe shows, using \index[persons]{Belgun, F. A.} Belgun's
  technique, that if $M$ admitted an LCB structure, then
  it admits an invariant LCB structure. Then he shows that
  this implies the vanishing of the Lee form,\index[terms]{form!Lee} and hence
  the structure would be balanced. However, a balanced\index[terms]{structure!balanced}
  structure on the product would project in a balanced,
  thus K\"ahler structure on the Inoue surface, that is 
  a contradiction. In particular, as in the LCK case, the
  LCB class is not stable at small deformations.
\item[(v)] Recall that a blow-up of an LCK manifold along a\index[terms]{blow-up}
submanifold is LCK only if the restriction of Lee form\index[terms]{form!Lee} to
the submanifold is exact (\ref{bdown}). By contrast, if $\tilde X\arrow
X$ is a proper modification of compact manifolds, then $X$ 
is LCB if and only if $\tilde X$ is LCB (\cite[Theorem 6.9]{shi},  \cite[Theorem 4.1]{shi2}).\footnote{Prior to these results, it was proven that the blow-up at points and the blow-up along compact submanifolds with induced GCB structure preserve the LCB class (\cite{lya}, \cite{y}).}
\end{description}

\example {(\cite{daniele})} Let $\mathfrak{g}$ be a 6-dimensional 
nilpotent Lie algebra whose generators $e^1,\ldots, e^6$ obey the relations:
$$de^1=\cdots=de^5=0,\quad de^6=e^{12}:=e^1\wedge e^2.$$
Define the complex structure $I$ by
$$Ie^{2k-1}=e^{2k}, \quad Ie^{2k}=-e^{2k-1}, \quad k=1,2,3,$$
and the two-form
$$\omega:=e^{12}+e^{34}+e^{56}.$$
Let $G$ be the corresponding simply connected nilpotent Lie group.
Then we have:
\begin{equation*}
\begin{split}
\frac 12\omega^2&=e^{1234}+e^{1256}+e^{3456},\\
d\omega^2&=e^{12345}=\omega^2\wedge e^5.
\end{split}
\end{equation*}
The pair $(I, \omega)$ defines
a left-invariant  Hermitian
structure on $G$, clearly integrable
by Newlander--Nirenberg theorem.
Since  $de^5=0$, the manifold $(G,I,\omega)$ is LCB with $\theta=e^5.$
Moreover, as the relations are rational, by \index[persons]{Maltsev, A.} Maltsev there
exists a cocompact lattice \index[terms]{lattice!cocompact}$\Gamma$ such that the quotient
$G/\Gamma$ is a compact LCB manifold.

\hfill

The analogue of Vaisman's \ref{vailcknotk} (and which also implies  \ref{vai_gen_au}) is:

\hfill

\theorem {(\cite{au})} Let $(M,I)$ be a $2n$-dimensional compact, complex  manifold. Suppose that the natural map 
\begin{equation}\label{bcinj}
H^{n-1,n}_{BC}(M)\arrow H^{n-1,n}_{\bar\6}(M),
\end{equation}
 induced by the identity, is injective. Then any LCB structure on $(M,I)$ is GCB.

\hfill

\proof  By Serre duality, from \eqref{bcinj} one obtains:\index[terms]{duality!Serre}
\begin{equation}\label{aesurj}
H^{1,0}_{\bar\6}(M)\rightarrow H^{1,0}_{\Ae}(M),
\end{equation}
where $H^{*,*}_{\Ae}$ represents the  Aeppli cohomology.\index[terms]{cohomology!Aeppli} Clearly, if \eqref{bcinj} is injective, then \eqref{aesurj} is surjective. By dimensional reason, if \eqref{aesurj} is surjective, it is in fact an isomorphism.

Now, if $\omega$ is LCB with Lee form\index[terms]{form!Lee} $\theta$, then $\theta\wedge\omega^{n-1}=d\omega^{n-1}$. If \eqref{aesurj} is an isomorphism, then $\theta$ can be chosen to be $\6\bar\6$-closed and  holomorphic (because representatives of  $H^{1,0}_{\bar \6}(M)$ are holomorphic). Then:
$\6\bar \6 \omega^{n-1}=\theta^{1,0}\wedge\theta^{0,1}\wedge \omega^{n-1}$ is exact and positive, and thus $\theta=0$. \endproof

\section{Locally conformally parallel 
$G_2$, $\Spin(7)$ and $\Spin(9)$ structures}

Recall  the definitions of the special holonomies
$G_2$ (in dimension $7$) and $\Spin(7)$ (in dimension
$8$). We follow the conventions in \cite{ipp},
\cite{oppv}, \cite{_V:intrinsic_volume_}.

Let $e_1,\ldots,e_7$ be a oriented orthonormal basis of
$\R^7$, $e^1,\ldots, e^7$ its dual basis and denote, as
above, by $e^{i_1\ldots i_k}:=e^{i_1}\wedge\cdots\wedge
e^{i_k}$. Then $G_2$ can be defined as the compact Lie
subgroup of $\SO(7)$ that fixes the three-form
$\phi$ on $\R^7$, identified with the imaginary part of the octonion algebra,
$\phi(x, y, z) = g(z, x \cdot y - y \cdot x)$,
where $x\cdot y$ denotes the octonion product, and $g$ is
the invariant
scalar product on the octonion algebra. In the appropriate coordinates,
the form $\phi$ can be written as 
$$\phi=e^{145}+e^{167}+e^{246}+e^{257}+e^{347}-e^{356}+e^{123}.$$
A {\bf $G_2$-structure} on a  Riemannian manifold $(M^7,g)$ is 
a three-form $\f$ that is  given
at each tangent space by the above expression in
an appropriate orthonormal frame. 
Equivalently, the structure group of 
the tangent bundle $TM$ is reduced to $G_2$. 

Similarly, on $\R^8=\R\oplus \R^7$ with oriented 
orthonormal basis $e_0,\ldots,e_7$, consider the four-form 
$$\Phi=e^0\wedge\f+*_7\f,$$
where $*_7$ is the Hodge operator on $\R^7$. Then
$\Spin(7)$ can be defined as the compact Lie subgroup of
$\SO(8)$ that fixes $\Phi$. {\bf  A $\Spin(7)$-structure}
on a Riemannian manifold $(M^8,g)$ is a 
four-form $\Phi$ that is  given at each tangent space by the above
expression. Equivalently, the structure group of the
tangent bundle $TM$ is reduced to $\Spin(7)$.

Similarly, $\Spin(9)$ can be viewed as a subgroup of
$\SO(16)$ fixing a certain eight-form whose expression is
not important here. {\bf A $\Spin(9)$-structure} on a
Riemannian manifold $(M^{16},g)$ can be defined as
reduction of the structural group to $\Spin(9)$ or,
equivalently (see \cite{fr}), can be given in terms of
existence of $9$ self-adjoint endomorphisms $I_\al$ of the
tangent bundle satisfying the identities:
$$I_\al^2=\id,\quad I_\al I_\beta=-I_\beta I_\al \, (\al\neq\be).$$

Since the groups $G_2$, $\Spin(7)$ and $\Spin(9)$ are
compact, the $G_2$-, $\Spin(7)$- and $\Spin(9)$-structures,
as defined above, are all compatible with a Riemannian
metric. In fact there is a canonical choice of a
Riemannian metric explicitly expressed through the
corresponding 3-, 4- and 8-forms. We say that
a $G_2$-, $\Spin(7)$- or $\Spin(9)$-manifold is
{\bf holonomy $G_2$-, $\Spin(7)$- or $\Spin(9)$-manifold}
if this 3-, 4- and 8-form is preserved by the
corresponding Levi--Civita
connection.\index[terms]{group!holonomy}

Originally the \index[persons]{Berger, M.} Berger's list included all groups
that can act transitively on a sphere (\cite{besse}),
and this includes $\Spin(9)$. Later, D. \index[persons]{Alekseevsky, D. V.} Alekseevsky
proved that all manifolds with holonomy $\Spin(9)$ are
locally symmetric (\cite{_Alekseevsky:Spin9_});
this is why the $\Spin(9)$ holonomy manifolds are omitted
from the contemporary version of the Berger's list of
holonomies.

Note that the manifolds with holonomies  
$G_2$, $\Spin(7)$ and $\Spin(9)$  
are Einstein, \cite{besse}, \cite{fr}.\index[terms]{manifold!Einstein}

Now, if $G$ is any of the above three groups, we can
consider  {\bf a locally conformally parallel $G$-structure}
on a given manifold to be a Riemannian metric locally
conformal with (local) metrics with holonomy contained in
$G$, see \cite{oppv}. We remark that this definition also
applies to LCK  (with $G=\U(n)$), LCHK (with
$G=\Sp(n)$),  but not to LCB structures.

In each of these three cases, a locally conformally
parallel $G$ structure gives rise to a closed Lee form\index[terms]{form!Lee}
$\theta$ such that the specific differential form of the
$G$-structure is $d_\theta$-closed.

Let $\psi:\; S \arrow S$ be an isometry of a Riemannian
manifold. Recall that {\bf the mapping torus} of $\psi$
is a quotient of a product $[0,1]\times S$ by 
the equivalence relation $(0, s) \sim (1, \psi(s))$.
This is a manifold locally isometric to the product
of $S$ and an interval, and equipped with a natural
projection to $S^1$, with each fibre isometric to $S$.

Similar to Sasakian manifold, we define
{\bf a nearly K\"ahler structure} on a manifold $S$
as a holonomy $G_2$ structure on its Riemannian cone,\index[terms]{cone!Riemannian}
defined in such a way that the natural homothety
of the Riemannian cone multiplies the structure
3-form by a constant. We define {\bf a nearly parallel
  $G_2$ structure} on $S$ as a holonomy $\Spin(7)$\index[terms]{structure!nearly parallel $G_2$}
structure on its Riemannian cone,
defined in such a way that the natural homothety
of the Riemannian cone multiplies the structure
4-form by a constant.

The structure theorem for compact Vaisman manifolds is the
model for the following results:\index[terms]{manifold!Vaisman}

\hfill

\theorem\label{spini} 
\begin{description}
	\item[(i)]{(\cite{ipp}, \cite{_V:intrinsic_volume_})} 
A compact Riemannian 7-manifold $M$ admits a locally conformally parallel $G_2$ 
structure if and only if it is isometric to a mapping
torus of an automorphism of a nearly K\"ahler manifold $S$.
Moreover, the natural $\Z$-covering\index[terms]{cover!K\"ahler $\Z$-} of $M$
is conformal to the Riemannian cone\index[terms]{cone!Riemannian} of $S$.

\item[(ii)]{ (\cite{ipp})} 
A compact Riemannian 8-manifold $M$ admits a locally conformally parallel
$\Spin(7)$ structure if and only if it is isometric to a mapping
torus of an automorphism of a nearly parallel
$G_2$-manifold $S$. Moreover, the natural $\Z$-covering of $M$
is conformal to the Riemannian cone of $S$.

\item[(iii)]{(\cite{oppv})} 
A compact Riemannian 16-manifold $M$ admits a locally conformally parallel
$\Spin(9)$ structure if and only if $M$ is isometric to a mapping torus
of an automorphism of 
a finite quotient $S^{15}/\Gamma$ of a 15-sphere. 
Moreover, the universal cover of $M$ is
conformally equivalent to $\R^{16}\setminus \{0\}$ and is
finitely covered by $S^{15}\times\R$.
\end{description}

\hfill

\remark Note that (i) and (ii) above implicitly provide
examples of locally conformally $G_2$ and $\Spin(9)$
structures. A unified discussion about the mapping tori
that appear in (i) and (ii) of \ref{spini}, in the context of
spinor geometry, can be found in \cite{af}, where locally
conformally parallel $\Spin(9)$ structures appear as
``$\Spin(9)$-structures of vectorial type'', and an
example is constructed on $S^{15}\times S^1$. Other examples of  locally conformally parallel $\Spin(9)$ structures are provided in \cite{oppv} by finding finite subgroups of $\Spin(9)$ acting freely on $S^{15}$. 

\section{Notes} \label{_SKT_Notes_}
\begin{enumerate}


\item In almost contact geometry\index[terms]{geometry!almost contact}, a great deal of work was 
done on  {\bf locally conformally cosymplectic}
structures, see \cite[Section 17.6]{do}, \cite{matzeu}
(where LCS manifolds are related to Einstein--Weyl\index[terms]{manifold!LCS} almost\index[terms]{manifold!Einstein--Weyl!almost contact}
contact manifolds), and \cite{falc} for some curvature
properties etc.  This structure gains importance due to
recent results on the geometry and topology of\index[terms]{manifold!cosymplectic}
cosymplectic (sometimes called ``co-K\"ahler'') manifolds
by G. \index[persons]{Bazzoni, G.} Bazzoni, J. C. \index[persons]{Marrero, J. C.} Marrero, J. \index[persons]{Oprea, J.} Oprea \textit{et. al.} On the
contrary,  a notion of locally conformally Sasakian
structure can not exist, see \cite{va_contact}.\index[terms]{manifold!locally conformally cosymplectic}

\item A version of LCK structure\index[terms]{structure!LCK} is also studied in
generalized complex geometry, see \cite{aisa},
\cite{va_gc}. 
\end{enumerate}

\chapter{Open questions}

\epigraph{\it What is important is to spread confusion, not eliminate it.}{\sc\scriptsize Salvador Dal\'i}

\epigraph{\it Why do people always expect authors to answer questions? I am an author because I want to ask questions. If I had answers, I'd be a politician.} {\sc\scriptsize Eug\`ene Ionesco}


In the sequel, ``strict LCK structure'' means
an LCK structure \index[terms]{structure!LCK}that is  not globally conformally K\"ahler,
that is, with non-trivial weight bundle. Most of the questions
asked about the Vaisman manifolds\index[terms]{manifold!Vaisman} are referred to compact
Vaisman manifolds; sometimes we omit ``compact''. 
An LCK manifold with potential\index[terms]{manifold!LCK!with potential} is always assumed
to be compact.


\section{Existence of LCK structures}\index[terms]{structure!LCK}


\subsection{LCK structures on complex manifolds}
 
Let $(M,I)$ be an almost complex manifold, $\dim_\C M >2$.
S.-T. \index[persons]{Yau, S.-T.} Yau conjectured
(\cite[Problem 2018001]{_Liu_Xu:open_problems_})\index[terms]{conjecture!Yau} 
that $I$ can be deformed to an integrable complex structure. 
In dimension 2, this is clearly false: there are
almost complex compact 4-dimensional manifolds
not admitting complex structures (\cite{_Van_de_Ven:surfaces_}),
even the ones with the tangent bundle 
that is  topologically trivial (\cite{_Brotherton:parallelizable_}).

Every year there are several publications 
appearing in \texttt{\scriptsize arXiv.org}
aimed to prove the existence or non-existence of a complex structure
on $S^6$; obviously, this is a question of great importance and
difficulty.

For LCK manifolds, it seems that even the following
question is open. 

\hfill

\question
In \cite{em}, \index[persons]{Eliashberg, Y.} Eliashberg and \index[persons]{Murphy, E.} Murphy proved that any
compact almost complex manifold $M$ with $b_1(M)>0$ admits a LCS structure (see also \cite{_Bertelson_Meigniez_} for a stronger result). Suppose that $M$ is a compact manifold,
 $\dim_\R M > 4$, $b_1(M)>0$ and
$M$ is almost complex.
Will $M$ always admit an LCK structure?\index[terms]{structure!LCK}

\hfill

One can ask many questions about existence or non-existence
of LCK metrics on a given complex manifold.

\hfill

\remark
Non-K\"ahler complex structures on compact tori were 
obtained by A. \index[persons]{Sommese, A. J.} Sommese in 
\cite[52, Remark III, E)]{_Sommese:Q_}; 
see also \cite{_Catanese:types_}. Let $M$ be
the $4n$-dimensional torus $T^{4n}$ equipped
with a standard hyperk\"ahler structure, $M = {\Bbb H}^n/\Z^{4n}$,
and $\Tw(M)$ its {\em twistor space}. Recall that {\bf the
twistor space} of a hypercomplex manifold $(M, I, J, K)$
is a product $M \times S^2$, equipped with a 
complex structure as follows. Given $m\in M, s\in \C P^1=S^2$,
let $L= aI + bJ + c K$ be the complex structure associated
with $s= (a, b, c) \in S^2$, and let ${\cal I}$ be
a complex structure operator acting on\index[terms]{twistor space}
$T_{(m, s)}M \times S^2= T_m M \times T_s \C P^1$ by $L$ on $T_m M$
and as the standard complex structure operator on $T_s \C P^1$.
The almost complex structure ${\cal I}$
constructed this way is integrable
(\cite{_Salamon_}, \cite[Theorem 14.68]{besse},
\cite{_Kaledin:twistor_}). Now let  $\sigma:\; E \arrow \C P^1$
be a 2-sheeted ramified covering, with $E$ being an elliptic curve,
and $X:= E\times_{\C P^1} \Tw(M)$
be a fibred product of $E$ and $\Tw(M)$ over $\C P^1$.
Since the definition of the fibred product can be given in
the category of topological spaces,
and topologically, $\Tw(M)= \C P^1 \times T^{4n}$,
the space $X$ is diffeomorphic to the product 
$E \times T^{4n}\cong T^{4n+2}$. However,
$X$ is not K\"ahler. To see that,
consider the product Riemannian form
on $\Tw(M) = M \times \C P^1$. It is Hermitian with respect to the
complex structure on $\Tw(M)$; let
$\omega_H$ be the corresponding Hermitian form.\index[terms]{form!Hermitian}
Denote by $\pi$ the standard holomorphic projection
$\pi:\; \Tw(M) \arrow \C P^1$, and let $\sigma:\; X \arrow \Tw(M)$
be the two-sheeted ramified covering constructed above, 
that is  also holomorphic.
In \cite[(8.12)]{_NHYM_}, it was shown
that $dd^c\omega_H= \omega_H \wedge \pi^* \omega_{\C P^1}$,
where $\omega_{\C P^1}$ is the Fubini--Study\index[terms]{form!Fubini--Study} form on $\C P^1$.
Then $\sigma^*(dd^c\omega_H)$ is an non-zero exact positive
(2,2)-form on $X$; existence of such a form is incompatible
with existence of a K\"ahler structure.

\hfill

\question
Let $M$ be a compact manifold diffeomorphic 
to a torus, and equip\-ped with a non-K\"ahler
complex structure. Prove that $M$ does not admit a strict LCK metric,
or find one.

\hfill

This question is a special case of the following,
due to I. Vaisman.

\hfill

\question (I. \index[persons]{Vaisman, I.} Vaisman, \cite[p. 283]{_Vaisman:_PAMS_79_}) 
Any known compact LCK manifold which
is homotopy equivalent to a compact K\"ahler manifold
is in fact globally conformally K\"ahler. 
Can a complex manifold that is  not of K\"ahler type
but homotopy equivalent to a K\"ahler one admit
a strict LCK structure?\index[terms]{structure!LCK}

\hfill

\question
Another question about exotic complex structures on
a compact torus was asked by F. \index[persons]{Catanese, F.} Catanese in 
\cite[question after Remark 4.9]{_Catanese:types_}.
In all known examples of non-K\"ahler complex structures on
a compact torus, the canonical bundle is always
non-trivial. \index[persons]{Catanese, F.} Catanese asked whether this is always
true, but the question looks difficult.
Will it be simplified if we also assume that
the complex structure is LCK? Note that 
Vaisman manifolds\index[terms]{manifold!Vaisman} with trivial canonical
bundle are classified, see Exercise\index[terms]{bundle!vector bundle!canonical}
\ref{_trivial_canonical_then_QR_Exercise_}.

\hfill

\question
Can a product 
of two compact complex manifolds admit
a strict LCK metric? 
Several results are now
available which give a negative answer when one of the
factors has special properties (see Subsection
\ref{noprod}).

\hfill

\question A manifold is called {\bf \index[persons]{Kobayashi, S.} Kobayashi hyperbolic} if
its Kobayashi pseudometric (Exercise \ref{_Kobayashi_pseudometric_Exercise_})
is non-degenerate. Kobayashi has conjectured that any Kobayashi\index[terms]{conjecture!Kobayashi (hyperbolic manifolds)}
hyperbolic manifold has ample canonical bundle\index[terms]{manifold!Kobayashi hyperbolic} 
(\cite{_Lang:hyperbolic_,_Chen_Yang:negative_}); this would
imply that it is projective. This conjecture seems to be difficult.
In \cite{_Diverio:Kobayashi_MO_}, S. \index[persons]{Diverio, S.} Diverio asked whether a
compact Kobayashi hyperbolic manifold is K\"ahler.\index[terms]{conjecture!GSS}
The answer is positive for surfaces, if one assumes
the GSS conjecture. We expect that no strict LCK manifold
is Kobayashi hyperbolic. All known examples of LCK manifolds either
contain a curve of genus $\leq 1$, or have non-zero holomorphic
vector fields; such manifolds cannot be \index[persons]{Kobayashi, S.} Kobayashi hyperbolic.

\hfill

\question 
Using earlier results of
I. \index[persons]{Vaisman, I.} Vaisman, F. \index[persons]{Belgun, F. A.} Belgun has shown that some of the Inoue surfaces of class
$S^+$ cannot be LCK; this implies that the LCK condition is not open
(Section \ref{_LCK_pot_Intro_Section_}).
Is the LCK condition closed in the deformation families? 
In other words, let $D$ be a continuous family
of complex structures on $M$, and $D_0\subset D$
the set of $I\in D$ that admit an LCK metric.\index[terms]{surface!Inoue}
Is $D_0$ closed? 

\hfill

\question\label{_LCK_form_exact_Question_}
Most examples of LCS manifolds which occur in symplectic
geometry have Morse--Novikov exact LCS form $\omega$. 
The LCS structures obtained by \index[persons]{Eliashberg, Y.} Eliashberg and \index[persons]{Murphy, E.} Murphy in \cite{em} have\index[terms]{manifold!LCS}\index[terms]{form!LCS}
Morse--Novikov exact LCS form; the LCS structures obtained by 
\index[persons]{Bertelson, M.} Bertelson and \index[persons]{Meigniez, G.} Meigniez in \cite{_Bertelson_Meigniez_}
have arbitrary LCS form, but this includes Morse--Novikov exact as well.
A pullback of a Morse--Novikov exact LCS form $\omega$ to the universal
covering of the base manifold $M$ is exact; however, there
are non-Morse--Novikov-exact LCS forms that are  pulled back to 
Morse--Novikov exact
forms on the universal cover. Such are the LCK forms on
the OT-manifolds and the \index[terms]{surface!Inoue} Inoue surfaces (\cite[Proposition
4.14, Proposition 4.16]{oti2}). \index[terms]{manifold!Oeljeklaus--Toma (OT)}However, on \index[persons]{Kato, Ma.} Kato manifolds
and the blow-ups of LCK manifolds, the LCK form is never
Morse--Novikov exact on the universal covering.
Indeed, if a manifold $M$ contains
a rational curve, the LCK form can never be Morse--Novikov exact.
We are interested in LCK manifolds $(M, \theta, \omega)$
such that the corresponding K\"ahler 
form $\tilde \omega$ is exact on the universal cover of
$M$.  Clearly, the LCK manifolds with potential \index[terms]{manifold!LCK!with potential}
(and Vaisman manifolds)\index[terms]{manifold!Vaisman} are of this type; however,
the OT-manifolds and Inoue surfaces\index[terms]{surface!Inoue} are also of this
type, and they are not LCK with potential.
Suppose that $M$ is a compact LCK manifold without rational
curves. Does $M$ always admit an LCK form $\omega$
such that $\tilde \omega$ is exact on the universal cover?

\subsection{Existence of LCK potential and Vaisman structures}\index[terms]{potential!LCK}

\question
Let $(M, \theta, \omega)$ be a compact LCK manifold
such that $\omega= d_\theta d^c_\theta \phi$,
that is, $\omega$ has a $d_\theta d^c_\theta$-potential.
In \ref{_main_positive_Theorem_} we prove that
$M$ admits an LCK metric $\omega'$ with LCK potential.\index[terms]{potential!LCK}
It is not clear whether $\omega$ itself has a positive
$d_\theta d^c_\theta$-potential; we expect that there
exist counterexamples.

\hfill

\question \label{_MN_vanishing-then-pote_Question_}
Let $(M,\omega,\theta)$ be a compact LCK
manifold and $L$ its weight local system (Section 
\ref{_MN_BC_Intro_Section_}). \index[terms]{class!Bott--Chern!twisted}
If the twisted Bott--Chern class of $\omega$ vanishes,
the manifold $(M,\omega,\theta)$ is LCK with potential\index[terms]{manifold!LCK!with potential}
(\ref{vanishing_BC_class_implies_potential_Theorem_}).
This condition clearly implies the vanishing of the Morse--Novikov class \index[terms]{class!Morse--Novikov}
$[\omega]_{M\!N}\in H^2(M,L)$, where $H^*(M,L)$ denotes the
cohomology of $L$, considered to be  a local system (see Chapter \ref{_MN_and_BC_Chapter_}).\index[terms]{Mores-Novikov class}
Is it true that the vanishing of $[\omega]_{M\!N}\in H^2(M,L)$ implies
the existence of an LCK potential\index[terms]{potential!LCK} on $M$?\index[terms]{bundle!weight}

\hfill

\question Let $M$ be 
a nilmanifold\index[terms]{nilmanifold} admitting an LCK structure\index[terms]{structure!LCK}
(not necessarily locally invariant, see 
Chapter \ref{_nil_and_solv_}).\index[terms]{nilmanifold}
The only examples known are diffeomorphic
to a mapping torus\index[terms]{mapping torus} over the circle
with fibre the Sasakian Heisenberg
group\index[terms]{group!Heisenberg}.
It is still unknown if other examples exist
(\cite{uga}). Since the Morse--Novikov cohomology
of any non-trivial local system on $M$ vanishes
(\cite{_Dixmier_}), the Morse--Novikov class\index[terms]{class!Morse--Novikov} $[\omega]\in H^2(M, L)$
of $M$ vanishes. If the conjecture stated in  
\ref{_MN_vanishing-then-pote_Question_}
is true, this should imply that $M$ is LCK with potential.\index[terms]{manifold!LCK!with potential}
We expect that in this case, $M$ is a Heisenberg
nilmanifold (\ref{_local_homo_nilma_LCK_pot_Question_}).

\hfill

\question
Let $M$ be a compact homogeneous complex manifold
admitting an LCK metric. Will $M$\index[terms]{manifold!complex!homogeneous}
admit a $G$-invariant LCK structure\index[terms]{structure!LCK}, for
some Lie group $G$ acting on $M$ holomorphically
and transitively? Note that by \ref{_homo_admits_Vais_Proposition_},
any homogeneous LCK manifold is of Vaisman type.

\hfill
 
%

\question
Let $G$ be a finite group acting freely and holomorphically on a compact LCK
manifold $M$. When $M$ is Vaisman, the quotient $M/G$
is also Vaisman; if $M$ is LCK with potential, $M/G$
is LCK with potential\index[terms]{manifold!LCK!with potential} (Exercise \ref{_quotients_Vaisman_Extercise_}). 
For a general LCK manifold, the answer
is not immediately clear. It is not hard to see that $M/G$ admits an LCK structure\index[terms]{class!Lee}\index[terms]{structure!LCK}
if $G$ preserves the Lee class.\index[terms]{class!Lee} Does the quotient $M/G$ admit an LCK
structure without this assumption? 

\hfill

\question
In \cite{_angella_pediconi_}, \index[persons]{Angella, D.} Angella and \index[persons]{Pediconi, F.} Pediconi
consider an invariant LCK metric on certain cohomogeneity
one manifolds (originally constructed by \index[persons]{B\'erard-Bergery, L.} B\'erard-Bergery,
\cite{_bergery_}, who has found a cohomogeneity 1
Einstein metric on some of these manifolds).\index[terms]{manifold!cohomogeneity one}
They prove that all compact B\'erard-Bergery
cohomogeneity 1 manifolds without singular orbits are
strict LCK. These manifolds are $S^1$-bundles over symmetric
spaces; the LCK metric is known to be non-Vaisman
for non-homogeneous B\'erard-Bergery manifolds 
 (\cite[Theorem A]{_angella_pediconi_}).
Is it always, or in some cases, LCK  with potential?\index[terms]{manifold!LCK!with potential}
We expect the answer to be positive.\index[terms]{manifold!homogeneous}


\section{Complex geometry of LCK manifolds}


\subsection{Hodge theory on LCK manifolds}

\question
The Fr\"olicher spectral sequence degenerates in $E^{*,*}_1$ page
for almost all known examples of LCK manifolds 
\index[terms]{spectral sequence!Fr\"olicher}
(\ref{_OT_degenerate_Theorem_}, \ref{_LCK_pot_degenerate_Corollary_}). 
Are there  LCK manifolds for which it {\em does not} degenerate?
A particular case of interest
is Kato manifolds (\ref{_Kato_mflds_Remark_}, \cite{iop}), for which 
this seems to be still unknown.

\hfill

\question
Let $(M, \theta, \omega)$ be a compact LCK manifold, and
$H_d^{1,0}(M)$ the space of closed holomorphic 1-forms. Consider
the natural map
\begin{equation}\label{_b_1_via_holo_Equation_}
H_d^{1,0}(M) \oplus \overline{H_d^{1,0}(M)} \oplus \langle \theta \rangle
\stackrel \rho \arrow H^1(M, \C);
\end{equation}
it is injective by \ref{_H^1_holo_LCK_Lemma_}.
This map is an isomorphism for LCK manifolds with potential\index[terms]{manifold!LCK!with potential}
(\ref{_LCK_pot_Hodge_decompo_Theorem_}), and for complex surfaces
(\ref{_H^1_odd_Theorem_}) and not an isomorphism 
for the LCK  OT-manifolds of dimension $> 2$ (\ref{LCKbun}).
 Are there any other examples of LCK\index[terms]{manifold!Oeljeklaus--Toma (OT)}
manifolds such that \index[terms]{surface!Inoue}
$H_d^{1,0}(M) \oplus \overline{H_d^{1,0}(M)} \oplus \langle \theta \rangle\neq H^1(M,\C)$? What are the geometric implications of this condition?
The isomorphism \eqref{_b_1_via_holo_Equation_}, that is  
true for complex surfaces, was one of the motivations
for the conjecture that $b_1(M)$ is odd for all\index[terms]{conjecture!Vaisman}
LCK manifolds (\ref{_Vaisman_conjecture_disproved_}).

\hfill

\question Let $(M, \theta, \omega)$ be a compact\index[terms]{connection!Chern}
LCK manifold, $\omega_g$ a 
Gauduchon metric on it,  and $\Theta_L$ the curvature 
of the Chern connection on its weight bundle $L$. Recall that the number
$\deg_\omega L:=\int_M \Theta_L\wedge \omega_g^{n-1}$\index[terms]{metric!Gauduchon}
is independent on the choice of the metric on the holomorphic 
bundle $L$ (Subsection \ref{_BC_degree_Subsection_}). For
a Vaisman manifold $(M, \omega, \theta)$, the curvature
of $L$ is $\omega_0$ \index[terms]{manifold!Vaisman}
(\ref{_Subva_Vaisman_Theorem_}); this implies that   the degree $\deg_\omega L$
is positive. The same is true for any other LCK manifold:
indeed, the curvature of $L$ is $d^c\theta$, but
\[
\int_M  d^c\theta\wedge \omega^{n-1} = \int_M \theta \wedge d^c(\omega^{n-1})=
\int_M \theta \wedge \theta^c \wedge \omega^{n-1} >0.\footnote{We are grateful to 
Nicolina Istrati \index[persons]{Istrati, N.} for this observation.}
\]
However, for Vaisman manifolds, $\deg_{\omega_1} L>0$ for any
Gauduchon metric $\omega_1$, because the curvature of $L$ is non-negative.
It seems that $\deg_{\omega_1} L>0$ also for OT LCK manifolds.

Is it true that \index[terms]{bundle!weight}
$\deg_{\omega_1} L >0$ for any LCK manifold $(M, \omega, \theta)$,
and any Gauduchon metric $\omega_1$ on $M$?

\hfill

\question
Let $M$ be a compact complex $n$-manifold.
Recall that {\bf the  Gau\-du\-chon cone}  of $M$ 
(\cite{_Popovici:Gauduchon_};\index[terms]{cone!Gauduchon}
see also Exercise \ref{_Gauduchon_cone_Exercise_})
is the set ${\cal G}\subset H^{n-1, n-1}_\Ae(M)$ of all 
{\bf \index[persons]{Gauduchon, P.} Gauduchon classes}\index[terms]{class!Gauduchon}
$[\omega^{n-1}]\in H^{n-1, n-1}_\Ae(M)$\index[terms]{cohomology!Aeppli}
in the  cohomology of $M$ represented
by $\omega^{n-1}$, where $\omega$ is a Gauduchon
form. The {\bf pseudo-effective cone}\index[terms]{cone!pseudo-effective}
is the set of all classes in the Bott--Chern cohomology
$H^{1,1}_{BC}(M)$ represented by positive
closed (1,1)-currents\index[terms]{current}.\index[terms]{cohomology!Bott--Chern}
In \cite[Corollary 3.6]{_Teleman:cone_}, A. Teleman
has computed the Gauduchon cone\index[terms]{cone!Gauduchon} of a complex surface
and proved that it is dual to the pseudo-effective cone,
under the Poincar\'e pairing \index[terms]{pairing!Poincar\'e}
$H^{1,1}_{BC}(M,\R)\times H^{n-1, n-1}_\Ae(M,\R)\arrow \R$.
The same result was proven in arbitrary dimension
by A. \index[persons]{Lamari, A.} Lamari, \cite[Lemma 3.3]{_Lamari_}.
Suppose now that a complex manifold $M$ admits an LCK structure,\index[terms]{structure!LCK}
and let ${\cal G}_{LCK}\subset {\cal G}$ be the set of all
Gauduchon class associated with LCK forms. Can
${\cal G}_{LCK}$ be strictly smaller than ${\cal G}$?
What is ${\cal G}_{LCK}$ for a Vaisman manifold?
For an LCK manifold with\index[terms]{manifold!LCK!with potential} potential?\index[terms]{manifold!Vaisman}

\subsection{Bimeromorphic geometry of LCK manifolds}

\question
Let $M$ be a Hopf manifold, and $Z\subset M$
a subvariety, possibly singular. 
\begin{enumerate}\index[terms]{manifold!Hopf}
\item Suppose that $\dim_\C Z =2$.
Prove that the normalization of $Z$ is a smooth
LCK surface.\footnote{We proved a special case of this result in
\ref{_surface_non-diagonal_Proposition_}.}\index[terms]{normalization}
\item More generally, suppose that the normalization
  $\tilde Z$
of $Z$ is smooth. Will $\tilde Z$ always admit an LCK structure?\index[terms]{structure!LCK}
\item Suppose that $\dim_\C Z >2$, and the normalization
of $Z$ is also singular. It is unclear if $Z$ can be
bimeromorphic to an LCK manifold in this case.
We expect that there are complex subvarieties of a Hopf
manifold that are  not bimeromorphic to an LCK manifold.
\end{enumerate}

\definition
Recall that {\bf the algebraic dimension} $a(M)$ of a compact complex manifold
$M$ is the transcendental dimension of its field of global
meromorphic functions. Moishezon has shown that
$a(M) \leq \dim(M)$, and the equality is realized if
and only if $M$ is bimeromorphic\index[terms]{dimension!algebraic}
to a projective manifold (\cite{moi}). Such\index[terms]{manifold!projective}
manifolds are now called {\bf Moishezon}.\index[terms]{manifold!Moishezon}
 By \ref{_dd^c_then_non-LCK_Theorem_}, a Moishezon manifold
does not admit a strict LCK structure\index[terms]{structure!LCK}; then  $a(M) <\dim (M)$
when $M$ is LCK. 

\hfill

\question
In Exercise \ref{_alge_dime_Hopf_Exercise_}, we have
established that the algebraic dimension of an $n$-dimensional Vaisman
Hopf manifold can take any values from $0$ to $n-1$.
Let $M$ be a $n$-dimensional LCK manifold with potential.\index[terms]{manifold!LCK!with potential}
What is the maximal algebraic dimension of $M$?\index[terms]{dimension!algebraic}

\hfill

\question
LCK OT-manifolds, such as \index[terms]{surface!Inoue} Inoue surfaces, have\index[terms]{surface!Inoue}\index[terms]{manifold!Oeljeklaus--Toma (OT)}
no complex subvarieties, and hence  their algebraic dimension
is 0. Vaisman manifolds\index[terms]{manifold!Vaisman} can have any algebraic dimension
between $0$ and $\dim_\C M-1$. The Kato manifolds, defined
by \index[persons]{Istrati, N.} Istrati, \index[persons]{Otiman, A.} Otiman and \index[persons]{Pontecorvo, M.} Pontecorvo (\cite{iop}, 
\ref{_Kato_mflds_Remark_}),\index[terms]{manifold!Oeljeklaus--Toma (OT)}
can have algebraic dimension  $\dim_\C M-2$ 
(\cite[Proposition 7.1, Theorem 8.1]{iop}).
Let $M$ be an LCK manifold that is  not bimeromorphic
to an LCK manifold with potential\index[terms]{manifold!LCK!with potential} or a Kato manifold. 
Can it have non-zero algebraic dimension?\index[terms]{manifold!Kato}

\hfill

\question
Let $M=\frac{\tilde M}{\langle A\rangle}$ be an LCK manifold with potential,
$\tilde M$ its open algebraic cone,\index[terms]{cone!algebraic} and ${\cal G}_A$ the complex algebraic group
obtained as the Zariski closure of $\langle A\rangle$ acting on the
algebraic cone. In  Exercise \ref{_orbit_space_and_a(H)_Exercise_},
it was shown that a meromorphic function on $\C^n$ is\index[terms]{cone!algebraic}
 $\langle A\rangle$-invariant if and only if it is\index[terms]{Zariski closure}
${\cal G}_A$-invariant. This implies that (roughly speaking)
meromorphic functions on a Hopf manifold\index[terms]{manifold!Hopf}
$\frac{\C^n\backslash 0}{\langle A\rangle}$ 
are meromorphic functions on the orbit space
of ${\cal G}_A$. The orbit space of ${\cal G}_A$
on $\C^n\backslash 0$ is quasi-projective,
if one removes the unstable orbits, as follows from 
the geometric invariant theory; then  its algebraic
dimension is equal to its dimension. \index[terms]{GIT}
When $M$ is a Hopf manifold, it seems\index[terms]{variety!quasi-projective}
that $a(M)$ is equal to $\dim_\C M - d$, where
 $d$ is the dimension of the general orbit of ${\cal G}_A$.
Is it true that the algebraic dimension of $M$
is $\dim_\C M - d$? We expect this to be true,
and the argument above should be sufficient
to prove it for Hopf manifolds.\index[terms]{dimension!algebraic}

\hfill

\question
{\bf A minimal model} of a complex variety $Z$
is a variety $Z'$, bimeromorphic to $Z$, which \index[terms]{minimal model}
does not admit holomorphic bimeromorphic contractions.
By \ref{bdown}, the blow-up of $M$ in a subvariety $Z$
that is  not of induced K\"ahler type cannot admit
an LCK structure.\index[terms]{structure!LCK}  It seems that this implies that
$M$ has a unique minimal model\index[terms]{minimal model}, is it true?

\hfill

\question
Let $M$ be an LCK manifold that does not admit
a non-trivial holomorphic bimeromorphic map to a complex manifold.
Examples of such $M$ are OT, Kato and LCK with potential.\index[terms]{manifold!LCK!with potential}
 It would be interesting\index[terms]{manifold!Kato}\index[terms]{manifold!Oeljeklaus--Toma (OT)}
to classify the minimal LCK manifolds and to
see whether the end result of a sequence of  holomorphic
 bimeromorphic contractions remains LCK (or ``singular LCK'', 
whatever it is).

\hfill

A compact complex manifold is called {\bf
  Fujiki class C} if it is bimeromorphic to a K\"ahler
manifold. A {\bf K\"ahler current} $\Theta$ is a positive, closed\index[terms]{current!K\"ahler}
(1,1)-current that satisfies $\Theta - \omega >0$\index[terms]{Fujiki class C}
for some Hermitian metric $\omega$ on $M$. \index[persons]{Demailly, J.-P.} Demailly--P\u aun
proved in \cite{dp} that a manifold is Fujiki class C
if and only if it admits a K\"ahler current.\index[terms]{current!K\"ahler}

\hfill

\question
An {\bf LCK Fujiki class C} manifold 
(Subsection \ref{_LCK_class_C_Subsection_})
is a compact complex manifold\index[terms]{current!LCK}
bimeromorphic to an LCK manifold. Let $\theta$ be a 
closed real 1-form on a complex manifold, and $d_\theta:=
d+ \theta \wedge$ be the standard Morse--Novikov differential.
 An {\bf LCK current} is a positive, $d_\theta$-closed (1,1)-current $\Theta$
that satisfies $\Theta - \omega >0$\index[terms]{theorem!Demailly--P\u aun}
for some Hermitian metric $\omega$ on $M$.
Is the LCK Fujiki class C condition equivalent to
the existence of an LCK current? In other words,
is the LCK version of the \index[persons]{Demailly, J.-P.} Demailly--P\u aun theorem valid?
We expect that it is false.

\hfill

\question
Let $M$ be a compact LCK manifold, and
$\phi:\; M \arrow M'$ a holomorphic, bimeromorphic map.
Clearly, $M'$ is equipped with an LCK current.\index[terms]{current!LCK} 
Is $M'$ also LCK? It should not always be LCK, otherwise
the ``LCK Fujiki class C'' condition makes little sense.

\hfill

\question
{\bf The \index[persons]{Lelong, P.} Lelong numbers} of positive, closed (1,1)-currents
are an important notion of complex geometry\index[terms]{Lelong number}
(\cite{_Demailly:analytic_}). The {\bf Lelong set}\index[terms]{Lelong set}
$Z_c$ of a positive, closed (1,1)-current\index[terms]{current} is the set
where the Lelong numbers are greater than a certain positive number $c\in \R^{>0}$. Y.-T. \index[persons]{Siu, Y.-T.} Siu has proven that the Lelong sets are
complex analytic (\cite{_Siu:1974_}).
For LCK currents the notion of Lelong number does not make
sense, because it is not conformally invariant;
however, the union of all Lelong sets of such a current
is well-defined, because it is locally conformally
equivalent to a closed, positive current. By \index[persons]{Siu, Y.-T.} Siu's theorem,
then, the Lelong set of an LCK current is a countable 
union of complex analytic subvarieties. \index[terms]{current!LCK}
Are all these subvarieties of\index[terms]{theorem!Siu}
IGCK\footnote{Subvarieties of IGCK type are 
  introduced in \ref{_IGCK_Definition_}.} type?

%

\hfill

\question
Let $\tilde M$ be the universal cover of a compact K\"ahler
manifold $M$. Shafarevich conjecture on uniformization
(\cite{_Eyssidieux:lectures_}) predicts that $\tilde M$
is holomorphically convex, that is, admits a proper
holomorphic map to a Stein variety. It is proven
when $\pi_1(M)$ admits an exact representation
(\cite{_EKRP:linear_Shafarevich_}). We\index[terms]{conjecture!Shafarevich} suggest  
the following LCK version of Shafarevich conjecture.
Let $M$ be an LCK manifold, and $\tilde M$ its
universal cover. Then $\tilde M$ admits 
a holomorphic map $\sigma:\; \tilde M\arrow S$ to a Stein
variety, with the following properties:
$\sigma$ is bijective in a general point, has compact fibres,
and $S \backslash (\im \sigma)$ is a complex subvariety in $S$.
For Vaisman manifolds,\index[terms]{manifold!Vaisman} this conjecture seems to follow
from the Shafarevich conjecture. \index[terms]{conjecture!Shafarevich}

\hfill

\question
Suppose  the 
``LCK Shafarevich conjecture'' stated above 
is true for Vaisman manifolds.\index[terms]{conjecture!Shafarevich (``LCK")}
Will it imply the Shafarevich 
conjecture for projective manifolds?

\hfill

\question
Any LCK manifold with potential\index[terms]{manifold!LCK!with potential} is a small deformation
of a Vaisman manifold, that has the same algebraic cone\index[terms]{cone!algebraic}
(\ref{_Vaisman_defo_same_cone_Remark_}). By Exercise
\ref{_trivial_canonical_then_QR_Exercise_}, 
a $\Z$-covering\index[terms]{cover!K\"ahler $\Z$-} of a Vaisman manifold \index[terms]{manifold!Vaisman}with trivial canonical
bundle is an elliptic fibration over a Calabi--Yau orbifold.
Since the group of automorphisms of a polarized\index[terms]{orbifold!Calabi--Yau}
Calabi--Yau orbifold is finite, this should imply that
any LCK manifold with potential with trivial canonical
bundle is biholomorphic to a quasi-regular Vaisman manifold.\index[terms]{manifold!Vaisman!quasi-regular}

Assume that $M$ is an LCK manifold with potential,
and the canonical bundle $K_M$ is holomorphically trivial.
Does $M$ admit a Vaisman metric? 

\hfill

\remark
The uniruled varieties\index[terms]{variety!uniruled}
(\ref{_uniruled_Definition_}) have Kodaira dimension $\kappa(M) =-\infty$,
because they contain a rational curve with non-negative normal bundle
(\cite{_Amerik:uniruled_}). However, a projective manifold with an
effective action of an algebraic group that preserves the
polarization is uniruled (\cite[Lemma 1.4]{_Andreatta:actions_}).
Let $M$ be a Vaisman manifold\index[terms]{manifold!Vaisman} such that its open algebraic cone\index[terms]{cone!algebraic} $\tilde M$
is a $\C^*$-bundle over a uniruled projective orbifold $X$.
By Exercise \ref{_Kodaira_dimension_Vaisman_}, $M$
has Kodaira dimension $\kappa(M)=-\infty$. Therefore, 
whenever $\kappa(M) \geq 0$, the group of automorphisms
of the orbifold $X$ that can be extended to the cone 
$\tilde M$ is finite. Then 
any Vaisman $\Z$-quotient of this cone is quasi-regular.
This implies that any Vaisman manifold with non-negative \index[terms]{dimension!Kodaira} Kodaira 
dimension must be quasi-regular.\index[terms]{manifold!Vaisman!quasi-regular}

\hfill

\question
Let $M$ be an LCK manifold with potential\index[terms]{manifold!LCK!with potential}. Assume that 
the Kodaira dimension of $M$ is non-negative.\index[terms]{dimension!Kodaira}
Will it follow that $M$ is Vaisman and quasi-regular?
We expect that the answer is positive.

\hfill

\question
Are LCK manifold of non-negative
 Kodaira dimension Vaisman?

\subsection{Complex subvarieties in LCK manifolds}

\question All LCK manifolds with potential \index[terms]{manifold!LCK!with potential}have 
complex subvarieties by \ref{flag}.  On the other hand, 
LCK OT manifolds have no proper complex subvarieties whatsoever (\ref{no_subv}). \index[terms]{manifold!Campana simple}
Generally, we have no easy way of characterizing
LCK manifolds with few or many complex subvarieties.
Recall that a complex manifold $M$ is called \index[terms]{dimension!algebraic}
{\bf \index[persons]{Campana, F.} Campana simple} if the set of all proper
complex subvarieties is contained in a countable union
of nowhere dense closed subsets.\footnote{A set that is  contained in a countable union \index[terms]{surface!K3}
of nowhere dense closed sets is sometimes called {\bf meager}, 
or {\bf Baire category 1}.} Clearly, complex manifolds\index[terms]{Baire category}
of positive algebraic dimension cannot be Campana simple.
However, a manifold of zero algebraic dimension is not
always \index[persons]{Campana, F.} Campana simple; for example, one can take a product
of two non-algebraic K3 surfaces. The products of LCK
manifolds conjecturally do not admit an LCK metric,
so in the LCK case this kind of situation is impossible. 
Are LCK manifolds of algebraic dimension 0 always \index[persons]{Campana, F.} Campana simple?

\hfill

\question
Recall that a complex manifold $M$ is called {\bf uniruled}\index[terms]{manifold!uniruled}
if for every point $p$ of $M$, there exists a rational curve on $M$ passing
through $p$. From \ref{_non-uniruled_Corollary_} it follows that a strict
LCK manifold cannot be uniruled when the Barlet space of the curves
that cover the whole $M$ is compact. Is the compactness assumption 
necessary?\index[terms]{space!Barlet}

\hfill

\question
Let $C\subset M$ be a rational curve in a compact LCK manifold, and
$Z$ the union of all deformations of $C$. 
\begin{enumerate}
\item  Is $Z$ always compact? We expect that it is.
\item Assume that $Z$ is compact;
it is possible to show, using Remmert's proper mapping theorem,
that $Z$ is a complex subvariety in $M$. \index[terms]{theorem!Remmert's proper mapping}
Suppose that $C$ is smooth. Will it follow that $Z$ is
{\bf contractible}, that is, there is a bimeromorphic holomorphic
map $M \arrow M'$ contracting $Z$? It is true in many cases
we could check (for the manifolds bimeromorphic
to LCK with potential, OT and solvmanifolds\index[terms]{solvmanifold})\index[terms]{manifold!LCK!with potential}.
\end{enumerate}

\question
Currently, the only known LCK manifolds that have families
of complex curves passing through every point are Vaisman.
Suppose that an LCK manifold $M$ is obtained as the union
of a collection of complex curves in the same deformation family.
Will it imply that all these curves are elliptic?\index[terms]{curve!elliptic}
Will it imply that $M$ is Vaisman?

\hfill

\question
Let $M$ be a compact complex manifold of K\"ahler
type. By Bi\-shop's compactness theorem \index[terms]{theorem!Bishop's compactness}
(\cite{_Bishop:conditions_}), any connected component
of the \index[persons]{Barlet, D.} Barlet space of complex subvarieties of $M$ is compact.
The K\"ahler condition is often not necessary.
\index[persons]{Gromov, M.} Gromov's compactness theorem \index[terms]{space!Barlet}
(\cite{_Gromov:pseudoholomorphic_})\index[terms]{form!symplectic!tamed}
implies that any connected component $Z_i$ of 
the Barlet space of complex curves is compact
if the complex structure on $M$ is tamed\index[terms]{theorem!Gromov's compactness}
by a symplectic structure. In \cite{_Verbitsky:twistor_},
it was shown that $Z_i$ is compact when
the Hermitian form $\omega$ is pluriclosed, or when
$dd^c\omega$ is a negative $(2,2)$-form.\index[terms]{form!pluriclosed (SKT)}
Suppose now that $M$ is a compact complex manifold
that admits an LCK structure\index[terms]{structure!LCK}. Will it imply that
any connected component of 
the \index[persons]{Barlet, D.} Barlet space of complex curves in $M$ is compact?
This is true in all examples we know.


\section{Sasakian and Vaisman manifolds}\index[terms]{manifold!Vaisman}


\subsection{Vaisman manifolds}

\question\label{_local_homo_nilma_LCK_pot_Question_}
Locally homogeneous LCK structures\index[terms]{structure!LCK} on nilmanifolds
are classified by \index[persons]{Sawai, H.} Sawai (\ref{inv_nil_vai}), who proved
that all locally homogeneous LCK nilmanifolds 
are Vaisman and Heisenberg type.\index[terms]{nilmanifold}
Let $M$ be a compact solvmanifold\index[terms]{solvmanifold} with locally
invariant LCK structure. Is it always Vaisman?
If it is Vaisman, will it follow that $M$
\index[terms]{manifold!Vaisman!Heisenberg type}
is a Heisenberg-type nilmanifold?
\index[terms]{manifold!LCK!locally homogeneous}

\hfill

\question
(see Subsection \ref{_qr_chapter_Intro_Section_}).
Let $X$ be a projective orbifold, and $L$
an ample line bundle, such that $\Tot^\circ(L)$\index[terms]{bundle!line!ample}
(the space of all non-zero vectors in $L$) is smooth.
In this case, the orbifold singularities of $X$
are abelian, that is, each orbifold chart on\index[terms]{orbifold!projective}
$X$ is modeled on $B/G$, where $B$ is an open ball
in $\C^n$, and $G$ is an abelian group. Let $X$ be a projective orbifold
with abelian quotient singularities. What are
the conditions which guarantee that $\Tot^\circ(L)$
is smooth?

%

\hfill

\definition
Let $M_1, M_2$ be LCK manifolds with proper potential.\index[terms]{manifold!LCK!with potential}
We say that $M_1$ is {\bf cone equivalent} to $M_2$ if
their open algebraic cones\index[terms]{cone!algebraic} are biholomorphic \index[terms]{manifold!LCK!cone-equivalent}
(by \ref{_same_algebra_structure_on_cone_Theorem_}, 
they are consequently equivalent as algebraic varieties).

\hfill

\question
Let $M$ be a Vaisman manifold. Is $M$ always cone equivalent
to a regular Vaisman manifold?\index[terms]{manifold!Vaisman!regular} For a Hopf manifold this is clearly true.

\hfill

\question
Let $M_1$ and $M_2$ be two cone equivalent
regular Vaisman manifolds, and 
$X_1$, $X_2$ the leaf spaces of their canonical foliations.\index[terms]{foliation!canonical}
Are they always homeomorphic?\footnote{When $M_i$ are quasi-regular,\index[terms]{manifold!Vaisman!quasi-regular}
the leaf space orbifolds $X_1, X_2$
are non-homeomorphic in many known
cases. Indeed, any weighted projective
space can be obtained as the orbit space
of an appropriate $\C^*$-action on $\C^n\backslash 0$.}
For Vaisman surfaces this 
follows directly from their classification
(\ref{_Vaisman_is_Hopf_or_elli_Proposition_}).

This question seems to be related
to \index[persons]{Milnor, J.} Milnor's construction of exotic
smooth manifolds homeomorphic to 
odd-dimensional spheres (that are  known to be 
always Sasakian, \cite{bog}). 
For toric Vaisman manifolds,\index[terms]{manifold!Vaisman}
a similar question is considered below,
\ref{_toric_cone_equivalent_Question_}.

\hfill

\question
Let $M$ be a Vaisman manifold, $\tilde M$\index[terms]{Zariski closure}
its K\"ahler $\Z$-cover\index[terms]{cover!K\"ahler $\Z$-}, and ${\cal G}$ the Zariski closure
of the $\Z$-action on $\tilde M$, acting on $M$.
A priori, the same complex manifold might have
two different Vaisman structures of LCK rank 1\index[terms]{rank!LCK},
giving two different K\"ahler $\Z$-covers.
Is the ${\cal G}$-action independent
from the choice of $\Z$-cover?

\subsection{Sasakian manifolds}

\question
In \cite{_MSY:volumes_Sasaki-Einstein_},
\index[persons]{Martelli, D.} Martelli, Sparks and \index[persons]{Yau, S.-T.} Yau proved\index[terms]{manifold!Sasaki--Einstein}
that the volume $a(S)$ of a Sasaki--Einstein manifold $S$, relative to
that of the round sphere with the same Einstein constant, 
is always an algebraic number.\index[terms]{manifold!Einstein--Weyl}
When $S$ is quasi-regular, this number is rational.\index[terms]{manifold!Sasaki!quasi-regular}
They defined the {\bf rank} of a Sasaki manifold as
the dimension of the closure if the generic Reeb orbit.%
\footnote{By \ref{_limits_torus_Proposition_}, 
the dimension of an orbit of a torus acting
by isometries is semicontinuous, and hence 
outside of a closed, nowhere dense set, the closure
of the Reeb orbit has maximal dimension.}
 Martelli, Sparks and  Yau conjecture that the
degree of the algebraic number $a(S)$ is
equal to the rank of $S$. Since 2006,\index[terms]{conjecture!Martelli--Sparks--Yau}
this question remains open (see \cite[Proposition 2.21,
  Lemma 4.2]{_Donaldson_Sun:limits_2_} 
for some developments in this direction). 
We expect that a similar result would
hold for a Vaisman manifold\index[terms]{manifold!Vaisman}, and propose
a conjecture in the same vein. Let $M$ be a
Vaisman manifold; it is called {\bf Einstein--Weyl}
if its K\"ahler cover is Ricci-flat\index[terms]{manifold!Vaisman!Einstein--Weyl} 
(\ref{ewcy}). Locally, such a manifold
is a product of a line tangent to the Lee
field and a Sasaki--Einstein manifold.\index[terms]{manifold!Sasaki--Einstein}
Let us normalize the metric on the 
Vaisman  Einstein--Weyl manifold in such
a way that $|\theta|=1$. Then the volume
of $M$, relative to the volume of $S^{2n-1}\times S^1$
of the same transversal Ricci curvature, should be 
algebraic; this statement should probably follow
from the same argument as the  Martelli--Sparks--Yau
theorem. The analogue of the Martelli--Sparks--Yau
conjecture would state the degree of this
algebraic number is equal to the dimension
of the closure of the generic orbit of the 
anti-Lee field\index[terms]{Lee field!anti-} $I(\theta^\sharp)$.\index[terms]{theorem!Martelli--Sparks--Yau}

\hfill

\question
The Reeb dynamics on 
Sasakian manifolds is well understood by now 
(\cite{_Rukimbira_}, Chapter 
\ref{_Elliptic_curves:Chapter_}), 
because the Sasakian manifold is transversally K\"ahler.
Is it possible to prove similar results on 
the Reeb dynamics for the contact manifolds
obtained as level sets of LCK potentials\index[terms]{potential!LCK}
on compact LCK manifolds with potential?\index[terms]{manifold!LCK!with potential}

\hfill

\question
Let $S$ be a compact Sasakian manifold,
and $(B\subset TS, I\in \End(B))$ its
CR-structure.\index[terms]{structure!CR} As follows from 
\ref{_Reeb_fie_from_CR_Theorem_},
the Sasakian metrics compatible
with this CR-structure are in \index[terms]{vector field!CR-holomorphic}
bijective correspondence with \index[terms]{cone!Sasakian}
CR-ho\-lo\-mor\-phic vector fields which\index[terms]{cone!K\"ahler}
are transversal to $B$ and positively oriented
with respect to the standard orientation on $TS/B$.
The set ${\goth t}^+$ of such vector fields is called
{\bf the Sasakian cone}; it was introduced
in \cite{_BGS:cone_,_BGS:canonical_}
and studied in the series of papers by\index[terms]{manifold!Fano}
Ch. P. Boyer\index[persons]{Boyer, C. P.} and his collaborators, e.g., 
\cite{_Boyer:maximal_tori_,_Boyer_TF:positivity_}.
In \cite{_BGS:canonical_} Boyer\index[persons]{Boyer, C. P.}, \index[persons]{Galicki, K.} 
Galicki and \index[persons]{Simanca, S.} Simanca
suggested that ${\goth t}^+$ plays
 the same role for the Sasakian manifolds as
the K\"ahler cone plays for the K\"ahler manifolds.
For a K\"ahler manifold $M$, the convex geometry of the
K\"ahler cone reflects the algebraic geometry \index[terms]{geometry!algebraic}of $M$;
for example, the effective cone of curves on a 
Fano manifold is expressed through the K\"ahler
cone and the duality
(\cite[Theorem 6.1]{_Debarre_}). From
the idea of \index[persons]{Boyer, C. P.}Boyer, \index[persons]{Galicki, K.}Galicki and \index[persons]{Simanca, S.} Simanca it 
would follow that the Sasaki cone should
play an equally important role.

Clearly, ${\goth t}^+$ is an open convex cone
in the Lie algebra ${\goth t}$ of CR-holomorphic
automorphisms of $S$. It is not hard to compute
${\goth t}^+$  in some special cases.\index[terms]{resonance}
Let, for example, $S\subset \C^n$ be a strictly
pseudoconvex submanifold that is  realized\index[terms]{cone!polyhedral}
as a pseudoconvex shell,\index[terms]{pseudoconvex shell} and $H= \frac{\C^n \backslash 0}{\Z}$ the 
corresponding Hopf manifold, \ref{_linear_LCK_pot_Corollary_}.
Assume that $H$ is Vaisman and has no resonance.
It is possible to show that in this case,
all holomorphic vector fields on $H$
are linear (\ref{_global_vector_on_Hopf_Proposition_}),
thus ${\goth t}$ is a subalgebra in $\goth{gl}(n,\C)$.
The Lie group
of linear CR-holomorphic automorphisms of $S$ is compact,
because it preserves a Finsler metric on $\C^n$,
hence ${\goth t}$ is compatible with a
unitary structure in $\goth{gl}(n,\C)$,
giving ${\goth t}\subset \goth{u}(n)$.
The transversality and positive orientation
condition for $X\in {\goth t}\subset \goth{u}(n)$
is equivalent to the Hermitian matrix 
$IX\in -\1 \goth{u}(n)$ defining a contraction of $\C^n$.
This is the same as $IX$ having all eigenvalues $\alpha_i< 0$.
When $\goth t$ is commutative, all 
$X\in {\goth t}\subset \goth{u}(n)$ are diagonal in the
same basis, and the cone ${\goth t}^+$ is polyhedral.
When $S$ is a round sphere, 
${\goth t}=\goth{u}(n)$, and ${\goth t}^+$
can be identified with the set of all Hermitian
forms on $\C^n$ with non-negative eigenvalues.
In this case, ${\goth t}^+$ is clearly non-polyhedral.
It seems that the above argument implies that 
${\goth t}^+$ is always polyhedral when $\goth t$
is commutative. Can it be polyhedral when
the algebra $\goth t$ is not commutative?
How does the shape of the cone ${\goth t}^+$ 
reflect the geometry of the Sasakian manifold $S$?


\section{LCK manifolds with potential}\index[terms]{manifold!LCK!with potential}


%



\question 
Sasakian manifolds can be defined as level sets of LCK
potential on Vaisman manifolds.\index[terms]{manifold!Vaisman} However, they also have
and intrinsic definition (\cite{bl, bog}).
Is it possible to have an intrinsic definition
for the metric contact structure that appears
on the level sets of LCK potentials?\index[terms]{potential!LCK} \index[terms]{structure!metric contact}

\hfill

%

\question
Let $S$ be the level set of
the K\"ahler potential on the K\"ahler cover
of an LCK manifold $(M, \omega)$ with potential,\index[terms]{manifold!LCK!with potential}
considered to be  a hypersurface in $M$.
When $\omega$ is Vaisman, $S$ is totally geodesic in $M$. Can we modify the LCK
metric with potential in its conformal class, or by adding
a $d_\theta d^c_\theta$-exact term, in such a way that all
level sets are minimal submanifolds in $(M, \omega)$?\index[terms]{submanifold!minimal}
\index[terms]{submanifold!totally geodesic}

\hfill

\question
Let $M$ be an LCK manifold, and $D$ its Kuranishi space
of complex deformations (\ref{_Kuranishi_Definition_}), 
understood as a germ of a complex variety.\index[terms]{space!Kuranishi}
When $M=H_A$ is a linear, non-resonant Hopf manifold, $H_A = \frac{\C^n
  \backslash 0}{\langle A \rangle}$,
we identify $D$ with the space of all $E \in \End(\C^n)$
commuting with $A$ (\ref{_global_vector_on_Hopf_Proposition_}).
A general
deformation $H_{A'}$ of\index[terms]{resonance}
$H_A$ has diagonal $A'$ with all eigenvalues distinct,
hence it is Vaisman. Moreover, the points of $D$
corresponding to non-Vaisman manifolds\index[terms]{manifold!Vaisman} are a meager
subset of measure zero in $D$.
\begin{enumerate}
\item   For a general LCK manifold $M$\index[terms]{manifold!LCK!with potential}
with potential, is the set of Vaisman deformations of $M$
always dense in $D$? Is the set of non-Vaisman 
deformations of $M$ always meager and of measure zero in $D$?
\item
Let $D$ be the  Kuranishi\index[terms]{space!Kuranishi} space of a Vaisman manifold
$M$, and $D_0$ the set of points $x\in D$ corresponding
to Vaisman manifolds. Can $D_0$ be closed in $D$?
Can it be closed and 0-dimensional? Neither of these
can happen for Hopf manifolds.\index[terms]{manifold!Hopf}
\end{enumerate}

\question 
In complex dimension $\geq 3$, all LCK manifolds with
potential can be deformed to Vaisman (\ref{def_lckpot2Vai})
and can be embedded to a Hopf manifold (\ref{embedding}).
The proof uses  Andreotti--Siu and  Rossi theorem,
that is  false in dimension 2. However, it is
possible to deduce the embedding and the deformation
result from the Global Spherical Shell Conjecture
(\ref{_sphe_implies_dim2_Theorem_}).\index[terms]{conjecture!GSS} \index[terms]{theorem!Rossi, Andreotti--Siu}
Since the Global Spherical Shell Conjecture 
seems to be out at reach an the moment, we would
like to have a direct proof of these two results
in dimension 2.


\section{Extremal metrics on LCK manifolds}\index[terms]{metric!extremal}


\question
The \index[persons]{Futaki, A.} Futaki invariant is a character of the
Lie group of conformal isometries of a 
compact complex manifold defined in \index[terms]{Futaki invariant}
\cite{_Futaki:complex_}, see also \ref{_Futaki_Ricci_Definition_}.
It is easy to produce an example of an LCK
manifold without holomorphic vector fields.
As shown already in \cite{ot}, an OT manifold
has no holomorphic vector fields; the same
is true for its special case, an \index[terms]{surface!Inoue} Inoue surface.
If one blows up a sufficient number of points,
a given LCK manifold can be transformed to 
a manifold that has no holomorphic vector fields
(and no holomorphic automorphisms at all).
It is not hard to see that general Kato surfaces
have no holomorphic vector fields.\index[terms]{manifold!toric Kato}
On the other hand, the toric Kato manifolds
(\cite{_IOPR:toric_Kato_}) might admit 
holomorphic vector fields (\cite{iop}). Any LCK manifold with 
potential admits a holomorphic vector field 
(\ref{_S^1_action_exists_Corollary_}).
When $M$ is Vaisman, its \index[persons]{Futaki, A.} Futaki invariant
vanishes by \ref{_Futaki_vanishes_Vaisman_Theorem_}.
Will it vanish for all LCK manifolds with potential?
We expect that the answer is positive. Will it vanish
for toric Kato manifolds?\index[terms]{manifold!toric Kato}

\hfill

\question Let $(M,I,g)$ be a compact K\"ahler manifold and
$\Scal_g$ its scalar curvature. The functional $g\mapsto
\mathcal{R}(g):=\int_M\Scal_g^2 \vol_g$ was introduced by
\index[persons]{Calabi, E.} Calabi (\cite{_Calabi:extremal_}, see also
\cite{_Gauduchon:book_}).
 Fix a K\"ahler class\index[terms]{class!K\"ahler} $[\omega]\in H^{1,1}(M)$,
and let $\Omega$ be the set of all K\"ahler metric
with this K\"ahler class. A metric $g\in \Omega$ is a critical 
 point of $\mathcal{R}\restrict{\Omega}$ if and only if  
 the gradient of its scalar curvature is a real holomorphic 
 vector field. In particular, K\"ahler--Einstein metrics are 
 critical points of $\mathcal{R}\restrict{\Omega}$ 
 (and indeed its absolute minima, 
\cite{_Gauduchon:book_}). \index[terms]{metric!K\"ahler--Einstein}
Now let $(M,I,g,\omega)$ be a compact LCK manifold 
of Vaisman type. Can we build a similar functional 
on the space of LCK metrics with potential\index[terms]{metric!LCK!with potential} and the
preferred gauge, with the critical points corresponding to
 Einstein--Weyl Vaisman metrics?\index[terms]{metric!Vaisman!Einstein--Weyl}
Can we build a functional on the space of conformal
classes of LCK metrics such that its critical points
are conformal classes of Einstein--Weyl type?

\hfill

\question 
Let $X$ be a projective orbifold and $M$ a quasi-regular
Vaisman manifold \index[terms]{manifold!Vaisman!quasi-regular}$M$ equipped with a principal elliptic
fibration $f:\; M \arrow X$.
Using the adjunction formula, one can compute \index[terms]{adjunction formula}
the canonical bundle $K_M$ as follows:
$K_M=f^* K_X$, where $K_X$ is the canonical bundle of $X$.
When $K_M$ is trivial, $M$ is quasi-regular, and $X$ is a Calabi--Yau orbifold by
Exercise \ref{_trivial_canonical_then_QR_Exercise_}.\index[terms]{orbifold!Calabi--Yau}
\begin{enumerate}
\item 
Let $M$ 
be a compact LCK manifold with trivial canonical bundle.
Fix the Lee class $[\theta]$\index[terms]{class!Lee} on $M$. By \ref{_CY_Vaisman_by_volume_Theorem_}, the Vaisman metric
with Lee class $[\theta]$,
inducing a flat metric on $K_M$, is unique up to a
constant multiplier.\index[terms]{metric!Calabi--Yau}
By the adjunction formula $K_M=f^* K_X$, the corresponding
K\"ahler metric on $X$ is Calabi--Yau. Let $M$ 
be non-Vaisman. Will it admit an LCK metric such that
a holomorphic section of $K_M$ has constant length?
This would give an analogue of Calabi--Yau theorem for
LCK manifolds with trivial canonical bundle.
Such a metric is Ricci-flat in the sense of
\ref{_Futaki_Ricci_Definition_}.
\item
Is a Ricci-flat metric unique if we fix\index[terms]{class!Lee}
the Lee class and the twisted Bott--Chern class \index[terms]{class!Bott--Chern}
(\ref{_twi_BC_class_Definition_}) of an LCK manifold?

\item  This statement would be implied by 
the following LCK version of the Ca\-la\-bi-Yau theorem.
Let $\omega_1, \omega_2$ be Gauduchon metrics\index[terms]{metric!Gauduchon}
on an LCK manifold. Assume that $\omega_1$ and $\omega_2$
have the same Lee class\index[terms]{class!Lee} and the same twisted Bott--Chern
class. Will $\omega_1^n= \omega_2^n$ imply that 
$\omega_1=\omega_2$?\footnote{The existence of a Gauduchon
metric with prescribed Ricci form\index[terms]{form!Ricci} is proven 
in \cite{_STW:Gauduchon_}.} 
\end{enumerate}

\remark
Let $S$ be a Sasaki--Einstein manifold.
Then, for an appropriate choice of\index[terms]{manifold!Sasaki--Einstein}
the cone metric, its Riemannian cone\index[terms]{cone!Riemannian} is a Ricci-flat
K\"ahler--Einstein manifold (\cite[Lemma 11.1.5]{bog}), in particular,\index[terms]{manifold!Einstein--Weyl}
it is Einstein--Weyl, and all its Vaisman $\Z$-quo\-tients
are Einstein--Weyl. These manifolds have\index[terms]{manifold!K\"ahler--Einstein}
LCK rank 1 \index[terms]{rank!LCK}, and, indeed, from \index[persons]{Gauduchon, P.} Gauduchon's
classification of Einstein--Weyl LCK
manifolds, the structure theorem\index[terms]{manifold!LCK!Einstein--Weyl}
(\ref{halfstr}) and from \ref{betti_1} it follows that any
Einstein--Weyl LCK manifold is obtained
this way. Therefore, the Sasaki--Einstein\index[terms]{manifold!Sasaki--Einstein}
geometry is roughly equivalent to the Einstein--Weyl 
geometry of Vaisman manifolds.\index[terms]{manifold!Vaisman} However, the 
Vaisman geometry gives a smooth passage between
different Sasakian manifolds realized as
pseudoconvex shells in the same LCK manifold.
 All Hermitian Einstein--Weyl manifolds admit
a Vaisman metric, but\index[terms]{pseudoconvex shell}
not all Vaisman-type complex manifolds admit
an Einstein--Weyl structure. Indeed, for such a 
manifold the weight bundle is proportional 
to the canonical bundle, and hence  the canonical
bundle must be flat. This already puts some
restriction to the geometry of $M$.\index[terms]{bundle!vector bundle!canonical}
Clearly, any compact Einstein--Weyl \index[terms]{manifold!Einstein--Weyl}
manifold $M$ can be deformed to a quasi-regular
Vaisman manifold $M_1$,\index[terms]{manifold!Vaisman!quasi-regular} not necessarily Einstein--Weyl
anymore. Still, $c_1(M)=c_1(M_1)=0$.
On the other hand, $M_1$ is the total space
of an elliptic bundle $\tau:\; M_1\arrow X$ on a projective
orbifold $X$, and this elliptic bundle
is obtained as a $\Z$-quotient of $\Tot^\circ(L)$
of non-zero vectors in the total space of a negative
line bundle (\ref{_Structure_of_quasi_regular_Vasman:Theorem_}).
By the adjunction formula, $\tau^*(K_X) = K_{M_1}$;
this implies $\tau^*(c_1(X))=0$.\index[terms]{adjunction formula}
Since the kernel of $\tau^*$ on $H^2(X)$
is generated by $c_1(L)$ (Exercise \ref{_kernel_pullback_H^2_Exercise_}), 
this implies that
$c_1(L)$ is proportional to $c_1(X)$.
The proportionality coefficient is negative,
because all Sasaki--Einstein manifolds\index[terms]{manifold!Sasaki--Einstein}
have positive Einstein constant.\footnote{Indeed, all
Sasaki--Einstein manifolds have Killing spinors,\index[terms]{spinor!Killing} 
\cite{_Bar:killing_}, and all manifolds with non-zero Killing spinors
are Einstein of positive Ricci curvature, \cite{bfgk}.}
This implies that $X$ is a Fano orbifold.\index[terms]{orbifold!Fano}

The correspondence between K\"ahler--Einstein metrics
on a Fano manifold and conical Calabi--Yau metrics
on its algebraic cone\index[terms]{cone!algebraic} is called {\bf the Calabi ansatz} 
(\cite{_Conlon_Hein:Asymptotical_}).\index[terms]{metric!Calabi--Yau}\index[terms]{Calabi ansatz}

\hfill

\question
Let $\tilde M$ be an algebraic cone\index[terms]{cone!algebraic} constructed from
a Fano orbifold. What are the conditions sufficient and
necessary for $\tilde M$ to admit a Ricci-flat
K\"ahler metric with conical symmetry? Equivalently,
what are the  sufficient and\index[terms]{cone!algebraic}
necessary conditions for $\frac{\tilde M}{\Z}$ to admit
an Einstein--Weyl Vaisman structure?\index[terms]{manifold!Vaisman!Einstein--Weyl}

\hfill

We are going to make a more precise version of this
question now. A quasi-regular Sasaki--Einstein manifold can always \index[terms]{manifold!Sasaki--Einstein!quasi-regular}
be obtained as a circle bundle over a K\"ahler--Einstein
Fano orbifold (\cite{bog}). However, not all Sasaki--Einstein manifolds\index[terms]{manifold!Sasaki--Einstein}
are quasi-regular. In one of the first papers written on the 
subject, \index[persons]{Cheeger, J.} Cheeger and \index[persons]{Tian, G.} Tian conjectured that 
all Sasaki--Einstein manifolds\index[terms]{manifold!Fano}
are quasi-regular (\cite{_Cheeger_Tian:cone_}), but some years
later, physicists found the\index[terms]{conjecture!Cheeger-Tian} now\index[terms]{manifold!Sasaki--Einstein}
famous set of examples\index[terms]{manifold!K\"ahler--Einstein}
(\cite{_GMSW:S^2xS^3_,_GMSW:new_,_Martelli_Sparks:2005_})
of irregular Sasaki--Einstein manifolds.

\hfill

\question
Let $\tilde M= \Tot^\circ(L)$ be an algebraic cone\index[terms]{cone!algebraic} over a Fano
\index[terms]{cone!algebraic}
orbifold $X$. By \cite[Theorem 1.1]{_Ross_Thomas:stability+CSCK_},
and also \cite{_CS:Sasaki-Einstein_}, the
existence of a Sasaki--Einstein pseudoconvex shell $S$
in $\tilde M$ is equivalent to the orbifold K-stability of
$X$, provided that $S$ is quasi-regular. The existence of a
Sasaki--Einstein pseudoconvex shell\index[terms]{pseudoconvex shell} in $\tilde M =
\Tot^\circ(L)$\index[terms]{manifold!Sasaki--Einstein} implies\index[terms]{manifold!Sasaki!quasi-regular}
that $X$ is a K\"ahler--Einstein orbifold, by \cite[Theorem 11.1.11]{bog}.
Ross--Thomas theorem
\cite[Theorem 1.1]{_Ross_Thomas:stability+CSCK_} is parallel to a celebrated
theorem of \index[persons]{Chen, X.} Chen--\index[persons]{Donaldson, S. K.}Donaldson--\index[persons]{Sun, S.}Sun and \index[persons]{Tian, G.} Tian,
who showed that K-stability on Fano manifolds is equivalent
to existence of K\"ahler--Einstein metrics
(\cite{_Chen_Donaldson_Sun_1_,_Chen_Donaldson_Sun_2_,_Chen_Donaldson_Sun_3_,_Tian:K-stability_}).
Suppose that $X$ does not admit a
K\"ahler--Einstein metric. There are cases when
 $\tilde M$ still contains a pseudoconvex shell 
admitting\index[terms]{pseudoconvex shell}
a Sasaki--Einstein metric compatible with its\index[terms]{structure!CR}\index[terms]{manifold!Sasaki--Einstein}
CR-structure. By \index[persons]{Lichnerowicz, A.} Lichnerowicz-\index[persons]{Matsushima, Y.}Matsushima \index[terms]{theorem!Lichnerowicz}
(see \cite[Chapter 19]{_moroianu:book_} and \index[terms]{theorem!Matsushima}
\cite{_Shu:CSCK_}), a K\"ahler--Einstein
(and even constant scalar curvature K\"ahler, a. k. a. CSCK) manifold
has reductive group of holomorphic automorphisms;
therefore, the space $\C P^2$ blown up in 2 points 
does not admit CSCK metrics. By \cite[Corollary 1.3]{_FOW:toric_},
the cone of a toric manifold always admits a 
Sasaki--Einstein pseudoconvex shell.\index[terms]{pseudoconvex shell}
In the proof of \cite[Corollary 1.3]{_FOW:toric_}
this Sasaki--Einstein structure is constructed explicitly;
by its construction it is apparent that it is irregular.
Existence of Sasaki--Einstein  pseudoconvex shell in an
open algebraic cone\index[terms]{cone!algebraic}\index[terms]{manifold!Sasaki--Einstein}
(or, equivalently, a cone-type  Calabi--Yau metric on the cone)
is an interesting condition, not yet understood.\index[terms]{metric!Calabi--Yau}
It is apparent from \cite[Corollary 1.3]{_FOW:toric_} 
that it is strictly weaker than the existence of a 
K\"ahler--Einstein metric. However, it is not 
clear yet if it is weaker than the existence
of a CSCK metric. In fact, we do not know any
Fano manifold $X$ such that its 
algebraic cone $\Tot^\circ(K_X)$
is known not to admit a cone-type Calabi--Yau metric.
Equivalently, we do not know of any Sasakian
manifold $S$ with positive transversal $c_1(S)$
such that its cone does not have a 
pseudoconvex shell that is  Sasaki--Einstein\index[terms]{manifold!Sasaki--Einstein}
(irregular, most likely). We are interested
to know if there are Fano orbifolds $X$ for which
the algebraic cone\index[terms]{cone!algebraic} $\Tot^\circ(K_X)$ does not
contain a Sasaki--Einstein pseudoconvex shell\index[terms]{pseudoconvex shell},
and how they can be characterized.

\hfill

\remark
Let $\tilde M=\Tot^\circ(K_X)$ be a cone over a Fano orbifold
$X$, and $\xi$ a Reeb field of a Sasakian pseudoconvex shell\index[terms]{pseudoconvex shell}
in $\tilde M$. One can interpret $\xi$ as a holomorphic\index[terms]{orbifold!Fano}
vector field on $\tilde M$.  In \cite{_CS:irregular_Sasakian_},
Collins and \index[persons]{Sz\'ekelyhidi, G.} Sz\'ekelyhidi redefined the notion of orbifold K-stability
to include $\xi$ as an automorphism of the test configuration.
They prove that the existence of \index[terms]{K-stability}
constant scalar curvature Sasakian metric on this
(possibly irregular) Sasakian pseudoconvex shell\index[terms]{pseudoconvex shell} 
with the Reeb field $\xi$ implies the orbifold 
K-stability of the pair $(X, \xi)$.

\hfill

\question
The existence of a Sasaki--Einstein pseudoconvex shell\index[terms]{pseudoconvex shell} in the
cone $\tilde M=\Tot^\circ(K_X)$ over a Fano manifold
might lead to a new notion of ``extremal metric''\index[terms]{manifold!Sasaki--Einstein}
condition for Fano manifolds. Will this condition
imply existence of the constant scalar curvature
metrics? Of the extremal metrics, in the sense
of \index[persons]{Calabi, E.} Calabi (\cite{_Calabi:extremal_})? \index[terms]{metric!extremal}
The papers of Ross-Thomas and Collins--\index[persons]{Sz\'ekelyhidi, G.}Sz\'ekelyhidi
seem to be relevant here: 
\cite{_Ross_Thomas:stability+CSCK_,_CS:Sasaki-Einstein_,_CS:irregular_Sasakian_}.

\hfill

\question
The celebrated theorem of \index[persons]{Bando, S.} Bando--\index[persons]{Mabuchi, T.}Mabuchi
\cite{_Bando_Mabuchi:transitivity_} implies that
the connected component of the\index[terms]{theorem!Bando--Mabuchi}
group of automorphisms of a K\"ahler manifold
acts transitively on the set of K\"ahler--Einstein metrics.
For Sasaki--Einstein manifolds\index[terms]{manifold!Sasaki--Einstein}, a similar result was proven
by \index[persons]{Nitta, Y.} Nitta and \index[persons]{Sekiya, K.} Sekiya (\cite{_Nitta_Sekiya_}):
the connected component of the group of CR-automorphisms
of a Sasakian manifold acts transitively on the set of
Sasaki--Einstein metrics compatible with this CR-structure.\index[terms]{structure!CR}
The same question can be asked about Einstein--Weyl
LCK manifolds (that are  a posteriori Vaisman by
\index[persons]{Gauduchon, P.} Gauduchon's \index[terms]{theorem!Gauduchon} \ref{lchk_vai}). 
Let ${\cal U}$ be the set of all Einstein--Weyl\index[terms]{manifold!Einstein--Weyl}
metrics of total volume 1 on a complex manifold $M$,
and ${\cal G}$ a connected component of its group
of holomorphic automorphisms. Does ${\cal G}$ act
transitively on ${\cal U}$? We expect that the
answer is positive, and follows from 
\cite{_Nitta_Sekiya_}.

\hfill

\question
Let $(M, \omega, \theta)$ be a Vaisman manifold
with trivial canonical bundle. From \ref{_CY_Vaisman_Theorem_} \index[terms]{manifold!Vaisman}
it follows that any volume form on $M$ can be realized
as the volume form of a Vaisman metric with  $(M, \omega', \theta')$
with the Lee class \index[terms]{class!Lee}$[\theta']$ proportional to $[\theta]$.
The volume form on a complex manifold determines its Ricci curvature
(\ref{_Futaki_Ricci_Definition_}). Then \ref{_CY_Vaisman_Theorem_} 
implies that a Vaisman manifold with trivial canonical
bundle admits a unique Vaisman metric with a given Lee class
and Ricci curvature equal to any given exact $(1,1)$-form.
The K\"ahler analogue of this statement is the Calabi--Yau theorem
for $c_1(M)=0$. What about $c_1(M)$ negative? Let $L$ be the
weight bundle of a Vaisman manifold,\index[terms]{manifold!Vaisman} and suppose 
$K_M= L^{\otimes \alpha}$, where $\alpha \in \R^{<0}$.
\begin{enumerate}
\item Let $\omega_0 = d^c\theta$ be the standard transversally
K\"ahler form\index[terms]{form!K\"ahler!transversal} Will $M$ admit a Vaisman metric $\omega$ such that
the Ricci curvature of $\omega$ is $- \alpha\omega_0$?
We expect that the answer is affirmative:
the same as for K\"ahler--Einstein manifolds
with anti-ample canonical bundle, there exists
a unique Vaisman metric $\omega$ on $M$ that satisfies
$\Ric_\omega=- \alpha\omega_0$. It should
follow from a transversal version of the\index[terms]{theorem!Aubin--Calabi--Yau}
 Aubin--Calabi--Yau theorem (\cite{_Kacimi_}) and the same 
argument as used to prove \ref{_CY_Vaisman_Theorem_}.

\item
Would a similar result hold for other LCK manifolds?
Suppose  $(M, \omega, \theta)$ is 
an LCK manifold with potential\index[terms]{manifold!LCK!with potential} that satisfies
 $\Ric_\omega=\alpha d^c\theta$, $\alpha \in \R^{<0}$. Is it Vaisman?

\item Let $(M, \omega, \theta)$ be an LCK manifold,
and $L$ its weight bundle. Suppose that $L^{-\alpha}= K_M$,
when $\alpha\in \R^{>0}$. Will it imply that $M$ is LCK with potential?

\item Let $(M, \omega, \theta)$ be an LCK manifold,
and $L$ its weight bundle. Suppose that $L^{-\alpha}= K_M$,
when $\alpha\in \R^{>0}$, and $M$ is LCK with potential. 
 Will it imply that $M$ is
LCK Vaisman?
\end{enumerate}


\section[LCHK and holomorphic symplectic struc\-tures]{LCHK and holomorphic symplectic\\ struc\-tures}


\subsection{LCHK structures}\index[terms]{manifold!LCHK}

Recall that a locally conformally hyperk\"ahler (LCHK)
manifold is a compact Riemannian manifold
equipped with 3 Vaisman structures $(I, \theta)$,
$(J, \theta)$ and $(K, \theta)$, with the same
Lee field\index[terms]{Lee field}, that satisfy the quaternionic
relations $IJ = - JI = K$.
An LCHK manifold $M$ is obtained as
a mapping torus over a 3-Sasakian manifold $S$
(\ref{_LCHK_mapping_torus_Theorem_}).
We write this $M= S\times \R^{>0}/\Z$,
where the $\Z$-action is generated by
$(s, t) \arrow (\phi(s), \lambda t)$
where $\phi$ is a 3-Sasakian automorphism.
Moreover, this correspondence is canonical
(\cite{_Verbitsky:Vanishing_LCHK_}).
Therefore, the LCHK geometry is essentially
the same as 3-Sasakian geometry\index[terms]{geometry!3-Sasaki}.\index[terms]{manifold!3-Sasakian}

For details about 3-Sasakian geometry, 
see \cite{bog} and Chapter \ref{other}.
Recall that {\bf a 3-Sasakian structure}
on a Riemannian manifold $(S,g)$ of dimension
$4n-1$ is defined by a triple of Killing vector fields\index[terms]{vector field!Killing}
$\xi_1, \xi_2, \xi_3$ satisfying the
$\goth{su}(2)$-relations
\[ [\xi_1, \xi_2]= 2\xi_3, [\xi_2, \xi_3]= 2\xi_1,
[\xi_3, \xi_1]= 2\xi_2.
\]
and the CR-structures\index[terms]{structure!CR} $I_1,I_2,I_3$ on the orthogonal complement $\xi_i^\bot$.
We say that $(g,\xi_1, \xi_2, \xi_3, I_1,I_2,I_3)$
is {\bf a 3-Sasakian structure if} each of these vector fields $\xi_i$
is the Reeb field of a Sasakian structure $(S,g, \xi_i^\bot,I_i)$ on 
$(S,g)$, and the complex structures $I_1,I_2,I_3$
preserve the $4n-4$-dimensional distribution
$\xi_1^\bot\cap \xi_2^\bot \cap \xi_3^\bot\subset TS$
and satisfy the quaternionic relations.
Equivalently, a 3-Sasakian structure on $S$
is the same as the hyperk\"ahler structure
on its Riemannian cone\index[terms]{cone!Riemannian} $(S \times \R^{>0})$,
such that the vector field $\frac d {dt}$
acts by homotheties preserving three complex structures
(\ref{_3_Sasaki_Einstein_}). Every 3-Sasakian manifold
is equipped with a 3-dimensional foliation
generated by $\xi_1, \xi_2, \xi_3$; the 
leaves of this foliation are compact by 
Meyers theorem, \index[terms]{theorem!Myers}
because they are locally isometric to $S^3$.
By \index[persons]{Molino, P.} Molino (\cite{_Molino_}), the leaf space
is a Riemannian orbifold. This orbifold
is quaternionic-K\"ahler (\ref{_qK_Remark_}).
Conversely, let $Q$ be a\index[terms]{manifold!quaternionic-K\"ahler} quaternionic-K\"ahler\index[terms]{orbifold!quaternionic-K\"ahler}
orbifold of positive Ricci curvature\footnote{Quaternionic-K\"ahler
manifolds are Einstein, and thus  the Ricci curvature
is constant; if it vanishes, the manifold is
locally hyperk\"ahler, and this case is usually
omitted from the definition of quaternionic-K\"ahler
manifolds. We call a  quaternionic-K\"ahler
manifold {\bf positive} if its Einstein
constant is positive.}. 
Denote by ${\cal B}$ the $\SU(2)$-principal bundle
over $Q$, associated with the set of quaternionic triples
$I,J,K$ (the fibres of ${\cal B}$ freely and transitively
act on the set of quaternionic triples in each tangent space 
$T_q Q$). This is an orbifold bundle; since $Q$ is locally
a quotient of $\R^{4n}$ by a finite group action, the
total space $\Tot({\cal B})$ is also an orbifold.
By \cite[Theorem 13.2.14]{bog}, it is 3-Sasakian.

\hfill

The foremost open question is a classical 
conjecture of quaternionic-K\"ahler geometry.\index[terms]{geometry!quaternionic-K\"ahler}

\hfill

\question
(\index[persons]{LeBrun, C.}LeBrun-\index[persons]{Salamon, S.}Salamon conjecture, \cite{_LeBrun_Salamon:qK_})\\
Let $Q$ be a complete (hence, compact)\index[terms]{conjecture!LeBrun-Salamon}
quaternionic-K\"ahler manifold of positive
Ricci curvature. Then $Q$ is a symmetric space. \index[terms]{manifold!symmetric}

\hfill

\remark
The \index[persons]{LeBrun, C.} LeBrun-Salamon conjecture is clearly false for quaternionic-K\"ahler orbifolds.
In fact, we have many non-homogeneous compact 3-Sasakian manifolds; the corresponding quaternionic-K\"ahler orbifolds
are singular.

\hfill

\remark
From the classification of
symmetric spaces it follows that 
positive quaternionic-K\"ahler symmetric spaces
are Wolf spaces.\index[terms]{space!Wolf}
Recall that {\bf a Wolf space} (\cite{_Wolf:Wolf_}) is
a compact symmetric space of form $G/H$, where
$G$ is a simple group and $H=\Sp(1)\cdot K$,
where $K$ is a compact subgroup which might
intersect $\Sp(1)$ in the central subgroup $\{\pm1\}$.
For each compact simple Lie group, such a subgroup
$H$ exists and is unique to an automorphism.

\hfill

\remark
For \index[persons]{LeBrun, C.} LeBrun-\index[persons]{Salamon, S.}Salamon 
conjecture, it suffices to show 
that all compact\index[terms]{manifold!homogeneous}
quaternionic-K\"ahler manifolds \index[terms]{conjecture!LeBrun-Salamon}
are homogeneous; then they are
symmetric by a result of \index[persons]{Alekseevsky, D. V.} Alekseevsky
(\cite{_Alekseevskii:qK_}),
who has shown that any positive homogeneous
quaternionic-K\"ahler manifold
is a Wolf space.\index[terms]{space!Wolf}

\hfill

Consider a 3-Sasakian manifold
$S$, and let $Q$ be the corresponding qua\-ter\-ni\-o\-nic-K\"ahler 
orbifold. Every Reeb field $r \in \goth{su}(2)$ 
is quasi-regular and induces an $S^1$-action,\index[terms]{action!$S^1$-}
because the corresponding 1-parametric 
subgroup of $\SU(2)$ is compact. Therefore,
its leaf space is
a holomorphic Fano contact orbifold $X_r$, called 
{\bf a twistor space}\index[terms]{twistor space} of $Q$. It is a Fano contact 
orbifold, by \cite{_LeBrun_Salamon:qK_}.\index[terms]{orbifold!contact Fano}
By \cite[Theorem 13.3.1]{bog}, the spaces $X_r$ 
are biholomorphic for different choices of $r$.
As usual, we identify the sphere $S^2$ with 
$\C P^1$ and with the set of all quaternions
$L\in {\Bbb H}$ such that $L^2=-1$.
Let $M$ be the corresponding LCHK manifold, $M=S \times S^1$,
and $\Tw(M)$ its twistor space, fibred with the fiber
$(M,L)$ over each $L\in \C P^1\subset {\Bbb H}$.
Each space $(M, L)$ is a principal elliptic
fibration over $X_r$, where $L\in \C P^1\subset {\Bbb H}$
is the complex structure associated with $r\in \goth{su}(2)$.
This elliptic fibration is generated by the
action of the Lee field\index[terms]{Lee field} on the hyperk\"ahler
cone of the 3-Sasakian manifold. The complexification
of the Lee field action defines the holomorphic action of $T^2$ on $M$
and on $\Tw(M)$. The quotient space $\Tw(M)/T^2$ is 
locally trivially fibred over $\C P^1$ with the fibre $X_r$.
Each fibre of the twistor projection $\Tw(M)\arrow \C P^1$
is a Vaisman manifold.\index[terms]{manifold!Vaisman} Since this fibration is
 locally trivial in the holomorphic category,
the canonical foliations on each Vaisman fibre can be glued
together to give the 1-dimensional foliation on $\Tw(M)$,
also called {\bf the canonical foliation}.\index[terms]{foliation!canonical}
The space $\Tw(M)/T^2$ can be considered to be  the space of leaves
of the canonical foliation on $\Tw(M)$. 

\hfill

\question
Let $M$ be an LCHK manifold.
\begin{enumerate}
\item Is the fibration $\Tw(M)/T^2\arrow \C P^1$ 
with the fibre $X_r$ always trivial, as a complex fibration?
\item 
The twistor fibration $\Tw(M)\arrow \C P^1$ 
of an LCHK manifold is never trivial.
Indeed, by \cite{_Tomberg_}, the twistor space
of a hypercomplex manifold is balanced. Were
$\Tw(M)$ trivial, we would obtain a smooth
holomorphic submersion $\Tw(M) \stackrel p \arrow M$
with fibre $\C P^1$.  Let $\alpha$ be a closed,
strictly positive $(d-1,d-1)$-form on $\Tw(M)$ defining the
balanced metric. Then $p_*(\alpha)$ is a\index[terms]{metric!balanced}
strictly positive $(d-2,d-2)$-form on $M$,
thus $M$ is balanced; this is impossible,
since $M$ is Vaisman (Exercise 
\ref{_Vaisman_non_balanced_Exercise_}).
Let $M$ be an LCHK manifold obtained as a quotient
of the cone $C(S)$ of a 3-Sasakian manifold $S$ by 
the standard $\Z$-action. \index[terms]{action!$\Z$-}Then 
the twistor fibration $\Tw(M)\arrow \C P^1$
is locally trivial as a complex fibration. 
Indeed, the group $\SU(2)$ acts on $C(S)$
rotating the complex structures (\cite{_Swann:hk_}).
However, locally trivial fibrations over  $\C P^1$ are determined by
$H^1(\C P^1, {\cal F})$, where ${\cal F}$ is the sheaf of fibrewise
automorphisms of the fibration. We want to compute this
cocycle for $\Tw(M)$. The automorphism group of 
$\Tw(M)$ contains $\C^*$, and $H^1(\C P^1, \calo^*)=\Z$.
Is the twistor space $\Tw(M)$, as a complex fibration,
isomorphic to the locally trivial fibration
associated with this cocycle?
\end{enumerate}

\question
Let $(M,I,J,K)$ be a hypercomplex manifold.
A {\bf trianalytic subvariety} of $M$
 is a closed subset $Z\subset M$ that is  complex
analytic with respect to $I, J, K$. For any hyperk\"ahler
manifold $(M,I,J,K)$ (not necessarily compact), and
all complex structures of form $L=aI+bJ+cK$, $(a,b,c)\in S^2$
outside of a countable set $R\subset S^2$,\index[terms]{subvariety!trianalytic}
all compact complex analytic subvarieties $Z\subset (M,L)$
are trianalytic (\cite{_V:non-compact_triana_}).
This is false for LCHK manifolds.\index[terms]{manifold!hypercomplex}
Indeed, let $M$ be an LCHK manifold obtained as a quotient
of the cone $C(S)$ of a 3-Sasakian manifold $S$ by 
the standard $\Z$-action.\index[terms]{action!$\Z$-} Then $(M,L)$ is a quasi-regular
Vaisman manifold\index[terms]{manifold!Vaisman!quasi-regular} for all $L=aI+bJ+cK$, and hence  $(M,L)$ 
contains an elliptic curve for any such $L$.
Such subvarieties can be easily constructed
from (for example) embedded Hopf hypercomplex manifolds.
Singularities of trianalytic subvarieties were
classified in \cite{_V:hypercomplex_}. It was shown
that a normalization of a trianalytic subvariety
is always smooth and hypercomplex. However, in all
examples, irreducible trianalytic subvarieties
in hypercomplex manifolds are smooth. Suppose
that $Z\subset M$ is an irreducible trianalytic 
subvariety of an LCHK manifold. Is it smooth?

\hfill

\question
The twistor space of a positive quaternionic-K\"ahler manifold
is a K\"ahler--Einstein holomorphically contact Fano manifold (\cite{bog}, 
\cite{_LeBrun:qK_contact_}). This correspondence 
can be reversed: in \cite[Theorem A]{_LeBrun:qK_contact_},
C. \index[persons]{LeBrun, C.} LeBrun shows that any compact\index[terms]{twistor space}
K\"ahler--Einstein holomorphically contact Fano manifold \index[terms]{manifold!contact Fano}
$X$ is the twistor space of a positive quaternionic-K\"ahler manifold.
His argument uses the  Calabi--Yau metric on the open\index[terms]{metric!Calabi--Yau}
algebraic cone\index[terms]{cone!algebraic} $\tilde M:= \Tot^\circ(K_X)$: he proves that the
holomorphic symplectic form  on $\tilde M$,
induced by the holomorphically contact structure on $X$,
is parallel with respect to the Levi--Civita\index[terms]{structure!holomorphic contact} connection.\index[terms]{connection!Levi--Civita}
The proof uses the stability of vector bundles 
and the  Kobayashi-Hitchin correspondence between\index[terms]{correspondence!Kobayashi-Hitchin}
stable holomorphic structures and the Hermitian\index[terms]{metric!Yang--Mills}
Yang--Mills connections. It seems that all the ingredients
of this proof carry on to the orbifold case. If so,
it should imply that any K\"ahler--Einstein 
holomorphically contact Fano orbifold with\index[terms]{orbifold!contact Fano} 
the smooth Calabi--Yau cone is a twistor\index[terms]{twistor space}
space over a quaternionic-K\"ahler orbifold, and
the corresponding Sasaki-Ein\-stein pseudoconvex\index[terms]{manifold!Sasaki--Einstein}
shell is 3-Sasakian. By \index[persons]{Chen, X.} Chen--\index[persons]{Donaldson, S. K.}Donaldson--\index[persons]{Sun, S.}Sun and \index[persons]{Tian, G.} Tian
theorem\index[terms]{theorem!Chen--Donaldson--Sun, Tian} (\cite{_Chen_Donaldson_Sun_1_,_Chen_Donaldson_Sun_2_,_Chen_Donaldson_Sun_3_,_Tian:K-stability_}) and its orbifold generalization
(\cite{_Ross_Thomas:stability+CSCK_}), 
a holomorphically contact Fano orbifold
admits a K\"ahler-Ein\-stein metric if and
only if it is K-stable. We do not know of any
example of a compact holomorphically contact Fano orbifold
that is  {\em not} K-stable. Are there any? If the\index[terms]{K-stability}
K-stability for holomorphically contact Fano orbifolds
is not automatic, is this notion optimal? We expect
that some weaker notions, such as the \index[persons]{Futaki, A.} Futaki\index[terms]{Futaki invariant}
invariant, might be sufficient in this case.
Note that \index[persons]{LeBrun, C.} LeBrun-\index[persons]{Salamon, S.}Salamon conjecture implies\index[terms]{conjecture!LeBrun-Salamon}
that any {\em smooth} holomorphically contact Fano
manifold is actually a symmetric space.\index[terms]{manifold!symmetric}

\hfill

\question
A $2n+1$-dimensional
contact manifold is {\bf toric of Reeb type}
if it is equipped with an effective action of $T^{n+1}$
by contact automorphisms, tangent to the Reeb field.
All compact toric contact manifolds of Reeb type
are Sasakian (\cite[Corollary 8.5.7]{bog});
these manifolds are called {\bf toric Sasakian}.
Toric contact manifolds of Reeb type
(and hence, Sasakian manifolds)\index[terms]{manifold!toric of Reeb type}
are classified  in terms of \index[persons]{Delzant, T.} Delzant polytopes 
(\cite{_Lerman:toric_contact_}).\index[terms]{Delzant polytope}
A similar description exists for 3-Sasakian
manifolds with a maximal possible torus action.

The notion of a {\bf hyperk\"ahler toric manifold}
was developed by \index[persons]{Bielawski, R.} Bi\-e\-law\-ski and \index[persons]{Dancer, A.} Dancer in 
\cite{_BD:hypertoric_}. These are manifolds
obtained by taking the hyperk\"ahler reduction of
 quaternionic Hermitian vector spaces $V={\Bbb H}^n$ 
over a compact torus acting by linear\index[terms]{reduction!hyperk\"ahler}
quaternionic isometries with $n$-dimensional orbits.
Similarly to the toric manifolds, that are 
classified in terms of \index[persons]{Delzant, T.} Delzant polytopes,
hyperk\"ahler toric manifolds are classified
in terms of the arrangements of affine subspaces.
When such manifolds have an action of $\R^{>0}$
acting by hyperk\"ahler homotheties, the
quotient manifold is compact and 3-Sasakian;
this 3-Sasakian manifold is called {\bf a
toric 3-Sasakian manifold}. Toric 3-Sasakian
manifolds can also be defined as 3-Sasakian
manifolds of real dimension $4n+3$ equipped
with an effective action of a torus $T^{n+1}$
by 3-Sasakian automorphisms 
(\cite[Theorem 6]{_Bielawski:hk_toric_}). 
In \cite[Theorem 3]{_Bielawski:hk_toric_},
\index[persons]{Bielawski, R.} Bielawski classified the toric  3-Sasakian
manifolds in terms of explicit hyperplane arrangements.
By \index[persons]{LeBrun, C.} LeBrun (\cite{_LeBrun:qK_contact_}),
3-Sasakian manifolds are in one-to-one
correspondence with K\"ahler--Einstein
holomorphically contact orbifolds
with smooth symplectization.
A $(2n-1)$-dimensional \index[terms]{orbifold!holomorphically contact}
holomorphically contact orbifold
is called {\bf toric} if it is equipped
with an effective $(\C^*)^{n-1}$-action.
We expect that these manifolds all
admit a K\"ahler--Einstein metric;
then the classification of toric
holomorphically contact orbifolds
with smooth symplectization
would bring the same answers as 
in \cite[Theorem 3]{_Bielawski:hk_toric_}.
Otherwise there exists some new and possibly notable
holomorphically contact orbifolds not admitting
K\"ahler--Einstein metrics.

\subsection[Locally conformally holomorphic symplectic
  structures]{Locally conformally\\ holomorphic symplectic
  structures}
\label{_LCHS_Subsection_}

\question
A {\bf conformally closed $2$-form} is a
2-form $\Phi$ that satisfies $d\Phi=\Phi \wedge \theta$,
where $\theta$ is closed.
Let $(M,I)$ be a complex manifold, 
$L$ a holomorphic line bundle on $M$, and $\Omega\in \Omega^2(M)\otimes L$\index[terms]{manifold!locally conformally holomorphic symplectic}
a non-degenerate $L$-valued 2-form on $M$. 
We say that $(M,I, L, \Omega)$ is {\bf locally
conformally holomorphic symplectic} (LCHS) 
if for any local trivialization $\phi:\; L\arrow \calo_M$,
the corresponding 2-form $\phi(\Omega)\in \Omega^2(M)$ is
conformally closed. Any complex surface is tautologically
LCHS, but for greater dimension this condition 
becomes restrictive. Any LCHK manifold is 
LCHS. Are there any projective LCHS manifolds
in complex dimension $>2$ that are  not holomorphically symplectic?

\hfill

\question
Let $L$ be a holomorphic line bundle
on a complex manifold, and $\Omega$ a non-degenerate
$L$-valued holomorphic symplectic form.
In \cite{_Istrati:twisted_}, such a form
was called {\bf a twisted holomorphic
symplectic form}. \index[persons]{Istrati, N.} Istrati proves that\index[terms]{form!twisted holomorphic symplectic}
on a compact K\"ahler manifold $M$, any twisted holomorphic
symplectic form takes values in a flat unitary
line bundle. This implies that any twisted holomorphic
symplectic form is LCHS.
Is the same true without the K\"ahler assumption?

\hfill

\question
The ``locally conformally holomorphic symplectic
structure'' can be considered to be  a generalization
of the notion of ``complex locally conformally symplectic
structure'' suggested by \index[persons]{Lichnerowicz, A.} Lichnerowicz (\cite[Section 5]{lic1}). 
A complex LCS structure\index[terms]{structure!complex locally conformally symplectic}
is a non-degenerate (2,0)-form $\Omega$ that satisfies
$d\Omega= \Omega\wedge \theta$, where $\theta$
is a closed 1-form, not necessarily real.
This notion is a special case of locally 
conformally holomorphically symplectic
structures: if $L$ is a trivial line bundle,
and $\nabla= \nabla_0 + \theta$, 
the condition $d_\nabla(\Omega)=0$
becomes $d(\Omega)= \Omega\wedge \theta$,
thus $\Omega$ can be considered
as a closed $(L, \nabla)$-valued form. 
The connection $\nabla$ defines a holomorphic
structure operator $\nabla^{0,1}$ on $L$;
since $d_\nabla(\Omega)=0$, the form
$\Omega$, considered to be  an $L$-valued
form, is actually holomorphic.
The LCHS structure with 
the line bundle $L$ equipped with a flat connection,
is the same as a complex locally conformally symplectic
structure of \index[persons]{Lichnerowicz, A.} Lichnerowicz. 
This implies that the three conditions
mentioned above (complex LCS,
LCHS, twisted holomorphically symplectic) are organized
in an hierarchy: complex LCS
is stronger than LCHS, that is  stronger than
twisted holomorphically symplectic. It seems
that locally one can show that this hierarchy
is strict, but on a compact complex manifold 
of dimension $\geq 4$, the strictness is not entirely evident.
Our main example (LCHK manifolds) belongs to all three classes.

Since all complex surfaces are 
LCHS, in dimension 2 the LCHS condition is clearly distinct
from the \index[persons]{Lichnerowicz, A.} Lichnerowicz's ``complex LCS
structures''. In higher dimension, this is not clear.
Indeed, \index[persons]{Istrati, N.} Istrati proves that any twisted 
holomorphic symplectic form on a compact K\"ahler
manifold takes values in a flat unitary line bundle,
and for compact K\"ahler manifolds, these
three classes (twisted holomorphically symplectic,
LCHS, Lichnerowicz) coincide. 
Are these inclusions strict for compact complex manifolds?
Are there any compact twisted  holomorphic symplectic 
manifolds or LCHS manifolds
with  $c_1(L)\neq 0$ in  complex dimension $\geq 4$? 

\hfill

\remark
In dimension 2, there exists only one Vaisman manifold
with trivial canonical bundle, the \index[terms]{surface!Kodaira} Kodaira surface.
It is holomorphically symplectic.
A Vaisman manifold\index[terms]{manifold!Vaisman} $M$ with trivial canonical\index[terms]{manifold!Vaisman!quasi-regular}
bundle is quasi-regular, and the space of leaves
of its canonical foliation\index[terms]{foliation!canonical} is a  Calabi--Yau orbifold\index[terms]{orbifold!Calabi--Yau} $X$
(Exercise \ref{_trivial_canonical_then_QR_Exercise_}).
By Exercise \ref{_Vaisman_never_holo_symple_Exercise_}, it cannot be holomorphically symplectic.

\hfill

\question
Can an LCK manifold with potential\index[terms]{manifold!LCK!with potential} be holomorphically 
symplectic? In view of  Exercise \index[terms]{manifold!holomorphically symplectic}
\ref{_Vaisman_never_holo_symple_Exercise_},
we expect that the answer is negative.
What about other LCK manifolds, can an LCK structure\index[terms]{structure!LCK}
occur on a compact holomorphically symplectic manifold?

\subsection{Holomorphic Poisson structures}\index[terms]{structure!holomorphic Poisson}

\question\label{_holo_Poisson_Question_}
Let $M$ be a compact LCK manifold.\index[terms]{structure!holomorphic Poisson}
We are interested in the existence of holomorphic Poisson
structures that are  of maximal rank somewhere of $M$.
Such Poisson structures are called {\bf non-degenerate}.
Unlike the holomorphic symplectic LCK manifolds, that are  in
short supply (Exercise \ref{_Vaisman_never_holo_symple_Exercise_}), the\index[terms]{manifold!holomorphic symplectic}
holomorphic Poisson structures exist on (primary)\index[terms]{structure!Poisson} 
Hopf surfaces (\cite{_Pym:survey_}).\index[terms]{surface!Hopf!primary}
However, it is not clear which of the Hopf manifolds
admit a non-degenerate Poisson structure in complex dimension $> 2$.
\begin{enumerate}
\item The classical Hopf manifold
$\frac{\C^{2n}\backslash 0}{\lambda \Id}$, $|\lambda| >1$
admits a non-degenerate Poisson structure 
$\sum_{i,j} p_i q_j \frac d {dp_i} \wedge \frac d{dq_j}$,
where $p_i, q_i$ are standard coordinates on $\C^{2n}$. Let 
$H:=\frac{\C^{2n}\backslash 0}{A}$ be a linear Hopf manifold,
with $A$ acting on $\C^{2n}$ by a symplectic homothety.
Will $H$ be holomorphic Poisson? \index[terms]{structure!Poisson}
\item
Let $M$ be a quasi-regular Vaisman manifold\index[terms]{manifold!Vaisman!quasi-regular} 
equipped with a holomorphic Poisson structure,
and $X$ the space of leaves of its canonical foliation,
that is  a projective orbifold \index[terms]{foliation!canonical}
(\ref{_Structure_of_quasi_regular_Vasman:Theorem_}).
Denote by $\pi:\; M \arrow X$ the standard projection.
By \ref{_holo_tensor_on_Vaisman_Lee_invariant_Theorem_}, 
the Poisson bivector $\sigma$ is Lee-invariant.
This should imply that 
$\sigma= \pi^* \sigma_1 + \rho \wedge \theta^\sharp$, 
where  $\sigma_1$ is a holomorphic Poisson bivector on\index[terms]{Poisson bivector}
$X$, $\theta^\sharp$ the Lee field\index[terms]{Lee field}, and $\rho\in TX$
a holomorphic vector field satisfying $\Lie_\rho\sigma_1=0$.
Is this true? Can we classify quasi-regular Vaisman
manifolds admitting non-degenerate holomorphic
Poisson structures in terms of the holomorphic
Poisson structures\index[terms]{structure!Poisson} on K\"ahler orbifolds?
Suppose that $\sigma$ is non-degenerate.
Does this imply that $X$ is a weak Fano\footnote{A 
variety $M$ is called\index[terms]{variety!weak Fano}
{\bf weak Fano} if its anti-canonical bundle $K_M^{-1}$ is
nef and big, that is, has non-negative degree on all
curves, and the dimension $\dim H^0(M, K_M^{-N})$ satisfies
 $\lim\sup_{N\to \infty}\frac{\dim H^0\left(M,
  K_M^{-N}\right)}{N^{\dim_\C M}} >0$.}  orbifold? Note that holomorphic symplectic structures
are incom\-pa\-ti\-ble with Vaisman structures, as
explained in Subsection \ref{_LCHS_Subsection_},
hence the Poisson structure on a Vaisman\index[terms]{structure!holomorphic Poisson} manifold\index[terms]{manifold!Vaisman}
necessarily has singularities.
\end{enumerate}

\question 
Let $M$ be a Vaisman manifold equipped with a non-degenerate
Poisson structure. Can $M$ be irregular?

\hfill

\question
Can a non-Vaisman LCK manifold with potential\index[terms]{manifold!LCK!with potential}
support a non-degenerate holomorphic Poisson structure?\index[terms]{structure!holomorphic Poisson}
What about other LCK manifolds?


\section{Riemannian geometry of LCK manifolds}


\subsection{Curvature of Vaisman manifolds}\index[terms]{manifold!Vaisman}
\label{_Vaisman_curvature_Subsection_}

Let $R\in \Sym^2(\Lambda^2 M)$ 
 denote the Riemann curvature tensor. The holomorphic bisectional
curvature $R(X, IX, X, IX)$ of a Vaisman manifold
vanishes whenever $X$ is\index[terms]{holomorphic bisectional curvature}
tangent to the canonical foliation,\index[terms]{foliation!canonical} because it is
flat and totally geodesic. Using the O'Neill\index[terms]{O'Neill formulas}
formulas, or computing the curvature explicitly
as in Exercise \ref{_holo_curva_Heisenberg_Exercise_},
it is actually possible to show that the curvature
of the Sasakian manifold and of the Vaisman manifold
in the transversal direction is always smaller
than the curvature of the corresponding K\"ahler leaf
space. Therefore, the curvature bounds on Sasaki and
Vaisman manifold should have
the same or stronger consequences as
for K\"ahler manifolds. For K\"ahler manifolds,
curvature bounds bring bountiful fruits, such as the celebrated Frankel
conjecture, solved by \index[persons]{Siu, Y.-T.} Siu-\index[persons]{Yau, S.-T.}Yau for strictly\index[terms]{conjecture!Frankel}
positive curvature and \index[persons]{Mok, N.} Mok for semi-positive, \cite{_Mok:Frankel_}. 
A work in this direction was already started for Sasakian
manifolds (\cite{_He_Sun:Sasaki_Frenkel_}),
but for Vaisman manifolds\index[terms]{manifold!Vaisman} nothing was done yet.

Since Frankel formulated his famous conjecture
\cite{_Frankel:conjecture_}, there were many results
relating (semi-)positive sectional and bisectional
holomorphic curvature with algebraic geometry\index[terms]{geometry!algebraic}.\index[terms]{holomorphic sectional curvature}
Frankel conjectured that any compact K\"ahler
manifold with positive bisectional curvature is
biholomorphic to $\C P^n$; this was proven by
Siu and \index[persons]{Yau, S.-T.} Yau in \cite{_Siu_Yau:Frankel_}.
Then, \index[persons]{Mok, N.} Mok proved a stronger version of 
Frankel conjecture, and showed that any
compact K\"ahler manifold with non-negative
holomorphic bisectional curvature is
(up to a finite covering) biholomorphic
to a product of tori and compact \index[terms]{manifold!symmetric}
K\"ahler symmetric spaces (\cite{_Mok:Frankel_}). For Sasakian manifolds,
an analogue of this result was proven in 
\cite{_He_Sun:Sasaki_Frenkel_}, by Weiyong He and Song
Sun. They have shown that any simply connected compact
Sasakian manifold with non-negative transversal holomorphic bisectional
curvature is either irregular and is isomorphic
to a sphere with the Sasakian metric, or is
quasi-regular, and has a product of locally
symmetric K\"ahler orbifolds as its Reeb orbit space.
The analogue of this result for Vaisman manifolds
with non-negative transversal holomorphic bisectional
curvature should be also true.\index[terms]{manifold!Vaisman}

\hfill

\question
Let $M$ be a compact Vaisman manifold with holomorphic bisectional
curvature that is  positive in the\index[terms]{manifold!Hopf}
transversal direction. Is $M$ biholomorphic to a Hopf manifold?
Conversely, let $M$ be a Hopf manifold of Vaisman type. Will it admit a
Vaisman metric of non-negative holomorphic bisectional
curvature, positive in the transversal direction? 

\hfill

\question
Let $M$ be a compact Vaisman manifold\index[terms]{manifold!Vaisman} with  non-negative holomorphic bisectional
curvature. Suppose that $M$ is not a finite quotient of a
Hopf manifold.  Will it follow that $M$ is a
suspension over a quasi-regular Sasakian manifold with
non-negative transversal holomorphic bisectional
curvature, such as described in\index[terms]{manifold!Sasaki!quasi-regular} \cite{_He_Sun:Sasaki_Frenkel_}?

\hfill
 
Another question relating the complex geometry and the
curvature is due to S.-T. \index[persons]{Yau, S.-T.} Yau, \cite{_Yau:problems_82_}.
Recall that a projective manifold $M$ is called {\bf
rationally connected} if there is a compact family\index[terms]{manifold!rationally connected}\index[terms]{manifold!Fano!weak}
${\goth S}$ of rational curves such that for any $x, y \in M$
there exists a curve $C \in {\goth S}$ passing through $x$
and $y$. It was long conjectured that $M$ is rationally
connected if and only if $M$ is birationally equivalent to
a weak Fano manifold (see \ref{_holo_Poisson_Question_} 
for the definition of a weak Fano).
In \cite{_Yau:problems_82_}, Yau conjectured
that any compact K\"ahler manifold with positive
holomorphic sectional curvature is rationally connected.
This conjecture was proven (for manifolds with a semi-positive
holomorphic sectional curvature strictly positive at some point)
in \cite{_Heier_Wong:positive_RC_} for projective and in 
\cite{_Yang:RC-pos_} for compact K\"ahler manifolds.
In \cite{_Matsumura:positive_}, \index[persons]{Matsumura, S.-I.} Matsumura proved a version
of this  conjecture without the assumption of strict
positivity: he has shown that any projective manifold with a semi-positive
 holomorphic sectional curvature admits a holomorphic
map to a compact torus with rationally connected fibres.
He also proved that a projective manifold with a (semi-)positive
 holomorphic sectional curvature is rationally connected,
if there is a point $x\in M$ such that no tangent vector 
$v\in T_x M$ is truly flat.\footnote{A tangent vector 
$v\in T_x M$ is {\bf truly flat} if $R(v, \cdot, \cdot,
  \cdot)=0$, where $R$ is the curvature tensor.}
Note that the Lee field\index[terms]{Lee field} on a Vaisman manifold\index[terms]{manifold!Vaisman}
is ``truly flat'' in this sense.\index[terms]{vector field!truly flat}

\hfill

\question
Let $M$ be a Vaisman manifold with 
semi-positive holomorphic sectional curvature.
Assume that there is a point $x\in M$
such that the Lee field\index[terms]{Lee field} $\theta^\sharp\restrict{T_x M}$ 
is the only truly flat direction in $x$.
Can we prove that $M$ is a deformation of
an elliptic fibration over a rationally 
connected orbifold?

\hfill

I. \index[persons]{Vaisman, I.} Vaisman has shown that the positivity 
of holomorphic sectional curvature on an LCK manifold can be 
quite restrictive (\cite{va_tr2}).
Let $M$ be a complete LCK manifold
with positive holomorphic sectional curvature
that satisfies \[ \frac{R(X, IX, X, IX)}{|X|^4} \geq
\frac{|\theta|^2}{4}+\delta,\]
where $\delta >0$. Then $M$ is compact, simply connected,
and thus globally conformally K\"ahler.
Non-positive holomorphic sectional curvature
is also restrictive: in \cite{_CCN:holo_sec_},
Chen, Chen and Nie \index[persons]{Chen, H.}\index[persons]{Chen, L.}\index[persons]{Nie, X.}prove that a compact LCK manifold with
constant negative holomorphic sectional curvature 
is K\"ahler.

\hfill

\question
The Heisenberg-type\index[terms]{nilmanifold}
Vaisman nilmanifold has vanishing  holomorphic sectional
curvature (Exercise \ref{_holo_curva_Heisenberg_Exercise_}).
Let $M$ be a compact Vaisman manifold,\index[terms]{manifold!Vaisman} and $\Sigma$
its canonical foliation.\index[terms]{foliation!canonical} As usual, we consider $\Sigma$
as a transversal K\"ahler foliation and treat the
leaf space of $\Sigma$ as a K\"ahler manifold.
Assume that the
holomorphic sectional curvature of the transversal K\"ahler
structure on this leaf space vanishes.
Can we prove that $M$ is biholomorphic
to a finite quotient of the Heisenberg-type
Vaisman nilmanifold?\index[terms]{holomorphic sectional curvature}

\hfill

\question
In \cite{_Gray_Vanhecke_},
\index[persons]{Gray, A.} Gray and \index[persons]{Vanhecke, L.} Vanhecke treated almost complex Hermitian 
manifolds with the holomorphic sectional curvature
$\frac{R(X, IX, X, IX)}{|X|^4}$ constant pointwise. The examples
include $\C P^n$ and $S^6$ with the standard (non-integrable)
almost complex structure. In Riemannian geometry,
if the sectional curvature\index[terms]{holomorphic sectional curvature!pointwise constant}
$\frac{R(X, Y, X, Y)}{|X\wedge Y|^2}$
is pointwise constant, the manifold is a space form,
and the sectional curvature is constant on the manifold.
The same is true for the holomorphic sectional curvature
on a K\"ahler manifold, by a result of \index[persons]{Bochner, S.} Bochner
(\cite[Vol. II, Theorem 7.5]{KoNFD2}, 
\cite[Theorem 5]{_Bochner:curvature_Hermitian_}).
On the contrary, for almost complex Hermitian 
manifolds, a pointwise constant holomorphic sectional curvature
can be non-constant.
\index[persons]{Gray, A.} Gray and \index[persons]{Vanhecke, L.} Vanhecke found examples of complex
Hermitian structures with the holomorphic sectional curvature 
constant pointwise, but not constant on the manifold.
For a Vaisman manifold,\index[terms]{manifold!Vaisman} the holomorphic sectional curvature
vanishes in the Lee field\index[terms]{Lee field} direction, and hence 
the constancy of the holomorphic sectional curvature
implies that it vanishes, and (by polarization)
the holomorphic bisectional curvature also vanishes.
Vanishing of the holomorphic bisectional curvature 
for K\"ahler manifold implies vanishing of the curvature,
by the same result of \index[persons]{Bochner, S.} Bochner. Are there any Vaisman
manifolds with constant (necessarily vanishing)
holomorphic bisectional curvature? From 
O'Neill formulas\index[terms]{O'Neill formulas} it should follow that the
corresponding transversal K\"ahler manifold
is negatively curved.

\hfill

\question  
Let $M$ be a Vaisman manifold.\index[terms]{manifold!Vaisman} Then the
holomorphic sectional curvature vanishes
in the Lee direction. Since an LCK manifold
with potential can be obtained as a small
deformation of a Vaisman, the holomorphic 
sectional curvature in the Lee direction
should oscillate close to zero. Let $M$
be a compact LCK manifold with potential.\index[terms]{manifold!LCK!with potential} 
 Can the holomorphic 
sectional curvature be positive or negative
everywhere on $M$ for some choice of
the conformal gauge?
When the conformal gauge of
$M$ is Gauduchon, we expect that the
answer is negative.

\hfill

\remark
By Exercise \ref{_holo_curva_Heisenberg_Exercise_},
the sectional curvature of the Sasakian metric on the
Heisenberg manifold $W$ is non-positive.\index[terms]{sectional curvature}
This implies the non-positivity of the\index[terms]{manifold!Heisenberg}
sectional curvature of the Kodaira surface $M$ with the 
corresponding Vaisman metric. Indeed, $M$ is locally
isometric to $W\times \R$. However, it is impossible
to find a metric of strictly negative sectional \index[terms]{surface!Kodaira}
curvature on $M$, or on any other Vaisman manifold.\index[terms]{manifold!Vaisman}
Indeed, the universal cover of a manifold
of negative sectional curvature is contractible
by the \index[persons]{Cartan, E.} Cartan--Hadamard theorem, and hence  
a Sasakian manifold $W_1$ with finite fundamental
group does not give a Vaisman manifold $M_1$ with
negative sectional curvature: its universal cover
$\tilde M_1$ is homeomorphic to a compact manifold times $\R$, and 
thus $\tilde M_1$ cannot be contractible. For all 
other Sasakian manifolds $W_1$, the Vaisman manifold
is homeomorphic to a mapping torus of an isometry
$\phi$ of $W_1$. Since $\phi$ is an isometry,
a power $\phi^n$ is isotopic to identity for 
some $n\in \Z^{>0}$. Passing to a finite cover of $M_1$,
we may assume that $\phi$ is isotopic to identity, and
$M_1$ is homeomorphic to a product of $S^1$ and $W_1$.
This implies that $M_1$ contains a torus $T^2=S^1\times S^1$
such that $\Z^2=\pi_1(T^2)$ is a subgroup in $\pi_1(M)$.
Then, $M_1$ cannot admit a metric of negative sectional
curvature (\cite{_Preissmann_}, \index[terms]{theorem!Preissmann}
\cite[Theorem 9.3.3]{_Burago_Burago_Ivanov_},
\cite[Corollaire 7.3]{_CDP:geometrie_groupes_}).
The reason is that a manifold of negative sectional
curvature is \index[persons]{Gromov, M.} Gromov hyperbolic, Gromov hyperbolic
manifolds have Gromov hyperbolic fundamental groups\index[terms]{fundamental group},
and a Gromov hyperbolic group cannot contain $\Z^2$.
The fundamental group of the OT manifold also
contains $\Z^2$, and hence  it cannot be\index[terms]{manifold!Gromov hyperbolic}
Gromov hyperbolic either; the same\index[terms]{fundamental group!Gromov hyperbolic}
is true for Seifert 3-manifolds, which implies that  
they cannot be hyperbolic.

\hfill

\question
Can a  manifold admit both a strict LCK structure\index[terms]{structure!LCK} and a
metric of negative sectional curvature?\index[terms]{sectional curvature}
Based on all examples considered, we expect that it is impossible.
Note that there are many K\"ahler manifolds of negative \index[terms]{Riemann surface}
sectional curvature, such as the compact Riemann surfaces of genus $>1$.

\subsection{Special Hermitian metrics on LCK manifolds}

\question (Dan \index[persons]{Popovici, D.} Popovici)\label{_Popovici_sGauduchon_} Let
$(M,I,g,\omega)$ be a complex Hermitian manifold of
complex dimension $n$. In \cite{_Popovici:SG_}, the metric $g$ is called {\bf strongly
Gauduchon}  if $\6\omega^{n-1}$ is $\bar\6$-exact.  
For manifolds that have $dd^c$-lemma, such as Moishezon and
Fujiki class C manifold, the Gauduchon metric is \index[terms]{manifold!Moishezon}
strongly Gauduchon, but $dd^c$-lemma does not\index[terms]{metric!strongly Gauduchon}\index[terms]{manifold!Fujiki class C}
hold for strict LCK manifolds (\ref{_theta_not_d^c_closed_Corollary_}). 
Moreover, the Vaisman manifolds\index[terms]{manifold!Vaisman} do not admit strongly Gauduchon \index[terms]{metric!strongly Gauduchon}
metrics (Exercise \ref{_Vaisman_all_are_not_sG_Exercise_}). 
Can a compact strict LCK manifold admit a strongly Gauduchon metric?

\hfill

\question\label{_SKT_balanced_LCK_Question_}
Let $(M, \omega)$ be a compact 
Hermitian manifold, $\dim_\C M >2$, 
not admitting a K\"ahler structure.\index[terms]{metric!pluriclosed (SKT)}\index[terms]{metric!balanced}
Recall that it is called {\bf SKT} (from ``strongly K\"ahler torsion'')
or {\bf pluriclosed} if $dd^c \omega=0$, and {\bf  balanced}
if $d(\omega^{\dim_\C M -1})=0$. It is known that
these conditions are complementary, when $\dim_\C M >2$, i.e. the same
metric cannot satisfy two of them, unless it is K\"ahler
(\cite{uga,_FPS:6dim_,_Alexandrov_Ivanov_}; Exercises 
\ref{_LCK_balanced_Kahler_Exercise_},
\ref{_LCK_SKT_Kahler_Exercise_},
\ref{_balanced_SKT_Kahler_Exercise_}). 
A locally
homogeneous nilmanifold\index[terms]{nilmanifold}
with a locally homogeneous SKT metric\index[terms]{metric!SKT}
cannot admit a locally homogeneous balanced metric
(\cite{_Fino_Verzzoni:SKT_}).\footnote{The proof given in
\cite{_Fino_Verzzoni:SKT_} was faulty, and applies only
to 2-step nilpotent manifolds 
(\cite{_Fino_Vezzoni:correction_}). However,
as shown in \cite{_Arroyo_Nicolini_}, any SKT nilmanifold is 
2-step nilpotent, and hence  the result of \cite{_Fino_Verzzoni:SKT_}
is valid. We are grateful to Alexandra Otiman for this observation.}
\begin{enumerate}
\item 
Can an LCK manifold admit an SKT metric? 
Can an LCK manifold admit a balanced metric?
For LCK manifolds with potential,\index[terms]{manifold!LCK!with potential} 
this is impossible by Exercise \ref{_LCK+pot_not_balanced_Exercise_}.
Also, some of the toric Kato manifolds cannot admit
a pluriclosed Hermitian metric, as shown in 
\cite{_IOPR:toric_Kato_}.\index[terms]{manifold!toric Kato}

\item
The pluriclosed (a. k. a. SKT) condition can be 
generalized: let's call a Hermitian form $\omega$ {\bf $p$-pluriclosed}\index[terms]{form!$p$-pluriclosed}
if $dd^c(\omega^p)=0$. An $(n-2)$-pluriclosed $n$-manifold is 
called {\bf astheno-K\"ahler}. It is not hard to see that
a metric that is  LCK, $p$-pluriclosed and $q$-pluriclosed
at the same time, $0 <p< q < n-1$ is also K\"ahler
(Exercise \ref{_LCK_p_pluriclo_Exercise_}). 
However, the $p$-pluriclosedness is not exclusive
with other special Hermitian conditions; for example,
the same nilmanifold can admit a balanced and an astheno-K\"ahler metric
(\cite{_Fino_Gra_Vezzoni:Astheno_,_Lattore_Ugarte_}).
An LCK metric cannot be 1-pluriclosed (that is, SKT), as shown 
in \cite{_Ivanov_Papadopoulos_,_Alexandrov_Ivanov_}.
It is interesting to find examples of LCK manifolds
admitting $p$-pluriclosed Hermitian structures, or to 
prove that a certain classes of LCK manifolds
do not admit a $p$-pluriclosed Hermitian structures.
A Vaisman manifold does not admit $p$-pluriclosed \index[terms]{form!$p$-pluriclosed}
metric (Exercise \ref{_Vaisman_no_dd^c_closed_pos_p_p_Exercise_}; 
see also \cite{_Angella_Otiman_}).
\end{enumerate}

\question
By \cite{_Angella_Otiman_}, 
a Vaisman manifold $M$ with $\dim_\C M >2$ \index[terms]{manifold!Vaisman}cannot admit an SKT metric.
Can an LCK manifold  with potential $M$,  $\dim_\C M >2$,\index[terms]{manifold!LCK!with potential} admit an SKT metric? \index[terms]{metric!SKT}

\hfill

\question
All known examples of LCK manifolds with non-irreducible
Weyl holonomy have flat Weyl holonomy 
(\cite{_Madani_Moroianu_Pilca:Einstein_,_Madani_Moroianu_Pilca:Weyl_}). Are there any LCK\index[terms]{holonomy!Weyl}
manifolds with non-flat Weyl connection and reducible
holonomy?

\hfill

\question
A Riemannian manifold is called ``geometrically formal'' if
the wedge product of two harmonic\index[terms]{manifold!geometrically formal}
forms is harmonic (\cite{dieter}). It is proven in \cite{Kas} that the LCK
metric constructed by \index[persons]{Oeljeklaus, K.} Oeljeklaus and \index[persons]{Toma, M.} Toma on a\index[terms]{manifold!formal}
3-dimensional OT manifold of type (2,1) is geometrically
formal, and thus this particular OT manifold is also
formal in the sense of rational homotopy theory
(that is, its de Rham algebra is quasi-isomorphic
to its cohomology algebra).  In \cite{opi} it was shown that
Vaisman manifolds\index[terms]{manifold!Vaisman} are geometrically formal if and only if
they are cohomological Hopf manifolds. Are there any other
geometrically formal LCK manifolds?

\hfill 

\question The relation between LCK geometry and spin\index[terms]{geometry!spin}
geometry has not yet been investigated. Is there any kind
of spinorial equation which forces a Hermitian metric to
be Vaisman? This seems natural, taking into account that
Sasakian manifolds carry a canonical $\Spin^c$ structure
with a Killing spinor (see \cite[Theorem 4.1, Corollary 4.2]{_moroianu:spin_}), and\index[terms]{spinor!Killing}
Sasakian-Einstein structures correspond to real Killing spinors
(\cite{bfgk}). See also \cite{nakad} where different
classes of CR-manifolds are characterized in terms of
$\Spin^c$ structures.

\hfill

\question
Recall that a real (1,1)-form on a 
complex manifold is called {\bf pseu\-do-K\"ahler}
if it is nowhere degenerate and closed. In this case, it defines a pseudo-Riemannian\index[terms]{form!pseudo-K\"ahler}
metric, and its signature is constant. In \cite{_Petean:surfaces_}, J. \index[persons]{Petean, J.} Petean obtained a classification of
pseudo-K\"ahler metrics on complex surfaces. In particular, he has constructed
a locally invariant pseudo-K\"ahler metric on a Kodaira surface (that is  Vaisman), see Exercise \index[terms]{metric!pseudo-K\"ahler}
\ref{_nilma_Heisenberg_forms_Exercise_}.\index[terms]{surface!Kodaira}
It is not hard to see that a manifold with $H^{1,1}(M)=0$ does not admit
a pseudo-K\"ahler metric; in particular, a Hopf manifold is never pseudo-K\"ahler.
What other Vaisman and LCK manifolds admit (and do not admit) pseudo-K\"ahler metrics?
Can a compact Vaisman manifold\index[terms]{manifold!Vaisman} $M$, $\dim_\C M >2$ admit a pseudo-K\"ahler metric?\index[terms]{metric!pseudo-K\"ahler}
We do not seem to have an example. The Heisenberg-type Vaisman nilmanifold\index[terms]{nilmanifold}
of dimension $>2$ does not admit a pseudo-K\"ahler metric by Exercise 
\ref{_nilma_Heisenberg_forms_Exercise_}.

\hfill

\question
Can an OT manifold admit a pseudo-K\"ahler metric? 
We expect the answer to be negative.\index[terms]{manifold!Oeljeklaus--Toma (OT)}

\hfill

\question
Recall that
an LCK manifold $(M,J,g,\theta)$ is called {\bf pluricanonical}
if its Lee form \index[terms]{form!Lee} satisfies  $(\nabla \theta)^{1,1}=0$,\index[terms]{manifold!pluricanonical}
where $\nabla$ is the Levi Civita connection of $g$. This
is a local condition. A compact pluricanonical LCK
manifold is Vaisman, see \ref{_pluri_vais_}. Is there any
example of a pluricanonical 
non-Vaisman metric on a non-compact complex manifold?

\hfill

\question  
Let $(M,I,g)$ be a Hermitian manifold and $\nabla$ a\index[terms]{connection!K\"ahler-like}\index[terms]{connection!Strominger K\"ahler-like (SKL)}
metric connection with curvature tensor $R^\nabla$. In 
\cite{_Angella_Otal_Ugarte_Villacampa_}, the connection $\nabla$
is called {\bf K\"ahler-like} if it satisfies the first
Bianchi identity, $\sum_{cycl}R^\nabla(x,y)z=0$, and is of
type $(1,1)$, i.e. $R^\nabla(Ix,Iy)=R^\nabla(x,y)$, for all $x,y,z\in TM$. In \cite{_Yau_Zhao_Zheng_}, a Hermitian manifold whose \index[persons]{Bismut, J.-M.} Bismut
connection is K\"ahler-like is called {\bf \index[persons]{Strominger, A.} Strominger
	K\"ahler-like} (SKL for short). Let $(M,I,g)$ be a
SKL manifold of dimension $\geq 3$. \index[persons]{Yau, S.-T.} Yau, Zhao and Zheng\index[terms]{conjecture!Yau, Zhao and Zheng}
conjectured that it cannot admit LCK metrics. They already
proved in \cite{_Yau_Zhao_Zheng_} that it cannot admit
Vaisman metrics. One can also ask the same question about
 the other Hermitian connections in the 1-parameter family constructed by
\index[persons]{Gauduchon, P.} Gauduchon (see Subsection \ref{_Bismut_connection_Note_}).
Suppose that one of these connections on an LCK manifold $M$
is K\"ahler-like, will it imply that $M$ is K\"ahler?\index[terms]{connection!Bismut}



\section{LCK reduction and LCS geometry }\index[terms]{geometry!LCS}


\subsection{The Lee cone of taming LCS structures}

\remark
Let $\omega$ be a 2-form on a complex or almost
complex manifold $(M,I)$. We say that
{\bf $\omega$ tames $I$}, or {\bf $I$ tames $\omega$}
if $\omega(x,Ix)>0$ for any non-zero tangent vector.
A complex manifold tamed by a symplectic form
is called {\bf Hermitian symplectic};
this condition is equivalent to having
a positive definite Hermitian form $\omega^{1,1}$ which 
is (1,1)-part of a closed\index[terms]{form!pluriclosed (SKT)} one.\index[terms]{manifold!Hermitian symplectic}
The form $\omega^{1,1}$ is pluriclosed, but
pluriclosedness is a weaker condition.\index[terms]{form!pluriclosed}
For instance, all complex surfaces admit\index[terms]{metric!pluriclosed (SKT)}
pluriclosed metrics (for complex surfaces,
pluriclosed is the same as Gauduchon),
but the only surfaces admitting 
Hermitian symplectic structures are
K\"ahler. This is easy to see because
all non-K\"ahler surfaces admit an
exact, positive (1,1)-current\index[terms]{current!positive} $\Theta$
(\ref{_Lamari-current-intro_Theorem_}). Then the integral
$\int \Theta\wedge \omega= \int \Theta\wedge \omega^{1,1}$ is positive
whenever the form $\omega^{1,1}$ is Hermitian, but
$\int \Theta\wedge \omega=0$, because $\omega$ is closed
and $\Theta$ is exact.
Streets and \index[persons]{Tian, G.} Tian (\cite{_Streets_Tian_}) conjectured that in 
$\dim_\C M >2$, Hermitian symplectic\index[terms]{conjecture!Streets-Tian}
implies K\"ahler; this is unknown even in dimension 3
(\cite{_Verbitsky:twistor_}).
For complex nilmanifold,\index[terms]{nilmanifold} this question was considered
in \cite{_Enrietti_Grancharov_Fino_} and \cite{_Enrietti_Fino_Vezzoni_}.
In \cite{_Enrietti_Fino_Vezzoni_}
it was shown that complex nilmanifolds cannot admit
Hermitian symplectic metrics, with the exception of a torus.
However, the pluriclosed metrics exist 
on many complex nilmanifolds.\index[terms]{nilmanifold}

\hfill

\question
A similar question to the one asked by Streets-\index[persons]{Tian, G.}Tian
makes sense for LCK manifolds.
Let $\omega$ be an LCS form taming a complex\index[terms]{form!LCS}
structure on a complex manifold $(M,I)$. Is
$(M,I)$ an LCK manifold? For complex surfaces, this question
was solved in negative in \cite{_Apostolov_Dloussky_},
who showed that any  complex surface admits an LCS form taming
the complex structure. For greater dimension, the 
question is still open: we know that any
almost complex manifold with $b_1>0$
admits an LCS form (\cite{em}),\index[terms]{form!LCS}
but we do not know if this LCS form is taming.

\hfill

\question
By \cite{_Bertelson_Meigniez_}, a compact manifold
$M$ with $b_1(M) >0$ that admits a non-degenerate 2-form 
carries an LCS form $\omega$ with any non-zero Lee class\index[terms]{class!Lee}
$[\theta]\in H^1(M)$ and any real Morse--Novikov class\index[terms]{class!Morse--Novikov}
$[\omega]\in H^2_\theta(M)$. 
Let $(M,I)$ be a complex manifold tamed by an LCS form.
What are the restrictions on the cohomology invariants
of this form? \index[persons]{Apostolov, V.} Apostolov and \index[persons]{Dloussky, G.} Dloussky partially solved this problem
for complex surfaces (\cite{ad2,ad3}), showing that any complex surface
is tamed by an LCS form, and computed the set of Lee
classes of taming LCS structures for all\index[terms]{structure!LCS} compact\index[terms]{class!Morse--Novikov}
complex surfaces with $b_1=1$.\index[terms]{class!Lee}

\hfill

\question
For Vaisman manifolds \index[terms]{manifold!Vaisman}
and LCK manifolds with potential\index[terms]{manifold!LCK!with potential}
the set of Lee classes of LCK structures
(``the Lee cone of LCK structures'')\index[terms]{structure!LCK}
is a half-space in $H^1(M,\R)$\index[terms]{cone!Lee}
(\ref{_Lee_cone_on_Vaisman_Theorem_}, 
\ref{_Lee_cone_on_LCK-pot_Theorem_}).\index[terms]{cone!LCS!taming}
By Exercise \ref{_taming_Lee_cone_Vaisman_Exercise_}, the 
set of Lee classes\index[terms]{class!Lee} of LCS structures
taming the complex structure of Vaisman type
(``the Lee cone of
taming LCS structures'') is equal to
the Lee cone of LCK structures. 
For an LCK Inoue surface, the Lee cone
of LCK structures\index[terms]{structure!LCK} is a point. As shown 
by \index[persons]{Apostolov, V.} Apostolov-Dloussky, the taming LCS
cone is also a point, the same as for
LCK cone. Does the Lee cone of LCK
structures always coincide with the Lee cone of
taming LCS structures?\index[terms]{surface!Inoue}

%
%

%

\subsection{Twisted Hamiltonian action, LCS reduction and
  toric Vaisman geometry}\index[terms]{action!twisted Hamiltonian}\index[terms]{reduction!LCS}

In this short expos\'e we follow  \cite{hr1}.
Let $(M, \omega, \theta)$ be an LCS manifold,\index[terms]{manifold!LCS}
and $\Symp_0(M)$ the connected component of the 
group of conformally symplectic diffeomorphisms of $M$.
We consider $\omega$ as a closed, non-degenerate 2-form
taking values in the flat line bundle $L$, called
{\bf the weight bundle}.\index[terms]{bundle!weight}

The Lie algebra $\goth{symp}(M)$ of $\Symp_0(M)$ is the algebra
of conformally symplectic vector fields.
We identify it with the space of vector
fields that are  lifted to symplectic
homotheties of the symplectic cover 
$(\tilde M, \tilde \omega)$ of $M$.\index[terms]{Lee homomorphism!extended}
The homothety coefficient defines a character
$\underline{\chi}:\; \goth{symp}(M) \arrow \R$.
The {\bf extended Lee homomorphism} 
is the character $\chi:\; \Symp_0(M)\arrow \R^{>0}$
obtained by the integration of the
character $\underline \chi$ to the Lie group. 
By \cite[Lemma 4.7]{hr1}, the kernel
$\ker \chi$ is arc connected; we denote it by
$\Symp_\chi(M)$.
A vector field $v$ in the Lie algebra of
$\Symp_\chi(M)$ acts on $(\tilde M, \tilde \omega)$
by symplectomorphisms, thus  corresponding 
to a closed 1-form $\tilde \eta_v:=i_v(\tilde \omega)$ on $\tilde M$.
Then $\tilde \eta_v$ is the pullback of the closed $L$-valued 1-form on $M$,
$\eta_v=i_v(\omega)$.
We denote the Lie algebra of such vector fields by
$\goth{symp}_\chi(M)\subset TM$.

For any $u_1 \in \Symp_\chi(M)$, consider
a path $u_t\in \ker \chi$, connecting $u_1$
to $u_0 =\Id_M$. By the Leibniz-Newton formula, $u_1$
can be obtained by integrating the vector fields
$\dot u_t\in \goth{symp}_\chi(M)$.
The {\bf twisted flux homomorphism} $F$
is the map from the universal cover\index[terms]{flux homomorphism!twisted}
$\widetilde \Symp_\chi(M)$ of $\Symp_\chi(M)$ to 
the Morse--Novikov cohomology group
$H^1_\theta(M)$ taking $u_1$
to $\int_0^1 [i_{\dot u_t}(\omega)] dt$,
where $i_{\dot u_t}(\omega)$ is 
the closed $L$-valued 1-form
constructed above and $[i_{\dot u_t}(\omega)] $
its cohomology class in $H^1_\theta(M)$.

When $\theta$ is exact, this construction
gives the standard flux homomorphism,
\cite{_McDuff:flux_,_Ono:flux_}.\index[terms]{flux homomorphism}

Recall that a vector field $v\in \goth{symp}_\chi(M)$
is called {\bf twisted Hamiltonian} if\index[terms]{vector field!twisted Hamiltonian}
the cohomology class $[i_v(\omega)]\in H^1_\theta(M)$
vanishes. {\bf The group of twisted Hamiltonian
symplectomorphisms} is the group generated by 
the flows of time-de\-pen\-dent twisted Hamiltonian
vector fields. Its Lie algebra is equal to 
the kernel of the twisted flux homomorphism
 $F:\; \goth{symp}_\chi(M)\arrow H^1_\theta(M)$. Therefore,
the group of twisted Hamiltonian
symplectomorphisms is generated by 
the image of $\ker F\subset \widetilde \Symp_\chi(M)$. 
By \cite[Corollary 7.5]{hr1},
this group is simple.

\hfill

\question 
Let $N$ be a compact symplectic manifold, 
$\widetilde \Symp_0(N)$ the u\-ni\-versal cover of the group $\Symp_0(N)$
of symplectic isotopies, and $\widetilde \Symp_0(N) \stackrel
{F_N} \arrow H^1(N)$ the flux map. 
The flux conjecture, stated by\index[terms]{conjecture!McDuff's flux}
D. \index[persons]{McDuff, D.} McDuff in \cite{_McDuff:flux_}
and proven by K. Ono in \cite{_Ono:flux_},
claims that the image of $\ker F_N$ in
$\Symp_0(N)$ is closed, and coincides with the group of
Hamiltonian diffeomorphisms. Its analogue for LCS manifolds is still
not proven, though the preliminary steps
were all taken in \cite{hr1}.
Let $M$ be a compact LCS manifold,\index[terms]{manifold!LCS}
and  $\widetilde \Symp_\chi(M) \stackrel {F} \arrow H^1(M)$ 
the twisted flux map. Is the image
of $\ker F$ in $\Symp_\chi(M)$ closed?

\hfill

\question
A holomorphic conformally symplectic vector field
on a compact LCK manifold preserves the Gauduchon
form, and, conversely, every holomorphic vector field
preserving the Gauduchon metric is conformally symplectic.
Therefore, the Lie algebra of automorphisms of a compact LCK
manifold is the algebra of holomorphic Killing vector fields (with respect to ,
 the Gauduchon metric).\index[terms]{vector field!Killing}\index[terms]{vector field!holomorphic}
Let $X$ be a Killing holomorphic vector field
on a compact K\"ahler manifold. Then $X$ 
is Hamiltonian if and only if it has a fixed point
(Exercise \ref{_holo_Killing_fixed_point_Exercise_}).
\begin{enumerate}
\item Is there a similar way of identifying the twisted
Hamiltonian holomorphic vector fields in the Lie algebra
of holomorphic Killing vector fields
on a compact LCK manifold?\index[terms]{vector field!twisted Hamiltonian}
\item The lift of $X$ to a K\"ahler cover\index[terms]{cover!K\"ahler} of $M$
acts by holomorphic homotheties. If it acts by non-isometric
holomorphic homotheties, $M$ is LCK with potential\index[terms]{manifold!LCK!with potential}
(Exercise \ref{_non-LCK+pot_isometry_on_cover_Exercise_}).
If, in addition, the K\"ahler form\index[terms]{form!K\"ahler} $\tilde \omega$
on the covering $\tilde M$ is exact, then any
symplectic vector field on $\tilde M$
is Hamiltonian, and $X$ is twisted Hamiltonian.
Suppose now that $M$ does not admit an LCK metric 
with potential. Are all Killing holomorphic vector fields
on $M$ twisted Hamiltonian?\index[terms]{vector field!twisted Hamiltonian}\index[terms]{vector field!Killing}\index[terms]{vector field!holomorphic}
\end{enumerate}

\smallskip

\remark
Let $X$ be
a holomorphic Killing vector field\index[terms]{vector field!holomorphic}\index[terms]{vector field!Killing}
on a Vaisman manifold $M$.\index[terms]{manifold!Vaisman} Assume that
$X$ belongs to the group $\goth{symp}_\chi(M)$,
that is, lifts to an isometry on a K\"ahler cover
$(\tilde M, \tilde \omega)$ of $M$. Since $\tilde \omega$
is exact, $X$ is automatically twisted Hamiltonian.
By Exercise 
\ref{_decompo_holo_on_Vaisman_Exercise_}, any Killing holomorphic
vector field on a Vaisman manifold
belongs to the space generated by the Lee field,\index[terms]{Lee field} 
the anti-Lee field\index[terms]{Lee field!anti-} and a $\Sigma$-basic Killing 
holomorphic vector field; the latter two belong
to $\goth{symp}_\chi(M)$, and hence  they are twisted
Hamiltonian.

\hfill

\question\label{_Delzant_LCS_Question_}
A toric symplectic manifold is uniquely determined
by its {\bf Del\-zant polytope}. A \index[persons]{Delzant, T.} Delzant polytope\index[terms]{Delzant polytope}
is a convex polytope in $\R^n$ such that the slopes of its edges\index[terms]{manifold!symplectic!toric}
define a basis in $\Z^n \subset \R^n$. The image of the momentum
map of a toric symplectic manifold \index[terms]{momentum polytope}
(its {\bf momentum polytope}) is a \index[persons]{Delzant, T.} Delzant polytope, \index[terms]{theorem!Delzant}
and by Delzant theorem a toric symplectic manifold is uniquely
determined by its momentum polytope (\cite{_Cannas_da_Silva_,_Delzant:moment_}).
The notions of toric LCS manifold and\index[terms]{manifold!LCS!toric} 
toric LCK manifold are well studied, and the image of the momentum\index[terms]{manifold!LCK!toric}
map is known to be convex and polyhedral (\cite{bgp}, \cite{nico3}), 
if the action is of Lee type, that is, when the anti-Lee field 
is induced by the toric action. What about the toric action that is  not necessarily of Lee type?
Can we obtain a complete description of toric LCS manifolds in 
terms of the image of the momentum map?

\hfill

\question
It is well known that all  toric symplectic manifolds are K\"ahler type,
and the notions of a toric symplectic manifolds and a toric K\"ahler manifolds
coincide. In \cite[Example 6.4]{nico1}, N. \index[persons]{Istrati, N.} Istrati has constructed a toric LCS
structure\index[terms]{structure!LCS!toric} on the 4-dimensional torus $T^4$. Since
any complex structure on $T^4$ is K\"ahler (\ref{_B-L-intro_Theorem_}), 
this manifold\index[terms]{manifold!toric}
cannot be strict LCK by Vaisman theorem. This gives an
example of a toric LCS manifold which 
does not admit a toric LCK\index[terms]{manifold!LCS} structure.\index[terms]{structure!LCK!toric}
This is the only example known. How big is the gap between
the toric LCS and toric LCK manifolds? Can we characterize
the toric LCK manifolds among the toric LCS manifolds in
terms of the combinatorics of the momentum polytope?

\hfill

\question
Let $M$ be a quasi-regular toric Vaisman\index[terms]{manifold!Vaisman!quasi-regular toric} manifold,\index[terms]{manifold!Vaisman}
and $X$ the leaf space of its canonical foliation\index[terms]{foliation!canonical},
that is   a 
projective orbifold by \ref{_Structure_of_quasi_regular_Vasman:Theorem_}.
Since the torus acts on $X$ with an open orbit,
it is a toric orbifold. In \cite{_Lerman_Tolman:orbifold_},
\index[persons]{Lerman, E.} Lerman and \index[persons]{Tolman, S.} Tolman have generalized the \index[persons]{Delzant, T.} Delzant
polytope construction, classifying the toric\index[terms]{Delzant polytope}
orbifolds by their \index[persons]{Delzant, T.} Delzant polytopes.\index[terms]{orbifold!toric}
Can we characterize the toric orbifolds
obtained from Vaisman manifolds by their
combinatorial data? \index[terms]{manifold!Vaisman}

\hfill

\question\label{_toric_cone_equivalent_Question_}
 Let $M_i= \frac{\tilde M}{\langle A_i \rangle}$, $i=1, 2$,
be two Vaisman manifolds with the same algebraic cone, and
$\rho:\; (\C^*)^n \arrow \Aut(\tilde M)$\index[terms]{cone!algebraic}
an action of the torus with an open orbit,
commuting with $A_1, A_2$.\footnote{This would
actually imply that the image of $\rho$ contains
$A_1$ and $A_2$, as follows from 
\ref{_Toric_Vaisman_Nicolina_Theorem_},
Step 1.} Let $M_i$ be  quasi-regular,\index[terms]{manifold!Vaisman!quasi-regular}
and let $X_1, X_2$ be the leaf spaces of the
canonical foliations,\index[terms]{foliation!canonical} that are  a posteriori toric
orbifolds. These manifolds are not necessarily
isomorphic, or even homeomorphic. Even when 
$\tilde M = \C^n$, the different Vaisman quotients
will give all weighted projective spaces,
that are  clearly non-isomorphic and non-homeomorphic.
Can we relate the \index[persons]{Delzant, T.} Delzant-Lerman-\index[persons]{Tolman, S.}Tolman
combinatorial data of these two different
orbifold quotients? Can we enumerate all 
toric orbifold quotients obtained 
from a given algebraic cone?\index[terms]{cone!algebraic}
This question seems to be related to 
\ref{_Delzant_LCS_Question_}.

\hfill

\question\label{_Donaldson_scal_as_moment_Question_}
In \cite{_Donaldson:Gauge_} 
(see also \cite{_Donaldson:Conjectures_}), Donaldson 
gives a momentum map interpretation of constant scalar
curvature K\"ahler metrics. He shows that the scalar
curvature is a momentum map for the group of
Hamiltonian diffeomorphisms acting on the space
of K\"ahler metrics. A similar result for twisted
Hamiltonian diffeomorphisms acting on the space\index[terms]{momentum map}
of LCK metrics is still unknown.\index[terms]{metric!K\"ahler!csc}

\hfill

\question
In \cite{_Angella_Calamai_Pediconi_Spotti_},
\ref{_Donaldson_scal_as_moment_Question_} 
was answered in the special case
when the anti-Lee field\index[terms]{Lee field!anti-} $V:=I(\theta^\sharp)$,
is tangent to a compact Lie group action
acting on $M$ by twisted Hamiltonian
automorphisms. Note that $V$ is always twisted
Hamiltonian by \ref{_Lee_twisted_Ham_Example_}.
For a Vaisman manifold,\index[terms]{manifold!Vaisman} $V$ is actually Killing,\index[terms]{vector field!Killing}
hence it is tangent to a compact Lie group
(the group of isometries of a Riemannian manifold
is compact). It seems that for LCK manifolds
of OT-type and LCK with potential\index[terms]{manifold!LCK!with potential}, the vector field\index[terms]{manifold!Oeljeklaus--Toma (OT)}
$V$ is not tangent to a compact Lie group. Are there 
any other examples of LCK manifolds with $V$ 
tangent to a compact Lie group?


\section{Hopf manifolds}\index[terms]{manifold!Hopf}


\question
Let $B$ be a \index[persons]{Mall, D.} Mall bundle on a Hopf manifold $M$.\index[terms]{bundle!vector bundle!Mall}
If $B$ is non-resonant, $B$ admits a flat holomorphic connection
(\ref{_Mall_flat_connection_Theorem_}). However, the non-resonance
condition is not always necessary. For example, if $B$ is a line
bundle, it also admits a flat holomorphic connection
(\ref{_Picard_Hopf_Proposition_}). In fact we have
no example of a Mall bundle that does not admit\index[terms]{connection!flat}
a holomorphic connection. Do all Mall bundles admit
a holomorphic connection?\index[terms]{connection!holomorphic}

\hfill

\question
A non-resonant \index[persons]{Mall, D.} Mall bundle on a Hopf manifold has a unique
holomorphic connection (\ref{_conne_unique_Remark_}), that is  a posteriori flat.
A resonant Mall bundle may admit a non-trivial
family of holomorphic connections. Does this family
always contain a flat connection?

\hfill

\question
Recall that {\bf a hypercomplex structure}\index[terms]{structure!hypercomplex}
on a manifold $M$ is a triple of complex structure operators
$I,J,K\in \End(TM)$ that are  integrable and satisfy the
quaternionic relation $IJ = -JI = K$. The {\bf \index[persons]{Obata, M.} Obata
  connection} on $(M, I, J, K)$ is a\index[terms]{connection!torsion-free} torsion-free\index[terms]{connection!Obata}
connection $\nabla$ on $M$ such that $\nabla(I)= \nabla(J)=\nabla(K)=0$.
The Obata connection exists and is unique for any
hypercomplex manifold. Let $(H,I)$ be a Hopf
manifold admitting a hypercomplex structure $(I,J,K)$.
Is $(H,I)$ a linear Hopf manifold? This would 
follow if the Obata connection is flat, by
\ref{_flat_Hopf_is_linear_Theorem_}.

\hfill

\question
Can a Hopf manifold of dimension $n>2$ 
admit a pluriclosed (SKT) metric? See also
\ref{_SKT_balanced_LCK_Question_}. Can it have
a Hermitian metric with \index[terms]{metric!pluriclosed (SKT)}
$dd^c(\omega^p)=0$, where $1 \leq p \leq n-2$? It is not
hard to see that the Hopf manifolds do not admit balanced 
metrics (Exercise \ref{_balanced_exact_no_divisors_Exercise_}).\index[terms]{metric!balanced}

\hfill

\question
A holomorphic bundle $B$ on a Hopf surface\index[terms]{surface!Hopf}
is called {\bf filtrable} if $B$ admits a filtration
by coherent sheaves, with successive subquotients \index[terms]{bundle!vector bundle!extensible}
coherent sheaves of rank $\leq 1$, and {\bf extensible}
if the pullback of $B$ to $\C^2 \backslash 0$\index[terms]{bundle!vector bundle!filtrable}
can be extended to a coherent sheaf on $\C^2$.
Is filtrability equivalent to extensibility for any
vector bundle? 
For diagonal Hopf manifolds of dimension $>2$, filtrability follows
from \cite{_V:filtrable_}, and extensibility from \index[persons]{Siu, Y.-T.} Siu 
theorem (\ref{_Siu_extension_over_Z_Theorem_}),\index[terms]{theorem!Siu extension}
see also \ref{_extension_coherent_Theorem_}.

\hfill

\question
Let $H$ be a non-diagonal Hopf \index[terms]{manifold!Hopf!non-diagonal}
manifold, $\dim_\C H > 2$. Is any torsion-free
coherent sheaf over $H$\index[terms]{sheaf!coherent!torsion-free} filtrable?\index[terms]{sheaf!coherent}
We expect that the answer is affirmative.
Is this true when $H$ is an arbitrary
LCK manifold with potential\index[terms]{manifold!LCK!with potential}, $\dim_\C H > 2$?\index[terms]{manifold!Vaisman!quasi-regular}
A holomorphic vector bundle over a quasi-regular Vaisman
manifold is always filtrable by \cite{_Verbitsky:Sta_Elli_}. 
\index[terms]{bundle!vector bundle!filtrable}


\section{Foliations on LCK manifolds}


We are grateful to Jorge Vit\'orio \index[persons]{Pereira, J. V.} Pereira for a discussion
that led us to the questions raised in this section.

In \cite{_Demailly:Frobenius_}, J.-P. Demailly proved a
general integrability theorem for holomorphic
distributions. Let $L$ be a holomorphic line bundle on a
complex manifold $M$, and \index[terms]{distribution!holomorphic}
$\Omega\in H^0(\Omega^p(M) \otimes L^{-p})$ a holomorphic
$p$-form such that the sheaf $D:=\{v\in TM \ \ |\ \ i_v \Omega=0\}$
has rank $(n-p)$, that is, $\Omega$ is locally identified
with a transversal holomorphic volume form of the
holomorphic distribution $D\subset TM$. Assume that $L$ is
{\bf pseudo-effective}, that is, admits a \index[terms]{bundle!line!pseudo-effective}
singular Hermitian metric with curvature equal to a 
positive, closed current.\index[terms]{current!positive} Suppose that $M$ is compact and K\"ahler.
Then the distribution $D$ is integrable.

\index[persons]{Brunella, M.} Brunella's theorem claims that the cotangent sheaf
of a rank 1 foliation $F\subset TM$ on a complex manifold
is pseudo-effective, unless the closure of any leaf of $F$ 
is a rational curve. Since strict LCK manifolds 
cannot be covered by compact families\index[terms]{curve!rational}
of rational curves (\ref{_non-uniruled_Corollary_}), 
Brunella's theorem implies that a 1-dimensional  \index[terms]{theorem!Brunella} 
holomorphic foliation on a strict LCK manifold has
pseudoeffective cotangent sheaf. \index[terms]{foliation!holomorphic}

Consider a decomposition of holomorphic vector bundles
$TM= A \oplus B$, where $B$ is a rank 1 bundle.
We consider $B$ as an integrable rank 1 distribution on $M$.
Assume that $M$ is not uniruled  and K\"ahler.\footnote{A complex manifold
is called {\bf uniruled} if there exists a family of
rational curves passing through each point of $M$.}
By \index[persons]{Brunella, M.} Brunella, $B^*$ is pseudo-effective. 
Applying \index[persons]{Demailly, J.-P.} Demailly's theorem to a section of
the canonical bundle $K_M =\Lambda^p(A^*) \otimes B^*$,
we obtain that the distribution $A$ is integrable,
unless $H^0(K_M \otimes (B^*)^{\otimes k})=0$ for all $k$.

Using a similar argument based on Brunella's theorem,
Demailly proved in \cite{_Demailly:Frobenius_} that a compact K\"ahler manifold with 
a holomorphic contact structure is uniruled. \index[terms]{structure!holomorphic contact} \index[terms]{manifold!uniruled}

\index[persons]{Brunella, M.} Brunella's theorem is valid on any compact complex manifold,
however, \index[persons]{Demailly, J.-P.} Demailly's argument requires $M$ to be K\"ahler.
We are interested in how far we can go in that direction on
LCK manifolds.

\hfill

\question
Let $M$ be a compact LCK manifold, 
$L$ a pseudo-effective\index[terms]{bundle!line!pseudo-effective}
line bundle, and $\Omega\in H^0(\Omega^p(M) \otimes L^{-p})$ a holomorphic
$p$-form. Will it follow that the distribution
$D:=\{v\in TM \ \ |\ \ i_v \Omega=0\}$ is integrable?

\hfill

\question 
Note that compact LCK manifolds cannot be uniruled \index[terms]{manifold!uniruled}
(\ref{_non-uniruled_Corollary_}).
\index[persons]{Brunella, M.} Brunella's theorem (\cite{_Brunella:pseudoeff_})\index[terms]{theorem!Brunella} 
implies that any rank 1 subsheaf $L \subset TM$ is
pseudo-effective. Will it follow that  any compact 
LCK manifold $M$ does not admit a holomorphic contact structure?

\hfill

\remark
Let $M$ be a Vaisman manifold,\index[terms]{manifold!Vaisman} and $\Sigma \subset TM$
the canonical foliation\index[terms]{foliation!canonical}. Denote by $B$ the quotient
holomorphic bundle, $B:= TM/\Sigma$. The exact sequence
\begin{equation} \label{_cano_exact_Equation_}
0 \arrow \Sigma \arrow TM \arrow B \arrow 0
\end{equation}
splits when $M$ is a linear Hopf manifold, because
$TM$ is a flat bundle with monodromy\index[terms]{monodromy} $\Z$, 
acting diagonally, and $TM$ splits as a flat vector
bundle. On the other hand, on the Kodaira surface\index[terms]{surface!Kodaira} the
exact sequence \eqref{_cano_exact_Equation_} does not
split, because for a Kodaira surface $M$, the bundle $B$
is trivial, and hence  splitting of
\eqref{_cano_exact_Equation_} 
would imply that $M$ is parallelizable
and the group of its holomorphic
automorphisms is 2-dimensional. It is 1-dimensional
by \ref{_toric_LCK_surfaces_Proposition_}, Step 2,
thus \eqref{_cano_exact_Equation_} does not split.
The same argument implies that
\eqref{_cano_exact_Equation_} does not split
when $M$ is a Vaisman manifold\index[terms]{manifold!Vaisman!Heisenberg type} of Heisenberg type.

\hfill

\question
How can one characterize Vaisman manifolds for which
the exact sequence \eqref{_cano_exact_Equation_} splits?
Can \eqref{_cano_exact_Equation_} split when $M$
is an elliptic fibration over a projective 
manifold with ample canonical bundle?

\hfill

%
%
%
%

\question
Recall that a complex manifold
is called {\bf complex parallelizable}\index[terms]{manifold!complex parallelizable} 
(\cite{_Winkelmann:parallelizable_book_}) if
its tangent bundle is holomorphically trivial.
In Exercise \ref{_parallelizable_Exercise_}
it is shown that a compact Vaisman manifold\index[terms]{manifold!Vaisman} cannot
be complex parallelizable. We expect that an LCK manifold with
potential is never complex parallelizable. It is unknown whether
a compact complex parallelizable manifold can admit
an LCK structure.\index[terms]{structure!LCK}

\hfill

\question Let $M$ be an LCK manifold that has a 1-dimensional 
transversally K\"ahler holomorphic foliation. Is it necessarily Vaisman?

\hfill

\question Let $M$ be an LCK manifold with Lee form\index[terms]{form!Lee}
$\theta$ and $\Sigma$ the distribution generated by
$\theta^\sharp$ and $I\theta^\sharp$. On a Vaisman
manifold  $\Sigma$  is integrable, holomorphic  and has
totally geodesic leaves.  
\begin{enumerate}
\item Characterize the (compact)
LCK manifolds on which $\Sigma$ is integrable and
holomorphic. 
\item
 Characterize the (compact) LCK manifolds
on which $\Sigma$ is integrable and has totally geodesic
leaves.  
\end{enumerate}

\question\label{_BB_compact_leaf_Question_}
Using the Baum--Bott theorem,\index[terms]{theorem!Baum--Bott}
it is not hard to deduce that any 
smooth rank 1 holomorphic foliation $S$ on a classical
Hopf manifold has a compact leaf 
(Exercise \ref{_foliation_on_classical_Hopf_compact_leaf_Exercise_}).
Is the same true for all Hopf manifolds?
Suppose that $S$ is a rank 1 foliation on a Hopf manifold,
not necessarily smooth. Will it follow that $S$
has a leaf with compact closure?

\hfill

\question
The argument proving the existence of compact leaves
of rank 1 foliations on classical Hopf manifolds
(Exercise \ref{_foliation_on_classical_Hopf_compact_leaf_Exercise_})
uses the Baum--Bott formula and the non-existence of smooth foliations 
on complex projective spaces. A Heisenberg-type Vaisman nilmanifold $M$
projects to a compact torus that has non-singular rank 1 foliations; consequently, the Baum--Bott argument cannot be applied. However,\index[terms]{nilmanifold} 
it seems that a rank 1 smooth holomorphic foliation on $M$ is
unique (the canonical foliation)\index[terms]{foliation!canonical}. For Kodaira 
surfaces, this is proven in \index[terms]{foliation!holomorphic}
\cite[Remark 2.5]{_Ballico:foliations_}. 


\section[Logarithmic foliations on LCK manifolds
  with potential]{Logarithmic foliations\\ on LCK manifolds
  with potential}\index[terms]{manifold!LCK!with potential}\index[terms]{foliation!logarithmic}


\question
Let $M$ be an LCK manifold with potential,
and $\tilde M$ its K\"ahler $\Z$-covering\index[terms]{cover!K\"ahler $\Z$-}. Assume that
the $\Z$-action on $\tilde M$ admits {\bf a logarithm $\vec r$,}
that is, a holomorphic, $\Z$-invariant vector
field on $\tilde M$ which 
satisfies $\gamma= e^{\vec r}$, where $\gamma$
is the generator of the $\Z$-action.\index[terms]{action!$\Z$-} Then $\langle \vec r, I(\vec r)\rangle$
defines a holomorphic foliation on $M$, that we call
{\bf a logarithmic foliation}.\index[terms]{foliation!logarithmic} \index[terms]{foliation!taut}
The notion of a taut foliation was defined by D. \index[persons]{Sullivan, D.} Sullivan in
\cite{_Sullivan:taut_}. An involutive distribution $F \subset TM$
is {\bf taut} if there exists a metric on $M$ such that
all leaves of $F$ are minimal, that is, have vanishing 
mean curvature. In 
\cite{_Masa:duality_} it is shown that
$F$ is taut if and only if the basic cohomology \index[terms]{cohomology!basic}
group $H^q_b(M)$ is non-zero, where $q= \codim F$. In this case, $H^q_b(M)$
is 1-dimensional, and the multiplication in basic cohomology defines 
a Poincar\'e-type duality on $H^*_b(M)$ (Subsection \ref{_basic_taut_Subsection_}).
\index[terms]{duality!Poincar\'e}
Since an LCK manifold with potential\index[terms]{manifold!LCK!with potential} can be obtained as a
limit of Vaisman manifolds,\index[terms]{manifold!Vaisman} and the basic cohomology\index[terms]{cohomology!basic}
``seems to be'' semicontinuous, it makes sense
to ask whether the logarithmic foliation\index[terms]{foliation!logarithmic} on an 
LCK manifold with potential is always taut.

\hfill

\question Can a logarithmic foliation\index[terms]{foliation!logarithmic} on a non-Vaisman
LCK manifold with potential admit a transversally\index[terms]{foliation!transversally K\"ahler}
K\"ahler structure? We expect that the answer is negative.

\hfill

\question
Recall that a foliation on a Riemannian manifold\index[terms]{foliation!totally geodesic}
is {\bf totally geo\-desic} if its leaves are totally
geodesic. The canonical foliation\index[terms]{foliation!canonical} on a Vaisman manifold\index[terms]{manifold!Vaisman}
is totally geodesic. Let $\Xi$ be a logarithmic foliation\index[terms]{foliation!logarithmic} on a non-Vaisman
LCK manifold $M$ with potential.\index[terms]{manifold!LCK!with potential} Can $M$ admit an LCK
metric such that $\Xi$ is totally geodesic?
We expect that the answer is negative.

\hfill

\question
Let $M$ be a Vaisman manifold of LCK rank 1\index[terms]{rank!LCK}, 
$M = \tilde M/\langle \Z\rangle$, and $\Sigma$ the canonical
foliation. Assume that the generator of the $\Z$-action
on $\tilde M$ has a logarithm (this is not given,
see \ref{_only_for_gamma^k_log_exists_Example_}), 
and let $\Xi$ be the corresponding
logarithmic foliation.\index[terms]{foliation!logarithmic} Suppose that $\Xi\neq \Sigma$.
Will $\Xi$ always admit a transversally K\"ahler structure?
Is $\Xi$ always taut?\index[terms]{foliation!transversally K\"ahler}

\hfill

\question
Consider a smooth 1-dimensional holomorphic
foliation on a compact Vaisman manifold $M$.\index[terms]{manifold!Vaisman}
The canonical foliation\index[terms]{foliation!canonical} on $M$ always has 
at least $\dim_\C M$ compact leaves by
\ref{_Vaisman_bound_from_Lfsch_Corollary_}.
The logarithmic foliation\index[terms]{foliation!logarithmic} is by construction
induced from a foliation on the ambient
Hopf manifold $H \supset M$; from this
it is possible to see that it has at
least one compact leaf. On the other hand,
the OT manifold has many 1-dimensional
foliations without compact leaves.
Consider a Heisenberg type Vaisman manifold,\index[terms]{manifold!Vaisman!Heisenberg type}
obtained as the total space of $T^2$-bundle
$\pi:\; M \arrow X$ over an abelian manifold $X$,\footnote{A
compact complex torus is called {\bf an abelian
manifold}, if it is projective.} $\dim_\C X >1$.
Choose an irrational 1-dimensional foliation $R_X$ 
on $X$ that is  invariant under the parallel\index[terms]{parallel transport}
transport. Clearly, this foliation has no compact
leaves. Then 
$R:= \{ v\in TM \ \ \ |\ \ \  D\pi(v) \in R_X\}$
is a 2-dimensional holomorphic foliation that has no
compact leaves. Can a Vaisman manifold\index[terms]{manifold!Vaisman} admit
a smooth 1-dimensional holomorphic
foliation without any compact leaves?
What about LCK with potential?\index[terms]{manifold!LCK!with potential}
Can a Hopf manifold admit
a smooth 1-dimensional holomorphic
foliation without any compact leaves? (See also \ref{_BB_compact_leaf_Question_}).


\section{Flat affine structures on LCK manifolds}\index[terms]{structure!flat affine}


Many important examples of LCK manifolds come equipped
with a flat, torsion-free connection preserving the
complex structure. Manifolds with a flat, torsion-free connection preserving the\index[terms]{connection!flat}\index[terms]{connection!torsion-free}
complex structure are called {\bf complex
flat affine manifolds} (see also \ref{_flat_affine_Definition_}).
A flat affine manifold is called\index[terms]{manifold!complex!flat affine}
{\bf complete} if it can be obtained as a quotient
of a vector space by a free affine action of a discrete
group (\ref{_flat_affine_complete_Definition_}
and \ref{_dev_for_complete_affine_Theorem_}). 
Note that the completeness does not follow\index[terms]{manifold!flat affine!complete}
from compactness, indeed, many examples of compact flat affine
manifolds, such as the Hopf manifold, are not complete.
It is interesting that for many examples of flat affine LCK manifolds
the monodromy representation\index[terms]{monodromy!representation} on the tangent space 
splits onto a direct sum of 1-dimensi\-onal representations;
equivalently, the tangent bundle is a direct sum of
flat holomorphic line bundles.

\hfill

\question
It is not hard to show that the tangent bundle of an OT
manifold is a direct sum of holomorphic line bundles, $TM= \bigoplus L_i$.
We call such manifolds {\bf manifolds with split tangent bundle}.\index[terms]{manifold!with split tangent bundle}
For OT manifolds, the  bundles $L_i$ admit\index[terms]{manifold!Oeljeklaus--Toma (OT)}
a flat connection (\cite{_OVV:OT_subvarieties_}). 
This is also true for diagonal Hopf manifolds.
A projective complex surface \index[terms]{connection!flat}
is called {\bf a Beauville surface} if it is \index[terms]{surface!Beauville}
the quotient of a product $S\times S'$ of two Riemann surfaces\index[terms]{Riemann surface} 
by a finite group freely acting on $S\times S'$. As Beauville has
shown in \cite[Theorem C]{_Beauville:surface_}, a surface with 
a split tangent bundle 
is a quotient of a product of two Riemann surfaces
by a group acting freely (this case includes the \index[terms]{surface!Inoue} Inoue\index[terms]{surface!Hopf}
surfaces),  or a Hopf surface. When a projective
complex surface has split tangent bundle, it
is a \index[persons]{Beauville, A.} Beauville surface.
For arbitrary dimension, some results towards 
classifying projective manifolds with split tangent
bundle are obtained in \cite{_Brunella_Pereira_Touzet_}
and in \cite{_Catanese_Franciosi_,_Catanese_di_Scala_}. In the
first of these papers, Brunella, Pereira and \index[persons]{Touzet, F.} Touzet have shown that
any K\"ahler manifold with split tangent bundle
is a quotient of a product of Riemann surfaces.\index[terms]{theorem!Beauville}
However, at the moment, there is no analogue of Beauville 
or \index[persons]{Brunella, M.} Brunella--\index[persons]{Pereira, J. V.}Pereira--Touzet theorem for\index[terms]{theorem!Brunella--Pereira--Touzet}
non-K\"ahler manifolds with split tangent bundle 
and $\dim>2$. Are there any other LCK manifolds with split
tangent bundle, $TM= \bigoplus L_i$?
What happens when all $L_i$ admit a flat
holomorphic connection?

\hfill

\question
All known examples of complex, complete affine manifolds
admitting an LCK metric are in fact Vaisman nilmanifolds,
that is, nilmanifolds of Heisenberg type (\ref{kod_nil}).
Are all complex, complete affine manifolds admitting an\index[terms]{nilmanifold}
LCK metric this type?\index[terms]{manifold!Vaisman!Heisenberg type}
The OT-manifolds are flat affine complex manifolds, 
but they are never complete: their universal covering
is $\C^r \times {\Bbb H}^k$.\index[terms]{manifold!Oeljeklaus--Toma (OT)}

\hfill

\question
Let $M$ be a complex, not necessarily complete affine manifold
admitting an LCK structure.\index[terms]{structure!LCK} The examples include
 Vaisman nilmanifolds\index[terms]{nilmanifold} of Heisenberg type, Hopf
manifolds and OT manifolds. In all these cases, the monodromy\index[terms]{monodromy} of the flat
connection  is abelian. 
Let $M$ be a complex, not necessarily complete, affine manifold
admitting an LCK structure\index[terms]{structure!LCK}. Is its monodromy  abelian?

\hfill

\hfill

{
\renewcommand {\epigraphflush}{center}

\epigraph{\it In conclusion, there is no conclusion. Things go on as
they always have, getting weirder all the
time.}{\sc\scriptsize The Principia Discordia, 4th edition.}

\epigraph{\it Math is hard, let's go shopping.}{\sc\scriptsize  attributed to Teen Talk Barbie, by Mattel Inc. (1992) }
}


{\scriptsize
\printindex[terms]
\printindex[persons]
}

\end{document}